\documentclass{article}
\usepackage{amsmath}
\usepackage{graphicx}
\usepackage{amssymb}
\usepackage{float}
\usepackage{tabularx}
\usepackage{here}
\usepackage{amsfonts}
\newcommand{\field}[1]{\mathbb{#1}}
\newcommand{\C}{\field{C}}
\newcommand{\R}{\field{R}}
\def\gR{$\mathfrak{R}$}
\def\gr{$\mathfrak{r}$}
\def\gD{$\mathfrak{D}$}
\def\gh{$\mathfrak{h}$}
\def\gs{$\mathfrak{s}$}
\def\go{$\mathfrak{o}$}
\def\gp{$\mathfrak{p}$}
\def\gu{$\mathfrak{u}$}
\def\ge{$\mathfrak{e}$}
\def\gg{$\mathfrak{g}$}
\def\gd{$\mathfrak{d}$}
\def\gf{$\mathfrak{f}$}
\def\gP{$\mathfrak{P}$}
\def\gJ{$\mathfrak{J}$}
\def\gA{$\mathfrak{A}$}
\def\ga{$\mathfrak{a}$}
\def\gb{$\mathfrak{b}$}
\def\gx{$\mathfrak{x}$}
\def\gy{$\mathfrak{y}$}
\def\gk{$\mathfrak{k}$}
\def\gi{$\mathfrak{i}$}
\def\gz{$\mathfrak{z}$}
\def\gw{$\mathfrak{w}$}
\def\gl{$\mathfrak{l}$}
\def\gC{$\mathfrak{C}$}

\def\gl{$\mathfrak{l}$}
\def\ggC{\mathfrak{C}}
\def\ggD{\mathfrak{D}}

\def\ggR{\mathfrak{R}}
\def\gS{$\mathfrak{S}$}
\def\gM{$\mathfrak{M}$}
\def\gm{$\mathfrak{m}$}
\def\gn{$\mathfrak{n}$}
\def\gt{$\mathfrak{t}$}
\def\gF{$\mathfrak{F}$}

\begin{document}

\title{Digitalization of exceptional simple Lie algebras into matrices over complex numbers}
\author{Takao IMAI}
\maketitle


\begin{abstract}
We give the images of the adjoint representations of exceptional simple Lie algebras by matrices over complex numbers. 
Next, we digitalize these matrices by the use of Maxima, which is a computer algebra system.
These digitalized matrices are provided by using Maxima's function: \emph{matrix}.

We prove that these digitalized matrices are closed by Lie bracket operations and make up simple Lie algebras. 
 Moreover, to prove the type of the exceptional simple Lie algebra, we calculate the root system using Maxima for Lie bracket operations as matrix calculations.

We show some examples of classical Lie subalgebras of these digitalized matrices.
\end{abstract}

\section*{Introduction}
\ \ \ \ 
H.Freudenthal\cite{Freudenthal1},\cite{Freudenthal2} constructed exceptional Lie algebras over real numbers by algebraic product in the Cayley algebra and the Jordan algebra. For these algebraic product see \emph{Section 1}.
Afterward 
I.Yokota, O.Sukuzawa, O.Yasukura, and T.Imai reconstructed exceptional simple Lie algebras to realize exceptional simple, connected, simply connected, compact Lie groups.
Especially, 
in T.Imai and I.Yokota\cite{YokotaImai2}, the exceptional simple Lie algebra $e^{\C}_{8}$ of type $E_{8}$ is given by a $\C$-linear space

\ \ \ \ \ \ \ \ \ \ \ \ \ \ \ \ \ \ \ge$_{8}^{\C}$ = \ge%
$_{7}^{\C} \oplus $\gP$^{\C} \oplus $\gP$^{\C} \oplus \C\oplus \C\oplus \C,$

\noindent
where \ge$_{7}^{\C}$ is the exceptional simple Lie algebra of type $E_{7}$,
\gP$^{\C}= $\gJ$^{\C} \oplus $\gJ$^{\C} \oplus \C\oplus \C$, and \gJ \ is the exceptional Jordan algebra.

Since the Cayley algebra and the Jordan algebra are hard to calculate, we provide digitalization of exceptional simple 
Lie algebras into matrices over complex numbers. So we can calculate on a computer
as a matrix calculation using Maxima, which is a computer algebra system.

To digitalize exceptional simple Lie algebras we have 3 steps.

At first (\emph{Sections:2,3,4,5}), we break \ge$_{7}^{\C}$ finely as a following  vector space:

\ \ \ \ \ \ \ \ \ \ \ge$_{7}^{\C}=$\gD$_{4}^{\C} \oplus $\gA$^{\C} \oplus $\gJ$_{0}^{\C} \oplus $\gJ$^{\C} \oplus $\gJ$^{\C} \oplus \C $.

\noindent
 Then  we define the $\R$-vector space \gR$_{4}$ for the real compact, simple exceptional Lie algebra \gf$_{4}$ (see $\emph{Proposition 2.2}$) of type $F_{4}$. And we build a real Lie algebraic structure on \gR$_{4}$, by showing that it is isomorphic to \gf$_{4}$.
We define $\R$-vector space \gR$_{6}$ by extension of \gR$_{4}$. Similarly, we define $\R$-vectore spaces  \gR$_{7}$ by extension of \gR$_{6}$ and \gR$_{8}$ by extension of \gR$_{7}$.
And we build a real Lie algebraic structure on \gR$_{6}$, \gR$_{7}$, and \gR$_{8}$,by showing that these are isomorphic to 
the real exceptional simple Lie algebra \ge$_{6,1}$(see $\emph{Proposition 3.6}$) of type $E_{6(-26)}$, 
\ge$_{7,1}$(see $\emph{Proposition 4.9}$) of type $E_{7(-25)}$, and \ge$_{8,1}$(see $\emph{Proposition}$ $\emph{ 5.5}$) of type $E_{8(-24)}$ respectively. 

\gR$_{4}=\{(D,M) \mid D\in $\gD$_{4},M\in $\gA $\}\cong $ \gf$_{4}$,

\gR$_{6}=\{(D,M,T) \mid D\in $\gD$_{4},M\in $\gA$,T\in $\gJ$_{0} \}\cong $ \ge$_{6,1}$,

\gR$_{7}=\{(D,M,T,A,B,\rho) \mid D\in $\gD$_{4},M\in $\gA$,T\in $\gJ$%
_{0},A,B\in $\gJ$,\rho \in \R\}\cong $ \ge$_{7,1}$,

\gR$_{8}=\{(D,M,T,A,B,\rho,X,Y,\xi,\eta,Z,W,\zeta,\omega,r,s,u) \mid D\in $\gD$_{4},M\in $\gA,

\ \ \ \ \ \ \ \ \ \ $T\in $\gJ$_{0},A,B,X,Y,Z,W\in $\gJ$,\rho,\xi,\eta,\zeta,\omega,r,s,u \in \R\}\cong$ \ge$_{8,1}$,

 \noindent
 where $(a_{1},a_{2},\cdot \cdot \cdot ,a_{n}),(a_{i} \in A_{i})$ means an element of a vectore space 
 $A_{1}\oplus A_{2}\oplus \cdot \cdot \cdot \oplus A_{n}$.
 For those who are not good at calculating the Cayley algebra and the Jordan algebra, we describe the calculation of the Cayley algebra  and the Jordan algebra in detail.
 
Second (\emph{Sections:6,7,8,9,10}), we define some matrices for use in the adjoint representations(see \emph{Section 6}), and we give the matrix images over real numbers of the adjoint representations $ad($\gR$_{4})$ of \gR$_{4}$,  $ad($\gR$_{6})$ of \gR$_{6}$, $ad($\gR$_{7})$ of \gR$_{7}$, and $ad($\gR$_{8})$ of \gR$_{8}$ respectively ( see \emph{Theorems:7.2, 8.3, 9.2, and 10.7}).
Let we put the $248 \times 248$ matrix image of $ad($\gR$_{8}^{\C})$ over complex numbers by \gr$_{8}^{\C}$.
Since \gR$_{4}$ is the subalgebra of \gR$_{8}$, let we put the subset of \gr$_{8}^{\C}$ by \gr$_{4}^{\C}$ 
which is restricted \gR$_{8}^{\C}$ to \gR$_{4}^{\C}$.
Similarly, let we put the subsets of \gr$_{8}^{\C}$ by \gr$_{6}^{\C}$ and  \gr$_{7}^{\C}$ 
which are restricted \gR$_{8}^{\C}$ to \gR$_{6}^{\C}$ and \gR$_{7}^{\C}$ respectively,
such that \gr$_{4}^{\C} \subset $\gr$_{6}^{\C} \subset $\gr$_{7}^{\C}\subset $\gr$_{8}^{\C}$.

Finally (\emph{Sections:11,12}) , we digitalize \gr$_{4}^{\C}$,  \gr$_{6}^{\C}$,  \gr$_{7}^{\C}$, and \gr$_{8}^{\C}$.
The digitalized matrix of \gr$_{8}^{C}$(see \emph{Section 11 and 12}) are provided by using Maxima's function:\emph{matrix}.
We prove that these digitalized matrices of  \gr$_{4}^{\C}$,  \gr$_{6}^{\C}$,  \gr$_{7}^{\C}$, and \gr$_{8}^{\C}$ are closed by Lie bracket operations and make up simple Lie algebras. 

Moreover in \emph{Sections:13,14,15,16,17, and 18} we prove the types of exceptional simple Lie algebras of these digitalized matrices. 
For this purpose, we caluculate the root systems to draw the Dynkin diagrams.
In I.Yokota\cite{Yokota1} the roots and root vectors were described explicitly but the root vectors were partially omitted.
Therefore we describe all the root vectors explicitly. Lie bracket operations are calculated by matrix calculations using the data of these digitalized matrices.

In \emph{Sections: 14, 15, and 18}, we also digitalize the exceptional simple Lie algebra of type $F_{4}$ as matrices in $M(27 \times 27,\R)$, type $E_{6}$ as matrices in $M(27 \times 27,\C)$, and type $G_{2}$ as matrices in $M(8 \times 8,\R)$ respectively.

In \emph{Section 19}, we show some examples of classical Lie subalgebras on the digitalized matrices of \gr$_{8}^{\C}$.

In Appendix A, we fix some printing mistakes in I.Yokota\cite{Yokota1} and M.Sato and T.Kimura\cite{SatoKimura1}, 
which we found in the process of calculating the root systems by using these digitalized matrices.

In Appendix B, we denote elements of the Lie groups as matrices in $M(27 \times 27,\R)$ corresponding to the bases of the Lie algebra \gf$_{4}$, 
 as matrices in $M(27 \times 27,\C)$ corresponding to the bases of the Lie algebra  \ge$_{6}$, and as matrices in $M(8 \times 8,\R)$ corresponding to the bases of the Lie algebra \gg$_{2}$.

\bigskip

We prepared the digitalized 248 basis matrices of \gr$_{8}^{\C}$ by using Maxima's function:\emph{matrix}
and programming codes to calculate the bracket operation.
These data can be obtained from the following \emph{GitHub} repository.

\ \ \ \ \ \emph{https://github.com/TakaoIMAI/An-example-of-e8}

\noindent
The programming codes have following 6 contens.
With the programming codes, readers can simulate the process of digitalization and verify the proof of this paper.

\noindent
(1)Operations (see \emph{Section 1}) of the Cayley algebra and the Jordan algebra by using Maxima.

\noindent
(2)Lie bracket operations of \gR$_{4}^{}$, \gR$_{6}^{}$, \gR$_{7}^{}$, and \gR$_{8}^{}$ by using Maxima.

\noindent
(3)Convertes from $ad($\gR$_{8}^{\C})$ over the Cayley algebra to the matrix of \gr$_{8}^{\C}$,  
by using Maxima for performing mathematical operation and using the programming language Python for string conversion. 

\noindent
(4)Check the Lie bracket operations as matrix calculations of \gr$_{4}^{\C}$, \gr$_{6}^{\C}$, \gr$_{7}^{\C}$, and \gr$_{8}^{\C}$ and confirming these root systems.

\noindent
(5)Decompose an element of \gr$_{8}^{\C}$ into a linear combination of the 248 bases matrices.

\noindent
(6)Check the Lie bracket operations of Lie subalgebras of \gr$_{8}^{\C}$.

\noindent
Using these digitalized data, Lie bracket operations of exceptional Lie algebras are performed as matrix calculations in the same way as classical Lie algebras.
Computer algebra systems are useful for explicitly calculating high-dimensional Lie algebras.
So readers can find these digitalized data are useful. 

\bigskip

Deep appreciation to Professor Osami Yasukura for his kind advice during the preparation of this paper.

\bigskip

\tableofcontents

\part{Digitalization of exceptional Lie algebras}
\section{Preliminaries}

\ \ \ \ \S 1,\S 2, and \S 3 are quoted from I.Yokota\cite{Yokota1}.

\bigskip

\S 1 Cayley algebra \gC

\ We denote the division Cayley algebra by \gC. We now explain this
algebra.

\bigskip

\emph{Definition 1.1. } We consider an 8-dimensional $\R$%
-vector space with bases

\noindent
$\{e_{0}=1,e_{1},e_{2},e_{3},e_{4},e_{5},e_{6},e_{7}\}$ and define a
multiplication between them as follows in the figure below.

\ \ \ \ \ \ \ \ \ \ \ \ \ \ \ \ $e_{1}e_{2}=e_{3}, e_{2}e_{3}=e_{1}, e_{3}e_{1}=e_{2},$

\ \ \ \ \ \ \ \ \ \ \ \ \ \ \ \ $e_{1}e_{4}=e_{5},e_{4}e_{5}=e_{1},e_{5}e_{1}=e_{4},$

\ \ \ \ \ \ \ \ \ \ \ \ \ \ \ \ $e_{3}e_{4}=e_{7},e_{4}e_{7}=e_{3},e_{7}e_{3}=e_{4},$

\ \ \ \ \ \ \ \ \ \ \ \ \ \ \ \ $e_{3}e_{5}=e_{6},e_{5}e_{6}=e_{3},e_{6}e_{3}=e_{5},$

\ \ \ \ \ \ \ \ \ \ \ \ \ \ \ \ $e_{6}e_{4}=e_{2},e_{4}e_{2}=e_{6},e_{2}e_{6}=e_{4},$

\ \ \ \ \ \ \ \ \ \ \ \ \ \ \ \ $e_{6}e_{7}=e_{1},e_{7}e_{1}=e_{6},e_{1}e_{6}=e_{7},$

\ \ \ \ \ \ \ \ \ \ \ \ \ \ \ \ $e_{2}e_{5}=e_{7},e_{5}e_{7}=e_{2},e_{7}e_{2}=e_{5}.$

$e_{0}=1$ is the unit of the multiplication and assume

\ \ \ \ \ \ \ \ \ \ \ \ \ \ $e_{i}^{2}=-1,i\neq 0,\ \
e_{i}e_{j}=-e_{j}e_{i},i\neq j,i\neq 0,j\neq 0,$

\noindent
and the distributive law. Thus \gC\  has a multiplication. $x1\  (x\in \R)$
is briefly denoted by $x$

\begin{figure}[H]
\centering
\includegraphics[width=5cm, height=5cm]{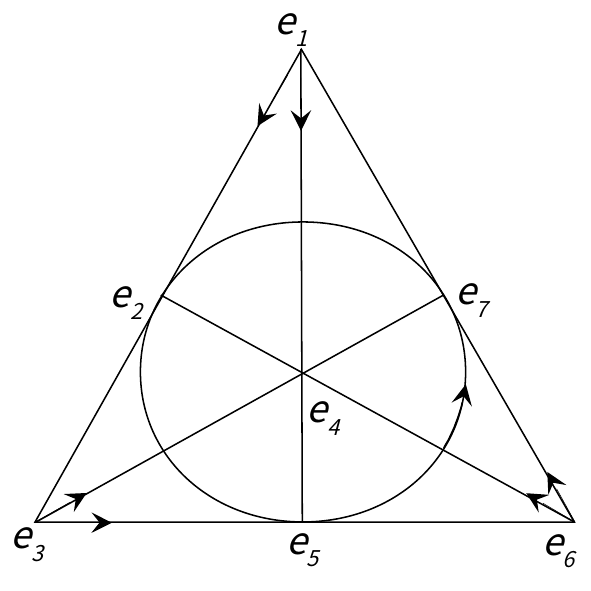}
\end{figure}

\emph{Definition 1.2. } In \gC, the conjugate $\overline{x}$%
, an inner product
$(x,y)$, the length $|x|$ and the real part $\R(x)$ are defined respectively by

\ \ \ \ \ \ $\overline{x_{0}+\sum_{i=1}^{7}x_{i}e_{i}}=x%
_{0}-\sum_{i=1}^{7}x_{i}e_{i}$ , $(\sum_{i=0}^{7}x_{i}e_{i},\sum_{i=0}^{7}y%
_{i}e_{i})=\sum_{i=0}^{7}x_{i}y_{i}$,

$\ \ \ \ \ \ \ \ \ \ \ \ |x|=\sqrt{(x,x)}$ \ , \ \ $\R(x_{0}+\sum_{i=1}^{7}x%
_{i}e_{i})=x_{0}.$

\bigskip

\S 2 Exceptional Jordan algebra \gJ

\ Let \gJ$=J(3,$\gC $)$ denote all 3 $\times $ 3 Hermitian
matrices with entries in the Cayley algebra \gC:

\ \ \ \ \ \ \ \ \ \ \ \ \gJ$=\{X\in M(3\times 3,$\gC $) \mid X^{\ast}=X\}$,

\noindent
where $X^{\ast }$= $^{t}\overline{X},\ ^{t}$ means transpose with matrix. 

Any element $X\in $\gJ\  is of the form

$\ \ \ \ \ \ \ \ \ \ \ \ X=X(\chi ,x)=\left( 
\begin{array}{ccc}
\chi _{1} & x_{3} & \overline{x_{2}} \\ 
\overline{x_{3}} & \chi _{2} & x_{1} \\ 
x_{2} & \overline{x_{1}} & \chi _{3}%
\end{array}%
\right) $, \ \ $\chi_{i}\in \R,x_{i}\in $\gC .

\noindent
\gJ\  is a 27-dimensional $\R$-vector space.

\bigskip

\emph{Definition 1.3. } In \gJ, the multiplication $X \circ Y$ , called the Jordan

\noindent
multiplication, is defined by

$\ \ \ \ \ \ \ \ \ \ \ \ \ X \circ Y=\frac{1}{2}(XY+YX).$

In \gJ, we define the trace $tr(X)$ and an inner product $(X,Y)$
respectively by

$\ \ \ \ \ \ \ \ \ \ \ \ \ tr(X)=\chi _{1}+\chi _{2}+\chi _{3},$ $X=X(\chi
,x)$,

$\ \ \ \ \ \ \ \ \ \ \ \ \ (X,Y)=tr(X \circ Y).$

\bigskip

\emph{Definition 1.4. } In \gJ, we define a multiplication $%
X\times Y$, called the

\noindent
Freudenthal multiplication,by

$\ \ \ \ \ \ \ \ \ \ \ \ \ X\times Y=\frac{1}{2}(2X \circ Y-tr(X)Y-tr(Y)X+(tr(X)tr(Y)-(X,Y))E),$

(where $E$ is the unit matrix) and a trilinear form $(X,Y,Z)$ and the
determinant
detX respectively by

$\ \ \ \ \ \ \ \ \ \ \ (X,Y,Z)=(X,Y\times Z)$, $\ \ \ \ detX=\frac{1}{3}%
(X,X,X).$

\bigskip

\ In \gJ, we adopt the following notations:

\ \ \ \ \ $E_{1}=\left( 
\begin{array}{ccc}
1 & 0 & 0 \\ 
0 & 0 & 0 \\ 
0 & 0 & 0%
\end{array}%
\right) ,$ \ \ \ \ $E_{2}=\left( 
\begin{array}{ccc}
0 & 0 & 0 \\ 
0 & 1 & 0 \\ 
0 & 0 & 0%
\end{array}%
\right) ,$ \ \ \ \ \ $E_{3}=\left( 
\begin{array}{ccc}
0 & 0 & 0 \\ 
0 & 0 & 0 \\ 
0 & 0 & 1%
\end{array}%
\right) ,$

$\ F_{1}(x)=\left( 
\begin{array}{ccc}
0 & 0 & 0 \\ 
0 & 0 & x \\ 
0 & \overline{x} & 0%
\end{array}%
\right) ,$ \ $F_{2}(x)=\left( 
\begin{array}{ccc}
0 & 0 & \overline{x} \\ 
0 & 0 & 0 \\ 
x & 0 & 0%
\end{array}%
\right) ,$ \ $F_{3}(x)=\left( 
\begin{array}{ccc}
0 & x & 0 \\ 
\overline{x} & 0 & 0 \\ 
0 & 0 & 0%
\end{array}%
\right) .$

For $X=X(\chi ,x)$ and $Y=Y(\gamma ,y)$ $\in $\gJ, the explicit
forms in the term of their entries are as follows.

$\ \ \ \ \ \ \ \ \ \ \ \ \ \ \ X \circ Y=(\chi _{1}\gamma
_{1}+(x_{2},y_{2})+(x_{3},y_{3}))E_{1}$

$\ \ \ \ \ \ \ \ \ \ \ \ \ \ \ \ \ \ \ \ \ \ \ \ +(\chi _{2}\gamma
_{2}+(x_{3},y_{3})+(x_{1},y_{1}))E_{2}$

$\ \ \ \ \ \ \ \ \ \ \ \ \ \ \ \ \ \ \ \ \ \ \ \ +(\chi _{3}\gamma
_{3}+(x_{1},y_{1})+(x_{2},y_{2}))E_{3}$

\ \ \ \ \ \ \ \ \ \ \ \ \ \ \ \ \ \ \ \ \ \ \ $\ +F_{1}(\frac{1}{2}((\chi _{2}+\chi
_{3})y_{1}+(\gamma _{2}+\gamma _{3})x_{1}+\overline{y_{2}x_{3}}+\overline{%
x_{2}y_{3}}))$

$\ \ \ \ \ \ \ \ \ \ \ \ \ \ \ \ \ \ \ \ \ \ \ \ +F_{2}(\frac{1}{2}((\chi _{3}+\chi
_{1})y_{2}+(\gamma _{3}+\gamma _{1})x_{2}+\overline{y_{3}x_{1}}+\overline{%
x_{3}y_{1}}))$

$\ \ \ \ \ \ \ \ \ \ \ \ \ \ \ \ \ \ \ \ \ \ \ \ +F_{3}(\frac{1}{2}((\chi _{1}+\chi
_{2})y_{3}+(\gamma _{1}+\gamma _{2})x_{3}+\overline{y_{1}x_{2}}+\overline{%
x_{1}y_{2}})),$

$\ \ \ \ \ \ \ \ \ \ \ \ \ \ (X,Y)=\sum_{i=1}^{3}(\chi _{i}\gamma
_{i}+2(x_{i},y_{i})),$

\ \ \ \ \ \ \ \ \ \ \ \ \ $X\times Y=(\frac{1}{2}(\chi _{2}\gamma _{3}+\chi
_{3}\gamma _{2})-(x_{1},y_{1}))E_{1}$

\ \ \ \ \ \ \ \ \ \ \ \ \ \ \ \ \ \ \ \ \ \ $\ +(\frac{1}{2}(\chi _{3}\gamma
_{1}+\chi _{1}\gamma _{3})-(x_{2},y_{2}))E_{2}$

\ \ \ \ \ \ \ \ \ \ \ \ \ \ \ \ \ \ \ \ \ \ $\ +(\frac{1}{2}(\chi _{1}\gamma
_{2}+\chi _{2}\gamma _{1})-(x_{3},y_{3}))E_{3}$

\ \ \ \ \ \ \ \ \ \ \ \ \ \ \ \ \ \ \ \ \ \ $\ +F_{1}(\frac{1}{2}(-\gamma
_{1}x_{1}-\chi _{1}y_{1}+\overline{x_{2}y_{3}}+\overline{y_{2}x_{3}}))$

\ \ \ \ \ \ \ \ \ \ \ \ \ \ \ \ \ \ \ \ \ \ $\ +F_{2}(\frac{1}{2}(-\gamma
_{2}x_{2}-\chi _{2}y_{2}+\overline{x_{3}y_{1}}+\overline{y_{3}x_{1}}))$

\ \ \ \ \ \ \ \ \ \ \ \ \ \ \ \ \ \ \ \ \ \ $\ +F_{3}(\frac{1}{2}(-\gamma
_{3}x_{3}-\chi _{3}y_{3}+\overline{x_{1}y_{2}}+\overline{y_{1}x_{2}})),$

\bigskip

\S 3 Lie algebras \gd \ and \gd$_{4}$

For an $K$-vector space \textbf{V} $(K = \R , \C )$, $Iso_{K}($\textbf{V}%
$)$ denotes all $K$-linear isomorphisms of \textbf{V} .

For an $K$-vector spaces \textbf{V},\textbf{W} $(K = \R , \C )$, $Hom_{K}%
($\textbf{V}$,$\textbf{W}$)$ denotes all $K$-homomorphisms
$\textrm{f} : \textbf{V} \to \textbf{W}$. $Hom_{K}($\textbf{V}$,$\textbf{V%
}$)$ is briefly denoted by $Hom_{K}($\textbf{V}$)$.

For an $\R$-vector space \textbf{V},we denote an complexification \textbf{V}$%
\otimes \C$ of \textbf{V} by \textbf{V}$^{\C}$.

We extend linear operations on \textbf{V} to on \textbf{V}$^{\C}$as complex
linear operations.

For the Lie group $G\subset Iso_{K}($\textbf{V}$)$ , we consider the
Lie algebra \gg\  of $G$ as follows.

\ \ \ \ \ \ \ \ \ \ \ \ \ \ \ \gg$=\{X\in
Hom_{K}($\textbf{V}$)$ $|$ $e^{tX}\in G,(^{\forall }t\in $K$)\}$.

\ We consider the Lie algebras

\ \ \ \ \ \ \ \ \ \ \ \ \ \ \gD$_{4}=$\gs \go$(8)=\{D\in
M(8\times 8,R) \mid D+^{t}D=0\},$

\ \ \ \ \ \ \ \ \ \ \ \ \ \ \gd$=\{D\in Hom_{\R}($\gC $)|(Dx,y)+(x,Dy)=0\}$ .

\bigskip

\emph{Definition 1.5. }\ We define $\R$-linear mappings $\ G_{ij}:$%
\gC $\to$ \gC $,i,j=0,1,\cdot \cdot \cdot ,7,i\neq j$
\ satisfying

$\ \ \ \ \ \ \ \ \ \ \ \ \ \ \ \ \ \ \ \ \ \ \ \ \
G_{ij}e_{j}=e_{i},G_{ij}e_{i}=-e_{j},G_{ij}e_{k}=0,k\neq i,j.$

In other words, we have

\ \ \ \ \ \ \ \ \ \ \ \ \ \ \ \ \ \ \ \ \ \ \ \ $%
G_{ij}x=e_{i}(e_{j},x)-e_{j}(e_{i},x),x\in $\gC.

\noindent
Furthermore we define $R$-linear mappings $F_{ij}:$\gC$ \to$ \gC$
,i,j=0,1,\cdot \cdot \cdot ,7,i\neq j$ by

$\ \ \ \ \ \ \ \ \ \ \ \ \ \ \ \ \ \ \ \ \ \ \ \ F_{ij}x=\frac{1}{2}e_{i}(%
\overline{e_{j}}x),x\in $\gC .

\bigskip

\emph{\ Proposition 1.6. }(I.Yokota\cite[\emph{Lemma 1.3.1.}]{Yokota1}) \ For $i,j=0,1,\cdot
\cdot \cdot ,7,i\neq j$, we have $F_{ij} \in$ \gd,
and when $i<j,F_{ij}$ is expressed in terms of $G_{ij}$ as follows.

{\fontsize{8pt}{10pt} \selectfont%
$\left[ 
\begin{array}{c}
2F_{01}=G_{01}+G_{23}+G_{45}+G_{67} \\ 
2F_{23}=G_{01}+G_{23}-G_{45}-G_{67} \\ 
2F_{45}=G_{01}-G_{23}+G_{45}-G_{67}\  \\ 
2F_{67}=G_{01}-G_{23}-G_{45}+G_{67}%
\end{array}%
\right] $ $\ \left[ 
\begin{array}{c}
2F_{02}=G_{02}-G_{13}-G_{46}+G_{57} \\ 
2F_{13}=-G_{02}+G_{13}-G_{46}+G_{57} \\ 
2F_{46}=-G_{02}-G_{13}+G_{46}+G_{57} \\ 
2F_{57}=G_{02}+G_{13}+G_{46}+G_{57}%
\end{array}%
\right] $

$\left[ 
\begin{array}{c}
2F_{03}=G_{03}+G_{12}+G_{47}+G_{56} \\ 
2F_{12}=G_{03}+G_{12}-G_{47}-G_{56} \\ 
2F_{47}=G_{03}-G_{12}+G_{47}-G_{56}\  \\ 
2F_{56}=G_{03}-G_{12}-G_{47}+G_{56}%
\end{array}%
\right] $ $\ \left[ 
\begin{array}{c}
2F_{04}=G_{04}-G_{15}+G_{26}-G_{37} \\ 
2F_{15}=-G_{04}+G_{15}+G_{26}-G_{37} \\ 
2F_{26}=G_{04}+G_{15}+G_{26}+G_{37} \\ 
2F_{37}=-G_{04}-G_{15}+G_{26}+G_{37}%
\end{array}%
\right] $

$\left[ 
\begin{array}{c}
2F_{05}=G_{05}+G_{14}-G_{27}-G_{36}\\ 
2F_{14}=G_{05}+G_{14}+G_{27}+G_{36}\\ 
2F_{27}=-G_{05}+G_{14}+G_{27}-G_{36}\\ 
2F_{36}=-G_{05}+G_{14}-G_{27}+G_{36}%
\end{array}%
\right] \ \left[ 
\begin{array}{c}
2F_{06}=G_{06}-G_{17}-G_{24}+G_{35} \\ 
2F_{17}=-G_{06}+G_{17}-G_{24}+G_{35} \\ 
2F_{24}=-G_{06}-G_{17}+G_{24}+G_{35} \\ 
2F_{35}=G_{06}+G_{17}+G_{24}+G_{35}%
\end{array}%
\right] $

$\left[ 
\begin{array}{c}
2F_{07}=G_{07}+G_{16}+G_{25}+G_{34} \\ 
2F_{16}=G_{07+}G_{16}-G_{25}-G_{34} \\ 
2F_{25}=G_{07}-G_{16}+G_{25}-G_{34} \\ 
2F_{34}=G_{07}-G_{16}-G_{25}+G_{34}%
\end{array}%
\right] $
}

In particular, $\{F_{ij} \mid 0\leq i<j\leq 7\}$ forms an $\R$-bases of \gd.

\bigskip

\emph{Definition 1.7.} \ We define $\R$-linear mappings $\kappa,%
\pi ,\nu :$\gd $\rightarrow $ \gd\  respectively by

$\ \ \ \ \ \ \ \ \ \ \ \ \ \ \ (\kappa D)x=\overline{D\overline{x}},x\in $
\gC ,

$\ \ \ \ \ \ \ \ \ \ \ \ \ \ \ \pi (G_{ij})=F_{ij},i,j=0,1,\cdot \cdot
\cdot ,7,i\neq j$,

$\ \ \ \ \ \ \ \ \ \ \ \ \ \ \ \nu =\pi \kappa $.

\bigskip

\emph{Proposition 1.8. }(I.Yokota\cite[\emph{Lemma1.3.2.}]{Yokota1}) \ The mappings $%
\kappa ,\pi ,\nu $ are automorphisms of
the Lie algebra \gd:

\ \ \ \ \ \ \ \ \ \ \ \ \ \ \ \ \ \ \ \ \ \ \ \ \ \ \ \ \ \ \ \ \ \ \ $%
\kappa ,\pi ,\nu \in Aut($\gd$)$.

\bigskip

\emph{Proposition 1.9.}  (I.Yokota\cite[\emph{Theorem1.3.5.}]{Yokota1}) In the
automorphism group 
\noindent
$Aut($\gd$) $\ of \gd,
the subgroup \gS$_{3}$ generated by $\kappa $ and $\pi 
$ is isomorphic to the symmetric group S$_{3}$ of
degree 3. \ Furthermore, $\kappa ,\pi ,\nu $ have the following
relations.

$\ \ \ \ \ \ \ \ \ \ \ \ \ \ \ \ \ \ \ \ \ \ \ \ \ \ \ \ \ \ \ \ \ \kappa
^{2}=1,\pi ^{2}=1,\nu ^{3}=1,\nu =\pi \kappa .$

\bigskip

\emph{Corollary 1.9.1. } By \emph{Proposition 1.9} , we have

\ \ \ \ \ \ \ \ \ \ \ \ \ \ $\kappa \pi =\kappa \pi \nu ^{3}=\kappa \pi (\pi
\kappa \pi \kappa \pi \kappa )=(\kappa \pi \pi \kappa )\pi \kappa \pi \kappa
=\pi \kappa \pi \kappa =\nu ^{2}.$

\bigskip

\emph{Lemma 1.10.} \ By definition of $\nu ,$ we have

\ \ \ $(1)\nu G_{ij}=F_{ij},1\leq i<j\leq 7,$

\ \ \ $(2)\nu G_{0j}=-F_{0j},1\leq j\leq 7.$

\bigskip

\emph{Proof.} $\ $For $x\in $ \gC, we have

$\ \ \ \ \ \ \ \ \ \ \ \ \ \ \ \ \ \ \ \ \ \ \ \ \ \ \ \ \ \kappa (G_{ij}x)=%
\overline{G_{ij}\overline{x}},$

\ \ \ \ \ \ \ \ \ \ \ \ \ \ \ \ \ \ \ \ \ \ \ \ \ \ \ \ \ \ \ \ \ \ \ \ \ \
\ \ $=\overline{e_{i}(e_{j},\overline{x})-e_{j}(e_{i},\overline{x})},$

\ \ \ \ \ \ \ \ \ \ \ \ \ \ \ \ \ \ \ \ \ \ \ \ \ \ \ \ \ \ \ \ \ \ \ \ \ \
\ \ =$\overline{e_{i}}(\overline{e_{j}},x)-\overline{e_{j}}(\overline{e_{i}}%
,x).$

$(1)$ For $1\leq i<j\leq 7,$ we have $\kappa (G_{ij}x)=G_{ij}x.$

Hence $\nu G_{ij}=\pi \kappa G_{ij}=\pi G_{ij}=F_{ij}.$

$(2)$ For $1\leq j\leq 7,$ we have $\kappa (G_{0j}x)=-G_{0j}x.$

Hence $\nu G_{0j}=\pi \kappa G_{0j}=-\pi G_{0j}=-F_{0j}.$\ \ \ \ \emph{%
Q.E.D.}

\bigskip

\emph{Proposition 1.11. } (I.Yokota\cite[\emph{Theorem1.3.6.}]{Yokota1})

For any $D_{1}\in $\gd, there exist $D_{2},D_{3}\in $\gd \ such
that

\ \ \ \ \ \ \ \ \ \ \ \ \ \ $\ (D_{1}x)y+x(D_{2}y)=D_{3}(xy),x,y\in $ \gC .

Also such $D_{2},D_{3}$ are uniquely determined for $D_{1}$ and we have

\ \ \ \ \ \ \ \ \ \ \ \ \ \ \ $D_{2}=\nu D_{1},D_{3}=\pi D_{1}.$

\bigskip

\emph{Definition 1.12. }
We denote $D_{ij}\in $\gD$_{4}$ $(i,j=0,1,\cdot \cdot \cdot ,7,i\neq
j)$ by

\ \ \ \ \ \ \ \ \ \ \ \ \ \ \ $D_{ij}=\left( 
\begin{array}{cccccccc}
0 & 0 & 0 & 0 & 0 & 0 & 0 & 0 \\ 
0 & 0 & 0 & 0 & 0 & 0 & 0 & 0 \\ 
0 & 0 & 0 & 0 & 1 & 0 & 0 & 0 \\ 
0 & 0 & 0 & 0 & 0 & 0 & 0 & 0 \\ 
0 & 0 & -1 & 0 & 0 & 0 & 0 & 0 \\ 
0 & 0 & 0 & 0 & 0 & 0 & 0 & 0 \\ 
0 & 0 & 0 & 0 & 0 & 0 & 0 & 0 \\ 
0 & 0 & 0 & 0 & 0 & 0 & 0 & 0%
\end{array}%
\right) ,$

\ \ \ \ \ \ \ \ \ \ \ \ \ \ \ \ \ \ \ \ \ \ \ element of $%
(i+1)th\ row,(j+1)th\ column$ is $1$,

\ \ \ \ \ \ \ \ \ \ \ \ \ \ \ \ \ \ \ \ \ \ \ element of $%
(j+1)th\ row,(i+1)th\ column$ is $-1$.

Then $\{D_{ij} \mid 0\leq i<j\leq 7\}$ forms an $\R$-bases of \gD$_{4}$.

\bigskip

\emph{Definition 1.13. } We define $\R$-linear mappings $\textrm{g}_{d}$:\gD$_{4} \rightarrow $ \gd\  and\ $\textrm{d}_{g}$:\gd 
$\rightarrow $ \gD$_{4}$ 

\noindent
respectively by

\ \ \ \ \ \ \ \ \ \ \ \ \ \ \ \ \ \ $\textrm{g}_{d}$\ :\ \gD$_{4}$ $\ \rightarrow $
\ \gd \ $\ \mathbf{;}\ \ \textrm{g}_{d}(D_{ij})=G_{ij},$

\ \ \ \ \ \ \ \ \ \ \ \ \ \ \ \ \ \ $\textrm{d}_{g}$\ :\ \gd\ \ $\ \rightarrow $
\ \gD$_{4}$ $\mathbf{;}\ \ \textrm{d}_{g}(G_{ij})=D_{ij}.$

\bigskip

\section{Definition of the vector space \gR$_{4}$ and it's Lie bracket}

\bigskip 

\ \ \ \emph{\ Proposition 2.1.} (I.Yokota\cite[\emph{Theorem 2.10.2}]{Yokota1}) 

\ \ \ \ \ \ $F_{4}=\{\alpha \in Iso_{\R}($\gJ$) \mid \alpha (X \circ Y)=\alpha X\circ 
\alpha Y\}$ 

\noindent
is a simply connected compact Lie group of type $F_{4}$.

\bigskip

\emph{Proposition 2.2. }(I.Yokota\cite[\emph{Theorem2.3.2.}]{Yokota1})

The Lie algebra  \gf$_{4}$ of the Lie group $F_{4}$ is given by

\ \ \ \ \ \ \ \ \gf$_{4}=\{\delta \in Hom_{\R}($%
\gJ$) \mid \delta (X\circ Y)=\delta X\circ Y+X\circ \delta Y$
\}.

\bigskip

\emph{Proposition 2.3. }(I.Yokota\cite[\emph{Theorem2.3.7.}]{Yokota1}) The Lie subalgebra 
\gd$_{4}$ of \gf$_{4}$:

\ \ \ \ \ \ \ \ \ \ \ \ \ \ \gd$_{4}=\{\delta \in $\gf$_{4} \mid %
\delta E_{i}=0,i=1,2,3\}$

\noindent
is isomorphic to the Lie algebra \gd\  under the correspondence

\ \ \ \ \ \ \ \ \ \ \ \ \ \ $\delta _{d}:$\gd$\ni D_{1}\rightarrow 
\delta _{d}(D_{1})\in $\gd$_{4}$

\noindent
given by

\ \ \ \ \ \ \ \ \ \ \ \ $\ \delta _{d}(D_{1})\left( 
\begin{array}{ccc}
\chi _{1} & x_{3} & \overline{x_{2}} \\ 
\overline{x_{3}} & \chi _{2} & x_{1} \\ 
x_{2} & \overline{x_{1}} & \chi _{3}%
\end{array}%
\right) =\left( 
\begin{array}{ccc}
0 & D_{3}x_{3} & \overline{D_{2}x_{2}} \\ 
\overline{D_{3}x_{3}} & 0 & D_{1}x_{1} \\ 
D_{2}x_{2} & \overline{D_{1}x_{1}} & 0%
\end{array}%
\right) ,$

\noindent
where $D_{2}$ and $D_{3}$ are elements of \gd\  which are determined by $D%
_{1}$ from the principle of triality :

$\ \ \ \ \ \ \ \ \ \ \ \ \ \ (D_{1}x)y+x(D_{2}y)=\overline{D_{3}(\overline{xy})},x,y\in $\gC ,

\ \ \ \ \ \ \ \ \ \ \ \ \ \ \ \ \ \ \ \ \ \ \ \ \ \ \ \ \ \ \ \ \ \ \ \ \ \ $%
=\kappa D_{3}(xy).$

\bigskip

\emph{Lemma 2.4. } We have

$\ \ \ \ \ \ \ \ \ \ \ \ \delta _{d}(D_{1})\left( 
\begin{array}{ccc}
\chi _{1} & x_{3} & \overline{x_{2}} \\ 
\overline{x_{3}} & \chi _{2} & x_{1} \\ 
x_{2} & \overline{x_{1}} & \chi _{3}%
\end{array}%
\right) =\left( 
\begin{array}{ccc}
0 & \nu ^{2}D_{1}x_{3} & \overline{\nu D_{1}x_{2}} \\ 
\overline{\nu ^{2}D_{1}x_{3}} & 0 & D_{1}x_{1} \\ 
\nu D_{1}x_{2} & \overline{D_{1}x_{1}} & 0%
\end{array}%
\right) .$

\bigskip

\emph{Proof.} \ Using \emph{Proposition 1.11 } we have the above $\kappa
D_{3}=\pi D_{1},$then

\ \ \ \ \ \ \ \ \ \ \ \ \ $D_{3}=\kappa \kappa D_{3}=\kappa (\kappa
D_{3})=\kappa (\pi D_{1})=\kappa \pi D_{1}=\nu ^{2}D_{1}.$ \ \ \ \ \ \emph{Q.E.D.}

\bigskip

\emph{Lemma 2.5. } For $D_{ij},D_{kl}\in $\gD$%
_{4},G_{ij},G_{kl},D_{1},D_{2}\in $\gd, We have

$\ \ \ \ \ \ \ \ \ \ [\textrm{g}_{d}(D_{ij}),\textrm{g}_{d}(D_{kl})]=\textrm{g}_{d}([D_{ij},D_{kl}]),$

$\ \ \ \ \ \ \ \ \ \ [\textrm{d}_{g}(G_{ij}),\textrm{d}_{g}(G_{kl})]=\textrm{d}_{g}([G_{ij},G_{kl}]),$

$\ \ \ \ \ \ \ \ \ \ [\delta _{d}(D_{1}),\delta _{d}(D_{2})]=\delta
_{d}([D_{1},D_{2}]).$

\bigskip

\emph{Proof.} Firstly, we have the followings by calculation.

\ \ \ \ \ \ \ $\ \ [D_{ij},D_{jk}]=D_{ik},\ \ \ \ i,j,k$ are distinct$,$

\ \ \ \ \ \ \ $\ \ [D_{ij},D_{kl}]=0,\ \ \ \ \ \ \ i,j,k,l$ are distinct,

\ \ \ \ \ \ \ $\ \ [G_{ij},G_{jk}]=G_{ik},\ \ \ \ \ i,j,k$ are distinct,

\ \ \ \ \ \ \ $\ \ [G_{ij},G_{kl}]=0,\ \ \ \ \ \ \ i,j,k,l$ are distinct.

For distinct $i,j,k$ , also we have

$\ \ \ \ \
[\textrm{g}_{d}(D_{ij}),\textrm{g}_{d}(D_{jk})]=[G_{ij},G_{jk}]=G_{ik}=\textrm{g}_{d}(D_{ik})=\textrm{g}_{d}([D_{ij},D_{jk}]), 
$

$\ \ \ \ \
[\textrm{d}_{g}(G_{ij}),\textrm{d}_{g}(G_{jk})]=[D_{ij},D_{jk}]=D_{ik}=\textrm{d}_{g}(G_{ik})=\textrm{d}_{g}([G_{ij},G_{jk}]). 
$

\noindent
Hence we have

$\ \ \ \ \ [\textrm{g}_{d}(D_{ij}),\textrm{g}_{d}(D_{kl})]=\textrm{g}_{d}([D_{ij},D_{kl}]),$

$\ \ \ \ \ [\textrm{d}_{g}(G_{ij}),\textrm{d}_{g}(G_{kl})]=\textrm{d}_{g}([G_{ij},G_{kl}]).$

For $X=X(\chi ,x)\in $\gJ$,$we have

$[\delta _{d}(D_{1}),\delta _{d}(D_{2})]X=\delta _{d}(D_{1})\delta
_{d}(D_{2})X-\delta _{d}(D_{2})\delta _{d}(D_{1})X$

\ \ \ =$\left( 
\begin{array}{ccc}
0 & \nu ^{2}D_{1}\nu ^{2}D_{2}x_{3} & \overline{\nu D_{1}\nu D_{2}x_{2}} \\ 
\overline{\nu ^{2}D_{1}\nu ^{2}D_{2}x_{3}} & 0 & D_{1}D_{2}x_{1} \\ 
\nu D_{1}\nu D_{2}x_{2} & \overline{D_{1}D_{2}x_{1}} & 0%
\end{array}%
\right) $

$\ -\left( 
\begin{array}{ccc}
0 & \nu ^{2}D_{2}\nu ^{2}D_{1}x_{3} & \overline{\nu D_{2}\nu D_{1}x_{2}} \\ 
\overline{\nu ^{2}D_{2}\nu ^{2}D_{1}x_{3}} & 0 & D_{2}D_{1}x_{1} \\ 
\nu D_{2}\nu D_{1}x_{2} & \overline{D_{2}D_{1}x_{1}} & 0%
\end{array}%
\right) ,$

\ \ =$\left( 
\begin{array}{ccc}
0 & \nu ^{2}[D_{1},D_{2}]x_{3} & \overline{\nu \lbrack D_{1},D_{2}]x_{2}} \\ 
\overline{\nu ^{2}[D_{1},D_{2}]x_{3}} & 0 & [D_{1},D_{2}]x_{1} \\ 
\nu \lbrack D_{1},D_{2}]x_{2} & \overline{[D_{1},D_{2}]x_{1}} & 0%
\end{array}%
\right) $ $($using $\nu \in Aut($\gd$)),$

\ \ =$\delta _{d}([D_{1},D_{2}]).$\ \ \ \ \emph{Q.E.D.}

\bigskip

\emph{Lemma 2.6. \ }For $x\in $\gC $,0\leq i<j\leq 7,$We have

\ \ \ \ \ \ \ $\ (1)\nu G_{ij}x=\frac{1}{4}(\overline{e_{i}}(e_{j}x)-%
\overline{e_{j}}(e_{i}x)),$

\ \ \ \ \ \ \ $\ (2)\nu ^{2}G_{ij}x=\frac{1}{4}((xe_{j})\overline{e_{i}}%
-(xe_{i})\overline{e_{j}}).$

\bigskip

\emph{Proof. }$(1)$ By \emph{Definition 1.5} ,we have

\ \ \ \ \ \ \ \ \ \ \ \ \ \ \ $\nu G_{ij}x=\nu
(e_{i}(e_{j},x)-e_{j}(e_{i},x))$

$\ \ \ \ \ \ \ \ \ \ \ \ \ \ \ \ \ \ \ \ \ \ =\pi \overline{(e_{i}(e_{j},%
\overline{x})-e_{j}(e_{i},\overline{x}))}$ (using $\nu =\pi \kappa )$

\ \ \ \ \ \ \ \ \ \ \ \ \ \ \ \ \ \ \ \ \ \ \ \ \ $=\pi (\overline{e_{i}}(%
\overline{e_{j}},x)-\overline{e_{j}}(\overline{e_{i}},x))$

\ \ \ \ \ \ \ \ \ \ \ \ \ \ \ \ \ \ \ \ \ \ \ \ \ $=\frac{1}{2}(\overline{%
e_{i}}(e_{j}x))$ \ \ (using \emph{Definition 1.7)}

$($i$)$If $i\neq 0,j\neq 0,$ we have

\ \ \ \ \ \ \ \ \ \ \ \ \ \ $\frac{1}{2}(\overline{e_{i}}(e_{j}x))=-\frac{1}{%
2}(e_{i}(e_{j}x))$

On the other hand, on the right side of $(1)$

\ \ \ \ \ \ \ \ \ \ \ \ \ $\frac{1}{4}(\overline{e_{i}}(e_{j}x)-\overline{%
e_{j}}(e_{i}x))=\frac{1}{4}(-e_{i}(e_{j}x)+e_{j}(e_{i}x))$

\ \ \ \ \ \ \ \ \ \ \ \ \ \ \ \ \ \ \ \ \ \ \ \ \ \ \ \ \ \ \ \ \ \ \ \ \ \ $%
\ \ \ =\frac{1}{4}(-e_{i}(e_{j}x)-e_{i}(e_{j}x))$

\ \ \ \ \ \ \ \ \ \ \ \ \ \ \ \ \ \ \ \ \ \ \ \ \ \ \ \ \ \ \ \ \ \ \ \ \ \
\ \ $\ =-\frac{1}{2}(e_{i}(e_{j}x))$

$($ii$)$If $i=0,j\neq 0,$we have

\ \ \ \ \ \ \ \ \ \ \ \ \ \ $\frac{1}{2}(\overline{e_{i}}(e_{j}x))=\frac{1}{2%
}(e_{j}x)$

On the other hand, on the right side of $(1)$

\ \ \ \ \ \ \ \ \ \ \ \ \ $\frac{1}{4}(\overline{e_{i}}(e_{j}x)-\overline{%
e_{j}}(e_{i}x))=\frac{1}{4}(e_{j}x+e_{j}x))$

$\ \ \ \ \ \ \ \ \ \ \ \ \ \ \ \ \ \ \ \ \ \ \ \ \ \ \ \ \ \ \ \ \ \ \ \ \ \
\ \ \ =\frac{1}{2}(e_{j}x)$

Thus we have the equation $(1)$.

$(2)$ By \emph{Definition 1.5} $,$we have

\ \ \ \ \ \ \ \ \ \ \ \ \ $\nu ^{2}G_{ij}x=\nu
^{2}(e_{i}(e_{j},x)-e_{j}(e_{i},x))$

$\ \ \ \ \ \ \ \ \ \ \ \ \ \ \ \ \ \ \ \ \ =\kappa (\frac{1}{2}e_{i}(%
\overline{e_{j}}x))$ (using $\nu ^{2}=\kappa \pi )$

\ \ \ \ \ \ \ \ \ \ \ \ \ \ \ \ \ \ \ \ \ \ \ =$\frac{1}{2}\overline{(%
\overline{e_{j}}\overline{x})}\overline{e_{i}}$

\ \ \ \ \ \ \ \ \ \ \ \ \ \ \ \ \ \ \ \ \ \ \ =$\frac{1}{2}(xe_{j})\overline{%
e_{i}}$

In the same way as proof $(1),$ we have the equation $(2).$\ \ \ \emph{\
Q.E.D.}

\bigskip

The Lie algebras \gD$_{4}$,\gd, and \gd$_{4}$ are
isomorphic to each other under the correspondences

\ \ \ \ \ \ \ \ \ \ \ \ \ \ \gD$_{4}\ni D_{ij}\rightarrow $ \gd$%
\ni \textrm{g}_{d}(D_{ij})=G_{ij}\rightarrow $ \gd$_{4}\ni \delta _{d}(G_{ij}).$

\bigskip

\emph{Definition 2.7.} \ We define an $\R$-vector space \gM$^{-}$
by

\ \ \ \ \ \ \ \ \ \ \ \ \ \ \ \ \ \ \ \gM$^{-}$ $=\{A\in M(3,$\gC $
) \mid A^{\ast }=-A\}.$

\noindent
Any element $A\in $\gM$^{-}$ induces an $\R$-linear mapping $\widehat{A}$%
: \gJ $\rightarrow$ \gJ\  defined by

\ \ \ \ \ \ \ \ \ \ \ \ \ \ \ \ \ \ \ $\widehat{A}X=\frac{1}{2}[A,X],X\in $%
\gJ.

\bigskip

In \gM$^{-}$, we adopt the following notations:

{\fontsize{8pt}{10pt} \selectfont%
$A_{1}(a)=\left( 
\begin{array}{ccc}
0 & 0 & 0 \\ 
0 & 0 & a \\ 
0 & -\overline{a} & 0%
\end{array}%
\right) ,$ $A_{2}(a)=\left( 
\begin{array}{ccc}
0 & 0 & -\overline{a} \\ 
0 & 0 & 0 \\ 
a & 0 & 0%
\end{array}%
\right) ,$ $A_{3}(a)=\left( 
\begin{array}{ccc}
0 & a & 0 \\ 
-\overline{a} & 0 & 0 \\ 
0 & 0 & 0%
\end{array}%
\right) $.
}

For $X(\chi ,x)\in $\gJ$,A=\sum\limits_{i=1}^{3}A_{i}(a_{i}),$the
explicit form in the term of $\widehat{A}X\in $\gJ\  entries is as
follows.

\ \ \ \ \ \ \ \ \ \ \ \ \ $\widehat{A}X=((a_{3},x_{3})-(a_{2},x_{2}))E_{1}$

$\ \ \ \ \ \ \ \ \ \ \ \ \ \ \ \ \ \ \ \ +((a_{1},x_{1})-(a_{3},x_{3}))E_{2}$

$\ \ \ \ \ \ \ \ \ \ \ \ \ \ \ \ \ \ \ \ +((a_{2},x_{2})-(a_{1},x_{1}))E_{3}$

\ \ \ \ \ \ \ \ \ \ \ \ \ \ \ \ \ \ \ $\ +\frac{1}{2}F_{1}((\chi _{3}-\chi
_{2})a_{1}-\overline{x_{2}a_{3}}+\overline{a_{2}x_{3}})$

$\ \ \ \ \ \ \ \ \ \ \ \ \ \ \ \ \ \ \ \ +\frac{1}{2}F_{2}((\chi _{1}-\chi
_{3})a_{2}-\overline{x_{3}a_{1}}+\overline{a_{3}x_{1}})$

$\ \ \ \ \ \ \ \ \ \ \ \ \ \ \ \ \ \ \ \ +\frac{1}{2}F_{3}((\chi _{2}-\chi
_{1})a_{3}-\overline{x_{1}a_{2}}+\overline{a_{1}x_{2}}).$

Let \gA \ be an $\R$-vector space defined by

\ \ \ \ \ \ \ \ \ \ \ \ \ \ \gA$=\{A\in $\gM$^{-} \mid diagA=0\}.$

\noindent
where $diagA=0$ means that all diagonal elements $a_{ii}$ of $A$ are 0.

\bigskip

\emph{Definition 2.8.} \ For $T\in $ \gJ, we define a
$\R$-linear mapping $\widetilde{T}:$ \gJ $\rightarrow$ \gJ by

\ \ \ \ \ \ \ \ \ \ \ \ \ $\widetilde{T}X=T\circ X,X\in $ \gJ.

\bigskip

\emph{Proposition 2.9.}(I.Yokota\cite[\emph{Theorem2.3.8.}]{Yokota1}) Any element $\delta
\in $\gf$_{4}$ is uniquely expressed by

\ \ \ \ \ \ \ \ \ \ \ \ \ $\delta =D+\widehat{A},D\in $\gd$_{4},A\in 
$\gA$,$

\bigskip

\emph{Proposition 2.10.} (I.Yokota\cite[\emph{Lemma2.4.5}]{Yokota1}) For $\delta \in $ \gf%
$_{4}^{\C}$ and $A,B\in $\gJ$^{\C}$, we have

\ \ \ \ \ \ \ \ \ \ \ $\ [\delta ,[\widetilde{A},\widetilde{B}]]=[%
\widetilde{\delta A},\widetilde{B}]+[\widetilde{A},\widetilde{\delta B}].$

\bigskip

\emph{Lemma 2.11. }For $a,b\in $\gC $,0\leq i<j\leq 7,$We have

\ \ \ \ \ \ \ \ \ \ $\ (1)[\widetilde{F}_{1}(e_{i}),\widetilde{F}%
_{1}(e_{j})]=\delta _{d}(G_{ij}),$

$\ \ \ \ \ \ \ \ \ \ \ (2)[\widetilde{F}_{2}(e_{i}),\widetilde{F}%
_{2}(e_{j})]=\delta _{d}(\nu ^{2}G_{ij}),$

$\ \ \ \ \ \ \ \ \ \ \ (3)[\widetilde{F}_{3}(e_{i}),\widetilde{F}%
_{3}(e_{j})]=\delta _{d}(\nu G_{ij}),$

$\ \ \ \ \ \ \ \ \ \ \ (4)[\widehat{A}_{i}(a),\widehat{A}_{i}(b)]=-[%
\widetilde{F}_{i}(a),\widetilde{F}_{i}(b)],$

\ \ \ \ $\ \ \ \ \ \ \ (5)[\widehat{A}_{i}(a),\widehat{A}_{i+1}(b)]=-\frac{1%
}{2}\widehat{A}_{i+2}(\overline{ab}),$

$\ \ \ \ \ \ \ \ \ \ \ (6)[\widehat{A}_{i}(a),\widehat{A}_{i+2}(b)]=\frac{1}{%
2}\widehat{A}_{i+1}(\overline{ba}),$

$\ \ \ \ \ \ \ \ \ \ \ (7)[\widetilde{F}_{i}(a),\widetilde{F}_{i+1}(b)]=-%
\frac{1}{2}\widehat{A}_{i+2}(\overline{ab}),$

$\ \ \ \ \ \ \ \ \ \ \ (8)[\widetilde{F}_{i}(a),\widetilde{F}_{i+2}(b)]=%
\frac{1}{2}\widehat{A}_{i+1}(\overline{ba}),$

$\ \ \ \ \ \ \ \ \ \ \ (9)[\widehat{A}_{i}(a),\widetilde{F}%
_{i}(b)]=(a,b)(E_{i+1}-E_{i+2})^{\widetilde{}},$

\ \ \ \ \ \ \ \ \ \ $\ (10)[\widehat{A}_{i}(a),\widetilde{F}_{i+1}(b)]=\frac{%
1}{2}\widetilde{F}_{i+2}(\overline{ab}),$

\ \ \ \ \ \ \ \ \ \ \ $(11)[\widehat{A}_{i}(a),\widetilde{F}_{i+2}(b)]=-%
\frac{1}{2}\widetilde{F}_{i+1}(\overline{ba}),$

$\ \ \ \ \ \ \ \ \ \ \ (12)[(E_{i}-E_{i+1})^{\widetilde{}},\widetilde{F}%
_{i}(a)]=-\frac{1}{2}\widehat{A}_{i}(a),$

$\ \ \ \ \ \ \ \ \ \ \ (13)[(E_{i}-E_{i+1})^{\widetilde{}},\widetilde{F}%
_{i+1}(a)]=-\frac{1}{2}\widehat{A}_{i+1}(a),$

$\ \ \ \ \ \ \ \ \ \ \ (14)[(E_{i}-E_{i+1})^{\widetilde{}},\widetilde{F}%
_{i+2}(a)]=\widehat{A}_{i+2}(a),$

\ \ \ \ \ \ \ $\ \ \ \ (15)[(E_{i}-E_{i+1})^{\widetilde{}},\widehat{A}%
_{i}(a)]=-\frac{1}{2}\widetilde{F}_{i}(a),$

\ \ \ \ \ \ $\ \ \ \ \ (16)[(E_{i}-E_{i+1})^{\widetilde{}},\widehat{A}%
_{i+1}(a)]=-\frac{1}{2}\widetilde{F}_{i+1}(a),$

\ \ \ \ \ \ \ \ $\ \ \ (17)[(E_{i}-E_{i+1})^{\widetilde{}},\widehat{A}%
_{i+2}(a)]=\widetilde{F}_{i+2}(a),$

\ \ \ \ \ \ \ \ $\ \ \ (18)[[\widetilde{F}_{1}(ei),\widetilde{F}_{1}(e_{j})],%
\widehat{A}_{1}(a)]=\widehat{A}_{1}(G_{ij}a),$

\ \ \ \ \ \ \ $\ \ \ \ (19)[[\widetilde{F}_{1}(ei),\widetilde{F}_{1}(e_{j})],%
\widehat{A}_{2}(a)]=\widehat{A}_{2}(\nu G_{ij}a),$

\ \ \ \ \ \ \ \ \ $\ \ (20)[[\widetilde{F}_{1}(ei),\widetilde{F}_{1}(e_{j})],%
\widehat{A}_{3}(a)]=\widehat{A}_{3}(\nu ^{2}G_{ij}a),$

\ \ \ \ \ \ \ \ $\ \ \ (21)[[\widetilde{F}_{2}(e_{i}),\widetilde{F}%
_{2}(e_{j})],\widehat{A}_{1}(a)]=\widehat{A}_{1}(\nu ^{2}G_{ij}a),$

\ \ \ \ \ \ \ \ \ $\ \ (22)[[\widetilde{F}_{2}(e_{i}),\widetilde{F}%
_{2}(e_{j})],\widehat{A}_{2}(a)]=\widehat{A}_{2}(G_{ij}a),$

\ \ \ \ \ \ $\ \ \ \ \ (23)[[\widetilde{F}_{2}(e_{i}),\widetilde{F}%
_{2}(e_{j})],\widehat{A}_{3}(a)]=\widehat{A}_{3}(\nu G_{ij}a),$

\ \ \ \ \ \ \ \ \ $\ \ (24)[[\widetilde{F}_{3}(e_{i}),\widetilde{F}%
_{3}(e_{j})],\widehat{A}_{1}(a)]=\widehat{A}_{1}(\nu G_{ij}a),$

\ \ \ \ \ \ \ \ $\ \ \ (25)[[\widetilde{F}_{3}(e_{i}),\widetilde{F}%
_{3}(e_{j})],\widehat{A}_{2}(a)]=\widehat{A}_{2}(\nu ^{2}G_{ij}a),$

\ \ \ \ \ \ $\ \ \ \ \ (26)[[\widetilde{F}_{3}(e_{i}),\widetilde{F}%
_{3}(e_{j})],\widehat{A}_{3}(a)]=\widehat{A}_{3}(G_{ij}a),$

\ \ \ \ \ \ \ $\ \ \ \ (27)[(E_{k}-E_{k+1})^{\widetilde{}},[\widetilde{F}%
_{1}(e_{i}),\widetilde{F}_{1}(e_{j})]]=0,$

\ \ \ \ \ \ \ \ \ \ \ $(28)[\widetilde{E}_{i},\widetilde{F}_{i}(a)]=0,$

\ \ \ \ \ \ \ \ \ \ \ $(29)[\widetilde{E}_{i},\widehat{A}_{i}(a)]=0.$

\bigskip

\emph{Proof. }For $X=X(\chi ,x)\in $\gJ ,we show that %
(the left side of equation)$X$ and

(the right side of equation)$X$ are equal.

$(1)$ $[\widetilde{F}_{1}(e_{i}),\widetilde{F}_{1}(e_{j})]X=F_{1}(e_{i})%
\circ (F_{1}(e_{j})\circ X)-F_{1}(e_{j})\circ (F_{1}(e_{i})\circ
X)$

\ \ \ \ \ \ \ \ \ \ \ \ \ \ \ \ \ \ \ \ \ \ \ \ \ \ $=\frac{1}{4}%
(F_{1}(e_{i})(F_{1}(e_{j})X+XF_{1}(e_{j}))+(F_{1}(e_{j})X+XF_{1}(e_{j}))F_{1}(e_{i})) 
$

\ \ \ \ \ \ \ \ \ \ \ \ \ \ \ \ \ \ \ \ \ \ \ \ $\ -\frac{1}{4}%
(F_{1}(e_{j})(F_{1}(e_{i})X+XF_{1}(e_{i}))+(F_{1}(e_{i})X+XF_{1}(e_{i}))F_{1}(e_{j})) 
$

\ \ \ \ \ \ \ \ \ \ \ \ \ \ \ \ \ \ \ \ $%
=F_{1}(e_{i}(e_{j},x_{1})-e_{j}(e_{i},x_{1}))+F_{2}(\frac{1}{4}(\overline{%
e_{i}}(e_{j}x_{2})-\overline{e_{j}}(e_{i}x_{2})))$

$\ \ \ \ \ \ \ \ \ \ \ \ \ \ \ \ \ \ \ \ \ +F_{3}(\frac{1}{4}%
((x_{3}e_{j})\overline{e_{i}}-(x_{3}e_{i})\overline{e_{j}}))$

\ \ \ \ \ \ \ \ \ \ \ \ \ \ \ \ \ \ \ \ $%
=F_{1}(G_{ij}x_{1})+F_{2}(\nu G_{ij}x_{2})+F_{3}(\nu ^{2}G_{ij}x_{3})$ (by 
\emph{Lemma 2.6})

\ \ \ \ \ \ \ \ \ \ \ \ \ \ \ \ \ \ \ \ $=\delta _{d}(G_{ij})X$

$(2)$ $[\widetilde{F}_{2}(e_{i}),\widetilde{F}_{2}(e_{j})]X=F_{1}(\frac{1}{4}%
((x_{1}e_{j})\overline{e_{i}}-(x_{1}e_{i})\overline{e_{j}}%
))+F_{2}(e_{i}(e_{j},x_{2})-e_{j}(e_{i},x_{2}))$

$\ \ \ \ \ \ \ \ \ \ \ \ \ \ \ \ \ \ \ \ \ \ \ \ \ \ +F_{3}(\frac{1}{4}(%
\overline{e_{i}}(e_{j}x_{3})-\overline{e_{j}}(e_{i}x_{3})))$

\ \ \ \ \ \ \ \ \ \ \ \ \ \ \ \ \ \ \ \ \ \ $=F_{1}(\nu
^{2}G_{ij}x_{1})+F_{2}(G_{ij}x_{2})+F_{3}(\nu G_{ij}x_{3})$ (by \emph{Lemma 2.6})

\ \ \ \ \ \ \ \ \ \ \ \ \ \ \ \ \ \ \ \ \ \ $=\delta _{d}(\nu
^{2}G_{ij})X$

$(3)$ $[\widetilde{F}_{3}(e_{i}),\widetilde{F}_{3}(e_{j})]X=F_{1}(\frac{1}{4}%
(\overline{e_{i}}(e_{j}x_{1})-\overline{e_{j}}(e_{i}x_{1})))+F_{2}(\frac{1}{4%
}((x_{2}e_{j})\overline{e_{i}}-(x_{2}e_{i})\overline{e_{j}}))$

$\ \ \ \ \ \ \ \ \ \ \ \ \ \ \ \ \ \ \ \ \ \ \ \ \ \ \
+F_{3}(e_{i}(e_{j},x_{3})-e_{j}(e_{i},x_{3}))$

\ \ \ \ \ \ \ \ \ \ \ \ \ \ \ \ \ \ \ \ \ \ \ $=F_{1}(\nu
G_{ij}x_{1})+F_{2}(\nu ^{2}G_{ij}x_{2})+F_{3}(G_{ij}x_{3})$ (by \emph{Lemma 2.6})

\ \ \ \ \ \ \ \ \ \ \ \ \ \ \ \ \ \ \ \ \ \ $=\delta _{d}(\nu
G_{ij})X$ ,

$(4)$ $[\widehat{A}_{1}(a),\widehat{A}%
_{1}(b)]X=-F_{1}(a(b,x_{1})-b(a,x_{1}))-F_{2}(\frac{1}{4}(\overline{a}%
(bx_{2})-\overline{b}(ax_{2})))$

$\ \ \ \ \ \ \ \ \ \ \ \ \ \ \ \ \ \ \ \ \ \ \ \ \ \ -F_{3}(\frac{1}{4}%
((x_{3}b)\overline{a}-(x_{3}a)\overline{b}))$

$\ \ \ \ \ [\widetilde{F}_{1}(a),\widetilde{F}%
_{1}(b)]X=F_{1}(a(b,x_{1})-b(a,x_{1}))+F_{2}(\frac{1}{4}(\overline{a}%
(bx_{2})-\overline{b}(ax_{2})))$

$\ \ \ \ \ \ \ \ \ \ \ \ \ \ \ \ \ \ \ \ \ \ \ \ \ \ +F_{3}(\frac{1}{4}((x_{3}b)%
\overline{a}-(x_{3}a)\overline{b}))$

Hence [$\widehat{A}_{1}(a),\widehat{A}_{1}(b)]=-[\widetilde{F}_{1}(a),%
\widetilde{F}_{1}(b)]$

Similarly we have [$\widehat{A}_{2}(a),\widehat{A}_{2}(b)]=-[\widetilde{F}%
_{2}(a),\widetilde{F}_{2}(b)]$ 

\ \ \ \ \ \ \ \ \ \ \ \ \ \ \ \ \ \ and [$\widehat{A}_{3}(a),\widehat{A}%
_{3}(b)]=-[\widetilde{F}_{3}(a),\widetilde{F}_{3}(b)].$

$(5)$ [$\widehat{A}_{1}(a),\widehat{A}_{2}(b)]X=-\frac{1}{2}(b,\overline{%
x_{3}a})E_{1}+\frac{1}{2}(a,\overline{bx_{3}})E_{2}+\frac{1}{2}((b,\overline{%
x_{3}a})-(a,\overline{bx_{3}}))E_{3}$

\ \ \ \ \ \ \ \ \ \ \ \ \ \ \ \ \ \ \ \ \ \ \ $+\frac{1}{2}F_{1}((b,x_{2})a-%
\frac{1}{2}((ax_{2})\overline{b})+\frac{1}{2}F_{2}(\frac{1}{2}\overline{a}%
(x_{1}b)-b(a,x_{1}))$

\ \ \ \ \ \ \ \ \ \ \ \ \ \ \ \ \ \ \ \ \ \ \ $+\frac{1}{2}F_{3}(\frac{1}{2}(\chi _{1}-\chi _{2})%
\overline{ab})$

(using $(b,\overline{a}x_{3})=(ab,x_{3}),(b,x_{2})a=\frac{1}{2}(ab)\overline{%
x_{2}}+\frac{1}{2}(ax_{2})\overline{b},$

\ \ \ \ \ \ \ \ \ $b(a,x_{1})=\frac{1}{2}\overline{x_{1}}%
(ab)+\frac{1}{2}\overline{a}(x_{1}b)$ )

\ \ \ \ \ \ \ \ \ \ \ \ \ \ \ \ \ \ \ \ \ \ \ \ $=-\frac{1}{2}(\overline{ab}%
,x_{3})E_{1}+\frac{1}{2}(\overline{ab},x_{3})E_{2}$

$\ \ \ \ \ \ \ \ \ \ \ \ \ \ \ \ \ \ \ \ \ \ \ \ +\frac{1}{2}F_{1}(\frac{1}{2%
}(ab)\overline{x_{2}})+\frac{1}{2}F_{2}(-\frac{1}{2}\overline{x_{1}}(ab))+%
\frac{1}{2}F_{3}(\frac{1}{2}(\chi _{1}-\chi _{2})\overline{ab})$

On the other hand,

\ \ \ \ \ \ \ \ \ \ \ \ $\widehat{A}_{3}(\overline{ab})X=(\overline{ab},x_{3})E_{1}-(\overline{ab}%
,x_{3})E_{2}$

\ \ \ \ \ \ \ \ \ \ \ \ \ \ \ \ \ \ \ \ \ \ \ \ $-F_{1}(\frac{1}{2}((ab)\overline{x_{2}})+F_{2}(\frac{1}{2}%
\overline{x_{1}}(ab))+F_{3}(\frac{1}{2}(\chi _{2}-\chi _{1})\overline{ab}))$

Hence [$\widehat{A}_{1}(a),\widehat{A}_{2}(b)]X=-\frac{1}{2}\widehat{A}_{3}(%
\overline{ab})X$

Similarly we have [$\widehat{A}_{2}(a),\widehat{A}_{3}(b)]X=-\frac{1}{2}%
\widehat{A}_{1}(\overline{ab})X$ 

\ \ \ \ \ \ \ \ \ \ \ \ \ \ \ \ \ \ and [$\widehat{A}_{3}(a),\widehat{A}%
_{1}(b)]X=-\frac{1}{2}\widehat{A}_{2}(\overline{ab})X.$

$(6)$ $[\widehat{A}_{i}(a),\widehat{A}_{i+2}(b)]=-[\widehat{A}_{i+2}(b),%
\widehat{A}_{i}(a)]$

\ \ \ \ \ \ \ \ \ \ \ \ \ \ \ \ \ \ \ \ \ \ \ \ \ \ =$-(-\frac{1}{2}\widehat{%
A}_{i+1}(\overline{ba}))$ (using $(5)$)

\ \ \ \ \ \ \ \ \ \ \ \ \ \ \ \ \ \ \ \ \ \ \ \ \ \ $=\frac{1}{2}\widehat{A}%
_{i+1}(\overline{ba})$

$(7) [\widetilde{F}_{1}(a),\widetilde{F}_{2}(b)]X=-\frac{1}{2}(b,\overline{%
x_{3}a})E_{1}+\frac{1}{2}(a,\overline{bx_{3}})E_{2}+\frac{1}{2}(-(b,%
\overline{x_{3}a})+(a,\overline{bx_{3}}))E_{3}$

\ \ \ \ \ \ \ \ \ \ \ \ \ \ \ \ \ \ \ \ \ \ \ \ \ \ \ \ $+\frac{1}{2}F_{1}((b,x_{2})a-%
\frac{1}{2}((ax_{2})\overline{b})+\frac{1}{2}F_{2}(\frac{1}{2}\overline{a}%
(x_{1}b)-b(a,x_{1}))$

\ \ \ \ \ \ \ \ \ \ \ \ \ \ \ \ \ \ \ \ \ \ \ \ \ \ \ \ $+\frac{1}{2}F_{3}(\frac{1}{2}(\chi _{1}-\chi _{2})%
\overline{ab})$

\ \ \ \ \ \ \ \ \ \ \ \ \ \ \ \ \ \ \ \ \ \ \ \ \ \ $=[\widehat{A}_{1}(a),\widehat{A}_{2}(b)]X ($ using $(5))$

\ \ \ \ \ \ \ \ \ \ \ \ \ \ \ \ \ \ \ \ \ \ \ \ \ \ $=-\frac{1}{2}\widehat{A}_{3}(\overline{ab})$

Similarly we have $[\widetilde{F}_{2}(a),\widetilde{F}_{3}(b)]X=-\frac{1}{2}%
\widehat{A}_{1}(\overline{ab})X$ 

\ \ \ \ \ \ \ \ \ \ \ \ \ \ \ \ \ \ and $[\widetilde{F}_{3}(a),\widetilde{F}%
_{1}(b)]X=-\frac{1}{2}\widehat{A}_{2}(\overline{ab})X$

$(8)$ $[\widetilde{F}_{i}(a),\widetilde{F}_{i+2}(b)]=-[\widetilde{F}%
_{i+2}(b),\widetilde{F}_{i}(a)]$

\ \ \ \ \ \ \ \ \ \ \ \ \ \ \ \ \ \ \ \ \ \ \ $\ =-(-\frac{1}{2}\widehat{A}%
_{i+1}(\overline{ba}))$ (using $(7)$)

\ \ \ \ \ \ \ \ \ \ \ \ \ \ \ \ \ \ \ \ \ \ \ $\ =\frac{1}{2}\widehat{A}%
_{i+1}(\overline{ba}))$

$(9)$ $[\widehat{A}_{1}(a),\widetilde{F}_{1}(b)]X=(a,b)\chi
_{2}E_{2}-(a,b)\chi _{3}E_{3}+\frac{1}{4}F_{2}(-\overline{a}(bx_{2})-%
\overline{b}(ax_{2}))$

\ \ \ \ \ \ \ \ \ \ \ \ \ \ \ \ \ \ \ \ \ \ \ \ $+\frac{1}{4}F_{3}((x_{3}b)\overline{a}+(x_{3}a)\overline{b})$

\ \ \ \ \ \ \ \ \ \ \ \ \ \ \ \ \ \ \ \ \ \ \ \ \ $=(a,b)\chi
_{2}E_{2}-(a,b)\chi _{3}E_{3}+\frac{1}{2}F_{2}(-(a,b)x_{2})+\frac{1}{2}%
F_{3}(x_{3}(a,b))$

\ \ \ \ \ \ \ \ \ \ \ \ \ \ \ \ \ \ \ \ \ \ $\ =(a,b)(\chi _{2}E_{2}-\chi
_{3}E_{3}+\frac{1}{2}F_{2}(-x_{2})+\frac{1}{2}F_{3}(x_{3}))$

On the other hand,

\ $(E_{2}-E_{3})^{\widetilde{}}X=\chi _{2}E_{2}-\chi _{3}E_{3}+\frac{1}{2}%
F_{2}(-x_{2})+\frac{1}{2}F_{3}(x_{3}))$

Hence $[\widehat{A}_{1}(a),\widetilde{F}_{1}(b)]X=(a,b)(E_{2}-E_{3})^{%
\widetilde{}}X$

Similarly we have $[\widehat{A}_{2}(a),\widetilde{F}%
_{2}(b)]X=(a,b)(E_{3}-E_{1})^{\widetilde{}}X$ and

$\ \ \ \ \ \ \ \ \ \ \ \ \ \ \ \ \ \ \ \ \ \ \ \ [\widehat{A}_{3}(a),%
\widetilde{F}_{3}(b)]X=(a,b)(E_{1}-E_{2})^{\widetilde{}}X.$

$(10)$ $[\widehat{A}_{1}(a),\widetilde{F}_{2}(b)]X=\frac{1}{2}(b,\overline{%
x_{3}a})E_{1}+\frac{1}{2}(a,\overline{bx_{3}})E_{2}+\frac{1}{2}((b,\overline{%
x_{3}a})-(a,\overline{bx_{3}}))E_{3}$

\ \ \ \ \ \ \ \ \ \ \ \ \ \ \ \ \ \ \ \ \ \ \ \ \ \ \ \ $+\frac{1}{2}F_{1}(a(b,x_{2})-\frac{1}{2}%
(ax_{2})\overline{b})+\frac{1}{2}F_{2}(-\frac{1}{2}\overline{a}%
(x_{1}b)+(a,x_{1})b)$

\ \ \ \ \ \ \ \ \ \ \ \ \ \ \ \ \ \ \ \ \ \ \ \ \ \ \ \ $+\frac{1}{4}F_{3}((\chi _{1}+\chi _{2})\overline{ab})$

\ \ \ \ \ \ \ \ \ \ \ \ \ \ \ \ \ \ \ \ \ \ \ \ \ \ \ \ =$\frac{1}{2}(\overline{b},x_{3}a)E_{1}+\frac{1%
}{2}(\overline{a},bx_{3})E_{2}+\frac{1}{2}((\overline{b},x_{3}a)-(\overline{a%
},bx_{3}))E_{3}$

\ \ \ \ \ \ \ \ \ \ \ \ \ \ \ \ \ \ \ \ \ \ \ \ \ \ \ \ $+\frac{1}{2}F_{1}((b,x_{2})a-\frac{1}{2}%
(ax_{2})\overline{b})+\frac{1}{2}F_{2}(-\frac{1}{2}\overline{a}%
(x_{1}b)+b(a,x_{1}))$

\ \ \ \ \ \ \ \ \ \ \ \ \ \ \ \ \ \ \ \ \ \ \ \ \ \ \ \ $+\frac{1}{4}F_{3}((\chi _{1}+\chi _{2})\overline{ab})$

(using $(\overline{b},x_{3}a)=(\overline{b}\overline{a},x_{3}),(b,x_{2})a=%
\frac{1}{2}(ab)\overline{x_{2}}+\frac{1}{2}(ax_{2})\overline{b},$

\ \ \ \ \ \ \ \ \ $b(a,x_{1})=%
\frac{1}{2}\overline{x_{1}}(ab)+\frac{1}{2}\overline{a}(x_{1}b)$ )

\ \ \ \ \ \ \ \ \ \ \ \ \ \ \ =$\frac{1}{2}(\overline{ab},x_{3})E_{1}+\frac{1%
}{2}(\overline{ab},x_{3})E_{2}+\frac{1}{2}F_{1}(\frac{1}{2}(ab)\overline{%
x_{2}})+\frac{1}{2}F_{2}(\frac{1}{2}\overline{x_{1}}(ab))$

\ \ \ \ \ \ \ \ \ \ \ \ \ \ \ $+\frac{1}{2}F_{3}(\frac{1}{2}(\chi _{1}+\chi _{2})\overline{ab})$

On the other hand,

\ $\widetilde{F}_{3}(\overline{ab})X=(\overline{ab},x_{3})E_{1}+(\overline{ab%
},x_{3})E_{2}+F_{1}(\frac{1}{2}(ab)\overline{x_{2}})+F_{2}(\frac{1}{2}%
\overline{x_{1}}(ab))$

\ \ \ \ \ \ \ \ \ \ \ \ \ \ $+F_{3}(\frac{1}{2}((\chi _{1}+\chi _{2})\overline{ab}))$

Hence $[\widehat{A}_{1}(a),\widetilde{F}_{2}(b)]X=\frac{1}{2}\widetilde{F}%
_{3}(\overline{ab})X$

Similarly we have $[\widehat{A}_{2}(a),\widetilde{F}_{3}(b)]X=\frac{1}{2}%
\widetilde{F}_{1}(\overline{ab})X$ 

\ \ \ \ \ \ \ \ \ \ \ \ \ \ \ \ \ \ and $[\widehat{A}_{3}(a),\widetilde{F}%
_{1}(b)]X=\frac{1}{2}\widetilde{F}_{2}(\overline{ab})X$

$(11)$ $[\widehat{A}_{1}(a),\widetilde{F}_{3}(b)]X=-\frac{1}{2}(b,\overline{%
ax_{2}})E_{1}+\frac{1}{2}(-(b,\overline{ax_{2}})+(a,\overline{x_{2}b}))E_{2}$

\ \ \ \ \ \ \ \ \ \ \ \ \ \ \ \ \ \ \ \ \ \ \ \ \ \ \ \ \ \ 
$+\frac{1}{2}(-(a,\overline{x_{2}b})E_{3}$+$\frac{1}{2}F_{1}(-a(b,x_{3})+\frac{1}{2}\overline{b}(x_{3}a))$

\ \ \ \ \ \ \ \ \ \ \ \ \ \ \ \ \ \ \ \ \ \ \ \ \ \ \ \ \ \ 
$+\frac{1}{4}F_{2}(-\overline{ba}(\chi_{1}+\chi _{3}))$
$+\frac{1}{2}F_{3}(\frac{1}{2}(bx_{1})\overline{a}-(a,x_{1})b)$

(using $(b,\overline{ax_{2}})=(\overline{ba},x_{2}),a(b,x_{3})=\frac{1}{2}(%
\overline{x_{3}}(ba)+\overline{b}(x_{3}a)),$

\ \ \ \ \ \ \ \ \ $(a,x_{1})b=\frac{1}{2}((ba)%
\overline{x_{1}}+(bx_{1})\overline{a})$ )

\ \ \ \ \ \ \ \ \ \ \ 
\ \ \ \ \ \ \ \ \ \ \ \ \ \ \ \ $=-\frac{1}{2}(\overline{ba},x_{2})E_{1}-%
\frac{1}{2}(\overline{ba},x_{2})E_{3}-\frac{1}{4}F_{1}(\overline{x_{3}}(ba))$

\ \ \ \ \ \ \ \ \ \ \ \ \ \ \ \ \ \ \ \ \ \ \ \ \ \ \ \ \ \ 
$-\frac{1}{4}F_{2}((\chi _{1}+\chi _{3})\overline{ba})-\frac{1}{4}F_{3}((ba)%
\overline{x_{1}})$

On the other hand,

$\widetilde{F}_{2}(\overline{ba})X=(\overline{ba},x_{2})E_{1}+(\overline{ba}%
,x_{2})E_{3}+\frac{1}{2}F_{1}(\overline{x_{3}}(ba))+\frac{1}{2}F_{2}((\chi
_{1}+\chi _{3})\overline{ba})$

\ \ \ \ \ \ \ \ \ \ \ \ \ \ $+\frac{1}{2}F_{3}((ba)\overline{x_{1}})$

Hence $[\widehat{A}_{1}(a),\widetilde{F}_{3}(b)]X=-\frac{1}{2}\widetilde{F}%
_{2}(\overline{ba})X$

Similarly we have $[\widehat{A}_{2}(a),\widetilde{F}_{1}(b)]X=-\frac{1}{2}%
\widetilde{F}_{3}(\overline{ba})X$ 

\ \ \ \ \ \ \ \ \ \ \ \ \ \ \ \ \ \ 
and $[\widehat{A}_{3}(a),\widetilde{F}%
_{2}(b)]X=-\frac{1}{2}\widetilde{F}_{1}(\overline{ba})X$

$(12)$ $[(E_{1}-E_{2})^{\widetilde{}},\widetilde{F}_{1}(a)]X=-\frac{1}{2}%
(a,x_{1})E_{2}+\frac{1}{2}(a,x_{1})E_{3}$

$\ \ \ \ \ \ \ \ \ \ \ \ \ \ \ \ \ \ \ \ \ \ \ \ \ \ \ \ \ \ \ \ \ \ +\frac{1%
}{2}F_{1}(\frac{1}{2}(\chi _{2}-\chi _{3})a)+\frac{1}{2}F_{2}(\frac{1}{2}%
\overline{x_{3}a})+\frac{1}{2}F_{3}(-\frac{1}{2}\overline{ax_{2}})$

On the other hand,

$\widehat{A}_{1}(a)X=(a,x_{1})E_{2}-(a,x_{1})E_{3}+\frac{1}{2}F_{1}((\chi
_{3}-\chi _{2})a)-\frac{1}{2}F_{2}(\overline{x_{3}a})+\frac{1}{2}F_{3}(%
\overline{ax_{2}})$

Hence $[(E_{1}-E_{2})^{\widetilde{}},\widetilde{F}_{1}(a)]X=-\frac{1}{2}%
\widehat{A}_{1}(a)X.$

Similarly we have $[(E_{2}-E_{3})^{\widetilde{}},\widetilde{F}_{2}(a)]X=-%
\frac{1}{2}\widehat{A}_{2}(a)X$ 

\ \ \ \ \ \ \ \ \ \ \ \ \ \ \ \ \ \ 
and $[(E_{3}-E_{1})^{\widetilde{}},%
\widetilde{F}_{3}(a)]X=-\frac{1}{2}\widehat{A}_{3}(a)X.$

$(13)$ $[(E_{1}-E_{2})^{\widetilde{}},\widetilde{F}_{2}(a)]X=\frac{1}{2}%
(a,x_{2})E_{1}-\frac{1}{2}(a,x_{2})E_{3}$

$\ \ \ \ \ \ \ \ \ \ \ \ \ \ \ \ \ \ \ \ \ \ \ \ \ \ \ \ \ \ \ \ +\frac{1}{2}%
F_{1}(-\frac{1}{2}\overline{ax_{3}})+\frac{1}{2}F_{2}(\frac{1}{2}(\chi
_{3}-\chi _{1})a)+\frac{1}{2}F_{3}(\frac{1}{2}\overline{x_{1}a})$

On the other hand,

$\widehat{A}_{2}(a)X=-(a,x_{2})E_{1}+(a,x_{2})E_{3}+F_{1}(\frac{1}{2}%
ax_{3})+F_{2}(\frac{1}{2}(\chi _{1}-\chi _{3})a)+F_{3}(-\frac{1}{2}\overline{%
x_{1}a})$

Hence $[(E_{1}-E_{2})^{\widetilde{}},\widetilde{F}_{2}(a)]X=-\frac{1}{2}%
\widehat{A}_{2}(a)X.$

Similarly we have $[(E_{2}-E_{3})^{\widetilde{}},\widetilde{F}_{3}(a)]X=-%
\frac{1}{2}\widehat{A}_{3}(a)X$ 

\ \ \ \ \ \ \ \ \ \ \ \ \ \ \ \ \ \ 
and $[(E_{3}-E_{1})^{\widetilde{}},%
\widetilde{F}_{1}(a)]X=-\frac{1}{2}\widehat{A}_{1}(a)X.$

$(14)$ $[(E_{1}-E_{2})^{\widetilde{}},\widetilde{F}%
_{3}(a)]X=(a,x_{3})E_{1}-(a,x_{3})E_{2}+\frac{1}{2}F_{1}(-\overline{x_{2}a})+%
\frac{1}{2}F_{2}(\overline{ax_{1}})$

\ \ \ \ \ \ \ \ \ \ \ \ \ \ \ \ \ \ \ \ \ \ \ \ \ \ \ \ \ \ \ \ \ \ \ 
$+\frac{1}{2}F_{3}((\chi _{2}-\chi _{1})a)$

On the other hand,

$\widehat{A}_{3}(a)X=(a,x_{3})E_{1}-(a,x_{3})E_{2}+\frac{1}{2}F_{1}(-%
\overline{x_{2}a})+\frac{1}{2}F_{2}(\overline{ax_{1}})+\frac{1}{2}F_{3}((\chi
_{2}-\chi _{1})a)$

Hence $[(E_{1}-E_{2})^{\widetilde{}},\widetilde{F}_{3}(a)]X=\widehat{A}%
_{3}(a)X.$

Similarly we have $[(E_{2}-E_{3})^{\widetilde{}},\widetilde{F}_{1}(a)]X=%
\widehat{A}_{1}(a)X$ 

\ \ \ \ \ \ \ \ \ \ \ \ \ \ \ \ \ \ 
and $[(E_{3}-E_{1})^{\widetilde{}},\widetilde{F}%
_{2}(a)]X=\widehat{A}_{2}(a)X.$

$(15)$ $[(E_{1}-E_{2})^{\widetilde{}},\widehat{A}_{1}(a)]X=-\frac{1}{2}%
(a,x_{1})E_{2}-\frac{1}{2}(a,x_{1})E_{3}$

$\ \ \ \ \ \ \ \ \ \ \ \ \ \ \ \ \ \ \ \ \ \ \ \ \ \ \ \ \ \ \ \ \ \ \ \ -%
\frac{1}{4}F_{1}((\chi _{2}+\chi _{3})a)-\frac{1}{4}F_{2}(\overline{x_{3}a})-%
\frac{1}{4}F_{3}(\overline{ax_{2})}$

On the other hand,

$\widetilde{F}_{1}(a)X=(a,x_{1})E_{2}+(a,x_{1})E_{3}+\frac{1}{2}F_{1}((\chi
_{2}+\chi _{3})a)+\frac{1}{2}F_{2}(\overline{x_{3}a})+\frac{1}{2}F_{3}(%
\overline{ax_{2}}).$

Hence $[(E_{1}-E_{2})^{\widetilde{}},\widehat{A}_{1}(a)]X=-\frac{1}{2}%
\widetilde{F}_{1}(a)X.$

Similarly we have $[(E_{2}-E_{3})^{\widetilde{}},\widehat{A}_{2}(a)]X=-\frac{%
1}{2}\widetilde{F}_{2}(a)X$ 

\ \ \ \ \ \ \ \ \ \ \ \ \ \ \ \ \ \ 
and $[(E_{3}-E_{1})^{\widetilde{}},\widehat{A}%
_{3}(a)]X=-\frac{1}{2}\widetilde{F}_{3}(a)X.$

$(16)$ $[(E_{1}-E_{2})^{\widetilde{}},\widehat{A}_{2}(a)]X=-\frac{1}{2}%
(a,x_{2})E_{1}-\frac{1}{2}(a,x_{2})E_{3}$

$\ \ \ \ \ \ \ \ \ \ \ \ \ \ \ \ \ \ \ \ \ \ \ \ \ \ \ \ \ \ \ \ \ \ \ \ \ -%
\frac{1}{4}F_{1}(\overline{ax_{3}})-\frac{1}{4}F_{2}((\chi _{1}+\chi _{3})a)-%
\frac{1}{4}F_{3}(\overline{x_{1}a})$

On the other hand,

$\widetilde{F}_{2}(a)X=(a,x_{2})E_{1}+(a,x_{2})E_{3}+\frac{1}{2}F_{1}(%
\overline{ax_{3}})+\frac{1}{2}F_{2}((\chi _{1}+\chi _{3})a)+\frac{1}{2}F_{3}(%
\overline{x_{1}a}).$

Hence $[(E_{1}-E_{2})^{\widetilde{}},\widehat{A}_{2}(a)]X=-\frac{1}{2}%
\widetilde{F}_{2}(a)X.$

Similarly we have $[(E_{2}-E_{3})^{\widetilde{}},\widehat{A}_{3}(a)]X=-\frac{%
1}{2}\widetilde{F}_{3}(a)X$ 

\ \ \ \ \ \ \ \ \ \ \ \ \ \ \ \ \ \ 
and $[(E_{3}-E_{1})^{\widetilde{}},\widehat{A}%
_{1}(a)]X=-\frac{1}{2}\widetilde{F}_{1}(a)X.$

$(17)$ $[(E_{1}-E_{2})^{\widetilde{}},\widehat{A}%
_{3}(a)]X=(a,x_{3})E_{1}+(a,x_{3})E_{2}+\frac{1}{2}F_{1}(\overline{x_{2}a})+%
\frac{1}{2}F_{2}(\overline{ax_{1}})$

\ \ \ \ \ \ \ \ \ \ \ \ \ \ \ \ \ \ \ \ \ \ \ \ \ \ \ \ \ \ \ \ \ \ \ 
$+\frac{1}{2}F_{3}((\chi _{1}+\chi _{2})a)$

On the other hand,

$\widetilde{F}_{3}(a)X=(a,x_{3})E_{1}+(a,x_{3})E_{2}+\frac{1}{2}F_{1}(%
\overline{x_{2}a})+\frac{1}{2}F_{2}(\overline{ax1})+\frac{1}{2}F_{3}((\chi
_{1}+\chi _{2})a)$

Hence $[(E_{1}-E_{2})^{\widetilde{}},\widehat{A}_{3}(a)]X=\widetilde{F}%
_{3}(a)X.$

Similarly we have $[(E_{2}-E_{3})^{\widetilde{}},\widehat{A}_{1}(a)]X=%
\widetilde{F}_{1}(a)X$ 

\ \ \ \ \ \ \ \ \ \ \ \ \ \ \ \ \ \ 
and $[(E_{3}-E_{1})^{\widetilde{}},\widehat{A}%
_{2}(a)]X=\widetilde{F}_{2}(a)X.$

$(18)$ By \emph{Proposition 2.10}, We have

\ \ $[[\widetilde{F}_{1}(ei),\widetilde{F}_{1}(e_{j})],\widehat{A}%
_{1}(a)]=-[(\widehat{A}_{1}(a)F_{1}(e_{i}))^{\widetilde{}},\widetilde{F}%
_{1}(e_{j})]-[\widetilde{F}_{1}(e_{i}),(\widehat{A}_{1}(a)F_{1}(e_{j}))^{%
\widetilde{}}]$

\ \ \ \ \ \ \ \ \ \ \ \ \ \ \ \ \ \ \ \ \ \ \ \ \ \ \ \ \ \ \ \ \ $\
=-[(a,e_{i})(E_{2}-E_{3})^{\widetilde{}},\widetilde{F}_{1}(e_{j})]-[%
\widetilde{F}_{1}(e_{i}),(a,e_{j})(E_{2}-E_{3})^{\widetilde{}}]$

\ \ \ \ \ \ \ \ \ \ \ \ \ \ \ \ \ \ \ \ \ \ \ \ \ \ \ \ \ \ \ \ \ \ $\
=-(a,e_{i})[(E_{2}-E_{3})^{\widetilde{}},\widetilde{F}%
_{1}(e_{j})]+(a,e_{j})[(E_{2}-E_{3})^{\widetilde{}},\widetilde{F}%
_{1}(e_{i})] $

\ \ \ \ \ \ \ \ \ \ \ \ \ \ \ \ \ \ \ \ \ \ \ \ \ \ \ $\
=-[(E_{2}-E_{3})^{\widetilde{}},\widetilde{F}%
_{1}(e_{j}(a,e_{i})-e_{i}(a,e_{j}))]$

\ \ \ \ \ \ \ \ \ \ \ \ \ \ \ \ \ \ \ \ \ \ \ \ \ \ $\ \
=[(E_{2}-E_{3})^{\widetilde{}},\widetilde{F}_{1}(G_{ij}a)]$

\ \ \ \ \ \ \ \ \ \ \ \ \ \ \ \ \ \ \ \ \ \ \ \ \ \ \ $\ =%
\widehat{A}1(G_{ij}(a))\ ($using$(14))$

$(19)$ By \emph{Proposition 2.9}, We have

\ \ $[[\widetilde{F}_{1}(ei),\widetilde{F}_{1}(e_{j})],\widehat{A}%
_{2}(a)]=-[(\widehat{A}_{2}(a)F_{1}(e_{i}))^{\widetilde{}},\widetilde{F}%
_{1}(e_{j})]-[\widetilde{F}_{1}(e_{i}),(\widehat{A}_{2}(a)F_{1}(e_{j}))^{%
\widetilde{}}]$

\ \ \ \ \ \ \ \ \ \ \ \ \ \ \ \ \ \ \ \ \ \ \ \ \ \ \ \ \ \ \ \ 
$=[\widetilde{F}_{3}(\frac{1}{2}\overline{e_{i}a}),\widetilde{F}_{1}(e_{j})]+[%
\widetilde{F}_{1}(e_{i}),\widetilde{F}_{3}(\frac{1}{2}\overline{e_{j}a})]$

\ \ \ \ \ \ \ \ \ \ \ \ \ \ \ \ \ \ \ \ \ \ \ \ \ \ \ \ \ \ \ \ $=-%
\widehat{A}_{2}(\frac{1}{4}(\overline{(\overline{e_{i}a})e_{j}}-\overline{(%
\overline{e_{j}a})e_{i}}))$ \ (using $(7)$)

\ \ \ \ \ \ \ \ \ \ \ \ \ \ \ \ \ \ \ \ \ \ \ \ \ \ \ \ \ \ \ 
\ $=\widehat{A}_{2}(\frac{1}{4}(\overline{e_{i}}(e_{j}a)-\overline{e_{j}}(e_{i}a)))$

\ \ \ \ \ \ \ \ \ \ \ \ \ \ \ \ \ \ \ \ \ \ \ \ \ \ \ \ \ \ \ \ 
$=\widehat{A}_{2}(\nu G_{ij}a)$ \ (by\ \emph{Proposition 1.6)}

$(20)$ By \emph{Lemma 2.6}, We have

 \ \ $[[\widetilde{F}_{1}(ei),\widetilde{F}_{1}(e_{j})],\widehat{A}%
_{3}(a)]=-[(\widehat{A}_{3}(a)F_{1}(e_{i}))^{\widetilde{}},\widetilde{F}%
_{1}(e_{j})]-[\widetilde{F}_{1}(e_{i}),(\widehat{A}_{3}(a)F_{1}(e_{j}))^{%
\widetilde{}}]$

\ \ \ \ \ \ \ \ \ \ \ \ \ \ \ \ \ \ \ \ \ \ \ \ \ \ \ \ \ \ \ \ $=-[%
\widetilde{F}_{2}(\frac{1}{2}\overline{ae_{i}}),\widetilde{F}_{1}(e_{j})]-[%
\widetilde{F}_{1}(e_{i}),\widetilde{F}_{2}(\frac{1}{2}\overline{ae_{j}})]$

\ \ \ \ \ \ \ \ \ \ \ \ \ \ \ \ \ \ \ \ \ \ \ \ \ \ \ \ \ \ \ \ 
$=\widehat{A}_{3}(\frac{1}{4}(-(\overline{e_{j}(\overline{ae_{i}}}))+\overline{%
e_{i}(\overline{ae_{j}})}))$ \ (using $(8)$)

\ \ \ \ \ \ \ \ \ \ \ \ \ \ \ \ \ \ \ \ \ \ \ \ \ \ \ \ \ \ \ \ $=%
\widehat{A}_{3}(\frac{1}{4}((ae_{j})\overline{e_{i}}-(ae_{i})\overline{e_{j}}%
))$

\ \ \ \ \ \ \ \ \ \ \ \ \ \ \ \ \ \ \ \ \ \ \ \ \ \ \ \ \ \ \ \ \ \ \ \ 
$=\widehat{A}_{3}(\nu ^{2}G_{ij}a)$ \ (by \emph{Lemma 2.6)}

$(21),(22),(23),(24),(25)$ and $(26)$

Calculate in the same way as above $(18)$,$(19)$ and $(20)$ ,we have the
equations.

$(27)$ $[(E_{1}-E_{2})^{\widetilde{}},[\widetilde{F}_{1}(e_{i}),\widetilde{F}%
_{1}(e_{j})]]X=(E_{1}-E_{2})^{\widetilde{}}([\widetilde{F}_{1}(e_{i}),%
\widetilde{F}_{1}(e_{j})]X)$

$\ \ \ \ \ \ \ \ \ \ \ \ \ \ \ \ \ \ \ \ \ \ \ \ \ \ \ \ \ \ \ \ \ \ \ \ \ \
\ \ \ \ \ \ -[\widetilde{F}_{1}(e_{i}),\widetilde{F}%
_{1}(e_{j})]((E_{1}-E_{2})^{\widetilde{}}$X$)$

\ \ \ \ \ \ \ \ \ \ \ \ \ \ \ \ \ $=(E_{1}-E_{2})^{%
\widetilde{}}(F_{1}(G_{ij}(x_{1}))+F_{2}(\nu G_{ij}(x_{2}))+F_{3}(\nu
^{2}G_{ij}(x_{3})))$

\ \ \ \ \ \ \ \ \ \ \ \ \ \ \ \ \ $\ -[\widetilde{F}%
_{1}(e_{i}),\widetilde{F}_{1}(e_{j})](\chi _{1}E_{1}-\chi _{2}E_{2}+F_{1}(-%
\frac{1}{2}x_{1})+F_{2}(\frac{1}{2}x_{2}))$

\ \ \ \ \ \ \ \ \ \ \ \ \ \ \ \ \ \ \ \ \ \ \ \ \ \ \ \  (using $(1)$)

\ \ \ \ \ \ \ \ \ \ \ \ \ \ \ \ \ \ \ \ \ \ $=F_{1}(-\frac{1}{2}%
G_{ij}x_{1})+F_{2}(\frac{1}{2}\nu G_{ij}x_{2})-(F_{1}(-\frac{1}{2}%
G_{ij}x_{1})+F_{2}(\frac{1}{2}\nu G_{ij}x_{2}))$

\ \ \ \ \ \ \ \ \ \ \ \ \ \ \ \ \ $=0$

Similarly we have $[(E_{2}-E_{3})^{\widetilde{}},[\widetilde{F}_{1}(e_{i}),%
\widetilde{F}_{1}(e_{j})]]X=0$ 

\ \ \ \ \ \ \ \ \ \ \ \ \ \ \ \ \ \ 
and $[(E_{3}-E_{1})^{\widetilde{}},[%
\widetilde{F}_{1}(e_{i}),\widetilde{F}_{1}(e_{j})]]X=0.$

$(28)[\widetilde{E}_{1},\widetilde{F}_{1}(a)]X=\frac{1}{2}\widetilde{E}%
_{1}(F_{1}(a)X+XF_{1}(a))-\frac{1}{2}\widetilde{F}_{1}(a)(E_{1}X+XE_{1})$

\ \ \ \ \ \ \ \ \ \ \ \ \ \ \ \ \ \ \ \ \ \ \ \ $=\frac{1}{2}\widetilde{E}%
_{1}(2(a,x_{1})E_{2}+2(a,x_{2})E_{3}+F_{1}((\chi _{2}+\chi _{3})a)+F_{2}(%
\overline{x_{3}a})$

\ \ \ \ \ \ \ \ \ \ \ \ \ \ \ \ \ \ \ \ \ \ \ \ \ 
$+F_{3}(\overline{ax_{2}}))-\frac{1}{2}\widetilde{F%
}_{1}(a)(2\chi _{1}E_{1}+F_{2}(x_{2})+F_{3}(x_{3}))$

\ \ \ \ \ \ \ \ \ \ \ \ \ \ \ \ \ \ \ \ \ $\ =\frac{1}{4}(F_{2}(\overline{%
x_{3}a})+F_{3}(\overline{ax_{2}}))-\frac{1}{4}(F_{2}(\overline{x_{3}a}%
)+F_{3}(\overline{ax_{2}}))$

\ \ \ \ \ \ \ \ \ \ \ \ \ \ \ \ \ \ \ \ \ $\ =0$

Similarly we have $[\widetilde{E}_{2},\widetilde{F}_{2}(a)]=0$ and $[%
\widetilde{E}_{3},\widetilde{F}_{3}(a)]=0$ .

\ $(29)[\widetilde{E}_{1},\widehat{A}_{1}(a)]X=\frac{1}{2}\widetilde{E}%
_{1}(2(a,x_{1})E_{2}-2(a,x_{1})E_{3}+F_{1}((\chi _{3}-\chi _{2})a)-F_{2}(%
\overline{x_{3}a})$

\ \ \ \ \ \ \ \ \ \ \ \ \ \ \ \ \ \ \ \ \ \ \ \ \ \ \ \ 
$+F_{3}(\overline{ax_{2}}))-\frac{1}{2}\widehat{A}%
_{1}(a)(2\chi _{1}E_{1}+F_{2}(x_{2})+F_{3}(x_{3}))$

\ \ \ \ \ \ \ \ \ \ \ \ \ \ \ \ \ \ \ \ \ \ \ $\ =\frac{1}{4}(-F_{2}(%
\overline{x_{3}a})+F_{3}(\overline{ax_{2}}))-\frac{1}{4}(-F_{2}(\overline{%
x_{3}a})+F_{3}(\overline{ax_{2}}))$

\ \ \ \ \ \ \ \ \ \ \ \ \ \ \ \ \ \ \ \ \ \ $\ \ =0$

Similarly we have $[\widetilde{E}_{2},\widehat{A}_{2}(a)]=0$ and $[%
\widetilde{E}_{3},\widehat{A}_{3}(a)]=0$ .

\ \ \ \ \emph{Q.E.D.}

\bigskip

\emph{Definition 2.12, } We define an $\R$-vector space \gR$_{4}$\textbf{\ }by

\ \ \ \ \ \ \ \ \ \ \ \ \gR$_{4}=\{(D,M) \mid D\in $\gD$%
_{4},M\in $\gA$\}.$

And we consider \gR$_{4}$ as a vector space \gD$%
_{4}\oplus $\gA$=\R^{28}\oplus \R^{8}\oplus \R^{8}\oplus \R^{8}$

\bigskip

\emph{Definition 2.13, } For $a=\sum_{i=0}^{7}a_{i}e_{i},b=%
\sum_{j=0}^{7}b_{j}e_{j},a_{i},b_{j}\in \R,$we define $\textrm{JD}(a,b)$ by

\ \ \ \ \ \ \ \ \ \ \ \ \ \ \ $\textrm{JD}(a,b)=\sum\limits_{i,j=0(i\neq
j)}^{7}a_{i}b_{j}D_{ij}=\sum\limits_{0 \leq i <j \leq 7}(a_{i}b_{j}-a_{j}b_{i})D_{ij}\in $\gD$_{4}.$

\bigskip

\emph{Definition 2.14, } For $\delta _{1}=(D_{1},M_{1}),\delta
_{2}=(D_{2},M_{2})\in $\gR$_{4}$\textbf{,}

$%
M_{1}=A_{1}(m_{11})+A_{2}(m_{12})+A_{3}(m_{13}),M_{2}=A_{1}(m_{21})+A_{2}(m_{22})+A_{3}(m_{23}), 
$

$\ m_{11},m_{12},m_{13},m_{21},m_{22},m_{23}\in $\gC ,

We define a bracket operation $[\delta _{1},\delta _{2}]_{4}$ by

\ \ \ \ \ \ \ \ \ \ \ $[\delta _{1},\delta _{2}]_{4}=([\delta _{1},\delta
_{2}]_{4D},[\delta _{1},\delta _{2}]_{4M})$ $\ ([\delta _{1},\delta
_{2}]_{4D}\in $\gD$_{4},[\delta _{1},\delta _{2}]_{4M}\in $\gA )%
$,$

where $[\delta _{1},\delta
_{2}]_{4D}=[D_{1},D_{2}]-\textrm{JD}(m_{11},m_{21})-\textrm{d}_{g}\nu
^{2}\textrm{g}_{d}(\textrm{JD}(m_{12},m_{22}))$

\ \ \ \ \ \ \ \ \ \ \ \ \ \ \ \ \ \ \ \ \ \ \ \ 
$-\textrm{d}_{g}\nu \textrm{g}_{d}(\textrm{JD}(m_{13},m_{23})),$

\ \ \ \ \ \ \ \ \ $[\delta _{1},\delta
_{2}]_{4M}=A_{1}(\textrm{g}_{d}(D_{1})m_{21}-\textrm{g}_{d}(D_{2})m_{11})$

\ \ \ \ \ \ \ \ \ \ \ \ \ \ \ \ \ \ \ \ \ \ \ 
$+A_{2}(\nu \textrm{g}_{d}(D_{1})m_{22}-\nu \textrm{g}_{d}(D_{2})m_{12})$

\ \ \ \ \ \ \ \ \ \ \ \ \ \ \ \ \ \ \ \ \ \ \ $+A_{3}(\nu
^{2}\textrm{g}_{d}(D_{1})m_{23}-\nu ^{2}\textrm{g}_{d}(D_{2})m_{13})$

\ \ \ \ \ \ \ \ \ \ \  \ \ \ \ \ \ \ \ \ \ \ \ \ \ \ \ \ \ \ \ $+A_{1}(-\frac{1}{2}%
\overline{m_{12}m_{23}}+\frac{1}{2}\overline{m_{22}m_{13}})+A_{2}(-\frac{1}{2%
}\overline{m_{13}m_{21}}+\frac{1}{2}\overline{m_{23}m_{11}})$

\ \ \ \ \ \ \ \ \ \ \ \ \ \ \ \ \ \ \ \ \ \ \ $+A_{3}(-\frac{1}{2}%
\overline{m_{11}m_{22}}+\frac{1}{2}\overline{m_{21}m_{12}}).$

\bigskip

\emph{Lemma 2.15, } \gR$_{4}$ is isomorphic to 
\gf$_{4}$ under the correspondence

$\ \ \ \ \ \ \ \ \ \ \ \ \ \ \ \ \ \ \ \ \ \textrm{f}_{y}:(D,M)\in $\gR%
$_{4}\rightarrow \textrm{f}_{y}(D,M)=\delta _{d}(\textrm{g}_{d}(D))+\widehat{M}\in $\gf$%
_{4}.$

\bigskip

\emph{Proof.} \ Since $\delta _{d}(\textrm{g}_{d}(D_{1})),\delta
_{d}(\textrm{g}_{d}(D_{2}))\in $\gd$_{4},\delta _{d}(\textrm{g}_{d}(D_{1}))+\widehat{M}%
_{1},\delta _{d}(\textrm{g}_{d}(D_{2}))+\widehat{M}_{2}\in $\gf$_{4},$
we have

[$\delta _{d}(\textrm{g}_{d}(D_{1}))+\widehat{M}_{1},\delta _{d}(\textrm{g}_{d}(D_{2}))+%
\widehat{M}_{2}]$

$=[\delta _{d}(\textrm{g}_{d}(D_{1})),\delta _{d}(\textrm{g}_{d}(D_{2}))]$
$+[\delta _{d}(\textrm{g}_{d}(D_{1})),\widehat{M}_{2}]-[\delta
_{d}(\textrm{g}_{d}(D_{2})),\widehat{M}_{1}]+[\widehat{M}_{1},\widehat{M}_{2}],$

$[\delta _{d}(\textrm{g}_{d}(D_{1})),\delta _{d}(\textrm{g}_{d}(D_{2}))]\in $\gd$%
_{4},[\delta _{d}(\textrm{g}_{d}(D_{1})),\widehat{M}_{2}],[\delta _{d}(\textrm{g}_{d}(D_{2})),%
\widehat{M}_{1}]\in $\gA$,$

$[\widehat{M}_{1},\widehat{M}_{2}]=\sum\limits_{k=1}^{3}[\widehat{A}%
_{k}(m_{1k}),\widehat{A}_{k}(m_{2k})]+\sum\limits_{k,l=1(k\neq l)}^{3}[%
\widehat{A}_{k}(m_{1k}),\widehat{A}_{l}(m_{2l})],$

$[\widehat{A}_{k}(m_{1k}),\widehat{A}_{k}(m_{2k})]\in $\gd$_{4},[%
\widehat{A}_{k}(m_{1k}),\widehat{A}_{l}(m_{2l})]_{(k\neq l)}\in $\gA$%
. $

On the other hand,

(By \emph{Lemma 2.5}) we have

$[\delta _{d}(\textrm{g}_{d}(D_{1})),\delta _{d}(\textrm{g}_{d}(D_{2}))]=\delta
_{d}(\textrm{g}_{d}([D1,D2]).$

(\emph{Lemma 2.11 (18),(19),(20)}) we have,

$[\delta _{d}(\textrm{g}_{d}(D_{1})),\widehat{M}_{2}]=\sum\limits_{k=1}^{3}[\delta
_{d}(\textrm{g}_{d}(D_{1})),\widehat{A}_{k}(m_{2k})],$

\ \ \ \ \ \ \ \ \ \ \ \ \ \ \ \ \ \ \ \ \ \ =$\widehat{A}%
_{1}(\textrm{g}_{d}(D_{1})m_{21})+\widehat{A}_{2}(\nu \textrm{g}_{d}(D_{1})m_{22})+\widehat{A}%
_{3}(\nu ^{2}\textrm{g}_{d}(D_{1})m_{23})$,

$[\delta _{d}(\textrm{g}_{d}(D_{2})),\widehat{M}_{1}]=\widehat{A}%
_{1}(\textrm{g}_{d}(D_{2})m_{11})+\widehat{A}_{2}(\nu \textrm{g}_{d}(D_{2})m_{12})+\widehat{A}%
_{3}(\nu ^{2}\textrm{g}_{d}(D_{2})m_{13}).$

(By \emph{Lemma 2.11 (1),(2),(3),(4)}) For $m_{kl}=\sum\limits_{i=0}^{7}$m$%
_{kli}e_{i},k=1,2,l=1,2,3,$ we have

$[\widehat{A}_{1}(m_{11}),\widehat{A}_{1}(m_{21})]=-\sum\limits_{i=0}^{7}$m$%
_{11i}\sum\limits_{j=0}^{7}m_{21j}\delta _{d}(G_{ij})$

$\ \ \ \ \ \ \ \ \ \ \ \ \ \ \ \ \ \ \ \ \ \ \ \ \ \ =-\delta
_{d}(\textrm{g}_{d}(\sum\limits_{i=0}^{7}$m$_{11i}\sum\limits_{j=0}^{7}$m$%
_{21j}D_{ij})),$

$\ \ \ \ \ \ \ \ \ \ \ \ \ \ \ \ \ \ \ \ \ \ \ \ \ \ =-\delta _{d}(\textrm{g}_{d}(\textrm{JD}($m$%
_{11},$m$_{21}))),$

$[\widehat{A}_{2}(m_{12}),\widehat{A}_{2}(m_{22})]=-\sum\limits_{i=0}^{7}$m$%
_{12i}\sum\limits_{j=0}^{7}$m$_{22j}\delta _{d}(\nu ^{2}G_{ij})$

$\ \ \ \ \ \ \ \ \ \ \ \ \ \ \ \ \ \ \ \ \ \ \ \ \ \ =-\delta
_{d}(\nu ^{2}\textrm{g}_{d}(\textrm{JD}(m_{12},m_{22}))),$

$\ \ \ \ \ \ \ \ \ \ \ \ \ \ \ \ \ \ \ \ \ \ \ \ \ \ =-\delta _{d}(\textrm{g}_{d}(\textrm{d}_{g}\nu
^{2}\textrm{g}_{d}(\textrm{JD}(m_{12},m_{22})))),$

$[\widehat{A}_{3}(m_{13}),\widehat{A}_{3}(m_{23})]=-\sum\limits_{i=0}^{7}$m$%
_{13i}\sum\limits_{j=0}^{7}$m$_{23j}\delta _{d}(\nu G_{ij})$

$\ \ \ \ \ \ \ \ \ \ \ \ \ \ \ \ \ \ \ \ \ \ \ \ \ \ =-\delta _{d}(\nu
\textrm{g}_{d}(\textrm{JD}(m_{13},m_{23}))),$

$\ \ \ \ \ \ \ \ \ \ \ \ \ \ \ \ \ \ \ \ \ \ \ \ \ \ =-\delta _{d}(\textrm{g}_{d}(\textrm{d}_{g}\nu
\textrm{g}_{d}(\textrm{JD}(m_{13},m_{23})))),$

(By \emph{Lemma 2.11 (5),(6)}) For $m_{kl}=\sum\limits_{i=0}^{7}$m$%
_{kli}e_{i},k=1,2,l=1,2,3,$ we have

$[\widehat{A}_{1}(m_{11}),\widehat{A}_{2}(m_{22})]=-\frac{1}{2}\widehat{A}%
_{3}(\overline{m_{11}m_{22}})$,

$[\widehat{A}_{1}(m_{11}),\widehat{A}_{3}(m_{23})]$
$=\ \ \frac{1}{2}\widehat{A}_{2}(\overline{m_{23}m_{11}}),$

$[\widehat{A}_{2}(m_{12}),\widehat{A}_{1}(m_{21})]=\ \ \frac{1}{2}\widehat{A}%
_{3}(\overline{m_{21}m_{13}})$,

$[\widehat{A}_{2}(m_{12}),\widehat{A}_{3}(m_{23})]$
$=-\frac{1}{2}\widehat{A}_{1}(\overline{m_{12}m_{23}}),$

$[\widehat{A}_{3}(m_{13}),\widehat{A}_{1}(m_{21})]=-\frac{1}{2}\widehat{A}%
_{2}(\overline{m_{13}m_{21}})$,

$[\widehat{A}_{3}(m_{13}),\widehat{A}_{2}(m_{22})]$
$=\ \ \frac{1}{2}\widehat{A}_{1}(\overline{m_{22}m_{13}}).$

Therfor we have,

[$\delta _{d}(\textrm{g}_{d}(D_{1}))+\widehat{M}_{1},\delta _{d}(\textrm{g}_{d}(D_{2}))+%
\widehat{M}_{2}]$

$=\delta _{d}(\textrm{g}_{d}([D_{1},D_{2}])$

\ \ +$\widehat{A}_{1}(\textrm{g}_{d}(D_{1})m_{21})+\widehat{A}_{2}(\nu
\textrm{g}_{d}(D_{1})m_{22})+\widehat{A}_{3}(\nu ^{2}\textrm{g}_{d}(D_{1})m_{23})$

\ \ -$\widehat{A}_{1}(\textrm{g}_{d}(D_{2})m_{11})-\widehat{A}_{2}(\nu
\textrm{g}_{d}(D_{2})m_{12})-\widehat{A}_{3}(\nu ^{2}\textrm{g}_{d}(D_{2})m_{13})$

\ \ -$\delta _{d}(\textrm{g}_{d}(\textrm{JD}(m_{11},m_{21})))-\delta _{d}(\textrm{g}_{d}(\textrm{d}_{g}\nu
^{2}\textrm{g}_{d}(\textrm{JD}(m_{12},m_{22}))))$

\ \ -$\delta _{d}(\textrm{g}_{d}(\textrm{d}_{g}\nu \textrm{g}_{d}(\textrm{JD}(m_{13},m_{23}))))$

\ \ $-\frac{1}{2}A_{3}(\overline{m_{11}m_{22}})+\frac{1}{2}A_{2}(\overline{%
m_{23}m_{11}})+\frac{1}{2}A_{3}(\overline{m_{21}m_{12}})$

\ \ $-\frac{1}{2}A_{1}(%
\overline{m_{12}m_{23}})-\frac{1}{2}A_{2}(\overline{m_{13}m_{21}})+\frac{1}{2%
}A_{1}(\overline{m_{22}m_{13}}),$

=$\delta _{d}(\textrm{g}_{d}([D_{1},D_{2}]-\textrm{JD}(m_{11},m_{21})-\textrm{d}_{g}\nu
^{2}\textrm{g}_{d}(\textrm{JD}(m_{12},m_{22}))$

\ \ \ \ \ \ \ \ \ $-\textrm{d}_{g}\nu \textrm{g}_{d}(\textrm{JD}(m_{13},m_{23})))$

\ \ +$\widehat{A}_{1}(\textrm{g}_{d}(D_{1})m_{21}-\textrm{g}_{d}(D_{2})m_{11})+\widehat{A}_{2}(\nu
\textrm{g}_{d}(D_{1})m_{22}-\nu \textrm{g}_{d}(D_{2})m_{12})$

\ \ $+\widehat{A}_{3}(\nu ^{2}\textrm{g}_{d}(D_{1})m_{23}-\nu ^{2}\textrm{g}_{d}(D_{2})m_{13})$

\ \ $+\widehat{A}_{1}(-\frac{1}{2}\overline{m_{12}m_{23}}+\frac{1}{2}\overline{%
m_{22}m_{13}})+\widehat{A}_{2}(-\frac{1}{2}\overline{m_{13}m_{21}}+\frac{1}{2%
}\overline{m_{23}m_{11}})$

\ \ $+\widehat{A}_{3}(-\frac{1}{2}\overline{m_{11}m_{22}}+\frac{1}{2}\overline{%
m_{21}m_{12}}).$

Let we put

$[\delta _{1},\delta _{2}]_{4D}=[D1,D2]-\textrm{JD}(m_{11},m_{21})-\textrm{d}_{g}\nu
^{2}\textrm{g}_{d}(\textrm{JD}(m_{12},m_{22}))$

\ \ \ \ \ \ \ \ \ \ \ \ \ \ \ 
$-\textrm{d}_{g}\nu \textrm{g}_{d}(\textrm{JD}(m_{13},m_{23})),$

$[\delta _{1},\delta
_{2}]_{4M}=A_{1}(\textrm{g}_{d}(D_{1})m_{21}-\textrm{g}_{d}(D_{2})m_{11})+A_{2}(\nu
\textrm{g}_{d}(D_{1})m_{22}-\nu \textrm{g}_{d}(D_{2})m_{12})$

$\ \ \ \ \ \ \ \ \ \ \ \ \ \ +A_{3}(\nu ^{2}\textrm{g}_{d}(D_{1})m_{23}-\nu
^{2}\textrm{g}_{d}(D_{2})m_{13})$

\ \ \ \ \ \ \ \ \ \ \ \ \ \ $+A_{1}(-\frac{1}{2}\overline{m_{12}m_{23}}+%
\frac{1}{2}\overline{m_{22}m_{13}})+A_{2}(-\frac{1}{2}\overline{m_{13}m_{21}}%
+\frac{1}{2}\overline{m_{23}m_{11}})$

\ \ \ \ \ \ \ \ \ \ \ \ \ \ $+A_{3}(-\frac{1}{2}\overline{m_{11}m_{22}}+%
\frac{1}{2}\overline{m_{21}m_{12}}),$

\noindent
then $[\delta _{1},\delta _{2}]_{4}=([\delta _{1},\delta _{2}]_{4D},[\delta
_{1},\delta _{2}]_{4M})$ is a Lie bracket of \gR$_{4}$%
\textbf{\ }and also \gR$_{4}$\textbf{\ }is

\noindent
isomorphic to \gf$_{4}$.\ \ \ \ \emph{Q.E.D.}

\bigskip

\section{Definition of the vector space \gR$_{6}$ and it's Lie bracket}

\bigskip 

\ \ \ \emph{\ Definition 3.1.} \ \ For an $\R$-vector space \textbf{V}, we
denote complex conjugation of v$\in $\textbf{V}$^{\C}$
by \ \ $\iota $v .

\bigskip

Let \gJ$^{\C}$ $=\{X_{1}+iX_{2}|X_{1},X_{2}\in $\gJ\} be the
complexification of the Jordan algebra \gJ.

\bigskip

\emph{Definition 3.2. \ }For $X,Y\in $\gJ$^{\C},$We define a
positive definite Hermitian inner product
$\langle X,Y \rangle$ in \gJ$^{\C}$ by

\ \ \ \ \ \ \ \ \ \ \ \ \ \ \ \ \ \ \ \ \ \ \ \ \ \ \ \ \ \ \ \ \ \ $%
\langle X,Y \rangle =(\iota X,Y).$

\bigskip

\emph{Definition 3.3.} \ We define the groups $E_{6}^{\C}$ and $E_{6}$
respectively by

$\ \ \ \ \ \ \ \ \ \ \ E_{6}^{\C}=\{\alpha \in Iso_{\C}($\gJ$^{\C}) \mid det(%
\alpha X)=detX\},$

\ \ \ \ \ \ \ \ \ \ \ $E_{6}=\{\alpha \in Iso_{\C}($\gJ$%
^{\C}) \mid det(\alpha X)=detX, \langle \alpha X,\alpha Y \rangle = \langle X,Y \rangle \}$,

\ \ \ \ \ \ \ \ \ \ \ $\ \ \ \ =\{\alpha \in Iso_{\C}($\gJ$%
^{\C}) \mid (\alpha X, \alpha Y, \alpha Z)=(X, Y, Z),\langle \alpha X,\alpha Y \rangle= \langle X,Y \rangle \}$,

\ \ \ \ \ \ \ \ \ \ \ \ \ \ $\ \ \ \ =\{\alpha \in Iso_{\C}($\gJ$%
^{\C}) \mid \alpha X \times \alpha Y=(\alpha ^{t})^{-1}(X \times Y),\langle \alpha X,\alpha Y \rangle= \langle X,Y \rangle \}$,

\ \ \ \ \ \ \ \ $\ \ \ \ =\{\alpha \in Iso_{\C}($\gJ$%
^{\C}) \mid \alpha X \times \alpha Y=\iota \alpha \iota (X \times Y),\langle \alpha X,\alpha Y \rangle= \langle X,Y \rangle \}$.

\bigskip

\emph{Proposition 3.4. } (I.Yokota\cite[\emph{Theorem 3.9.2.}]{Yokota1})

$E_{6}$ is a simply connected compact Lie group of type $E_{6}$.

\bigskip

\emph{Proposition 3.5. } (I.Yokota\cite[\emph{Theorem 3.14.1.}]{Yokota1}) 

The polar decomposition of the Lie group $E_{6}^{\C}$
is given by

\ \ \ \ \ \ \ \ \ \ \ \ \ \ \ \ \ \ \ \ \ \ \ \ \ \ \ $E_{6}^{\C}\simeq
E_{6}\times \R^{78}.$

In particular, $E_{6}^{\C}$ is a simply connected complex Lie group of type $%
E_{6}$.

\bigskip

\emph{Proposition 3.6.  \ }

(I.Yokota\cite[\emph{Theorem3.2.1.}]{Yokota1}, T.Imai and I.Yokota\cite[\emph{Proposition 1.}]{YokotaImai5})

$(1)$ The Lie algebra \ge$_{6}^{\C}$ of the Lie group $E_{6}^{\C}$ is
given by

\ \ \ \ \ \ \ \ \ \ \ \ \ \ \ge$_{6}^{\C}=\{\phi \in Hom_{\C}($\gJ%
$^{\C}) \mid (\phi X,X,X)=0\}.$

$(2)$ The Lie algebra \ge$_{6,1}$ of the Lie group $E_{6(-26)}$ is
given by

\ \ \ \ \ \ \ \ \ \ \ \ \ \ \ge$_{6,1}=\{\phi \in Hom_{\R}($\gJ%
$) \mid (\phi X,X,X)=0\}.$

$(3)$ Any element $\phi \in $\ge$_{6}^{\C}$ is uniquely expressed by

$\ \ \ \ \ \ \ \ \ \ \ \ \ \phi =\delta +\widetilde{T},\delta \in $\gf%
$_{4}^{\C},T\in $\gJ$_{0}^{\C},$

\noindent
where \gJ$_{0}^{\C}=\{X\in $\gJ$^{\C} \mid tr(X)=0\}.$

\bigskip

\emph{Proposition 3.7.} (I.Yokota\cite[\emph{  Theorem3.2.2.}]{Yokota1})

The Lie bracket $[\phi _{1},\phi _{2}]$ in \ge$_{6}^{\C}$ is given by

\ \ \ \ $[\delta _{1}+\widetilde{T}_{1},\delta _{2}+\widetilde{T}_{2}]=([%
\delta _{1},\delta _{2}]+[\widetilde{T}_{1},\widetilde{T}_{2}])+(%
\widetilde{\delta _{1}T}_{2}-\widetilde{\delta _{2}T}_{1}),$

\noindent
where $\phi _{i}=\delta _{i}+\widetilde{T}_{i},\delta _{i}\in $\gf$%
_{4}^{\C},T_{i}\in $\gJ$_{0}^{\C}.$

\bigskip

\emph{Corollary 3.7.1. \ }The definition of the Lie bracket $[\phi
_{1},\phi _{2}]$ in \emph{Proposition 3.7 } can be
extended to $T_{i}\in $\gJ$^{\C}$ $(i=1,2).$

\bigskip

\emph{Proof.} \ For $\phi _{i}=\delta _{i}+\widetilde{T}_{i},\delta
_{i}\in $\gf$_{4}^{\C},T_{i},X\in $\gJ$^{\C},$ we have

$[\delta _{1}+\widetilde{T}_{1},\delta _{2}+\widetilde{T}_{2}]=([\delta
_{1},\delta _{2}]+[\widetilde{T}_{1},\widetilde{T}_{2}])+([\delta _{1},%
\widetilde{T}_{2}]-[\delta _{2},\widetilde{T}_{1}])$

By \emph{Lemma 2.11} $ (1),(2),(3),(7),(8),(12),(13),(14)$ and $(28)$, 

\noindent
we have $[%
\widetilde{T}_{1},\widetilde{T}_{2}]\in $\gA$^{\C}.$
On the otherhand, we have

$[\delta _{1},\widetilde{T}_{2}]X=\delta _{1}(T_{2}\circ X)-T_{2}\circ 
\delta _{1}X$

\ \ \ \ \ \ \ \ \ \ \ \ $=\delta _{1}T_{2}\circ X+T_{2}\circ \delta _{1}X-T_{2}\circ 
\delta _{1}X$ (by \emph{Proposition 2.2}),

\ \ \ \ \ \ \ \ \ \ \ \ =$\widetilde{\delta _{1}T_{2}}X.$

\noindent
Hence we have

$[\delta _{1}+\widetilde{T}_{1},\delta _{2}+\widetilde{T}_{2}]=([\delta
_{1},\delta _{2}]+[\widetilde{T}_{1},\widetilde{T}_{2}])+(\widetilde{\delta
_{1}T}_{2}-\widetilde{\delta _{2}T}_{1})$ . \ \ \ \ \emph{Q.E.D.}

\bigskip

\emph{Proposition 3.8. } (I.Yokota\cite[\emph{Lemma 3.2.3.}]{Yokota1}) 

For $\phi =\delta +%
\widetilde{T}\in $\ge$_{6}^{\C},\delta \in $\gf$_{%
4}^{\C},T\in $\gJ$_{0}^{\C}$,we have

\ \ \ \ \ \ \ \ \ \ \ \ \ \ \ \ \ \ \ \ \ \ \ \ \ \ \ \ \ \ \ $-$ $\phi^{t}
=-$ $(\delta +\widetilde{T})^{t}=\delta -\widetilde{T}$.

\noindent
In particular,$ -\phi^{t}  \in $ \ge$_{6}^{\C}$, where upper right $^{t}$ means transpose with inner product.

\bigskip

\emph{Corollary 3.8.1. \ }For $T\in $\gJ$^{\C}$,we have

\ \ \ \ \ \ \ \ \ \ \ \ \ \ \ \ \ \ \ \ \ \ \ \ \ \ \ \ \ \ $\widetilde{T}^{t}=\widetilde{T}$

\bigskip

\emph{Proof.} \ For $X,Y\in $\gJ$^{\C},$we have

\ \ $(\widetilde{T}^{t}X,Y)=(X,\widetilde{T}Y)=(X,T\circ Y)=(X\circ
T,Y)=(T\circ X,Y)=(\widetilde{T}X,Y).$ \ \ \ \ \emph{Q.E.D.}

\bigskip

\emph{Definition 3.9.} We define an $\R$-vector space \gR$_{6}$\textbf{\ }by

\ \ \ \ \ \ \ \ \ \ \ \ \gR$_{6}=\{(D,M,T) \mid D\in $\gD%
$_{4},M\in $\gA$,T\in $\gJ$_{0}\}.$

\noindent
And we consider \gR$_{6}$ as a vector space \gD$%
_{4}\oplus $\gA$\oplus $\gJ$_{0}=\R^{28}\oplus
\R^{8}\oplus \R^{8}\oplus \R^{8}\oplus $

\noindent
$\R\oplus \R\oplus \R^{8}\oplus \R^{8}\oplus \R^{8}.$

\bigskip

\emph{Definition 3.10. } For $\phi _{1}=(D_{1},M_{1},T1),\phi
_{2}=(D_{2},M_{2},T2)\in $\gR$_{6}$\textbf{,}

$%
M_{1}=A_{1}(m_{11})+A_{2}(m_{12})+A_{3}(m_{13}),M_{2}=A_{1}(m_{21})+A_{2}(m_{22})+A_{3}(m_{23}), 
$

$\ m_{11},m_{12},m_{13},m_{21},m_{22},m_{23}\in $\gC ,

T$_{1}$=$\tau _{11}E_{1}+\tau _{12}E_{2}+(-\tau _{11}-\tau
_{12})E_{3}+F_{1}(t_{11})+F_{2}(t_{12})+F_{3}(t_{13}),$

T$_{2}$=$\tau _{21}E_{1}+\tau _{22}E_{2}+(-\tau _{21}-\tau
_{22})E_{3}+F_{1}(t_{21})+F_{2}(t_{22})+F_{3}(t_{23}),$

$t_{11},t_{12},t_{13},t_{21},t_{22},t_{23}\in $\gC $,\tau _{11},\tau
_{12},\tau _{21},\tau _{22}\in \R,$

We define a bracket operation $[\phi _{1},\phi _{2}]_{6}$ by

\ $[\phi _{1},\phi _{2}]_{6}=([\phi _{1},\phi _{2}]_{6D},[\phi _{1},\phi
_{2}]_{6M},[\phi _{1},\phi _{2}]_{6T}),$ \ 

$\ \ \ \ \ \ \ \ \ \ \ \ \ \ \ \ \ \
\ \ ([\phi _{1},\phi _{2}]_{6D}\in $\gD$_{4},[\phi _{1},\phi
_{2}]_{6M}\in $\gA$,[\phi _{1},\phi _{2}]_{6T}\in $\gJ$%
_{0}).$

\noindent
where

$[\phi _{1},\phi
_{2}]_{6D}=[(D_{1},M_{1}),(D_{2},M_{2})]_{4D}+([T_{1},T_{2}])_{D},$

$[\phi _{1},\phi
_{2}]_{6M}=[(D_{1},M_{1}),(D_{2},M_{2})]_{4M}+([T_{1},T_{2}])_{M},$

$[\phi _{1},\phi
_{2}]_{6T}=([(D_{1},M_{1}),T_{2}])_{T}-([(D_{2},M_{2}),T_{1}])_{T}.$

$([T_{1},T_{2}])_{D}=\textrm{JD}(t_{11},t_{21})+\textrm{d}_{g}\nu
^{2}\textrm{g}_{d}(\textrm{JD}(t_{12},t_{22}))+\textrm{d}_{g}\nu \textrm{g}_{d}(\textrm{JD}(t_{13},t_{23})),$

$([T_{1},T_{2}])_{M}=A_{1}(\frac{1}{2}((\tau _{11}+2\tau _{12})t_{21}-(\tau
_{21}+2\tau _{22})t_{11}-\overline{t_{12}t_{23}}+\overline{t_{22}t_{13}}))$

$\ \ \ \ \ \ \ \ \ \ \ \ \ \ \ \ +A_{2}(\frac{1}{2}((-2\tau _{11}-\tau
_{12})t_{22}-(-2\tau _{21}-\tau _{22})t_{12}+\overline{t_{23}t_{11}}-%
\overline{t_{13}t_{21}}))$

$\ \ \ \ \ \ \ \ \ \ \ \ \ \ \ \ +A_{3}(\frac{1}{2}((\tau _{11}-\tau
_{12})t_{23}-(\tau _{21}-\tau _{22})t_{13}-\overline{t_{11}t_{22}}+\overline{%
t_{21}t_{12}})),$

$([(D_{1},M_{1}),T_{2}])_{T}=F_{1}(\textrm{g}_{d}(D_{1})t_{21})+F_{2}(\nu
\textrm{g}_{d}(D_{1})t_{22})+F_{3}(\nu ^{2}\textrm{g}_{d}(D_{1})t_{23})$

\ \ \ \ \ \ \ \ \ \ \ \ \ \ \ \ \ \ $\
+(-(m_{12},t_{22})+(m_{13},t_{23}))E_{1}+(-(m_{13},t_{23})+(m_{11},t_{21}))E_{2} 
$

$\ \ \ \ \ \ \ \ \ \ \ \ \ \ \ \ \ \ \
+(-(m_{11},t_{21})+(m_{12},t_{22}))E_{3}$

$\ \ \ \ \ \ \ \ \ \ \ \ \ \ \ \ \ \ \ +F_{1}(\frac{1}{2}(\overline{%
m_{12}t_{23}}-\overline{t_{22}m_{13}}+(-\tau _{21}-2\tau _{22})m_{11}))$

\ \ \ \ \ \ \ \ \ \ \ \ \ \ \ \ \ \ $\ +F_{2}(\frac{1}{2}(\overline{%
m_{13}t_{21}}-\overline{t_{23}m_{11}}+(2\tau _{21}+\tau _{22})m_{12}))$

\ \ \ \ \ \ \ \ \ \ \ \ \ \ \ \ \ \ $\ +F_{3}(\frac{1}{2}(\overline{%
m_{11}t_{22}}-\overline{t_{21}m_{12}}+(-\tau _{21}+\tau _{22})m_{13})),$

$([(D_{2},M_{2}),T_{1}])_{T}=F_{1}(\textrm{g}_{d}(D_{2})t_{11})+F_{2}(\nu
\textrm{g}_{d}(D_{2})t_{12})+F_{3}(\nu ^{2}\textrm{g}_{d}(D_{2})t_{13})$

\ \ \ \ \ \ \ \ \ \ \ \ \ \ \ \ \ $\
+(-(m_{22},t_{12})+(m_{23},t_{13}))E_{1}+(-(m_{23},t_{13})+(m_{21},t_{11}))E_{2} 
$

$\ \ \ \ \ \ \ \ \ \ \ \ \ \ \ \ \ \ +(-(m_{21},t_{11})+(m_{22},t_{12}))E_{3}$

$\ \ \ \ \ \ \ \ \ \ \ \ \ \ \ \ \ \ +F_{1}(\frac{1}{2}(\overline{m_{22}t_{13}}%
-\overline{t_{12}m_{23}}+(-\tau _{11}-2\tau _{12})m_{21}))$

\ \ \ \ \ \ \ \ \ \ \ \ \ \ \ \ \ $\ +F_{2}(\frac{1}{2}(\overline{m_{23}t_{11}}%
-\overline{t_{13}m_{21}}+(2\tau _{11}+\tau _{12})m_{22}))$

\ \ \ \ \ \ \ \ \ \ \ \ \ \ \ \ \ $\ +F_{3}(\frac{1}{2}(\overline{m_{21}t_{12}}%
-\overline{t_{11}m_{22}}+(-\tau _{11}+\tau _{12})m_{23})).$

\bigskip

\emph{Lemma 3.11.} \ \gR$_{6}$ is isomorphic to 
\ge$_{6,1}$ under the correspondence

$\ \ \ \ \ \ \ \ \ \ \ \ \ \ \ \ \ \ \ \ \ \textrm{f}_{y}:(D,M,T)\in $\gR$_{%
6}\rightarrow \textrm{f}_{y}(D,M,T)=\delta _{d}(\textrm{g}_{d}(D))+\widehat{M}+%
\widetilde{T}\in $\ge$_{6,1}$.

\bigskip

\emph{Proof.} We have $\delta _{d}(\textrm{g}_{d}(D_{1})),\delta
_{d}(\textrm{g}_{d}(D_{2}))\in $\gd$_{4},$

$\ \ \ \ \ \ \ \ \ \ \ \ \ \ \ \ \ \ \ \ \ \ \delta _{d}(\textrm{g}_{d}(D_{1}))+%
\widehat{M}_{1},\delta _{d}(\textrm{g}_{d}(D_{2}))+\widehat{M}_{2}\in $\gf$%
_{4},$

\ \ \ \ \ \ \ \ \ \ \ \ \ \ \ \ \ \ \ \ $\ \ \delta _{d}(\textrm{g}_{d}(D_{1}))+%
\widehat{M}_{1}+\widetilde{T}_{1},\delta _{d}(\textrm{g}_{d}(D_{2}))+\widehat{M}_{2}+%
\widetilde{T}_{2}\in $\ge$_{6,1},$

$\ \ \ \ \ \ \ \ \ \ \ \ \ \ \ \ \ \ \ \ \ (\delta _{d}(\textrm{g}_{d}(D_{1}))+\widehat{%
M}_{1})T_{2},(\delta _{d}(\textrm{g}_{d}(D_{2}))+\widehat{M}_{2})T_{1}\in $\gJ$%
_{0}.$

\noindent
Hence we have by \emph{Proposition 3.7 },

[$\delta _{d}(\textrm{g}_{d}(D_{1}))+\widehat{M}_{1}+\widetilde{T}_{1},\delta
_{d}(\textrm{g}_{d}(D_{2}))+\widehat{M}_{2}+\widetilde{T}_{2}]$

$=[\delta _{d}(\textrm{g}_{d}(D_{1})),\delta _{d}(\textrm{g}_{d}(D_{2}))]$

$+[\delta _{d}(\textrm{g}_{d}(D_{1})),\widehat{M}_{2}]-[\delta _{d}(\textrm{g}_{d}(D_{2})),%
\widehat{M}_{1}]+[\widehat{M}_{1},\widehat{M}_{2}]$

+$[\widetilde{T}_{1},\widetilde{T}_{2}]+\widetilde{\delta
_{d}(\textrm{g}_{d}(D_{1}))T_{2}}+\widetilde{\widehat{M}_{1}T_{2}}-\widetilde{\delta
_{d}(\textrm{g}_{d}(D_{2})T_{1}}-\widetilde{\widehat{M}_{2}T_{1}},$

\noindent
On the other hand,

$[\widetilde{T}_{1},\widetilde{T}_{2}]$

$=[\frac{1}{3}(\tau _{11}-\tau _{12})(E_{1}-E_{2})^{\widetilde{}}+\frac{1}{3}%
(\tau _{11}+2\tau _{12})(E_{2}-E_{3})^{\widetilde{}}+\frac{1}{3}(-2\tau
_{11}-\tau _{12})(E_{3}-E_{1})^{\widetilde{}}$

$\ \ \ \ +\widetilde{F}_{1}(t_{11})+\widetilde{F}_{2}(t_{12})+\widetilde{F}%
_{3}(t_{13}),$

$\ \ \ \frac{1}{3}(\tau _{21}-\tau _{22})(E_{1}-E_{2})^{\widetilde{}}+\frac{1%
}{3}(\tau _{21}+2\tau _{22})(E_{2}-E_{3})^{\widetilde{}}+\frac{1}{3}(-2\tau
_{21}-\tau _{22})(E_{3}-E_{1})^{\widetilde{}}$

$\ \ \ \ \ +\widetilde{F}_{1}(t_{21})+\widetilde{F}_{2}(t_{22})+\widetilde{F}%
_{3}(t_{23})],$

$=\frac{1}{3}(\tau _{11}-\tau _{12})[(E_{1}-E_{2})^{\widetilde{}},\widetilde{%
F}_{1}(t_{21})+\widetilde{F}_{2}(t_{22})+\widetilde{F}_{3}(t_{23})]$

$-\frac{1}{3}(\tau _{21}-\tau _{22})[(E_{1}-E_{2})^{\widetilde{}},\widetilde{%
F}_{1}(t_{11})+\widetilde{F}_{2}(t_{12})+\widetilde{F}_{3}(t_{13})]$

$+\frac{1}{3}(\tau _{11}+2\tau _{12})[(E_{2}-E_{3})^{\widetilde{}},%
\widetilde{F}_{1}(t_{21})+\widetilde{F}_{2}(t_{22})+\widetilde{F}%
_{3}(t_{23})]$

$-\frac{1}{3}(\tau _{21}+2\tau _{22})[(E_{2}-E_{3})^{\widetilde{}},%
\widetilde{F}_{1}(t_{11})+\widetilde{F}_{2}(t_{12})+\widetilde{F}%
_{3}(t_{13})]$

$+\frac{1}{3}(-2\tau _{11}-\tau _{12})[(E_{3}-E_{1})^{\widetilde{}},%
\widetilde{F}_{1}(t_{21})+\widetilde{F}_{2}(t_{22})+\widetilde{F}%
_{3}(t_{23})]$

$-\frac{1}{3}(-2\tau _{21}-\tau _{22})[(E_{3}-E_{1})^{\widetilde{}},%
\widetilde{F}_{1}(t_{11})+\widetilde{F}_{2}(t_{12})+\widetilde{F}%
_{3}(t_{13})]$

$+[\widetilde{F}_{1}(t_{11}),\widetilde{F}_{1}(t_{21})]+[\widetilde{F}%
_{1}(t_{11}),\widetilde{F}_{2}(t_{22})]+[\widetilde{F}_{1}(t_{11}),%
\widetilde{F}_{3}(t_{23})]$

$+[\widetilde{F}_{2}(t_{12}),\widetilde{F}_{1}(t_{21})]+[\widetilde{F}%
_{2}(t_{12}),\widetilde{F}_{2}(t_{22})]+[\widetilde{F}_{2}(t_{12}),%
\widetilde{F}_{3}(t_{23})]$

$+[\widetilde{F}_{3}(t_{13}),\widetilde{F}_{1}(t_{21})]+[\widetilde{F}%
_{3}(t_{13}),\widetilde{F}_{2}(t_{22})]+[\widetilde{F}_{3}(t_{13}),%
\widetilde{F}_{3}(t_{23})]$, \ 

(By \emph{Lemma 2.11 } $(1),(2),(3),(7),(8),(12),(13)$ and $(14)$)

$=\frac{1}{6}(\tau _{11}-\tau _{12})(-\widehat{A}_{1}(t_{21})-\widehat{A}%
_{2}(t_{22})+2\widehat{A}_{3}(t_{23}))$

\ $-\frac{1}{6}(\tau _{21}-\tau _{22})(-%
\widehat{A}_{1}(t_{11})-\widehat{A}_{2}(t_{12})+2\widehat{A}_{3}(t_{13}))$

\ $+\frac{1}{6}(\tau _{11}+2\tau _{12})(2\widehat{A}_{1}(t_{21})-\widehat{A}%
_{2}(t_{22})-\widehat{A}_{3}(t_{23}))$

\ $-\frac{1}{6}(\tau _{21}+2\tau _{22})(2%
\widehat{A}_{1}(t_{11})-\widehat{A}_{2}(t_{12})-\widehat{A}_{3}(t_{13}))$

\ $+\frac{1}{6}(-2\tau _{11}-\tau _{12})(-\widehat{A}_{1}(t_{21})+2\widehat{A}%
_{2}(t_{22})-\widehat{A}_{3}(t_{23}))$

\ $-\frac{1}{6}(-2\tau _{21}-\tau _{22})(-%
\widehat{A}_{1}(t_{11})+2\widehat{A}_{2}(t_{12})-\widehat{A}_{3}(t_{13}))$

\ $+\delta _{d}(\textrm{g}_{d}(\textrm{JD}(t_{11},t_{21})))-\frac{1}{2}\widehat{A}_{3}(\overline{%
t_{11}t_{22}})+\frac{1}{2}\widehat{A}_{2}(\overline{t_{23}t_{11}})$

\ $+\frac{1}{2}\widehat{A}_{3}(\overline{t_{21}t_{12}})+\delta
_{d}(\textrm{g}_{d}(\textrm{d}_{g}\nu ^{2}\textrm{g}_{d}(\textrm{JD}(t_{12},t_{22}))))-\frac{1}{2}\widehat{A}%
_{1}(\overline{t_{12}t_{23}})$

\ $-\frac{1}{2}\widehat{A}_{2}(\overline{t_{13}t_{21}})+\frac{1}{2}\widehat{A}%
_{1}(\overline{t_{22}t_{13}})+\delta _{d}(\textrm{g}_{d}(\textrm{d}_{g}\nu
\textrm{g}_{d}(\textrm{JD}(t_{13},t_{23})))),$

$=\delta _{d}(\textrm{g}_{d}(\textrm{JD}(t_{11},t_{21})+\textrm{d}_{g}\nu
^{2}\textrm{g}_{d}(\textrm{JD}(t_{12},t_{22}))+\textrm{d}_{g}\nu \textrm{g}_{d}(\textrm{JD}(t_{13},t_{23}))))$

\ $+\widehat{A}_{1}(\frac{1}{2}((\tau _{11}+2\tau _{12})t_{21}-(\tau
_{21}+2\tau _{22})t_{11}-\overline{t_{12}t_{23}}+\overline{t_{22}t_{13}}))$

\ $+\widehat{A}_{2}(\frac{1}{2}((-2\tau _{11}-\tau _{12})t_{22}-(-2\tau
_{21}-\tau _{22})t_{12}+\overline{t_{23}t_{11}}-\overline{t_{13}t_{21}}))$

\ $+\widehat{A}_{3}(\frac{1}{2}((\tau _{11}-\tau _{12})t_{23}-(\tau _{21}-\tau
_{22})t_{13}-\overline{t_{11}t_{22}}+\overline{t_{21}t_{12}})),$

$\widetilde{\delta _{d}(\textrm{g}_{d}(D_{1}))T_{2}}=\widetilde{F}%
_{1}(\textrm{g}_{d}(D_{1})t_{21})+\widetilde{F}_{2}(\nu \textrm{g}_{d}(D_{1})t_{22})+%
\widetilde{F}_{3}(\nu ^{2}\textrm{g}_{d}(D_{1})t_{23}),$

$\widetilde{\delta _{d}(\textrm{g}_{d}(D_{2}))T_{1}}=\widetilde{F}%
_{1}(\textrm{g}_{d}(D_{2})t_{11})+\widetilde{F}_{2}(\nu \textrm{g}_{d}(D_{2})t_{12})+%
\widetilde{F}_{3}(\nu ^{2}\textrm{g}_{d}(D_{2})t_{13}),$

$\widehat{M}_{1}T_{2}=(\widehat{A}_{1}(m_{11})+\widehat{A}_{2}(m_{12})+%
\widehat{A}_{3}(m_{13}))(\tau _{21}E_{1}+\tau _{22}E_{2}-(\tau _{21}+\tau
_{22})E_{3}$

\ \ \ \ \ \ \ \ \ \ \ \ \ \ \ \ \ \ \ \ \ \ \ \ \ \ \ \ \ \ \ \ \ \ \ \ \ \ \ \ 
\ \ \ \ \ \ \ \ \ \ \ \ \ \ \ \ 
$+F_{1}(t_{21})+F_{2}(t_{22})+F_{3}(t_{23})),$

\ \ \ \ \ \ \ \ $=(m_{11},t_{21})(E_{2}-E_{3})+F_{1}(\frac{1}{2}(-\tau
_{21}-2\tau _{22})m_{11})+F_{2}(-\frac{1}{2}\overline{t_{23}m_{11}})$

\ \ \ \ \ \ \ \ \ $+F_{3}(\frac{1}{2}\overline{m_{11}t_{22}})
+(m_{12},t_{22})(E_{3}-E_{1})+F_{1}(\frac{1}{2}\overline{%
m_{12}t_{23}})$

\ \ \ \ \ \ \ \ \ $+F_{2}(\frac{1}{2}(2\tau _{21}+\tau _{22})m_{12})+F_{3}(-\frac{%
1}{2}\overline{t_{21}m_{12}})+(m_{13},t_{23})(E_{1}-E_{2})$

\ \ \ \ \ \ \ \ \ $+F_{1}(-\frac{1}{2}\overline{%
t_{22}m_{13}})+F_{2}(\frac{1}{2}\overline{m_{13}t_{21}})+F_{3}(\frac{1}{2}%
(\tau _{22}-\tau _{21})m_{13}),$

\ \ \ \ \ \ \ $\
=(-(m_{12},t_{22})+(m_{13},t_{23}))E_{1}+(-(m_{13},t_{23})+(m_{11},t_{21}))E_{2}$

\ \ \ \ \ \ \ \ \ $+(-(m_{11},t_{21})+(m_{12},t_{22}))E_{3}$
 
\ \ \ \ \ \ \ \ \ $+F_{1}(\frac{1}{2}(\overline{m_{12}t_{23}}-\overline{%
t_{22}m_{13}}+(-\tau _{21}-2\tau _{22})m_{11}))$

\ \ \ \ \ \ \ \ $\ +F_{2}(\frac{1}{2}(\overline{m_{13}t_{21}}-\overline{%
t_{23}m_{11}}+(2\tau _{21}+\tau _{22})m_{12}))$

\ \ \ \ \ \ \ \ $\ +F_{3}(\frac{1}{2}(\overline{m_{11}t_{22}}-\overline{%
t_{21}m_{12}}+(-\tau _{21}+\tau _{22})m_{13})).$

$\widehat{M}_{2}T_{1}=(\widehat{A}_{1}(m_{21})+\widehat{A}_{2}(m_{22})+%
\widehat{A}_{3}(m_{23}))(\tau _{11}E_{1}+\tau _{12}E_{2}-(\tau _{11}+\tau
_{12})E_{3}$

\ \ \ \ \ \ \ \ \ \ \ \ \ \ \ \ \ \ \ \ \ \ \ \ \ \ \ \ \ \ \ \ \ \ \ \ \ \ \ \ 
\ \ \ \ \ \ \ \ \ \ \ \ \ \ \ \ 
$+\widetilde{F}_{1}(t_{11})+\widetilde{F}_{2}(t_{12})+\widetilde{F}_{3}(t_{13})),$

\ \ \ \ \ \ \ \ $%
=(-(m_{22},t_{12})+(m_{23},t_{13}))E_{1}+(-(m_{23},t_{13})+(m_{21},t_{11}))E_{2}$

\ \ \ \ \ \ \ \ \ $+(-(m_{21},t_{11})+(m_{22},t_{12}))E_{3} 
$

$\ \ \ \ \ \ \ \ \ +F_{1}(\frac{1}{2}(\overline{m_{22}t_{13}}-\overline{%
t_{12}m_{23}}+(-\tau _{11}-2\tau _{12})m_{21}))$

\ \ \ \ \ \ \ \ $\ +F_{2}(\frac{1}{2}(\overline{m_{23}t_{11}}-\overline{%
t_{13}m_{21}}+(2\tau _{11}+\tau _{12})m_{22}))$

\ \ \ \ \ \ \ \ $\ +F_{3}(\frac{1}{2}(\overline{m_{21}t_{12}}-\overline{%
t_{11}m_{22}}+(-\tau _{11}+\tau _{12})m_{23})).$

Let we put

$[\phi _{1},\phi
_{2}]_{6D}=[(D_{1},M_{1}),(D_{2},M_{2})]_{4D}+([T_{1},T_{2}])_{D},$

$[\phi _{1},\phi
_{2}]_{6M}=[(D_{1},M_{1}),(D_{2},M_{2})]_{4M}+([T_{1},T_{2}])_{M},$

$[\phi _{1},\phi
_{2}]_{6T}=([(D_{1},M_{1}),T_{2}])_{T}-([(D_{2},M_{2}),T_{1}])_{T},$

$([T_{1},T_{2}])_{D}=\textrm{JD}(t_{11},t_{21})+\textrm{d}_{g}\nu
^{2}\textrm{g}_{d}(\textrm{JD}(t_{12},t_{22}))+\textrm{d}_{g}\nu \textrm{g}_{d}(\textrm{JD}(t_{13},t_{23})),$

$([T_{1},T_{2}])_{M}=A_{1}(\frac{1}{2}(\tau _{11}+2\tau _{12})t_{21}-(\tau
_{21}+2\tau _{22})t_{11}-\overline{t_{12}t_{23}}+\overline{t_{22}t_{13}}))$

\ \ $\ \ \ \ \ \ \ \ \ \ \ \ \ \ +A_{2}(\frac{1}{2}((-2\tau _{11}-\tau
_{12})t_{22}-(-2\tau _{21}-\tau _{22})t_{12}+\overline{t_{23}t_{11}}-%
\overline{t_{13}t_{21}}))$

\ \ $\ \ \ \ \ \ \ \ \ \ \ \ \ \ +A_{3}(\frac{1}{2}((\tau _{11}-\tau
_{12})t_{23}-(\tau _{21}-\tau _{22})t_{13}-\overline{t_{11}t_{22}}+\overline{%
t_{21}t_{12}})),$

$([(D_{1},M_{1}),T_{2}])_{T}=F_{1}(\textrm{g}_{d}(D_{1})t_{21})+F_{2}(\nu
\textrm{g}_{d}(D_{1})t_{22})+F_{3}(\nu ^{2}\textrm{g}_{d}(D_{1})t_{23})$

\ \ \ \ \ \ \ \ \ \ \ \ \ \ \ \ \ \ \ \ $\
+(-(m_{12},t_{22})+(m_{13},t_{23}))E_{1}+(-(m_{13},t_{23})+(m_{11},t_{21}))E_{2} 
$

$\ \ \ \ \ \ \ \ \ \ \ \ \ \ \ \ \ \ \ \
+(-(m_{11},t_{21})+(m_{12},t_{22}))E_{3}$

$\ \ \ \ \ \ \ \ \ \ \ \ \ \ \ \ \ \ \ \ +F_{1}(\frac{1}{2}(\overline{%
m_{12}t_{23}}-\overline{t_{22}m_{13}}+(-\tau _{21}-2\tau _{22})m_{11}))$

\ \ \ \ \ \ \ \ \ \ \ \ \ \ \ \ \ \ \ $\ +F_{2}(\frac{1}{2}(\overline{%
m_{13}t_{21}}-\overline{t_{23}m_{11}}+(2\tau _{21}+\tau _{22})m_{12}))$

\ \ \ \ \ \ \ \ \ \ \ \ \ \ \ \ \ \ \ $\ +F_{3}(\frac{1}{2}(\overline{%
m_{11}t_{22}}-\overline{t_{21}m_{12}}+(-\tau _{21}+\tau _{22})m_{13})),$

$([(D_{2},M_{2}),T_{1}])_{T}=F_{1}(\textrm{g}_{d}(D_{2})t_{11})+F_{2}(\nu
\textrm{g}_{d}(D_{2})t_{12})+F_{3}(\nu ^{2}\textrm{g}_{d}(D_{2})t_{13})$

\ \ \ \ \ \ \ \ \ \ \ \ \ \ \ \ \ \ \ \ $\
+(-(m_{22},t_{12})+(m_{23},t_{13}))E_{1}+(-(m_{23},t_{13})+(m_{21},t_{11}))E_{2} 
$

$\ \ \ \ \ \ \ \ \ \ \ \ \ \ \ \ \ \ \ \
+(-(m_{21},t_{11})+(m_{22},t_{12}))E_{3}$

$\ \ \ \ \ \ \ \ \ \ \ \ \ \ \ \ \ \ \ \ +F_{1}(\frac{1}{2}(\overline{%
m_{22}t_{13}}-\overline{t_{12}m_{23}}+(-\tau _{11}-2\tau _{12})m_{21}))$

\ \ \ \ \ \ \ \ \ \ \ \ \ \ \ \ \ \ \ $\ +F_{2}(\frac{1}{2}(\overline{%
m_{23}t_{11}}-\overline{t_{13}m_{21}}+(2\tau _{11}+\tau _{12})m_{22}))$

\ \ \ \ \ \ \ \ \ \ \ \ \ \ \ \ \ \ \ $\ +F_{3}(\frac{1}{2}(\overline{%
m_{21}t_{12}}-\overline{t_{11}m_{22}}+(-\tau _{11}+\tau _{12})m_{23})),$

\noindent
then $[\phi _{1},\phi _{2}]_{6}=([\phi _{1},\phi _{2}]_{6D},[\phi _{1},\phi
_{2}]_{6M},[\phi _{1},\phi _{2}]_{6T})$ is a Lie bracket of \gR$_{6}$\ 
and also \gR$_{6}$ is isomorphic to 
\ge$_{6,1}$.\ \ \ \ \emph{Q.E.D.}

\bigskip

\emph{Corollary 3.11.1. \ }The definition of the Lie bracket $[\phi
_{1},\phi _{2}]_{6}$ in \emph{Definition 3.10} can be
extended to $T_{i}\in $\gJ$(i=1,2).$

\bigskip

\emph{Proof.} \ Evidently by \emph{Corollary 3.7.1 }and\emph{\
Corollary 3.8.1} \ \ \ \ \emph{Q.E.D.}

\bigskip 

For $\phi \in $\gR$_{6},(\phi )_{D},(\phi )_{M}$ and $%
(\phi )_{T}$ means \gD$_{4}$ element of $\phi ,$\gA
\ element of $\phi $ and \gJ$_{0}$ element
of $\phi $ respectively. And $\phi ^{6}\in e_{6,1},(\phi
^{6})_{D},(\phi ^{6})_{M}$ and $(\phi ^{6})_{T}$ means \gD$_{4}$
element of $\textrm{f}_{y}^{-1}(\phi ^{6}),$
\gA \ element of $\textrm{f}_{y}^{-1}(\phi ^{6})$ and \gJ$_{0}$
element of $\textrm{f}_{y}^{-1}(\phi ^{6})$ respectively.

\bigskip

\emph{Lemma 3.12. \ }$[\phi _{1},\phi _{2}]_{6D},[\phi _{1},\phi
_{2}]_{6M},[\phi _{1},\phi _{2}]_{6T}$ are expressed by the following
expression respectively.

$[\phi _{1},\phi _{2}]_{6D}=[D1,D2]-\textrm{JD}(m_{11},m_{21})$
$-\textrm{d}_{g}\nu^{2}\textrm{g}_{d}(\textrm{JD}(m_{12},m_{22}))$

\ \ \ \ \ \ \ \ \ \ \ \ \ \ \ 
$-\textrm{d}_{g}\nu \textrm{g}_{d}(\textrm{JD}(m_{13},m_{23}))$
$\ +\textrm{JD}(t_{11},t_{21})+\textrm{d}_{g}\nu^{2}\textrm{g}_{d}(\textrm{JD}(t_{12},t_{22}))$

\ \ \ \ \ \ \ \ \ \ \ \ \ \ \ 
$+\textrm{d}_{g}\nu \textrm{g}_{d}(\textrm{JD}(t_{13},t_{23}))$

$[\phi _{1},\phi_{2}]_{6M}=A_{1}(\textrm{g}_{d}(D_{1})m_{21}-\textrm{g}_{d}(D_{2})m_{11})+A_{2}(\nu
\textrm{g}_{d}(D_{1})m_{22}-\nu \textrm{g}_{d}(D_{2})m_{12})$

$\ \ \ \ \ \ \ \ \ \ \ \ \ \ \ \ +A_{3}(\nu ^{2}\textrm{g}_{d}(D_{1})m_{23}-\nu
^{2}\textrm{g}_{d}(D_{2})m_{13})$

\ \ \ \ \ \ \ \ \ \ \ \ \ \ \ $\ +A_{1}(-\frac{1}{2}\overline{m_{12}m_{23}}+%
\frac{1}{2}\overline{m_{22}m_{13}})+A_{2}(-\frac{1}{2}\overline{m_{13}m_{21}}%
+\frac{1}{2}\overline{m_{23}m_{11}})$

\ \ \ \ \ \ \ \ \ \ \ \ \ \ \ $\ +A_{3}(-\frac{1}{2}\overline{m_{11}m_{22}}+%
\frac{1}{2}\overline{m_{21}m_{12}})$

\ \ \ \ \ \ \ \ \ \ \ \ \ \ \ $\ +A_{1}(\frac{1}{2}((\tau _{11}+2\tau
_{12})t_{21}-(\tau _{21}+2\tau _{22})t_{11}-\overline{t_{12}t_{23}}+%
\overline{t_{22}t_{13}})$

$\ \ \ \ \ \ \ \ \ \ \ \ \ \ \ \ +A_{2}(\frac{1}{2}((-2\tau _{11}-\tau
_{12})t_{22}-(-2\tau _{21}-\tau _{22})t_{12}+\overline{t_{23}t_{11}}-%
\overline{t_{13}t_{21}})$

$\ \ \ \ \ \ \ \ \ \ \ \ \ \ \ \ +A_{3}(\frac{1}{2}((\tau _{11}-\tau
_{12})t_{23}-(\tau _{21}-\tau _{22})t_{13}-\overline{t_{11}t_{22}}+\overline{%
t_{21}t_{12}})$

$[\phi _{1},\phi _{2}]_{6T}=$\ $F_{1}(\textrm{g}_{d}(D_{1})t_{21})+F_{2}(\nu
\textrm{g}_{d}(D_{1})t_{22})+F_{3}(\nu ^{2}\textrm{g}_{d}(D_{1})t_{23})$

\ \ \ \ \ \ \ \ \ \ \ \ \ \ $\
+(-(m_{12},t_{22})+(m_{13},t_{23}))E_{1}+(-(m_{13},t_{23})+(m_{11},t_{21}))E_{2} 
$

$\ \ \ \ \ \ \ \ \ \ \ \ \ \ \ +(-(m_{11},t_{21})+(m_{12},t_{22}))E_{3}$

$\ \ \ \ \ \ \ \ \ \ \ \ \ \ \ +F_{1}(\frac{1}{2}(\overline{m_{12}t_{23}}-%
\overline{t_{22}m_{13}}+(-\tau _{21}-2\tau _{22})m_{11}))$

\ \ \ \ \ \ \ \ \ \ \ \ \ \ $\ +F_{2}(\frac{1}{2}(\overline{m_{13}t_{21}}-%
\overline{t_{23}m_{11}}+(2\tau _{21}+\tau _{22})m_{12}))$

\ \ \ \ \ \ \ \ \ \ \ \ \ \ $\ +F_{3}(\frac{1}{2}(\overline{m_{11}t_{22}}-%
\overline{t_{21}m_{12}}+(-\tau _{21}+\tau _{22})m_{13}))$

$\ \ \ \ \ \ \ \ \ \ \ \ \ \ -(F_{1}(\textrm{g}_{d}(D_{2})t_{11})+F_{2}(\nu
\textrm{g}_{d}(D_{2})t_{12})+F_{3}(\nu ^{2}\textrm{g}_{d}(D_{2})t_{13})$

\ \ \ \ \ \ \ \ \ \ \ \ \ \ $\
+(-(m_{22},t_{12})+(m_{23},t_{13}))E_{1}+(-(m_{23},t_{13})+(m_{21},t_{11}))E_{2} 
$

$\ \ \ \ \ \ \ \ \ \ \ \ \ \ \ +(-(m_{21},t_{11})+(m_{22},t_{12}))E_{3}$

$\ \ \ \ \ \ \ \ \ \ \ \ \ \ \ +F_{1}(\frac{1}{2}(\overline{m_{22}t_{13}}-%
\overline{t_{12}m_{23}}+(-\tau _{11}-2\tau _{12})m_{21}))$

\ \ \ \ \ \ \ \ \ \ \ \ \ \ $\ +F_{2}(\frac{1}{2}(\overline{m_{23}t_{11}}-%
\overline{t_{13}m_{21}}+(2\tau _{11}+\tau _{12})m_{22}))$

\ \ \ \ \ \ \ \ \ \ \ \ \ \ $\ +F_{3}(\frac{1}{2}(\overline{m_{21}t_{12}}-%
\overline{t_{11}m_{22}}+(-\tau _{11}+\tau _{12})m_{23}))$ $\ )$

\bigskip

\emph{Proof}. By \emph{Lemma 2.15} amd \emph{Lemma 3.11} , we have
the above expression.\ \ \ \ \emph{Q.E.D.}

\bigskip

\emph{Definition 3.13.}  For $A,B\in $\gJ$^{\C}$, we define an
operator $A\vee B$ by

$\ \ \ \ \ \ \ \ \ \ \ \ \ \ \ \ \ \ \ \ \ \ \ \ \ A\vee B=[\widetilde{A},%
\widetilde{B}]+(A\circ B-\frac{1}{3}(A,B)E)^{\widetilde{}}.$

\bigskip

\emph{Proposition 3.14.} (I.Yokota\cite[\emph{Proposition2.4.1.}]{Yokota1}) 

$(1)$ For $A\in $%
\gJ$^{\C},tr(A)=0$, we have $\widetilde{A}$ $\in$ \ge$%
_{6}^{\C}$.

$(2)$ For $A,B\in $\gJ$^{\C}$, we have $[\widetilde{A},\widetilde{B}%
]\in $\gf$_{4}^{\C}$.

\bigskip

Since $(A\circ B-\frac{1}{3}(A,B)E)\in $\gJ$_{0}^{\C}$ and $[%
\widetilde{A},\widetilde{B}]\in $\gf$_{4}^{\C},$ by \emph{Proposition 3.14 }%
we have $A\vee B\in $\ge$_{6}^{\C}$.

\bigskip

\emph{Lemma 3.15. \ }For $A(\alpha ,a),B(\beta ,b)\in $\gJ$%
^{\C},$ $A\vee B$ is expressed by the following expression.

$A\vee B=\delta _{d}(\textrm{g}_{d}((A\vee B)_{D}))+(A\vee B)_{M}^{\widehat{}}+(A\vee
B)_{T}^{\widetilde{}},$

$(A\vee B)_{D}=\textrm{JD}(a_{1},b_{1})+\textrm{d}_{g}\nu ^{2}\textrm{g}_{d}(\textrm{JD}(a_{2},b_{2}))+\textrm{d}_{g}\nu
\textrm{g}_{d}(\textrm{JD}(a_{3},b_{3})),$

$(A\vee B)_{M}=\frac{1}{2}A_{1}((\alpha _{2}-\alpha _{3})b_{1}-(\beta
_{2}-\beta _{3})a_{1}-\overline{a_{2}b_{3}}+\overline{b_{2}a_{3}})$

$\ \ \ \ \ \ \ \ \ \ \ \ \ \ \ +\frac{1}{2}A_{2}((\alpha _{3}-\alpha
_{1})b_{2}-(\beta _{3}-\beta _{1})a_{2}-\overline{a_{3}b_{1}}+\overline{%
b_{3}a_{1}})$

$\ \ \ \ \ \ \ \ \ \ \ \ \ \ \ +\frac{1}{2}A_{3}((\alpha _{1}-\alpha
_{2})b_{3}-(\beta _{1}-\beta _{2})a_{3}-\overline{a_{1}b_{2}}+\overline{%
b_{1}a_{2}}),$

$(A\vee B)_{T}=\frac{1}{3}(\alpha _{1}\beta _{1}-\alpha _{2}\beta
_{2}-(a_{1},b_{1})+(a_{2},b_{2}))(E_{1}-E_{2})$

$\ \ \ \ \ \ \ \ \ \ \ \ \ \ +\frac{1}{3}(\alpha _{2}\beta _{2}-\alpha
_{3}\beta _{3}-(a_{2},b_{2})+(a_{3},b_{3}))(E_{2}-E_{3})$

$\ \ \ \ \ \ \ \ \ \ \ \ \ \ +\frac{1}{3}(\alpha _{3}\beta _{3}-\alpha
_{1}\beta _{1}-(a_{3},b_{3})+(a_{1},b_{1}))(E_{3}-E_{1})$

$\ \ \ \ \ \ \ \ \ \ \ \ \ \ +\frac{1}{2}F_{1}((\alpha _{2}+\alpha
_{3})b_{1}+(\beta _{2}+\beta _{3})a_{1}+\overline{b_{2}a_{3}}+\overline{%
a_{2}b_{3}})$

$\ \ \ \ \ \ \ \ \ \ \ \ \ \ +\frac{1}{2}F_{2}((\alpha _{3}+\alpha
_{1})b_{2}+(\beta _{3}+\beta _{1})a_{2}+\overline{b_{3}a_{1}}+\overline{%
a_{3}b_{1}})$

$\ \ \ \ \ \ \ \ \ \ \ \ \ \ +\frac{1}{2}F_{3}((\alpha _{1}+\alpha
_{2})b_{3}+(\beta _{1}+\beta _{2})a_{3}+\overline{b_{1}a_{2}}+\overline{%
a_{1}b_{2}}).$

\bigskip

\emph{Proof. }For $A(\alpha ,a),B(\beta ,b)\in $\gJ$^{\C}$,by 
\emph{Lemma 3.11}. we have

$\ \ \ \ \ \ \ \ \ \ \ \ \ \ \ \ A\vee B=\delta _{d}(\textrm{g}_{d}((A\vee
B)_{D}))+(A\vee B)_{M}^{\widehat{}}+(A\vee B)_{T}^{\widetilde{}}$ .

\noindent
Then, by \emph{Definition  3.13} of $A\vee B$ and \emph{Proposition 3.14} , we have

$\ \ \ \ \ \ \ \ \ \ \ \ \ \ \delta _{d}(\textrm{g}_{d}((A\vee B)_{D}))+(A\vee
B)_{M}^{\widehat{}}=[\widetilde{A},\widetilde{B}]$ ,

\ \ \ \ \ \ \ \ \ \ \ \ \ $(A\vee B)_{D}=([\widetilde{A},\widetilde{B}])_{D} 
$ ,

$\ \ \ \ \ \ \ \ \ \ \ \ \ \ (A\vee B)_{M}=([\widetilde{A},\widetilde{B}%
])_{M}$ ,

$\ \ \ \ \ \ \ \ \ \ \ \ \ \ (A\vee B)_{T}^{\widetilde{}}=(A\circ B-\frac{%
1}{3}(A,B)E)^{\widetilde{}}$.

\noindent
By calculation , we have$.$

$[\widetilde{A}-\frac{1}{3}tr(A)\widetilde{E},\widetilde{B}-\frac{1}{3}tr(B)%
\widetilde{E}]=[\widetilde{A},\widetilde{B}]-\frac{1}{3}tr(A)[\widetilde{E},%
\widetilde{B}]-\frac{1}{3}tr(B)[\widetilde{A},\widetilde{E}]$

\ \ \ \ \ \ \ \ \ \ \ \ \ \ \ \ \ \ \ \ \ \ \ \ \ \ \ \ \ \ \ \ \ \ \ \ \ \ 
\ \ \ $+\frac{1}{9}tr(A)tr(B)[\widetilde{E},\widetilde{E}],$

\ \ \ \ \ \ \ \ \ \ \ \ \ \ \ \ \ \ \ \ \ \ \ \ \ \ \ \ \ \ \ \ \ \ \ \ \ \
\ $=[\widetilde{A},\widetilde{B}].$

\noindent
Since $(A-\frac{1}{3}tr(A)E),(B-\frac{1}{3}tr(B)E)\in $\gJ$_{0}^{\C},$%
we have

$([\widetilde{A},\widetilde{B}])_{D}=([A-\frac{1}{3}tr(A)E,B-\frac{1}{3}%
tr(B)E])_{D}$

$([\widetilde{A},\widetilde{B}])_{M}=([A-\frac{1}{3}tr(A)E,B-\frac{1}{3}%
tr(B)E])_{M}$ .

\noindent
By the definitions $([T_{1},T_{2}])_{D}$ and $([T_{1},T_{2}])_{M}$ \ of 
\emph{Lemma 3.11} ,we have

$(A\vee B)_{D}=([A-\frac{1}{3}tr(A)E,B-\frac{1}{3}tr(B)E])_{D}$

$=([\frac{1}{3}((2\alpha _{1}-\alpha _{2}-\alpha _{3})E_{1}+(-\alpha
_{1}+2\alpha _{2}-\alpha _{3})E_{2}+(-\alpha _{1}-\alpha _{2}+2\alpha
_{3})E_{3}$

$\ \ \ \ +F_{1}(a_{1})+F_{2}(a_{2})+F_{3}(a_{3}),$

\ \ \ $\frac{1}{3}((2\beta _{1}-\beta _{2}-\beta _{3})E_{1}+(-\beta
_{1}+2\beta _{2}-\beta _{3})E_{2}+(-\beta _{1}-\beta _{2}+2\beta _{3})E_{3}$

\ $\ +F_{1}(b_{1})+F_{2}(b_{2})+F_{3}(b_{3})])_{D}$

$=\textrm{JD}(a_{1},b_{1})+\textrm{d}_{g}\nu ^{2}\textrm{g}_{d}(\textrm{JD}(a_{2},b_{2}))+\textrm{d}_{g}\nu
\textrm{g}_{d}(\textrm{JD}(a_{3},b_{3})),$

$(A\vee B)_{M}=([A-\frac{1}{3}tr(A)E,B-\frac{1}{3}tr(B)E])_{M}$

$=\frac{1}{2}A_{1}((\alpha _{2}-\alpha _{3})b_{1}-(\beta _{2}-\beta
_{3})a_{1}-\overline{a_{2}b_{3}}+\overline{b_{2}a_{3}})$

$+\frac{1}{2}A_{2}((\alpha _{3}-\alpha _{1})b_{2}-(\beta _{3}-\beta
_{1})a_{2}-\overline{a_{3}b_{1}}+\overline{b_{3}a_{1}})$

$+\frac{1}{2}A_{3}((\alpha _{1}-\alpha _{2})b_{3}-(\beta _{1}-\beta
_{2})a_{3}-\overline{a_{1}b_{2}}+\overline{b_{1}a_{2}})$.

\noindent
Furthermore by calculation of the following formula,

$(A\vee B)_{T}=A\circ B-\frac{1}{3}(A,B)E,$

\noindent
we have the expression of \emph{Lemma 3.15 }.\ \ \ \ \emph{Q.E.D.}

\bigskip

\section{Definition of the vector space \gR$_{7}$ and it's Lie bracket}

\bigskip 

\ \ \ \emph{Definition 4.1. } We define an $\R$-vector space \gP%
, called the Freudenthal $\R$-vector space by

\ \ \ \ \ \ \ \ \ \ \ \ \ \ \ \ \ \ \ \ \ \ \ \ \ \ \ \ \ \ \gP\  $= $ \gJ\  $\oplus$ \gJ\  $\oplus\  \R\oplus \R.$

An element of \gP\ is often denoted by $(X,Y,\xi ,\eta )$.

For $P=(X,Y,\xi ,\eta )$ , $(P)_{X}$ , $(P)_{Y}$ , $(P)_{\xi }$ and $%
(P)_{\eta }$ means $X$ element of $P,Y$ element of $P,$
$\xi $ element of $P$ and $\eta $ element of $P$ respectively.

\bigskip

\emph{Definition 4.2.} \ In \gP, we define an inner
product $(P,Q)$ and a skew-symmetric inner
product $\{P,Q\}$ respectively by

$\ \ \ \ \ \ \ \ \ \ \ \ \ \ \ \ \ \ \ \ \ \ \ \ \ \ \ \ (P,Q)=(X,Z)+(Y,W)+%
\xi \zeta +\eta \omega ,$

$\ \ \ \ \ \ \ \ \ \ \ \ \ \ \ \ \ \ \ \ \ \ \ \ \ \ \ \ \{P,Q\}=(X,W)-(Z,Y)+%
\xi \omega -\zeta \eta ,$

\noindent
where $P=(X,Y,\xi ,\eta ),Q=(Z,W,\zeta ,\omega )\in $\gP.

\bigskip

\emph{Definition 4.3.} \ For $\phi \in $\ge$_{6,1},A,B\in $%
\gJ$,\rho \in \R$, we define a $\R$-linear mapping $\Phi (\phi
,A,B,\rho )$ :
\gP\ $\rightarrow$ \gP\ by

$\Phi (\phi ,A,B,\rho )\left( 
\begin{array}{c}
X \\ 
Y \\ 
\xi \\ 
\eta%
\end{array}%
\right) =\left( 
\begin{array}{cccc}
\phi \text{-}\frac{1}{3}\rho & 2B & 0 & A \\ 
2A & \text{-}\phi^{t} \text{+}\frac{1}{3}\rho & B & 0 \\ 
0 & A & \rho & 0 \\ 
B & 0 & 0 & -\rho%
\end{array}%
\right) \left( 
\begin{array}{c}
X \\ 
Y \\ 
\xi \\ 
\eta%
\end{array}%
\right) ,$

\ \ \ \ \ \ \ \ \ \ \ \ \ \ \ \ \ \ \ \ \ \ \ \ \ \ \ \ \ $=\left( 
\begin{array}{c}
\phi X\text{-}\frac{1}{3}\rho X+2B\times Y+\eta A \\ 
2A\times X\text{-}\phi^{t} Y\text{+}\frac{1}{3}\rho Y+\xi B \\ 
(A,Y)+\rho \xi \\ 
(B,X)-\rho \eta%
\end{array}%
\right) ,$

\noindent
where $(\phi^{t} X,Y)=(X,\phi Y)$.

\bigskip

\emph{Definition 4.4. } For $P=(X,Y,\xi ,\eta ),Q=(Z,W,\zeta
,\omega )\in $\gP, we define a $\R$-linear
mapping $P\times Q:\ $\gP\ $\rightarrow $\gP\ \ by

$\ \ \ \ \ \ \ \ \ \ \ \ \ \ \ \ \ P\times Q=\Phi (\phi ,A,B,\rho ),$

\ \ \ \ \ \ \ \ \ \ \ \ \ \ \ \ \ \ \ \ \ \ \ \ \ \ \ \ \ $\phi =-\frac{1}{2}%
(X\vee W+Z\vee Y)$

\ \ \ \ \ \ \ \ \ \ \ \ \ \ \ \ \ \ \ \ \ \ \ \ \ \ \ \ \ $A=-\frac{1}{4}%
(2Y\times W-\xi Z-\zeta X)$

\ \ \ \ \ \ \ \ \ \ \ \ \ \ \ \ \ \ \ \ \ \ \ \ \ \ \ \ \ $B=\ \ \frac{1}{4}%
(2X\times Z-\eta W-\omega Y)$

\ \ \ \ \ \ \ \ \ \ \ \ \ \ \ \ \ \ \ \ \ \ \ \ \ \ \ \ $\ \rho =\ \ \ \frac{%
1}{8}((X,W)+(Z,Y)-3(\xi \omega +\zeta \eta )).$

\bigskip

\emph{Proposition 4.5.}\ (I.Yokota\cite[\emph{Theorem 4.9.2.}]{Yokota1})

$\ \ \ \ \ E_{7}=\{\alpha \in Iso_{\C}($\gP%
$^{\C})|\alpha (P\times Q)\alpha ^{-1}=\alpha P\times \alpha Q,<\alpha P,%
\alpha Q>=<P,Q>\}$

\noindent
is a simply connected compact Lie group of type E$_{7}$.

\bigskip

\emph{Definition 4.6.} \ We define the Lie Group $E_{7}^{\C}$ by

$\ \ \ \ \ \ \ \ \ \ \ \ \ \ \ \ \ \ \ \ \ E_{7}^{\C}=\{\alpha \in Iso_{%
\C}($\gP$^{\C}) \mid \alpha (P\times Q)\alpha ^{-1}=\alpha P\times 
\alpha Q\}$.

\bigskip

\emph{Proposition 4.7.} (I.Yokota\cite[\emph{Theorem 4.14.1.}]{Yokota1}) 

The polar decomposition of the Lie group $E_{7}^{\C}$
is given by

\ \ \ \ \ \ \ \ \ \ \ \ \ \ \ \ \ \ \ \ \ \ \ \ \ \ \ \ \ \ $E_{7}^{\C}\simeq
E_{7}\times \R^{133}.$

\noindent
In particular, $E_{7}^{\C}$ is a simply connected complex Lie group of type $%
E_{7}$.

\bigskip

\emph{Lemma 4.8. }For $P=(X(\chi ,x),Y(\gamma ,y),\xi ,\eta
),Q=(Z(\mu ,z),W(\psi ,w),\zeta ,\omega )$

$\in $\gP,\ $P\times Q$ is expressed by the following expression.

$P\times Q=\Phi (\phi ,A,B,\rho ),$

$\phi =\delta _{d}(\textrm{g}_{d}((P\times Q)_{D}))+(P\times Q)_{M}^{\widehat{}%
}+(P\times Q)_{T}^{\widetilde{}},$

$A=(P\times Q)_{A},B=(P\times Q)_{B},\rho
=(P\times Q)_{\rho },$

\noindent
where

$(P\times Q)_{D}=-\frac{1}{2}(\textrm{JD}(x_{1},w_{1})+\textrm{d}_{g}\nu
^{2}\textrm{g}_{d}(\textrm{JD}(x_{2},w_{2}))+\textrm{d}_{g}\nu \textrm{g}_{d}(\textrm{JD}(x_{3},w_{3}))$

$\ \ \ \ \ \ \ \ \ \ \ \ \ \ \ \ \ +\textrm{JD}(z_{1},y_{1})+\textrm{d}_{g}\nu
^{2}\textrm{g}_{d}(\textrm{JD}(z_{2},y_{2}))+\textrm{d}_{g}\nu \textrm{g}_{d}(\textrm{JD}(z_{3},y_{3}))),$

$(P\times Q)_{M}=$
$-\frac{1}{4}A_{1}((\chi _{2}-\chi _{3})w_{1}-(\psi _{2}-\psi _{3})x_{1}-%
\overline{x_{2}w_{3}}+\overline{w_{2}x_{3}}$

\ \ \ \ \ \ \ \ \ \ \ \ \ \ \ \ 
$+(\mu _{2}-\mu _{3})y_{1}-(\gamma
_{2}-\gamma _{3})z_{1}-\overline{z_{2}y_{3}}+\overline{y_{2}z_{3}})$

\ \ \ \ \ \ \ \ \ \ \ \ \ \ \ \ 
$-\frac{1}{4}A_{2}((\chi _{3}-\chi _{1})w_{2}-(\psi _{3}-\psi _{1})x_{2}-%
\overline{x_{3}w_{1}}+\overline{w_{3}x_{1}}$

\ \ \ \ \ \ \ \ \ \ \ \ \ \ \ \ 
$+(\mu _{3}-\mu _{1})y_{2}-(\gamma
_{3}-\gamma _{1})z_{2}-\overline{z_{3}y_{1}}+\overline{y_{3}z_{1}})$

\ \ \ \ \ \ \ \ \ \ \ \ \ \ \ \ 
$-\frac{1}{4}A_{3}((\chi _{1}-\chi _{2})w_{3}-(\psi _{1}-\psi _{2})x_{3}-%
\overline{x_{1}w_{2}}+\overline{w_{1}x_{2}}$

\ \ \ \ \ \ \ \ \ \ \ \ \ \ \ \ 
$+(\mu _{1}-\mu _{2})y_{3}-(\gamma
_{1}-\gamma _{2})z_{3}-\overline{z_{1}y_{2}}+\overline{y_{1}z_{2}})$

$(P\times Q)_{T}=$
$-\frac{1}{6}(\chi _{1}\psi _{1}-\chi _{2}\psi
_{2}-(x_{1},w_{1})+(x_{2},w_{2})+\mu _{1}\gamma _{1}-\mu _{2}\gamma
_{2}$

\ \ \ \ \ \ \ \ \ \ \ \ \ \ \ \ 
$-(z_{1},y_{1})+(z_{2},y_{2}))(E_{1}-E_{2})$

\ \ \ \ \ \ \ \ \ \ \ \ \ \ \ \ 
$-\frac{1}{6}(\chi _{2}\psi _{2}-\chi _{3}\psi
_{3}-(x_{2},w_{2})+(x_{3},w_{3})+\mu _{2}\gamma _{2}-\mu _{3}\gamma
_{3}$

\ \ \ \ \ \ \ \ \ \ \ \ \ \ \ \ 
$-(z_{2},y_{2})+(z_{3},y_{3}))(E_{2}-E_{3})$

\ \ \ \ \ \ \ \ \ \ \ \ \ \ \ \ 
$-\frac{1}{6}(\chi _{3}\psi _{3}-\chi _{1}\psi
_{1}-(x_{3},w_{3})+(x_{1},w_{1})+\mu _{3}\gamma _{3}-\mu _{1}\gamma
_{1}$

\ \ \ \ \ \ \ \ \ \ \ \ \ \ \ \ 
$-(z_{3},y_{3})+(z_{1},y_{1}))(E_{3}-E_{1})$

\ \ \ \ \ \ \ \ \ \ \ \ \ \ \ \ 
$-\frac{1}{4}F_{1}((\chi _{2}+\chi _{3})w_{1}+(\psi _{2}+\psi _{3})x_{1}+%
\overline{w_{2}x_{3}}+\overline{x_{2}w_{3}}$

\ \ \ \ \ \ \ \ \ \ \ \ \ \ \ \ 
$+(\mu _{2}+\mu _{3})y_{1}+(\gamma
_{2}+\gamma _{3})z_{1}+\overline{y_{2}z_{3}}+\overline{z_{2}y_{3}})$

\ \ \ \ \ \ \ \ \ \ \ \ \ \ \ \ 
$-\frac{1}{4}F_{2}((\chi _{3}+\chi _{1})w_{2}+(\psi _{3}+\psi _{1})x_{2}+%
\overline{w_{3}x_{1}}+\overline{x_{3}w_{1}}$

\ \ \ \ \ \ \ \ \ \ \ \ \ \ \ \ 
$+(\mu _{3}+\mu _{1})y_{2}+(\gamma
_{3}+\gamma _{1})z_{2}+\overline{y_{3}z_{1}}+\overline{z_{3}y_{1}})$

\ \ \ \ \ \ \ \ \ \ \ \ \ \ \ \ 
$-\frac{1}{4}F_{3}((\chi _{1}+\chi _{2})w_{3}+(\psi _{1}+\psi _{2})x_{3}+%
\overline{w_{1}x_{2}}+\overline{x_{1}w_{2}}$

\ \ \ \ \ \ \ \ \ \ \ \ \ \ \ \ 
$+(\mu _{1}+\mu _{2})y_{3}+(\gamma
_{1}+\gamma _{2})z_{3}+\overline{y_{1}z_{2}}+\overline{z_{1}y_{2}}).$

$(P\times Q)_{A}=-\frac{1}{4}(((\gamma _{2}\psi _{3}+\gamma _{3}\psi
_{2})-2(y_{1},w_{1})-\xi \mu _{1}-\zeta \chi _{1})E_{1}$

$\ \ \ \ \ \ \ \ \ \ \ \ \ \ \ +((\gamma _{3}\psi _{1}+\gamma _{1}\psi
_{3})-2(y_{2},w_{2})-\xi \mu _{2}-\zeta \chi _{2})E_{2}$

\ $\ \ \ \ \ \ \ \ \ \ \ \ \ \ +((\gamma _{1}\psi _{2}+\gamma _{2}\psi
_{1})-2(y_{3},w_{3})-\xi \mu _{3}-\zeta \chi _{3})E_{3}$

\ \ \ \ \ \ \ \ \ \ \ \ \ \ $\ +F_{1}((-\psi _{1}y_{1}-\gamma _{1}w_{1}+%
\overline{y_{2}w_{3}}+\overline{w_{2}y_{3}})-\xi z_{1}-\zeta x_{1})$

\ \ \ \ \ \ \ \ \ \ \ \ \ \ $\ +F_{2}((-\psi _{2}y_{2}-\gamma _{2}w_{2}+%
\overline{y_{3}w_{1}}+\overline{w_{3}y_{1}})-\xi z_{2}-\zeta x_{2})$

\ \ \ \ \ \ \ \ \ \ \ \ \ \ $\ +F_{3}((-\psi _{3}y_{3}-\gamma _{3}w_{3}+%
\overline{y_{1}w_{2}}+\overline{w_{1}y_{2}})-\xi z_{3}-\zeta x_{3})$ $),$

$(P\times Q)_{B}=$ $\ \ \frac{1}{4}(((\chi _{2}\mu _{3}+\chi _{3}\mu
_{2})-2(x_{1},z_{1})-\eta \psi _{1}-\omega \gamma _{1})E_{1}$

$\ \ \ \ \ \ \ \ \ \ \ \ \ \ \ +((\chi _{3}\mu _{1}+\chi _{1}\mu
_{3})-2(x_{2},z_{2})-\eta \psi _{2}-\omega \gamma _{2})E_{2}$

\ $\ \ \ \ \ \ \ \ \ \ \ \ \ \ +((\chi _{1}\mu _{2}+\chi _{2}\mu
_{1})-2(x_{3},z_{3})-\eta \psi _{3}-\omega \gamma _{3})E_{3}$

\ \ \ \ \ \ \ \ \ \ \ \ \ \ $\ +F_{1}((-\mu _{1}x_{1}-\chi _{1}z_{1}+\overline{%
x_{2}z_{3}}+\overline{z_{2}x_{3}})-\eta w_{1}-\omega y_{1})$

\ \ \ \ \ \ \ \ \ \ \ \ \ \ $\ +F_{2}((-\mu _{2}x_{2}-\chi _{2}z_{2}+\overline{%
x_{3}z_{1}}+\overline{z_{3}x_{1}})-\eta w_{2}-\omega y_{2})$

\ \ \ \ \ \ \ \ \ \ \ \ \ \ $\ +F_{3}((-\mu _{3}x_{3}-\chi _{3}z_{3}+\overline{%
x_{1}z_{2}}+\overline{z_{1}x_{2}})-\eta w_{3}-\omega y_{3})$ $),$

$(P\times Q)_{\rho }=$\ \ \ $\frac{1}{8}(\sum_{i=1}^{3}(\chi _{i}\psi
_{i}+2(x_{i},w_{i})+\mu _{i}\gamma _{i}+2(z_{i},y_{i}))-3(\xi \omega +%
\zeta \eta )$ $)$,

\bigskip

\emph{Proof.} \ By \emph{Definition 4.4}, We have

\ \ \ \ \ \ \ \ $(P\times Q)_{D}=(\phi )_{D}=(-\frac{1}{2}(X\vee W+Z\vee
Y))_{D},$

\ \ \ \ \ \ \ \ \ \ \ \ \ \ \ \ \ \ \ \ \ \ \ \ \ \ \ \ \ \ $\ \ =-\frac{1}{2%
}((X\vee W)_{D}+(Z\vee Y)_{D}),$

\ \ \ \ \ \ \ \ $(P\times Q)_{M}=(\phi )_{M}=(-\frac{1}{2}(X\vee W+Z\vee
Y))_{M},$

\ \ \ \ \ \ \ \ \ \ \ \ \ \ \ \ \ \ \ \ \ \ \ \ \ \ \ \ \ \ $\ \ =-\frac{1}{2%
}((X\vee W)_{M}+(Z\vee Y)_{M}),$

\ \ \ \ \ \ \ \ \ $(P\times Q)_{T}=(\phi )_{T}=(-\frac{1}{2}(X\vee W+Z\vee
Y))_{T},$

\ \ \ \ \ \ \ \ \ \ \ \ \ \ \ \ \ \ \ \ \ \ \ \ \ \ \ \ \ \ $\ \ =-\frac{1}{2%
}((X\vee W)_{T}+(Z\vee Y)_{T}),$

\ \ \ \ \ \ \ \ \ $(P\times Q)_{A}=-\frac{1}{4}(2Y\times W-\xi Z-\zeta
X),$

\ \ \ \ \ \ \ \ \ \ $(P\times Q)_{B}=\ \frac{1}{4}(2X\times Z-\eta W-%
\omega Y),$

\ \ \ \ \ \ \ \ \ $(P\times Q)_{\rho }=\ \frac{1}{8}((X,W)+(Z,Y)-3(\xi 
\omega +\zeta \eta )).$

\noindent
By \emph{Lemma 3.15}, we have the expression of $(P\times Q)_{D}$ , $%
(P\times Q)_{M}$ and $(P\times Q)_{T}.$

\noindent
By calculation we have the expression of $(P\times Q)_{A}$ , $(P\times
Q)_{B} $ and $(P\times Q)_{\rho }$ .\ \ \ \ \ \emph{Q.E.D.}

\bigskip

\emph{Proposition 4.9.}

(I.Yokota\cite[\emph{Theorem4.3.1.}]{Yokota1},\ T.Imai and I.Yokota\cite[\emph{Theorem 2.}]{YokotaImai5}) 

(1)The Lie algebra 
\ge$_{7}^{\C}$ of the Lie group $E_{7}^{\C}$
is given by

\ \ \ \ \ \ \ \ \ \ \ \ \ \ge$_{7}^{\C}=\{\Phi (\phi ,A,B,\rho )\in
Hom_{\C}($\gP$^{\C})|\phi \in $\ge$_{6}^{\C},A,B\in $\gJ$%
^{\C},\rho \in \C\}.$

\noindent
The Lie bracket $[\Phi _{1},\Phi _{2}]$ in \ge$_{7}^{\C}$ is
given by

$\ \ \ \ \ \ \ \ \ \ \ \ \ \ \ \ [\Phi (\phi
_{1},A_{1},B_{1},\rho _{1}),\Phi (\phi _{2},A_{2},B_{2},\rho _{2})]=%
\Phi (\phi ,A,B,\rho ),$

\noindent
where

$\ \ \ \ \ \ \ \ \ \ \ \ \ \ \ \ \ \ \phi =[\phi _{1},\phi _{2}]+2A_{1}\vee
B_{2}-2A_{2}\vee B_{1},$

$\ \ \ \ \ \ \ \ \ \ \ \ \ \ \ \ \ \ A=(\phi _{1}+\frac{2}{3}\rho
_{1})A_{2}-(\phi _{2}+\frac{2}{3}\rho _{2})A_{1},$

$\ \ \ \ \ \ \ \ \ \ \ \ \ \ \ \ \ \ B=-(\phi^{t} _{1}+\frac{2}{3}\rho
_{1})B_{2}+(\phi^{t} _{2}+\frac{2}{3}\rho _{2})B_{1},$

$\ \ \ \ \ \ \ \ \ \ \ \ \ \ \ \ \ \ \rho =(A_{1},B_{2})-(B_{1},A_{2}).$

(2)The Lie algebra 
\ge$_{7,1}$ of the Lie group $E_{7(-25)}$
is given by

\ \ \ \ \ \ \ \ \ \ \ \ \ \ge$_{7,1}=\{\Phi (\phi ,A,B,\rho )\in
Hom_{\R}($\gP$) \mid \phi \in $\ge$_{6,1},A,B\in $\gJ$,\rho \in \R\}.$

\bigskip

\emph{Definition 4.10. } We define an $\R$-vector space \gR$_{%
7}$\textbf{\ }by

\ \ \ \ \ \ \ \ \gR$_{7}=\{(D,M,T,A,B,\rho
) \mid D\in $\gD$_{4},M\in $\gA$,T\in $\gJ$%
_{0},A,B\in $\gJ$,\rho \in \R\}.$

\noindent
And we consider \ \gR$_{7}$ as a vector space \gD
$_{4}\oplus $\gA$\oplus $\gJ$_{0}\oplus $\gJ$%
\oplus $\gJ$\oplus \R$

$=\R^{28}\oplus \R^{8}\oplus \R^{8}\oplus \R^{8}\oplus \R\oplus \R\oplus
\R^{8}\oplus \R^{8}\oplus \R^{8}\oplus \R\oplus \R\oplus \R\oplus \R^{8}\oplus
\R^{8}\oplus \R^{8}\oplus $

$\R\oplus \R\oplus \R\oplus \R^{8}\oplus \R^{8}\oplus \R^{8}\oplus \R.$

\bigskip

\emph{Definition 4.11. } For $\Phi
_{1}=(D_{1},M_{1},T_{1},A_{1},B_{1},\rho _{1}),$

$\Phi
_{2}=(D_{2},M_{2},T_{2},A_{2},B_{2},\rho _{2})\in $\gR$_{7},$
We define a bracket operation 

\noindent
$[\Phi _{1},\Phi _{2}]_{7}$ by

$[\Phi _{1},\Phi _{2}]_{7}$

\ \ \ \ $=([\Phi _{1},\Phi _{2}]_{7D},[\Phi
_{1},\Phi _{2}]_{7M},[\Phi _{1},\Phi _{2}]_{7T},[\Phi _{1},%
\Phi _{2}]_{7A},[\Phi _{1},\Phi _{2}]_{7B},[\Phi _{1},\Phi
_{2}]_{7\rho }),$

\noindent
where 

\ \ \ \ \ \ \ \ \ $[\Phi _{1},\Phi _{2}]_{7D}=[\phi _{1},\phi _{2}]_{6D}+2([A_{1}-%
\frac{1}{3}tr(A_{1})E,B_{2}-\frac{1}{3}tr(B_{2})E])_{D}$

$\ \ \ \ \ \ \ \ \ \ \ \ \ \ \ \ \ \ 
\ \ \ \ \ \ -2([A_{2}-\frac{1}{3}tr(A_{2})E,B_{1}-\frac{1}{3}%
tr(B_{1})E])_{D},$

\ \ \ \ \ \ \ \ \ $[\Phi _{1},\Phi _{2}]_{7M}=[\phi _{1},\phi
_{2}]_{6M}+2([A_{1}-\frac{1}{3}tr(A_{1})E,B_{2}-\frac{1}{3}tr(B_{2})E])_{M}$

$\ \ \ \ \ \ \ \ \ \ \ \ \ \ \ \ \ \ 
\ \ \ \ \ \ -2([A_{2}-\frac{1}{3}tr(A_{2})E,B_{1}-\frac{1}{3}%
tr(B_{1})E])_{M},$

\ \ \ \ \ \ \ \ \ $[\Phi _{1},\Phi _{2}]_{7T}=[\phi _{1},\phi
_{2}]_{6T}+2(A_{1}\circ B_{2}-\frac{1}{3}(A_{1},B_{2})E)$

\ \ \ \ \ \ \ \ \ \ \ \ \ \ \ \ \ \ \ \ \ \ \ \ 
$-2(A_{2}\circ
B_{1}-\frac{1}{3}(A_{2},B_{1})E),$

\ \ \ \ \ \ \ \ $[\Phi _{1},\Phi _{2}]_{7A}=(\delta _{d}\textrm{g}_{d}(D_{1})+%
\widehat{M}_{1}+\widetilde{T}_{1}+\frac{2}{3}\rho _{1})A_{2}$

\ \ \ \ \ \ \ \ \ \ \ \ \ \ \ \ \ \ \ \ \ \ \ $-(\delta
_{d}\textrm{g}_{d}(D_{2})+\widehat{M}_{2}+\widetilde{T}_{2}+\frac{2}{3}\rho
_{2})A_{1},$

\ \ \ \ \ \ \ \ $[\Phi _{1},\Phi _{2}]_{7B}=(\delta _{d}\textrm{g}_{d}(D_{1})+%
\widehat{M}_{1}-\widetilde{T}_{1}+\frac{2}{3}\rho _{1})B_{2}$

\ \ \ \ \ \ \ \ \ \ \ \ \ \ \ \ \ \ \ \ \ \ \ $-(\delta
_{d}\textrm{g}_{d}(D_{2})+\widehat{M}_{2}-\widetilde{T}_{2}+\frac{2}{3}\rho
_{2})B_{1},$

\ \ \ \ \ \ \ \ $[\Phi _{1},\Phi _{2}]_{7\rho
}=(A_{1},B_{2})-(B_{1},A_{2}),$

\ \ \ \ \ \ \ \ \ \ \ \ \ \ \ \ \ \  
$\phi _{1}=(D_{1},M_{1},T_{1}),\phi _{2}=(D_{2},M_{2},T_{2}).$

\bigskip

\emph{Lemma 4.12.} \ \gR$_{7}$ is isomorphic to 
\ge$_{7,1}$ under the correspondence

$\textrm{f}_{y}:(D,M,T,A,B,\rho )\in $\gR$_{7}^{\C}\rightarrow \textrm{f}_{y}(D,M,T,A,B,\rho )$

\ \ \ \ \ \ \ \ \ \ \ \ \ \ \ \ \ \ \ \ \ \ \ \ \ \ \ \ \ \ \ \ \ \ \ \ \ \ \ \ \ \ \ 
$=\Phi (\delta _{d}(\textrm{g}_{d}(D))+%
\widehat{M}+\widetilde{T},$
$A,B,\rho )\in $\ge$_{7,1}$.

\bigskip

\emph{Proof.} For $\Phi _{1}=(D_{1},M_{1},T_{1},A_{1},B_{1},\rho
_{1}),\Phi _{2}=(D_{2},M_{2},T_{2},A_{2},B_{2},\rho _{2}),$

\noindent
We have $\phi _{1}=\delta _{d}\textrm{g}_{d}(D_{1})+\widehat{M}_{1}+\widetilde{T}%
_{1}\in \mathbf{e}_{6,1}$ $,\phi _{2}=\delta _{d}\textrm{g}_{d}(D_{2})+\widehat{M}%
_{2}+\widetilde{T}_{2}\in $\ge$_{6,1}.$

\noindent
Hence we have\ $[\Phi (\phi _{1},A_{1},B_{1},\rho _{1}),\Phi (\phi
_{2},A_{2},B_{2},\rho _{2})]=\Phi (\phi ,A,B,\rho ),$where

$\phi =[\phi _{1},\phi _{2}]+2A_{1}\vee B_{2}-2A_{2}\vee B_{1}$ $($by \emph{%
Lemma 4.4}$)$

\ \ \ $= \delta _{d}\textrm{g}_{d}([\phi _{1},\phi _{2}]_{6D})+([\phi _{1},\phi
_{2}]_{6M})^{\widehat{}}+([\phi _{1},\phi _{2}]_{6T})^{\widetilde{}}$

\ \ \ \ $+2\delta _{d}\textrm{g}_{d}([A_{1}-\frac{1}{3}tr(A_{1})E,B_{2}-\frac{1}{3}%
tr(B_{2})E]_{D})$

\ \ \ \ $+2([A_{1}-\frac{1}{3}tr(A_{1})E,B_{2}-\frac{1}{3}%
tr(B_{2})E])_{M}^{\widehat{}}$

\ \ \ \ $+2(A_{1}\circ B_{2}-\frac{1}{3}(A_{1},B_{2})E)^{\widetilde{}}$

\ \ \ \ $-2\delta _{d}\textrm{g}_{d}(([A_{2}-\frac{1}{3}tr(A_{2})E,B_{1}-\frac{1}{3}%
tr(B_{1})E])_{D}))$

\ \ \ \ $-2([A_{2}-\frac{1}{3}tr(A_{2})E,B_{1}-\frac{1}{3}%
tr(B_{1})E])_{M}^{\widehat{}}$

\ \ \ \ $-2(A_{2}\circ B_{1}-\frac{1}{3}(A_{2},B_{1})E)^{\widetilde{}}$ $($by 
\emph{Lemma 3.9}$)$ ,

$A=(\phi _{1}+\frac{2}{3}\rho _{1})A_{2}-(\phi _{2}+\frac{2}{3}\rho
_{2})A_{1},$

$\ \ \ =(\delta _{d}\textrm{g}_{d}(D_{1})+\widehat{M}_{1}+\widetilde{T}_{1}+\frac{2}{3}%
\rho _{1})A_{2}-(\delta _{d}\textrm{g}_{d}(D_{2})+\widehat{M}_{2}+\widetilde{T}_{2}+%
\frac{2}{3}\rho _{2})A_{1},$

$B=-(\phi^{t} _{1}+\frac{2}{3}\rho _{1})B_{2}+(\phi^{t} _{2}+\frac{2}{3}\rho
_{2})B_{1},$

$\ \ \ =(\delta _{d}\textrm{g}_{d}(D_{1})+\widehat{M}_{1}-\widetilde{T}_{1}+\frac{2}{3}%
\rho _{1})B_{2}-(\delta _{d}\textrm{g}_{d}(D_{2})+\widehat{M}_{2}-\widetilde{T}_{2}+%
\frac{2}{3}\rho _{2})B_{1}$ 

\ \ \ \ \ \ $($ by \emph{Proposition 3.5}$)$,

$\rho =(A_{1},B_{2})-(B_{1},A_{2}),$

$\phi _{1}=(D_{1},M_{1},T_{1}),\phi _{2}=(D_{2},M_{2},T_{2}),$

Let we put

$[\Phi _{1},\Phi _{2}]_{7D}=[\phi _{1},\phi _{2}]_{6D}+2([A_{1}-\frac{1%
}{3}tr(A_{1})E,B_{2}-\frac{1}{3}tr(B_{2})E])_{D}$

$\ \ \ \ \ \ \ \ \ \ \ \ \ \ \ \ \ \ \ \ \ \ \ \ \ \ \ \ \ \ \
-2([A_{2}-\frac{1}{3}tr(A_{2})E,B_{1}-\frac{1}{3}%
tr(B_{1})E])_{D},$

$[\Phi _{1},\Phi _{2}]_{7M}=[\phi _{1},\phi _{2}]_{6M}+2([A_{1}-\frac{1%
}{3}tr(A_{1})E,B_{2}-\frac{1}{3}tr(B_{2})E])_{M}$

$\ \ \ \ \ \ \ \ \ \ \ \ \ \ \ \ \ \ \ \ \ \ \ \ \ \ \ \ \ \ \ \
-2([A_{2}-\frac{1}{3}tr(A_{2})E,B_{1}-\frac{1}{3}%
tr(B_{1})E])_{M},$

$[\Phi _{1},\Phi _{2}]_{7T}=[\phi _{1},\phi _{2}]_{6T}+2(A_{1}\circ
B_{2}-\frac{1}{3}(A_{1},B_{2})E)-2(A_{2}\circ B_{1}-\frac{1}{3}%
(A_{2},B_{1})E),$

$[\Phi _{1},\Phi _{2}]_{7A}=(\delta _{d}\textrm{g}_{d}(D_{1})+\widehat{M}_{1}+%
\widetilde{T}_{1}+\frac{2}{3}\rho _{1})A_{2}-(\delta _{d}\textrm{g}_{d}(D_{2})+%
\widehat{M}_{2}+\widetilde{T}_{2}+\frac{2}{3}\rho _{2})A_{1},$

$[\Phi _{1},\Phi _{2}]_{7B}=(\delta _{d}\textrm{g}_{d}(D_{1})+\widehat{M}_{1}-%
\widetilde{T}_{1}+\frac{2}{3}\rho _{1})B_{2}-(\delta _{d}\textrm{g}_{d}(D_{2})+%
\widehat{M}_{2}-\widetilde{T}_{2}+\frac{2}{3}\rho _{2})B_{1},$

$[\Phi _{1},\Phi _{2}]_{7\rho }=(A_{1},B_{2})-(B_{1},A_{2}),$

$\phi _{1}=(D_{1},M_{1},T_{1}),\phi _{2}=(D_{2},M_{2},T_{2}),$

\noindent
then

\ \ $[\Phi _{1},\Phi _{2}]_{7}=([\Phi _{1},\Phi _{2}]_{7D},[\Phi
_{1},\Phi _{2}]_{7M},[\Phi _{1},\Phi _{2}]_{7T},[\Phi _{1},%
\Phi _{2}]_{7A},[\Phi _{1},\Phi _{2}]_{7B},$

\ \ \ \ \ \ \ \ \ \ \ \ \ \ \ \ \ \ \ 
$[\Phi _{1},\Phi_{2}]_{7\rho })$

\noindent
is a Lie bracket of \gR$_{7}$\textbf{\ }and also \textbf{%
\gR}$_{7}$\textbf{\ }is\textbf{\ }isomorphic to \ge$%
_{7,1}$. \ \ \ \ \ \ \ \ \emph{Q.E.D}.

\bigskip

\emph{Corollary 4.12.1. \ }The definition of the Lie bracket $[\Phi
_{1},\Phi _{2}]_{7}$ in \emph{Definition 4.11} can be
extended to $T_{i}\in $\gJ $(i=1,2).$

\bigskip

\emph{Proof.} \ Evidently by \emph{Corollary 3.11.1} \ \ \ \ \emph{%
Q.E.D.}

\bigskip

\emph{Lemma 4.13. } $[\Phi _{1},\Phi _{2}]_{7D},[\Phi _{1},%
\Phi _{2}]_{7M},[\Phi _{1},\Phi _{2}]_{7T},[\Phi _{1},\Phi
_{2}]_{7A},[\Phi _{1},\Phi _{2}]_{7B}$ 

\noindent
and $[\Phi _{1},\Phi
_{2}]_{7\rho }$
are expressed by the following expression respectively, 

\noindent
where, $A(\alpha
,a),B(\beta ,b)\in $\gJ.

$[\Phi _{1},\Phi _{2}]_{7D}=[D1,D2]-\textrm{JD}(m_{11},m_{21})-\textrm{d}_{g}\nu
^{2}\textrm{g}_{d}(\textrm{JD}(m_{12},m_{22}))$

\ \ \ \ \ \ \ \ \ \ \ \ \ \ \ \ \ 
$-\textrm{d}_{g}\nu \textrm{g}_{d}(\textrm{JD}(m_{13},m_{23}))$

\ \ \ \ \ \ \ \ \ \ \ \ \ \ \ \ $\ +\textrm{JD}(t_{11},t_{21})+\textrm{d}_{g}\nu
^{2}\textrm{g}_{d}(\textrm{JD}(t_{12},t_{22}))+\textrm{d}_{g}\nu \textrm{g}_{d}(\textrm{JD}(t_{13},t_{23}))$

$\ \ \ \ \ \ \ \ \ \ \ \ \ \ \ \ \ \ +2(\textrm{JD}(a_{11},b_{21})+\textrm{d}_{g}\nu
^{2}\textrm{g}_{d}(\textrm{JD}(a_{12},b_{22}))+\textrm{d}_{g}\nu \textrm{g}_{d}(\textrm{JD}(a_{13},b_{23})))$

\ \ \ \ \ \ \ \ \ \ \ \ \ \ \ \ \ $\ -2(\textrm{JD}(a_{21},b_{11})+\textrm{d}_{g}\nu
^{2}\textrm{g}_{d}(\textrm{JD}(a_{22},b_{12}))+\textrm{d}_{g}\nu \textrm{g}_{d}(\textrm{JD}(a_{23},b_{13}))),$

$[\Phi _{1},\Phi
_{2}]_{7M}=A_{1}(\textrm{g}_{d}(D_{1})m_{21}-\textrm{g}_{d}(D_{2})m_{11})+A_{2}(\nu
\textrm{g}_{d}(D_{1})m_{22}-\nu \textrm{g}_{d}(D_{2})m_{12})$

$\ \ \ \ \ \ \ \ \ \ \ \ \ \ \ \ +A_{3}(\nu ^{2}\textrm{g}_{d}(D_{1})m_{23}-\nu
^{2}\textrm{g}_{d}(D_{2})m_{13})$

\ \ \ \ \ \ \ \ \ \ \ \ \ \ \ $\ +\frac{1}{2}A_{1}(-\overline{m_{12}m_{23}}+%
\overline{m_{22}m_{13}})+\frac{1}{2}A_{2}(-\overline{m_{13}m_{21}}+\overline{%
m_{23}m_{11}})$

\ \ \ \ \ \ \ \ \ \ \ \ \ \ \ $\ +\frac{1}{2}A_{3}(-\overline{m_{11}m_{22}}+%
\overline{m_{21}m_{12}})$

\ \ \ \ \ \ \ \ \ \ \ \ \ \ \ $\ +\frac{1}{2}A_{1}((\tau _{11}+2\tau
_{12})t_{21}-(\tau _{21}+2\tau _{22})t_{11}-\overline{t_{12}t_{23}}+%
\overline{t_{22}t_{13}})$

$\ \ \ \ \ \ \ \ \ \ \ \ \ \ \ \ +\frac{1}{2}A_{2}((-2\tau _{11}-\tau
_{12})t_{22}-(-2\tau _{21}-\tau _{22})t_{12}-\overline{t_{23}t_{11}}+%
\overline{t_{13}t_{21}})$

$\ \ \ \ \ \ \ \ \ \ \ \ \ \ \ \ +\frac{1}{2}A_{3}((\tau _{11}-\tau
_{12})t_{23}-(\tau _{21}-\tau _{22})t_{13}-\overline{t_{11}t_{22}}+\overline{%
t_{21}t_{12}})$

$\ \ \ \ \ \ \ \ \ \ \ \ \ \ \ \ +A_{1}((\alpha _{12}-\alpha
_{13})b_{21}-(\beta _{22}-\beta _{23})a_{11}-\overline{a_{12}b_{23}}+%
\overline{b_{22}a_{13}}$

$\ \ \ \ \ \ \ \ \ \ \ \ \ \ \ \ \ \ \ \ -(\alpha _{22}-\alpha
_{23})b_{11}+(\beta _{12}-\beta _{13})a_{21}+\overline{a_{22}b_{13}}-%
\overline{b_{12}a_{23}})$

$\ \ \ \ \ \ \ \ \ \ \ \ \ \ \ \ +A_{2}((\alpha _{13}-\alpha
_{11})b_{22}-(\beta _{23}-\beta _{21})a_{12}-\overline{a_{13}b_{21}}+%
\overline{b_{23}a_{11}})$

$\ \ \ \ \ \ \ \ \ \ \ \ \ \ \ \ \ \ \ \ -(\alpha _{23}-\alpha
_{21})b_{12}+(\beta _{13}-\beta _{11})a_{22}+\overline{a_{23}b_{11}}-%
\overline{b_{13}a_{21}})$

$\ \ \ \ \ \ \ \ \ \ \ \ \ \ \ \ +A_{3}((\alpha _{11}-\alpha
_{12})b_{23}-(\beta _{21}-\beta _{22})a_{13}-\overline{a_{11}b_{22}}+%
\overline{b_{21}a_{12}}$

$\ \ \ \ \ \ \ \ \ \ \ \ \ \ \ \ \ \ \ \ -(\alpha _{21}-\alpha
_{22})b_{13}+(\beta _{11}-\beta _{12})a_{23}+\overline{a_{21}b_{12}}-%
\overline{b_{11}a_{22}}),$

$[\Phi _{1},\Phi _{2}]_{7T}=F_{1}(\textrm{g}_{d}(D_{1})t_{21})+F_{2}(\nu
\textrm{g}_{d}(D_{1})t_{22})+F_{3}(\nu ^{2}\textrm{g}_{d}(D_{1})t_{23})$

$\ \ \ \ \ \ \ \ \ \ \ \ \ \ \ \ -(F_{1}(\textrm{g}_{d}(D_{2})t_{11})-F_{2}(\nu
\textrm{g}_{d}(D_{2})t_{12})-F_{3}(\nu ^{2}\textrm{g}_{d}(D_{2})t_{13})$

\ \ \ \ \ \ \ \ \ \ \ \ \ \ $\
+(-(m_{12},t_{22})+(m_{13},t_{23})+(m_{22},t_{12})-(m_{23},t_{13}))E_{1}$

$\ \ \ \ \ \ \ \ \ \ \ \ \ \ \
+(-(m_{13},t_{23})+(m_{11},t_{21})+(m_{23},t_{13})-(m_{21},t_{11}))E_{2}$

$\ \ \ \ \ \ \ \ \ \ \ \ \ \ \
+(-(m_{11},t_{21})+(m_{12},t_{22})+(m_{21},t_{11})-(m_{22},t_{12}))E_{3}$

$\ \ \ \ \ \ \ \ \ \ \ \ \ \ \ +\frac{1}{3}(4\alpha _{11}\beta _{21}-2\alpha
_{12}\beta _{22}-2\alpha _{13}\beta
_{23}-4(a_{11},b_{21})+2(a_{12},b_{22})$

$\ \ \ \ \ \ \ \ \ \ \ \ \ \ \ \ \ \ \ +2(a_{13},b_{23}))E_{1}$

$\ \ \ \ \ \ \ \ \ \ \ \ \ \ -\frac{1}{3}(4\alpha _{21}\beta _{11}-2\alpha
_{22}\beta _{12}-2\alpha _{23}\beta
_{13}-4(a_{21},b_{11})+2(a_{22},b_{12})$

$\ \ \ \ \ \ \ \ \ \ \ \ \ \ \ \ \ \ \ +2(a_{23},b_{13}))E_{1}$

\ \ \ \ \ \ \ \ \ \ \ \ \ \ $\ +\frac{1}{3}(-2\alpha _{11}\beta _{21}+4\alpha
_{12}\beta _{22}-2\alpha _{13}\beta
_{23}+2(a_{11},b_{21})-4(a_{12},b_{22})$

$\ \ \ \ \ \ \ \ \ \ \ \ \ \ \ \ \ \ \ +2(a_{13},b_{23}))E_{2}$

$\ \ \ \ \ \ \ \ \ \ \ \ \ \ \ -\frac{1}{3}(-2\alpha _{21}\beta _{11}+4\alpha
_{22}\beta _{21}-2\alpha _{23}\beta
_{13}+2(a_{21},b_{11})-4(a_{22},b_{12})$

$\ \ \ \ \ \ \ \ \ \ \ \ \ \ \ \ \ \ \ +2(a_{23},b_{13}))E_{2}$

$\ \ \ \ \ \ \ \ \ \ \ \ \ \ \ +\frac{1}{3}(-2\alpha _{11}\beta _{21}-2\alpha
_{12}\beta _{22}+4\alpha _{13}\beta
_{23}+2(a_{11},b_{21})+2(a_{12},b_{22})$

$\ \ \ \ \ \ \ \ \ \ \ \ \ \ \ \ \ \ \ -4(a_{13},b_{23}))E_{3}$

$\ \ \ \ \ \ \ \ \ \ \ \ \ \ \ -\frac{1}{3}(-2\alpha _{21}\beta _{11}-2\alpha
_{22}\beta _{12}+4\alpha _{23}\beta
_{13}+2(a_{21},b_{11})+2(a_{22},b_{12})$

$\ \ \ \ \ \ \ \ \ \ \ \ \ \ \ \ \ \ \ -4(a_{23},b_{13}))E_{3}$

$\ \ \ \ \ \ \ \ \ \ \ \ \ \ \ +F_{1}(\frac{1}{2}(\overline{m_{12}t_{23}}-%
\overline{t_{22}m_{13}}+(-\tau _{21}-2\tau _{22})m_{11}$

\ \ \ \ \ \ \ \ \ \ \ \ \ \ \ \ \ \ \ \ \ \  
$-\overline{m_{22}t_{13}}+\overline{t_{12}m_{23}}-(-\tau _{11}-2\tau _{12})m_{21}))$

\ \ \ \ \ \ \ \ \ \ \ \ \ \ $\ +F_{2}(\frac{1}{2}(\overline{m_{13}t_{21}}-%
\overline{t_{23}m_{11}}+(2\tau _{21}+\tau _{22})m_{12}$

\ \ \ \ \ \ \ \ \ \ \ \ \ \ \ \ \ \ \ \ \ \  
$-\overline{m_{23}t_{11}}+\overline{t_{13}m_{21}}-(2\tau _{11}+\tau _{12})m_{22}))$

\ \ \ \ \ \ \ \ \ \ \ \ \ \ $\ +F_{3}(\frac{1}{2}(\overline{m_{11}t_{22}}-%
\overline{t_{21}m_{12}}+(-\tau _{21}+\tau _{22})m_{13}$

\ \ \ \ \ \ \ \ \ \ \ \ \ \ \ \ \ \ \ 
$-\overline{m_{21}t_{12}}+\overline{t_{11}m_{22}}%
-(-\tau _{11}+\tau _{12})m_{23}))$

$\ \ \ \ \ \ \ \ \ \ \ \ \ \ \ +F_{1}((\alpha _{12}+\alpha _{13})b_{21}+(\beta
_{22}+\beta _{23})a_{11}+\overline{b_{22}a_{13}}+\overline{a_{12}b_{23}}$

$\ \ \ \ \ \ \ \ \ \ \ \ \ \ \ \ \ \ \ -(\alpha _{22}+\alpha
_{23})b_{11}-(\beta _{12}+\beta _{13})a_{21}-\overline{b_{12}a_{23}}-%
\overline{a_{22}b_{13}})$

$\ \ \ \ \ \ \ \ \ \ \ \ \ \ \ +F_{2}((\alpha _{13}+\alpha
_{11})b_{22}+(\beta _{23}+\beta _{21})a_{12}+\overline{b_{23}a_{11}}+%
\overline{a_{13}b_{21}}$

$\ \ \ \ \ \ \ \ \ \ \ \ \ \ \ \ \ \ \ \ -(\alpha _{23}+\alpha
_{21})b_{12}-(\beta _{13}+\beta _{11})a_{22}-\overline{b_{13}a_{21}}-%
\overline{a_{23}b_{11}})$

$\ \ \ \ \ \ \ \ \ \ \ \ \ \ \ +F_{3}((\alpha _{11}+\alpha
_{12})b_{23}+(\beta _{21}+\beta _{22})a_{13}+\overline{b_{21}a_{12}}+%
\overline{a_{11}b_{22}}$

$\ \ \ \ \ \ \ \ \ \ \ \ \ \ \ \ \ \ \ \ -(\alpha _{21}+\alpha
_{22})b_{13}-(\beta _{11}+\beta _{12})a_{23}-\overline{b_{11}a_{22}}-%
\overline{a_{21}b_{12}}),$

$[\Phi _{1},\Phi _{2}]_{7A}=F_{1}(\textrm{g}_{d}(D_{1})a_{21})+F_{2}(\nu
\textrm{g}_{d}(D_{1})a_{22})+F_{3}(\nu ^{2}\textrm{g}_{d}(D_{1})a_{23})$

\ \ \ \ \ \ \ \ \ \ \ \ $\ \ \ -F_{1}(\textrm{g}_{d}(D_{2})a_{11})-F_{2}(\nu
\textrm{g}_{d}(D_{2})a_{12})-F_{3}(\nu ^{2}\textrm{g}_{d}(D_{2})a_{13})$

\ \ \ \ \ \ \ \ \ \ \ \ \ \ \ \ \ $\
+((m_{13},a_{23})-(m_{12},a_{22})-(m_{23},a_{13})+(m_{22},a_{12}))E_{1}$

$\ \ \ \ \ \ \ \ \ \ \ \ \ \ \ \ \ \
+((m_{11},a_{21})-(m_{13},a_{23})-(m_{21},a_{11})+(m_{23},a_{13}))E_{2}$

$\ \ \ \ \ \ \ \ \ \ \ \ \ \ \ \ \ \
+((m_{12},a_{22})-(m_{11},a_{21})-(m_{22},a_{12})+(m_{21},a_{11}))E_{3}$

\ \ \ \ \ \ \ \ \ \ \ \ \ \ \ \ \ $\ +\frac{1}{2}F_{1}((\alpha _{23}-\alpha
_{22})m_{11}-\overline{a_{22}m_{13}}+\overline{m_{12}a_{23}}$

\ \ \ \ \ \ \ \ \ \ \ \ \ \ \ \ \ \ \ \ \ \ \ \ \ 
$-(\alpha_{13}-\alpha _{12})m_{21}+\overline{a_{12}m_{23}}-\overline{m_{22}a_{13}})$

$\ \ \ \ \ \ \ \ \ \ \ \ \ \ \ \ \ \ +\frac{1}{2}F_{2}((\alpha _{21}-\alpha
_{23})m_{12}-\overline{a_{23}m_{11}}+\overline{m_{13}a_{21}}$

\ \ \ \ \ \ \ \ \ \ \ \ \ \ \ \ \ \ \ \ \ \ \ \ \ 
$-(\alpha_{11}-\alpha _{13})m_{22}+\overline{a_{13}m_{21}}-\overline{m_{23}a_{11}})$

$\ \ \ \ \ \ \ \ \ \ \ \ \ \ \ \ \ \ +\frac{1}{2}F_{3}((\alpha _{22}-\alpha
_{21})m_{13}-\overline{a_{21}m_{12}}+\overline{m_{11}a_{22}}$

\ \ \ \ \ \ \ \ \ \ \ \ \ \ \ \ \ \ \ \ \ \ \ \ \ 
$-(\alpha_{12}-\alpha _{11})m_{23}+\overline{a_{11}m_{22}}-\overline{m_{21}a_{12}})$

$\ \ \ \ \ \ \ \ \ \ \ \ \ \ \ \ \ \ \ \ \ \  +(\tau _{11}\alpha
_{21}+(t_{12},a_{22})+(t_{13},a_{23})-\tau _{21}\alpha
_{11}-(t_{22},a_{12})-(t_{23},a_{13})$

\ \ \ \ \ \ \ \ \ \ \ \ \ \ \ \ \ \ \ \   
$+\frac{2}{3}\rho _{1}\alpha _{21}-\frac{2}{3}\rho _{2}\alpha _{11})E_{1}$

$\ \ \ \ \ \ \ \ \ \ \ \ \ \  \ \ +(\tau _{12}\alpha
_{22}+(t_{13},a_{23})+(t_{11},a_{21})-\tau _{22}\alpha
_{12}-(t_{23},a_{13})-(t_{21},a_{11})$

\ \ \ \ \ \ \ \ \ \ \ \ \ \ \ \ \ \ \ \   
$+\frac{2}{3}\rho _{1}\alpha _{22}-\frac{2}{3}\rho _{2}\alpha _{12})E_{2}$

\ \ \ \ \ \ \ \ \ \ \ \ \ \  \ \ 
$+((-\tau _{11}-\tau _{12})\alpha
_{23}+(t_{11},a_{21})+(t_{12},a_{22})$

\ \ \ \ \ \ \ \ \ \ \ \ \ \ \  \ \ 
$-(-\tau _{21}-\tau _{22})\alpha
_{13}-(t_{21},a_{11})-(t_{22},a_{12})$

\ \ \ \ \ \ \ \ \ \ \ \ \ \ \ \ \ \ \ \ 
$+\frac{2}{3}\rho _{1}\alpha _{23}-\frac{2}{3}\rho
_{2}\alpha _{13})E_{3}$

\ \ \ \ \ \ \ \ \ \ \ \ \ \ \ $\ +F_{1}(\frac{1}{2}(-\tau
_{11}a_{21}+(\alpha _{22}+\alpha _{23})t_{11}+\overline{a_{22}t_{13}}+%
\overline{t_{12}a_{23}})+\frac{2}{3}\rho _{1}a_{21}$

\ \ \ \ \ \ \ \ \ \ \ \ \ \ \ \ \ \ $\ \ -\frac{1}{2}(-\tau
_{21}a_{11}+(\alpha _{12}+\alpha _{13})t_{21}+\overline{a_{12}t_{23}}+%
\overline{t_{22}a_{13}})-\frac{2}{3}\rho _{2}a_{11})$

$\ \ \ \ \ \ \ \ \ \ \ \ \ \ \ \ +F_{2}(\frac{1}{2}(-\tau
_{12}a_{22}+(\alpha _{23}+\alpha _{21})t_{12}+\overline{a_{23}t_{11}}+%
\overline{t_{13}a_{21}})+\frac{2}{3}\rho _{1}a_{22}$

$\ \ \ \ \ \ \ \ \ \ \ \ \ \ \ \ \ \ \ \ -\frac{1}{2}(-\tau
_{22}a_{12}+(\alpha _{13}+\alpha _{11})t_{22}+\overline{a_{13}t_{21}}+%
\overline{t_{23}a_{11}})-\frac{2}{3}\rho _{2}a_{12})$

$\ \ \ \ \ \ \ \ \ \ \ \ \ \ \ \ \ \ \ +F_{3}(\frac{1}{2}((\tau _{11}+\tau
_{12})a_{23}+(\alpha _{21}+\alpha _{22})t_{13}+\overline{a_{21}t_{12}}+%
\overline{t_{11}a_{22}})+\frac{2}{3}\rho _{1}a_{23}$

\ \ \ \ \ \ \ \ \ \ \ \ \ \ \ \ \ \ \ \ \ \ \ \ \ \ 
$-\frac{1}{2}((\tau _{21}+\tau
_{22})a_{13}+(\alpha _{11}+\alpha _{12})t_{23}+\overline{a_{11}t_{22}}+%
\overline{t_{21}a_{12}})-\frac{2}{3}\rho _{2}a_{13}),$

$[\Phi _{1},\Phi _{2}]_{7B}=F_{1}(\textrm{g}_{d}(D_{1})b_{21})+F_{2}(\nu
\textrm{g}_{d}(D_{1})b_{22})+F_{3}(\nu ^{2}\textrm{g}_{d}(D_{1})b_{23})$

\ \ \ \ \ \ \ \ \ \ \ \ $\ \ \ -F_{1}(\textrm{g}_{d}(D_{2})b_{11})-F_{2}(\nu
\textrm{g}_{d}(D_{2})b_{12})-F_{3}(\nu ^{2}\textrm{g}_{d}(D_{2})b_{13})$

\ \ \ \ \ \ \ \ \ \ \ \ \ \ \ \ \ $\
+((m_{13},b_{23})-(m_{12},b_{22})-(m_{23},b_{13})+(m_{22},b_{12}))E_{1}$

$\ \ \ \ \ \ \ \ \ \ \ \ \ \ \ \ \ \
+((m_{11},b_{21})-(m_{13},b_{23})-(m_{21},b_{11})+(m_{23},b_{13}))E_{2}$

$\ \ \ \ \ \ \ \ \ \ \ \ \ \ \ \ \ \
+((m_{12},b_{22})-(m_{11},b_{21})-(m_{22},b_{12})+(m_{21},b_{11}))E_{3}$

\ \ \ \ \ \ \ \ \ \ \ \ \ \ \ \ $\ +\frac{1}{2}F_{1}((\beta _{23}-\beta
_{22})m_{11}-\overline{b_{22}m_{13}}+\overline{m_{12}b_{23}}$

\ \ \ \ \ \ \ \ \ \ \ \ \ \ \ \ \ \ \ \ \ \ \ \ \ 
$-(\beta
_{13}-\beta _{12})m_{21}+\overline{b_{12}m_{23}}-\overline{m_{22}b_{13}})$

$\ \ \ \ \ \ \ \ \ \ \ \ \ \ \ \ \ +\frac{1}{2}F_{2}((\beta _{21}-\beta
_{23})m_{12}-\overline{b_{23}m_{11}}+\overline{m_{13}b_{21}}$

\ \ \ \ \ \ \ \ \ \ \ \ \ \ \ \ \ \ \ \ \ \ \ \ \ 
$-(\beta
_{11}-\beta _{13})m_{22}+\overline{b_{13}m_{21}}-\overline{m_{23}b_{11}})$

$\ \ \ \ \ \ \ \ \ \ \ \ \ \ \ \ \ +\frac{1}{2}F_{3}((\beta _{22}-\beta
_{21})m_{13}-\overline{b_{21}m_{12}}+\overline{m_{11}b_{22}}$

\ \ \ \ \ \ \ \ \ \ \ \ \ \ \ \ \ \ \ \ \ \ \ \ \ 
$-(\beta
_{12}-\beta _{11})m_{23}+\overline{b_{11}m_{22}}-\overline{m_{21}b_{12}})$

$\ \ \ \ \ \ \ \ \ \ \ \ \ \ \ \ \ \ -(\tau _{11}\beta
_{21}+(t_{12},b_{22})+(t_{13},b_{23})-\tau _{21}\beta
_{11}-(t_{22},b_{12})-(t_{23},b_{13})$

\ \ \ \ \ \ \ \ \ \ \ \ \ \ \ \ \ \ \ \   
$+\frac{2}{3}\rho _{1}\beta _{21}-\frac{2}{3}\rho _{2}\beta _{11})E_{1}$

$\ \ \ \ \ \ \ \ \ \ \ \ \ \ \ \ \ \ -(\tau _{12}\beta
_{22}+(t_{13},b_{23})+(t_{11},b_{21})-\tau _{22}\beta
_{12}-(t_{23},b_{13})-(t_{21},b_{11})$

\ \ \ \ \ \ \ \ \ \ \ \ \ \ \ \ \ \ \ \   
$+\frac{2}{3}\rho _{1}\beta _{22}-\frac{2}{3}\rho _{2}\beta _{12})E_{2}$

$\ \ \ \ \ \ \ \ \ \ \ \ \ \ \ \ \ -((-\tau _{11}-\tau _{12})\beta
_{23}+(t_{11},b_{21})+(t_{12},b_{22})-(-\tau _{21}-\tau _{22})\beta
_{13}$

\ \ \ \ \ \ \ \ \ \ \ \ \ \ \ \ \ \ \ \   
$-(t_{21},b_{11})-(t_{22},b_{12})$
$+\frac{2}{3}\rho _{1}\beta _{23}-\frac{2}{3}\rho
_{2}\beta _{13})E_{3}$

\ \ \ \ \ \ \ \ \ \ \ \ \ \ \ \ $\ -F_{1}(\frac{1}{2}(-\tau
_{11}b_{21}+(\beta _{22}+\beta _{23})t_{11}+\overline{b_{22}t_{13}}+%
\overline{t_{12}b_{23}})+\frac{2}{3}\rho _{1}b_{21}$

\ \ \ \ \ \ \ \ \ \ \ \ \ \ \ \ \ \ \ $\ \ -\frac{1}{2}(-\tau
_{21}b_{11}+(\beta _{12}+\beta _{13})t_{21}+\overline{b_{12}t_{23}}+%
\overline{t_{22}b_{13}})-\frac{2}{3}\rho _{2}b_{11})$

$\ \ \ \ \ \ \ \ \ \ \ \ \ \ \ \ \ -F_{2}(\frac{1}{2}(-\tau
_{12}b_{22}+(\beta _{23}+\beta _{21})t_{12}+\overline{b_{23}t_{11}}+%
\overline{t_{13}b_{21}})+\frac{2}{3}\rho _{1}b_{22}$

$\ \ \ \ \ \ \ \ \ \ \ \ \ \ \ \ \ \ \ \ \ \ -\frac{1}{2}(-\tau
_{22}b_{12}+(\beta _{13}+\beta _{11})t_{22}+\overline{b_{13}t_{21}}+%
\overline{t_{23}b_{11}})-\frac{2}{3}\rho _{2}b_{12})$

$\ \ \ \ \ \ \ \ \ \ \ \ \ \ \ \ \ \ \ -F_{3}(\frac{1}{2}((\tau _{11}+\tau
_{12})b_{23}+(\beta _{21}+\beta _{22})t_{13}+\overline{b_{21}t_{12}}+%
\overline{t_{11}b_{22}})+\frac{2}{3}\rho _{1}b_{23}$

$\ \ \ \ \ \ \ \ \ \ \ \ \ \ \ \ \ \ \ \ \ \ \ \ -\frac{1}{2}((\tau _{21}+\tau
_{22})b_{13}+(\beta _{11}+\beta _{12})t_{23}+\overline{b_{11}t_{22}}+%
\overline{t_{21}b_{12}})-\frac{2}{3}\rho _{2}b_{13}),$

$[\Phi _{1},\Phi _{2}]_{7\rho }=\alpha _{11}\beta _{21}+\alpha
_{12}\beta _{22}+\alpha _{13}\beta _{23}-\alpha _{21}\beta _{11}-\alpha
_{22}\beta _{12}-\alpha _{23}\beta _{13}$

$\ \ \ \ \ \ \ \ \ \ \ \ \ \ \ \ \ \ 
+2(a_{11},b_{21})+2(a_{12},b_{22})+2(a_{13},b_{23})-2(a_{21},b_{11})-2(a_{22},b_{12})$

$\ \ \ \ \ \ \ \ \ \ \ \ \ \ \ \ \ -2(a_{23},b_{13}). $

\bigskip

\emph{Proof.} \ By \emph{Lemma 4.12} and \emph{Lemma 3.15} , we have the
expression of $[\Phi _{1},\Phi _{2}]_{7D}$ ,
$[\Phi _{1},\Phi _{2}]_{7M}$ and $[\Phi _{1},\Phi _{2}]_{7T}$ .

\noindent
By \emph{Lemma 4.12} and furthermore calculation, we have the expression of 

\noindent
$[\Phi _{1},\Phi _{2}]_{7A},$
$[\Phi _{1},\Phi _{2}]_{7B}$ and $[\Phi _{1},\Phi _{2}]_{7\rho }$
.\ \ \ \ \ \ \emph{Q.E.D.}

\bigskip

\section{Definition of the vector space \gR$_{8}$ and it's Lie bracket}

\bigskip 

\ \ \ \ \emph{Proposition 5.1.} (I.Yokota\cite[\emph{Theorem 5.1.1.}]{Yokota1})

 In a $133+56\times
2+3=248$-dimensional $\C$-vector space

\ \ \ \ \ \ \ \ \ \ \ \ \ \ \ \ \ \ \ \ \ \ \ge$_{8}^{\C}=$%
\ge$_{7}^{\C}\oplus $\gP$^{\C}\oplus $\gP$%
^{\C}\oplus \C\oplus \C\oplus \C$,

\noindent
if we define a Lie bracket $[R_{1}^{8},R_{2}^{8}]$ by

$\ \ \ \ \ \ \ \ \ \ \ \ [(\Phi _{1}^{7},P_{1},Q_{1},r_{1},s_{1},u_{1}),(%
\Phi _{2}^{7},P_{2},Q_{2},r_{2},s_{2},u_{2})]=(\Phi ^{7},P,Q,r,s,u),$

\noindent
where

$\ \ \ \ \ \ \ \ \ \ \ \ \Phi ^{7}=[\Phi _{1}^{7},\Phi
_{2}^{7}]+P_{1}\times Q_{2}-P_{2}\times Q_{1},$

$\ \ \ \ \ \ \ \ \ \ \ \ P=\Phi _{1}^{7}P_{2}-\Phi
_{2}^{7}P_{1}+r_{1}P_{2}-r_{2}P_{1}+s_{1}Q_{2}-s_{2}Q_{1},$

$\ \ \ \ \ \ \ \ \ \ \ \ Q=\Phi _{1}^{7}Q_{2}-\Phi
_{2}^{7}Q_{1}-r_{1}Q_{2}+r_{2}Q_{1}+u_{1}P_{2}-u_{2}P_{1},$

$\ \ \ \ \ \ \ \ \ \ \ \ \ r=-\frac{1}{8}\{P_{1},Q_{2}\}+\frac{1}{8}%
\{P_{2},Q_{1}\}+s_{1}u_{2}-s_{2}u_{1},$

$\ \ \ \ \ \ \ \ \ \ \ \ \ s=\frac{1}{4}\{P_{1},P_{2}%
\}+2r_{1}s_{2}-2r_{2}s_{1},$

$\ \ \ \ \ \ \ \ \ \ \ \ \ u=-\frac{1}{4}\{Q_{1},Q_{2}%
\}-2r_{1}u_{2}+2r_{2}u_{1},$

\noindent
then \ge$_{8}^{\C}$ is a $\C$-Lie algebra.

\bigskip

\emph{Proposition 5.2.} (I.Yokota\cite[\emph{Theorem 5.2.1.}]{Yokota1}) 

The $\C$-Lie algebra \ge$_{8}^{\C}$ is simple.

\bigskip

\emph{Definition 5.3.} \ The group $E_{8}^{\C}$ is defined to be the
automorphism group of the Lie algebra \ge$_{8}^{\C}$:

\ \ \ \ \ \ \ \ \ \ \ \ \ \ \ \ \ \ \ \ \ \ \ \ $E_{8}^{\C}=\{\alpha \in
Iso_{\C}($\ge$_{8}^{\C}) \mid \alpha \lbrack R1,R2]=[\alpha R1,%
\alpha R2]\}.$

\bigskip

\emph{Proposition 5.4.} (I.Yokota\cite[\emph{Theorem 5.4.2.}]{Yokota1}) 

$E_{8}^{\C}$ is a simply connected complex Lie group of type $E_{8}$.

\bigskip

\emph{Proposition 5.5.} (T.Imai and I.Yokota\cite[\emph{Proposition 9.}]{YokotaImai5}) 

\ge$_{8,1}$ is a simple Lie algebra of type $E_{8(-24)}$,

\noindent
where \ge$_{8,1}=$\ \ge$_{7,1}\oplus$\gP$\oplus$\gP$\oplus \R \oplus \R \oplus \R$.

\bigskip

\emph{Definition 5.6. } We define an $\R$-vector space \gR$_{8}$\textbf{\ }by

\ \ \ \ \ \ \ \ \ \ \ \ \gR$_{8}=$ \gR$_{7}\oplus $\gP$\oplus $\gP$\oplus \R\oplus \R\oplus \R$

\bigskip

For $R=(\Phi ,P,Q,r,s,u)\in $\gR$_{8},(R)_{\Phi
},(R)_{P}$,$(R)_{Q},(R)_{r},(R)_{s}$ and $(R)_{u}$ means \gR$_{7}$
element of $R,$
$P$ element of $R$ $,$ $Q$ element of $R,r$ element of $R,s$ element of $R$
and $u$ element of $R$ respectively.

\bigskip

\emph{Definition 5.7. } For $R_{1}=(\Phi
_{1},P_{1},Q_{1},r_{1},s_{1},u_{1}),R_{2}=(\Phi
_{2},P_{2},Q_{2},r_{2},s_{2},u_{2})$

$\in $\gR$_{8},$
$\Phi _{1}=(D_{1},M_{1},T_{1},A_{1},B_{1},\rho _{1}),\Phi
_{2}=(D_{2},M_{2},T_{2},A_{2},B_{2},\rho _{2})\in $\gR$_{7},$

$\phi _{1}=(D_{1},M_{1},T_{1}),\phi _{2}=(D_{2},M_{2},T_{2})\in $\gR$%
_{6},$
$P_{1}=(X_{1},Y_{1},\xi _{1},\eta _{1}),$

$Q_{1}=(Z_{1},W_{1},\zeta
_{1},\omega _{1}),P_{2}=(X_{2},Y_{2},\xi _{2},\eta
_{2}),Q_{2}=(Z_{2},W_{2},\zeta _{2},\omega _{2})\in $\gP,

\noindent
We define a bracket operation $[R_{1},R_{2}]_{8}$ by

\ \ \ \ $[R_{1},R_{2}]_{8}$

\ \ \ \ $=([R_{1},R_{2}]_{8\Phi
},[R_{1},R_{2}]_{8P},[R_{1},R_{2}]_{8Q},[R_{1},R_{2}]_{8r},[R_{1},R_{2}]_{8s},[R_{1},R_{2}]_{8u}), 
$

\noindent
where

$[R_{1},R_{2}]_{8\Phi }=([\Phi _{1},\Phi _{2}]_{7D}+(P_{1}\times
Q_{2})_{D}-(P_{2}\times Q_{1})_{D},$

$\ \ \ \ \ \ \ \ \ \ \ \ \ \ \ \ \ \ \ [\Phi _{1},\Phi
_{2}]_{7M}+(P_{1}\times Q_{2})_{M}-(P_{2}\times Q_{1})_{M},$

$\ \ \ \ \ \ \ \ \ \ \ \ \ \ \ \ \ \ \ [\Phi _{1},\Phi
_{2}]_{7T}+(P_{1}\times Q_{2})_{T}-(P_{2}\times Q_{1})_{T},$

$\ \ \ \ \ \ \ \ \ \ \ \ \ \ \ \ \ \ \ [\Phi _{1},\Phi
_{2}]_{7A}+(P_{1}\times Q_{2})_{A}-(P_{2}\times Q_{1})_{A},$

\ \ \ \ \ \ \ \ \ \ \ \ \ \ \ \ \ \ \ $[\Phi _{1},\Phi
_{2}]_{7B}+(P_{1}\times Q_{2})_{B}-(P_{2}\times Q_{1})_{B},$

$\ \ \ \ \ \ \ \ \ \ \ \ \ \ \ \ \ \ \ [\Phi _{1},\Phi _{2}]_{7\rho
}+(P_{1}\times Q_{2})_{\rho }-(P_{2}\times Q_{1})_{\rho }),$

$[R_{1},R_{2}]_{8P}=\left( 
\begin{array}{c}
(\delta _{d}\textrm{g}_{d}(D_{1})+\widehat{M}_{1}+\widetilde{T}_{1})X_{2}\text{-}%
\frac{1}{3}\rho _{1}X_{2}+2B_{1}\times Y_{2}+\eta _{2}A_{1} \\ 
2A_{1}\times X_{2}+(\delta _{d}\textrm{g}_{d}(D_{1})+\widehat{M}_{1}-\widetilde{T}%
_{1})Y_{2}\text{+}\frac{1}{3}\rho _{1}Y_{2}+\xi _{2}B_{1} \\ 
(A_{1},Y_{2})+\rho _{1}\xi _{2} \\ 
(B_{1},X_{2})-\rho _{1}\eta _{2}%
\end{array}%
\right) $

$\ \ \ \ \ \ \ \ \ \ \ \ \ \ \ -\left( 
\begin{array}{c}
(\delta _{d}\textrm{g}_{d}(D_{2})+\widehat{M}_{2}+\widetilde{T}_{2})X_{1}\text{-}%
\frac{1}{3}\rho _{2}X_{1}+2B_{2}\times Y_{1}+\eta _{1}A_{2} \\ 
2A_{2}\times X_{1}+(\delta _{d}\textrm{g}_{d}(D_{2})+\widehat{M}_{2}-\widetilde{T}%
_{2})Y_{1}\text{+}\frac{1}{3}\rho _{2}Y_{1}+\xi _{1}B_{2} \\ 
(A_{2},Y_{1})+\rho _{2}\xi _{1} \\ 
(B_{2},X_{1})-\rho _{2}\eta _{1}%
\end{array}%
\right) $

$\ \ \ \ \ \ \ \ \ \ \ \ \ \ \ \ \
+r_{1}P_{2}-r_{2}P_{1}+s_{1}Q_{2}-s_{2}Q_{1},$

$[R_{1},R_{2}]_{8Q}=\left( 
\begin{array}{c}
(\delta _{d}\textrm{g}_{d}(D_{1})+\widehat{M}_{1}+\widetilde{T}_{1})Z_{2}\text{-}%
\frac{1}{3}\rho _{1}Z_{2}+2B_{1}\times W_{2}+\omega _{2}A_{1} \\ 
2A_{1}\times Z_{2}+(\delta _{d}\textrm{g}_{d}(D_{1})+\widehat{M}_{1}-\widetilde{T}%
_{1})W_{2}\text{+}\frac{1}{3}\rho _{1}W_{2}+\zeta _{2}B_{1} \\ 
(A_{1},W_{2})+\rho _{1}\zeta _{2} \\ 
(B_{1},Z_{2})-\rho _{1}\omega _{2}%
\end{array}%
\right) $

$\ \ \ \ \ \ \ \ \ \ \ \ \ \ \ -\left( 
\begin{array}{c}
(\delta _{d}\textrm{g}_{d}(D_{2})+\widehat{M}_{2}+\widetilde{T}_{2})Z_{1}\text{-}%
\frac{1}{3}\rho _{2}Z_{1}+2B_{2}\times W_{1}+\omega _{1}A_{2} \\ 
2A_{2}\times Z_{1}\text{+}(\delta _{d}\textrm{g}_{d}(D_{2})+\widehat{M}_{2}-%
\widetilde{T}_{2})W_{1}\text{+}\frac{1}{3}\rho _{2}W_{1}+\zeta _{1}B_{2}
\\ 
(A_{2},W_{1})+\rho _{2}\zeta _{1} \\ 
(B_{2},Z_{1})-\rho _{2}\omega _{1}%
\end{array}%
\right) $

$\ \ \ \ \ \ \ \ \ \ \ \ \ \ \ \ \ \ \
-r_{1}Q_{2}+r_{2}Q_{1}+u_{1}P_{2}-u_{2}P_{1},$

$[R_{1},R_{2}]_{8r}=-\frac{1}{8}\{P_{1},Q_{2}\}+\frac{1}{8}%
\{P_{2},Q_{1}\}+s_{1}u_{2}-s_{2}u_{1},$

$[R_{1},R_{2}]_{8s}=\frac{1}{4}\{P_{1},P_{2}\}+2r_{1}s_{2}-2r_{2}s_{1},$

$[R_{1},R_{2}]_{8u}=-\frac{1}{4}\{Q_{1},Q_{2}\}-2r_{1}u_{2}+2r_{2}u_{1}$.

\bigskip

\emph{Lemma 5.8.}\ \ \gR$_{8}$ is isomorphic to 
\ge$_{8,1}$ under the correspondence

$\ \textrm{f}_{y}:R=(\Phi ,P,Q,r,s,u)\in $\gR$_{8}\rightarrow \
\textrm{f}_{y}(R)=(\textrm{f}_{y}(\Phi ),P,Q,r,s,u)\in $\ge$_{8,1}$.

\bigskip

\emph{Proof.} For $R_{1}=(\Phi
_{1},P_{1},Q_{1},r_{1},s_{1},u_{1}),R_{2}=(\Phi
_{2},P_{2},Q_{2},r_{2},s_{2},u_{2})\in $\gR$_{8},$

$\Phi _{1}=(D_{1},M_{1},T_{1},A_{1},B_{1},\rho _{1}),\Phi
_{2}=(D_{2},M_{2},T_{2},A_{2},B_{2},\rho _{2})\in $\gR$_{7},$

$\phi _{1}=(D_{1},M_{1},T_{1}),\phi _{2}=(D_{2},M_{2},T_{2})\in $\gR$%
_{6},$

$P_{1}=(X_{1},Y_{1},\xi _{1},\eta _{1}),Q_{1}=(Z_{1},W_{1},\zeta
_{1},\omega _{1}),P_{2}=(X_{2},Y_{2},\xi _{2},\eta
_{2}),$

$Q_{2}=(Z_{2},W_{2},\zeta _{2},\omega _{2})\in $\gP,

\noindent
We have

\noindent
\ \ \  $\Phi (\delta _{d}\textrm{g}_{d}(D_{1})+\widehat{M}_{1}+\widetilde{T}%
_{1},A_{1},B_{1},\rho _{1}),\Phi (\delta _{d}\textrm{g}_{d}(D_{2})+\widehat{M}%
_{2}+\widetilde{T}_{2},A_{2},B_{2},\rho _{2})\in $ \ge $_{7,1}$ .

\noindent
Let we put these as $\Phi _{1}^{7},\Phi _{2}^{7}$ respectively, we have
R$_{1}^{8}$=$(\Phi _{1}^{7},P_{1},Q_{1},r_{1},s_{1},u_{1}),$

$R_{2}^{8}=(\Phi _{2}^{7},P_{2},Q_{2},r_{2},s_{2},u_{2})\in $\ge$%
_{8,1}.$

\noindent
Hence by \emph{Proposition 5.1 }we have

$[R_{1}^{8},R_{2}^{8}]=(\Phi ^{7},P,Q,r,s,u)$,

\noindent
where

$\Phi ^{7}=[\Phi (\delta _{d}\textrm{g}_{d}(D_{1})+%
\widehat{M}_{1}+\widetilde{T}_{1},A_{1},B_{1},\rho _{1}),\Phi (\delta
_{d}\textrm{g}_{d}(D_{2})+\widehat{M}_{2}+\widetilde{T}_{2},A_{2},B_{2},\rho _{2})]$

$\ \ \ \ \ \ +P_{1}\times Q_{2}-P_{2}\times Q_{1},$

$ P=\Phi (\delta _{d}\textrm{g}_{d}(D_{1})+\widehat{M}_{1}+%
\widetilde{T}_{1},A_{1},B_{1},\rho _{1})P_{2}$

\ \ \ \ $-\Phi (\delta
_{d}\textrm{g}_{d}(D_{2})+\widehat{M}_{2}+\widetilde{T}_{2},A_{2},B_{2},\rho
_{2})P_{1}$
$+r_{1}P_{2}-r_{2}P_{1}+s_{1}Q_{2}-s_{2}Q_{1},$

$ Q=\Phi (\delta _{d}\textrm{g}_{d}(D_{1})+\widehat{M}_{1}+%
\widetilde{T}_{1},A_{1},B_{1},\rho _{1})Q_{2}$

\ \ \ \ $-\Phi (\delta
_{d}\textrm{g}_{d}(D_{2})+\widehat{M}_{2}+\widetilde{T}_{2},A_{2},B_{2},\rho
_{2})Q_{1}$
$-r_{1}Q_{2}+r_{2}Q_{1}+u_{1}P_{2}-u_{2}P_{1},$

$ r=-\frac{1}{8}\{P_{1},Q_{2}\}+\frac{1}{8}%
\{P_{2},Q_{1}\}+s_{1}u_{2}-s_{2}u_{1},$

$ s=\frac{1}{4}\{P_{1},P_{2}%
\}+2r_{1}s_{2}-2r_{2}s_{1},$

$ u=-\frac{1}{4}\{Q_{1},Q_{2}%
\}-2r_{1}u_{2}+2r_{2}u_{1}$.

\noindent
By \emph{Definition 4.3}. and \emph{Proposition 5.1} , we have 
$P=[R_{1},R_{2}]_{8P},$

$Q=[R_{1},R_{2}]_{8Q},r=[R_{1},R_{2}]_{8r},$
$s=[R_{1},R_{2}]_{8s},u=[R_{1},R_{2}]_{8u}$ . 

\noindent
On the other hand, by \emph{%
Lemma 4.12 }$\Phi ^{7}$ is isomorphic to $[R_{1},R_{2}]_{8\Phi }$ . 

\noindent
So\ $[R_{1},R_{2}]_{8}=([R_{1},R_{2}]_{8\Phi
},[R_{1},R_{2}]_{8P},[R_{1},R_{2}]_{8Q},[R_{1},R_{2}]_{8r},[R_{1},R_{2}]_{8s},$

\noindent
$[R_{1},R_{2}]_{8u}) $ is
a Lie bracket of \gR$_{8}$\textbf{\ }and also \gR$%
_{8}$\textbf{\ }is\textbf{\ }isomorphic to \ge$_{8,1}$%
. \ \ \ \ \emph{Q.E.D.}

\bigskip

\emph{Corollary 5.8.1. \ }The definition of the Lie bracket $%
[R_{1},R_{2}]_{8}$ in \emph{Definition 5.7} can be
extended to $T_{i}\in $\gJ $(i=1,2).$

\bigskip

\emph{Proof.} \ Evidently by \emph{Corollary 4.12.1} \ \ \ \ \emph{%
Q.E.D.}

\bigskip

\emph{Lemma 5.9.} $\ [R_{1},R_{2}]_{8\Phi
},[R_{1},R_{2}]_{8P},[R_{1},R_{2}]_{8Q},[R_{1},R_{2}]_{8r},[R_{1},R_{2}]_{8s},$

\noindent
$[R_{1},R_{2}]_{8u}$
are expressed by the following expression respectively.

For $A(\alpha ,a),B(\beta ,b),X(\chi ,x),Y(\gamma ,y),Z(\mu ,z),W(\psi
,w)\in $\gJ ,

$P(X(\chi ,x),Y(\gamma ,y),\xi ,\eta ),Q(Z(\mu ,z),W(\psi ,w),\zeta ,\omega
)\in $\gP,

$[R_{1},R_{2}]_{8\Phi }=(([R_{1},R_{2}]_{8\Phi })_{D},([R_{1},R_{2}]_{8%
\Phi })_{M},([R_{1},R_{2}]_{8\Phi })_{T},([R_{1},R_{2}]_{8\Phi
})_{A},$

\ \ \ \ \ \ \ \ \ \ \ \ \ \ \ \ \ \ \ 
$([R_{1},R_{2}]_{8\Phi })_{B},$
$([R_{1},R_{2}]_{8\Phi })_{\rho}),$

($[R_{1},R_{2}]_{8\Phi })_{D}$

$\ =[D1,D2]-\textrm{JD}(m_{11},m_{21})-\textrm{d}_{g}\nu ^{2}\textrm{g}_{d}(\textrm{JD}(m_{12},m_{22}))-\textrm{d}_{g}\nu
\textrm{g}_{d}(\textrm{JD}(m_{13},m_{23}))$

$\ +\textrm{JD}(t_{11},t_{21})+\textrm{d}_{g}\nu ^{2}\textrm{g}_{d}(\textrm{JD}(t_{12},t_{22}))+\textrm{d}_{g}\nu
\textrm{g}_{d}(\textrm{JD}(t_{13},t_{23}))$

$\ +2(\textrm{JD}(a_{11},b_{21})+\textrm{d}_{g}\nu ^{2}\textrm{g}_{d}(\textrm{JD}(a_{12},b_{22}))+\textrm{d}_{g}\nu
\textrm{g}_{d}(\textrm{JD}(a_{13},b_{23})))$

$\ -2(\textrm{JD}(a_{21},b_{11})+\textrm{d}_{g}\nu ^{2}\textrm{g}_{d}(\textrm{JD}(a_{22},b_{12}))+\textrm{d}_{g}\nu
\textrm{g}_{d}(\textrm{JD}(a_{23},b_{13})))$

$\ -\frac{1}{2}(\textrm{JD}(x_{11},w_{21})+\textrm{d}_{g}\nu
^{2}\textrm{g}_{d}(\textrm{JD}(x_{12},w_{22}))+\textrm{d}_{g}\nu \textrm{g}_{d}(\textrm{JD}(x_{13},w_{23}))$

$\ +\textrm{JD}(z_{21},y_{11})+\textrm{d}_{g}\nu ^{2}\textrm{g}_{d}(\textrm{JD}(z_{22},y_{12}))+\textrm{d}_{g}\nu
\textrm{g}_{d}(\textrm{JD}(z_{23},y_{13}))),$

$\ +\frac{1}{2}(\textrm{JD}(x_{21},w_{11})+\textrm{d}_{g}\nu
^{2}\textrm{g}_{d}(\textrm{JD}(x_{22},w_{12}))+\textrm{d}_{g}\nu \textrm{g}_{d}(\textrm{JD}(x_{23},w_{13}))$

$\ +\textrm{JD}(z_{11},y_{21})+\textrm{d}_{g}\nu ^{2}\textrm{g}_{d}(\textrm{JD}(z_{12},y_{22}))+\textrm{d}_{g}\nu
\textrm{g}_{d}(\textrm{JD}(z_{13},y_{23}))),$

\bigskip

$([R_{1},R_{2}]_{8\Phi })_{M}$

$\ =A_{1}(\textrm{g}_{d}(D_{1})m_{21}-\textrm{g}_{d}(D_{2})m_{11})+A_{2}(\nu
\textrm{g}_{d}(D_{1})m_{22}-\nu \textrm{g}_{d}(D_{2})m_{12})$

$\ +A_{3}(\nu ^{2}\textrm{g}_{d}(D_{1})m_{23}-\nu ^{2}\textrm{g}_{d}(D_{2})m_{13})$

$\ +\frac{1}{2}A_{1}(-\overline{m_{12}m_{23}}+\overline{m_{22}m_{13}})+\frac{%
1}{2}A_{2}(-\overline{m_{13}m_{21}}+\overline{m_{23}m_{11}})$

$\ +\frac{1}{2}A_{3}(-\overline{m_{11}m_{22}}+\overline{m_{21}m_{12}})$

$\ +\frac{1}{2}A_{1}((\tau _{11}+2\tau _{12})t_{21}-(\tau _{21}+2\tau
_{22})t_{11}-\overline{t_{12}t_{23}}+\overline{t_{22}t_{13}})$

$\ +\frac{1}{2}A_{2}((-2\tau _{11}-\tau _{12})t_{22}-(-2\tau _{21}-\tau
_{22})t_{12}-\overline{t_{23}t_{11}}+\overline{t_{13}t_{21}})$

$\ +\frac{1}{2}A_{3}((\tau _{11}-\tau _{12})t_{23}-(\tau _{21}-\tau
_{22})t_{13}-\overline{t_{11}t_{22}}+\overline{t_{21}t_{12}})$

$\ +A_{1}((\alpha _{12}-\alpha _{13})b_{21}-(\beta _{22}-\beta _{23})a_{11}-%
\overline{a_{12}b_{23}}+\overline{b_{22}a_{13}}$

$\ \ \ \ \ -(\alpha _{22}-\alpha _{23})b_{11}+(\beta _{12}-\beta
_{13})a_{21}+\overline{a_{22}b_{13}}-\overline{b_{12}a_{23}})$

$\ +A_{2}((\alpha _{13}-\alpha _{11})b_{22}-(\beta _{23}-\beta _{21})a_{12}-%
\overline{a_{13}b_{21}}+\overline{b_{23}a_{11}})$

$\ \ \ \ \ -(\alpha _{23}-\alpha _{21})b_{12}+(\beta _{13}-\beta
_{11})a_{22}+\overline{a_{23}b_{11}}-\overline{b_{13}a_{21}})$

$\ +A_{3}((\alpha _{11}-\alpha _{12})b_{23}-(\beta _{21}-\beta _{22})a_{13}-%
\overline{a_{11}b_{22}}+\overline{b_{21}a_{12}}$

$\ \ \ \ \ -(\alpha _{21}-\alpha _{22})b_{13}+(\beta _{11}-\beta
_{12})a_{23}+\overline{a_{21}b_{12}}-\overline{b_{11}a_{22}})$

$\ -\frac{1}{4}A_{1}((\chi _{12}-\chi _{13})w_{21}-(\psi _{22}-\psi
_{23})x_{11}-\overline{x_{12}w_{23}}+\overline{w_{22}x_{13}}$

$\ \ \ \ \ \ \ \ \ +(\mu _{22}-\mu _{23})y_{11}-(\gamma _{12}-\gamma
_{13})z_{21}-\overline{z_{22}y_{13}}+\overline{y_{12}z_{23}})$

\ $-\frac{1}{4}A_{2}((\chi _{13}-\chi _{11})w_{22}-(\psi _{23}-\psi
_{21})x_{12}-\overline{x_{13}w_{21}}+\overline{w_{23}x_{11}}$

$\ \ \ \ \ \ \ +(\mu _{23}-\mu _{21})y_{12}-(\gamma _{13}-\gamma
_{11})z_{22}-\overline{z_{23}y_{11}}+\overline{y_{13}z_{21}})$

$\ -\frac{1}{4}A_{3}((\chi _{11}-\chi _{12})w_{23}-(\psi _{21}-\psi
_{22})x_{13}-\overline{x_{11}w_{22}}+\overline{w_{21}x_{12}}$

$\ \ \ \ \ \ \ \ +(\mu _{21}-\mu _{22})y_{13}-(\gamma _{11}-\gamma
_{12})z_{23}-\overline{z_{21}y_{12}}+\overline{y_{11}z_{22}})$

$\ +\frac{1}{4}A_{1}((\chi _{22}-\chi _{23})w_{11}-(\psi _{12}-\psi
_{13})x_{21}-\overline{x_{22}w_{13}}+\overline{w_{12}x_{23}}$

$\ \ \ \ \ \ \ \ +(\mu _{12}-\mu _{13})y_{21}-\gamma _{22}-\gamma
_{23})z_{11}-\overline{z_{12}y_{23}}+\overline{y_{22}z_{13}})$

$\ +\frac{1}{4}A_{2}((\chi _{23}-\chi _{21})w_{12}-(\psi _{13}-\psi
_{11})x_{22}-\overline{x_{23}w_{11}}+\overline{w_{13}x_{21}}$

$\ \ \ \ \ \ \ \ +(\mu _{13}-\mu _{11})y_{22}-(\gamma _{23}-\gamma
_{21})z_{12}-\overline{z_{13}y_{21}}+\overline{y_{23}z_{11}})$

$\ +\frac{1}{4}A_{3}((\chi _{21}-\chi _{22})w_{13}-(\psi _{11}-\psi
_{12})x_{23}-\overline{x_{21}w_{12}}+\overline{w_{11}x_{22}}$

$\ \ \ \ \ \ \ \ +(\mu _{11}-\mu _{12})y_{23}-(\gamma _{21}-\gamma
_{22})z_{13}-\overline{z_{11}y_{22}}+\overline{y_{21}z_{12}}),$

\bigskip

$([R_{1},R_{2}]_{8\Phi })_{T}$

=$F_{1}(\textrm{g}_{d}(D_{1})t_{21})+F_{2}(\nu \textrm{g}_{d}(D_{1})t_{22})+F_{3}(\nu
^{2}\textrm{g}_{d}(D_{1})t_{23})$

$-(F_{1}(\textrm{g}_{d}(D_{2})t_{11})-F_{2}(\nu \textrm{g}_{d}(D_{2})t_{12})-F_{3}(\nu
^{2}\textrm{g}_{d}(D_{2})t_{13})$

$\ +(-(m_{12},t_{22})+(m_{13},t_{23})+(m_{22},t_{12})-(m_{23},t_{13}))E_{1}$

$\ +(-(m_{13},t_{23})+(m_{11},t_{21})+(m_{23},t_{13})-(m_{21},t_{11}))E_{2}$

$\ +(-(m_{11},t_{21})+(m_{12},t_{22})+(m_{21},t_{11})-(m_{22},t_{12}))E_{3}$

$+\frac{1}{3}(4\alpha _{11}\beta _{21}-2\alpha _{12}\beta _{22}-2\alpha
_{13}\beta _{23}-4(a_{11},b_{21})+2(a_{12},b_{22})+2(a_{13},b_{23}))E_{1}$

$-\frac{1}{3}(4\alpha _{21}\beta _{11}-2\alpha _{22}\beta _{12}-2\alpha
_{23}\beta _{13}-4(a_{21},b_{11})+2(a_{22},b_{12})+2(a_{23},b_{13}))E_{1}$

$+\frac{1}{3}(-2\alpha _{11}\beta _{21}+4\alpha _{12}\beta _{22}-2\alpha
_{13}\beta _{23}+2(a_{11},b_{21})-4(a_{12},b_{22})+2(a_{13},b_{23}))E_{2}$

$-\frac{1}{3}(-2\alpha _{21}\beta _{11}+4\alpha _{22}\beta _{21}-2\alpha
_{23}\beta _{13}+2(a_{21},b_{11})-4(a_{22},b_{12})+2(a_{23},b_{13}))E_{2}$

$+\frac{1}{3}(-2\alpha _{11}\beta _{21}-2\alpha _{12}\beta _{22}+4\alpha
_{13}\beta _{23}+2(a_{11},b_{21})+2(a_{12},b_{22})-4(a_{13},b_{23}))E_{3}$

$-\frac{1}{3}(-2\alpha _{21}\beta _{11}-2\alpha _{22}\beta _{12}+4\alpha
_{23}\beta _{13}+2(a_{21},b_{11})+2(a_{22},b_{12})-4(a_{23},b_{13}))E_{3}$

$\ +F_{1}(\frac{1}{2}(\overline{m_{12}t_{23}}-\overline{t_{22}m_{13}}+(-\tau
_{21}-2\tau _{22})m_{11}$

\ \ \ \ \ \ \ \ \ $-\overline{m_{22}t_{13}}+\overline{t_{12}m_{23}}%
-(-\tau _{11}-2\tau _{12})m_{21}))$

$\ +F_{2}(\frac{1}{2}(\overline{m_{13}t_{21}}-\overline{t_{23}m_{11}}+(2\tau
_{21}+\tau _{22})m_{12}$

\ \ \ \ \ \ \ \ \ $-\overline{m_{23}t_{11}}+\overline{t_{13}m_{21}}%
-(2\tau _{11}+\tau _{12})m_{22}))$

$\ +F_{3}(\frac{1}{2}(\overline{m_{11}t_{22}}-\overline{t_{21}m_{12}}+(-\tau
_{21}+\tau _{22})m_{13}$

\ \ \ \ \ \ \ \ \ $-\overline{m_{21}t_{12}}+\overline{t_{11}m_{22}}%
-(-\tau _{11}+\tau _{12})m_{23}))$

$\ +F_{1}((\alpha _{12}+\alpha _{13})b_{21}+(\beta _{22}+\beta _{23})a_{11}+%
\overline{b_{22}a_{13}}+\overline{a_{12}b_{23}}$

$\ \ \ \ \ -(\alpha _{22}+\alpha _{23})b_{11}-(\beta _{12}+\beta
_{13})a_{21}-\overline{b_{12}a_{23}}-\overline{a_{22}b_{13}})$

$\ +F_{2}((\alpha _{13}+\alpha _{11})b_{22}+(\beta _{23}+\beta _{21})a_{12}+%
\overline{b_{23}a_{11}}+\overline{a_{13}b_{21}}$

$\ \ \ \ \ -(\alpha _{23}+\alpha _{21})b_{12}-(\beta _{13}+\beta
_{11})a_{22}-\overline{b_{13}a_{21}}-\overline{a_{23}b_{11}})$

$\ +F_{3}((\alpha _{11}+\alpha _{12})b_{23}+(\beta _{21}+\beta _{22})a_{13}+%
\overline{b_{21}a_{12}}+\overline{a_{11}b_{22}}$

$\ \ \ \ \ -(\alpha _{21}+\alpha _{22})b_{13}-(\beta _{11}+\beta
_{12})a_{23}-\overline{b_{11}a_{22}}-\overline{a_{21}b_{12}})$

$-\frac{1}{6}(\chi _{11}\psi _{21}-\chi _{12}\psi
_{22}-(x_{11},w_{21})+(x_{12},w_{22})$

$\ \ \ \ \ +\mu _{21}\gamma _{11}-\mu _{22}\gamma
_{12}-(z_{21},y_{11})+(z_{22},y_{12}))(E_{1}-E_{2})$

$-\frac{1}{6}(\chi _{12}\psi _{22}-\chi _{13}\psi
_{23}-(x_{12},w_{22})+(x_{13},w_{23})$

$\ \ \ \ +\mu _{22}\gamma _{12}-\mu _{23}\gamma
_{13}-(z_{22},y_{12})+(z_{23},y_{13}))(E_{2}-E_{3})$

$-\frac{1}{6}(\chi _{13}\psi _{23}-\chi _{11}\psi
_{21}-(x_{13},w_{23})+(x_{11},w_{21})$

$\ \ \ +\mu _{23}\gamma _{13}-\mu _{21}\gamma
_{11}-(z_{23},y_{13})+(z_{21},y_{11}))(E_{3}-E_{1})$

$-\frac{1}{4}F_{1}((\chi _{12}+\chi _{13})w_{21}+(\psi _{22}+\psi
_{23})x_{11}+\overline{w_{22}x_{13}}+\overline{x_{12}w_{23}}$

$\ \ \ \ \ \ +(\mu _{22}+\mu _{23})y_{11}+(\gamma _{12}+\gamma _{13})z_{21}+%
\overline{y_{12}z_{23}}+\overline{z_{22}y_{13}})$

$-\frac{1}{4}F_{2}((\chi _{13}+\chi _{11})w_{22}+(\psi _{23}+\psi
_{21})x_{12}+\overline{w_{23}x_{11}}+\overline{x_{13}w_{21}}$

$\ \ \ \ \ \ +(\mu _{23}+\mu _{21})y_{12}+(\gamma _{13}+\gamma _{11})z_{22}+%
\overline{y_{13}z_{21}}+\overline{z_{23}y_{11}})$

$-\frac{1}{4}F_{3}((\chi _{11}+\chi _{12})w_{23}+(\psi _{21}+\psi
_{22})x_{13}+\overline{w_{21}x_{12}}+\overline{x_{11}w_{22}}$

$\ \ \ \ \ \ +(\mu _{21}+\mu _{22})y_{13}+(\gamma _{11}+\gamma _{12})z_{23}+%
\overline{y_{11}z_{22}}+\overline{z_{21}y_{12}})$

$+\frac{1}{6}(\chi _{21}\psi _{11}-\chi _{22}\psi
_{12}-(x_{21},w_{11})+(x_{22},w_{12})$

$\ \ \ +\mu _{11}\gamma _{21}-\mu _{12}\gamma
_{22}-(z_{11},y_{21})+(z_{12},y_{22}))(E_{1}-E_{2})$

$+\frac{1}{6}(\chi _{22}\psi _{12}-\chi _{23}\psi
_{13}-(x_{22},w_{12})+(x_{23},w_{13})$

$\ \ \ +\mu _{12}\gamma _{22}-\mu _{23}\gamma
_{23}-(z_{12},y_{22})+(z_{13},y_{23}))(E_{2}-E_{3})$

$+\frac{1}{6}(\chi _{23}\psi _{13}-\chi _{21}\psi
_{11}-(x_{23},w_{13})+(x_{21},w_{11})$

$\ \ +\mu _{13}\gamma _{23}-\mu _{11}\gamma
_{21}-(z_{13},y_{23})+(z_{11},y_{21}))(E_{3}-E_{1})$

$+\frac{1}{4}F_{1}((\chi _{22}+\chi _{23})w_{11}+(\psi _{12}+\psi
_{13})x_{21}+\overline{w_{12}x_{23}}+\overline{x_{22}w_{13}}$

$\ \ \ \ \ \ +(\mu _{12}+\mu _{13})y_{21}+(\gamma _{22}+\gamma _{23})z_{11}+%
\overline{y_{22}z_{13}}+\overline{z_{12}y_{23}})$

$+\frac{1}{4}F_{2}((\chi _{23}+\chi _{21})w_{12}+(\psi _{13}+\psi
_{11})x_{22}+\overline{w_{13}x_{21}}+\overline{x_{23}w_{11}}$

$\ \ \ \ \ \ +(\mu _{13}+\mu _{11})y_{22}+(\gamma _{23}+\gamma _{21})z_{12}+%
\overline{y_{23}z_{11}}+\overline{z_{13}y_{21}})$

$+\frac{1}{4}F_{3}((\chi _{21}+\chi _{22})w_{13}+(\psi _{11}+\psi
_{12})x_{23}+\overline{w_{11}x_{22}}+\overline{x_{21}w_{12}}$

$\ \ \ \ \ +(\mu _{11}+\mu _{12})y_{23}+(\gamma _{21}+\gamma _{22})z_{13}+%
\overline{y_{21}z_{12}}+\overline{z_{11}y_{22}}),$

\bigskip

$([R_{1},R_{2}]_{8\Phi })_{A}$

=$F_{1}(\textrm{g}_{d}(D_{1})a_{21})+F_{2}(\nu \textrm{g}_{d}(D_{1})a_{22})+F_{3}(\nu
^{2}\textrm{g}_{d}(D_{1})a_{23})$

$-F_{1}(\textrm{g}_{d}(D_{2})a_{11})-F_{2}(\nu \textrm{g}_{d}(D_{2})a_{12})-F_{3}(\nu
^{2}\textrm{g}_{d}(D_{2})a_{13})$

$+((m_{13},a_{23})-(m_{12},a_{22})-(m_{23},a_{13})+(m_{22},a_{12}))E_{1}$

$+((m_{11},a_{21})-(m_{13},a_{23})-(m_{21},a_{11})+(m_{23},a_{13}))E_{2}$

$+((m_{12},a_{22})-(m_{11},a_{21})-(m_{22},a_{12})+(m_{21},a_{11}))E_{3}$

$+\frac{1}{2}F_{1}((\alpha _{23}-\alpha _{22})m_{11}-\overline{a_{22}m_{13}}+%
\overline{m_{12}a_{23}}-(\alpha _{13}-\alpha _{12})m_{21}+\overline{%
a_{12}m_{23}}-\overline{m_{22}a_{13}})$

$+\frac{1}{2}F_{2}((\alpha _{21}-\alpha _{23})m_{12}-\overline{a_{23}m_{11}}+%
\overline{m_{13}a_{21}}-(\alpha _{11}-\alpha _{13})m_{22}+\overline{%
a_{13}m_{21}}-\overline{m_{23}a_{11}})$

$+\frac{1}{2}F_{3}((\alpha _{22}-\alpha _{21})m_{13}-\overline{a_{21}m_{12}}+%
\overline{m_{11}a_{22}}-(\alpha _{12}-\alpha _{11})m_{23}+\overline{%
a_{11}m_{22}}-\overline{m_{21}a_{12}})$

$+(\tau _{11}\alpha _{21}+(t_{12},a_{22})+(t_{13},a_{23})-\tau _{21}\alpha
_{11}-(t_{22},a_{12})-(t_{23},a_{13})+\frac{2}{3}\rho _{1}\alpha _{21}$

$-\frac{2}{3}\rho _{2}\alpha _{11})E_{1}$

$+(\tau _{12}\alpha _{22}+(t_{13},a_{23})+(t_{11},a_{21})-\tau _{22}\alpha
_{12}-(t_{23},a_{13})-(t_{21},a_{11})+\frac{2}{3}\rho _{1}\alpha _{22}$

$-\frac{2}{3}\rho _{2}\alpha _{12})E_{2}$

$+((-\tau _{11}-\tau _{12})\alpha
_{23}+(t_{11},a_{21})+(t_{12},a_{22})-(-\tau _{21}-\tau _{22})\alpha
_{13}-(t_{21},a_{11})$

$-(t_{22},a_{12})$
$+\frac{2}{3}\rho _{1}\alpha _{23}-\frac{2}{3}\rho _{2}\alpha _{13})E_{3}$

$+F_{1}(\frac{1}{2}(-\tau _{11}a_{21}+(\alpha _{22}+\alpha _{23})t_{11}+%
\overline{a_{22}t_{13}}+\overline{t_{12}a_{23}})+\frac{2}{3}\rho _{1}a_{21}$

$\ \ -\frac{1}{2}(-\tau _{21}a_{11}+(\alpha _{12}+\alpha _{13})t_{21}+%
\overline{a_{12}t_{23}}+\overline{t_{22}a_{13}})-\frac{2}{3}\rho _{2}a_{11})$

$+F_{2}(\frac{1}{2}(-\tau _{12}a_{22}+(\alpha _{23}+\alpha _{21})t_{12}+%
\overline{a_{23}t_{11}}+\overline{t_{13}a_{21}})+\frac{2}{3}\rho _{1}a_{22}$

$\ \ -\frac{1}{2}(-\tau _{22}a_{12}+(\alpha _{13}+\alpha _{11})t_{22}+%
\overline{a_{13}t_{21}}+\overline{t_{23}a_{11}})-\frac{2}{3}\rho _{2}a_{12})$

$+F_{3}(\frac{1}{2}((\tau _{11}+\tau _{12})a_{23}+(\alpha _{21}+\alpha
_{22})t_{13}+\overline{a_{21}t_{12}}+\overline{t_{11}a_{22}})+\frac{2}{3}%
\rho _{1}a_{23}$

$\ \ \ -\frac{1}{2}((\tau _{21}+\tau _{22})a_{13}+(\alpha _{11}+\alpha
_{12})t_{23}+\overline{a_{11}t_{22}}+\overline{t_{21}a_{12}})-\frac{2}{3}%
\rho _{2}a_{13})$

$-\frac{1}{4}(((\gamma _{12}\psi _{23}+\gamma _{13}\psi
_{22})-2(y_{11},w_{21})-\xi _{1}\mu _{21}-\zeta _{2}\chi _{11})E_{1}$

$\ \ +((\gamma _{13}\psi _{21}+\gamma _{11}\psi _{23})-2(y_{12},w_{22})-\xi
_{1}\mu _{22}-\zeta _{2}\chi _{12})E_{2}$

$\ \ +((\gamma _{11}\psi _{22}+\gamma _{12}\psi _{21})-2(y_{13},w_{23})-\xi
_{1}\mu _{23}-\zeta _{2}\chi _{13})E_{3}$

$\ +F_{1}((-\psi _{21}y_{11}-\gamma _{11}w_{21}+\overline{y_{12}w_{23}}+%
\overline{w_{22}y_{13}})-\xi _{1}z_{21}-\zeta _{2}x_{11})$

$\ +F_{2}((-\psi _{22}y_{12}-\gamma _{12}w_{22}+\overline{y_{13}w_{21}}+%
\overline{w_{23}y_{11}})-\xi _{1}z_{22}-\zeta _{2}x_{12})$

$\ +F_{3}((-\psi _{23}y_{13}-\gamma _{13}w_{23}+\overline{y_{11}w_{22}}+%
\overline{w_{21}y_{12}})-\xi _{1}z_{23}-\zeta _{2}x_{13})$ $)$

$+\frac{1}{4}(((\gamma _{22}\psi _{13}+\gamma _{23}\psi
_{12})-2(y_{21},w_{11})-\xi _{2}\mu _{11}-\zeta _{1}\chi _{21})E_{1}$

$\ \ +((\gamma _{23}\psi _{11}+\gamma _{21}\psi _{13})-2(y_{22},w_{12})-\xi
_{2}\mu _{12}-\zeta _{1}\chi _{22})E_{2}$

$\ \ +((\gamma _{21}\psi _{12}+\gamma _{22}\psi _{11})-2(y_{23},w_{13})-\xi
_{2}\mu _{13}-\zeta _{1}\chi _{23})E_{3}$

$\ +F_{1}((-\psi _{11}y_{21}-\gamma _{21}w_{11}+\overline{y_{22}w_{13}}+%
\overline{w_{12}y_{23}})-\xi _{2}z_{11}-\zeta _{1}x_{21})$

$\ +F_{2}((-\psi _{12}y_{22}-\gamma _{22}w_{12}+\overline{y_{23}w_{11}}+%
\overline{w_{13}y_{21}})-\xi _{2}z_{12}-\zeta _{1}x_{22})$

$\ +F_{3}((-\psi _{13}y_{23}-\gamma _{23}w_{13}+\overline{y_{21}w_{12}}+%
\overline{w_{11}y_{22}})-\xi _{2}z_{13}-\zeta _{1}x_{23})$ $),$

\bigskip

$([R_{1},R_{2}]_{8\Phi })_{B}$

=$F_{1}(\textrm{g}_{d}(D_{1})b_{21})+F_{2}(\nu \textrm{g}_{d}(D_{1})b_{22})+F_{3}(\nu
^{2}\textrm{g}_{d}(D_{1})b_{23})$

$-F_{1}(\textrm{g}_{d}(D_{2})b_{11})-F_{2}(\nu \textrm{g}_{d}(D_{2})b_{12})-F_{3}(\nu
^{2}\textrm{g}_{d}(D_{2})b_{13})$

$+((m_{13},b_{23})-(m_{12},b_{22})-(m_{23},b_{13})+(m_{22},b_{12}))E_{1}$

$+((m_{11},b_{21})-(m_{13},b_{23})-(m_{21},b_{11})+(m_{23},b_{13}))E_{2}$

$+((m_{12},b_{22})-(m_{11},b_{21})-(m_{22},b_{12})+(m_{21},b_{11}))E_{3}$

$+\frac{1}{2}F_{1}((\beta _{23}-\beta _{22})m_{11}-\overline{b_{22}m_{13}}+%
\overline{m_{12}b_{23}}-(\beta _{13}-\beta _{12})m_{21}+\overline{%
b_{12}m_{23}}-\overline{m_{22}b_{13}})$

$+\frac{1}{2}F_{2}((\beta _{21}-\beta _{23})m_{12}-\overline{b_{23}m_{11}}+%
\overline{m_{13}b_{21}}-(\beta _{11}-\beta _{13})m_{22}+\overline{%
b_{13}m_{21}}-\overline{m_{23}b_{11}})$

$+\frac{1}{2}F_{3}((\beta _{22}-\beta _{21})m_{13}-\overline{b_{21}m_{12}}+%
\overline{m_{11}b_{22}}-(\beta _{12}-\beta _{11})m_{23}+\overline{%
b_{11}m_{22}}-\overline{m_{21}b_{12}})$

$-(\tau _{11}\beta _{21}+(t_{12},b_{22})+(t_{13},b_{23})-\tau _{21}\beta
_{11}-(t_{22},b_{12})-(t_{23},b_{13})+\frac{2}{3}\rho _{1}\beta _{21}$

$-\frac{2}{3}\rho _{2}\beta _{11})E_{1}$

$-(\tau _{12}\beta _{22}+(t_{13},b_{23})+(t_{11},b_{21})-\tau _{22}\beta
_{12}-(t_{23},b_{13})-(t_{21},b_{11})+\frac{2}{3}\rho _{1}\beta _{22}$

$-\frac{2}{3}\rho _{2}\beta _{12})E_{2}$

$-((-\tau _{11}-\tau _{12})\beta
_{23}+(t_{11},b_{21})+(t_{12},b_{22})-(-\tau _{21}-\tau _{22})\beta
_{13}-(t_{21},b_{11})$

$-(t_{22},b_{12})$
$\ +\frac{2}{3}\rho _{1}\beta _{23}-\frac{2}{3}\rho _{2}\beta _{13})E_{3}$

$-F_{1}(\frac{1}{2}(-\tau _{11}b_{21}+(\beta _{22}+\beta _{23})t_{11}+%
\overline{b_{22}t_{13}}+\overline{t_{12}b_{23}})+\frac{2}{3}\rho _{1}b_{21}$

$\ -\frac{1}{2}(-\tau _{21}b_{11}+(\beta _{12}+\beta _{13})t_{21}+\overline{%
b_{12}t_{23}}+\overline{t_{22}b_{13}})-\frac{2}{3}\rho _{2}b_{11})$

$-F_{2}(\frac{1}{2}(-\tau _{12}b_{22}+(\beta _{23}+\beta _{21})t_{12}+%
\overline{b_{23}t_{11}}+\overline{t_{13}b_{21}})+\frac{2}{3}\rho _{1}b_{22}$

$\ -\frac{1}{2}(-\tau _{22}b_{12}+(\beta _{13}+\beta _{11})t_{22}+\overline{%
b_{13}t_{21}}+\overline{t_{23}b_{11}})-\frac{2}{3}\rho _{2}b_{12})$

$-F_{3}(\frac{1}{2}((\tau _{11}+\tau _{12})b_{23}+(\beta _{21}+\beta
_{22})t_{13}+\overline{b_{21}t_{12}}+\overline{t_{11}b_{22}})+\frac{2}{3}%
\rho _{1}b_{23}$

$\ \ -\frac{1}{2}((\tau _{21}+\tau _{22})b_{13}+(\beta _{11}+\beta
_{12})t_{23}+\overline{b_{11}t_{22}}+\overline{t_{21}b_{12}})-\frac{2}{3}%
\rho _{2}b_{13})$

+$\frac{1}{4}(((\chi _{12}\mu _{23}+\chi _{13}\mu
_{22})-2(x_{11},z_{21})-\eta _{1}\psi _{21}-\omega _{2}\gamma _{11})E_{1}$

$\ \ +((\chi _{13}\mu _{21}+\chi _{11}\mu _{23})-2(x_{12},z_{22})-\eta
_{1}\psi _{22}-\omega _{2}\gamma _{12})E_{2}$

\ $\ +((\chi _{11}\mu _{22}+\chi _{12}\mu _{21})-2(x_{13},z_{23})-\eta
_{1}\psi _{23}-\omega _{2}\gamma _{13})E_{3}$

$\ +F_{1}((-\mu _{21}x_{11}-\chi _{11}z_{21}+\overline{x_{12}z_{23}}+%
\overline{z_{22}x_{13}})-\eta _{1}w_{21}-\omega _{2}y_{11})$

$\ +F_{2}((-\mu _{22}x_{12}-\chi _{12}z_{22}+\overline{x_{13}z_{21}}+%
\overline{z_{23}x_{11}})-\eta _{1}w_{22}-\omega _{2}y_{12})$

$\ +F_{3}((-\mu _{23}x_{13}-\chi _{13}z_{23}+\overline{x_{11}z_{22}}+%
\overline{z_{21}x_{12}})-\eta _{1}w_{23}-\omega _{2}y_{13})$ $)$

$-\frac{1}{4}(((\chi _{22}\mu _{13}+\chi _{23}\mu
_{12})-2(x_{21},z_{11})-\eta _{2}\psi _{11}-\omega _{1}\gamma _{21})E_{1}$

$\ \ +((\chi _{23}\mu _{11}+\chi _{21}\mu _{13})-2(x_{22},z_{12})-\eta
_{2}\psi _{12}-\omega _{1}\gamma _{22})E_{2}$

\ $\ +((\chi _{21}\mu _{12}+\chi _{22}\mu _{11})-2(x_{23},z_{13})-\eta
_{2}\psi _{13}-\omega _{1}\gamma _{23})E_{3}$

$\ +F_{1}((-\mu _{11}x_{21}-\chi _{21}z_{11}+\overline{x_{22}z_{13}}+%
\overline{z_{12}x_{23}})-\eta _{2}w_{11}-\omega _{1}y_{21})$

$\ +F_{2}((-\mu _{12}x_{22}-\chi _{22}z_{12}+\overline{x_{23}z_{11}}+%
\overline{z_{13}x_{21}})-\eta _{2}w_{12}-\omega _{1}y_{22})$

$\ +F_{3}((-\mu _{13}x_{23}-\chi _{23}z_{13}+\overline{x_{21}z_{12}}+%
\overline{z_{11}x_{22}})-\eta _{2}w_{13}-\omega _{1}y_{23})$ $),$

\bigskip

$([R_{1},R_{2}]_{8\Phi })_{\rho }$

=$\alpha _{11}\beta _{21}+\alpha _{12}\beta _{22}+\alpha _{13}\beta
_{23}-\alpha _{21}\beta _{11}-\alpha _{22}\beta _{12}-\alpha _{23}\beta
_{13} $

$%
+2(a_{11},b_{21})+2(a_{12},b_{22})+2(a_{13},b_{23})-2(a_{21},b_{11})-2(a_{22},b_{12})-2(a_{23},b_{13}) 
$

+$\frac{1}{8}(\sum_{i=1}^{3}(\chi _{1i}\psi _{2i}+2(x_{1i},w_{2i})+\mu
_{2i}\gamma _{1i}+2(z_{2i},y_{1i}))-3(\xi _{1}\omega _{2}+\zeta _{2}%
\eta _{1})$ $)$

\ -$\frac{1}{8}(\sum_{i=1}^{3}(\chi _{2i}\psi _{1i}+2(x_{2i},w_{1i})+\mu
_{1i}\gamma _{2i}+2(z_{1i},y_{2i}))-3(\xi _{2}\omega _{1}+\zeta _{1}%
\eta _{2})$ $),$

\bigskip

$%
[R_{1},R_{2}]_{8P}=(([R_{1},R_{2}]_{8P})_{X},([R_{1},R_{2}]_{8P})_{Y},([R_{1},R_{2}]_{8P})_{\xi },([R_{1},R_{2}]_{8P})_{\eta }), 
$

$([R_{1},R_{2}]_{8P})_{X}$

$\ =F_{1}(\textrm{g}_{d}(D_{1})x_{21})+F_{2}(\nu \textrm{g}_{d}(D_{1})x_{22})+F_{3}(\nu
^{2}\textrm{g}_{d}(D_{1})x_{23})$

$\ -F_{1}(\textrm{g}_{d}(D_{2})x_{11})-F_{2}(\nu \textrm{g}_{d}(D_{2})x_{12})-F_{3}(\nu
^{2}\textrm{g}_{d}(D_{2})x_{13})$

$\ +(-(m_{12},x_{22})+(m_{13},x_{23}))E_{1}+(-(m_{13},x_{23})+(m_{11},x_{21}))E_{2}$

$\ +(-(m_{11},x_{21})+(m_{12},x_{22}))E_{3}$

$\
-(-(m_{22},x_{12})+(m_{23},x_{13}))E_{1}-(-(m_{23},x_{13})+(m_{21},x_{11}))E_{2}$

$\ -(-(m_{21},x_{11})+(m_{22},x_{12}))E_{3}$

$\ +F_{1}(\frac{1}{2}(\overline{m_{12}x_{23}}-\overline{x_{22}m_{13}}+(-\chi
_{22}+\chi _{23})m_{11}))$

$\ +F_{2}(\frac{1}{2}(\overline{m_{13}x_{21}}-\overline{x_{23}m_{11}}+(-\chi
_{23}+\chi _{21})m_{12}))$

$\ +F_{3}(\frac{1}{2}(\overline{m_{11}x_{22}}-\overline{x_{21}m_{12}}+(-\chi
_{21}+\chi _{22})m_{13}))$

$\ -F_{1}(\frac{1}{2}(\overline{m_{22}x_{13}}-\overline{x_{12}m_{23}}+(-\chi
_{12}+\chi _{13})m_{21}))$

$\ -F_{2}(\frac{1}{2}(\overline{m_{23}x_{11}}-\overline{x_{13}m_{21}}+(-\chi
_{13}+\chi _{11})m_{22}))$

$\ -F_{3}(\frac{1}{2}(\overline{m_{21}x_{12}}-\overline{x_{11}m_{22}}+(-\chi
_{11}+\chi _{12})m_{23}))$

$\ +(\tau _{11}\chi _{21}+(t_{12},x_{22})+(t_{13},x_{23}))E_{1}$

$\ +(\tau _{12}\chi _{22}+(t_{13},x_{23})+(t_{11},x_{21}))E_{2}$

$\ +((-\tau _{11}-\tau _{12})\chi
_{23}+(t_{11},x_{21})+(t_{12},x_{22}))E_{3} $

$\ -(\tau _{21}\chi _{11}+(t_{22},x_{12})+(t_{23},x_{13}))E_{1}$

$\ -(\tau _{22}\chi _{12}+(t_{23},x_{13})+(t_{21},x_{11}))E_{2}$

$\ -((-\tau _{21}-\tau _{22})\chi
_{13}+(t_{21},x_{11})+(t_{22},x_{12}))E_{3} $

$\ +F_{1}(\frac{1}{2}(-\tau _{11}x_{21}+(\chi _{22}+\chi _{23})t_{11}+%
\overline{x_{22}t_{13}}+\overline{t_{12}x_{23}}))$

$\ +F_{2}(\frac{1}{2}(-\tau _{12}x_{22}+(\chi _{23}+\chi _{21})t_{12}+%
\overline{x_{23}t_{11}}+\overline{t_{13}x_{21}}))$

$\ +F_{3}(\frac{1}{2}((\tau _{11}+\tau _{12})x_{23}+(\chi _{21}+\chi
_{22})t_{13}+\overline{x_{21}t_{12}}+\overline{t_{11}x_{22}}))$

$\ -F_{1}(\frac{1}{2}(-\tau _{21}x_{11}+(\chi _{12}+\chi _{13})t_{21}+%
\overline{x_{12}t_{23}}+\overline{t_{22}x_{13}}))$

$\ -F_{2}(\frac{1}{2}(-\tau _{22}x_{12}+(\chi _{13}+\chi _{11})t_{22}+%
\overline{x_{13}t_{21}}+\overline{t_{23}x_{11}}))$

$\ -F_{3}(\frac{1}{2}((\tau _{21}+\tau _{22})x_{13}+(\chi _{11}+\chi
_{12})t_{23}+\overline{x_{11}t_{22}}+\overline{t_{21}x_{12}}))$

$\ -\frac{1}{3}\rho _{1}(\chi _{21}E_{1}+\chi _{22}E_{2}+\chi
_{23}E_{3}+F_{1}(x_{21})+F_{2}(x_{22})+F_{3}(x_{23}))$

$\ +\frac{1}{3}\rho _{2}(\chi _{11}E_{1}+\chi _{12}E_{2}+\chi
_{13}E_{3}+F_{1}(x_{11})+F_{2}(x_{12})+F_{3}(x_{13}))$

$\ +((\beta _{12}\gamma _{23}+\beta _{13}\gamma
_{22})-2(b_{11},y_{21}))E_{1} $

$\ +((\beta _{13}\gamma _{21}+\beta _{11}\gamma
_{23})-2(b_{12},y_{22}))E_{2} $

$\ +((\beta _{11}\gamma _{22}+\beta _{12}\gamma
_{21})-2(b_{13},y_{23}))E_{3} $

$\ -((\beta _{22}\gamma _{13}+\beta _{23}\gamma
_{12})-2(b_{21},y_{11}))E_{1} $

$\ -((\beta _{23}\gamma _{11}+\beta _{21}\gamma
_{13})-2(b_{22},y_{12}))E_{2} $

$\ -((\beta _{21}\gamma _{12}+\beta _{22}\gamma
_{11})-2(b_{23},y_{13}))E_{3} $

$\ +F_{1}(-\gamma _{21}b_{11}-\beta _{11}y_{21}+\overline{b_{12}y_{23}}+%
\overline{y_{22}b_{13}})$

$\ +F_{2}(-\gamma _{22}b_{12}-\beta _{12}y_{22}+\overline{b_{13}y_{21}}+%
\overline{y_{23}b_{11}})$

$\ +F_{3}(-\gamma _{23}b_{13}-\beta _{13}y_{23}+\overline{b_{11}y_{22}}+%
\overline{y_{21}b_{12}})$

$\ -F_{1}(-\gamma _{11}b_{21}-\beta _{21}y_{11}+\overline{b_{22}y_{13}}+%
\overline{y_{12}b_{23}})$

$\ -F_{2}(-\gamma _{12}b_{22}-\beta _{22}y_{12}+\overline{b_{23}y_{11}}+%
\overline{y_{13}b_{21}})$

$\ -F_{3}(-\gamma _{13}b_{23}-\beta _{23}y_{13}+\overline{b_{21}y_{12}}+%
\overline{y_{11}b_{22}})$

$\ +\eta _{2}(\alpha _{11}E_{1}+\alpha _{12}E_{2}+\alpha
_{13}E_{3}+F_{1}(a_{11})+F_{2}(a_{12})+F_{3}(a_{13}))$

$\ -\eta _{1}(\alpha _{21}E_{1}+\alpha _{22}E_{2}+\alpha
_{23}E_{3}+F_{1}(a_{21})+F_{2}(a_{22})+F_{3}(a_{23}))$

$\ +r_{1}(\chi _{21}E_{1}+\chi _{22}E_{2}+\chi
_{23}E_{3}+F1(x_{21})+F_{2}(x_{22})+F_{3}(x_{23}))$

$\ -r_{2}(\chi _{11}E_{1}+\chi _{12}E_{2}+\chi
_{13}E_{3}+F1(x_{11})+F_{2}(x_{12})+F_{3}(x_{13}))$

$\ +s_{1}(\mu _{21}E_{1}+\mu _{22}E_{2}+\mu
_{23}E_{3}+F1(z_{21})+F_{2}(z_{22})+F_{3}(z_{23}))$

$\ -s_{2}(\mu _{11}E_{1}+\mu _{12}E_{2}+\mu
_{13}E_{3}+F1(z_{11})+F_{2}(z_{12})+F_{3}(z_{13})),$

$([R_{1},R_{2}]_{8P})_{Y}$

$\ =F_{1}(\textrm{g}_{d}(D_{1})y_{21})+F_{2}(\nu \textrm{g}_{d}(D_{1})y_{22})+F_{3}(\nu
^{2}\textrm{g}_{d}(D_{1})y_{23})$

$\ -F_{1}(\textrm{g}_{d}(D_{2})y_{11})-F_{2}(\nu \textrm{g}_{d}(D_{2})y_{12})-F_{3}(\nu
^{2}\textrm{g}_{d}(D_{2})y_{13})$

$\
+(-(m_{12},y_{22})+(m_{13},y_{23}))E_{1}+(-(m_{13},y_{23})+(m_{11},y_{21}))E_{2}$

$\ +(-(m_{11},y_{21})+(m_{12},y_{22}))E_{3} $

$\
-(-(m_{22},y_{12})+(m_{23},y_{13}))E_{1}-(-(m_{23},y_{13})+(m_{21},y_{11}))E_{2}$

$\ -(-(m_{21},y_{11})+(m_{22},y_{12}))E_{3} $

$\ +F_{1}(\frac{1}{2}(\overline{m_{12}y_{23}}-\overline{y_{22}m_{13}}%
+(-\gamma _{22}+\gamma _{23})m_{11}))$

$\ +F_{2}(\frac{1}{2}(\overline{m_{13}y_{21}}-\overline{y_{23}m_{11}}%
+(-\gamma _{23}+\gamma _{21})m_{12}))$

$\ +F_{3}(\frac{1}{2}(\overline{m_{11}y_{22}}-\overline{y_{21}m_{12}}%
+(-\gamma _{21}+\gamma _{22})m_{13}))$

$\ -F_{1}(\frac{1}{2}(\overline{m_{22}y_{13}}-\overline{y_{12}m_{23}}%
+(-\gamma _{12}+\gamma _{13})m_{21}))$

$\ -F_{2}(\frac{1}{2}(\overline{m_{23}y_{11}}-\overline{y_{13}m_{21}}%
+(-\gamma _{13}+\gamma _{11})m_{22}))$

$\ -F_{3}(\frac{1}{2}(\overline{m_{21}y_{12}}-\overline{y_{11}m_{22}}%
+(-\gamma _{11}+\gamma _{12})m_{23}))$

$\ -(\tau _{11}\gamma _{21}+(t_{12},y_{22})+(t_{13},y_{23}))E_{1}$

$\ -(\tau _{12}\gamma _{22}+(t_{13},y_{23})+(t_{11},y_{21}))E_{2}$

$\ -((-\tau _{11}-\tau _{12})\gamma
_{23}+(t_{11},y_{21})+(t_{12},y_{22}))E_{3}$

$\ +(\tau _{21}\gamma _{11}+(t_{22},y_{12})+(t_{23},y_{13}))E_{1}$

$\ +(\tau _{22}\gamma _{12}+(t_{23},y_{13})+(t_{21},y_{11}))E_{2}$

$\ +((-\tau _{21}-\tau _{22})\gamma
_{13}+(t_{21},y_{11})+(t_{22},y_{12}))E_{3}$

$\ -F_{1}(\frac{1}{2}(-\tau _{11}y_{21}+(\gamma _{22}+\gamma _{23})t_{11}+%
\overline{y_{22}t_{13}}+\overline{t_{12}y_{23}}))$

$\ -F_{2}(\frac{1}{2}(-\tau _{12}y_{22}+(\gamma _{23}+\gamma _{21})t_{12}+%
\overline{y_{23}t_{11}}+\overline{t_{13}y_{21}}))$

$\ -F_{3}(\frac{1}{2}((\tau _{11}+\tau _{12})y_{23}+(\gamma _{21}+\gamma
_{22})t_{13}+\overline{y_{21}t_{12}}+\overline{t_{11}y_{22}}))$

$\ +F_{1}(\frac{1}{2}(-\tau _{21}y_{11}+(\gamma _{12}+\gamma _{13})t_{21}+%
\overline{y_{12}t_{23}}+\overline{t_{22}y_{13}}))$

$\ +F_{2}(\frac{1}{2}(-\tau _{22}y_{12}+(\gamma _{13}+\gamma _{11})t_{22}+%
\overline{y_{13}t_{21}}+\overline{t_{23}y_{11}}))$

$\ +F_{3}(\frac{1}{2}((\tau _{21}+\tau _{22})y_{13}+(\gamma _{11}+\gamma
_{12})t_{23}+\overline{y_{11}t_{22}}+\overline{t_{21}y_{12}}))$

$\ +\frac{1}{3}\rho _{1}(\gamma _{21}E_{1}+\gamma _{22}E_{2}+\gamma
_{23}E_{3}+F_{1}(y_{21})+F_{2}(y_{22})+F_{3}(y_{23}))$

$\ -\frac{1}{3}\rho _{2}(\gamma _{11}E_{1}+\gamma _{12}E_{2}+\gamma
_{13}E_{3}+F_{1}(y_{11})+F_{2}(y_{12})+F_{3}(y_{13}))$

$\ +((\alpha _{12}\chi _{23}+\alpha _{13}\chi _{22})-2(a_{11},x_{21}))E_{1}$

$\ +((\alpha _{13}\chi _{21}+\alpha _{11}\chi _{23})-2(a_{12},x_{22}))E_{2}$

$\ +((\alpha _{11}\chi _{22}+\alpha _{12}\chi _{21})-2(a_{13},x_{23}))E_{3}$

$\ -((\alpha _{22}\chi _{13}+\alpha _{23}\chi _{12})-2(a_{21},x_{11}))E_{1}$

$\ -((\alpha _{23}\chi _{11}+\alpha _{21}\chi _{13})-2(a_{22},x_{12}))E_{2}$

$\ -((\alpha _{21}\chi _{12}+\alpha _{22}\chi _{11})-2(a_{23},x_{13}))E_{3}$

$\ +F_{1}(-\chi _{21}a_{11}-\alpha _{11}x_{21}+\overline{a_{12}x_{23}}+%
\overline{x_{22}a_{13}})$

$\ +F_{2}(-\chi _{22}a_{12}-\alpha _{12}x_{22}+\overline{a_{13}x_{21}}+%
\overline{x_{23}a_{11}})$

$\ +F_{3}(-\chi _{23}a_{13}-\alpha _{13}x_{23}+\overline{a_{11}x_{22}}+%
\overline{x_{21}a_{12}})$

$\ -F_{1}(-\chi _{11}a_{21}-\alpha _{21}x_{11}+\overline{a_{22}x_{13}}+%
\overline{x_{12}a_{23}})$

$\ -F_{2}(-\chi _{12}a_{22}-\alpha _{22}x_{12}+\overline{a_{23}x_{11}}+%
\overline{x_{13}a_{21}})$

$\ -F_{3}(-\chi _{13}a_{23}-\alpha _{23}x_{13}+\overline{a_{21}x_{12}}+%
\overline{x_{11}a_{22}})$

$\ +\xi _{2}(\beta _{11}E_{1}+\beta _{12}E_{2}+\beta
_{13}E_{3}+F_{1}(b_{11})+F_{2}(b_{12})+F_{3}(b_{13}))$

$\ -\xi _{1}(\beta _{21}E_{1}+\beta _{22}E_{2}+\beta
_{23}E_{3}+F_{1}(b_{21})+F_{2}(b_{22})+F_{3}(b_{23}))$

$\ +r_{1}(\gamma _{21}E_{1}+\gamma _{22}E_{2}+\gamma
_{23}E_{3}+F1(y_{21})+F_{2}(y_{22})+F_{3}(y_{23}))$

$\ -r_{2}(\gamma _{11}E_{1}+\gamma _{12}E_{2}+\gamma
_{13}E_{3}+F1(y_{11})+F_{2}(y_{12})+F_{3}(y_{13}))$

$\ +s_{1}(\psi _{21}E_{1}+\psi _{22}E_{2}+\psi
_{23}E_{3}+F1(w_{21})+F_{2}(w_{22})+F_{3}(w_{23}))$

$\ -s_{2}(\psi _{11}E_{1}+\psi _{12}E_{2}+\psi
_{13}E_{3}+F1(w_{11})+F_{2}(w_{12})+F_{3}(w_{13})),$

$([R_{1},R_{2}]_{8P})_{\xi }$

$=\sum_{i=1}^{3}(\alpha _{1i}\gamma
_{2i}+2(a_{1i},y_{2i}))-\sum_{i=1}^{3}(\alpha _{2i}\gamma
_{1i}+2(a_{2i},y_{1i}))$

$\ \ +\rho _{1}\xi _{2}-\rho _{2}\xi _{1}+r_{1}\xi _{2}-r_{2}\xi
_{1}+s_{1}\zeta _{2}-s_{2}\zeta _{1},$

$([R_{1},R_{2}]_{8P})_{\eta }$

=$\sum_{i=1}^{3}(\beta _{1i}\chi
_{2i}+2(b_{1i},x_{2i}))-\sum_{i=1}^{3}(\beta _{2i}\chi
_{1i}+2(b_{2i},x_{1i}))$

$\ \ -\rho _{1}\eta _{2}+\rho _{2}\eta _{1}+r_{1}\eta _{2}-r_{2}\eta
_{1}+s_{1}\omega _{2}-s_{2}\omega _{1},$

\bigskip

$%
[R_{1},R_{2}]_{8Q}=(([R_{1},R_{2}]_{8Q})_{Z},([R_{1},R_{2}]_{8Q})_{W},([R_{1},R_{2}]_{8Q})_{\zeta },([R_{1},R_{2}]_{8Q})_{\omega }), 
$

$([R_{1},R_{2}]_{8Q})_{Z}$

$\ =F_{1}(\textrm{g}_{d}(D_{1})z_{21})+F_{2}(\nu \textrm{g}_{d}(D_{1})z_{22})+F_{3}(\nu
^{2}\textrm{g}_{d}(D_{1})z_{23})$

$\ -F_{1}(\textrm{g}_{d}(D_{2})z_{11})-F_{2}(\nu \textrm{g}_{d}(D_{2})z_{12})-F_{3}(\nu
^{2}\textrm{g}_{d}(D_{2})z_{13})$

$\
+(-(m_{12},z_{22})+(m_{13},z_{23}))E_{1}+(-(m_{13},z_{23})+(m_{11},z_{21}))E_{2}$

\ $+(-(m_{11},z_{21})+(m_{12},z_{22}))E_{3}$

$\
-(-(m_{22},z_{12})+(m_{23},z_{13}))E_{1}-(-(m_{23},z_{13})+(m_{21},z_{11}))E_{2}$

\ $-(-(m_{21},z_{11})+(m_{22},z_{12}))E_{3} $

$\ +F_{1}(\frac{1}{2}(\overline{m_{12}z_{23}}-\overline{z_{22}m_{13}}+(-\mu
_{22}+\mu _{23})m_{11}))$

$\ +F_{2}(\frac{1}{2}(\overline{m_{13}z_{21}}-\overline{z_{23}m_{11}}+(-\mu
_{23}+\mu _{21})m_{12}))$

$\ +F_{3}(\frac{1}{2}(\overline{m_{11}z_{22}}-\overline{z_{21}m_{12}}+(-\mu
_{21}+\mu _{22})m_{13}))$

$\ -F_{1}(\frac{1}{2}(\overline{m_{22}z_{13}}-\overline{z_{12}m_{23}}+(-\mu
_{12}+\mu _{13})m_{21}))$

$\ -F_{2}(\frac{1}{2}(\overline{m_{23}z_{11}}-\overline{z_{13}m_{21}}+(-\mu
_{13}+\mu _{11})m_{22}))$

$\ -F_{3}(\frac{1}{2}(\overline{m_{21}z_{12}}-\overline{z_{11}m_{22}}+(-\mu
_{11}+\mu _{12})m_{23}))$

$\ +(\tau _{11}\mu _{21}+(t_{12},z_{22})+(t_{13},z_{23}))E_{1}$

$\ +(\tau _{12}\mu _{22}+(t_{13},z_{23})+(t_{11},z_{21}))E_{2}$

$\ +((-\tau _{11}-\tau _{12})\mu _{23}+(t_{11},z_{21})+(t_{12},z_{22}))E_{3}$

$\ -(\tau _{21}\mu _{11}+(t_{22},z_{12})+(t_{23},z_{13}))E_{1}$

$\ -(\tau _{22}\mu _{12}+(t_{23},z_{13})+(t_{21},z_{11}))E_{2}$

$\ -((-\tau _{21}-\tau _{22})\mu _{13}+(t_{21},z_{11})+(t_{22},z_{12}))E_{3}$

$\ +F_{1}(\frac{1}{2}(-\tau _{11}z_{21}+(\mu _{22}+\mu _{23})t_{11}+%
\overline{z_{22}t_{13}}+\overline{t_{12}z_{23}}))$

$\ +F_{2}(\frac{1}{2}(-\tau _{12}z_{22}+(\mu _{23}+\mu _{21})t_{12}+%
\overline{z_{23}t_{11}}+\overline{t_{13}z_{21}}))$

$\ +F_{3}(\frac{1}{2}((\tau _{11}+\tau _{12})z_{23}+(\mu _{21}+\mu
_{22})t_{13}+\overline{z_{21}t_{12}}+\overline{t_{11}z_{22}}))$

$\ -F_{1}(\frac{1}{2}(-\tau _{21}z_{11}+(\mu _{12}+\mu _{13})t_{21}+%
\overline{z_{12}t_{23}}+\overline{t_{22}z_{13}}))$

$\ -F_{2}(\frac{1}{2}(-\tau _{22}z_{12}+(\mu _{13}+\mu _{11})t_{22}+%
\overline{z_{13}t_{21}}+\overline{t_{23}z_{11}}))$

$\ -F_{3}(\frac{1}{2}((\tau _{21}+\tau _{22})z_{13}+(\mu _{11}+\mu
_{12})t_{23}+\overline{z_{11}t_{22}}+\overline{t_{21}z_{12}}))$

$\ -\frac{1}{3}\rho _{1}(\mu _{21}E_{1}+\mu _{22}E_{2}+\mu
_{23}E_{3}+F_{1}(z_{21})+F_{2}(z_{22})+F_{3}(z_{23}))$

$\ +\frac{1}{3}\rho _{2}(\mu _{11}E_{1}+\mu _{12}E_{2}+\mu
_{13}E_{3}+F_{1}(z_{11})+F_{2}(z_{12})+F_{3}(z_{13}))$

$\ +((\beta _{12}\psi _{23}+\beta _{13}\psi _{22})-2(b_{11},w_{21}))E_{1}$

$\ +((\beta _{13}\psi _{21}+\beta _{11}\psi _{23})-2(b_{12},w_{22}))E_{2}$

$\ +((\beta _{11}\psi _{22}+\beta _{12}\psi _{21})-2(b_{13},w_{23}))E_{3}$

$\ -((\beta _{22}\psi _{13}+\beta _{23}\psi _{12})-2(b_{21},w_{11}))E_{1}$

$\ -((\beta _{23}\psi _{11}+\beta _{21}\psi _{13})-2(b_{22},w_{12}))E_{2}$

$\ -((\beta _{21}\psi _{12}+\beta _{22}\psi _{11})-2(b_{23},w_{13}))E_{3}$

$\ +F_{1}(-\psi _{21}b_{11}-\beta _{11}w_{21}+\overline{b_{12}w_{23}}+%
\overline{w_{22}b_{13}})$

$\ +F_{2}(-\psi _{22}b_{12}-\beta _{12}w_{22}+\overline{b_{13}w_{21}}+%
\overline{w_{23}b_{11}})$

$\ +F_{3}(-\psi _{23}b_{13}-\beta _{13}w_{23}+\overline{b_{11}w_{22}}+%
\overline{w_{21}b_{12}})$

$\ -F_{1}(-\psi _{11}b_{21}-\beta _{21}w_{11}+\overline{b_{22}w_{13}}+%
\overline{w_{12}b_{23}})$

$\ -F_{2}(-\psi _{12}b_{22}-\beta _{22}w_{12}+\overline{b_{23}w_{11}}+%
\overline{w_{13}b_{21}})$

$\ -F_{3}(-\psi _{13}b_{23}-\beta _{23}w_{13}+\overline{b_{21}w_{12}}+%
\overline{w_{11}b_{22}})$

$\ +\omega _{2}(\alpha _{11}E_{1}+\alpha _{12}E_{2}+\alpha
_{13}E_{3}+F_{1}(a_{11})+F_{2}(a_{12})+F_{3}(a_{13}))$

$\ -\omega _{1}(\alpha _{21}E_{1}+\alpha _{22}E_{2}+\alpha
_{23}E_{3}+F_{1}(a_{21})+F_{2}(a_{22})+F_{3}(a_{23}))$

$\ -r_{1}(\mu _{21}E_{1}+\mu _{22}E_{2}+\mu
_{23}E_{3}+F1(z_{21})+F_{2}(z_{22})+F_{3}(z_{23}))$

$\ +r_{2}(\mu _{11}E_{1}+\mu _{12}E_{2}+\mu
_{13}E_{3}+F1(z_{11})+F_{2}(z_{12})+F_{3}(z_{13}))$

$\ +u_{1}(\chi _{21}E_{1}+\chi _{22}E_{2}+\chi
_{23}E_{3}+F1(x_{21})+F_{2}(x_{22})+F_{3}(x_{23}))$

$\ -u_{2}(\chi _{11}E_{1}+\chi _{12}E_{2}+\chi
_{13}E_{3}+F1(x_{11})+F_{2}(x_{12})+F_{3}(x_{13})),$

$([R_{1},R_{2}]_{8Q})_{W}$

$\ =F_{1}(\textrm{g}_{d}(D_{1})w_{21})+F_{2}(\nu \textrm{g}_{d}(D_{1})w_{22})+F_{3}(\nu
^{2}\textrm{g}_{d}(D_{1})w_{23})$

$\ -F_{1}(\textrm{g}_{d}(D_{2})w_{11})-F_{2}(\nu \textrm{g}_{d}(D_{2})w_{12})-F_{3}(\nu
^{2}\textrm{g}_{d}(D_{2})w_{13})$

$\
+(-(m_{12},w_{22})+(m_{13},w_{23}))E_{1}+(-(m_{13},w_{23})+(m_{11},w_{21}))E_{2}$

\ $+(-(m_{11},w_{21})+(m_{12},w_{22}))E_{3} $

$\
-(-(m_{22},w_{12})+(m_{23},w_{13}))E_{1}-(-(m_{23},w_{13})+(m_{21},w_{11}))E_{2}$

\ $-(-(m_{21},w_{11})+(m_{22},w_{12}))E_{3} $

$\ +F_{1}(\frac{1}{2}(\overline{m_{12}w_{23}}-\overline{w_{22}m_{13}}+(-\psi
_{22}+\psi _{23})m_{11}))$

$\ +F_{2}(\frac{1}{2}(\overline{m_{13}w_{21}}-\overline{w_{23}m_{11}}+(-\psi
_{23}+\psi _{21})m_{12}))$

$\ +F_{3}(\frac{1}{2}(\overline{m_{11}w_{22}}-\overline{w_{21}m_{12}}+(-\psi
_{21}+\psi _{22})m_{13}))$

$\ -F_{1}(\frac{1}{2}(\overline{m_{22}w_{13}}-\overline{w_{12}m_{23}}+(-\psi
_{12}+\psi _{13})m_{21}))$

$\ -F_{2}(\frac{1}{2}(\overline{m_{23}w_{11}}-\overline{w_{13}m_{21}}+(-\psi
_{13}+\psi _{11})m_{22}))$

$\ -F_{3}(\frac{1}{2}(\overline{m_{21}w_{12}}-\overline{w_{11}m_{22}}+(-\psi
_{11}+\psi _{12})m_{23}))$

$\ -(\tau _{11}\psi _{21}+(t_{12},w_{22})+(t_{13},w_{23}))E_{1}$

$\ -(\tau _{12}\psi _{22}+(t_{13},w_{23})+(t_{11},w_{21}))E_{2}$

$\ -((-\tau _{11}-\tau _{12})\psi
_{23}+(t_{11},w_{21})+(t_{12},w_{22}))E_{3} $

$\ +(\tau _{21}\psi _{11}+(t_{22},w_{12})+(t_{23},w_{13}))E_{1}$

$\ +(\tau _{22}\psi _{12}+(t_{23},w_{13})+(t_{21},w_{11}))E_{2}$

$\ +((-\tau _{21}-\tau _{22})\psi
_{13}+(t_{21},w_{11})+(t_{22},w_{12}))E_{3} $

$\ -F_{1}(\frac{1}{2}(-\tau _{11}w_{21}+(\psi _{22}+\psi _{23})t_{11}+%
\overline{w_{22}t_{13}}+\overline{t_{12}w_{23}}))$

$\ -F_{2}(\frac{1}{2}(-\tau _{12}w_{22}+(\psi _{23}+\psi _{21})t_{12}+%
\overline{w_{23}t_{11}}+\overline{t_{13}w_{21}}))$

$\ -F_{3}(\frac{1}{2}((\tau _{11}+\tau _{12})w_{23}+(\psi _{21}+\psi
_{22})t_{13}+\overline{w_{21}t_{12}}+\overline{t_{11}w_{22}}))$

$\ +F_{1}(\frac{1}{2}(-\tau _{21}w_{11}+(\psi _{12}+\psi _{13})t_{21}+%
\overline{w_{12}t_{23}}+\overline{t_{22}w_{13}}))$

$\ +F_{2}(\frac{1}{2}(-\tau _{22}w_{12}+(\psi _{13}+\psi _{11})t_{22}+%
\overline{w_{13}t_{21}}+\overline{t_{23}w_{11}}))$

$\ +F_{3}(\frac{1}{2}((\tau _{21}+\tau _{22})w_{13}+(\psi _{11}+\psi
_{12})t_{23}+\overline{w_{11}t_{22}}+\overline{t_{21}w_{12}}))$

$\ +\frac{1}{3}\rho _{1}(\psi _{21}E_{1}+\psi _{22}E_{2}+\psi
_{23}E_{3}+F_{1}(w_{21})+F_{2}(w_{22})+F_{3}(w_{23}))$

$\ -\frac{1}{3}\rho _{2}(\psi _{11}E_{1}+\psi _{12}E_{2}+\psi
_{13}E_{3}+F_{1}(w_{11})+F_{2}(w_{12})+F_{3}(w_{13}))$

$\ +((\alpha _{12}\mu _{23}+\alpha _{13}\mu _{22})-2(a_{11},z_{21}))E_{1}$

$\ +((\alpha _{13}\mu _{21}+\alpha _{11}\mu _{23})-2(a_{12},z_{22}))E_{2}$

$\ +((\alpha _{11}\mu _{22}+\alpha _{12}\mu _{21})-2(a_{13},z_{23}))E_{3}$

$\ -((\alpha _{22}\mu _{13}+\alpha _{23}\mu _{12})-2(a_{21},z_{11}))E_{1}$

$\ -((\alpha _{23}\mu _{11}+\alpha _{21}\mu _{13})-2(a_{22},z_{12}))E_{2}$

$\ -((\alpha _{21}\mu _{12}+\alpha _{22}\mu _{11})-2(a_{23},z_{13}))E_{3}$

$\ +F_{1}(-\mu _{21}a_{11}-\alpha _{11}z_{21}+\overline{a_{12}z_{23}}+%
\overline{z_{22}a_{13}})$

$\ +F_{2}(-\mu _{22}a_{12}-\alpha _{12}z_{22}+\overline{a_{13}z_{21}}+%
\overline{z_{23}a_{11}})$

$\ +F_{3}(-\mu _{23}a_{13}-\alpha _{13}z_{23}+\overline{a_{11}z_{22}}+%
\overline{z_{21}a_{12}})$

$\ -F_{1}(-\mu _{11}a_{21}-\alpha _{21}z_{11}+\overline{a_{22}z_{13}}+%
\overline{z_{12}a_{23}})$

$\ -F_{2}(-\mu _{12}a_{22}-\alpha _{22}z_{12}+\overline{a_{23}z_{11}}+%
\overline{z_{13}a_{21}})$

$\ -F_{3}(-\mu _{13}a_{23}-\alpha _{23}z_{13}+\overline{a_{21}z_{12}}+%
\overline{z_{11}a_{22}})$

$\ +\zeta _{2}(\beta _{11}E_{1}+\beta _{12}E_{2}+\beta
_{13}E_{3}+F_{1}(b_{11})+F_{2}(b_{12})+F_{3}(b_{13}))$

$\ -\zeta _{1}(\beta _{21}E_{1}+\beta _{22}E_{2}+\beta
_{23}E_{3}+F_{1}(b_{21})+F_{2}(b_{22})+F_{3}(b_{23}))$

$\ -r_{1}(\psi _{21}E_{1}+\psi _{22}E_{2}+\psi
_{23}E_{3}+F1(w_{21})+F_{2}(w_{22})+F_{3}(w_{23}))$

$\ +r_{2}(\psi _{11}E_{1}+\psi _{12}E_{2}+\psi
_{13}E_{3}+F1(w_{11})+F_{2}(w_{12})+F_{3}(w_{13}))$

$\ +u_{1}(\gamma _{21}E_{1}+\gamma _{22}E_{2}+\gamma
_{23}E_{3}+F1(y_{21})+F_{2}(y_{22})+F_{3}(y_{23}))$

$\ -u_{2}(\gamma _{11}E_{1}+\gamma _{12}E_{2}+\gamma
_{13}E_{3}+F1(y_{11})+F_{2}(y_{12})+F_{3}(y_{13})),$

$([R_{1},R_{2}]_{8Q})_{\zeta }$

$=\sum_{i=1}^{3}(\alpha _{1i}\psi
_{2i}+2(a_{1i},w_{2i}))-\sum_{i=1}^{3}(\alpha _{2i}\psi
_{1i}+2(a_{2i},w_{1i}))$

$\ \ +\rho _{1}\xi _{2}-\rho _{2}\xi _{1}-r_{1}\zeta _{2}+r_{2}\zeta
_{1}+u_{1}\xi _{2}-u_{2}\xi _{1},$

$([R_{1},R_{2}]_{8Q})_{\omega }$

=$\sum_{i=1}^{3}(\beta _{1i}\mu _{2i}+2(b_{1i},z_{2i}))-\sum_{i=1}^{3}(\beta
_{2i}\mu _{1i}+2(b_{2i},z_{1i}))$

$\ \ -\rho _{1}\omega _{2}+\rho _{2}\omega _{1}-r_{1}\omega _{2}+r_{2}\omega
_{1}+u_{1}\eta _{2}-u_{2}\eta _{1},$

\bigskip

$[R_{1},R_{2}]_{8r}$

=$-\frac{1}{8}(\sum_{i=1}^{3}(\chi _{1i}\psi
_{2i}+2(x_{1i},w_{2i}))-\sum_{i=1}^{3}(\gamma _{1i}\mu
_{2i}+2(y_{1i},z_{2i}))+\xi _{1}\omega _{2}-\zeta _{2}\eta _{1})$

$\ +\frac{1}{8}(\sum_{i=1}^{3}(\psi _{1i}\chi
_{2i}+2(w_{1i},x_{2i}))-\sum_{i=1}^{3}(\mu _{1i}\gamma
_{2i}+2(z_{1i},y_{2i}))+\xi _{2}\omega _{1}-\zeta _{1}\eta
_{2})$

\ $+s_{1}u_{2}-s_{2}u_{1},$

\bigskip

$[R_{1},R_{2}]_{8s}$

=$\frac{1}{4}(\sum_{i=1}^{3}(\chi _{1i}\gamma
_{2i}+2(x_{1i},y_{2i}))-\sum_{i=1}^{3}(\gamma _{1i}\chi
_{2i}+2(y_{1i},x_{2i}))+\xi _{1}\eta _{2}-\xi _{2}\eta
_{1})$

\ $+2r_{1}s_{2}-2r_{2}s_{1},$

\bigskip

$[R_{1},R_{2}]_{8u}$

=$-\frac{1}{4}(\sum_{i=1}^{3}(\mu _{1i}\psi
_{2i}+2(z_{1i},w_{2i}))-\sum_{i=1}^{3}(\psi _{1i}\mu
_{2i}+2(w_{1i},z_{2i}))+\zeta _{1}\omega _{2}-\zeta _{2}\omega
_{1})$

\ $-2r_{1}u_{2}+2r_{2}u_{1}$.

\bigskip

\emph{Proof}. By \emph{Definition 5.6} and \emph{Lemma 5.8} , we
have

($[R_{1},R_{2}]_{8\Phi })_{D}=[\Phi _{1},\Phi _{2}]_{7D}+(P_{1}%
\times Q_{2})_{D}-(P_{2}\times Q_{1})_{D},$

$([R_{1},R_{2}]_{8\Phi })_{M}=[\Phi _{1},\Phi _{2}]_{7M}+(P_{1}%
\times Q_{2})_{M}-(P_{2}\times Q_{1})_{M},$

$([R_{1},R_{2}]_{8\Phi })_{T}=[\Phi _{1},\Phi _{2}]_{7T}+(P_{1}%
\times Q_{2})_{T}-(P_{2}\times Q_{1})_{T},$

$([R_{1},R_{2}]_{8\Phi })_{A}=[\Phi _{1},\Phi _{2}]_{7A}+(P_{1}%
\times Q_{2})_{A}-(P_{2}\times Q_{1})_{A},$

$([R_{1},R_{2}]_{8\Phi })_{B}=[\Phi _{1},\Phi _{2}]_{7B}+(P_{1}%
\times Q_{2})_{B}-(P_{2}\times Q_{1})_{B},$

$([R_{1},R_{2}]_{8\Phi })_{\rho }=[\Phi _{1},\Phi _{2}]_{7\rho
}+(P_{1}\times Q_{2})_{\rho }-(P_{2}\times Q_{1})_{\rho }),$

$%
[R_{1},R_{2}]_{8P}=(([R_{1},R_{2}]_{8P})_{X},([R_{1},R_{2}]_{8P})_{Y},([R_{1},R_{2}]_{8P})_{\xi },([R_{1},R_{2}]_{8P})_{\eta }), 
$

$([R_{1},R_{2}]_{8P})_{X}=(\delta _{d}\textrm{g}_{d}(D_{1})+\widehat{M}_{1}+%
\widetilde{T}_{1})X_{2}$-$\frac{1}{3}\rho _{1}X_{2}+2B_{1}\times Y_{2}+\eta
_{2}A_{1}$

$\ \ \ \ \ \ \ \ \ \ \ \ \ \ \ \ \ \ \ \ \ \ \ -(\delta _{d}\textrm{g}_{d}(D_{2})+%
\widehat{M}_{2}+\widetilde{T}_{2})X_{1}$+$\frac{1}{3}\rho
_{2}X_{1}-2B_{2}\times Y_{1}-\eta _{1}A_{2}$

$\ \ \ \ \ \ \ \ \ \ \ \ \ \ \ \ \ \ \ \ \ \
+r_{1}X_{2}-r_{2}X_{1}+s_{1}Z_{2}-s_{2}Z_{1},$

$([R_{1},R_{2}]_{8P})_{Y}=2A_{1}\times X_{2}+(\delta _{d}\textrm{g}_{d}(D_{1})+%
\widehat{M}_{1}-\widetilde{T}_{1})Y_{2}+\frac{1}{3}\rho _{1}Y_{2}+\xi
_{2}B_{1}$

$\ \ \ \ \ \ \ \ \ \ \ \ \ \ \ \ \ \ -2A_{2}\times X_{1}-(\delta
_{d}\textrm{g}_{d}(D_{2})+\widehat{M}_{2}-\widetilde{T}_{2})Y_{1}-\frac{1}{3}\rho
_{2}Y_{1}-\xi _{1}B_{2}$

\ \ \ \ \ \ \ \ \ \ \ \ \ \ \ \ \ \ \ \ $%
+r_{1}Y_{2}-r_{2}Y_{1}+s_{1}W_{2}-s_{2}W_{1},$

$([R_{1},R_{2}]_{8P})_{\xi }=(A_{1},Y_{2})+\rho _{1}\xi
_{2}-(A_{2},Y_{1})-\rho _{2}\xi _{1}+r_{1}\xi _{2}-r_{2}\xi
_{1}+s_{1}\zeta _{2}-s_{2}\zeta _{1},$

$([R_{1},R_{2}]_{8P})_{\eta }=(B_{1},X_{2})-\rho _{1}\eta
_{2}-(B_{2},X_{1})+\rho _{2}\eta _{1}+r_{1}\eta _{2}-r_{2}\eta
_{1}+s_{1}\omega _{2}-s_{2}\omega _{1},$

$%
[R_{1},R_{2}]_{8Q}=(([R_{1},R_{2}]_{8Q})_{Z},([R_{1},R_{2}]_{8Q})_{W},([R_{1},R_{2}]_{8Q})_{\zeta },([R_{1},R_{2}]_{8Q})_{\omega }), 
$

$([R_{1},R_{2}]_{8Q})_{Z}=(\delta _{d}\textrm{g}_{d}(D_{1})+\widehat{M}_{1}+%
\widetilde{T}_{1})Z_{2}-\frac{1}{3}\rho _{1}Z_{2}+2B_{1}\times W_{2}+\omega
_{2}A_{1}$

$\ \ \ \ \ \ \ \ \ \ \ \ \ \ \ \ \ \ \ \ \ \ \ -(\delta _{d}\textrm{g}_{d}(D_{2})+%
\widehat{M}_{2}+\widetilde{T}_{2})Z_{1}+\frac{1}{3}\rho
_{2}Z_{1}-2B_{2}\times W_{1}-\omega _{1}A_{2}$

$\ \ \ \ \ \ \ \ \ \ \ \ \ \ \ \ \ \ \ \ \ \
-r_{1}Z_{2}+r_{2}Z_{1}+u_{1}X_{2}-u_{2}X_{1},$

$([R_{1},R_{2}]_{8Q})_{W}=2A_{1}\times Z_{2}+(\delta _{d}\textrm{g}_{d}(D_{1})+%
\widehat{M}_{1}-\widetilde{T}_{1})W_{2}+\frac{1}{3}\rho _{1}W_{2}+\zeta
_{2}B_{1}$

$\ \ \ \ \ \ \ \ \ \ \ \ \ \ \ \ \ \ -2A_{2}\times Z_{1}-(\delta
_{d}\textrm{g}_{d}(D_{2})+\widehat{M}_{2}-\widetilde{T}_{2})W_{1}-\frac{1}{3}\rho
_{2}W_{1}-\zeta _{1}B_{2}$

\ \ \ \ \ \ \ \ \ \ \ \ \ \ \ \ \ $\
-r_{1}W_{2}+r_{2}W_{1}+u_{1}Y_{2}-u_{2}Y_{1},$

$([R_{1},R_{2}]_{8Q})_{\zeta }=(A_{1},W_{2})+\rho _{1}\zeta
_{2}-(A_{2},W_{1})-\rho _{2}\zeta _{1}-r_{1}\zeta _{2}+r_{2}\zeta
_{1}+u_{1}\xi _{2}-u_{2}\xi _{1},$

$([R_{1},R_{2}]_{8Q})_{\omega }=(B_{1},Z_{2})-\rho _{1}\omega
_{2}-(B_{2},Z_{1})+\rho _{2}\omega _{1}-r_{1}\omega _{2}+r_{2}\omega
_{1}+u_{1}\eta _{2}-u_{2}\eta _{1},$

$[R_{1},R_{2}]_{8r}=-\frac{1}{8}\{P_{1},Q_{2}\}+\frac{1}{8}%
\{P_{2},Q_{1}\}+s_{1}u_{2}-s_{2}u_{1},$

\ \ \ \ \ \ \ \ \ \ \ \ \ $\ \ =-\frac{1}{8}((X_{1},W_{2})-(Z_{2},Y_{1})+\xi
_{1}\omega _{2}-\zeta _{2}\eta _{1})+\frac{1}{8}%
((X_{2},W_{1})-(Z_{1},Y_{2})$

\ \ \ \ \ \ \ \ \ \ \ \ \ \ \ \ $+\xi _{2}\omega _{1}-\zeta _{1}\eta _{2})$
$+s_{1}u_{2}-s_{2}u_{1},$

$[R_{1},R_{2}]_{8s}=\frac{1}{4}\{P_{1},P_{2}\}+2r_{1}s_{2}-2r_{2}s_{1},$

\ \ \ \ \ \ \ \ \ \ \ \ \ \ $\ =\frac{1}{4}((X_{1},Y_{2})-(X_{2},Y_{1})+\xi
_{1}\eta _{2}-\xi _{2}\eta _{1})+2r_{1}s_{2}-2r_{2}s_{1},$

$[R_{1},R_{2}]_{8u}=-\frac{1}{4}\{Q_{1},Q_{2}\}-2r_{1}u_{2}+2r_{2}u_{1}$,

\ \ \ \ \ \ \ \ \ \ \ \ \ \ $\ =-\frac{1}{4}((Z_{1},W_{2})-(Z_{2},W_{1})+%
\zeta _{1}\omega _{2}-\zeta _{2}\omega _{1})-2r_{1}u_{2}+2r_{2}u_{1}$.

\noindent
By\emph{\ Lemma 4.13} and \emph{Lemma 4.8} , we have the
expression with respect to ($[R_{1},R_{2}]_{8\Phi })_{D}$ ,
$([R_{1},R_{2}]_{8\Phi })_{M}$ , $([R_{1},R_{2}]_{8\Phi })_{T}$ , $%
([R_{1},R_{2}]_{8\Phi })_{A}$ , $([R_{1},R_{2}]_{8\Phi })_{B}$ and $%
([R_{1},R_{2}]_{8\Phi })_{\rho }$ .

\noindent
Furthermore calculation by \emph{Definition 1.4} and \emph{Definition 4.2}
, we have the 

\noindent
expression
with respect to  $([R_{1},R_{2}]_{8P})_{X}$ ,  $([R_{1},R_{2}]_{8P})_{Y}$ ,  $%
([R_{1},R_{2}]_{8P})_{\xi }$ ,

\noindent
$([R_{1},R_{2}]_{8P})_{\eta }$ ,\ $([R_{1},R_{2}]_{8Q})_{Z}$ ,
 \ $([R_{1},R_{2}]_{8Q})_{W}$ ,\   $([R_{1},R_{2}]_{8Q})_{\zeta }$ ,
 
 \noindent
$([R_{1},R_{2}]_{8Q})_{\omega }$ ,\   $[R_{1},R_{2}]_{8r}$ ,
$[R_{1},R_{2}]_{8s}$ and $[R_{1},R_{2}]_{8u}$ . \ \ \ \ \emph{Q.E.D.}

\bigskip\ 

\section{Preliminaries for the adjoint representation}

\bigskip 

\ \ \ \ \ Let \gg\ be a Lie algebra over complex numbers.
We consider that $ad$ is the representation of \gg\ to \gg \gl(\gg),where \gg \gl(\gg) is the Lie algebra 
of the general linear group of the vector space \gg.

\ \ \ \ \ \ \ \ \ \ \ \ \ \ \ \ \ \ \ \ \ $ad: $\gg\ $\rightarrow $\gg\gl(\gg)

\bigskip

\emph{Definition 6.1.} \ We denote $D\in $\gD$_{4}^{}$ with
bases $\{D_{ij}|0\leq i<j\leq 7\}$ as follows.

\ \ \ \ \ \ \ \ \ \ \ \ \ \ \ \ \ \ \ \ \ $\ D=\sum\limits_{0\leq i<j\leq 7}%
d_{ij}D_{ij} ,\ d_{ij}\in \R$,

\noindent
And we define an $\R$-linear mapping $\ \textrm{fv}:$ \gD$_{4}^{}\rightarrow
\R^{28}$ by

$D=\left( 
\begin{smallmatrix}
0 & d_{01} & d_{02} & d_{03} & d_{04} & d_{05} & d_{06} & d_{07} \\ 
-d_{01} & 0 & d_{12} & d_{13} & d_{14} & d_{15} & d_{16} & d_{17} \\ 
-d_{02} & -d_{12} & 0 & d_{23} & d_{24} & d_{25} & d_{26} & d_{27} \\ 
-d_{03} & -d_{13} & -d_{23} & 0 & d_{34} & d_{35} & d_{36} & d_{37} \\ 
-d_{04} & -d_{14} & -d_{24} & -d_{34} & 0 & d_{45} & d_{46} & d_{47} \\ 
-d_{05} & -d_{15} & -d_{25} & -d_{35} & -d_{45} & 0 & d_{56} & d_{57} \\ 
-d_{06} & -d_{16} & -d_{26} & -d_{36} & -d_{46} & -d_{56} & 0 & d_{67} \\ 
-d_{07} & -d_{17} & -d_{27} & -d_{37} & -d_{47} & -d_{57} & -d_{67} & 0%
\end{smallmatrix}%
\right) \rightarrow \textrm{fv}(D)=\left( 
\begin{smallmatrix}
d_{01} \\ 
d_{02} \\ 
d_{03} \\ 
d_{04} \\ 
d_{05} \\ 
d_{06} \\ 
d_{07} \\ 
d_{12} \\ 
d_{13} \\ 
d_{14} \\ 
d_{15} \\ 
d_{16} \\ 
d_{17} \\ 
d_{23} \\ 
d_{24} \\ 
d_{25} \\ 
d_{26} \\ 
d_{27} \\ 
d_{34} \\ 
d_{35} \\ 
d_{36} \\ 
d_{37} \\ 
d_{45} \\ 
d_{46} \\ 
d_{47} \\ 
d_{56} \\ 
d_{57} \\ 
d_{67}%
\end{smallmatrix}%
\right) $.

\noindent
Then, we identify $D$ with $\textrm{fv}(D)$.

\bigskip

\emph{Definition 6.2.} We denote $m\in $\gC$^{}$ with bases $%
\{e_{i}|0\leq i\leq 7\}$ as follows.

\ \ \ \ \ \ \ \ \ \ \ \ \ \ \ \ \ \ \ \ \ \ $m=\sum\limits_{0\leq i\leq 7}$m$%
_{i}e_{i}$ ,\ m$_{i}\in \R$,

\noindent
And we define $\R$-linear mapping $\ \textrm{fv}:$ \gC$^{}\rightarrow \R^{8}$ by

\ \ \ \ \ \ \ \ \ \ \ \ \ \ \ \ \ \ \ \ \ \ \ \ 
$m=\sum\limits_{0\leq i\leq 7}$m$_{i}e_{i}\rightarrow \textrm{fv}(m)=\left( 
\begin{smallmatrix}
m_{0} \\ 
m_{1} \\ 
m_{2} \\ 
m_{3} \\ 
m_{4} \\ 
m_{5} \\ 
m_{6} \\ 
m_{7}%
\end{smallmatrix}%
\right) .$

\noindent
Then, we identify $m$ with $\textrm{fv}(m)$.

\noindent
For $x=(x_{0},\cdot\cdot\cdot,x_{7}) \in \R^{8}$, also $\overline{x}$ means $(x_{0},-x_{1},\cdot\cdot\cdot,-x_{7})$.

\bigskip

\emph{Remark 6.3.} 
The Caylay multiplication C$_{m}$\ is expressed on $\R^{8}$ as following:

For $x=(x_{0},x_{1},x_{2},x_{3},x_{4},x_{5},x_{6},x_{7}),y=(y_{0}%
,y_{1},y_{2},y_{3},y_{4},y_{5},y_{6},y_{7}) \in \R^{8},$

$\textrm{C}_{m}(x, y)=$

$(-x_{7}y_{7}-x_{6}y_{6}-x_{5}y_{5}-x_{4}y_{4}-x_{3}y_{3}-x_{2}y_{2}%
-x_{1}y_{1}+x_{0}y_{0},$

$x_{6}y_{7}-x_{7}y_{6}+x_{4}y_{5}-x_{5}y_{4}+x_{2}y_{3}-x_{3}y_{2}+x_{0}%
y_{1}+x_{1}y_{0},$

$x_{5}y_{7}-x_{4}y_{6}-x_{7}y_{5}+x_{6}y_{4}-x_{1}y_{3}+x_{0}y_{2}+x_{3}%
y_{1}+x_{2}y_{0},$

$x_{4}y_{7}+x_{5}y_{6}-x_{6}y_{5}-x_{7}y_{4}+x_{0}y_{3}+x_{1}y_{2}-x_{2}%
y_{1}+x_{3}y_{0},$

$-x_{3}y_{7}+x_{2}y_{6}-x_{1}y_{5}+x_{0}y_{4}+x_{7}y_{3}-x_{6}y_{2}+x_{5}%
y_{1}+x_{4}y_{0},$

$-x_{2}y_{7}-x_{3}y_{6}+x_{0}y_{5}+x_{1}y_{4}+x_{6}y_{3}+x_{7}y_{2}-x_{4}%
y_{1}+x_{5}y_{0},$

$-x_{1}y_{7}+x_{0}y_{6}+x_{3}y_{5}-x_{2}y_{4}-x_{5}y_{3}+x_{4}y_{2}+x_{7}%
y_{1}+x_{6}y_{0},$

$x_{0}y_{7}+x_{1}y_{6}+x_{2}y_{5}+x_{3}y_{4}-x_{4}y_{3}-x_{5}y_{2}-x_{6}%
y_{1}+x_{7}y_{0}).$

C$_{m}$ is naturally extended on $\C^{8}$.

\bigskip

\emph{Definition 6.4.} \ We denote $X(\chi ,x)\in $\gJ$^{}$
with bases $\left\{ E_{1},E_{2},E_{3},F_{1},F_{2},F_{3}\right\} $ as follows.
\begin{multline*} 
\ \ \ \ \ \ X=\chi _{1}E_{1}+\chi _{2}E_{2}+\chi _{3}E_{3}+F_{1}(x_{1})+F_{2}(x_{2})+F_{3}(x_{3}),\\
\chi _{1},\chi _{2},\chi _{3}\in \R,x_{1},x_{2},x_{3}\in  \ggC^{}.
\end{multline*} 

\noindent
And we define an $\R$-linear mapping $\ \textrm{fv}:$ \gJ$^{}\rightarrow \R\oplus
\R\oplus \R\oplus \R^{8}\oplus \R^{8}\oplus \R^{8}$ by

\ \ \ \ \ \ \ \ \ \ \ \ \ \ \ \ \ \ \ \ \ \ \ \ 
$ X(\chi ,x)\rightarrow \textrm{fv}(X)=\left( 
\begin{smallmatrix}
\chi _{1} \\ 
\chi _{2} \\ 
\chi _{3} \\ 
\textrm{fv}(x_{1}) \\ 
\textrm{fv}(x_{2}) \\ 
\textrm{fv}(x_{3})%
\end{smallmatrix}%
\right) .$

\bigskip

\emph{Remark 6.5.}  The Jordan multiplication J$_{m}$ and the Freudenthal multiplication F$_{m}$ are expressed on $\R^{27}$ as followings:

For $X=(\chi_{1},\chi_{2},\chi_{3},x_{1},x_{2},x_{3}), Y=(\gamma_{1},\gamma_{2},\gamma_{3},y_{1},y_{2},y_{3}) \in \R^{27}$,
$\chi_{i},\gamma_{i} \in \R, x_{i},y_{i} \in \R^{8}$,

J$_{m}(X, Y)=(\chi_{1}\gamma_{1}+(x_{2},y_{2})+(x_{3},y_{3})$,

\ \ \ \ \ \ \ \ \ \ \ \ \ \ \ \ \ \  $\chi_{2}\gamma_{2}+(x_{3},y_{3})+(x_{1},y_{1}),$

\ \ \ \ \ \ \ \ \ \ \ \ \ \ \ \ \ \  $\chi_{3}\gamma_{3}+(x_{1},y_{1})+(x_{2},y_{2}),$

\ \ \ \ \ \ \ \ \ \ \ \ \ \ \ \ $((\chi_{2}+\chi_{3})y_{1}+(\gamma_{2}+\gamma_{3})x_{1}+\textrm{C}_{m}(\overline{x_{3}},\overline{y_{2}})%
+\textrm{C}_{m}(\overline{y_{3}},\overline{x_{2}}))/2$,

\ \ \ \ \ \ \ \ \ \ \ \ \ \ \ \ $((\chi_{3}+\chi_{1})y_{2}+(\gamma_{3}+\gamma_{1})x_{2}+\textrm{C}_{m}(\overline{x_{1}},\overline{y_{3}})%
+\textrm{C}_{m}(\overline{y_{1}},\overline{x_{3}}))/2$,

\ \ \ \ \ \ \ \ \ \ \ \ \ \ \ \ $((\chi_{1}+\chi_{2})y_{3}+(\gamma_{1}+\gamma_{2})x_{3}+\textrm{C}_{m}(\overline{x_{2}},\overline{y_{1}})%
+\textrm{C}_{m}(\overline{y_{2}},\overline{x_{1}}))/2 \ )$.

F$_{m}(X, Y)=((\chi_{2}\gamma_{3}+\chi_{3}\gamma_{2})/2-(x_{1},y_{1})$,

\ \ \ \ \ \ \ \ \ \ \ \ \ \ \ \ \ \  $(\chi_{3}\gamma_{1}+\chi_{1}\gamma_{1})/2-(x_{2},y_{2})$,

\ \ \ \ \ \ \ \ \ \ \ \ \ \ \ \ \ \  $(\chi_{1}\gamma_{2}+\chi_{2}\gamma_{1})/2-(x_{3},y_{3})$,

\ \ \ \ \ \ \ \ \ \ \ \ \ \ \ \ $(-\gamma_{1}x_{1}-\chi_{1}y_{1}+\textrm{C}_{m}(\overline{x_{3}},\overline{y_{2}})%
+\textrm{C}_{m}(\overline{y_{3}},\overline{x_{2}}))/2$,

\ \ \ \ \ \ \ \ \ \ \ \ \ \ \ \ $(-\gamma_{2}x_{2}-\chi_{2}y_{2}+\textrm{C}_{m}(\overline{x_{1}},\overline{y_{3}})%
+\textrm{C}_{m}(\overline{y_{1}},\overline{x_{3}}))/2$,

\ \ \ \ \ \ \ \ \ \ \ \ \ \ \ \ $(-\gamma_{3}x_{3}-\chi_{3}y_{3}+\textrm{C}_{m}(\overline{x_{2}},\overline{y_{1}})%
+\textrm{C}_{m}(\overline{y_{2}},\overline{x_{1}}))/2 \ )$.

J$_{m}$ and F$_{m}$ are naturally extended on $\C^{27}$.

\bigskip

\emph{Lemma 6.6.}  We define the inner product on $\R^{27}$ in the sane way as on \gJ.

\noindent
For $A=(a_{ij}) \in M(27 \times 27, \R)$, the transpose with the inner product is expressed
by

 $\ \ \ \ \ \ \ \ \ \ \ \ \ \ \ \ \
\begin{array}
[c]{cccccc}%
\text{ }1\text{ \ } & \text{ \ \ \ \ }1\text{ } & \text{ \ \ \ \ \ }1\text{ } & \text{
\ \ \ \ \ }8\text{ } & \text{ \ \ }8\text{ } & \text{ \ \ \ \ }8
\end{array}
$

$A^{t}=%
\begin{array}
[c]{c}%
1\\
1\\
1\\
8\\
8\\
8
\end{array}
\left(
\begin{tabular}
[c]{ccc|ccc}%
$a_{11}$ & $a_{21}$ & $a_{31}$ & $2a_{l1}$ & $2a_{m1}$ & $2a_{n1}$\\
$a_{12}$ & $a_{22}$ & $a_{32}$ & $2a_{l2}$ & $2a_{m2}$ & $2a_{n2}$\\
$a_{13}$ & $a_{23}$ & $a_{33}$ & $2a_{l3}$ & $2a_{m3}$ & $2a_{n3}$\\\hline
$a_{1l}/2$ & $a_{2l}/2$ & $a_{3l}/2$ &  &  & \\
$a_{1m}/2$ & $a_{2m}/2$ & $a_{3m}/2$ &  & $a_{ji}$ & \\
$a_{1n}/2$ & $a_{2n}/2$ & $a_{3n}/2$ &  &  &
\end{tabular}
\right)  $

\bigskip

\emph{Proof.}  For $X, Y \in \R^{27}$, by calculating using Maxima, we have $(AX, Y)=(X, A^{t}Y)$. \ \ \ \ \ \emph{Q.E.D.}

\bigskip

\emph{\ Definition 6.7.} \ We consider \gR$_{4}^{}=\{(D,M)|D\in $\gD$_{4}^{},M\in $\gA$^{}\}$ is a vector
space of
$\underbrace{\R^{28}}_{\mathfrak{D}_{4}}\oplus \underbrace{\R^{8}\oplus \R^{8}\oplus \R^{8}}_{\mathfrak{A}}.$

\noindent
We identify $(D,M)\in $ \gR$_{4}^{}$ as an element of a
vector space by

\gR$_{\mathbf{4}}^{}\ni $ $(D,M)\rightarrow
\textrm{fv}(D,M)=\left( 
\begin{smallmatrix}
\textrm{fv}(D) \\ 
\textrm{fv}(m_{1}) \\ 
\textrm{fv}(m_{2}) \\ 
\textrm{fv}(m_{3})%
\end{smallmatrix}%
\right) \in \R^{28}\oplus \R^{8}\oplus \R^{8}\oplus \R^{8}=\R^{28}\oplus \R^{24},$

\noindent
where $M=A_{1}(m_{1})+A_{2}(m_{2})+A_{3}(m_{3}),m_{1},m_{2},m_{3}\in $\gC$^{}.$

\bigskip

\emph{\ Definition 6.8. }We consider \gR$_{6}^{}=\{(D,M,T)|D\in $\gD$_{4}^{},M\in $\gA$^{},T\in $%
\gJ$_{0}^{}\}$ is a vector
space of $\underbrace{\R^{28}}_{\mathfrak{D}_{4}}\oplus \underbrace{\R^{8}\oplus \R^{8}\oplus \R^{8}%
}_{\mathfrak{A}}\oplus \underbrace{\R\oplus \R\oplus \R^{8}\oplus \R^{8}\oplus \R^{8}}_{\mathfrak{J}_{0}}.$

\noindent
We identify $(D,M,T)\in $\gR$_{6}^{}$ as an element of a
vector space by

\bigskip

\ \ \gR$_{6}^{}\ni (D,M,T)\rightarrow
\textrm{fv}(D,M,T)=\left( 
\begin{smallmatrix}
\textrm{fv}(D) \\ 
\textrm{fv}(m_{1}) \\ 
\textrm{fv}(m_{2}) \\ 
\textrm{fv}(m_{3}) \\ 
\tau _{1} \\ 
\tau _{2} \\ 
\textrm{fv}(t_{1}) \\ 
\textrm{fv}(t_{2}) \\ 
\textrm{fv}(t_{3})%
\end{smallmatrix}%
\right) \in 
\begin{array}{l}
\R^{28} \\ 
\oplus \R^{24} \\ 
\oplus \R^{26} \\ 
=\R^{78}%
\end{array}%
,$

\bigskip

\noindent
where $T=\tau _{1}E_{1}+\tau _{2}E_{2}+(-\tau _{1}-\tau
_{2})E_{3}+F_{1}(t_{1})+F_{2}(t_{2})+F_{3}(t_{3}),\tau _{1},\tau _{2}\in \R,$

\noindent
$t_{1},t_{2},t_{3}\in $\gC$^{}.$

\bigskip

\emph{\ Definition 6.9.} \ We consider
\gR$_{7}^{}=\{(D,M,T,A,B,\rho )|D\in $\gD$%
_{4}^{},M\in $\gA$^{},$

\noindent
$T\in $\gJ$_{0}^{},A,B\in $\gJ$%
^{},\rho \in \R\}$ is a vector space of

$\underbrace{\R^{28}}_{\mathfrak{D}_{4}}\oplus \underbrace{\R^{8}\oplus \R^{8}\oplus \R^{8}}_{\mathfrak{A}}\oplus 
\underbrace{\R\oplus \R\oplus \R^{8}\oplus \R^{8}\oplus \R^{8}}_{\mathfrak{J}_{0}}$

$\oplus \underbrace{\R\oplus \R\oplus \R\oplus \R^{8}\oplus \R^{8}\oplus \R^{8}}_{\mathfrak{J}}%
\oplus \underbrace{\R\oplus \R\oplus \R\oplus \R^{8}\oplus \R^{8}\oplus \R^{8}}_{\mathfrak{J}}%
\oplus \R.$

\noindent
We identify $(D,M,T,A,B,\rho )\in $\gR$_{7}^{}$ as an
element of a vector space by 

\gR$_{7}^{}\ni (D,M,T,A,B,\rho )\rightarrow
\textrm{fv}(D,M,T,A,B,\rho )=\left( 
\begin{smallmatrix}
\textrm{fv}(D) \\ 
\textrm{fv}(m_{1}) \\ 
\textrm{fv}(m_{2}) \\ 
\textrm{fv}(m_{3}) \\ 
\tau _{1} \\ 
\tau _{2} \\ 
\textrm{fv}(t_{1}) \\ 
\textrm{fv}(t_{2}) \\ 
\textrm{fv}(t_{3}) \\ 
\textrm{fv}(A) \\ 
\textrm{fv}(B) \\ 
\rho%
\end{smallmatrix}%
\right) \in 
\begin{array}{l}
\R^{78} \\ 
\oplus \R^{27} \\ 
\oplus \R^{27} \\ 
\oplus \R \\ 
=\R^{133}%
\end{array}%
,$

\noindent
where
\begin{align*}
A &=\alpha _{1}E_{1}+\alpha _{2}E_{2}+\alpha _{3}E_{3}+F_{1}(a_{1})+F_{2}(a_{2})+F_{3}(a_{3}), \\
B &=\beta _{1}E_{1}+\beta _{2}E_{2}+\beta _{3}E_{3} +F_{1}(b_{1})+F_{2}(b_{2})+F_{3}(b_{3}), \\
&\alpha _{1},\alpha _{2},\alpha _{3},\beta _{1},\beta _{2},\beta _{3}\in \R,a_{1},a_{2},a_{3},b_{1},b_{2},b_{3}\in \ggC^{}.
\end{align*} 

\bigskip

\emph{Definition 6.10. }\ \ We consider \gR$_{8}^{}=$ \gR$_{7}^{}\oplus $\gP$^{}\oplus $\gP%
$^{}\oplus \R\oplus \R\oplus \R$ is a vector space of

\ \ \ \ $\underbrace{\R^{28}}_{\mathfrak{D}_{4}}\oplus \underbrace{\R^{8}\oplus \R^{8}\oplus \R^{8}}_{\mathfrak{A}}%
\oplus \underbrace{\R\oplus \R\oplus \R^{8}\oplus \R^{8}\oplus \R^{8}}_{\mathfrak{J}_{0}}.$

$\ \ \oplus \underbrace{\R\oplus \R\oplus \R\oplus \R^{8}\oplus \R^{8}\oplus \R^{8}%
}_{\mathfrak{J}}\oplus \underbrace{\R\oplus \R\oplus \R\oplus \R^{8}\oplus \R^{8}\oplus \R^{8}}_{\mathfrak{J}}%
\oplus \R$

$\ \ \oplus \underbrace{\underbrace{\R\oplus \R\oplus \R\oplus \R^{8}\oplus \R^{8}\oplus \R^{8}%
}_{\mathfrak{J}}\oplus \underbrace{\R\oplus \R\oplus \R\oplus \R^{8}\oplus \R^{8}\oplus \R^{8}}_{\mathfrak{J}}%
\oplus \R\oplus \R}_{\mathfrak{P}}$

$\ \ \oplus \underbrace{\underbrace{\R\oplus \R\oplus \R\oplus \R^{8}\oplus \R^{8}\oplus \R^{8}%
}_{\mathfrak{J}}\oplus \underbrace{\R\oplus \R\oplus \R\oplus \R^{8}\oplus \R^{8}\oplus \R^{8}}_{\mathfrak{J}}%
\oplus \R\oplus \R}_{\mathfrak{P}}$

$\ \ \oplus \R\oplus \R\oplus \R.$

We identify $(\Phi ,P,Q,r,s,u)\in $\gR$_{8}^{}$ $(\Phi
\in $\gR$_{7}^{})$ as an element of a vector space by

\gR$_{8}^{}\ni (\Phi ,P,Q,r,s,u)\rightarrow \textrm{fv}(\Phi
,P,Q,r,s,u)=\left( 
\begin{smallmatrix}
\textrm{fv}(\Phi ) \\ 
\textrm{fv}(X) \\ 
\textrm{fv}(Y) \\ 
\xi \\ 
\eta \\ 
\textrm{fv}(Z) \\ 
\textrm{fv}(W) \\ 
\zeta \\ 
\omega \\ 
r \\ 
s \\ 
u%
\end{smallmatrix}%
\right) \in 
\begin{array}{l}
\R^{133} \\ 
\oplus \R^{27} \\ 
\oplus \R^{27} \\ 
\oplus \R\oplus \R \\ 
\oplus \R^{27} \\ 
\oplus \R^{27} \\ 
\oplus \R\oplus \R \\ 
\oplus \R\oplus \R\oplus \R \\ 
=\R^{248}%
\end{array}%
,$

\noindent
where
\begin{align*}
P&=(X,Y,\xi ,\eta ),Q=(Z,W,\zeta ,\omega ),\\
X&=\chi _{1}E_{1}+\chi _{2}E_{2}+\chi _{3}E_{3}+F_{1}(x_{1})+F_{2}(x_{2})+F_{3}(x_{3}),\\
&\chi _{1},\chi _{2},\chi _{3}\in \R,x_{1},x_{2},x_{3}\in \ggC^{},\\
Y&=\gamma _{1}E_{1}+\gamma _{2}E_{2}+\gamma
_{3}E_{3}+F_{1}(y_{1})+F_{2}(y_{2})+F_{3}(y_{3}),\\
&\gamma _{1},\gamma _{2},\gamma _{3}\in
\R,y_{1},y_{2},y_{3}\in \ggC^{},\\
Z&=\mu _{1}E_{1}+\mu _{2}E_{2}+\mu _{3}E_{3}+F_{1}(z_{1})+F_{2}(z_{2})+F_{3}(z_{3}),\\
&\mu _{1},\mu _{2},\mu _{3}\in \R,z_{1},z_{2},z_{3}\in \ggC^{},\\
W&=\psi _{1}E_{1}+\psi _{2}E_{2}+\psi _{3}E_{3}+F_{1}(w_{1})+F_{2}(w_{2})+F_{3}(w_{3}),\\
&\psi _{1},\psi _{2},\psi _{3}\in \R,w_{1},w_{2},w_{3}\in \ggC^{},\ \xi ,\eta ,\zeta ,\omega ,r,s,u\in \R.
\end{align*}

\bigskip

\emph{Definition 6.11. } For $D,D_{1},D_{2}\in $\gD$_{4}^{}$,$m,m_{1},m_{2}\in $\gJ$^{}$, $\alpha ,\beta \in \R,$\ We define
$\R$-linear mappings as follows:
\begin{align*}
\textrm{f}_{DDD}&:\ggD_{4}^{}\times \ggD_{4}^{}\rightarrow \ggD_{4}^{} %
&\textrm{f}_{DDD}(D_{1},D_{2})&=[D_{1},D_{2}],\\
& \textrm{f}_{DDD}(D_{1},?):\ggD_{4}^{}\rightarrow \ggD_{4}^{} \ \ %
&\textrm{f}_{DDD}(D_{1},?)D_{2}&=[D_{1},D_{2}],\\
\textrm{f}_{DJJ}&:\ggD_{4}^{}\times \ggC%
^{}\rightarrow \ggC^{}
&\textrm{f}_{DJJ}(D,m)&=\textrm{g}_{d}(D)m,\\
&\textrm{f}_{DJJ}(D,?):\ggC%
^{}\rightarrow \ggC^{}%
&\textrm{f}_{DJJ}(D,?)m&=\textrm{g}_{d}(D)m,\\
&\textrm{f}_{DJJ}(?,m):\ggD_{4}^{}\rightarrow \ggC^{}%
&\textrm{f}_{DJJ}(?,m)D&=\textrm{g}_{d}(D)m,\\
\textrm{f}_{JJD}&:\ggC^{}\times \ggC^{}\rightarrow \ggD_{4}^{}%
&\textrm{f}_{JJD}(m_{1},m_{2})&=\textrm{JD}(m_{1},m_{2}),\\
&\textrm{f}_{JJD}(m_{1},?):\ggC%
^{}\rightarrow \ggD_{4}^{}%
&\textrm{f}_{JJD}(m_{1},?)m_{2}&=\textrm{JD}(m_{1},m_{2}),\\
&\textrm{f}_{JJD}(?,m_{2}):\ggC%
^{}\rightarrow \ggD_{4}^{}%
&\textrm{f}_{JJD}(?,m_{2})m_{1}&=\textrm{JD}(m_{1},m_{2}),\\
\textrm{f}_{JJJ}&:\ggC^{}\times \ggC^{}\rightarrow \ggC^{}%
&\textrm{f}_{JJJ}(m_{1},m_{2})&=m_{1}m_{2},\\
&\textrm{f}_{JJJ}(m_{1},?):\ggC%
^{}\rightarrow \ggC^{}%
&\textrm{f}_{JJJ}(m_{1},?)m_{2}&=m_{1}m_{2},\\
&\textrm{f}_{JJJ}(?,m_{2}):\ggC%
^{}\rightarrow \ggC^{}%
&\textrm{f}_{JJJ}(?,m_{2})m_{1}&=m_{1}m_{2},\\
\textrm{f}_{iJJ}&:\ggC^{}\times \ggC^{}\rightarrow \ggC^{}%
&\textrm{f}_{iJJ}(m_{1},m_{2})&=\overline{m_{1}m_{2}},\\
&\textrm{f}_{iJJ}(m_{1},?):\ggC%
^{}\rightarrow \ggC^{}%
&\textrm{f}_{iJJ}(m_{1},?)m_{2}&=\overline{m_{1}m_{2}},\\
&\textrm{f}_{iJJ}(?,m_{2}):\ggC%
^{}\rightarrow \ggC^{}%
&\textrm{f}_{iJJ}(?,m_{2})m_{1}&=\overline{m_{1}m_{2}},
\end{align*}
\begin{align*}
\textrm{f}_{JJC}&:\ggC^{}\times \ggC^{}\rightarrow \R
&\textrm{f}_{JJC}(m_{1},m_{2})&=(m_{1},m_{2}),\\
&\textrm{f}_{JJC}(m_{1},?):\ggC%
^{}\rightarrow \R  %
&\textrm{f}_{JJC}(m_{1},?)m_{2}&=(m_{1},m_{2}),\\
&\textrm{f}_{JJC}(?,m_{2}):\ggC%
^{}\rightarrow \R %
&\textrm{f}_{JJC}(?,m_{2})m_{1}&=(m_{1},m_{2}),\\
\textrm{f}_{CJJ}&:\R\times \ggC^{}\rightarrow \ggC^{}%
&\textrm{f}_{CJJ}(\alpha ,m)&=\alpha \ast m,\\
&\textrm{f}_{CJJ}(\alpha ,?):\ggC%
^{}\rightarrow \ggC^{} %
&\textrm{f}_{CJJ}(\alpha,?)m&=\alpha \ast m,\\
&\textrm{f}_{CJJ}(?,m):\R\rightarrow \ggC^{}
&\textrm{f}_{CJJ}(?,m)\alpha& =\alpha \ast m,\\
\textrm{f}_{CCC}&:\R\times \R\rightarrow \R%
&\textrm{f}_{CCC}(\alpha ,\beta )&=\alpha \ast \beta ,\\
&\textrm{f}_{CCC}(\alpha ,?):\R\rightarrow \R%
&\textrm{f}_{CCC}(\alpha ,?)\beta &=\alpha \ast \beta ,\\
&\textrm{f}_{CCC}(?,\beta ):\R\rightarrow \R%
&\textrm{f}_{CCC}(?,\beta )\alpha& =\alpha \ast \beta .
\end{align*}

\noindent
where $\ast $ means commutative multiplication operator.

\bigskip

Now let's express $\textrm{f}_{DDD}(D,?)$ as a matrix $\textrm{LD}(D)\in M(28\times 28,\R)$.

\bigskip

\emph{Lemma 6.12. }For $D=\sum\limits_{0\leq i<j\leq 7}d_{ij}D_{ij} (
d_{ij}\in \R),$ we have $\textrm{f}_{DDD}(D,?)$ is expressed by

$\textrm{LD}(D)=$

\begin{flushleft}
{\fontsize{4pt}{8pt} \selectfont%
$\left( 
\begin{tabular}{@{}c@{}c@{}c@{}c@{}c@{}c@{}c@{}c@{}c@{}c@{}c@{}c@{}c@{}c@{}c@{}c@{}c@{}c@{}c@{}c@{}c@{}c@{}c@{}c@{}c@{}c@{}c@{}c@{}}
0 & $d_{\scalebox{0.8}{12}}$ & $d_{\scalebox{0.8}{13}}$ & $d_{\scalebox{0.8}{14}}$ & $d_{\scalebox{0.8}{15}}$ & $d_{\scalebox{0.8}{16}}$ & $d_{\scalebox{0.8}{17}}$ & $%
\mathchar`-d_{\scalebox{0.8}{02}}$ & $\mathchar`-d_{\scalebox{0.8}{03}}$ & $\mathchar`-d_{\scalebox{0.8}{04}}$ & $\mathchar`-d_{\scalebox{0.8}{05}}$ & $\mathchar`-d_{\scalebox{0.8}{06}}$ & $\mathchar`-d_{\scalebox{0.8}{07}}$ & 0 & 0
& 0 & 0 & 0 & 0 & 0 & 0 & 0 & 0 & 0 & 0 & 0 & 0 & 0 \\ 
$\mathchar`-d_{\scalebox{0.8}{12}}$ & 0 & $d_{\scalebox{0.8}{23}}$ & $d_{\scalebox{0.8}{24}}$ & $d_{\scalebox{0.8}{25}}$ & $d_{\scalebox{0.8}{26}}$ & $d_{\scalebox{0.8}{27}}$ & $%
d_{\scalebox{0.8}{01}}$ & 0 & 0 & 0 & 0 & 0 & $\mathchar`-d_{\scalebox{0.8}{03}}$ & $\mathchar`-d_{\scalebox{0.8}{04}}$ & $\mathchar`-d_{\scalebox{0.8}{05}}$ & $\mathchar`-d_{\scalebox{0.8}{06}}$
& $\mathchar`-d_{\scalebox{0.8}{07}}$ & 0 & 0 & 0 & 0 & 0 & 0 & 0 & 0 & 0 & 0 \\ 
$\mathchar`-d_{\scalebox{0.8}{13}}$ & $\mathchar`-d_{\scalebox{0.8}{23}}$ & 0 & $d_{\scalebox{0.8}{34}}$ & $d_{\scalebox{0.8}{35}}$ & $d_{\scalebox{0.8}{36}}$ & $d_{\scalebox{0.8}{37}}$ & 0 & $%
d_{\scalebox{0.8}{01}}$ & 0 & 0 & 0 & 0 & $d_{\scalebox{0.8}{02}}$ & 0 & 0 & 0 & 0 & $\mathchar`-d_{\scalebox{0.8}{04}}$ & $\mathchar`-d_{\scalebox{0.8}{05}}$ & 
$\mathchar`-d_{\scalebox{0.8}{06}}$ & $\mathchar`-d_{\scalebox{0.8}{07}}$ & 0 & 0 & 0 & 0 & 0 & 0 \\ 
$\mathchar`-d_{\scalebox{0.8}{14}}$ & $\mathchar`-d_{\scalebox{0.8}{24}}$ & $\mathchar`-d_{\scalebox{0.8}{34}}$ & 0 & $d_{\scalebox{0.8}{45}}$ & $d_{\scalebox{0.8}{46}}$ & $d_{\scalebox{0.8}{47}}$ & 0 & 
00 & $d_{\scalebox{0.8}{01}}$ & 0 & 0 & 0 & 0 & $d_{\scalebox{0.8}{02}}$ & 0 & 0 & 0 & $d_{\scalebox{0.8}{03}}$ & 0 & 0 & 0
& $\mathchar`-d_{\scalebox{0.8}{05}}$ & $\mathchar`-d_{\scalebox{0.8}{06}}$ & $\mathchar`-d_{\scalebox{0.8}{07}}$ & 0 & 0 & 0 \\ 
$\mathchar`-d_{\scalebox{0.8}{15}}$ & $\mathchar`-d_{\scalebox{0.8}{25}}$ & $\mathchar`-d_{\scalebox{0.8}{35}}$ & $\mathchar`-d_{\scalebox{0.8}{45}}$ & 0 & $d_{\scalebox{0.8}{56}}$ & $d_{\scalebox{0.8}{57}}$ & 0
& 0 & 0 & $d_{\scalebox{0.8}{01}}$ & 0 & 0 & 0 & 0 & $d_{\scalebox{0.8}{02}}$ & 0 & 0 & 0 & $d_{\scalebox{0.8}{03}}$ & 0 & 0
& $d_{\scalebox{0.8}{04}}$ & 0 & 0 & $\mathchar`-d_{\scalebox{0.8}{06}}$ & $\mathchar`-d_{\scalebox{0.8}{07}}$ & 0 \\ 
$\mathchar`-d_{\scalebox{0.8}{16}}$ & $\mathchar`-d_{\scalebox{0.8}{26}}$ & $\mathchar`-d_{\scalebox{0.8}{36}}$ & $\mathchar`-d_{\scalebox{0.8}{46}}$ & $\mathchar`-d_{\scalebox{0.8}{56}}$ & 0 & $d_{\scalebox{0.8}{67}}$ & 0
& 0 & 0 & 0 & $d_{\scalebox{0.8}{01}}$ & 0 & 0 & 0 & 0 & $d_{\scalebox{0.8}{02}}$ & 0 & 0 & 0 & $d_{\scalebox{0.8}{03}}$ & 0
& 0 & $d_{\scalebox{0.8}{04}}$ & 0 & $d_{\scalebox{0.8}{05}}$ & 0 & $\mathchar`-d_{\scalebox{0.8}{07}}$ \\ 
$\mathchar`-d_{\scalebox{0.8}{17}}$ & $\mathchar`-d_{\scalebox{0.8}{27}}$ & $\mathchar`-d_{\scalebox{0.8}{37}}$ & $\mathchar`-d_{\scalebox{0.8}{47}}$ & $\mathchar`-d_{\scalebox{0.8}{57}}$ & $\mathchar`-d_{\scalebox{0.8}{67}}$ & 0 & 0
& 0 & 0 & 0 & 0 & $d_{\scalebox{0.8}{01}}$ & 0 & 0 & 0 & 0 & $d_{\scalebox{0.8}{02}}$ & 0 & 0 & 0 & $d_{\scalebox{0.8}{03}}$
& 0 & 0 & $d_{\scalebox{0.8}{04}}$ & 0 & $d_{\scalebox{0.8}{05}}$ & $d_{\scalebox{0.8}{06}}$ \\ 
$d_{\scalebox{0.8}{02}}$ & $\mathchar`-d_{\scalebox{0.8}{01}}$ & 0 & 0 & 0 & 0 & 0 & 0 & $d_{\scalebox{0.8}{23}}$ & $d_{\scalebox{0.8}{24}}$ & $d_{\scalebox{0.8}{25}}$
& $d_{\scalebox{0.8}{26}}$ & $d_{\scalebox{0.8}{27}}$ & $\mathchar`-d_{\scalebox{0.8}{13}}$ & $\mathchar`-d_{\scalebox{0.8}{14}}$ & $\mathchar`-d_{\scalebox{0.8}{15}}$ & $\mathchar`-d_{\scalebox{0.8}{16}}$ & $%
\mathchar`-d_{\scalebox{0.8}{17}}$ & 0 & 0 & 0 & 0 & 0 & 0 & 0 & 0 & 0 & 0 \\ 
$d_{\scalebox{0.8}{03}}$ & 0 & $\mathchar`-d_{\scalebox{0.8}{01}}$ & 0 & 0 & 0 & 0 & $\mathchar`-d_{\scalebox{0.8}{23}}$ & 0 & $d_{\scalebox{0.8}{34}}$ & $d_{\scalebox{0.8}{35}}
$ & $d_{\scalebox{0.8}{36}}$ & $d_{\scalebox{0.8}{37}}$ & $d_{\scalebox{0.8}{12}}$ & 0 & 0 & 0 & 0 & $\mathchar`-d_{\scalebox{0.8}{14}}$ & $\mathchar`-d_{\scalebox{0.8}{15}}$ & 
$\mathchar`-d_{\scalebox{0.8}{16}}$ & $\mathchar`-d_{\scalebox{0.8}{17}}$ & 0 & 0 & 0 & 0 & 0 & 0 \\ 
$d_{\scalebox{0.8}{04}}$ & 0 & 0 & $\mathchar`-d_{\scalebox{0.8}{01}}$ & 0 & 0 & 0 & $\mathchar`-d_{\scalebox{0.8}{24}}$ & $\mathchar`-d_{\scalebox{0.8}{34}}$ & 0 & $%
d_{\scalebox{0.8}{45}}$ & $d_{\scalebox{0.8}{46}}$ & $d_{\scalebox{0.8}{47}}$ & 0 & $d_{\scalebox{0.8}{12}}$ & 0 & 0 & 0 & $d_{\scalebox{0.8}{13}}$ & 0 & 0
& 0 & $\mathchar`-d_{\scalebox{0.8}{15}}$ & $\mathchar`-d_{\scalebox{0.8}{16}}$ & $\mathchar`-d_{\scalebox{0.8}{17}}$ & 0 & 0 & 0 \\ 
$d_{\scalebox{0.8}{05}}$ & 0 & 0 & 0 & $\mathchar`-d_{\scalebox{0.8}{01}}$ & 0 & 0 & $\mathchar`-d_{\scalebox{0.8}{25}}$ & $\mathchar`-d_{\scalebox{0.8}{35}}$ & $\mathchar`-d_{\scalebox{0.8}{45}}$
& 0 & $d_{\scalebox{0.8}{56}}$ & $d_{\scalebox{0.8}{57}}$ & 0 & 0 & $d_{\scalebox{0.8}{12}}$ & 0 & 0 & 0 & $d_{\scalebox{0.8}{13}}$ & 0 & 0
& $d_{\scalebox{0.8}{14}}$ & 0 & 0 & $\mathchar`-d_{\scalebox{0.8}{16}}$ & $\mathchar`-d_{\scalebox{0.8}{17}}$ & 0 \\ 
$d_{\scalebox{0.8}{06}}$ & 0 & 0 & 0 & 0 & $\mathchar`-d_{\scalebox{0.8}{01}}$ & 0 & $\mathchar`-d_{\scalebox{0.8}{26}}$ & $\mathchar`-d_{\scalebox{0.8}{36}}$ & $\mathchar`-d_{\scalebox{0.8}{46}}$
& $\mathchar`-d_{\scalebox{0.8}{56}}$ & 0 & $d_{\scalebox{0.8}{67}}$ & 0 & 0 & 0 & $d_{\scalebox{0.8}{12}}$ & 0 & 0 & 0 & $d_{\scalebox{0.8}{13}}$ & 0
& 0 & $d_{\scalebox{0.8}{14}}$ & 0 & $d_{\scalebox{0.8}{15}}$ & 0 & $\mathchar`-d_{\scalebox{0.8}{17}}$ \\ 
$d_{\scalebox{0.8}{07}}$ & 0 & 0 & 0 & 0 & 0 & $\mathchar`-d_{\scalebox{0.8}{01}}$ & $\mathchar`-d_{\scalebox{0.8}{27}}$ & $\mathchar`-d_{\scalebox{0.8}{37}}$ & $\mathchar`-d_{\scalebox{0.8}{47}}$
& $\mathchar`-d_{\scalebox{0.8}{57}}$ & $\mathchar`-d_{\scalebox{0.8}{67}}$ & 0 & 0 & 0 & 0 & 0 & $d_{\scalebox{0.8}{12}}$ & 0 & 0 & 0 & $d_{\scalebox{0.8}{13}}$
& 0 & 0 & $d_{\scalebox{0.8}{14}}$ & 0 & $d_{\scalebox{0.8}{15}}$ & $d_{\scalebox{0.8}{16}}$ \\ 
0 & $d_{\scalebox{0.8}{03}}$ & $\mathchar`-d_{\scalebox{0.8}{02}}$ & 0 & 0 & 0 & 0 & $d_{\scalebox{0.8}{13}}$ & $\mathchar`-d_{\scalebox{0.8}{12}}$ & 0 & 0 & 0
& 0 & 0 & $d_{\scalebox{0.8}{34}}$ & $d_{\scalebox{0.8}{35}}$ & $d_{\scalebox{0.8}{36}}$ & $d_{\scalebox{0.8}{37}}$ & $\mathchar`-d_{\scalebox{0.8}{24}}$ & $\mathchar`-d_{\scalebox{0.8}{25}}$
& $\mathchar`-d_{\scalebox{0.8}{26}}$ & $\mathchar`-d_{\scalebox{0.8}{27}}$ & 0 & 0 & 0 & 0 & 0 & 0 \\ 
0 & $d_{\scalebox{0.8}{04}}$ & 0 & $\mathchar`-d_{\scalebox{0.8}{02}}$ & 0 & 0 & 0 & $d_{\scalebox{0.8}{14}}$ & 0 & $\mathchar`-d_{\scalebox{0.8}{12}}$ & 0 & 0
& 0 & $\mathchar`-d_{\scalebox{0.8}{34}}$ & 0 & $d_{\scalebox{0.8}{45}}$ & $d_{\scalebox{0.8}{46}}$ & $d_{\scalebox{0.8}{47}}$ & $d_{\scalebox{0.8}{23}}$ & 0 & 0 & 0
& $\mathchar`-d_{\scalebox{0.8}{25}}$ & $\mathchar`-d_{\scalebox{0.8}{26}}$ & $\mathchar`-d_{\scalebox{0.8}{27}}$ & 0 & 0 & 0 \\ 
0 & $d_{\scalebox{0.8}{05}}$ & 0 & 0 & $\mathchar`-d_{\scalebox{0.8}{02}}$ & 0 & 0 & $d_{\scalebox{0.8}{15}}$ & 0 & 0 & $\mathchar`-d_{\scalebox{0.8}{12}}$ & 0
& 0 & $\mathchar`-d_{\scalebox{0.8}{35}}$ & $\mathchar`-d_{\scalebox{0.8}{45}}$ & 0 & $d_{\scalebox{0.8}{56}}$ & $d_{\scalebox{0.8}{57}}$ & 0 & $d_{\scalebox{0.8}{23}}$ & 0 & 0
& $d_{\scalebox{0.8}{24}}$ & 0 & 0 & $\mathchar`-d_{\scalebox{0.8}{26}}$ & $\mathchar`-d_{\scalebox{0.8}{27}}$ & 0 \\ 
0 & $d_{\scalebox{0.8}{06}}$ & 0 & 0 & 0 & $\mathchar`-d_{\scalebox{0.8}{02}}$ & 0 & $d_{\scalebox{0.8}{16}}$ & 0 & 0 & 0 & $\mathchar`-d_{\scalebox{0.8}{12}}$
& 0 & $\mathchar`-d_{\scalebox{0.8}{36}}$ & $\mathchar`-d_{\scalebox{0.8}{46}}$ & $\mathchar`-d_{\scalebox{0.8}{56}}$ & 0 & $d_{\scalebox{0.8}{67}}$ & 0 & 0 & $d_{\scalebox{0.8}{23}}$ & 0
& 0 & $d_{\scalebox{0.8}{24}}$ & 0 & $d_{\scalebox{0.8}{25}}$ & 0 & $\mathchar`-d_{\scalebox{0.8}{27}}$ \\ 
0 & $d_{\scalebox{0.8}{07}}$ & 0 & 0 & 0 & 0 & $\mathchar`-d_{\scalebox{0.8}{02}}$ & $d_{\scalebox{0.8}{17}}$ & 0 & 0 & 0 & 0 & $%
\mathchar`-d_{\scalebox{0.8}{12}}$ & $\mathchar`-d_{\scalebox{0.8}{37}}$ & $\mathchar`-d_{\scalebox{0.8}{47}}$ & $\mathchar`-d_{\scalebox{0.8}{57}}$ & $\mathchar`-d_{\scalebox{0.8}{67}}$ & 0 & 0 & 0 & 0 & $%
d_{\scalebox{0.8}{23}}$ & 0 & 0 & $d_{\scalebox{0.8}{24}}$ & 0 & $d_{\scalebox{0.8}{25}}$ & $d_{\scalebox{0.8}{26}}$ \\ 
0 & 0 & $d_{\scalebox{0.8}{04}}$ & $\mathchar`-d_{\scalebox{0.8}{03}}$ & 0 & 0 & 0 & 0 & $d_{\scalebox{0.8}{14}}$ & $\mathchar`-d_{\scalebox{0.8}{13}}$ & 0 & 0
& 0 & $d_{\scalebox{0.8}{24}}$ & $\mathchar`-d_{\scalebox{0.8}{23}}$ & 0 & 0 & 0 & 0 & $d_{\scalebox{0.8}{45}}$ & $d_{\scalebox{0.8}{46}}$ & $d_{\scalebox{0.8}{47}}$
& $\mathchar`-d_{\scalebox{0.8}{35}}$ & $\mathchar`-d_{\scalebox{0.8}{36}}$ & $\mathchar`-d_{\scalebox{0.8}{37}}$ & 0 & 0 & 0 \\ 
0 & 0 & $d_{\scalebox{0.8}{05}}$ & 0 & $\mathchar`-d_{\scalebox{0.8}{03}}$ & 0 & 0 & 0 & $d_{\scalebox{0.8}{15}}$ & 0 & $\mathchar`-d_{\scalebox{0.8}{13}}$ & 0
& 0 & $d_{\scalebox{0.8}{25}}$ & 0 & $\mathchar`-d_{\scalebox{0.8}{23}}$ & 0 & 0 & $\mathchar`-d_{\scalebox{0.8}{45}}$ & 0 & $d_{\scalebox{0.8}{56}}$ & $d_{\scalebox{0.8}{57}}$
& $d_{\scalebox{0.8}{34}}$ & 0 & 0 & $\mathchar`-d_{\scalebox{0.8}{36}}$ & $\mathchar`-d_{\scalebox{0.8}{37}}$ & 0 \\ 
0 & 0 & $d_{\scalebox{0.8}{06}}$ & 0 & 0 & $\mathchar`-d_{\scalebox{0.8}{03}}$ & 0 & 0 & $d_{\scalebox{0.8}{16}}$ & 0 & 0 & $\mathchar`-d_{\scalebox{0.8}{13}}$
& 0 & $d_{\scalebox{0.8}{26}}$ & 0 & 0 & $\mathchar`-d_{\scalebox{0.8}{23}}$ & 0 & $\mathchar`-d_{\scalebox{0.8}{46}}$ & $\mathchar`-d_{\scalebox{0.8}{56}}$ & 0 & $d_{\scalebox{0.8}{67}}$
& 0 & $d_{\scalebox{0.8}{34}}$ & 0 & $d_{\scalebox{0.8}{35}}$ & 0 & $\mathchar`-d_{\scalebox{0.8}{37}}$ \\ 
0 & 0 & $d_{\scalebox{0.8}{07}}$ & 0 & 0 & 0 & $\mathchar`-d_{\scalebox{0.8}{03}}$ & 0 & $d_{\scalebox{0.8}{17}}$ & 0 & 0 & 0 & $%
\mathchar`-d_{\scalebox{0.8}{13}}$ & $d_{\scalebox{0.8}{27}}$ & 0 & 0 & 0 & $\mathchar`-d_{\scalebox{0.8}{23}}$ & $\mathchar`-d_{\scalebox{0.8}{47}}$ & $\mathchar`-d_{\scalebox{0.8}{57}}$ & $%
\mathchar`-d_{\scalebox{0.8}{67}}$ & 0 & 0 & 0 & $d_{\scalebox{0.8}{34}}$ & 0 & $d_{\scalebox{0.8}{35}}$ & $d_{\scalebox{0.8}{36}}$ \\ 
0 & 0 & 0 & $d_{\scalebox{0.8}{05}}$ & $\mathchar`-d_{\scalebox{0.8}{04}}$ & 0 & 0 & 0 & 0 & $d_{\scalebox{0.8}{15}}$ & $\mathchar`-d_{\scalebox{0.8}{14}}$ & 0
& 0 & 0 & $d_{\scalebox{0.8}{25}}$ & $\mathchar`-d_{\scalebox{0.8}{24}}$ & 0 & 0 & $d_{\scalebox{0.8}{35}}$ & $\mathchar`-d_{\scalebox{0.8}{34}}$ & 0 & 0 & 0 & $%
d_{\scalebox{0.8}{56}}$ & $d_{\scalebox{0.8}{57}}$ & $\mathchar`-d_{\scalebox{0.8}{46}}$ & $\mathchar`-d_{\scalebox{0.8}{47}}$ & 0 \\ 
0 & 0 & 0 & $d_{\scalebox{0.8}{06}}$ & 0 & $\mathchar`-d_{\scalebox{0.8}{04}}$ & 0 & 0 & 0 & $d_{\scalebox{0.8}{16}}$ & 0 & $\mathchar`-d_{\scalebox{0.8}{14}}$
& 0 & 0 & $d_{\scalebox{0.8}{26}}$ & 0 & $\mathchar`-d_{\scalebox{0.8}{24}}$ & 0 & $d_{\scalebox{0.8}{36}}$ & 0 & $\mathchar`-d_{\scalebox{0.8}{34}}$ & 0 & $%
\mathchar`-d_{\scalebox{0.8}{56}}$ & 0 & $d_{\scalebox{0.8}{67}}$ & $d_{\scalebox{0.8}{45}}$ & 0 & $\mathchar`-d_{\scalebox{0.8}{47}}$ \\ 
0 & 0 & 0 & $d_{\scalebox{0.8}{07}}$ & 0 & 0 & $\mathchar`-d_{\scalebox{0.8}{04}}$ & 0 & 0 & $d_{\scalebox{0.8}{17}}$ & 0 & 0 & $%
\mathchar`-d_{\scalebox{0.8}{14}}$ & 0 & $d_{\scalebox{0.8}{27}}$ & 0 & 0 & $\mathchar`-d_{\scalebox{0.8}{24}}$ & $d_{\scalebox{0.8}{37}}$ & 0 & 0 & $\mathchar`-d_{\scalebox{0.8}{34}}$ & 
$\mathchar`-d_{\scalebox{0.8}{57}}$ & $\mathchar`-d_{\scalebox{0.8}{67}}$ & 0 & 0 & $d_{\scalebox{0.8}{45}}$ & $d_{\scalebox{0.8}{46}}$ \\ 
0 & 0 & 0 & 0 & $d_{\scalebox{0.8}{06}}$ & $\mathchar`-d_{\scalebox{0.8}{05}}$ & 0 & 0 & 0 & 0 & $d_{\scalebox{0.8}{16}}$ & $\mathchar`-d_{\scalebox{0.8}{15}}$
& 0 & 0 & 0 & $d_{\scalebox{0.8}{26}}$ & $\mathchar`-d_{\scalebox{0.8}{25}}$ & 0 & 0 & $d_{\scalebox{0.8}{36}}$ & $\mathchar`-d_{\scalebox{0.8}{35}}$ & 0 & $%
d_{\scalebox{0.8}{46}}$ & 0 & 0 & 0 & $d_{\scalebox{0.8}{67}}$ & $\mathchar`-d_{\scalebox{0.8}{57}}$ \\ 
0 & 0 & 0 & 0 & $d_{\scalebox{0.8}{07}}$ & 0 & $\mathchar`-d_{\scalebox{0.8}{05}}$ & 0 & 0 & 0 & $d_{\scalebox{0.8}{17}}$ & 0 & $%
\mathchar`-d_{\scalebox{0.8}{15}}$ & 0 & 0 & $d_{\scalebox{0.8}{27}}$ & 0 & $\mathchar`-d_{\scalebox{0.8}{25}}$ & 0 & $d_{\scalebox{0.8}{37}}$ & 0 & $\mathchar`-d_{\scalebox{0.8}{35}}$ & 
$d_{\scalebox{0.8}{47}}$ & 0 & $\mathchar`-d_{\scalebox{0.8}{45}}$ & $\mathchar`-d_{\scalebox{0.8}{67}}$ & 0 & $d_{\scalebox{0.8}{56}}$ \\ 
0 & 0 & 0 & 0 & 0 & $d_{\scalebox{0.8}{07}}$ & $\mathchar`-d_{\scalebox{0.8}{06}}$ & 0 & 0 & 0 & 0 & $d_{\scalebox{0.8}{17}}$ & $%
\mathchar`-d_{\scalebox{0.8}{16}}$ & 0 & 0 & 0 & $d_{\scalebox{0.8}{27}}$ & $\mathchar`-d_{\scalebox{0.8}{26}}$ & 0 & 0 & $d_{\scalebox{0.8}{37}}$ & $\mathchar`-d_{\scalebox{0.8}{36}}$ & 
0 & $d_{\scalebox{0.8}{47}}$ & $\mathchar`-d_{\scalebox{0.8}{46}}$ & $d_{\scalebox{0.8}{57}}$ & $\mathchar`-d_{\scalebox{0.8}{56}}$ & 0%
\end{tabular}%
\right) $ 
}
\end{flushleft}
\bigskip

\emph{Proof. }Using\emph{\ Lemma 2.4 \ }we calculate\emph{\ }as
follows,

$[D_{0j},D_{jk}]=D_{0k},0<j\neq k\leq 7,$

$[D_{1j},D_{jk}]=D_{1k},1<j\leq 7,0\leq k\leq 7,j\neq k$ ,

$[D_{2j},D_{jk}]=D_{2k},2<j\leq 7,0\leq k\leq 7,j\neq k$ ,

$[D_{3j},D_{jk}]=D_{3k},3<j\leq 7,0\leq k\leq 7,j\neq k$ ,

$[D_{4j},D_{jk}]=D_{4k},4<j\leq 7,0\leq k\leq 7,j\neq k$ ,

$[D_{5j},D_{jk}]=D_{5k},5<j\leq 7,0\leq k\leq 7,j\neq k$ ,

$[D_{67},D_{7k}]=D_{6k},0\leq k\leq 5,$

\noindent
Hence we have the above expression. \ \ \ \ \emph{Q.E.D.}

\bigskip

\emph{Lemma 6.13. \ }For $m=\sum\limits_{i=0}^{7}m_{i}e_{i} (m%
_{i}\in \R)$ $,$ let's express $\textrm{f}_{DJJ}(?,m)$ as a matrix
$\textrm{MJr}(m)\in M(8\times 28,\R)$. We have $\textrm{f}_{DJJ}(?,m)$ is expressed by

$\textrm{MJr}(m)=$
\begin{flushleft}
{\fontsize{4pt}{10pt} \selectfont%
$\left( 
\begin{tabular}{@{}c@{}c@{}c@{}c@{}c@{}c@{}c@{}c@{}c@{}c@{}c@{}c@{}c@{}c@{}c@{}c@{}c@{}c@{}c@{}c@{}c@{}c@{}c@{}c@{}c@{}c@{}c@{}c@{}}
$m_{\scalebox{0.8}{1}}$ & $m_{\scalebox{0.8}{2}}$ & $m_{\scalebox{0.8}{3}}$ & $m_{\scalebox{0.8}{4}}$ & $m_{\scalebox{0.8}{5}}$ & $m_{\scalebox{0.8}{6}}$ & $m_{\scalebox{0.8}{7}}$ & $0$ & $%
0$ & $0$ & $0$ & $0$ & $0$ & $0$ & $0$ & $0$ & $0$ & $0$ & $0$ & $0$ & $0$ & 
$0$ & $0$ & $0$ & $0$ & $0$ & $0$ & $0$ \\ 
$\mathchar`-m_{\scalebox{0.8}{0}}$ & $0$ & $0$ & $0$ & $0$ & $0$ & $0$ & $m_{\scalebox{0.8}{2}}$ & $m_{\scalebox{0.8}{3}}$ & $m_{\scalebox{0.8}{4}}$ & 
$m_{\scalebox{0.8}{5}}$ & $m_{\scalebox{0.8}{6}}$ & $m_{\scalebox{0.8}{7}}$ & $0$ & $0$ & $0$ & $0$ & $0$ & $0$ & $0$ & $0$
& $0$ & $0$ & $0$ & $0$ & $0$ & $0$ & $0$ \\ 
$0$ & $\mathchar`-m_{\scalebox{0.8}{0}}$ & $0$ & $0$ & $0$ & $0$ & $0$ & $\mathchar`-m_{\scalebox{0.8}{1}}$ & $0$ & $0$ & $0$ & $%
0$ & $0$ & $m_{\scalebox{0.8}{3}}$ & $m_{\scalebox{0.8}{4}}$ & $m_{\scalebox{0.8}{5}}$ & $m_{\scalebox{0.8}{6}}$ & $m_{\scalebox{0.8}{7}}$ & $0$ & $0$ & $0$
& $0$ & $0$ & $0$ & $0$ & $0$ & $0$ & $0$ \\ 
$0$ & $0$ & $\mathchar`-m_{\scalebox{0.8}{0}}$ & $0$ & $0$ & $0$ & $0$ & $0$ & $\mathchar`-m_{\scalebox{0.8}{1}}$ & $0$ & $0$ & $%
0$ & $0$ & $\mathchar`-m_{\scalebox{0.8}{2}}$ & $0$ & $0$ & $0$ & $0$ & $m_{\scalebox{0.8}{4}}$ & $m_{\scalebox{0.8}{5}}$ & $m_{\scalebox{0.8}{6}}$ & $%
m_{\scalebox{0.8}{7}}$ & $0$ & $0$ & $0$ & $0$ & $0$ & $0$ \\ 
$0$ & $0$ & $0$ & $\mathchar`-m_{\scalebox{0.8}{0}}$ & $0$ & $0$ & $0$ & $0$ & $0$ & $\mathchar`-m_{\scalebox{0.8}{1}}$ & $0$ & $%
0$ & $0$ & $0$ & $\mathchar`-m_{\scalebox{0.8}{2}}$ & $0$ & $0$ & $0$ & $\mathchar`-m_{\scalebox{0.8}{3}}$ & $0$ & $0$ & $0$ & $%
m_{\scalebox{0.8}{5}}$ & $m_{\scalebox{0.8}{6}}$ & $m_{\scalebox{0.8}{7}}$ & $0$ & $0$ & $0$ \\ 
$0$ & $0$ & $0$ & $0$ & $\mathchar`-m_{\scalebox{0.8}{0}}$ & $0$ & $0$ & $0$ & $0$ & $0$ & $\mathchar`-m_{\scalebox{0.8}{1}}$ & $%
0$ & $0$ & $0$ & $0$ & $\mathchar`-m_{\scalebox{0.8}{2}}$ & $0$ & $0$ & $0$ & $\mathchar`-m_{\scalebox{0.8}{3}}$ & $0$ & $0$ & $%
\mathchar`-m_{\scalebox{0.8}{4}}$ & $0$ & $0$ & $m_{\scalebox{0.8}{6}}$ & $m_{\scalebox{0.8}{7}}$ & $0$ \\ 
$0$ & $0$ & $0$ & $0$ & $0$ & $\mathchar`-m_{\scalebox{0.8}{0}}$ & $0$ & $0$ & $0$ & $0$ & $0$ & $%
\mathchar`-m_{\scalebox{0.8}{1}}$ & $0$ & $0$ & $0$ & $0$ & $\mathchar`-m_{\scalebox{0.8}{2}}$ & $0$ & $0$ & $0$ & $\mathchar`-m_{\scalebox{0.8}{3}}$ & $0$
& $0$ & $\mathchar`-m_{\scalebox{0.8}{4}}$ & $0$ & $\mathchar`-m_{\scalebox{0.8}{5}}$ & $0$ & $m_{\scalebox{0.8}{7}}$ \\ 
$0$ & $0$ & $0$ & $0$ & $0$ & $0$ & $\mathchar`-m_{\scalebox{0.8}{0}}$ & $0$ & $0$ & $0$ & $0$ & $0$ & 
$\mathchar`-m_{\scalebox{0.8}{1}}$ & $0$ & $0$ & $0$ & $0$ & $\mathchar`-m_{\scalebox{0.8}{2}}$ & $0$ & $0$ & $0$ & $\mathchar`-m_{\scalebox{0.8}{3}}$ & $0
$ & $0$ & $\mathchar`-m_{\scalebox{0.8}{4}}$ & $0$ & $\mathchar`-m_{\scalebox{0.8}{5}}$ & $\mathchar`-m_{\scalebox{0.8}{6}}$%
\end{tabular}%
\right) $
}
\end{flushleft}

\bigskip

\emph{Proof. }By calculation of $\textrm{g}_{d}(D)m,$we have the above
expression. \ \ \ \emph{Q.E.D.}

\bigskip

\emph{Lemma 6.14. \ }For $D=\sum\limits_{0\leq i<j\leq 7}d_{ij}D_{ij}$
$(d_{ij}\in \R)$ $,$let's express $\textrm{f}_{DJJ}(D,?)$ as a matrix
$\textrm{MJl}(D)\in M(8\times 8,\R)$. We have $\textrm{f}_{DJJ}(D,?)$ is expressed by

$\textrm{MJl}(D)=\left( 
\begin{smallmatrix}
0 & d_{01} & d_{02} & d_{03} & d_{04} & d_{05} & d_{06} & d_{07} \\ 
-d_{01} & 0 & d_{12} & d_{13} & d_{14} & d_{15} & d_{16} & d_{17} \\ 
-d_{02} & -d_{12} & 0 & d_{23} & d_{24} & d_{25} & d_{26} & d_{27} \\ 
-d_{03} & -d_{13} & -d_{23} & 0 & d_{34} & d_{35} & d_{36} & d_{37} \\ 
-d_{04} & -d_{14} & -d_{24} & -d_{34} & 0 & d_{45} & d_{46} & d_{47} \\ 
-d_{05} & -d_{15} & -d_{25} & -d_{35} & -d_{45} & 0 & d_{56} & d_{57} \\ 
-d_{06} & -d_{16} & -d_{26} & -d_{36} & -d_{46} & -d_{56} & 0 & d_{67} \\ 
-d_{07} & -d_{17} & -d_{27} & -d_{37} & -d_{47} & -d_{57} & -d_{67} & 0%
\end{smallmatrix}%
\right)$

\bigskip

\emph{Proof. }By calculation of $\textrm{g}_{d}(D)m,$we have the above
expression. \ \ \emph{Q.E.D.}

\bigskip

\emph{Lemma 6.15. \ }For $m=\sum\limits_{i=0}^{7}m_{i}e_{i} (m%
_{i}\in \R)$ $,$let's express $\textrm{f}_{JJD}(m,?)$ as a matrix
$\textrm{MDl}(m)\in M(28\times 8,\R)$. We have $\textrm{f}_{JJD}(m,?)$ is expressed by

$\textrm{MDl}(m)=$
{\fontsize{8pt}{8pt} \selectfont%
$\left( 
\begin{tabular}{@{}c@{}c@{}c@{}c@{}c@{}c@{}c@{}c@{}}
$\mathchar`-m_{\scalebox{0.8}{1}}$ & $m_{\scalebox{0.8}{0}}$ & $0$ & $0$ & $0$ & $0$ & $0$ & $0$ \\ 
$\mathchar`-m_{\scalebox{0.8}{2}}$ & $0$ & $m_{\scalebox{0.8}{0}}$ & $0$ & $0$ & $0$ & $0$ & $0$ \\ 
$\mathchar`-m_{\scalebox{0.8}{3}}$ & $0$ & $0$ & $m_{\scalebox{0.8}{0}}$ & $0$ & $0$ & $0$ & $0$ \\ 
$\mathchar`-m_{\scalebox{0.8}{4}}$ & $0$ & $0$ & $0$ & $m_{\scalebox{0.8}{0}}$ & $0$ & $0$ & $0$ \\ 
$\mathchar`-m_{\scalebox{0.8}{5}}$ & $0$ & $0$ & $0$ & $0$ & $m_{\scalebox{0.8}{0}}$ & $0$ & $0$ \\ 
$\mathchar`-m_{\scalebox{0.8}{6}}$ & $0$ & $0$ & $0$ & $0$ & $0$ & $m_{\scalebox{0.8}{0}}$ & $0$ \\ 
$\mathchar`-m_{\scalebox{0.8}{7}}$ & $0$ & $0$ & $0$ & $0$ & $0$ & $0$ & $m_{\scalebox{0.8}{0}}$ \\ 
$0$ & $\mathchar`-m_{\scalebox{0.8}{2}}$ & $m_{\scalebox{0.8}{1}}$ & $0$ & $0$ & $0$ & $0$ & $0$ \\ 
$0$ & $\mathchar`-m_{\scalebox{0.8}{3}}$ & $0$ & $m_{\scalebox{0.8}{1}}$ & $0$ & $0$ & $0$ & $0$ \\ 
$0$ & $\mathchar`-m_{\scalebox{0.8}{4}}$ & $0$ & $0$ & $m_{\scalebox{0.8}{1}}$ & $0$ & $0$ & $0$ \\ 
$0$ & $\mathchar`-m_{\scalebox{0.8}{5}}$ & $0$ & $0$ & $0$ & $m_{\scalebox{0.8}{1}}$ & $0$ & $0$ \\ 
$0$ & $\mathchar`-m_{\scalebox{0.8}{6}}$ & $0$ & $0$ & $0$ & $0$ & $m_{\scalebox{0.8}{1}}$ & $0$ \\ 
$0$ & $\mathchar`-m_{\scalebox{0.8}{7}}$ & $0$ & $0$ & $0$ & $0$ & $0$ & $m_{\scalebox{0.8}{1}}$ \\ 
$0$ & $0$ & $\mathchar`-m_{\scalebox{0.8}{3}}$ & $m_{\scalebox{0.8}{2}}$ & $0$ & $0$ & $0$ & $0$ \\ 
$0$ & $0$ & $\mathchar`-m_{\scalebox{0.8}{4}}$ & $0$ & $m_{\scalebox{0.8}{2}}$ & $0$ & $0$ & $0$ \\ 
$0$ & $0$ & $\mathchar`-m_{\scalebox{0.8}{5}}$ & $0$ & $0$ & $m_{\scalebox{0.8}{2}}$ & $0$ & $0$ \\ 
$0$ & $0$ & $\mathchar`-m_{\scalebox{0.8}{6}}$ & $0$ & $0$ & $0$ & $m_{\scalebox{0.8}{2}}$ & $0$ \\ 
$0$ & $0$ & $\mathchar`-m_{\scalebox{0.8}{7}}$ & $0$ & $0$ & $0$ & $0$ & $m_{\scalebox{0.8}{2}}$ \\ 
$0$ & $0$ & $0$ & $\mathchar`-m_{\scalebox{0.8}{4}}$ & $m_{\scalebox{0.8}{3}}$ & $0$ & $0$ & $0$ \\ 
$0$ & $0$ & $0$ & $\mathchar`-m_{\scalebox{0.8}{5}}$ & $0$ & $m_{\scalebox{0.8}{3}}$ & $0$ & $0$ \\ 
$0$ & $0$ & $0$ & $\mathchar`-m_{\scalebox{0.8}{6}}$ & $0$ & $0$ & $m_{\scalebox{0.8}{3}}$ & $0$ \\ 
$0$ & $0$ & $0$ & $\mathchar`-m_{\scalebox{0.8}{7}}$ & $0$ & $0$ & $0$ & $m_{\scalebox{0.8}{3}}$ \\ 
$0$ & $0$ & $0$ & $0$ & $\mathchar`-m_{\scalebox{0.8}{5}}$ & $m_{\scalebox{0.8}{4}}$ & $0$ & $0$ \\ 
$0$ & $0$ & $0$ & $0$ & $\mathchar`-m_{\scalebox{0.8}{6}}$ & $0$ & $m_{\scalebox{0.8}{4}}$ & $0$ \\ 
$0$ & $0$ & $0$ & $0$ & $\mathchar`-m_{\scalebox{0.8}{7}}$ & $0$ & $0$ & $m_{\scalebox{0.8}{4}}$ \\ 
$0$ & $0$ & $0$ & $0$ & $0$ & $\mathchar`-m_{\scalebox{0.8}{6}}$ & $m_{\scalebox{0.8}{5}}$ & $0$ \\ 
$0$ & $0$ & $0$ & $0$ & $0$ & $\mathchar`-m_{\scalebox{0.8}{7}}$ & $0$ & $m_{\scalebox{0.8}{5}}$ \\ 
$0$ & $0$ & $0$ & $0$ & $0$ & $0$ & $\mathchar`-m_{\scalebox{0.8}{7}}$ & $m_{\scalebox{0.8}{6}}$%
\end{tabular}%
\right) $
}

\bigskip

\emph{Proof. }By calculation of $\textrm{JD}(m_{1},m_{2}),$we have the above
expression. \ \ \emph{Q.E.D.}

\bigskip

\emph{Lemma 6.16. \ }For $m=\sum\limits_{i=0}^{7}m_{i}e_{i} \ (m%
_{i}\in \R)$ $,$ let's express $\textrm{f}_{JJD}(?,m)$ as a matrix
$\textrm{MDr}(m)\in M(28\times 8,\R)$. We have $\textrm{f}_{JJD}(?,m)$ is expressed by

$\textrm{MDr}(m)=-\textrm{MDl}(m)=$ $^{t}\textrm{MJr}(m)$

\bigskip

\emph{Proof. }By calculation of $\textrm{JD}(m_{1},m_{2}),$we have the above
equation. \ \ \emph{Q.E.D.}

\bigskip

\emph{Lemma 6.17. }For $m=\sum\limits_{i=0}^{7}m_{i}e_{i}  \ (m_{i}\in
\R)$ $,$ let's express $\textrm{f}_{iJJ}(?,m)$ as a matrix
$\textrm{MIr}(m)\in M(8\times 8,\R)$.\ We have $\textrm{f}_{iJJ}(?,m)$ is expressed by

$\textrm{MIr}(m)=$
{\fontsize{10pt}{10pt} \selectfont%
$\left( 
\begin{tabular}{@{}c@{}c@{}c@{}c@{}c@{}c@{}c@{}c@{}}
$m_{\scalebox{0.8}{0}}$ & $\mathchar`-m_{\scalebox{0.8}{1}}$ & $\mathchar`-m_{\scalebox{0.8}{2}}$ & $\mathchar`-m_{\scalebox{0.8}{3}}$ & $\mathchar`-m_{\scalebox{0.8}{4}}$ & $\mathchar`-m_{\scalebox{0.8}{5}}$ & $\mathchar`-m_{\scalebox{0.8}{6}}$ & $%
\mathchar`-m_{\scalebox{0.8}{7}}$ \\ 
$\mathchar`-m_{\scalebox{0.8}{1}}$ & $\mathchar`-m_{\scalebox{0.8}{0}}$ & $\mathchar`-m_{\scalebox{0.8}{3}}$ & $m_{\scalebox{0.8}{2}}$ & $\mathchar`-m_{\scalebox{0.8}{5}}$ & $m_{\scalebox{0.8}{4}}$ & $\mathchar`-m_{\scalebox{0.8}{7}}$ & $%
m_{\scalebox{0.8}{6}}$ \\ 
$\mathchar`-m_{\scalebox{0.8}{2}}$ & $m_{\scalebox{0.8}{3}}$ & $\mathchar`-m_{\scalebox{0.8}{0}}$ & $\mathchar`-m_{\scalebox{0.8}{1}}$ & $m_{\scalebox{0.8}{6}}$ & $\mathchar`-m_{\scalebox{0.8}{7}}$ & $\mathchar`-m_{\scalebox{0.8}{4}}$ & $%
m_{\scalebox{0.8}{5}}$ \\ 
$\mathchar`-m_{\scalebox{0.8}{3}}$ & $\mathchar`-m_{\scalebox{0.8}{2}}$ & $m_{\scalebox{0.8}{1}}$ & $\mathchar`-m_{\scalebox{0.8}{0}}$ & $\mathchar`-m_{\scalebox{0.8}{7}}$ & $\mathchar`-m_{\scalebox{0.8}{6}}$ & $m_{\scalebox{0.8}{5}}$ & $%
m_{\scalebox{0.8}{4}}$ \\ 
$\mathchar`-m_{\scalebox{0.8}{4}}$ & $m_{\scalebox{0.8}{5}}$ & $\mathchar`-m_{\scalebox{0.8}{6}}$ & $m_{\scalebox{0.8}{7}}$ & $\mathchar`-m_{\scalebox{0.8}{0}}$ & $\mathchar`-m_{\scalebox{0.8}{1}}$ & $m_{\scalebox{0.8}{2}}$ & $%
\mathchar`-m_{\scalebox{0.8}{3}}$ \\ 
$\mathchar`-m_{\scalebox{0.8}{5}}$ & $\mathchar`-m_{\scalebox{0.8}{4}}$ & $m_{\scalebox{0.8}{7}}$ & $m_{\scalebox{0.8}{6}}$ & $m_{\scalebox{0.8}{1}}$ & $\mathchar`-m_{\scalebox{0.8}{0}}$ & $\mathchar`-m_{\scalebox{0.8}{3}}$ & $%
\mathchar`-m_{\scalebox{0.8}{2}}$ \\ 
$\mathchar`-m_{\scalebox{0.8}{6}}$ & $m_{\scalebox{0.8}{7}}$ & $m_{\scalebox{0.8}{4}}$ & $\mathchar`-m_{\scalebox{0.8}{5}}$ & $\mathchar`-m_{\scalebox{0.8}{2}}$ & $m_{\scalebox{0.8}{3}}$ & $\mathchar`-m_{\scalebox{0.8}{0}}$ & $%
\mathchar`-m_{\scalebox{0.8}{1}}$ \\ 
$\mathchar`-m_{\scalebox{0.8}{7}}$ & $\mathchar`-m_{\scalebox{0.8}{6}}$ & $\mathchar`-m_{\scalebox{0.8}{5}}$ & $\mathchar`-m_{\scalebox{0.8}{4}}$ & $m_{\scalebox{0.8}{3}}$ & $m_{\scalebox{0.8}{2}}$ & $m_{\scalebox{0.8}{1}}$ & $%
\mathchar`-m_{\scalebox{0.8}{0}}$%
\end{tabular}%
\right) $
}

\bigskip

\emph{Proof. }By direct calculation of $\overline{m_{1}m_{2}},$we
have the above expression evidently. 

\emph{Q.E.D.}

\bigskip

\emph{Lemma 6.18. }For $m=\sum\limits_{i=0}^{7}m_{i}e_{i} (m_{i}\in \R)$ $%
,$let's express $\textrm{f}_{iJJ}(m,?)$ as a matrix

\noindent
$\textrm{MIl}(m)\in M(8\times 8,\R)$. We have $\textrm{f}_{iJJ}(m,?)$ is expressed by

$\textrm{MIl}(m)=$ $^{t}\textrm{MIr}(m)=$
{\fontsize{10pt}{10pt} \selectfont%
$\left( 
\begin{tabular}{@{}c@{}c@{}c@{}c@{}c@{}c@{}c@{}c@{}}
$m_{\scalebox{0.8}{0}}$ & $\mathchar`-m_{\scalebox{0.8}{1}}$ & $\mathchar`-m_{\scalebox{0.8}{2}}$ & $\mathchar`-m_{\scalebox{0.8}{3}}$ & $\mathchar`-m_{\scalebox{0.8}{4}}$ & $\mathchar`-m_{\scalebox{0.8}{5}}$ & $\mathchar`-m_{\scalebox{0.8}{6}}$ & -%
$m_{\scalebox{0.8}{7}}$ \\ 
$\mathchar`-m_{\scalebox{0.8}{1}}$ & $\mathchar`-m_{\scalebox{0.8}{0}}$ & $m_{\scalebox{0.8}{3}}$ & $\mathchar`-m_{\scalebox{0.8}{2}}$ & $m_{\scalebox{0.8}{5}}$ & $\mathchar`-m_{\scalebox{0.8}{4}}$ & $m_{\scalebox{0.8}{7}}$ & -$%
m_{\scalebox{0.8}{6}}$ \\ 
$\mathchar`-m_{\scalebox{0.8}{2}}$ & $\mathchar`-m_{\scalebox{0.8}{3}}$ & $\mathchar`-m_{\scalebox{0.8}{0}}$ & $m_{\scalebox{0.8}{1}}$ & $\mathchar`-m_{\scalebox{0.8}{6}}$ & $m_{\scalebox{0.8}{7}}$ & $m_{\scalebox{0.8}{4}}$ & -$%
m_{\scalebox{0.8}{5}}$ \\ 
$\mathchar`-m_{\scalebox{0.8}{3}}$ & $m_{\scalebox{0.8}{2}}$ & $\mathchar`-m_{\scalebox{0.8}{1}}$ & $\mathchar`-m_{\scalebox{0.8}{0}}$ & $m_{\scalebox{0.8}{7}}$ & $m_{\scalebox{0.8}{6}}$ & $\mathchar`-m_{\scalebox{0.8}{5}}$ & -$%
m_{\scalebox{0.8}{4}}$ \\ 
$\mathchar`-m_{\scalebox{0.8}{4}}$ & $\mathchar`-m_{\scalebox{0.8}{5}}$ & $m_{\scalebox{0.8}{6}}$ & $\mathchar`-m_{\scalebox{0.8}{7}}$ & $\mathchar`-m_{\scalebox{0.8}{0}}$ & $m_{\scalebox{0.8}{1}}$ & $\mathchar`-m_{\scalebox{0.8}{2}}$ & $%
m_{\scalebox{0.8}{3}}$ \\ 
$\mathchar`-m_{\scalebox{0.8}{5}}$ & $m_{\scalebox{0.8}{4}}$ & $\mathchar`-m_{\scalebox{0.8}{7}}$ & $\mathchar`-m_{\scalebox{0.8}{6}}$ & $\mathchar`-m_{\scalebox{0.8}{1}}$ & $\mathchar`-m_{\scalebox{0.8}{0}}$ & $m_{\scalebox{0.8}{3}}$ & $%
m_{\scalebox{0.8}{2}}$ \\ 
$\mathchar`-m_{\scalebox{0.8}{6}}$ & $\mathchar`-m_{\scalebox{0.8}{7}}$ & $\mathchar`-m_{\scalebox{0.8}{4}}$ & $m_{\scalebox{0.8}{5}}$ & $m_{\scalebox{0.8}{2}}$ & $\mathchar`-m_{\scalebox{0.8}{3}}$ & $\mathchar`-m_{\scalebox{0.8}{0}}$ & $%
m_{\scalebox{0.8}{1}}$ \\ 
$\mathchar`-m_{\scalebox{0.8}{7}}$ & $m_{\scalebox{0.8}{6}}$ & $m_{\scalebox{0.8}{5}}$ & $m_{\scalebox{0.8}{4}}$ & $\mathchar`-m_{\scalebox{0.8}{3}}$ & $\mathchar`-m_{\scalebox{0.8}{2}}$ & $\mathchar`-m_{\scalebox{0.8}{1}}$ & -$%
m_{\scalebox{0.8}{0}}$%
\end{tabular}%
\right) $
}
\bigskip

\emph{Proof. }By direct calculation of $\overline{m_{1}m_{2}},$we
have the above expression. 

\emph{Q.E.D.}

\bigskip

\emph{Lemma 6.19. \ }For $m=\sum\limits_{i=0}^{7}m_{i}e_{i} \ (m%
_{i}\in \R)$ $,$let's express $\textrm{f}_{JJC}(m,?)$ as a matrix
$\textrm{MC}(m)\in M(1\times 8,\R)$. We have $\textrm{f}_{JJC}(m,?)$ is expressed by

$\textrm{MC}(m)=\left( 
\begin{array}{cccccccc}
m_{0} & m_{1} & m_{2} & m_{3} & m_{4} & m_{5} & m_{6} & m_{7}%
\end{array}%
\right) $

\bigskip

\emph{Proof. }By calculation of $(m_{1},m_{2}),$we have the above
expression. \ \ \ \ \emph{Q.E.D.}

\bigskip

\emph{Lemma 6.20. }For $\tau \in \R$ $,$let's express $\textrm{f}_{CJJ}(\tau
,?)$ as a matrix
$\textrm{ME}(\tau )\in M(8\times 8,\R)$. We have $\textrm{f}_{CJJ}(\tau ,?)$ is expressed by

$\textrm{ME}(\tau )=\left( 
\begin{array}{cccccccc}
\tau &  &  &  &  &  &  &  \\ 
& \tau &  &  &  &  &  &  \\ 
&  & \tau &  &  &  &  &  \\ 
&  &  & \tau &  &  &  &  \\ 
&  &  &  & \tau &  &  &  \\ 
&  &  &  &  & \tau &  &  \\ 
&  &  &  &  &  & \tau &  \\ 
&  &  &  &  &  &  & \tau%
\end{array}%
\right) =\tau E$

\noindent
where blank elements means 0.

\bigskip

\emph{Proof. } Evidently. \ \ \ \emph{\ Q.E.D.}

\bigskip

\emph{Lemma 6.21. \ }For $m=\sum\limits_{i=0}^{7}m_{i}e_{i} \ (m%
_{i}\in \R)$ $,$let's express $\textrm{f}_{CJJ}(?,m)$ as a matrix
of $M(8\times 1,\R)$. We have $\textrm{f}_{CJJ}(?,m)$ is expressed by

$\ \ \ \ \ \ \ \ \ \ \ \ \ \left( 
\begin{array}{c}
m_{0} \\ 
m_{1} \\ 
m_{2} \\ 
m_{3} \\ 
m_{4} \\ 
m_{5} \\ 
m_{6} \\ 
m_{7}%
\end{array}%
\right) =$ $^{t}\textrm{MC}(m)$

\bigskip

\emph{Proof. }By \emph{Definition 6.11}$,$we have the above
expression. \ \ \ \ \emph{Q.E.D.}

\bigskip

Let's express $\textrm{d}_{g}\nu \textrm{g}_{d}$ as a matrix $Mv\in M(28\times 28,\R)$.

\bigskip

\emph{Lemma 6.22. \ }We have

$2Mv=$

{\fontsize{6pt}{8pt} \selectfont%
$\left( 
\begin{smallmatrix}
-1 & 0 & 0 & 0 & 0 & 0 & 0 & 0 & 0 & 0 & 0 & 0 & 0 & 1 & 0 & 0 & 0 & 0 & 0 & 
0 & 0 & 0 & 1 & 0 & 0 & 0 & 0 & 1 \\ 
0 & -1 & 0 & 0 & 0 & 0 & 0 & 0 & -1 & 0 & 0 & 0 & 0 & 0 & 0 & 0 & 0 & 0 & 0
& 0 & 0 & 0 & 0 & -1 & 0 & 0 & 1 & 0 \\ 
0 & 0 & -1 & 0 & 0 & 0 & 0 & 1 & 0 & 0 & 0 & 0 & 0 & 0 & 0 & 0 & 0 & 0 & 0 & 
0 & 0 & 0 & 0 & 0 & 1 & 1 & 0 & 0 \\ 
0 & 0 & 0 & -1 & 0 & 0 & 0 & 0 & 0 & 0 & -1 & 0 & 0 & 0 & 0 & 0 & 1& 0 & 0 &
 0 & 0 & -1 & 0 & 0 & 0 & 0 & 0 & 0 \\ 
0 & 0 & 0 & 0 & -1 & 0 & 0 & 0 & 0 & 1 & 0 & 0 & 0 & 0 & 0 & 0 & 0 & 
-1 & 0 & 0 & -1 & 0 & 0 & 0 & 0 & 0 & 0 & 0 \\ 
0 & 0 & 0 & 0 & 0 & -1 & 0 & 0 & 0 & 0 & 0 & 0 & -1 & 0 & -1 & 0 & 0 & 0 & 0
& 1 & 0 & 0 & 0 & 0 & 0 & 0 & 0 & 0 \\ 
0 & 0 & 0 & 0 & 0 & 0 & -1 & 0 & 0 & 0 & 0 & 1 & 0 & 0 & 0 & 1 & 0 & 0 & 1 & 
0 & 0 & 0 & 0 & 0 & 0 & 0 & 0 & 0 \\ 
0 & 0 & -1 & 0 & 0 & 0 & 0 & 1 & 0 & 0 & 0 & 0 & 0 & 0 & 0 & 0 & 0 & 0 & 0 & 0
& 0 & 0 & 0 & 0 & -1 & -1 & 0 & 0 \\ 
0 & 1 & 0 & 0 & 0 & 0 & 0 & 0 & 1 & 0 & 0 & 0 & 0 & 0 & 0 & 0 & 0 & 0 & 0 & 0
& 0 & 0 & 0 & -1 & 0 & 0 & 1 & 0 \\ 
0 & 0 & 0 & 0 & -1 & 0 & 0 & 0 & 0 & 1 & 0 & 0 & 0 & 0 & 0 & 0 & 0 & 1 & 0 & 0
& 1 & 0 & 0 & 0 & 0 & 0 & 0 & 0 \\ 
0 & 0 & 0 & 1 & 0 & 0 & 0 & 0 & 0 & 0 & 1 & 0 & 0 & 0 & 0 & 0 & 1 & 0 & 0 & 0
& 0 & -1 & 0 & 0 & 0 & 0 & 0 & 0 \\ 
0 & 0 & 0 & 0 & 0 & 0 & -1 & 0 & 0 & 0 & 0 & 1 & 0 & 0 & 0 & -1 & 0 & 0 & -1
& 0 & 0 & 0 & 0 & 0 & 0 & 0 & 0 & 0 \\ 
0 & 0 & 0 & 0 & 0 & 1 & 0 & 0 & 0 & 0 & 0 & 0 & 1 & 0 & -1 & 0 & 0 & 0 & 0 & 
1 & 0 & 0 & 0 & 0 & 0 & 0 & 0 & 0 \\ 
-1 & 0 & 0 & 0 & 0 & 0 & 0 & 0 & 0 & 0 & 0 & 0 & 0 & 1 & 0 & 0 & 0 & 0 & 0 & 
0 & 0 & 0 & -1 & 0 & 0 & 0 & 0 & -1 \\ 
0 & 0 & 0 & 0 & 0 & 1 & 0 & 0 & 0 & 0 & 0 & 0 & -1 & 0 & 1 & 0 & 0 & 0 & 0 & 
1 & 0 & 0 & 0 & 0 & 0 & 0 & 0 & 0 \\ 
0 & 0 & 0 & 0 & 0 & 0 & -1 & 0 & 0 & 0 & 0 & -1 & 0 & 0 & 0 & 1 & 0 & 0 & -1
& 0 & 0 & 0 & 0 & 0 & 0 & 0 & 0 & 0 \\ 
0 & 0 & 0 & -1 & 0 & 0 & 0 & 0 & 0 & 0 & 1 & 0 & 0 & 0 & 0 & 0 & 1 & 0 & 0 & 
0 & 0 & 1 & 0 & 0 & 0 & 0 & 0 & 0 \\ 
0 & 0 & 0 & 0 & 1 & 0 & 0 & 0 & 0 & 1 & 0 & 0 & 0 & 0 & 0 & 0 & 0 & 1 & 0 & 0
& -1 & 0 & 0 & 0 & 0 & 0 & 0 & 0 \\ 
0 & 0 & 0 & 0 & 0 & 0 & -1 & 0 & 0 & 0 & 0 & -1 & 0 & 0 & 0 & -1 & 0 & 0 & 1
& 0 & 0 & 0 & 0 & 0 & 0 & 0 & 0 & 0 \\ 
0 & 0 & 0 & 0 & 0 & -1 & 0 & 0 & 0 & 0 & 0 & 0 & 1 & 0 & 1 & 0 & 0 & 0 & 0 & 
1 & 0 & 0 & 0 & 0 & 0 & 0 & 0 & 0 \\ 
0 & 0 & 0 & 0 & 1 & 0 & 0 & 0 & 0 & 1 & 0 & 0 & 0 & 0 & 0 & 0 & 0 & -1 & 0 & 
0 & 1 & 0 & 0 & 0 & 0 & 0 & 0 & 0 \\ 
0 & 0 & 0 & 1 & 0 & 0 & 0 & 0 & 0 & 0 & -1 & 0 & 0 & 0 & 0 & 0 & 1 & 0 & 0 & 
0 & 0 & 1 & 0 & 0 & 0 & 0 & 0 & 0 \\ 
-1 & 0 & 0 & 0 & 0 & 0 & 0 & 0 & 0 & 0 & 0 & 0 & 0 & -1 & 0 & 0 & 0 & 0 & 0
& 0 & 0 & 0 & 1 & 0 & 0 & 0 & 0 & -1 \\ 
0 & 1 & 0 & 0 & 0 & 0 & 0 & 0 & -1 & 0 & 0 & 0 & 0 & 0 & 0 & 0 & 0 & 0 & 0 & 
0 & 0 & 0 & 0 & 1 & 0 & 0 & 1 & 0 \\ 
0 & 0 & -1 & 0 & 0 & 0 & 0 & -1 & 0 & 0 & 0 & 0 & 0 & 0 & 0 & 0 & 0 & 0 & 0
& 0 & 0 & 0 & 0 & 0 & 1 & -1 & 0 & 0 \\ 
0 & 0 & -1 & 0 & 0 & 0 & 0 & -1 & 0 & 0 & 0 & 0 & 0 & 0 & 0 & 0 & 0 & 0 & 0
& 0 & 0 & 0 & 0 & 0 & -1 & 1 & 0 & 0 \\ 
0 & -1 & 0 & 0 & 0 & 0 & 0 & 0 & 1 & 0 & 0 & 0 & 0 & 0 & 0 & 0 & 0 & 0 & 0 & 
0 & 0 & 0 & 0 & 1 & 0 & 0 & 1 & 0 \\ 
-1 & 0 & 0 & 0 & 0 & 0 & 0 & 0 & 0 & 0 & 0 & 0 & 0 & -1 & 0 & 0 & 0 & 0 & 0
& 0 & 0 & 0 & -1 & 0 & 0 & 0 & 0 & 1%
\end{smallmatrix}%
\right). $ }
\bigskip

\emph{Proof. }By \emph{Definition 1.5, 1.7, 1.13 }and \emph{Proposition 1.9}, 
we have the above expression.
\ \ \ \ \emph{Q.E.D.}

\bigskip 

\emph{Corollary 6.22.1.} \ $M\nu ^{2}=$ $^{t}M\nu $

\bigskip

\emph{Proof. }By calculation we have $M\nu ^{2}=$ $^{t}M\nu $. \ \
\ \ \emph{Q.E.D.}

\bigskip

\emph{Lemma 6.23. \ }For $D=\sum\limits_{0\leq i<j\leq 7}d_{ij}D_{ij}$
$(d_{ij}\in \R)$ ,let's express 

\noindent
$\textrm{f}_{DJJ}(D,?)\textrm{d}_{g}\nu \textrm{g}_{d}$ and
$\textrm{f}_{DJJ}(D,?)\textrm{d}_{g}\nu ^{2}\textrm{g}_{d}$ as a matrix $\textrm{MJl}(MvD)$  and

\noindent
$\textrm{MJl}(Mv^{2}D)\in
M(8\times 8,\R)$ respectively.

\noindent
We have $\textrm{f}_{DJJ}(D,?)\textrm{d}_{g}\nu \textrm{g}_{d}$ and $\textrm{f}_{DJJ}(D,?)\textrm{d}_{g}\nu ^{2}\textrm{g}_{d}$
are expressed by

$2\textrm{MJl}(MvD)=$

\begin{flushleft}
{\fontsize{6pt}{8pt} \selectfont%
$\left( 
\begin{tabular}{@{}ccccccc@{}}
$0$ & $d_{\scalebox{0.8}{67}}+d_{\scalebox{0.8}{45}}+d_{\scalebox{0.8}{23}}$-$d_{\scalebox{0.8}{01}}$ & $d_{\scalebox{0.8}{57}}$-$d_{\scalebox{0.8}{46}}$-$d_{\scalebox{0.8}{13}}$-$d_{\scalebox{0.8}{02}}$ & $%
d_{\scalebox{0.8}{56}}+d_{\scalebox{0.8}{47}}+d_{\scalebox{0.8}{12}}$-$d_{\scalebox{0.8}{03}}$ \\ 
-$d_{\scalebox{0.8}{67}}$-$d_{\scalebox{0.8}{45}}$-$d_{\scalebox{0.8}{23}}+d_{\scalebox{0.8}{01}}$ & $0$ & -$d_{\scalebox{0.8}{56}}$-$d_{\scalebox{0.8}{47}}+d_{\scalebox{0.8}{12}}$-$d_{\scalebox{0.8}{03}}$ & $%
d_{\scalebox{0.8}{57}}$-$d_{\scalebox{0.8}{46}}+d_{\scalebox{0.8}{13}}+d_{\scalebox{0.8}{02}}$ \\ 
-$d_{\scalebox{0.8}{57}}+d_{\scalebox{0.8}{46}}+d_{\scalebox{0.8}{13}}+d_{\scalebox{0.8}{02}}$ & $d_{\scalebox{0.8}{56}}+d_{\scalebox{0.8}{47}}$-$d_{\scalebox{0.8}{12}}+d_{\scalebox{0.8}{03}}$ & $0$ & $%
$-$d_{\scalebox{0.8}{67}}$-$d_{\scalebox{0.8}{45}}+d_{\scalebox{0.8}{23}}$-$d_{\scalebox{0.8}{01}}$ \\ 
-$d_{\scalebox{0.8}{56}}$-$d_{\scalebox{0.8}{47}}$-$d_{\scalebox{0.8}{12}}+d_{\scalebox{0.8}{03}}$ & -$d_{\scalebox{0.8}{57}}+d_{\scalebox{0.8}{46}}$-$d_{\scalebox{0.8}{13}}$-$d_{\scalebox{0.8}{02}}$ & $%
d_{\scalebox{0.8}{67}}+d_{\scalebox{0.8}{45}}$-$d_{\scalebox{0.8}{23}}+d_{\scalebox{0.8}{01}}$ & $0$ \\ 
$d_{\scalebox{0.8}{37}}$-$d_{\scalebox{0.8}{26}}+d_{\scalebox{0.8}{15}}+d_{\scalebox{0.8}{04}}$ & -$d_{\scalebox{0.8}{36}}$-$d_{\scalebox{0.8}{27}}$-$d_{\scalebox{0.8}{14}}+d_{\scalebox{0.8}{05}}$ & $%
$-$d_{\scalebox{0.8}{35}}$-$d_{\scalebox{0.8}{24}}+d_{\scalebox{0.8}{17}}$-$d_{\scalebox{0.8}{06}}$ & -$d_{\scalebox{0.8}{34}}+d_{\scalebox{0.8}{25}}+d_{\scalebox{0.8}{16}}+d_{\scalebox{0.8}{07}}$ \\
\end{tabular}%
\right. $

\bigskip

$\left. 
\begin{tabular}{@{}cccccccc@{}}
-$d_{\scalebox{0.8}{37}}+d_{\scalebox{0.8}{26}}$-$d_{\scalebox{0.8}{15}}$-$d_{\scalebox{0.8}{04}}$ & -$d_{\scalebox{0.8}{36}}$-$d_{\scalebox{0.8}{27}}+d_{\scalebox{0.8}{14}}$-$d_{\scalebox{0.8}{05}}$ & $%
d_{\scalebox{0.8}{35}}$-$d_{\scalebox{0.8}{24}}$-$d_{\scalebox{0.8}{17}}$-$d_{\scalebox{0.8}{06}}$ & $d_{\scalebox{0.8}{34}}+d_{\scalebox{0.8}{25}}+d_{\scalebox{0.8}{16}}$-$d_{\scalebox{0.8}{07}}$ \\ 
$d_{\scalebox{0.8}{36}}+d_{\scalebox{0.8}{27}}+d_{\scalebox{0.8}{14}}$-$d_{\scalebox{0.8}{05}}$ & -$d_{\scalebox{0.8}{37}}+d_{\scalebox{0.8}{26}}+d_{\scalebox{0.8}{15}}+d_{\scalebox{0.8}{04}}$ & $%
$-$d_{\scalebox{0.8}{34}}$-$d_{\scalebox{0.8}{25}}+d_{\scalebox{0.8}{16}}$-$d_{\scalebox{0.8}{07}}$ & $d_{\scalebox{0.8}{35}}$-$d_{\scalebox{0.8}{24}}+d_{\scalebox{0.8}{17}}+d_{\scalebox{0.8}{06}}$ \\ 
$d_{\scalebox{0.8}{35}}+d_{\scalebox{0.8}{24}}$-$d_{\scalebox{0.8}{17}}+d_{\scalebox{0.8}{06}}$ & -$d_{\scalebox{0.8}{34}}+d_{\scalebox{0.8}{25}}$-$d_{\scalebox{0.8}{16}}$-$d_{\scalebox{0.8}{07}}$ & $%
d_{\scalebox{0.8}{37}}+d_{\scalebox{0.8}{26}}+d_{\scalebox{0.8}{15}}$-$d_{\scalebox{0.8}{04}}$ & -$d_{\scalebox{0.8}{36}}+d_{\scalebox{0.8}{27}}+d_{\scalebox{0.8}{14}}+d_{\scalebox{0.8}{05}}$ \\ 
$d_{\scalebox{0.8}{34}}$-$d_{\scalebox{0.8}{25}}$-$d_{\scalebox{0.8}{16}}$-$d_{\scalebox{0.8}{07}}$ & $d_{\scalebox{0.8}{35}}+d_{\scalebox{0.8}{24}}+d_{\scalebox{0.8}{17}}$-$d_{\scalebox{0.8}{06}}$ & $%
d_{\scalebox{0.8}{36}}$-$d_{\scalebox{0.8}{27}}+d_{\scalebox{0.8}{14}}+d_{\scalebox{0.8}{05}}$ & $d_{\scalebox{0.8}{37}}+d_{\scalebox{0.8}{26}}$-$d_{\scalebox{0.8}{15}}+d_{\scalebox{0.8}{04}}$ \\ 
$0$ & -$d_{\scalebox{0.8}{67}}+d_{\scalebox{0.8}{45}}$-$d_{\scalebox{0.8}{23}}$-$d_{\scalebox{0.8}{01}}$ & $d_{\scalebox{0.8}{57}}+d_{\scalebox{0.8}{46}}$-$d_{\scalebox{0.8}{13}}+d_{\scalebox{0.8}{02}}$ & $%
$-$d_{\scalebox{0.8}{56}}+d_{\scalebox{0.8}{47}}$-$d_{\scalebox{0.8}{12}}$-$d_{\scalebox{0.8}{03}}$ \\
\end{tabular}%
\right) $
}
\end{flushleft}

\bigskip

$2\textrm{MJl}(Mv^{2}D)=$

{\fontsize{6pt}{8pt} \selectfont%
$\left( 
\begin{tabular}{@{}cccccccc@{}}
$0$ & -$d_{\scalebox{0.8}{67}}$-$d_{\scalebox{0.8}{45}}$-$d_{\scalebox{0.8}{23}}$-$d_{\scalebox{0.8}{01}}$ & -$d_{\scalebox{0.8}{57}}+d_{\scalebox{0.8}{46}}+d_{\scalebox{0.8}{13}}$-$d_{\scalebox{0.8}{02}}$ & $%
$-$d_{\scalebox{0.8}{56}}$-$d_{\scalebox{0.8}{47}}$-$d_{\scalebox{0.8}{12}}$-$d_{\scalebox{0.8}{03}}$ \\ 
$d_{\scalebox{0.8}{67}}+d_{\scalebox{0.8}{45}}+d_{\scalebox{0.8}{23}}+d_{\scalebox{0.8}{01}}$ & $0$ & -$d_{\scalebox{0.8}{56}}$-$d_{\scalebox{0.8}{47}}+d_{\scalebox{0.8}{12}}+d_{\scalebox{0.8}{03}}$ & $%
d_{\scalebox{0.8}{57}}$-$d_{\scalebox{0.8}{46}}+d_{\scalebox{0.8}{13}}$-$d_{\scalebox{0.8}{02}}$ \\ 
$d_{\scalebox{0.8}{57}}$-$d_{\scalebox{0.8}{46}}$-$d_{\scalebox{0.8}{13}}+d_{\scalebox{0.8}{02}}$ & $d_{\scalebox{0.8}{56}}+d_{\scalebox{0.8}{47}}$-$d_{\scalebox{0.8}{12}}$-$d_{\scalebox{0.8}{03}}$ & $0$ & $%
$-$d_{\scalebox{0.8}{67}}$-$d_{\scalebox{0.8}{45}}+d_{\scalebox{0.8}{23}}+d_{\scalebox{0.8}{01}}$ \\ 
$d_{\scalebox{0.8}{56}}+d_{\scalebox{0.8}{47}}+d_{\scalebox{0.8}{12}}+d_{\scalebox{0.8}{03}}$ & -$d_{\scalebox{0.8}{57}}+d_{\scalebox{0.8}{46}}$-$d_{\scalebox{0.8}{13}}+d_{\scalebox{0.8}{02}}$ & $%
d_{\scalebox{0.8}{67}}+d_{\scalebox{0.8}{45}}$-$d_{\scalebox{0.8}{23}}$-$d_{\scalebox{0.8}{01}}$ & $0$ \\ 
-$d_{\scalebox{0.8}{37}}+d_{\scalebox{0.8}{26}}$-$d_{\scalebox{0.8}{15}}+d_{\scalebox{0.8}{04}}$ & -$d_{\scalebox{0.8}{36}}$-$d_{\scalebox{0.8}{27}}$-$d_{\scalebox{0.8}{14}}$-$d_{\scalebox{0.8}{05}}$ & $%
$-$d_{\scalebox{0.8}{35}}$-$d_{\scalebox{0.8}{24}}+d_{\scalebox{0.8}{17}}+d_{\scalebox{0.8}{06}}$ & -$d_{\scalebox{0.8}{34}}+d_{\scalebox{0.8}{25}}+d_{\scalebox{0.8}{16}}$-$d_{\scalebox{0.8}{07}}$ \\
\end{tabular}%
\right. $

\bigskip 

$\left. 
\begin{tabular}{@{}cccccccc@{}}
$d_{\scalebox{0.8}{37}}$-$d_{\scalebox{0.8}{26}}+d_{\scalebox{0.8}{15}}$-$d_{\scalebox{0.8}{04}}$ & $d_{\scalebox{0.8}{36}}+d_{\scalebox{0.8}{27}}$-$d_{\scalebox{0.8}{14}}$-$d_{\scalebox{0.8}{05}}$ & $%
$-$d_{\scalebox{0.8}{35}}+d_{\scalebox{0.8}{24}}+d_{\scalebox{0.8}{17}}$-$d_{\scalebox{0.8}{06}}$ & -$d_{\scalebox{0.8}{34}}$-$d_{\scalebox{0.8}{25}}$-$d_{\scalebox{0.8}{16}}$-$d_{\scalebox{0.8}{07}}$ \\ 
$d_{\scalebox{0.8}{36}}+d_{\scalebox{0.8}{27}}+d_{\scalebox{0.8}{14}}+d_{\scalebox{0.8}{05}}$ & -$d_{\scalebox{0.8}{37}}+d_{\scalebox{0.8}{26}}+d_{\scalebox{0.8}{15}}$-$d_{\scalebox{0.8}{04}}$ & $%
$-$d_{\scalebox{0.8}{34}}$-$d_{\scalebox{0.8}{25}}+d_{\scalebox{0.8}{16}}+d_{\scalebox{0.8}{07}}$ & $d_{\scalebox{0.8}{35}}$-$d_{\scalebox{0.8}{24}}+d_{\scalebox{0.8}{17}}$-$d_{\scalebox{0.8}{06}}$ \\ 
$d_{\scalebox{0.8}{35}}+d_{\scalebox{0.8}{24}}$-$d_{\scalebox{0.8}{17}}$-$d_{\scalebox{0.8}{06}}$ & -$d_{\scalebox{0.8}{34}}+d_{\scalebox{0.8}{25}}$-$d_{\scalebox{0.8}{16}}+d_{\scalebox{0.8}{07}}$ & $%
d_{\scalebox{0.8}{37}}+d_{\scalebox{0.8}{26}}+d_{\scalebox{0.8}{15}}+d_{\scalebox{0.8}{04}}$ & -$d_{\scalebox{0.8}{36}}+d_{\scalebox{0.8}{27}}+d_{\scalebox{0.8}{14}}$-$d_{\scalebox{0.8}{05}}$ \\ 
$d_{\scalebox{0.8}{34}}$-$d_{\scalebox{0.8}{25}}$-$d_{\scalebox{0.8}{16}}+d_{\scalebox{0.8}{07}}$ & $d_{\scalebox{0.8}{35}}+d_{\scalebox{0.8}{24}}+d_{\scalebox{0.8}{17}}+d_{\scalebox{0.8}{06}}$ & $%
d_{\scalebox{0.8}{36}}$-$d_{\scalebox{0.8}{27}}+d_{\scalebox{0.8}{14}}$-$d_{\scalebox{0.8}{05}}$ & $d_{\scalebox{0.8}{37}}+d_{\scalebox{0.8}{26}}$-$d_{\scalebox{0.8}{15}}$-$d_{\scalebox{0.8}{04}}$ \\ 
$0$ & -$d_{\scalebox{0.8}{67}}+d_{\scalebox{0.8}{45}}$-$d_{\scalebox{0.8}{23}}+d_{\scalebox{0.8}{01}}$ & $d_{\scalebox{0.8}{57}}+d_{\scalebox{0.8}{46}}$-$d_{\scalebox{0.8}{13}}$-$d_{\scalebox{0.8}{02}}$ & $%
$-$d_{\scalebox{0.8}{56}}+d_{\scalebox{0.8}{47}}$-$d_{\scalebox{0.8}{12}}+d_{\scalebox{0.8}{03}}$ \\
\end{tabular}%
\right) $
}

\bigskip

\emph{Proof. }By \emph{Lemma 6.14 }and \emph{Lemma 6.12} $,$we
have the above expression. 

\emph{Q.E.D.}

\bigskip

\emph{\ Lemma 6.24.} We have

\ $2Mv\textrm{MDl}(m)=$
{\fontsize{8pt}{10pt} \selectfont%
$\left( 
\begin{tabular}{@{}c@{}c@{}c@{}c@{}c@{}c@{}c@{}c@{}}
$m_{\scalebox{0.8}{1}}$ & $-m_{\scalebox{0.8}{0}}$ & $-m_{\scalebox{0.8}{3}}$ & $m_{\scalebox{0.8}{2}}$ & $-m_{\scalebox{0.8}{5}}$ & $m_{\scalebox{0.8}{4}}$ & $-m_{\scalebox{0.8}{7}}$ & 
$m_{\scalebox{0.8}{6}}$ \\ 
$m_{\scalebox{0.8}{2}}$ & $m_{\scalebox{0.8}{3}}$ & $-m_{\scalebox{0.8}{0}}$ & $-m_{\scalebox{0.8}{1}}$ & $m_{\scalebox{0.8}{6}}$ & $-m_{\scalebox{0.8}{7}}$ & $-m_{\scalebox{0.8}{4}}$ & $%
m_{\scalebox{0.8}{5}}$ \\ 
$m_{\scalebox{0.8}{3}}$ & $-m_{\scalebox{0.8}{2}}$ & $m_{\scalebox{0.8}{1}}$ & $-m_{\scalebox{0.8}{0}}$ & $-m_{\scalebox{0.8}{7}}$ & $-m_{\scalebox{0.8}{6}}$ & $m_{\scalebox{0.8}{5}}$ & $%
m_{\scalebox{0.8}{4}}$ \\ 
$m_{\scalebox{0.8}{4}}$ & $m_{\scalebox{0.8}{5}}$ & $-m_{\scalebox{0.8}{6}}$ & $m_{\scalebox{0.8}{7}}$ & $-m_{\scalebox{0.8}{0}}$ & $-m_{\scalebox{0.8}{1}}$ & $m_{\scalebox{0.8}{2}}$ & $%
-m_{\scalebox{0.8}{3}}$ \\ 
$m_{\scalebox{0.8}{5}}$ & $-m_{\scalebox{0.8}{4}}$ & $m_{\scalebox{0.8}{7}}$ & $m_{\scalebox{0.8}{6}}$ & $m_{\scalebox{0.8}{1}}$ & $-m_{\scalebox{0.8}{0}}$ & $-m_{\scalebox{0.8}{3}}$ & $%
-m_{\scalebox{0.8}{2}}$ \\ 
$m_{\scalebox{0.8}{6}}$ & $m_{\scalebox{0.8}{7}}$ & $m_{\scalebox{0.8}{4}}$ & $-m_{\scalebox{0.8}{5}}$ & $-m_{\scalebox{0.8}{2}}$ & $m_{\scalebox{0.8}{3}}$ & $-m_{\scalebox{0.8}{0}}$ & $%
-m_{\scalebox{0.8}{1}}$ \\ 
$m_{\scalebox{0.8}{7}}$ & $-m_{\scalebox{0.8}{6}}$ & $-m_{\scalebox{0.8}{5}}$ & $-m_{\scalebox{0.8}{4}}$ & $m_{\scalebox{0.8}{3}}$ & $m_{\scalebox{0.8}{2}}$ & $m_{\scalebox{0.8}{1}}$ & $%
-m_{\scalebox{0.8}{0}}$ \\ 
$m_{\scalebox{0.8}{3}}$ & $-m_{\scalebox{0.8}{2}}$ & $m_{\scalebox{0.8}{1}}$ & $-m_{\scalebox{0.8}{0}}$ & $m_{\scalebox{0.8}{7}}$ & $m_{\scalebox{0.8}{6}}$ & $-m_{\scalebox{0.8}{5}}$ & $%
-m_{\scalebox{0.8}{4}}$ \\
\end{tabular}%
\right) $
}

\bigskip 

\ $2Mv^{2}\textrm{MDl}(m)=$
{\fontsize{8pt}{10pt} \selectfont%
$\left( 
\begin{tabular}{@{}c@{}c@{}c@{}c@{}c@{}c@{}c@{}c@{}}
$m_{\scalebox{0.8}{1}}$ & $-m_{\scalebox{0.8}{0}}$ & $m_{\scalebox{0.8}{3}}$ & $-m_{\scalebox{0.8}{2}}$ & $m_{\scalebox{0.8}{5}}$ & $-m_{\scalebox{0.8}{4}}$ & $m_{\scalebox{0.8}{7}}$ & $%
-m_{\scalebox{0.8}{6}}$ \\ 
$m_{\scalebox{0.8}{2}}$ & $-m_{\scalebox{0.8}{3}}$ & $-m_{\scalebox{0.8}{0}}$ & $m_{\scalebox{0.8}{1}}$ & $-m_{\scalebox{0.8}{6}}$ & $m_{\scalebox{0.8}{7}}$ & $m_{\scalebox{0.8}{4}}$ & $%
-m_{\scalebox{0.8}{5}}$ \\ 
$m_{\scalebox{0.8}{3}}$ & $m_{\scalebox{0.8}{2}}$ & $-m_{\scalebox{0.8}{1}}$ & $-m_{\scalebox{0.8}{0}}$ & $m_{\scalebox{0.8}{7}}$ & $m_{\scalebox{0.8}{6}}$ & $-m_{\scalebox{0.8}{5}}$ & $%
-m_{\scalebox{0.8}{4}}$ \\ 
$m_{\scalebox{0.8}{4}}$ & $-m_{\scalebox{0.8}{5}}$ & $m_{\scalebox{0.8}{6}}$ & $-m_{\scalebox{0.8}{7}}$ & $-m_{\scalebox{0.8}{0}}$ & $m_{\scalebox{0.8}{1}}$ & $-m_{\scalebox{0.8}{2}}$ & $%
m_{\scalebox{0.8}{3}}$ \\ 
$m_{\scalebox{0.8}{5}}$ & $m_{\scalebox{0.8}{4}}$ & $-m_{\scalebox{0.8}{7}}$ & $-m_{\scalebox{0.8}{6}}$ & $-m_{\scalebox{0.8}{1}}$ & $-m_{\scalebox{0.8}{0}}$ & $m_{\scalebox{0.8}{3}}$ & $%
m_{\scalebox{0.8}{2}}$ \\ 
$m_{\scalebox{0.8}{6}}$ & $-m_{\scalebox{0.8}{7}}$ & $-m_{\scalebox{0.8}{4}}$ & $m_{\scalebox{0.8}{5}}$ & $m_{\scalebox{0.8}{2}}$ & $-m_{\scalebox{0.8}{3}}$ & $-m_{\scalebox{0.8}{0}}$ & $%
m_{\scalebox{0.8}{1}}$ \\ 
$m_{\scalebox{0.8}{7}}$ & $m_{\scalebox{0.8}{6}}$ & $m_{\scalebox{0.8}{5}}$ & $m_{\scalebox{0.8}{4}}$ & $-m_{\scalebox{0.8}{3}}$ & $-m_{\scalebox{0.8}{2}}$ & $-m_{\scalebox{0.8}{1}}$ & $%
-m_{\scalebox{0.8}{0}}$ \\ 
$-m_{\scalebox{0.8}{3}}$ & $-m_{\scalebox{0.8}{2}}$ & $m_{\scalebox{0.8}{1}}$ & $m_{\scalebox{0.8}{0}}$ & $m_{\scalebox{0.8}{7}}$ & $m_{\scalebox{0.8}{6}}$ & $-m_{\scalebox{0.8}{5}}$ & $%
-m_{\scalebox{0.8}{4}}$ \\
\end{tabular}%
\right) $
}

\bigskip

\emph{Proof.} \ By calculation we have the above expressions. \ \ \
\ \emph{Q.E.D.}

\bigskip

\emph{Lemma 6.25.} We have

$2\textrm{MJr}(m)Mv=$

{\fontsize{8pt}{8pt} \selectfont%
$\left( 
\begin{tabular}{@{}c@{}c@{}c@{}c@{}c@{}c@{}c@{}c@{}c@{}c@{}c@{}c@{}c@{}c@{}}
$-m_{\scalebox{0.8}{1}}$ & $-m_{\scalebox{0.8}{2}}$ & $-m_{\scalebox{0.8}{3}}$ & $-m_{\scalebox{0.8}{4}}$ & $-m_{\scalebox{0.8}{5}}$ & $-m_{\scalebox{0.8}{6}}$ & $-m_{\scalebox{0.8}{7}}$ & 
$m_{\scalebox{0.8}{3}}$ & $-m_{\scalebox{0.8}{2}}$ & $m_{\scalebox{0.8}{5}}$ & $-m_{\scalebox{0.8}{4}}$ & $m_{\scalebox{0.8}{7}}$ & $-m_{\scalebox{0.8}{6}}$ & $m_{\scalebox{0.8}{1}}$ \\ 
$m_{\scalebox{0.8}{0}}$ & $m_{\scalebox{0.8}{3}}$ & $-m_{\scalebox{0.8}{2}}$ & $m_{\scalebox{0.8}{5}}$ & $-m_{\scalebox{0.8}{4}}$ & $m_{\scalebox{0.8}{7}}$ & $-m_{\scalebox{0.8}{6}}$ & $%
m_{\scalebox{0.8}{2}}$ & $m_{\scalebox{0.8}{3}}$ & $m_{\scalebox{0.8}{4}}$ & $m_{\scalebox{0.8}{5}}$ & $m_{\scalebox{0.8}{6}}$ & $m_{\scalebox{0.8}{7}}$ & $-m_{\scalebox{0.8}{0}}$ \\ 
$-m_{\scalebox{0.8}{3}}$ & $m_{\scalebox{0.8}{0}}$ & $m_{\scalebox{0.8}{1}}$ & $-m_{\scalebox{0.8}{6}}$ & $m_{\scalebox{0.8}{7}}$ & $m_{\scalebox{0.8}{4}}$ & $-m_{\scalebox{0.8}{5}}$ & $%
-m_{\scalebox{0.8}{1}}$ & $m_{\scalebox{0.8}{0}}$ & $m_{\scalebox{0.8}{7}}$ & $m_{\scalebox{0.8}{6}}$ & $-m_{\scalebox{0.8}{5}}$ & $-m_{\scalebox{0.8}{4}}$ & $m_{\scalebox{0.8}{3}}$ \\ 
$m_{\scalebox{0.8}{2}}$ & $-m_{\scalebox{0.8}{1}}$ & $m_{\scalebox{0.8}{0}}$ & $m_{\scalebox{0.8}{7}}$ & $m_{\scalebox{0.8}{6}}$ & $-m_{\scalebox{0.8}{5}}$ & $-m_{\scalebox{0.8}{4}}$ & $%
-m_{\scalebox{0.8}{0}}$ & $-m_{\scalebox{0.8}{1}}$ & $m_{\scalebox{0.8}{6}}$ & $-m_{\scalebox{0.8}{7}}$ & $-m_{\scalebox{0.8}{4}}$ & $m_{\scalebox{0.8}{5}}$ & $-m_{\scalebox{0.8}{2}}$ \\ 
$-m_{\scalebox{0.8}{5}}$ & $m_{\scalebox{0.8}{6}}$ & $-m_{\scalebox{0.8}{7}}$ & $m_{\scalebox{0.8}{0}}$ & $m_{\scalebox{0.8}{1}}$ & $-m_{\scalebox{0.8}{2}}$ & $m_{\scalebox{0.8}{3}}$ & $%
-m_{\scalebox{0.8}{7}}$ & $-m_{\scalebox{0.8}{6}}$ & $-m_{\scalebox{0.8}{1}}$ & $m_{\scalebox{0.8}{0}}$ & $m_{\scalebox{0.8}{3}}$ & $m_{\scalebox{0.8}{2}}$ & $-m_{\scalebox{0.8}{5}}$ \\
\end{tabular}%
\right. $

\bigskip

$\left. 
\begin{tabular}{@{}c@{}c@{}c@{}c@{}c@{}c@{}c@{}c@{}c@{}c@{}c@{}c@{}c@{}c@{}}
$-m_{\scalebox{0.8}{6}}$ & $m_{\scalebox{0.8}{7}}$ & $m_{\scalebox{0.8}{4}}$ & $-m_{\scalebox{0.8}{5}}$ & $m_{\scalebox{0.8}{7}}$ & $m_{\scalebox{0.8}{6}}$ & $-m_{\scalebox{0.8}{5}}$ & $%
-m_{\scalebox{0.8}{4}}$ & $m_{\scalebox{0.8}{1}}$ & $-m_{\scalebox{0.8}{2}}$ & $m_{\scalebox{0.8}{3}}$ & $m_{\scalebox{0.8}{3}}$ & $m_{\scalebox{0.8}{2}}$ & $m_{\scalebox{0.8}{1}}$ \\ 
$-m_{\scalebox{0.8}{7}}$ & $-m_{\scalebox{0.8}{6}}$ & $m_{\scalebox{0.8}{5}}$ & $m_{\scalebox{0.8}{4}}$ & $-m_{\scalebox{0.8}{6}}$ & $m_{\scalebox{0.8}{7}}$ & $m_{\scalebox{0.8}{4}}$ & $%
-m_{\scalebox{0.8}{5}}$ & $-m_{\scalebox{0.8}{0}}$ & $-m_{\scalebox{0.8}{3}}$ & $-m_{\scalebox{0.8}{2}}$ & $-m_{\scalebox{0.8}{2}}$ & $m_{\scalebox{0.8}{3}}$ & $-m_{\scalebox{0.8}{0}}$ \\ 
$m_{\scalebox{0.8}{4}}$ & $m_{\scalebox{0.8}{5}}$ & $m_{\scalebox{0.8}{6}}$ & $m_{\scalebox{0.8}{7}}$ & $-m_{\scalebox{0.8}{5}}$ & $m_{\scalebox{0.8}{4}}$ & $-m_{\scalebox{0.8}{7}}$ & $%
m_{\scalebox{0.8}{6}}$ & $-m_{\scalebox{0.8}{3}}$ & $m_{\scalebox{0.8}{0}}$ & $m_{\scalebox{0.8}{1}}$ & $m_{\scalebox{0.8}{1}}$ & $-m_{\scalebox{0.8}{0}}$ & $-m_{\scalebox{0.8}{3}}$ \\ 
$m_{\scalebox{0.8}{5}}$ & $-m_{\scalebox{0.8}{4}}$ & $m_{\scalebox{0.8}{7}}$ & $-m_{\scalebox{0.8}{6}}$ & $m_{\scalebox{0.8}{4}}$ & $m_{\scalebox{0.8}{5}}$ & $m_{\scalebox{0.8}{6}}$ & $%
m_{\scalebox{0.8}{7}}$ & $m_{\scalebox{0.8}{2}}$ & $m_{\scalebox{0.8}{1}}$ & $-m_{\scalebox{0.8}{0}}$ & $-m_{\scalebox{0.8}{0}}$ & $-m_{\scalebox{0.8}{1}}$ & $m_{\scalebox{0.8}{2}}$ \\ 
$-m_{\scalebox{0.8}{2}}$ & $m_{\scalebox{0.8}{3}}$ & $-m_{\scalebox{0.8}{0}}$ & $-m_{\scalebox{0.8}{1}}$ & $-m_{\scalebox{0.8}{3}}$ & $-m_{\scalebox{0.8}{2}}$ & $-m_{\scalebox{0.8}{1}}$ & $%
m_{\scalebox{0.8}{0}}$ & $m_{\scalebox{0.8}{5}}$ & $m_{\scalebox{0.8}{6}}$ & $m_{\scalebox{0.8}{7}}$ & $-m_{\scalebox{0.8}{7}}$ & $m_{\scalebox{0.8}{6}}$ & $-m_{\scalebox{0.8}{5}}$ \\
\end{tabular}%
\right) $
}

\bigskip

$2\textrm{MJr}(m)Mv^{2}=$

{\fontsize{8pt}{8pt} \selectfont%
$\left( 
\begin{tabular}{@{}c@{}c@{}c@{}c@{}c@{}c@{}c@{}c@{}c@{}c@{}c@{}c@{}c@{}c@{}}
$-m_{\scalebox{0.8}{1}}$ & $-m_{\scalebox{0.8}{2}}$ & $-m_{\scalebox{0.8}{3}}$ & $-m_{\scalebox{0.8}{4}}$ & $-m_{\scalebox{0.8}{5}}$ & $-m_{\scalebox{0.8}{6}}$ & $-m_{\scalebox{0.8}{7}}$ & 
$-m_{\scalebox{0.8}{3}}$ & $m_{\scalebox{0.8}{2}}$ & $-m_{\scalebox{0.8}{5}}$ & $m_{\scalebox{0.8}{4}}$ & $-m_{\scalebox{0.8}{7}}$ & $m_{\scalebox{0.8}{6}}$ & $-m_{\scalebox{0.8}{1}}$ \\ 
$m_{\scalebox{0.8}{0}}$ & $-m_{\scalebox{0.8}{3}}$ & $m_{\scalebox{0.8}{2}}$ & $-m_{\scalebox{0.8}{5}}$ & $m_{\scalebox{0.8}{4}}$ & $-m_{\scalebox{0.8}{7}}$ & $m_{\scalebox{0.8}{6}}$ & $%
m_{\scalebox{0.8}{2}}$ & $m_{\scalebox{0.8}{3}}$ & $m_{\scalebox{0.8}{4}}$ & $m_{\scalebox{0.8}{5}}$ & $m_{\scalebox{0.8}{6}}$ & $m_{\scalebox{0.8}{7}}$ & $m_{\scalebox{0.8}{0}}$ \\ 
$m_{\scalebox{0.8}{3}}$ & $m_{\scalebox{0.8}{0}}$ & $-m_{\scalebox{0.8}{1}}$ & $m_{\scalebox{0.8}{6}}$ & $-m_{\scalebox{0.8}{7}}$ & $-m_{\scalebox{0.8}{4}}$ & $m_{\scalebox{0.8}{5}}$ & -$%
m_{\scalebox{0.8}{1}}$ & $-m_{\scalebox{0.8}{0}}$ & $m_{\scalebox{0.8}{7}}$ & $m_{\scalebox{0.8}{6}}$ & $-m_{\scalebox{0.8}{5}}$ & $-m_{\scalebox{0.8}{4}}$ & $m_{\scalebox{0.8}{3}}$ \\ 
$-m_{\scalebox{0.8}{2}}$ & $m_{\scalebox{0.8}{1}}$ & $m_{\scalebox{0.8}{0}}$ & $-m_{\scalebox{0.8}{7}}$ & $-m_{\scalebox{0.8}{6}}$ & $m_{\scalebox{0.8}{5}}$ & $m_{\scalebox{0.8}{4}}$ & $%
m_{\scalebox{0.8}{0}}$ & $-m_{\scalebox{0.8}{1}}$ & $m_{\scalebox{0.8}{6}}$ & $-m_{\scalebox{0.8}{7}}$ & $-m_{\scalebox{0.8}{4}}$ & $m_{\scalebox{0.8}{5}}$ & $-m_{\scalebox{0.8}{2}}$ \\ 
$m_{\scalebox{0.8}{5}}$ & $-m_{\scalebox{0.8}{6}}$ & $m_{\scalebox{0.8}{7}}$ & $m_{\scalebox{0.8}{0}}$ & $-m_{\scalebox{0.8}{1}}$ & $m_{\scalebox{0.8}{2}}$ & $-m_{\scalebox{0.8}{3}}$ & -$%
m_{\scalebox{0.8}{7}}$ & $-m_{\scalebox{0.8}{6}}$ & $-m_{\scalebox{0.8}{1}}$ & $-m_{\scalebox{0.8}{0}}$ & $m_{\scalebox{0.8}{3}}$ & $m_{\scalebox{0.8}{2}}$ & $-m_{\scalebox{0.8}{5}}$ \\
\end{tabular}%
\right. $

\bigskip 

$\left. 
\begin{tabular}{@{}c@{}c@{}c@{}c@{}c@{}c@{}c@{}c@{}c@{}c@{}c@{}c@{}c@{}c@{}}
$m_{\scalebox{0.8}{6}}$ & $-m_{\scalebox{0.8}{7}}$ & $-m_{\scalebox{0.8}{4}}$ & $m_{\scalebox{0.8}{5}}$ & $-m_{\scalebox{0.8}{7}}$ & $-m_{\scalebox{0.8}{6}}$ & $m_{\scalebox{0.8}{5}}$ & $%
m_{\scalebox{0.8}{4}}$ & $-m_{\scalebox{0.8}{1}}$ & $m_{\scalebox{0.8}{2}}$ & $-m_{\scalebox{0.8}{3}}$ & $-m_{\scalebox{0.8}{3}}$ & $-m_{\scalebox{0.8}{2}}$ & $-m_{\scalebox{0.8}{1}}$ \\ 
$-m_{\scalebox{0.8}{7}}$ & $-m_{\scalebox{0.8}{6}}$ & $m_{\scalebox{0.8}{5}}$ & $m_{\scalebox{0.8}{4}}$ & $-m_{\scalebox{0.8}{6}}$ & $m_{\scalebox{0.8}{7}}$ & $m_{\scalebox{0.8}{4}}$ & -$%
m_{\scalebox{0.8}{5}}$ & $m_{\scalebox{0.8}{0}}$ & $-m_{\scalebox{0.8}{3}}$ & $-m_{\scalebox{0.8}{2}}$ & $-m_{\scalebox{0.8}{2}}$ & $m_{\scalebox{0.8}{3}}$ & $m_{\scalebox{0.8}{0}}$ \\ 
$m_{\scalebox{0.8}{4}}$ & $m_{\scalebox{0.8}{5}}$ & $m_{\scalebox{0.8}{6}}$ & $m_{\scalebox{0.8}{7}}$ & $-m_{\scalebox{0.8}{5}}$ & $m_{\scalebox{0.8}{4}}$ & $-m_{\scalebox{0.8}{7}}$ & $%
m_{\scalebox{0.8}{6}}$ & $-m_{\scalebox{0.8}{3}}$ & $-m_{\scalebox{0.8}{0}}$ & $m_{\scalebox{0.8}{1}}$ & $m_{\scalebox{0.8}{1}}$ & $m_{\scalebox{0.8}{0}}$ & $-m_{\scalebox{0.8}{3}}$ \\ 
$m_{\scalebox{0.8}{5}}$ & $-m_{\scalebox{0.8}{4}}$ & $m_{\scalebox{0.8}{7}}$ & $-m_{\scalebox{0.8}{6}}$ & $m_{\scalebox{0.8}{4}}$ & $m_{\scalebox{0.8}{5}}$ & $m_{\scalebox{0.8}{6}}$ & $%
m_{\scalebox{0.8}{7}}$ & $m_{\scalebox{0.8}{2}}$ & $m_{\scalebox{0.8}{1}}$ & $m_{\scalebox{0.8}{0}}$ & $m_{\scalebox{0.8}{0}}$ & $-m_{\scalebox{0.8}{1}}$ & $m_{\scalebox{0.8}{2}}$ \\ 
$-m_{\scalebox{0.8}{2}}$ & $m_{\scalebox{0.8}{3}}$ & $m_{\scalebox{0.8}{0}}$ & $-m_{\scalebox{0.8}{1}}$ & $-m_{\scalebox{0.8}{3}}$ & $-m_{\scalebox{0.8}{2}}$ & $-m_{\scalebox{0.8}{1}}$ & -$%
m_{\scalebox{0.8}{0}}$ & $m_{\scalebox{0.8}{5}}$ & $m_{\scalebox{0.8}{6}}$ & $m_{\scalebox{0.8}{7}}$ & $-m_{\scalebox{0.8}{7}}$ & $m_{\scalebox{0.8}{6}}$ & $-m_{\scalebox{0.8}{5}}$ \\
\end{tabular}%
\right) $
}

\bigskip

\emph{\ Proof.} \ By calculation we have the above expressions. \ \ \
\ \emph{Q.E.D.}

\bigskip

\emph{Definition 6.26.} \ For $X(\chi ,x)\in $\gJ%
$^{},x_{i}=\sum\limits_{j=0}^{7}x_{ij}e_{j}$ $(x_{i}\in $\gC%
$^{},x_{ij},\chi _{i}\in \R),$

\noindent
we define $%
\textrm{M}_{D}\textrm{Jl},\textrm{M}_{D}(X),\textrm{MCC}(X)\textrm{MCCC}(X),\textrm{M}_{D}\textrm{C}(X),\textrm{M}_{R}\textrm{Dl}(X),$

\noindent
$\textrm{M}_{C}\textrm{Jr}(X),\textrm{MI}(X),\textrm{M}^{+}\textrm{I}(X),\textrm{M}_{D}\textrm{E}(X)$,
$\textrm{M}_{D}^{-}\textrm{E}(X),\textrm{M}_{D}^{+}\textrm{E}(X),\textrm{MJC}(X),$

\noindent
$\textrm{M}^{+}\textrm{JC}(X),\textrm{M}_{C}\textrm{CJ}(X),\textrm{M}_{R}\textrm{JC}(X), \textrm{M}_{C}\textrm{CCC}(X)$ respectively as follows.

$\textrm{M}_{D}\textrm{Jl}(D)=\left( 
\begin{array}{ccc}
\textrm{MJl}(D) & 0 & 0 \\ 
0 & \textrm{MJl}(MvD) & 0 \\ 
0 & 0 & \textrm{MJl}(Mv^{2}D)%
\end{array}%
\right) \in M(24\times 24,\R),$

$\textrm{M}_{D}(X)=\left( 
\begin{array}{ccc}
\chi _{1} & 0 & 0 \\ 
0 & \chi _{2} & 0 \\ 
0 & 0 & \chi _{3}%
\end{array}%
\right) \in M(3\times 3,\R),$

$M_{D3}(X)=\left( 
\begin{array}{ccc}
\chi _{1} & 0 & 0 \\ 
0 & \chi _{2} & 0%
\end{array}%
\right) \in M(2\times 3,\R),$

$\textrm{M}_{D}^{3}(X)=\left( 
\begin{array}{cc}
\chi _{1} & 0 \\ 
0 & \chi _{2} \\ 
-\chi _{3} & -\chi _{3}%
\end{array}%
\right) \in M(3\times 2,\R),$

$\textrm{MCC}(X)=\left( 
\begin{array}{ccc}
0 & \chi _{3} & \chi _{2} \\ 
\chi _{3} & 0 & \chi _{1} \\ 
\chi _{2} & \chi _{1} & 0%
\end{array}%
\right) \in M(3\times 3,\R),$

$\textrm{MCCC}(X)=\left( 
\begin{array}{ccc}
0 & \chi _{2} & \chi _{3} \\ 
\chi _{1} & 0 & \chi _{3} \\ 
\chi _{1} & \chi _{2} & 0%
\end{array}%
\right) \in M(3\times 3,\R),$

$\textrm{MCCC}_{3}(X)=\left( 
\begin{array}{ccc}
0 & \chi _{2} & \chi _{3} \\ 
\chi _{1} & 0 & \chi _{3}%
\end{array}%
\right) \in M(2\times 3,\R),$

$\textrm{M}_{D}\textrm{C}(X)=\left( 
\begin{array}{ccc}
\textrm{MC}(x_{1}) & 0 & 0 \\ 
0 & \textrm{MC}(x_{2}) & 0 \\ 
0 & 0 & \textrm{MC}(x_{3})%
\end{array}%
\right) \in M(3\times 24,\R)$

$=$
{\fontsize{8pt}{10pt} \selectfont%
$\left( 
\begin{tabular}{@{}c@{}c@{}c@{}c@{}c@{}c@{}c@{}c@{}c@{}c@{}c@{}c@{}c@{}c@{}c@{}c@{}c@{}c@{}c@{}c@{}c@{}c@{}c@{}c@{}}
$x_{\scalebox{0.8}{10}}$ & $x_{\scalebox{0.8}{11}}$ & $x_{\scalebox{0.8}{12}}$ & $x_{\scalebox{0.8}{13}}$ & $x_{\scalebox{0.8}{14}}$ & $x_{\scalebox{0.8}{15}}$ & $x_{\scalebox{0.8}{16}}$ & 
$x_{\scalebox{0.8}{17}}$ & $0$ & $0$ & $0$ & $0$ & $0$ & $0$ & $0$ & $0$ & $0$ & $0$ & $0$ & 
$0$ & $0$ & $0$ & $0$ & $0$ \\ 
$0$ & $0$ & $0$ & $0$ & $0$ & $0$ & $0$ & $0$ & $x_{\scalebox{0.8}{20}}$ & $x_{\scalebox{0.8}{21}}$ & $x_{\scalebox{0.8}{22}}
$ & $x_{\scalebox{0.8}{23}}$ & $x_{\scalebox{0.8}{24}}$ & $x_{\scalebox{0.8}{25}}$ & $x_{\scalebox{0.8}{26}}$ & $x_{\scalebox{0.8}{27}}$ & $0$ & $0$ & $0$ & 
$0$ & $0$ & $0$ & $0$ & $0$ \\ 
$0$ & $0$ & $0$ & $0$ & $0$ & $0$ & $0$ & $0$ & $0$ & $0$ & $0$ & $0$ & $0$
& $0$ & $0$ & $0$ & $x_{\scalebox{0.8}{30}}$ & $x_{\scalebox{0.8}{31}}$ & $x_{\scalebox{0.8}{32}}$ & $x_{\scalebox{0.8}{33}}$ & $x_{\scalebox{0.8}{34}}$ & $%
x_{\scalebox{0.8}{35}}$ & $x_{\scalebox{0.8}{36}}$ & $x_{\scalebox{0.8}{37}}$%
\end{tabular}%
\right) $
}

\bigskip

$\textrm{M}_{D}\textrm{C}_{3}(X)=\left( 
\begin{array}{ccc}
\textrm{MC}(x_{1}) & 0 & 0 \\ 
0 & \textrm{MC}(x_{2}) & 0%
\end{array}%
\right) \in M(2\times 24,\R),$

$2\textrm{M}_{R}\textrm{Dl}(X)$

$=2\left( 
\begin{array}{ccc}
\textrm{MDl}(x_{1}) & Mv^{2}\textrm{MDl}(x_{2}) & Mv\textrm{MDl}(x_{3})%
\end{array}%
\right) \in M(28\times 24,\R)$

$=$

\begin{flushleft}
{\fontsize{4pt}{8pt} \selectfont%
$\left( 
\begin{tabular}{@{}c@{}c@{}c@{}c@{}c@{}c@{}c@{}c@{}c@{}c@{}c@{}c@{}c@{}c@{}c@{}c@{}c@{}c@{}c@{}c@{}c@{}c@{}c@{}c@{}}
-$2x_{\scalebox{0.8}{11}}$ & $2x_{\scalebox{0.8}{10}}$ & $0$ & $0$ & $0$ & $0$ & $0$ & $0$ & $x_{\scalebox{0.8}{21}}$ & $%
-x_{\scalebox{0.8}{20}}$ & $x_{\scalebox{0.8}{23}}$ & -$x_{\scalebox{0.8}{22}}$ & $x_{\scalebox{0.8}{25}}$ & -$x_{\scalebox{0.8}{24}}$ & $x_{\scalebox{0.8}{27}}$ & -$x_{\scalebox{0.8}{26}}$
& $x_{\scalebox{0.8}{31}}$ & -$x_{\scalebox{0.8}{30}}$ & -$x_{\scalebox{0.8}{33}}$ & $x_{\scalebox{0.8}{32}}$ & -$x_{\scalebox{0.8}{35}}$ & $x_{\scalebox{0.8}{34}}$ & %
-$x_{\scalebox{0.8}{37}}$ & $x_{\scalebox{0.8}{36}}$ \\ 
-$2x_{\scalebox{0.8}{12}}$ & $0$ & $2x_{\scalebox{0.8}{10}}$ & $0$ & $0$ & $0$ & $0$ & $0$ & $x_{\scalebox{0.8}{22}}$ & $%
-x_{\scalebox{0.8}{23}}$ & -$x_{\scalebox{0.8}{20}}$ & $x_{\scalebox{0.8}{21}}$ & -$x_{\scalebox{0.8}{26}}$ & $x_{\scalebox{0.8}{27}}$ & $x_{\scalebox{0.8}{24}}$ & -$x_{\scalebox{0.8}{25}}$
& $x_{\scalebox{0.8}{32}}$ & $x_{\scalebox{0.8}{33}}$ & -$x_{\scalebox{0.8}{30}}$ & -$x_{\scalebox{0.8}{31}}$ & $x_{\scalebox{0.8}{36}}$ & -$x_{\scalebox{0.8}{37}}$ & %
-$x_{\scalebox{0.8}{34}}$ & $x_{\scalebox{0.8}{35}}$ \\ 
-$2x_{\scalebox{0.8}{13}}$ & $0$ & $0$ & $2x_{\scalebox{0.8}{10}}$ & $0$ & $0$ & $0$ & $0$ & $x_{\scalebox{0.8}{23}}$ & $%
x_{\scalebox{0.8}{22}}$ & -$x_{\scalebox{0.8}{21}}$ & -$x_{\scalebox{0.8}{20}}$ & $x_{\scalebox{0.8}{27}}$ & $x_{\scalebox{0.8}{26}}$ & -$x_{\scalebox{0.8}{25}}$ & -$x_{\scalebox{0.8}{24}}$
& $x_{\scalebox{0.8}{33}}$ & -$x_{\scalebox{0.8}{32}}$ & $x_{\scalebox{0.8}{31}}$ & -$x_{\scalebox{0.8}{30}}$ & -$x_{\scalebox{0.8}{37}}$ & -$x_{\scalebox{0.8}{36}}$ & $%
x_{\scalebox{0.8}{35}}$ & $x_{\scalebox{0.8}{34}}$ \\
\end{tabular}%
\right) $
}
\end{flushleft}

\bigskip

$\textrm{M}_{C}\textrm{Jr}(X)=\left( 
\begin{array}{c}
\textrm{MJr}(x_{1}) \\ 
\textrm{MJr}(x_{2})Mv \\ 
\textrm{MJr}(x_{3})Mv^{2}%
\end{array}%
\right) =-\ ^{t}\textrm{M}_{R}\textrm{Dl}(X)\in M(24\times 28,\R)$

\bigskip

$\textrm{MI}(X)=\left( 
\begin{array}{ccc}
0 & \textrm{MIr}(x_{3}) & -\textrm{MIl}(x_{2}) \\ 
-\textrm{MIl}(x_{3}) & 0 & \textrm{MIr}(x_{1}) \\ 
\textrm{MIr}(x_{2}) & -\textrm{MIl}(x_{1}) & 0%
\end{array}%
\right) \in M(24\times 24,\R)$

$=$

\begin{flushleft}
{\fontsize{6pt}{8pt} \selectfont%
$\left( 
\begin{tabular}{@{}c@{}c@{}c@{}c@{}c@{}c@{}c@{}c@{}c@{}c@{}c@{}c@{}c@{}c@{}c@{}c@{}c@{}c@{}c@{}c@{}c@{}c@{}c@{}c@{}}
$0$ & $0$ & $0$ & $0$ & $0$ & $0$ & $0$ & $0$ & $x_{\scalebox{0.8}{30}}$ & $x_{\scalebox{0.8}{31}}$ & $x_{\scalebox{0.8}{32}}
$ & $x_{\scalebox{0.8}{33}}$ & $x_{\scalebox{0.8}{34}}$ & $x_{\scalebox{0.8}{35}}$ & $x_{\scalebox{0.8}{36}}$ & $x_{\scalebox{0.8}{37}}$ & -$x_{\scalebox{0.8}{20}}$ & %
-$x_{\scalebox{0.8}{21}}$ & -$x_{\scalebox{0.8}{22}}$ & -$x_{\scalebox{0.8}{23}}$ & -$x_{\scalebox{0.8}{24}}$ & -$x_{\scalebox{0.8}{25}}$ & -$x_{\scalebox{0.8}{26}}$ & %
-$x_{\scalebox{0.8}{27}}$ \\ 
$0$ & $0$ & $0$ & $0$ & $0$ & $0$ & $0$ & $0$ & -$x_{\scalebox{0.8}{31}}$ & $x_{\scalebox{0.8}{30}}$ & %
--$x_{\scalebox{0.8}{33}}$ & $x_{\scalebox{0.8}{32}}$ & -$x_{\scalebox{0.8}{35}}$ & $x_{\scalebox{0.8}{34}}$ & -$x_{\scalebox{0.8}{37}}$ & $x_{\scalebox{0.8}{36}}$ & $x_{\scalebox{0.8}{21}}$
& -$x_{\scalebox{0.8}{20}}$ & -$x_{\scalebox{0.8}{23}}$ & $x_{\scalebox{0.8}{22}}$ & -$x_{\scalebox{0.8}{25}}$ & $x_{\scalebox{0.8}{24}}$ & -$x_{\scalebox{0.8}{27}}$ & $%
x_{\scalebox{0.8}{26}}$ \\ 
$0$ & $0$ & $0$ & $0$ & $0$ & $0$ & $0$ & $0$ & -$x_{\scalebox{0.8}{32}}$ & $x_{\scalebox{0.8}{33}}$ & $%
x_{\scalebox{0.8}{30}}$ & -$x_{\scalebox{0.8}{31}}$ & $x_{\scalebox{0.8}{36}}$ & -$x_{\scalebox{0.8}{37}}$ & -$x_{\scalebox{0.8}{34}}$ & $x_{\scalebox{0.8}{35}}$ & $x_{\scalebox{0.8}{22}}$
& $x_{\scalebox{0.8}{23}}$ & -$x_{\scalebox{0.8}{20}}$ & -$x_{\scalebox{0.8}{21}}$ & $x_{\scalebox{0.8}{26}}$ & -$x_{\scalebox{0.8}{27}}$ & -$x_{\scalebox{0.8}{24}}$ & $%
x_{\scalebox{0.8}{25}}$ \\
\end{tabular}%
\right) $
}
\end{flushleft}
\bigskip

$\textrm{M}^{+}\textrm{I}(X)=\left( 
\begin{array}{ccc}
0 & \textrm{MIr}(x_{3}) & \textrm{MIl}(x_{2}) \\ 
\textrm{MIl}(x_{3}) & 0 & \textrm{MIr}(x_{1}) \\ 
\textrm{MIr}(x_{2}) & \textrm{MIl}(x_{1}) & 0%
\end{array}%
\right) \in M(24\times 24,\R)$

$=$

\begin{flushleft}
{\fontsize{6pt}{8pt} \selectfont%
$\left( 
\begin{tabular}{@{}c@{}c@{}c@{}c@{}c@{}c@{}c@{}c@{}c@{}c@{}c@{}c@{}c@{}c@{}c@{}c@{}c@{}c@{}c@{}c@{}c@{}c@{}c@{}c@{}}
$0$ & $0$ & $0$ & $0$ & $0$ & $0$ & $0$ & $0$ & $x_{\scalebox{0.8}{30}}$ & $x_{\scalebox{0.8}{31}}$ & $x_{\scalebox{0.8}{32}}
$ & $x_{\scalebox{0.8}{33}}$ & $x_{\scalebox{0.8}{34}}$ & $x_{\scalebox{0.8}{35}}$ & $x_{\scalebox{0.8}{36}}$ & $x_{\scalebox{0.8}{37}}$ & $x_{\scalebox{0.8}{20}}$ & $x_{\scalebox{0.8}{21}}
$ & $x_{\scalebox{0.8}{22}}$ & $x_{\scalebox{0.8}{23}}$ & $x_{\scalebox{0.8}{24}}$ & $x_{\scalebox{0.8}{25}}$ & $x_{\scalebox{0.8}{26}}$ & $x_{\scalebox{0.8}{27}}$ \\ 
$0$ & $0$ & $0$ & $0$ & $0$ & $0$ & $0$ & $0$ & -$x_{\scalebox{0.8}{31}}$ & $x_{\scalebox{0.8}{30}}$ & %
-$x_{\scalebox{0.8}{33}}$ & $x_{\scalebox{0.8}{32}}$ & -$x_{\scalebox{0.8}{35}}$ & $x_{\scalebox{0.8}{34}}$ & -$x_{\scalebox{0.8}{37}}$ & $x_{\scalebox{0.8}{36}}$ & -$x_{\scalebox{0.8}{21}}$
& $x_{\scalebox{0.8}{20}}$ & $x_{\scalebox{0.8}{23}}$ & -$x_{\scalebox{0.8}{22}}$ & $x_{\scalebox{0.8}{25}}$ & -$x_{\scalebox{0.8}{24}}$ & $x_{\scalebox{0.8}{27}}$ & %
-$x_{\scalebox{0.8}{26}}$ \\ 
$0$ & $0$ & $0$ & $0$ & $0$ & $0$ & $0$ & $0$ & -$x_{\scalebox{0.8}{32}}$ & $x_{\scalebox{0.8}{33}}$ & $%
x_{\scalebox{0.8}{30}}$ & -$x_{\scalebox{0.8}{31}}$ & $x_{\scalebox{0.8}{36}}$ & -$x_{\scalebox{0.8}{37}}$ & -$x_{\scalebox{0.8}{34}}$ & $x_{\scalebox{0.8}{35}}$ & -$x_{\scalebox{0.8}{22}}$
& -$x_{\scalebox{0.8}{23}}$ & $x_{\scalebox{0.8}{20}}$ & $x_{\scalebox{0.8}{21}}$ & -$x_{\scalebox{0.8}{26}}$ & $x_{\scalebox{0.8}{27}}$ & $x_{\scalebox{0.8}{24}}$ & %
-$x_{\scalebox{0.8}{25}}$ \\
\end{tabular}%
\right) $
}
\end{flushleft}

\bigskip

$\textrm{M}_{D}\textrm{E}(X)=\left( 
\begin{array}{ccc}
\textrm{ME}(\chi _{1}) & 0 & 0 \\ 
0 & \textrm{ME}(\chi _{2}) & 0 \\ 
0 & 0 & \textrm{ME}(\chi _{3})%
\end{array}%
\right) $

$\ \ \ \ \ \ \ \ \ \ \ \ \ \ =\left( 
\begin{array}{ccc}
\chi _{1}E & 0 & 0 \\ 
0 & \chi _{2}E & 0 \\ 
0 & 0 & \chi _{3}E%
\end{array}%
\right) \in M(24\times 24,\R),$

$\textrm{M}_{D}^{-}\textrm{E}(X)=\left( 
\begin{array}{ccc}
\textrm{ME}(\chi _{2}-\chi _{3}) & 0 & 0 \\ 
0 & \textrm{ME}(\chi _{3}-\chi _{1}) & 0 \\ 
0 & 0 & \textrm{ME}(\chi _{1}-\chi _{2})%
\end{array}%
\right) $

$\ \ \ \ \ \ \ \ \ \ \ \ \ \ =\left( 
\begin{array}{ccc}
(\chi _{2}-\chi _{3})E & 0 & 0 \\ 
0 & (\chi _{3}-\chi _{1})E & 0 \\ 
0 & 0 & (\chi _{1}-\chi _{2})E%
\end{array}%
\right) \in M(24\times 24,\R),$

$\textrm{M}_{D}^{+}\textrm{E}(X)=\left( 
\begin{array}{ccc}
\textrm{ME}(\chi _{2}+\chi _{3}) & 0 & 0 \\ 
0 & \textrm{ME}(\chi _{3}+\chi _{1}) & 0 \\ 
0 & 0 & \textrm{ME}(\chi _{1}+\chi _{2})%
\end{array}%
\right) $

$\ \ \ \ \ \ \ \ \ \ \ \ \ \ =\left( 
\begin{array}{ccc}
(\chi _{2}+\chi _{3})E & 0 & 0 \\ 
0 & (\chi _{3}+\chi _{1})E & 0 \\ 
0 & 0 & (\chi _{1}+\chi _{2})E%
\end{array}%
\right) \in M(24\times 24,\R),$

$\textrm{MJC}(X)=\left( 
\begin{array}{ccc}
0 & \textrm{MC}(x_{2}) & -\textrm{MC}(x_{3}) \\ 
-\textrm{MC}(x_{1}) & 0 & \textrm{MC}(x_{3}) \\ 
\textrm{MC}(x_{1}) & -\textrm{MC}(x_{2}) & 0%
\end{array}%
\right) \in M(3\times 24,\R)$

$=$
\begin{flushleft}
{\fontsize{6pt}{10pt} \selectfont%
$\left( 
\begin{tabular}{@{}c@{}c@{}c@{}c@{}c@{}c@{}c@{}c@{}c@{}c@{}c@{}c@{}c@{}c@{}c@{}c@{}c@{}c@{}c@{}c@{}c@{}c@{}c@{}c@{}}
$0$ & $0$ & $0$ & $0$ & $0$ & $0$ & $0$ & $0$ & $x_{\scalebox{0.8}{20}}$ & $x_{\scalebox{0.8}{21}}$ & $x_{\scalebox{0.8}{22}}
$ & $x_{\scalebox{0.8}{23}}$ & $x_{\scalebox{0.8}{24}}$ & $x_{\scalebox{0.8}{25}}$ & $x_{\scalebox{0.8}{26}}$ & $x_{\scalebox{0.8}{27}}$ & -$x_{\scalebox{0.8}{30}}$ & %
-$x_{\scalebox{0.8}{31}}$ & -$x_{\scalebox{0.8}{32}}$ & -$x_{\scalebox{0.8}{33}}$ & -$x_{\scalebox{0.8}{34}}$ & -$x_{\scalebox{0.8}{35}}$ & -$x_{\scalebox{0.8}{36}}$ & %
-$x_{\scalebox{0.8}{37}}$ \\ 
-$x_{\scalebox{0.8}{10}}$ & -$x_{\scalebox{0.8}{11}}$ & -$x_{\scalebox{0.8}{12}}$ & -$x_{\scalebox{0.8}{13}}$ & -$x_{\scalebox{0.8}{14}}$ & -$x_{\scalebox{0.8}{15}}$ & %
-$x_{\scalebox{0.8}{16}}$ & -$x_{\scalebox{0.8}{17}}$ & $0$ & $0$ & $0$ & $0$ & $0$ & $0$ & $0$ & $0$ & $%
x_{\scalebox{0.8}{30}}$ & $x_{\scalebox{0.8}{31}}$ & $x_{\scalebox{0.8}{32}}$ & $x_{\scalebox{0.8}{33}}$ & $x_{\scalebox{0.8}{34}}$ & $x_{\scalebox{0.8}{35}}$ & $x_{\scalebox{0.8}{36}}$ & $%
x_{\scalebox{0.8}{37}}$ \\ 
$x_{\scalebox{0.8}{10}}$ & $x_{\scalebox{0.8}{11}}$ & $x_{\scalebox{0.8}{12}}$ & $x_{\scalebox{0.8}{13}}$ & $x_{\scalebox{0.8}{14}}$ & $x_{\scalebox{0.8}{15}}$ & $x_{\scalebox{0.8}{16}}$ & 
$x_{\scalebox{0.8}{17}}$ & -$x_{\scalebox{0.8}{20}}$ & -$x_{\scalebox{0.8}{21}}$ & -$x_{\scalebox{0.8}{22}}$ & -$x_{\scalebox{0.8}{23}}$ & -$x_{\scalebox{0.8}{24}}$ & %
-$x_{\scalebox{0.8}{25}}$ & -$x_{\scalebox{0.8}{26}}$ & -$x_{\scalebox{0.8}{27}}$ & $0$ & $0$ & $0$ & $0$ & $0$ & $0$ & $0$ & 
$0$%
\end{tabular}%
\right) $
}
\end{flushleft}

\bigskip

$\textrm{MJC}_{3}(X)=\left( 
\begin{array}{ccc}
0 & \textrm{MC}(x_{2}) & -\textrm{MC}(x_{3}) \\ 
-\textrm{MC}(x_{1}) & 0 & \textrm{MC}(x_{3})%
\end{array}%
\right) \in M(2\times 24,\R)$

$=$

\begin{flushleft}
{\fontsize{6pt}{10pt} \selectfont%
$\left( 
\begin{tabular}{@{}c@{}c@{}c@{}c@{}c@{}c@{}c@{}c@{}c@{}c@{}c@{}c@{}c@{}c@{}c@{}c@{}c@{}c@{}c@{}c@{}c@{}c@{}c@{}c@{}}
$0$ & $0$ & $0$ & $0$ & $0$ & $0$ & $0$ & $0$ & $x_{\scalebox{0.8}{20}}$ & $x_{\scalebox{0.8}{21}}$ & $x_{\scalebox{0.8}{22}}
$ & $x_{\scalebox{0.8}{23}}$ & $x_{\scalebox{0.8}{24}}$ & $x_{\scalebox{0.8}{25}}$ & $x_{\scalebox{0.8}{26}}$ & $x_{\scalebox{0.8}{27}}$ & -$x_{\scalebox{0.8}{30}}$ & %
-$x_{\scalebox{0.8}{31}}$ & -$x_{\scalebox{0.8}{32}}$ & -$x_{\scalebox{0.8}{33}}$ & -$x_{\scalebox{0.8}{34}}$ & -$x_{\scalebox{0.8}{35}}$ & -$x_{\scalebox{0.8}{36}}$ & %
-$x_{\scalebox{0.8}{37}}$ \\ 
-$x_{\scalebox{0.8}{10}}$ & -$x_{\scalebox{0.8}{11}}$ & -$x_{\scalebox{0.8}{12}}$ & -$x_{\scalebox{0.8}{13}}$ & -$x_{\scalebox{0.8}{14}}$ & -$x_{\scalebox{0.8}{15}}$ & %
-$x_{\scalebox{0.8}{16}}$ & -$x_{\scalebox{0.8}{17}}$ & $0$ & $0$ & $0$ & $0$ & $0$ & $0$ & $0$ & $0$ & $%
x_{\scalebox{0.8}{30}}$ & $x_{\scalebox{0.8}{31}}$ & $x_{\scalebox{0.8}{32}}$ & $x_{\scalebox{0.8}{33}}$ & $x_{\scalebox{0.8}{34}}$ & $x_{\scalebox{0.8}{35}}$ & $x_{\scalebox{0.8}{36}}$ & $%
x_{\scalebox{0.8}{37}}$%
\end{tabular}%
\right) $
}
\end{flushleft}

\bigskip

$^{t}\textrm{MJC}^{3}(X)$=$\left( 
\begin{array}{cc}
-^{t}\textrm{MC}(x_{1}) & -2^{t}\textrm{MC}(x_{1}) \\ 
2^{t}\textrm{MC}(x_{2}) & ^{t}\textrm{MC}(x_{2}) \\ 
-^{t}\textrm{MC}(x_{3}) & ^{t}\textrm{MC}(x_{3})%
\end{array}%
\right) \in M(24\times 2,\R)$

$=$
{\fontsize{8pt}{8pt} \selectfont%
$\left( 
\begin{tabular}{@{}c@{}c@{}}
$-x_{\scalebox{0.8}{10}}$ & $-2x_{\scalebox{0.8}{10}}$ \\ 
$-x_{\scalebox{0.8}{11}}$ & $-2x_{\scalebox{0.8}{11}}$ \\ 
$-x_{\scalebox{0.8}{12}}$ & $-2x_{\scalebox{0.8}{12}}$ \\ 
$-x_{\scalebox{0.8}{13}}$ & $-2x_{\scalebox{0.8}{13}}$ \\ 
$-x_{\scalebox{0.8}{14}}$ & $-2x_{\scalebox{0.8}{14}}$ \\ 
$-x_{\scalebox{0.8}{15}}$ & $-2x_{\scalebox{0.8}{15}}$ \\ 
$-x_{\scalebox{0.8}{16}}$ & $-2x_{\scalebox{0.8}{16}}$ \\ 
$-x_{\scalebox{0.8}{17}}$ & $-2x_{\scalebox{0.8}{17}}$ \\ 
$2x_{\scalebox{0.8}{20}}$ & $x_{\scalebox{0.8}{20}}$ \\ 
$2x_{\scalebox{0.8}{21}}$ & $x_{\scalebox{0.8}{21}}$ \\ 
$2x_{\scalebox{0.8}{22}}$ & $x_{\scalebox{0.8}{22}}$ \\ 
$2x_{\scalebox{0.8}{23}}$ & $x_{\scalebox{0.8}{23}}$ \\ 
$2x_{\scalebox{0.8}{24}}$ & $x_{\scalebox{0.8}{24}}$ \\ 
$2x_{\scalebox{0.8}{25}}$ & $x_{\scalebox{0.8}{25}}$ \\ 
$2x_{\scalebox{0.8}{26}}$ & $x_{\scalebox{0.8}{26}}$ \\ 
$2x_{\scalebox{0.8}{27}}$ & $x_{\scalebox{0.8}{27}}$ \\ 
$-x_{\scalebox{0.8}{30}}$ & $x_{\scalebox{0.8}{30}}$ \\ 
$-x_{\scalebox{0.8}{31}}$ & $x_{\scalebox{0.8}{31}}$ \\ 
$-x_{\scalebox{0.8}{32}}$ & $x_{\scalebox{0.8}{32}}$ \\ 
$-x_{\scalebox{0.8}{33}}$ & $x_{\scalebox{0.8}{33}}$ \\ 
$-x_{\scalebox{0.8}{34}}$ & $x_{\scalebox{0.8}{34}}$ \\ 
$-x_{\scalebox{0.8}{35}}$ & $x_{\scalebox{0.8}{35}}$ \\ 
$-x_{\scalebox{0.8}{36}}$ & $x_{\scalebox{0.8}{36}}$ \\ 
$-x_{\scalebox{0.8}{37}}$ & $x_{\scalebox{0.8}{37}}$%
\end{tabular}%
\right) $
}

\bigskip

$\textrm{M}^{+}\textrm{JC}(X)=\left( 
\begin{array}{ccc}
0 & \textrm{MC}(x_{2}) & \textrm{MC}(x_{3}) \\ 
\textrm{MC}(x_{1}) & 0 & \textrm{MC}(x_{3}) \\ 
\textrm{MC}(x_{1}) & \textrm{MC}(x_{2}) & 0%
\end{array}%
\right) \in M(3\times 24,\R)$

$=$

\begin{flushleft}
{\fontsize{8pt}{10pt} \selectfont%
$\left( 
\begin{tabular}{@{}c@{}c@{}c@{}c@{}c@{}c@{}c@{}c@{}c@{}c@{}c@{}c@{}c@{}c@{}c@{}c@{}c@{}c@{}c@{}c@{}c@{}c@{}c@{}c@{}}
$0$ & $0$ & $0$ & $0$ & $0$ & $0$ & $0$ & $0$ & $x_{\scalebox{0.8}{20}}$ & $x_{\scalebox{0.8}{21}}$ & $x_{\scalebox{0.8}{22}}
$ & $x_{\scalebox{0.8}{23}}$ & $x_{\scalebox{0.8}{24}}$ & $x_{\scalebox{0.8}{25}}$ & $x_{\scalebox{0.8}{26}}$ & $x_{\scalebox{0.8}{27}}$ & $x_{\scalebox{0.8}{30}}$ & $x_{\scalebox{0.8}{31}}
$ & $x_{\scalebox{0.8}{32}}$ & $x_{\scalebox{0.8}{33}}$ & $x_{\scalebox{0.8}{34}}$ & $x_{\scalebox{0.8}{35}}$ & $x_{\scalebox{0.8}{36}}$ & $x_{\scalebox{0.8}{37}}$ \\ 
$x_{\scalebox{0.8}{10}}$ & $x_{\scalebox{0.8}{11}}$ & $x_{\scalebox{0.8}{12}}$ & $x_{\scalebox{0.8}{13}}$ & $x_{\scalebox{0.8}{14}}$ & $x_{\scalebox{0.8}{15}}$ & $x_{\scalebox{0.8}{16}}$ & 
$x_{\scalebox{0.8}{17}}$ & $0$ & $0$ & $0$ & $0$ & $0$ & $0$ & $0$ & $0$ & $x_{\scalebox{0.8}{30}}$ & $x_{\scalebox{0.8}{31}}
$ & $x_{\scalebox{0.8}{32}}$ & $x_{\scalebox{0.8}{33}}$ & $x_{\scalebox{0.8}{34}}$ & $x_{\scalebox{0.8}{35}}$ & $x_{\scalebox{0.8}{36}}$ & $x_{\scalebox{0.8}{37}}$ \\ 
$x_{\scalebox{0.8}{10}}$ & $x_{\scalebox{0.8}{11}}$ & $x_{\scalebox{0.8}{12}}$ & $x_{\scalebox{0.8}{13}}$ & $x_{\scalebox{0.8}{14}}$ & $x_{\scalebox{0.8}{15}}$ & $x_{\scalebox{0.8}{16}}$ & 
$x_{\scalebox{0.8}{17}}$ & $x_{\scalebox{0.8}{20}}$ & $x_{\scalebox{0.8}{21}}$ & $x_{\scalebox{0.8}{22}}$ & $x_{\scalebox{0.8}{23}}$ & $x_{\scalebox{0.8}{24}}$ & $x_{\scalebox{0.8}{25}}$ & 
$x_{\scalebox{0.8}{26}}$ & $x_{\scalebox{0.8}{27}}$ & $0$ & $0$ & $0$ & $0$ & $0$ & $0$ & $0$ & $0$%
\end{tabular}%
\right) $
}
\end{flushleft}

\bigskip

$\textrm{M}^{+}\textrm{JC}_{3}(X)=\left( 
\begin{array}{ccc}
0 & \textrm{MC}(x_{2}) & \textrm{MC}(x_{3}) \\ 
\textrm{MC}(x_{1}) & 0 & \textrm{MC}(x_{3})%
\end{array}%
\right) \in M(2\times 24,\R),$

$^{t}\textrm{M}^{+}\textrm{JC}^{3}(X)=\left( 
\begin{array}{cc}
-\text{ }^{t}\textrm{MC}(x_{1}) & 0 \\ 
\text{0} & -\text{ }^{\text{t}}\textrm{MC}(x_{2}) \\ 
^{t}\textrm{MC}(x_{3}) & ^{t}\textrm{MC}(x_{3})%
\end{array}%
\right) \in M(24\times 2,\R)$

$\textrm{M}_{R}\textrm{JC}(X)=\left( 
\begin{array}{ccc}
\textrm{MC}(x_{1}) & \textrm{MC}(x_{2}) & \textrm{MC}(x_{3})%
\end{array}%
\right) \in M(1\times 24,\R)$

$=$

\begin{flushleft}
{\fontsize{8pt}{8pt} \selectfont%
$\left( 
\begin{tabular}{@{}c@{}c@{}c@{}c@{}c@{}c@{}c@{}c@{}c@{}c@{}c@{}c@{}c@{}c@{}c@{}c@{}c@{}c@{}c@{}c@{}c@{}c@{}c@{}c@{}}
$x_{\scalebox{0.8}{10}}$ & $x_{\scalebox{0.8}{11}}$ & $x_{\scalebox{0.8}{12}}$ & $x_{\scalebox{0.8}{13}}$ & $x_{\scalebox{0.8}{14}}$ & $x_{\scalebox{0.8}{15}}$ & $x_{\scalebox{0.8}{16}}$ & 
$x_{\scalebox{0.8}{17}}$ & $x_{\scalebox{0.8}{20}}$ & $x_{\scalebox{0.8}{21}}$ & $x_{\scalebox{0.8}{22}}$ & $x_{\scalebox{0.8}{23}}$ & $x_{\scalebox{0.8}{24}}$ & $x_{\scalebox{0.8}{25}}$ & 
$x_{\scalebox{0.8}{26}}$ & $x_{\scalebox{0.8}{27}}$ & $x_{\scalebox{0.8}{30}}$ & $x_{\scalebox{0.8}{31}}$ & $x_{\scalebox{0.8}{32}}$ & $x_{\scalebox{0.8}{33}}$ & $x_{\scalebox{0.8}{34}}$ & 
$x_{\scalebox{0.8}{35}}$ & $x_{\scalebox{0.8}{36}}$ & $x_{\scalebox{0.8}{37}}$%
\end{tabular}%
\right) $
}
\end{flushleft}

\bigskip

$\textrm{M}_{C}\textrm{CJ}(X)=\left( 
\begin{array}{c}
^{t}\textrm{MC}(x_{1}) \\ 
^{t}\textrm{MC}(x_{2}) \\ 
^{t}\textrm{MC}(x_{3})%
\end{array}%
\right) =$ $^{t}\textrm{M}_{R}\textrm{JC}(X)\in M(24\times 1,\R),$

$\textrm{M}_{C}\textrm{CCC}(X)=\left( 
\begin{array}{c}
\chi _{1} \\ 
\chi _{2} \\ 
\chi _{3}%
\end{array}%
\right) \in M(3\times 1,\R).$

\bigskip

\section{The adjoint representation of \gR$_{4}^{}$}

\bigskip 

\ \ \ \ \emph{Lemma 7.1. } For $(D_{1},M_{1}),(D_{2},M_{2})\in \ggR%
_{4}^{},M_{1}=A_{1}(m_{11})+A_{2}(m_{12})$

$+A_{3}(m_{13}),M_{2}=A_{1}(m_{21})+A_{2}(m_{22})+A_{3}(m_{23}),
m_{11},m_{12},m_{13},m_{21},m_{22},$

$m_{23}\in \ggC ^{},$

\noindent
Lie bracket $[(D_{1},M_{1}),(D_{2},M_{2})]_{4}$ is expressed by the
following expression.

$[(D_{1},M_{1}),(D_{2},M_{2})]_{4}=$
\begin{flushleft}
$\left( 
%
%
\right) .$
}

\bigskip

\emph{Proof.} By \emph{Definition 2.14} and mappings $%
\textrm{f}_{DDD},\textrm{f}_{JJD},\textrm{f}_{DJJ},\textrm{f}_{iJJ},$we have the above 
expression. \ \ \ \emph{Q.E.D.}

\bigskip

\emph{Theorem 7.2. } The image of the adjoit representation  $ad($\gR$_{4})$ of \textbf{\gR%
}$_{4}^{}$ is expressed by

$ad(D,M)=ad(\textrm{fv}(D,M))=\left( 
\begin{tabular}{@{}|c|c|@{}}
\hline
$\textrm{LD}(D)$ & $-\textrm{M}_{R}\textrm{Dl}(M)$ \\ \hline
$\text{ }^{t}\textrm{M}_{R}\textrm{Dl}(M)$ & $\textrm{M}_{D}\textrm{Jl}(D)+\frac{1}{2}\textrm{MI}(M)$ \\ \hline
\end{tabular}%
\right) \in M(52\times 52,\R).$

\noindent
where $D=\sum\limits_{0\leq i<j\leq 7}d_{ij}D_{ij} \ (d_{ij}\in
\R),M=A_{1}(m_{1})+A_{2}(m_{2})+A_{3}(m_{3}),$

\ \ \ $m_{i}=\sum\limits_{j=0}^{7}$m$%
_{ij}e_{j}$ $($m$_{ij}\in \R).$

\noindent
Since \gR$_{4}^{}$ and \gf$_{4}^{}$ are isomorphic, 
$ad(D,M)$ is also the image of the adjoit representation of \gf$_{4}^{}.$

\bigskip

\emph{Proof.} \ Using \emph{Lemma
7.1, 6.10, 6.11, 6.12, 6.13, 6.14, 6.15, 6.16, 6.17, 6.18, 6.19, 6.20, }and %
\emph{Lemma 6.24,} we have

$ad(\textrm{fv}(D,M))=$
\begin{flushleft}
$\left(
%
%
$ \\ \hline
\end{tabular}%
\right) .$
\end{flushleft}

\noindent
Hence we have the above expression. \ \ \ \ \emph{Q.E.D.}

\bigskip

\section{The adjoint representation of \gR$_{6}^{}$}

\bigskip 

\ \ \ \ \emph{Lemma 8.1.} \ For $\phi _{1}=(D_{1},M_{1},T1),\phi
_{2}=(D_{2},M_{2},T2)\in $\gR$_{6}^{}$\textbf{,}

$%
M_{1}=A_{1}(m_{11})+A_{2}(m_{12})+A_{3}(m_{13}),M_{2}=A_{1}(m_{21})+A_{2}(m_{22})+A_{3}(m_{23}), 
$

$\ m_{11},m_{12},m_{13},m_{21},m_{22},m_{23}\in $\gC $^{},$

$T_{1}=\tau _{11}E_{1}+\tau _{12}E_{2}+(-\tau _{11}-\tau
_{12})E_{3}+F_{1}(t_{11})+F_{2}(t_{12})+F_{3}(t_{13}),$

$T_{2}=\tau _{21}E_{1}+\tau _{22}E_{2}+(-\tau _{21}-\tau
_{22})E_{3}+F_{1}(t_{21})+F_{2}(t_{22})+F_{3}(t_{23}),$

$t_{11},t_{12},t_{13},t_{21},t_{22},t_{23}\in $\gC $^{},\tau _{11},\tau
_{12},\tau _{21},\tau _{22}\in \R,$

\noindent
Lie bracket $[\phi _{1},\phi _{2}]_{6}$ is expressed by the following
expression.

$[\phi _{1},\phi _{2}]_{6}=\left( 
\begin{array}{cccc}
F_{11} & F_{12} & 0 & F_{14} \\ 
F_{21} & F_{22} & F_{23} & F_{24} \\ 
0 & F_{32} & 0 & F_{34} \\ 
F_{41} & F_{42} & F_{43} & F_{44}%
\end{array}%
\right) \left( 
\begin{tabular}{c}
$D_{2}$ \\ \hline
$%
\begin{array}{c}
m_{21} \\ 
m_{22} \\ 
m_{23}%
\end{array}%
$ \\ \hline
$%
\begin{array}{c}
\tau _{21} \\ 
\tau _{22} \\ 
-\tau _{21-}\tau _{22}%
\end{array}%
$ \\ \hline
$%
\begin{array}{c}
t_{21} \\ 
t_{22} \\ 
t_{23}%
\end{array}%
$ \\ 
\end{tabular}%
\right) $

\ \ \ \ \ \ \ \ \ \ $\ \ =\left( 
\begin{tabular}{c}
$\lbrack \phi _{1},\phi _{2}]_{6D}$ \\ \hline
$%
\begin{array}{c}
\lbrack \phi _{1},\phi _{2}]_{6M\_A_{1}} \\ 
\lbrack \phi _{1},\phi _{2}]_{6M\_A_{2}} \\ 
\lbrack \phi _{1},\phi _{2}]_{6M\_A_{3}}%
\end{array}%
$ \\ \hline
$%
\begin{array}{c}
\lbrack \phi _{1},\phi _{2}]_{6T\_E_{1}} \\ 
\lbrack \phi _{1},\phi _{2}]_{6T\_E_{2}} \\ 
\lbrack \phi _{1},\phi _{2}]_{6T\_E_{3}}%
\end{array}%
$ \\ \hline
$%
\begin{array}{c}
\lbrack \phi _{1},\phi _{2}]_{6T\_F_{1}} \\ 
\lbrack \phi _{1},\phi _{2}]_{6T\_F_{2}} \\ 
\lbrack \phi _{1},\phi _{2}]_{6T\_F_{3}}%
\end{array}%
$ \\ 
\end{tabular}%
\right) ,$

\noindent
where

$\left( 
\begin{array}{cc}
F_{11} & F_{12} \\ 
F_{21} & F_{22}%
\end{array}%
\right) =$

{\scriptsize
\begin{flushleft}
$\left( 
\begin{tabular}{@{}c|c@{}}
$\textrm{f}_{DDD}(D_{1},?)$ & $%
\begin{array}{@{}ccc@{}}
-\textrm{f}_{JJD}(m_{11},?) & -\textrm{d}_{g}\nu ^{2}\textrm{g}_{d}\textrm{f}_{JJD}(m_{12},?) & -\textrm{d}_{g}\nu
\textrm{g}_{d}\textrm{f}_{JJD}(m_{13},?)%
\end{array}%
$ \\ \hline
$%
\begin{array}{@{}c@{}}
-\textrm{f}_{DJJ}(?,m_{11}) \\ 
-\textrm{f}_{DJJ}(?,m_{12})\textrm{d}_{g}\nu \textrm{g}_{d} \\ 
-\textrm{f}_{DJJ}(?,m_{13})\textrm{d}_{g}\nu ^{2}\textrm{g}_{d}%
\end{array}%
$ & $%
\begin{array}{@{}ccc@{}}
\textrm{f}_{DJJ}(D_{1},?) & \frac{1}{2}\textrm{f}_{iJJ}(?,m_{13}) & -\frac{1}{2}%
\textrm{f}_{iJJ}(m_{12},?) \\ 
-\frac{1}{2}\textrm{f}_{iJJ}(m_{13},?) & \textrm{f}_{DJJ}(\textrm{d}_{g}\nu \textrm{g}_{d}D_{1},?) & \frac{1}{2}%
\textrm{f}_{iJJ}(?,m_{11}) \\ 
\frac{1}{2}\textrm{f}_{iJJ}(?,m_{12}) & -\frac{1}{2}\textrm{f}_{iJJ}(m_{11},?) & 
\textrm{f}_{DJJ}(\textrm{d}_{g}\nu ^{2}\textrm{g}_{d}D_{1},?)%
\end{array}%
$ \\ 
\end{tabular}%
\right) ,$
\end{flushleft}
}
\begin{flushleft}
$F_{14}=\left( 
\begin{array}{ccc}
\textrm{f}_{JJD}(t_{11},?) & \textrm{d}_{g}\nu ^{2}\textrm{g}_{d}\textrm{f}_{JJD}(t_{12},?) & \textrm{d}_{g}\nu
\textrm{g}_{d}\textrm{f}_{JJD}(t_{13},?)%
\end{array}%
\right) ,$

$F_{23}=\left( 
\begin{array}{ccc}
0 & -\frac{1}{2}\textrm{f}_{CJJ}(?,t_{11}) & \frac{1}{2}\textrm{f}_{CJJ}(?,t_{11}) \\ 
\frac{1}{2}\textrm{f}_{CJJ}(?,t_{12}) & 0 & -\frac{1}{2}\textrm{f}_{CJJ}(?,t_{12}) \\ 
-\frac{1}{2}\textrm{f}_{CJJ}(?,t_{13}) & \frac{1}{2}\textrm{f}_{CJJ}(?,t_{13}) & 0%
\end{array}%
\right) ,$

$F_{24}=\left( 
\begin{array}{ccc}
0 & \frac{1}{2}\textrm{f}_{iJJ}(?,t_{13}) & -\frac{1}{2}\textrm{f}_{iJJ}(t_{12},?) \\ 
-\frac{1}{2}\textrm{f}_{iJJ}(t_{13},?) & 0 & \frac{1}{2}\textrm{f}_{iJJ}(?,t_{11}) \\
\frac{1}{2}\textrm{f}_{iJJ}(?,t_{12}) & -\frac{1}{2}\textrm{f}_{iJJ}(t_{11},?) &0 \\ 
\end{array}%
\right) $

$\ \ \ \ \ \ +\left( 
\begin{array}{@{}c@{}c@{}c@{}}
\frac{1}{2}\textrm{f}_{CJJ}(\tau _{11}+2\tau _{12},?) & 0 & 0 \\ 
0 & \frac{1}{2}\textrm{f}_{CJJ}(-2\tau _{11}-\tau _{12},?) & 0 \\ 
0 & 0 & \frac{1}{2}\textrm{f}_{CJJ}(\tau _{11}-\tau _{12},?) \\ 
\end{array}%
\right) ,$

$F_{32}=\left( 
\begin{array}{ccc}
0 & \textrm{f}_{JJC}(t_{12},?) & -\textrm{f}_{JJC}(t_{13},?) \\ 
-\textrm{f}_{JJC}(t_{11},?) & 0 & \textrm{f}_{JJC}(t_{13},?) \\ 
\textrm{f}_{JJC}(t_{11},?) & -\textrm{f}_{JJC}(t_{12},?) & 0%
\end{array}%
\right) ,$

$F_{34}=\left( 
\begin{array}{ccc}
0 & -\textrm{f}_{JJC}(m_{12},?) & \textrm{f}_{JJC}(m_{13},?) \\ 
\textrm{f}_{JJC}(m_{11},?) & 0 & -\textrm{f}_{JJC}(m_{13},?) \\ 
-\textrm{f}_{JJC}(m_{11},?) & \textrm{f}_{JJC}(m_{12},?) & 0%
\end{array}%
\right) ,$

$F_{41}=\left( 
\begin{array}{c}
-\textrm{f}_{DJJ}(?,t_{11}) \\ 
-\textrm{f}_{DJJ}(?,t_{12})\textrm{d}_{g}\nu \textrm{g}_{d} \\ 
-\textrm{f}_{DJJ}(?,t_{13})\textrm{d}_{g}\nu ^{2}\textrm{g}_{d}%
\end{array}%
\right) ,$

$F_{42}=\left( 
\begin{array}{ccc}
0 & -\frac{1}{2}\textrm{f}_{iJJ}(?,t_{13}) & \frac{1}{2}\textrm{f}_{iJJ}(t_{12},?) \\ 
\frac{1}{2}\textrm{f}_{iJJ}(t_{13},?) & 0 & -\frac{1}{2}\textrm{f}_{iJJ}(?,t_{11}) \\ 
-\frac{1}{2}\textrm{f}_{iJJ}(?,t_{12}) & \frac{1}{2}\textrm{f}_{iJJ}(t_{11},?) & 0 \\ 
\end{array}%
\right) $

$\ \ \ \ \ \ +\left( 
\begin{array}{@{}c@{}c@{}c@{}}
\frac{1}{2}\textrm{f}_{CJJ}(\tau _{11}+2\tau _{12},?) & 0 & 0 \\ 
0 & \frac{1}{2}\textrm{f}_{CJJ}(-2\tau _{11}-\tau _{12},?) & 0 \\ 
0 & 0 & \frac{1}{2}\textrm{f}_{CJJ}(\tau _{11}-\tau _{12},?) \\ 
\end{array}%
\right) ,$

$F_{43}=\left( 
\begin{array}{ccc}
0 & -\frac{1}{2}\textrm{f}_{CJJ}(?,m_{11}) & \frac{1}{2}\textrm{f}_{CJJ}(?,m_{11}) \\ 
\frac{1}{2}\textrm{f}_{CJJ}(?,m_{12}) & 0 & -\frac{1}{2}\textrm{f}_{CJJ}(?,m_{12}) \\ 
-\frac{1}{2}\textrm{f}_{CJJ}(?,m_{13}) & \frac{1}{2}\textrm{f}_{CJJ}(?,m_{13}) & 0%
\end{array}%
\right) $

$F_{44}=\left( 
\begin{array}{ccc}
0 & -\frac{1}{2}\textrm{f}_{iJJ}(?,m_{13}) & \frac{1}{2}\textrm{f}_{iJJ}(m_{12},?) \\ 
\frac{1}{2}\textrm{f}_{iJJ}(m_{13},?) & 0 & -\frac{1}{2}\textrm{f}_{iJJ}(?,m_{11}) \\ 
-\frac{1}{2}\textrm{f}_{iJJ}(?,m_{12}) & \frac{1}{2}\textrm{f}_{iJJ}(m_{11},?) & 0 \\ 
\end{array}%
\right) $

$\ \ \ \ \ \ +\left( 
\begin{array}{ccc}
\textrm{f}_{DJJ}(D_{1},?) & 0 & 0 \\ 
0 & \textrm{f}_{DJJ}(\textrm{d}_{g}\nu \textrm{g}_{d}D_{1},?) & 0 \\ 
0 & 0 & \textrm{f}_{DJJ}(\textrm{d}_{g}\nu ^{2}\textrm{g}_{d}D_{1},?) \\ 
\end{array}%
\right) ,$

\end{flushleft}

\begin{align*}
[\phi _{1},\phi _{2}]_{6D}&=[D1,D2]-\textrm{JD}(m_{11},m_{21})-\textrm{d}_{g}\nu
^{2}\textrm{g}_{d}(\textrm{JD}(m_{12},m_{22}))\\
&-\textrm{d}_{g}\nu \textrm{g}_{d}(\textrm{JD}(m_{13},m_{23}))+\textrm{JD}(t_{11},t_{21})+\textrm{d}_{g}\nu
^{2}\textrm{g}_{d}(\textrm{JD}(t_{12},t_{22}))\\
&+\textrm{d}_{g}\nu \textrm{g}_{d}(\textrm{JD}(t_{13},t_{23})),\\
[\phi _{1},\phi _{2}]_{6M\_A_{1}}&=\textrm{g}_{d}(D_{1})m_{21}-\textrm{g}_{d}(D_{2})m_{11}-\frac{1%
}{2}\overline{m_{12}m_{23}}+\frac{1}{2}\overline{m_{22}m_{13}}\\
&+\frac{1}{2}((\tau _{11}+2\tau
_{12})t_{21}-(\tau _{21}+2\tau _{22})t_{11}-\overline{t_{12}t_{23}}+%
\overline{t_{22}t_{13}}),\\
[\phi _{1},\phi _{2}]_{6M\_A_{2}}&=\nu \textrm{g}_{d}(D_{1})m_{22}-\nu
\textrm{g}_{d}(D_{2})m_{12}-\frac{1}{2}\overline{m_{13}m_{21}}+\frac{1}{2}\overline{%
m_{23}m_{11}}\\
&+\frac{1}{2}((-2\tau _{11}-\tau
_{12})t_{22}-(-2\tau _{21}-\tau _{22})t_{12}+\overline{t_{23}t_{11}}-%
\overline{t_{13}t_{21}}),\\
[\phi _{1},\phi _{2}]_{6M\_A_{3}}&=\nu ^{2}\textrm{g}_{d}(D_{1})m_{23}-\nu
^{2}\textrm{g}_{d}(D_{2})m_{13}-\frac{1}{2}\overline{m_{11}m_{22}}+\frac{1}{2}%
\overline{m_{21}m_{12}}\\
&+\frac{1}{2}((\tau _{11}-\tau
_{12})t_{23}-(\tau _{21}-\tau _{22})t_{13}-\overline{t_{11}t_{22}}+\overline{%
t_{21}t_{12}}),\\
[\phi _{1},\phi
_{2}]_{6T%
\_E_{1}}&=-(m_{12},t_{22})+(m_{13},t_{23})+(m_{22},t_{12})-(m_{23},t_{13}),\\
[\phi _{1},\phi
_{2}]_{6T%
\_E_{2}}&=-(m_{13},t_{23})+(m_{11},t_{21})+(m_{23},t_{13})-(m_{21},t_{11}),\\
[\phi _{1},\phi
_{2}]_{6T%
\_E_{3}}&=-(m_{11},t_{21})+(m_{12},t_{22})+(m_{21},t_{11})-(m_{22},t_{12}),\\
[\phi _{1},\phi _{2}]_{6T\_F_{1}}&=\textrm{g}_{d}(D_{1})t_{21}-\textrm{g}_{d}(D_{2})t_{11}\\
&+\frac{1%
}{2}(\overline{m_{12}t_{23}}-\overline{t_{22}m_{13}}+(-\tau _{21}-2\tau
_{22})m_{11})\\
&-\frac{1}{2}(\overline{%
m_{22}t_{13}}-\overline{t_{12}m_{23}}+(-\tau _{11}-2\tau _{12})m_{21}),\\
[\phi _{1},\phi _{2}]_{6T\_F_{2}}&=\nu \textrm{g}_{d}(D_{1})t_{22}-\nu
\textrm{g}_{d}(D_{2})t_{12}\\
&+\frac{1}{2}(\overline{m_{13}t_{21}}-\overline{t_{23}m_{11}%
}+(2\tau _{21}+\tau _{22})m_{12})\\
& -\frac{1}{2}(\overline{%
m_{23}t_{11}}-\overline{t_{13}m_{21}}+(2\tau _{11}+\tau _{12})m_{22}),\\
[\phi _{1},\phi _{2}]_{6T\_F_{3}}&=\nu ^{2}\textrm{g}_{d}(D_{1})t_{23}-\nu
^{2}\textrm{g}_{d}(D_{2})t_{13}\\
&+\frac{1}{2}(\overline{m_{11}t_{22}}-\overline{%
t_{21}m_{12}}+(-\tau _{21}+\tau _{22})m_{13})\\
&-\frac{1}{2}(\overline{%
m_{21}t_{12}}-\overline{t_{11}m_{22}}+(-\tau _{11}+\tau _{12})m_{23}).
\end{align*}

\emph{Proof.} Using \emph{Lemma 3.11} and definition of mappings $%
\textrm{f}_{DDD},\textrm{f}_{JJD},\textrm{f}_{DJJ},$
$\textrm{f}_{iJJ},\textrm{f}_{JJC},\textrm{f}_{CJJ},$
we have the above expression. \ \ \ \ \emph{Q.E.D.}

\bigskip

\emph{Lemma 8.2.\ }When we put $T\in $\gJ$_{0}^{}$ of \emph{%
Lemma 8.1} as 

$\tau _{1}E_{1}+\tau _{2}E_{2}+F_{1}(t_{1})+F_{2}(t_{2})+F_{3}(t_{3}),$ we
have

$\left( 
\begin{tabular}{cccc}
$F_{11}$ & $F_{12}$ & $0$ & $F_{14}$ \\ \cline{3-3}
$F_{21}$ & $F_{22}$ & \multicolumn{1}{|c}{$F_{23}^{\prime }$} & \multicolumn{1}{|c}{$%
F_{24}$} \\ \cline{2-4}
$0$ & \multicolumn{1}{|c}{$F_{32}^{\prime }$} & \multicolumn{1}{|c}{$0$} & 
\multicolumn{1}{|c|}{$F_{34}^{\prime }$} \\ \cline{2-4}
$F_{41}$ & $F_{42}$ & \multicolumn{1}{|c}{$F_{43}^{\prime }$} & \multicolumn{1}{|c}{$%
F_{44}$} \\ \cline{3-3}
\end{tabular}%
\right) \left( 
\begin{tabular}{c}
$D_{2}$ \\ \hline
$%
\begin{array}{c}
m_{21} \\ 
m_{22} \\ 
m_{23}%
\end{array}%
$ \\ \hline
$%
\begin{tabular}{c}
$\tau _{21}$ \\ 
$\tau _{22}$%
\end{tabular}%
$ \\ \hline
$%
\begin{array}{c}
t_{21} \\ 
t_{22} \\ 
t_{23}%
\end{array}%
$ \\ 
\end{tabular}%
\right) =\left( 
\begin{tabular}{c}
$\lbrack \phi _{1},\phi _{2}]_{6D}$ \\ \hline
$%
\begin{array}{c}
\lbrack \phi _{1},\phi _{2}]_{6M\_A_{1}} \\ 
\lbrack \phi _{1},\phi _{2}]_{6M\_A_{2}} \\ 
\lbrack \phi _{1},\phi _{2}]_{6M\_A_{3}}%
\end{array}%
$ \\ \hline
$%
\begin{tabular}{c}
$\lbrack \phi _{1},\phi _{2}]_{6T\_E_{1}}$ \\ 
$\lbrack \phi _{1},\phi _{2}]_{6T\_E_{2}}$%
\end{tabular}%
$ \\ \hline
$%
\begin{array}{c}
\lbrack \phi _{1},\phi _{2}]_{6T\_F_{1}} \\ 
\lbrack \phi _{1},\phi _{2}]_{6T\_F_{2}} \\ 
\lbrack \phi _{1},\phi _{2}]_{6T\_F_{3}}%
\end{array}%
$ \\ 
\end{tabular}%
\right) ,$

\noindent
where

$F_{23}^{\prime }=\left( 
\begin{array}{cc}
-\frac{1}{2}\textrm{f}_{CJJ}(?,t_{11}) & -\textrm{f}_{CJJ}(?,t_{11}) \\ 
\textrm{f}_{CJJ}(?,t_{12}) & \frac{1}{2}\textrm{f}_{CJJ}(?,t_{12}) \\ 
-\frac{1}{2}\textrm{f}_{CJJ}(?,t_{13}) & \frac{1}{2}\textrm{f}_{CJJ}(?,t_{13})%
\end{array}%
\right) ,$

$F_{32}^{\prime }=\left( 
\begin{array}{ccc}
0 & \textrm{f}_{JJC}(t_{12},?) & -\textrm{f}_{JJC}(t_{13},?) \\ 
-\textrm{f}_{JJC}(t_{11},?) & 0 & \textrm{f}_{JJC}(t_{13},?)%
\end{array}%
\right) ,$

$F_{34}^{\prime }=\left( 
\begin{array}{ccc}
0 & -\textrm{f}_{JJC}(m_{12},?) & \textrm{f}_{JJC}(m_{13},?) \\ 
\textrm{f}_{JJC}(m_{11},?) & 0 & -\textrm{f}_{JJC}(m_{13},?)%
\end{array}%
\right) $

$F_{43}^{\prime }=\left( 
\begin{array}{cc}
-\frac{1}{2}\textrm{f}_{CJJ}(?,m_{11}) & -\textrm{f}_{CJJ}(?,m_{11}) \\ 
\textrm{f}_{CJJ}(?,m_{12}) & \frac{1}{2}\textrm{f}_{CJJ}(?,m_{12}) \\ 
-\frac{1}{2}\textrm{f}_{CJJ}(?,m_{13}) & \frac{1}{2}\textrm{f}_{CJJ}(?,m_{13})%
\end{array}%
\right) .$

\bigskip

\emph{Proof.} By \emph{Lemma 8.1} and direct calculation, we have
the above expression.

\ \ \ \ \emph{Q.E.D.}

\bigskip

\emph{Theorem 8.3.} The image of the adjoit representation $ad($\gR$_{6})$ of \textbf{%
\gR}$_{6}^{}$ is expressed by

$ad(D,M,T)=ad(\textrm{fv}(D,M,T))=$

\begin{flushleft}
$\left( 
\begin{array}{@{}c@{}c|@{}c@{}c@{}}
\textrm{LD}(D) & -\textrm{M}_{R}\textrm{Dl}(M) & 0 & \textrm{M}_{R}\textrm{Dl}(T) \\
\text{ }^{t}\textrm{M}_{R}\textrm{Dl}(M) & \ \textrm{M}_{D}\textrm{Jl}(D)+\frac{1}{2}\textrm{MI}(M)%
&\frac{1}{2}^{\text{t}}\textrm{MJC}^{3}(T) & \ \frac{1}{2}%
\textrm{M}_{D}\textrm{E}(T)+\frac{1}{2}\textrm{MI}(T) \\ \hline
0 & \textrm{MJC}_{3}(T) & 0 & -MJC_{3}(M) \\ 
\text{ }^{t}\textrm{M}_{R}\textrm{Dl}(T) & \ \frac{1}{2}\textrm{M}_{D}\textrm{E}(T)-\frac{1}{2}\textrm{MI}(T) & \frac{%
1}{2}^{\text{t}}\textrm{MJC}^{3}(M) &\ \textrm{M}_{D}\textrm{Jl}(D)-\frac{1}{2}\textrm{MI}(M) %
\end{array}%
\right) $
\end{flushleft}

$\ \in M(78\times 78,\R),$

\noindent
where $D=\sum\limits_{0\leq i<j\leq 7}d_{ij}D_{ij} \ (d_{ij}\in
\R),M=A_{1}(m_{1})+A_{2}(m_{2})+A_{3}(m_{3})$,

\noindent
$m_{i}=\sum\limits_{j=0}^{7}$m$%
_{ij}e_{j}$ $(m_{ij}\in \R),$
$T=\tau _{1}E_{1}+\tau _{2}E_{2}+(-\tau _{1}-\tau
_{2})E_{3}+F_{1}(t_{1})+F_{2}(t_{2})+F_{3}(t_{3}),t_{i}=\sum%
\limits_{j=0}^{7} $t$_{ij}e_{j}(t_{ij}\in \R),\tau _{1},\tau _{2}\in \R.$

\noindent
Since \gR$_{6}^{}$ and \ge$_{6,1}^{}$ are isomorphic, 
$ad(D,M,T)$ is also the image of the adjoit representation of \ge$_{6,1}^{}. $

\bigskip

\emph{Proof.} Using \emph{Lemma 8.1, 8.2, Lemma
6.10, 6.11, 6.12, 6.13, 6.14, 6.15},

\noindent
\emph{6.16, 6.17, 6.18, 6.19, 6.20,} and \emph{Lemma 6.24}, we have

$ad(\textrm{fv}(D,M,T))=$

$\left( 
\begin{tabular}{cccc}
\cline{1-2}
 \multicolumn{1}{|c}{$M_{11}$} & \multicolumn{1}{|c}{$M_{12}$} & \multicolumn{1}{|c}{$0$} & $M_{14}$ \\ 
\cline{1-2}
 \multicolumn{1}{|c}{$M_{21}$} & \multicolumn{1}{|c}{$M_{22}$} & \multicolumn{1}{|c}{$M_{23}$} & $M_{24}$ \\ 
\cline{1-2}
$0$ & $M_{32}$ & $0$ & $M_{34}$ \\ 
$M_{41}$ & $M_{42}$ & $M_{43}$ & $M_{44}$ \\ 
\end{tabular}%
\right) \in M(78\times 78,\R),$

$\left( 
\begin{array}{cc}
M_{11} & M_{12} \\ 
M_{21} & M_{22}%
\end{array}%
\right) =ad(D,M)\in M(52\times 52,\R),$

$M_{14}=\left( 
\begin{array}{ccc}
\textrm{MDl}(t_{1}) & Mv^{2}\textrm{MDl}(t_{2}) & Mv\textrm{MDl}(t_{3})%
\end{array}%
\right) \in M(28\times 24,\R),$

$M_{23}=\left( 
\begin{array}{cc}
-\frac{1}{2}^{t}\textrm{MC}(t_{1}) & -^{t}\textrm{MC}(t_{1}) \\ 
^{t}\textrm{MC}(t_{2}) & \frac{1}{2}^{t}\textrm{MC}(t_{2}) \\ 
-\frac{1}{2}^{t}\textrm{MC}(t_{3}) & \frac{1}{2}^{t}\textrm{MC}(t_{3})%
\end{array}%
\right) \in M(24\times 2,\R),$

$M_{24}=\left( 
\begin{tabular}{ccc}
\cline{1-1}
\multicolumn{1}{|c}{$\frac{1}{2}\textrm{ME}(\tau _{1}+2\tau _{2})$} & 
\multicolumn{1}{|c}{$\frac{1}{2}\textrm{MIr}(t_{3})$} & $-\frac{1}{2}\textrm{MIl}(t_{2})$ \\ 
\cline{1-2}
$-\frac{1}{2}\textrm{MIl}(t_{3})$ & \multicolumn{1}{|c}{$\frac{1}{2}\textrm{ME}(-2\tau
_{1}-\tau _{2})$} & \multicolumn{1}{|c}{$\frac{1}{2}\textrm{MIr}(t_{1})$} \\ 
\cline{2-3}
$\frac{1}{2}\textrm{MIr}(t_{2})$ & $-\frac{1}{2}\textrm{MIl}(t_{1})$ & \multicolumn{1}{|c|}{$%
\frac{1}{2}\textrm{ME}(\tau _{1}-\tau _{2})$} \\ \cline{3-3}
\end{tabular}%
\right)$ 

$\ \ \ \ \ \ \ \ \in M(24\times 24,\R),$

$M_{32}=\left( 
\begin{array}{ccc}
0 & \textrm{MC}(t_{2}) & -\textrm{MC}(t_{3}) \\ 
-\textrm{MC}(t_{1}) & 0 & \textrm{MC}(t_{3})%
\end{array}%
\right) \in M(2\times 24,\R),$

$M_{34}=\left( 
\begin{array}{ccc}
0 & -\textrm{MC}(m_{2}) & \textrm{MC}(m_{3}) \\ 
\textrm{MC}(m_{1}) & 0 & -\textrm{MC}(m_{3})%
\end{array}%
\right) \in M(2\times 24,\R),$

$M_{41}=\left( 
\begin{array}{c}
^{t}\textrm{MDl}(t_{1}) \\ 
^{t}\textrm{MDl}(t_{2})Mv \\ 
^{t}\textrm{MDl}(t_{3})Mv^{2}%
\end{array}%
\right) \in M(24\times 28,\R),$

$M_{42}=\left( 
\begin{tabular}{ccc}
\cline{1-1}
\multicolumn{1}{|c}{$\frac{1}{2}\textrm{ME}(\tau _{1}+2\tau _{2})$} & 
\multicolumn{1}{|c}{$-\frac{1}{2}\textrm{MIr}(t_{3})$} & $\frac{1}{2}\textrm{MIl}(t_{2})$ \\ 
\cline{1-2}
$\frac{1}{2}\textrm{MIl}(t_{3})$ & \multicolumn{1}{|c}{$\frac{1}{2}\textrm{ME}(-2\tau
_{1}-\tau _{2})$} & \multicolumn{1}{|c}{$-\frac{1}{2}\textrm{MIr}(t_{1})$} \\ 
\cline{2-3}
$-\frac{1}{2}\textrm{MIr}(t_{2})$ & $\frac{1}{2}\textrm{MIl}(t_{1})$ & \multicolumn{1}{|c|}{$%
\frac{1}{2}\textrm{ME}(\tau _{1}-\tau _{2})$} \\ \cline{3-3}
\end{tabular}%
\right)$

$\ \ \ \ \ \ \ \ \in M(24\times 24,\R),$

$M_{43}=\left( 
\begin{array}{cc}
-\frac{1}{2}^{t}\textrm{MC}(m_{1}) & -^{t}\textrm{MC}(m_{1}) \\ 
^{t}\textrm{MC}(m_{2}) & \frac{1}{2}^{t}\textrm{MC}(m_{2}) \\ 
-\frac{1}{2}^{t}\textrm{MC}(m_{3}) & \frac{1}{2}^{t}\textrm{MC}(m_{3})%
\end{array}%
\right)\in M(24\times 2,\R),$

$M_{44}=\left( 
\begin{tabular}{ccc}
\cline{1-1}
\multicolumn{1}{|c}{$\textrm{MJl}(D)$} & \multicolumn{1}{|c}{$-\frac{1}{2}\textrm{MIr}(m_{3})$}
& $\frac{1}{2}\textrm{MIl}(m_{2})$ \\ \cline{1-2}
$\frac{1}{2}\textrm{MIl}(m_{3})$ & \multicolumn{1}{|c}{$\textrm{MJl}(MvD)$} & 
\multicolumn{1}{|c}{$-\frac{1}{2}\textrm{MIr}(m_{1})$} \\ \cline{2-3}
$-\frac{1}{2}\textrm{MIr}(m_{2})$ & $\frac{1}{2}\textrm{MIl}(m_{1})$ & \multicolumn{1}{|c|}{$%
\textrm{MJl}(Mv^{2}D)$} \\ \cline{3-3}
\end{tabular}%
\right)$

$\ \ \ \ \ \ \ \ \in M(24\times 24,\R),$

\noindent
Hence we have the above expression. \ \ \ \ \emph{Q.E.D.}

\bigskip

\section{The adjoint representation of \gR$_{7}^{}$}

\bigskip 

\ \ \ \ \emph{Lemma 9.1.} \ For $\Phi
_{1}=(D_{1},M_{1},T_{1},A_{1},B_{1},\rho _{1})$,

$\Phi_{2}=(D_{2},M_{2},T_{2},A_{2},B_{2},\rho _{2})\in $\gR$_{7}^{},$

$%
M_{1}=A_{1}(m_{11})+A_{2}(m_{12})+A_{3}(m_{13}),M_{2}=A_{1}(m_{21})+A_{2}(m_{22})+A_{3}(m_{23}), 
$

$\ m_{11},m_{12},m_{13},m_{21},m_{22},m_{23}\in $\gC $^{},$

$T_{1}=\tau _{11}E_{1}+\tau _{12}E_{2}+(-\tau _{11}-\tau
_{12})E_{3}+F_{1}(t_{11})+F_{2}(t_{12})+F_{3}(t_{13}),$

$T_{2}=\tau _{21}E_{1}+\tau _{22}E_{2}+(-\tau _{21}-\tau
_{22})E_{3}+F_{1}(t_{21})+F_{2}(t_{22})+F_{3}(t_{23}),$

$t_{11},t_{12},t_{13},t_{21},t_{22},t_{23}\in $\gC $^{},\tau _{11},\tau
_{12},\tau _{21},\tau _{22}\in \R,$

$A_{1}=\alpha _{11}E_{1}+\alpha _{12}E_{2}+\alpha
_{13}E_{3}+F_{1}(a_{11})+F_{2}(a_{12})+F_{3}(a_{13}),$

$A_{2}=\alpha _{21}E_{1}+\alpha _{22}E_{2}+\alpha
_{23}E_{3}+F_{1}(a_{21})+F_{2}(a_{22})+F_{3}(a_{23}),$

$B_{1}=\beta _{11}E_{1}+\beta _{12}E_{2}+\beta
_{13}E_{3}+F_{1}(b_{11})+F_{2}(b_{12})+F_{3}(b_{13}),$

$B_{2}=\beta _{21}E_{1}+\beta _{22}E_{2}+\beta
_{23}E_{3}+F_{1}(b_{21})+F_{2}(b_{22})+F_{3}(b_{23}),$

$a_{11},a_{12},a_{13},a_{21},a_{22},a_{23}\in $\gC $^{},\alpha
_{11},\alpha _{12},\alpha _{13},\alpha _{21},\alpha _{22},\alpha _{23}\in \R,$

$b_{11},b_{12},b_{13},b_{21},b_{22},b_{23}\in $\gC $^{},\beta
_{11},\beta _{12},\beta _{13},\beta _{21},\beta _{22},\beta _{23}\in \R,$

$\rho _{1},\rho _{2}\in \R,$

\noindent
Lie bracket $[\Phi _{1},\Phi _{2}]_{7}$ is expressed by the following
expression.

{\fontsize{8pt}{10pt} \selectfont%
\begin{tabular}
[c]{@{}c@{}c@{}}
& $D$ \  $m$ \  $\tau$ \ \ $t$ \  $\alpha$ \ \ $a$
\ \  $\beta$ \ \ \ $b$\ \ \  $\rho$\\%
\begin{tabular}
[c]{@{}c@{}}%
$D$\\
$m$\\
$\tau$\\
$t$\\
$\alpha$\\
$a$\\
$\beta$\\
$b$\\
$\rho$%
\end{tabular}
& $\left(
\begin{tabular}
[c]{@{}c@{}c@{}c@{}c@{}c@{}c@{}c@{}c@{}c@{}}%
$F_{11}$ & $F_{12}$ & $0$ & $F_{14}$ & $0$ & $F_{16}$ & $0$ & $F_{18}$ & $0$\\
$F_{21}$ & $F_{22}$ & $F_{23}$ & $F_{24}$ & $F_{25}$ & $F_{26}$ & $F_{27}$ &
$F_{28}$ & $0$\\
$0$ & $F_{32}$ & $0$ & $F_{34}$ & $F_{35}$ & $F_{36}$ & $F_{37}$ & $F_{38}$ &
$0$\\
$F_{41}$ & $F_{42}$ & $F_{43}$ & $F_{44}$ & $F_{45}$ & $F_{46}$ & $F_{47}$ &
$F_{48}$ & $0$\\
$0$ & $F_{52}$ & $F_{53}$ & $F_{54}$ & $F_{55}$ & $F_{56}$ & $0$ & $0$ &
$F_{59}$\\
$F_{61}$ & $F_{62}$ & $F_{63}$ & $F_{64}$ & $F_{65}$ & $F_{66}$ & $0$ & $0$ &
$F_{69}$\\
$0$ & $F_{72}$ & $F_{73}$ & $F_{74}$ & $0$ & $0$ & $F_{77}$ & $F_{78}$ &
$F_{79}$\\
$F_{81}$ & $F_{82}$ & $F_{83}$ & $F_{84}$ & $0$ & $0$ & $F_{87}$ & $F_{88}$ &
$F_{89}$\\
$0$ & $0$ & $0$ & $0$ & $F_{95}$ & $F_{96}$ & $F_{97}$ & $F_{98}$ & $0$%
\end{tabular}
\right)  $%
\end{tabular}
$\left(
\begin{tabular}
[c]{@{}c@{}}%
$D_{2}$\\
$m_{21}$\\
$m_{22}$\\
$m_{23}$\\
$\tau_{21}$\\
$\tau_{22}$\\
$-\tau_{21}-\tau_{22}$\\
$t_{21}$\\
$t_{22}$\\
$t_{23}$\\
$\alpha_{21}$\\
$\alpha_{22}$\\
$\alpha_{23}$\\
$a_{21}$\\
$a_{22}$\\
$a_{23}$\\
$\beta_{21}$\\
$\beta_{22}$\\
$\beta_{23}$\\
$b_{21}$\\
$b_{22}$\\
$b_{23}$\\
$\rho_{2}$%
\end{tabular}
\right)  =\left(
\begin{tabular}
[c]{@{}c@{}}%
$\lbrack\Phi_{1},\Phi_{2}]_{7D}$\\
$\lbrack\Phi_{1},\Phi_{2}]_{7M\_A_{1}}$\\
$\lbrack\Phi_{1},\Phi_{2}]_{7M\_A_{2}}$\\
$\lbrack\Phi_{1},\Phi_{2}]_{7M\_A_{3}}$\\
$\lbrack\Phi_{1},\Phi_{2}]_{7T\_E_{1}}$\\
$\lbrack\Phi_{1},\Phi_{2}]_{7T\_E_{2}}$\\
$\lbrack\Phi_{1},\Phi_{2}]_{7T\_E_{3}}$\\
$\lbrack\Phi_{1},\Phi_{2}]_{7T\_F_{1}}$\\
$\lbrack\Phi_{1},\Phi_{2}]_{7T\_F_{2}}$\\
$\lbrack\Phi_{1},\Phi_{2}]_{7T\_F_{3}}$\\
$\lbrack\Phi_{1},\Phi_{2}]_{7A\_E_{1}}$\\
$\lbrack\Phi_{1},\Phi_{2}]_{7A\_E_{2}}$\\
$\lbrack\Phi_{1},\Phi_{2}]_{7A\_E_{3}}$\\
$\lbrack\Phi_{1},\Phi_{2}]_{7A\_F_{1}}$\\
$\lbrack\Phi_{1},\Phi_{2}]_{7A\_F_{2}}$\\
$\lbrack\Phi_{1},\Phi_{2}]_{7A\_F_{3}}$\\
$\lbrack\Phi_{1},\Phi_{2}]_{7B\_E_{1}}$\\
$\lbrack\Phi_{1},\Phi_{2}]_{7B\_E_{2}}$\\
$\lbrack\Phi_{1},\Phi_{2}]_{7B\_E_{3}}$\\
$\lbrack\Phi_{1},\Phi_{2}]_{7B\_F_{1}}$\\
$\lbrack\Phi_{1},\Phi_{2}]_{7B\_F_{2}}$\\
$\lbrack\Phi_{1},\Phi_{2}]_{7B\_F_{3}}$\\
$\lbrack\Phi_{1},\Phi_{2}]_{7\rho}$%
\end{tabular}
\right)  $
}

\noindent
where

$F_{16}=\left( 
\begin{array}{ccc}
-2\textrm{f}_{JJD}(?,b_{11}) & -2\textrm{d}_{g}\nu ^{2}\textrm{g}_{d}\textrm{f}_{JJD}(?,b_{12}) & -2\textrm{d}_{g}\nu
\textrm{g}_{d}\textrm{f}_{JJD}(?,b_{13})%
\end{array}%
\right) ,$

$F_{18}=\left( 
\begin{array}{ccc}
2\textrm{f}_{JJD}(a_{11},?) & 2\textrm{d}_{g}\nu ^{2}\textrm{g}_{d}\textrm{f}_{JJD}(a_{12},?) & 2\textrm{d}_{g}\nu
\textrm{g}_{d}\textrm{f}_{JJD}(a_{13},?)%
\end{array}%
\right) ,$

$F_{25}=\left( 
\begin{array}{ccc}
0 & -\textrm{f}_{CJJ}(?,b_{11}) & \textrm{f}_{CJJ}(?,b_{11}) \\ 
\textrm{f}_{CJJ}(?,b_{12}) & 0 & -\textrm{f}_{CJJ}(?,b_{12}) \\ 
-\textrm{f}_{CJJ}(?,b_{13}) & \textrm{f}_{CJJ}(?,b_{13}) & 0%
\end{array}%
\right) ,$

$F_{26}=\left( 
\begin{array}{ccc}
\textrm{f}_{CJJ}(\beta _{12}-\beta _{13},?) &
0 & 0 \\ 
0 & \textrm{f}_{CJJ}(\beta _{13}-\beta
_{11},?) & 0 \\ 
0 & 0 & %
\textrm{f}_{CJJ}(\beta _{11}-\beta _{12},?) 
\end{array}%
\right) ,$

$\ \ \ \ \ \ +\left( 
\begin{array}{ccc}
0 &
\textrm{f}_{iJJ}(?,b_{13}) & -\textrm{f}_{iJJ}(b_{12},?) \\ 
-\textrm{f}_{iJJ}(b_{13},?) & 0 & \textrm{f}_{iJJ}(?,b_{11}) \\ 
\textrm{f}_{iJJ}(?,b_{12}) & -\textrm{f}_{iJJ}(b_{11},?) & %
0 
\end{array}%
\right) ,$

$F_{27}=\left( 
\begin{array}{ccc}
0 & -\textrm{f}_{CJJ}(?,a_{11}) & \textrm{f}_{CJJ}(?,a_{11}) \\ 
\textrm{f}_{CJJ}(?,a_{12}) & 0 & -\textrm{f}_{CJJ}(?,a_{12}) \\ 
-\textrm{f}_{CJJ}(?,a_{13}) & \textrm{f}_{CJJ}(?,a_{13}) & 0%
\end{array}%
\right) ,$

$F_{28}=\left( 
\begin{array}{@{}ccc@{}}
\textrm{f}_{CJJ}(\alpha _{12}-\alpha _{13},?) & 
\textrm{f}_{iJJ}(?,a_{13}) & -\textrm{f}_{iJJ}(a_{12},?) \\ 
-\textrm{f}_{iJJ}(a_{13},?) & \textrm{f}_{CJJ}(\alpha _{13}-\alpha
_{11},?) & \textrm{f}_{iJJ}(?,a_{11}) \\ 
\textrm{f}_{iJJ}(?,a_{12}) &-\textrm{f}_{iJJ}(a_{11},?) & 
\textrm{f}_{CJJ}(\alpha _{11}-\alpha _{12},?) 
\end{array}%
\right) ,$

$F_{35}=\left( 
\begin{array}{ccc}
-\frac{4}{3}\textrm{f}_{CCC}(\beta _{11},?) & \frac{2}{3}\textrm{f}_{CCC}(\beta _{12},?) & 
\frac{2}{3}\textrm{f}_{CCC}(\beta _{13},?) \\ 
\frac{2}{3}\textrm{f}_{CCC}(\beta _{11},?) & -\frac{4}{3}\textrm{f}_{CCC}(\beta _{12},?) & 
\frac{2}{3}\textrm{f}_{CCC}(\beta _{13},?) \\ 
\frac{2}{3}\textrm{f}_{CCC}(\beta _{11},?) & \frac{2}{3}\textrm{f}_{CCC}(\beta _{12},?) & -%
\frac{4}{3}\textrm{f}_{CCC}(\beta _{13},?)%
\end{array}%
\right) ,$

$F_{36}=\left( 
\begin{array}{ccc}
\frac{4}{3}\textrm{f}_{JJC}(b_{11},?) & -\frac{2}{3}\textrm{f}_{JJC}(b_{12},?) & -\frac{2}{3}%
\textrm{f}_{JJC}(b_{13},?) \\ 
-\frac{2}{3}\textrm{f}_{JJC}(b_{11},?) & \frac{4}{3}\textrm{f}_{JJC}(b_{12},?) & -\frac{2}{3}%
\textrm{f}_{JJC}(b_{13},?) \\ 
-\frac{2}{3}\textrm{f}_{JJC}(b_{11},?) & -\frac{2}{3}\textrm{f}_{JJC}(b_{12},?) & \frac{4}{3}%
\textrm{f}_{JJC}(b_{13},?)%
\end{array}%
\right) ,$

$F_{37}=\left( 
\begin{array}{ccc}
\frac{4}{3}\textrm{f}_{CCC}(\alpha _{11},?) & -\frac{2}{3}\textrm{f}_{CCC}(\alpha _{12},?) & -%
\frac{2}{3}\textrm{f}_{CCC}(\alpha _{13},?) \\ 
-\frac{2}{3}\textrm{f}_{CCC}(\alpha _{11},?) & \frac{4}{3}\textrm{f}_{CCC}(\alpha _{12},?) & -%
\frac{2}{3}\textrm{f}_{CCC}(\alpha _{13},?) \\ 
-\frac{2}{3}\textrm{f}_{CCC}(\alpha _{11},?) & -\frac{2}{3}\textrm{f}_{CCC}(\alpha _{12},?) & 
\frac{4}{3}\textrm{f}_{CCC}(\alpha _{13},?)%
\end{array}%
\right) ,$

$F_{38}=\left( 
\begin{array}{ccc}
-\frac{4}{3}\textrm{f}_{JJC}(a_{11},?) & \frac{2}{3}\textrm{f}_{JJC}(a_{12},?) & \frac{2}{3}%
\textrm{f}_{JJC}(a_{13},?) \\ 
\frac{2}{3}\textrm{f}_{JJC}(a_{11},?) & -\frac{4}{3}\textrm{f}_{JJC}(a_{12},?) & \frac{2}{3}%
\textrm{f}_{JJC}(a_{13},?) \\ 
\frac{2}{3}\textrm{f}_{JJC}(a_{11},?) & \frac{2}{3}\textrm{f}_{JJC}(a_{12},?) & -\frac{4}{3}%
\textrm{f}_{JJC}(a_{13},?)%
\end{array}%
\right) ,$

$F_{45}=\left( 
\begin{array}{ccc}
0 & -\textrm{f}_{CJJ}(?,b_{11}) & -\textrm{f}_{CJJ}(?,b_{11}) \\ 
-\textrm{f}_{CJJ}(?,b_{12}) & 0 & -\textrm{f}_{CJJ}(?,b_{12}) \\ 
-\textrm{f}_{CJJ}(?,b_{13}) & -\textrm{f}_{CJJ}(?,b_{13}) & 0%
\end{array}%
\right) ,$

$F_{46}=\left( 
\begin{array}{@{}c@{}c@{}c@{}}
-\textrm{f}_{CJJ}(\beta _{12}+\beta _{13},?) & 
0 & 0 \\ 
0 & -\textrm{f}_{CJJ}(\beta _{13}+\beta
_{11},?) &0 \\ 
0 & 0 & %
-\textrm{f}_{CJJ}(\beta _{11}+\beta _{12},?)
\end{array}%
\right) ,$

$\ \ \ \ \ \ +\left( 
\begin{array}{ccc}
0 & 
-\textrm{f}_{iJJ}(?,b_{13}) & -\textrm{f}_{iJJ}(b_{12},?) \\ 
-\textrm{f}_{iJJ}(b_{13},?) & 0 & -\textrm{f}_{iJJ}(?,b_{11}) \\ 
-\textrm{f}_{iJJ}(?,b_{12}) & -\textrm{f}_{iJJ}(b_{11},?) & 0
\end{array}%
\right) ,$

$F_{47}=\left( 
\begin{array}{ccc}
0 & \textrm{f}_{CJJ}(?,a_{11}) & \textrm{f}_{CJJ}(?,a_{11}) \\ 
\textrm{f}_{CJJ}(?,a_{12}) & 0 & \textrm{f}_{CJJ}(?,a_{12}) \\ 
\textrm{f}_{CJJ}(?,a_{13}) & \textrm{f}_{CJJ}(?,a_{13}) & 0%
\end{array}%
\right) ,$

$F_{48}=\left( 
\begin{array}{@{}ccc@{}}
\textrm{f}_{CJJ}(\alpha _{12}+\alpha _{13},?) & 
0 & 0 \\ 
0 & \textrm{f}_{CJJ}(\alpha _{13}+\alpha
_{11},?) & 0 \\ 
0 & 0 & %
\textrm{f}_{CJJ}(\alpha _{11}+\alpha _{12},?) 
\end{array}%
\right) ,$

$\ \ \ \ \ \ +\left( 
\begin{array}{@{}ccc@{}}
0 & 
\textrm{f}_{iJJ}(?,a_{13}) & \textrm{f}_{iJJ}(a_{12},?) \\ 
\textrm{f}_{iJJ}(a_{13},?) & 0 & \textrm{f}_{iJJ}(?,a_{11}) \\ 
\textrm{f}_{iJJ}(?,a_{12}) & \textrm{f}_{iJJ}(a_{11},?) & 0 
\end{array}%
\right) ,$

$F_{52}=\left( 
\begin{array}{ccc}
0 & \textrm{f}_{JJC}(a_{12},?) & -\textrm{f}_{JJC}(a_{13},?) \\ 
-\textrm{f}_{JJC}(a_{11},?) & 0 & \textrm{f}_{JJC}(a_{13},?) \\ 
\textrm{f}_{JJC}(a_{11},?) & -\textrm{f}_{JJC}(a_{12},?) & 0%
\end{array}%
\right) ,$

$F_{53}=\left( 
\begin{array}{ccc}
-fccc(\alpha _{11},?) & 0 & 0 \\ 
0 & -fccc(\alpha _{12},?) & 0 \\ 
0 & 0 & -fccc(\alpha _{13},?)%
\end{array}%
\right) ,$

$F_{54}=\left( 
\begin{array}{ccc}
0 & -\textrm{f}_{JJC}(a_{12},?) & -\textrm{f}_{JJC}(a_{13},?) \\ 
-\textrm{f}_{JJC}(a_{11},?) & 0 & -\textrm{f}_{JJC}(a_{13},?) \\ 
-\textrm{f}_{JJC}(a_{11},?) & -\textrm{f}_{JJC}(a_{12},?) & 0%
\end{array}%
\right) ,$

$F_{55}=\left( 
\begin{array}{ccc}
\frac{2}{3}\textrm{f}_{CCC}((\rho _{1},?) & 0 & 0 \\ 
0 & \frac{2}{3}\textrm{f}_{CCC}((\rho _{1},?) & 0 \\ 
0 & 0 & 
\frac{2}{3}\textrm{f}_{CCC}((\rho _{1},?)%
\end{array}%
\right) $

$\ \ \ \ \ \ +\left( 
\begin{array}{ccc}
fccc(\tau _{11},?) & 0 & 0 \\ 
0 & fccc(\tau _{12},?) & 0 \\ 
0 & 0 & 
fccc(-\tau _{11}-\tau _{12},?)%
\end{array}%
\right)$ ,

$F_{56}=\left( 
\begin{array}{ccc}
0 & -\textrm{f}_{JJC}(m_{12},?) & 
\textrm{f}_{JJC}(m_{13},?) \\ 
\textrm{f}_{JJC}(m_{11},?) & 0 & 
-\textrm{f}_{JJC}(m_{13},?) \\ 
-\textrm{f}_{JJC}(m_{11},?) & \textrm{f}_{JJC}(m_{12},?) & 
0%
\end{array}%
\right) ,$

$\ \ \ \ \ \ +\left( 
\begin{array}{ccc}
0 & \textrm{f}_{JJC}(t_{12},?) & 
\textrm{f}_{JJC}(t_{13},?) \\ 
\textrm{f}_{JJC}(t_{11},?) & 0 & 
\textrm{f}_{JJC}(t_{13},?) \\ 
\textrm{f}_{JJC}(t_{11},?) & \textrm{f}_{JJC}(t_{12},?) & 
0%
\end{array}%
\right) ,$

$F_{59}=\left( 
\begin{array}{c}
-\frac{2}{3}\textrm{f}_{CCC}(\alpha _{11},?) \\ 
-\frac{2}{3}\textrm{f}_{CCC}(\alpha _{12},?) \\ 
-\frac{2}{3}\textrm{f}_{CCC}(\alpha _{13},?)%
\end{array}%
\right) ,$

$F_{61}=\left( 
\begin{array}{c}
-\textrm{f}_{DJJ}(?,a_{11}) \\ 
-\textrm{f}_{DJJ}(?,a_{12})\textrm{d}_{g}\nu \textrm{g}_{d} \\ 
-\textrm{f}_{DJJ}(?,a_{13})\textrm{d}_{g}\nu ^{2}\textrm{g}_{d}%
\end{array}%
\right) ,$

$F_{62}=\left( 
\begin{array}{@{}c@{}c@{}c@{}}
\frac{1}{2}\textrm{f}_{CJJ}(\alpha _{12}-\alpha _{13},?) & 0 & 0\\ 
0 & \frac{1}{2}\textrm{f}_{CJJ}(\alpha _{13}-\alpha _{11},?) & 0 \\
0 & 0 & \frac{1}{2}\textrm{f}_{CJJ}(\alpha _{11}-\alpha _{12},?) \\ 
\end{array}%
\right) $

$\ \ \ \ \ \ +\left( 
\begin{array}{@{}c@{}c@{}c@{}}
0 & -\frac{1}{2}\textrm{f}_{iJJ}(?,a_{13}) & \frac{1}{2}%
\textrm{f}_{iJJ}(a_{12},?)\\ 
\frac{1}{2}\textrm{f}_{iJJ}(a_{13},?) & 0 & -\frac{1}{2}%
\textrm{f}_{iJJ}(?,a_{11}) \\
-\frac{1}{2}\textrm{f}_{iJJ}(?,a_{12}) & \frac{1}{2}\textrm{f}_{iJJ}(a_{11},?) & 0 \\ 
\end{array}%
\right) ,$

$F_{63}=\left( 
\begin{array}{ccc}
0 & -\frac{1}{2}\textrm{f}_{CJJ}(?,a_{11}) & -\frac{1}{2}\textrm{f}_{CJJ}(?,a_{11}) \\ 
-\frac{1}{2}\textrm{f}_{CJJ}(?,a_{12}) & 0 & -\frac{1}{2}\textrm{f}_{CJJ}(?,a_{12}) \\ 
-\frac{1}{2}\textrm{f}_{CJJ}(?,a_{13}) & -\frac{1}{2}\textrm{f}_{CJJ}(?,a_{13}) & 0%
\end{array}%
\right) ,$

$F_{64}=\left( 
\begin{array}{ccc}
0 & -\frac{1}{2}\textrm{f}_{iJJ}(?,a_{13}) & -\frac{1}{2}\textrm{f}_{iJJ}(a_{12},?) \\ 
-\frac{1}{2}\textrm{f}_{iJJ}(a_{13},?) & 0 & -\frac{1}{2}\textrm{f}_{iJJ}(?,a_{11}) \\ 
-\frac{1}{2}\textrm{f}_{iJJ}(?,a_{12}) & -\frac{1}{2}\textrm{f}_{iJJ}(a_{11},?) & 0 \\ 
\end{array}%
\right) $

$\ \ +\left( 
\begin{array}{@{}c@{}c@{}c@{}}
-\frac{1}{2}\textrm{f}_{CJJ}(\alpha _{12}+\alpha _{13},?) & 0 & 0 \\ 
0 & -\frac{1}{2}\textrm{f}_{CJJ}(\alpha _{13}+\alpha _{11},?) & 0 \\ 
0 & 0 & -\frac{1}{2}\textrm{f}_{CJJ}(\alpha _{11}+\alpha _{12},?) \\ 
\end{array}%
\right) $

$F_{65}=\left( 
\begin{array}{ccc}
0 & -\frac{1}{2}\textrm{f}_{CJJ}(?,m_{11}) & \frac{1}{2}\textrm{f}_{CJJ}(?,m_{11}) \\ 
\frac{1}{2}\textrm{f}_{CJJ}(?,m_{12}) & 0 & -\frac{1}{2}\textrm{f}_{CJJ}(?,m_{12}) \\ 
-\frac{1}{2}\textrm{f}_{CJJ}(?,m_{13}) & \frac{1}{2}\textrm{f}_{CJJ}(?,m_{13}) & 0%
\end{array}%
\right) $

$\ \ \ \ \ \ +\left( 
\begin{array}{ccc}
0 & \frac{1}{2}\textrm{f}_{CJJ}(?,t_{11}) & \frac{1}{2}\textrm{f}_{CJJ}(?,t_{11}) \\ 
\frac{1}{2}\textrm{f}_{CJJ}(?,t_{12}) & 0 & \frac{1}{2}\textrm{f}_{CJJ}(?,t_{12}) \\ 
\frac{1}{2}\textrm{f}_{CJJ}(?,t_{13}) & \frac{1}{2}\textrm{f}_{CJJ}(?,t_{13}) & 0%
\end{array}%
\right) ,$

$F_{66}=\left( 
\begin{tabular}{ccc}
\cline{1-1}
\multicolumn{1}{|c}{$\textrm{f}_{DJJ}(D_{1},?)$} & \multicolumn{1}{|c}{$-\frac{1}{2}%
\textrm{f}_{iJJ}(?,m_{13})$} & $\frac{1}{2}\textrm{f}_{iJJ}(m_{12},?)$ \\ \cline{1-2}
$\frac{1}{2}\textrm{f}_{iJJ}(m_{13},?)$ & \multicolumn{1}{|c}{$\textrm{f}_{DJJ}(\textrm{d}_{g}\nu
\textrm{g}_{d}fD_{1},?)$} & \multicolumn{1}{|c}{$-\frac{1}{2}\textrm{f}_{iJJ}(?,m_{11})$} \\ 
\cline{2-3}
$-\frac{1}{2}\textrm{f}_{iJJ}(?,m_{12})$ & $\frac{1}{2}\textrm{f}_{iJJ}(m_{11},?)$ & 
\multicolumn{1}{|c|}{$\textrm{f}_{DJJ}(\textrm{d}_{g}\nu ^{2}\textrm{g}_{d}D_{1},?)$} \\ \cline{3-3}
\end{tabular}%
\right) $

$\ \ \ \ \ \ +\left( 
\begin{array}{ccc}
\frac{1}{2}\textrm{f}_{CJJ}(-\tau _{11},?) & \frac{1}{2}\textrm{f}_{iJJ}(?,t_{13}) & \frac{1}{2%
}\textrm{f}_{iJJ}(t_{12},?) \\ 
\frac{1}{2}\textrm{f}_{iJJ}(t_{13},?) & \frac{1}{2}\textrm{f}_{CJJ}(-\tau _{12},?) & \frac{1}{2%
}\textrm{f}_{iJJ}(?,t_{11}) \\ 
\frac{1}{2}\textrm{f}_{iJJ}(?,t_{12}) & \frac{1}{2}\textrm{f}_{iJJ}(t_{11},?) & \frac{1}{2}%
\textrm{f}_{CJJ}(\tau _{11}+\tau _{12},?)%
\end{array}%
\right) $

\ \ \ \ \ \ \ \ $\ +\left( 
\begin{array}{ccc}
\frac{2}{3}\textrm{f}_{CJJ}(\rho _{1},?) & 0 & 0 \\ 
0 & \frac{2}{3}\textrm{f}_{CJJ}(\rho _{1},?) & 0 \\ 
0 & 0 & \frac{2}{3}\textrm{f}_{CJJ}(\rho _{1},?)%
\end{array}%
\right) ,$

$F_{69}=\left( 
\begin{array}{c}
-\frac{2}{3}\textrm{f}_{CJJ}(?,a_{11}) \\ 
-\frac{2}{3}\textrm{f}_{CJJ}(?,a_{12}) \\ 
-\frac{2}{3}\textrm{f}_{CJJ}(?,a_{13})%
\end{array}%
\right) ,$

$F_{72}=\left( 
\begin{array}{ccc}
0 & \textrm{f}_{JJC}(b_{12},?) & -\textrm{f}_{JJC}(b_{13},?) \\ 
-\textrm{f}_{JJC}(b_{11},?) & 0 & \textrm{f}_{JJC}(b_{13},?) \\ 
\textrm{f}_{JJC}(b_{11},?) & -\textrm{f}_{JJC}(b_{12},?) & 0%
\end{array}%
\right) ,$

$F_{73}=\left( 
\begin{array}{ccc}
\textrm{f}_{CCC}(\beta _{11},?) & 0 & 0 \\ 
0 & \textrm{f}_{CCC}(\beta _{12},?) & 0 \\ 
0 & 0 & \textrm{f}_{CCC}(\beta _{13},?)%
\end{array}%
\right) ,$

$F_{74}=\left( 
\begin{array}{ccc}
0 & \textrm{f}_{JJC}(b_{12},?) & \textrm{f}_{JJC}(b_{13},?) \\ 
\textrm{f}_{JJC}(b_{11},?) & 0 & \textrm{f}_{JJC}(b_{13},?) \\ 
\textrm{f}_{JJC}(b_{11},?) & \textrm{f}_{JJC}(b_{12},?) & 0%
\end{array}%
\right) ,$

$F_{77}=-F55,$

$F_{78}=\left( 
\begin{array}{ccc}
0 & -\textrm{f}_{JJC}(m_{12},?) & \textrm{f}_{JJC}(m_{13},?) \\ 
\textrm{f}_{JJC}(m_{11},?) & 0 & -\textrm{f}_{JJC}(m_{13},?) \\ 
-\textrm{f}_{JJC}(m_{11},?) & \textrm{f}_{JJC}(m_{12},?) & 0
\end{array}%
\right) $

$\ \ \ \ \ \ +\left( 
\begin{array}{ccc}
0 & -\textrm{f}_{JJC}(t_{12},?) & -\textrm{f}_{JJC}(t_{13},?) \\ 
-\textrm{f}_{JJC}(t_{11},?) & 0 & -\textrm{f}_{JJC}(t_{13},?) \\ 
-\textrm{f}_{JJC}(t_{11},?) & -\textrm{f}_{JJC}(t_{12},?) & 0
\end{array}%
\right) ,$

$F_{79}=\left( 
\begin{array}{c}
-\frac{2}{3}\textrm{f}_{CCC}(\beta _{11},?) \\ 
-\frac{2}{3}\textrm{f}_{CCC}(\beta _{12},?) \\ 
-\frac{2}{3}\textrm{f}_{CCC}(\beta _{13},?)%
\end{array}%
\right) ,$

$F_{81}=\left( 
\begin{array}{c}
-\textrm{f}_{DJJ}(?,b_{11}) \\ 
-\textrm{f}_{DJJ}(?,b_{12})\textrm{d}_{g}\nu \textrm{g}_{d} \\ 
-\textrm{f}_{DJJ}(?,b_{13})\textrm{d}_{g}\nu ^{2}\textrm{g}_{d}%
\end{array}%
\right) ,$

$F_{82}=\left( 
\begin{array}{ccc}
0 & -\frac{1}{2}\textrm{f}_{iJJ}(?,b_{13}) & \frac{1}{2}\textrm{f}_{iJJ}(b_{12},?) \\ 
\frac{1}{2}\textrm{f}_{iJJ}(b_{13},?) &0 & -\frac{1}{2}\textrm{f}_{iJJ}(?,b_{11}) \\ 
-\frac{1}{2}\textrm{f}_{iJJ}(?,b_{12}) & \frac{1}{2}\textrm{f}_{iJJ}(b_{11},?) & 0 \\ 
\end{array}%
\right) $

$\ \ +\left( 
\begin{array}{@{}c@{}c@{}c@{}}
\frac{1}{2}\textrm{f}_{CJJ}(\beta _{12}-\beta _{13},?) & 0 & 0 \\ 
0 & \frac{1}{2}\textrm{f}_{CJJ}(\beta _{13}-\beta _{11},?) & 0 \\ 
0 &0 & \frac{1}{2}\textrm{f}_{CJJ}(\beta _{11}-\beta _{12},?) \\ 
\end{array}%
\right) ,$

$F_{83}=\left( 
\begin{array}{ccc}
0 & \frac{1}{2}\textrm{f}_{CJJ}(?,b_{11}) & \frac{1}{2}\textrm{f}_{CJJ}(?,b_{11}) \\ 
\frac{1}{2}\textrm{f}_{CJJ}(?,b_{12}) & 0 & \frac{1}{2}\textrm{f}_{CJJ}(?,b_{12}) \\ 
\frac{1}{2}\textrm{f}_{CJJ}(?,b_{13}) & \frac{1}{2}\textrm{f}_{CJJ}(?,b_{13}) & 0%
\end{array}%
\right) ,$

$F_{84}=\left( 
\begin{array}{ccc}
0 & \frac{1}{2}\textrm{f}_{iJJ}(?,b_{13}) & \frac{1}{2}\textrm{f}_{iJJ}(b_{12},?) \\ 
\frac{1}{2}\textrm{f}_{iJJ}(b_{13},?) & 0 & \frac{1}{2}\textrm{f}_{iJJ}(?,b_{11}) \\ 
\frac{1}{2}\textrm{f}_{iJJ}(?,b_{12}) & \frac{1}{2}\textrm{f}_{iJJ}(b_{11},?) & 0 \\ 
\end{array}%
\right) $

$\ \ +\left( 
\begin{array}{ccc}
\frac{1}{2}\textrm{f}_{CJJ}(\beta _{12}+\beta _{13},?) & 0 & 0 \\ 
0 & \frac{1}{2}\textrm{f}_{CJJ}(\beta _{13}+\beta _{11},?) &  \\ 
0 & 0 & \frac{1}{2}\textrm{f}_{CJJ}(\beta _{11}+\beta _{12},?) \\ 
\end{array}%
\right) ,$

$F_{87}=\left( 
\begin{array}{ccc}
0 & -\frac{1}{2}\textrm{f}_{CJJ}(?,m_{11}) & \frac{1}{2}\textrm{f}_{CJJ}(?,m_{11}) \\ 
\frac{1}{2}\textrm{f}_{CJJ}(?,m_{12}) & 0 & -\frac{1}{2}\textrm{f}_{CJJ}(?,m_{12}) \\ 
-\frac{1}{2}\textrm{f}_{CJJ}(?,m_{13}) & \frac{1}{2}\textrm{f}_{CJJ}(?,m_{13}) & 0%
\end{array}%
\right) $

$\ \ \ \ \ \ -\left( 
\begin{array}{ccc}
0 & \frac{1}{2}\textrm{f}_{CJJ}(?,t_{11}) & \frac{1}{2}\textrm{f}_{CJJ}(?,t_{11}) \\ 
\frac{1}{2}\textrm{f}_{CJJ}(?,t_{12}) & 0 & \frac{1}{2}\textrm{f}_{CJJ}(?,t_{12}) \\ 
\frac{1}{2}\textrm{f}_{CJJ}(?,t_{13}) & \frac{1}{2}\textrm{f}_{CJJ}(?,t_{13}) & 0%
\end{array}%
\right) ,$

$F_{88}=\left( 
\begin{tabular}{ccc}
\cline{1-1}
\multicolumn{1}{|c}{$\textrm{f}_{DJJ}(D_{1},?)$} & \multicolumn{1}{|c}{$-\frac{1}{2}%
\textrm{f}_{iJJ}(?,m_{13})$} & $\frac{1}{2}\textrm{f}_{iJJ}(m_{12},?)$ \\ \cline{1-2}
$\frac{1}{2}\textrm{f}_{iJJ}(m_{13},?)$ & \multicolumn{1}{|c}{$\textrm{f}_{DJJ}(\textrm{d}_{g}\nu
\textrm{g}_{d}D_{1},?)$} & \multicolumn{1}{|c}{$-\frac{1}{2}\textrm{f}_{iJJ}(?,m_{11})$} \\ 
\cline{2-3}
$-\frac{1}{2}\textrm{f}_{iJJ}(?,m_{12})$ & $\frac{1}{2}\textrm{f}_{iJJ}(m_{11},?)$ & 
\multicolumn{1}{|c|}{$\textrm{f}_{DJJ}(\textrm{d}_{g}\nu ^{2}\textrm{g}_{d}D_{1},?)$} \\ \cline{3-3}
\end{tabular}%
\right) $

$\ \ \ \ \ \ -\left( 
\begin{tabular}{ccc}
\cline{1-1}
\multicolumn{1}{|c}{$\frac{1}{2}\textrm{f}_{CJJ}(-\tau _{11},?)$} & 
\multicolumn{1}{|c}{$\frac{1}{2}\textrm{f}_{iJJ}(?,t_{13})$} & $\frac{1}{2}%
\textrm{f}_{iJJ}(t_{12},?)$ \\ \cline{1-2}
$\frac{1}{2}\textrm{f}_{iJJ}(t_{13},?)$ & \multicolumn{1}{|c}{$\frac{1}{2}%
\textrm{f}_{CJJ}(-\tau _{12},?)$} & \multicolumn{1}{|c}{$\frac{1}{2}\textrm{f}_{iJJ}(?,t_{11})$%
} \\ \cline{2-3}
$\frac{1}{2}\textrm{f}_{iJJ}(?,t_{12})$ & $\frac{1}{2}\textrm{f}_{iJJ}(t_{11},?)$ & 
\multicolumn{1}{|c|}{$\frac{1}{2}\textrm{f}_{CJJ}(\tau _{11}+\tau _{12},?)$} \\ 
\cline{3-3}
\end{tabular}%
\right) $

\ \ \ \ \ $\ -\left( 
\begin{array}{ccc}
\frac{2}{3}\textrm{f}_{CJJ}(\rho _{1},?) & 0 & 0 \\ 
0 & \frac{2}{3}\textrm{f}_{CJJ}(\rho _{1},?) & 0 \\ 
0 & 0 & \frac{2}{3}\textrm{f}_{CJJ}(\rho _{1},?)%
\end{array}%
\right) ,$

$F_{89}=\left( 
\begin{array}{c}
\frac{2}{3}\textrm{f}_{CJJ}(?,b_{11}) \\ 
\frac{2}{3}\textrm{f}_{CJJ}(?,b_{12}) \\ 
\frac{2}{3}\textrm{f}_{CJJ}(?,b_{13})%
\end{array}%
\right) ,$

$F_{95}=\left( 
\begin{array}{ccc}
-\textrm{f}_{CCC}(\beta _{11},?) & -\textrm{f}_{CCC}(\beta _{12},?) & -\textrm{f}_{CCC}(\beta _{13},?)%
\end{array}%
\right) ,$

$F_{96}=\left( 
\begin{array}{ccc}
-2\textrm{f}_{JJC}(b_{11},?) & -2\textrm{f}_{JJC}(b_{12},?) & -2\textrm{f}_{JJC}(b_{13},?)%
\end{array}%
\right) ,$

$F_{97}=\left( 
\begin{array}{ccc}
\textrm{f}_{CCC}(\alpha _{11},?) & \textrm{f}_{CCC}(\alpha _{12},?) & \textrm{f}_{CCC}(\alpha _{13},?)%
\end{array}%
\right) ,$

$F_{98}=\left( 
\begin{array}{ccc}
2\textrm{f}_{JJC}(a_{11},?) & 2\textrm{f}_{JJC}(a_{12},?) & 2\textrm{f}_{JJC}(a_{13},?)%
\end{array}%
\right) ,$

$[\Phi _{1},\Phi _{2}]_{7D}=[D1,D2]$

\ \ \ \ \ \ \ \ \ \ \ \ \ \ \ \ \ $-\textrm{JD}(m_{11},m_{21})-\textrm{d}_{g}\nu
^{2}\textrm{g}_{d}(\textrm{JD}(m_{12},m_{22}))-\textrm{d}_{g}\nu \textrm{g}_{d}(\textrm{JD}(m_{13},m_{23}))$

\ \ \ \ \ \ \ \ \ \ \ \ \ $\ +\textrm{JD}(t_{11},t_{21})+\textrm{d}_{g}\nu
^{2}\textrm{g}_{d}(\textrm{JD}(t_{12},t_{22}))+\textrm{d}_{g}\nu \textrm{g}_{d}(\textrm{JD}(t_{13},t_{23}))$

$\ \ \ \ \ \ \ \ \ \ \ \ \ \ +2(\textrm{JD}(a_{11},b_{21})+\textrm{d}_{g}\nu
^{2}\textrm{g}_{d}(\textrm{JD}(a_{12},b_{22}))+\textrm{d}_{g}\nu \textrm{g}_{d}(\textrm{JD}(a_{13},b_{23})))$

\ \ \ \ \ \ \ \ \ \ \ \ \ $\ -2(\textrm{JD}(a_{21},b_{11})+\textrm{d}_{g}\nu
^{2}\textrm{g}_{d}(\textrm{JD}(a_{22},b_{12}))+\textrm{d}_{g}\nu \textrm{g}_{d}(\textrm{JD}(a_{23},b_{13}))),$

$[\Phi _{1},\Phi _{2}]_{7M\_A_{1}}=\textrm{g}_{d}(D_{1})m_{21}-\textrm{g}_{d}(D_{2})m_{11}+%
\frac{1}{2}(-\overline{m_{12}m_{23}}+\overline{m_{22}m_{13}})$

$\ \ \ \ \ \ \ \ \ \ \ \ \ \ \ \ \ \ \ \ +\frac{1}{2}((\tau _{11}+2\tau
_{12})t_{21}-(\tau _{21}+2\tau _{22})t_{11}-\overline{t_{12}t_{23}}+%
\overline{t_{22}t_{13}})$

$\ \ \ \ \ \ \ \ \ \ \ \ \ \ \ \ \ \ \ \ +(\alpha _{12}-\alpha
_{13})b_{21}-(\beta _{22}-\beta _{23})a_{11}-\overline{a_{12}b_{23}}+%
\overline{b_{22}a_{13}}$

$\ \ \ \ \ \ \ \ \ \ \ \ \ \ \ \ \ \ \ \ -(\alpha _{22}-\alpha
_{23})b_{11}+(\beta _{12}-\beta _{13})a_{21}+\overline{a_{22}b_{13}}-%
\overline{b_{12}a_{23}},$

$[\Phi _{1},\Phi _{2}]_{7M\_A_{2}}=\nu \textrm{g}_{d}(D_{1})m_{22}-\nu
\textrm{g}_{d}(D_{2})m_{12}+\frac{1}{2}(-\overline{m_{13}m_{21}}+\overline{%
m_{23}m_{11}})$

$\ \ \ \ \ \ \ \ \ \ \ \ \ \ \ \ \ \ \ \ +\frac{1}{2}((-2\tau _{11}-\tau
_{12})t_{22}-(-2\tau _{21}-\tau _{22})t_{12}-\overline{t_{23}t_{11}}+%
\overline{t_{13}t_{21}})$

$\ \ \ \ \ \ \ \ \ \ \ \ \ \ \ \ \ \ \ \ +(\alpha _{13}-\alpha
_{11})b_{22}-(\beta _{23}-\beta _{21})a_{12}-\overline{a_{13}b_{21}}+%
\overline{b_{23}a_{11}})$

$\ \ \ \ \ \ \ \ \ \ \ \ \ \ \ \ \ \ \ \ -(\alpha _{23}-\alpha
_{21})b_{12}+(\beta _{13}-\beta _{11})a_{22}+\overline{a_{23}b_{11}}-%
\overline{b_{13}a_{21}},$

$[\Phi _{1},\Phi _{2}]_{7M\_A_{3}}=\nu ^{2}\textrm{g}_{d}(D_{1})m_{23}-\nu
^{2}\textrm{g}_{d}(D_{2})m_{13}+\frac{1}{2}(-\overline{m_{11}m_{22}}+\overline{%
m_{21}m_{12}})$

$\ \ \ \ \ \ \ \ \ \ \ \ \ \ \ \ \ \ \ \ +\frac{1}{2}((\tau _{11}-\tau
_{12})t_{23}-(\tau _{21}-\tau _{22})t_{13}-\overline{t_{11}t_{22}}+\overline{%
t_{21}t_{12}})$

$\ \ \ \ \ \ \ \ \ \ \ \ \ \ \ \ \ \ \ \ +(\alpha _{11}-\alpha
_{12})b_{23}-(\beta _{21}-\beta _{22})a_{13}-\overline{a_{11}b_{22}}+%
\overline{b_{21}a_{12}}$

$\ \ \ \ \ \ \ \ \ \ \ \ \ \ \ \ \ \ \ \ -(\alpha _{21}-\alpha
_{22})b_{13}+(\beta _{11}-\beta _{12})a_{23}+\overline{a_{21}b_{12}}-%
\overline{b_{11}a_{22}},$

$[\Phi _{1},\Phi _{2}]_{7T%
\_E_{1}}=-(m_{12},t_{22})+(m_{13},t_{23})+(m_{22},t_{12})-(m_{23},t_{13})$

$\ \ \ \ \ +\frac{1}{3}(4\alpha _{11}\beta
_{21}-2\alpha _{12}\beta _{22}-2\alpha _{13}\beta
_{23}-4(a_{11},b_{21})+2(a_{12},b_{22})+2(a_{13},b_{23}))$

$\ \ \ \ \ -\frac{1}{3}(4\alpha _{21}\beta
_{11}-2\alpha _{22}\beta _{12}-2\alpha _{23}\beta
_{13}-4(a_{21},b_{11})+2(a_{22},b_{12})+2(a_{23},b_{13})),$

$[\Phi _{1},\Phi _{2}]_{7T%
\_E_{2}}=-(m_{13},t_{23})+(m_{11},t_{21})+(m_{23},t_{13})-(m_{21},t_{11})$

$\ \ \ \ \ +\frac{1}{3}(-2\alpha _{11}\beta
_{21}+4\alpha _{12}\beta _{22}-2\alpha _{13}\beta
_{23}+2(a_{11},b_{21})-4(a_{12},b_{22})+2(a_{13},b_{23}))$

$\ \ \ \ \ -\frac{1}{3}(-2\alpha _{21}\beta
_{11}+4\alpha _{22}\beta _{21}-2\alpha _{23}\beta
_{13}+2(a_{21},b_{11})-4(a_{22},b_{12})+2(a_{23},b_{13})),$

$[\Phi _{1},\Phi _{2}]_{7T%
\_E_{3}}=-(m_{11},t_{21})+(m_{12},t_{22})+(m_{21},t_{11})-(m_{22},t_{12}),$

$\ \ \ \ \ +\frac{1}{3}(-2\alpha _{11}\beta
_{21}-2\alpha _{12}\beta _{22}+4\alpha _{13}\beta
_{23}+2(a_{11},b_{21})+2(a_{12},b_{22})-4(a_{13},b_{23}))$

$\ \ \ \ \ -\frac{1}{3}(-2\alpha _{21}\beta
_{11}-2\alpha _{22}\beta _{12}+4\alpha _{23}\beta
_{13}+2(a_{21},b_{11})+2(a_{22},b_{12})-4(a_{23},b_{13})),$

$[\Phi _{1},\Phi _{2}]_{7T\_F_{1}}=\textrm{g}_{d}(D_{1})t_{21}-\textrm{g}_{d}(D_{2})t_{11}$

$\ \ \ \ \ \ \ \ \ \ \ \ \ \ \ \ \ \ \ \ +\frac{1}{2}(\overline{m_{12}t_{23}}%
-\overline{t_{22}m_{13}}+(-\tau _{21}-2\tau _{22})m_{11}$

$\ \ \ \ \ \ \ \ \ \ \ \ \ \ \ \ \ \ \ \ \ \ \ -\overline{%
m_{22}t_{13}}+\overline{t_{12}m_{23}}-(-\tau _{11}-2\tau _{12})m_{21})$

$\ \ \ \ \ \ \ \ \ \ \ \ \ \ \ \ \ \ \ \ +(\alpha _{12}+\alpha
_{13})b_{21}+(\beta _{22}+\beta _{23})a_{11}+\overline{b_{22}a_{13}}+%
\overline{a_{12}b_{23}}$

$\ \ \ \ \ \ \ \ \ \ \ \ \ \ \ \ \ \ \ \ -(\alpha _{22}+\alpha
_{23})b_{11}-(\beta _{12}+\beta _{13})a_{21}-\overline{b_{12}a_{23}}-%
\overline{a_{22}b_{13}},$

$[\Phi _{1},\Phi _{2}]_{7T\_F_{2}}=\nu \textrm{g}_{d}(D_{1})t_{22}-\nu
\textrm{g}_{d}(D_{2})t_{12}$

$\ \ \ \ \ \ \ \ \ \ \ \ \ \ \ \ \ \ \ \ +\frac{1}{2}(\overline{m_{13}t_{21}}%
-\overline{t_{23}m_{11}}+(2\tau _{21}+\tau _{22})m_{12}$

$\ \ \ \ \ \ \ \ \ \ \ \ \ \ \ \ \ \ \ \ \ \ \ -\overline{%
m_{23}t_{11}}+\overline{t_{13}m_{21}}-(2\tau _{11}+\tau _{12})m_{22})$

\ \ \ \ \ \ \ \ \ \ \ \ \ \ \ \ \ \ \ \ \ +$(\alpha _{13}+\alpha
_{11})b_{22}+(\beta _{23}+\beta _{21})a_{12}+\overline{b_{23}a_{11}}+%
\overline{a_{13}b_{21}}$

$\ \ \ \ \ \ \ \ \ \ \ \ \ \ \ \ \ \ \ \ \ -(\alpha _{23}+\alpha
_{21})b_{12}-(\beta _{13}+\beta _{11})a_{22}-\overline{b_{13}a_{21}}-%
\overline{a_{23}b_{11}},$

$[\Phi _{1},\Phi _{2}]_{7T\_F_{3}}=\nu ^{2}\textrm{g}_{d}(D_{1})t_{23}-\nu
^{2}\textrm{g}_{d}(D_{2})t_{13}$

$\ \ \ \ \ \ \ \ \ \ \ \ \ \ \ \ \ \ \ \ +\frac{1}{2}(\overline{m_{11}t_{22}}%
-\overline{t_{21}m_{12}}+(-\tau _{21}+\tau _{22})m_{13}$

$\ \ \ \ \ \ \ \ \ \ \ \ \ \ \ \ \ \ \ \ \ \ \ -\overline{%
m_{21}t_{12}}+\overline{t_{11}m_{22}}-(-\tau _{11}+\tau _{12})m_{23})$

\ \ \ \ \ \ \ \ \ \ \ \ \ \ \ \ \ \ \ \ \ +$(\alpha _{11}+\alpha
_{12})b_{23}+(\beta _{21}+\beta _{22})a_{13}+\overline{b_{21}a_{12}}+%
\overline{a_{11}b_{22}}$

$\ \ \ \ \ \ \ \ \ \ \ \ \ \ \ \ \ \ \ \ \ -(\alpha _{21}+\alpha
_{22})b_{13}-(\beta _{11}+\beta _{12})a_{23}-\overline{b_{11}a_{22}}-%
\overline{a_{21}b_{12}},$

$[\Phi _{1},\Phi _{2}]_{7A%
\_E_{1}}=(m_{13},a_{23})-(m_{12},a_{22})-(m_{23},a_{13})+(m_{22},a_{12})$

\ \ \ \ \ \ \ \ \ \ \ \ \ \ \ \ \ \ \ \ \ \ $+\tau _{11}\alpha
_{21}+(t_{12},a_{22})+(t_{13},a_{23})-\tau _{21}\alpha
_{11}-(t_{22},a_{12})-(t_{23},a_{13})$

$\ \ \ \ \ \ \ \ \ \ \ \ \ \ \ \ \ \ \ +\frac{2}{3}\rho _{1}\alpha _{21}-%
\frac{2}{3}\rho _{2}\alpha _{11},$

$[\Phi _{1},\Phi _{2}]_{7A%
\_E_{2}}=(m_{11},a_{21})-(m_{13},a_{23})-(m_{21},a_{11})+(m_{23},a_{13})$

\ \ \ \ \ \ \ \ \ \ \ \ \ \ \ \ \ \ \ \ \ \ $+\tau _{12}\alpha
_{22}+(t_{13},a_{23})+(t_{11},a_{21})-\tau _{22}\alpha
_{12}-(t_{23},a_{13})-(t_{21},a_{11})$

$\ \ \ \ \ \ \ \ \ \ \ \ \ \ \ \ \ \ \ +\frac{2}{3}\rho _{1}\alpha _{22}-%
\frac{2}{3}\rho _{2}\alpha _{12},$

$[\Phi _{1},\Phi _{2}]_{7A%
\_E_{3}}=(m_{12},a_{22})-(m_{11},a_{21})-(m_{22},a_{12})+(m_{21},a_{11})$

\ \ \ \ \ \ \ \ \ \ \ \ \ \ \ \ \ \ \ $+(-\tau _{11}-\tau _{12})\alpha
_{23}+(t_{11},a_{21})+(t_{12},a_{22})$

$\ \ \ \ \ \ \ \ \ \ \ \ \ \ \ \ \ \ \ -(-\tau _{21}-\tau _{22})\alpha
_{13}-(t_{21},a_{11})-(t_{22},a_{12})$

$\ \ \ \ \ \ \ \ \ \ \ \ \ \ \ \ \ \ \ +\frac{2}{3}\rho _{1}\alpha _{23}-%
\frac{2}{3}\rho _{2}\alpha _{13}$

$[\Phi _{1},\Phi _{2}]_{7A\_F_{1}}=\textrm{g}_{d}(D_{1})a_{21}-\textrm{g}_{d}(D_{2})a_{11}$

$\ \ \ \ \ \ \ \ \ \ \ \ \ \ \ \ \ \ \ \ \ +\frac{1}{2}((\alpha _{23}-\alpha
_{22})m_{11}-\overline{a_{22}m_{13}}+\overline{m_{12}a_{23}}$

$\ \ \ \ \ \ \ \ \ \ \ \ \ \ \ \ \ \ \ \ \ \ \ \ -(\alpha
_{13}-\alpha _{12})m_{21}+\overline{a_{12}m_{23}}-\overline{m_{22}a_{13}})$

\ \ \ \ \ \ \ \ \ \ \ \ \ \ \ \ \ \ \ \ \ $+\frac{1}{2}(-\tau
_{11}a_{21}+(\alpha _{22}+\alpha _{23})t_{11}+\overline{a_{22}t_{13}}+%
\overline{t_{12}a_{23}})+\frac{2}{3}\rho _{1}a_{21}$

\ \ \ \ \ \ \ \ \ \ \ \ \ \ \ \ \ \ \ $\ \ -\frac{1}{2}(-\tau
_{21}a_{11}+(\alpha _{12}+\alpha _{13})t_{21}+\overline{a_{12}t_{23}}+%
\overline{t_{22}a_{13}})-\frac{2}{3}\rho _{2}a_{11},$

$[\Phi _{1},\Phi _{2}]_{7A\_F_{2}}=\nu \textrm{g}_{d}(D_{1})a_{22}-\nu
\textrm{g}_{d}(D_{2})a_{12}$

$\ \ \ \ \ \ \ \ \ \ \ \ \ \ \ \ \ \ \ \ +\frac{1}{2}((\alpha _{21}-\alpha
_{23})m_{12}-\overline{a_{23}m_{11}}+\overline{m_{13}a_{21}}$

\ \ \ \ \ \ \ \ \ \ \ \ \ \ \ \ \ \ \ \ \ \ \ $-(\alpha
_{11}-\alpha _{13})m_{22}+\overline{a_{13}m_{21}}-\overline{m_{23}a_{11}})$

\ \ \ \ \ \ \ \ \ \ \ \ \ \ \ \ \ \ \ \ \ $+\frac{1}{2}(-\tau
_{12}a_{22}+(\alpha _{23}+\alpha _{21})t_{12}+\overline{a_{23}t_{11}}+%
\overline{t_{13}a_{21}})+\frac{2}{3}\rho _{1}a_{22}$

$\ \ \ \ \ \ \ \ \ \ \ \ \ \ \ \ \ \ \ \ \ -\frac{1}{2}(-\tau
_{22}a_{12}+(\alpha _{13}+\alpha _{11})t_{22}+\overline{a_{13}t_{21}}+%
\overline{t_{23}a_{11}})-\frac{2}{3}\rho _{2}a_{12},$

$[\Phi _{1},\Phi _{2}]_{7A\_F_{3}}=\nu ^{2}\textrm{g}_{d}(D_{1})a_{23}-\nu
^{2}\textrm{g}_{d}(D_{2})a_{13}$

$\ \ \ \ \ \ \ \ \ \ \ \ \ \ \ \ \ \ +\frac{1}{2}((\alpha _{22}-\alpha
_{21})m_{13}-\overline{a_{21}m_{12}}+\overline{m_{11}a_{22}}$

\ \ \ \ \ \ \ \ \ \ \ \ \ \ \ \ \ \ \ \ \ $-(\alpha
_{12}-\alpha _{11})m_{23}+\overline{a_{11}m_{22}}-\overline{m_{21}a_{12}})$

\ \ \ \ \ \ \ \ \ \ \ \ \ \ \ \ \ \ \ \ $+\frac{1}{2}((\tau _{11}+\tau
_{12})a_{23}+(\alpha _{21}+\alpha _{22})t_{13}+\overline{a_{21}t_{12}}+%
\overline{t_{11}a_{22}})+\frac{2}{3}\rho _{1}a_{23}$

$\ \ \ \ \ \ \ \ \ \ \ \ \ \ \ \ \ \ \ \ \ -\frac{1}{2}((\tau _{21}+\tau
_{22})a_{13}+(\alpha _{11}+\alpha _{12})t_{23}+\overline{a_{11}t_{22}}+%
\overline{t_{21}a_{12}})-\frac{2}{3}\rho _{2}a_{13},$

$[\Phi _{1},\Phi _{2}]_{7B%
\_E_{1}}=(m_{13},b_{23})-(m_{12},b_{22})-(m_{23},b_{13})+(m_{22},b_{12})$

$\ \ \ \ \ \ \ \ \ \ \ \ \ \ \ \ \ \ \ \ -\tau _{11}\beta
_{21}-(t_{12},b_{22})-(t_{13},b_{23})+\tau _{21}\beta
_{11}+(t_{22},b_{12})+(t_{23},b_{13})$

$\ \ \ \ \ \ \ \ \ \ \ \ \ \ \ \ \ \ -\frac{2}{3}\rho _{1}\beta _{21}+%
\frac{2}{3}\rho _{2}\beta _{11},$

$[\Phi _{1},\Phi _{2}]_{7B%
\_E_{2}}=(m_{11},b_{21})-(m_{13},b_{23})-(m_{21},b_{11})+(m_{23},b_{13})$

\ \ \ \ \ \ \ \ \ \ \ \ \ \ \ \ \ \ $\ \ -\tau _{12}\beta
_{22}-(t_{13},b_{23})-(t_{11},b_{21})+\tau _{22}\beta
_{12}+(t_{23},b_{13})+(t_{21},b_{11})$

$\ \ \ \ \ \ \ \ \ \ \ \ \ \ \ \ \ \ -\frac{2}{3}\rho _{1}\beta _{22}+%
\frac{2}{3}\rho _{2}\beta _{12},$

$[\Phi _{1},\Phi _{2}]_{7B%
\_E_{3}}=(m_{12},b_{22})-(m_{11},b_{21})-(m_{22},b_{12})+(m_{21},b_{11})$

\ \ \ \ \ \ \ \ \ \ \ \ \ \ \ \ \ \ \ \ \ \ $-(-\tau _{11}-\tau _{12})\beta
_{23}-(t_{11},b_{21})-(t_{12},b_{22})$

$\ \ \ \ \ \ \ \ \ \ \ \ \ \ \ \ \ \ \ \ \ \ +(-\tau _{21}-\tau _{22})\beta
_{13}+(t_{21},b_{11})+(t_{22},b_{12})$

$\ \ \ \ \ \ \ \ \ \ \ \ \ \ \ \ \ \ \ \ \ -\frac{2}{3}\rho _{1}\beta _{23}+%
\frac{2}{3}\rho _{2}\beta _{13},$

$[\Phi _{1},\Phi _{2}]_{7B\_F_{1}}=\textrm{g}_{d}(D_{1})b_{21}-\textrm{g}_{d}(D_{2})b_{11}$

$\ \ \ \ \ \ \ \ \ \ \ \ \ \ \ \ \ \ \ \ \ +\frac{1}{2}((\beta _{23}-\beta
_{22})m_{11}-\overline{b_{22}m_{13}}+\overline{m_{12}b_{23}}$

\ \ \ \ \ \ \ \ \ \ \ \ \ \ \ \ \ \ \ \ \ \ \ \ $-(\beta
_{13}-\beta _{12})m_{21}+\overline{b_{12}m_{23}}-\overline{m_{22}b_{13}})$

$\ \ \ \ \ \ \ \ \ \ \ \ \ \ \ \ \ \ \ \ \ -\frac{1}{2}(-\tau
_{11}b_{21}+(\beta _{22}+\beta _{23})t_{11}+\overline{b_{22}t_{13}}+%
\overline{t_{12}b_{23}})-\frac{2}{3}\rho _{1}b_{21}$

\ \ \ \ \ \ \ \ \ \ \ \ \ \ \ \ \ \ \ $\ \ +\frac{1}{2}(-\tau
_{21}b_{11}+(\beta _{12}+\beta _{13})t_{21}+\overline{b_{12}t_{23}}+%
\overline{t_{22}b_{13}})+\frac{2}{3}\rho _{2}b_{11},$

$[\Phi _{1},\Phi _{2}]_{7B\_F_{2}}=\nu \textrm{g}_{d}(D_{1})b_{22}-\nu
\textrm{g}_{d}(D_{2})b_{12}$

$\ \ \ \ \ \ \ \ \ \ \ \ \ \ \ \ \ \ \ \ \ +\frac{1}{2}((\beta _{21}-\beta
_{23})m_{12}-\overline{b_{23}m_{11}}+\overline{m_{13}b_{21}}$

\ \ \ \ \ \ \ \ \ \ \ \ \ \ \ \ \ \ \ \ \ \ \ \ $-(\beta
_{11}-\beta _{13})m_{22}+\overline{b_{13}m_{21}}-\overline{m_{23}b_{11}})$

$\ \ \ \ \ \ \ \ \ \ \ \ \ \ \ \ \ \ \ \ \ -\frac{1}{2}(-\tau
_{12}b_{22}+(\beta _{23}+\beta _{21})t_{12}+\overline{b_{23}t_{11}}+%
\overline{t_{13}b_{21}})-\frac{2}{3}\rho _{1}b_{22}$

$\ \ \ \ \ \ \ \ \ \ \ \ \ \ \ \ \ \ \ \ \ +\frac{1}{2}(-\tau
_{22}b_{12}+(\beta _{13}+\beta _{11})t_{22}+\overline{b_{13}t_{21}}+%
\overline{t_{23}b_{11}})+\frac{2}{3}\rho _{2}b_{12},$

$[\Phi _{1},\Phi _{2}]_{7B\_F_{3}}=\nu ^{2}\textrm{g}_{d}(D_{1})b_{23}-\nu
^{2}\textrm{g}_{d}(D_{2})b_{13}$

$\ \ \ \ \ \ \ \ \ \ \ \ \ \ \ \ \ \ \ \ \ +\frac{1}{2}((\beta _{22}-\beta
_{21})m_{13}-\overline{b_{21}m_{12}}+\overline{m_{11}b_{22}}$

\ \ \ \ \ \ \ \ \ \ \ \ \ \ \ \ \ \ \ \ \ \ \ \ $-(\beta
_{12}-\beta _{11})m_{23}+\overline{b_{11}m_{22}}-\overline{m_{21}b_{12}})$

$\ \ \ \ \ \ \ \ \ \ \ \ \ \ \ \ \ \ \ \ \ -\frac{1}{2}((\tau _{11}+\tau
_{12})b_{23}+(\beta _{21}+\beta _{22})t_{13}+\overline{b_{21}t_{12}}+%
\overline{t_{11}b_{22}})-\frac{2}{3}\rho _{1}b_{23}$

$\ \ \ \ \ \ \ \ \ \ \ \ \ \ \ \ \ \ \ \ \ +\frac{1}{2}((\tau _{21}+\tau
_{22})b_{13}+(\beta _{11}+\beta _{12})t_{23}+\overline{b_{11}t_{22}}+%
\overline{t_{21}b_{12}})+\frac{2}{3}\rho _{2}b_{13},$

$[\Phi _{1},\Phi _{2}]_{7\rho }=\alpha _{11}\beta _{21}+\alpha
_{12}\beta _{22}+\alpha _{13}\beta _{23}-\alpha _{21}\beta _{11}-\alpha
_{22}\beta _{12}-\alpha _{23}\beta _{13}$

$\ \ \ \ \ \ \ \ \ \ \ \ \ \ \ \ \ \
+2(a_{11},b_{21})+2(a_{12},b_{22})+2(a_{13},b_{23})$

\ \ \ \ \ \ \ \ \ \ \ \ \ \ \ \ \ \ $-2(a_{21},b_{11})-2(a_{22},b_{12})-2(a_{23},b_{13}). $

\bigskip

\emph{Proof.} Using \emph{Proposition 4.9} and \emph{Definition 6.1} of
mappings $\textrm{f}_{DDD},\textrm{f}_{JJD},$

\noindent
$\textrm{f}_{DJJ},\textrm{f}_{iJJ},\textrm{f}_{JJC},\textrm{f}_{CJJ},\textrm{f}_{CCC},$
we have the above expression. \ \ \ \ \emph{Q.E.D.}

\bigskip

\emph{Theorem 9.2. } The image of the adjoit representation $ad($\gR$_{7})$
of \gR$_{7}^{}$ is expressed by

$ad(D,M,T,A,B,\rho )=ad(\textrm{fv}(D,M,T,A,B,\rho ))$

$=\left( 
\begin{array}{ccc}
M_{\phi \phi } & M_{\phi AB} & 0 \\ 
M_{AB\phi } & M_{ABAB} & M_{AB\rho } \\ 
0 & M_{\rho AB} & 0%
\end{array}%
\right) \in M(133\times 133,\R),$

\noindent
where

$M_{\phi \phi }=ad(D,M,T)\in M(78\times 78,\R),$

$M_{\phi AB}=%
\begin{array}{@{}l@{}l@{}}
& \text{  \ \ \ \ \ \ \ \ \ \ \  \ \ \  }\ \ \ \ \ 3\text{\ \ \ \ \  }\ \ \ \
\ \ \ \ \ \ \ \ \ \ \ \ \ \ \ \ \ \ \ \ \ \ \ \ \ \ \text{24 } \\ 
\begin{array}{@{}l@{}}
28 \\ 
24 \\ 
2 \\ 
24%
\end{array}
& \left( 
\begin{array}{@{}cc@{}}
0 & -2\textrm{M}_{R}\textrm{Dl}(B) \\ 
^{t}\textrm{MJC}(B) & \textrm{M}_{D}^{-}\textrm{E}(B)+\textrm{MI}(B) \\ 
-\frac{4}{3}M_{D3}(B)+\frac{2}{3}\textrm{MCCC}_{3}(B) & \frac{4}{3}\textrm{M}_{D}\textrm{C}_{3}(B)-%
\frac{2}{3}\textrm{M}^{+}\textrm{JC}_{3}(B) \\ 
-\text{ }^{t}\textrm{MJC}(B) & -\textrm{M}_{D}^{+}\textrm{E}(A)+\textrm{M}^{+}\textrm{I}(A)%
\end{array}%
\right. 
\end{array}%
$

$\ \ \ \ \ \ \ \ \ \ \ \ \ \ \ \ \ \ \ \ \ \ 
\begin{array}{l@{}}
\text{ \ \ \ \ \ \ \ \ \ \  \ \  }\ \ \ \ \ 3\text{\ \ \ \ \  }\ \ \ \ \ \
\ \ \ \ \ \ \ \ \ \ \ \ \ \ \ \ \ \ \ \ \ \ \ \ \ \text{24 } \\ 
\left. 
\begin{array}{cc}
0 & 2\textrm{M}_{R}\textrm{Dl}(A) \\ 
^{t}\textrm{MJC}(A) & \textrm{M}_{D}^{-}\textrm{E}(A)+\textrm{MI}(A) \\ 
\frac{4}{3}M_{D3}(A)-\frac{2}{3}\textrm{MCCC}_{3}(A) & -\frac{4}{3}\textrm{M}_{D}\textrm{C}_{3}(A)+%
\frac{2}{3}\textrm{M}^{+}\textrm{JC}_{3}(A) \\ 
^{t}\textrm{M}^{+}\textrm{JC}(A) & \textrm{M}_{D}^{+}\textrm{E}(A)+\textrm{M}^{+}\textrm{I}(A)%
\end{array}%
\right) 
\end{array}$%

$\in M(78\times 54,\R),$

$M_{AB\phi }=%
\begin{array}{@{}l@{}l@{}}
& \text{ \ \ \ \ \ \ \ \ \ \ 28 \ \ \ \ \ \ \ \ \ \ \ \ \ \ \ \ \ \ \ \ \ \ \
\ \ \ 24} \\ 
\begin{array}{@{}l@{}}
3 \\ 
24 \\ 
3 \\ 
24%
\end{array}
& \left( 
\begin{array}{@{}cc@{}}
0 & \textrm{MJC}(A) \\ 
-\text{ }^{t}\textrm{M}_{R}\textrm{Dl}(A) & \frac{1}{2}\textrm{M}_{D}^{-}\textrm{E}(A)-\frac{1}{2}\textrm{MI}(A) \\ 
0 & \textrm{MJC}(B) \\ 
-\text{ }^{t}\textrm{M}_{R}\textrm{Dl}(B) & \frac{1}{2}\textrm{M}_{D}^{-}\textrm{E}(B)-\frac{1}{2}\textrm{MI}(B)%
\end{array}%
\right. 
\end{array}%
$

$\ \ \ \ \ \ \ \ \ \ \ \ \ \ \ \ 
\begin{array}{l}
\text{ \ \ \ \ \ \ \ \ \ 2 \ \ \ \ \ \ \ \ \ \ \ \ \ \ \ \ \ \ \ \ \ \ \ \
\ \ 24} \\ 
\left. 
\begin{array}{cc}
-\textrm{M}_{D}^{3}(A) & -\textrm{MJC}(A) \\ 
-\frac{1}{2}^{t}\textrm{M}^{+}\textrm{JC}^{3}(A) & -\frac{1}{2}\textrm{M}_{D}^{+}\textrm{E}(A)-\frac{1}{2}%
\textrm{M}^{+}\textrm{I}(A) \\ 
\textrm{M}_{D}^{3}(B) & \textrm{M}^{+}\textrm{JC}(B) \\ 
^{t}\textrm{M}^{+}\textrm{JC}^{3}(B) & \frac{1}{2}\textrm{M}_{D}^{+}\textrm{E}(B)+\frac{1}{2}\textrm{M}^{+}\textrm{I}(B)%
\end{array}%
\right) 
\end{array}$%

$\in M(54\times 78,\R),$

$M_{ABAB}=%
\begin{array}{@{}l@{}l@{}}
& \text{ \ \ \ \ \ \ \ \ \ \ \ 3 \ \ \ \ \ \ \ \ \ \ \ \ \ \ \ \ \ \ \ \ \
\ 24} \\ 
\begin{array}{@{}c@{}}
3 \\ 
24 \\ 
3 \\ 
24%
\end{array}
& \left( 
\begin{array}{cc}
0 & -\textrm{MJC}(M) \\ 
-\frac{1}{2}^{t}\textrm{MJC}(M) & \textrm{M}_{D}\textrm{Jl}(D)-\frac{1}{2}\textrm{MI}(M) \\ 
0 & 0 \\ 
0 & 0%
\end{array}%
\right. 
\end{array}%
$

$\ \ \ \ \ \ \ \ \ \ \ \ \ \ \ \ \ \ \ 
\begin{array}{l}
\text{\ \ \ \ \ \ \ \ \ \ 3 \ \ \ \ \ \ \ \ \ \ \ \ \ \ \ \ \ \ \ \ \ \ \ 24}
\\ 
\left. 
\begin{array}{cc}
0 & 0 \\ 
0 & 0 \\ 
0 & -\textrm{MJC}(M) \\ 
-\frac{1}{2}^{t}\textrm{MJC}(M) & \textrm{M}_{D}\textrm{Jl}(D)-\frac{1}{2}\textrm{MI}(M)%
\end{array}%
\right) 
\end{array}%
$

$\ \ \ \ \ \ \ \ \ \ \ \ +%
\begin{array}{l@{}l@{}}
& \text{ \ \ \ \ \ \ \ \ \ \ \ \ 3 \ \ \ \ \ \ \ \ \ \ \ \ \ \ \ \ \ \ \ \ \
\ \ \ \ \ \ \ \ \ 24} \\ 
\begin{array}{@{}c@{}}
3 \\ 
24 \\ 
3 \\ 
24%
\end{array}
& \left( 
\begin{array}{cc}
\frac{2}{3}\rho E+\textrm{M}_{D}(T) & \textrm{M}^{+}\textrm{JC}(T) \\ 
\frac{1}{2}^{t}\textrm{M}^{+}\textrm{JC}(T) & \frac{1}{2}\textrm{M}_{D}^{+}\textrm{E}(T)+\frac{1}{2}\textrm{M}^{+}\textrm{I}(T)+%
\frac{2}{3}\rho E \\ 
0 & 0 \\ 
0 & 0%
\end{array}%
\right. 
\end{array}%
$

$\ \ \ \ \ \ \ \ \ \ \ \ \ \ \ \ \ \ \ \ 
\begin{array}{l}
\text{ \ \ \ \ \ \ \ \ \ \ \ 3 \ \ \ \ \ \ \ \ \ \ \ \ \ \ \ \ \ \ \ \ \ \
\ \ \ \ \ \ \ \ \ \ 24} \\ 
\left. 
\begin{array}{cc}
0 & 0 \\ 
0 & 0 \\ 
-\frac{2}{3}\rho E-\textrm{M}_{D}(T) & -\textrm{M}^{+}\textrm{JC}(T) \\ 
-\frac{1}{2}^{t}\textrm{M}^{+}\textrm{JC}(T) & -\frac{1}{2}\textrm{M}_{D}^{+}\textrm{E}(T)-\frac{1}{2}\textrm{M}^{+}\textrm{I}(T)-%
\frac{2}{3}\rho E%
\end{array}%
\right) 
\end{array}%
$

$\in M(54\times 54,\R),$

$M_{AB\rho }=%
\begin{array}{@{}l@{}l@{}}
\begin{array}{@{}c@{}}
3 \\ 
24 \\ 
3 \\ 
24%
\end{array}
& \left( 
\begin{array}{c}
-\frac{2}{3}\textrm{M}_{C}\textrm{CCC}(A) \\ 
-\frac{2}{3}\textrm{M}_{C}\textrm{CJ}(A) \\ 
-\frac{2}{3}\textrm{M}_{C}\textrm{CCC}(B) \\ 
\frac{2}{3}\textrm{M}_{C}\textrm{CJ}(B)%
\end{array}%
\right)%
\end{array}%
\in M(54\times 1,\R),$

$M_{\rho AB}=%
\begin{array}[b]{@{}l@{}}
\begin{array}{llll@{}}
\text{ \ \ \ \ \ \ \ \ \ \ \ 3} & \text{ \ \ \ \ \ \ \ \ \ \ \ \ \ \ 24} & 
\text{ \ \ \ \ \ \ \ \ \ \ \ \ \ \ 3} & \text{ \ \ \ \ \ \ \ \ \ \ \ \ \ 24}%
\end{array}
\\ 
\left( 
\begin{array}{@{}cccc@{}}
-\text{ }^{t}\textrm{M}_{C}\textrm{CCC}(B) & -2^{t}\textrm{M}_{C}\textrm{CJ}(B) & \text{ }^{t}\textrm{M}_{C}\textrm{CCC}(A) & 
2^{t}\textrm{M}_{C}\textrm{CJ}(A)%
\end{array}%
\right)%
\end{array}$%

$\in M(1\times 54,\R),$

$D=\sum\limits_{0\leq i<j\leq 7}d_{ij}D_{ij} \ (d_{ij}\in
\R),M=A_{1}(m_{1})+A_{2}(m_{2})+A_{3}(m_{3}),$

$m_{i}=\sum\limits_{j=0}^{7}m_{ij}e_{j}(m_{ij}\in \R),$

$T=\tau _{1}E_{1}+\tau
_{2}E_{2}+F_{1}(t_{1})+F_{2}(t_{2})+F_{3}(t_{3}),t_{i}=\sum\limits_{j=0}^{7}$%
t$_{ij}e_{j}(t_{ij}\in \R),\tau _{1},\tau _{2}\in \R,$

$A=\alpha _{1}E_{1}+\alpha _{2}E_{2}+\alpha
_{3}E_{3}+F_{1}(a_{1})+F_{2}(a_{2})+F_{3}(a_{3}),$

$B=\beta _{1}E_{1}+\beta _{2}E_{2}+\beta
_{3}E_{3}+F_{1}(b_{1})+F_{2}(b_{2})+F_{3}(b_{3}),$

$a_{1},a_{2},a_{3},b_{1},b_{2},b_{3}\in $\gC $^{},\alpha _{1},\alpha
_{2},\alpha _{3},\beta _{1},\beta _{2},\beta _{3},\rho \in \R.$

\noindent
Since \gR$_{7}^{}$ and \ge$_{7,1}^{}$ are isomorphic, 
$ad(D,M,T,A,B,\rho ) $ is also the image of the adjoit representation of \ge$_{7,1}^{}.$

\bigskip

\emph{Proof.} Using \emph{Lemma 9.1, Lemma
6.10, 6.11, 6.12, 6.13, 6.14, 6.15, 6.16,}

\noindent
\emph{6.17, 6.18, 6.19, 6.20,} and \emph{Lemma 6.24,} we have

$ad(D,M,T,A,B,\rho )=$

$\left( 
\begin{tabular}{ccccccccc}
\cline{1-4}
\multicolumn{1}{|c}{$M_{11}$} & \multicolumn{1}{|c}{$M_{12}$} & 
\multicolumn{1}{|c}{$0$} & \multicolumn{1}{|c}{$M_{14}$} & \multicolumn{1}{|c}{0
} & $M_{16}$ & 0 & $M_{18}$ & 0 \\ \cline{1-4}
\multicolumn{1}{|c}{$M_{21}$} & \multicolumn{1}{|c}{$M_{22}$} & 
\multicolumn{1}{|c}{$M_{23}$} & \multicolumn{1}{|c}{$M_{24}$} & 
\multicolumn{1}{|c}{$M_{25}$} & $M_{26}$ & $M_{27}$ & $M_{28}$ & 0 \\ \cline{1-4}
\multicolumn{1}{|c}{0} & \multicolumn{1}{|c}{$M_{32}$} & \multicolumn{1}{|c}{0}
& \multicolumn{1}{|c}{$M_{34}$} & \multicolumn{1}{|c}{$M_{35}$} & $M_{36}$ & $M_{37}$ & $%
M_{38}$ & 0 \\ \cline{1-4}
\multicolumn{1}{|c}{$M_{41}$} & \multicolumn{1}{|c}{$M_{42}$} & 
\multicolumn{1}{|c}{$M_{43}$} & \multicolumn{1}{|c}{$M_{44}$} & 
\multicolumn{1}{|c}{$M_{45}$} & $M_{46}$ & $M_{47}$ & $M_{48}$ & 0 \\ \cline{1-4}
0 & $M_{52}$ & $M_{53}$ & $M_{54}$ & $M_{55}$ & $M_{56}$ & 0 & 0 & $M_{59}$ \\ 
$M_{61}$ & $M_{62}$ & $M_{63}$ & $M_{64}$ & $M_{65}$ & $M_{66}$ & 0 & 0 & $M_{69}$ \\ 
0 & $M_{72}$ & $M_{73}$ & $M_{74}$ & 0 & 0 & $M_{77}$ & $M_{78}$ & $M_{79}$ \\ 
$M_{81}$ & $M_{82}$ & $M_{83}$ & $M_{84}$ & 0 & 0 & $M_{87}$ & $M_{88}$ & $M_{89}$ \\ 
0 & 0 & 0 & 0 & $M_{95}$ & $M_{96}$ & $M_{97}$ & $M_{98}$ & 0%
\end{tabular}%
\right) $

$\in M(133\times 133,\R),$

$\left( 
\begin{array}{cccc}
M_{11} & M_{12} & 0 & M_{14} \\ 
M_{21} & M_{22} & M_{23} & M_{24} \\ 
0 & M_{32} & 0 & M_{34} \\ 
M_{41} & M_{42} & M_{43} & M_{44}%
\end{array}%
\right) =ad(D,M,T)\in M(78\times 78,\R),$

$M_{16}=\left( 
\begin{array}{ccc}
-2\textrm{MDr}(b_{1}) & -2Mv^{2}\textrm{MDr}(b_{2}) & -2Mv\textrm{MDr}(b_{3})%
\end{array}%
\right) $

$\in M(28\times 24,\R),$

$M_{18}=\left( 
\begin{array}{ccc}
2\textrm{MDl}(a_{1}) & 2Mv^{2}\textrm{MDl}(a_{2}) & 2Mv\textrm{MDl}(a_{3})%
\end{array}%
\right) \in M(28\times 24,\R),$

$M_{25}=\left( 
\begin{array}{ccc}
0 & -^{t}\textrm{MC}(b_{1}) & ^{t}\textrm{MC}(b_{1}) \\ 
^{t}\textrm{MC}(b_{2}) & 0 & -^{t}\textrm{MC}(b_{2}) \\ 
-^{t}\textrm{MC}(b_{3}) & ^{t}\textrm{MC}(b_{3}) & 0%
\end{array}%
\right) \in M(24\times 3,\R),$

$M_{26}=\left( 
\begin{tabular}{ccc}
\cline{1-1}
\multicolumn{1}{|c}{$\textrm{ME}(\beta _{2}-\beta _{3})$} & \multicolumn{1}{|c}{$%
\textrm{MIr}(b_{3})$} & $-\textrm{MIl}(b_{2})$ \\ \cline{1-2}
$-\textrm{MIl}(b_{3})$ & \multicolumn{1}{|c}{$\textrm{ME}(\beta _{3}-\beta _{1})$} & 
\multicolumn{1}{|c}{$\textrm{MIr}(b_{1})$} \\ \cline{2-3}
$\textrm{MIr}(b_{2})$ & $-\textrm{MIl}(b_{1})$ & \multicolumn{1}{|c|}{$\textrm{ME}(\beta _{1}-\beta
_{2})$} \\ \cline{3-3}
\end{tabular}%
\right) $

$\in M(24\times 24,\R),$

$M_{27}=\left( 
\begin{array}{ccc}
0 & -^{t}\textrm{MC}(a_{1}) & ^{t}\textrm{MC}(a_{1}) \\ 
^{t}\textrm{MC}(a_{2}) & 0 & -^{t}\textrm{MC}(a_{2}) \\ 
-^{t}\textrm{MC}(a_{3}) & ^{t}\textrm{MC}(a_{3}) & 0%
\end{array}%
\right) \in M(24\times 3,\R),$

$M_{28}=\left( 
\begin{tabular}{ccc}
\cline{1-1}
\multicolumn{1}{|c}{$\textrm{ME}(\alpha _{2}-\alpha _{3})$} & \multicolumn{1}{|c}{$%
\textrm{MIr}(a_{3})$} & $-\textrm{MIl}(a_{2})$ \\ \cline{1-2}
$-\textrm{MIl}(a_{3})$ & \multicolumn{1}{|c}{$\textrm{ME}(\alpha _{3}-\alpha _{1})$} & 
\multicolumn{1}{|c}{$\textrm{MIr}(a_{1})$} \\ \cline{2-3}
$\textrm{MIr}(a_{2})$ & $-\textrm{MIl}(a_{1})$ & \multicolumn{1}{|c|}{$\textrm{ME}(\alpha _{1}-\alpha
_{2})$} \\ \cline{3-3}
\end{tabular}%
\right) $

$\in M(24\times 24,\R),$

$M_{35}=\left( 
\begin{array}{ccc}
-\frac{4}{3}\beta _{1} & \frac{2}{3}\beta _{2} & \frac{2}{3}\beta _{3} \\ 
\frac{2}{3}\beta _{1} & -\frac{4}{3}\beta _{2} & \frac{2}{3}\beta _{3}%
\end{array}%
\right) \in M(2\times 3,\R),$

$M_{36}=\left( 
\begin{array}{ccc}
\frac{4}{3}\textrm{MC}(b_{1}) & \text{-}\frac{2}{3}\textrm{MC}(b_{2}) & \text{-}\frac{2}{3}%
\textrm{MC}(b_{3}) \\ 
\text{-}\frac{2}{3}\textrm{MC}(b_{1}) & \frac{4}{3}\textrm{MC}(b_{2}) & \text{-}\frac{2}{3}%
\textrm{MC}(b_{3})%
\end{array}%
\right) \in M(2\times 24,\R),$

$M_{37}=\left( 
\begin{array}{ccc}
\frac{4}{3}\mathbf{\alpha }_{1} & \mathbf{-}\frac{2}{3}\mathbf{\alpha }_{2}
& \mathbf{-}\frac{2}{3}\mathbf{\alpha }_{3} \\ 
\mathbf{-}\frac{2}{3}\mathbf{\alpha }_{1} & \frac{4}{3}\mathbf{\alpha }_{2}
& \mathbf{-}\frac{2}{3}\mathbf{\alpha }_{3}%
\end{array}%
\right) \in M(2\times 3,\R),$

$M_{38}=\left( 
\begin{array}{ccc}
-\frac{4}{3}\textrm{MC}(a_{1}) & \frac{2}{3}\textrm{MC}(a_{2}) & \frac{2}{3}\textrm{MC}(a_{3}) \\ 
\frac{2}{3}\textrm{MC}(a_{1}) & -\frac{4}{3}\textrm{MC}(a_{2}) & \frac{2}{3}\textrm{MC}(a_{3})%
\end{array}%
\right) \in M(2\times 24,\R),$

$M_{45}=\left( 
\begin{array}{ccc}
0 & -^{t}\textrm{MC}(b_{1}) & -^{t}\textrm{MC}(b_{1}) \\ 
-^{t}\textrm{MC}(b_{2}) & 0 & -^{t}\textrm{MC}(b_{2}) \\ 
-^{t}\textrm{MC}(b_{3}) & -^{t}\textrm{MC}(b_{3}) & 0%
\end{array}%
\right) \in M(24\times 3,\R),$

$M_{46}=\left( 
\begin{tabular}{ccc}
\cline{1-1}
\multicolumn{1}{|c}{$\mathbf{-}ME\mathbf{(}\beta _{2}\mathbf{+}\beta _{3}%
\mathbf{)}$} & \multicolumn{1}{|c}{$\mathbf{-}\textrm{MIr}(b_{3})$} & $\mathbf{-}%
\textrm{MIl}(b_{2})$ \\ \cline{1-2}
$\mathbf{-}\textrm{MIl}(b_{3})$ & \multicolumn{1}{|c}{$\mathbf{-}\textrm{ME}(\beta _{3}+\beta
_{1})$} & \multicolumn{1}{|c}{$\mathbf{-}\textrm{MIr}(b_{1})$} \\ \cline{2-3}
$\mathbf{-}\textrm{MIr}(b_{2})$ & $\mathbf{-}\textrm{MIl}(b_{1})$ & \multicolumn{1}{|c|}{$%
\mathbf{-}\textrm{ME}(\beta _{1}+\beta _{2})$} \\ \cline{3-3}
\end{tabular}%
\right) $

$\in M(24\times 24,\R),$

$M_{47}=\left( 
\begin{array}{ccc}
0 & ^{t}\textrm{MC}(a_{1}) & ^{t}\textrm{MC}(a_{1}) \\ 
^{t}\textrm{MC}(a_{2}) & 0 & ^{t}\textrm{MC}(a_{2}) \\ 
^{t}\textrm{MC}(a_{3}) & ^{t}\textrm{MC}(a_{3}) & 0%
\end{array}%
\right) \in M(24\times 3,\R),$

$M_{48}=\left( 
\begin{tabular}{ccc}
\cline{1-1}
\multicolumn{1}{|c}{$\textrm{ME}(\alpha _{2}+\alpha _{3})$} & \multicolumn{1}{|c}{$%
\textrm{MIr}(a_{3})$} & $\textrm{MIl}(a_{2})$ \\ \cline{1-2}
$\textrm{MIl}(a_{3})$ & \multicolumn{1}{|c}{$\textrm{ME}(\alpha _{3}+\alpha _{1})$} & 
\multicolumn{1}{|c}{$\textrm{MIr}(a_{1})$} \\ \cline{2-3}
$\textrm{MIr}(a_{2})$ & $\textrm{MIl}(a_{1})$ & \multicolumn{1}{|c|}{$\textrm{ME}(\alpha _{1}+\alpha
_{2})$} \\ \cline{3-3}
\end{tabular}%
\right) $

$\in M(24\times 24,\R),$

$M_{52}=\left( 
\begin{array}{ccc}
0 & \textrm{MC}(a_{2}) & -\textrm{MC}(a_{3}) \\ 
-\textrm{MC}(a_{1}) & 0 & \textrm{MC}(a_{3}) \\ 
\textrm{MC}(a_{1}) & -\textrm{MC}(a_{2}) & 0%
\end{array}%
\right) \in M(3\times 24,\R),$

$M_{53}=\left( 
\begin{array}{cc}
-\alpha _{1} & 0 \\ 
0 & -\alpha _{2} \\ 
\alpha _{3} & \alpha _{3}%
\end{array}%
\right) \in M(3\times 2,\R),$

$M_{54}=\left( 
\begin{array}{ccc}
0 & -\textrm{MC}(a_{2}) & -\textrm{MC}(a_{3}) \\ 
-\textrm{MC}(a_{1}) & 0 & -\textrm{MC}(a_{3}) \\ 
-\textrm{MC}(a_{1}) & -\textrm{MC}(a_{2}) & 0%
\end{array}%
\right) \in M(3\times 24,\R),$

$M_{55}=\left( 
\begin{array}{ccc}
\frac{2}{3}\rho +\tau _{1} & 0 & 0 \\ 
0 & \frac{2}{3}\rho +\tau _{2} & 0 \\ 
0 & 0 & \frac{2}{3}\rho +(-\tau _{11}-\tau _{12})%
\end{array}%
\right) \in M(3\times 3,\R),$

$M_{56}=\left( 
\begin{array}{@{}ccc@{}}
0 & -\textrm{MC}(m_{2}) & \textrm{MC}(m_{3}) \\ 
\textrm{MC}(m_{1}) & 0 & -\textrm{MC}(m_{3}) \\ 
-\textrm{MC}(m_{1}) & \textrm{MC}(m_{2}) & 0%
\end{array}%
\right) $

\ \ \ \ \ \ \ \ $+\left( 
\begin{array}{@{}ccc@{}}
0 & \textrm{MC}(t_{2}) & \textrm{MC}(t_{3}) \\ 
\textrm{MC}(t_{1}) & 0 & \textrm{MC}(t_{3}) \\ 
\textrm{MC}(t_{1}) & \textrm{MC}(t_{2}) & 0%
\end{array}%
\right) $

$\in M(3\times 24,C),$

$M_{59}=\left( 
\begin{array}{c}
-\frac{2}{3}\alpha _{1} \\ 
-\frac{2}{3}\alpha _{2} \\ 
-\frac{2}{3}\alpha _{3}%
\end{array}%
\right) \in M(3\times 1,\R),$

$M_{61}=\left( 
\begin{array}{c}
-^{t}\textrm{MDr}(a_{1}) \\ 
-^{t}\textrm{MDr}(a_{2})Mv \\ 
-^{t}\textrm{MDr}(a_{3})Mv^{2}%
\end{array}%
\right) \in M(24\times 28,\R),$

$M_{62}=\left( 
\begin{tabular}{ccc}
\cline{1-1}
\multicolumn{1}{|c}{$\frac{1}{2}\textrm{ME}(\alpha _{2}-\alpha _{3})$} & 
\multicolumn{1}{|c}{$-\frac{1}{2}\textrm{MIr}(a_{3})$} & $\frac{1}{2}\textrm{MIl}(a_{2})$ \\ 
\cline{1-2}
$\frac{1}{2}\textrm{MIl}(a_{3})$ & \multicolumn{1}{|c}{$\frac{1}{2}\textrm{ME}(\alpha
_{3}-\alpha _{1})$} & \multicolumn{1}{|c}{$-\frac{1}{2}\textrm{MIr}(a_{1})$} \\ 
\cline{2-3}
$-\frac{1}{2}\textrm{MIr}(a_{2})$ & $\frac{1}{2}\textrm{MIl}(a_{1})$ & \multicolumn{1}{|c|}{$%
\frac{1}{2}\textrm{ME}(\alpha _{1}-\alpha _{2})$} \\ \cline{3-3}
\end{tabular}%
\right) $

$\in M(24\times 24,\R),$

$M_{63}=\left( 
\begin{array}{cc}
\frac{1}{2}^{t}\textrm{MC}(a_{1}) & 0 \\ 
0 & \frac{1}{2}^{t}\textrm{MC}(a_{2}) \\ 
-\frac{1}{2}^{t}\textrm{MC}(a_{3}) & -\frac{1}{2}^{t}\textrm{MC}(a_{3})%
\end{array}%
\right) \in M(24\times 2,\R),$

$M_{64}=\left( 
\begin{tabular}{ccc}
\cline{1-1}
\multicolumn{1}{|c}{$-\frac{1}{2}\textrm{ME}(\alpha _{2}+\alpha _{3})$} & 
\multicolumn{1}{|c}{$-\frac{1}{2}\textrm{MIr}(a_{3})$} & $-\frac{1}{2}\textrm{MIl}(a_{2})$ \\ 
\cline{1-2}
$-\frac{1}{2}\textrm{MIl}(a_{3})$ & \multicolumn{1}{|c}{$-\frac{1}{2}\textrm{ME}(\alpha
_{3}+\alpha _{1})$} & \multicolumn{1}{|c}{$-\frac{1}{2}\textrm{MIr}(a_{1})$} \\ 
\cline{2-3}
$-\frac{1}{2}\textrm{MIr}(a_{2})$ & $-\frac{1}{2}\textrm{MIl}(a_{1})$ & \multicolumn{1}{|c|}{$-%
\frac{1}{2}\textrm{ME}(\alpha _{1}+\alpha _{2})$} \\ \cline{3-3}
\end{tabular}%
\right) $

$\in M(24\times 24,\R),$

$M_{65}=\left( 
\begin{array}{ccc}
0 & -\frac{1}{2}^{t}\textrm{MC}(m_{1}) & \frac{1}{2}^{t}\textrm{MC}(m_{1}) \\ 
\frac{1}{2}^{t}\textrm{MC}(m_{2}) & 0 & -\frac{1}{2}^{t}\textrm{MC}(m_{2}) \\ 
-\frac{1}{2}^{t}\textrm{MC}(m_{3}) & \frac{1}{2}^{t}\textrm{MC}(m_{3}) & 0%
\end{array}%
\right) $

$\ \ \ \ \ \ \ +\left( 
\begin{array}{ccc}
0 & \frac{1}{2}^{t}\textrm{MC}(t_{1}) & \frac{1}{2}^{t}\textrm{MC}(t_{1}) \\ 
\frac{1}{2}^{t}\textrm{MC}(t_{2}) & 0 & \frac{1}{2}^{t}\textrm{MC}(t_{2}) \\ 
\frac{1}{2}^{t}\textrm{MC}(t_{3}) & \frac{1}{2}^{t}\textrm{MC}(t_{3}) & 0%
\end{array}%
\right) \in M(24\times 3,\R),$

$M_{66}=M_{44}+\left( 
\begin{tabular}{ccc}
\cline{1-1}
\multicolumn{1}{|c}{$\frac{1}{2}\textrm{ME}(-\tau _{1})$} & \multicolumn{1}{|c}{$%
\frac{1}{2}\textrm{MIr}(t_{3})$} & $\frac{1}{2}\textrm{MIl}(t_{2})$ \\ \cline{1-2}
$\frac{1}{2}\textrm{MIl}(t_{3})$ & \multicolumn{1}{|c}{$\frac{1}{2}\textrm{ME}(-\tau _{2})$} & 
\multicolumn{1}{|c}{$\frac{1}{2}\textrm{MIr}(t_{1})$} \\ \cline{2-3}
$\frac{1}{2}\textrm{MIr}(t_{2})$ & $\frac{1}{2}\textrm{MIl}(t_{1})$ & \multicolumn{1}{|c|}{$%
\frac{1}{2}\textrm{ME}(\tau _{1}+\tau _{2})$} \\ \cline{3-3}
\end{tabular}%
\right) $

\ \ \ \ \ \ \ \ \ \ \ \ \ \ \ \ \ $\ +\left( 
\begin{array}{ccc}
\frac{2}{3}\textrm{ME}(\rho ) & 0 & 0 \\ 
0 & \frac{2}{3}\textrm{ME}(\rho ) & 0 \\ 
0 & 0 & \frac{2}{3}\textrm{ME}(\rho )%
\end{array}%
\right) \in M(24\times 24,\R),$

$M_{69}=\left( 
\begin{array}{c}
-\frac{2}{3}^{t}\textrm{MC}(a_{1}) \\ 
-\frac{2}{3}^{t}\textrm{MC}(a_{2}) \\ 
-\frac{2}{3}^{t}\textrm{MC}(a_{3})%
\end{array}%
\right) \in M(24\times 1,\R),$

$M_{72}=\left( 
\begin{array}{ccc}
0 & \textrm{MC}(b_{2}) & -\textrm{MC}(b_{3}) \\ 
-\textrm{MC}(b_{1}) & 0 & \textrm{MC}(b_{3}) \\ 
\textrm{MC}(b_{1}) & -\textrm{MC}(b_{2}) & 0%
\end{array}%
\right) \in M(3\times 24,\R),$

$M_{73}=\left( 
\begin{array}{cc}
\beta _{1} & 0 \\ 
0 & \beta _{2} \\ 
-\beta _{3} & -\beta _{3}%
\end{array}%
\right) \in M(3\times 2,\R),$

$M_{74}=\left( 
\begin{array}{ccc}
0 & \textrm{MC}(b_{2}) & \textrm{MC}(b_{3}) \\ 
\textrm{MC}(b_{1}) & 0 & \textrm{MC}(b_{3}) \\ 
\textrm{MC}(b_{1}) & \textrm{MC}(b_{2}) & 0%
\end{array}%
\right) \in M(3\times 24,\R),$

$M_{77}=-M_{55}\in M(3\times 3,\R),$

$M_{78}=\left( 
\begin{array}{ccc}
0 & -\textrm{MC}(m_{2}) & \textrm{MC}(m_{3}) \\ 
\textrm{MC}(m_{1}) & 0 & -\textrm{MC}(m_{3}) \\ 
-\textrm{MC}(m_{1}) & \textrm{MC}(m_{2}) & 0%
\end{array}%
\right) $

\ \ \ \ \ \ \ $+\left( 
\begin{array}{ccc}
0 & -\textrm{MC}(t_{2}) & -\textrm{MC}(t_{3}) \\ 
-\textrm{MC}(t_{1}) & 0 & -\textrm{MC}(t_{3}) \\ 
-\textrm{MC}(t_{1}) & -\textrm{MC}(t_{2}) & 0%
\end{array}%
\right) $

$\in M(3\times 24,\R),$

$M_{79}=\left( 
\begin{array}{c}
-\frac{2}{3}\beta _{1} \\ 
-\frac{2}{3}\beta _{2} \\ 
-\frac{2}{3}\beta _{3}%
\end{array}%
\right) \in M(3\times 1,\R),$

$M_{81}=\left( 
\begin{array}{c}
^{t}\textrm{MDl}(b_{1}) \\ 
^{t}\textrm{MDl}(b_{2})Mv \\ 
^{t}\textrm{MDl}(b_{3})Mv^{2}%
\end{array}%
\right) \in M(24\times 28,\R),$

$M_{82}=\left( 
\begin{tabular}{ccc}
\cline{1-1}
\multicolumn{1}{|c}{$\frac{1}{2}\textrm{ME}(\beta _{2}-\beta _{3})$} & 
\multicolumn{1}{|c}{$-\frac{1}{2}\textrm{MIr}(b_{3})$} & $\frac{1}{2}\textrm{MIl}(b_{2})$ \\ 
\cline{1-2}
$\frac{1}{2}\textrm{MIl}(b_{3})$ & \multicolumn{1}{|c}{$\frac{1}{2}\textrm{ME}(\beta
_{3}-\beta _{1})$} & \multicolumn{1}{|c}{$-\frac{1}{2}\textrm{MIr}(b_{1})$} \\ 
\cline{2-3}
$-\frac{1}{2}\textrm{MIr}(b_{2})$ & $\frac{1}{2}\textrm{MIl}(b_{1})$ & \multicolumn{1}{|c|}{$%
\frac{1}{2}\textrm{ME}(\beta _{1}-\beta _{2})$} \\ \cline{3-3}
\end{tabular}%
\right) $

$\in M(24\times 24,\R),$

$M_{83}=\left( 
\begin{array}{cc}
-\frac{1}{2}^{t}\textrm{MC}(b_{1}) & 0 \\ 
0 & -\frac{1}{2}^{t}\textrm{MC}(b_{2}) \\ 
\frac{1}{2}^{t}\textrm{MC}(b_{3}) & \frac{1}{2}^{t}\textrm{MC}(b_{3})%
\end{array}%
\right) \in M(24\times 2,\R),$

$M_{84}=\left( 
\begin{tabular}{ccc}
\cline{1-1}
\multicolumn{1}{|c}{$\frac{1}{2}\textrm{ME}(\beta _{2}+\beta _{3})$} & 
\multicolumn{1}{|c}{$\frac{1}{2}\textrm{MIr}(b_{3})$} & $\frac{1}{2}\textrm{MIl}(b_{2})$ \\ 
\cline{1-2}
$\frac{1}{2}\textrm{MIl}(b_{3})$ & \multicolumn{1}{|c}{$\frac{1}{2}\textrm{ME}(\beta
_{3}+\beta _{1})$} & \multicolumn{1}{|c}{$\frac{1}{2}\textrm{MIr}(b_{1})$} \\ 
\cline{2-3}
$\frac{1}{2}\textrm{MIr}(b_{2})$ & $\frac{1}{2}\textrm{MIl}(b_{1})$ & \multicolumn{1}{|c|}{$%
\frac{1}{2}\textrm{ME}(\beta _{1}+\beta _{2})$} \\ \cline{3-3}
\end{tabular}%
\right) $

$\in M(24\times 24,C),$

$M_{87}=\left( 
\begin{array}{ccc}
0 & -\frac{1}{2}^{t}\textrm{MC}(m_{1}) & \frac{1}{2}^{t}\textrm{MC}(m_{1}) \\ 
\frac{1}{2}^{t}\textrm{MC}(m_{2}) & 0 & -\frac{1}{2}^{t}\textrm{MC}(m_{2}) \\ 
-\frac{1}{2}^{t}\textrm{MC}(m_{3}) & \frac{1}{2}^{t}\textrm{MC}(m_{3}) & 0%
\end{array}%
\right) $

$\ \ \ \ \ \ \ -\left( 
\begin{array}{ccc}
0 & \frac{1}{2}^{t}\textrm{MC}(t_{1}) & \frac{1}{2}^{t}\textrm{MC}(t_{1}) \\ 
\frac{1}{2}^{t}\textrm{MC}(t_{2}) & 0 & \frac{1}{2}^{t}\textrm{MC}(t_{2}) \\ 
\frac{1}{2}^{t}\textrm{MC}(t_{3}) & \frac{1}{2}^{t}\textrm{MC}(t_{3}) & 0%
\end{array}%
\right) \in M(24\times 3,\R),$

$M_{88}=M_{44}-\left( 
\begin{tabular}{ccc}
\cline{1-1}
\multicolumn{1}{|c}{$\frac{1}{2}\textrm{ME}(-\tau _{1})$} & \multicolumn{1}{|c}{$%
\frac{1}{2}\textrm{MIr}(t_{3})$} & $\frac{1}{2}\textrm{MIl}(t_{2})$ \\ \cline{1-2}
$\frac{1}{2}\textrm{MIl}(t_{3})$ & \multicolumn{1}{|c}{$\frac{1}{2}\textrm{ME}(-\tau _{2})$} & 
\multicolumn{1}{|c}{$\frac{1}{2}\textrm{MIr}(t_{1})$} \\ \cline{2-3}
$\frac{1}{2}\textrm{MIr}(t_{2})$ & $\frac{1}{2}\textrm{MIl}(t_{1})$ & \multicolumn{1}{|c|}{$%
\frac{1}{2}\textrm{ME}(\tau _{1}+\tau _{2})$} \\ \cline{3-3}
\end{tabular}%
\right) $

\ \ \ \ \ \ \ \ \ $\ -\left( 
\begin{array}{ccc}
\frac{2}{3}\textrm{ME}(\rho ) & 0 & 0 \\ 
0 & \frac{2}{3}\textrm{ME}(\rho ) & 0 \\ 
0 & 0 & \frac{2}{3}\textrm{ME}(\rho )%
\end{array}%
\right) \in M(24\times 24,\R),$

$M_{89}=\left( 
\begin{array}{c}
\frac{2}{3}^{t}\textrm{MC}(b_{1}) \\ 
\frac{2}{3}^{t}\textrm{MC}(b_{2}) \\ 
\frac{2}{3}^{t}\textrm{MC}(b_{3})%
\end{array}%
\right) \in M(24\times 1,\R),$

$M95=\left( 
\begin{array}{ccc}
-\beta _{1} & -\beta _{2} & -\beta _{3}%
\end{array}%
\right) \in M(1\times 3,\R),$

$M_{96}=\left( 
\begin{array}{ccc}
-2\textrm{MC}(b_{1}) & -2\textrm{MC}(b_{2}) & -2\textrm{MC}(b_{3})%
\end{array}%
\right) \in M(1\times 24,\R),$

$M_{97}=\left( 
\begin{array}{ccc}
\alpha _{1} & \alpha _{2} & \alpha _{3}%
\end{array}%
\right) \in M(1\times 3,\R),$

$M_{98}=\left( 
\begin{array}{ccc}
2\textrm{MC}(a_{1}) & 2\textrm{MC}(a_{2}) & 2\textrm{MC}(a_{3})%
\end{array}%
\right) \in M(1\times 24,\R).$

Let we put

$M_{\phi \phi }=\left( 
\begin{array}{cccc}
M_{11} & M_{12} & 0 & M_{14} \\ 
M_{21} & M_{22} & M_{23} & M_{24} \\ 
0 & M_{32} & 0 & M_{34} \\ 
M_{41} & M_{42} & M_{43} & M_{44}%
\end{array}%
\right) ,$

$\ M_{\phi AB}=\left( 
\begin{array}{cccc}
0 & M_{16} & 0 & M_{18} \\ 
M_{25} & M_{26} & M_{27} & M_{28} \\ 
M_{35} & M_{36} & M_{37} & M_{38} \\ 
M_{45} & M_{46} & M_{47} & M_{48}%
\end{array}%
\right) ,$

$M_{AB\phi }=\left( 
\begin{array}{cccc}
0 & M_{52} & M_{53} & M_{54} \\ 
M_{61} & M_{62} & M_{63} & M_{64} \\ 
0 & M_{72} & M_{73} & M_{74} \\ 
M_{81} & M_{82} & M_{83} & M_{84}%
\end{array}%
\right) ,$ 

$\ M_{ABAB}=\left( 
\begin{array}{cccc}
M_{55} & M_{56} & 0 & 0 \\ 
M_{65} & M_{66} & 0 & 0 \\ 
0 & 0 & M_{77} & M_{78} \\ 
0 & 0 & M_{87} & M_{88}%
\end{array}%
\right) ,$

$M_{AB\rho }=\left( 
\begin{array}{c}
M_{59} \\ 
M_{69} \\ 
M_{79} \\ 
M_{89}%
\end{array}%
\right) ,$ $\ \ M_{\rho AB}=\left( 
\begin{array}{cccc}
M_{95} & M_{96} & M_{97} & M_{98}%
\end{array}%
\right) ,$

\noindent
then we have the above expression. \ \ \ \ \emph{Q.E.D.}

\bigskip

\section{The adjoint representation of \gR$_{8}^{}$}

\bigskip 

\ \ \ \emph{\ }For $R_{1}=(\Phi _{1},P_{1},Q_{1},r_{1},s_{1},u_{1}),R_{2}=(%
\Phi _{2},P_{2},Q_{2},r_{2},s_{2},u_{2})\in $\gR$_{8}^{},
$

$\Phi _{1}=(D_{1},M_{1},T_{1},A_{1},B_{1},\rho _{1}),\Phi
_{2}=(D_{2},M_{2},T_{2},A_{2},B_{2},\rho _{2})\in $\gR$_{7}^{},$

$\phi _{1}=(D_{1},M_{1},T_{1}),\phi _{2}=(D_{2},M_{2},T_{2})\in $\gR$%
_{6}^{},$

$P_{1}=(X_{1},Y_{1},\xi _{1},\eta _{1}),Q_{1}=(Z_{1},W_{1},\zeta
_{1},\omega _{1}),$

$P_{2}=(X_{2},Y_{2},\xi _{2},\eta
_{2}),Q_{2}=(Z_{2},W_{2},\zeta _{2},\omega _{2})\in $\gP$^{},$

$%
M_{1}=A_{1}(m_{11})+A_{2}(m_{12})+A_{3}(m_{13}),M_{2}=A_{1}(m_{21})+A_{2}(m_{22})+A_{3}(m_{23}), 
$

$\ m_{11},m_{12},m_{13},m_{21},m_{22},m_{23}\in $\gC $^{},$

$T_{1}=\tau _{11}E_{1}+\tau _{12}E_{2}+(-\tau _{11}-\tau
_{12})E_{3}+F_{1}(t_{11})+F_{2}(t_{12})+F_{3}(t_{13}),$

$T_{2}=\tau _{21}E_{1}+\tau _{22}E_{2}+(-\tau _{21}-\tau
_{22})E_{3}+F_{1}(t_{21})+F_{2}(t_{22})+F_{3}(t_{23}),$

$t_{11},t_{12},t_{13},t_{21},t_{22},t_{23}\in $\gC $^{},\tau _{11},\tau
_{12},\tau _{21},\tau _{22}\in \R,$

$A_{1}=\alpha _{11}E_{1}+\alpha _{12}E_{2}+\alpha
_{13}E_{3}+F_{1}(a_{11})+F_{2}(a_{12})+F_{3}(a_{13}),$

$A_{2}=\alpha _{21}E_{1}+\alpha _{22}E_{2}+\alpha
_{23}E_{3}+F_{1}(a_{21})+F_{2}(a_{22})+F_{3}(a_{23}),$

$B_{1}=\beta _{11}E_{1}+\beta _{12}E_{2}+\beta
_{13}E_{3}+F_{1}(b_{11})+F_{2}(b_{12})+F_{3}(b_{13}),$

$B_{2}=\beta _{21}E_{1}+\beta _{22}E_{2}+\beta
_{23}E_{3}+F_{1}(b_{21})+F_{2}(b_{22})+F_{3}(b_{23}),$

$a_{11},a_{12},a_{13},a_{21},a_{22},a_{23}\in $\gC $^{},\alpha
_{11},\alpha _{12},\alpha _{13},\alpha _{21},\alpha _{22},\alpha _{23}\in \R,$

$b_{11},b_{12},b_{13},b_{21},b_{22},b_{23}\in $\gC $^{},\beta
_{11},\beta _{12},\beta _{13},\beta _{21},\beta _{22},\beta _{23}\in \R,$

$\rho _{1},\rho _{2}\in \R,$

$X_{1}=\chi _{11}E_{1}+\chi _{12}E_{2}+\chi
_{13}E_{3}+F_{1}(x_{11})+F_{2}(x_{12})+F_{3}(x_{13}),$

$X_{2}=\chi _{21}E_{1}+\chi _{22}E_{2}+\chi
_{23}E_{3}+F_{1}(x_{21})+F_{2}(x_{22})+F_{3}(x_{23}),$

$Y_{1}=\gamma _{11}E_{1}+\gamma _{12}E_{2}+\gamma
_{13}E_{3}+F_{1}(y_{11})+F_{2}(y_{12})+F_{3}(y_{13}),$

$Y_{2}=\gamma _{21}E_{1}+\gamma _{22}E_{2}+\gamma
_{23}E_{3}+F_{1}(y_{21})+F_{2}(y_{22})+F_{3}(y_{23}),$

$Z_{1}=\mu _{11}E_{1}+\mu _{12}E_{2}+\mu
_{13}E_{3}+F_{1}(z_{11})+F_{2}(z_{12})+F_{3}(z_{13}),$

$Z_{2}=\mu _{21}E_{1}+\mu _{22}E_{2}+\mu
_{23}E_{3}+F_{1}(z_{21})+F_{2}(z_{22})+F_{3}(z_{23}),$

$W_{1}=\psi _{11}E_{1}+\psi _{12}E_{2}+\psi
_{13}E_{3}+F_{1}(w_{11})+F_{2}(w_{12})+F_{3}(w_{13}),$

$W_{2}=\psi _{21}E_{1}+\psi _{22}E_{2}+\psi
_{23}E_{3}+F_{1}(w_{21})+F_{2}(w_{22})+F_{3}(w_{23}),$

$x_{11},x_{12},x_{13},x_{21},x_{22},x_{23}\in $\gC $^{},\chi _{11},\chi
_{12},\chi _{13},\chi _{21},\chi _{22},\chi _{23}\in \R,$

$y_{11},y_{12},y_{13},y_{21},y_{22},y_{23}\in $\gC $^{},\gamma
_{11},\gamma _{12},\gamma _{13},\gamma _{21},\gamma _{22},\gamma _{23}\in \R,$

$z_{11},z_{12},z_{13},z_{21},z_{22},z_{23}\in $\gC $^{},\mu _{11},\mu
_{12},\mu _{13},\mu _{21},\mu _{22},\mu _{23}\in \R,$

$w_{11},w_{12},w_{13},w_{21},w_{22},w_{23}\in $\gC $^{},\psi _{11},\psi
_{12},\psi _{13},\psi _{21},\psi _{22},\psi _{23}\in \R,$

$\xi _{1},\eta _{1},\zeta _{1},\omega _{1},\xi _{2},\eta
_{2},\zeta _{2},\omega _{2}\in \R,r_{1},r_{2},s_{1},s_{2},u_{1},u_{2}\in
\R,$

\noindent
Let Lie bracket $[R_{1},R_{2}]_{8}$ be the following expression.

\ \ $\left( 
\begin{array}{cccccc}
F_{\Phi \Phi } & F_{\Phi P} & F_{\Phi Q} & 0 & 0 & 0 \\ 
F_{P\Phi } & F_{PP} & F_{PQ} & F_{P\text{\textit{r}}} & F_{Ps} & 0 \\ 
F_{Q\Phi } & F_{QP} & F_{QQ} & F_{Qr} & 0 & F_{Qu} \\ 
0 & F_{rP} & F_{rQ} & 0 & F_{rs} & F_{ru} \\ 
0 & F_{sP} & 0 & F_{sr} & F_{ss} & 0 \\ 
0 & 0 & F_{uQ} & F_{ur} & 0 & F_{uu}%
\end{array}%
\right) \left( 
\begin{array}{c}
\Phi _{2} \\ 
P_{2} \\ 
Q_{2} \\ 
r_{2} \\ 
s_{2} \\ 
u_{2}%
\end{array}%
\right) =\left( 
\begin{array}{c}
\lbrack R_{1},R_{2}]_{8\Phi } \\ 
\lbrack R_{1},R_{2}]_{8P} \\ 
\lbrack R_{1},R_{2}]_{8Q} \\ 
\lbrack R_{1},R_{2}]_{8r} \\ 
\lbrack R_{1},R_{2}]_{8s} \\ 
\lbrack R_{1},R_{2}]_{8u}%
\end{array}%
\right) .$

\bigskip

\emph{Lemma 10.1. }

$\left( 
\begin{array}{cccccc}
F_{\Phi \Phi } & F_{\Phi P} & F_{\Phi Q} & 0 & 0 & 0%
\end{array}%
\right) \left( 
\begin{array}{c}
\Phi _{2} \\ 
P_{2} \\ 
Q_{2} \\ 
r_{2} \\ 
s_{2} \\ 
u_{2}%
\end{array}%
\right) =\left( [R_{1},R_{2}]_{8\Phi }\right) ,$

\noindent
where

$F_{\Phi \Phi }=\left( 
\begin{array}{ccccccccc}
F_{11} & F_{12} & 0 & F_{14} & 0 & F_{16} & 0 & F_{18} & 0 \\ 
F_{21} & F_{22} & F_{23} & F_{24} & F_{25} & F_{26} & F_{27} & F_{28} & 0 \\ 
0 & F_{32} & 0 & F_{34} & F_{35} & F_{36} & F_{37} & F_{38} & 0 \\ 
F_{41} & F_{42} & F_{43} & F_{44} & F_{45} & F_{46} & F_{47} & F_{48} & 0 \\ 
0 & F_{52} & F_{53} & F_{54} & F_{55} & F_{56} & 0 & 0 & F_{59} \\ 
F_{61} & F_{62} & F_{63} & F_{64} & F_{65} & F_{66} & 0 & 0 & F_{69} \\ 
0 & F_{72} & F_{73} & F_{74} & 0 & 0 & F_{77} & F_{78} & F_{79} \\ 
F_{81} & F_{82} & F_{83} & F_{84} & 0 & 0 & F_{87} & F_{88} & F_{89} \\ 
0 & 0 & 0 & 0 & F_{95} & F_{96} & F_{97} & F_{98} & 0%
\end{array}%
\right) ,$

\bigskip

$F_{\Phi P}=\left( 
\begin{array}{cccccc}
0 & F_{\Phi P12} & 0 & F_{\Phi P14} & 0 & 0 \\ 
F_{\Phi P21} & F_{\Phi P22} & F_{\Phi P23} & F_{\Phi P24} & 0 & 0
\\ 
F_{\Phi P31} & F_{\Phi P32} & F_{\Phi P33} & F_{\Phi P34} & 0 & 0
\\ 
F_{\Phi P41} & F_{\Phi P42} & F_{\Phi P43} & F_{\Phi P44} & 0 & 0
\\ 
F_{\Phi P51} & 0 & F_{\Phi P53} & F_{\Phi P54} & F_{\Phi P55} & 0
\\ 
0 & F_{\Phi P62} & F_{\Phi P63} & F_{\Phi P64} & F_{\Phi P65} & 0
\\ 
F_{\Phi P71} & F_{\Phi P72} & F_{\Phi P73} & 0 & 0 & F_{\Phi P76}
\\ 
F_{\Phi P81} & F_{\Phi P82} & 0 & F_{\Phi P84} & 0 & F_{\Phi P86}
\\ 
F_{\Phi P91} & F_{\Phi P92} & F_{\Phi P93} & F_{\Phi P94} & F_{%
\Phi P95} & F_{\Phi P96}%
\end{array}%
\right) ,$

$F_{\Phi P12}=\left( 
\begin{array}{ccc}
\frac{1}{2}\textrm{f}_{JJD}(?,w_{11}) & \frac{1}{2}\textrm{d}_{g}\nu ^{2}\textrm{g}_{d}\textrm{f}_{JJD}(?,w_{12})
& \frac{1}{2}\textrm{d}_{g}\nu \textrm{g}_{d}\textrm{f}_{JJD}(?,w_{13})%
\end{array}%
\right) ,$

$F_{\Phi P14}=\left( 
\begin{array}{ccc}
\frac{1}{2}\textrm{f}_{JJD}(z_{11},?) & \frac{1}{2}\textrm{d}_{g}\nu ^{2}\textrm{g}_{d}\textrm{f}_{JJD}(z_{12},?)
& \frac{1}{2}\textrm{d}_{g}\nu \textrm{g}_{d}\textrm{f}_{JJD}(z_{13},?)%
\end{array}%
\right) ,$

$F_{\Phi P21}=\left( 
\begin{array}{ccc}
0 & \frac{1}{4}\textrm{f}_{CJJ}(?,w_{11}) & -\frac{1}{4}\textrm{f}_{CJJ}(?,w_{11}) \\ 
-\frac{1}{4}\textrm{f}_{CJJ}(?,w_{12}) & 0 & \frac{1}{4}\textrm{f}_{CJJ}(?,w_{12}) \\ 
\frac{1}{4}\textrm{f}_{CJJ}(?,w_{13}) & -\frac{1}{4}\textrm{f}_{CJJ}(?,w_{13}) & 0%
\end{array}%
\right) ,$

$F_{\Phi P22}=$

$\left( 
\begin{array}{@{}ccc@{}}
-\frac{1}{4}\textrm{f}_{CJJ}(\psi _{12}-\psi _{13},?) & -%
\frac{1}{4}\textrm{f}_{iJJ}(?,w_{13}) & \frac{1}{4}\textrm{f}_{iJJ}(w_{12},?) \\ 
\frac{1}{4}\textrm{f}_{iJJ}(w_{13},?) & -\frac{1}{4}\textrm{f}_{CJJ}(\psi _{13}-\psi _{11},?) & -\frac{1}{4}\textrm{f}_{iJJ}(?,w_{11}) \\ 
-\frac{1}{4}\textrm{f}_{iJJ}(?,w_{12}) & \frac{1}{4}\textrm{f}_{iJJ}(w_{11},?) & -\frac{1}{4}%
\textrm{f}_{CJJ}(\psi _{11}-\psi _{12},?)%
\end{array}%
\right) $

$F_{\Phi P23}=\left( 
\begin{array}{ccc}
0 & -\frac{1}{4}\textrm{f}_{CJJ}(?,z_{11}) & \frac{1}{4}\textrm{f}_{CJJ}(?,z_{11}) \\ 
\frac{1}{4}\textrm{f}_{CJJ}(?,z_{12}) & 0 & -\frac{1}{4}\textrm{f}_{CJJ}(?,z_{12}) \\ 
-\frac{1}{4}\textrm{f}_{CJJ}(?,z_{13}) & \frac{1}{4}\textrm{f}_{CJJ}(?,z_{13}) & 0%
\end{array}%
\right) ,$

$F_{\Phi P24}=$

$\left( 
\begin{array}{@{}ccc@{}}
\frac{1}{4}\textrm{f}_{CJJ}(\mu _{12}-\mu _{13},?) & 
\frac{1}{4}\textrm{f}_{iJJ}(?,z_{13}) & -\frac{1}{4}\textrm{f}_{iJJ}(z_{12},?) \\ 
-\frac{1}{4}\textrm{f}_{iJJ}(z_{13},?) & \frac{1}{4}\textrm{f}_{CJJ}(\mu _{13}-\mu _{11},?) & \frac{1}{4}\textrm{f}_{iJJ}(?,z_{11}) \\ 
\frac{1}{4}\textrm{f}_{iJJ}(?,z_{12}) & -\frac{1}{4}\textrm{f}_{iJJ}(z_{11},?) & \frac{1}{4}%
\textrm{f}_{CJJ}(\mu _{11}-\mu _{12},?)%
\end{array}%
\right) ,$

$F_{\Phi P31}=\left( 
\begin{array}{ccc}
\frac{1}{3}\textrm{f}_{CCC}(\psi _{11},?) & -\frac{1}{6}\textrm{f}_{CCC}(\psi _{12},?) & -%
\frac{1}{6}\textrm{f}_{CCC}(\psi _{13},?) \\ 
-\frac{1}{6}\textrm{f}_{CCC}(\psi _{11},?) & \frac{1}{3}\textrm{f}_{CCC}(\psi _{12},?) & -%
\frac{1}{6}\textrm{f}_{CCC}(\psi _{13},?) \\ 
-\frac{1}{6}\textrm{f}_{CCC}(\psi _{11},?) & -\frac{1}{6}\textrm{f}_{CCC}(\psi _{12},?) & 
\frac{1}{3}\textrm{f}_{CCC}(\psi _{13},?)%
\end{array}%
\right) ,$

$F_{\Phi P32}=\left( 
\begin{array}{ccc}
-\frac{1}{3}\textrm{f}_{JJC}(?,w_{11}) & \frac{1}{6}\textrm{f}_{JJC}(?,w_{12}) & \frac{1}{6}%
\textrm{f}_{JJC}(?,w_{13}) \\ 
\frac{1}{6}\textrm{f}_{JJC}(?,w_{11}) & -\frac{1}{3}\textrm{f}_{JJC}(?,w_{12}) & \frac{1}{6}%
\textrm{f}_{JJC}(?,w_{13}) \\ 
\frac{1}{6}\textrm{f}_{JJC}(?,w_{11}) & \frac{1}{6}\textrm{f}_{JJC}(?,w_{12}) & -\frac{1}{3}%
\textrm{f}_{JJC}(?,w_{13})%
\end{array}%
\right) ,$

$F_{\Phi P33}=\left( 
\begin{array}{ccc}
\frac{1}{3}\textrm{f}_{CCC}(\mu _{11},?) & -\frac{1}{6}\textrm{f}_{CCC}(\mu _{12},?) & -\frac{1%
}{6}\textrm{f}_{CCC}(\mu _{13},?) \\ 
-\frac{1}{6}\textrm{f}_{CCC}(\mu _{11},?) & \frac{1}{3}\textrm{f}_{CCC}(\mu _{12},?) & -\frac{1%
}{6}\textrm{f}_{CCC}(\mu _{13},?) \\ 
-\frac{1}{6}\textrm{f}_{CCC}(\mu _{11},?) & -\frac{1}{6}\textrm{f}_{CCC}(\mu _{12},?) & \frac{1%
}{3}\textrm{f}_{CCC}(\mu _{13},?)%
\end{array}%
\right) ,$

$F_{\Phi P34}=\left( 
\begin{array}{ccc}
-\frac{1}{3}\textrm{f}_{JJC}(?,z_{11}) & \frac{1}{6}\textrm{f}_{JJC}(?,z_{12}) & \frac{1}{6}%
\textrm{f}_{JJC}(?,z_{13}) \\ 
\frac{1}{6}\textrm{f}_{JJC}(?,z_{11}) & -\frac{1}{3}\textrm{f}_{JJC}(?,z_{12}) & \frac{1}{6}%
\textrm{f}_{JJC}(?,z_{13}) \\ 
\frac{1}{6}\textrm{f}_{JJC}(?,z_{11}) & \frac{1}{6}\textrm{f}_{JJC}(?,z_{12}) & -\frac{1}{3}%
\textrm{f}_{JJC}(?,z_{13})%
\end{array}%
\right) ,$

$F_{\Phi P41}=\left( 
\begin{array}{ccc}
0 & \frac{1}{4}\textrm{f}_{CJJ}(?,w_{11}) & \frac{1}{4}\textrm{f}_{CJJ}(?,w_{11}) \\ 
\frac{1}{4}\textrm{f}_{CJJ}(?,w_{12}) & 0 & \frac{1}{4}\textrm{f}_{CJJ}(?,w_{12}) \\ 
\frac{1}{4}\textrm{f}_{CJJ}(?,w_{13}) & \frac{1}{4}\textrm{f}_{CJJ}(?,w_{13}) & 0%
\end{array}%
\right) ,$

$F_{\Phi P42}=$

$\left( 
\begin{array}{@{}ccc@{}}
\frac{1}{4}\textrm{f}_{CJJ}(\psi _{12}+\psi _{13},?) & 
\frac{1}{4}\textrm{f}_{iJJ}(?,w_{13}) & \frac{1}{4}\textrm{f}_{iJJ}(w_{12},?) \\ 
\frac{1}{4}\textrm{f}_{iJJ}(w_{13},?) & \frac{1}{4}\textrm{f}_{CJJ}(\psi _{13}+\psi _{11},?) & \frac{1}{4}\textrm{f}_{iJJ}(?,w_{11}) \\ 
\frac{1}{4}\textrm{f}_{iJJ}(?,w_{12}) & \frac{1}{4}\textrm{f}_{iJJ}(w_{11},?) & \frac{1}{4}%
\textrm{f}_{CJJ}(\psi _{11}+\psi _{12},?)%
\end{array}%
\right) ,$

$F_{\Phi P43}=\left( 
\begin{array}{ccc}
0 & \frac{1}{4}\textrm{f}_{CJJ}(?,z_{11}) & \frac{1}{4}\textrm{f}_{CJJ}(?,z_{11}) \\ 
\frac{1}{4}\textrm{f}_{CJJ}(?,z_{12}) & 0 & \frac{1}{4}\textrm{f}_{CJJ}(?,z_{12}) \\ 
\frac{1}{4}\textrm{f}_{CJJ}(?,z_{13}) & \frac{1}{4}\textrm{f}_{CJJ}(?,z_{13}) & 0%
\end{array}%
\right) ,$

$F_{\Phi P44}=$

$\left( 
\begin{array}{@{}ccc@{}}
\frac{1}{4}\textrm{f}_{CJJ}(\mu _{12}+\mu _{13},?) & 
\frac{1}{4}\textrm{f}_{iJJ}(?,z_{13}) & \frac{1}{4}\textrm{f}_{iJJ}(z_{12},?) \\ 
\frac{1}{4}\textrm{f}_{iJJ}(z_{13},?) & \frac{1}{4}\textrm{f}_{CJJ}(\mu _{13}+\mu _{11},?) & \frac{1}{4}\textrm{f}_{iJJ}(?,z_{11}) \\ 
\frac{1}{4}\textrm{f}_{iJJ}(?,z_{12}) & \frac{1}{4}\textrm{f}_{iJJ}(z_{11},?) & \frac{1}{4}%
\textrm{f}_{CJJ}(\mu _{11}+\mu _{12},?)%
\end{array}%
\right) ,$

$F_{\Phi P51}=\left( 
\begin{array}{ccc}
-\textrm{f}_{CCC}(\zeta _{1},?) & 0 & 0 \\ 
0 & -\textrm{f}_{CCC}(\zeta _{1},?) & 0 \\ 
0 & 0 & -\textrm{f}_{CCC}(\zeta _{1},?)%
\end{array}%
\right) ,$

$F_{\Phi P53}=\left( 
\begin{array}{ccc}
0 & \frac{1}{4}\textrm{f}_{CCC}(\psi _{13},?) & \frac{1}{4}\textrm{f}_{CCC}(\psi _{12},?) \\ 
\frac{1}{4}\textrm{f}_{CCC}(\psi _{13},?) & 0 & \frac{1}{4}\textrm{f}_{CCC}(\psi _{11},?) \\ 
\frac{1}{4}\textrm{f}_{CCC}(\psi _{12},?) & \frac{1}{4}\textrm{f}_{CCC}(\psi _{11},?) & 0%
\end{array}%
\right) ,$

$F_{\Phi P54}=\left( 
\begin{array}{ccc}
-\frac{1}{2}\textrm{f}_{JJC}(?,w_{11}) & 0 & 0 \\ 
0 & -\frac{1}{2}\textrm{f}_{JJC}(?,w_{12}) & 0 \\ 
0 & 0 & -\frac{1}{2}\textrm{f}_{JJC}(?,w_{13})%
\end{array}%
\right) ,$

$F_{\Phi P55}=\left( 
\begin{array}{c}
-\frac{1}{4}\textrm{f}_{CCC}(\mu _{11},?) \\ 
-\frac{1}{4}\textrm{f}_{CCC}(\mu _{12},?) \\ 
-\frac{1}{4}\textrm{f}_{CCC}(\mu _{13},?)%
\end{array}%
\right) ,$

$F_{\Phi P62}=\left( 
\begin{array}{ccc}
-\frac{1}{4}\textrm{f}_{CJJ}(\zeta _{1},?) & 0 & 0 \\ 
0 & -\frac{1}{4}\textrm{f}_{CJJ}(\zeta _{1},?) & 0 \\ 
0 & 0 & -\frac{1}{4}\textrm{f}_{CJJ}(\zeta _{1},?)%
\end{array}%
\right) ,$

$F_{\Phi P63}=\left( 
\begin{array}{ccc}
-\frac{1}{4}\textrm{f}_{CJJ}(?,w_{11}) & 0 & 0 \\ 
0 & -\frac{1}{4}\textrm{f}_{CJJ}(?,w_{12}) & 0 \\ 
0 & 0 & -\frac{1}{4}\textrm{f}_{CJJ}(?,w_{13})%
\end{array}%
\right) ,$

$F_{\Phi P64}=\left( 
\begin{array}{ccc}
-\frac{1}{4}\textrm{f}_{CJJ}(\psi _{11},?) & \frac{1}{4}\textrm{f}_{iJJ}(?,w_{13}) & \frac{1}{4%
}\textrm{f}_{iJJ}(w_{12},?) \\ 
\frac{1}{4}\textrm{f}_{iJJ}(w_{13},?) & -\frac{1}{4}\textrm{f}_{CJJ}(\psi _{12},?) & \frac{1}{4%
}\textrm{f}_{iJJ}(?,w_{11}) \\ 
\frac{1}{4}\textrm{f}_{iJJ}(?,w_{12}) & \frac{1}{4}\textrm{f}_{iJJ}(w_{11},?) & -\frac{1}{4}%
\textrm{f}_{CJJ}(\psi _{13},?)%
\end{array}%
\right) ,$

$F_{\Phi P65}=\left( 
\begin{array}{c}
-\frac{1}{4}\textrm{f}_{CJJ}(?,z_{11}) \\ 
-\frac{1}{4}\textrm{f}_{CJJ}(?,z_{12}) \\ 
-\frac{1}{4}\textrm{f}_{CJJ}(?,z_{13})%
\end{array}%
\right) ,$

$F_{\Phi P71}=\left( 
\begin{array}{ccc}
0 & -\frac{1}{4}\textrm{f}_{CCC}(\mu _{13},?) & -\frac{1}{4}\textrm{f}_{CCC}(\mu _{12},?) \\ 
-\frac{1}{4}\textrm{f}_{CCC}(\mu _{13},?) & 0 & -\frac{1}{4}\textrm{f}_{CCC}(\mu _{11},?) \\ 
-\frac{1}{4}\textrm{f}_{CCC}(\mu _{12},?) & -\frac{1}{4}\textrm{f}_{CCC}(\mu _{11},?) & 0%
\end{array}%
\right) ,$

$F_{\Phi P72}=\left( 
\begin{array}{ccc}
\frac{1}{2}\textrm{f}_{JJC}(?,z_{11}) & 0 & 0 \\ 
0 & \frac{1}{2}\textrm{f}_{JJC}(?,z_{12}) & 0 \\ 
0 & 0 & \frac{1}{2}\textrm{f}_{JJC}(?,z_{13})%
\end{array}%
\right) ,$

$F_{\Phi P73}=\left( 
\begin{array}{ccc}
\frac{1}{4}\textrm{f}_{CCC}(\omega _{1},?) & 0 & 0 \\ 
0 & \frac{1}{4}\textrm{f}_{CCC}(\omega _{1},?) & 0 \\ 
0 & 0 & \frac{1}{4}\textrm{f}_{CCC}(\omega _{1},?)%
\end{array}%
\right) ,$

$F_{\Phi P76}=\left( 
\begin{array}{c}
\frac{1}{4}\textrm{f}_{CCC}(\psi _{11},?) \\ 
\frac{1}{4}\textrm{f}_{CCC}(\psi _{12},?) \\ 
\frac{1}{4}\textrm{f}_{CCC}(\psi _{13},?)%
\end{array}%
\right) ,$

$F_{\Phi P81}=\left( 
\begin{array}{ccc}
\frac{1}{4}\textrm{f}_{CJJ}(?,z_{11}) & 0 & 0 \\ 
0 & \frac{1}{4}\textrm{f}_{CJJ}(?,z_{12}) & 0 \\ 
0 & 0 & \frac{1}{4}\textrm{f}_{CJJ}(?,z_{13})%
\end{array}%
\right) ,$

$F_{\Phi P82}=\left( 
\begin{array}{ccc}
\frac{1}{4}\textrm{f}_{CJJ}(\mu _{11},?) & -\frac{1}{4}\textrm{f}_{iJJ}(?,z_{13}) & -\frac{1}{4%
}\textrm{f}_{iJJ}(z_{12},?) \\ 
-\frac{1}{4}\textrm{f}_{iJJ}(z_{13},?) & \frac{1}{4}\textrm{f}_{CJJ}(\mu _{12},?) & -\frac{1}{4%
}\textrm{f}_{iJJ}(?,z_{11}) \\ 
-\frac{1}{4}\textrm{f}_{iJJ}(?,z_{12}) & -\frac{1}{4}\textrm{f}_{iJJ}(z_{11},?) & \frac{1}{4}%
\textrm{f}_{CJJ}(\mu _{13},?)%
\end{array}%
%
\right) ,$

$([R_{1},R_{2}]_{8\Phi })_{D}=[\Phi _{1},\Phi
_{2}]_{7D}+([R_{1},R_{2}]_{8\Phi })_{oD}$

$([R_{1},R_{2}]_{8\Phi })_{oD}=$

\ \ \ \ \ \ \ \ \ \ \ \ $-\frac{1}{2}(\ \textrm{JD}(x_{11},w_{21})+\textrm{d}_{g}\nu
^{2}\textrm{g}_{d}(\textrm{JD}(x_{12},w_{22}))+\textrm{d}_{g}\nu \textrm{g}_{d}(\textrm{JD}(x_{13},w_{23}))$

$\ \ \ \ \ \ \ \ \ \ \ \ \ \ \ \ 
+\textrm{JD}(z_{21},y_{11})+\textrm{d}_{g}\nu ^{2}\textrm{g}_{d}(\textrm{JD}(z_{22},y_{12}))+\textrm{d}_{g}\nu
\textrm{g}_{d}(\textrm{JD}(z_{23},y_{13})))$

$\ \ \ \ \ \ \ \ \ \ \ \ \ +\frac{1}{2}%
(\textrm{JD}(x_{21},w_{11})+\textrm{d}_{g}\nu ^{2}\textrm{g}_{d}(\textrm{JD}(x_{22},w_{12}))+\textrm{d}_{g}\nu
\textrm{g}_{d}(\textrm{JD}(x_{23},w_{13}))$

$\ \ \ \ \ \ \ \ \ \ \ \ \ \ \ \ 
+\textrm{JD}(z_{11},y_{21})+\textrm{d}_{g}\nu ^{2}\textrm{g}_{d}(\textrm{JD}(z_{12},y_{22}))+\textrm{d}_{g}\nu
\textrm{g}_{d}(\textrm{JD}(z_{13},y_{23}))),$

$([R_{1},R_{2}]_{8\Phi })_{M}=\left( 
\begin{array}{c}
([R_{1},R_{2}]_{8\Phi })_{M\_A_{1}} \\ 
([R_{1},R_{2}]_{8\Phi })_{M\_A_{2}} \\ 
([R_{1},R_{2}]_{8\Phi })_{M\_A_{3}}%
\end{array}%
\right) $

$\ \ \ \ \ \ \ \ \ \ \ \ \ \ \ \ \ \ \ =\left( 
\begin{array}{c}
\lbrack \Phi _{1},\Phi _{2}]_{7M\_A_{1}} \\ 
\lbrack \Phi _{1},\Phi _{2}]_{7M\_A_{2}} \\ 
\lbrack \Phi _{1},\Phi _{2}]_{7M\_A_{3}}%
\end{array}%
\right) +\left( 
\begin{array}{c}
([R_{1},R_{2}]_{8\Phi })_{oM\_A_{1}} \\ 
([R_{1},R_{2}]_{8\Phi })_{oM\_A_{2}} \\ 
([R_{1},R_{2}]_{8\Phi })_{oM\_A_{3}}%
\end{array}%
\right) ,$

$([R_{1},R_{2}]_{8\Phi })_{oM\_A_{1}}=-\frac{1}{4}((\chi _{12}-\chi
_{13})w_{21}-(\psi _{22}-\psi _{23})x_{11}-\overline{x_{12}w_{23}}+\overline{%
w_{22}x_{13}}$

$\ \ \ \ \ \ \ \ \ \ \ \ \ \ \ \ \ \ \ \ \ \ \ \ \ \ \ \ \ \ \ \ \ +(\mu
_{22}-\mu _{23})y_{11}-(\gamma _{12}-\gamma _{13})z_{21}-\overline{%
z_{22}y_{13}}+\overline{y_{12}z_{23}})$

$\ \ \ \ \ \ \ \ \ \ \ \ \ \ \ \ \ \ \ \ \ \ \ \ \ \ \ \ +\frac{1}{4}((\chi
_{22}-\chi _{23})w_{11}-(\psi _{12}-\psi _{13})x_{21}-\overline{x_{22}w_{13}}%
+\overline{w_{12}x_{23}}$

$\ \ \ \ \ \ \ \ \ \ \ \ \ \ \ \ \ \ \ \ \ \ \ \ \ \ \ \ \ \ \ \ \ \ +(\mu
_{12}-\mu _{13})y_{21}-(\gamma _{22}-\gamma _{23})z_{11}-\overline{%
z_{12}y_{23}}+\overline{y_{22}z_{13}}),$

$([R_{1},R_{2}]_{8\Phi })_{oM\_A_{2}}=-\frac{1}{4}((\chi _{13}-\chi
_{11})w_{22}-(\psi _{23}-\psi _{21})x_{12}-\overline{x_{13}w_{21}}+\overline{%
w_{23}x_{11}}$

$\ \ \ \ \ \ \ \ \ \ \ \ \ \ \ \ \ \ \ \ \ \ \ \ \ \ \ \ \ \ \ \ \ +(\mu
_{23}-\mu _{21})y_{12}-(\gamma _{13}-\gamma _{11})z_{22}-\overline{%
z_{23}y_{11}}+\overline{y_{13}z_{21}})$

$\ \ \ \ \ \ \ \ \ \ \ \ \ \ \ \ \ \ \ \ \ \ \ \ \ \ \ \ +\frac{1}{4}((\chi
_{23}-\chi _{21})w_{12}-(\psi _{13}-\psi _{11})x_{22}-\overline{x_{23}w_{11}}%
+\overline{w_{13}x_{21}}$

$\ \ \ \ \ \ \ \ \ \ \ \ \ \ \ \ \ \ \ \ \ \ \ \ \ \ \ \ \ \ \ \ \ \ +(\mu
_{13}-\mu _{11})y_{22}-(\gamma _{23}-\gamma _{21})z_{12}-\overline{%
z_{13}y_{21}}+\overline{y_{23}z_{11}}),$

$([R_{1},R_{2}]_{8\Phi })_{oM\_A_{3}}=-\frac{1}{4}((\chi _{11}-\chi
_{12})w_{23}-(\psi _{21}-\psi _{22})x_{13}-\overline{x_{11}w_{22}}+\overline{%
w_{21}x_{12}}$

$\ \ \ \ \ \ \ \ \ \ \ \ \ \ \ \ \ \ \ \ \ \ \ \ \ \ \ \ \ \ \ \ \ \ +(\mu
_{21}-\mu _{22})y_{13}-(\gamma _{11}-\gamma _{12})z_{23}-\overline{%
z_{21}y_{12}}+\overline{y_{11}z_{22}})$

$\ \ \ \ \ \ \ \ \ \ \ \ \ \ \ \ \ \ \ \ \ \ \ \ \ \ \ \ +\frac{1}{4}((\chi
_{21}-\chi _{22})w_{13}-(\psi _{11}-\psi _{12})x_{23}-\overline{x_{21}w_{12}}%
+\overline{w_{11}x_{22}}$

$\ \ \ \ \ \ \ \ \ \ \ \ \ \ \ \ \ \ \ \ \ \ \ \ \ \ \ \ \ \ \ \ \ \ \ +(\mu
_{11}-\mu _{12})y_{23}-(\gamma _{21}-\gamma _{22})z_{13}-\overline{%
z_{11}y_{22}}+\overline{y_{21}z_{12}}),$

$([R_{1},R_{2}]_{8\Phi })_{T}=\left( 
\begin{array}{c}
([R_{1},R_{2}]_{8\Phi })_{T\_E_{1}} \\ 
([R_{1},R_{2}]_{8\Phi })_{T\_E_{2}} \\ 
([R_{1},R_{2}]_{8\Phi })_{T\_E_{3}} \\ 
([R_{1},R_{2}]_{8\Phi })_{T\_F_{1}} \\ 
([R_{1},R_{2}]_{8\Phi })_{T\_F_{2}} \\ 
([R_{1},R_{2}]_{8\Phi })_{T\_F_{3}}%
\end{array}%
\right) $

$\ \ \ \ \ \ \ \ \ \ \ \ \ \ \ \ \ \ \ =\left( 
\begin{array}{c}
\lbrack \Phi _{1},\Phi _{2}]_{7T\_E_{1}} \\ 
\lbrack \Phi _{1},\Phi _{2}]_{7T\_E_{2}} \\ 
\lbrack \Phi _{1},\Phi _{2}]_{7T\_E_{3}} \\ 
\lbrack \Phi _{1},\Phi _{2}]_{7T\_F_{1}} \\ 
\lbrack \Phi _{1},\Phi _{2}]_{7T\_F_{2}} \\ 
\lbrack \Phi _{1},\Phi _{2}]_{7T\_F_{3}}%
\end{array}%
\right) +\left( 
\begin{array}{c}
([R_{1},R_{2}]_{8\Phi })_{oT\_E_{1}} \\ 
([R_{1},R_{2}]_{8\Phi })_{oT\_E_{2}} \\ 
([R_{1},R_{2}]_{8\Phi })_{oT\_E_{3}} \\ 
([R_{1},R_{2}]_{8\Phi })_{oT\_F_{1}} \\ 
([R_{1},R_{2}]_{8\Phi })_{oT\_F_{2}} \\ 
([R_{1},R_{2}]_{8\Phi })_{oT\_F_{3}}%
\end{array}%
\right) ,$

$([R_{1},R_{2}]_{8\Phi })_{oT\_E_{1}}=-\frac{1}{6}(\chi _{11}\psi _{21}-\chi
_{12}\psi _{22}-(x_{11},w_{21})+(x_{12},w_{22})$

$\ \ \ \ \ \ \ \ \ \ \ \ \ \ \ \ \ \ \ \ \ \ \ \ \ \ \ \ \ \ \ \ +\mu
_{21}\gamma _{11}-\mu _{22}\gamma _{12}-(z_{21},y_{11})+(z_{22},y_{12}))$

$\ \ \ \ \ \ \ \ \ \ \ \ \ \ \ \ \ \ \ \ \ \ \ \ \ \ \ +\frac{1}{6}(\chi
_{13}\psi _{23}-\chi _{11}\psi _{21}-(x_{13},w_{23})+(x_{11},w_{21})$

$\ \ \ \ \ \ \ \ \ \ \ \ \ \ \ \ \ \ \ \ \ \ \ \ \ \ \ \ \ \ \ \ +\mu
_{23}\gamma _{13}-\mu _{21}\gamma _{11}-(z_{23},y_{13})+(z_{21},y_{11}))$

\ \ \ $\ \ \ \ \ \ \ \ \ \ \ \ \ \ \ \ \ \ \ \ \ \ \ +\frac{1}{6}(\chi
_{21}\psi _{11}-\chi _{22}\psi _{12}-(x_{21},w_{11})+(x_{22},w_{12})$

$\ \ \ \ \ \ \ \ \ \ \ \ \ \ \ \ \ \ \ \ \ \ \ \ \ \ \ \ \ \ \ \ +\mu
_{11}\gamma _{21}-\mu _{12}\gamma _{22}-(z_{11},y_{21})+(z_{12},y_{22}))$

$\ \ \ \ \ \ \ \ \ \ \ \ \ \ \ \ \ \ \ \ \ \ \ \ \ \ -\frac{1}{6}(\chi
_{23}\psi _{13}-\chi _{21}\psi _{11}-(x_{23},w_{13})+(x_{21},w_{11})$

$\ \ \ \ \ \ \ \ \ \ \ \ \ \ \ \ \ \ \ \ \ \ \ \ \ \ \ \ \ \ \ \ \ +\mu
_{13}\gamma _{23}-\mu _{11}\gamma _{21}-(z_{13},y_{23})+(z_{11},y_{21})),$

$([R_{1},R_{2}]_{8\Phi })_{oT\_E_{2}}=-\frac{1}{6}(\chi _{12}\psi _{22}-\chi
_{13}\psi _{23}-(x_{12},w_{22})+(x_{13},w_{23})$

$\ \ \ \ \ \ \ \ \ \ \ \ \ \ \ \ \ \ \ \ \ \ \ \ \ \ \ \ \ \ \ \ +\mu
_{22}\gamma _{12}-\mu _{23}\gamma _{13}-(z_{22},y_{12})+(z_{23},y_{13}))$

$\ \ \ \ \ \ \ \ \ \ \ \ \ \ \ \ \ \ \ \ \ \ \ \ \ +\frac{1}{6}(\chi
_{11}\psi _{21}-\chi _{12}\psi _{22}-(x_{11},w_{21})+(x_{12},w_{22})$

$\ \ \ \ \ \ \ \ \ \ \ \ \ \ \ \ \ \ \ \ \ \ \ \ \ \ \ \ \ \ \ \ \ +\mu
_{21}\gamma _{11}-\mu _{22}\gamma _{12}-(z_{21},y_{11})+(z_{22},y_{12}))$

$\ \ \ \ \ \ \ \ \ \ \ \ \ \ \ \ \ \ \ \ \ \ \ \ \ +\frac{1}{6}(\chi
_{22}\psi _{12}-\chi _{23}\psi _{13}-(x_{22},w_{12})+(x_{23},w_{13})$

$\ \ \ \ \ \ \ \ \ \ \ \ \ \ \ \ \ \ \ \ \ \ \ \ \ \ \ \ \ \ \ \ +\mu
_{12}\gamma _{22}-\mu _{23}\gamma _{23}-(z_{12},y_{22})+(z_{13},y_{23}))$

$\ \ \ \ \ \ \ \ \ \ \ \ \ \ \ \ \ \ \ \ \ \ \ \ \ -\frac{1}{6}(\chi
_{21}\psi _{11}-\chi _{22}\psi _{12}-(x_{21},w_{11})+(x_{22},w_{12})$

$\ \ \ \ \ \ \ \ \ \ \ \ \ \ \ \ \ \ \ \ \ \ \ \ \ \ \ \ \ \ \ \ +\mu
_{11}\gamma _{21}-\mu _{12}\gamma _{22}-(z_{11},y_{21})+(z_{12},y_{22})),$

$([R_{1},R_{2}]_{8\Phi })_{T\_E_{3}}=-\frac{1}{6}(\chi _{13}\psi _{23}-\chi
_{11}\psi _{21}-(x_{13},w_{23})+(x_{11},w_{21})$

$\ \ \ \ \ \ \ \ \ \ \ \ \ \ \ \ \ \ \ \ \ \ \ \ \ \ \ \ \ \ \ \ \ +\mu
_{23}\gamma _{13}-\mu _{21}\gamma _{11}-(z_{23},y_{13})+(z_{21},y_{11}))$

$\ \ \ \ \ \ \ \ \ \ \ \ \ \ \ \ \ \ \ \ \ \ \ \ \ \ +\frac{1}{6}(\chi
_{12}\psi _{22}-\chi _{13}\psi _{23}-(x_{12},w_{22})+(x_{13},w_{23})$

$\ \ \ \ \ \ \ \ \ \ \ \ \ \ \ \ \ \ \ \ \ \ \ \ \ \ \ \ \ \ \ +\mu
_{22}\gamma _{12}-\mu _{23}\gamma _{13}-(z_{22},y_{12})+(z_{23},y_{13}))$

$\ \ \ \ \ \ \ \ \ \ \ \ \ \ \ \ \ \ \ \ \ \ \ \ \ \ +\frac{1}{6}(\chi
_{23}\psi _{13}-\chi _{21}\psi _{11}-(x_{23},w_{13})+(x_{21},w_{11})$

$\ \ \ \ \ \ \ \ \ \ \ \ \ \ \ \ \ \ \ \ \ \ \ \ \ \ \ \ \ \ \ \ \ +\mu
_{13}\gamma _{23}-\mu _{11}\gamma _{21}-(z_{13},y_{23})+(z_{11},y_{21}))$

$\ \ \ \ \ \ \ \ \ \ \ \ \ \ \ \ \ \ \ \ \ \ \ \ \ -\frac{1}{6}(\chi
_{22}\psi _{12}-\chi _{23}\psi _{13}-(x_{22},w_{12})+(x_{23},w_{13})$

$\ \ \ \ \ \ \ \ \ \ \ \ \ \ \ \ \ \ \ \ \ \ \ \ \ \ \ \ \ \ \ \ \ +\mu
_{12}\gamma _{22}-\mu _{23}\gamma _{23}-(z_{12},y_{22})+(z_{13},y_{23})),$

$([R_{1},R_{2}]_{8\Phi })_{oT\_F_{1}}=-\frac{1}{4}((\chi _{12}+\chi
_{13})w_{21}+(\psi _{22}+\psi _{23})x_{11}+\overline{w_{22}x_{13}}+\overline{%
x_{12}w_{23}}$

$\ \ \ \ \ \ \ \ \ \ \ \ \ \ \ \ \ \ \ \ \ \ \ \ \ \ \ \ \ \ +(\mu _{22}+\mu
_{23})y_{11}+(\gamma _{12}+\gamma _{13})z_{21}+\overline{y_{12}z_{23}}+%
\overline{z_{22}y_{13}})$

$\ \ \ \ \ \ \ \ \ \ \ \ \ \ \ \ \ \ \ \ \ \ \ \ \ \ \ +\frac{1}{4}((\chi
_{22}+\chi _{23})w_{11}+(\psi _{12}+\psi _{13})x_{21}+\overline{w_{12}x_{23}}%
+\overline{x_{22}w_{13}}$

$\ \ \ \ \ \ \ \ \ \ \ \ \ \ \ \ \ \ \ \ \ \ \ \ \ \ \ \ \ \ +(\mu _{12}+\mu
_{13})y_{21}+(\gamma _{22}+\gamma _{23})z_{11}+\overline{y_{22}z_{13}}+%
\overline{z_{12}y_{23}}),$

$([R_{1},R_{2}]_{8\Phi })_{oT\_F_{2}}=-\frac{1}{4}((\chi _{13}+\chi
_{11})w_{22}+(\psi _{23}+\psi _{21})x_{12}+\overline{w_{23}x_{11}}+\overline{%
x_{13}w_{21}}$

$\ \ \ \ \ \ \ \ \ \ \ \ \ \ \ \ \ \ \ \ \ \ \ \ \ \ \ \ \ \ +(\mu _{23}+\mu
_{21})y_{12}+(\gamma _{13}+\gamma _{11})z_{22}+\overline{y_{13}z_{21}}+%
\overline{z_{23}y_{11}})$

$\ \ \ \ \ \ \ \ \ \ \ \ \ \ \ \ \ \ \ \ \ \ \ \ \ \ \ +\frac{1}{4}((\chi
_{23}+\chi _{21})w_{12}+(\psi _{13}+\psi _{11})x_{22}+\overline{w_{13}x_{21}}%
+\overline{x_{23}w_{11}}$

$\ \ \ \ \ \ \ \ \ \ \ \ \ \ \ \ \ \ \ \ \ \ \ \ \ \ \ \ \ \ \ +(\mu
_{13}+\mu _{11})y_{22}+(\gamma _{23}+\gamma _{21})z_{12}+\overline{%
y_{23}z_{11}}+\overline{z_{13}y_{21}}),$

$([R_{1},R_{2}]_{8\Phi })_{oT\_F_{3}}=-\frac{1}{4}((\chi _{11}+\chi
_{12})w_{23}+(\psi _{21}+\psi _{22})x_{13}+\overline{w_{21}x_{12}}+\overline{%
x_{11}w_{22}}$

$\ \ \ \ \ \ \ \ \ \ \ \ \ \ \ \ \ \ \ \ \ \ \ \ \ \ \ \ \ \ +(\mu _{21}+\mu
_{22})y_{13}+(\gamma _{11}+\gamma _{12})z_{23}+\overline{y_{11}z_{22}}+%
\overline{z_{21}y_{12}})$

$\ \ \ \ \ \ \ \ \ \ \ \ \ \ \ \ \ \ \ \ \ \ \ \ \ \ \ +\frac{1}{4}((\chi
_{21}+\chi _{22})w_{13}+(\psi _{11}+\psi _{12})x_{23}+\overline{w_{11}x_{22}}%
+\overline{x_{21}w_{12}}$

$\ \ \ \ \ \ \ \ \ \ \ \ \ \ \ \ \ \ \ \ \ \ \ \ \ \ \ \ \ +(\mu _{11}+\mu
_{12})y_{23}+(\gamma _{21}+\gamma _{22})z_{13}+\overline{y_{21}z_{12}}+%
\overline{z_{11}y_{22}}),$

$([R_{1},R_{2}]_{8\Phi })_{A}=\left( 
\begin{array}{c}
([R_{1},R_{2}]_{8\Phi })_{A\_E_{1}} \\ 
([R_{1},R_{2}]_{8\Phi })_{A\_E_{2}} \\ 
([R_{1},R_{2}]_{8\Phi })_{A\_E_{3}} \\ 
([R_{1},R_{2}]_{8\Phi })_{A\_F_{1}} \\ 
([R_{1},R_{2}]_{8\Phi })_{A\_F_{2}} \\ 
([R_{1},R_{2}]_{8\Phi })_{A\_F_{3}}%
\end{array}%
\right) $

$\ \ \ \ \ \ \ \ \ \ \ \ \ \ \ \ \ \ \ =\left( 
\begin{array}{c}
\lbrack \Phi _{1},\Phi _{2}]_{7A\_E_{1}} \\ 
\lbrack \Phi _{1},\Phi _{2}]_{7A\_E_{2}} \\ 
\lbrack \Phi _{1},\Phi _{2}]_{7A\_E_{3}} \\ 
\lbrack \Phi _{1},\Phi _{2}]_{7A\_F_{1}} \\ 
\lbrack \Phi _{1},\Phi _{2}]_{7A\_F_{2}} \\ 
\lbrack \Phi _{1},\Phi _{2}]_{7A\_F_{3}}%
\end{array}%
\right) +\left( 
\begin{array}{c}
([R_{1},R_{2}]_{8\Phi })_{oA\_E_{1}} \\ 
([R_{1},R_{2}]_{8\Phi })_{oA\_E_{2}} \\ 
([R_{1},R_{2}]_{8\Phi })_{oA\_E_{3}} \\ 
([R_{1},R_{2}]_{8\Phi })_{oA\_F_{1}} \\ 
([R_{1},R_{2}]_{8\Phi })_{oA\_F_{2}} \\ 
([R_{1},R_{2}]_{8\Phi })_{oA\_F_{3}}%
\end{array}%
\right) ,$

$([R_{1},R_{2}]_{8\Phi })_{oA\_E_{1}}=-\frac{1}{4}((\gamma _{12}\psi
_{23}+\gamma _{13}\psi _{22})-2(y_{11},w_{21})-\xi _{1}\mu _{21}-\zeta
_{2}\chi _{11})$

$\ \ \ \ \ \ \ \ \ \ \ \ \ \ \ \ \ \ \ \ \ \ \ \ \ \ +\frac{1}{4}((\gamma
_{22}\psi _{13}+\gamma _{23}\psi _{12})-2(y_{21},w_{11})-\xi _{2}\mu
_{11}-\zeta _{1}\chi _{21}),$

$([R_{1},R_{2}]_{8\Phi })_{oA\_E_{2}}=-\frac{1}{4}((\gamma _{13}\psi
_{21}+\gamma _{11}\psi _{23})-2(y_{12},w_{22})-\xi _{1}\mu _{22}-\zeta
_{2}\chi _{12})$

$\ \ \ \ \ \ \ \ \ \ \ \ \ \ \ \ \ \ \ \ \ \ \ \ \ \ +\frac{1}{4}((\gamma
_{23}\psi _{11}+\gamma _{21}\psi _{13})-2(y_{22},w_{12})-\xi _{2}\mu
_{12}-\zeta _{1}\chi _{22}),$

$([R_{1},R_{2}]_{8\Phi })_{oA\_E_{3}}=-\frac{1}{4}((\gamma _{11}\psi
_{22}+\gamma _{12}\psi _{21})-2(y_{13},w_{23})-\xi _{1}\mu _{23}-\zeta
_{2}\chi _{13})$

$\ \ \ \ \ \ \ \ \ \ \ \ \ \ \ \ \ \ \ \ \ \ \ \ \ \ +\frac{1}{4}((\gamma
_{21}\psi _{12}+\gamma _{22}\psi _{11})-2(y_{23},w_{13})-\xi _{2}\mu
_{13}-\zeta _{1}\chi _{23}),$

$([R_{1},R_{2}]_{8\Phi })_{oA\_F_{1}}=-\frac{1}{4}(-\psi _{21}y_{11}-\gamma
_{11}w_{21}+\overline{y_{12}w_{23}}+\overline{w_{22}y_{13}}-\xi
_{1}z_{21}-\zeta _{2}x_{11})$

$\ \ \ \ \ \ \ \ \ \ \ \ \ \ \ \ \ \ \ \ \ \ \ \ \ \ +\frac{1}{4}(-\psi
_{11}y_{21}-\gamma _{21}w_{11}+\overline{y_{22}w_{13}}+\overline{w_{12}y_{23}%
}-\xi _{2}z_{11}-\zeta _{1}x_{21}),$

$([R_{1},R_{2}]_{8\Phi })_{oA\_F_{2}}=-\frac{1}{4}(-\psi _{22}y_{12}-\gamma
_{12}w_{22}+\overline{y_{13}w_{21}}+\overline{w_{23}y_{11}}-\xi
_{1}z_{22}-\zeta _{2}x_{12})$

$\ \ \ \ \ \ \ \ \ \ \ \ \ \ \ \ \ \ \ \ \ \ \ \ \ \ +\frac{1}{4}(-\psi
_{12}y_{22}-\gamma _{22}w_{12}+\overline{y_{23}w_{11}}+\overline{w_{13}y_{21}%
}-\xi _{2}z_{12}-\zeta _{1}x_{22}),$

$([R_{1},R_{2}]_{8\Phi })_{oA\_F_{3}}=-\frac{1}{4}(-\psi _{23}y_{13}-\gamma
_{13}w_{23}+\overline{y_{11}w_{22}}+\overline{w_{21}y_{12}}-\xi
_{1}z_{23}-\zeta _{2}x_{13})$

$\ \ \ \ \ \ \ \ \ \ \ \ \ \ \ \ \ \ \ \ \ \ \ \ \ \ +\frac{1}{4}(-\psi
_{13}y_{23}-\gamma _{23}w_{13}+\overline{y_{21}w_{12}}+\overline{w_{11}y_{22}%
}-\xi _{2}z_{13}-\zeta _{1}x_{23})$,

$([R_{1},R_{2}]_{8\Phi })_{B}=\left( 
\begin{array}{c}
([R_{1},R_{2}]_{8\Phi })_{B\_E_{1}} \\ 
([R_{1},R_{2}]_{8\Phi })_{B\_E_{2}} \\ 
([R_{1},R_{2}]_{8\Phi })_{B\_E_{3}} \\ 
([R_{1},R_{2}]_{8\Phi })_{B\_F_{1}} \\ 
([R_{1},R_{2}]_{8\Phi })_{B\_F_{2}} \\ 
([R_{1},R_{2}]_{8\Phi })_{B\_F_{3}}%
\end{array}%
\right) $

$\ \ \ \ \ \ \ \ \ \ \ \ \ \ \ \ \ \ \ =\left( 
\begin{array}{c}
\lbrack \Phi _{1},\Phi _{2}]_{7B\_E_{1}} \\ 
\lbrack \Phi _{1},\Phi _{2}]_{7B\_E_{2}} \\ 
\lbrack \Phi _{1},\Phi _{2}]_{7B\_E_{3}} \\ 
\lbrack \Phi _{1},\Phi _{2}]_{7B\_F_{1}} \\ 
\lbrack \Phi _{1},\Phi _{2}]_{7B\_F_{2}} \\ 
\lbrack \Phi _{1},\Phi _{2}]_{7B\_F_{3}}%
\end{array}%
\right) +\left( 
\begin{array}{c}
([R_{1},R_{2}]_{8\Phi })_{oB\_E_{1}} \\ 
([R_{1},R_{2}]_{8\Phi })_{oB\_E_{2}} \\ 
([R_{1},R_{2}]_{8\Phi })_{oB\_E_{3}} \\ 
([R_{1},R_{2}]_{8\Phi })_{oB\_F_{1}} \\ 
([R_{1},R_{2}]_{8\Phi })_{oB\_F_{2}} \\ 
([R_{1},R_{2}]_{8\Phi })_{oB\_F_{3}}%
\end{array}%
\right) ,$

$([R_{1},R_{2}]_{8\Phi })_{oB\_E_{1}}=\frac{1}{4}((\chi _{12}\mu _{23}+\chi
_{13}\mu _{22})-2(x_{11},z_{21})-\eta _{1}\psi _{21}-\omega _{2}\gamma
_{11}) $

$\ \ \ \ \ \ \ \ \ \ \ \ \ \ \ \ \ \ \ \ \ \ \ \ -\frac{1}{4}((\chi _{22}\mu
_{13}+\chi _{23}\mu _{12})-2(x_{21},z_{11})-\eta _{2}\psi _{11}-\omega
_{1}\gamma _{21}),$

$([R_{1},R_{2}]_{8\Phi })_{oB\_E_{2}}=\frac{1}{4}((\chi _{13}\mu _{21}+\chi
_{11}\mu _{23})-2(x_{12},z_{22})-\eta _{1}\psi _{22}-\omega _{2}\gamma
_{12}) $

$\ \ \ \ \ \ \ \ \ \ \ \ \ \ \ \ \ \ \ \ \ \ \ \ -\frac{1}{4}((\chi _{23}\mu
_{11}+\chi _{21}\mu _{13})-2(x_{22},z_{12})-\eta _{2}\psi _{12}-\omega
_{1}\gamma _{22}),$

$([R_{1},R_{2}]_{8\Phi })_{oB\_E_{3}}=\frac{1}{4}((\chi _{11}\mu _{22}+\chi
_{12}\mu _{21})-2(x_{13},z_{23})-\eta _{1}\psi _{23}-\omega _{2}\gamma
_{13}) $

$\ \ \ \ \ \ \ \ \ \ \ \ \ \ \ \ \ \ \ \ \ \ \ \ -\frac{1}{4}((\chi _{21}\mu
_{12}+\chi _{22}\mu _{11})-2(x_{23},z_{13})-\eta _{2}\psi _{13}-\omega
_{1}\gamma _{23}),$

$([R_{1},R_{2}]_{8\Phi })_{oB\_F_{1}}=\frac{1}{4}(-\mu _{21}x_{11}-\chi
_{11}z_{21}+\overline{x_{12}z_{23}}+\overline{z_{22}x_{13}})-\eta
_{1}w_{21}-\omega _{2}y_{11})$

$\ \ \ \ \ \ \ \ \ \ \ \ \ \ \ \ \ \ \ \ \ \ \ \ \ -\frac{1}{4}(-\mu
_{11}x_{21}-\chi _{21}z_{11}+\overline{x_{22}z_{13}}+\overline{z_{12}x_{23}}%
)-\eta _{2}w_{11}-\omega _{1}y_{21}),$

$([R_{1},R_{2}]_{8\Phi })_{oB\_F_{2}}=\frac{1}{4}(-\mu _{22}x_{12}-\chi
_{12}z_{22}+\overline{x_{13}z_{21}}+\overline{z_{23}x_{11}})-\eta
_{1}w_{22}-\omega _{2}y_{12})$

$\ \ \ \ \ \ \ \ \ \ \ \ \ \ \ \ \ \ \ \ \ \ \ \ \ -\frac{1}{4}(-\mu
_{12}x_{22}-\chi _{22}z_{12}+\overline{x_{23}z_{11}}+\overline{z_{13}x_{21}}%
)-\eta _{2}w_{12}-\omega _{1}y_{22}),$

$([R_{1},R_{2}]_{8\Phi })_{oB\_F_{3}}=\frac{1}{4}(-\mu _{23}x_{13}-\chi
_{13}z_{23}+\overline{x_{11}z_{22}}+\overline{z_{21}x_{12}})-\eta
_{1}w_{23}-\omega _{2}y_{13})$

$\ \ \ \ \ \ \ \ \ \ \ \ \ \ \ \ \ \ \ \ \ \ \ \ \ -\frac{1}{4}(-\mu
_{13}x_{23}-\chi _{23}z_{13}+\overline{x_{21}z_{12}}+\overline{z_{11}x_{22}}%
)-\eta _{2}w_{13}-\omega _{1}y_{23}),$

$([R_{1},R_{2}]_{8\Phi })_{\rho }=[\Phi _{1},\Phi _{2}]_{7\rho
}+([R_{1},R_{2}]_{8\Phi })_{o\rho }$

$([R_{1},R_{2}]_{8\Phi })_{o\rho }=$

\ \ \ \ \ \ \ \ \ \ \ \ \ \ \ \ \ \ \ $\frac{1}{8}(\sum_{i=1}^{3}(\chi
_{1i}\psi _{2i}+2(x_{1i},w_{2i})+\mu _{2i}\gamma _{1i}+2(z_{2i},y_{1i}))-3(%
\xi _{1}\omega _{2}+\zeta _{2}\eta _{1}))$

\ \ \ \ \ \ \ \ \ \ \ \ \ \ \ \ \ \ \ -$\frac{1}{8}%
(\sum_{i=1}^{3}(\chi _{2i}\psi _{1i}+2(x_{2i},w_{1i})+\mu _{1i}\gamma
_{2i}+2(z_{1i},y_{2i}))-3(\xi _{2}\omega _{1}+\zeta _{1}\eta _{2})).$

\bigskip

\emph{\ Proof. }Using \emph{Lemma 5.9} and definition of mappings $%
\textrm{f}_{DDD},\textrm{f}_{JJD},\textrm{f}_{DJJ}$,

\noindent
$\textrm{f}_{iJJ},\textrm{f}_{JJC},\textrm{f}_{CJJ},\textrm{f}_{CCC},$
we have the above expression. \ \ \ \ \emph{Q.E.D.}

\bigskip

\emph{Lemma 10.2. }

$\left( 
\begin{array}{cccccc}
F_{P\Phi } & F_{PP} & F_{PQ} & F_{P\text{\textit{r}}} & F_{Ps} & 0%
\end{array}%
%
\right) ,$

$([R_{1},R_{2}]_{8P})_{X%
\_E_{1}}=-(m_{12},x_{22})+(m_{13},x_{23})+(m_{22},x_{12})-(m_{23},x_{13})$

$\ +\tau _{11}\chi _{21}+(t_{12},x_{22})+(t_{13},x_{23})-\tau _{21}\chi
_{11}-(t_{22},x_{12})-(t_{23},x_{13})$

$-\frac{1}{3}\rho _{1}\chi _{21}+\frac{1%
}{3}\rho _{2}\chi _{11}$

$\ +((\beta _{12}\gamma _{23}+\beta _{13}\gamma
_{22})-2(b_{11},y_{21})-((\beta _{22}\gamma _{13}+\beta _{23}\gamma
_{12})-2(b_{21},y_{11}))$

$\ +\eta _{2}\alpha _{11}-\eta _{1}\alpha _{21}+r_{1}\chi _{21}-r_{2}\chi
_{11}+s_{1}\mu _{21}-s_{2}\mu _{11},$

$([R_{1},R_{2}]_{8P})_{X%
\_E_{2}}=-(m_{13},x_{23})+(m_{11},x_{21})+(m_{23},x_{13})-(m_{21},x_{11})$

$\ +\tau _{12}\chi _{22}+(t_{13},x_{23})+(t_{11},x_{21})-\tau _{22}\chi
_{12}-(t_{23},x_{13})-(t_{21},x_{11})$

$-\frac{1}{3}\rho _{1}\chi _{22}+\frac{1%
}{3}\rho _{2}\chi _{12}$

$\ +((\beta _{13}\gamma _{21}+\beta _{11}\gamma
_{23})-2(b_{12},y_{22}))-((\beta _{23}\gamma _{11}+\beta _{21}\gamma
_{13})-2(b_{22},y_{12}))$

$\ +\eta _{2}\alpha _{12}-\eta _{1}\alpha _{22}+r_{1}\chi _{22}-r_{2}\chi
_{12}+s_{1}\mu _{22}-s_{2}\mu _{12},$

$([R_{1},R_{2}]_{8P})_{X%
\_E_{3}}=-(m_{11},x_{21})+(m_{12},x_{22})+(m_{21},x_{11})-(m_{22},x_{12})$

$\ +(-\tau _{11}-\tau _{12})\chi
_{23}+(t_{11},x_{21})+(t_{12},x_{22})$

$-(-\tau _{21}-\tau _{22})\chi
_{13}-(t_{21},x_{11})-(t_{22},x_{12})-\frac{1}{3}\rho _{1}\chi _{23}+\frac{1%
}{3}\rho _{2}\chi _{13}$

$\ +(\beta _{11}\gamma _{22}+\beta _{12}\gamma
_{21})-2(b_{13},y_{23}))-((\beta _{21}\gamma _{12}+\beta _{22}\gamma
_{11})-2(b_{23},y_{13}))$

$\ +\eta _{2}\alpha _{13}-\eta _{1}\alpha _{23}+r_{1}\chi _{23}-r_{2}\chi
_{13}+s_{1}\mu _{23}-s_{2}\mu _{13},$

$([R_{1},R_{2}]_{8P})_{X%
\_F_{1}}=\textrm{d}_{g}\textrm{g}_{d}(D_{1})x_{21}-\textrm{d}_{g}\textrm{g}_{d}(D_{2})x_{11} $

$\ +\frac{1}{2}(\overline{m_{12}x_{23}}-\overline{x_{22}m_{13}}+(-\chi
_{22}+\chi _{23})m_{11})$

$-\frac{1}{2}(\overline{m_{22}x_{13}}-\overline{%
x_{12}m_{23}}+(-\chi _{12}+\chi _{13})m_{21})$

$\ +\frac{1}{2}(-\tau _{11}x_{21}+(\chi _{22}+\chi _{23})t_{11}+\overline{%
x_{22}t_{13}}+\overline{t_{12}x_{23}})$

$\ -\frac{1}{2}(-\tau _{21}x_{11}+(\chi _{12}+\chi _{13})t_{21}+\overline{%
x_{12}t_{23}}+\overline{t_{22}x_{13}})-\frac{1}{3}\rho _{1}x_{21}+\frac{1}{3}%
\rho _{2}x_{11}$

$\ +(-\gamma _{21}b_{11}-\beta _{11}y_{21}+\overline{b_{12}y_{23}}+\overline{%
y_{22}b_{13}})-(-\gamma _{11}b_{21}-\beta _{21}y_{11}+\overline{b_{22}y_{13}}%
+\overline{y_{12}b_{23}})$

$\ +\eta _{2}a_{11}-\eta
_{1}a_{21}+r_{1}x_{21}-r_{2}x_{11}+s_{1}z_{21}-s_{2}z_{11},$

$([R_{1},R_{2}]_{8P})_{X\_F_{2}}=\textrm{d}_{g}\nu \textrm{g}_{d}(D_{1})x_{22}-\textrm{d}_{g}\nu
\textrm{g}_{d}(D_{2})x_{12}$

$\ +\frac{1}{2}(\overline{m_{13}x_{21}}-\overline{x_{23}m_{11}}+(-\chi
_{23}+\chi _{21})m_{12})$

$-\frac{1}{2}(\overline{m_{23}x_{11}}-\overline{%
x_{13}m_{21}}+(-\chi _{13}+\chi _{11})m_{22})$

$\ +\frac{1}{2}(-\tau _{12}x_{22}+(\chi _{23}+\chi _{21})t_{12}+\overline{%
x_{23}t_{11}}+\overline{t_{13}x_{21}})$

$\ -\frac{1}{2}(-\tau _{22}x_{12}+(\chi _{13}+\chi _{11})t_{22}+\overline{%
x_{13}t_{21}}+\overline{t_{23}x_{11}})-\frac{1}{3}\rho _{1}x_{22}+\frac{1}{3}%
\rho _{2}x_{12}$

$\ +(-\gamma _{22}b_{12}-\beta _{12}y_{22}+\overline{b_{13}y_{21}}+\overline{%
y_{23}b_{11}})-(-\gamma _{12}b_{22}-\beta _{22}y_{12}+\overline{b_{23}y_{11}}%
+\overline{y_{13}b_{21}})$

$\ +\eta _{2}a_{12}-\eta
_{1}a_{22}+r_{1}x_{22}-r_{2}x_{12}+s_{1}z_{22}-s_{2}z_{12},$

$([R_{1},R_{2}]_{8P})_{X\_F_{3}}=\textrm{d}_{g}\nu ^{2}\textrm{g}_{d}(D_{1})x_{23}-\textrm{d}_{g}\nu
^{2}\textrm{g}_{d}(D_{2})x_{13}$

$\ +\frac{1}{2}(\overline{m_{11}x_{22}}-\overline{x_{21}m_{12}}+(-\chi
_{21}+\chi _{22})m_{13})$

$-\frac{1}{2}(\overline{m_{21}x_{12}}-\overline{%
x_{11}m_{22}}+(-\chi _{11}+\chi _{12})m_{23})$

$\ +\frac{1}{2}((\tau _{11}+\tau _{12})x_{23}+(\chi _{21}+\chi _{22})t_{13}+%
\overline{x_{21}t_{12}}+\overline{t_{11}x_{22}})$

$\ -\frac{1}{2}((\tau _{21}+\tau _{22})x_{13}+(\chi _{11}+\chi _{12})t_{23}+%
\overline{x_{11}t_{22}}+\overline{t_{21}x_{12}})-\frac{1}{3}\rho _{1}x_{23}+%
\frac{1}{3}\rho _{2}x_{13}$

$\ +(-\gamma _{23}b_{13}-\beta _{13}y_{23}+\overline{b_{11}y_{22}}+\overline{%
y_{21}b_{12}})-(-\gamma _{13}b_{23}-\beta _{23}y_{13}+\overline{b_{21}y_{12}}%
+\overline{y_{11}b_{22}})$

$\ +\eta _{2}a_{13}-\eta
_{1}a_{23}+r_{1}x_{23}-r_{2}x_{13}+s_{1}z_{23}-s_{2}z_{13},$

$([R_{1},R_{2}]_{8P})_{Y}=\left( 
\begin{array}{c}
([R_{1},R_{2}]_{8P})_{Y\_E_{1}} \\ 
([R_{1},R_{2}]_{8P})_{Y\_E_{2}} \\ 
([R_{1},R_{2}]_{8P})_{Y\_E_{3}} \\ 
([R_{1},R_{2}]_{8P})_{Y\_F_{1}} \\ 
([R_{1},R_{2}]_{8P})_{Y\_F_{2}} \\ 
([R_{1},R_{2}]_{8P})_{Y\_F_{3}}%
\end{array}%
\right) ,$

$([R_{1},R_{2}]_{8P})_{Y%
\_E_{1}}=(-(m_{12},y_{22})+(m_{13},y_{23}))-(-(m_{22},y_{12})+(m_{23},y_{13}))$

$\ -(\tau _{11}\gamma _{21}+(t_{12},y_{22})+(t_{13},y_{23}))+(\tau
_{21}\gamma _{11}+(t_{22},y_{12})+(t_{23},y_{13}))$

$+\frac{1}{3}\rho
_{1}\gamma _{21}-\frac{1}{3}\rho _{2}\gamma _{11}$

$\ +((\alpha _{12}\chi _{23}+\alpha _{13}\chi
_{22})-2(a_{11},x_{21}))-((\alpha _{22}\chi _{13}+\alpha _{23}\chi
_{12})-2(a_{21},x_{11}))$

$\ +\xi _{2}\beta _{11}-\xi _{1}\beta _{21}+r_{1}\gamma _{21}-r_{2}\gamma
_{11}+s_{1}\psi _{21}-s_{2}\psi _{11},$

$([R_{1},R_{2}]_{8P})_{Y%
\_E_{2}}=(-(m_{13},y_{23})+(m_{11},y_{21}))-(-(m_{23},y_{13})+(m_{21},y_{11}))$

$\ -(\tau _{12}\gamma _{22}+(t_{13},y_{23})+(t_{11},y_{21}))+(\tau
_{22}\gamma _{12}+(t_{23},y_{13})+(t_{21},y_{11}))$

$+\frac{1}{3}\rho
_{1}\gamma _{22}-\frac{1}{3}\rho _{2}\gamma _{12}$

$\ +((\alpha _{13}\chi _{21}+\alpha _{11}\chi
_{23})-2(a_{12},x_{22}))-((\alpha _{23}\chi _{11}+\alpha _{21}\chi
_{13})-2(a_{22},x_{12}))$

$\ +\xi _{2}\beta _{12}-\xi _{1}\beta _{22}+r_{1}\gamma _{22}-r_{2}\gamma
_{12}+s_{1}\psi _{22}-s_{2}\psi _{12,}$

$([R_{1},R_{2}]_{8P})_{Y%
\_E_{3}}=(-(m_{11},y_{21})+(m_{12},y_{22}))-(-(m_{21},y_{11})+(m_{22},y_{12}))$

$\ -((-\tau _{11}-\tau _{12})\gamma
_{23}+(t_{11},y_{21})+(t_{12},y_{22}))$

$+((-\tau _{21}-\tau _{22})\gamma
_{13}+(t_{21},y_{11})+(t_{22},y_{12}))$

$\ \ \ \ +\frac{1}{3}\rho _{1}\gamma _{23}-\frac{1}{3}\rho _{2}\gamma _{13}$

$\ +((\alpha _{11}\chi _{22}+\alpha _{12}\chi
_{21})-2(a_{13},x_{23}))-((\alpha _{21}\chi _{12}+\alpha _{22}\chi
_{11})-2(a_{23},x_{13}))$

$\ +\xi _{2}\beta _{13}-\xi _{1}\beta _{23}+r_{1}\gamma _{23}-r_{2}\gamma
_{13}+s_{1}\psi _{23}-s_{2}\psi _{13},$

$([R_{1},R_{2}]_{8P})_{Y%
\_F_{1}}=\textrm{d}_{g}\textrm{g}_{d}(D_{1})y_{21}-\textrm{d}_{g}\textrm{g}_{d}(D_{2})y_{11} $

$\ +(\frac{1}{2}(\overline{m_{12}y_{23}}-\overline{y_{22}m_{13}}+(-\gamma
_{22}+\gamma _{23})m_{11}))$

$-(\frac{1}{2}(\overline{m_{22}y_{13}}-\overline{%
y_{12}m_{23}}+(-\gamma _{12}+\gamma _{13})m_{21}))$

$\ -(\frac{1}{2}(-\tau _{11}y_{21}+(\gamma _{22}+\gamma _{23})t_{11}+%
\overline{y_{22}t_{13}}+\overline{t_{12}y_{23}}))$

$\ +(\frac{1}{2}(-\tau _{21}y_{11}+(\gamma _{12}+\gamma _{13})t_{21}+%
\overline{y_{12}t_{23}}+\overline{t_{22}y_{13}}))+\frac{1}{3}\rho _{1}y_{21}-%
\frac{1}{3}\rho _{2}y_{11}$

$\ +(-\chi _{21}a_{11}-\alpha _{11}x_{21}+\overline{a_{12}x_{23}}+\overline{%
x_{22}a_{13}})-(-\chi _{11}a_{21}-\alpha _{21}x_{11}+\overline{a_{22}x_{13}}+%
\overline{x_{12}a_{23}})$

$\ +\xi _{2}b_{11}-\xi
_{1}b_{21}+r_{1}y_{21}-r_{2}y_{11}+s_{1}w_{21}-s_{2}w_{11},$

$([R_{1},R_{2}]_{8P})_{Y\_F_{2}}=\textrm{d}_{g}\nu \textrm{g}_{d}(D_{1})y_{22}-\textrm{d}_{g}\nu
\textrm{g}_{d}(D_{2})y_{12}$

$\ +\frac{1}{2}(\overline{m_{13}y_{21}}-\overline{y_{23}m_{11}}+(-\gamma
_{23}+\gamma _{21})m_{12})$

$-\frac{1}{2}(\overline{m_{23}y_{11}}-\overline{%
y_{13}m_{21}}+(-\gamma _{13}+\gamma _{11})m_{22})$

$\ -\frac{1}{2}(-\tau _{12}y_{22}+(\gamma _{23}+\gamma _{21})t_{12}+%
\overline{y_{23}t_{11}}+\overline{t_{13}y_{21}})$

$\ +\frac{1}{2}(-\tau _{22}y_{12}+(\gamma _{13}+\gamma _{11})t_{22}+%
\overline{y_{13}t_{21}}+\overline{t_{23}y_{11}})+\frac{1}{3}\rho _{1}y_{22}-%
\frac{1}{3}\rho _{2}y_{12}$

$\ +(-\chi _{22}a_{12}-\alpha _{12}x_{22}+\overline{a_{13}x_{21}}+\overline{%
x_{23}a_{11}})-(-\chi _{12}a_{22}-\alpha _{22}x_{12}+\overline{a_{23}x_{11}}+%
\overline{x_{13}a_{21}})$

$\ +\xi _{2}b_{12}-\xi
_{1}b_{22}+r_{1}y_{22}-r_{2}y_{12}+s_{1}w_{22}-s_{2}w_{12},$

$([R_{1},R_{2}]_{8P})_{Y\_F_{3}}=\textrm{d}_{g}\nu ^{2}\textrm{g}_{d}(D_{1})y_{23}-\textrm{d}_{g}\nu
^{2}\textrm{g}_{d}(D_{2})y_{13}$

$\ +\frac{1}{2}(\overline{m_{11}y_{22}}-\overline{y_{21}m_{12}}+(-\gamma
_{21}+\gamma _{22})m_{13})$

$-\frac{1}{2}(\overline{m_{21}y_{12}}-\overline{%
y_{11}m_{22}}+(-\gamma _{11}+\gamma _{12})m_{23})$

$\ -\frac{1}{2}((\tau _{11}+\tau _{12})y_{23}+(\gamma _{21}+\gamma
_{22})t_{13}+\overline{y_{21}t_{12}}+\overline{t_{11}y_{22}})$

$\ +\frac{1}{2}((\tau _{21}+\tau _{22})y_{13}+(\gamma _{11}+\gamma
_{12})t_{23}+\overline{y_{11}t_{22}}+\overline{t_{21}y_{12}})+\frac{1}{3}%
\rho _{1}y_{23}-\frac{1}{3}\rho _{2}y_{13}$

$\ +(-\chi _{23}a_{13}-\alpha _{13}x_{23}+\overline{a_{11}x_{22}}+\overline{%
x_{21}a_{12}})-(-\chi _{13}a_{23}-\alpha _{23}x_{13}+\overline{a_{21}x_{12}}+%
\overline{x_{11}a_{22}})$

$\ +\xi _{2}b_{13}-\xi
_{1}b_{23}+r_{1}y_{23}-r_{2}y_{13}+s_{1}w_{23}-s_{2}w_{13},$

$([R_{1},R_{2}]_{8P})_{\xi }=\sum_{i=1}^{3}(\alpha _{1i}\gamma
_{2i}+2(a_{1i},y_{2i}))-\sum_{i=1}^{3}(\alpha _{2i}\gamma
_{1i}+2(a_{2i},y_{1i}))$

$\ \ +\rho _{1}\xi _{2}-\rho _{2}\xi _{1}+r_{1}\xi _{2}-r_{2}\xi
_{1}+s_{1}\zeta _{2}-s_{2}\zeta _{1},$

$([R_{1},R_{2}]_{8P})_{\eta }=\sum_{i=1}^{3}(\beta _{1i}\chi
_{2i}+2(b_{1i},x_{2i}))-\sum_{i=1}^{3}(\beta _{2i}\chi
_{1i}+2(b_{2i},x_{1i}))$

$\ \ -\rho _{1}\eta _{2}+\rho _{2}\eta _{1}+r_{1}\eta _{2}-r_{2}\eta
_{1}+s_{1}\omega _{2}-s_{2}\omega _{1}.$

\bigskip

\emph{\ Proof. }Using \emph{Lemma 5.9} and definition of mappings $%
\textrm{f}_{DDD},\textrm{f}_{JJD},\textrm{f}_{DJJ}$,

\noindent
$\textrm{f}_{iJJ},\textrm{f}_{JJC},\textrm{f}_{CJJ},\textrm{f}_{CCC},$
we have the above expression. \ \ \ \ \emph{Q.E.D.}

\bigskip

\emph{Lemma 10.3. }

$\left( 
%
\right) ,$

$([R_{1},R_{2}]_{8Q})_{Z%
\_E_{1}}=(-(m_{12},z_{22})+(m_{13},z_{23}))-(-(m_{22},z_{12})+(m_{23},z_{13}))$

$\ +(\tau _{11}\mu _{21}+(t_{12},z_{22})+(t_{13},z_{23}))-(\tau _{21}\mu
_{11}+(t_{22},z_{12})+(t_{23},z_{13}))$

$-\frac{1}{3}\rho _{1}\mu _{21}+\frac{1%
}{3}\rho _{2}\mu _{11}$

$\ +((\beta _{12}\psi _{23}+\beta _{13}\psi _{22})-2(b_{11},w_{21}))-((\beta
_{22}\psi _{13}+\beta _{23}\psi _{12})-2(b_{21},w_{11}))$

$\ +\omega _{2}\alpha _{11}-\omega _{1}\alpha _{21}-r_{1}\mu _{21}+r_{2}\mu
_{11}+u_{1}\chi _{21}-u_{2}\chi _{11},$

$([R_{1},R_{2}]_{8Q})_{Z%
\_E_{2}}=(-(m_{13},z_{23})+(m_{11},z_{21}))-(-(m_{23},z_{13})+(m_{21},z_{11}))$

$\ +(\tau _{12}\mu _{22}+(t_{13},z_{23})+(t_{11},z_{21}))-(\tau _{22}\mu
_{12}+(t_{23},z_{13})+(t_{21},z_{11}))$

$\ -\frac{1}{3}\rho _{1}\mu _{22}+\frac{1%
}{3}\rho _{2}\mu _{12}$

$\ +((\beta _{13}\psi _{21}+\beta _{11}\psi _{23})-2(b_{12},w_{22}))-((\beta
_{23}\psi _{11}+\beta _{21}\psi _{13})-2(b_{22},w_{12}))$

$\ +\omega _{2}\alpha _{12}-\omega _{1}\alpha _{22}-r_{1}\mu _{22}+r_{2}\mu
_{12}+u_{1}\chi _{22}-u_{2}\chi _{12},$

$([R_{1},R_{2}]_{8Q})_{Z%
\_E_{3}}=(-(m_{11},z_{21})+(m_{12},z_{22}))-(-(m_{21},z_{11})+(m_{22},z_{12}))$

$\ +((-\tau _{11}-\tau _{12})\mu
_{23}+(t_{11},z_{21})+(t_{12},z_{22}))$

$\ -((-\tau _{21}-\tau _{22})\mu
_{13}+(t_{21},z_{11})+(t_{22},z_{12}))$

$\ -\frac{1}{3}\rho _{1}\mu _{23}+\frac{1}{3}\rho _{2}\mu _{13}$

$\ +((\beta _{11}\psi _{22}+\beta _{12}\psi _{21})-2(b_{13},w_{23}))-((\beta
_{21}\psi _{12}+\beta _{22}\psi _{11})-2(b_{23},w_{13}))$

$\ +\omega _{2}\alpha _{13}-\omega _{1}\alpha _{23}-r_{1}\mu _{23}+r_{2}\mu
_{13}+u_{1}\chi _{23}-u_{2}\chi _{13},$

$([R_{1},R_{2}]_{8Q})_{Z%
\_F_{1}}=\textrm{d}_{g}\textrm{g}_{d}(D_{1})z_{21}-\textrm{d}_{g}\textrm{g}_{d}(D_{2})z_{11} $

$\ +\frac{1}{2}(\overline{m_{12}z_{23}}-\overline{z_{22}m_{13}}+(-\mu
_{22}+\mu _{23})m_{11})$

$\ -\frac{1}{2}(\overline{m_{22}z_{13}}-\overline{%
z_{12}m_{23}}+(-\mu _{12}+\mu _{13})m_{21})$

$\ +\frac{1}{2}(-\tau _{11}z_{21}+(\mu _{22}+\mu _{23})t_{11}+\overline{%
z_{22}t_{13}}+\overline{t_{12}z_{23}})$

$\ -\frac{1}{2}(-\tau _{21}z_{11}+(\mu _{12}+\mu _{13})t_{21}+\overline{%
z_{12}t_{23}}+\overline{t_{22}z_{13}})-\frac{1}{3}\rho _{1}z_{21}+\frac{1}{3}%
\rho _{2}z_{11}$

$\ +(-\psi _{21}b_{11}-\beta _{11}w_{21}+\overline{b_{12}w_{23}}+\overline{%
w_{22}b_{13}})\ -(-\psi _{11}b_{21}-\beta _{21}w_{11}+\overline{b_{22}w_{13}}%
+\overline{w_{12}b_{23}})$

$\ +\omega _{2}a_{11}-\omega
_{1}a_{21}-r_{1}z_{21}+r_{2}z_{11}+u_{1}x_{21}-u_{2}x_{11},$

$([R_{1},R_{2}]_{8Q})_{Z\_F_{2}}=\textrm{d}_{g}\nu \textrm{g}_{d}(D_{1})z_{22}-\textrm{d}_{g}\nu
\textrm{g}_{d}(D_{2})z_{12}$

$\ +\frac{1}{2}(\overline{m_{13}z_{21}}-\overline{z_{23}m_{11}}+(-\mu
_{23}+\mu _{21})m_{12})$

$\ -\frac{1}{2}(\overline{m_{23}z_{11}}-\overline{%
z_{13}m_{21}}+(-\mu _{13}+\mu _{11})m_{22})$

$\ +\frac{1}{2}(-\tau _{12}z_{22}+(\mu _{23}+\mu _{21})t_{12}+\overline{%
z_{23}t_{11}}+\overline{t_{13}z_{21}})$

$\ -\frac{1}{2}(-\tau _{22}z_{12}+(\mu _{13}+\mu _{11})t_{22}+\overline{%
z_{13}t_{21}}+\overline{t_{23}z_{11}})-\frac{1}{3}\rho _{1}z_{22}+\frac{1}{3}%
\rho _{2}z_{12}$

$\ +(-\psi _{22}b_{12}-\beta _{12}w_{22}+\overline{b_{13}w_{21}}+\overline{%
w_{23}b_{11}})-(-\psi _{12}b_{22}-\beta _{22}w_{12}+\overline{b_{23}w_{11}}+%
\overline{w_{13}b_{21}})$

$\ +\omega _{2}a_{12}-\omega
_{1}a_{22}-r_{1}z_{22}+r_{2}z_{12}+u_{1}x_{22}-u_{2}x_{12},$

$([R_{1},R_{2}]_{8Q})_{Z\_F_{3}}=\textrm{d}_{g}\nu ^{2}\textrm{g}_{d}(D_{1})z_{23}-\textrm{d}_{g}\nu
^{2}\textrm{g}_{d}(D_{2})z_{13}$

$\ +\frac{1}{2}(\overline{m_{11}z_{22}}-\overline{z_{21}m_{12}}+(-\mu
_{21}+\mu _{22})m_{13})$

$\ -\frac{1}{2}(\overline{m_{21}z_{12}}-\overline{%
z_{11}m_{22}}+(-\mu _{11}+\mu _{12})m_{23})$

$\ +\frac{1}{2}((\tau _{11}+\tau _{12})z_{23}+(\mu _{21}+\mu _{22})t_{13}+%
\overline{z_{21}t_{12}}+\overline{t_{11}z_{22}})$

$\ -\frac{1}{2}((\tau _{21}+\tau _{22})z_{13}+(\mu _{11}+\mu _{12})t_{23}+%
\overline{z_{11}t_{22}}+\overline{t_{21}z_{12}})-\frac{1}{3}\rho _{1}z_{23}+%
\frac{1}{3}\rho _{2}z_{13}$

$\ +(-\psi _{23}b_{13}-\beta _{13}w_{23}+\overline{b_{11}w_{22}}+\overline{%
w_{21}b_{12}})-(-\psi _{13}b_{23}-\beta _{23}w_{13}+\overline{b_{21}w_{12}}+%
\overline{w_{11}b_{22}})$

$\ +\omega _{2}a_{13}-\omega
_{1}a_{23}-r_{1}z_{23}+r_{2}z_{13}+u_{1}x_{23}-u_{2}x_{13},$

$([R_{1},R_{2}]_{8Q})_{W}=\left( 
\begin{array}{c}
([R_{1},R_{2}]_{8Q})_{W\_E_{1}} \\ 
([R_{1},R_{2}]_{8Q})_{W\_E_{2}} \\ 
([R_{1},R_{2}]_{8Q})_{W\_E_{3}} \\ 
([R_{1},R_{2}]_{8Q})_{W\_F_{1}} \\ 
([R_{1},R_{2}]_{8Q})_{W\_F_{2}} \\ 
([R_{1},R_{2}]_{8Q})_{W\_F_{3}}%
\end{array}%
\right) ,$

$([R_{1},R_{2}]_{8Q})_{W%
\_E_{1}}=(-(m_{12},w_{22})+(m_{13},w_{23}))-(-(m_{22},w_{12})+(m_{23},w_{13}))$

$\ -(\tau _{11}\psi _{21}+(t_{12},w_{22})+(t_{13},w_{23}))+(\tau _{21}\psi
_{11}+(t_{22},w_{12})+(t_{23},w_{13}))$

$\ +\frac{1}{3}\rho _{1}\psi _{21}-\frac{1%
}{3}\rho _{2}\psi _{11}$

$\ +((\alpha _{12}\mu _{23}+\alpha _{13}\mu
_{22})-2(a_{11},z_{21}))-((\alpha _{22}\mu _{13}+\alpha _{23}\mu
_{12})-2(a_{21},z_{11}))$

$\ +\zeta _{2}\beta _{11}-\zeta _{1}\beta _{21}-r_{1}\psi _{21}+r_{2}\psi
_{11}+u_{1}\gamma _{21}-u_{2}\gamma _{11},$

$([R_{1},R_{2}]_{8Q})_{W%
\_E_{2}}=(-(m_{13},w_{23})+(m_{11},w_{21}))-(-(m_{23},w_{13})+(m_{21},w_{11}))$

$\ -(\tau _{12}\psi _{22}+(t_{13},w_{23})+(t_{11},w_{21}))+(\tau _{22}\psi
_{12}+(t_{23},w_{13})+(t_{21},w_{11}))$

$\ +\frac{1}{3}\rho _{1}\psi _{22}-\frac{1%
}{3}\rho _{2}\psi _{12}$

$\ +((\alpha _{13}\mu _{21}+\alpha _{11}\mu
_{23})-2(a_{12},z_{22}))-((\alpha _{23}\mu _{11}+\alpha _{21}\mu
_{13})-2(a_{22},z_{12}))$

$\ +\zeta _{2}\beta _{12}-\zeta _{1}\beta _{22}-r_{1}\psi _{22}+r_{2}\psi
_{12}+u_{1}\gamma _{22}-u_{2}\gamma _{12},$

$([R_{1},R_{2}]_{8Q})_{W%
\_E_{3}}=(-(m_{11},w_{21})+(m_{12},w_{22}))-(-(m_{21},w_{11})+(m_{22},w_{12}))$

$\ -((-\tau _{11}-\tau _{12})\psi
_{23}+(t_{11},w_{21})+(t_{12},w_{22}))$

$\ +((-\tau _{21}-\tau _{22})\psi
_{13}+(t_{21},w_{11})+(t_{22},w_{12}))$

$\ +\frac{1}{3}\rho _{1}\psi _{23}-\frac{1}{3}\rho _{2}\psi _{13}$

$\ +((\alpha _{11}\mu _{22}+\alpha _{12}\mu
_{21})-2(a_{13},z_{23}))-((\alpha _{21}\mu _{12}+\alpha _{22}\mu
_{11})-2(a_{23},z_{13}))$

$\ +\zeta _{2}\beta _{13}-\zeta _{1}\beta _{23}-r_{1}\psi _{23}+r_{2}\psi
_{13}+u_{1}\gamma _{23}-u_{2}\gamma _{13},$

$([R_{1},R_{2}]_{8Q})_{W%
\_F_{1}}=\textrm{d}_{g}\textrm{g}_{d}(D_{1})w_{21}-\textrm{d}_{g}\textrm{g}_{d}(D_{2})w_{11} $

$\ +\frac{1}{2}(\overline{m_{12}w_{23}}-\overline{w_{22}m_{13}}+(-\psi
_{22}+\psi _{23})m_{11})$

$\ -\frac{1}{2}(\overline{m_{22}w_{13}}-\overline{%
w_{12}m_{23}}+(-\psi _{12}+\psi _{13})m_{21})$

$\ -\frac{1}{2}(-\tau _{11}w_{21}+(\psi _{22}+\psi _{23})t_{11}+\overline{%
w_{22}t_{13}}+\overline{t_{12}w_{23}})$

$\ +\frac{1}{2}(-\tau _{21}w_{11}+(\psi _{12}+\psi _{13})t_{21}+\overline{%
w_{12}t_{23}}+\overline{t_{22}w_{13}})+\frac{1}{3}\rho _{1}w_{21}-\frac{1}{3}%
\rho _{2}w_{11}$

$\ +(-\mu _{21}a_{11}-\alpha _{11}z_{21}+\overline{a_{12}z_{23}}+\overline{%
z_{22}a_{13}})-(-\mu _{11}a_{21}-\alpha _{21}z_{11}+\overline{a_{22}z_{13}}+%
\overline{z_{12}a_{23}})$

$\ +\zeta _{2}b_{11}-\zeta
_{1}b_{21}-r_{1}w_{21}+r_{2}w_{11}+u_{1}y_{21}-u_{2}y_{11},$

$([R_{1},R_{2}]_{8Q})_{W\_F_{2}}=\textrm{d}_{g}\nu \textrm{g}_{d}(D_{1})w_{22}-\textrm{d}_{g}\nu
\textrm{g}_{d}(D_{2})w_{12}$

$\ +\frac{1}{2}(\overline{m_{13}w_{21}}-\overline{w_{23}m_{11}}+(-\psi
_{23}+\psi _{21})m_{12})$

$\ -\frac{1}{2}(\overline{m_{23}w_{11}}-\overline{%
w_{13}m_{21}}+(-\psi _{13}+\psi _{11})m_{22})$

$\ -\frac{1}{2}(-\tau _{12}w_{22}+(\psi _{23}+\psi _{21})t_{12}+\overline{%
w_{23}t_{11}}+\overline{t_{13}w_{21}})$

$\ +\frac{1}{2}(-\tau _{22}w_{12}+(\psi _{13}+\psi _{11})t_{22}+\overline{%
w_{13}t_{21}}+\overline{t_{23}w_{11}})+\frac{1}{3}\rho _{1}w_{22}-\frac{1}{3}%
\rho _{2}w_{12}$

$\ +(-\mu _{22}a_{12}-\alpha _{12}z_{22}+\overline{a_{13}z_{21}}+\overline{%
z_{23}a_{11}})-(-\mu _{12}a_{22}-\alpha _{22}z_{12}+\overline{a_{23}z_{11}}+%
\overline{z_{13}a_{21}})$

$\ +\zeta _{2}b_{12}-\zeta
_{1}b_{22}-r_{1}w_{22}+r_{2}w_{12}+u_{1}y_{22}-u_{2}y_{12},$

$([R_{1},R_{2}]_{8Q})_{W\_F_{3}}=\textrm{d}_{g}\nu ^{2}\textrm{g}_{d}(D_{1})w_{23}-\textrm{d}_{g}\nu
^{2}\textrm{g}_{d}(D_{2})w_{13}$

$\ +\frac{1}{2}(\overline{m_{11}w_{22}}-\overline{w_{21}m_{12}}+(-\psi
_{21}+\psi _{22})m_{13})$

$\ -\frac{1}{2}(\overline{m_{21}w_{12}}-\overline{%
w_{11}m_{22}}+(-\psi _{11}+\psi _{12})m_{23})$

$\ -\frac{1}{2}((\tau _{11}+\tau _{12})w_{23}+(\psi _{21}+\psi _{22})t_{13}+%
\overline{w_{21}t_{12}}+\overline{t_{11}w_{22}})$

$\ +\frac{1}{2}((\tau _{21}+\tau _{22})w_{13}+(\psi _{11}+\psi _{12})t_{23}+%
\overline{w_{11}t_{22}}+\overline{t_{21}w_{12}})+\frac{1}{3}\rho _{1}w_{23}-%
\frac{1}{3}\rho _{2}w_{13}$

$\ +(-\mu _{23}a_{13}-\alpha _{13}z_{23}+\overline{a_{11}z_{22}}+\overline{%
z_{21}a_{12}})-(-\mu _{13}a_{23}-\alpha _{23}z_{13}+\overline{a_{21}z_{12}}+%
\overline{z_{11}a_{22}})$

$\ +\zeta _{2}b_{13}-\zeta
_{1}b_{23}-r_{1}w_{23}+r_{2}w_{13}+u_{1}y_{23}-u_{2}y_{13},$

$([R_{1},R_{2}]_{8Q})_{\zeta }=\sum_{i=1}^{3}(\alpha _{1i}\psi
_{2i}+2(a_{1i},w_{2i}))-\sum_{i=1}^{3}(\alpha _{2i}\psi
_{1i}+2(a_{2i},w_{1i}))$

$\ \ \ \ \ \ \ \ \ \ \ \ \ \ \ \ \ \ \ \ \ +\rho _{1}\xi _{2}-\rho _{2}%
\xi _{1}-r_{1}\zeta _{2}+r_{2}\zeta _{1}+u_{1}\xi _{2}-u_{2}\xi _{1},$

$([R_{1},R_{2}]_{8Q})_{\omega }=\sum_{i=1}^{3}(\beta _{1i}\mu
_{2i}+2(b_{1i},z_{2i}))-\sum_{i=1}^{3}(\beta _{2i}\mu
_{1i}+2(b_{2i},z_{1i})) $

$\ \ \ \ \ \ \ \ \ \ \ \ \ \ \ \ \ \ \ \ \ -\rho _{1}\omega _{2}+\rho
_{2}\omega _{1}-r_{1}\omega _{2}+r_{2}\omega _{1}+u_{1}\eta _{2}-u_{2}\eta
_{1}.$

\bigskip

\emph{\ Proof. }Using \emph{Lemma 5.9} and definition of mappings $%
\textrm{f}_{DDD},\textrm{f}_{JJD},\textrm{f}_{DJJ},$

\noindent
$\textrm{f}_{iJJ},\textrm{f}_{JJC},\textrm{f}_{CJJ},\textrm{f}_{CCC},$
we have the above expression. \ \ \ \ \emph{Q.E.D.}

\bigskip

\emph{Lemma 10.4.}

$\left( 
\begin{array}{cccccc}
0 & F_{rP} & F_{rQ} & 0 & F_{rs} & F_{ru}%
\end{array}%
\right) \left( 
\begin{array}{c}
\Phi _{2} \\ 
P_{2} \\ 
Q_{2} \\ 
r_{2} \\ 
s_{2} \\ 
u_{2}%
\end{array}%
\right) =\left( [R_{1},R_{2}]_{8r}\right) ,$

\noindent
where

$F_{rP}=\left( 
\begin{array}{cccccc}
F_{rP11} & F_{rP12} & F_{rP13} & F_{rP14} & F_{rP15} & F_{rP16}%
\end{array}%
\right) ,$

$F_{rP11}=\left( 
\begin{array}{ccc}
\frac{1}{8}\textrm{f}_{CCC}(\psi _{11},?) & \frac{1}{8}\textrm{f}_{CCC}(\psi _{12},?) & \frac{1%
}{8}\textrm{f}_{CCC}(\psi _{13},?)%
\end{array}%
\right) ,$

$F_{rP12}=\left( 
\begin{array}{ccc}
\frac{1}{4}\textrm{f}_{JJC}(?,w_{11}) & \frac{1}{4}\textrm{f}_{JJC}(?,w_{12}) & \frac{1}{4}%
\textrm{f}_{JJC}(?,w_{13})%
\end{array}%
\right) ,$

$F_{rP13}=\left( 
\begin{array}{ccc}
-\frac{1}{8}\textrm{f}_{CCC}(\mu _{11},?) & -\frac{1}{8}\textrm{f}_{CCC}(\mu _{12},?) & -\frac{%
1}{8}\textrm{f}_{CCC}(\mu _{13},?)%
\end{array}%
\right) ,$

$F_{rP14}=\left( 
\begin{array}{ccc}
-\frac{1}{4}\textrm{f}_{JJC}(?,z_{11}) & -\frac{1}{4}\textrm{f}_{JJC}(?,z_{12}) & -\frac{1}{4}%
\textrm{f}_{JJC}(?,z_{13})%
\end{array}%
\right) ,$

$F_{rP15}=\left( \frac{1}{8}\textrm{f}_{CCC}(\omega _{1},?)\right) ,$ \ \ \ \ $%
F_{rP16}=\left( -\frac{1}{8}\textrm{f}_{CCC}(\zeta _{1},?)\right) ,$

\bigskip

$F_{rQ}=\left( 
\begin{array}{cccccc}
F_{rQ11} & F_{rQ12} & F_{rQ13} & F_{rQ14} & F_{rQ15} & F_{rQ16}%
\end{array}%
\right) ,$

$F_{rQ11}=\left( 
\begin{array}{ccc}
\frac{1}{8}\textrm{f}_{CCC}(\gamma _{11},?) & \frac{1}{8}\textrm{f}_{CCC}(\gamma _{12},?) & 
\frac{1}{8}\textrm{f}_{CCC}(\gamma _{13},?)%
\end{array}%
\right) ,$

$F_{rQ12}=\left( 
\begin{array}{ccc}
\frac{1}{4}\textrm{f}_{JJC}(?,y_{11}) & \frac{1}{4}\textrm{f}_{JJC}(?,y_{12}) & \frac{1}{4}%
\textrm{f}_{JJC}(?,y_{13})%
\end{array}%
\right) ,$

$F_{rQ13}=\left( 
\begin{array}{ccc}
-\frac{1}{8}\textrm{f}_{CCC}(\chi _{11},?) & -\frac{1}{8}\textrm{f}_{CCC}(\chi _{12},?) & -%
\frac{1}{8}\textrm{f}_{CCC}(\chi _{13},?)%
\end{array}%
\right) ,$

$F_{rQ14}=\left( 
\begin{array}{ccc}
-\frac{1}{4}\textrm{f}_{JJC}(?,x_{11}) & -\frac{1}{4}\textrm{f}_{JJC}(?,x_{12}) & -\frac{1}{4}%
\textrm{f}_{JJC}(?,x_{13})%
\end{array}%
\right) ,$

$F_{rQ15}=\left( \frac{1}{8}\textrm{f}_{CCC}(\eta _{1},?)\right) ,$ \ \ \ \ \ $%
F_{rQ16}=\left( -\frac{1}{8}\textrm{f}_{CCC}(\xi _{1},?)\right) ,$

\bigskip

$F_{rs}=\left( -\textrm{f}_{CCC}(u_{1},?)\right) ,$

$F_{ru}=\left( \textrm{f}_{CCC}(s_{1},?)\right) ,$

\bigskip

$[R_{1},R_{2}]_{8r}$

=$-\frac{1}{8}(\sum_{i=1}^{3}(\chi _{1i}\psi
_{2i}+2(x_{1i},w_{2i}))-\sum_{i=1}^{3}(\gamma _{1i}\mu
_{2i}+2(y_{1i},z_{2i}))+\xi _{1}\omega _{2}-\zeta _{2}\eta _{1})$

$\ \ \ +\frac{1}{8}(\sum_{i=1}^{3}(\psi _{1i}\chi
_{2i}+2(w_{1i},x_{2i}))-\sum_{i=1}^{3}(\mu _{1i}\gamma
_{2i}+2(z_{1i},y_{2i}))+\xi _{2}\omega _{1}-\zeta _{1}\eta_{2})$

\ \ \ $+s_{1}u_{2}-s_{2}u_{1}.$

\bigskip

\emph{\ Proof. }Using \emph{Lemma 5.9} and definition of mappings $%
\textrm{f}_{DDD},\textrm{f}_{JJD},\textrm{f}_{DJJ},$

\noindent
$\textrm{f}_{iJJ},\textrm{f}_{JJC},\textrm{f}_{CJJ},\textrm{f}_{CCC},$
we have the above expression. \ \ \ \ \emph{Q.E.D.}

\bigskip

\emph{Lemma 10.5.}

$\left( 
\begin{array}{cccccc}
0 & F_{sP} & 0 & F_{sr} & F_{ss} & 0%
\end{array}%
\right) \left( 
\begin{array}{c}
\Phi _{2} \\ 
P_{2} \\ 
Q_{2} \\ 
r_{2} \\ 
s_{2} \\ 
u_{2}%
\end{array}%
\right) =\left( [R_{1},R_{2}]_{8s}\right) ,$

\noindent
where

$F_{sP}=\left( 
\begin{array}{cccccc}
F_{sP11} & F_{sP12} & F_{sP13} & F_{sP14} & F_{sP15} & F_{sP16}%
\end{array}%
\right) ,$

$F_{sP11}=\left( 
\begin{array}{ccc}
-\frac{1}{4}\textrm{f}_{CCC}(\gamma _{11},?) & -\frac{1}{4}\textrm{f}_{CCC}(\gamma _{12},?) & -%
\frac{1}{4}\textrm{f}_{CCC}(\gamma _{13},?)%
\end{array}%
\right) ,$

$F_{sP12}=\left( 
\begin{array}{ccc}
-\frac{1}{2}\textrm{f}_{JJC}(?,y_{11}) & -\frac{1}{2}\textrm{f}_{JJC}(?,y_{12}) & -\frac{1}{2}%
\textrm{f}_{JJC}(?,y_{13})%
\end{array}%
\right) ,$

$F_{sP13}=\left( 
\begin{array}{ccc}
\frac{1}{4}\textrm{f}_{CCC}(\chi _{11},?) & \frac{1}{4}\textrm{f}_{CCC}(\chi _{12},?) & \frac{1%
}{4}\textrm{f}_{CCC}(\chi _{13},?)%
\end{array}%
\right) ,$

$F_{sP14}=\left( 
\begin{array}{ccc}
\frac{1}{2}\textrm{f}_{JJC}(?,x_{11}) & \frac{1}{2}\textrm{f}_{JJC}(?,x_{12}) & \frac{1}{2}%
\textrm{f}_{JJC}(?,x_{13})%
\end{array}%
\right) ,$

$F_{sP15}=\left( -\frac{1}{4}\textrm{f}_{CCC}(\eta _{1},?)\right) ,$ \ \ \ \ \ $%
F_{sP16}=\left( \frac{1}{4}\textrm{f}_{CCC}(\xi _{1},?)\right) ,$

\bigskip

$F_{sr}=\left( -2\textrm{f}_{CCC}(s_{1},?)\right) ,$

$F_{ss}=\left( 2\textrm{f}_{CCC}(r_{1},?)\right) ,$

\bigskip

$[R_{1},R_{2}]_{8s}$

=$\frac{1}{4}(\sum_{i=1}^{3}(\chi _{1i}\gamma
_{2i}+2(x_{1i},y_{2i}))-\sum_{i=1}^{3}(\gamma _{1i}\chi
_{2i}+2(y_{1i},x_{2i}))+\xi _{1}\eta _{2}-\xi _{2}\eta_{1})$

\ \ $+2r_{1}s_{2}-2r_{2}s_{1}.$

\bigskip

\emph{\ Proof. }Using \emph{Lemma 5.9} and definition of mappings $%
\textrm{f}_{DDD},\textrm{f}_{JJD},\textrm{f}_{DJJ},$

\noindent
$\textrm{f}_{iJJ},\textrm{f}_{JJC},\textrm{f}_{CJJ},\textrm{f}_{CCC},$
we have the above expression. \ \ \ \ \emph{Q.E.D.}

\bigskip

\emph{Lemma 10.6.}

$\left( 
\begin{array}{cccccc}
0 & 0 & F_{uQ} & F_{ur} & 0 & F_{uu}%
\end{array}%
\right) \left( 
\begin{array}{c}
\Phi _{2} \\ 
P_{2} \\ 
Q_{2} \\ 
r_{2} \\ 
s_{2} \\ 
u_{2}%
\end{array}%
\right) =\left( [R_{1},R_{2}]_{8u}\right) ,$

\noindent
where

$F_{uQ}=\left( 
\begin{array}{cccccc}
F_{uQ11} & F_{uQ12} & F_{uQ13} & F_{uQ14} & F_{uQ15} & F_{uQ16}%
\end{array}%
\right) ,$

$F_{uQ11}=\left( 
\begin{array}{ccc}
\frac{1}{4}\textrm{f}_{CCC}(\psi _{11},?) & \frac{1}{4}\textrm{f}_{CCC}(\psi _{12},?) & \frac{1%
}{4}\textrm{f}_{CCC}(\psi _{13},?)%
\end{array}%
\right) ,$

$F_{uQ12}=\left( 
\begin{array}{ccc}
\frac{1}{2}\textrm{f}_{JJC}(?,w_{11}) & \frac{1}{2}\textrm{f}_{JJC}(?,w_{12}) & \frac{1}{2}%
\textrm{f}_{JJC}(?,w_{13})%
\end{array}%
\right) ,$

$F_{uQ13}=\left( 
\begin{array}{ccc}
-\frac{1}{4}\textrm{f}_{CCC}(\mu _{11},?) & -\frac{1}{4}\textrm{f}_{CCC}(\mu _{12},?) & -\frac{%
1}{4}\textrm{f}_{CCC}(\mu _{13},?)%
\end{array}%
\right) ,$

$F_{uQ14}=\left( 
\begin{array}{ccc}
-\frac{1}{2}\textrm{f}_{JJC}(?,z_{11}) & -\frac{1}{2}\textrm{f}_{JJC}(?,z_{12}) & -\frac{1}{2}%
\textrm{f}_{JJC}(?,z_{13})%
\end{array}%
\right) ,$

$F_{uQ15}=\left( \frac{1}{4}\textrm{f}_{CCC}(\omega _{1},?)\right) ,$ \ \ \ \ \ \ \ $%
F_{uQ16}=\left( -\frac{1}{4}\textrm{f}_{CCC}(\zeta _{1},?)\right) ,$

\bigskip

$F_{ur}=\left( 2\textrm{f}_{CCC}(u_{1},?)\right) ,$

$F_{uu}=\left( -2\textrm{f}_{CCC}(r_{1},?)\right) ,$

\bigskip

$[R_{1},R_{2}]_{8u}$

=$-\frac{1}{4}(\sum_{i=1}^{3}(\mu _{1i}\psi
_{2i}+2(z_{1i},w_{2i}))-\sum_{i=1}^{3}(\psi _{1i}\mu
_{2i}+2(w_{1i},z_{2i}))+\zeta _{1}\omega _{2}-\zeta _{2}\omega_{1})$

$\ \ -2r_{1}u_{2}+2r_{2}u_{1}$.

\bigskip

\emph{\ Proof. }Using \emph{Lemma 5.9} and definition of mappings $%
\textrm{f}_{DDD},\textrm{f}_{JJD},\textrm{f}_{DJJ},$

\noindent
$\ \textrm{f}_{iJJ},\textrm{f}_{JJC},\textrm{f}_{CJJ},\textrm{f}_{CCC},$
we have the above expression. \ \ \ \ \emph{Q.E.D.}

\bigskip

\emph{Theorem 10.7.} The image of the adjoit representation $ad($\gR$_{8})$ of \gR$_{8}^{}$
is expressed by

$ad(\Phi,P,Q,r,s,u )=ad(\textrm{fv}(\Phi,P,Q,r,s,u ))$

$\ \ \ \ \ \ \ \ =%
%
,$

$M_{ur}=2uE$

$M_{uu}=-2rE$

Since \gR$_{8}^{}$ and \ge$_{8,1}^{}$ are isomorphic, 
$ad(\Phi ,P,Q,r,s,u)$ is also the image of the adjoit representation of \ge$_{8,1}^{}.$
\ Let we put the $248 \times 248$ matrix images of $ad($\gR$_{8}^{\C})$ over comples numbers  by \gr$_{8}^{\C}$.

\bigskip

\emph{Proof. }Using \emph{Lemma
10.1, 10.2, 10.3, 10.4, 10.5, 10.6, Lemma 6.12,}

\noindent
\emph{6.11, 6.12, 6.13, 6.14, 6.15,}
\emph{6.16, 6.17, 6.18, 6.19, 6.20, }and 
\emph{Lemma 6.24}, we have the
following expression.

$M_{\Phi P}=\left( 
%
%
\right) ,$

$M_{\Phi P12}=\frac{1}{2}\textrm{M}_{R}\textrm{Dl}(W),$ 

$M_{\Phi P14}=\frac{1}{2}\textrm{M}_{R}\textrm{Dl}(Z)$

$M_{\Phi P21}=-\frac{1}{4}^{t}\textrm{MJC}(W),$ 

$M_{\Phi P22}=-\frac{1}{4}\textrm{M}_{D}^{-}\textrm{E}(W)-\frac{1}{4}\textrm{MI}(W)$

$M_{\Phi P23}=\frac{1}{4}^{t}\textrm{MJC}(Z),$ 

$M_{\Phi P24}=\frac{1}{4}\textrm{M}_{D}^{-}\textrm{E}(Z)+\frac{1}{4}\textrm{MI}(Z)$

$M_{\Phi P31}=\frac{1}{3}M_{D3}(W)-\frac{1}{6}\textrm{MCCC}_{3}(W),$

$M_{\Phi P32}=-\frac{1}{3}\textrm{M}_{D}\textrm{C}_{3}(W)+\frac{1}{6}\textrm{M}^{+}\textrm{JC}_{3}(W)$

$M_{\Phi P33}=\frac{1}{3}M_{D3}(Z)-\frac{1}{6}\textrm{MCCC}_{3}(Z),$ 

$M_{\Phi P34}=-\frac{1}{3}\textrm{M}_{D}\textrm{C}_{3}(Z)+\frac{1}{6}\textrm{M}^{+}\textrm{JC}_{3}(Z)$

$M_{\Phi P41}=\frac{1}{4}^{t}\textrm{M}^{+}\textrm{JC}(W),$ 

$M_{\Phi P42}=\frac{1}{4}\textrm{M}_{D}^{+}\textrm{E}(W)+\frac{1}{4}\textrm{M}^{+}\textrm{I}(W)$

$M_{\Phi P43}=\frac{1}{4}^{t}\textrm{M}^{+}\textrm{JC}(Z),$ 

$M_{\Phi P44}=\frac{1}{4}\textrm{M}_{D}^{+}\textrm{E}(Z)+\frac{1}{4}\textrm{M}^{+}\textrm{I}(Z)$

$M_{\Phi P51}=-\zeta E,$ 

$M_{\Phi P53}=\frac{1}{4}\textrm{MCC}(W),$ 

$M_{\Phi P54}=-\frac{1}{2}\textrm{M}_{D}\textrm{C}(W)$

$M_{\Phi P55}=-\frac{1}{4}\textrm{M}_{C}\textrm{CCC}(Z)$

$M_{\Phi P62}=-\frac{1}{4}\zeta E,$ 

$M_{\Phi P63}=-\frac{1}{4}^{t}\textrm{M}_{D}\textrm{C}(W),$ 

$M_{\Phi P64}=-\frac{1}{4}\textrm{M}_{D}\textrm{E}(W)+\frac{1}{4}\textrm{M}^{+}\textrm{I}(W)$

$M_{\Phi P65}=-\frac{1}{4}\textrm{M}_{C}\textrm{CJ}(Z)$

$M_{\Phi P71}=-\frac{1}{4}\textrm{MCC}(Z),$ 

$M_{\Phi P72}=\frac{1}{2}\textrm{M}_{D}\textrm{C}(Z),$ 

$M_{\Phi P73}=\frac{1}{4}\omega E,$ 

$M_{\Phi P76}=\frac{1}{4}\textrm{M}_{C}\textrm{CCC}(W)$

$M_{\Phi P81}=\frac{1}{4}^{t}\textrm{M}_{D}\textrm{C}(Z),$ 

$M_{\Phi P82}=\frac{1}{4}\textrm{M}_{D}\textrm{E}(Z)-\frac{1}{4}\textrm{M}^{+}\textrm{I}(Z),$ 

$M_{\Phi P84}=\frac{1}{4}\omega E,$

$M_{\Phi P86}=\frac{1}{4}\textrm{M}_{C}\textrm{CJ}(W)$

$M_{\Phi P91}=-\frac{1}{8}^{t}\textrm{M}_{C}\textrm{CCC}(W),$ 

$M_{\Phi P92}=-\frac{1}{4}^{t}\textrm{M}_{C}\textrm{CJ}(W),$ 

$M_{\Phi P93}=-\frac{1}{8}^{t}\textrm{M}_{C}\textrm{CCC}(Z)$

$M_{\Phi P94}=-\frac{1}{4}^{t}\textrm{M}_{C}\textrm{CJ}(Z),$ 

$M_{\Phi P95}=\frac{3}{8}\omega ,$ 

$M_{\Phi P96}=\frac{3}{8}\zeta $

\bigskip

$M_{\Phi Q}=\left( 
\begin{array}{cccccc}
0 & M_{\Phi Q12} & 0 & M_{\Phi Q14} & 0 & 0 \\ 
M_{\Phi Q21} & M_{\Phi Q22} & M_{\Phi Q23} & M_{\Phi Q24} & 0 & 0
\\ 
M_{\Phi Q31} & M_{\Phi Q32} & M_{\Phi Q33} & M_{\Phi Q34} & 0 & 0
\\ 
M_{\Phi Q41} & M_{\Phi Q42} & M_{\Phi Q43} & M_{\Phi Q44} & 0 & 0
\\ 
M_{\Phi Q51} & 0 & M_{\Phi Q53} & M_{\Phi Q54} & M_{\Phi Q55} & 0
\\ 
0 & M_{\Phi Q62} & M_{\Phi Q63} & M_{\Phi Q64} & M_{\Phi Q65} & 0
\\ 
M_{\Phi Q71} & M_{\Phi Q72} & M_{\Phi Q73} & 0 & 0 & M_{\Phi Q76}
\\ 
M_{\Phi Q81} & M_{\Phi Q82} & 0 & M_{\Phi Q84} & 0 & M_{\Phi Q86}
\\ 
M_{\Phi Q91} & M_{\Phi Q92} & M_{\Phi Q93} & M_{\Phi Q94} & M_{%
\Phi Q95} & M_{\Phi Q96}%
\end{array}%
\right) ,$

$M_{\Phi Q12}=-\frac{1}{2}\textrm{M}_{R}\textrm{Dl}(Y),$ \ \ 

$M_{\Phi Q14}=-\frac{1}{2}\textrm{M}_{R}\textrm{Dl}(X),$

$M_{\Phi Q21}=\frac{1}{4}^{t}\textrm{MJC}(Y),$ \ \ 

$M_{\Phi Q22}=\frac{1}{4}\textrm{M}_{D}^{-}\textrm{E}(Y)+\frac{1}{4}\textrm{MI}(Y),$

$M_{\Phi Q23}=\frac{1}{4}^{t}\textrm{MJC}(X),$ \ \ 

$M_{\Phi Q24}=-\frac{1}{4}\textrm{M}_{D}^{-}\textrm{E}(X)-\frac{1}{4}\textrm{MI}(X),$

$M_{\Phi Q31}=-\frac{1}{3}M_{D3}(Y)+\frac{1}{6}\textrm{MCCC}_{3}(Y),$ \ \ 

$M_{\Phi Q32}=\frac{1}{3}\textrm{M}_{D}\textrm{C}_{3}(Y)-\frac{1}{6}\textrm{M}^{+}\textrm{JC}_{3}(Y),$

$M_{\Phi Q33}=-\frac{1}{3}M_{D3}(X)+\frac{1}{6}\textrm{MCCC}_{3}(X),$ \ \ 

$M_{\Phi Q34}=\frac{1}{3}\textrm{M}_{D}\textrm{C}_{3}(X)-\frac{1}{6}\textrm{M}^{+}\textrm{JC}_{3}(X),$

$M_{\Phi Q41}=-\frac{1}{4}^{t}\textrm{M}^{+}\textrm{JC}(Y),$ \ \ 

$M_{\Phi Q42}=-\frac{1}{4}\textrm{M}_{D}^{+}\textrm{E}(Y)-\frac{1}{4}\textrm{M}^{+}\textrm{I}(Y),$

$M_{\Phi Q43}=-\frac{1}{4}^{t}\textrm{M}^{+}\textrm{JC}(X),$ \ \ 

$M_{\Phi Q44}=-\frac{1}{4}\textrm{M}_{D}^{+}\textrm{E}(X)-\frac{1}{4}\textrm{M}^{+}\textrm{I}(X),$

$M_{\Phi Q51}=\xi E,$ \ \ 

$M_{\Phi Q53}=-\frac{1}{4}\textrm{MCC}(Y),$ \ \ 
$M_{\Phi Q54}=\frac{1}{2}\textrm{M}_{D}\textrm{C}(Y),$ \ \ 

$M_{\Phi Q55}=\frac{1}{4}\textrm{M}_{C}\textrm{CCC}(X),$

$M_{\Phi Q62}=\frac{1}{4}\xi E,$ \ \ 

$M_{\Phi Q63}=\frac{1}{4}^{t}\textrm{M}_{D}\textrm{C}(Y),$ \ \ 

$M_{\Phi Q64}=\frac{1}{4}\textrm{M}_{D}\textrm{E}(Y)-\frac{1}{4}\textrm{M}^{+}\textrm{I}(Y),$

$M_{\Phi Q65}=\frac{1}{4}\textrm{M}_{C}\textrm{CJ}(X)$

$M_{\Phi Q71}=\frac{1}{4}\textrm{MCC}(X),$ \ \ 

$M_{\Phi Q72}=-\frac{1}{2}\textrm{M}_{D}\textrm{C}(X),$ \ 

$M_{\Phi Q73}=-\frac{1}{4}\eta E,$ \ \ 

$M_{\Phi Q76}=-\frac{1}{4}\textrm{M}_{C}\textrm{CCC}(Y)$

$M_{\Phi Q81}=-\frac{1}{4}^{t}\textrm{M}_{D}\textrm{C}(X),$ \ \ 

$M_{\Phi Q82}=-\frac{1}{4}\textrm{M}_{D}\textrm{E}(X)+\frac{1}{4}\textrm{M}^{+}\textrm{I}(X),$

$M_{\Phi Q84}=-\frac{1}{4}\eta E,$ \ \ 

$M_{\Phi Q86}=-\frac{1}{4}\textrm{M}_{C}\textrm{CJ}(Y),$

$M_{\Phi Q91}=\frac{1}{8}^{t}\textrm{M}_{C}\textrm{CCC}(Y),$ \ \ 

$M_{\Phi Q92}=\frac{1}{4}^{t}\textrm{M}_{C}\textrm{CJ}(Y),$ \ \ 

$M_{\Phi Q93}=\frac{1}{8}^{t}\textrm{M}_{C}\textrm{CCC}(X),$

$M_{\Phi Q94}=\frac{1}{4}^{t}\textrm{M}_{C}\textrm{CJ}(X),$ \ \ 

$M_{\Phi Q95}=-\frac{3}{8}\eta ,$ \ \ 

$M_{\Phi Q96}=-\frac{3}{8}\xi ,$

\bigskip

$M_{P\Phi }=\left( 
\begin{array}{cccc}
0 & M_{P\Phi 12} & M_{P\Phi 13} & M_{P\Phi 14} \\ 
M_{P\Phi 21} & M_{P\Phi 22} & M_{P\Phi 23} & M_{P\Phi 24} \\ 
0 & M_{P\Phi 32} & M_{P\Phi 33} & M_{P\Phi 34} \\ 
M_{P\Phi 41} & M_{P\Phi 42} & M_{P\Phi 43} & M_{P\Phi 44} \\ 
0 & 0 & 0 & 0 \\ 
0 & 0 & 0 & 0%
\end{array}%
\right.$

$\ \ \ \ \ \ \ \ \ \ \ \ \ \left. 
\begin{array}{ccccc}
M_{P\Phi 15} & 0 & M_{P\Phi 17} & M_{P\Phi 18} & M_{P\Phi 19} \\ 
0 & M_{P\Phi 26} & M_{P\Phi 27} & M_{P\Phi 28} & M_{P\Phi 29} \\ 
M_{P\Phi 35} & M_{P\Phi 36} & M_{P\Phi 37} & 0 & M_{P\Phi 39} \\ 
M_{P\Phi 45} & M_{P\Phi 46} & 0 & M_{P\Phi 48} & M_{P\Phi 49} \\ 
M_{P\Phi 55} & M_{P\Phi 56} & 0 & 0 & M_{P\Phi 59}\\ 
0 & 0 & M_{P\Phi 67} & M_{P\Phi 68} & M_{P\Phi 69}%
\end{array}%
\right) ,$

$M_{P\Phi 12}=\textrm{MJC}(X),$ \ \ 

$M_{P\Phi 13}=-\textrm{M}_{D}^{3}(X),$ \ \ 

$M_{P\Phi 14}=-\textrm{M}^{+}\textrm{JC}(X),$ \ \ 

$M_{P\Phi 15}=-\eta E,$

$M_{P\Phi 17}=-\textrm{MCC}(Y),$ \ \ 

$M_{P\Phi 18}=2\textrm{M}_{D}\textrm{C}(Y),$ \ \ 

$M_{P\Phi19}=\frac{1}{3}\textrm{M}_{C}\textrm{CCC}(X),$

$M_{P\Phi 21}=-\textrm{M}_{C}\textrm{Jr}(X),$ \ \ 

$M_{P\Phi 22}=\frac{1}{2}\textrm{M}_{D}\textrm{E}(X)-\frac{1}{2}\textrm{MI}(X),$

$M_{P\Phi 23}=-\frac{1}{2}^{t}\textrm{M}^{+}\textrm{JC}^{3}(X),$ \ \ 

$M_{P\Phi 24}=-\frac{1}{2}\textrm{M}_{D}^{+}\textrm{E}(X)-\frac{1}{2}\textrm{M}^{+}\textrm{I}(X),$ \ \ 

$M_{P\Phi 26}=-\eta E, $

$M_{P\Phi 27}=$ $^{t}\textrm{M}_{D}\textrm{C}(Y),$ \ \ 

$M_{P\Phi 28}=\textrm{M}_{D}\textrm{E}(Y)-\textrm{M}^{+}\textrm{I}(Y), $ \ \ 

$M_{P\Phi 29}=\frac{1}{3}\textrm{M}_{C}\textrm{CJ}(X),$

$M_{P\Phi 32}=\textrm{MJC}(Y),$ \ \ 

$M_{P\Phi 33}=\textrm{M}_{D}^{3}(Y),$ \ \ 

$M_{P\Phi 34}=\textrm{M}^{+}\textrm{JC}(Y),$ \ \ 

$M_{P\Phi 35}=-\textrm{MCC}(X),$

$M_{P\Phi 36}=2\textrm{M}_{D}\textrm{C}(X),$ \ \ 

$M_{P\Phi 37}=-\xi E,$ \ \ 

$M_{P\Phi 39}=-\frac{1}{3}\textrm{M}_{C}\textrm{CCC}(Y),$

$M_{P\Phi 41}=-\textrm{M}_{C}\textrm{Jr}(Y),$ \ \ 

$M_{P\Phi 42}=\frac{1}{2}\textrm{M}_{D}\textrm{E}(Y)-\frac{1}{2}\textrm{MI}(Y),$ \ \ 

$M_{P\Phi 43}=\frac{1}{2}^{t}\textrm{M}^{+}\textrm{JC}^{3}(X),$

$M_{P\Phi 44}=\frac{1}{2}\textrm{M}_{D}^{+}\textrm{E}(Y)+\frac{1}{2}\textrm{M}^{+}\textrm{I}(Y),$ \ \ 

$M_{P\Phi 45}=$ $^{t}\textrm{M}_{D}\textrm{C}(X),$ \ \ 

$M_{P\Phi 46}=\textrm{M}_{D}\textrm{E}(X)-\textrm{M}^{+}\textrm{I}(X),$ \ \ 

$M_{P\Phi 48}=-\xi E,$ \ \ 

$M_{P\Phi 49}=-\frac{1}{3}\textrm{M}_{C}\textrm{CJ}(Y),$

$M_{P\Phi 55}=-$ $^{t}\textrm{M}_{C}\textrm{CCC}(Y),$ \ \ 

$M_{P\Phi 56}=-2^{t}\textrm{M}_{C}\textrm{CJ}(Y), $ \ \ 

$M_{P\Phi 59}=-\xi ,$

$M_{P\Phi 67}=-$ $^{t}\textrm{M}_{C}\textrm{CCC}(X),$ \ \ 

$M_{P\Phi 68}=-2^{t}\textrm{M}_{C}\textrm{CJ}(X), $ \ \ 

$M_{P\Phi 69}=\eta ,$

\bigskip

$M_{PP}=\left( 
\begin{array}{cccccc}
M_{PP11} & M_{PP12} & M_{PP13} & M_{PP14} & 0 & M_{PP16} \\ 
M_{PP21} & M_{PP22} & M_{PP23} & M_{PP24} & 0 & M_{PP26} \\ 
M_{PP31} & M_{PP32} & M_{PP33} & M_{PP34} & M_{PP35} & 0 \\ 
M_{PP41} & M_{PP42} & M_{PP43} & M_{PP44} & M_{PP45} & 0 \\ 
0 & 0 & M_{PP53} & M_{PP54} & M_{PP55} & 0 \\ 
M_{PP61} & M_{PP62} & 0 & 0 & 0 & M_{PP66}%
\end{array}%
\right) ,$

$M_{PP11}=\textrm{M}_{D}(T)+(-\frac{1}{3}\rho +r)E,$ \ \ 

$M_{PP12}=M_{34}+\textrm{M}^{+}\textrm{JC}(T),$

$M_{PP13}=\textrm{MCC}(B),$ \ \ 

$M_{PP14}=-2\textrm{M}_{D}\textrm{C}(B),$ \ \ 

$M_{PP16}=\textrm{M}_{C}\textrm{CCC}(A),$

$M_{PP21}=2M_{43}+\ ^{t}\textrm{M}^{+}\textrm{JC}(T),$ \ \ 

$M_{PP22}=M_{44}+\frac{1}{2}\textrm{M}_{D}^{+}\textrm{E}(T)+\frac{1}{2}\textrm{M}^{+}\textrm{I}(T)+(-\frac{1}{3}%
\rho +r)E,$

$M_{PP23}=-\ ^{t}\textrm{M}_{D}\textrm{C}(B),$ \ \ 

$M_{PP24}=-\textrm{M}_{D}\textrm{E}(B)+\textrm{M}^{+}\textrm{I}(B),$ \ \ 

$M_{PP26}=\textrm{M}_{C}\textrm{CJ}(A),$

$M_{PP31}=\textrm{MCC}(A),$ \ \ 

$M_{PP32}=-2\textrm{M}_{D}\textrm{C}(A),$ \ \ 

$M_{PP33}=-\textrm{M}_{D}(T)+(\frac{1}{3}\rho +r)E,$

$M_{PP34}=M_{34}-\textrm{M}^{+}\textrm{JC}(T),$ \ \ 

$M_{PP35}=\textrm{M}_{C}\textrm{CCC}(B),$

$M_{PP41}=-\ ^{t}\textrm{M}_{D}\textrm{C}(A),$ \ \ 

$M_{PP42}=-\textrm{M}_{D}\textrm{E}(A)+\textrm{M}^{+}\textrm{I}(A),$

$M_{PP43}=M_{43}-\frac{1}{2}\ ^{t}\textrm{M}^{+}\textrm{JC}(T),$

$M_{PP44}=M_{44}-\frac{1}{2}\textrm{M}_{D}^{+}\textrm{E}(T)-\frac{1}{2}\textrm{M}^{+}\textrm{I}(T)+(\frac{1}{3}\rho
+r)E,$ \ \ 

$M_{PP45}=\textrm{M}_{C}\textrm{CJ}(B),$

$M_{PP53}=$ $^{t}\textrm{M}_{C}\textrm{CCC}(A),$ \ \ 

$M_{PP54}=2^{t}\textrm{M}_{C}\textrm{CJ}(A),$ \ \ $%
M_{PP55}=(\rho +r),$

$M_{PP61}=$ $^{t}\textrm{M}_{C}\textrm{CCC}(B),$ \ \ 

$M_{PP62}=2^{t}\textrm{M}_{C}\textrm{CJ}(B),$ \ \ 

$M_{PP66}=(-\rho +r),$

\bigskip

$M_{PQ}=\left( 
\begin{array}{cccccc}
sE & 0 & 0 & 0 & 0 & 0 \\ 
0 & sE & 0 & 0 & 0 & 0 \\ 
0 & 0 & sE & 0 & 0 & 0 \\ 
0 & 0 & 0 & sE & 0 & 0 \\ 
0 & 0 & 0 & 0 & s & 0 \\ 
0 & 0 & 0 & 0 & 0 & s%
\end{array}%
\right) ,$

\bigskip

$M_{Pr}=\left( 
\begin{array}{c}
M_{Pr11} \\ 
M_{Pr21} \\ 
M_{Pr31} \\ 
M_{Pr41} \\ 
M_{Pr51} \\ 
M_{Pr61}%
\end{array}%
\right) ,$

$M_{Pr11}=-\textrm{M}_{C}\textrm{CCC}(X),$ \ \ 

$M_{Pr21}=-\textrm{M}_{C}\textrm{CJ}(X),$ \ \ 

$M_{Pr31}=-\textrm{M}_{C}\textrm{CCC}(Y),$

$M_{Pr41}=-\textrm{M}_{C}\textrm{CJ}(Y),$\ \ 

$M_{Pr51}=-\xi ,$\ \

$M_{Pr61}=-\eta ,$

$\bigskip $

$M_{Ps}=\left( 
\begin{array}{c}
M_{Ps11} \\ 
M_{Ps21} \\ 
M_{Ps31} \\ 
M_{Ps41} \\ 
M_{Ps51} \\ 
M_{Ps61}%
\end{array}%
\right) ,$

$M_{Ps11}=-\textrm{M}_{C}\textrm{CCC}(Z),$ \ \ 

$M_{Ps21}=-\textrm{M}_{C}\textrm{CJ}(Z),$ \ \ 

$M_{Ps31}=-\textrm{M}_{C}\textrm{CCC}(W),$

$M_{Ps41}=-\textrm{M}_{C}\textrm{CJ}(W),$ \ \ 

$M_{Ps51}=-\zeta,$ \ \ 

$M_{Ps61}=-\omega ,$

\bigskip

$M_{Q\Phi }=\left( 
\begin{array}{cccc}
0 & M_{Q\Phi 12} & M_{Q\Phi 13} & M_{Q\Phi 14} \\ 
M_{Q\Phi 21} & M_{Q\Phi 22} & M_{Q\Phi 23} & M_{Q\Phi 24}  \\ 
0 & M_{Q\Phi 32} & M_{Q\Phi 33} & M_{Q\Phi 34}  \\ 
M_{Q\Phi 41} & M_{Q\Phi 42} & M_{Q\Phi 43} & M_{Q\Phi 44}  \\ 
0 & 0 & 0 & 0 \\ 
0 & 0 & 0 & 0%
\end{array}%
\right.$

\ \ \ \ \ \ \ \ \ \ \ \ \ $\left. 
\begin{array}{ccccc}
M_{Q\Phi 15} & 0 & M_{Q\Phi 17} & M_{Q\Phi 18} & M_{Q\Phi 19} \\ 
0 & M_{Q\Phi 26} & M_{Q\Phi 27} & M_{Q\Phi 28} & M_{Q\Phi 29} \\ 
M_{Q\Phi 35} & M_{Q\Phi 36} & M_{Q\Phi 37} & 0 & M_{Q\Phi 39} \\ 
M_{Q\Phi 45} & M_{Q\Phi 46} & 0 & M_{Q\Phi 48} & M_{Q\Phi 49} \\ 
M_{Q\Phi 55} & M_{Q\Phi 56} & 0 & 0 & M_{Q\Phi 59}\\ 
0 & 0 & M_{Q\Phi 67} & M_{Q\Phi 68} & M_{Q\Phi 69}%
\end{array}%
\right) ,$

$M_{Q\Phi 12}=\textrm{MJC}(Z),$ \ \ 

$M_{Q\Phi 13}=-\textrm{M}_{D}^{3}(Z),$ \ \ 

$M_{Q\Phi 14}=-\textrm{M}^{+}\textrm{JC}(Z),$ \ \ 

$M_{Q\Phi 15}=-\omega E,$

$M_{Q\Phi 17}=-\textrm{MCC}(W),$ \ \ 

$M_{Q\Phi 18}=2\textrm{M}_{D}\textrm{C}(W),$ \ \ 

$M_{Q\Phi 19}=\frac{1}{3}\textrm{M}_{C}\textrm{CCC}(Z),$

$M_{Q\Phi 21}=-\textrm{M}_{C}\textrm{Jr}(Z),$ \ \ 

$M_{Q\Phi 22}=\frac{1}{2}\textrm{M}_{D}^{-}\textrm{E}(Z)-\frac{1}{2}\textrm{MI}(Z),$ \ \ 

$M_{Q\Phi 23}=-\frac{1}{2}^{t}\textrm{M}^{+}\textrm{JC}^{3}(Z),$

$M_{Q\Phi 24}=-\frac{1}{2}\textrm{M}_{D}^{+}\textrm{E}(Z)-\frac{1}{2}\textrm{M}^{+}\textrm{I}(Z),$ \ \ 

$M_{Q\Phi 26}=-\omega E,$ \ \ 

$M_{Q\Phi 27}=$ $^{t}\textrm{M}_{D}\textrm{C}(W),$ \ \ 

$M_{Q\Phi 28}=\textrm{M}_{D}\textrm{E}(W)-\textrm{M}^{+}\textrm{I}(W),$ \ \ 

$M_{Q\Phi 29}=\frac{1}{3}\textrm{M}_{C}\textrm{CJ}(Z),$

$M_{Q\Phi 32}=\textrm{MJC}(W),$ \ \ 

$M_{Q\Phi 33}=\textrm{M}_{D}^{3}(W),$ \ \ 

$M_{Q\Phi 34}=\textrm{M}^{+}\textrm{JC}(W),$ \ \ 

$M_{Q\Phi 35}=-\textrm{MCC}(Z),$

$M_{Q\Phi 36}=2\textrm{M}_{D}\textrm{C}(Z),$ \ \ 

$M_{Q\Phi 37}=-\zeta E,$ \ \ 

$M_{Q\Phi 39}=-\frac{1}{3}\textrm{M}_{C}\textrm{CCC}(W),$

$M_{Q\Phi 41}=-\textrm{M}_{C}\textrm{Jr}(W),$ \ \ 

$M_{Q\Phi 42}=\frac{1}{2}\textrm{M}_{D}^{-}\textrm{E}(W)-\frac{1}{2}\textrm{MI}(W),$ \ \ 

$M_{Q\Phi 43}=\frac{1}{2}^{t}\textrm{M}^{+}\textrm{JC}^{3}(W),$

$M_{Q\Phi 44}=\frac{1}{2}\textrm{M}_{D}^{+}\textrm{E}(W)+\frac{1}{2}\textrm{M}^{+}\textrm{I}(W),$ \ \ 

$M_{Q\Phi 45}=$ $^{t}\textrm{M}_{D}\textrm{C}(Z),$

$M_{Q\Phi 46}=\textrm{M}_{D}\textrm{E}(Z)-\textrm{M}^{+}\textrm{I}(Z),$ \ \ 

$M_{Q\Phi 48}=-\zeta E,$ \ \ 

$M_{Q\Phi 49}=-\frac{1}{3}\textrm{M}_{C}\textrm{CJ}(W),$

$M_{Q\Phi 55}=-$ $^{t}\textrm{M}_{C}\textrm{CCC}(W),$ \ \ 

$M_{Q\Phi 56}=-2^{t}\textrm{M}_{C}\textrm{CJ}(W), $ \ \ 

$M_{Q\Phi 59}=-\zeta ,$

$M_{Q\Phi 67}=-$ $^{t}\textrm{M}_{C}\textrm{CCC}(Z),$ \ \ 

$M_{Q\Phi68}=-2^{t}\textrm{M}_{C}\textrm{CJ}(Z), $ \ \ 

$M_{Q\Phi 69}=\omega ,$

\bigskip

$M_{QP}=\left( 
\begin{array}{cccccc}
uE & 0 & 0 & 0 & 0 & 0 \\ 
0 & uE & 0 & 0 & 0 & 0 \\ 
0 & 0 & uE & 0 & 0 & 0 \\ 
0 & 0 & 0 & uE & 0 & 0 \\ 
0 & 0 & 0 & 0 & u & 0 \\ 
0 & 0 & 0 & 0 & 0 & u%
\end{array}%
\right) ,$

\bigskip

$M_{QQ}=\left( 
\begin{array}{cccccc}
M_{QQ11} & M_{QQ12} & M_{QQ13} & M_{QQ14} & 0 & M_{QQ16} \\ 
M_{QQ21} & M_{QQ22} & M_{QQ23} & M_{QQ24} & 0 & M_{QQ26} \\ 
M_{QQ31} & M_{QQ32} & M_{QQ33} & M_{QQ34} & M_{QQ35} & 0 \\ 
M_{QQ41} & M_{QQ42} & M_{QQ43} & M_{QQ44} & M_{QQ45} & 0 \\ 
0 & 0 & M_{QQ53} & M_{QQ54} & M_{QQ55} & 0 \\ 
M_{QQ61} & M_{QQ62} & 0 & 0 & 0 & M_{QQ66}%
\end{array}%
\right) ,$

$M_{QQ11}=\textrm{M}_{D}(T)+(-\frac{1}{3}\rho -r)E,$ \ \ 

$M_{QQ12}=M_{34}+\textrm{M}^{+}\textrm{JC}(T),$

$M_{QQ13}=\textrm{MCC}(B),$ \ \ 

$M_{QQ14}=-2\textrm{M}_{D}\textrm{C}(B),$ \ \ 

$M_{QQ16}=\textrm{M}_{C}\textrm{CCC}(A),$

$M_{QQ21}=2M_{43}+\ ^{t}\textrm{M}^{+}\textrm{JC}(T),$ \ \ 

$M_{QQ22}=M_{44}+\frac{1}{2}\textrm{M}_{D}^{+}\textrm{E}(T)+\frac{1}{2}\textrm{M}^{+}\textrm{I}(T)+(-\frac{1}{3}\rho -r)E$

$M_{QQ23}=-\ ^{t}\textrm{M}_{D}\textrm{C}(B),$\ \ 

$M_{QQ24}=-\textrm{M}_{D}\textrm{E}(B)+\textrm{M}^{+}\textrm{I}(B),$

$M_{QQ26}=\textrm{M}_{C}\textrm{CJ}(A),$

$M_{QQ31}=\textrm{MCC}(A),$ \ \ 

$M_{QQ32}=-2\textrm{M}_{D}\textrm{C}(A),$ \ \ 

$M_{QQ33}=-\textrm{M}_{D}(T)+(\frac{1}{3}\rho -r)E,$

$M_{QQ34}=M_{34}-\textrm{M}^{+}\textrm{JC}(T),$ \ \ 

$M_{QQ35}=\textrm{M}_{C}\textrm{CCC}(B),$

$M_{QQ41}=-\ ^{t}\textrm{M}_{D}\textrm{C}(A),$ \ \ 

$M_{QQ42}=-\textrm{M}_{D}\textrm{E}(A)+\textrm{M}^{+}\textrm{I}(A),$

$M_{QQ43}=M_{43}-\frac{1}{2}\ ^{t}\textrm{M}^{+}\textrm{JC}(T),$

$M_{QQ44}=M_{44}-\frac{1}{2}\textrm{M}_{D}^{+}\textrm{E}(T)-\frac{1}{2}\textrm{M}^{+}\textrm{I}(T)+(\frac{1}{3}\rho-r)E,$

$M_{QQ45}=\textrm{M}_{C}\textrm{CJ}(B),$

$M_{QQ53}=$ $^{t}\textrm{M}_{C}\textrm{CCC}(A),$ 

$M_{QQ54}=2^{t}\textrm{M}_{C}\textrm{CJ}(A),$ \ \ 

$M_{QQ55}=(\rho -r),$

$M_{QQ61}=$ $^{t}\textrm{M}_{C}\textrm{CCC}(B),$ \ \ 

$M_{QQ62}=2^{t}\textrm{M}_{C}\textrm{CJ}(B),$ \ \ 

$M_{QQ66}=(-\rho -r),$

\bigskip

$M_{Qr}=\left( 
\begin{array}{c}
M_{Qr11} \\ 
M_{Qr21} \\ 
M_{Qr31} \\ 
M_{Qr41} \\ 
M_{Qr51} \\ 
M_{Qr61}%
\end{array}%
\right) ,$

$M_{Qr11}=-\textrm{M}_{C}\textrm{CCC}(Z),$ \ \ 

$M_{Qr21}=-\textrm{M}_{C}\textrm{CJ}(Z),$ \ \ 

$M_{Qr31}=-\textrm{M}_{C}\textrm{CCC}(W),$

$M_{Qr41}=-\textrm{M}_{C}\textrm{CJ}(W),$ \ \ 

$M_{Qr51}=-\zeta,$ \ \ 

$M_{Qr61}=-\omega ,$

$\bigskip $

$M_{Qu}=\left( 
\begin{array}{c}
M_{Qu11} \\ 
M_{Qu21} \\ 
M_{Qu31} \\ 
M_{Qu41} \\ 
M_{Qu51} \\ 
M_{Qu61}%
\end{array}%
\right) ,$

$M_{Qu11}=-\textrm{M}_{C}\textrm{CCC}(X),$ \ \ 

$M_{Qu21}=-\textrm{M}_{C}\textrm{CJ}(X),$ \ \ 

$M_{Qu31}=-\textrm{M}_{C}\textrm{CCC}(Y),$

$M_{Qu41}=-\textrm{M}_{C}\textrm{CJ}(Y),$ \ \ 

$M_{Qu51}=-\xi ,$

$M_{Qu61}=-\eta ,$

\bigskip

$M_{rP}=\left( 
\begin{array}{cccccc}
M_{rP11} & M_{rP12} & M_{rP13} & M_{rP14} & M_{rP15} & M_{rP16}%
\end{array}%
\right) ,$

$M_{rP11}=\frac{1}{8}^{t}\textrm{M}_{C}\textrm{CCC}(W),$ \ \ 

$M_{rP12}=\frac{1}{4}^{t}\textrm{M}_{C}\textrm{CJ}(W),$ \ \ 

$M_{rP13}=-\frac{1}{8}^{t}\textrm{M}_{C}\textrm{CCC}(Z),$

$M_{rP14}=-\frac{1}{4}^{t}\textrm{M}_{C}\textrm{CJ}(Z),$ \ \ 

$M_{rP15}=\frac{1}{8}\omega ,$ \

$M_{rP16}=\frac{1}{8}\zeta ,$

\bigskip

$M_{rQ}=\left( 
\begin{array}{cccccc}
M_{rQ11} & M_{rQ12} & M_{rQ13} & M_{rQ14} & M_{rQ15} & M_{rQ16}%
\end{array}%
\right) ,$

$M_{rQ11}=\frac{1}{8}^{t}\textrm{M}_{C}\textrm{CCC}(Y),$ \ \ 

$M_{rQ12}=\frac{1}{4}^{t}\textrm{M}_{C}\textrm{CJ}(Y),$ \ \ 

$M_{rQ13}=-\frac{1}{8}^{t}\textrm{M}_{C}\textrm{CCC}(X),$

$M_{rQ14}=-\frac{1}{4}^{t}\textrm{M}_{C}\textrm{CJ}(X),$ \ \ 

$M_{rQ15}=\frac{1}{8}\eta E,$ \ \ 

$M_{rQ16}=-\frac{1}{8}\xi E,$

$M_{rs}=-u,$

$M_{ru}=s,$

\bigskip

$M_{sP}=\left( 
\begin{array}{cccccc}
M_{sP11} & M_{sP12} & M_{sP13} & M_{sP14} & M_{sP15} & M_{sP16}%
\end{array}%
\right) ,$

$M_{sP11}=-\frac{1}{4}^{t}\textrm{M}_{C}\textrm{CCC}(Y),$ \ \ 

$M_{sP12}=-\frac{1}{2}^{t}\textrm{M}_{C}\textrm{CJ}(Y),$ 

$M_{sP13}=\frac{1}{4}^{t}\textrm{M}_{C}\textrm{CCC}(X),$

$M_{sP14}=\frac{1}{2}^{t}\textrm{M}_{C}\textrm{CJ}(X),$ \ \ 

$M_{sP15}=-\frac{1}{4}\eta E,$ \ \ 

$M_{sP16}=\frac{1}{4}\xi E,$

$M_{sr}=-2sE,$

$M_{ss}=2rE,$

\bigskip

$M_{uQ}=\left( 
\begin{array}{cccccc}
M_{uQ11} & M_{uQ12} & M_{uQ13} & M_{uQ14} & M_{uQ15} & M_{uQ16}%
\end{array}%
\right) ,$

$M_{uQ11}=\frac{1}{4}^{t}\textrm{M}_{C}\textrm{CCC}(W),$ \ \ 

$M_{uQ12}=\frac{1}{2}^{t}\textrm{M}_{C}\textrm{CJ}(W),$ \ \ 

$M_{uQ13}=-\frac{1}{4}^{t}\textrm{M}_{C}\textrm{CCC}(Z),$

$M_{uQ14}=-\frac{1}{2}^{t}\textrm{M}_{C}\textrm{CJ}(Z),$ \ \ 

$M_{uQ15}=\frac{1}{4}\omega ,$ \

$M_{uQ16}=-\frac{1}{4}\zeta ,$

$M_{ur}=2u,$

$M_{uu}=-2r,$

\noindent
Hence we have the above expression of $adR$ . \ \ \ \ \emph{Q.E.D.}

\bigskip

\section{Elements of the digitalized matrix \gr$_{8}^{\C}$}

\bigskip

\ \ \ \ \ From \emph{Theorem 10.7} we digitalize \gr$_{8}^{\C}$. As a result,  \gr$_{8}^{\C}$ is expressed on the
matrix $M(248\times 248,\C)$ as follows.

\pagebreak

\begin{figure}[H]
\centering
\includegraphics[width=12cm, height=19.0cm,page=1 ]{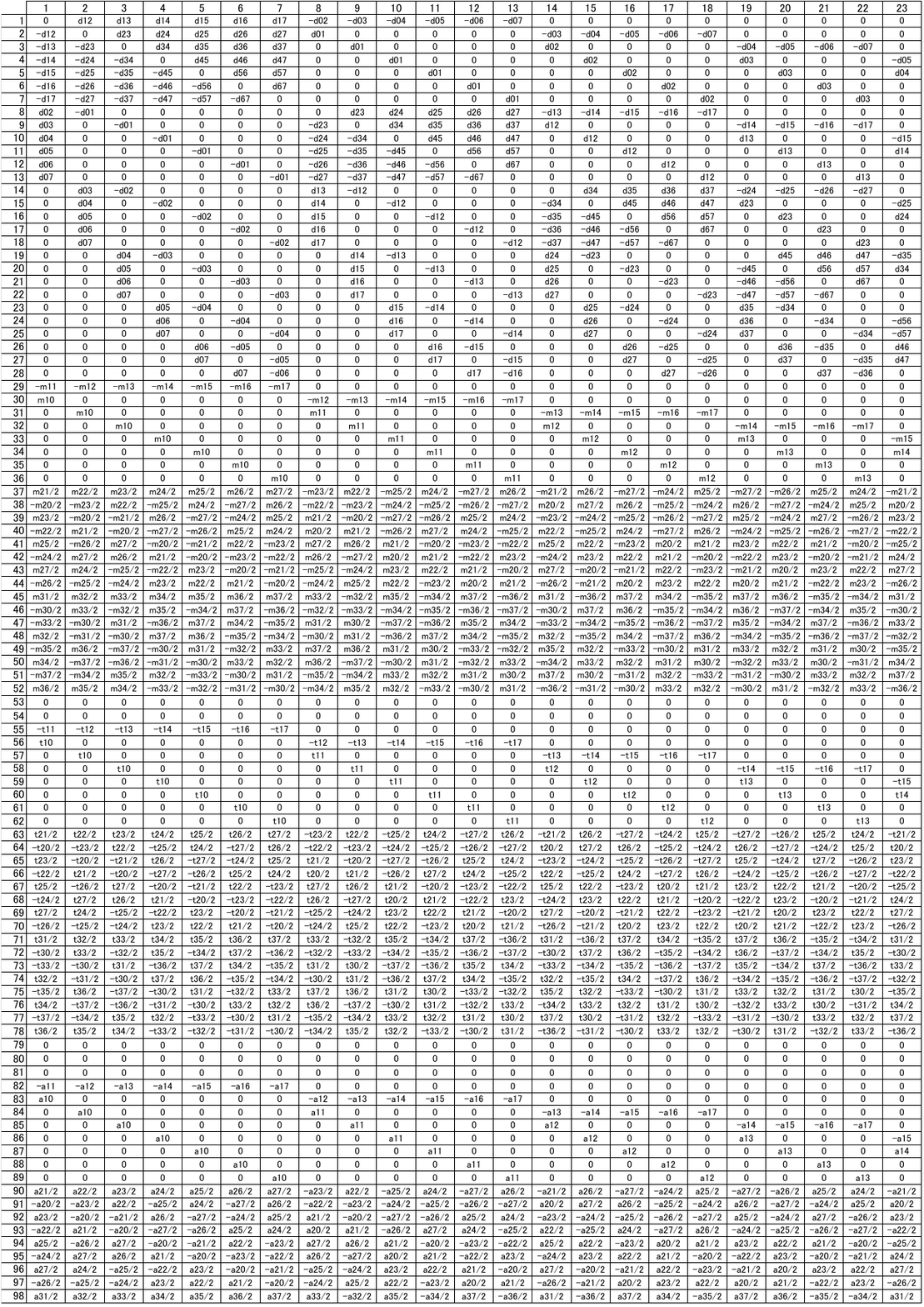}
\end{figure}

\begin{figure}[H]
\centering
\includegraphics[width=12cm, height=19.0cm,page=2 ]{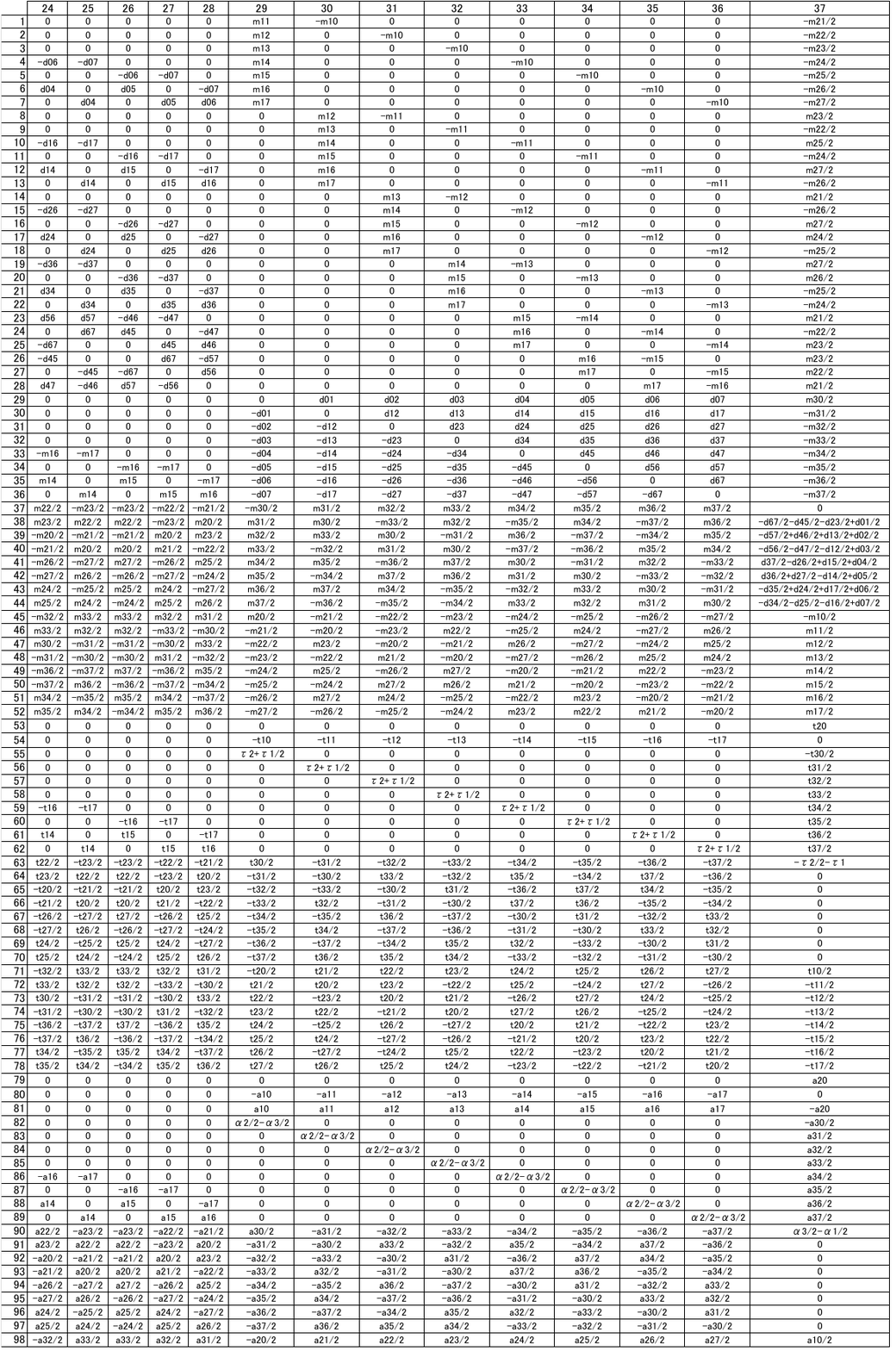}
\end{figure}

\begin{figure}[H]
\centering
\includegraphics[width=12cm, height=19.0cm,page=3 ]{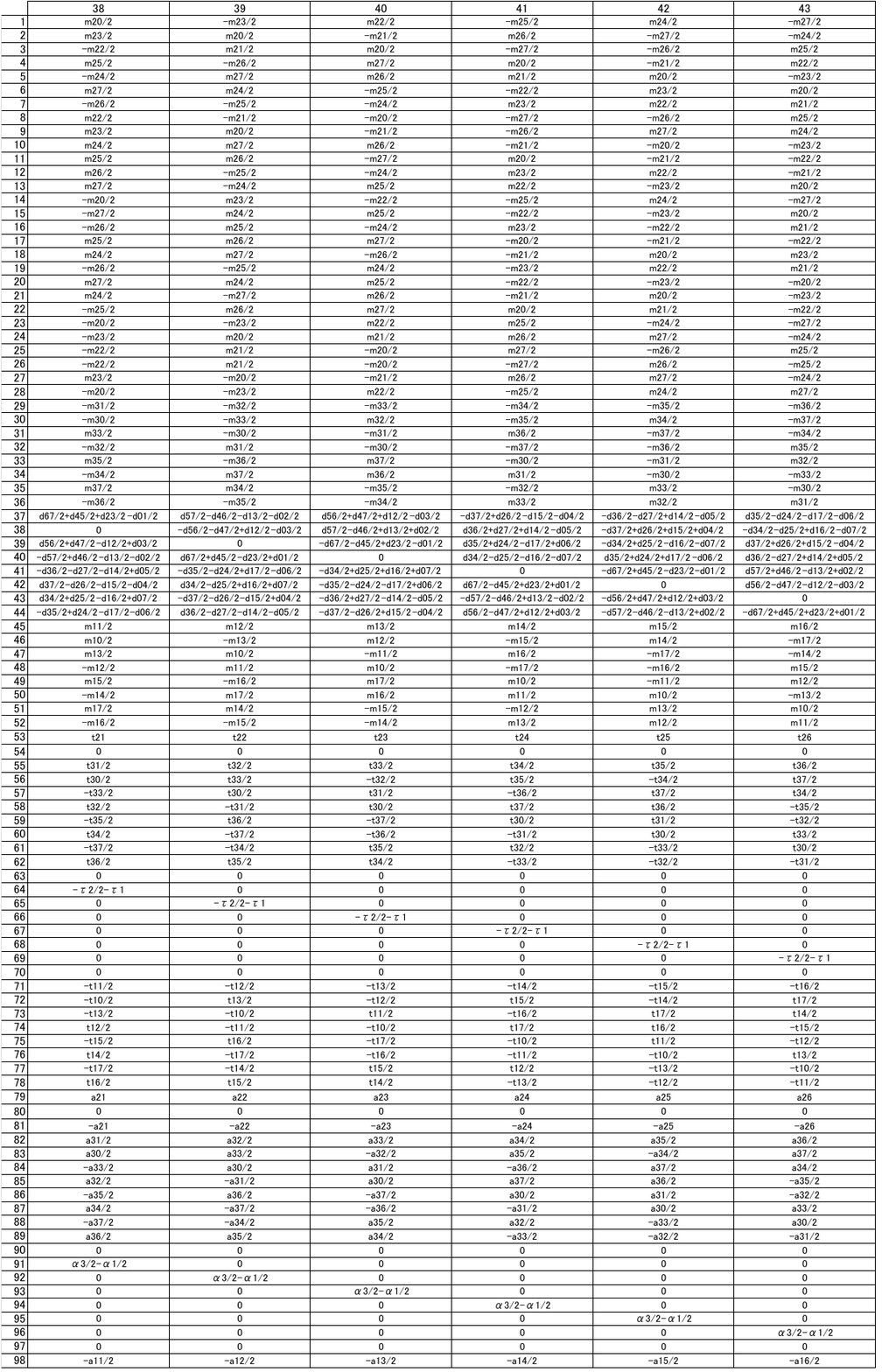}
\end{figure}

\begin{figure}[H]
\centering
\includegraphics[width=12cm, height=19.0cm,page=4 ]{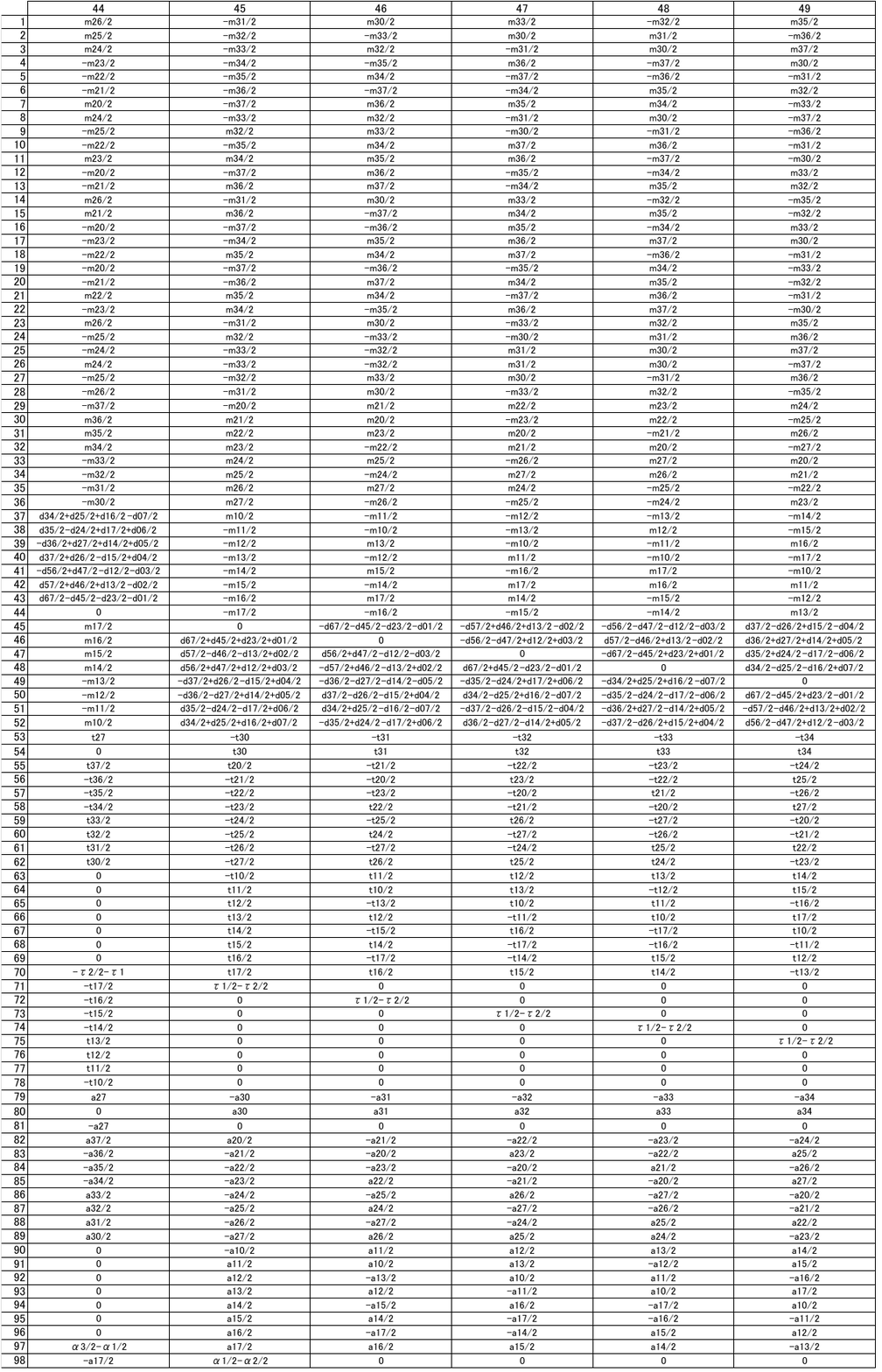}
\end{figure}

\begin{figure}[H]
\centering
\includegraphics[width=12cm, height=19.0cm,page=5 ]{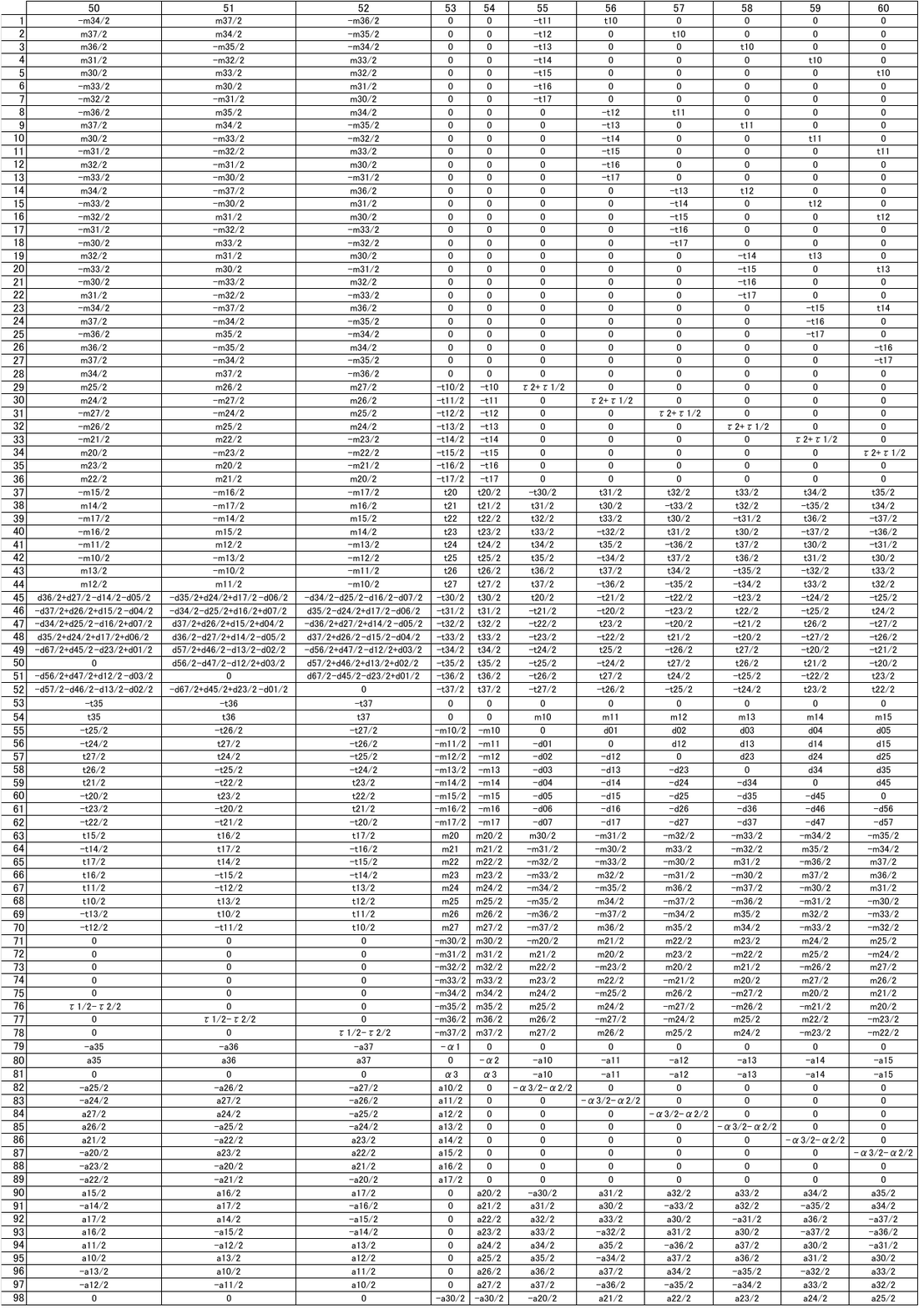}
\end{figure}

\begin{figure}[H]
\centering
\includegraphics[width=12cm, height=19.0cm,page=6 ]{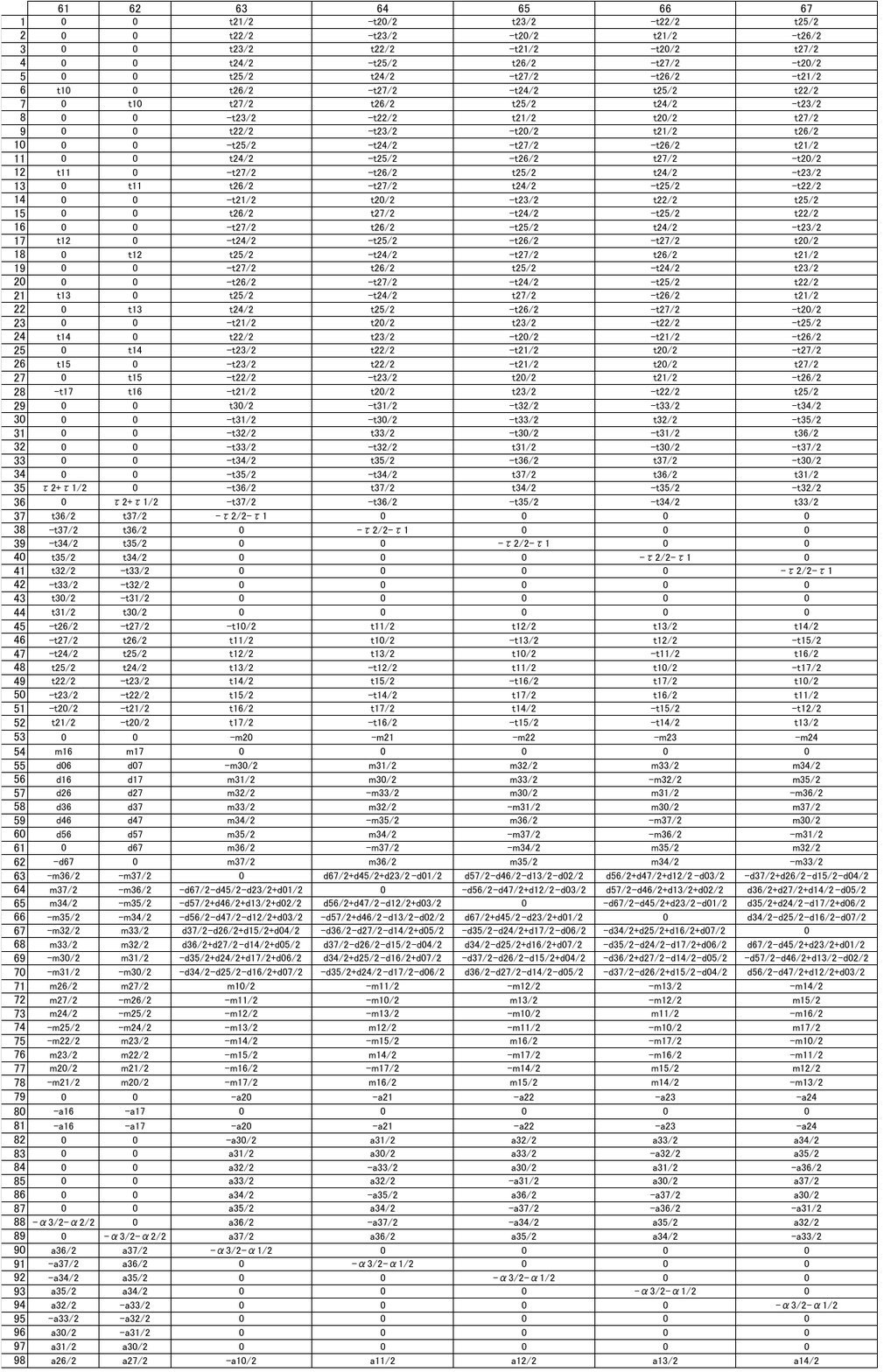}
\end{figure}

\begin{figure}[H]
\centering
\includegraphics[width=12cm, height=19.0cm,page= 7 ]{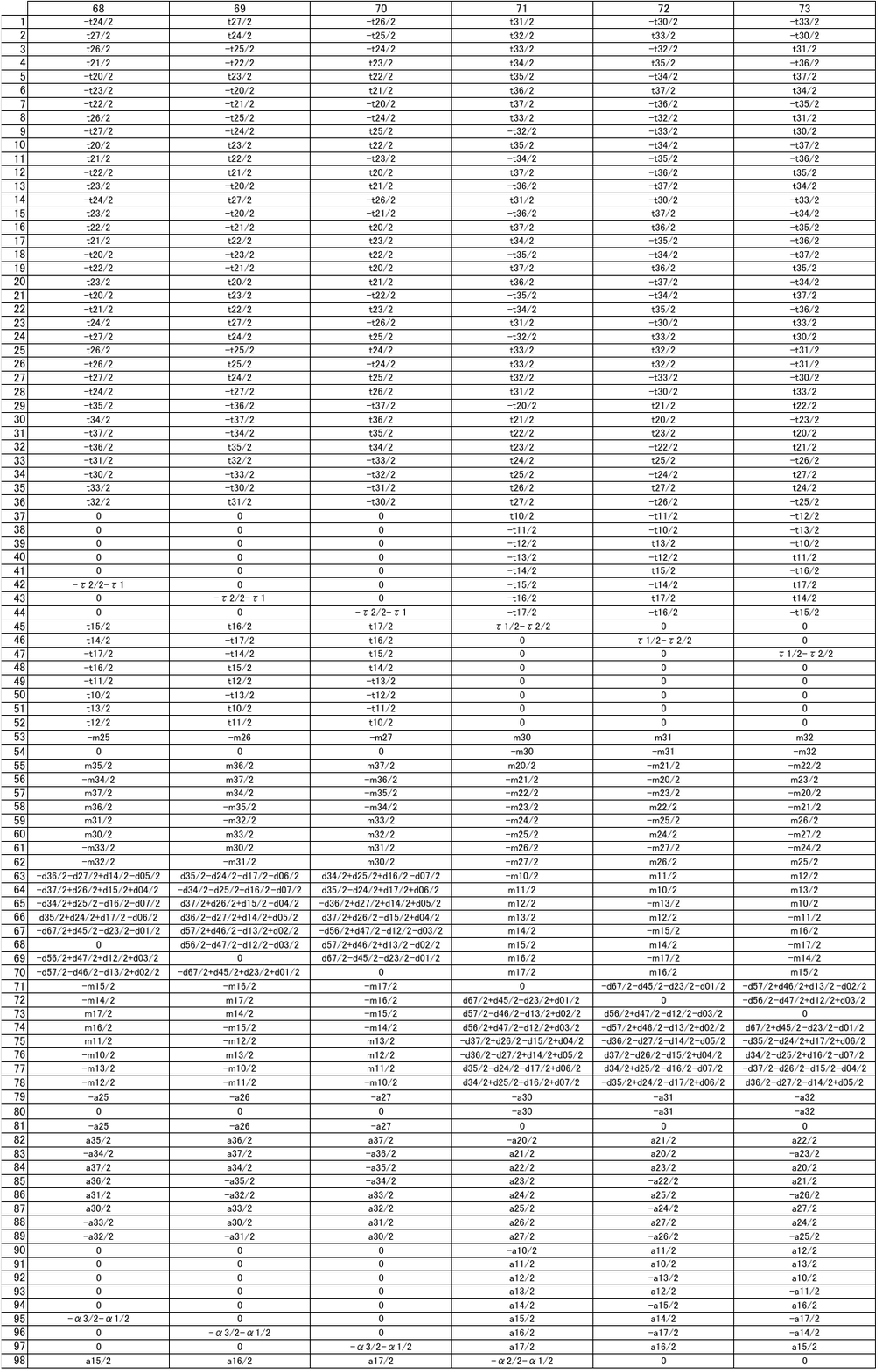}
\end{figure}

\begin{figure}[H]
\centering
\includegraphics[width=12cm, height=19.0cm,page= 8 ]{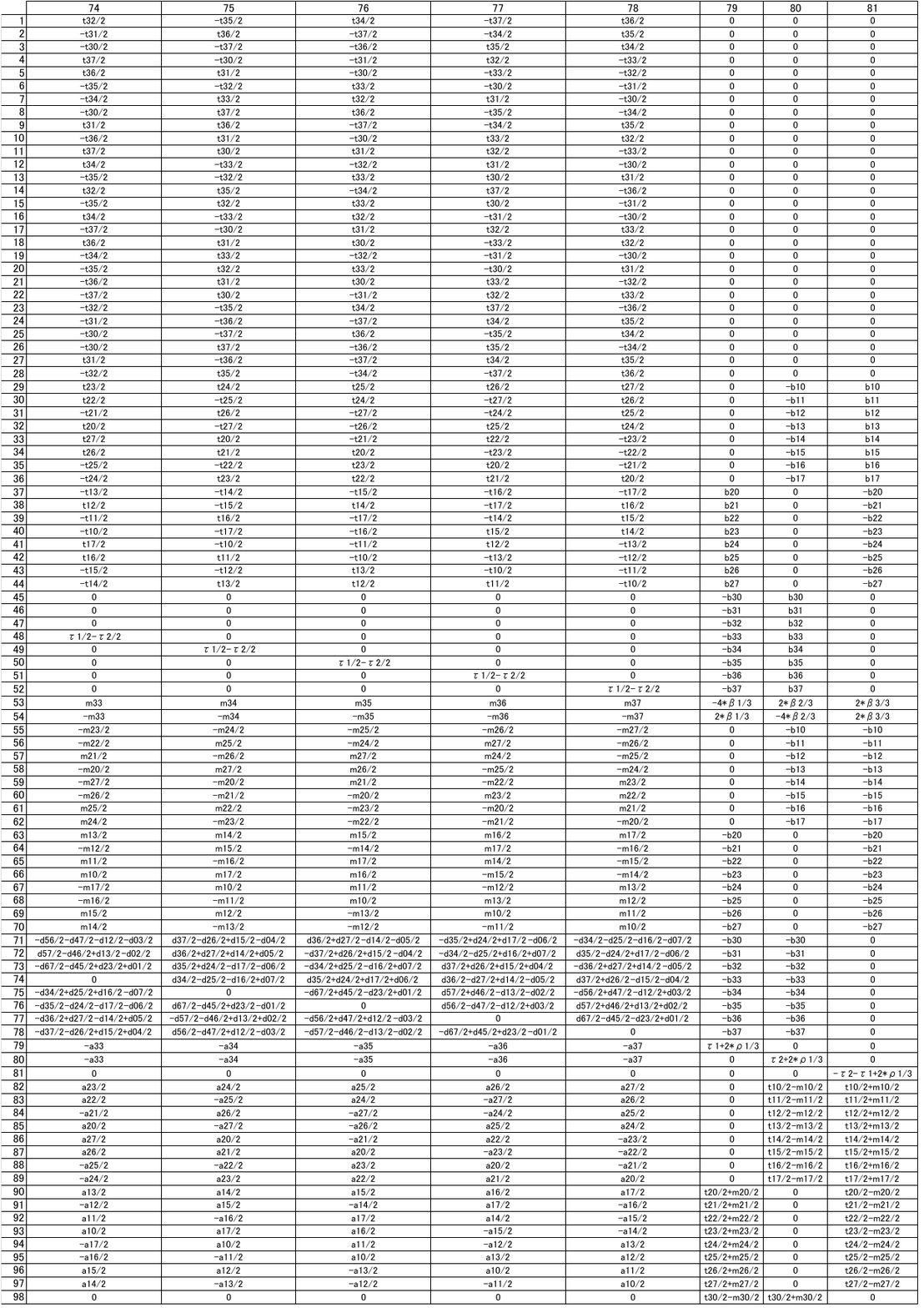}
\end{figure}

\begin{figure}[H]
\centering
\includegraphics[width=12cm, height=19.0cm,page= 9 ]{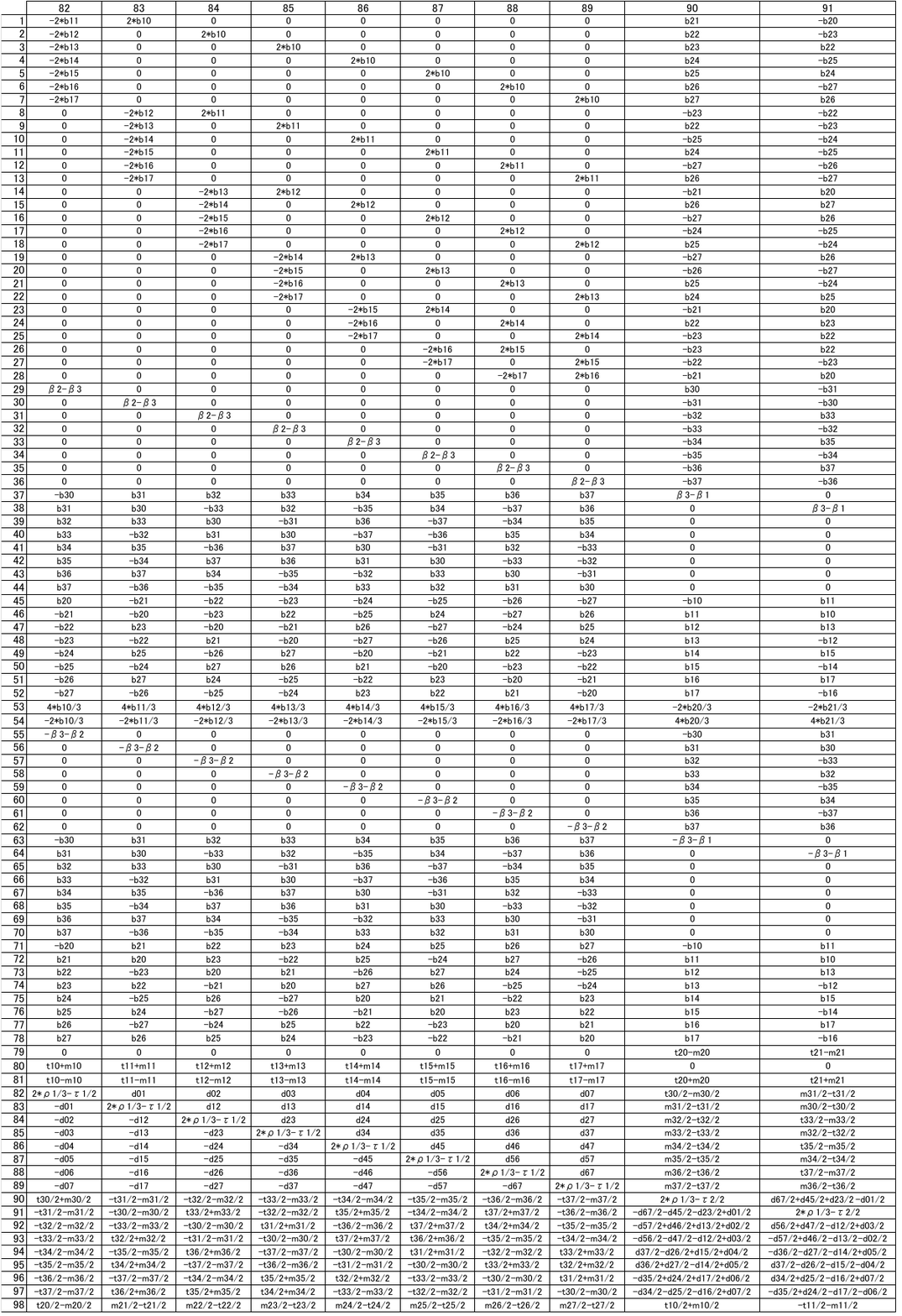}
\end{figure}

\begin{figure}[H]
\centering
\includegraphics[width=12cm, height=19.0cm,page= 10 ]{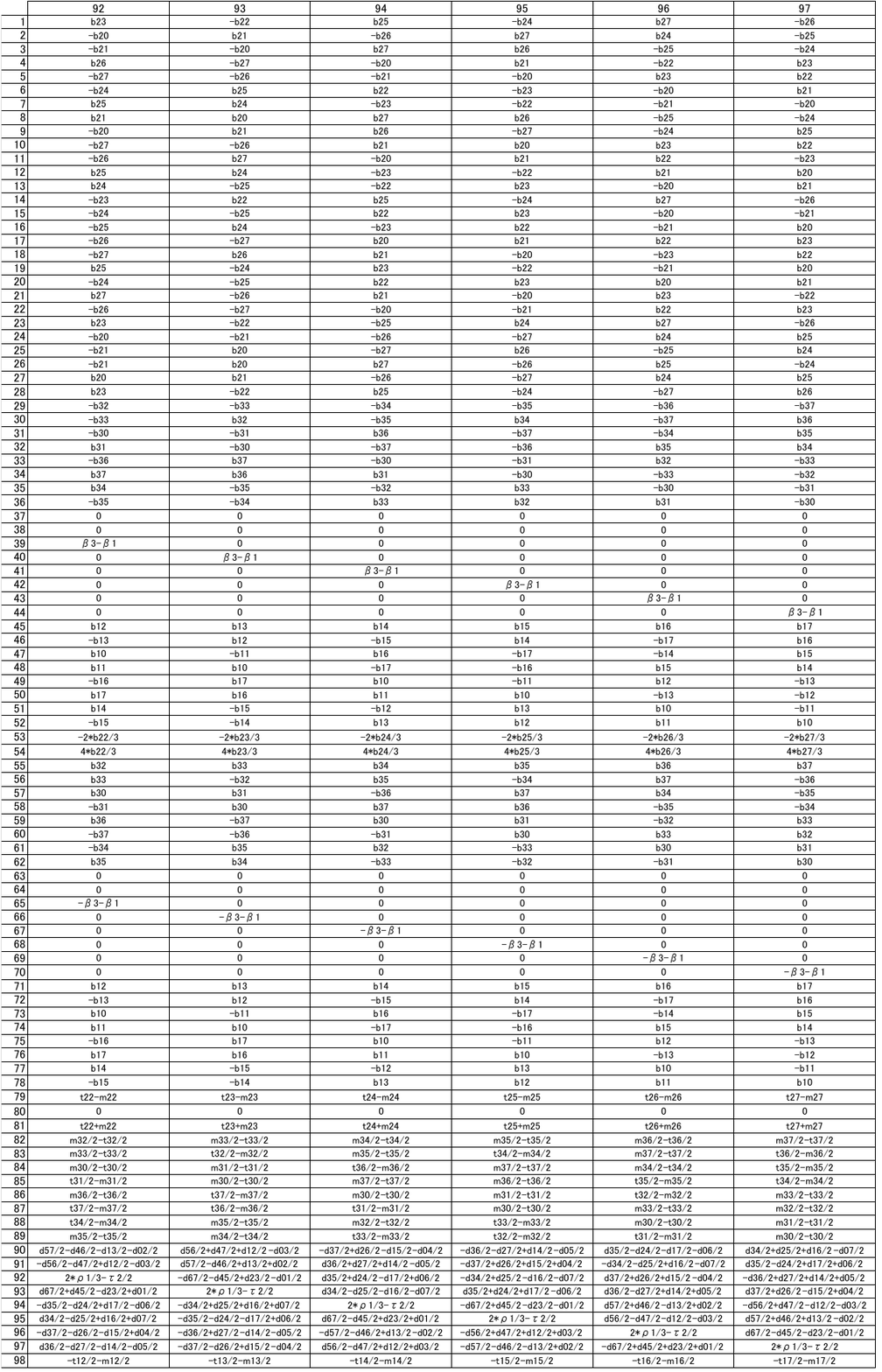}
\end{figure}

\begin{figure}[H]
\centering
\includegraphics[width=12cm, height=19.0cm,page= 11 ]{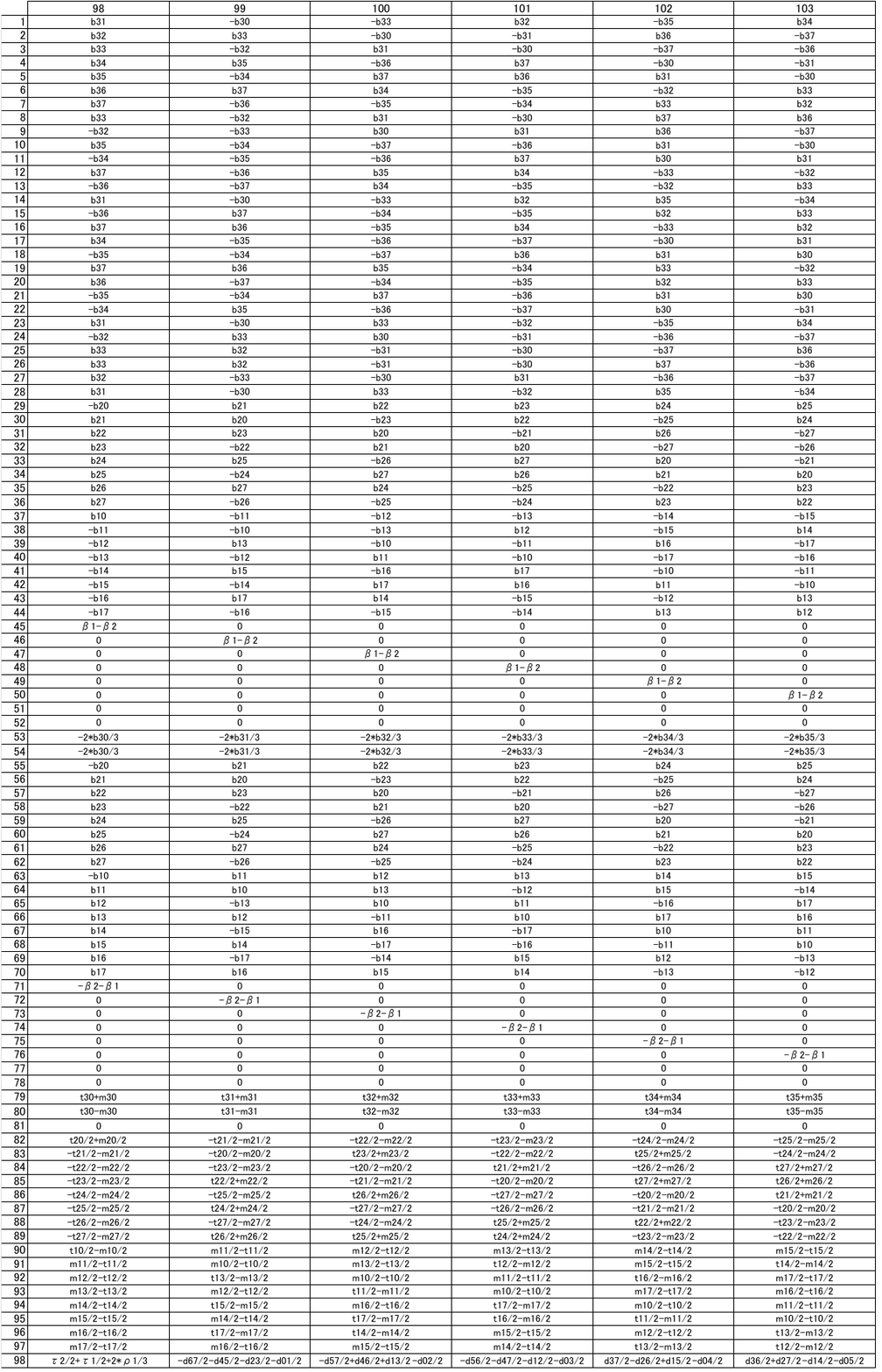}
\end{figure}

\begin{figure}[H]
\centering
\includegraphics[width=12cm, height=19.0cm,page= 12 ]{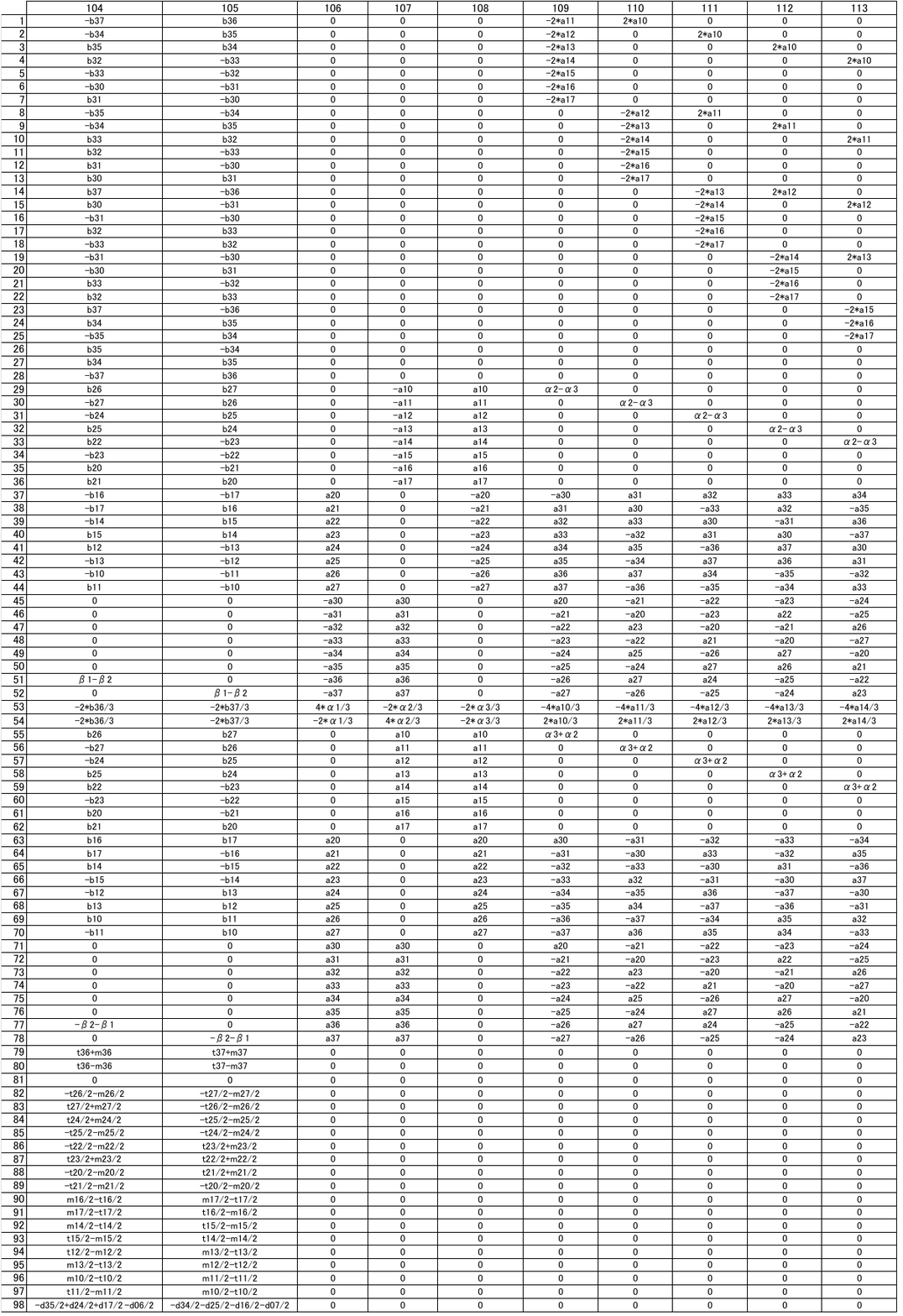}
\end{figure}

\begin{figure}[H]
\centering
\includegraphics[width=12cm, height=19.0cm,page= 13 ]{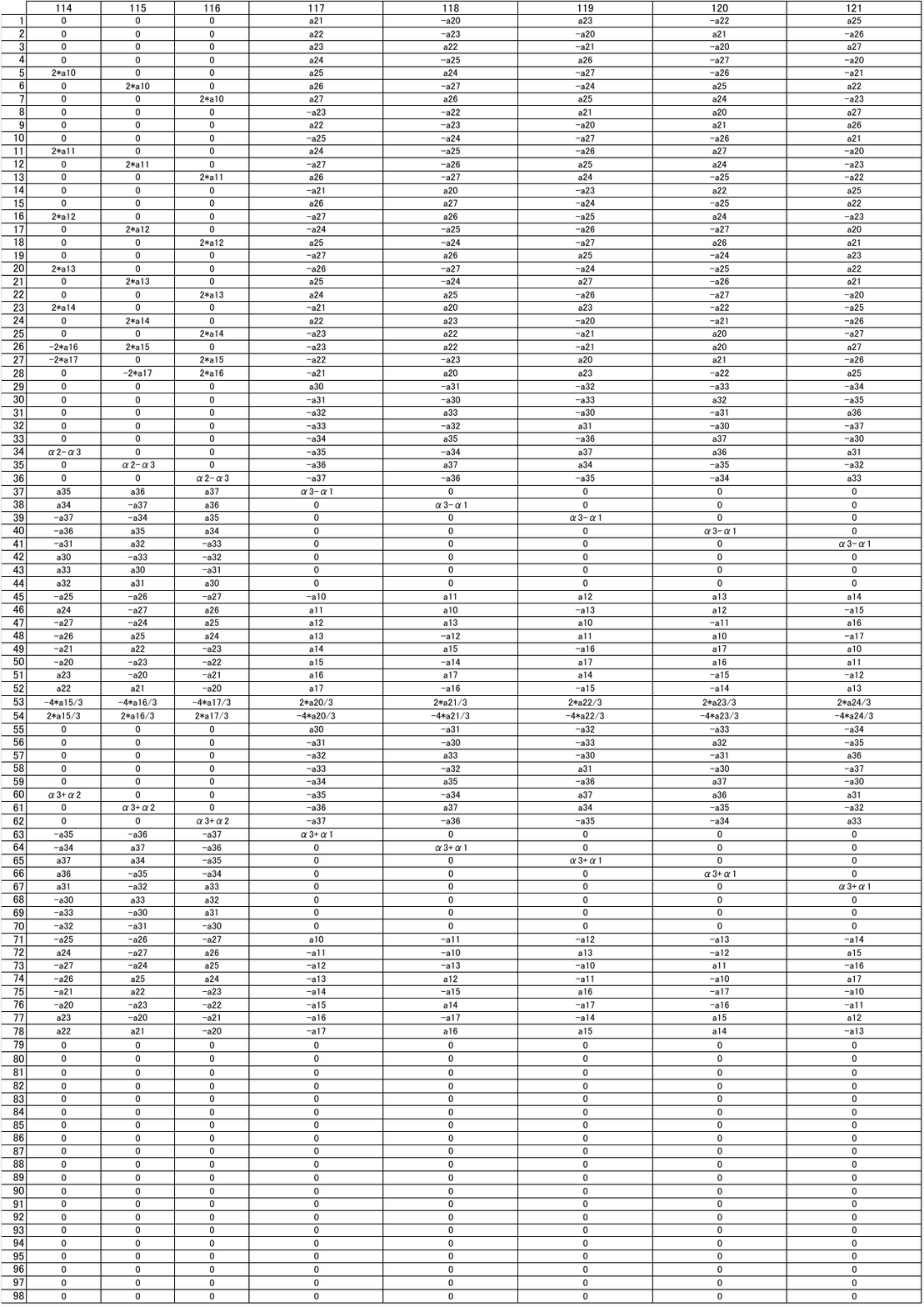}
\end{figure}

\begin{figure}[H]
\centering
\includegraphics[width=12cm, height=19.0cm,page= 14 ]{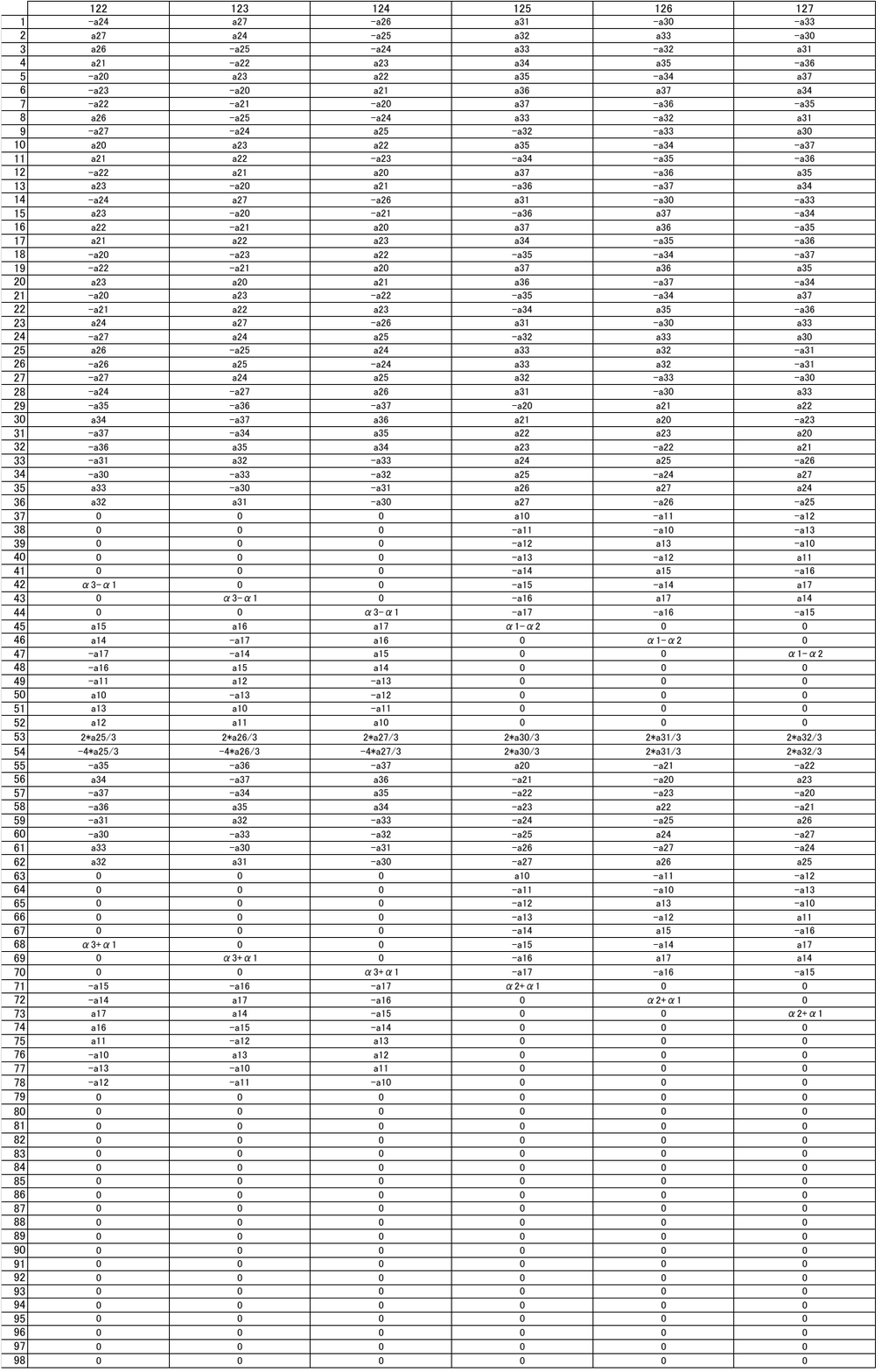}
\end{figure}

\begin{figure}[H]
\centering
\includegraphics[width=12cm, height=19.0cm,page= 15 ]{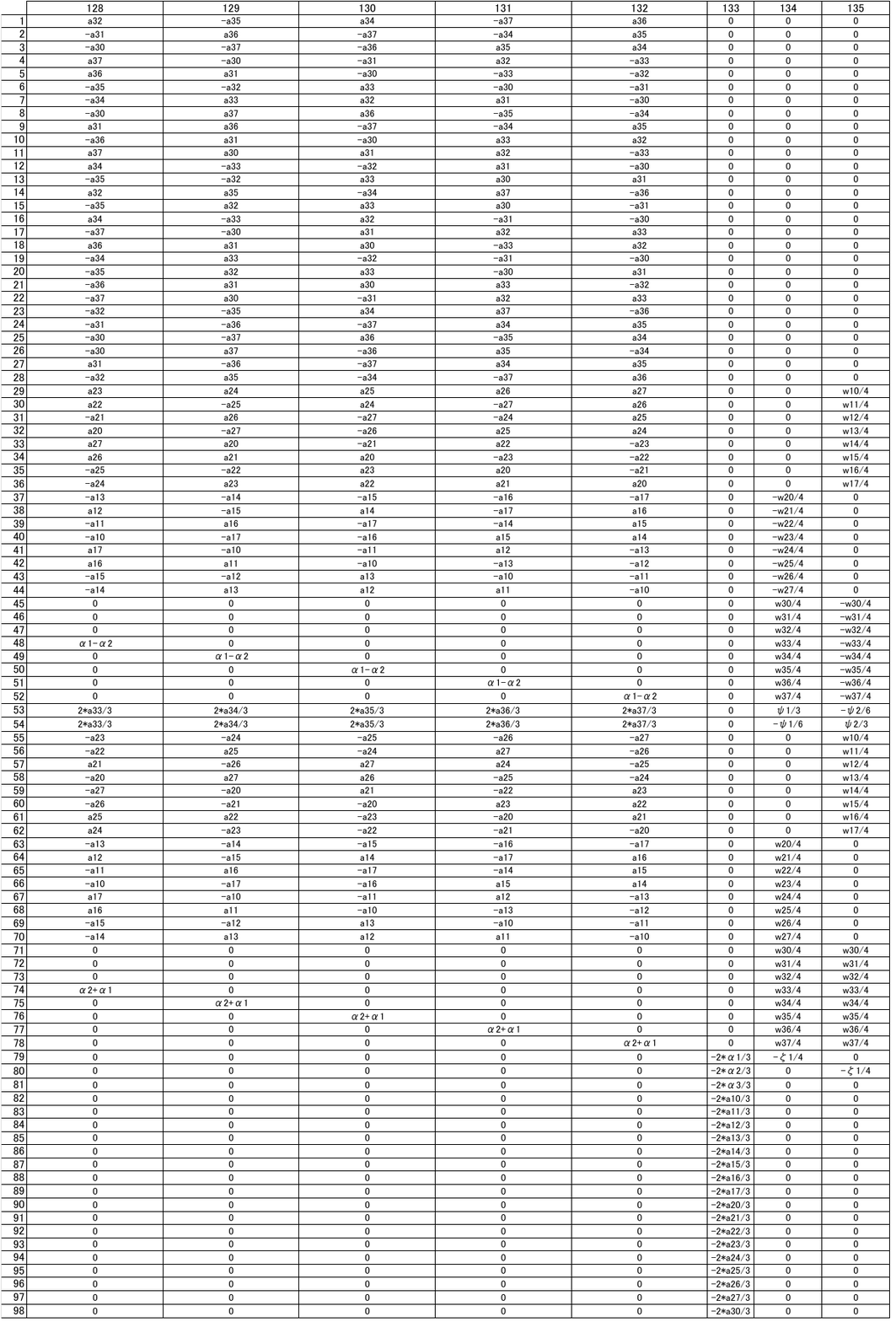}
\end{figure}

\begin{figure}[H]
\centering
\includegraphics[width=12cm, height=19.0cm,page= 16 ]{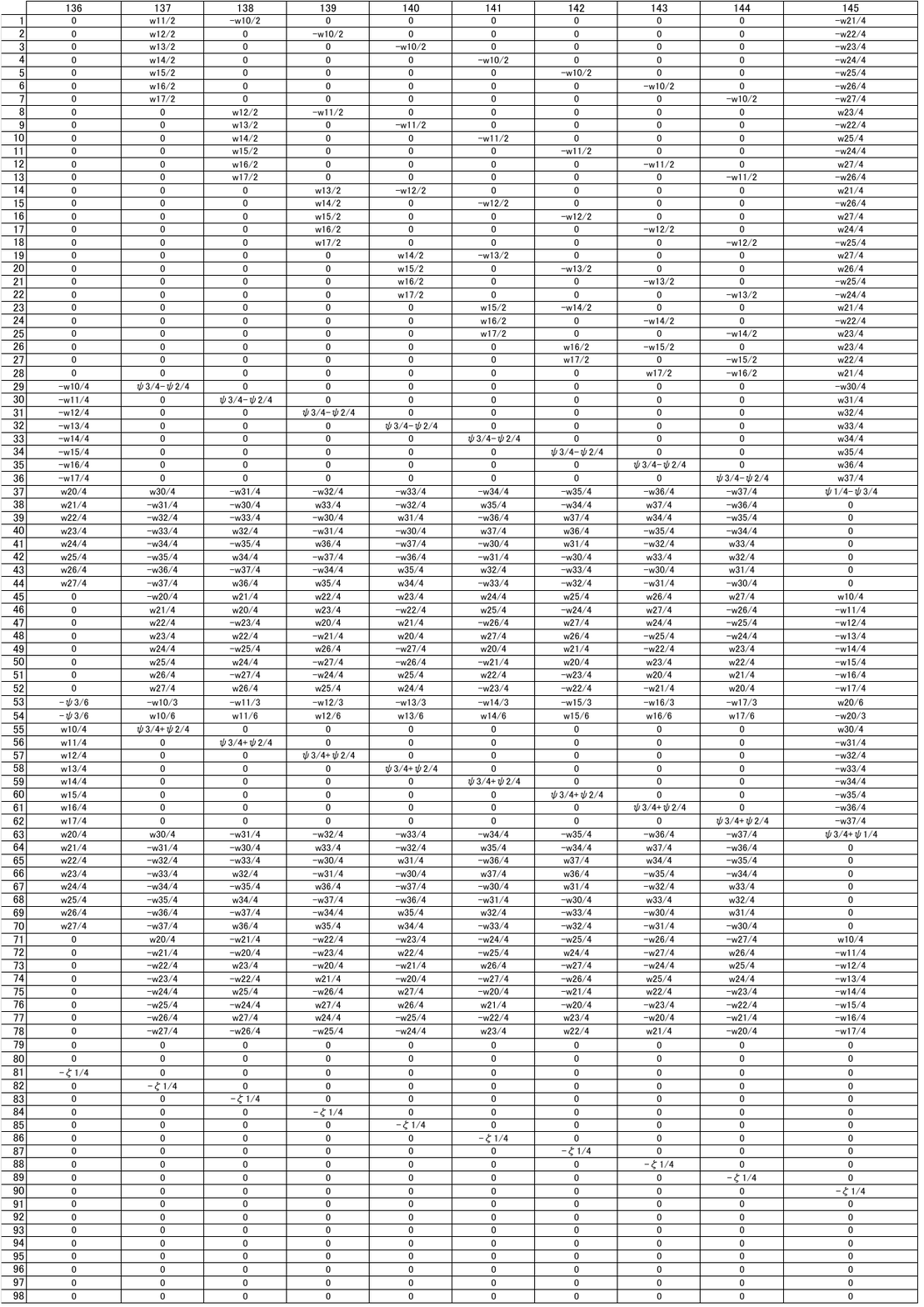}
\end{figure}

\begin{figure}[H]
\centering
\includegraphics[width=12cm, height=19.0cm,page= 17 ]{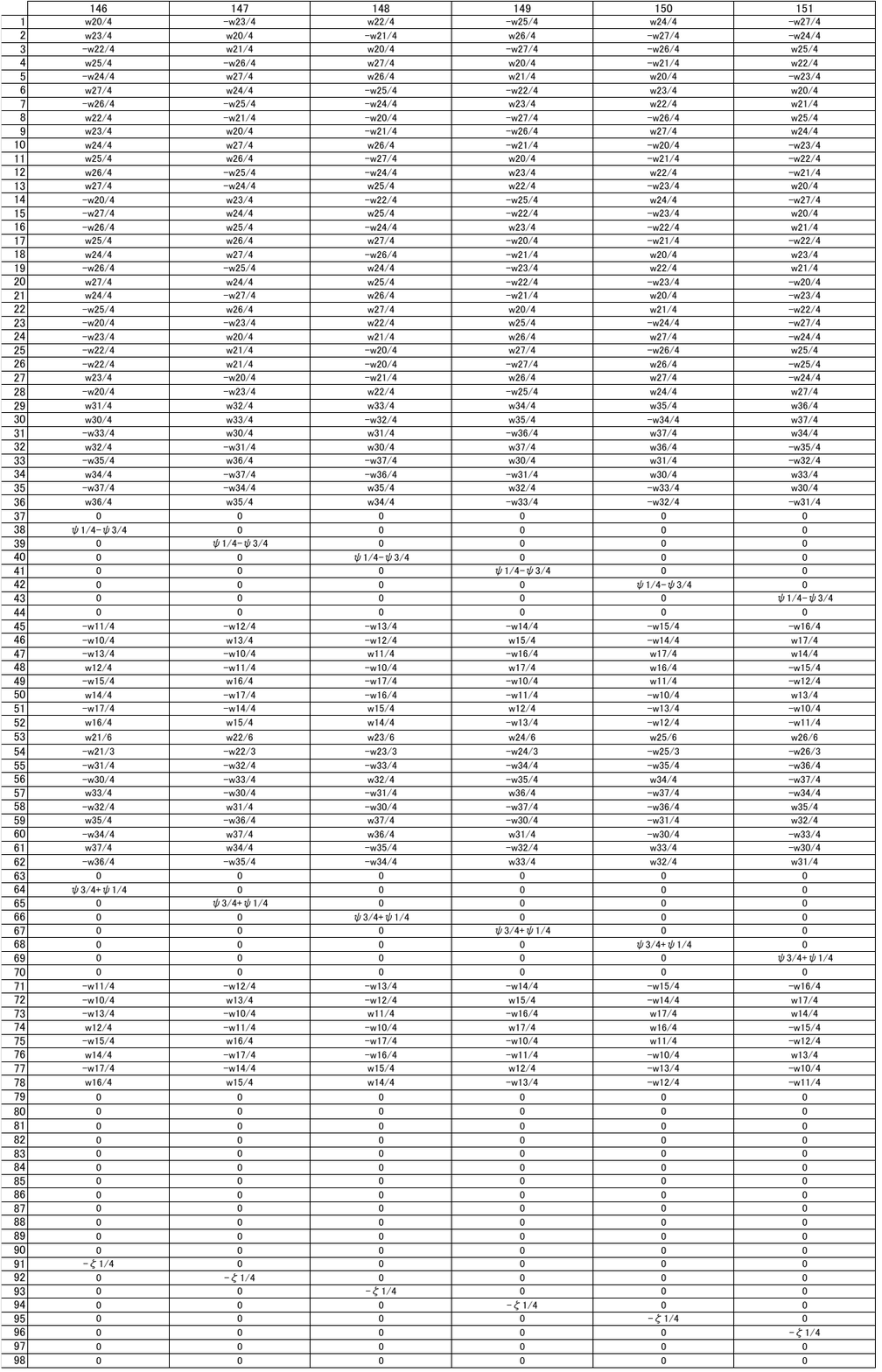}
\end{figure}

\begin{figure}[H]
\centering
\includegraphics[width=12cm, height=19.0cm,page= 18 ]{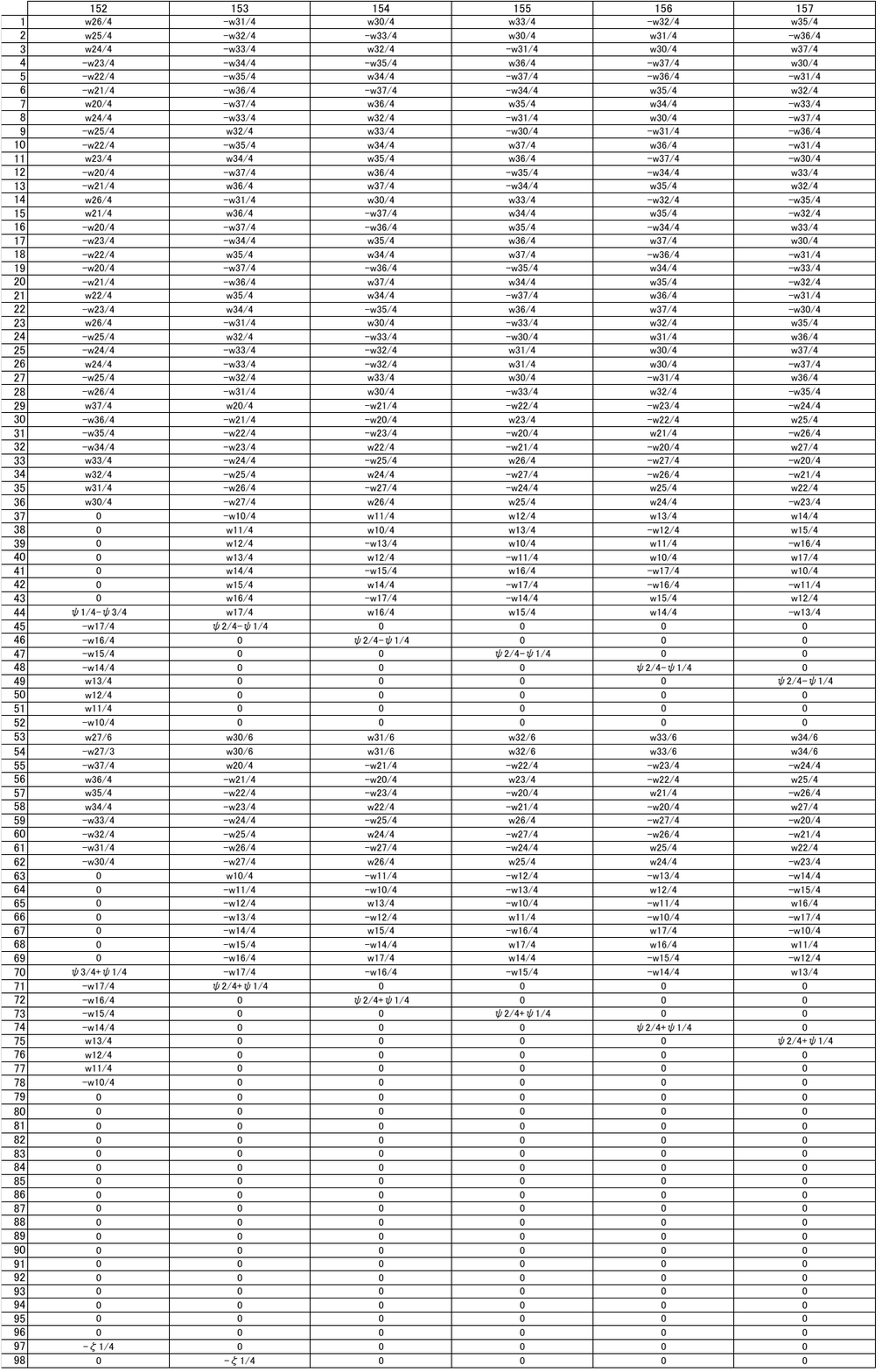}
\end{figure}

\begin{figure}[H]
\centering
\includegraphics[width=12cm, height=19.0cm,page= 19 ]{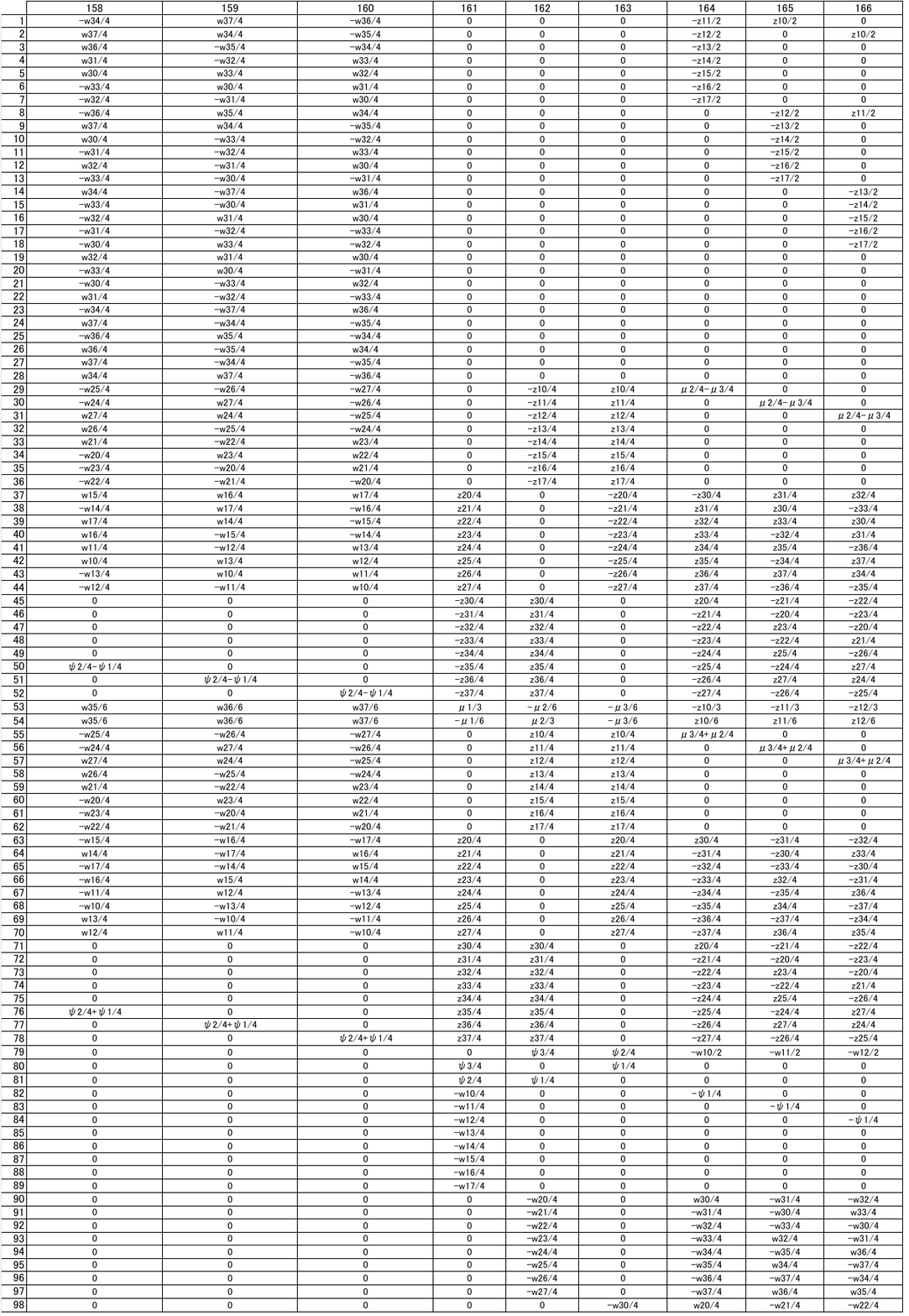}
\end{figure}

\begin{figure}[H]
\centering
\includegraphics[width=12cm, height=19.0cm,page= 20 ]{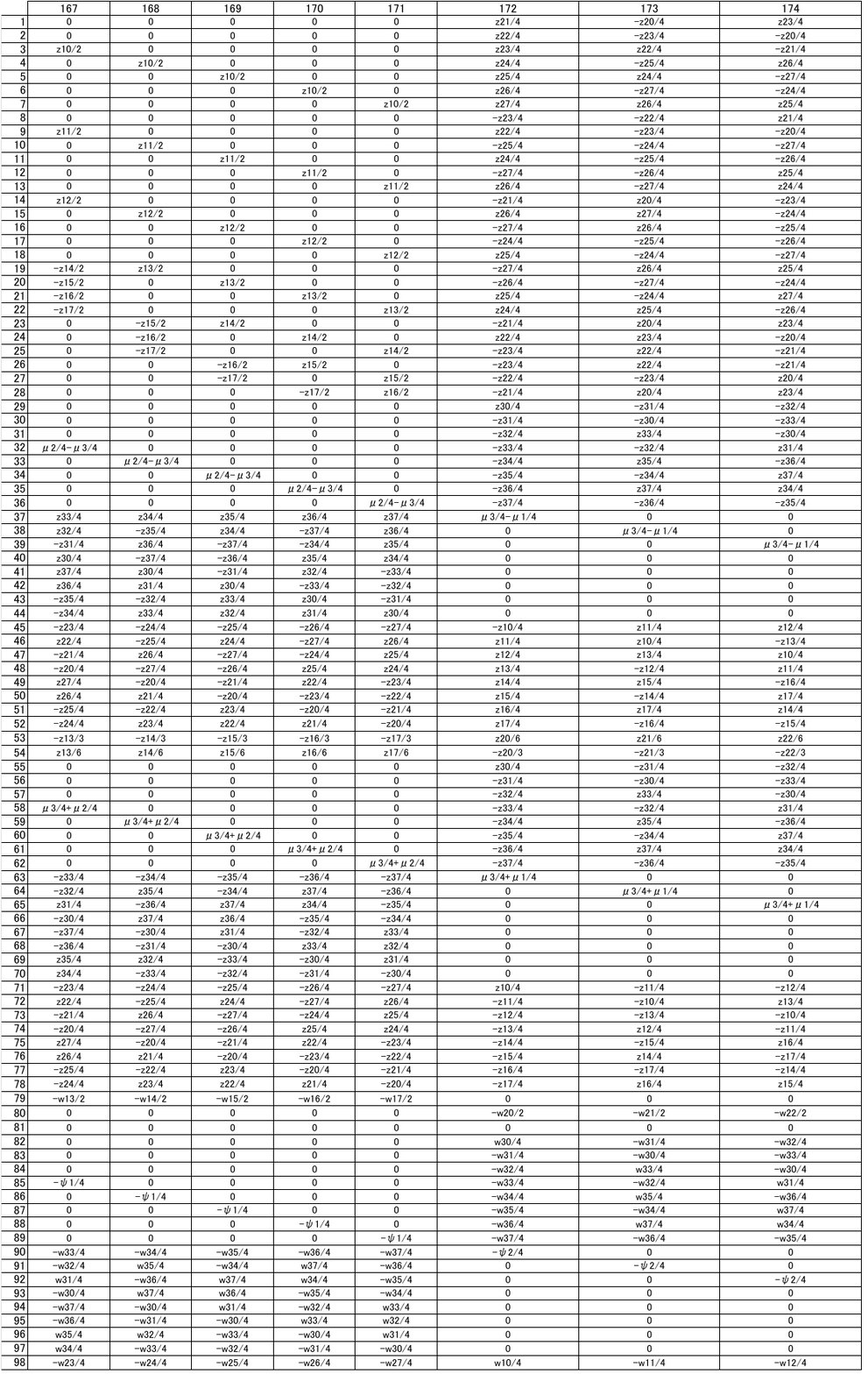}
\end{figure}

\begin{figure}[H]
\centering
\includegraphics[width=12cm, height=19.0cm,page= 21 ]{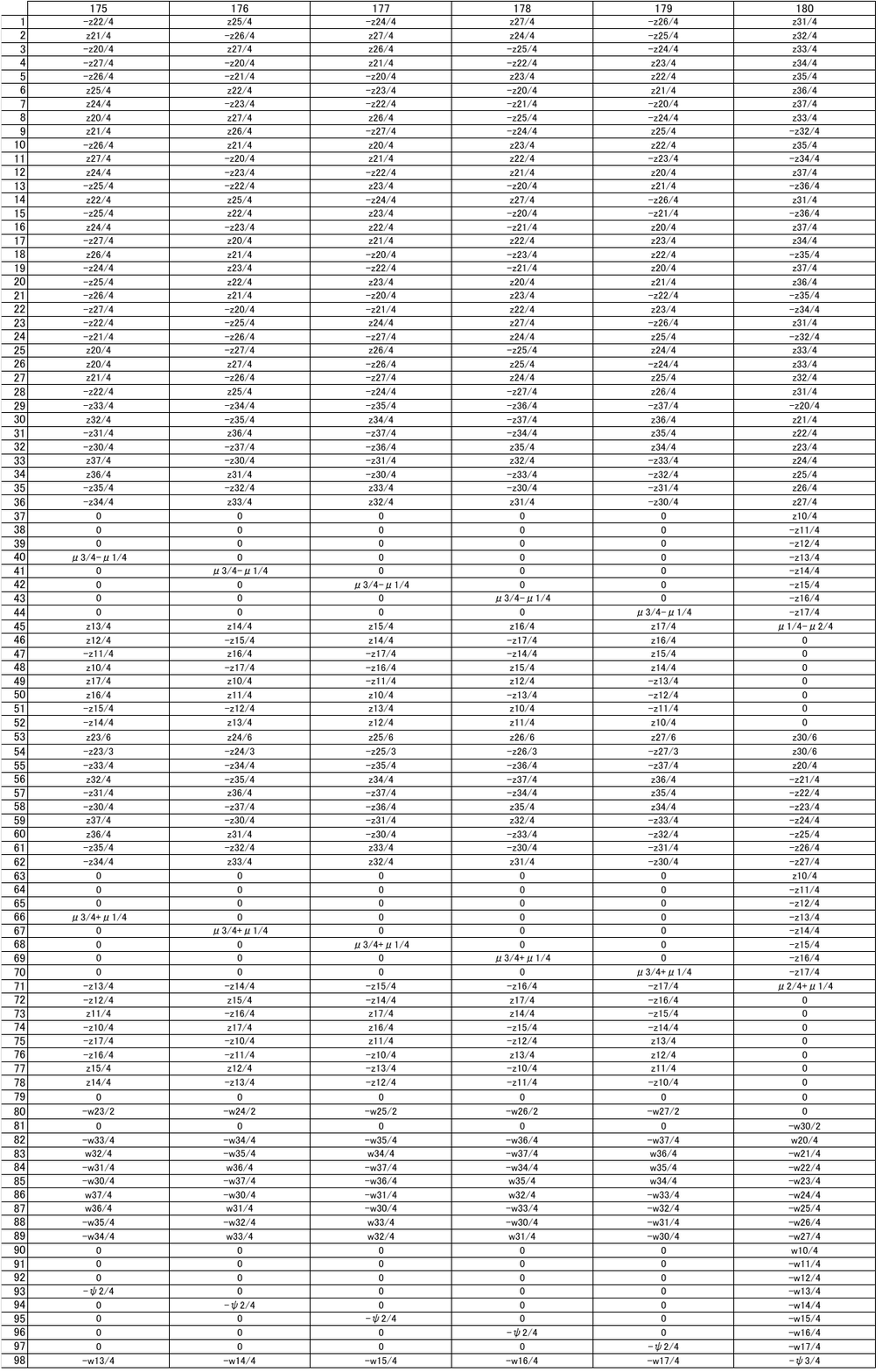}
\end{figure}

\begin{figure}[H]
\centering
\includegraphics[width=12cm, height=19.0cm,page= 22 ]{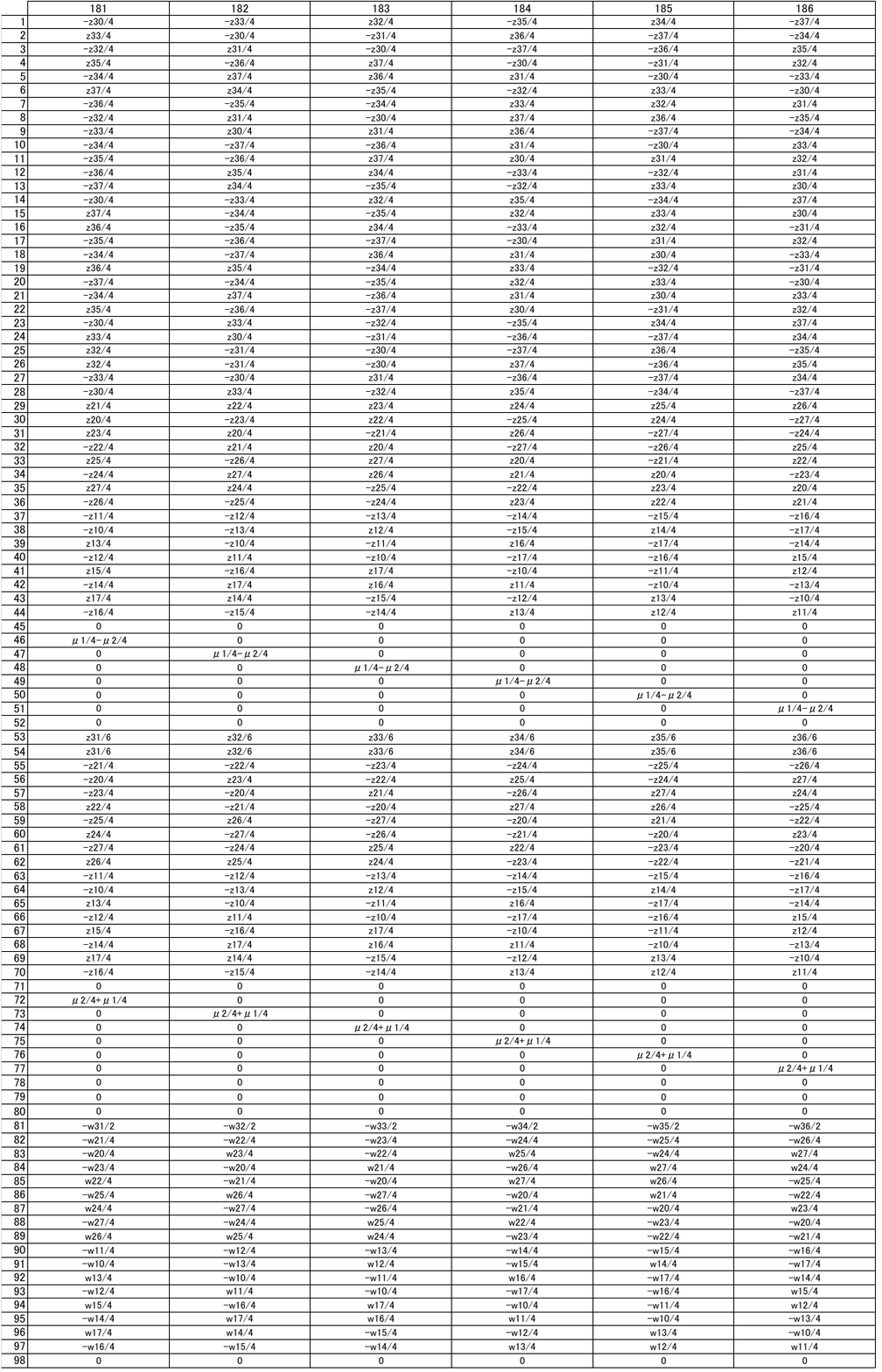}
\end{figure}

\begin{figure}[H]
\centering
\includegraphics[width=12cm, height=19.0cm,page= 23 ]{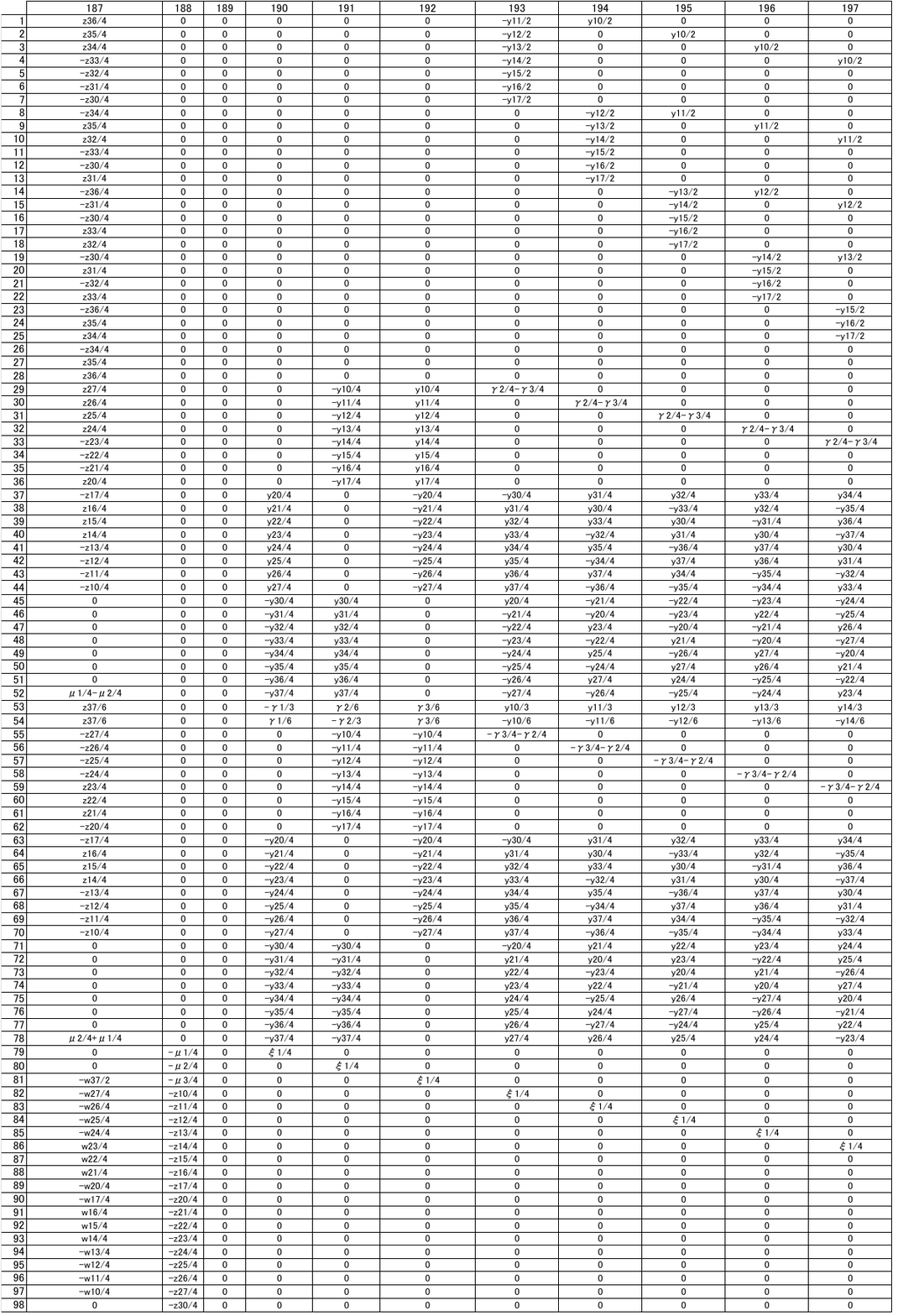}
\end{figure}

\begin{figure}[H]
\centering
\includegraphics[width=12cm, height=19.0cm,page= 24 ]{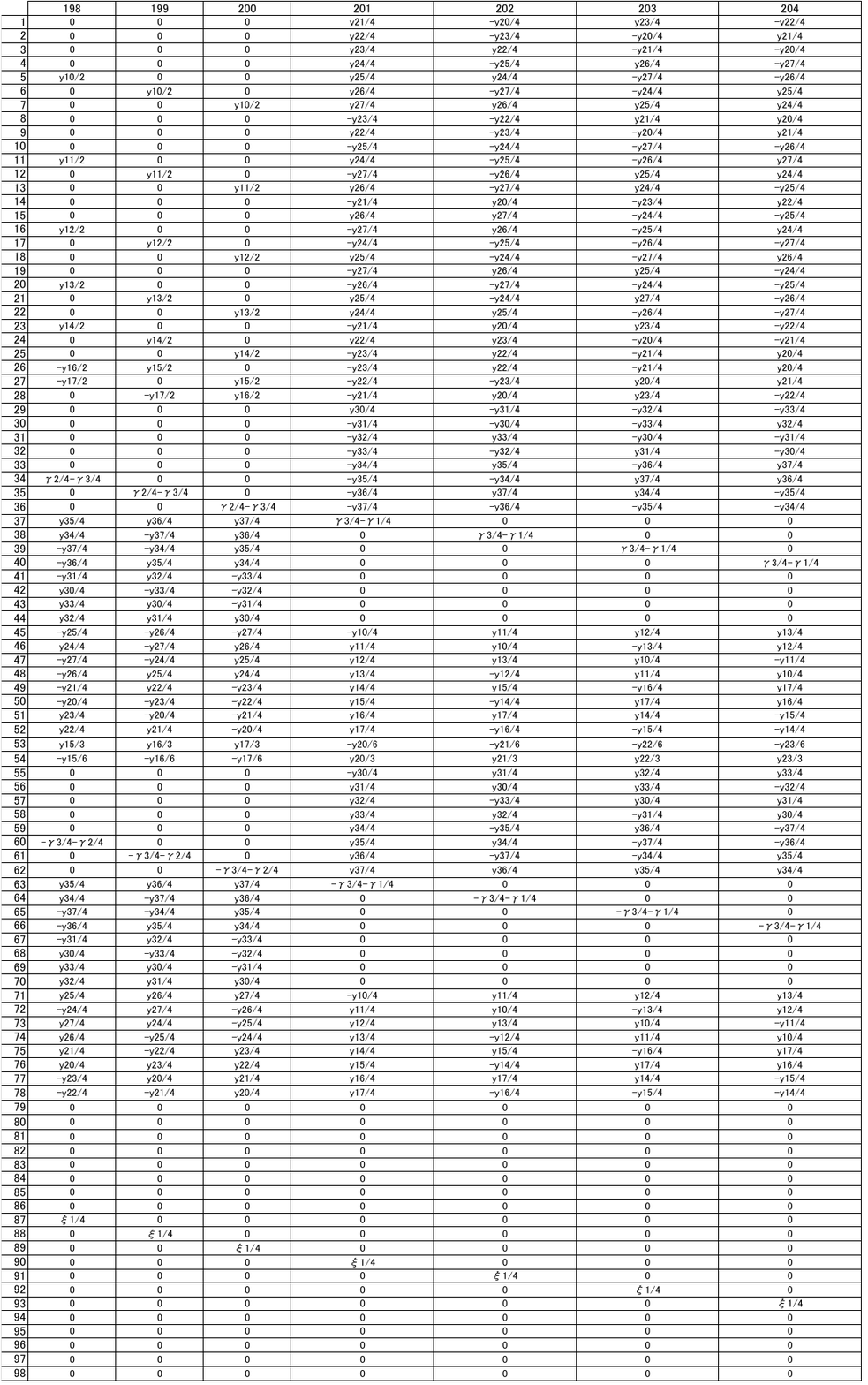}
\end{figure}

\begin{figure}[H]
\centering
\includegraphics[width=12cm, height=19.0cm,page= 25 ]{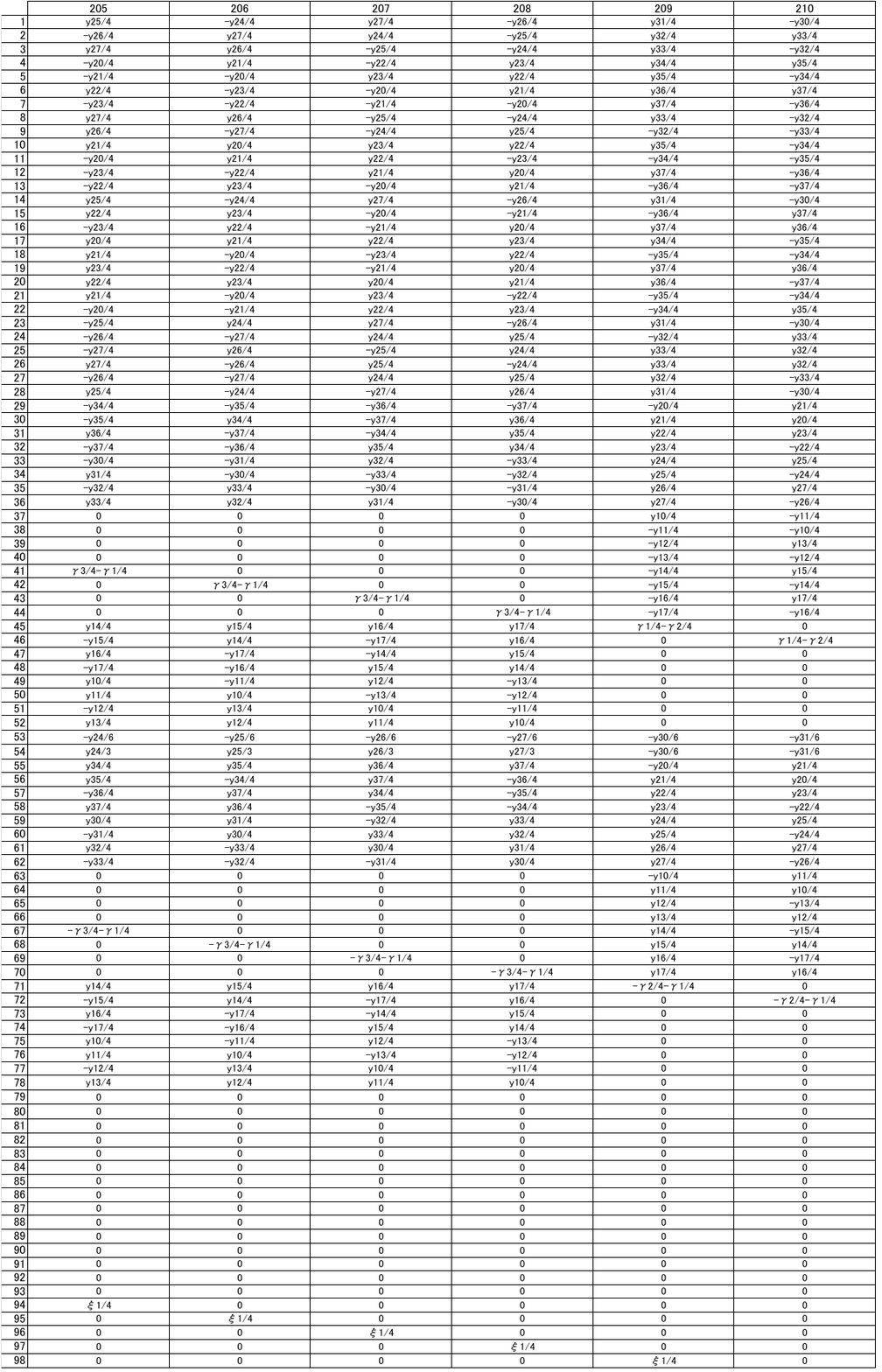}
\end{figure}

\begin{figure}[H]
\centering
\includegraphics[width=12cm, height=19.0cm,page= 26 ]{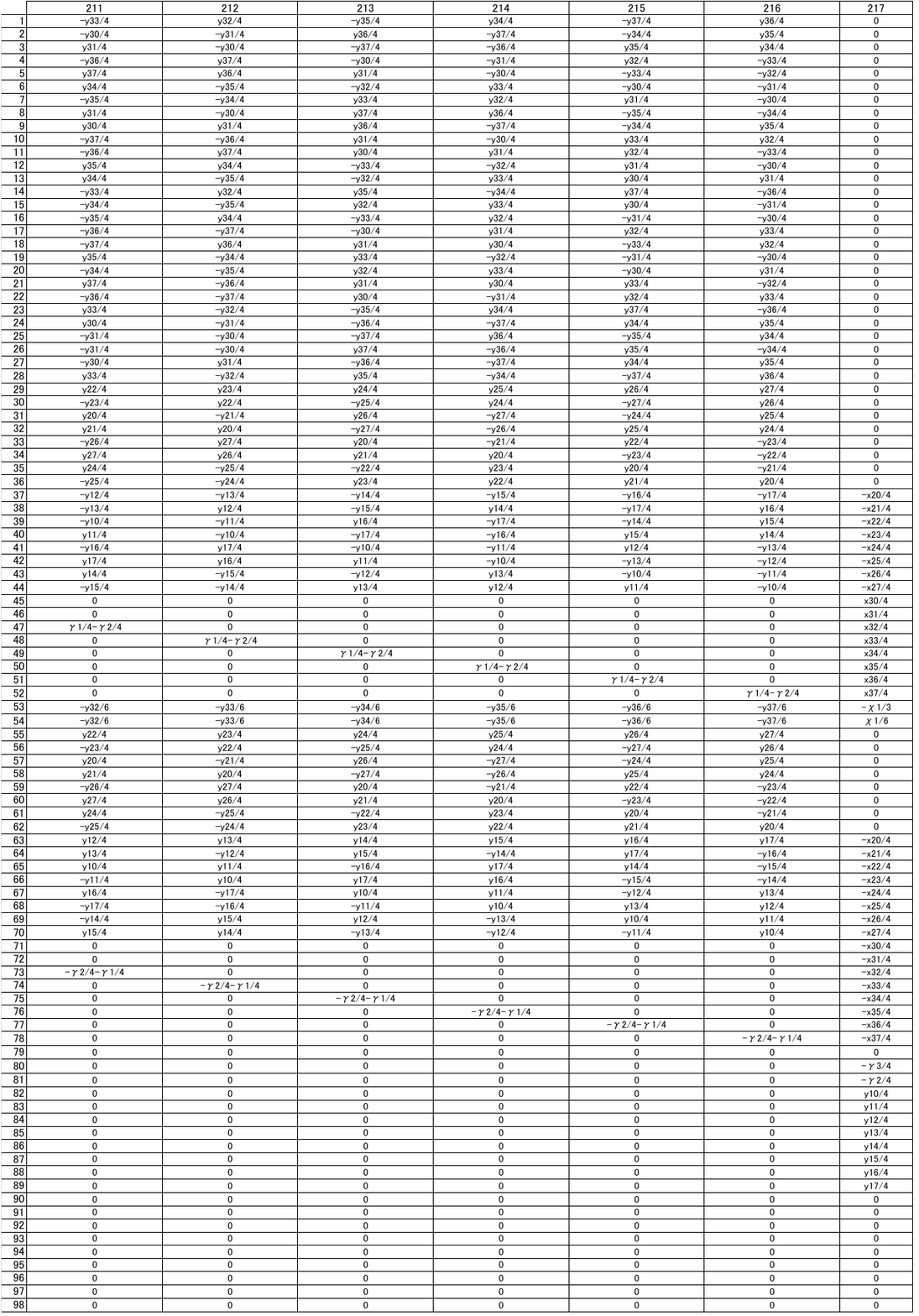}
\end{figure}

\begin{figure}[H]
\centering
\includegraphics[width=12cm, height=19.0cm,page= 27 ]{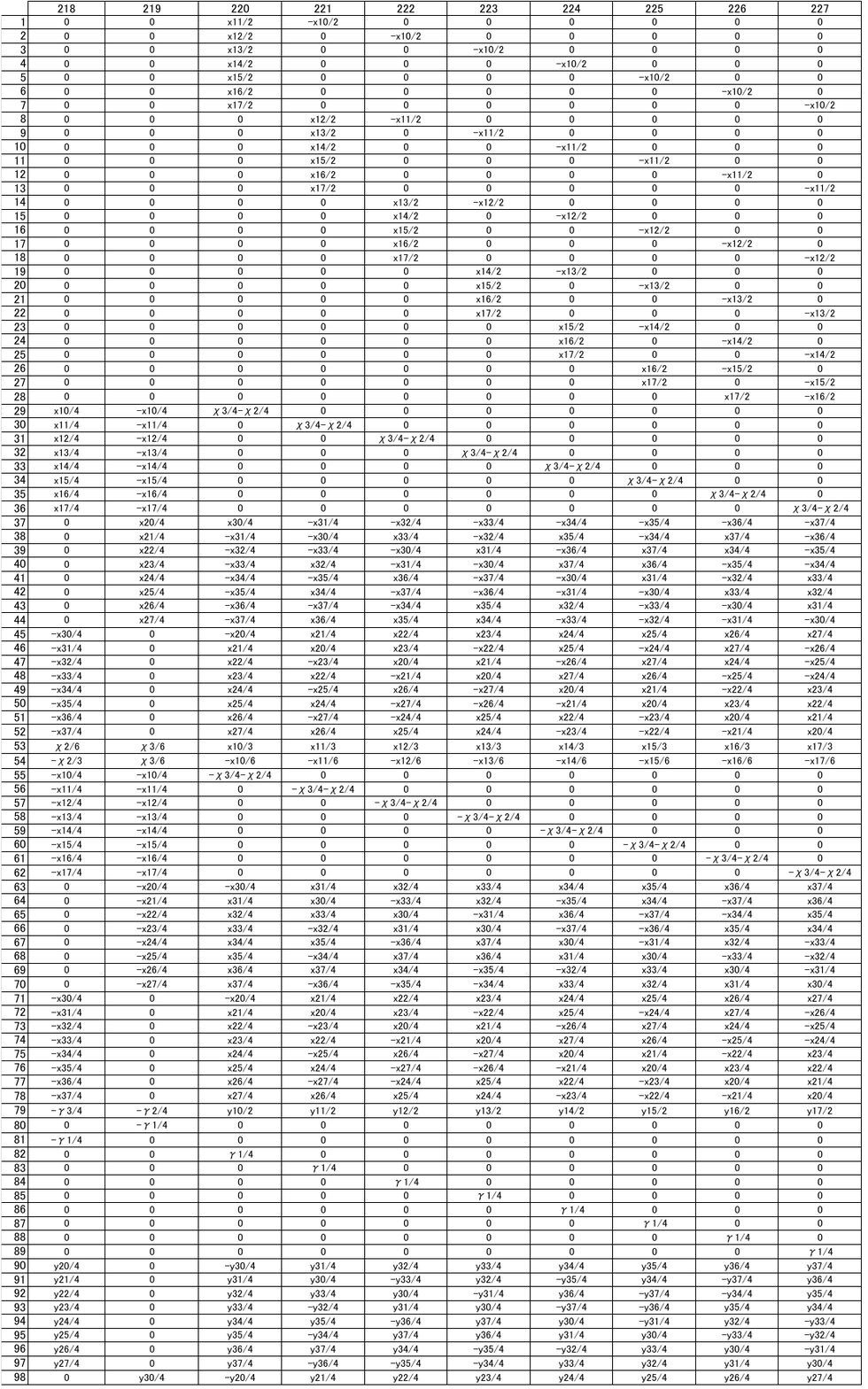}
\end{figure}

\begin{figure}[H]
\centering
\includegraphics[width=12cm, height=19.0cm,page= 28 ]{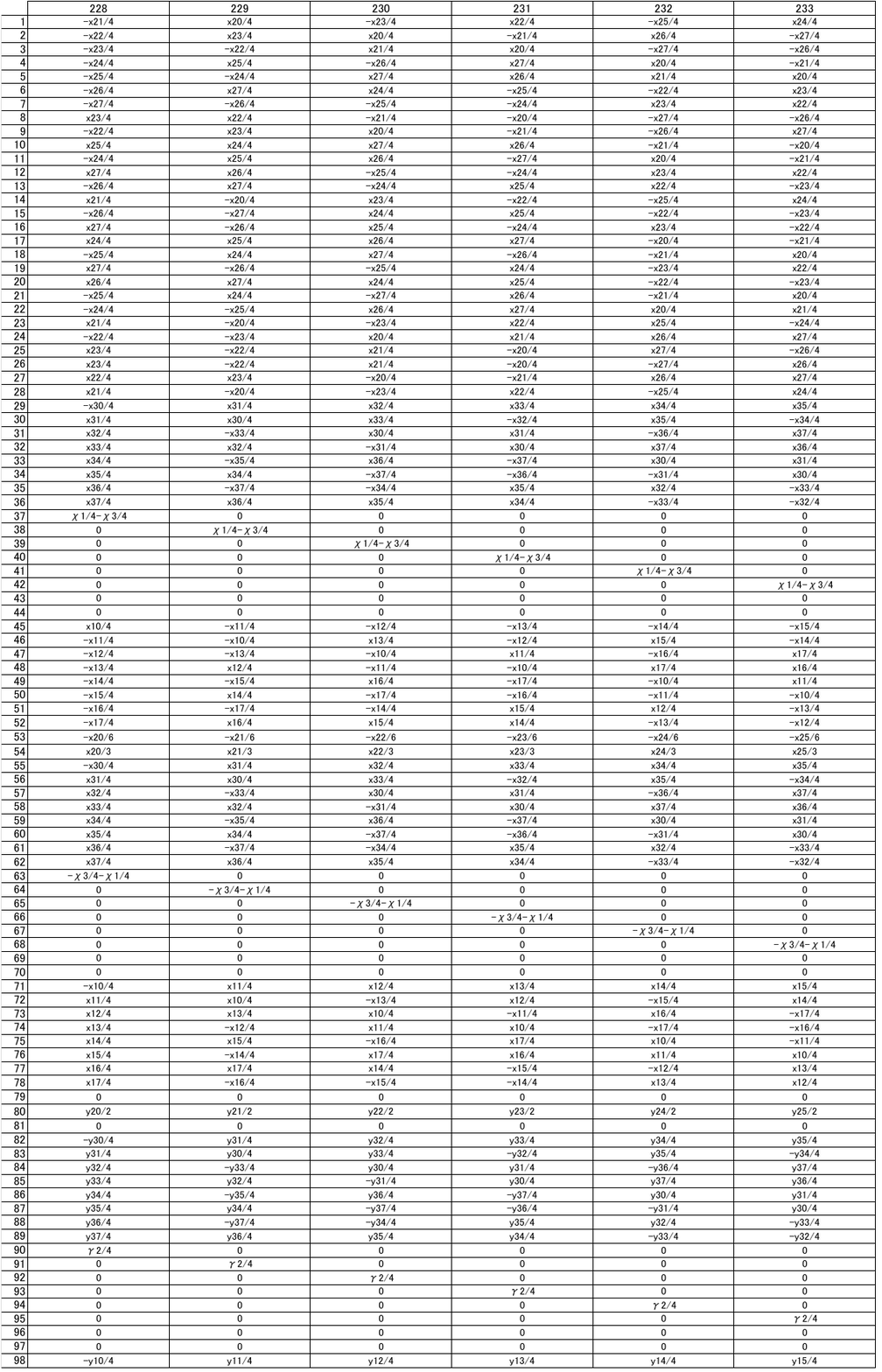}
\end{figure}

\begin{figure}[H]
\centering
\includegraphics[width=12cm, height=19.0cm,page= 29 ]{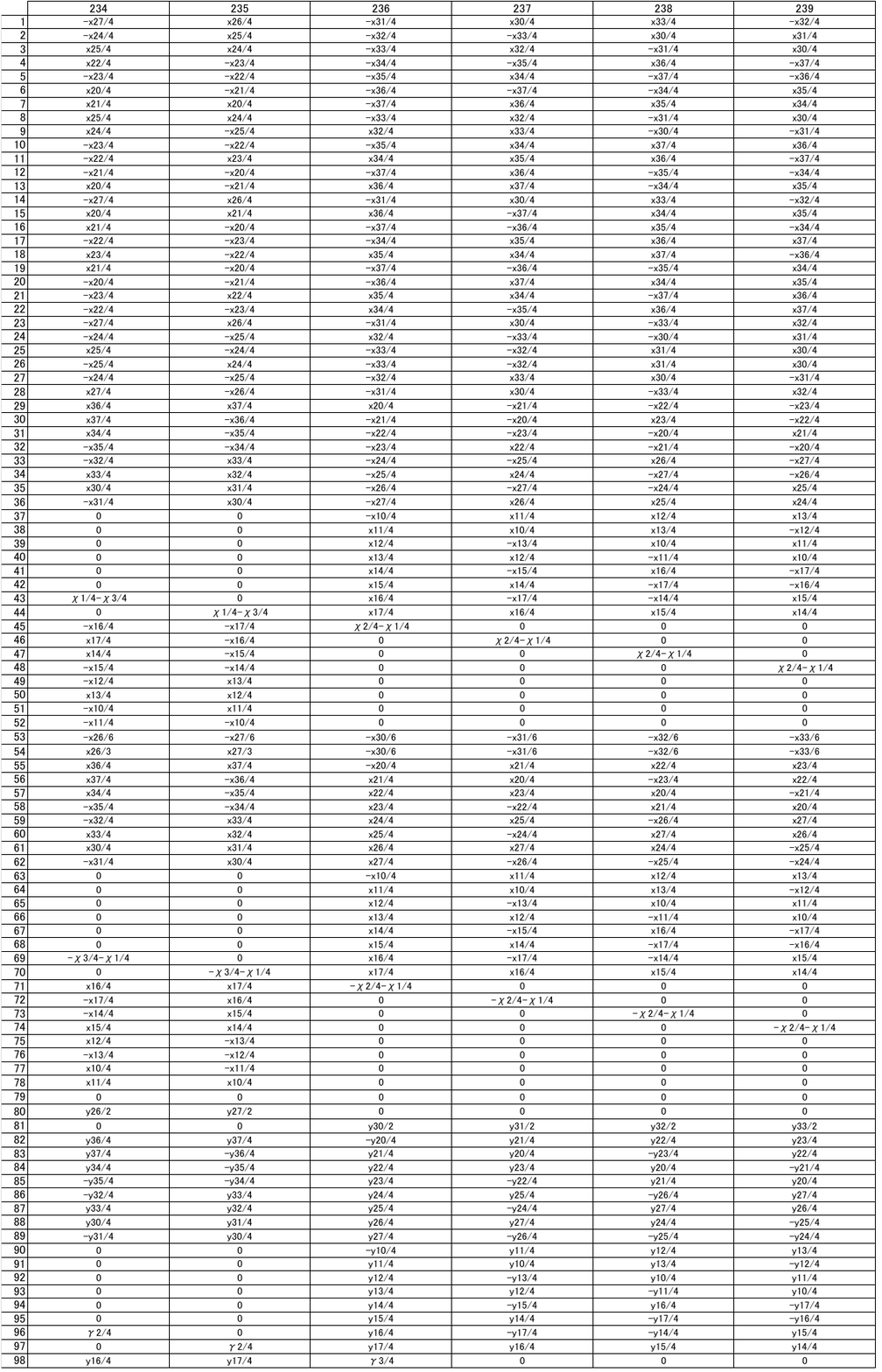}
\end{figure}

\begin{figure}[H]
\centering
\includegraphics[width=12cm, height=19.0cm,page= 30 ]{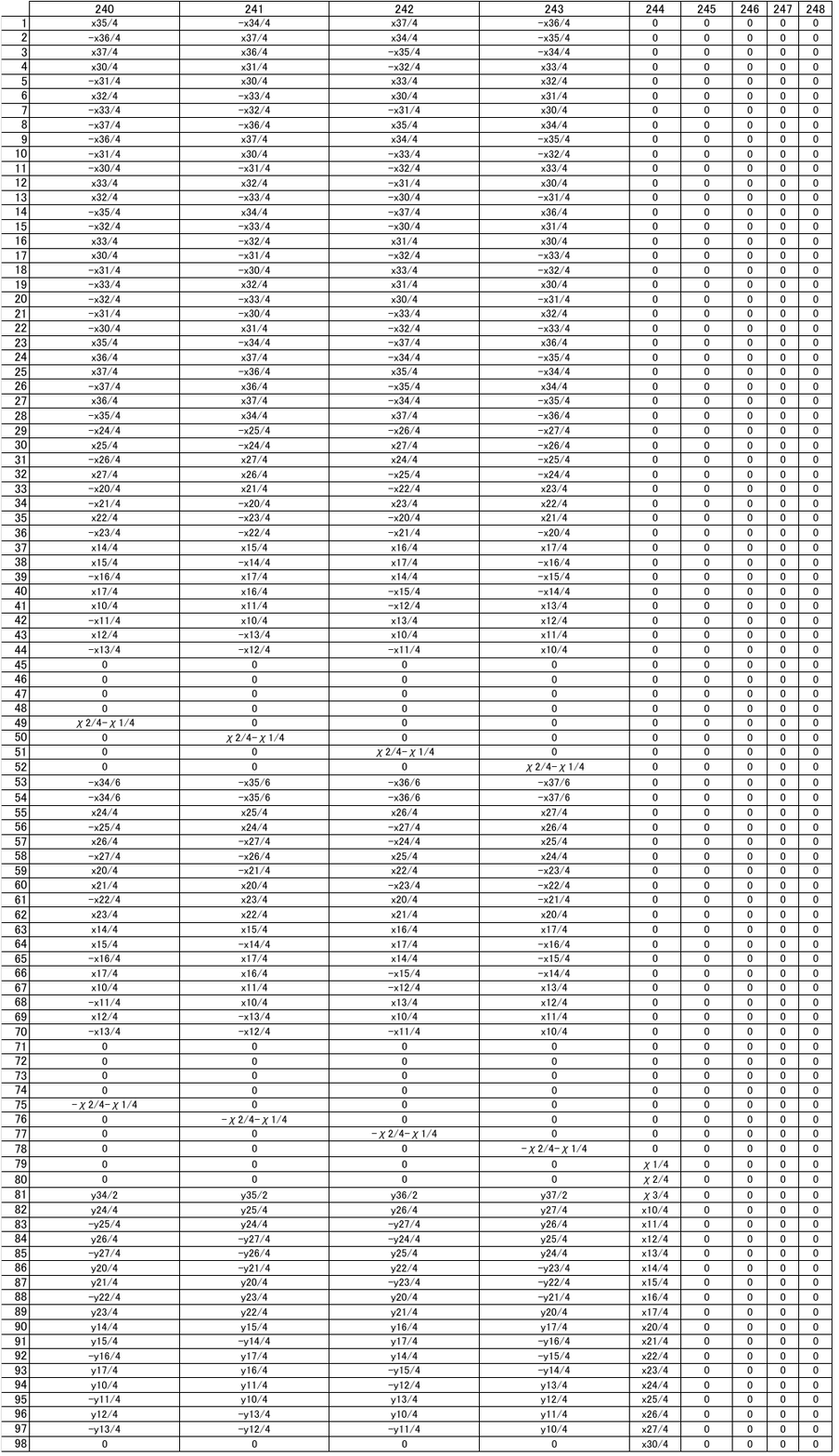}
\end{figure}

\begin{figure}[H]
\centering
\includegraphics[width=12cm, height=19.0cm,page= 31 ]{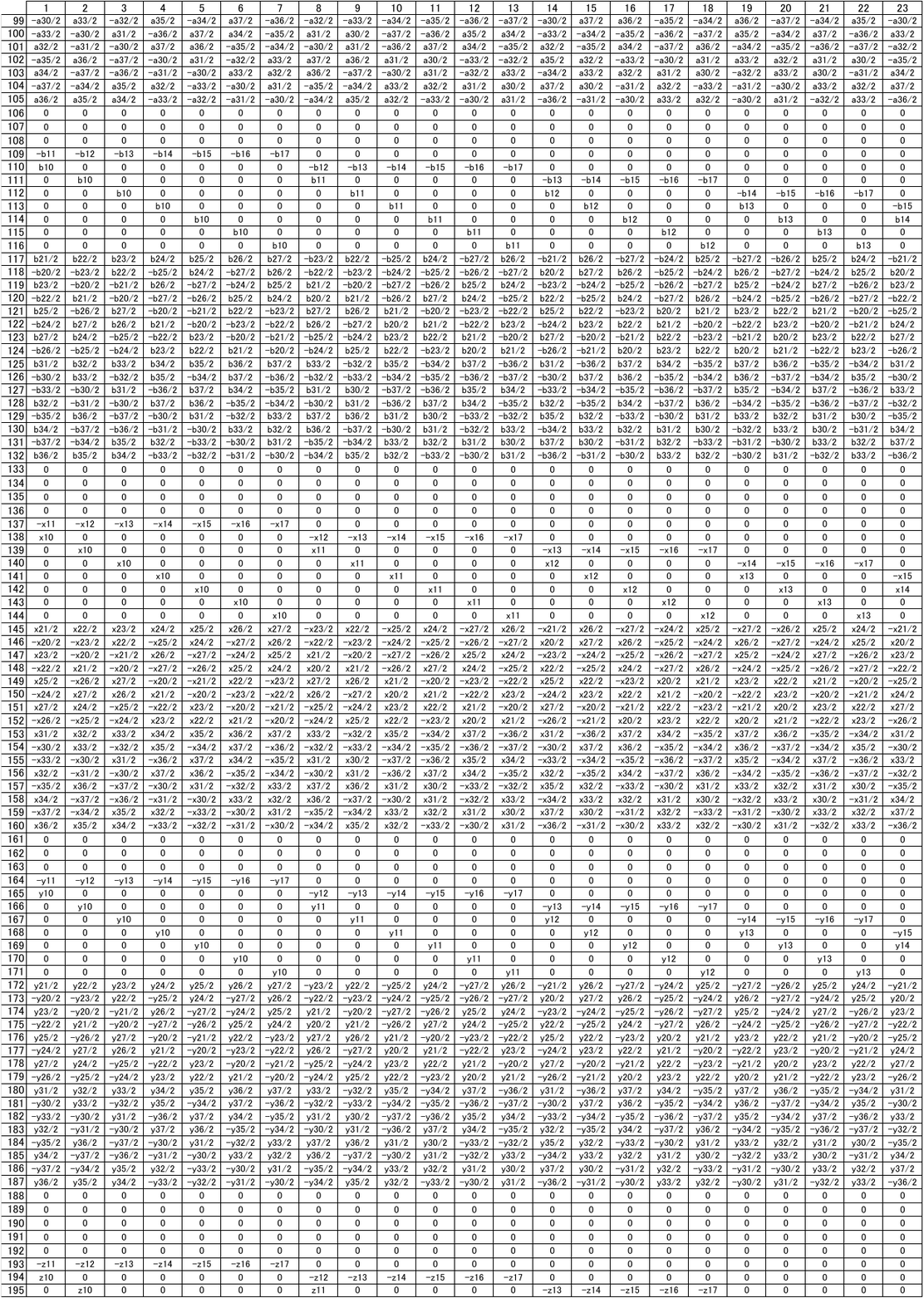}
\end{figure}

\begin{figure}[H]
\centering
\includegraphics[width=12cm, height=19.0cm,page= 32 ]{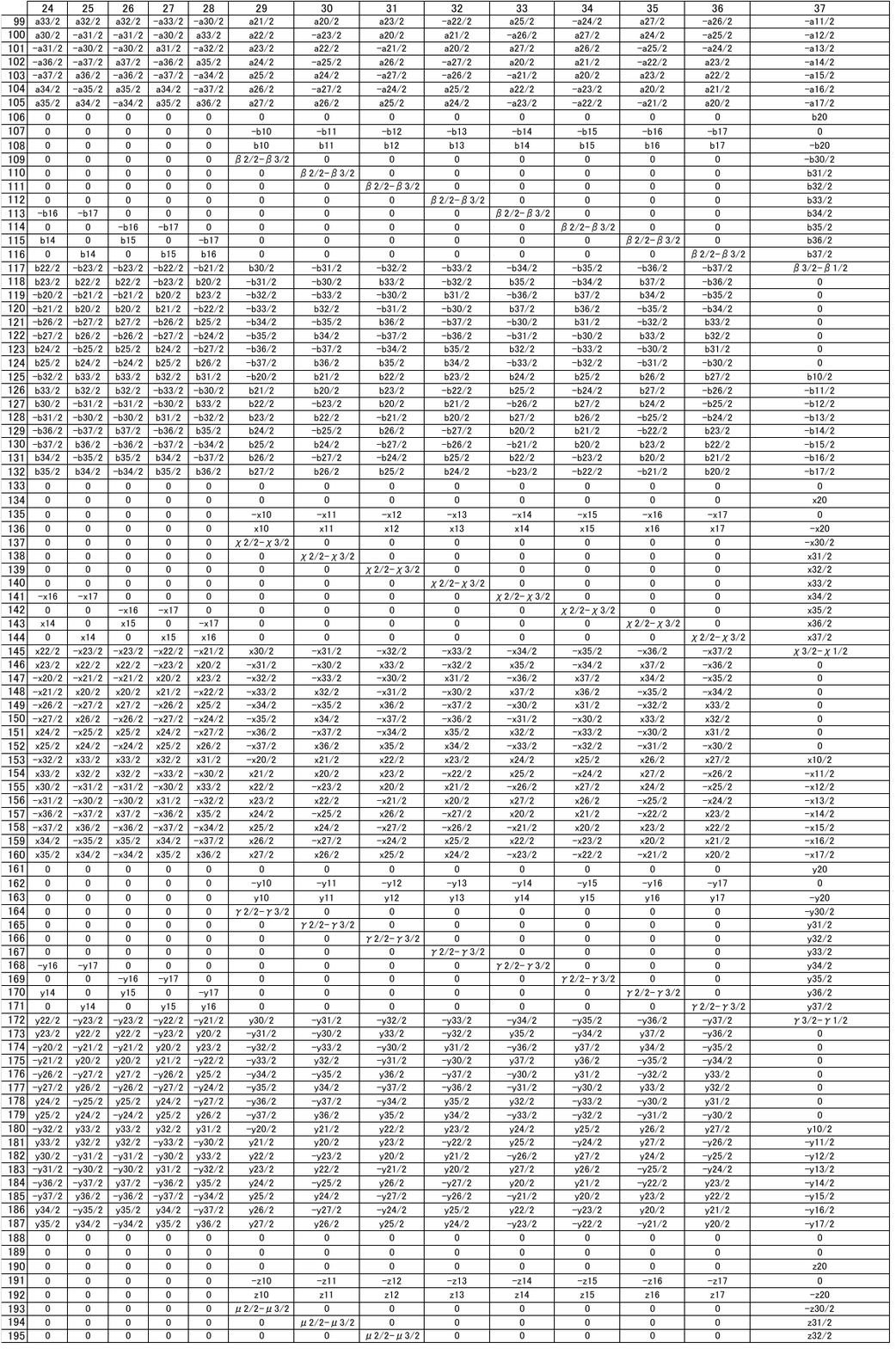}
\end{figure}

\begin{figure}[H]
\centering
\includegraphics[width=12cm, height=19.0cm,page= 33 ]{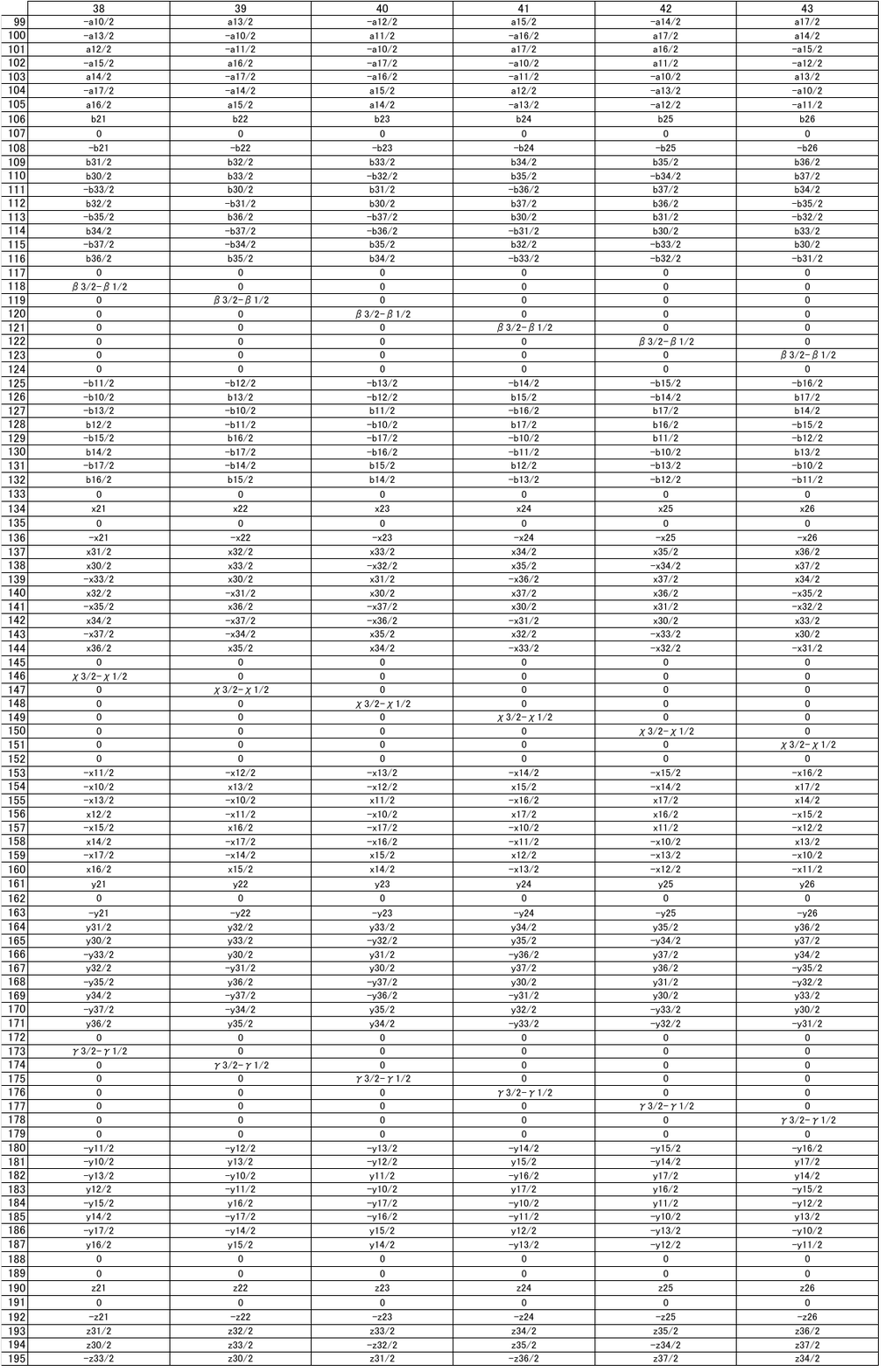}
\end{figure}

\begin{figure}[H]
\centering
\includegraphics[width=12cm, height=19.0cm,page= 34 ]{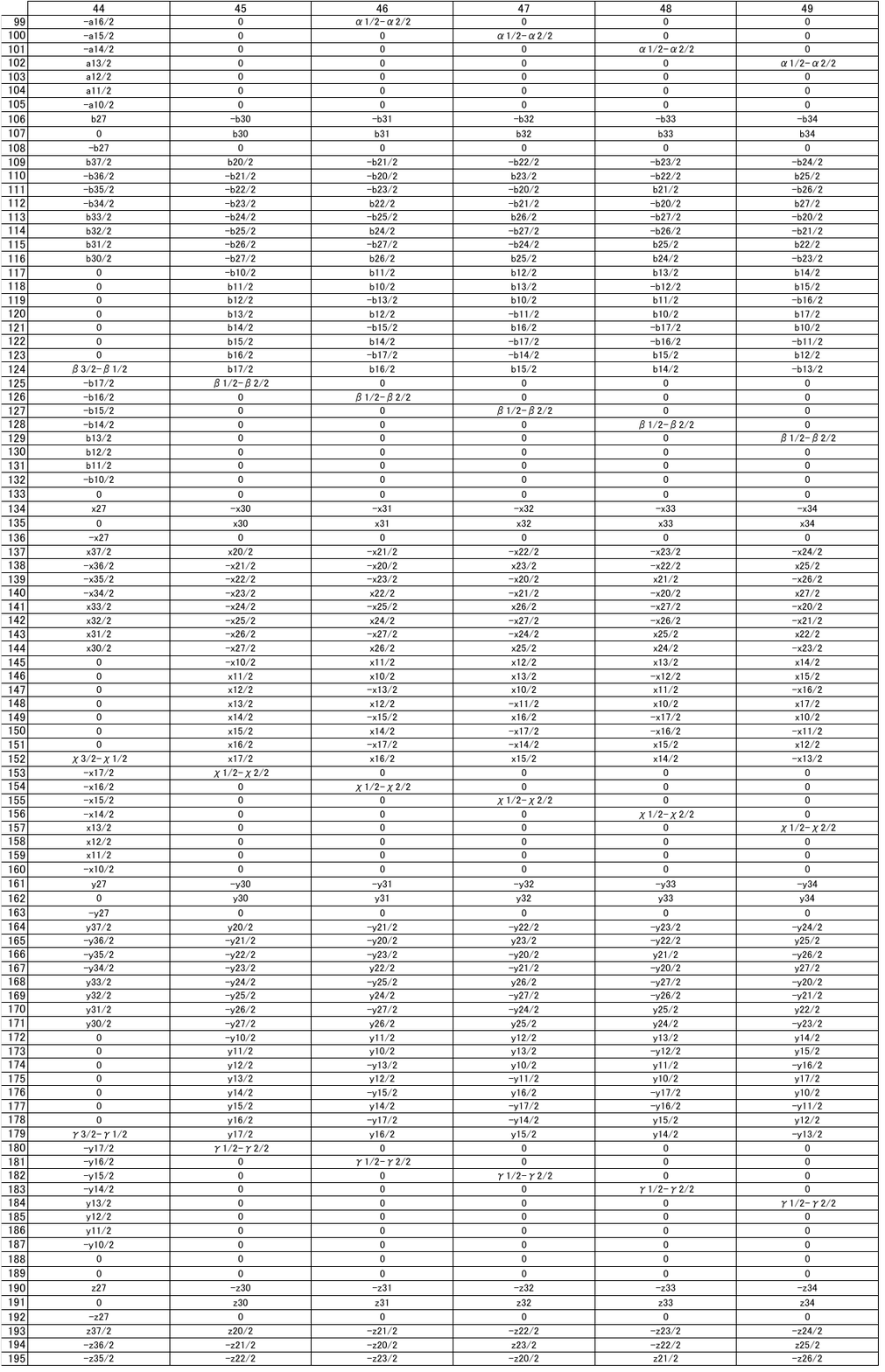}
\end{figure}

\begin{figure}[H]
\centering
\includegraphics[width=12cm, height=19.0cm,page= 35 ]{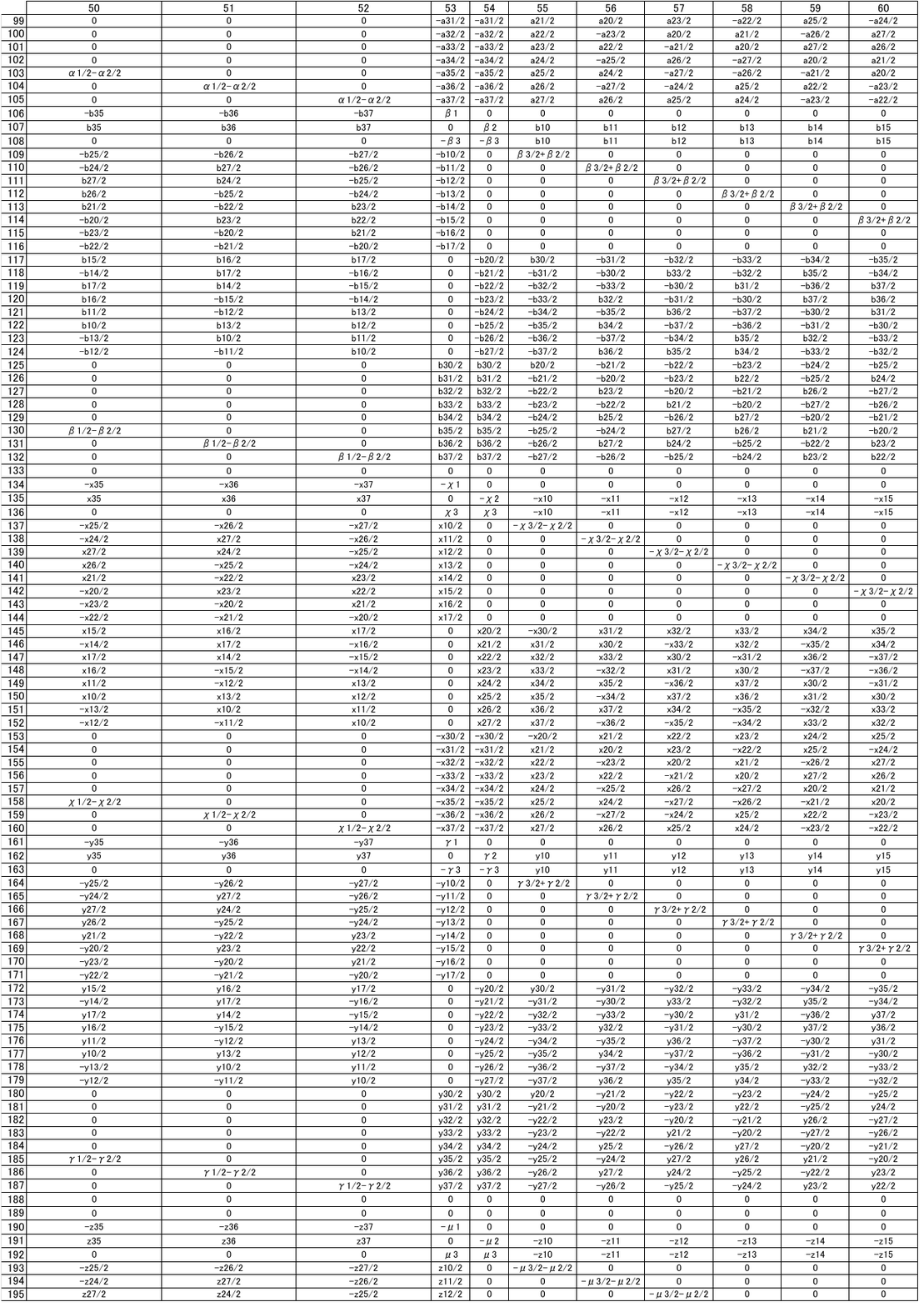}
\end{figure}

\begin{figure}[H]
\centering
\includegraphics[width=12cm, height=19.0cm,page= 36 ]{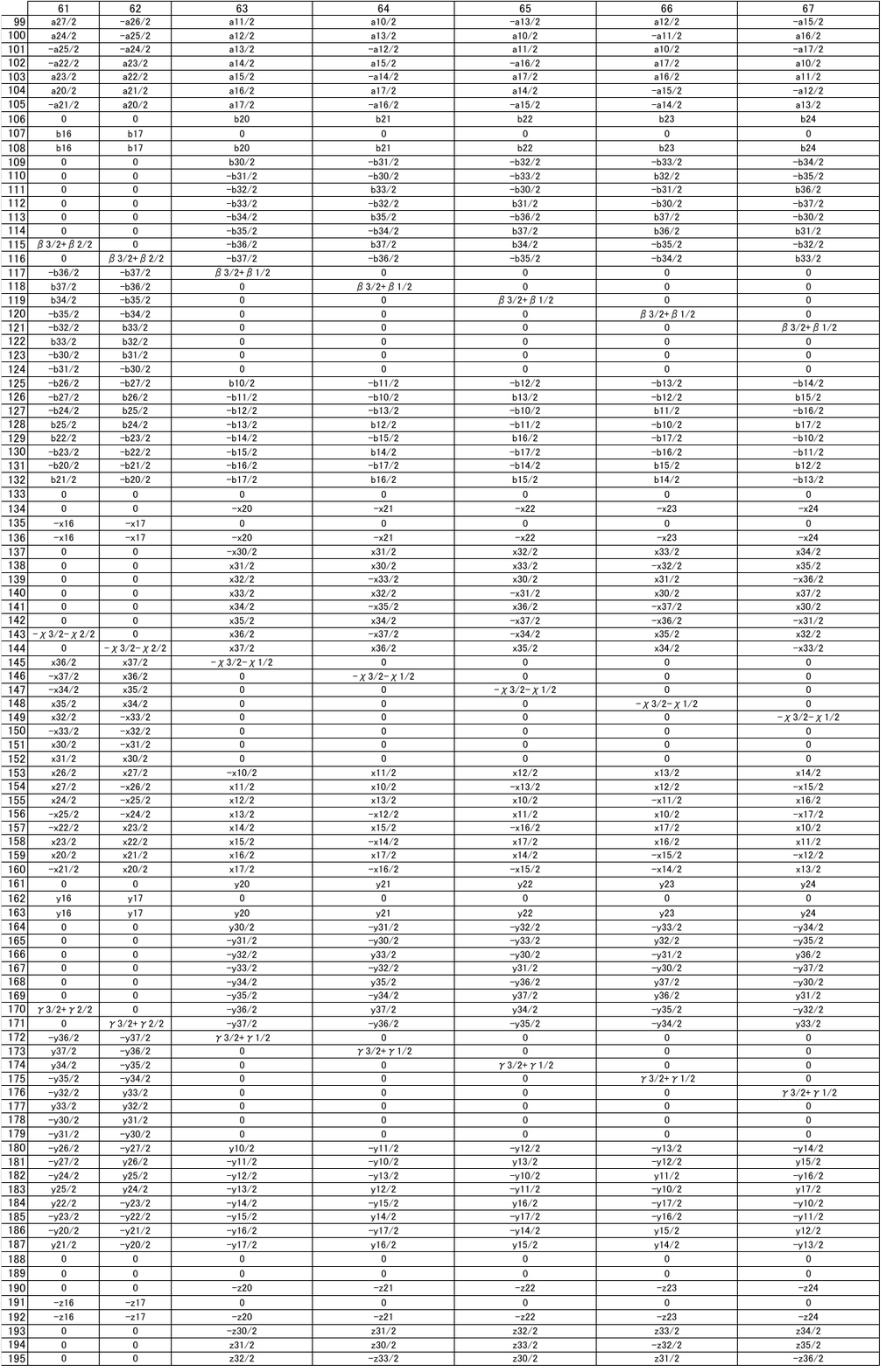}
\end{figure}

\begin{figure}[H]
\centering
\includegraphics[width=12cm, height=19.0cm,page= 37 ]{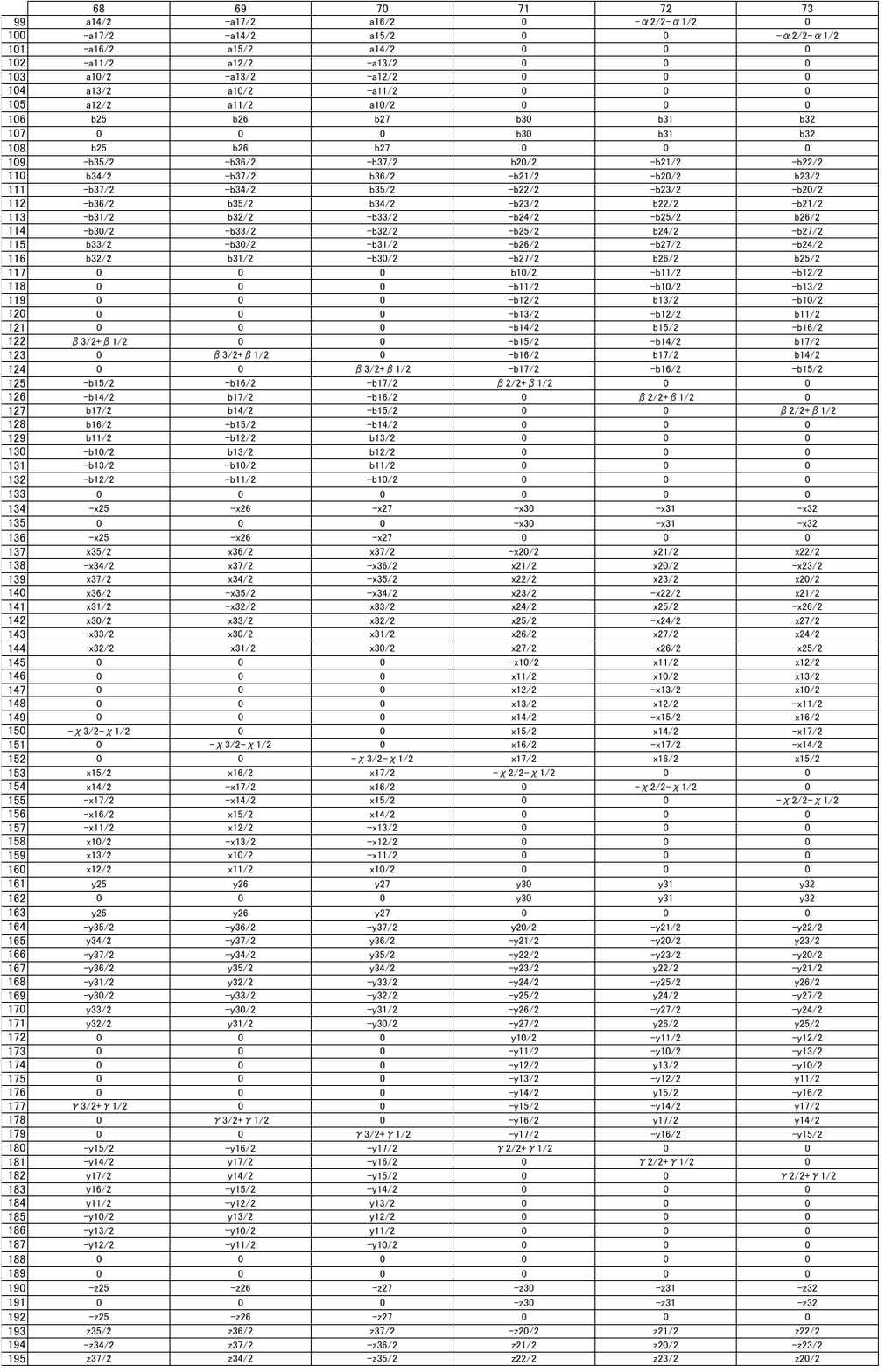}
\end{figure}

\begin{figure}[H]
\centering
\includegraphics[width=12cm, height=19.0cm,page= 38 ]{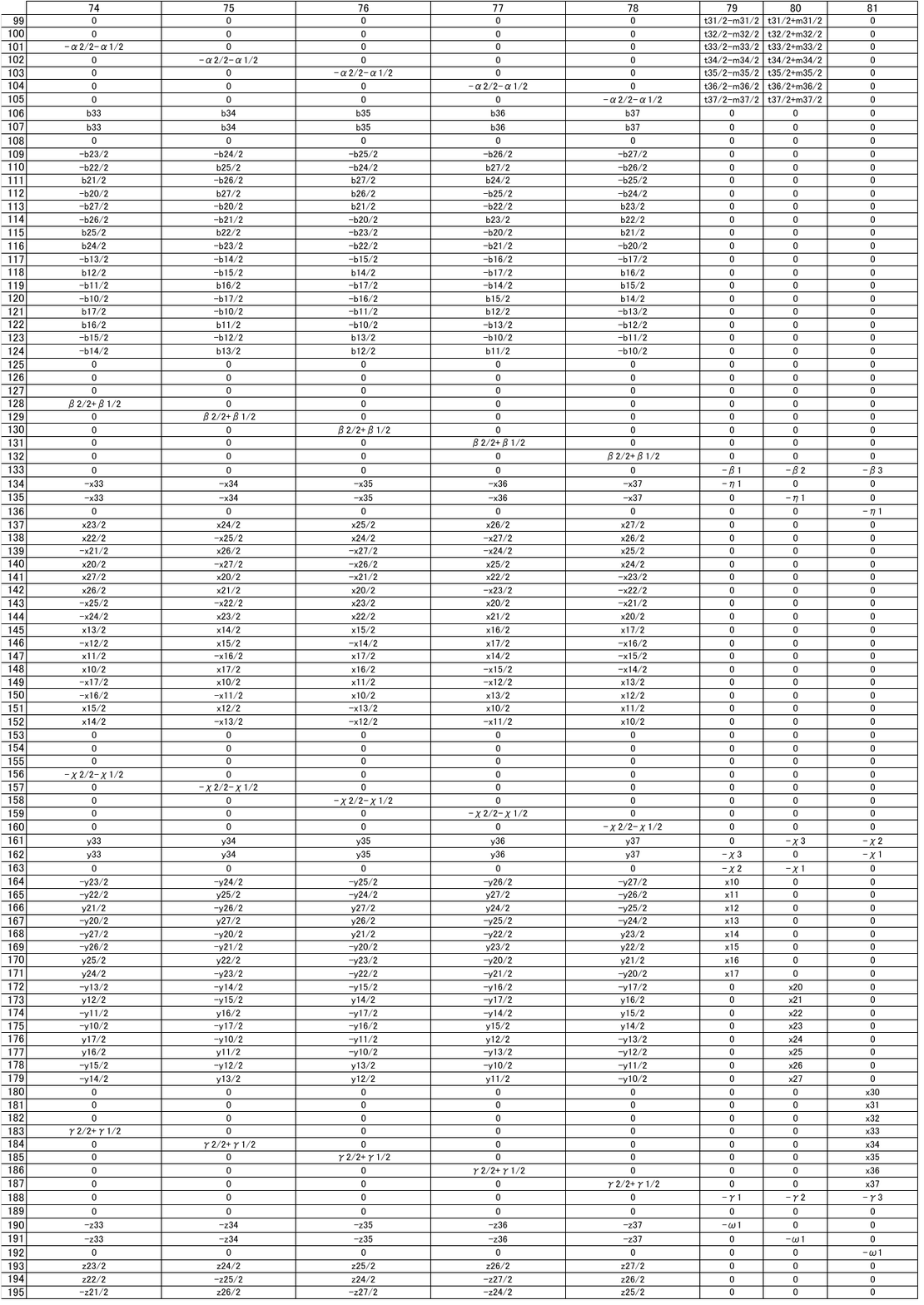}
\end{figure}

\begin{figure}[H]
\centering
\includegraphics[width=12cm, height=19.0cm,page= 39 ]{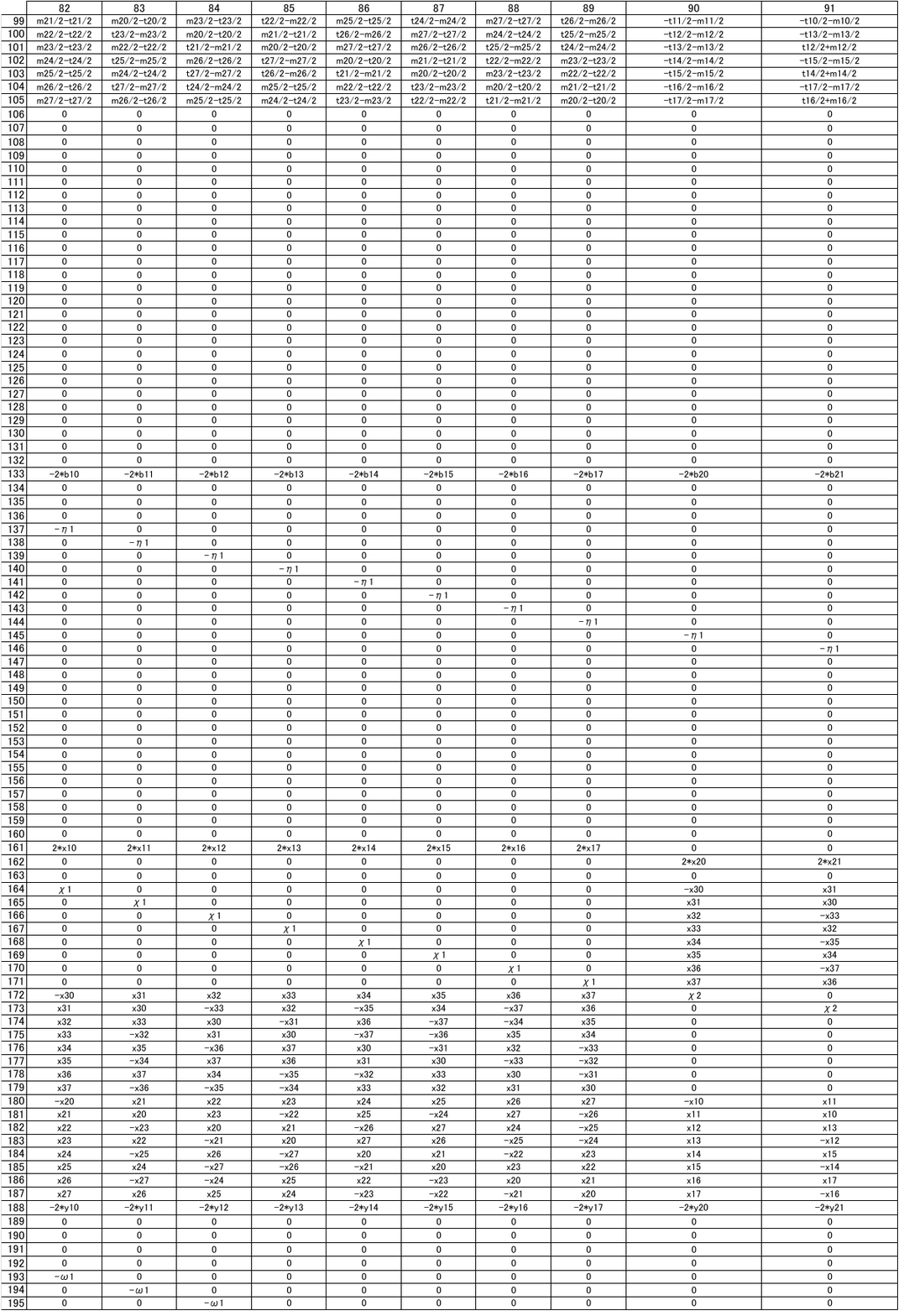}
\end{figure}

\begin{figure}[H]
\centering
\includegraphics[width=12cm, height=19.0cm,page= 40 ]{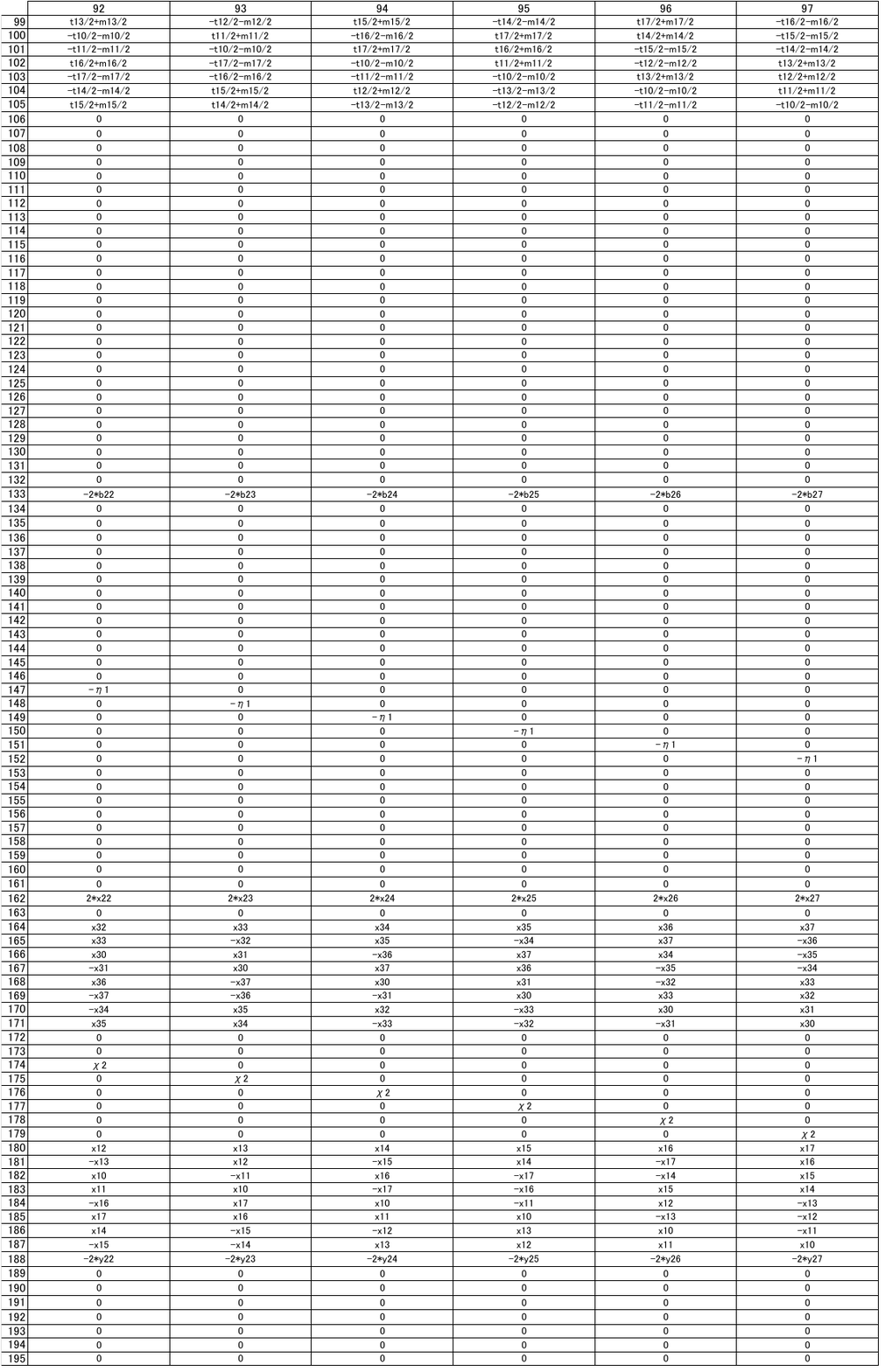}
\end{figure}

\begin{figure}[H]
\centering
\includegraphics[width=12cm, height=19.0cm,page= 41 ]{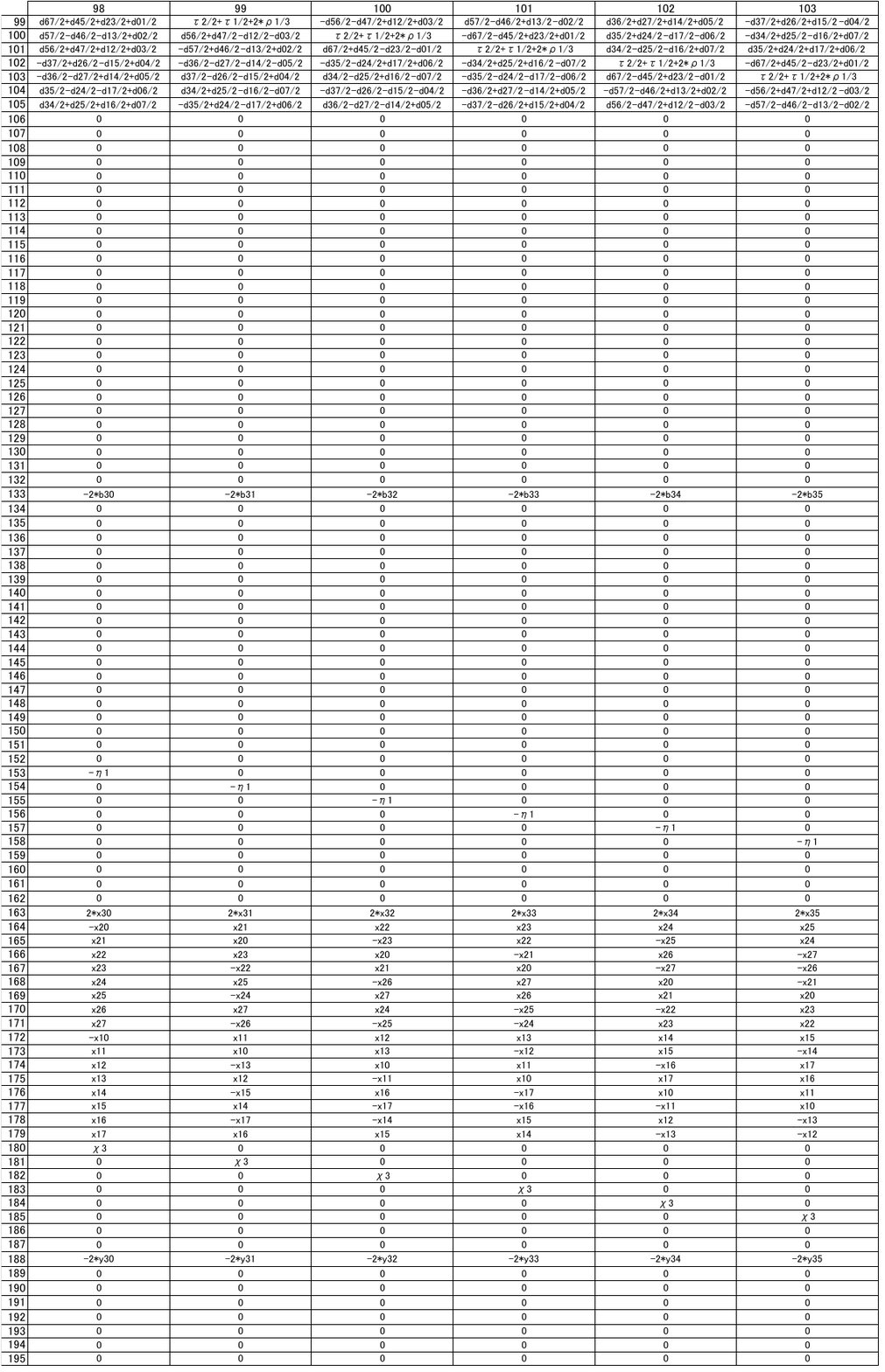}
\end{figure}

\begin{figure}[H]
\centering
\includegraphics[width=12cm, height=19.0cm,page= 42 ]{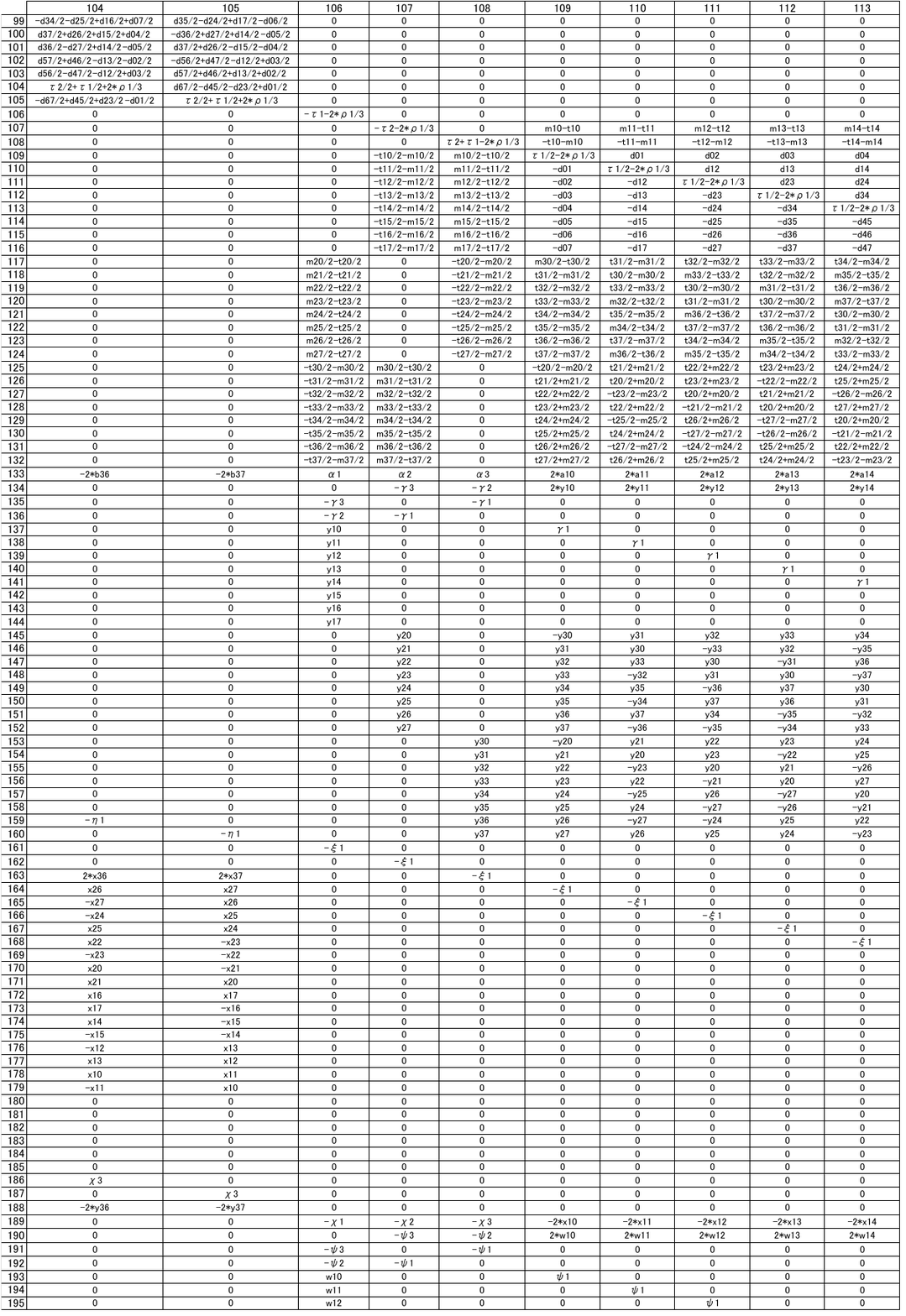}
\end{figure}

\begin{figure}[H]
\centering
\includegraphics[width=12cm, height=19.0cm,page= 43 ]{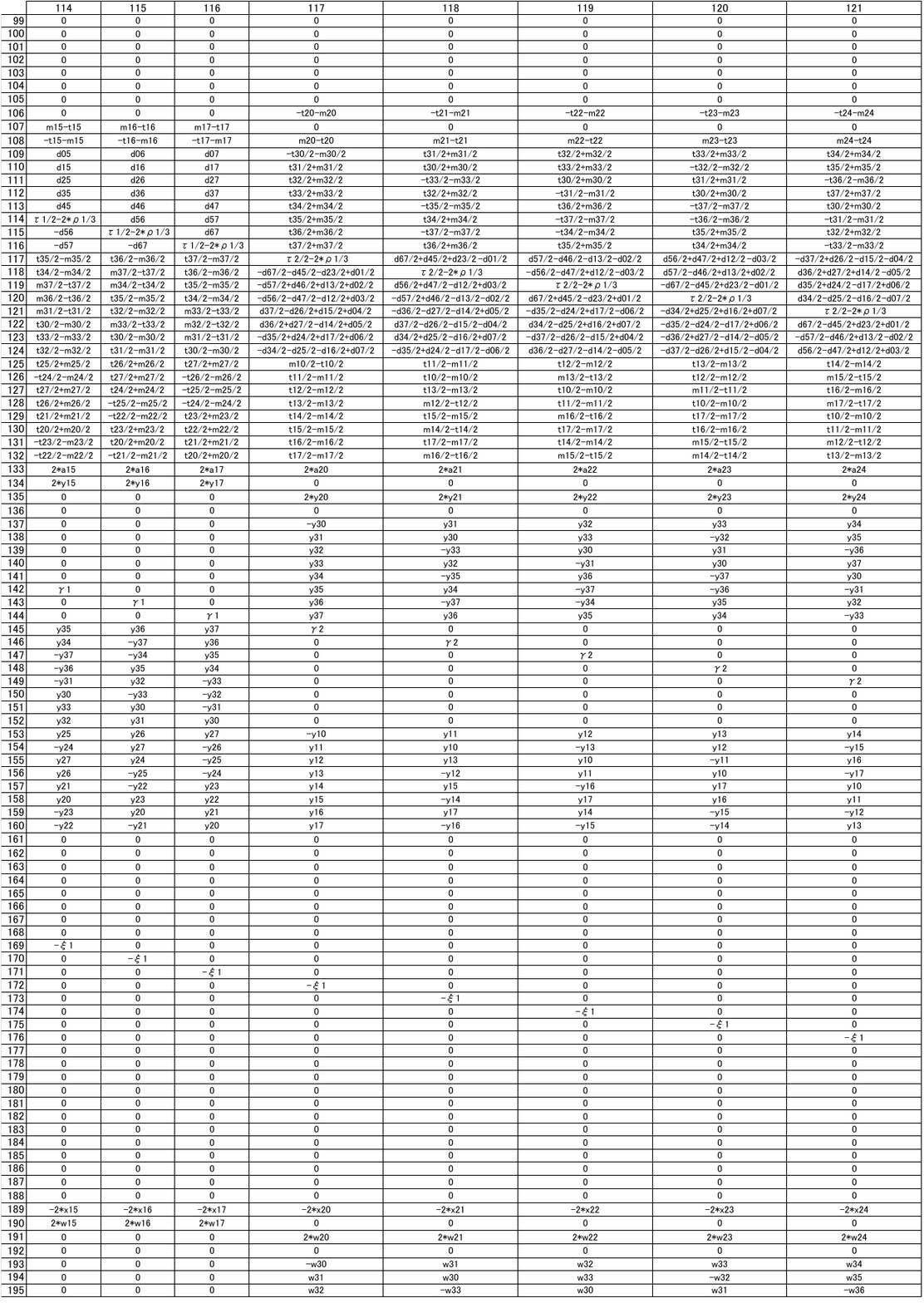}
\end{figure}

\begin{figure}[H]
\centering
\includegraphics[width=12cm, height=19.0cm,page= 44 ]{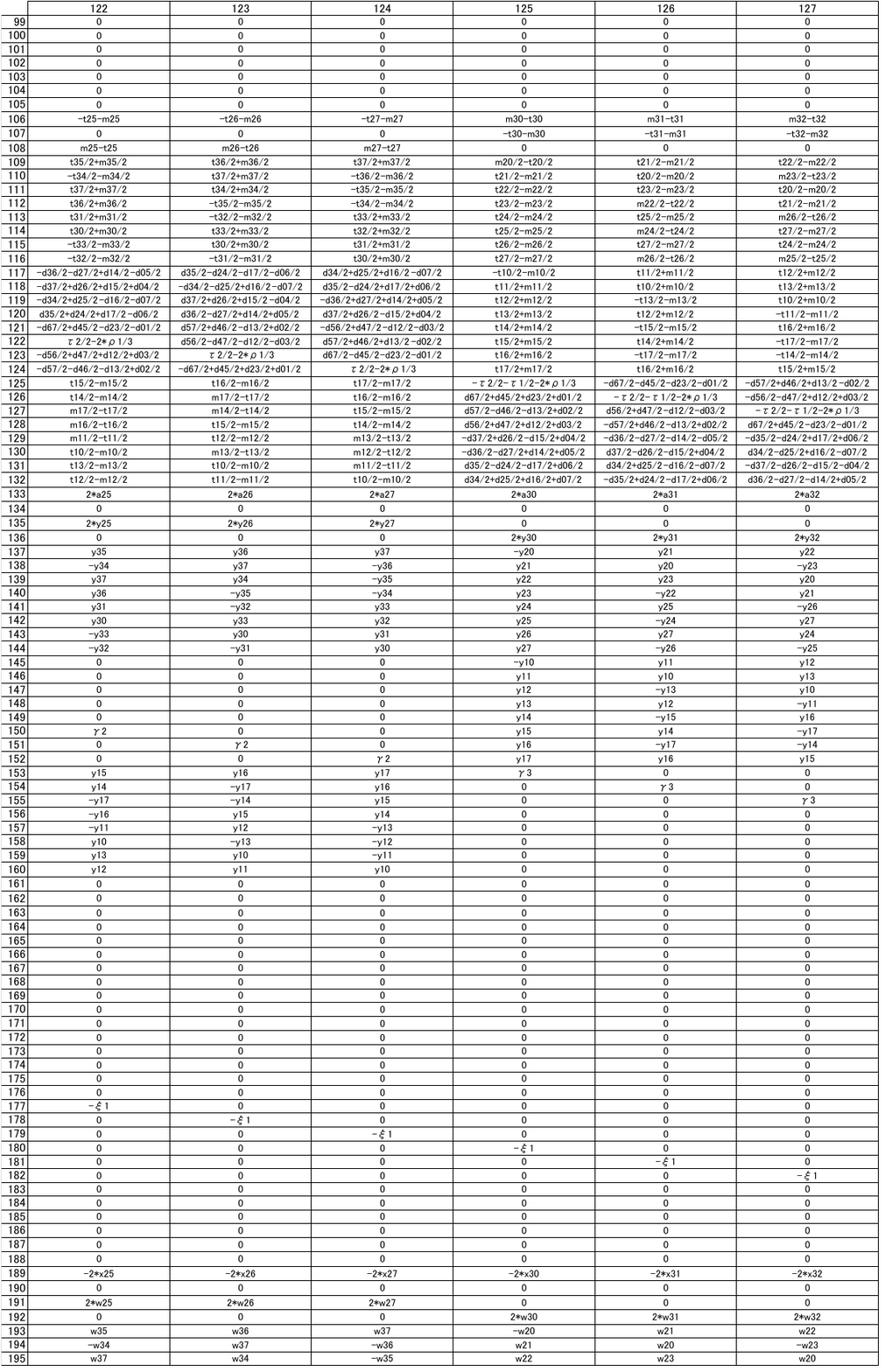}
\end{figure}

\begin{figure}[H]
\centering
\includegraphics[width=12cm, height=19.0cm,page= 45 ]{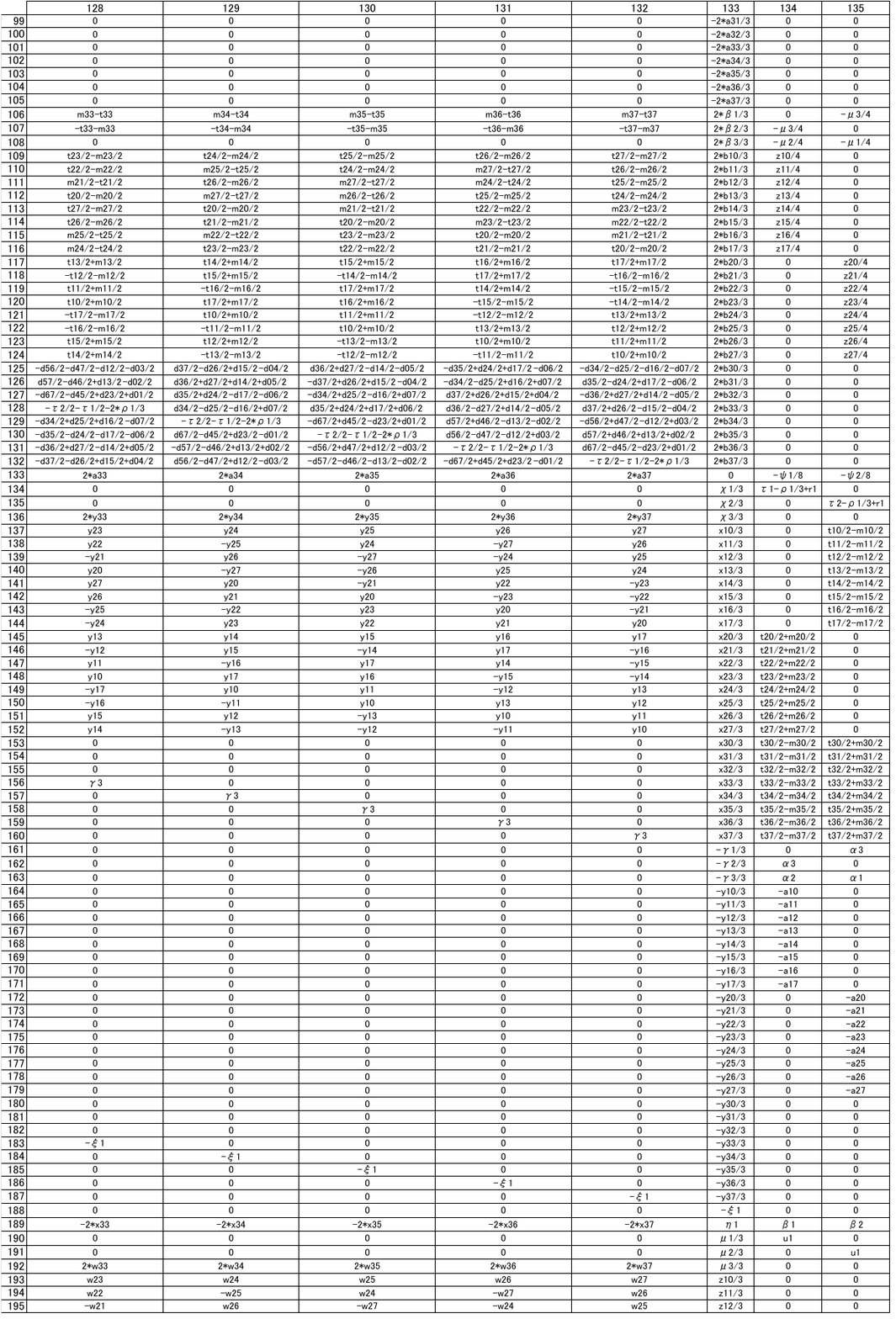}
\end{figure}

\begin{figure}[H]
\centering
\includegraphics[width=12cm, height=19.0cm,page= 46 ]{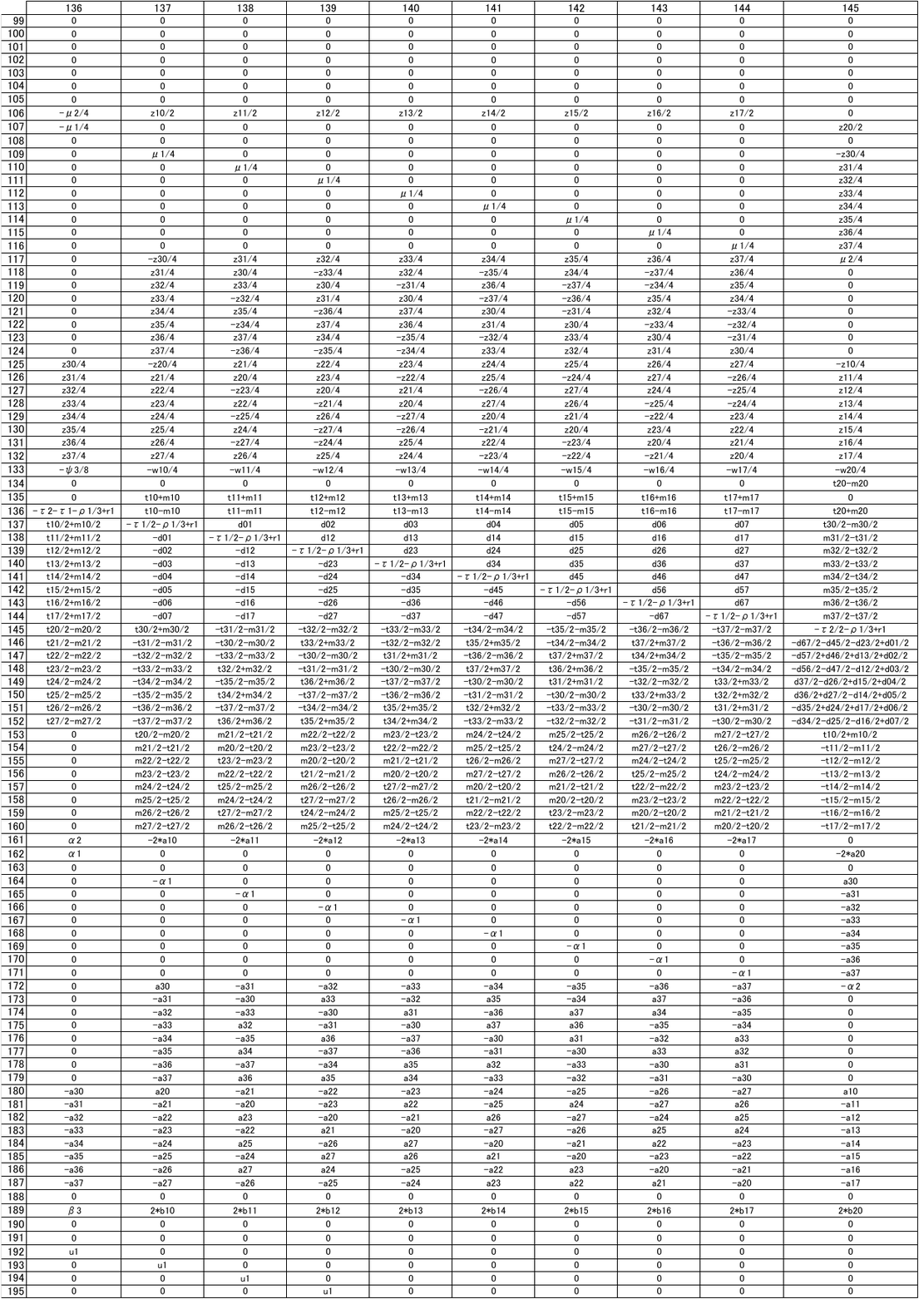}
\end{figure}

\begin{figure}[H]
\centering
\includegraphics[width=12cm, height=19.0cm,page= 47 ]{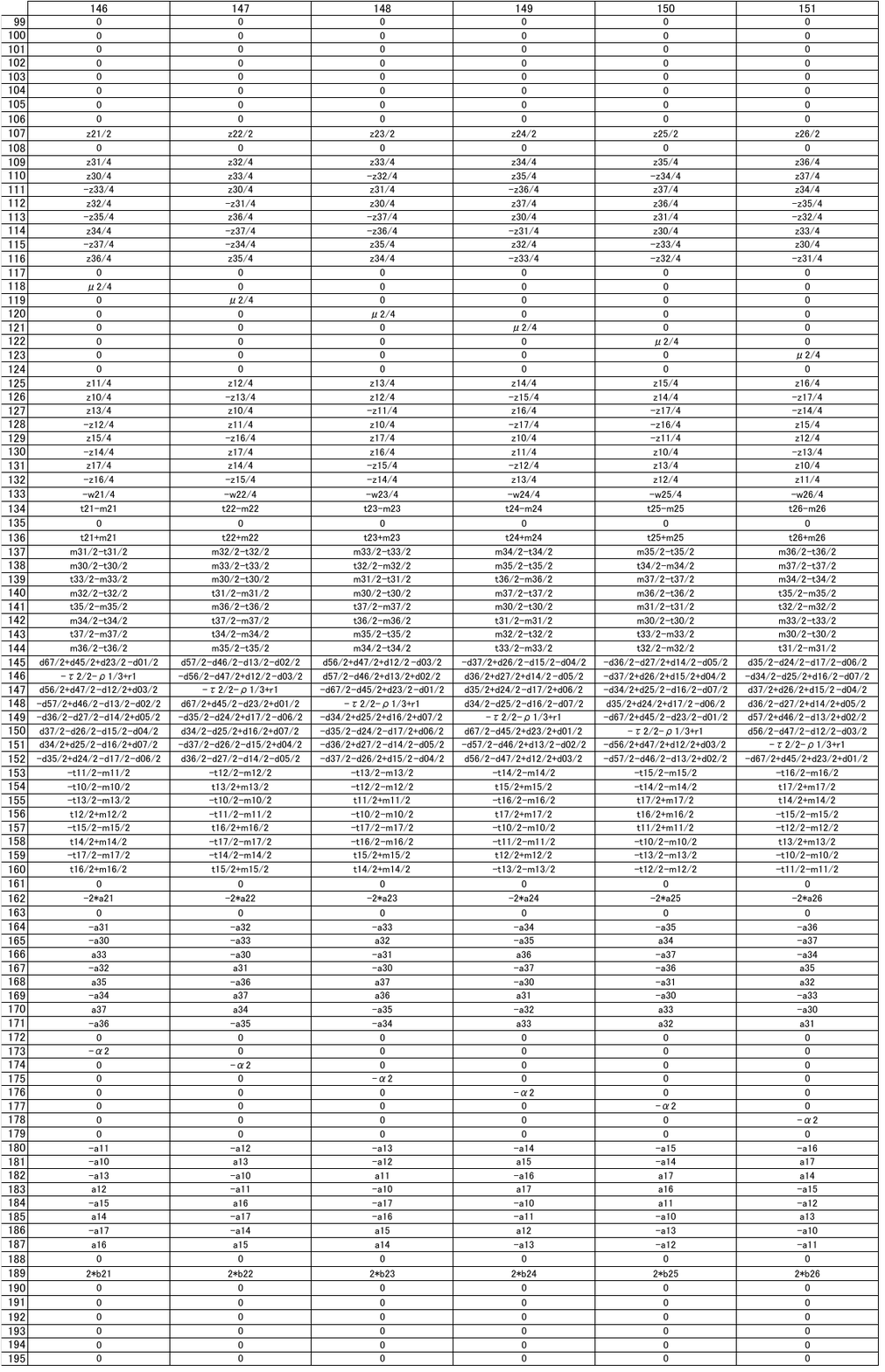}
\end{figure}

\begin{figure}[H]
\centering
\includegraphics[width=12cm, height=19.0cm,page= 48 ]{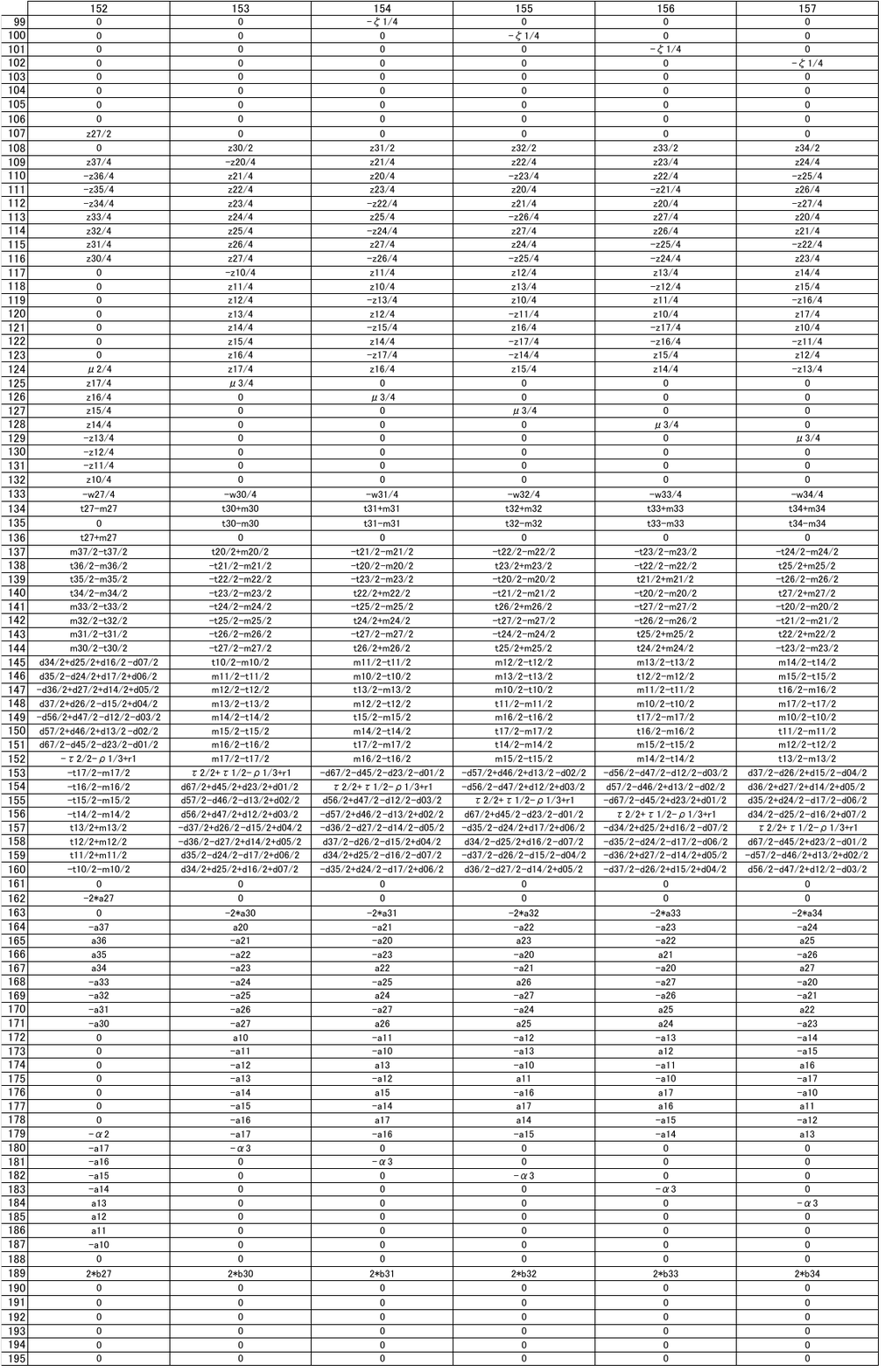}
\end{figure}

\begin{figure}[H]
\centering
\includegraphics[width=12cm, height=19.0cm,page= 49 ]{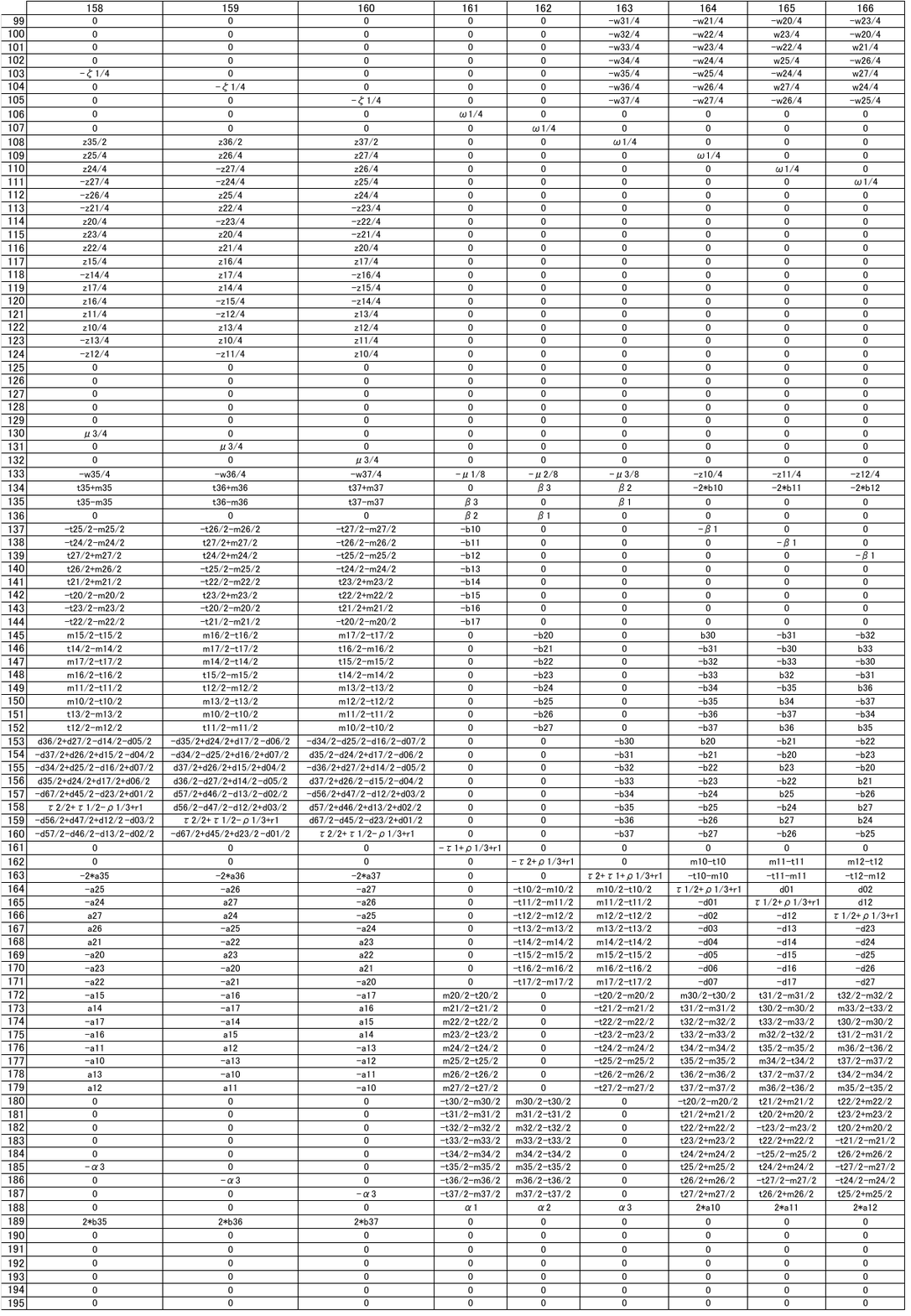}
\end{figure}

\begin{figure}[H]
\centering
\includegraphics[width=12cm, height=19.0cm,page= 50 ]{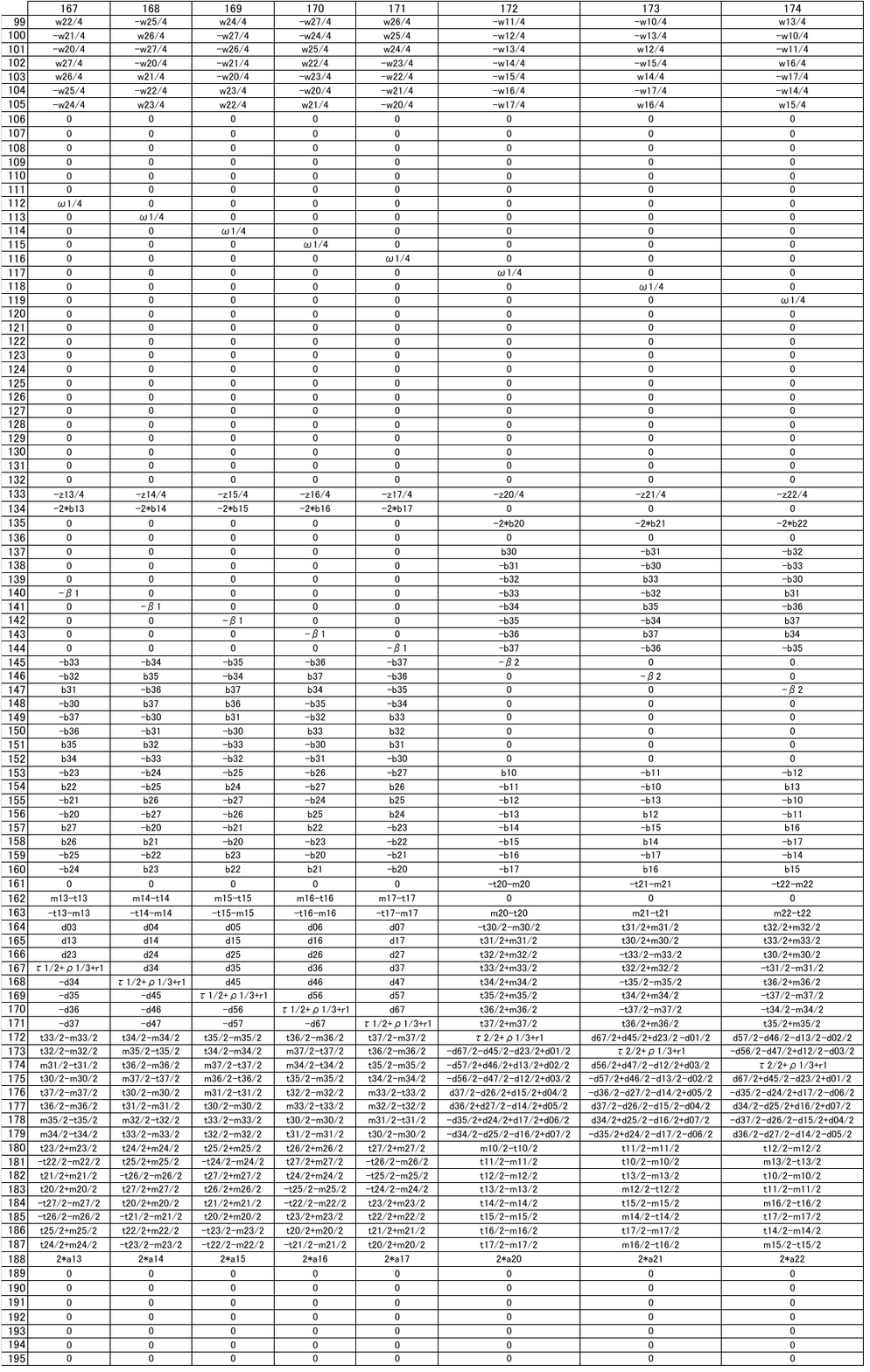}
\end{figure}

\begin{figure}[H]
\centering
\includegraphics[width=12cm, height=19.0cm,page= 51 ]{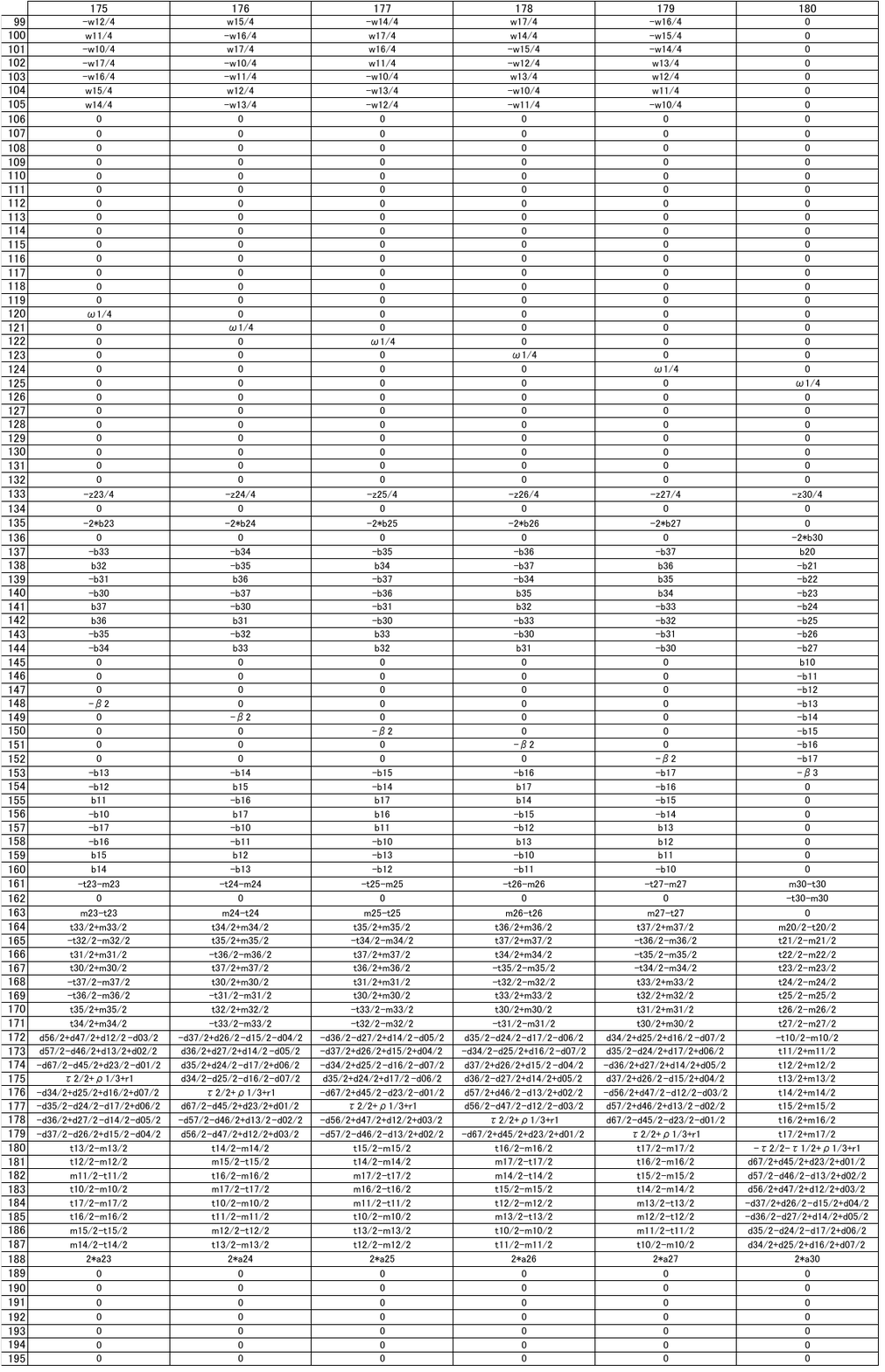}
\end{figure}

\begin{figure}[H]
\centering
\includegraphics[width=12cm, height=19.0cm,page= 52 ]{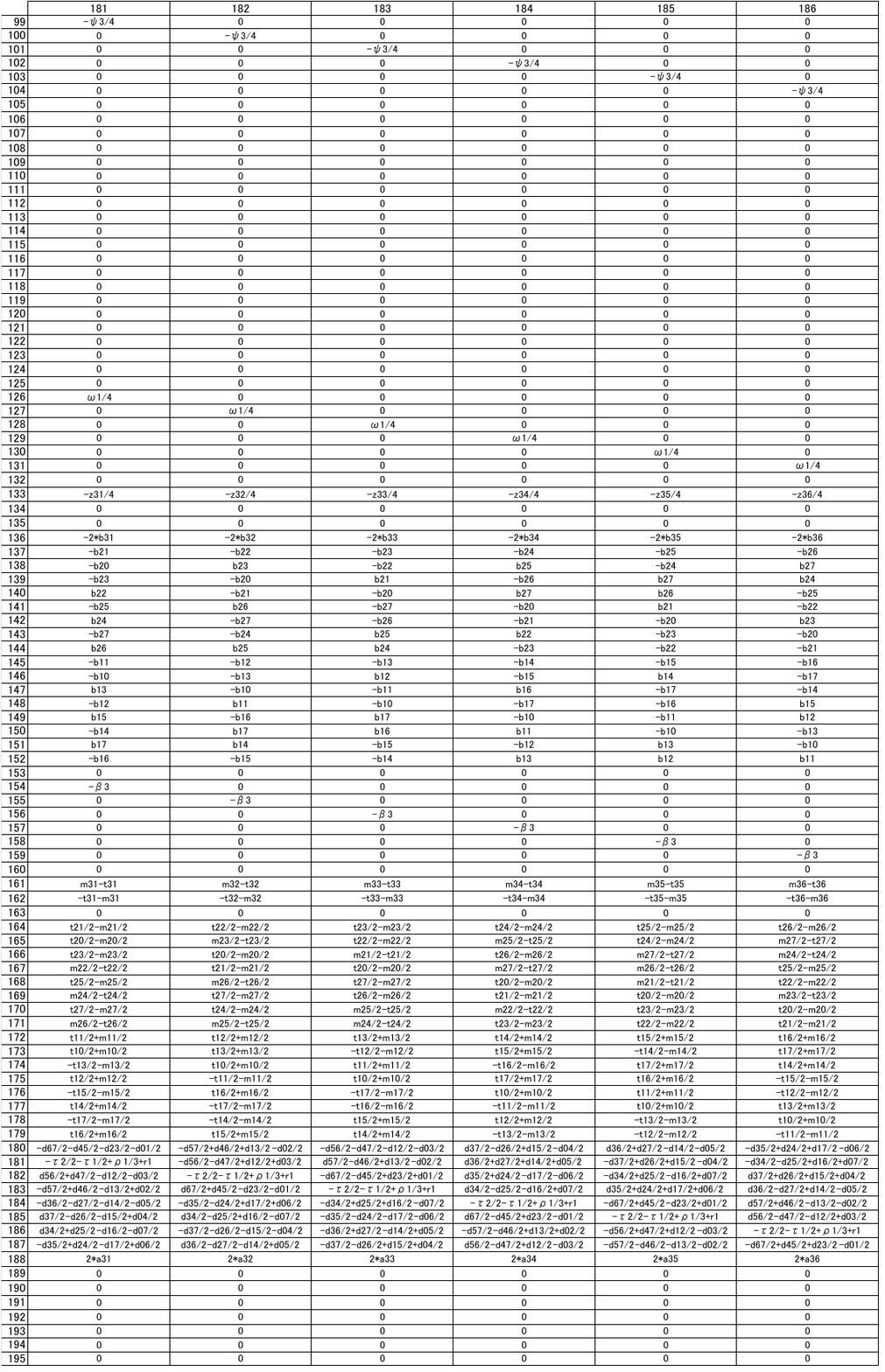}
\end{figure}

\begin{figure}[H]
\centering
\includegraphics[width=12cm, height=19.0cm,page= 53 ]{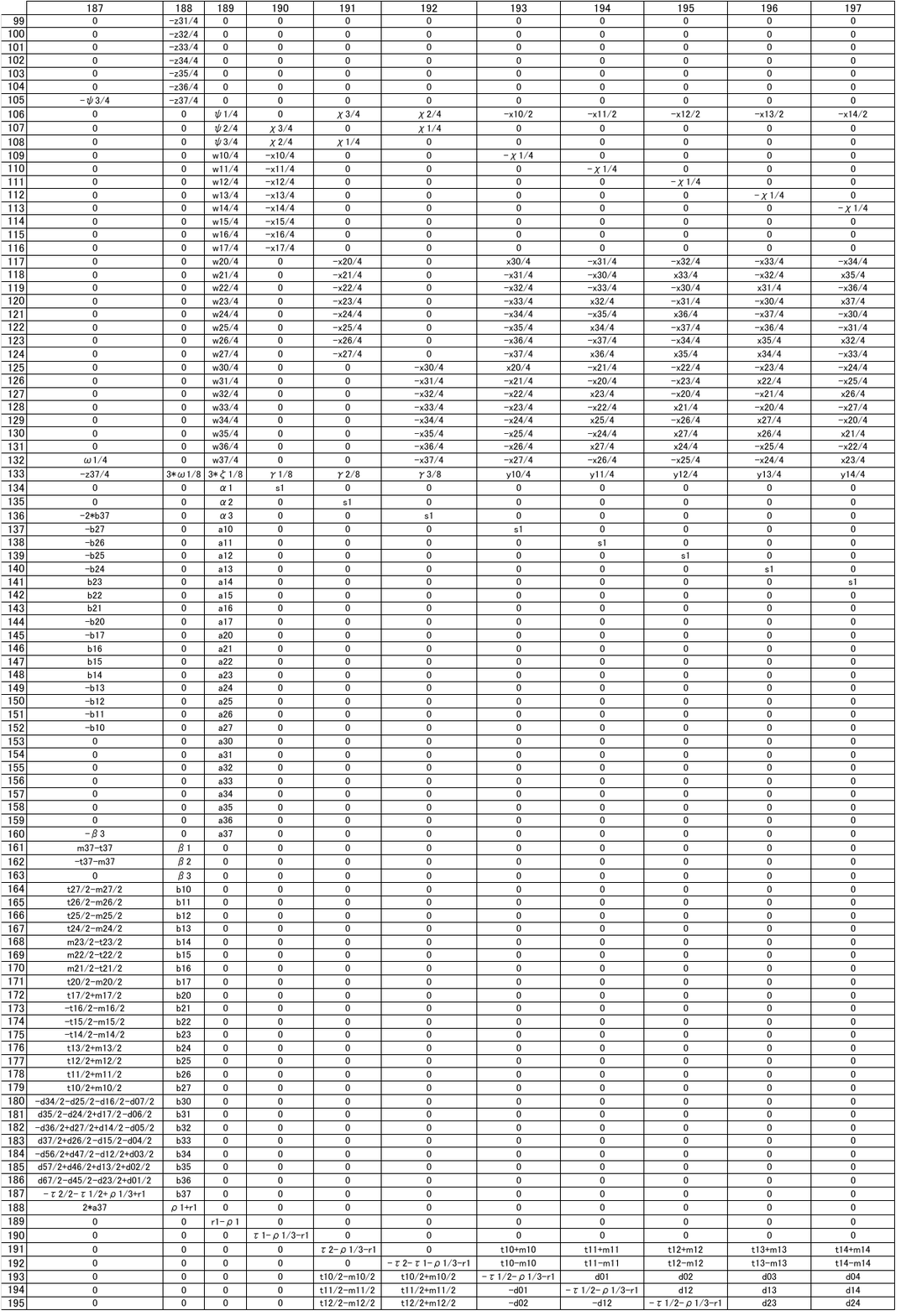}
\end{figure}

\begin{figure}[H]
\centering
\includegraphics[width=12cm, height=19.0cm,page= 54 ]{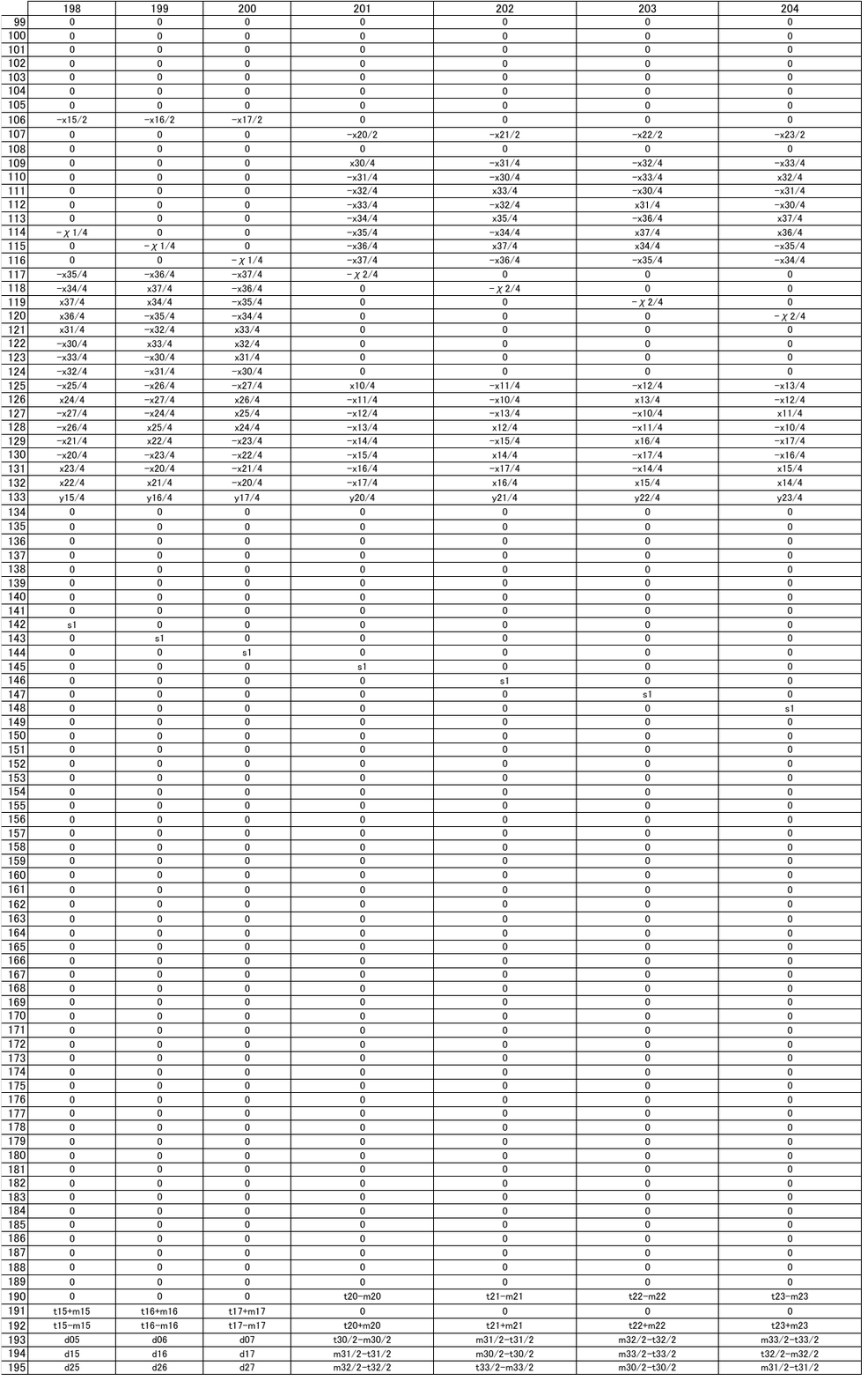}
\end{figure}

\begin{figure}[H]
\centering
\includegraphics[width=12cm, height=19.0cm,page= 55 ]{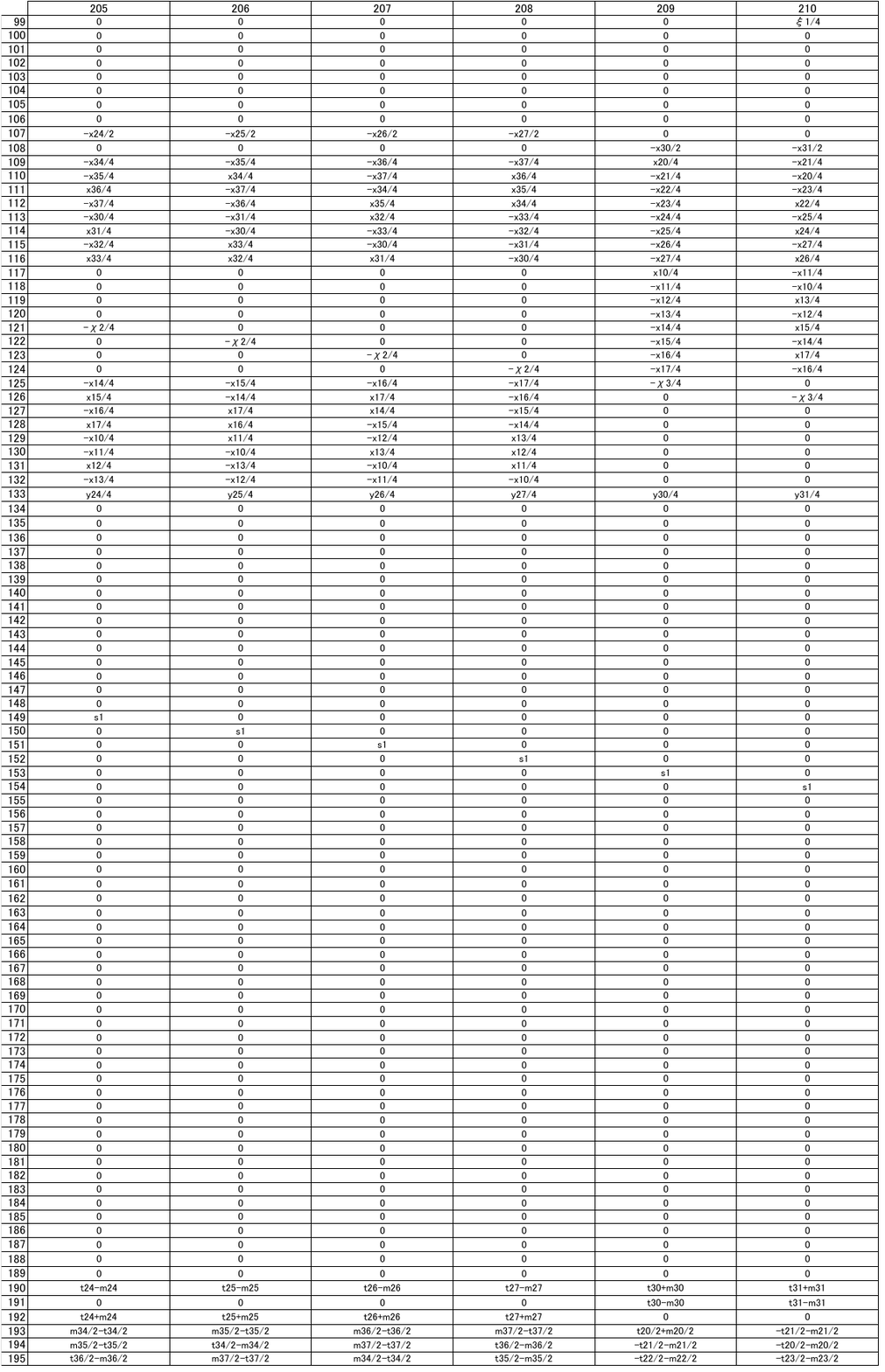}
\end{figure}

\begin{figure}[H]
\centering
\includegraphics[width=12cm, height=19.0cm,page= 56 ]{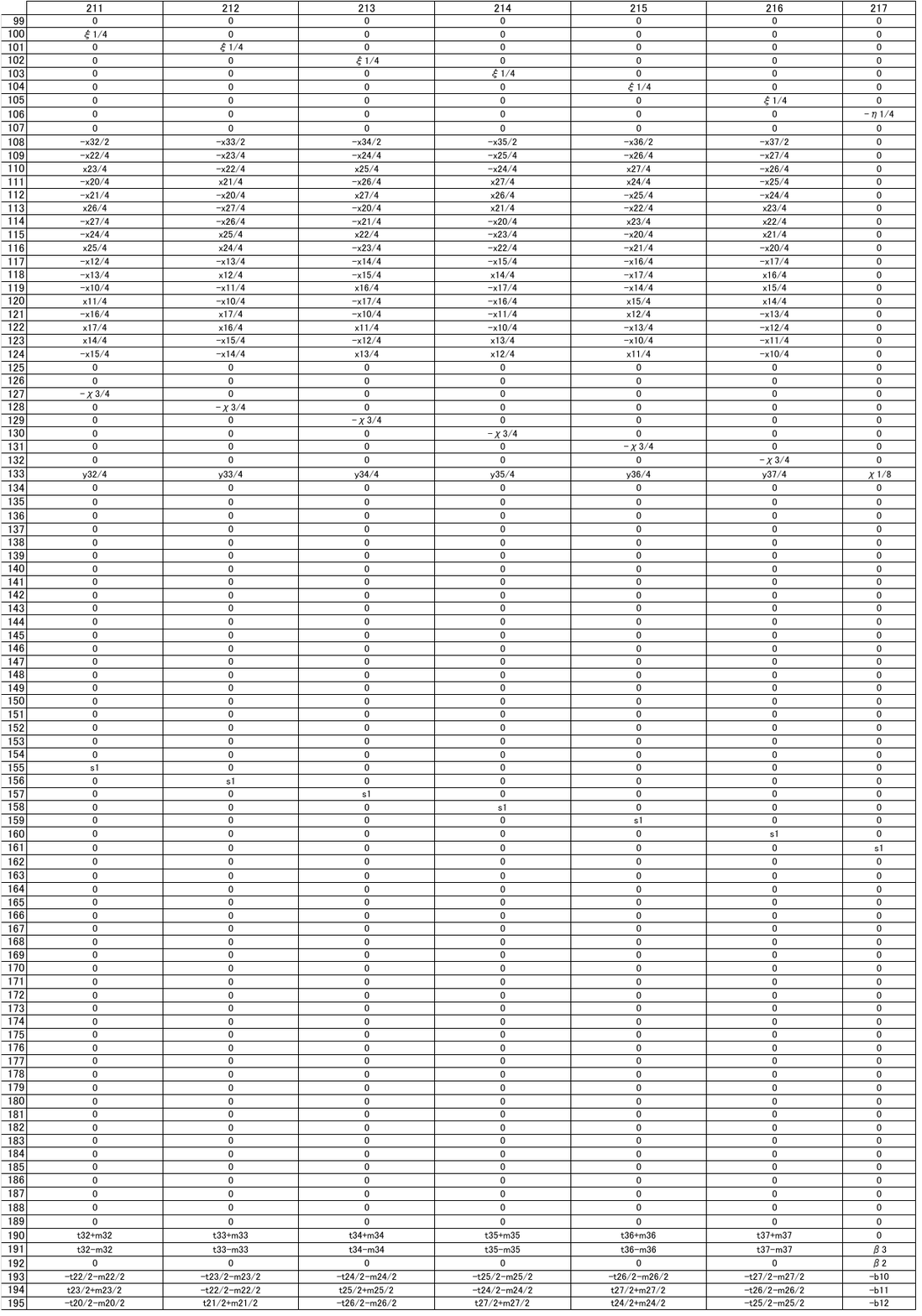}
\end{figure}

\begin{figure}[H]
\centering
\includegraphics[width=12cm, height=19.0cm,page= 57 ]{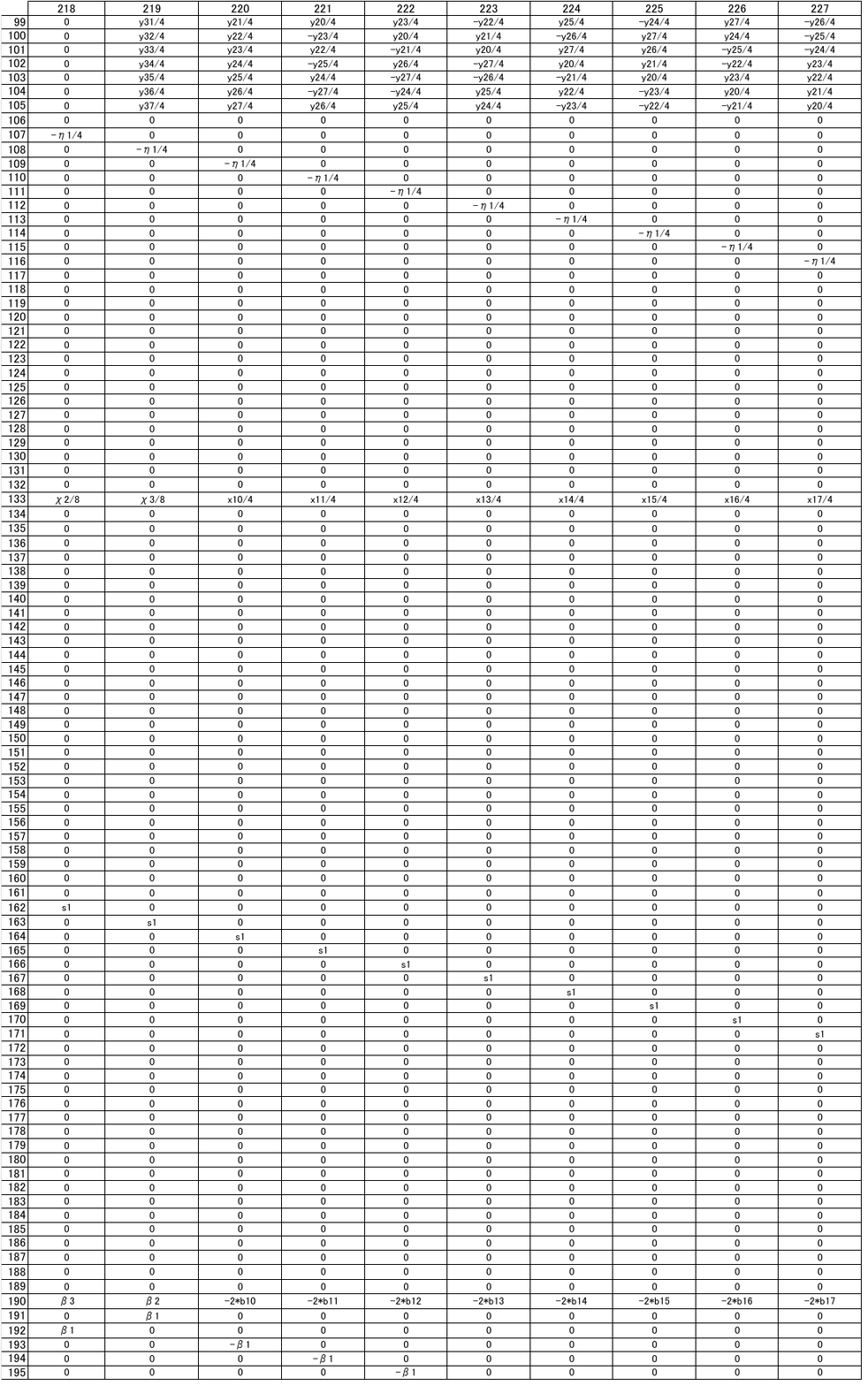}
\end{figure}

\begin{figure}[H]
\centering
\includegraphics[width=12cm, height=19.0cm,page= 58 ]{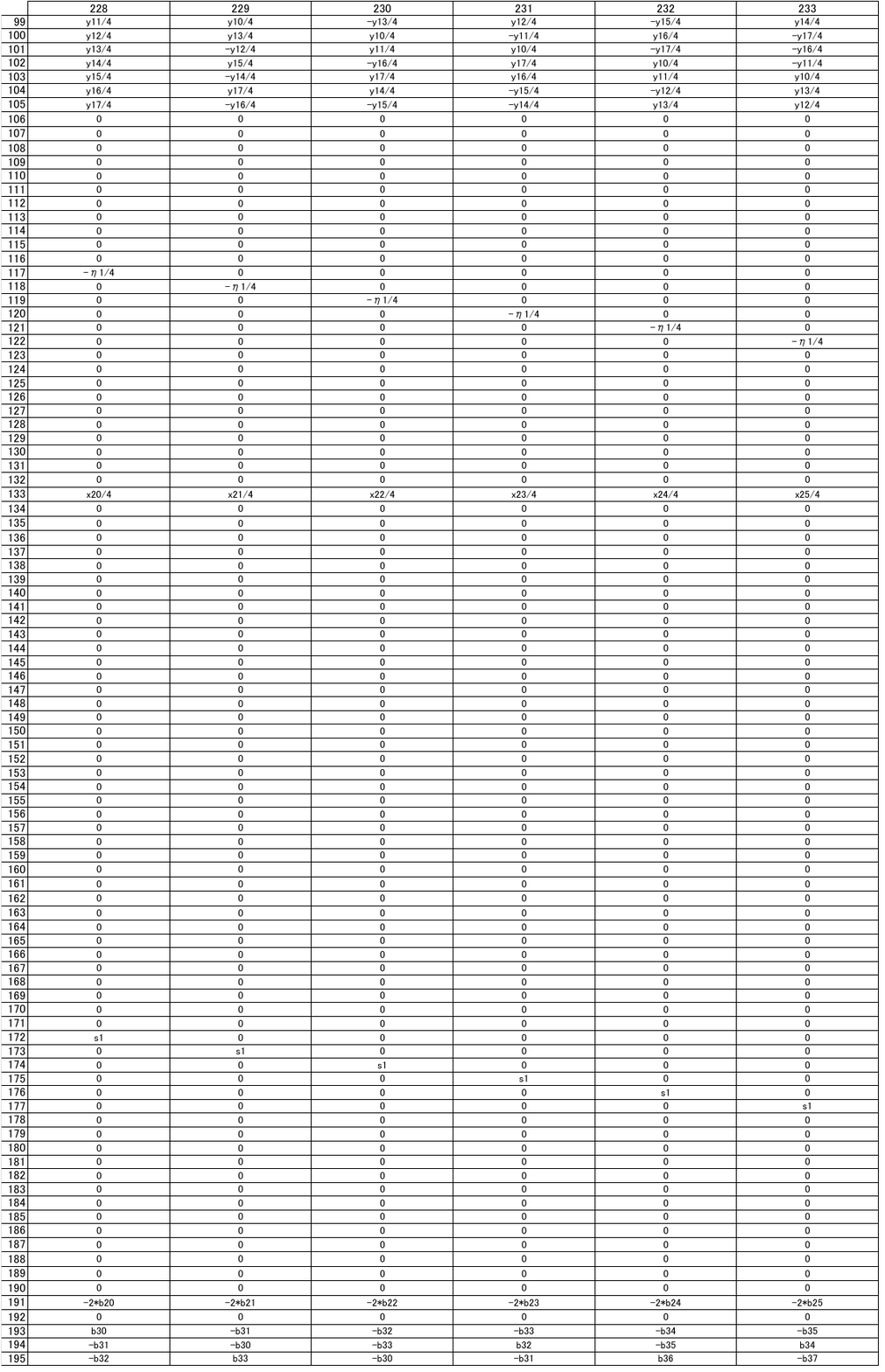}
\end{figure}

\begin{figure}[H]
\centering
\includegraphics[width=12cm, height=19.0cm,page= 59 ]{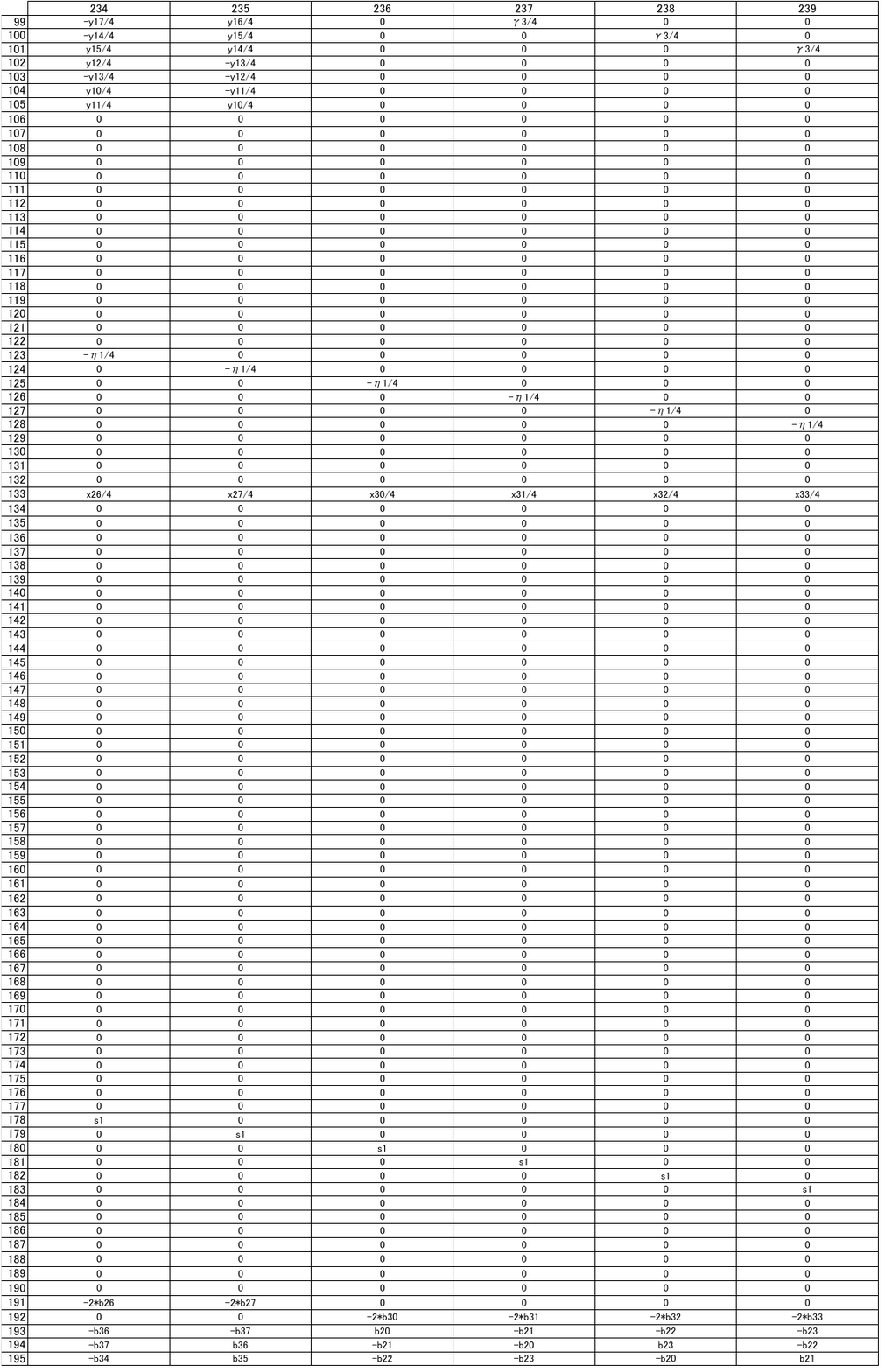}
\end{figure}

\begin{figure}[H]
\centering
\includegraphics[width=12cm, height=19.0cm,page= 60 ]{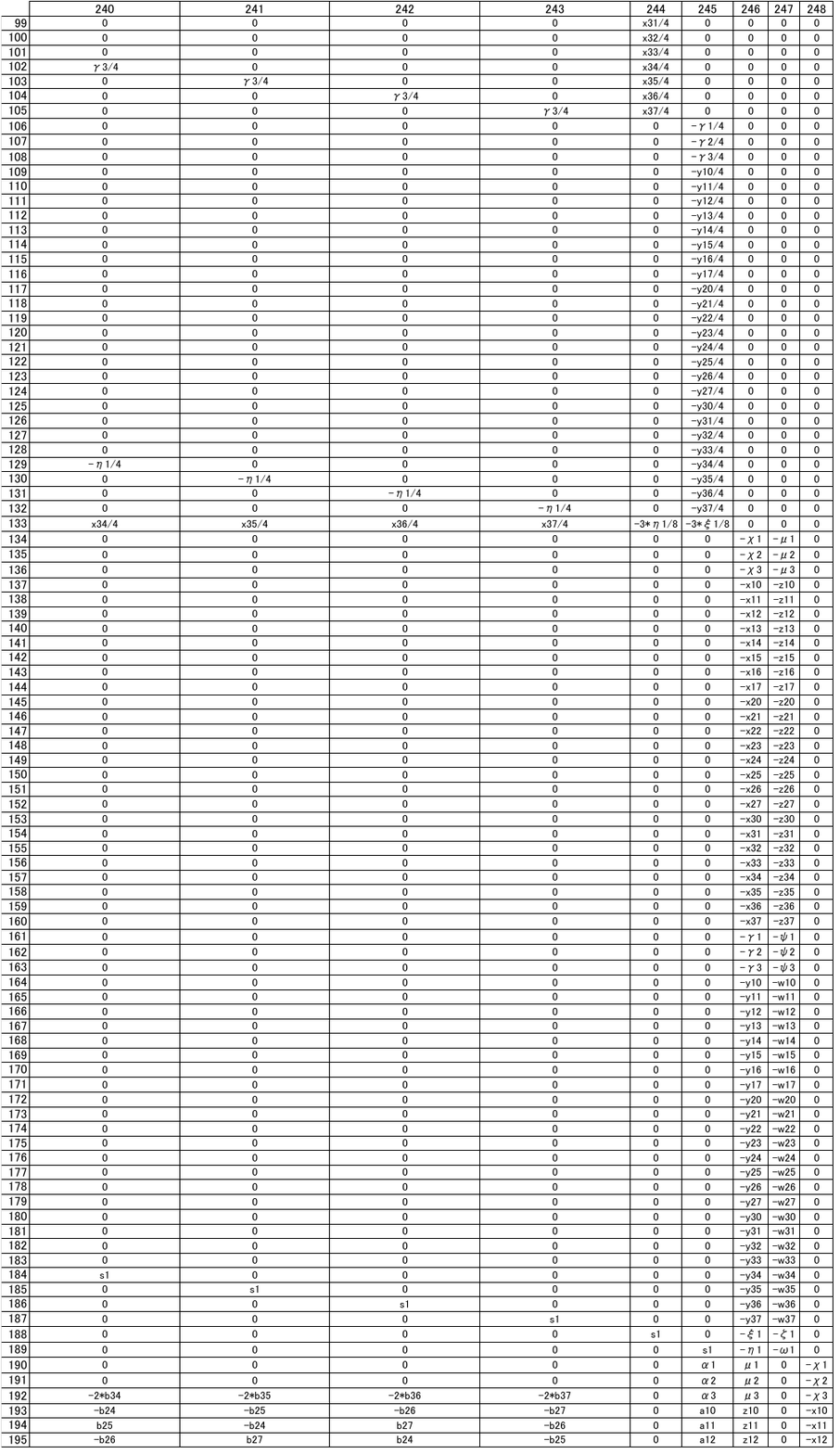}
\end{figure}

\begin{figure}[H]
\centering
\includegraphics[width=12cm, height=19.0cm,page= 61 ]{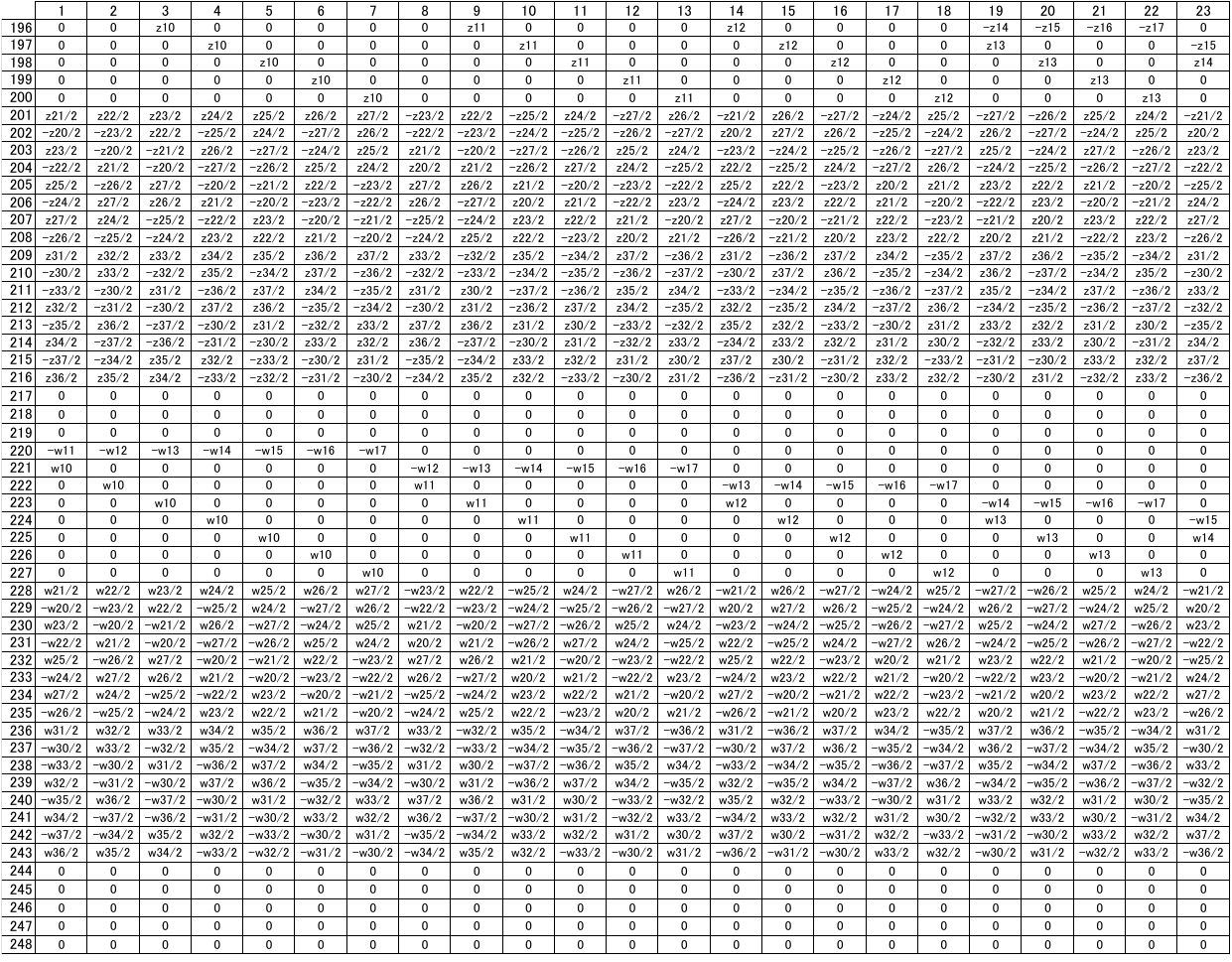}
\end{figure}

\begin{figure}[H]
\centering
\includegraphics[width=12cm, height=19.0cm,page= 62 ]{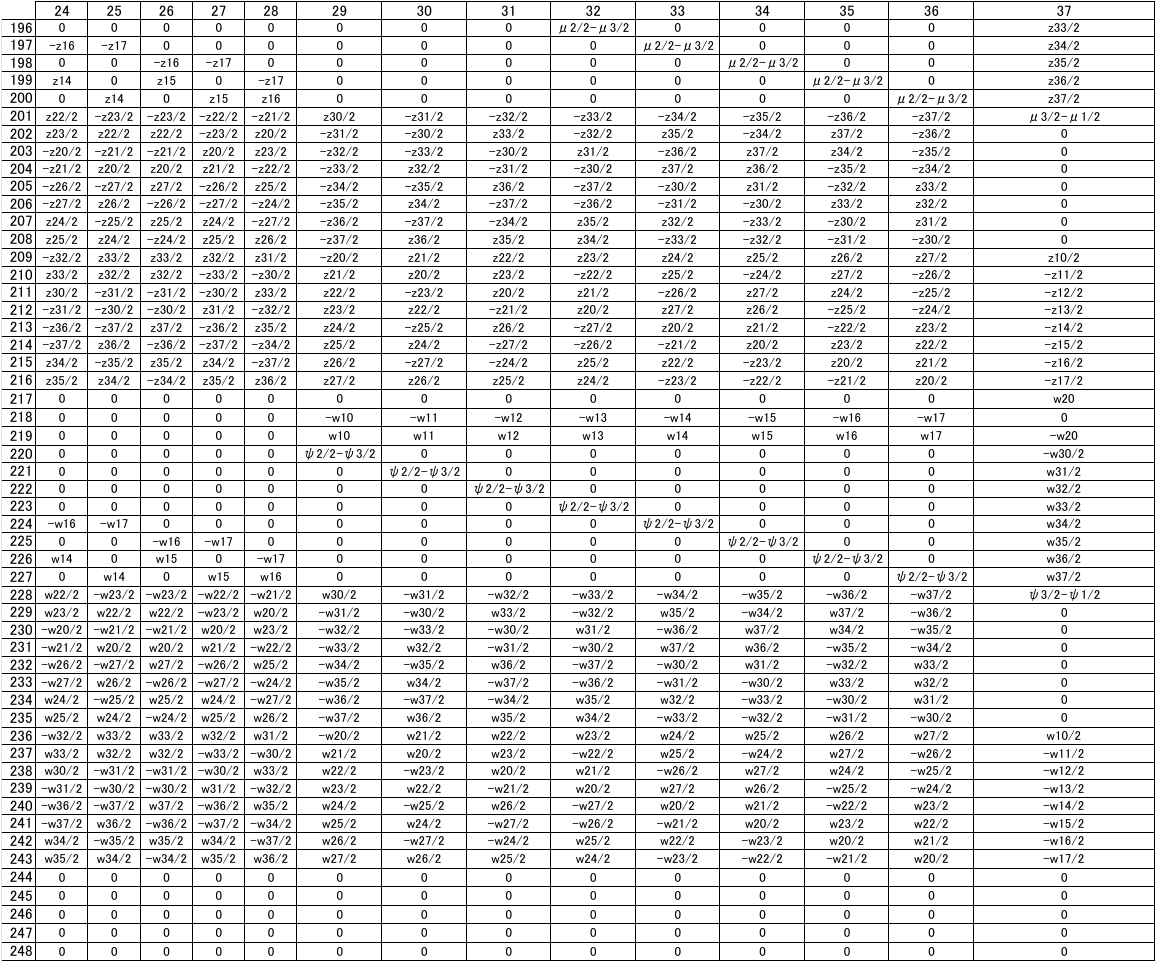}
\end{figure}

\begin{figure}[H]
\centering
\includegraphics[width=12cm, height=19.0cm,page= 63 ]{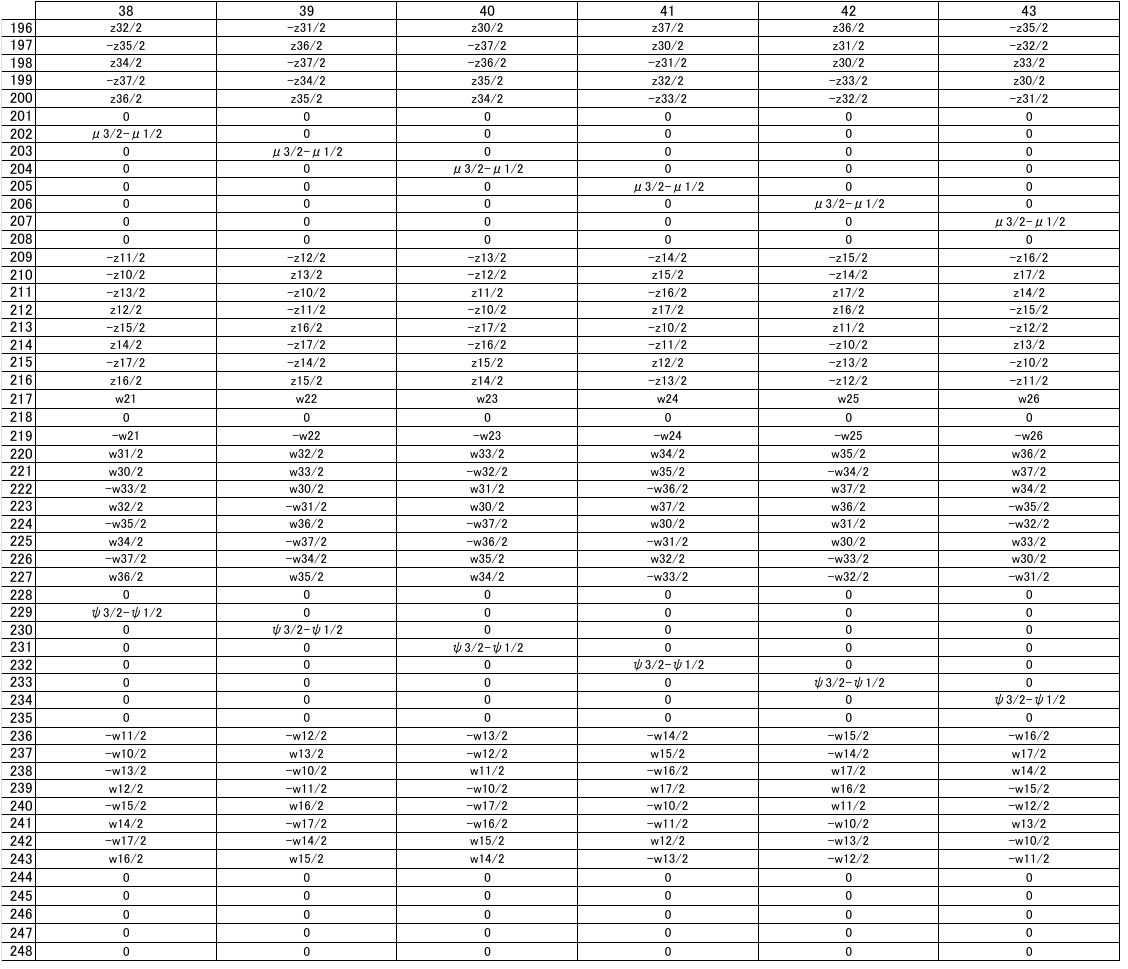}
\end{figure}

\begin{figure}[H]
\centering
\includegraphics[width=12cm, height=19.0cm,page= 64 ]{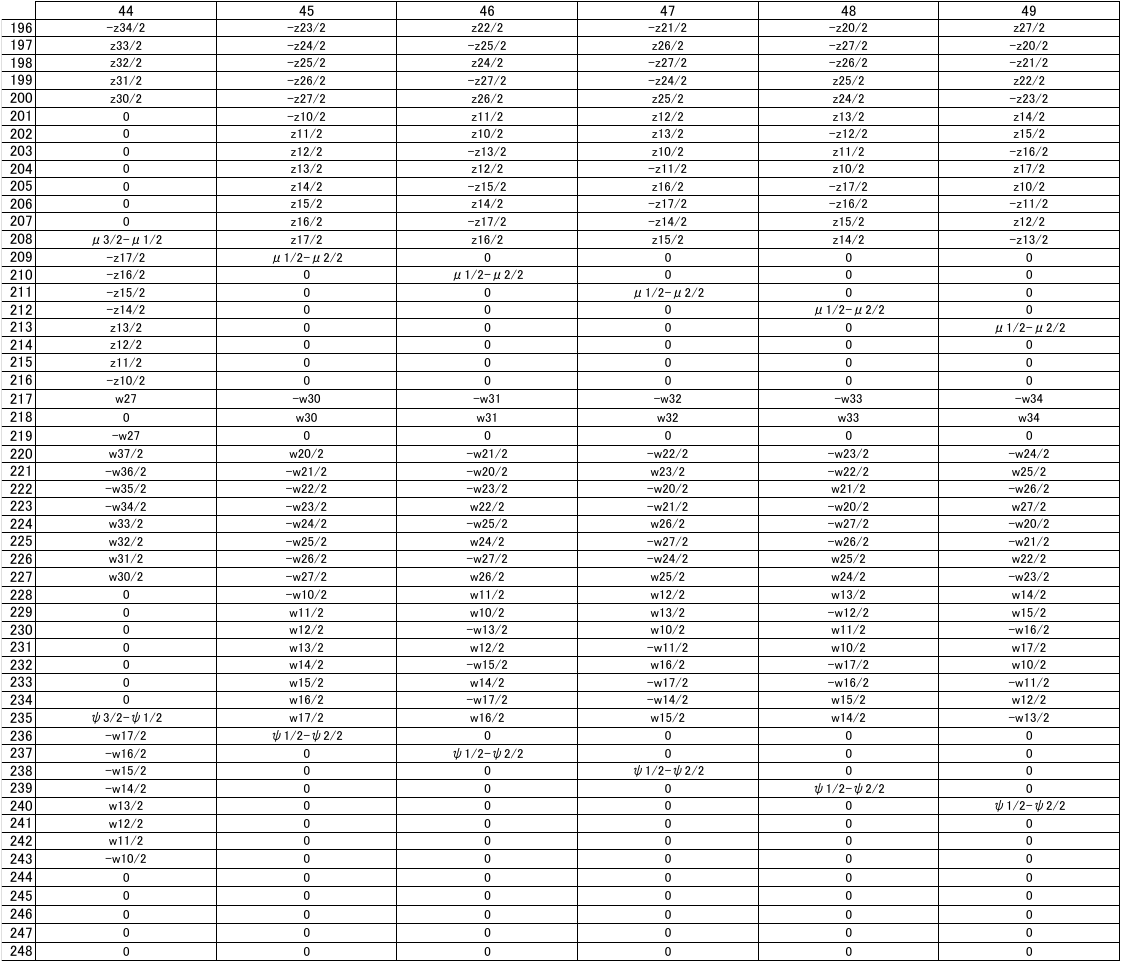}
\end{figure}

\begin{figure}[H]
\centering
\includegraphics[width=12cm, height=19.0cm,page= 65 ]{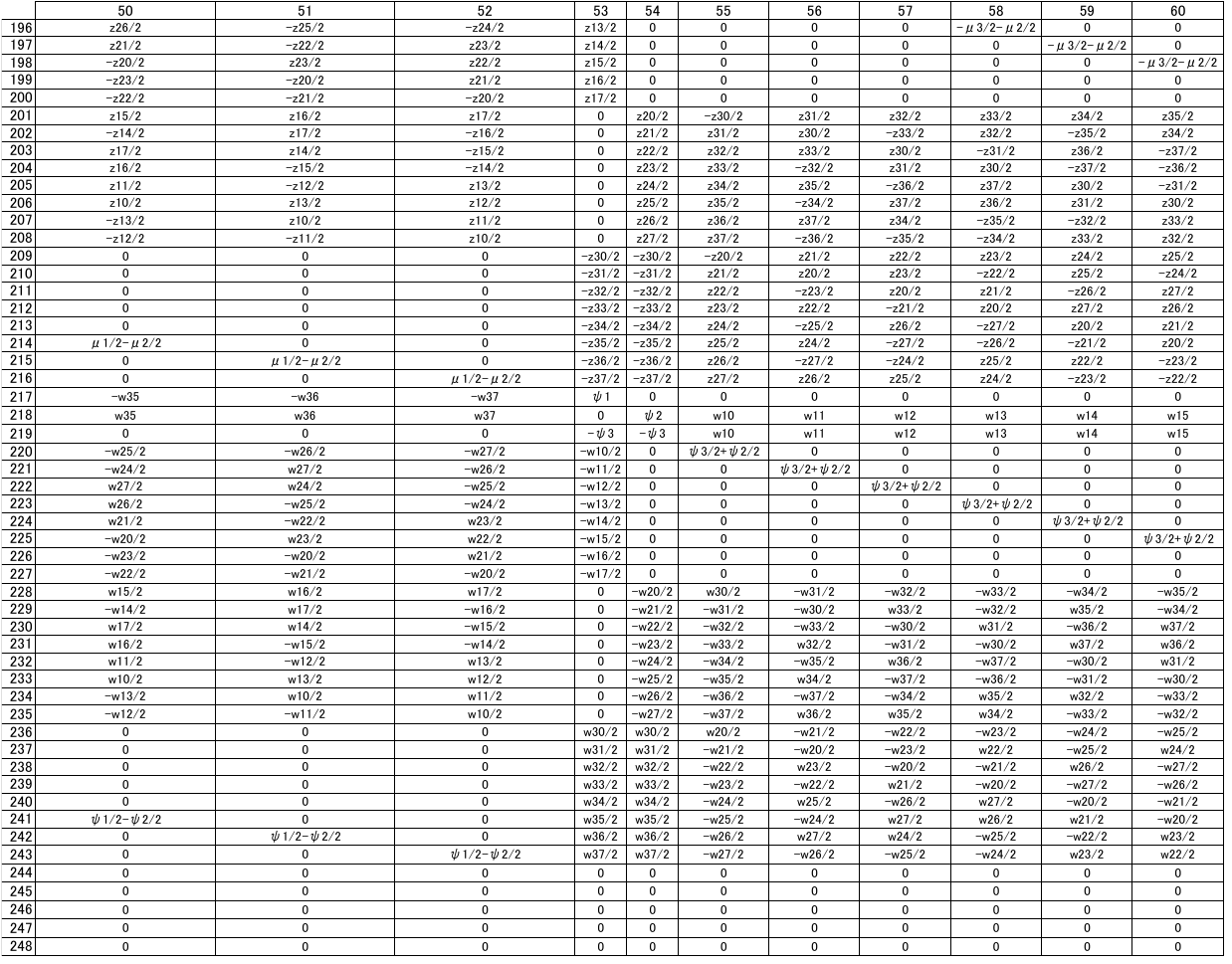}
\end{figure}

\begin{figure}[H]
\centering
\includegraphics[width=12cm, height=19.0cm,page= 66 ]{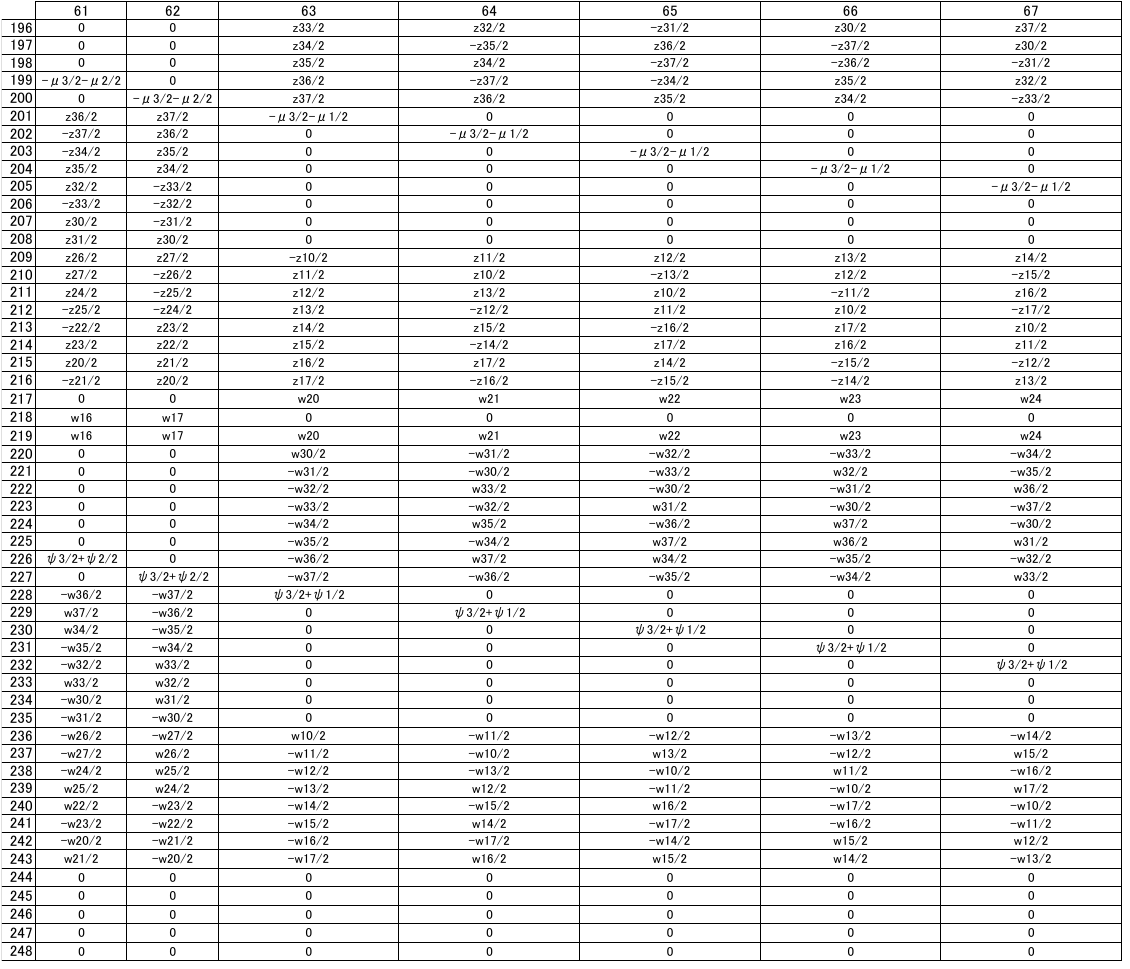}
\end{figure}

\begin{figure}[H]
\centering
\includegraphics[width=12cm, height=19.0cm,page= 67 ]{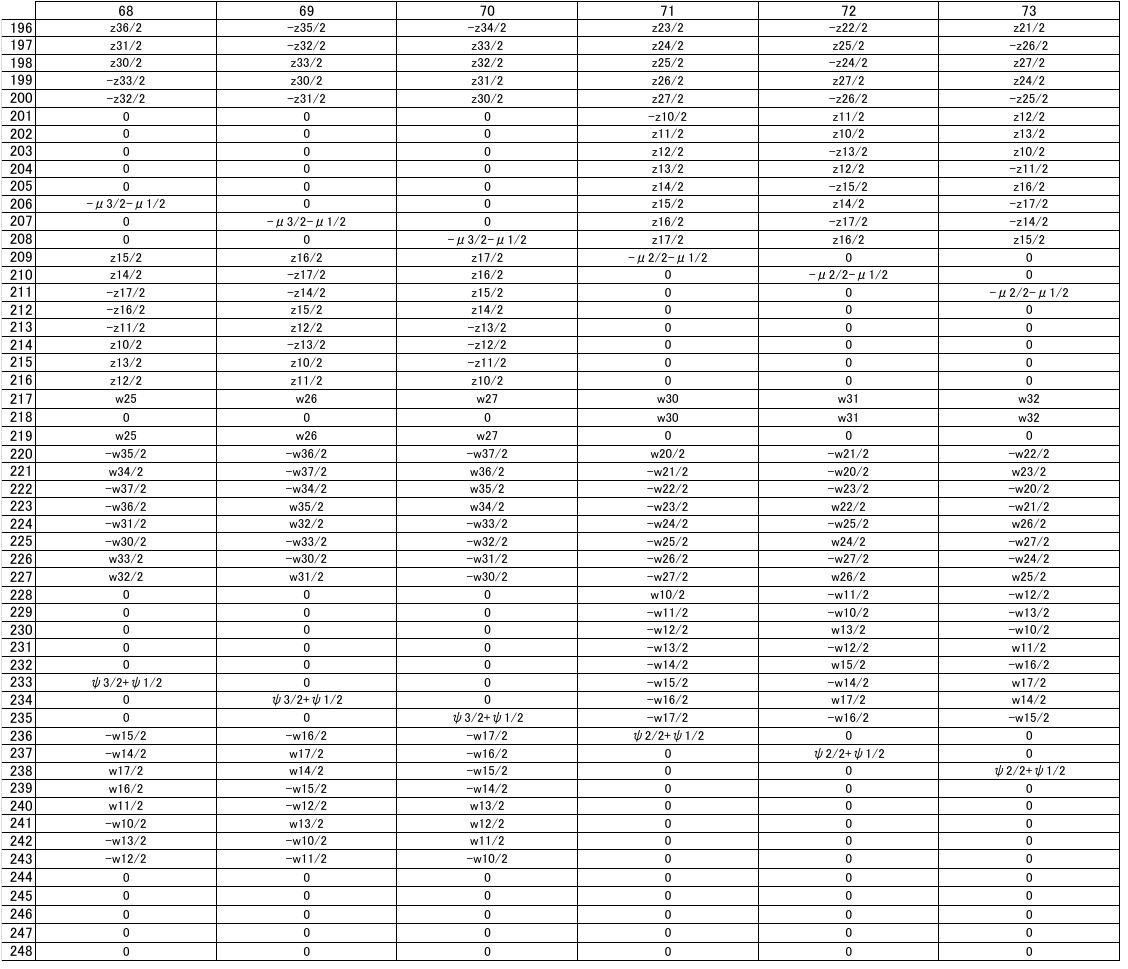}
\end{figure}

\begin{figure}[H]
\centering
\includegraphics[width=12cm, height=19.0cm,page= 68 ]{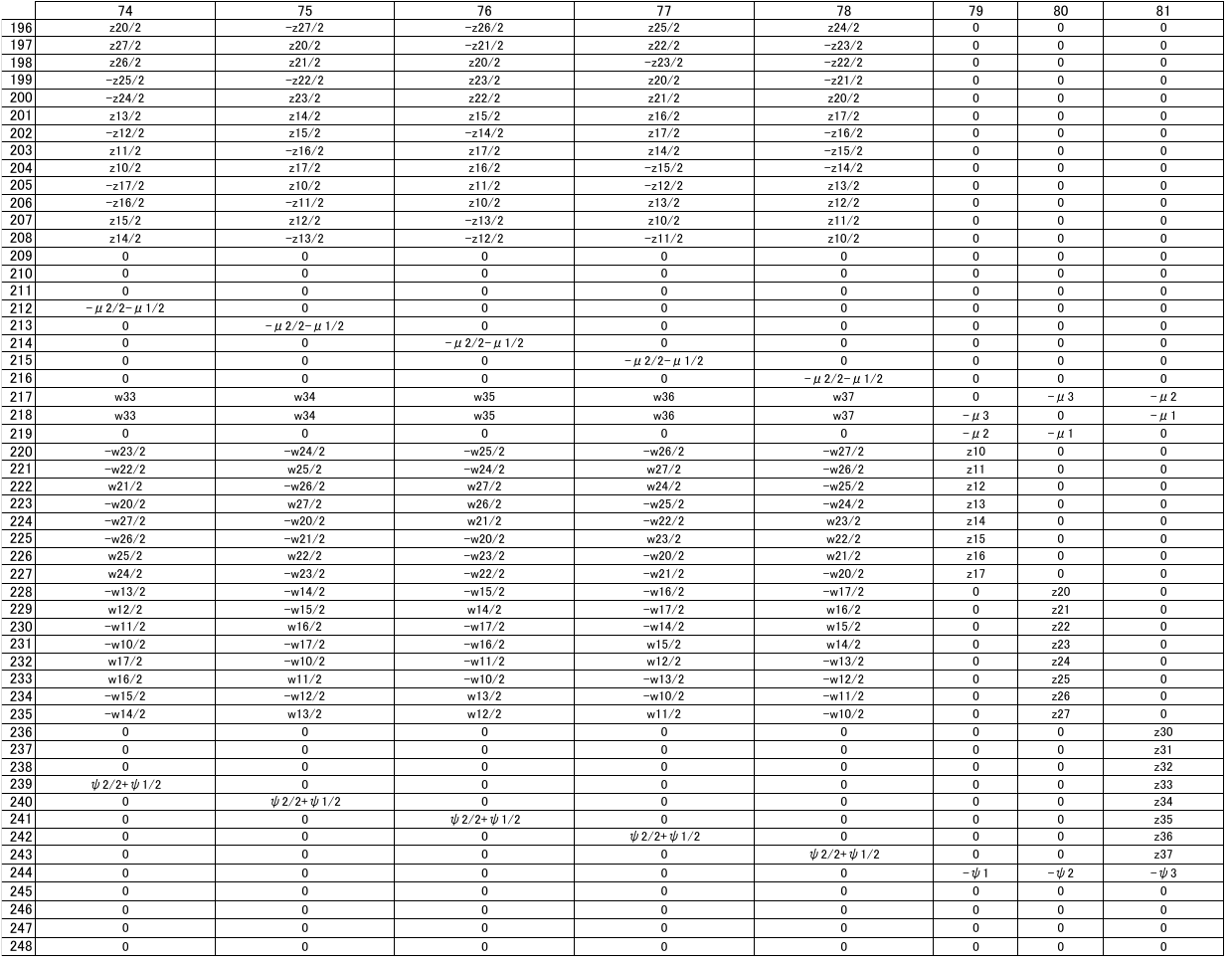}
\end{figure}

\begin{figure}[H]
\centering
\includegraphics[width=12cm, height=19.0cm,page= 69 ]{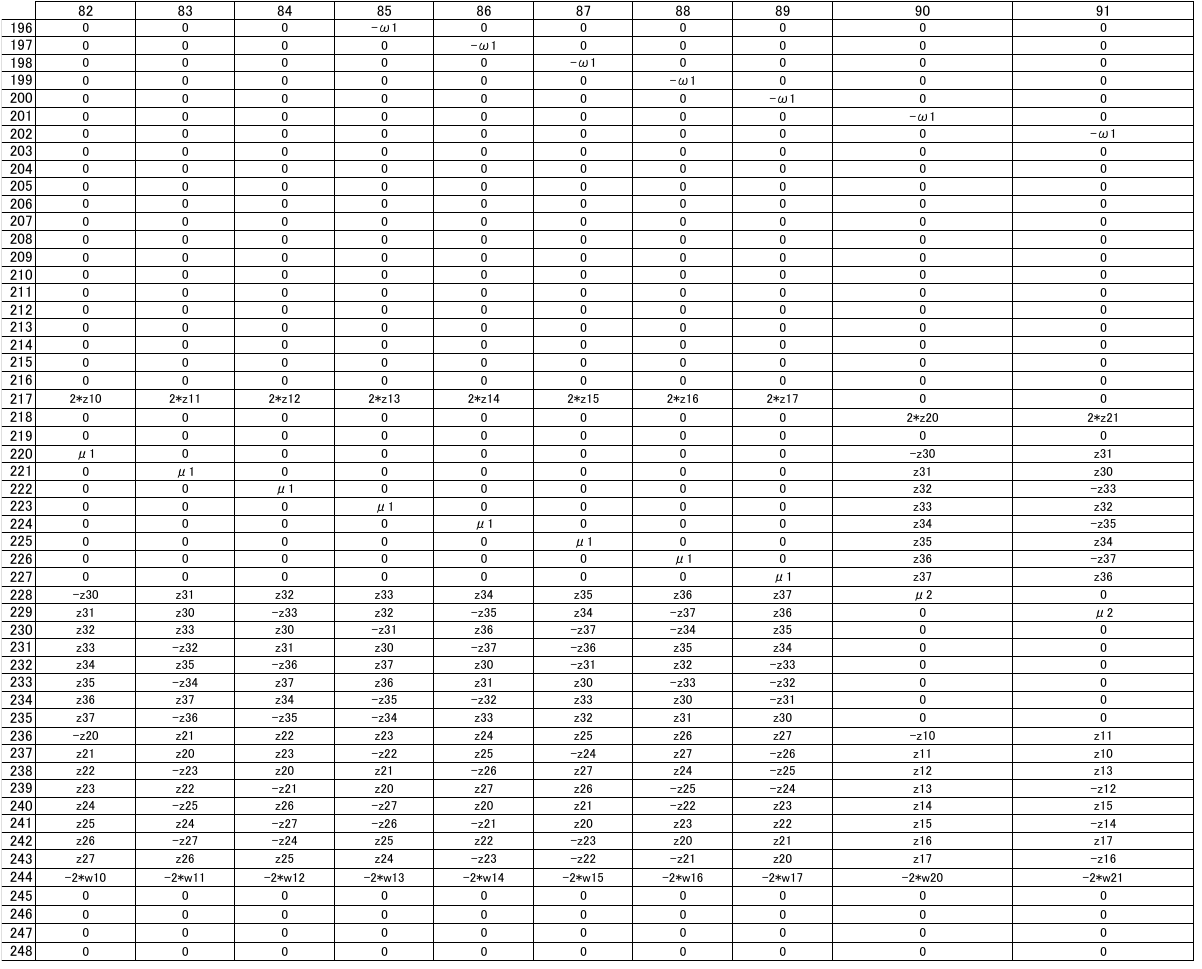}
\end{figure}

\begin{figure}[H]
\centering
\includegraphics[width=12cm, height=19.0cm,page= 70 ]{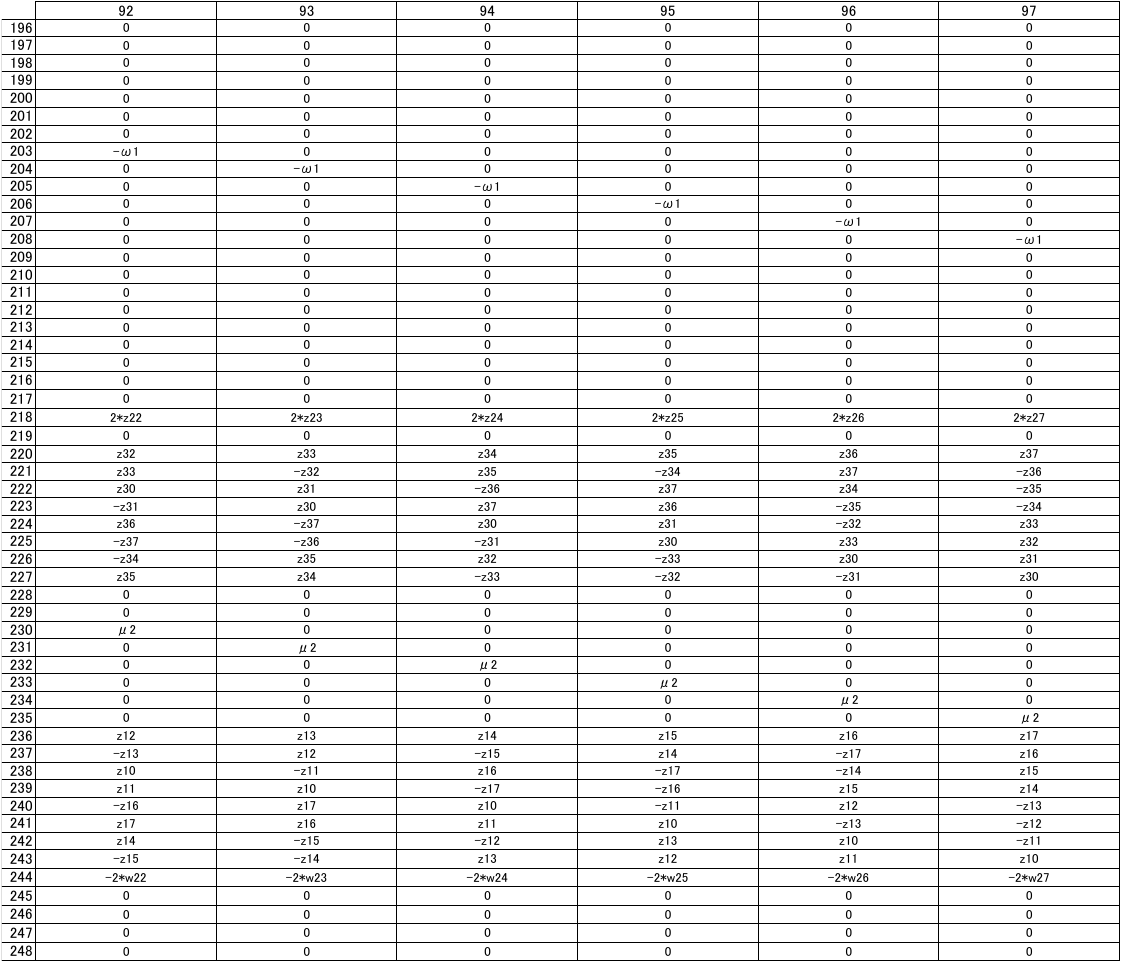}
\end{figure}

\begin{figure}[H]
\centering
\includegraphics[width=12cm, height=19.0cm,page= 71 ]{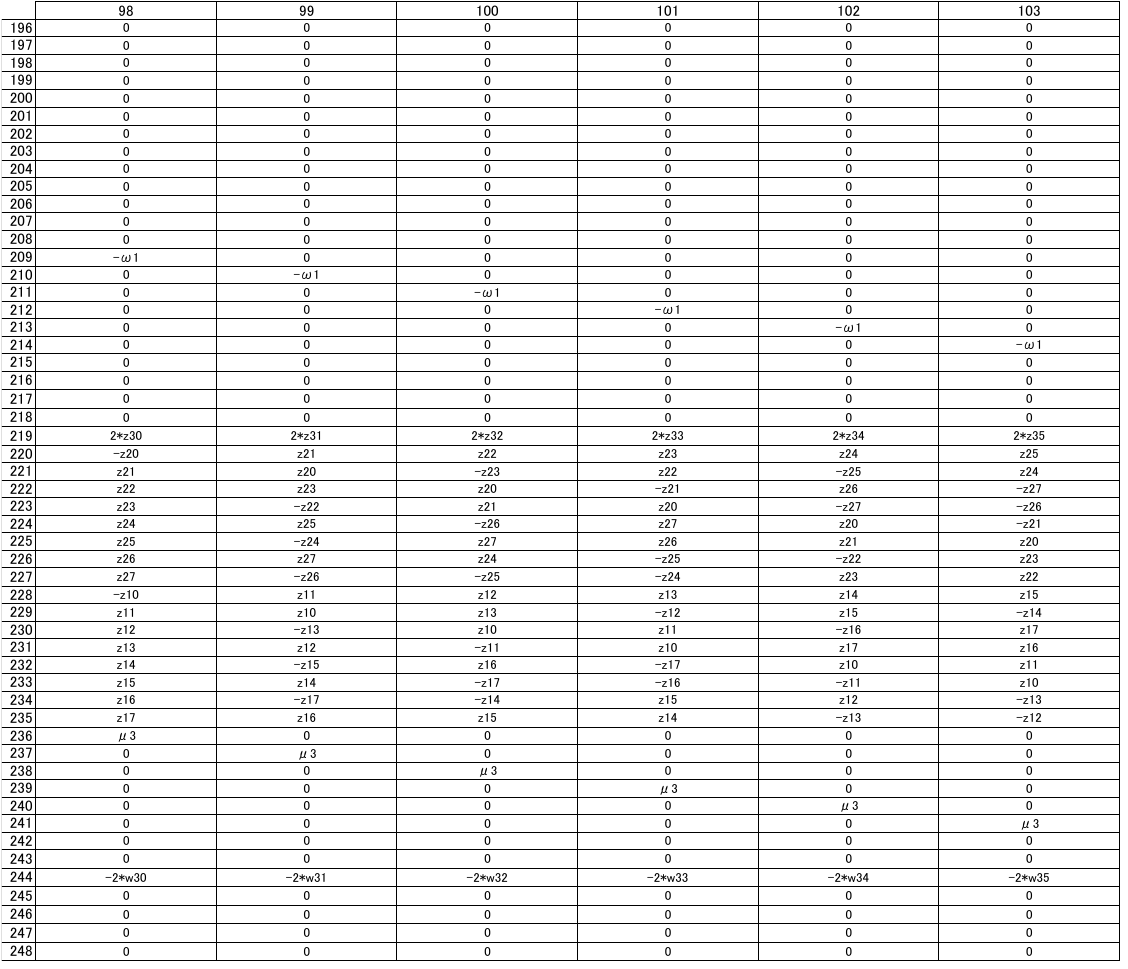}
\end{figure}

\begin{figure}[H]
\centering
\includegraphics[width=12cm, height=19.0cm,page= 72 ]{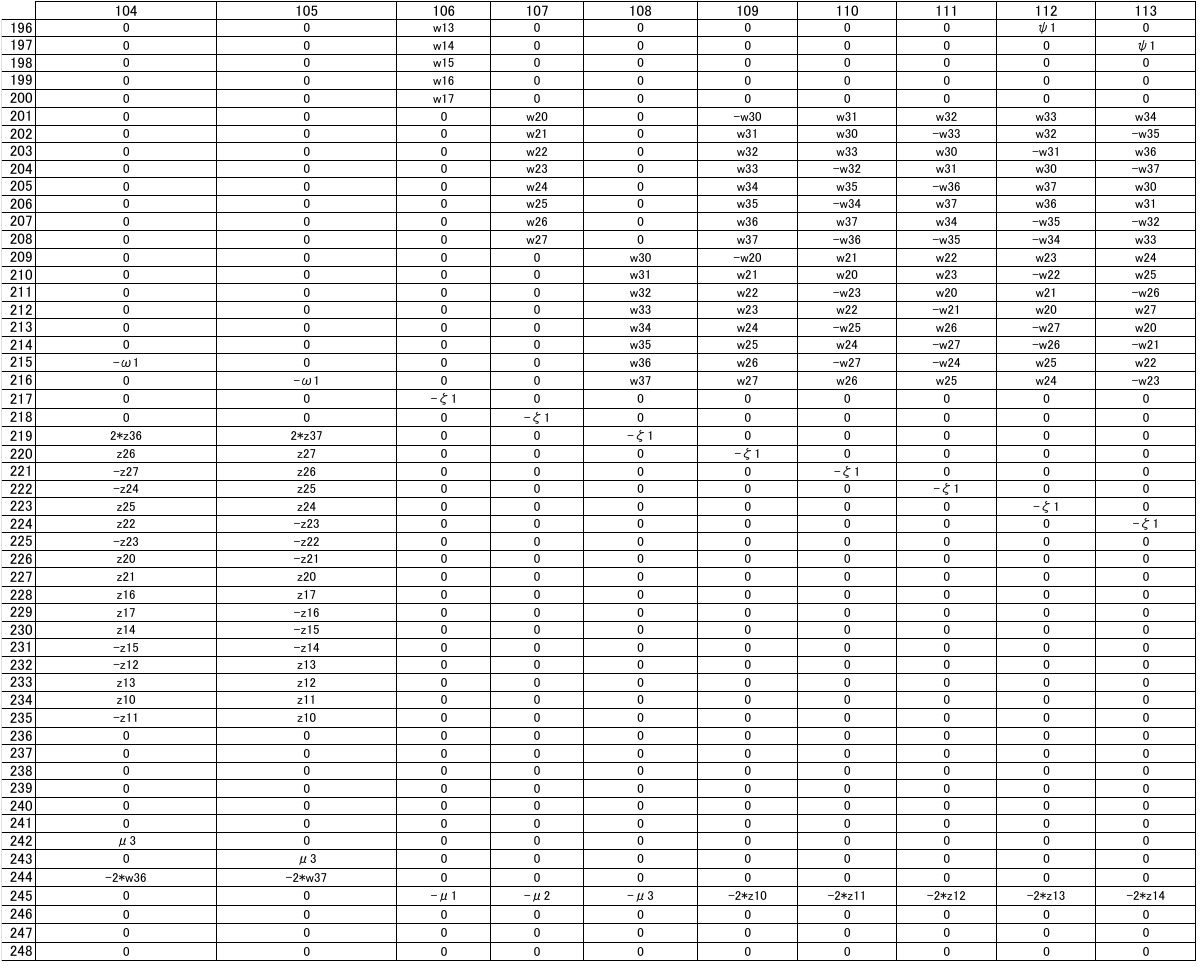}
\end{figure}

\begin{figure}[H]
\centering
\includegraphics[width=12cm, height=19.0cm,page= 73 ]{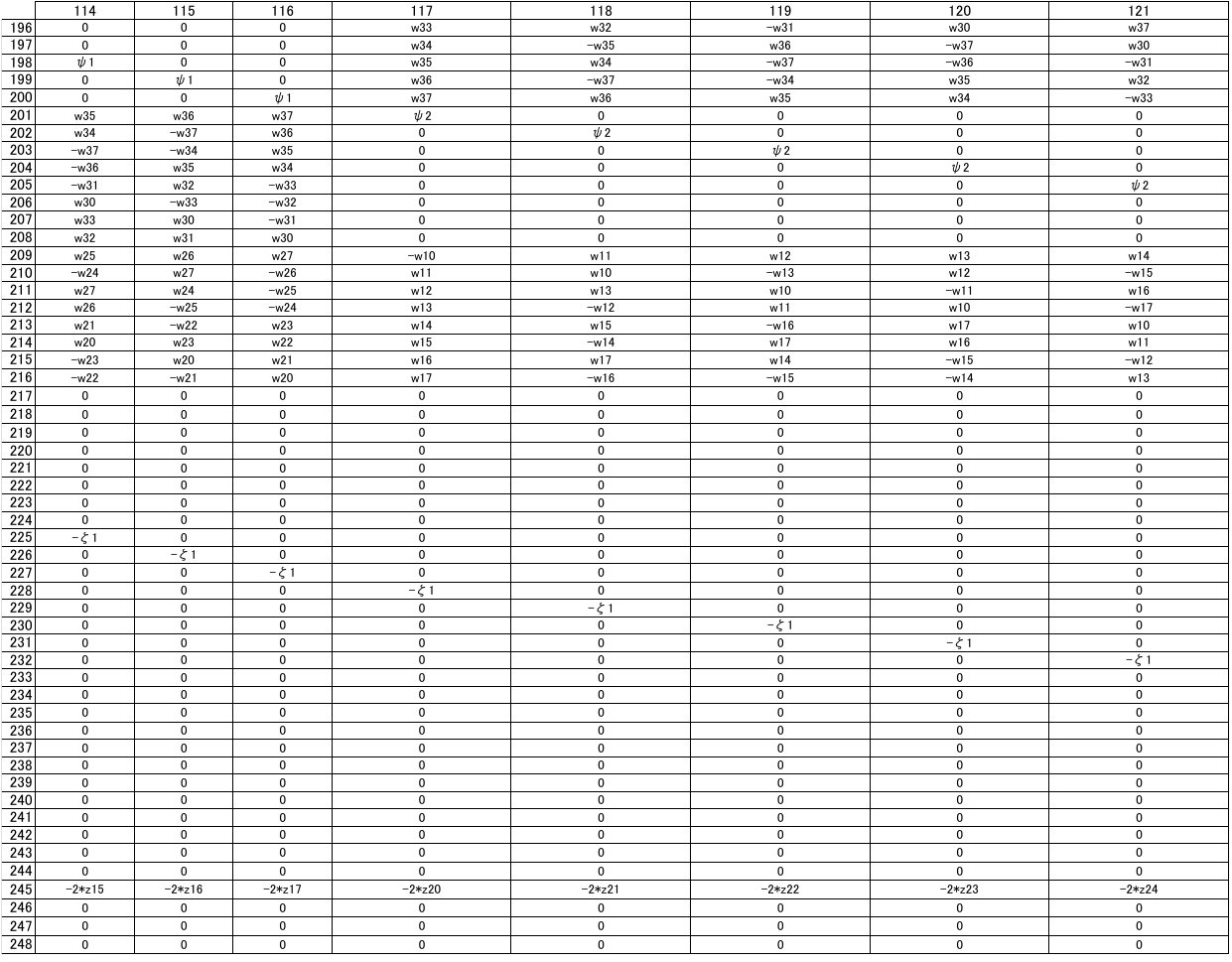}
\end{figure}

\begin{figure}[H]
\centering
\includegraphics[width=12cm, height=19.0cm,page= 74 ]{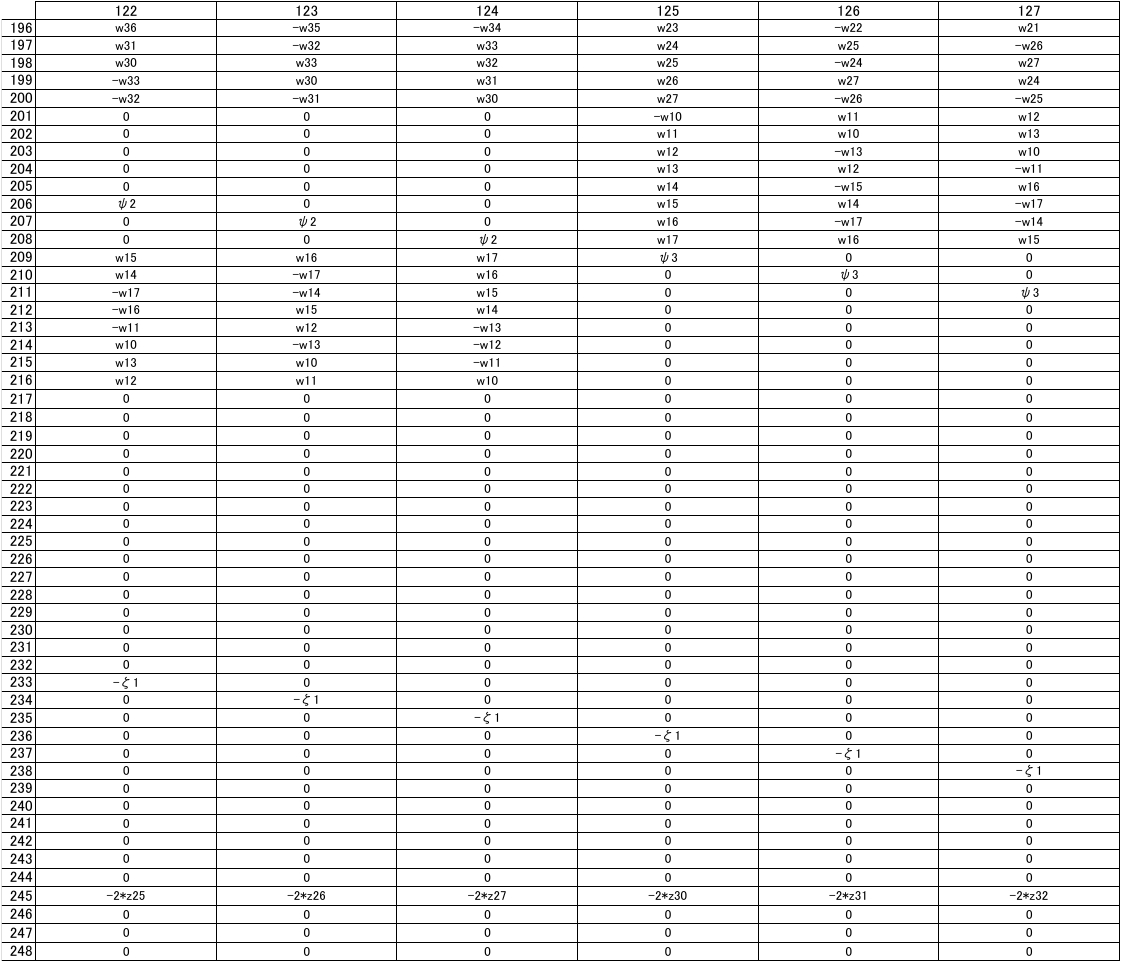}
\end{figure}

\begin{figure}[H]
\centering
\includegraphics[width=12cm, height=19.0cm,page= 75 ]{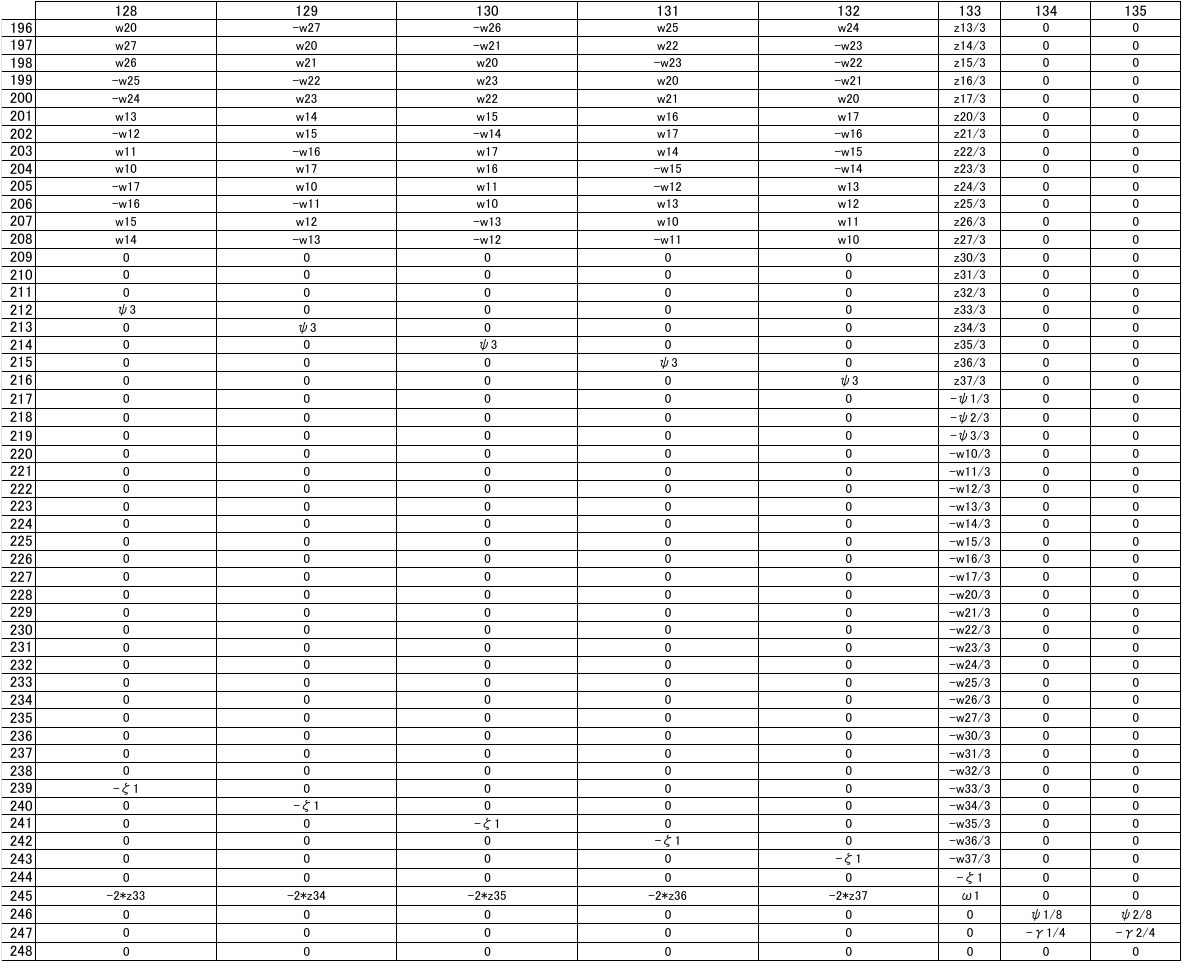}
\end{figure}

\begin{figure}[H]
\centering
\includegraphics[width=12cm, height=19.0cm,page= 76 ]{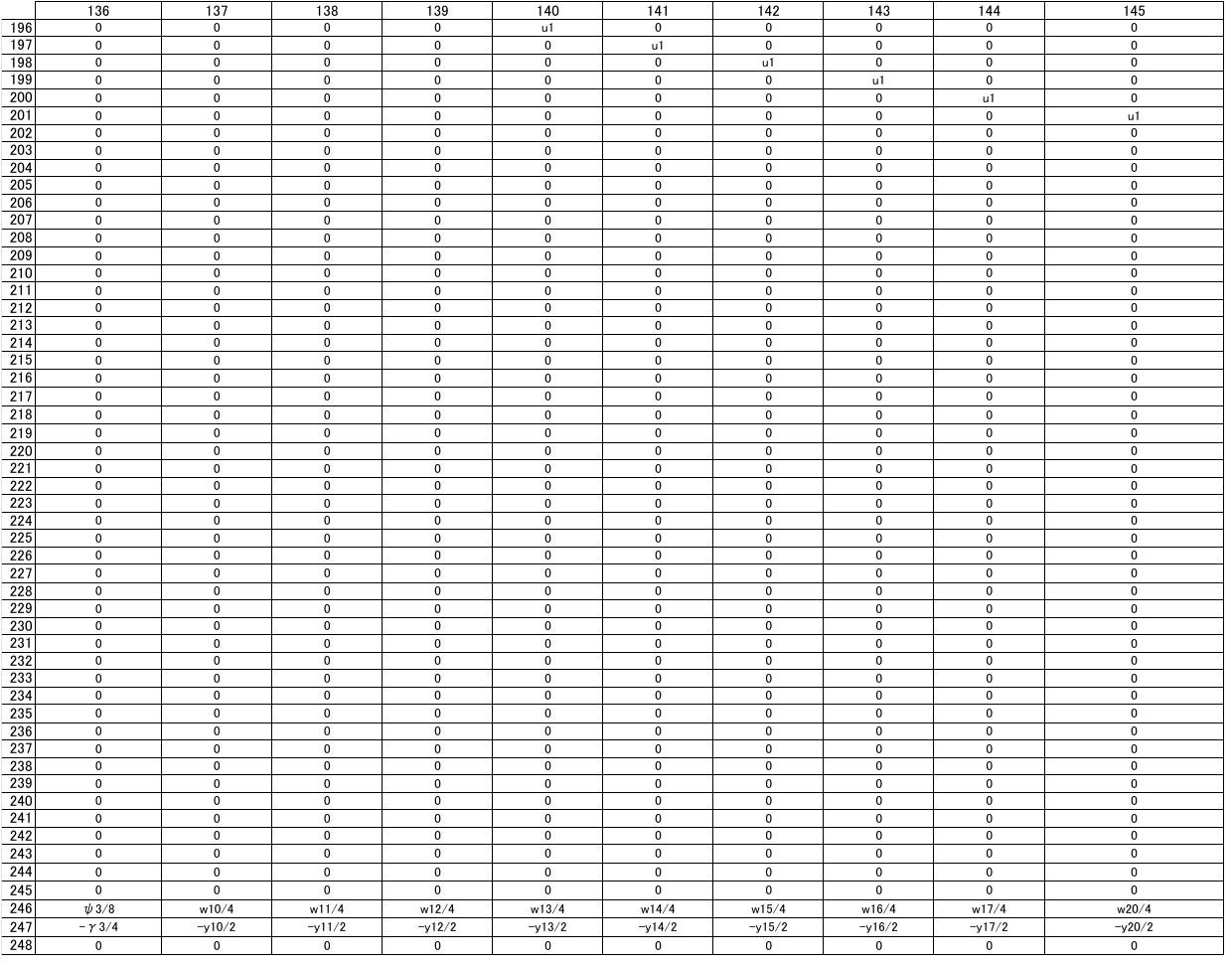}
\end{figure}

\begin{figure}[H]
\centering
\includegraphics[width=12cm, height=19.0cm,page= 77 ]{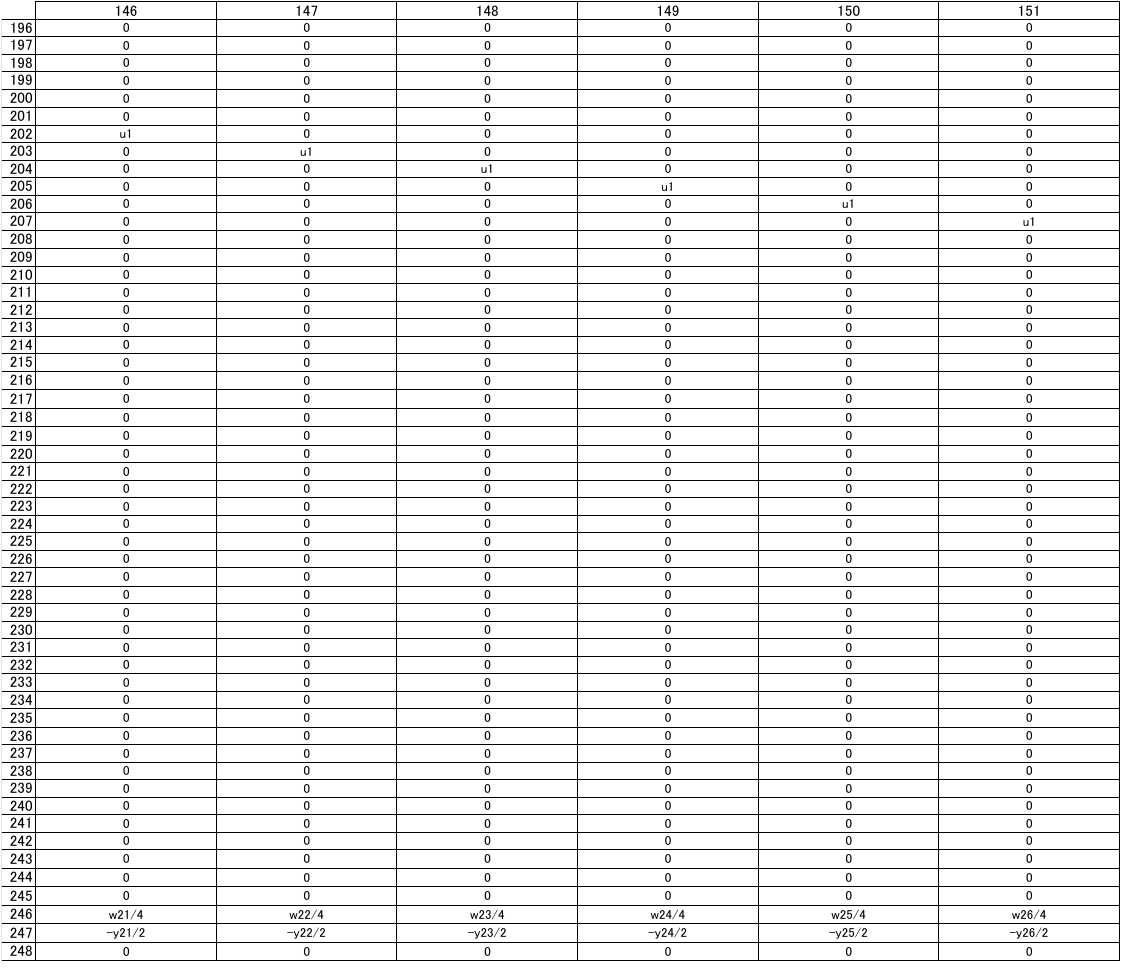}
\end{figure}

\begin{figure}[H]
\centering
\includegraphics[width=12cm, height=19.0cm,page= 78 ]{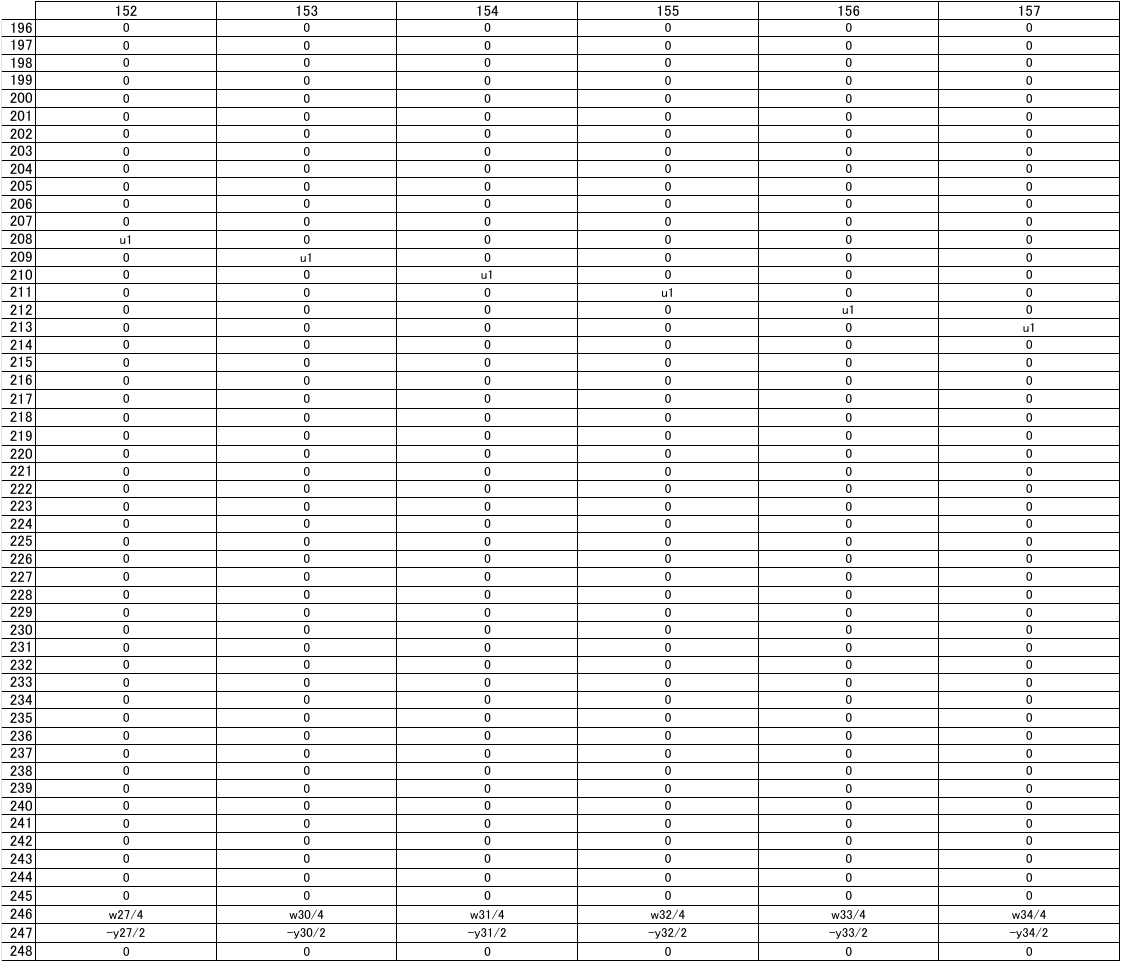}
\end{figure}

\begin{figure}[H]
\centering
\includegraphics[width=12cm, height=19.0cm,page= 79 ]{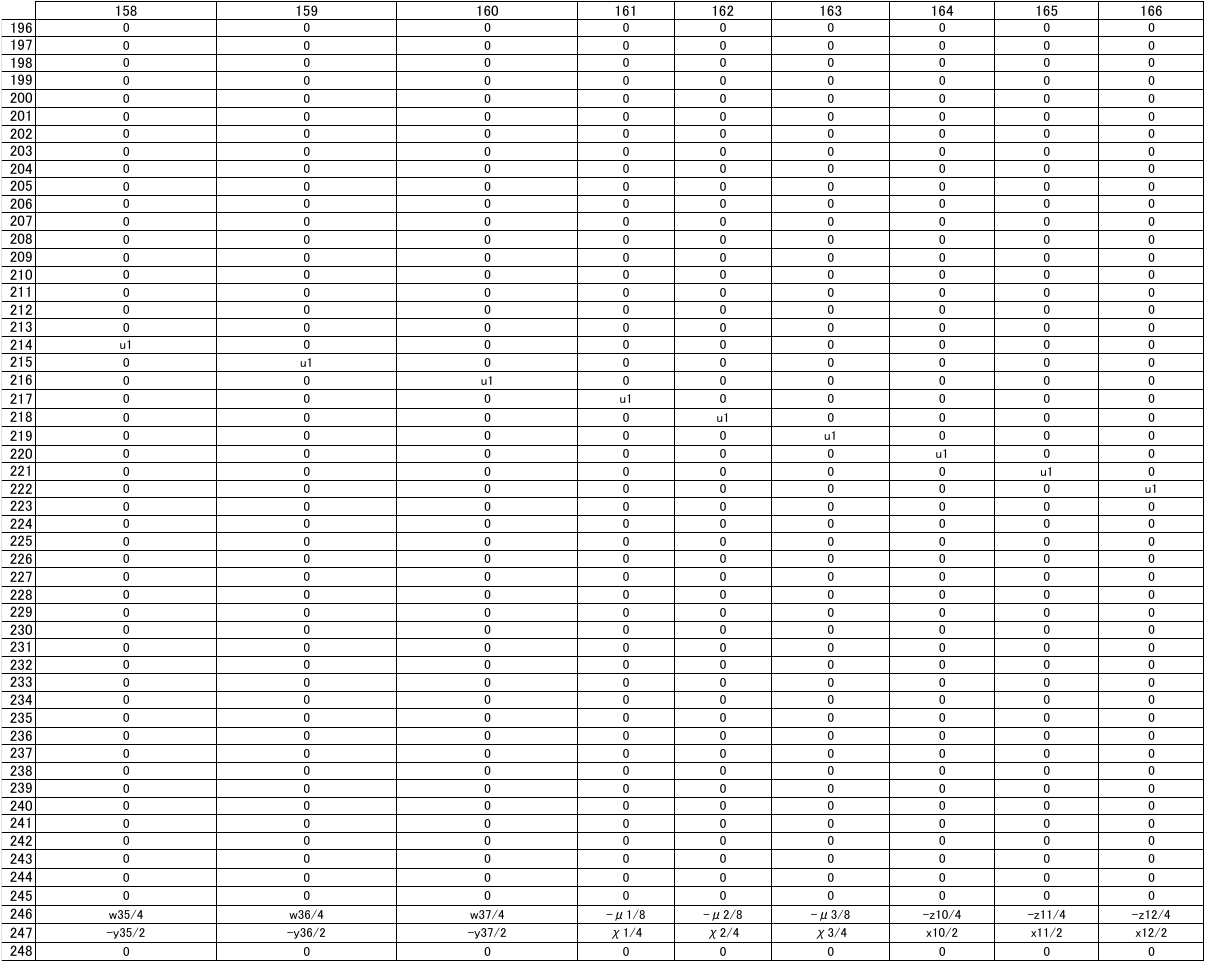}
\end{figure}

\begin{figure}[H]
\centering
\includegraphics[width=12cm, height=19.0cm,page= 80 ]{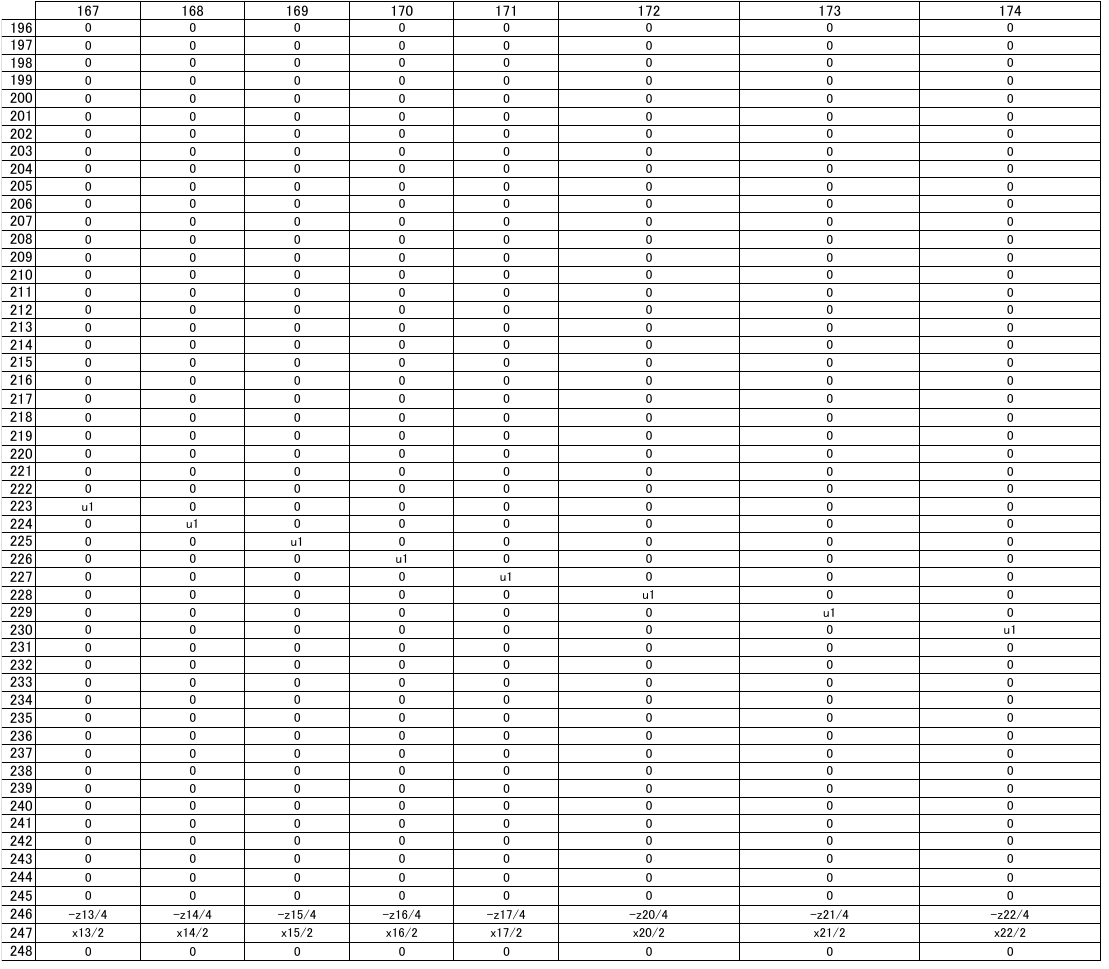}
\end{figure}

\begin{figure}[H]
\centering
\includegraphics[width=12cm, height=19.0cm,page= 81 ]{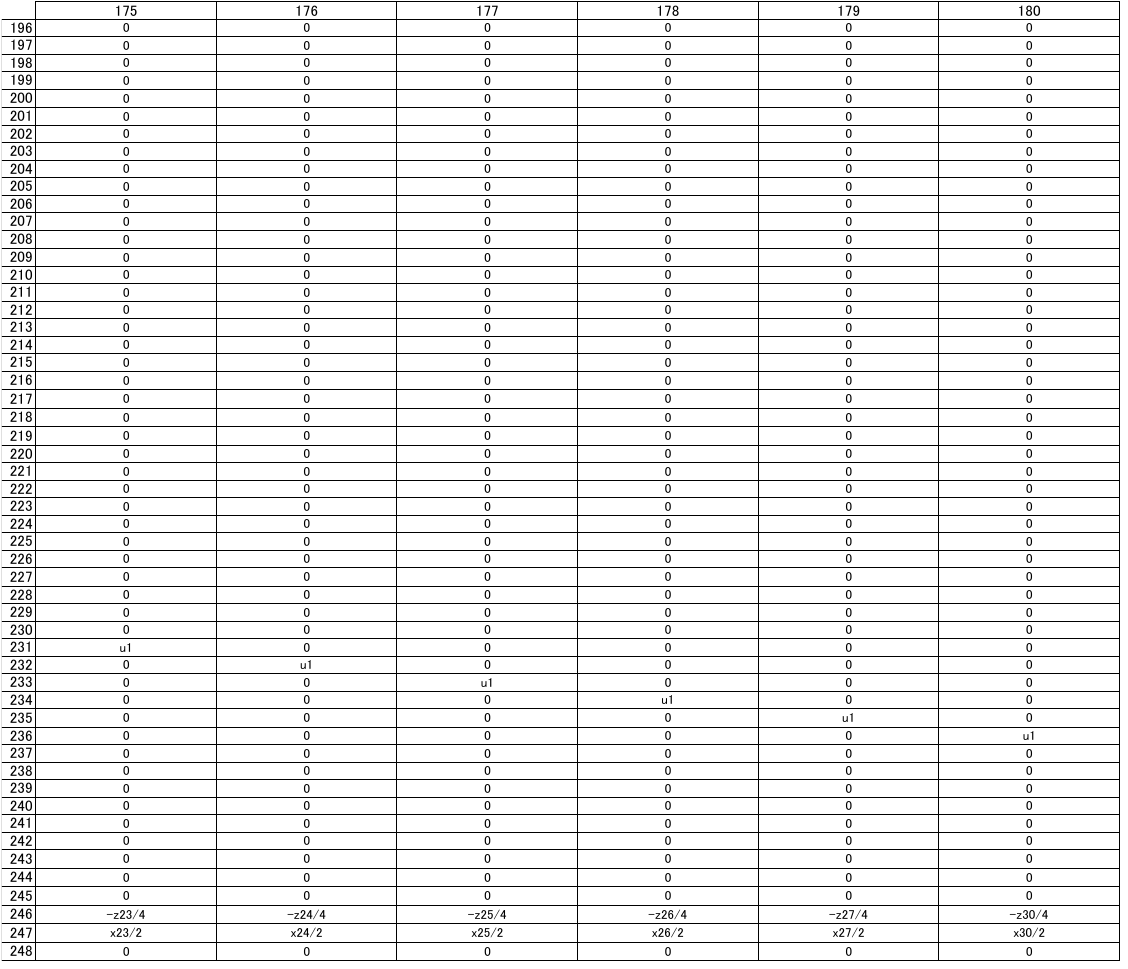}
\end{figure}

\begin{figure}[H]
\centering
\includegraphics[width=12cm, height=19.0cm,page= 82 ]{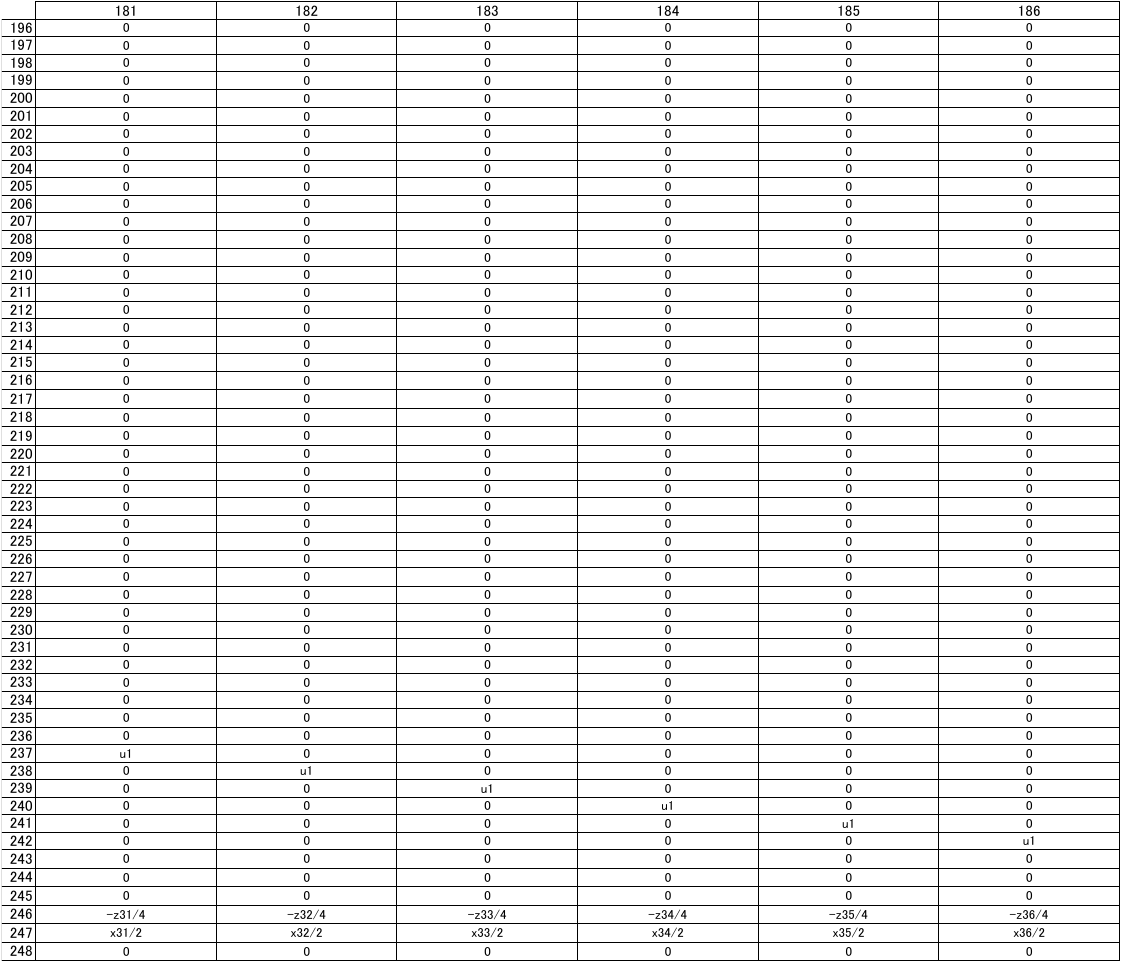}
\end{figure}

\begin{figure}[H]
\centering
\includegraphics[width=12cm, height=19.0cm,page= 83 ]{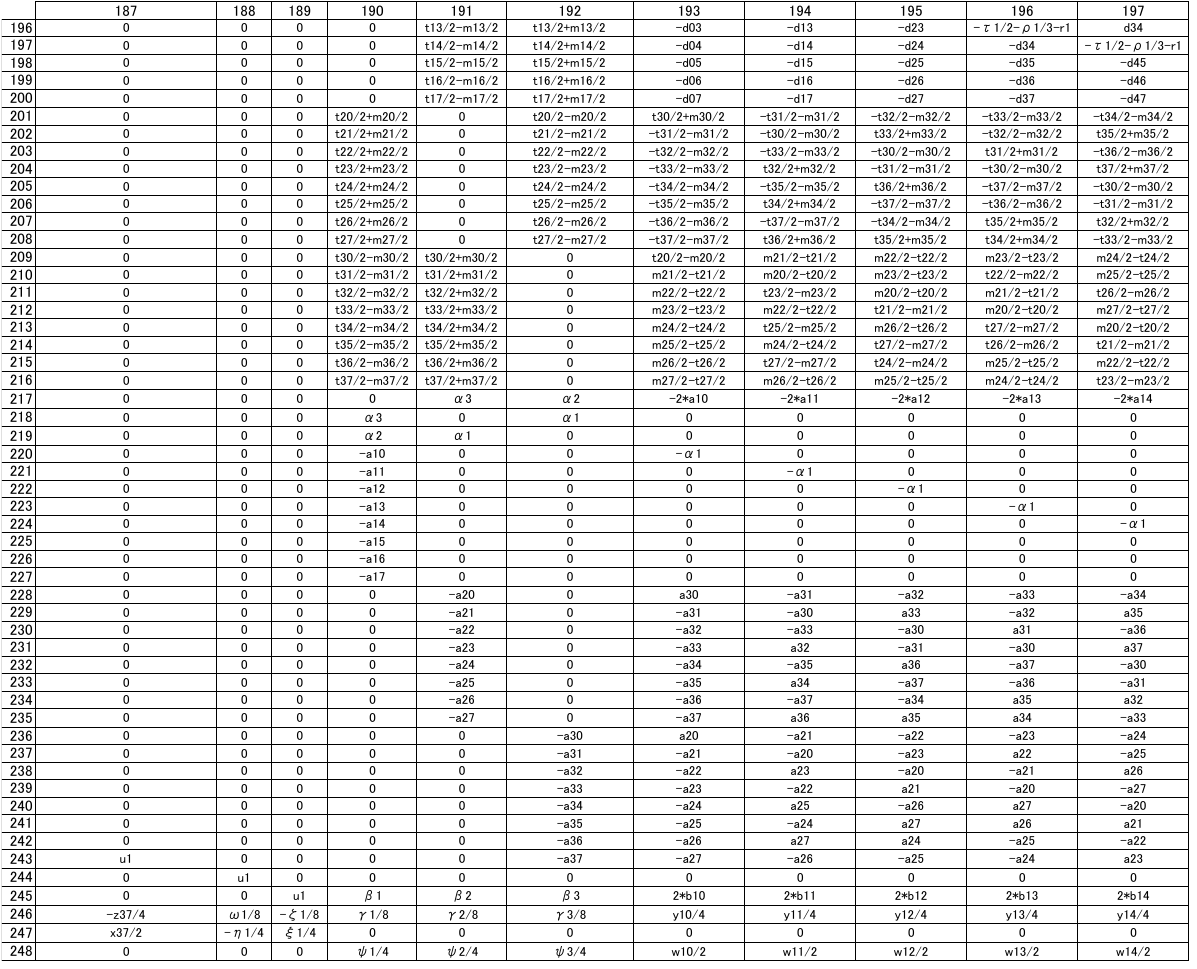}
\end{figure}

\begin{figure}[H]
\centering
\includegraphics[width=12cm, height=19.0cm,page= 84 ]{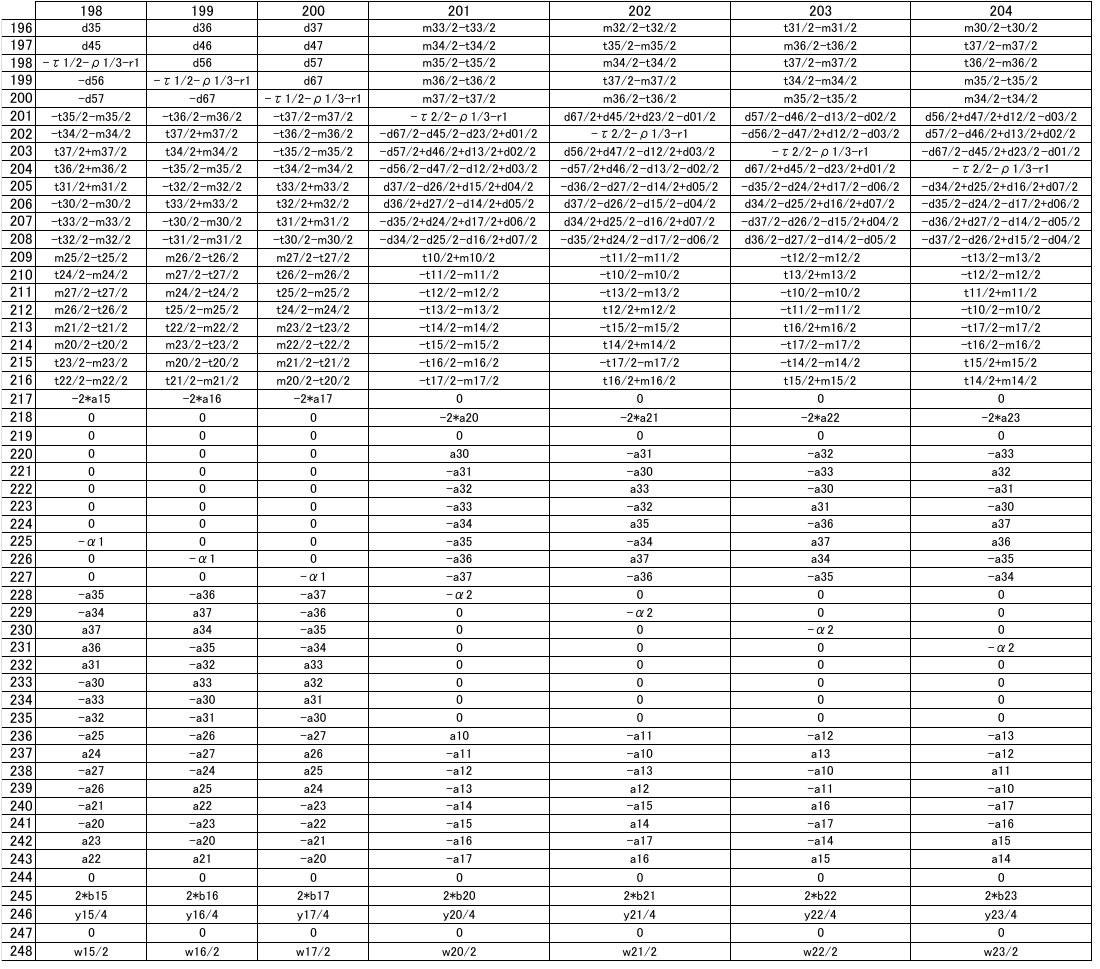}
\end{figure}

\begin{figure}[H]
\centering
\includegraphics[width=12cm, height=19.0cm,page= 85 ]{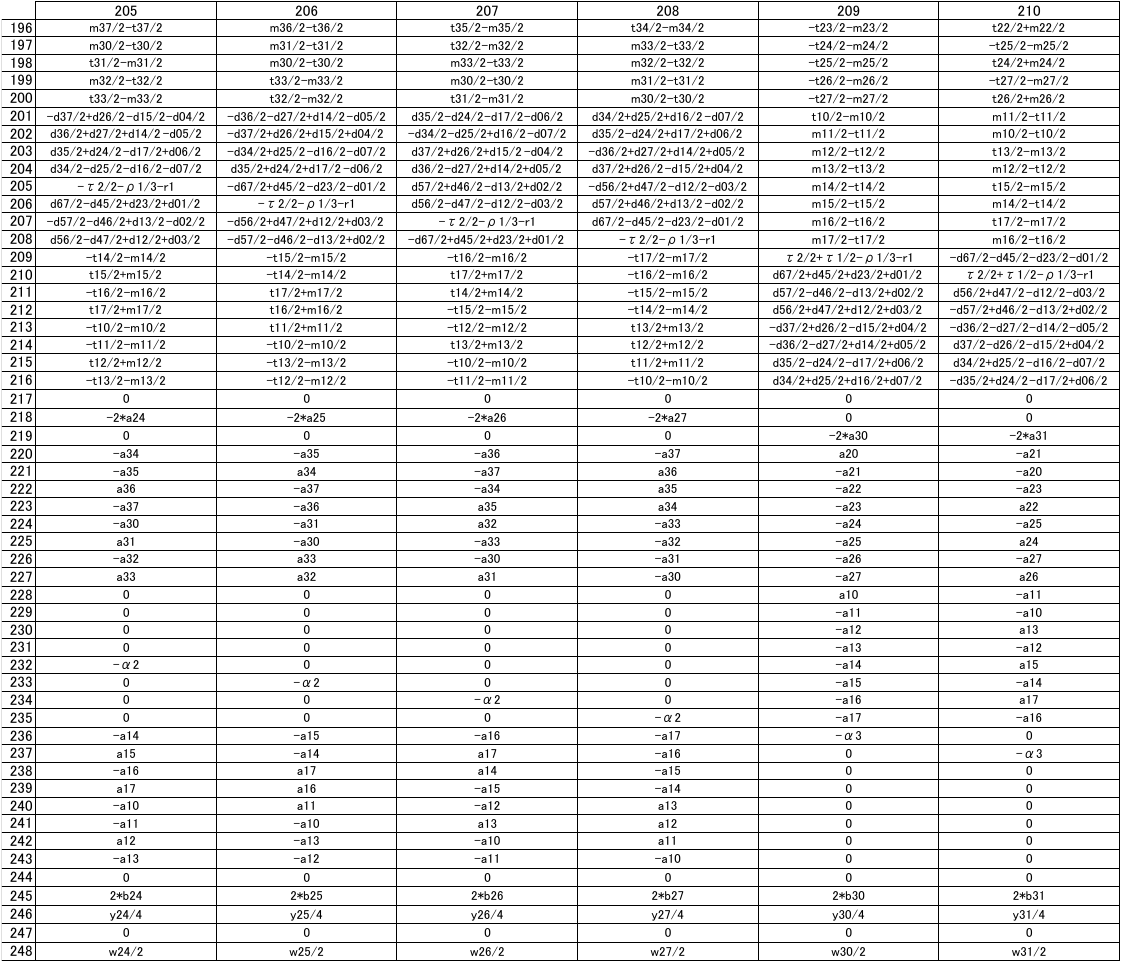}
\end{figure}

\begin{figure}[H]
\centering
\includegraphics[width=12cm, height=19.0cm,page= 86 ]{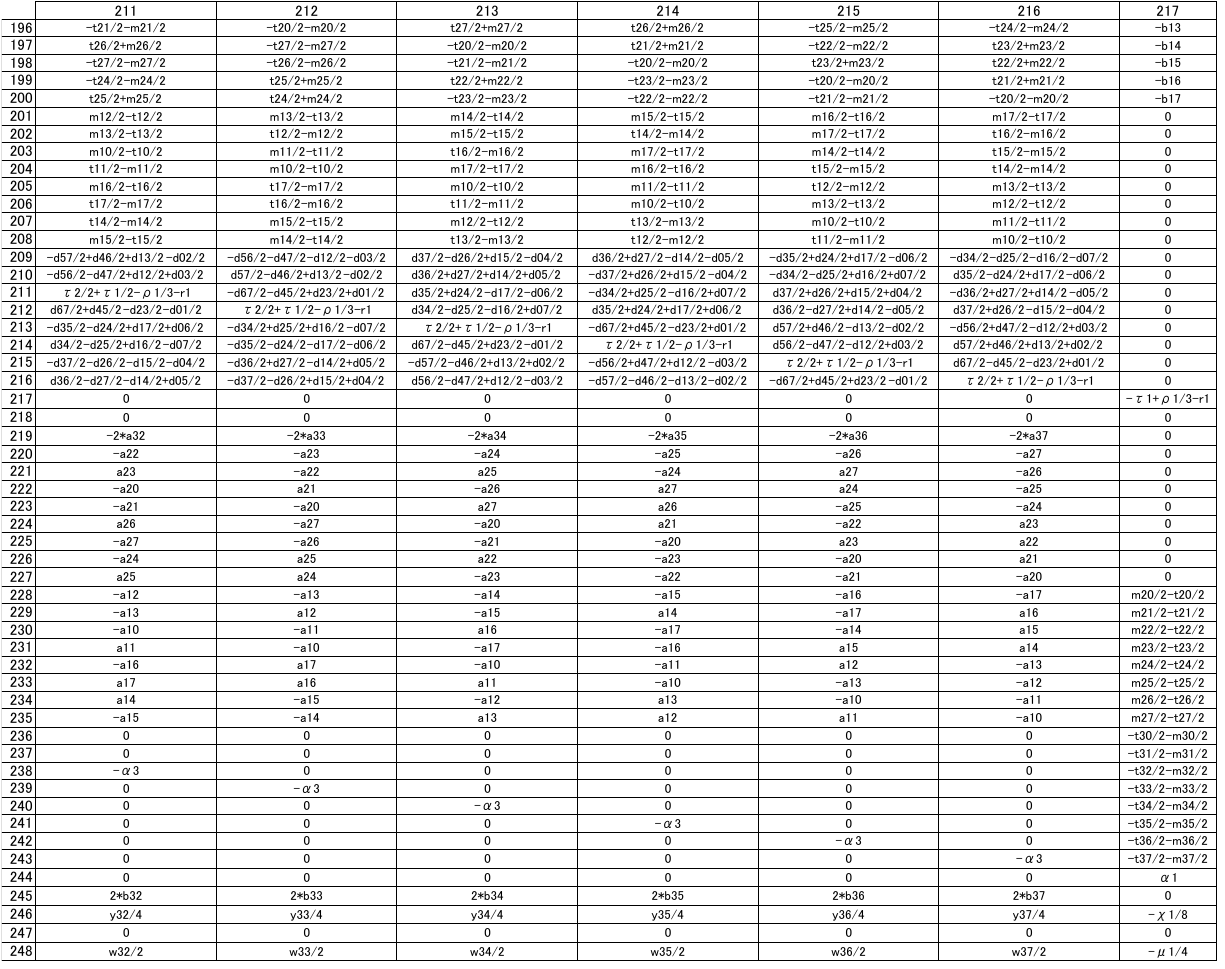}
\end{figure}

\begin{figure}[H]
\centering
\includegraphics[width=12cm, height=19.0cm,page= 87 ]{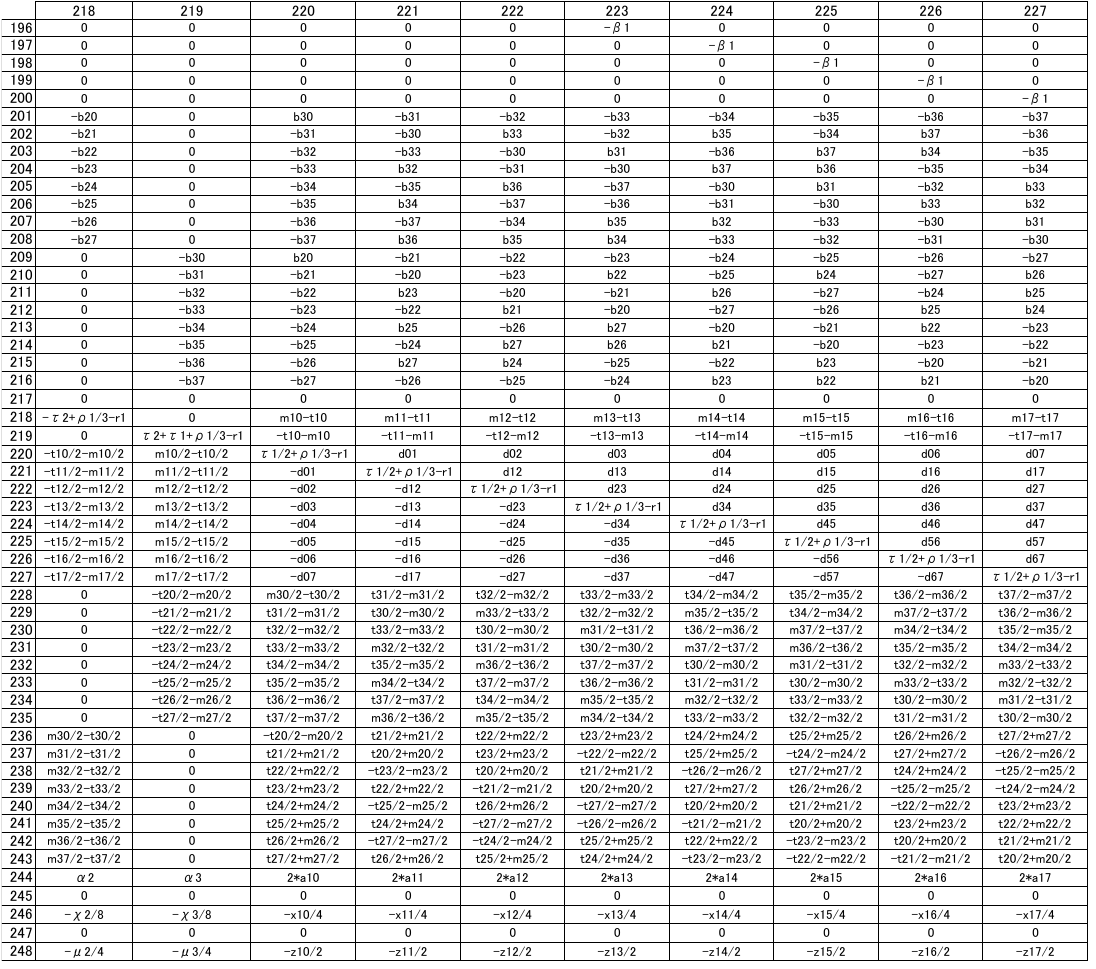}
\end{figure}

\begin{figure}[H]
\centering
\includegraphics[width=12cm, height=19.0cm,page= 88 ]{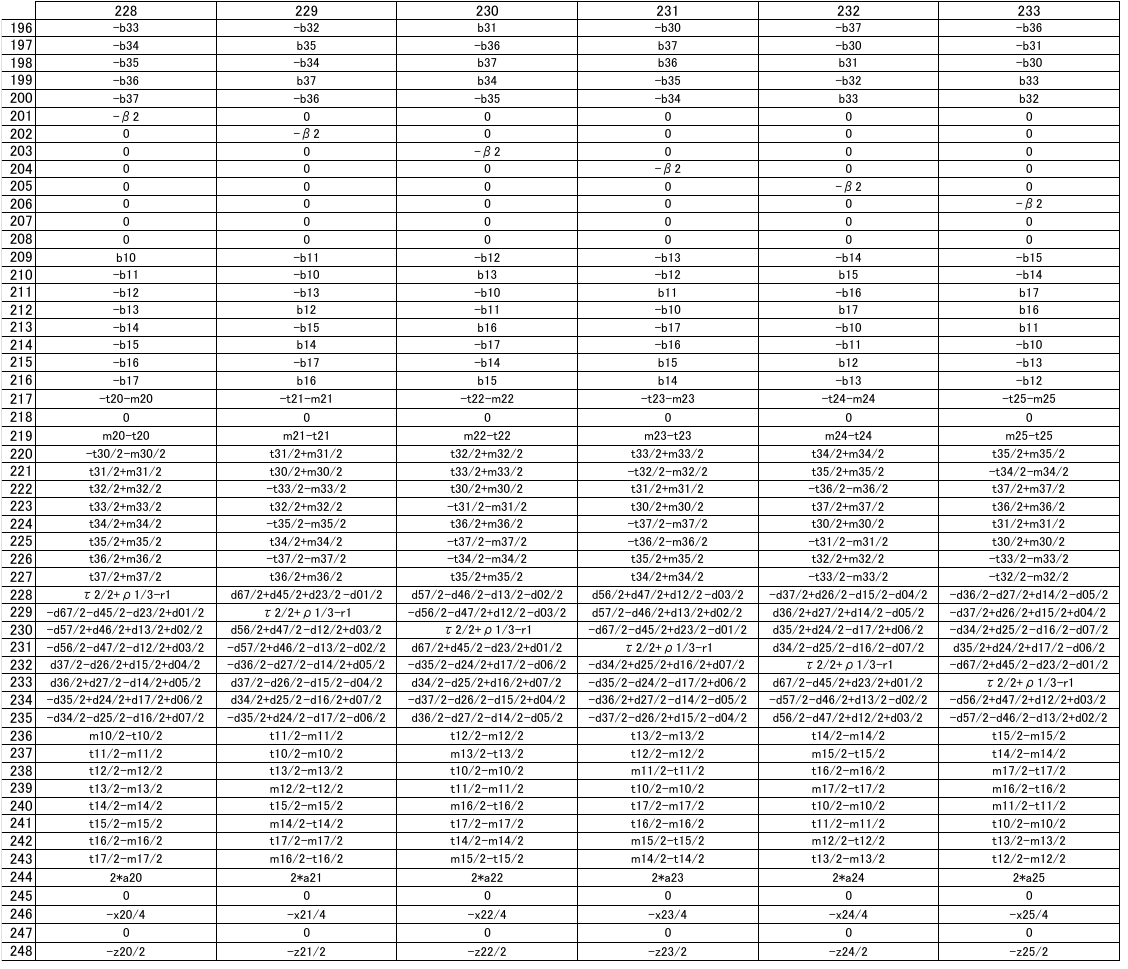}
\end{figure}

\begin{figure}[H]
\centering
\includegraphics[width=12cm, height=19.0cm,page= 89 ]{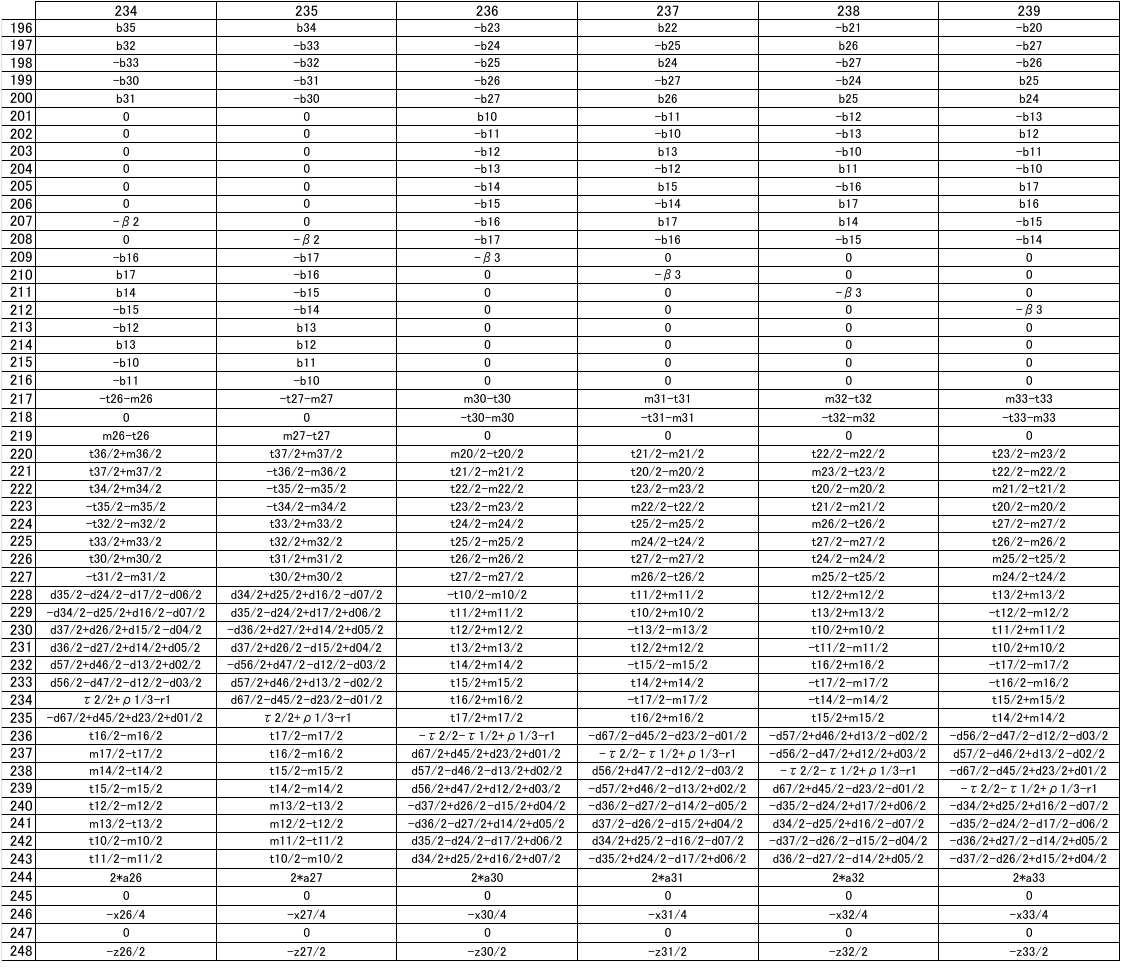}
\end{figure}

\begin{figure}[H]
\centering
\includegraphics[width=12cm, height=19.0cm,page= 90 ]{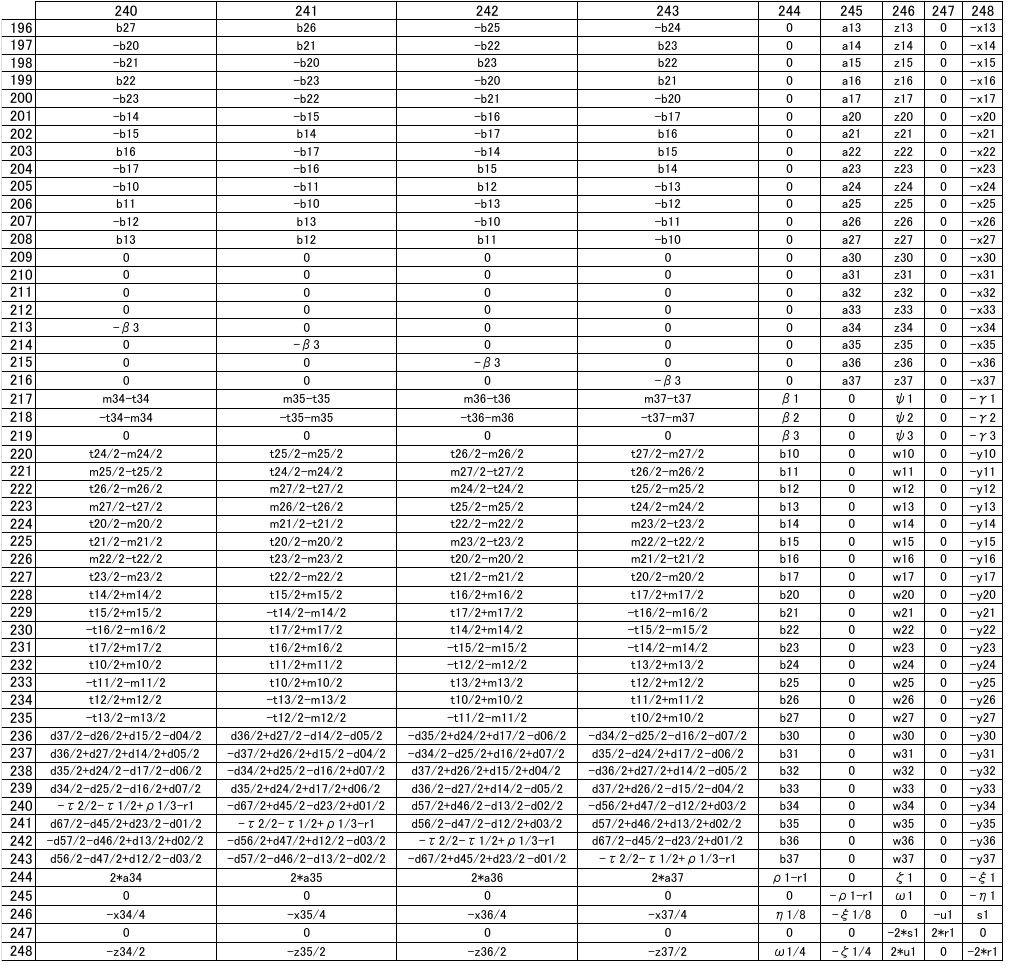}
\end{figure}

\bigskip

\part{Root systems}

\section{Bases of \gr$_{8}^{\C}$}

\ \ \ \ By using \emph{Theorem 10.7}, \gr$_{8}^{\C}$ is
expressed on the matrix $M(248\times 248,\C)$ as in Section 11.

\gr$_{8}^{\C}$ is an algebra made up of 248$\times 248$-dimensional matrices.

\bigskip

\ \emph{Definition 12.1.} We define the followings:

$Rd_{ij}=X\in $\gr$_{8}^{\C}$, $X$ is an element related to $d_{ij},(0\leq
i<j\leq 7),$

$Rd_{ji}=-Rd_{ij},(0\leq i<j\leq 7),Rd_{ii}=0,(0\leq i\leq 7)$

$Rm_{ij}=X\in $\gr$_{8}^{\C}$, $X$ is an element related to $m_{ij}$ $%
,(i=1,2,3,0\leq j\leq 7),$

$R\tau _{i}=X\in $\gr$_{8}^{\C}$, $X$ is an element related to $\tau
_{i},(i=1,2),$

$Rt_{ij}=X\in $\gr$_{8}^{\C}$ $,$ $X$ is an element related to $t_{ij}$ $%
,(i=1,2,3,0\leq j\leq 7),$

$R\alpha _{i}=X\in $\gr$_{8}^{\C}$ $,$ $X$ is an element related to $%
\alpha _{i},(i=1,2,3),$

$Ra_{ij}=X\in $\gr$_{8}^{\C}$ $,$ $X$ is an element related to $a_{ij}$ $%
,(i=1,2,3,0\leq j\leq 7),$

$R\beta _{i}=X\in $\gr$_{8}^{\C}$ $,$ $X$ is an element related to $%
\beta _{i},(i=1,2,3),$

$Rb_{ij}=X\in $\gr$_{8}^{\C}$ $,$ $X$ is an element related to $b_{ij}$ $%
,(i=1,2,3,0\leq j\leq 7),$

$R\rho _{1}=X\in $\gr$_{8}^{\C}$ $,$ $X$ is an element related to $\rho
_{1},$

$R\chi _{i}=X\in $\gr$_{8}^{\C}$ $,$ $X$ is an element related to $\chi
_{i},(i=1,2,3),$

$Rx_{ij}=X\in $\gr$_{8}^{\C}$ $,$ $X$ is an element related to $x_{ij}$ $%
,(i=1,2,3,0\leq j\leq 7),$

$R\gamma _{i}=X\in $,\gr$_{8}^{\C}$ $,$ $X$ is an element related to $%
\gamma _{i},(i=1,2,3),$

$Ry_{ij}=X\in $\gr$_{8}^{\C}$ $,$ $X$ is an element related to $y_{ij}$ $%
,(i=1,2,3,0\leq j\leq 7),$

$R\xi _{1}=X\in $\gr$_{8}^{\C}$ $,$ $X$ is an element related to $\xi _{1},$

$R\eta _{1}=X\in $\gr$_{8}^{\C}$ $,$ $X$ is an element related to $\eta
_{1},$

$R\mu _{i}=X\in $\gr$_{8}^{\C}$ $,$ $X$ is an element related to $\mu
_{i},(i=1,2,3),$

$Rz_{ij}=X\in $\gr$_{8}^{\C}$ $,$ $X$ is an element related to $z_{ij}$ $%
,(i=1,2,3,0\leq j\leq 7),$

$R\psi _{i}=X\in $\gr$_{8}^{\C}$ $,$ $X$ is an element related to $\psi
_{i},(i=1,2,3),$

$Rw_{ij}=X\in $,\gr$_{8}^{\C}$ $,$ $X$ is an element related to $w_{ij}$ $%
,(i=1,2,3,0\leq j\leq 7),$

$R\zeta _{1}=X\in $\gr$_{8}^{\C}$ $,$ $X$ is an element related to $%
\zeta _{1},$

$R\omega _{1}=X\in $\gr$_{8}^{\C}$ $,$ $X$ is an element related to $%
\omega _{1},$

$Rr_{1}=X\in $\gr$_{8}^{\C}$ $,$ $X$ is an element related to $r_{1},$

$Rs_{1}=X\in $\gr$_{8}^{\C}$ $,$ $X$ is an element related to $s_{1},$

$Ru_{1}=X\in $\gr$_{8}^{\C}$ $,$ $X$ is an element related to $u_{1},$

\gr\gd$=\{Rd_{ij}\in $\gr$_{8}^{\C} \mid d_{ij}\in \R,0\leq%
i<j\leq 7\},$ (dimension:28),

\gR\gm$=\{Rm_{ij}\in $\gr$_{8}^{\C} \mid m_{ij}\in \R,1\leq i\leq%
3,0\leq j\leq 7\},$ (dimension:24),

\gR\gt$=\{R\tau _{1},R\tau _{2},Rt_{ij}\in $\gr$_{8}^{\C} \mid \tau_{1},\tau_{2},t_{ij}\in \R,1\leq i\leq 3,0\leq j\leq 7\},$

(dimension:26),

\gR\ga$=\{R\alpha _{i},Ra_{ij}\in $\gr$_{8}^{\C} \mid \alpha _{i},a_{ij}\in \R,1\leq i\leq 3,0\leq j\leq 7\},$
(dimension:27),

\gR\gb$=\{R\beta _{i},Rb_{ij}\in $\gr$_{8}^{\C} \mid \beta _{i},b_{ij}\in \R,1\leq i\leq 3,0\leq j\leq 7\},$
(dimension:27),

\gR\gp$ =\{R\rho _{1}\in $\gr$_{8}^{\C} \mid \rho _{1}\in \R\}, $
(dimension:1),

\gR\gx$=\{R\chi _{i},Rx_{ij}\in $\gr$_{8}^{\C} \mid \chi _{i},x_{ij}\in \R,1\leq i\leq 3,0\leq j\leq 7\},$
(dimension:27),

\gR\gy$=\{R\gamma _{i},Ry_{ij}\in $\gr$_{8}^{\C} \mid \gamma _{i},y_{ij}\in \R,1\leq i\leq 3,0\leq j\leq 7\},$
(dimension:27),

\gR\gk$ =\{R\xi _{1}\in $\gr$_{8}^{\C} \mid \xi _{1}\in \R\},$
(dimension:1),

\gR\gi$ =\{R\eta _{1}\in $\gr$_{8}^{\C} \mid \eta _{1}\in \R\},$
(dimension:1),

\gR\gz$=\{R\mu _{i},Rz_{ij}\in $\gr$_{8}^{\C} \mid \mu _{i},z_{ij}\in \R,1\leq i\leq 3,0\leq j\leq 7\},$
(dimension:27),

\gR\gw$=\{R\psi _{i},Rw_{ij}\in $\gr$_{8}^{\C} \mid \psi _{i},w_{ij}\in \R,1\leq i\leq 3,0\leq j\leq 7\},$
(dimension:27),

\gR\gl$=\{R\zeta _{1}\in $\gr$_{8}^{\C} \mid \zeta _{1}\in \R\},$
(dimension:1),

\gR\go$ =\{R\omega _{1}\in $\gr$_{8}^{\C} \mid \omega _{1}\in \R\},$
(dimension:1),

\gR\gr$=\{Rr_{1}\in $\gr$_{8}^{\C} \mid r_{1}\in \R\},$
(dimension:1),

\gR\gs$=\{Rs_{1}\in $\gr$_{8}^{\C} \mid s_{1}\in \R\},$
(dimension:1),

\gR\gu$=\{Ru_{1}\in $\gr$_{8}^{\C} \mid u_{1}\in \R\}.$
(dimension:1).

$\{R*_{ij} \mid *_{ij} \in \C\}$ is simply represented as $\{R*_{ij} \}$, for example $\{Rx_{ij} \}$ means $\{Rx_{ij} \mid x_{ij} \in \C\}$.
Similarly $\{ [Rx_{ij},Ry_{kl}]\}$ means $\{ [Rx_{ij},Ry_{kl}] \mid x_{ij},y_{kl} \in \C \}$. 

\gr$_{4}^{\C}$=\ \gr\gd$^{\C} \bigoplus $\gR\gm$^{\C}\ ,$

\gr$_{6}^{\C}$=\ \gr$_{4}^{\C} \bigoplus $\gR\gt$^{\C}\ ,$

\gr$_{7}^{\C}$=\ \gr$_{6}^{\C} \bigoplus $\gR\ga$^{\C} \bigoplus $\gR\gb$^{\C} \bigoplus $\gR\gp$^{\C} \ ,$

\gr$_{8}^{\C}$=\ \gr$_{7}^{\C} \bigoplus $\gR\gx$^{\C} \bigoplus $\gR\gy$^{\C} \bigoplus $\gR\gk$^{\C} \bigoplus %
$\gR\gi$^{\C} \bigoplus $\gR\gz$^{\C} \bigoplus $\gR\gw$^{\C} \bigoplus $\gR\gl$^{\C} \bigoplus $\gR\go$^{\C}  $

\ \ \ \ \ \ \ \ \ $\bigoplus $\gR\gr$^{\C} \bigoplus $\gR\gs$^{\C} \bigoplus $\gR\gu$^{\C}\ .$

\gP$^{\C}$ in T.Imai and I.Yokota\cite{YokotaImai2} is following: 

\gP$^{\C}=$ \gR\gx$ \oplus $\gR\gy$ \oplus $\gR\gk$ \oplus $\gR\gi \ or \gR\gz$ \oplus $\gR\gw$ \oplus $\gR\gl$ \oplus $\gR\go.

\bigskip

\emph{Definition 12.2. \ }We define the mapping mtv

mtv$:M(248\times 248,\C)\ni x(i,j)\rightarrow C^{248\times 248}\ni v((i-1)\times 248+j)=x(i,j).$

\noindent
mtv is a mapping that converts a matrix into a column vector.

\bigskip

\emph{Lemma 12.3.}\ The following 248 matrices are bases of \gr$_{8}^{\C}$.

$Rd_{ij}(0\leq i<j\leq 7),Rm_{ij}(i=1,2,3,0\leq j\leq 7),$

$R\tau _{i}(i=1,2),Rt_{ij}(i=1,2,3,0\leq j\leq 7),$

$R\alpha _{i}(i=1,2,3),Ra_{ij}(i=1,2,3,0\leq j\leq 7),$

$R\beta _{i}(i=1,2,3),Rb_{ij}(i=1,2,3,0\leq j\leq 7),R\rho _{1},$

$R\chi _{i}(i=1,2,3),Rx_{ij}(i=1,2,3,0\leq j\leq 7),$

$R\gamma _{i}(i=1,2,3),Ry_{ij}(i=1,2,3,0\leq j\leq 7),R\xi _{1},R\eta _{1},$

$R\mu _{i}(i=1,2,3),Rz_{ij}(i=1,2,3,0\leq j\leq 7),$

$R\psi _{i}(i=1,2,3),Rw_{ij}(i=1,2,3,0\leq j\leq 7),R\zeta _{1},R\omega _{1},$

$Rr_{1},Rs_{1},Ru_{1}.$

\bigskip

\emph{Proof.} \ We consider the following 248 vectors.

$Vd_{ij}=\textrm{mtv}(Rd_{ij}),(0\leq i<j\leq 7),$(28vectors),

$Vm_{kj}=\textrm{mtv}(Rm_{kj}),(k=1,2,3,0\leq j\leq 7),$(24vectors),

$V\tau _{k}=\textrm{mtv}(R\tau _{k}),(k=1,2),$(2vectors)

$Vt_{kj}=\textrm{mtv}(Rt_{kj}),(k=1,2,3,0\leq j\leq7),$(24vectors),

$V\alpha _{k}=\textrm{mtv}(R\alpha _{k}),(k=1,2,3),$(3vectors)

$Va_{kj}=\textrm{mtv}(Ra_{kj}),(k=1,2,3,0\leq j\leq7),$(24vectors),

$V\beta _{k}=\textrm{mtv}(R\beta _{k}),(k=1,2,3),$(3vectors),

$Vb_{kj}=\textrm{mtv}(Rb_{kj}),(k=1,2,3,0\leq j\leq7),$(24vectors),

$V\rho _{1}=\textrm{mtv}(R\rho _{1}),$(1vectort),

$V\chi _{k}=\textrm{mtv}(R\chi _{k}),(k=1,2,3),$(3vectors),

$Vx_{kj}=\textrm{mtv}(Rx_{kj}),(k=1,2,3,0\leq j\leq7),$(24vectors),

$V\gamma _{k}=\textrm{mtv}(R\gamma _{k}),(k=1,2,3),$(3vectors),

$Vy_{kj}=\textrm{mtv}(Ry_{kj}),(k=1,2,3,0\leq j\leq7),$(24vectors),

$V\xi _{1}=\textrm{mtv}(R\xi _{1}),V\eta _{1}=\textrm{mtv}(R\eta _{1}),$(1+1=2vectors),

$V\mu _{k}=\textrm{mtv}(R\mu _{k}),(k=1,2,3),$(3vectors),

$Vz_{kj}=\textrm{mtv}(Rz_{kj}),(k=1,2,3,0\leq j\leq7),$(24vectors),

$V\psi _{k}=\textrm{mtv}(R\psi _{k}),(k=1,2,3),$(3vectors),

$Vw_{kj}=\textrm{mtv}(Rw_{kj}),(k=1,2,3,0\leq j\leq7),$(24vectors),

$V\zeta _{1}=\textrm{mtv}(R\zeta _{1}),V\omega _{1}=\textrm{mtv}(R\omega _{1}),$(1+1=2vectors),

$Vr_{1}=\textrm{\textrm{mtv}}(Rr_{1}),Vs_{1}=\textrm{mtv}(Rs_{1}),Vu_{1}=\textrm{mtv}(Ru_{1}),$(1+1+1=3vectors).

We prove that the above 248 vectors are independent.
For that, let we put 248-dimensional row vector $MR(l)$ as

$MR(l)=\{Vd_{ij}(l)/d_{ij},Vm_{kj}(l)/m_{kj},V\tau _{k}(l)/\tau _{k},Vt_{kj}(l)/t_{kj},V\alpha%
_{k}(l)/\alpha _{k},$

$Va_{kj}(l)/a_{kj},V\beta _{k}(l)/\beta _{k},Vb_{kj}(l)/b_{kj},$ $V\rho _{1}(l)/\rho _{1},V\chi _{k}(l)/\chi%
_{k},Vx_{kj}(l)/x_{kj},$

$V\gamma _{k}(l)/\gamma _{k},Vy_{kj}(l)/y_{kj},-V\xi _{1}(l)/\xi _{1},V\eta _{1}(l)/\eta _{1},V\mu _{k}(l)/\mu _{k},%
Vz_{kj}(l)/z_{kj},$

$V\psi _{k}(l)/\psi _{k},Vw_{kj}(l)/w_{kj},V\zeta _{1}(l)/\zeta _{1},V\omega _{1}(l)/\omega1,Vr_{1}(l)/r_{1},Vs_{1}(l)/s_{1},$

$Vu_{1}(l)/u_{1}\},$\ \ where $V(l)$ means\ $l$-th element of $V.$

\noindent
And let we put 248$\times$248-dimensional matrix $MR$ as

$MR= \{ MR(l) \mid l=$ 2,3,4,5,6,7,251,252,253,254,255,256,500,501,502,503,

749,750,751,998,999,1247,1517,1737,1985,2233,2481,2729,2977,7193,7448,

7702,7955,8207,8458,8708,9018,12933,12934,12935,12936,12937,12938,

12939,12940,12967,12968,12969,12970,12971,12972,12973,12974,13189,

13190,13191,13192,13193,13194,13195,13196,13421,13641,13896,14150,

14403,14655,14906,15156,15429,15677,15925,16173,16421,16669,16917,

17165,19381,19382,19383,19384,19385,19386,19387,19388,19423,19637,

19638,19639,19640,19641,19642,19643,19644,19681,19876,19893,20337,

20592,20846,21099,21351,21602,21852,26077,26078,26079,26080,26081,

26082,26083,26084,26093,26333,26334,26335,26336,26337,26338,26339,

26340,26342,26572,27033,27288,27542,27795,28047,28298,28548,32842,

32843,33021,33022,33023,33024,33025,33026,33027,33028,33118,33174,

33277,33278,33279,33280,33281,33282,33283,33284,33516,33533,33977,

34232,34486,34739,34991,35242,35492,37821,39717,39718,39719,39720,

39721,39722,39723,39724,39733,39786,39973,39974,39975,39976,39977,

39978,39979,39980,39982,40212,40506,40673,40928,41182,41435,41687,

41938,42188,42498,46564,46757,46909,46910,46911,46912,46913,46914,

46915,46916,47006,47165,47166,47167,47168,47169,47170,47171,47172,

47404,47421,47865,48120,48374,48627,48879,49130,49380,51709,53605,

53606,53607,53608,53609,53610,53611,53612,53621,53861,53862,53863,

53864,53865,53866,53867,53868,53870,54100,54394,54561,54816,55070,

55323,55575,55826,56076,56386,60510,60645\}

\noindent
Then by calculate the rank of $MR$ with the \emph{rank} package of Maxima, we have the rank of $MR$ is 248.\ \ \ \emph{\
Q.E.D.}

\bigskip

\emph{Remark 12.4.} In I.Yokota\cite[Section 2.1, 4.1, and 5.3]{Yokota1} ,
inner products $(X,Y)$ in \gJ,\ $(P,Q)$ in \gP$^{\C}$ and  $(R_{1} , R_{2} )_{8} \in $\ge$_{8}^{\C}$
are defined by respectively

$(X,Y)=tr(X \circ Y)$,
$(P,Q)=(X,Z)+(Y,W)+ \xi \zeta +\eta \omega $,

$(R_{1} , R_{2} )_{8} = (\Phi _{1} , \Phi _{2} )_{7}-\{Q_{1} , P_{2} \}+\{P_{1} , Q_{2} \}-8r_{1} r_{2}-4u_{1} s_{2}-s_{1} u_{2}$ ,

\noindent
where $R_{i} = (\Phi _{i} , P_{i} , Q_{i} , r_{i} , s_{i} , u_{i} )\in $\ge$_{8}^{\C}$ .

Transpose of a matrix and transpose with the inner products are different.

\noindent
So, we denote upper left ${}^{t}$ as transpose with a matrix and upper right ${}^{t}$ as transpose with the inner product. 

\bigskip

\section{The Lie algebra \gr\gd$^{\C}$ of type $D_{4}$}

\ \ \ \ \ \emph{Definition 13.1.} \ We define the following:

\ \ \ $CD=%
\{Rd_{01},Rd_{02},Rd_{03},Rd_{04},Rd_{05},Rd_{06},Rd_{07},Rd_{12},Rd_{13},Rd_{14},$

$Rd_{15},Rd_{16},Rd_{17},Rd_{23},Rd_{24},Rd_{25},Rd_{26},Rd_{27},Rd_{34},Rd_{35},Rd_{36},$

$Rd_{37},Rd_{45},Rd_{46},Rd_{47},Rd_{56},Rd_{57},Rd_{67}\}.$

$CD(i)$ means i-th element of $CD$. For example $CD(8)=Rd_{12}$.
\bigskip

\emph{Lemma 13.2.} \ \gr\gd\  is a Lie algebra.

\bigskip

\emph{Proof.} \ Let we put $c$ as

\ \ \ \ \ \ $c=[CD(i),CD(j)].$

\noindent
We have the following equation with calculations using Maxima.

$c=c(29,30)/d_{01}\times CD(1)+c(29,31)/d_{02}\mathbf{\times }%
CD(2)+c(29,32)/d_{03}\mathbf{\times }CD(3)$

$+c(29,33)/d_{04}\mathbf{\times }CD(4)+c(29,34)/d_{05}\mathbf{\times }%
CD(5)+c(29,35)/d_{06}\mathbf{\times }CD(6)$

$+c(29,36)/d_{07}\mathbf{\times }CD(7)+c(30,31)/d_{12}\mathbf{\times }%
CD(8)+c(30,32)/d_{13}\mathbf{\times }CD(9)$

$+c(30,33)/d_{14}\mathbf{\times }CD(10)+c(30,34)/d_{15}\mathbf{\times }%
CD(11)+c(30,35)/d_{16}\mathbf{\times }CD(12)$

$+c(30,36)/d_{17}\mathbf{\times }CD(13)+c(31,32)/d_{23}\mathbf{\times }%
CD(14)+c(31,33)/d_{24}\mathbf{\times }CD(15)$

$+c(31,34)/d_{25}\mathbf{\times }CD(16)+c(31,35)/d_{26}\mathbf{\times }%
CD(17)+c(31,36)/d_{27}\mathbf{\times }CD(18)$

$+c(32,33)/d_{34}\mathbf{\times }CD(19)+c(32,34)/d_{35}\mathbf{\times }%
CD(20)+c(32,35)/d_{36}\mathbf{\times }CD(21)$

$+c(32,36)/d_{37}\mathbf{\times }CD(22)+c(33,34)/d_{45}\mathbf{\times }%
CD(23)+c(33,35)/d_{46}\mathbf{\times }CD(24)$

$+c(33,36)/d_{47}\mathbf{\times }CD(25)+c(34,35)/d_{56}\mathbf{\times }%
CD(26)+c(34,36)/d_{57}\mathbf{\times }CD(27)$

$+c(35,36)/d_{67}\mathbf{\times }CD(28)$,

\noindent
where $c(k,l)$ means $k$th row and $l$th column of c.\ \ \ \ \emph{Q.E.D.}
\bigskip

\emph{Lemma 13.3.} We have

\ \ \ \ \ \ \ \ \ \ \ \ \ $^{t}Rd_{ij}+Rd_{ij}=0,Rd_{ij}\in $\gr\gd$,0\leq i,j\leq 7.$

\bigskip

\emph{Proof. \ }We have the above equation with calculations using Maxima.

\bigskip

\emph{Definition 13.4.} \ We define the followings for elements $Rd_{ij}\in $\gr\gd.

\ \ \ \ \ $\ \ \ \ \ \ \ \ \ \ \ \ \ Ud_{ij}=Rd_{ij}/d_{ij}\ \ ,(0\leq i\neq j\leq
7). $

Elements of $Ud_{ij}$ are 1 or 0, and the position of element 1 are related to  $d_{ij}$.

Subsequent $U*$  has the same meaning.

\bigskip

\emph{Lemma 13.5.} \ We have the following Lie bracket operations.

\ \ \ \ $[Ud_{ij},Ud_{jk}]=Ud_{ik}\ ,\ (0\leq i,j,k\leq 7),$

\ \ \ \ $[Ud_{ij},Ud_{kl}]=0$ \ \ \ \ \ , $(0\leq i,j,k,l\leq 7,i,j,k,l$ are distinct$).$

\bigskip

\emph{Proof. \ }We have the above Lie bracket operations with calculations
using Maxima.\ \ \ \ \emph{Q.E.D.}

\bigskip

\emph{Lemma 13.6.} \gr\gd$^{\C}$ is simple.

\bigskip

\emph{Proof}. \ Case 1 for $Rd_{01}\in$ \gr\gd$^{\C}$.

\noindent
By \emph{Lemma 13.5 }, we have

$(1)\ \ \ [Ud_{01},Ud_{1i}]=Ud_{0i},(2\leq i\leq 7).$

$(2)\ \ \ [Ud_{01},Ud_{0i}]=-Ud_{1i},(2\leq i\leq 7).$

$(3)\ \ \ [Ud_{02},Ud_{0i}]=-Ud_{2i},(3\leq i\leq 7).$

$(4)\ \ \ [Ud_{03},Ud_{0i}]=-Ud_{3i},(4\leq i\leq 7).$

$(5)\ \ \ [Ud_{04},Ud_{0i}]=-Ud_{4i},(5\leq i\leq 7).$

$(6)\ \ \ [Ud_{05},Ud_{0i}]=-Ud_{5i},(6\leq i\leq 7).$

$(7)\ \ \ [Ud_{06},Ud_{07}]=-Ud_{67}.$

\noindent
Then we have \ \ 

$(8)\ \ \{[Rd_{01},x],[[Rd_{01},x],y] \mid x,y \in$ \gr\gd$^{\C}\}=$\gr\gd$^{\C}.$

Case 2 for $Rd_{0i}\in$\gr\gd$^{\C},(2\leq i\leq 7)$.

\noindent
We have

$(9)\ \ [Ud_{0i},Ud_{1i}]=-Ud_{01}.$

\noindent
Then we have

$(10)$ $\ \{[Rd_{0i},x],[[Rd_{0i},x],y] \mid x,y \in $\gr\gd$^{\C}\}=$\gr\gd$^{\C},(2\leq
i\leq 7).$

Case 3 for $Rd_{ki}\in $\gr\gd$^{\C},(1\leq k<i\leq 7)$.

\noindent
We have

$(11)\ [Ud_{ki},Ud_{0i}]=Ud_{0k}.$

\noindent
Then we have

$(12)$ $\{[Rd_{ki},x],[[Rd_{ki},x],y] \mid x,y \in $\gr\gd$^{\C}\}=$\gr\gd$^{\C},(1\leq
k<i\leq 7).$

By $(8),(10)$ and $(12)$, \gr\gd$^{\C}$ is simple. \ \ \ \ \emph{Q.E.D.}

\bigskip

We consider the Lie algebra

\ \ \ \ \ \ \gD$_{4}=$\gs\go $(8)=\{D\in M(8\times 8,\R) 
\mid D+^{t}D=0\}.$

\bigskip

\emph{Definition 13.7.\ }\gD$_{4}$ is expressed by

{\fontsize{8pt}{10pt} \selectfont%
$\ \left\{ \left( 
\begin{array}{cccccccc}
0 & d_{01} & d_{02} & d_{03} & d_{04} & d_{05} & d_{06} & d_{07} \\ 
-d_{01} & 0 & d_{12} & d_{13} & d_{14} & d_{15} & d_{16} & d_{17} \\ 
-d_{02} & -d_{12} & 0 & d_{23} & d_{24} & d_{25} & d_{26} & d_{27} \\ 
-d_{03} & -d_{13} & -d_{23} & 0 & d_{34} & d_{35} & d_{36} & d_{37} \\ 
-d_{04} & -d_{14} & -d_{24} & -d_{34} & 0 & d_{45} & d_{46} & d_{47} \\ 
-d_{05} & -d_{15} & -d_{25} & -d_{35} & -d_{45} & 0 & d_{56} & d_{57} \\ 
-d_{06} & -d_{16} & -d_{26} & -d_{36} & -d_{46} & -d_{56} & 0 & d_{67} \\ 
-d_{07} & -d_{17} & -d_{27} & -d_{37} & -d_{47} & -d_{57} & -d_{67} & 0%
\end{array}%
\right) \middle| d_{ij}\in \R \right\} .$ }

\bigskip

\emph{Lemma 13.8.} \gr\gd \ is isomorphic to \gD$_{4}$.

\bigskip

\emph{Proof. \ }We consider $\R$-linear mapping

\ \ \ \ \ \ \ \ \ \ \ \ \ \ \ \ \ \ \ \ \ \ \ \ \ $\textrm{fd}:$\gr\gd$\ni Ud_{ij}\rightarrow
D_{ij}\in $\gD$_{4}.$

\noindent
By \emph{Lemma 13.5}, we have

\ \ \ \ $\textrm{fd}([Ud_{ij},U_{jk}])=\textrm{fd}(Ud_{ik})=D_{ik}=[D_{ij},D_{jk}]=[\textrm{fd}(Ud_{ij}),\textrm{fd}(Ud_{jk})],$

\hspace{10em}$(0\leq i<k\leq 7).$

\noindent
Then \gr\gd \ is isomorphic to \gD$_{4}$ under fd . \ \ \ \ \emph{%
Q.E.D.}

\bigskip

\section{The exceptional simle Lie algebra \gr$_{4}^{\C}$ of type $F_{4}$}

\bigskip

\ \ \ \ \emph{Definition 14.1.} \ We define the followings for elements $%
Rm_{ij}\in $\gR\gm$^{\C}.$

\ \ \ \ \ $\ \ \ \ \ \ \ \ \ \ \ \ \ Um_{ij}=Rm_{ij}/m_{ij}\ \ (1\leq i\leq 3,0\leq
j\leq 7).$

\bigskip

\emph{Definition 14.2.} \ We define the following matrices:

\ \ $\ Ca=\left( 
{\fontsize{8pt}{10pt} \selectfont%
\begin{array}{cccccccc}
0 & 1 & 2 & 3 & 4 & 5 & 6 & 7 \\ 
1 & 0 & 3 & 2 & 5 & 4 & 7 & 6 \\ 
2 & 3 & 0 & 1 & 6 & 7 & 4 & 5 \\ 
3 & 2 & 1 & 0 & 7 & 6 & 5 & 4 \\ 
4 & 5 & 6 & 7 & 0 & 1 & 2 & 3 \\ 
5 & 4 & 7 & 6 & 1 & 0 & 3 & 2 \\ 
6 & 7 & 4 & 5 & 2 & 3 & 0 & 1 \\ 
7 & 6 & 5 & 4 & 3 & 2 & 1 & 0%
\end{array} } %
\right) ,$

\ \ \ $Sn=\left( 
{\fontsize{8pt}{10pt} \selectfont%
\begin{array}{rrrrrrrr}
1 & -1 & -1 & -1 & -1 & -1 & -1 & -1 \\ 
-1 & -1 & -1 & 1 & -1 & 1 & -1 & 1 \\ 
-1 & 1 & -1 & -1 & 1 & -1 & -1 & 1 \\ 
-1 & -1 & 1 & -1 & 1 & -1 & 1 & 1 \\ 
-1 & 1 & -1 & 1 & -1 & -1 & 1 & -1 \\ 
-1 & -1 & 1 & 1 & 1 & -1 & -1 & -1 \\ 
-1 & 1 & 1 & -1 & -1 & 1 & -1 & -1 \\ 
-1 & -1 & -1 & -1 & 1 & 1 & 1 & -1%
\end{array} } %
\right) ,$

\ \ \ $Nu=\left( 
{\fontsize{8pt}{10pt} \selectfont%
\begin{array}{rrrrrrrr}
0 & 1 & 2 & 3 & 4 & 5 & 6 & 7 \\ 
1 & 0 & 8 & 9 & 10 & 11 & 12 & 13 \\ 
2 & 8 & 0 & 14 & 15 & 16 & 17 & 18 \\ 
3 & 9 & 14 & 0 & 19 & 20 & 21 & 22 \\ 
4 & 10 & 15 & 19 & 0 & 23 & 24 & 25 \\ 
5 & 11 & 16 & 20 & 23 & 0 & 26 & 27 \\ 
6 & 12 & 17 & 21 & 24 & 26 & 0 & 28 \\ 
7 & 13 & 18 & 22 & 25 & 27 & 28 & 0%
\end{array} } %
\right) $.

\bigskip

\emph{Lemma 14.3.} \ We have the following Lie bracket operations.

\noindent
$[Um_{1i},Um_{2j}]=-\frac{1}{2}Sn(i+1,j+1)Um_{3k}\ ,(k=Ca(i+1,j+1),\ 0\leq i,j\leq 7)$,

\noindent
$[Um_{2i},Um_{3j}]=-\frac{1}{2}Sn(i+1,j+1)Um_{1k}\ ,(k=Ca(i+1,j+1),\ 0\leq i,j\leq 7)$,

\noindent
$[Um_{3i},Um_{1j}]=-\frac{1}{2}Sn(i+1,j+1)Um_{2k}\ ,(k=Ca(i+1,j+1),\ 0\leq i,j\leq 7)$,

\noindent
$[Um_{1i},Um_{1j}]=-Ud_{ij}\ ,(0\leq i<j\leq 7)$,

\noindent
$[Um_{2i},Um_{2j}]=-\sum\limits_{0\leq n<l\leq 7}Mv(ki,kj)Ud_{nl}$,

\ \ \ \ \ \ \ \ \ \ $(ki=Nu(i+1,j+1),kj=Nu(n+1,l+1),0\leq i<j\leq 7)$,

\noindent
$[Um_{3i},Um_{3j}]=-\sum\limits_{0\leq n<l\leq 7}Mv^{2}(ki,kj)Ud_{nl}$,

\ \ \ \ \ \ \ \ \ \ $(ki=Nu(i+1,j+1),kj=Nu(n+1,l+1),0\leq i<j\leq 7)$.

\bigskip

\emph{Proof. \ }We have the above Lie bracket operations with
calculations using Maxima.\ \ \ \ \emph{Q.E.D.}

\bigskip

\emph{Lemma 14.4. }\ We have the following Lie bracket operations.
\begin{align*}
[Ud_{ij},Um_{1k}] & =-Um_{1j} (\text{in case of }k=i),\\
&=Um_{1i} (\text{in case of }k=j),\\
&=0(\text{in case of }k\neq i,j),(0\leq i<j\leq 7),\\
[Ud_{ij},Um_{2k}] &=-\sum\limits_{0\leq n<l\leq%
7}Mv^{2}(ki,kj)Um_{2l} (\text{where }k=n)\\
&\ \ \ +\sum\limits_{0\leq n<l\leq 7}Mv^{2}(ki,kj)Um_{2n} (\text{where }k=l),\\
&(ki=Nu(i+1,j+1),kj=Nu(n+1,l+1),0\leq i<j\leq 7),\\
[Ud_{ij},Um_{3k}] &=-\sum\limits_{0\leq n<l\leq%
7}Mv(ki,kj)Um_{3l} (\text{where }n=k)\\
&\ \ \ +\sum\limits_{0\leq n<l\leq 7}Mv(ki,kj)Um_{3n} (\text{where }l=k),\\
&(ki=Nu(i+1,j+1),kj=Nu(n+1,l+1),0\leq i<j\leq 7).
\end{align*}

\emph{Proof. \ }We have the above Lie bracket operations with
calculations using Maxima.\ \ \ \ \emph{Q.E.D.}

\bigskip

\emph{Lemma 14.5.} \ \gr$_{4}^{\C}$ is simple.

\bigskip

\emph{Proof. }\ By \emph{Lemma 14.3}, we have

$\ \ \ \ \{[Rm_{1i},Rm_{2j}] \mid 0\leq j\leq 7\}=\{Rm_{3k} 
\mid 0\leq k\leq 7\},(0\leq i\leq 7),$

$\ \ \ \ \{[Rm_{2i},Rm_{3j}] \mid 0\leq j\leq 7\}=\{Rm_{1k} 
\mid 0\leq k\leq 7\},(0\leq i\leq 7),$

$\ \ \ \ \{[Rm_{3i},Rm_{1j}] \mid 0\leq j\leq 7\}=\{Rm_{2k} 
\mid 0\leq k\leq 7\},(0\leq i\leq 7).$

\noindent
By \emph{Lemma 14.3}, we have

$\ \ \ \ \{[Rm_{1i},Rm_{1j}] \mid 0\leq j\leq
7\},\{[Rm_{2i},Rm_{2j}]|0\leq j\leq 7\},$

$\ \ \ \ \{[Rm_{3i},Rm_{3j}]\mid 0\leq j\leq 7\}\subset $\gr\gd$^{\C}$.

\noindent
Then we have by \emph{Lemma 13.6}

$\ \ \ \ \{[Rm_{1i},x],[[Rm_{1i},x],y] \mid  x,y \in \ $\gr\gd$^{\C} \bigoplus $\gR\gm$^{\C}\ \}=$ 
\gr$_{4}^{\C},(0\leq i\leq 7),$

$\ \ \ \ \{[Rm_{2i},x],[[Rm_{2i},x],y] \mid x,y\in \ $\gr\gd$^{\C} \bigoplus $\gR\gm$^{\C}\ \}=$ 
\gr$_{4}^{\C},(0\leq i\leq 7),$

$\ \ \ \ \{[Rm_{3i},x],[[Rm_{3i},x],y] \mid x,y\in \ $\gr\gd$^{\C} \bigoplus $\gR\gm$^{\C}\ \}=$ 
\gr$_{4}^{\C},(0\leq i\leq 7).$

By \emph{Lemma 14.4}, we have

$\ \ \ \ [Rd_{ij},Rm_{1i}]\in \{Rm_{1k} \mid 0\leq k\leq
7\},[Rd_{ij},Rm_{2i}]\in \{Rm_{2k} \mid 0\leq k\leq 7\},$

$\ \ \ \ [Rd_{ij},Rm_{3i}]\in \{Rm_{3k} \mid 0\leq k\leq 7\}.$

\noindent
Then \ we have

$\ \ \ \ \{[Rd_{ij},x],[[Rd_{ij},x],y] \mid x,y\in \ $\gr\gd$^{\C} \bigoplus $\gR\gm$^{\C}\ \}=$ 
\gr$_{4}^{\C},(0\leq i<j\leq 7).$

So \gr$_{4}^{\C}=\{$\gr\gd$^{\C}, $ \gR\gm$^{\C}\}$ is simple. \ \ \emph{Q.E.D.}

\bigskip

\emph{Lemma 14.6.} \ The Killing form $B_{4}$ of the Lie algebra 
\gr$_{4}^{\C}$ is given by

\ \ \ \ \ \ \ \ \ \ \ \ \ \ \ \ \ $B_{4}(R_{1},R_{2})=\frac{3}{10}%
tr(R_{1}R_{2}),R_{1},R_{2}\in $\gr$_{4}^{\C}.$

\bigskip

\emph{\ Proof.} \ Since \gr$_{4}^{\C}$\textbf{\ }is simple,there
exist $\kappa \in \C$ such that

\ \ \ \ \ \ \ \ \ \ \ \ \ \ \ \ $B_{4}(R_{1},R_{2})=\kappa tr(R_{1}R_{2}).$

\noindent
To determin this $\kappa $,we put $R=R_{1}=R_{2}=Um_{10}\in $\gR\gm$^{\C}.$

$(adR)^{2}$ is calculated as follows by \emph{Lemma 14.3} and \emph{Lemma 14.4} .

\ \ \ \ \ $\ [Um_{10},[Um_{10},Ud_{0i}]]=[Um_{10},Um_{1i}]=-Ud_{0i},1\leq i\leq 7,$

\ \ \ \ \ $\ [Um_{10},[Um_{10},Um_{1i}]]=[Um_{10},-Ud_{0i}]=-Um_{1i},1\leq i\leq 7,$

\ \ \ \ \ \ $[Um_{10},[Um_{10},Um_{20}]]=[Um_{10},-\frac{1}{2}Um_{30}]=-\frac{1}{4}Um_{20},$

\ \ \ \ \ \ $[Um_{10},[Um_{10},Um_{2i}]]=[Um_{10},\frac{1}{2}Um_{3i}]=-\frac{1}{4}%
Um_{2i},1\leq i\leq 7,$

\ \ \ \ \ \ $[Um_{10},[Um_{10},Um_{30}]]=[Um_{10},-\frac{1}{2}Um_{20}]=-\frac{1}{4}Um_{30},$

\ \ \ \ \ \ $[Um_{10},[Um_{10},Um_{3i}]]=[Um_{10},\frac{1}{2}Um_{2i}]=-\frac{1}{4}%
Um_{3i},1\leq i\leq 7,$

\ \ \ \ \ the others $=0$.

\noindent
Hence we have

\ \ \ \ $B_{4}(Um_{10},Um_{10})=tr((ad(Um_{10}))^{2})$

\ \ \ \ \ \ \ $=(-1)\times 7+(-1)\times 7+(-%
\frac{1}{4})\times 8+(-\frac{1}{4})\times 8=-18.$

On the other hand ,by calculate with Maxima we have

\ \ \ \ $tr((Um_{10})(Um_{10}))=-60.$

\noindent
Hence we have $\kappa =\frac{-18}{-60}=\frac{3}{10}.$ \ \ \ \ \emph{Q.E.D.}

\bigskip

\emph{Lemma 14.7.} The roots of \gr\gd$^{\C}$ relative to some Cartan
subalgebra of \gr\gd$^{\C}$ are given by

$\ \ \ \ \ \ \ \ \ \ \ \ \ \ \ \ \ \ \pm (\lambda _{k}-\lambda _{l}),\pm
(\lambda _{k}+\lambda _{l}),0\leq k<l\leq 3.$

\bigskip

\emph{Proof. \ }Since\emph{\ }the Lie algebra \gr\gd$^{\C}$ is isomorphic to 
\gD$_{4}^{\C},$\ if we let put $H_{k}=-iUd_{k4+k}$ for
$k=0,1,2,3$, then

\ \ \ \ \ \ \ \ \ \ \ \ \ \ \ \ \ \ \ \ \gh$=\{H=\sum%
\limits_{k=0}^{3}\lambda _{k}H_{k} \mid H_{k}=-iUd_{k4+k},%
\lambda _{k}\in \C\}.$

\noindent
is a Cartan subalgebra of \gr\gd$^{\C}.$

\noindent
By \emph{Lemma 13.5}, we have

\ \ \ \ \ \ \ \ \ \ [$H,(Ud_{kl}+Ud_{4+k4+l})-i(Ud_{k4+l}+Ud_{l4+k})]$

$\ \ \ \ \ \ \ \ \ \ =(\lambda _{k}-\lambda
_{l})((Ud_{kl}+Ud_{4+k4+l})-i(Ud_{k4+l}+Ud_{l4+k})),$

\ \ \ \ \ \ \ \ \ \ [$H,i(Ud_{kl}+Ud_{4+k4+l})-(Ud_{k4+l}+Ud_{l4+k})]$

$\ \ \ \ \ \ \ \ \ \ =(-\lambda _{k}+\lambda
_{l})(i(Ud_{kl}+Ud_{4+k4+l})-(Ud_{k4+l}+Ud_{l4+k})),$

\ \ \ \ \ \ \ \ \ \ [$H,(Ud_{kl}-Ud_{4+k4+l})+i(Ud_{k4+l}-Ud_{l4+k})]$

$\ \ \ \ \ \ \ \ \ \ =(\lambda _{k}+\lambda
_{l})((Ud_{kl}-Ud_{4+k4+l})+i(Ud_{k4+l}-Ud_{l4+k})),$

$\ \ \ \ \ \ \ \ \ \ [H,i(Ud_{kl}-Ud_{4+k4+l})+(Ud_{k4+l}-Ud_{l4+k})]$

$\ \ \ \ \ \ \ \ \ \ =(-\lambda _{k}-\lambda
_{l})(i(Ud_{kl}-Ud_{4+k4+l})+(Ud_{k4+l}-Ud_{l4+k})).$

\noindent
We do the same with calculations using Maxima.

\noindent
Hence $\pm (\lambda _{k}-\lambda _{l}),\pm (\lambda _{k}+\lambda _{l}),0\leq
k<l\leq 3$ are roots of \gr\gd$^{\C}$. \ \ \ \ \emph{Q.E.D.}

\bigskip

\emph{Lemma 14.8.} The rank of the Lie algebra \gr$_{4}^{\C}$
is 4. The roots of \gr$_{4}^{\C}$ relative to some Cartan
subalgebra of \gr$_{4}^{\C}$ are given by

$\ \ \ \ \ \ \ \ \ \ \ \ \ \ \ \ \ \ \pm (\lambda _{k}-\lambda _{l}),\pm
(\lambda _{k}+\lambda _{l}),0\leq k<l\leq 3,$

\ \ \ \ \ \ \ \ \ \ \ \ \ \ \ \ \ \ $\pm \lambda _{0},\pm \lambda _{1},\pm
\lambda _{2},\pm \lambda _{3},$

$\ \ \ \ \ \ \ \ \ \ \ \ \ \ \ \ \ \ \pm \frac{1}{2}(-\lambda _{0}-\lambda
_{1}+\lambda _{2}-\lambda _{3}),\pm \frac{1}{2}(\lambda _{0}+\lambda
_{1}+\lambda _{2}-\lambda _{3}),$

\ \ \ \ \ \ \ \ \ \ \ \ \ \ \ \ \ \ $\pm \frac{1}{2}(-\lambda _{0}+\lambda
_{1}+\lambda _{2}+\lambda _{3}),\pm \frac{1}{2}(\lambda _{0}-\lambda
_{1}+\lambda _{2}+\lambda _{3}),$

\ \ \ \ \ \ \ \ \ \ \ \ \ \ \ \ \ \ $\pm \frac{1}{2}(-\lambda _{0}+\lambda
_{1}-\lambda _{2}+\lambda _{3}),\pm \frac{1}{2}(-\lambda _{0}+\lambda
_{1}+\lambda _{2}-\lambda _{3}),$

\ \ \ \ \ \ \ \ \ \ \ \ \ \ \ \ \ \ $\pm \frac{1}{2}(\lambda _{0}+\lambda
_{1}+\lambda _{2}+\lambda _{3}),\pm \frac{1}{2}(-\lambda _{0}-\lambda
_{1}+\lambda _{2}+\lambda _{3}).$

\bigskip

\emph{Proof. \ }The Lie algebra \gr\gd$^{\C}$ is contained in \gr$%
_{4}^{\C}$.

Let \gh$=\{H=\sum\limits_{k=0}^{3}\lambda _{k}H_{k} \mid \lambda _{k}\in \C\}%
\subset $\gr\gd$^{\C}\subset $ \gr$_{4}^{\C}.$ Since 
\gh\ is a Cartan subalgebra of
\gr\gd$^{\C}$ (it will be also a Cartan subalgebra of \gr$_{4}^{\C}$ ),
the roots of \gr\gd$^{\C}$ :

$\ \ \ \ \ \ \ \ \ \ \ \ \ \ \ \ \ \ \pm (\lambda _{k}-\lambda _{l}),\pm
(\lambda _{k}+\lambda _{l}),0\leq k<l\leq 3$

\noindent
are also roots of \gr$_{4}^{\C}.$ 

\noindent
Furthermore, by \emph{Lemma 14.4} we have

$[H,Um_{1k}+iUm_{14+k}]=-i[\lambda _{k}Ud_{k4+k},Um_{1k}]+[\lambda
_{k}Ud_{k4+k},Um_{14+k}]$

\ \ \ \ \ \ \ \ \ \ \ \ \ \ \ \ \ \ \ \ \ \ \ \ \ \ \ \ $=\lambda
_{k}(iUm_{14+k})+\lambda _{k}(Um_{1k})=\lambda _{k}(Um_{1k}+iUm_{14+k}),$

\noindent
we see that $\lambda _{k}$ are roots of \gr$_{4}^{\C}$ and $%
Um_{1k}+iUm_{14+k}$ are associated root vectors for
$0\leq k\leq 3.$ Similarly,$-\lambda _{k}$ for $0\leq k\leq 3$ are roots of 
\gr$_{4}^{\C}$ and $Um_{1k}-iUm_{14+k}$ are associated root vectors. 

Next, by \emph{Lemma 14.4 }we have

$\ \ \ \ \ \ \ \ \ [H,Um_{20}+iUm_{24}]=\frac{1}{2}(-\lambda _{0}-\lambda
_{1}+\lambda _{2}-\lambda _{3})(Um_{20}+iUm_{24}),$

\ \ \ \ \ \ \ \ \ $[H,Um_{21}+iUm_{25}]=\frac{1}{2}(\lambda _{0}+\lambda
_{1}+\lambda _{2}-\lambda _{3})(Um_{21}+iUm_{25}),$

\ \ \ \ \ \ \ \ \ $[H,Um_{22}+iUm_{26}]=\frac{1}{2}(-\lambda _{0}+\lambda
_{1}+\lambda _{2}+\lambda _{3})(Um_{22}+iUm_{26}),$

\ \ \ \ \ \ \ \ \ $[H,Um_{23}+iUm_{27}]=\frac{1}{2}(\lambda _{0}-\lambda
_{1}+\lambda _{2}+\lambda _{3})(Um_{23}+iUm_{27}),$

$\ \ \ \ \ \ \ \ \ [H,Um_{20}-iUm_{24}]=-\frac{1}{2}(-\lambda _{0}-\lambda
_{1}+\lambda _{2}-\lambda _{3})(Um_{20}-iUm_{24}),$

\ \ \ \ \ \ \ \ \ $[H,Um_{21}-iUm_{25}]=-\frac{1}{2}(\lambda _{0}+\lambda
_{1}+\lambda _{2}-\lambda _{3})(Um_{21}-iUm_{25}),$

\ \ \ \ \ \ \ \ \ $[H,Um_{22}-iUm_{26}]=-\frac{1}{2}(-\lambda _{0}+\lambda
_{1}+\lambda _{2}+\lambda _{3})(Um_{22}-iUm_{26}),$

\ \ \ \ \ \ \ \ \ $[H,Um_{23}-iUm_{27}]=-\frac{1}{2}(\lambda _{0}-\lambda
_{1}+\lambda _{2}+\lambda _{3})(Um_{23}-iUm_{27}).$

\noindent
So we see $\frac{1}{2}(-\lambda _{0}-\lambda _{1}+\lambda _{2}-\lambda _{3}),%
\frac{1}{2}(\lambda _{0}+\lambda _{1}+\lambda _{2}-\lambda _{3}),\frac{1}{2}%
(-\lambda _{0}+\lambda _{1}+\lambda _{2}+\lambda _{3}),$

\noindent
$\frac{1}{2}(\lambda _{0}-\lambda _{1}+\lambda _{2}+\lambda _{3})$ are roots
of \gr$_{4}^{\C}$ and that $(Um_{2k}+iUm_{24+k})$ are associated root
vectors
for $0\leq k\leq 3.$ The roots above with negative sign are also roots of 
\gr$_{4}^{\C}$ and
$(Um_{2k}-iUm_{24+k})$ are associated root vectors. 

Similarly,by \emph{Lemma 14.4} we have

$\ \ \ \ \ \ \ \ \ [H,Um_{30}+iUm_{34}]=\frac{1}{2}(-\lambda _{0}+\lambda
_{1}-\lambda _{2}+\lambda _{3})(Um_{30}+iUm_{34}),$

\ \ \ \ \ \ \ \ \ $[H,Um_{31}+iUm_{35}]=\frac{1}{2}(-\lambda _{0}+\lambda
_{1}+\lambda _{2}-\lambda _{3})(Um_{31}+iUm_{35}),$

\ \ \ \ \ \ \ \ \ $[H,Um_{32}+iUm_{36}]=\frac{1}{2}(\lambda _{0}+\lambda
_{1}+\lambda _{2}+\lambda _{3})(Um_{32}+iUm_{36}),$

\ \ \ \ \ \ \ \ \ $[H,Um_{33}+iUm_{37}]=\frac{1}{2}(-\lambda _{0}-\lambda
_{1}+\lambda _{2}+\lambda _{3})(Um_{33}+iUm_{37}),$

$\ \ \ \ \ \ \ \ \ [H,Um_{30}-iUm_{34}]=-\frac{1}{2}(-\lambda _{0}+\lambda
_{1}-\lambda _{2}+\lambda _{3})(Um_{30}-iUm_{34}),$

\ \ \ \ \ \ \ \ \ $[H,Um_{31}-iUm_{35}]=-\frac{1}{2}(-\lambda _{0}+\lambda
_{1}+\lambda _{2}-\lambda _{3})(Um_{31}-iUm_{35}),$

\ \ \ \ \ \ \ \ \ $[H,Um_{32}-iUm_{36}]=-\frac{1}{2}(\lambda _{0}+\lambda
_{1}+\lambda _{2}+\lambda _{3})(Um_{32}-iUm_{36}),$

\ \ \ \ \ \ \ \ \ $[H,Um_{33}-iUm_{37}]=-\frac{1}{2}(-\lambda _{0}-\lambda
_{1}+\lambda _{2}+\lambda _{3})(Um_{33}-iUm_{37}).$

\noindent
So we see $\frac{1}{2}(-\lambda _{0}+\lambda _{1}-\lambda _{2}+\lambda _{3}),%
\frac{1}{2}(-\lambda _{0}+\lambda _{1}+\lambda _{2}-\lambda _{3}),\frac{1}{2}%
(\lambda _{0}+\lambda _{1}+\lambda _{2}+\lambda _{3}),$

\noindent
$\frac{1}{2}(-\lambda _{0}-\lambda _{1}+\lambda _{2}+\lambda _{3})$ are
roots of \gr$_{4}^{\C}$ and that $(Um_{3k}+iUm_{34+k})$ are associated
root vectors for $0\leq k\leq 3.$ The roots above with negative sign are also roots 
of  \gr$_{4}^{\C}$ and $(Um_{3k}-iUm_{34+k})$ are associated root vectors.

By direct calculations using Maxima, we also have the above. \ \ \ \ \emph{Q.E.D.}

\bigskip

\emph{Theorem 14.9.} \ In the root system of \emph{Lemma 14.7},

\ \ \ \ \ \ \ $\alpha _{1}=\lambda _{0}-\lambda _{1},\alpha _{2}=\lambda
_{1}-\lambda _{2},\alpha _{3}=\lambda _{2},\alpha _{4}=\frac{1}{2}(-\lambda
_{0}-\lambda _{1}-\lambda _{2}+\lambda _{3})$

\noindent
is a fundamental root sysyten of the Lie algebra \gr$_{4}^{\C}$ and

$\ \ \ \ \ \ \ \mu =2\alpha _{1}+3\alpha _{2}+4\alpha _{3}+2\alpha _{4}$

\noindent
is the highest root. The Dynkin diagram and the extended Dynkin diagram 
of  \gr$_{4}^{\C}$ are respectively given by

\begin{figure}[H]
\centering
\includegraphics[width=10cm, height=2.0cm]{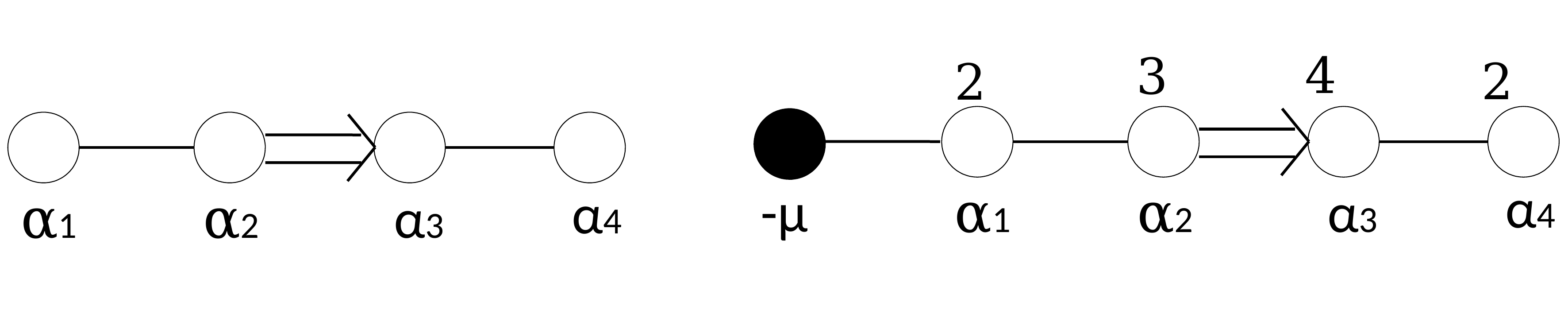}
\end{figure}

\emph{Proof.} In the following, the notation $n_{1}$ $n_{2}$ $n_{3}$ 
$n_{4}$ denotes the root $n_{1}\alpha _{1}+n_{2}\alpha _{2}+n_{3}\alpha
_{3}+n_{4}\alpha _{4}$ $.$

\noindent
For example, $\alpha =%
\begin{array}{@{}llll@{}}
n_{1} & n_{2} & n_{3} & n_{4}%
\end{array}%
$ means \ $\alpha =n_{1}\alpha _{1}+n_{2}\alpha _{2}+n_{3}\alpha
_{3}+n_{4}\alpha _{4}$ .

Now,all positive roots of \gr$_{4}^{\C}$ are represented by

\begin{align*}
\lambda _{0} & =  %
\begin{array}{llll}
1 & 1 & 1 & 0%
\end{array},\\
\lambda _{1} &= %
\begin{array}{llll}
0 & 1 & 1 & 0%
\end{array},\\
\lambda _{2} &= %
\begin{array}{llll}
0 & 0 & 1 & 0%
\end{array},\\
\lambda _{3}& = %
\begin{array}{llll}
1 & 2 & 3 & 2%
\end{array},\\
\lambda _{0}-\lambda _{1} &= %
\begin{array}{llll}
1 & 0 & 0 & 0%
\end{array},\\
\lambda _{0}-\lambda _{2} &= %
\begin{array}{llll}
1 & 1 & 0 & 0%
\end{array},\\
-\lambda_{0}+\lambda _{3} &= %
\begin{array}{llll}
0 & 1 & 2 & 2%
\end{array},\\
\lambda _{1}-\lambda _{2} &= %
\begin{array}{llll}
0 & 1 & 0 & 0%
\end{array},\\
 -\lambda_{1}+\lambda _{3} &= %
\begin{array}{llll}
1 & 1 & 2 & 2%
\end{array},\\
 -\lambda_{2}+\lambda _{3} &= %
\begin{array}{llll}
1 & 2 & 2 & 2%
\end{array}.\\
\lambda _{0}+\lambda _{1} &= %
\begin{array}{llll}
1 & 2 & 2 & 0%
\end{array},\\
\lambda _{0}+\lambda _{2} &= %
\begin{array}{llll}
1 & 1 & 2 & 0%
\end{array},\\
\lambda _{0}+\lambda _{3} &= %
\begin{array}{llll}
2 & 3 & 4 & 2%
\end{array},\\
\lambda _{1}+\lambda _{2} &= %
\begin{array}{llll}
0 & 1 & 2 & 0%
\end{array},\\
\lambda _{1}+\lambda _{3} &= %
\begin{array}{llll}
1 & 3 & 4 & 2%
\end{array},\\
\lambda _{2}+\lambda _{3} &= %
\begin{array}{llll}
1 & 2 & 4 & 2%
\end{array},\\
-\frac{1}{2}(-\lambda _{0}-\lambda _{1}+\lambda%
_{2}-\lambda _{3}) &= %
\begin{array}{llll}
1 & 2 & 2 & 1%
\end{array},\\
-\frac{1}{2}(\lambda _{0}+\lambda _{1}+\lambda%
_{2}-\lambda _{3}) &= %
\begin{array}{llll}
0 & 0 & 0 & 1%
\end{array},\\
\frac{1}{2}(-\lambda _{0}+\lambda%
_{1}+\lambda _{2}+\lambda _{3}) &= %
\begin{array}{llll}
0 & 1 & 2 & 1%
\end{array},\\
\end{align*}
\begin{align*}
\frac{1}{2}(\lambda _{0}-\lambda%
_{1}+\lambda _{2}+\lambda _{3}) &= %
\begin{array}{llll}
1 & 1 & 2 & 1%
\end{array},\\
\frac{1}{2}(-\lambda _{0}+\lambda%
_{1}-\lambda _{2}+\lambda _{3}) &= %
\begin{array}{llll}
0 & 1 & 1 & 1%
\end{array},\\
-\frac{1}{2}(-\lambda _{0}+\lambda _{1}+\lambda%
_{2}-\lambda _{3}) &= %
\begin{array}{llll}
1 & 1 & 1 & 1%
\end{array},\\
\frac{1}{2}(\lambda _{0}+\lambda%
_{1}+\lambda _{2}+\lambda _{3}) &= %
\begin{array}{llll}
1 & 2 & 3 & 1%
\end{array},\\
\frac{1}{2}(-\lambda _{0}-\lambda%
_{1}+\lambda _{2}+\lambda _{3}) &= %
\begin{array}{llll}
0 & 0 & 1 & 1%
\end{array}.
\end{align*}

\noindent
Hence $\Pi =\{\alpha _{1},\alpha _{2},\alpha _{3},\alpha _{4}\}$ is a
fundamental root system of \gr$_{4}^{\C}$ and

\noindent
$\mu =2\alpha _{1}+3\alpha _{2}+4\alpha _{3}+2\alpha _{4}(=\lambda
_{0}+\lambda _{3})$ is the highest root. The real part of \gh$_{R}$
of \gh\  is

\ \ \ \ \ \ \ \ \ \ \ \ \ \ \ \ \ \ \ \ \gh$_{R}=\{H=\sum%
\limits_{k=0}^{3}\lambda _{k}H_{k} \mid \lambda _{k}\in \R\}.$

The Killing form $B_{4}$ of \gr$_{4}^{\C}$ is $B_{4}(R_{1},R_{2})=\frac{3}{10%
}tr(R_{1}R_{2})$ (\emph{Lemma 14.6} ), so that

\ \ \ \ \ \ \ \ \ \ \ $B_{4}(H,H^{\prime} )=18\sum\limits_{k=0}^{3}\lambda
_{k}\lambda _{k}^{\prime },H=\sum\limits_{k=0}^{3}\lambda
_{k}H_{k},H^{\prime }=\sum\limits_{k=0}^{3}\lambda _{k}^{\prime
}H_{k}^{\prime }\in $\gh$_{R}.$

\noindent
Indeed,by calculate with Maxima we have

$\ \ \ \ \ \ \ \ \ \ \ \ \ \ B_{4}(H,H^{\prime} )=\frac{3}{10}tr(HH^{\prime} )=%
\frac{3}{10}(60\sum\limits_{k=0}^{3}\lambda _{k}\lambda _{k}^{\prime
})=18\sum\limits_{k=0}^{3}\lambda _{k}\lambda _{k}^{\prime }.$

Now,the canonical elements H$_{\alpha _{i}}\in $\gh$_{R}$ corresponding
to $\alpha _{i}(B_{4}(H_{\alpha _{i}},H)=\alpha _{i}(H),H\in $\gh$_{R})$
are determind by

$\ \ \ \ \ \ \ \ \ \ \ \ \ \ \ \ \ \ H_{\alpha 1}=\frac{1}{18}%
(H_{0}-H_{1}),H_{\alpha 2}=\frac{1}{18}(H_{1}-H_{2}),H_{\alpha 3}=\frac{1}{18%
}H_{2},$

$\ \ \ \ \ \ \ \ \ \ \ \ \ \ \ \ \ \ H_{\alpha 4}=\frac{1}{36}%
(-H_{0}-H_{1}-H_{2}+H_{3}).$

\noindent
Hence we have

$\ \ \ \ \ \ \ \ \ \ \ \ \ \ \ \ (\alpha _{1},\alpha _{1})=B_{4}(H_{\alpha
1},H_{\alpha 1})=18\frac{1}{18}\frac{1}{18}(1+1)=\frac{1}{9}$

\noindent
and the other inner products are similarly calculated. Hence,the inner
products induced by the Killing form $B_{4}$ between $\alpha _{1},\alpha _{2},\alpha _{3},\alpha
_{4}$ and $-\mu $ 
are given by

\ \ \ \ \ \ \ \ \ \ \ \ \ \ \ \ $(\alpha _{1},\alpha _{1})=(\alpha
_{2},\alpha _{2})=\frac{1}{9},(\alpha _{3},\alpha _{3})=(\alpha _{4},\alpha
_{4})=\frac{1}{18},$

\ \ \ \ \ \ \ \ \ \ \ \ \ \ \ $\ (\alpha _{1},\alpha _{2})=-\frac{1}{18}%
,(\alpha _{2},\alpha _{3})=-\frac{1}{18},(\alpha _{3},\alpha _{4})=-\frac{1}{%
36},$

\ \ \ \ \ \ \ \ \ \ \ \ \ \ \ \ $(\alpha _{1},\alpha _{3})=(\alpha
_{1},\alpha _{4})=(\alpha _{2},\alpha _{4})=0,$

\ \ \ \ \ \ \ \ \ \ \ \ \ \ \ \ $(-\mu ,-\mu )=\frac{1}{9},(-\mu ,\alpha
_{1})=-\frac{1}{18},(-\mu ,\alpha _{i})=0,i=2,3,4,$

\noindent
using them,we can draw the Dynkin diagram and the extended Dynkin 
diagram of  \gr$_{4}^{\C}.$ \ \ \ \ \emph{Q.E.D.}

\bigskip

\emph{Corollary 14.10.} For a $248\times 248$ matrix $X$,let $X|_{4}$ be the matrix in which the $52\times 52$
elements in the upper left corner are clipped from  $X$. Furthermore, let \gr\gd$|_{4}=\{Rd|_{4} \mid Rd\in $\gr\gd$\}$ and 
\gR\gm$|_{4}=\{Rm|_{4} \mid Rm\in $\gR\gm$\}$.
\emph{Theorem 14.9} holds for \gr$_{4}^{\C}|_{4}=$\gr\gd$^{\C}|_{4}\oplus $\gR\gm$^{\C}|_{4}$ as well. However, the
Killing form $B_{4}(R_{1},R_{2})=tr(R_{1}R_{2})\ \ (R_{1},R_{2} \in$ \gr$_{4}^{\C}|_{4}).$

\bigskip

\emph{Proof} We have the above with calculations using Maxima.\ \ \ \ \emph{Q.E.D.}

\bigskip 

\emph{Definition 14.11.} \ We define the real part of \gr\gd$^{\C}$ \ by \ge\gd,

\ge\gd$=\{Rd_{ij}\in $\gr\gd$^{\C} \mid d_{ij}\in \R,0\leq i<j\leq 7\}.$

\noindent
And we define the real part of \gR\gm$^{\C}$ \ by \ge\gm,

\ge\gm$=\{Rm_{ij}\in $\gR\gm$^{\C}  \mid m_{ij}\in \R,1\leq i\leq 3,0\leq j\leq 7\}.$

\noindent
Also we define as follows.

\ge\gd$_{4}=\{X|_{4} \mid X\in $\ge\gd$\},$

\ge\gm$_{4}=\{X|_{4} \mid X\in $\ge\gm$\}.$

\bigskip

\emph{Theorem 14.12.} Let we put

\gf$_{4}=$\ge\gd$_{4}\oplus $\ge\gm$_{4}$ . \ 

\noindent
Then \gf$_{4}$ is a compact exceptional simple Lie algebra of type $F_{4}.$

\bigskip

\emph{Proof.} \ By calculations using Maxima,we have as follows.

\noindent
For $X\in $\ge\gd$_{4},Y\in $\ge\gm$_{4},$

$B_{4}(X,X)$

$=-18(d_{67}^{2}+d_{57}^{2}+d_{56}^{2}+d_{47}^{2}+d_{46}%
^{2}+d_{45}^{2}+d_{37}^{2}$

\ \ \ \ \ \ $\ +d_{36}^{2}+d_{35}^{2}+d_{34}^{2}+d_{27}^{2}+d_{26}^{2}+d_{25}^{2}+d_{24}^{2}$

\ \ \ \ \ \ $\ +d_{23}^{2}+d_{17}^{2}+d_{16}^{2}+d_{15}^{2}+d_{14}^{2}+d_{13}^{2}+d_{12}^{2}$

\ \ \ \ \ \ $\ +d_{07}^{2}+d_{06}^{2}+d_{05}^{2}+d_{04}^{2}+d_{03}^{2}+d_{02}^{2}+d_{01}^{2}),$

$B_{4}(Y,Y)$

$=-18(m_{37}^{2}+m_{36}^{2}+m_{35}^{2}+m_{34}^{2}+m_{33}%
^{2}+m_{32}^{2}+m_{31}^{2}+m_{30}^{2}$

\ \ \ \ \ \ $\ +m_{27}^{2}+m_{26}^{2}+m_{25}^{2}+m_{24}^{2}+m_{23}^{2}+m_{22}^{2}+m_{21}^{2}+m_{20}^{2}$

\ \ \ \ \ \ $\ +m_{17}^{2}+m_{16}^{2}+m_{15}^{2}+m_{14}^{2}+m_{13}^{2}+m_{12}^{2}+m_{11}^{2}+m_{10}^{2}),$

$B_{4}(X,Y)=0.$

\noindent
Therefor we have $B_{4}(X,X)<0$,\ for $^{\forall }X\neq 0,X\in$ \gf$_{4}.$
Then \gf$_{4}$ is compact. \ \ \ \ \emph{Q.E.D.}

\bigskip

Next, let digitalize the exceptional simple Lie algebra of type $F_{4}$ as
matrices in $M(27 \times 27,\R)$.

\bigskip

\emph{Definition 14.13.}  We consider following elements in $M(27 \times 27,\R)$,

\noindent
 (See \emph{Theorem 7.2}).

For \ge\gd \ $\ni Rd_{ij} , (d_{ij} \in \R$),

{\fontsize{7pt}{8pt} \selectfont%
 \ \ \ \ \ \ \ \ \ \ \ \ \ \ \ \ \ \ \ \ \ \ \ $%

\right)  $
}

\noindent
See \emph{Lemma 6.17} for MIr(m), \emph{Lemma 6.19} for MC(m). 

We denote elemtns of  $Rd_{ij}^{`}, Rm_{1}^{`}, Rm_{2}^{`}$, and $Rm_{3}^{`}$. 

$Rd_{01}^{'}=d_{01}$
{\fontsize{6pt}{8pt} \selectfont%
$\left( 
%
%
\right)$ }.

We define the following Lie algebra:

\gf$^{`}_{4}=$\{ Lie algebra over real numbers generated by $Rd_{ij}^{`}(0 \leq i < j \leq 7)$,
$Rm_{ij}^{`}(1 \leq i \leq 3, 0 \leq j \leq 7) \}$

\bigskip

\emph{Lemma 14.14.}  The following correspondences are homomorphism.

\ \ \ \ \ \ \ \ \ \ \ \ \ \ \ \ \ \ \ge\gd$ \ni Rd_{ij} \rightarrow Rd_{ij}^{`} \in $ \gf$^{`}_{4}$

\ \ \ \ \ \ \ \ \ \ \ \ \ \ \ \ \ \ \ge\gm$ \ni Rm_{ij} \rightarrow Rm_{ij}^{`} \in $ \gf$^{`}_{4}$

\bigskip

\emph{Proof.}  Let's check these Lie bracket operations.

For $Rd_{ij}^{`}/d_{ij}$ and $Rm_{ij}^{`}/m_{ij}$, by calculation with Maxima we have same results of 
\emph{Lemma 13.5}, \emph{Lemma 14.3}, and \emph{Lemma 14.4}. \ \ \ \ \ \emph{Q.E.D.}

\bigskip

Therfore we have \emph{Theorem 14.15}

\bigskip

\emph{Theorem 14.15.}  \gf$^{`}_{4}$ is isomorphic to \gf$_{4}$. 

\bigskip

We identify \gf$_{4}$ with \gf$^{`}_{4}$.

\bigskip

\emph{Lemma 14.16.}  \gf$_{4}$ is expressed  by 

$\{ D \in Hom_{\R}(\R^{27}) \mid D(\textrm{J}_{m}(X, Y))=\textrm{J}_{m}(D(X), Y)+\textrm{J}_{m}(X, D(Y)), X,Y\in \R^{27} \}.$
And $\{ Rd_{ij}^{`}(0 \leq i < j \leq 7)$,
$Rm_{ij}^{`}(1 \leq i \leq 3, 0 \leq j \leq 7) \}$
 are orthogonal bases of \gf$_{4}$.

\bigskip

\emph{Proof} \ By \emph{Remark 6.5} and \emph{Proposition 2.2}, \gf$_{4}$ is expressed  by 

$\{ D \in Hom_{\R}(\R^{27}) \mid D(\textrm{J}_{m}(X, Y))=\textrm{J}_{m}(D(X), Y)+\textrm{J}_{m}(X, D(Y)), X,Y\in \R^{27} \}.$

\noindent
Let check each element of $\{ Rd_{ij}^{`}(0 \leq i < j \leq 7)$,
$Rm_{ij}^{`}(1 \leq i \leq 3, 0 \leq j \leq 7) \}$
satisfy the condition:

$D(\textrm{J}_{m}(X, Y))=\textrm{J}_{m}(D(X), Y)+\textrm{J}_{m}(X, D(Y))$,

$X=(\chi_{1},\chi_{2},\chi_{3},(x_{10},\cdot\cdot\cdot,x_{17}),(x_{20},\cdot\cdot\cdot,x_{27}),(x_{30},\cdot\cdot\cdot,x_{37})) \ \in \R^{27}$,

$Y=(\gamma_{1},\gamma_{2},\gamma_{3},(y_{10},\cdot\cdot\cdot,y_{17}),(y_{20},\cdot\cdot\cdot,y_{27}),(y_{30},%
\cdot\cdot\cdot,y_{37}))\ \in \R^{27}$.

In case of $D=Rd_{01}^{`}/d_{01}$:       ($A_{(j)}$ means jth row element of $A$)

\noindent
$D(\textrm{J}_{m}(X,Y))_{(1)}=0$

\noindent
$D(\textrm{J}_{m}(X,Y))_{(2)}=0$

\noindent
$D(\textrm{J}_{m}(X,Y))_{(3)}=0$

\noindent
$D(\textrm{J}_{m}(X,Y))_{(4)}=$
$-2 x_{26} y_{37}+2 x_{27} y_{36}-2 x_{24} y_{35}+2 x_{25} y_{34}-2 x_{22} y_{33}+2 x_{23} y_{32}$
$-2 x_{20} y_{31}-2 x_{21} y_{30}+2 x_{36} y_{27}-2 x_{37} y_{26}+2 x_{34} y_{25}-2 x_{35} y_{24}+2 x_{32} y_{23}-2 x_{33} y_{22}$
$-2 x_{30} y_{21}-2 x_{31} y_{20}+2 \chi_{3} y_{11}+2 \chi_{2} y_{11}+2 x_{11} \gamma_{3}+2 x_{11} \gamma_{2}$

\noindent
$D(\textrm{J}_{m}(X,Y))_{(5)}=$
$2 x_{27} y_{37}+2 x_{26} y_{36}+2 x_{25} y_{35}+2 x_{24} y_{34}+2 x_{23} y_{33}+2 x_{22} y_{32}$
$+2 x_{21} y_{31}-2 x_{20} y_{30}+2 x_{37} y_{27}+2 x_{36} y_{26}+2 x_{35} y_{25}+2 x_{34} y_{24}+2 x_{33} y_{23}+2 x_{32} y_{22}$
$+2 x_{31} y_{21}-2 x_{30} y_{20}-2 \chi_{3} y_{10}-2 \chi_{2} y_{10}-2 x_{10} \gamma_{3}-2 x_{10} \gamma_{2}$

\noindent
$D(\textrm{J}_{m}(X,Y))_{(6)}=0$

\noindent
$D(\textrm{J}_{m}(X,Y))_{(7)}=0$

\noindent
$D(\textrm{J}_{m}(X,Y))_{(8)}=0$

\noindent
$D(\textrm{J}_{m}(X,Y))_{(9)}=0$

\noindent
$D(\textrm{J}_{m}(X,Y))_{(10)}=0$

\noindent
$D(\textrm{J}_{m}(X,Y))_{(11)}=0$

\noindent
$D(\textrm{J}_{m}(X,Y))_{(12)}=%
-x_{16} y_{37}+x_{17} y_{36}-x_{14} y_{35}+x_{15} y_{34}-x_{12} y_{33}+x_{13} y_{32}+x_{10} y_{31}$
$+x_{11} y_{30}-\chi_{3} y_{21}-\chi_{1} y_{21}+x_{36} y_{17}-x_{37} y_{16}+x_{34} y_{15}-x_{35} y_{14}+x_{32} y_{13}-x_{33} y_{12}$
$+x_{30} y_{11}+x_{31} y_{10}-x_{21} \gamma_{3}-x_{21} \gamma_{1}$

\noindent
$D(\textrm{J}_{m}(X,Y))_{(13)}=%
-x_{17} y_{37}-x_{16} y_{36}-x_{15} y_{35}-x_{14} y_{34}-x_{13} y_{33}-x_{12} y_{32}-x_{11} y_{31}$
$+x_{10} y_{30}+\chi_{3} y_{20}+\chi_{1} y_{20}-x_{37} y_{17}-x_{36} y_{16}-x_{35} y_{15}-x_{34} y_{14}-x_{33} y_{13}-x_{32} y_{12}$
$-x_{31} y_{11}+x_{30} y_{10}+x_{20} \gamma_{3}+x_{20} \gamma_{1}$

\noindent
$D(\textrm{J}_{m}(X,Y))_{(14)}=%
-x_{14} y_{37}-x_{15} y_{36}+x_{16} y_{35}+x_{17} y_{34}+x_{10} y_{33}-x_{11} y_{32}+x_{12} y_{31}$
$+x_{13} y_{30}-\chi_{3} y_{23}-\chi_{1} y_{23}+x_{34} y_{17}+x_{35} y_{16}-x_{36} y_{15}-x_{37} y_{14}+x_{30} y_{13}+x_{31} y_{12}$
$-x_{32} y_{11}+x_{33} y_{10}-x_{23} \gamma_{3}-x_{23} \gamma_{1}$

\noindent
$D(\textrm{J}_{m}(X,Y))_{(15)}=%
x_{15} y_{37}-x_{14} y_{36}-x_{17} y_{35}+x_{16} y_{34}-x_{11} y_{33}-x_{10} y_{32}+x_{13} y_{31}$
$-x_{12} y_{30}+\chi_{3} y_{22}+\chi_{1} y_{22}-x_{35} y_{17}+x_{34} y_{16}+x_{37} y_{15}-x_{36} y_{14}+x_{31} y_{13}-x_{30} y_{12}$
$-x_{33} y_{11}-x_{32} y_{10}+x_{22} \gamma_{3}+x_{22} \gamma_{1}$

\noindent
$D(\textrm{J}_{m}(X,Y))_{(16)}=%
x_{12} y_{37}+x_{13} y_{36}+x_{10} y_{35}-x_{11} y_{34}-x_{16} y_{33}-x_{17} y_{32}+x_{14} y_{31}$
$+x_{15} y_{30}-\chi_{3} y_{25}-\chi_{1} y_{25}-x_{32} y_{17}-x_{33} y_{16}+x_{30} y_{15}+x_{31} y_{14}+x_{36} y_{13}+x_{37} y_{12}$
$-x_{34} y_{11}+x_{35} y_{10}-x_{25} \gamma_{3}-x_{25} \gamma_{1}$

\noindent
$D(\textrm{J}_{m}(X,Y))_{(17)}=%
-x_{13} y_{37}+x_{12} y_{36}-x_{11} y_{35}-x_{10} y_{34}+x_{17} y_{33}-x_{16} y_{32}+x_{15} y_{31}$
$-x_{14} y_{30}+\chi_{3} y_{24}+\chi_{1} y_{24}+x_{33} y_{17}-x_{32} y_{16}+x_{31} y_{15}-x_{30} y_{14}-x_{37} y_{13}+x_{36} y_{12}$
$-x_{35} y_{11}-x_{34} y_{10}+x_{24} \gamma_{3}+x_{24} \gamma_{1}$

\noindent
$D(\textrm{J}_{m}(X,Y))_{(18)}=%
x_{10} y_{37}-x_{11} y_{36}-x_{12} y_{35}-x_{13} y_{34}+x_{14} y_{33}+x_{15} y_{32}+x_{16} y_{31}$
$+x_{17} y_{30}-\chi_{3} y_{27}-\chi_{1} y_{27}+x_{30} y_{17}+x_{31} y_{16}+x_{32} y_{15}+x_{33} y_{14}-x_{34} y_{13}-x_{35} y_{12}$
$-x_{36} y_{11}+x_{37} y_{10}-x_{27} \gamma_{3}-x_{27} \gamma_{1}$

\noindent
$D(\textrm{J}_{m}(X,Y))_{(19)}=%
-x_{11} y_{37}-x_{10} y_{36}+x_{13} y_{35}-x_{12} y_{34}-x_{15} y_{33}+x_{14} y_{32}+x_{17} y_{31}$
$-x_{16} y_{30}+\chi_{3} y_{26}+\chi_{1} y_{26}+x_{31} y_{17}-x_{30} y_{16}-x_{33} y_{15}+x_{32} y_{14}+x_{35} y_{13}-x_{34} y_{12}$
$-x_{37} y_{11}-x_{36} y_{10}+x_{26} \gamma_{3}+x_{26} \gamma_{1}$

\noindent
$D(\textrm{J}_{m}(X,Y))_{(20)}=%
-\chi_{2} y_{31}-\chi_{1} y_{31}+x_{16} y_{27}-x_{17} y_{26}+x_{14} y_{25}-x_{15} y_{24}+x_{12} y_{23}$
$-x_{13} y_{22}+x_{10} y_{21}+x_{11} y_{20}-x_{26} y_{17}+x_{27} y_{16}-x_{24} y_{15}+x_{25} y_{14}-x_{22} y_{13}+x_{23} y_{12}$
$+x_{20} y_{11}+x_{21} y_{10}-x_{31} \gamma_{2}-x_{31} \gamma_{1}$

\noindent
$D(\textrm{J}_{m}(X,Y))_{(21)}=%
\chi_{2} y_{30}+\chi_{1} y_{30}-x_{17} y_{27}-x_{16} y_{26}-x_{15} y_{25}-x_{14} y_{24}-x_{13} y_{23}$
$-x_{12} y_{22}-x_{11} y_{21}+x_{10} y_{20}-x_{27} y_{17}-x_{26} y_{16}-x_{25} y_{15}-x_{24} y_{14}-x_{23} y_{13}-x_{22} y_{12}$
$-x_{21} y_{11}+x_{20} y_{10}+x_{30} \gamma_{2}+x_{30} \gamma_{1}$

\noindent
$D(\textrm{J}_{m}(X,Y))_{(22)}=%
\chi_{2} y_{33}+\chi_{1} y_{33}-x_{14} y_{27}-x_{15} y_{26}+x_{16} y_{25}+x_{17} y_{24}-x_{10} y_{23}$
$-x_{11} y_{22}+x_{12} y_{21}-x_{13} y_{20}+x_{24} y_{17}+x_{25} y_{16}-x_{26} y_{15}-x_{27} y_{14}-x_{20} y_{13}+x_{21} y_{12}$
$-x_{22} y_{11}-x_{23} y_{10}+x_{33} \gamma_{2}+x_{33} \gamma_{1}$

\noindent
$D(\textrm{J}_{m}(X,Y))_{(23)}=%
-\chi_{2} y_{32}-\chi_{1} y_{32}+x_{15} y_{27}-x_{14} y_{26}-x_{17} y_{25}+x_{16} y_{24}-x_{11} y_{23}$
$+x_{10} y_{22}+x_{13} y_{21}+x_{12} y_{20}-x_{25} y_{17}+x_{24} y_{16}+x_{27} y_{15}-x_{26} y_{14}+x_{21} y_{13}+x_{20} y_{12}$
$-x_{23} y_{11}+x_{22} y_{10}-x_{32} \gamma_{2}-x_{32} \gamma_{1}$

\noindent
$D(\textrm{J}_{m}(X,Y))_{(24)}=%
\chi_{2} y_{35}+\chi_{1} y_{35}+x_{12} y_{27}+x_{13} y_{26}-x_{10} y_{25}-x_{11} y_{24}-x_{16} y_{23}$
$-x_{17} y_{22}+x_{14} y_{21}-x_{15} y_{20}-x_{22} y_{17}-x_{23} y_{16}-x_{20} y_{15}+x_{21} y_{14}+x_{26} y_{13}+x_{27} y_{12}$
$-x_{24} y_{11}-x_{25} y_{10}+x_{35} \gamma_{2}+x_{35} \gamma_{1}$

\noindent
$D(\textrm{J}_{m}(X,Y))_{(25)}=%
-\chi_{2} y_{34}-\chi_{1} y_{34}-x_{13} y_{27}+x_{12} y_{26}-x_{11} y_{25}+x_{10} y_{24}+x_{17} y_{23}$
$-x_{16} y_{22}+x_{15} y_{21}+x_{14} y_{20}+x_{23} y_{17}-x_{22} y_{16}+x_{21} y_{15}+x_{20} y_{14}-x_{27} y_{13}+x_{26} y_{12}$
$-x_{25} y_{11}+x_{24} y_{10}-x_{34} \gamma_{2}-x_{34} \gamma_{1}$

\noindent
$D(\textrm{J}_{m}(X,Y))_{(26)}=%
\chi_{2} y_{37}+\chi_{1} y_{37}-x_{10} y_{27}-x_{11} y_{26}-x_{12} y_{25}-x_{13} y_{24}+x_{14} y_{23}$
$+x_{15} y_{22}+x_{16} y_{21}-x_{17} y_{20}-x_{20} y_{17}+x_{21} y_{16}+x_{22} y_{15}+x_{23} y_{14}-x_{24} y_{13}-x_{25} y_{12}$
$-x_{26} y_{11}-x_{27} y_{10}+x_{37} \gamma_{2}+x_{37} \gamma_{1}$

\noindent
$D(\textrm{J}_{m}(X,Y))_{(27)}=%
-\chi_{2} y_{36}-\chi_{1} y_{36}-x_{11} y_{27}+x_{10} y_{26}+x_{13} y_{25}-x_{12} y_{24}-x_{15} y_{23}$
$+x_{14} y_{22}+x_{17} y_{21}+x_{16} y_{20}+x_{21} y_{17}+x_{20} y_{16}-x_{23} y_{15}+x_{22} y_{14}+x_{25} y_{13}-x_{24} y_{12}$
$-x_{27} y_{11}+x_{26} y_{10}-x_{36} \gamma_{2}-x_{36} \gamma_{1}$

\bigskip

\noindent
$\textrm{J}_{m}(D(X),Y))_{(1)}=%
-2 x_{36} y_{37}+2 x_{37} y_{36}-2 x_{34} y_{35}+2 x_{35} y_{34}-2 x_{32} y_{33}+2 x_{33} y_{32}$
$+2 x_{30} y_{31}-2 x_{31} y_{30}+2 x_{26} y_{27}-2 x_{27} y_{26}+2 x_{24} y_{25}-2 x_{25} y_{24}+2 x_{22} y_{23}-2 x_{23} y_{22}$
$+2 x_{20} y_{21}-2 x_{21} y_{20}$

\noindent
$\textrm{J}_{m}(D(X),Y))_{(2)}=%
-2 x_{36} y_{37}+2 x_{37} y_{36}-2 x_{34} y_{35}+2 x_{35} y_{34}-2 x_{32} y_{33}+2 x_{33} y_{32}$
$+2 x_{30} y_{31}-2 x_{31} y_{30}-4 x_{10} y_{11}+4 x_{11} y_{10}$

\noindent
$\textrm{J}_{m}(D(X),Y))_{(3)}=%
-2 x_{36} y_{37}+2 x_{37} y_{36}-2 x_{34} y_{35}+2 x_{35} y_{34}-2 x_{32} y_{33}+2 x_{33} y_{32}$
$+2 x_{30} y_{31}-2 x_{31} y_{30}-4 x_{10} y_{11}+4 x_{11} y_{10}$

\noindent
$\textrm{J}_{m}(D(X),Y))_{(4)}=%
-x_{26} y_{37}+x_{27} y_{36}-x_{24} y_{35}+x_{25} y_{34}-x_{22} y_{33}+x_{23} y_{32}-x_{20} y_{31}$
$-x_{21} y_{30}+x_{36} y_{27}-x_{37} y_{26}+x_{34} y_{25}-x_{35} y_{24}+x_{32} y_{23}-x_{33} y_{22}-x_{30} y_{21}-x_{31} y_{20}$
$+2 x_{11} \gamma_{3}+2 x_{11} \gamma_{2}$

\noindent
$\textrm{J}_{m}(D(X),Y))_{(5)}=%
x_{27} y_{37}+x_{26} y_{36}+x_{25} y_{35}+x_{24} y_{34}+x_{23} y_{33}+x_{22} y_{32}+x_{21} y_{31}$
$-x_{20} y_{30}+x_{37} y_{27}+x_{36} y_{26}+x_{35} y_{25}+x_{34} y_{24}+x_{33} y_{23}+x_{32} y_{22}+x_{31} y_{21}-x_{30} y_{20}$
$-2 x_{10} \gamma_{3}-2 x_{10} \gamma_{2}$

\noindent
$\textrm{J}_{m}(D(X),Y))_{(6)}=%
-x_{24} y_{37}-x_{25} y_{36}+x_{26} y_{35}+x_{27} y_{34}+x_{20} y_{33}+x_{21} y_{32}-x_{22} y_{31}$
$+x_{23} y_{30}-x_{34} y_{27}-x_{35} y_{26}+x_{36} y_{25}+x_{37} y_{24}-x_{30} y_{23}+x_{31} y_{22}-x_{32} y_{21}-x_{33} y_{20}$

\noindent
$\textrm{J}_{m}(D(X),Y))_{(7)}=%
-x_{24} y_{37}-x_{25} y_{36}+x_{26} y_{35}+x_{27} y_{34}+x_{20} y_{33}+x_{21} y_{32}-x_{22} y_{31}$
$+x_{23} y_{30}-x_{34} y_{27}-x_{35} y_{26}+x_{36} y_{25}+x_{37} y_{24}-x_{30} y_{23}+x_{31} y_{22}-x_{32} y_{21}-x_{33} y_{20}$

\noindent
$\textrm{J}_{m}(D(X),Y))_{(8)}=%
-x_{24} y_{37}-x_{25} y_{36}+x_{26} y_{35}+x_{27} y_{34}+x_{20} y_{33}+x_{21} y_{32}-x_{22} y_{31}$
$+x_{23} y_{30}-x_{34} y_{27}-x_{35} y_{26}+x_{36} y_{25}+x_{37} y_{24}-x_{30} y_{23}+x_{31} y_{22}-x_{32} y_{21}-x_{33} y_{20}$

\noindent
$\textrm{J}_{m}(D(X),Y))_{(9)}=%
-x_{24} y_{37}-x_{25} y_{36}+x_{26} y_{35}+x_{27} y_{34}+x_{20} y_{33}+x_{21} y_{32}-x_{22} y_{31}$
$+x_{23} y_{30}-x_{34} y_{27}-x_{35} y_{26}+x_{36} y_{25}+x_{37} y_{24}-x_{30} y_{23}+x_{31} y_{22}-x_{32} y_{21}-x_{33} y_{20}$

\noindent
$\textrm{J}_{m}(D(X),Y))_{(10)}=%
-x_{24} y_{37}-x_{25} y_{36}+x_{26} y_{35}+x_{27} y_{34}+x_{20} y_{33}+x_{21} y_{32}-x_{22} y_{31}$
$+x_{23} y_{30}-x_{34} y_{27}-x_{35} y_{26}+x_{36} y_{25}+x_{37} y_{24}-x_{30} y_{23}+x_{31} y_{22}-x_{32} y_{21}-x_{33} y_{20}$

\noindent
$\textrm{J}_{m}(D(X),Y))_{(11)}=%
-x_{24} y_{37}-x_{25} y_{36}+x_{26} y_{35}+x_{27} y_{34}+x_{20} y_{33}+x_{21} y_{32}-x_{22} y_{31}$
$+x_{23} y_{30}-x_{34} y_{27}-x_{35} y_{26}+x_{36} y_{25}+x_{37} y_{24}-x_{30} y_{23}+x_{31} y_{22}-x_{32} y_{21}-x_{33} y_{20}$

\noindent
$\textrm{J}_{m}(D(X),Y))_{(12)}=%
-x_{24} y_{37}-x_{25} y_{36}+x_{26} y_{35}+x_{27} y_{34}+x_{20} y_{33}+x_{21} y_{32}-x_{22} y_{31}$
$+x_{23} y_{30}-x_{34} y_{27}-x_{35} y_{26}+x_{36} y_{25}+x_{37} y_{24}-x_{30} y_{23}+x_{31} y_{22}-x_{32} y_{21}-x_{33} y_{20}$

\noindent
$\textrm{J}_{m}(D(X),Y))_{(13)}=%
-x_{24} y_{37}-x_{25} y_{36}+x_{26} y_{35}+x_{27} y_{34}+x_{20} y_{33}+x_{21} y_{32}-x_{22} y_{31}$
$+x_{23} y_{30}-x_{34} y_{27}-x_{35} y_{26}+x_{36} y_{25}+x_{37} y_{24}-x_{30} y_{23}+x_{31} y_{22}-x_{32} y_{21}-x_{33} y_{20}$

\noindent
$\textrm{J}_{m}(D(X),Y))_{(14)}=%
-x_{24} y_{37}-x_{25} y_{36}+x_{26} y_{35}+x_{27} y_{34}+x_{20} y_{33}+x_{21} y_{32}-x_{22} y_{31}$
$+x_{23} y_{30}-x_{34} y_{27}-x_{35} y_{26}+x_{36} y_{25}+x_{37} y_{24}-x_{30} y_{23}+x_{31} y_{22}-x_{32} y_{21}-x_{33} y_{20}$

\noindent
$\textrm{J}_{m}(D(X),Y))_{(15)}=%
-x_{24} y_{37}-x_{25} y_{36}+x_{26} y_{35}+x_{27} y_{34}+x_{20} y_{33}+x_{21} y_{32}-x_{22} y_{31}$
$+x_{23} y_{30}-x_{34} y_{27}-x_{35} y_{26}+x_{36} y_{25}+x_{37} y_{24}-x_{30} y_{23}+x_{31} y_{22}-x_{32} y_{21}-x_{33} y_{20}$

\noindent
$\textrm{J}_{m}(D(X),Y))_{(16)}=%
-x_{24} y_{37}-x_{25} y_{36}+x_{26} y_{35}+x_{27} y_{34}+x_{20} y_{33}+x_{21} y_{32}-x_{22} y_{31}$
$+x_{23} y_{30}-x_{34} y_{27}-x_{35} y_{26}+x_{36} y_{25}+x_{37} y_{24}-x_{30} y_{23}+x_{31} y_{22}-x_{32} y_{21}-x_{33} y_{20}$

\noindent
$\textrm{J}_{m}(D(X),Y))_{(17)}=%
-x_{24} y_{37}-x_{25} y_{36}+x_{26} y_{35}+x_{27} y_{34}+x_{20} y_{33}+x_{21} y_{32}-x_{22} y_{31}$
$+x_{23} y_{30}-x_{34} y_{27}-x_{35} y_{26}+x_{36} y_{25}+x_{37} y_{24}-x_{30} y_{23}+x_{31} y_{22}-x_{32} y_{21}-x_{33} y_{20}$

\noindent
$\textrm{J}_{m}(D(X),Y))_{(18)}=%
-x_{24} y_{37}-x_{25} y_{36}+x_{26} y_{35}+x_{27} y_{34}+x_{20} y_{33}+x_{21} y_{32}-x_{22} y_{31}$
$+x_{23} y_{30}-x_{34} y_{27}-x_{35} y_{26}+x_{36} y_{25}+x_{37} y_{24}-x_{30} y_{23}+x_{31} y_{22}-x_{32} y_{21}-x_{33} y_{20}$

\noindent
$\textrm{J}_{m}(D(X),Y))_{(19)}=%
-x_{24} y_{37}-x_{25} y_{36}+x_{26} y_{35}+x_{27} y_{34}+x_{20} y_{33}+x_{21} y_{32}-x_{22} y_{31}$
$+x_{23} y_{30}-x_{34} y_{27}-x_{35} y_{26}+x_{36} y_{25}+x_{37} y_{24}-x_{30} y_{23}+x_{31} y_{22}-x_{32} y_{21}-x_{33} y_{20}$

\noindent
$\textrm{J}_{m}(D(X),Y))_{(20)}=%
-x_{24} y_{37}-x_{25} y_{36}+x_{26} y_{35}+x_{27} y_{34}+x_{20} y_{33}+x_{21} y_{32}-x_{22} y_{31}$
$+x_{23} y_{30}-x_{34} y_{27}-x_{35} y_{26}+x_{36} y_{25}+x_{37} y_{24}-x_{30} y_{23}+x_{31} y_{22}-x_{32} y_{21}-x_{33} y_{20}$

\noindent
$\textrm{J}_{m}(D(X),Y))_{(21)}=%
-x_{24} y_{37}-x_{25} y_{36}+x_{26} y_{35}+x_{27} y_{34}+x_{20} y_{33}+x_{21} y_{32}-x_{22} y_{31}$
$+x_{23} y_{30}-x_{34} y_{27}-x_{35} y_{26}+x_{36} y_{25}+x_{37} y_{24}-x_{30} y_{23}+x_{31} y_{22}-x_{32} y_{21}-x_{33} y_{20}$

\noindent
$\textrm{J}_{m}(D(X),Y))_{(22)}=%
-x_{24} y_{37}-x_{25} y_{36}+x_{26} y_{35}+x_{27} y_{34}+x_{20} y_{33}+x_{21} y_{32}-x_{22} y_{31}$
$+x_{23} y_{30}-x_{34} y_{27}-x_{35} y_{26}+x_{36} y_{25}+x_{37} y_{24}-x_{30} y_{23}+x_{31} y_{22}-x_{32} y_{21}-x_{33} y_{20}$

\noindent
$\textrm{J}_{m}(D(X),Y))_{(23)}=%
-x_{24} y_{37}-x_{25} y_{36}+x_{26} y_{35}+x_{27} y_{34}+x_{20} y_{33}+x_{21} y_{32}-x_{22} y_{31}$
$+x_{23} y_{30}-x_{34} y_{27}-x_{35} y_{26}+x_{36} y_{25}+x_{37} y_{24}-x_{30} y_{23}+x_{31} y_{22}-x_{32} y_{21}-x_{33} y_{20}$

\noindent
$\textrm{J}_{m}(D(X),Y))_{(24)}=%
-x_{24} y_{37}-x_{25} y_{36}+x_{26} y_{35}+x_{27} y_{34}+x_{20} y_{33}+x_{21} y_{32}-x_{22} y_{31}$
$+x_{23} y_{30}-x_{34} y_{27}-x_{35} y_{26}+x_{36} y_{25}+x_{37} y_{24}-x_{30} y_{23}+x_{31} y_{22}-x_{32} y_{21}-x_{33} y_{20}$

\noindent
$\textrm{J}_{m}(D(X),Y))_{(25)}=%
-x_{24} y_{37}-x_{25} y_{36}+x_{26} y_{35}+x_{27} y_{34}+x_{20} y_{33}+x_{21} y_{32}-x_{22} y_{31}$
$+x_{23} y_{30}-x_{34} y_{27}-x_{35} y_{26}+x_{36} y_{25}+x_{37} y_{24}-x_{30} y_{23}+x_{31} y_{22}-x_{32} y_{21}-x_{33} y_{20}$

\noindent
$\textrm{J}_{m}(D(X),Y))_{(26)}=%
-x_{24} y_{37}-x_{25} y_{36}+x_{26} y_{35}+x_{27} y_{34}+x_{20} y_{33}+x_{21} y_{32}-x_{22} y_{31}$
$+x_{23} y_{30}-x_{34} y_{27}-x_{35} y_{26}+x_{36} y_{25}+x_{37} y_{24}-x_{30} y_{23}+x_{31} y_{22}-x_{32} y_{21}-x_{33} y_{20}$

\noindent
$\textrm{J}_{m}(D(X),Y))_{(27)}=%
-x_{24} y_{37}-x_{25} y_{36}+x_{26} y_{35}+x_{27} y_{34}+x_{20} y_{33}+x_{21} y_{32}-x_{22} y_{31}$
$+x_{23} y_{30}-x_{34} y_{27}-x_{35} y_{26}+x_{36} y_{25}+x_{37} y_{24}-x_{30} y_{23}+x_{31} y_{22}-x_{32} y_{21}-x_{33} y_{20}$

\bigskip

\noindent
$\textrm{J}_{m}(X,D(Y))_{(1)}=%
2 x_{36} y_{37}-2 x_{37} y_{36}+2 x_{34} y_{35}-2 x_{35} y_{34}+2 x_{32} y_{33}-2 x_{33} y_{32}$
$-2 x_{30} y_{31}+2 x_{31} y_{30}-2 x_{26} y_{27}+2 x_{27} y_{26}-2 x_{24} y_{25}+2 x_{25} y_{24}-2 x_{22} y_{23}+2 x_{23} y_{22}$
$-2 x_{20} y_{21}+2 x_{21} y_{20}$

\noindent
$\textrm{J}_{m}(X,D(Y))_{(2)}=%
2 x_{36} y_{37}-2 x_{37} y_{36}+2 x_{34} y_{35}-2 x_{35} y_{34}+2 x_{32} y_{33}-2 x_{33} y_{32}$
$-2 x_{30} y_{31}+2 x_{31} y_{30}+4 x_{10} y_{11}-4 x_{11} y_{10}$

\noindent
$\textrm{J}_{m}(X,D(Y))_{(3)}=%
-2 x_{26} y_{27}+2 x_{27} y_{26}-2 x_{24} y_{25}+2 x_{25} y_{24}-2 x_{22} y_{23}+2 x_{23} y_{22}$
$-2 x_{20} y_{21}+2 x_{21} y_{20}+4 x_{10} y_{11}-4 x_{11} y_{10}$

\noindent
$\textrm{J}_{m}(X,D(Y))_{(4)}=%
-x_{26} y_{37}+x_{27} y_{36}-x_{24} y_{35}+x_{25} y_{34}-x_{22} y_{33}+x_{23} y_{32}-x_{20} y_{31}$
$-x_{21} y_{30}+x_{36} y_{27}-x_{37} y_{26}+x_{34} y_{25}-x_{35} y_{24}+x_{32} y_{23}-x_{33} y_{22}-x_{30} y_{21}-x_{31} y_{20}$
$+2 \chi_{3} y_{11}+2 \chi_{2} y_{11}$

\noindent
$\textrm{J}_{m}(X,D(Y))_{(5)}=%
x_{27} y_{37}+x_{26} y_{36}+x_{25} y_{35}+x_{24} y_{34}+x_{23} y_{33}+x_{22} y_{32}+x_{21} y_{31}$
$-x_{20} y_{30}+x_{37} y_{27}+x_{36} y_{26}+x_{35} y_{25}+x_{34} y_{24}+x_{33} y_{23}+x_{32} y_{22}+x_{31} y_{21}-x_{30} y_{20}$
$-2 \chi_{3} y_{10}-2 \chi_{2} y_{10}$

\noindent
$\textrm{J}_{m}(X,D(Y))_{(6)}=%
x_{24} y_{37}+x_{25} y_{36}-x_{26} y_{35}-x_{27} y_{34}-x_{20} y_{33}-x_{21} y_{32}+x_{22} y_{31}$
$-x_{23} y_{30}+x_{34} y_{27}+x_{35} y_{26}-x_{36} y_{25}-x_{37} y_{24}+x_{30} y_{23}-x_{31} y_{22}+x_{32} y_{21}+x_{33} y_{20}$

\noindent
$\textrm{J}_{m}(X,D(Y))_{(7)}=%
-x_{25} y_{37}+x_{24} y_{36}+x_{27} y_{35}-x_{26} y_{34}-x_{21} y_{33}+x_{20} y_{32}+x_{23} y_{31}$
$+x_{22} y_{30}-x_{35} y_{27}+x_{34} y_{26}+x_{37} y_{25}-x_{36} y_{24}-x_{31} y_{23}-x_{30} y_{22}+x_{33} y_{21}-x_{32} y_{20}$

\noindent
$\textrm{J}_{m}(X,D(Y))_{(8)}=%
-x_{22} y_{37}-x_{23} y_{36}-x_{20} y_{35}-x_{21} y_{34}+x_{26} y_{33}+x_{27} y_{32}+x_{24} y_{31}$
$-x_{25} y_{30}-x_{32} y_{27}-x_{33} y_{26}+x_{30} y_{25}-x_{31} y_{24}+x_{36} y_{23}+x_{37} y_{22}+x_{34} y_{21}+x_{35} y_{20}$

\noindent
$\textrm{J}_{m}(X,D(Y))_{(9)}=%
x_{23} y_{37}-x_{22} y_{36}-x_{21} y_{35}+x_{20} y_{34}-x_{27} y_{33}+x_{26} y_{32}+x_{25} y_{31}$
$+x_{24} y_{30}+x_{33} y_{27}-x_{32} y_{26}-x_{31} y_{25}-x_{30} y_{24}-x_{37} y_{23}+x_{36} y_{22}+x_{35} y_{21}-x_{34} y_{20}$

\noindent
$\textrm{J}_{m}(X,D(Y))_{10)}=%
-x_{20} y_{37}-x_{21} y_{36}+x_{22} y_{35}+x_{23} y_{34}-x_{24} y_{33}-x_{25} y_{32}+x_{26} y_{31}$
$-x_{27} y_{30}+x_{30} y_{27}-x_{31} y_{26}+x_{32} y_{25}+x_{33} y_{24}-x_{34} y_{23}-x_{35} y_{22}+x_{36} y_{21}+x_{37} y_{20}$

\noindent
$\textrm{J}_{m}(X,D(Y))_{(11)}=%
-x_{21} y_{37}+x_{20} y_{36}-x_{23} y_{35}+x_{22} y_{34}+x_{25} y_{33}-x_{24} y_{32}+x_{27} y_{31}$
$+x_{26} y_{30}-x_{31} y_{27}-x_{30} y_{26}-x_{33} y_{25}+x_{32} y_{24}+x_{35} y_{23}-x_{34} y_{22}+x_{37} y_{21}-x_{36} y_{20}$

\noindent
$\textrm{J}_{m}(X,D(Y))_{(12)}=%
-x_{16} y_{37}+x_{17} y_{36}-x_{14} y_{35}+x_{15} y_{34}-x_{12} y_{33}+x_{13} y_{32}-x_{10} y_{31}$
$-x_{11} y_{30}-\chi_{3} y_{21}-\chi_{1} y_{21}+2 x_{30} y_{11}+2 x_{31} y_{10}$

\noindent
$\textrm{J}_{m}(X,D(Y))_{(13)}=%
-x_{17} y_{37}-x_{16} y_{36}-x_{15} y_{35}-x_{14} y_{34}-x_{13} y_{33}-x_{12} y_{32}+x_{11} y_{31}$
$-x_{10} y_{30}+\chi_{3} y_{20}+\chi_{1} y_{20}-2 x_{31} y_{11}+2 x_{30} y_{10}$

\noindent
$\textrm{J}_{m}(X,D(Y))_{(14)}=%
-x_{14} y_{37}-x_{15} y_{36}+x_{16} y_{35}+x_{17} y_{34}-x_{10} y_{33}+x_{11} y_{32}+x_{12} y_{31}$
$+x_{13} y_{30}-\chi_{3} y_{23}-\chi_{1} y_{23}-2 x_{32} y_{11}+2 x_{33} y_{10}$

\noindent
$\textrm{J}_{m}(X,D(Y))_{(15)}=%
x_{15} y_{37}-x_{14} y_{36}-x_{17} y_{35}+x_{16} y_{34}+x_{11} y_{33}+x_{10} y_{32}+x_{13} y_{31}$
$-x_{12} y_{30}+\chi_{3} y_{22}+\chi_{1} y_{22}-2 x_{33} y_{11}-2 x_{32} y_{10}$

\noindent
$\textrm{J}_{m}(X,D(Y))_{(16)}=%
x_{12} y_{37}+x_{13} y_{36}-x_{10} y_{35}+x_{11} y_{34}-x_{16} y_{33}-x_{17} y_{32}+x_{14} y_{31}$
$+x_{15} y_{30}-\chi_{3} y_{25}-\chi_{1} y_{25}-2 x_{34} y_{11}+2 x_{35} y_{10}$

\noindent
$\textrm{J}_{m}(X,D(Y))_{(17)}=%
-x_{13} y_{37}+x_{12} y_{36}+x_{11} y_{35}+x_{10} y_{34}+x_{17} y_{33}-x_{16} y_{32}+x_{15} y_{31}$
$-x_{14} y_{30}+\chi_{3} y_{24}+\chi_{1} y_{24}-2 x_{35} y_{11}-2 x_{34} y_{10}$

\noindent
$\textrm{J}_{m}(X,D(Y))_{(18)}=%
-x_{10} y_{37}+x_{11} y_{36}-x_{12} y_{35}-x_{13} y_{34}+x_{14} y_{33}+x_{15} y_{32}+x_{16} y_{31}$
$+x_{17} y_{30}-\chi_{3} y_{27}-\chi_{1} y_{27}-2 x_{36} y_{11}+2 x_{37} y_{10}$

\noindent
$\textrm{J}_{m}(X,D(Y))_{(19)}=%
x_{11} y_{37}+x_{10} y_{36}+x_{13} y_{35}-x_{12} y_{34}-x_{15} y_{33}+x_{14} y_{32}+x_{17} y_{31}$
$-x_{16} y_{30}+\chi_{3} y_{26}+\chi_{1} y_{26}-2 x_{37} y_{11}-2 x_{36} y_{10}$

\noindent
$\textrm{J}_{m}(X,D(Y))_{(20)}=%
-\chi_{2} y_{31}-\chi_{1} y_{31}+x_{16} y_{27}-x_{17} y_{26}+x_{14} y_{25}-x_{15} y_{24}+x_{12} y_{23}$
$-x_{13} y_{22}-x_{10} y_{21}-x_{11} y_{20}+2 x_{20} y_{11}+2 x_{21} y_{10}$

\noindent
$\textrm{J}_{m}(X,D(Y))_{(21)}=%
\chi_{2} y_{30}+\chi_{1} y_{30}-x_{17} y_{27}-x_{16} y_{26}-x_{15} y_{25}-x_{14} y_{24}-x_{13} y_{23}$
$-x_{12} y_{22}+x_{11} y_{21}-x_{10} y_{20}-2 x_{21} y_{11}+2 x_{20} y_{10}$

\noindent
$\textrm{J}_{m}(X,D(Y))_{(22)}=%
\chi_{2} y_{33}+\chi_{1} y_{33}-x_{14} y_{27}-x_{15} y_{26}+x_{16} y_{25}+x_{17} y_{24}+x_{10} y_{23}$
$+x_{11} y_{22}+x_{12} y_{21}-x_{13} y_{20}-2 x_{22} y_{11}-2 x_{23} y_{10}$

\noindent
$\textrm{J}_{m}(X,D(Y))_{(23)}=%
-\chi_{2} y_{32}-\chi_{1} y_{32}+x_{15} y_{27}-x_{14} y_{26}-x_{17} y_{25}+x_{16} y_{24}+x_{11} y_{23}$
$-x_{10} y_{22}+x_{13} y_{21}+x_{12} y_{20}-2 x_{23} y_{11}+2 x_{22} y_{10}$

\noindent
$\textrm{J}_{m}(X,D(Y))_{(24)}=%
\chi_{2} y_{35}+\chi_{1} y_{35}+x_{12} y_{27}+x_{13} y_{26}+x_{10} y_{25}+x_{11} y_{24}-x_{16} y_{23}$
$-x_{17} y_{22}+x_{14} y_{21}-x_{15} y_{20}-2 x_{24} y_{11}-2 x_{25} y_{10}$

\noindent
$\textrm{J}_{m}(X,D(Y))_{(25)}=%
-\chi_{2} y_{34}-\chi_{1} y_{34}-x_{13} y_{27}+x_{12} y_{26}+x_{11} y_{25}-x_{10} y_{24}+x_{17} y_{23}$
$-x_{16} y_{22}+x_{15} y_{21}+x_{14} y_{20}-2 x_{25} y_{11}+2 x_{24} y_{10}$

\noindent
$\textrm{J}_{m}(X,D(Y))_{(26)}=%
\chi_{2} y_{37}+\chi_{1} y_{37}+x_{10} y_{27}+x_{11} y_{26}-x_{12} y_{25}-x_{13} y_{24}+x_{14} y_{23}$
$+x_{15} y_{22}+x_{16} y_{21}-x_{17} y_{20}-2 x_{26} y_{11}-2 x_{27} y_{10}$

\noindent
$\textrm{J}_{m}(X,D(Y))_{(27)}=%
-\chi_{2} y_{36}-\chi_{1} y_{36}+x_{11} y_{27}-x_{10} y_{26}+x_{13} y_{25}-x_{12} y_{24}-x_{15} y_{23}$
$+x_{14} y_{22}+x_{17} y_{21}+x_{16} y_{20}-2 x_{27} y_{11}+2 x_{26} y_{10}$

\noindent Therfore  we have $D(\textrm{J}_{m}(X, Y))=\textrm{J}_{m}(D(X), Y)+\textrm{J}_{m}(X, D(Y))$.

Other cases can be confirmed in the same way.

For $R_{1},R_{2} \in  \{ Rd_{ij}^{`}/d_{ij}(0 \leq i < j \leq 7),Rm_{ij}^{`}/m_{ij}(1 \leq i \leq 3, 0 \leq j \leq 7) \}$,
let's  calculate $tr(R_{1}.^{t}R_{2})$.

\ \ \ \ \ \ \ \ \ \ \ \ \ \ \ \ $tr(R_{1}.^{t}R_{2})=0 \ \ \ (R_{1} \neq R_{2})$,

\ \ \ \ \ \ \ \ \ \ \ \ \ \ \ \ $tr(R_{1}.^{t}R_{1})=6 \ \ \ (R_{1} \in  \{ Rd_{ij}^{`}/d_{ij}(0 \leq i < j \leq 7)\}$,

\ \ \ \ \ \ \ \ \ \ \ \ \ \ \ \ $tr(R_{2}.^{t}R_{2})=\frac{13}{2} \ \ (R_{2} \in \{Rm_{ij}^{`}/m_{ij}(1 \leq i \leq 3, 0 \leq j \leq 7) \}$.

Since dim \gf$_{4}=52$, then $\{ Rd_{ij}^{`}(0 \leq i < j \leq 7),Rm_{ij}^{`}(1 \leq i \leq 3, 0 \leq j \leq 7) \}$ are 
orthogonal bases of \gf$_{4}$.
\ \ \ \ \ \emph{Q.E.D.}

\bigskip

\emph{Lemma 14.17.}  The simply connected compact Lie group $F_{4}$ is expressed by

$F_{4}=\{\alpha \in Iso_{\R}(\R^{27}) \mid \alpha (\textrm{J}_{m}(X, Y))=\textrm{J}_{m}(\alpha X,\alpha Y) \}$ .

\bigskip

\emph{Proof.}  By \emph{Proposition 2.1} and \emph{Remark 6.5}, we have the above \emph{Lemma}.\ \ \ \ \ \emph{Q.E.D.}

\bigskip

\section{The exceptional simple Lie algebra \gr$_{6}^{\C}$ of type $E_{6}$}

\bigskip

\ \ \ \ \emph{Definition 15.1.} \ We define the followings for elements $Rt_{ij},R\tau
_{k}\in $\gR\gt$^{\C}.$
\begin{align*}
Ut_{ij}&=Rt_{ij}/t_{ij}\ \ (1\leq i\leq 3,0\leq j\leq 7),\\
U\tau _{k}&=R\tau _{k}/\tau _{k} \ (k=1,2).
\end{align*}

\bigskip

\emph{Lemma 15.2.} \ We have the following Lie bracket operations.

\noindent
$[Ut_{1i},Ut_{2j}]=-\frac{1}{2}Sn(i+1,j+1)Um_{3k}\,(k=Ca(i+1,j+1),\ 0\leq i,j\leq 7)$,

\noindent
$[Ut_{2i},Ut_{3j}]=-\frac{1}{2}Sn(i+1,j+1)Um_{1k}\,(k=Ca(i+1,j+1),\ 0\leq i,j\leq 7)$,

\noindent
$[Ut_{3i},Ut_{1j}]=-\frac{1}{2}Sn(i+1,j+1)Um_{2k}\,(k=Ca(i+1,j+1),\ 0\leq i,j\leq 7)$,

\noindent
$[Ut_{1i},Ut_{1j}]=Ud_{ij}\ ,(0\leq i<j\leq 7)$,

\noindent
$[Ut_{2i},Ut_{2j}]=\sum\limits_{0\leq n<l\leq 7}Mv(ki,kj)Ud_{nl}$,

\ \ \ \ \ \ \ \ \ \ $(ki=Nu(i+1,j+1),kj=Nu(n+1,l+1),0\leq i<j\leq 7)$,

\noindent
$[Ut_{3i},Ut_{3j}]=\sum\limits_{0\leq n<l\leq 7}Mv^{2}(ki,kj)Ud_{nl}$,

\ \ \ \ \ \ \ \ \ \ $(ki=Nu(i+1,j+1),kj=Nu(n+1,l+1),0\leq i<j\leq 7)$,

\noindent
$[U\tau _{1},Ut_{1i}]=\frac{1}{2}Um_{1i},(0\leq i\leq 7)$,

\noindent
$[U\tau _{2},Ut_{1i}]=Um_{1i},(0\leq i\leq 7)$,

\noindent
$[U\tau _{1},Ut_{2i}]=-Um_{2i},(0\leq i\leq 7)$,

\noindent
$[U\tau _{2},Ut_{2i}]=-\frac{1}{2}Um_{2i},(0\leq i\leq 7)$,

\noindent
$[U\tau _{1},Ut_{3i}]=\frac{1}{2}Um_{3i},(0\leq i\leq 7)$,

\noindent
$[U\tau _{2},Ut_{3i}]=-\frac{1}{2}Um_{3i},(0\leq i\leq 7)$,

\noindent
$[U\tau _{1},U\tau _{2}]=0$.

\bigskip

\emph{Proof. \ }We have the above Lie bracket operations with calculations
using Maxima.\ \ \ \ \emph{Q.E.D.}

\bigskip

\emph{Lemma 15.3. }\ We have the following Lie bracket operations.

\ \ \ \ $[Ud_{ij},Ut_{1k}]=-Ut_{1j}$ $($in case of $k=i),$

$\ \ \ \ \ \ \ \ \ \ \ \ \ \ \ \ \ \ \ =Ut_{1i}$ $($in case of $%
k=j),$

$\ \ \ \ \ \ \ \ \ \ \ \ \ \ \ \ \ \ \ =0$ $($in case of $%
k\neq i,j),(0\leq i<j\leq 7),$

\ \ \ \ $[Ud_{ij},Ut_{2k}]=-\sum\limits_{0\leq n<l\leq 7}Mv^{2}(ki,kj)Ut_{2l}$ $($where $k=n)$

\ \ \ \ \ \ \ \ \ \ \ \ \ \ \ \ \ \ \ \ \ \ \ $
+\sum\limits_{0\leq n<l\leq 7}Mv^{2}(ki,kj)Ut_{2n}$ $($where $k=l)$,

$\ \ \ \ \ \ \ \ \ \ \ \ \ \ \ \ \ \ \ \ \ \ \ \ \ \ \ \ \ \ \
(ki=Nu(i+1,j+1),kj=Nu(n+1,l+1),0\leq i<j\leq 7),$

\ \ \ \ $[Ud_{ij},Ut_{3k}]=-\sum\limits_{0\leq n<l\leq 7}Mv(ki,kj)Ut_{3l}$ $($where $n=k)$

\ \ \ \ \ \ \ \ \ \ \ \ \ \ \ \ \ \ \ \ \ \ \ $
+\sum\limits_{0\leq n<l\leq 7}Mv(ki,kj)Ut_{3n}$ $($where $l=k),$

$\ \ \ \ \ \ \ \ \ \ \ \ \ \ \ \ \ \ \ \ \ \ \ \ \ \ \ \ \ \ \
(ki=Nu(i+1,j+1),kj=Nu(n+1,l+1),0\leq i<j\leq 7),$

\ \ \ \ $[Ud_{ij},U\tau _{k}]=0,(0\leq i<j\leq 7,k=1,2).$

\bigskip

\emph{Proof. \ }We have the above Lie bracket operations with
calculations using Maxima.\ \ \ \ \emph{Q.E.D.}

\bigskip

\emph{Lemma 15.4. }\ We have the following Lie bracket operations.

\ \ \ \ $[Um_{1i},Ut_{1i}]=U\tau _{2},(0\leq i\leq 7),$

\ \ \ \ $[Um_{1i},Ut_{1j}]=0,(0\leq i,j\leq 7,i\neq j),$

\ \ \ \ $[Um_{1i},Ut_{2j}]=\frac{1}{2}Sn(i+1,j+1)Ut_{3k},(k=Ca(i+1,j+1),\ 0\leq i,j\leq 7),$

\ \ \ \ $[Um_{1i},Ut_{3j}]=-\frac{1}{2}Sn(j+1,i+1)Ut_{2k},(k=Ca(j+1,i+1),\ 0\leq i,j\leq 7),$

\ \ \ \ $[Um_{2i},Ut_{1j}]=-\frac{1}{2}Sn(j+1,i+1)Ut_{3k},(k=Ca(j+1,i+1),\ 0\leq i,j\leq 7),$

\ \ \ \ $[Um_{2i},Ut_{2i}]=-U\tau _{1},(0\leq i\leq 7),$

\ \ \ \ $[Um_{2i},Ut_{2j}]=0,(0\leq i,j\leq 7,i\neq j),$

\ \ \ \ $[Um_{2i},Ut_{3j}]=\frac{1}{2}Sn(i+1,j+1)Ut_{1k},(k=Ca(i+1,j+1),\ 0\leq i,j\leq 7),$

\ \ \ \ $[Um_{3i},Ut_{1j}]=\frac{1}{2}Sn(i+1,j+1)Ut_{2k},(k=Ca(i+1,j+1),\ 0\leq i,j\leq 7),$

\ \ \ \ $[Um_{3i},Ut_{2j}]=-\frac{1}{2}Sn(j+1,i+1)Ut_{1k},(k=Ca(j+1,i+1),\ 0\leq i,j\leq 7),$

\ \ \ \ $[Um_{3i},Ut_{3i}]=U\tau _{1}-U\tau _{2},(0\leq i\leq 7),$

\ \ \ \ $[Um_{3i},Ut_{3j}]=0,(0\leq i,j\leq 7,i\neq j),$

\ \ \ \ $[U\tau _{1},Um_{1i}]=\frac{1}{2}Ut_{1i},(0\leq i\leq 7),$

\ \ \ \ $[U\tau _{2},Um_{1i}]=Ut_{1i},(0\leq i\leq 7),$

\ \ \ \ $[U\tau _{1},Um_{2i}]=-Ut_{2i},(0\leq i\leq 7),$

\ \ \ \ $[U\tau _{2},Um_{2i}]=-\frac{1}{2}Ut_{2i},(0\leq i\leq 7),$

\ \ \ \ $[U\tau _{1},Um_{3i}]=\frac{1}{2}Ut_{3i},(0\leq i\leq 7),$

\ \ \ \ $[U\tau _{2},Um_{3i}]=-\frac{1}{2}Ut_{3i},(0\leq i\leq 7),$

\bigskip

\emph{Proof. \ }We have the above Lie bracket operations with
calculations using Maxima.\ \ \ \ \emph{Q.E.D.}

\bigskip

\emph{Lemma 15.5.} \ \gr$_{6}^{\C}$ is simple.

\bigskip

\emph{Proof.} \ By \emph{Lemma 15.2, 15.3 }and\emph{\ 15.4}, we have

\ \ \ \ $[Ut_{10},Ud_{0j}]=Ut_{1j},(0\leq j\leq 7),$

\ \ \ \ $[Ut_{10},Um_{10}]=-U\tau _{2},$

\ \ \ \ $[Ut_{10},Um_{2i}]=\frac{1}{2}Sn(1,i+1)Ut_{3i},(0\leq i\leq 7),$

\ \ \ \ $[Ut_{10},Um_{3i}]=-\frac{1}{2}Sn(i+1,1)Ut_{2i},(0\leq i\leq 7),$

\ \ \ \ $[Ut_{20},Um_{20}]=U\tau _{1},$

\ \ \ \ $[Ut_{1i},U\tau _{1}]=-\frac{1}{2}Um_{1i},$

\ \ \ \ $[Ut_{2i},U\tau _{1}]=Um_{2i},$

\ \ \ \ $[Ut_{3i},U\tau _{1}]=-\frac{1}{2}Um_{3i}.$

\noindent
Hence we have

$\{[Rt_{10},x],[[Rt_{10},x],y] \mid x,y\in \ $\gr$_{4}^{\C} \bigoplus $\gR\gt$^{\C}\ \}=$ 
\gr$_{6}^{\C}.$

\noindent
Furthermore,by \emph{Lemma 15.4} $\ [Ut_{1i},Ud_{0i}]=-Ut_{10},[Ut_{2i},Um_{30}]=$

$\frac{1}{2}Sn(i+1,1)Ut_{1i},$
$[Ut_{3i},Um_{20}]=-\frac{1}{2}Sn(1,i+1)Ut_{1i},(0\leq i\leq 7).$

\noindent
So we have

$\{[Rt_{ij},x],[[Rt_{ij},x],y] \mid x,y\in \ $\gr$%
_{4}^{\C} \bigoplus $\gR\gt$^{\C}\ \}=$ \gr$_{6}^{\C},(i=1,2,3,0\leq j\leq 7).$

On the other hand by \emph{Lemma 15.2} $[U\tau _{1},Um_{10}]=\frac{1}{2}%
Ut_{10},[U\tau _{2},Um_{10}]=Ut_{10}.$ 

\noindent
Therefore we have

$\{[R\tau _{k},x],[[R\tau _{k},x],y] \mid x,y\in \ $\gr$%
_{4}^{\C} \bigoplus $\gR\gt$^{\C}\ \}=$ \gr$_{6}^{\C},(0\leq k\leq 7).$

Then \gr$_{6}^{\C}$ is simple. \ \ \ \ \emph{Q.E.D.}

\bigskip

\emph{Lemma 15.6.} \ The Killing form B$_{6}$ of the Lie algebra 
\gr$_{6}^{\C}$ is given by

\ \ \ \ \ \ \ \ \ \ \ \ \ \ \ \ \ $B_{6}(R_{1},R_{2})=\frac{2}{5}tr(R_{1}R_{2}),R_{1},R_{2}\in 
$\gr$_{6}^{\C}.$

\bigskip

\emph{Proof}. \ Since \gr$_{6}^{\C}$ is simple ,there exist $%
\kappa \in \C$ such that

\ \ \ \ \ \ \ \ \ \ \ \ \ \ \ \ \ $B_{6}(R_{1},R_{2})=\kappa tr(R_{1}R_{2}).$

\noindent
To determine this $\kappa ,$let $R=U\tau _{1}.$ $(adR)^{2}$ is calculated as
follows.

\ \ \ \ \ \ \ \ \ \ $[U\tau _{1},[U\tau _{1},Um_{1i}]]=[U\tau _{1},\frac{1}{2}Ut_{1i}]=%
\frac{1}{4}Um_{1i},(0\leq i\leq 7),$

\ \ \ \ \ \ \ \ \ \ $[U\tau _{1},[U\tau _{1},Um_{2i}]]=[U\tau _{1},-Ut_{2i}]=Um_{2i},(0\leq
i\leq 7),$

\ \ \ \ \ \ \ \ \ \ $[U\tau _{1},[U\tau _{1},Um_{3i}]]=[U\tau _{1},\frac{1}{2}Ut_{3i}]=%
\frac{1}{4}Um_{3i},(0\leq i\leq 7),$

\ \ \ \ \ \ \ \ \ \ $[U\tau _{1},[U\tau _{1},Ut_{1i}]]=[U\tau _{1},\frac{1}{2}Um_{1i}]=%
\frac{1}{4}Ut_{1i},(0\leq i\leq 7),$

\ \ \ \ \ \ \ \ \ \ $[U\tau _{1},[U\tau _{1},Ut_{2i}]]=[U\tau _{1},-Um_{2i}]=Ut_{2i},(0\leq
i\leq 7),$

\ \ \ \ \ \ \ \ \ \ $[U\tau _{1},[U\tau _{1},Ut_{3i}]]=[U\tau _{1},\frac{1}{2}Um_{3i}]=%
\frac{1}{4}Ut_{3i},(0\leq i\leq 7),$

\ \ \ \ \ \ \ \ \ \ the others = 0.

\noindent
Hence we have

\ \ \ \ \ \ \ \ \ \ $B_{6}(U\tau _{1},U\tau _{1})=tr(ad(U\tau _{1}))^{2}=(\frac{1}{4}%
+1+\frac{1}{4}+\frac{1}{4}+1+\frac{1}{4})\times 8=24.$

\noindent
On the other hand,by calculate with Maxima we have

\ \ \ \ \ \ \ \ \ \ $tr((U\tau _{1})(U\tau _{1}))=60.$

\noindent
Hence we have $\kappa =\frac{24}{60}=\frac{2}{5}.$ \ \ \ \ \emph{Q.E.D.}

\bigskip

\emph{Lemma 15.7.} The rank of the Lie algebra \gr$_{6}^{\C}$
is 6. The roots of \gr$_{6}^{\C}$ relative to some Cartan
subalgebra of \gr$_{6}^{\C}$ are given by

$\ \ \ \ \ \ \ \ \ \ \ \ \ \ \ \ \ \pm (\lambda _{k}-\lambda _{l}),\pm
(\lambda _{k}+\lambda _{l}),0\leq k<l\leq 3,$

\ \ \ \ \ \ \ \ \ \ \ \ \ \ \ \ $\ \pm \lambda _{k}\pm \frac{1}{2}(\mu
_{1}+2\mu _{2}),0\leq k\leq 3,$

$\ \ \ \ \ \ \ \ \ \ \ \ \ \ \ \ \ \pm \frac{1}{2}(-\lambda _{0}-\lambda
_{1}+\lambda _{2}-\lambda _{3})\pm \frac{1}{2}(-2\mu _{1}-\mu _{2}),$

$\ \ \ \ \ \ \ \ \ \ \ \ \ \ \ \ \ \pm \frac{1}{2}(\lambda _{0}+\lambda
_{1}+\lambda _{2}-\lambda _{3})\pm \frac{1}{2}(-2\mu _{1}-\mu _{2}),$

\ \ \ \ \ \ \ \ \ \ \ \ \ \ \ \ \ $\pm \frac{1}{2}(-\lambda _{0}+\lambda
_{1}+\lambda _{2}+\lambda _{3})\pm \frac{1}{2}(-2\mu _{1}-\mu _{2}),$

$\ \ \ \ \ \ \ \ \ \ \ \ \ \ \ \ \ \pm \frac{1}{2}(\lambda _{0}-\lambda
_{1}+\lambda _{2}+\lambda _{3})\pm \frac{1}{2}(-2\mu _{1}-\mu _{2}),$

\ \ \ \ \ \ \ \ \ \ \ \ \ \ \ \ $\ \pm \frac{1}{2}(\lambda _{0}-\lambda
_{1}+\lambda _{2}-\lambda _{3})\pm \frac{1}{2}(\mu _{1}-\mu _{2}),$

$\ \ \ \ \ \ \ \ \ \ \ \ \ \ \ \ \ \pm \frac{1}{2}(-\lambda _{0}+\lambda
_{1}+\lambda _{2}-\lambda _{3})\pm \frac{1}{2}(\mu _{1}-\mu _{2}),$

\ \ \ \ \ \ \ \ \ \ \ \ \ \ \ \ $\ \pm \frac{1}{2}(\lambda _{0}+\lambda
_{1}+\lambda _{2}+\lambda _{3})\pm \frac{1}{2}(\mu _{1}-\mu _{2}),$

$\ \ \ \ \ \ \ \ \ \ \ \ \ \ \ \ \ \pm \frac{1}{2}(-\lambda _{0}-\lambda
_{1}+\lambda _{2}+\lambda _{3})\pm \frac{1}{2}(\mu _{1}-\mu _{2}).$

\bigskip

\emph{Proof. \ }Let

\ \ \ \ \ \ \gh$=\left\{ h=h_{\delta }+H\in \text{\gr}_{6}^{\C}%
\ \middle| %
\begin{array}{c}
\ h_{\delta }=\sum\limits_{k=0}^{3}\lambda
_{k}H_{k}=\sum\limits_{k=0}^{3}-\lambda _{k}iUd_{k4+k}, \\
H=\mu _{1}U\tau _{1}+\mu_{2}U\tau _{2}, \lambda _{k},\mu _{j},\nu \in \C%
\end{array}%
\right\} ,$

\noindent
then \gh\  is an abelian subalgebra of \gr$_{6}^{\C} \ ($it will
be a Cartan subalgebra of \gr$_{6}^{\C} )$. That \gh\  is
abelian is clear from

\ \ \ \ \ \ \ \ \ \ \ \ \ \ \ \ [$h_{\delta },h_{\delta }^{\prime }]=0,$

$\ \ \ \ \ \ \ \ \ \ \ \ \ \ \ \ [h_{\delta },H]=0$ $($\emph{Lemma 15.3}$),$

$\ \ \ \ \ \ \ \ \ \ \ \ \ \ \ \ [H,H^{\prime} ]=0$ $($\emph{Lemma 15.2}$).$

\ $\ (1)$ The roots $\pm \lambda _{k}\pm \lambda _{l}$ of \gr\gd $($ $\subset $%
\gr$_{4}^{\C}\subset $\gr$_{6}^{\C})$ are also roots of \gr%
$_{6}^{\C}$. Indeed,let $\alpha $ be a root of \gr\gd
\ and S$\in $\gr$_{4}^{\C}\subset $\gr$_{6}^{\C}$ be an associated
root vector. Then

\ \ \ \ \ \ \ \ \ \ \ \ \ \ \ \ $[h,S]=[h_{\delta }+H,S]=[h_{\delta
},S]+[H,S],$

\ \ \ \ \ \ \ \ \ \ \ \ \ \ \ \ \ \ \ \ \ \ \ \ \ $=\alpha (h_{\delta
})S+[H,S],$ $\ ($since $S\in $\gr\gd$,$ so $[H,S]=0),$

\ \ \ \ \ \ \ \ \ \ \ \ \ \ \ \ \ \ \ \ \ \ \ \ \ $=\alpha (h_{\delta
})S=(\pm \lambda _{k}\pm \lambda _{l})S.$ \ 

\noindent
Hence $\pm \lambda _{k}\pm \lambda _{l}$ are roots of \gr$_{6}^{\C}.$

\ \ $(2)$ Let we put $Sa:$

\ \ \ \ $\ \ \ \ \ \ \ \ \ \ \ \ \ \ \ \ \ Sa=\frac{1}{2}%
(Um_{1k}+iUm_{14+k}+Ut_{1k}+iUt_{14+k}),$

\noindent
then we have

\ \ \ \ \ \ \ \ \ \ \ \ \ \ \ \ $[h,Sa]=[h_{\delta },Sa]+[H,Sa]$

\ \ \ \ \ \ \ \ \ \ \ \ \ \ \ \ \ \ \ \ \ \ \ \ \ $=\frac{1}{2}\lambda _{k}$[%
$-iUd_{44+k},Um_{1k}+iUm_{14+k}+Ut_{1k}+iUt_{14+k}]$

\ \ \ \ \ \ \ \ \ \ \ \ \ \ \ \ \ \ \ \ \ \ \ \ \ \ $\ +\frac{1}{2}\mu
_{1}[U\tau _{1},Um_{1k}+iUm_{14+k}+Ut_{1k}+iUt_{14+k}]$

\ \ \ \ \ \ \ \ \ \ \ \ \ \ \ \ \ \ \ \ \ \ \ \ \ $\ \ +\frac{1}{2}\mu
_{2}[U\tau _{2},Um_{1k}+iUm_{14+k}+Ut_{1k}+iUt_{14+k}]$

\ \ \ \ \ \ \ \ \ \ \ \ \ \ \ \ \ \ \ \ \ \ \ \ $=\frac{1}{2}\lambda
_{k}(Um_{1k}+iUm_{14+k}+Ut_{1k}+iUt_{14+k})$

\ \ \ \ \ \ \ \ \ \ \ \ \ \ \ \ \ \ \ \ \ \ \ $\ \ \ \ +\frac{1}{4}\mu
_{1}(Um_{1k}+iUm_{14+k}+Ut_{1k}+iUt_{14+k})$

\ \ \ \ \ \ \ \ \ \ \ \ \ \ \ \ \ \ \ \ \ \ \ \ \ $\ \ +\frac{1}{2}\mu
_{2}(Um_{1k}+iUm_{14+k}+Ut_{1k}+iUt_{14+k})$

\ \ \ \ \ \ \ \ \ \ \ \ \ \ \ \ \ \ \ \ \ \ \ \ \ $=\lambda _{k}S+\frac{1}{2}%
(\mu _{1}+2\mu _{2})Sa=(\lambda _{k}+\frac{1}{2}(\mu _{1}+2\mu _{2}))Sa.$

\noindent
Hence $\lambda _{k}+\frac{1}{2}(\mu _{1}+2\mu _{2})$ are roots of \gr$%
_{6}^{\C}$ and $\frac{1}{2}(Um_{1k}+iUm_{14+k}+Ut_{1k}+iUt_{14+k})$ are their
root vectors.

\noindent
Let we put Sb:

\ \ \ \ \ \ \ \ \ \ \ \ \ \ \ \ \ \ \ \ \ \ $Sb=\frac{1}{2}%
(Um_{1k}-iUm_{14+k}+Ut_{1k}-iUt_{14+k}),$

\noindent
then similarly we have

\ \ \ \ \ \ \ \ \ \ \ \ \ \ \ \ \ $[h,Sb]=(-\lambda _{k}+\frac{1}{2}(\mu
_{1}+2\mu _{2}))Sb.$

\noindent
Hence $-\lambda _{k}+\frac{1}{2}(\mu _{1}+2\mu _{2})$ are roots of \gr%
$_{6}^{\C}$ and $\frac{1}{2}(Um_{1k}-iUm_{14+k}+Ut_{1k}-iUt_{14+k})$ are their
root vectors.

\noindent
Let we put Sc:

\ \ \ \ \ \ \ \ \ \ \ \ \ \ \ \ \ \ \ \ \ \ $Sc=\frac{1}{2}%
(Um_{1k}+iUm_{14+k}-Ut_{1k}-iUt_{14+k}),$

\noindent
then similarly we have

\ \ \ \ \ \ \ \ \ \ \ \ \ \ \ \ \ $[h,Sc]=(\lambda _{k}-\frac{1}{2}(\mu
_{1}+2\mu _{2}))Sc.$

\noindent
Hence $\lambda _{k}-\frac{1}{2}(\mu _{1}+2\mu _{2})$ are roots of \gr$%
_{6}^{\C}$ and $\frac{1}{2}(Um_{1k}+iUm_{14+k}-Ut_{1k}-iUt_{14+k})$ are their
root vectors.

\noindent
Let we put Sd:

\ \ \ \ \ \ \ \ \ \ \ \ \ \ \ \ \ \ \ \ \ \ $Sd=\frac{1}{2}%
(Um_{1k}-iUm_{14+k}-Ut_{1k}+iUt_{14+k}),$

\noindent
then similarly we have

\ \ \ \ \ \ \ \ \ \ \ \ \ \ \ \ \ $[h,Sd]=(-\lambda _{k}-\frac{1}{2}(\mu
_{1}+2\mu _{2}))Sd.$

\noindent
Hence $-\lambda _{k}-\frac{1}{2}(\mu _{1}+2\mu _{2})$ are roots of \gr%
$_{6}^{\C}$ and $\frac{1}{2}(Um_{1k}-iUm_{14+k}-Ut_{1k}+iUt_{14+k})$ are their
root vectors.

\ \ $(3)$ Let we put $Sa_{k}:$

\ \ \ \ $\ \ \ \ \ \ \ \ \ \ \ \ \ \ \ \ Sa_{k}=\frac{1}{2}%
(Um_{2k}+iUm_{24+k}+Ut_{2k}+iUt_{24+k}),0\leq k\leq 3,$

\noindent
then we have

\ \ \ \ \ \ \ \ \ \ \ \ \ \ \ $[h,Sa_{0}]=(\frac{1}{2}(-\lambda _{0}-\lambda
_{1}+\lambda _{2}-\lambda _{3})+\frac{1}{2}(-2\mu _{1}-\mu _{2}))Sa_{0}$

\ \ \ \ \ \ \ \ \ \ \ \ \ \ \ $[h,Sa_{1}]=(\frac{1}{2}(\lambda _{0}+\lambda
_{1}+\lambda _{2}-\lambda _{3})+\frac{1}{2}(-2\mu _{1}-\mu _{2}))Sa_{1},$

\ \ \ \ \ \ \ \ \ \ \ \ \ \ \ $[h,Sa_{2}]=(\frac{1}{2}(-\lambda _{0}+\lambda
_{1}+\lambda _{2}+\lambda _{3})+\frac{1}{2}(-2\mu _{1}-\mu _{2}))Sa_{2},$

\ \ \ \ \ \ \ \ \ \ \ \ \ \ \ $[h,Sa_{3}]=(\frac{1}{2}(\lambda _{0}-\lambda
_{1}+\lambda _{2}+\lambda _{3})+\frac{1}{2}(-2\mu _{1}-\mu _{2}))Sa_{3}.$

\noindent
So we see 

$\frac{1}{2}(-\lambda _{0}-\lambda _{1}+\lambda _{2}-\lambda _{3})+%
\frac{1}{2}(-2\mu _{1}-\mu _{2}),\frac{1}{2}(\lambda _{0}+\lambda
_{1}+\lambda _{2}-\lambda _{3})+\frac{1}{2}(-2\mu _{1}-\mu _{2}),$

$\frac{1}{2}(-\lambda _{0}+\lambda _{1}+\lambda _{2}+\lambda _{3})+\frac{1}{2%
}(-2\mu _{1}-\mu _{2}),\frac{1}{2}(\lambda _{0}-\lambda _{1}+\lambda
_{2}+\lambda _{3})+\frac{1}{2}(-2\mu _{1}-\mu _{2})$ are roots of \gr$%
_{6}^{\C}$ and
that $\frac{1}{2}(Um_{2k}+iUm_{24+k}+Ut_{2k}+iUt_{24+k})$ are associated root
vectors for $0\leq k\leq 3.$

\noindent
Let we put $Sb_{k}:$

\ \ \ \ $\ \ \ \ \ \ \ \ \ \ \ \ \ \ \ \ Sb_{k}=\frac{1}{2}%
(Um_{2k}-iUm_{24+k}+Ut_{2k}-iUt_{24+k}),0\leq k\leq 3,$

\noindent
then we have

\ \ \ \ \ \ \ \ \ \ \ \ \ \ \ $[h,Sb_{0}]=(-\frac{1}{2}(-\lambda
_{0}-\lambda _{1}+\lambda _{2}-\lambda _{3})+\frac{1}{2}(-2\mu _{1}-\mu
_{2}))Sb_{0}$

\ \ \ \ \ \ \ \ \ \ \ \ \ \ \ $[h,Sb_{1}]=(-\frac{1}{2}(\lambda _{0}+\lambda
_{1}+\lambda _{2}-\lambda _{3})+\frac{1}{2}(-2\mu _{1}-\mu _{2}))Sb_{1},$

\ \ \ \ \ \ \ \ \ \ \ \ \ \ \ $[h,Sb_{2}]=(-\frac{1}{2}(-\lambda
_{0}+\lambda _{1}+\lambda _{2}+\lambda _{3})+\frac{1}{2}(-2\mu _{1}-\mu
_{2}))Sb_{2},$

\ \ \ \ \ \ \ \ \ \ \ \ \ \ \ $[h,Sb_{3}]=(-\frac{1}{2}(\lambda _{0}-\lambda
_{1}+\lambda _{2}+\lambda _{3})+\frac{1}{2}(-2\mu _{1}-\mu _{2}))Sb_{3}.$

\noindent
So we see 

$-\frac{1}{2}(-\lambda _{0}-\lambda _{1}+\lambda _{2}-\lambda
_{3})+\frac{1}{2}(-2\mu _{1}-\mu _{2}),-\frac{1}{2}(\lambda _{0}+\lambda
_{1}+\lambda _{2}-\lambda _{3})+\frac{1}{2}(-2\mu _{1}-\mu _{2}),$

$-\frac{1}{2}(-\lambda _{0}+\lambda _{1}+\lambda _{2}+\lambda _{3})+\frac{1}{%
2}(-2\mu _{1}-\mu _{2}),-\frac{1}{2}(\lambda _{0}-\lambda _{1}+\lambda
_{2}+\lambda _{3})+\frac{1}{2}(-2\mu _{1}-\mu _{2})$ are roots of \gr$%
_{6}^{\C}$
and that $\frac{1}{2}(Um_{2k}-iUm_{24+k}+Ut_{2k}-iUt_{24+k})$ are associated
root vectors for $0\leq k\leq 3.$

\noindent
Let we put $Sc_{k}:$

\ \ \ \ $\ \ \ \ \ \ \ \ \ \ \ \ \ \ \ \ Sc_{k}=\frac{1}{2}%
(Um_{2k}+iUm_{24+k}-Ut_{2k}-iUt_{24+k}),0\leq k\leq 3,$

\noindent
then we have

\ \ \ \ \ \ \ \ \ \ \ \ \ \ \ $[h,Sc_{0}]=(\frac{1}{2}(-\lambda _{0}-\lambda
_{1}+\lambda _{2}-\lambda _{3})-\frac{1}{2}(-2\mu _{1}-\mu _{2}))Sc_{0}$

\ \ \ \ \ \ \ \ \ \ \ \ \ \ \ $[h,Sc_{1}]=(\frac{1}{2}(\lambda _{0}+\lambda
_{1}+\lambda _{2}-\lambda _{3})-\frac{1}{2}(-2\mu _{1}-\mu _{2}))Sc_{1},$

\ \ \ \ \ \ \ \ \ \ \ \ \ \ \ $[h,Sc_{2}]=(\frac{1}{2}(-\lambda _{0}+\lambda
_{1}+\lambda _{2}+\lambda _{3})-\frac{1}{2}(-2\mu _{1}-\mu _{2}))Sc_{2},$

\ \ \ \ \ \ \ \ \ \ \ \ \ \ \ $[h,Sc_{3}]=(\frac{1}{2}(\lambda _{0}-\lambda
_{1}+\lambda _{2}+\lambda _{3})-\frac{1}{2}(-2\mu _{1}-\mu _{2}))Sc_{3}.$

\noindent
So we see 

$\frac{1}{2}(-\lambda _{0}-\lambda _{1}+\lambda _{2}-\lambda _{3})-%
\frac{1}{2}(-2\mu _{1}-\mu _{2}),\frac{1}{2}(\lambda _{0}+\lambda
_{1}+\lambda _{2}-\lambda _{3})-\frac{1}{2}(-2\mu _{1}-\mu _{2}),$

$\frac{1}{2}(-\lambda _{0}+\lambda _{1}+\lambda _{2}+\lambda _{3})-\frac{1}{2%
}(-2\mu _{1}-\mu _{2}),\frac{1}{2}(\lambda _{0}-\lambda _{1}+\lambda
_{2}+\lambda _{3})-\frac{1}{2}(-2\mu _{1}-\mu _{2})$ are roots of \gr$%
_{6}^{\C}$
and that $\frac{1}{2}(Um_{2k}+iUm_{24+k}-Ut_{2k}-iUt_{24+k})$ are associated
root vectors for $0\leq k\leq 3.$

\noindent
Let we put $Sd_{k}:$

\ \ \ \ $\ \ \ \ \ \ \ \ \ \ \ \ \ \ \ \ Sd_{k}=\frac{1}{2}%
(Um_{2k}-iUm_{24+k}-Ut_{2k}+iUt_{24+k}),0\leq k\leq 3,$

\noindent
then we have

\ \ \ \ \ \ \ \ \ \ \ \ \ \ \ $[h,Sd_{0}]=(-\frac{1}{2}(-\lambda
_{0}-\lambda _{1}+\lambda _{2}-\lambda _{3})-\frac{1}{2}(-2\mu _{1}-\mu
_{2}))Sd_{0}$

\ \ \ \ \ \ \ \ \ \ \ \ \ \ \ $[h,Sd_{1}]=(-\frac{1}{2}(\lambda _{0}+\lambda
_{1}+\lambda _{2}-\lambda _{3})-\frac{1}{2}(-2\mu _{1}-\mu _{2}))Sd_{1},$

\ \ \ \ \ \ \ \ \ \ \ \ \ \ \ $[h,Sd_{2}]=(-\frac{1}{2}(-\lambda
_{0}+\lambda _{1}+\lambda _{2}+\lambda _{3})-\frac{1}{2}(-2\mu _{1}-\mu
_{2}))Sd_{2},$

\ \ \ \ \ \ \ \ \ \ \ \ \ \ \ $[h,Sd_{3}]=(-\frac{1}{2}(\lambda _{0}-\lambda
_{1}+\lambda _{2}+\lambda _{3})-\frac{1}{2}(-2\mu _{1}-\mu _{2}))Sd_{3}.$

\noindent
So we see 

$-\frac{1}{2}(-\lambda _{0}-\lambda _{1}+\lambda _{2}-\lambda
_{3})-\frac{1}{2}(-2\mu _{1}-\mu _{2}),-\frac{1}{2}(\lambda _{0}+\lambda
_{1}+\lambda _{2}-\lambda _{3})-\frac{1}{2}(-2\mu _{1}-\mu _{2}),$

$-\frac{1}{2}(-\lambda _{0}+\lambda _{1}+\lambda _{2}+\lambda _{3})-\frac{1}{%
2}(-2\mu _{1}-\mu _{2}),-\frac{1}{2}(\lambda _{0}-\lambda _{1}+\lambda
_{2}+\lambda _{3})-\frac{1}{2}(-2\mu _{1}-\mu _{2})$ are roots of \gr$%
_{6}^{\C}$
and that $\frac{1}{2}(Um_{2k}-iUm_{24+k}-Ut_{2k}+iUt_{24+k})$ are associated
root vectors for $0\leq k\leq 3.$

\ \ $(4)$ Let we put $Sa_{k}:$

\ \ \ \ $\ \ \ \ \ \ \ \ \ \ \ \ \ \ \ \ \ \ \ Sa_{k}=\frac{1}{2}%
(Um_{3k}+iUm_{34+k}+Ut_{3k}+iUt_{34+k}),0\leq k\leq 3,$

\noindent
then we have

\ \ \ \ \ \ \ \ \ \ \ \ \ \ \ $[h,Sa_{0}]=(\frac{1}{2}(-\lambda _{0}+\lambda
_{1}-\lambda _{2}+\lambda _{3})+\frac{1}{2}(\mu _{1}-\mu _{2}))Sa_{0},$

\ \ \ \ \ \ \ \ \ \ \ \ \ \ \ $[h,Sa_{1}]=(\frac{1}{2}(-\lambda _{0}+\lambda
_{1}+\lambda _{2}-\lambda _{3})+\frac{1}{2}(\mu _{1}-\mu _{2}))Sa_{1},$

\ \ \ \ \ \ \ \ \ \ \ \ \ \ \ $[h,Sa_{2}]=(\frac{1}{2}(\lambda _{0}+\lambda
_{1}+\lambda _{2}+\lambda _{3})+\frac{1}{2}(\mu _{1}-\mu _{2}))Sa_{2},$

\ \ \ \ \ \ \ \ \ \ \ \ \ \ \ $[h,Sa_{3}]=(\frac{1}{2}(-\lambda _{0}-\lambda
_{1}+\lambda _{2}+\lambda _{3})+\frac{1}{2}(\mu _{1}-\mu _{2}))Sa_{3}.$

\noindent
So we see 

$\frac{1}{2}(-\lambda _{0}+\lambda _{1}-\lambda _{2}+\lambda _{3})+%
\frac{1}{2}(\mu _{1}-\mu _{2}),\frac{1}{2}(-\lambda _{0}+\lambda
_{1}+\lambda _{2}-\lambda _{3})+\frac{1}{2}(\mu _{1}-\mu _{2}),$

$\frac{1}{2}(\lambda _{0}+\lambda _{1}+\lambda _{2}+\lambda _{3})+\frac{1}{2}%
(\mu _{1}-\mu _{2}),\frac{1}{2}(-\lambda _{0}-\lambda _{1}+\lambda
_{2}+\lambda _{3})+\frac{1}{2}(\mu _{1}-\mu _{2})$ are roots of \gr$%
_{6}^{\C}$ and
that $\frac{1}{2}(Um_{3k}+iUm_{34+k}+Ut_{3k}+iUt_{34+k})$ are associated root
vectors for $0\leq k\leq 3.$

\noindent
Let we put $Sb_{k}:$

\ \ \ \ $\ \ \ \ \ \ \ \ \ \ \ \ \ \ \ \ \ \ \ Sb_{k}=\frac{1}{2}%
(Um_{3k}-iUm_{34+k}+Ut_{3k}-iUt_{34+k}),0\leq k\leq 3,$

\noindent
then we have

\ \ \ \ \ \ \ \ \ \ \ \ \ \ \ $[h,Sb_{0}]=(-\frac{1}{2}(-\lambda
_{0}+\lambda _{1}-\lambda _{2}+\lambda _{3})+\frac{1}{2}(\mu _{1}-\mu
_{2}))Sb_{0},$

\ \ \ \ \ \ \ \ \ \ \ \ \ \ \ $[h,Sb_{1}]=(-\frac{1}{2}(-\lambda
_{0}+\lambda _{1}+\lambda _{2}-\lambda _{3})+\frac{1}{2}(\mu _{1}-\mu
_{2}))Sb_{1},$

\ \ \ \ \ \ \ \ \ \ \ \ \ \ \ $[h,Sb_{2}]=(-\frac{1}{2}(\lambda _{0}+\lambda
_{1}+\lambda _{2}+\lambda _{3})+\frac{1}{2}(\mu _{1}-\mu _{2}))Sb_{2},$

\ \ \ \ \ \ \ \ \ \ \ \ \ \ \ $[h,Sb_{3}]=(-\frac{1}{2}(-\lambda
_{0}-\lambda _{1}+\lambda _{2}+\lambda _{3})+\frac{1}{2}(\mu _{1}-\mu
_{2}))Sb_{3}.$

\noindent
So we see 

$-\frac{1}{2}(-\lambda _{0}+\lambda _{1}-\lambda _{2}+\lambda
_{3})+\frac{1}{2}(\mu _{1}-\mu _{2}),-\frac{1}{2}(-\lambda _{0}+\lambda
_{1}+\lambda _{2}-\lambda _{3})+\frac{1}{2}(\mu _{1}-\mu _{2}),$

$-\frac{1}{2}(\lambda _{0}+\lambda _{1}+\lambda _{2}+\lambda _{3})+\frac{1}{2%
}(\mu _{1}-\mu _{2}),-\frac{1}{2}(-\lambda _{0}-\lambda _{1}+\lambda
_{2}+\lambda _{3})+\frac{1}{2}(\mu _{1}-\mu _{2})$ are roots of \gr$%
_{6}^{\C}$ and
that $\frac{1}{2}(Um_{3k}-iUm_{34+k}+Ut_{3k}-iUt_{34+k})$ are associated root
vectors for $0\leq k\leq 3.$

\noindent
Let we put $Sc_{k}:$

\ \ \ \ $\ \ \ \ \ \ \ \ \ \ \ \ \ \ \ \ \ \ \ Sc_{k}=\frac{1}{2}%
(Um_{3k}+iUm_{34+k}-Ut_{3k}-iUt_{34+k}),0\leq k\leq 3,$

\noindent
then we have

\ \ \ \ \ \ \ \ \ \ \ \ \ \ \ $[h,Sc_{0}]=(\frac{1}{2}(-\lambda _{0}+\lambda
_{1}-\lambda _{2}+\lambda _{3})-\frac{1}{2}(\mu _{1}-\mu _{2}))Sc_{0},$

\ \ \ \ \ \ \ \ \ \ \ \ \ \ \ $[h,Sc_{1}]=(\frac{1}{2}(-\lambda _{0}+\lambda
_{1}+\lambda _{2}-\lambda _{3})-\frac{1}{2}(\mu _{1}-\mu _{2}))Sc_{1},$

\ \ \ \ \ \ \ \ \ \ \ \ \ \ \ $[h,Sc_{2}]=(\frac{1}{2}(\lambda _{0}+\lambda
_{1}+\lambda _{2}+\lambda _{3})-\frac{1}{2}(\mu _{1}-\mu _{2}))Sc_{2},$

\ \ \ \ \ \ \ \ \ \ \ \ \ \ \ $[h,Sc_{3}]=(\frac{1}{2}(-\lambda _{0}-\lambda
_{1}+\lambda _{2}+\lambda _{3})-\frac{1}{2}(\mu _{1}-\mu _{2}))Sc_{3}.$

\noindent
In the same way 

$\frac{1}{2}(-\lambda _{0}+\lambda _{1}-\lambda _{2}+\lambda
_{3})-\frac{1}{2}(\mu _{1}-\mu _{2}),\frac{1}{2}(-\lambda _{0}+\lambda
_{1}+\lambda _{2}-\lambda _{3})-\frac{1}{2}(\mu _{1}-\mu _{2}),$

$\frac{1}{2}(\lambda _{0}+\lambda _{1}+\lambda _{2}+\lambda _{3})-\frac{1}{2}%
(\mu _{1}-\mu _{2}),\frac{1}{2}(-\lambda _{0}-\lambda _{1}+\lambda
_{2}+\lambda _{3})-\frac{1}{2}(\mu _{1}-\mu _{2})$ 

\noindent
are roots of \gr$_{6}^{\C}$ and
that $\frac{1}{2}(Um_{3k}+iUm_{34+k}-Ut_{3k}-iUt_{34+k})$ are associated root
vectors for $0\leq k\leq 3.$

\noindent
Let we put $Sd_{k}:$

\ \ \ \ $\ \ \ \ \ \ \ \ \ \ \ \ \ \ \ \ \ \ \ Sd_{k}=\frac{1}{2}%
(Um_{3k}-iUm_{34+k}-Ut_{3k}+iUt_{34+k}),0\leq k\leq 3,$

\noindent
then we have

\ \ \ \ \ \ \ \ \ \ \ \ \ \ \ $[h,Sd_{0}]=(-\frac{1}{2}(-\lambda
_{0}+\lambda _{1}-\lambda _{2}+\lambda _{3})-\frac{1}{2}(\mu _{1}-\mu
_{2}))Sd_{0},$

\ \ \ \ \ \ \ \ \ \ \ \ \ \ \ $[h,Sd_{1}]=(-\frac{1}{2}(-\lambda
_{0}+\lambda _{1}+\lambda _{2}-\lambda _{3})-\frac{1}{2}(\mu _{1}-\mu
_{2}))Sd_{1},$

\ \ \ \ \ \ \ \ \ \ \ \ \ \ \ $[h,Sd_{2}]=(-\frac{1}{2}(\lambda _{0}+\lambda
_{1}+\lambda _{2}+\lambda _{3})-\frac{1}{2}(\mu _{1}-\mu _{2}))Sd_{2},$

\ \ \ \ \ \ \ \ \ \ \ \ \ \ \ $[h,Sd_{3}]=(-\frac{1}{2}(-\lambda
_{0}-\lambda _{1}+\lambda _{2}+\lambda _{3})-\frac{1}{2}(\mu _{1}-\mu
_{2}))Sd_{3}.$

\noindent
So we see 

$-\frac{1}{2}(-\lambda _{0}+\lambda _{1}-\lambda _{2}+\lambda
_{3})-\frac{1}{2}(\mu _{1}-\mu _{2}),-\frac{1}{2}(-\lambda _{0}+\lambda
_{1}+\lambda _{2}-\lambda _{3})-\frac{1}{2}(\mu _{1}-\mu _{2}),$

$-\frac{1}{2}(\lambda _{0}+\lambda _{1}+\lambda _{2}+\lambda _{3})-\frac{1}{2%
}(\mu _{1}-\mu _{2}),-\frac{1}{2}(-\lambda _{0}-\lambda _{1}+\lambda
_{2}+\lambda _{3})-\frac{1}{2}(\mu _{1}-\mu _{2})$ are roots of \gr$%
_{6}^{\C}$ and
that $\frac{1}{2}(Um_{3k}-iUm_{34+k}-Ut_{3k}+iUt_{34+k})$ are associated root
vectors for $0\leq k\leq 3.$

By direct calculations using Maxima, we also have $(1),(2),(3)$ and $(4)$. \ \ \ \ \emph{Q.E.D.}

\bigskip

\emph{Theorem 15.8.} \ In the root system of \emph{Lemma 15.7} ,

\ \ \ \ \ \ \ $\alpha _{1}=\lambda _{0}-\lambda _{1},\alpha _{2}=\lambda
_{1}-\lambda _{2},\alpha _{3}=\lambda _{2}-\lambda _{3},$

$\ \ \ \ \ \ \ \alpha _{4}=\lambda _{3}+\frac{1}{2}(\mu _{1}+2\mu _{2}),$

\ \ \ \ \ \ \ $\alpha _{5}=\frac{1}{2}(-\lambda _{0}-\lambda _{1}-\lambda
_{2}+\lambda _{3})+\frac{1}{2}(-2\mu _{1}-\mu _{2}),$

\ \ \ \ \ \ \ $\alpha _{6}=\frac{1}{2}(\lambda _{0}+\lambda _{1}+\lambda
_{2}+\lambda _{3})+\frac{1}{2}(\mu _{1}-\mu _{2}),$

\noindent
is a fundamental root sysyten of the Lie algebra \gr$_{6}^{\C}$ and

$\ \ \ \ \ \ \ \mu =\alpha _{1}+2\alpha _{2}+3\alpha _{3}+2\alpha
_{4}+2\alpha _{5}+\alpha _{6}$

\noindent
is the highest root. The Dynkin diagram and the extended Dynkin diagram 
of \gr$_{6}^{\C}$ are respectively given by

\begin{figure}[H]
\centering
\includegraphics[width=10cm, height=3cm]{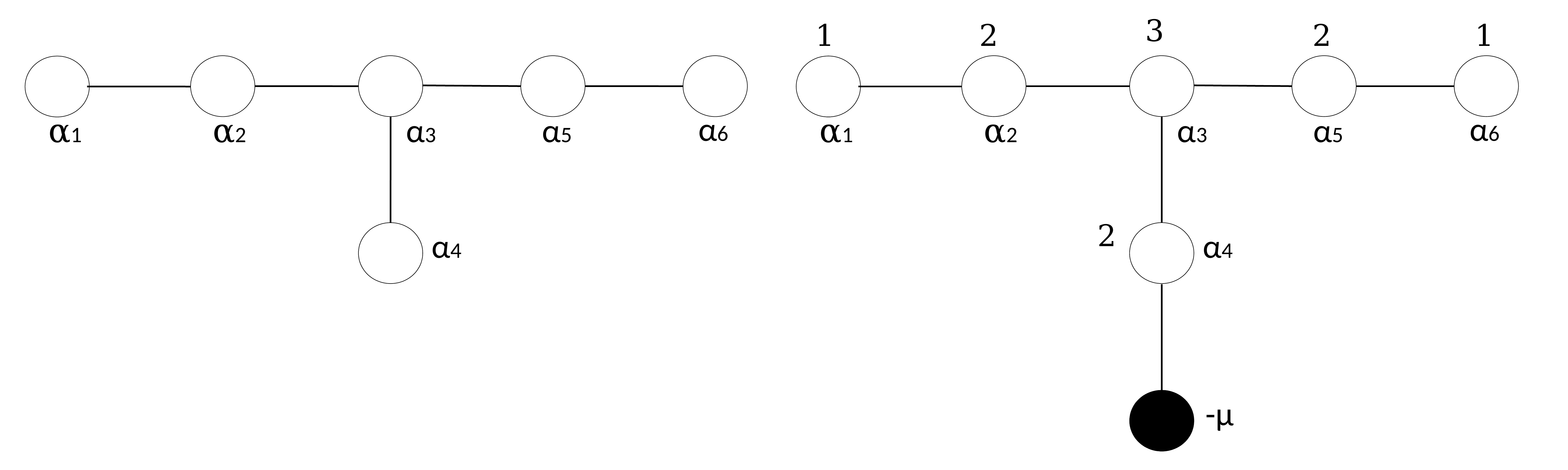}
\end{figure}

\emph{Proof}. \ In the following, the notation $n_{1}$ $n_{2}$ $%
n_{3} $ $n_{4}$ $n_{5}$ $n_{6}$ denotes the root

$n_{1}\alpha _{1}+n_{2}\alpha _{2}+n_{3}\alpha _{3}+n_{4}\alpha
_{4}+n_{5}\alpha _{5}+n_{6}\alpha _{6}.$

Now,all positive roots of \gr$_{6}^{\C}$ are represented by

\begin{flushright}
$\ \ \ \ \ \ \ \ \ \lambda _{0}-\lambda _{1}=%
\begin{array}{llllll}
1 & 0 & 0 & 0 & 0 & 0%
\end{array}%
,\lambda _{0}+\lambda _{1}=%
\begin{array}{llllll}
1 & 2 & 2 & 1 & 1 & 1%
\end{array}%
,$

$\ \ \ \ \ \ \ \ \ \lambda _{0}-\lambda _{2}=%
\begin{array}{llllll}
1 & 1 & 0 & 0 & 0 & 0%
\end{array}%
,\lambda _{0}+\lambda _{2}=%
\begin{array}{llllll}
1 & 1 & 2 & 1 & 1 & 1%
\end{array}%
,$

$\ \ \ \ \ \ \ \ \ \lambda _{0}-\lambda _{3}=%
\begin{array}{llllll}
1 & 1 & 1 & 0 & 0 & 0%
\end{array}%
,\lambda _{0}+\lambda _{3}=%
\begin{array}{llllll}
1 & 1 & 1 & 1 & 1 & 1%
\end{array}%
,$

$\ \ \ \ \ \ \ \ \ \lambda _{1}-\lambda _{2}=%
\begin{array}{llllll}
0 & 1 & 0 & 0 & 0 & 0%
\end{array}%
,\lambda _{1}+\lambda _{2}=%
\begin{array}{llllll}
0 & 1 & 2 & 1 & 1 & 1%
\end{array}%
,$

$\ \ \ \ \ \ \ \ \ \lambda _{1}-\lambda _{3}=%
\begin{array}{llllll}
0 & 1 & 1 & 0 & 0 & 0%
\end{array}%
,\lambda _{1}+\lambda _{3}=%
\begin{array}{llllll}
0 & 1 & 1 & 1 & 1 & 1%
\end{array}%
,$

$\ \ \ \ \ \ \ \ \ \lambda _{2}-\lambda _{3}=%
\begin{array}{llllll}
0 & 0 & 1 & 0 & 0 & 0%
\end{array}%
,\lambda _{2}+\lambda _{3}=%
\begin{array}{llllll}
0 & 0 & 1 & 1 & 1 & 1%
\end{array}%
,$

$\ \ \ \ \ \ \ \ \ \lambda _{0}+\frac{1}{2}(\mu _{1}+2\mu _{2})=%
\begin{array}{llllll}
1 & 1 & 1 & 1 & 0 & 0%
\end{array}%
,$

$\ \ \ \ \ \ \ \ \ \lambda _{1}+\frac{1}{2}(\mu _{1}+2\mu _{2})=%
\begin{array}{llllll}
0 & 1 & 1 & 1 & 0 & 0%
\end{array}%
,$

$\ \ \ \ \ \ \ \ \ \lambda _{2}+\frac{1}{2}(\mu _{1}+2\mu _{2})=%
\begin{array}{llllll}
0 & 0 & 1 & 1 & 0 & 0%
\end{array}%
,$

$\ \ \ \ \ \ \ \ \ \lambda _{3}+\frac{1}{2}(\mu _{1}+2\mu _{2})=%
\begin{array}{llllll}
0 & 0 & 0 & 1 & 0 & 0%
\end{array}%
,$

$\ \ \ \ \ \ \ \ \ \lambda _{0}-\frac{1}{2}(\mu _{1}+2\mu _{2})=%
\begin{array}{llllll}
1 & 1 & 1 & 0 & 1 & 1%
\end{array}%
,$

$\ \ \ \ \ \ \ \ \ \lambda _{1}-\frac{1}{2}(\mu _{1}+2\mu _{2})=%
\begin{array}{llllll}
0 & 1 & 1 & 0 & 1 & 1%
\end{array}%
,$

$\ \ \ \ \ \ \ \ \ \lambda _{2}-\frac{1}{2}(\mu _{1}+2\mu _{2})=%
\begin{array}{llllll}
0 & 0 & 1 & 0 & 1 & 1%
\end{array}%
,$

$\ \ \ \ \ \ \ \ \ \lambda _{3}-\frac{1}{2}(\mu _{1}+2\mu _{2})=%
\begin{array}{llllll}
0 & 0 & 0 & 0 & 1 & 1%
\end{array}%
,$

$\ \ \ \ \ \ \ \frac{1}{2}(-\lambda _{0}-\lambda _{1}+\lambda _{2}-\lambda
_{3})+\frac{1}{2}(-2\mu _{1}-\mu _{2})=%
\begin{array}{llllll}
0 & 0 & 1 & 0 & 1 & 0%
\end{array}%
,$

$\ \ \ \ \ \ \ \ \ \frac{1}{2}(\lambda _{0}+\lambda _{1}+\lambda
_{2}-\lambda _{3})+\frac{1}{2}(-2\mu _{1}-\mu _{2})=%
\begin{array}{llllll}
1 & 2 & 3 & 1 & 2 & 1%
\end{array}%
,$

$\ \ \ \ \ \ \ \frac{1}{2}(-\lambda _{0}+\lambda _{1}+\lambda _{2}+\lambda
_{3})+\frac{1}{2}(-2\mu _{1}-\mu _{2})=%
\begin{array}{llllll}
0 & 1 & 2 & 1 & 2 & 1%
\end{array}%
,$

$\ \ \ \ \ \ \ \ \ \frac{1}{2}(\lambda _{0}-\lambda _{1}+\lambda
_{2}+\lambda _{3})+\frac{1}{2}(-2\mu _{1}-\mu _{2})=%
\begin{array}{llllll}
1 & 1 & 2 & 1 & 2 & 1%
\end{array}%
,$

$\ \ \ \ \ \ \ \ \ \frac{1}{2}(\lambda _{0}+\lambda _{1}-\lambda
_{2}+\lambda _{3})+\frac{1}{2}(-2\mu _{1}-\mu _{2})=%
\begin{array}{llllll}
1 & 2 & 2 & 1 & 2 & 1%
\end{array}%
,$

$\ \ \ \ \ \ \ \frac{1}{2}(-\lambda _{0}-\lambda _{1}-\lambda _{2}+\lambda
_{3})+\frac{1}{2}(-2\mu _{1}-\mu _{2})=%
\begin{array}{llllll}
0 & 0 & 0 & 0 & 1 & 0%
\end{array}%
,$

$\ \ \ \ \ \ \ \ \frac{1}{2}(\lambda _{0}-\lambda _{1}-\lambda _{2}-\lambda
_{3})+\frac{1}{2}(-2\mu _{1}-\mu _{2})=%
\begin{array}{llllll}
1 & 1 & 1 & 0 & 1 & 0%
\end{array}%
,$

$\ \ \ \ \ \ \ \frac{1}{2}(-\lambda _{0}+\lambda _{1}-\lambda _{2}-\lambda
_{3})+\frac{1}{2}(-2\mu _{1}-\mu _{2})=%
\begin{array}{llllll}
0 & 1 & 1 & 0 & 1 & 0%
\end{array}%
,$

$\ \ \ \ \ \ \ \ \ \ \ \ \frac{1}{2}(\lambda _{0}-\lambda _{1}+\lambda
_{2}-\lambda _{3})-\frac{1}{2}(\mu _{1}-\mu _{2})=%
\begin{array}{llllll}
1 & 1 & 2 & 1 & 1 & 0%
\end{array}%
,$

$\ \ \ \ \ \ \ \ \ \ \ \ \frac{1}{2}(\lambda _{0}-\lambda _{1}-\lambda
_{2}+\lambda _{3})-\frac{1}{2}(\mu _{1}-\mu _{2})=%
\begin{array}{llllll}
1 & 1 & 1 & 1 & 1 & 0%
\end{array}%
,$

$\ \ \ \ \ \ \ \ \ \ \ \ \frac{1}{2}(\lambda _{0}+\lambda _{1}+\lambda
_{2}+\lambda _{3})+\frac{1}{2}(\mu _{1}-\mu _{2})=%
\begin{array}{llllll}
0 & 0 & 0 & 0 & 0 & 1%
\end{array}%
,$

$\ \ \ \ \ \ \ \ \ \ \ \ \frac{1}{2}(\lambda _{0}+\lambda _{1}-\lambda
_{2}-\lambda _{3})-\frac{1}{2}(\mu _{1}-\mu _{2})=%
\begin{array}{llllll}
1 & 2 & 2 & 1 & 1 & 0%
\end{array}%
,$

$\ \ \ \ \ \ \ \ \ \frac{1}{2}(-\lambda _{0}+\lambda _{1}-\lambda
_{2}+\lambda _{3})-\frac{1}{2}(\mu _{1}-\mu _{2})=%
\begin{array}{llllll}
0 & 1 & 1 & 1 & 1 & 0%
\end{array}%
,$

$\ \ \ \ \ \ \ \ \ \ \frac{1}{2}(-\lambda _{0}+\lambda _{1}+\lambda
_{2}-\lambda _{3})-\frac{1}{2}(\mu _{1}-\mu _{2})=%
\begin{array}{llllll}
0 & 1 & 2 & 1 & 1 & 0%
\end{array}%
,$

$\ \ \ \ \ \ \ \ \ \ \ \ \frac{1}{2}(\lambda _{0}+\lambda _{1}+\lambda
_{2}+\lambda _{3})-\frac{1}{2}(\mu _{1}-\mu _{2})=%
\begin{array}{llllll}
1 & 2 & 3 & 2 & 2 & 1%
\end{array}%
,$

$\ \ \ \ \ \ \ \ \ \ \frac{1}{2}(-\lambda _{0}-\lambda _{1}+\lambda
_{2}+\lambda _{3})-\frac{1}{2}(\mu _{1}-\mu _{2})=%
\begin{array}{llllll}
0 & 0 & 1 & 1 & 1 & 0%
\end{array}%
,$
\end{flushright}

\noindent
Hence $\Pi =\{\alpha _{1},\alpha _{2},\alpha _{3},\alpha _{4},\alpha
_{5},\alpha _{6}\}$ is a fundamental root system of \gr$_{6}^{\C}$ and

\noindent
$\mu =\alpha _{1}+2\alpha _{2}+3\alpha _{3}+2\alpha _{4}$+2$\alpha
_{5}+\alpha _{6}$ is the highest root. The real part of \gh$_{R}$ of 
\gh \ is

\ \ \ \ \ \ \gh$_{R}=\left\{ h=h_{\delta }+H\in \text{\gr}_{6}^{\C}%
\ \middle| %
\begin{array}{c}
\ h_{\delta }=\sum\limits_{k=0}^{3}\lambda
_{k}H_{k}=\sum\limits_{k=0}^{3}-\lambda _{k}iUd_{k4+k}, \\
H=\mu _{1}U\tau _{1}+\mu_{2}U\tau _{2}, \lambda _{k},\mu _{j},\nu \in \R%
\end{array}%
\right\} ,$

The Killing form $B_{6}$ of \gr$_{6}^{\C}$ is $B_{6}(R_{1},R_{2})=\frac{2}{5}%
tr(R_{1}R_{2})$ (\emph{Lemma 15.6} ), so that

$\ \ \ \ \ \ \ \ \ \ \ \ B_{6}(h,h^{\prime
})=12(2\sum\limits_{k=0}^{3}\lambda _{k}\lambda _{k}^{\prime }+2\mu _{1}\mu
_{1}^{\prime }+\mu _{1}\mu _{2}^{\prime }+\mu _{2}\mu _{1}^{\prime }+2\mu
_{2}\mu _{2}^{\prime }),$

\noindent
for $h=\sum\limits_{k=0}^{3}\lambda _{k}H_{k}+\mu _{1}U\tau _{1}+\mu _{2}u\tau
2,h^{\prime }=\sum\limits_{k=0}^{3}\lambda _{k}^{\prime }H_{k}+\mu
_{1}^{\prime }U\tau _{1}+\mu _{2}^{\prime }U\tau _{2}\in $\gh$_{R}.$

\noindent
Indeed,by calculate with Maxima we have

$\ \ \ \ \ B_{6}(h,h^{\prime })=\frac{2}{5}tr(hh^{\prime })=%
\frac{2}{5}(30(2\sum\limits_{k=0}^{3}\lambda _{k}\lambda _{k}^{\prime }+2\mu
_{1}\mu _{1}^{\prime }+\mu _{1}\mu _{2}^{\prime }+\mu _{2}\mu _{1}^{\prime
}+2\mu _{2}\mu _{2}^{\prime }))$

$\ \ \ \ \ \ \ \ \ \ \ \ \ \ \ \  \ 
=12(2\sum\limits_{k=0}^{3}\lambda _{k}\lambda _{k}^{\prime }+\mu _{1}(2\mu
_{1}^{\prime }+\mu _{2}^{\prime })+\mu _{2}(\mu _{1}^{\prime }+2\mu
_{2}^{\prime })).$

Now,the canonical elements H$_{\alpha _{i}}\in $\gh$_{R}$ corresponding
to $\alpha _{i}(B_{6}(H_{\alpha _{i}},H)=\alpha _{i}(H),H\in $\gh$_{R})$
are determind as follows.

$\ \ \ \ \ \ \ \ \ \ \ \ \ \ \ \ \ \ H_{\alpha 1}=\frac{1}{24}%
(H_{0}-H_{1}),H_{\alpha 2}=\frac{1}{24}(H_{1}-H_{2}),H_{\alpha 3}=\frac{1}{24%
}(H_{2}-H_{3}),$

$\ \ \ \ \ \ \ \ \ \ \ \ \ \ \ \ \ \ H_{\alpha 4}=\frac{1}{24}(H_{3}+u\tau
2),$

\ \ \ \ \ \ \ \ \ \ \ \ \ \ \ \ \ \ $H_{\alpha 5}=\frac{1}{48}%
((-H_{0}-H_{1}-H_{2}+H_{3})-2U\tau _{1}),$

\ \ \ \ \ \ \ \ \ \ \ \ \ \ \ \ \ \ $H_{\alpha 6}=\frac{1}{48}%
((H_{0}+H_{1}+H_{2}+H_{3})+2U\tau _{1}-2U\tau _{2}).$

\noindent
Therefore,by calculate with Maxima we have

\ \ \ \ \ \ \ \ \ \ \ \ \ \ \ $(\alpha _{1},\alpha _{1})=B_{6}(H_{\alpha
1},H_{\alpha 1})=\frac{1}{12},$

\ \ \ \ \ \ \ \ \ \ \ \ \ \ $\ (\alpha _{i},\alpha _{i})=\frac{1}{12}%
,(i=2,3,4,5,6),$

\ \ \ \ \ \ \ \ \ \ \ \ \ \ $\ (\alpha _{1},\alpha _{2})=(\alpha _{2},\alpha
_{3})=(\alpha _{3},\alpha _{4})=(\alpha _{3},\alpha _{5})=(\alpha
_{5},\alpha _{6})=-\frac{1}{24},$

\ \ \ \ \ \ \ \ \ \ \ \ \ \ $\ (\alpha _{i},\alpha _{j})=0$, otherwise,

\ \ \ \ \ \ \ \ \ \ \ \ \ \ $\ (-\mu ,-\mu )=\frac{1}{12},(-\mu ,\alpha
_{4})=-\frac{1}{24},(-\mu ,\alpha _{i})=0,(i=1,2,3,5,6),$

\noindent
using them,we can draw the Dynkin diagram and the extended Dynkin 
diagram of  \gr$_{6}^{\C}.$\ \ \ \ \emph{Q.E.D.}

\bigskip

\emph{Corollary 15.9.} For a $248\times 248$ matrix $X$,let $X|_{6}$ be the matrix in which the $78\times 78$
elements in the upper left corner are clipped from  $X$. Furthermore, let \gr\gd$|_{6}=\{Rd|_{6} \mid Rd\in $\gr\gd$\}$ , 
\gR\gm$|_{6}=\{Rm|_{6} \mid Rm\in $\gR\gm$\}$ and \gR\gt$|_{6}=\{Rt|_{6} \mid Rt\in $\gR\gt$\}$.
\emph{Theorem 15.8} holds for \gr$_{6}^{\C}|_{6}=$\gr\gd$^{\C}|_{6}\oplus $\gR\gm$^{\C}|_{6}\oplus $\gR\gt$^{\C}|_{6}$ as well. However, the
Killing form $B_{6}(R_{1},R_{2})=tr(R_{1}R_{2}) \ \ (R_{1},R_{2} \in $\gr$_{6}^{\C}|_{6}).$

\bigskip

\emph{Proof} We have the above with calculations using Maxima.\ \ \ \ \emph{Q.E.D.}

\bigskip

\emph{Definition 15.10.} \ We define the imaginary part of \gR\gt$^{\C}$ by \ge\gt,

\ge\gt$^{0}=\{i\tau_{1}(U\tau _{1}-U\tau_{2}),i\tau_{2}(U\tau_{1}+U\tau _{2})\in $\gR\gt$^{\C} \mid \tau _{1},\tau _{2}\in \R\},$

\ge\gt$^{1}=\{it_{kj}(Ut_{kj})\in $\gR\gt$^{\C} \mid t_{kj}\in \R,1\leq k\leq 3,0\leq j\leq 7\},$

\ge\gt$=$\ge\gt$^{0}\oplus $\ge\gt$^{1}$

\noindent
Also we define as follows.

\ge\gd$_{6}=\{X|_{6} \mid X\in $\ge\gd$\},$

\ge\gm$_{6}=\{X|_{6} \mid X\in $\ge\gm$\},$

\ge\gt$^{0}_{6}=\{X|_{6} \mid X\in $\ge\gt$^{0}\}, $ \ge\gt$^{1}_{6}=\{X|_{6} \mid X\in $\ge\gt$^{1}\}, $ \ge\gt$_{6}=\{X|_{6} \mid X\in $\ge\gt$\}$.

\bigskip

\emph{Theorem 15.11.} Let we put

\ge$_{6}=$\ge\gd$_{6}\oplus $\ge\gm$_{6}\oplus $\ge\gt$_{6}$ . \ 

\noindent
Then \ge$_{6}$ is a compact exceptional simple Lie algebra of type $E_{6}.$

\bigskip

\emph{Proof.} \ By \emph{Lemma 15.2,15.3}, and \emph{15.4}, \ge$_{6}$ is a Lie algebra.
And by calculations using Maxima,we have as follows.

For $X_{1}\in $\ge\gd$_{6},X_{2}\in $\ge\gm$_{6},X_{3}\in $\ge\gt$^{0}_{6},X_{4}\in $\ge\gt$^{1}_{6}$,

$B_{6}(X_{1},X_{1})$

$=-24(d_{67}^{2}+d_{57}^{2}+d_{56}^{2}+d_{47}^{2}+d_{46}%
^{2}+d_{45}^{2}+d_{37}^{2}$

\ \ \ \ \ \ $\ +d_{36}^{2}+d_{35}^{2}+d_{34}^{2}+d_{27}^{2}+d_{26}^{2}+d_{25}^{2}+d_{24}^{2}$

\ \ \ \ \ \ $\ +d_{23}^{2}+d_{17}^{2}+d_{16}^{2}+d_{15}^{2}+d_{14}^{2}+d_{13}^{2}+d_{12}^{2}$

\ \ \ \ \ \ $\ +d_{07}^{2}+d_{06}^{2}+d_{05}^{2}+d_{04}^{2}+d_{03}^{2}+d_{02}^{2}+d_{01}^{2}),$

$B_{6}(X_{2},X_{2})$

$=-24(m_{37}^{2}+m_{36}^{2}+m_{35}^{2}+m_{34}^{2}+m_{33}%
^{2}+m_{32}^{2}+m_{31}^{2}+m_{30}^{2}$

\ \ $\ +m_{27}^{2}+m_{26}^{2}+m_{25}^{2}+m_{24}^{2}+m_{23}%
^{2}+m_{22}^{2}+m_{21}^{2}+m_{20}^{2}$

\ \ $\ +m_{17}^{2}+m_{16}^{2}+m_{15}^{2}+m_{14}^{2}+m_{13}%
^{2}+m_{12}^{2}+m_{11}^{2}+m_{10}^{2}),$

$B_{6}(X_{3},X_{3})=-24(\tau _{1}^{2}+3\tau _{2}^{2})$

$B_{6}(X_{4},X_{4})$

$=-24(t_{37}^{2}+t_{36}^{2}+t_{35}^{2}+t_{34}^{2}+t_{33}%
^{2}+t_{32}^{2}+t_{31}^{2}+t_{30}^{2}$

\ \ \ $\ +t_{27}^{2}+t_{26}^{2}+t_{25}^{2}+t_{24}^{2}+t_{23}%
^{2}+t_{22}^{2}+t_{21}^{2}+t_{20}^{2}$

\ \ $\ \ $+$t_{17}^{2}+t_{16}^{2}+t_{15}^{2}+t_{14}^{2}+t_{13}%
^{2}+t_{12}^{2}+t_{11}^{2}+t_{10}^{2}),$

$B_{4}(X_{i},X_{j})=0,(i\neq j).$

\noindent
Therfor we have $B_{6}(X,X)<0$,for $^{\forall }X\neq 0,X\in $\ge$_{6}.$ Then \ge$_{6}$ is compact.\ \ \ \emph{Q.E.D.}

\bigskip

Next, let digitalize the exceptional simple Lie algebra of type $E_{6}$ as
matrices in $M(27 \times 27,\C)$.

\bigskip

\emph{Definition 15.12.}  We consider following elements in $M(27 \times 27,\R)$,

For \ge\gt\  $\ni iR\tau_{1}, iR\tau_{2}, (\tau_{1}, \tau_{2} \in \R$)

$R\tau_{1}^{`}=\tau_{1}$
{\fontsize{6pt}{8pt} \selectfont%
$\left( 

\right)  $
}

\noindent
See \emph{Lemma 6.17} for MIr(m), \emph{Lemma 6.19} for MC(m). 

We denote elements of $Rt_{1}^{`}, Rt_{2}^{`}$, and $Rt_{3}^{`}$. 

\noindent
$Rt_{10}^{'}=t_{10}$
{\fontsize{6pt}{8pt} \selectfont%
$\left( 
%
%
\right)$ }.

We define the following Lie algebra:

\ge$^{`}_{6}=$\{ Lie algebra over real numbers generated by $Rd_{ij}^{`}(0 \leq i < j \leq 7)$,
$Rm_{ij}^{`}(1 \leq i \leq 3, 0 \leq j \leq 7),  iR\tau_{1}^{`}, iR\tau_{2}^{`}, iRt_{ij}^{`}(1 \leq i \leq 3, 0 \leq j \leq 7)\}$

\bigskip

\emph{Lemma 15.13.}  The following correspondences are homomorphism.

\ \ \ \ \ \ \ \ \ \ \ \ \ \ \ \ \ \ \ge\gd$ \ni Rd_{ij} \rightarrow Rd_{ij}^{`} \in $ \ge$_{6}^{`}$

\ \ \ \ \ \ \ \ \ \ \ \ \ \ \ \ \ \ \ge\gm$ \ni Rm_{ij} \rightarrow Rm_{ij}^{`} \in $ \ge$_{6}^{`}$

\ \ \ \ \ \ \ \ \ \ \ \ \ \ \ \ \ \ \ge\gt$ \ni iR\tau_{i}\rightarrow iR\tau_{i}^{`} \in $ \ge$_{6}^{`}$

\ \ \ \ \ \ \ \ \ \ \ \ \ \ \ \ \ \ \ge\gt$ \ni iRt_{ij} \rightarrow iRt_{ij}^{`} \in $ \ge$_{6}^{`}$

\bigskip

\emph{Proof.}  Let's check these Lie bracket operations.

For $Rd_{ij}^{`}/d_{ij},Rm_{ij}^{`}/m_{ij}, R\tau_{i}^{`}/\tau_{i}$, and $Rt_{ij}^{`}/t_{ij}$, by calculation with Maxima we have same results of 
\emph{Lemma 15.2},\emph{Lemma 15.3}, and \emph{Lemma 15.4}. \ \ \ \ \ \emph{Q.E.D.}

\bigskip

Therfore we have \emph{Theorem 15.14}

\bigskip

\emph{Theorem 15.14.}  \ge$^{`}_{6}$ is isomorphic to \ge$_{6}$. 

\bigskip

We identify \ge$_{6}$ with \ge$_{6}^{`}$.

\bigskip

\emph{Lemma 15.15.}  \ge$_{6}^{\C}$ is expressed  by 

$\{ D \in Hom_{\C}(\C^{27}) \mid \textrm{C}_{m}(DX,\textrm{F}_{m}(X, X))=0, X\in \C^{27} \}.$
And $\{ Rd_{ij}^{`}(0 \leq i < j \leq 7)$,
$Rm_{ij}^{`}(1 \leq i \leq 3, 0 \leq j \leq 7), i(R\tau_{1}^{`}/\tau_{1}-R\tau_{2}^{`}/\tau_{2})\tau_{1}, i(R\tau_{1}^{`}/\tau_{1}+R\tau_{2}^{`}/\tau_{2})\tau_{2}, iRt_{ij}^{`}(1 \leq i \leq 3, 0 \leq j \leq 7) \}$
 are orthogonal bases of \ge$_{6}$.

\bigskip

\emph{Proof} \ By \emph{Remark 6.5} and \emph{Proposition 3.6}, \ge$_{6}^{\C}$ is expressed  by 

$\{ D \in Hom_{\C}(\C^{27}) \mid \textrm{C}_{m}(DX,\textrm{F}_{m}(X, X))=0, X\in \C^{27} \}.$

Let check each element of $Rd_{ij}^{`}(0 \leq i < j \leq 7)$,
$Rm_{ij}^{`}(1 \leq i \leq 3, 0 \leq j \leq 7), i(R\tau_{1}^{`}/\tau_{1}-R\tau_{2}^{`}/\tau_{2})\tau_{1}, i(R\tau_{1}^{`}/\tau_{1}+R\tau_{2}^{`}/\tau_{2})\tau_{2}, iRt_{ij}^{`}(1 \leq i \leq 3, 0 \leq j \leq 7)$
satisfy the condition:

$\textrm{C}_{m}(DX,\textrm{F}_{m}(X, X))=0$,

$X=(\chi_{1},\chi_{2},\chi_{3},(x_{10},\cdot\cdot\cdot,x_{17}),(x_{20},\cdot\cdot\cdot,x_{27}),(x_{30},\cdot\cdot\cdot,x_{37})) \ \in \C^{27}$.

By calculation, we have

F$_{m}(X,X)=%
(-x_{17}^{2}-x_{16}^{2}-x_{15}^{2}-x_{14}^{2}-x_{13}^{2}-x_{12}^{2}-x_{11}^{2}-x_{10}^{2}+\chi_{2}\chi_{3}$,

$-x_{27}^{2}-x_{26}^{2}-x_{25}^{2}-x_{24}^{2}-x_{23}^{2}-x_{22}^{2}-x_{21}^{2}-x_{20}^{2}+\chi_{1}\chi_{3}$,

$-x_{37}^{2}-x_{36}^{2}-x_{35}^{2}-x_{34}^{2}-x_{33}^{2}-x_{32}^{2}-x_{31}^{2}-x_{30}^{2}+\chi_{1}\chi_{2}$,

$(-x_{27}x_{37}-x_{26}x_{36}-x_{25}x_{35}-x_{24}x_{34}-x_{23}x_{33}-x_{22}x_{32}-x_{21}x_{31}+x_{20}x_{30}-\chi_{1}x_{10}$,

$-x_{26}x_{37}+x_{27}x_{36}-x_{24}x_{35}+x_{25}x_{34}-x_{22}x_{33}+x_{23}x_{32}-x_{20}x_{31}-x_{21}x_{30}-\chi_{1}x_{11}$,

$-x_{25}x_{37}+x_{24}x_{36}+x_{27}x_{35}-x_{26}x_{34}+x_{21}x_{33}-x_{20}x_{32}-x_{23}x_{31}-x_{22}x_{30}-\chi_{1}x_{12},$

$-x_{24}x_{37}-x_{25}x_{36}+x_{26}x_{35}+x_{27}x_{34}-x_{20}x_{33}-x_{21}x_{32}+x_{22}x_{31}-x_{23}x_{30}-\chi_{1}x_{13},$

$x_{23}x_{37}-x_{22}x_{36}+x_{21}x_{35}-x_{20}x_{34}-x_{27}x_{33}+x_{26}x_{32}-x_{25}x_{31}-x_{24}x_{30}-\chi_{1}x_{14},$

$x_{22}x_{37}+x_{23}x_{36}-x_{20}x_{35}-x_{21}x_{34}-x_{26}x_{33}-x_{27}x_{32}+x_{24}x_{31}-x_{25}x_{30}-\chi_{1}x_{15},$

$x_{21}x_{37}-x_{20}x_{36}-x_{23}x_{35}+x_{22}x_{34}+x_{25}x_{33}-x_{24}x_{32}-x_{27}x_{31}-x_{26}x_{30}-\chi_{1}x_{16},$

$-x_{20}x_{37}-x_{21}x_{36}-x_{22}x_{35}-x_{23}x_{34}+x_{24}x_{33}+x_{25}x_{32}+x_{26}x_{31}-x_{27}x_{30}-\chi_{1}x_{17})$,

($-x_{17}x_{37}-x_{16}x_{36}-x_{15}x_{35}-x_{14}x_{34}-x_{13}x_{33}-x_{12}x_{32}-x_{11}x_{31}+x_{10}x_{30}-\chi_{2}x_{20},$

$x_{16}x_{37}-x_{17}x_{36}+x_{14}x_{35}-x_{15}x_{34}+x_{12}x_{33}-x_{13}x_{32}-x_{10}x_{31}-x_{11}x_{30}-\chi_{2}x_{21},$

$x_{15}x_{37}-x_{14}x_{36}-x_{17}x_{35}+x_{16}x_{34}-x_{11}x_{33}-x_{10}x_{32}+x_{13}x_{31}-x_{12}x_{30}-\chi_{2}x_{22},$

$x_{14}x_{37}+x_{15}x_{36}-x_{16}x_{35}-x_{17}x_{34}-x_{10}x_{33}+x_{11}x_{32}-x_{12}x_{31}-x_{13}x_{30}-\chi_{2}x_{23},$

$-x_{13}x_{37}+x_{12}x_{36}-x_{11}x_{35}-x_{10}x_{34}+x_{17}x_{33}-x_{16}x_{32}+x_{15}x_{31}-x_{14}x_{30}-\chi_{2}x_{24},$

$-x_{12}x_{37}-x_{13}x_{36}-x_{10}x_{35}+x_{11}x_{34}+x_{16}x_{33}+x_{17}x_{32}-x_{14}x_{31}-x_{15}x_{30}-\chi_{2}x_{25},$

$-x_{11}x_{37}-x_{10}x_{36}+x_{13}x_{35}-x_{12}x_{34}-x_{15}x_{33}+x_{14}x_{32}+x_{17}x_{31}-x_{16}x_{30}-\chi_{2}x_{26},$

$-x_{10}x_{37}+x_{11}x_{36}+x_{12}x_{35}+x_{13}x_{34}-x_{14}x_{33}-x_{15}x_{32}-x_{16}x_{31}-x_{17}x_{30}-\chi_{2}x_{27}),$

$(-\chi_{3}x_{30}-x_{17}x_{27}-x_{16}x_{26}-x_{15}x_{25}-x_{14}x_{24}-x_{13}x_{23}-x_{12}x_{22}-x_{11}x_{21}+x_{10}x_{20},$

$-\chi_{3}x_{31}-x_{16}x_{27}+x_{17}x_{26}-x_{14}x_{25}+x_{15}x_{24}-x_{12}x_{23}+x_{13}x_{22}-x_{10}x_{21}-x_{11}x_{20},$

$-\chi_{3}x_{32}-x_{15}x_{27}+x_{14}x_{26}+x_{17}x_{25}-x_{16}x_{24}+x_{11}x_{23}-x_{10}x_{22}-x_{13}x_{21}-x_{12}x_{20},$

$-\chi_{3}x_{33}-x_{14}x_{27}-x_{15}x_{26}+x_{16}x_{25}+x_{17}x_{24}-x_{10}x_{23}-x_{11}x_{22}+x_{12}x_{21}-x_{13}x_{20},$

$-\chi_{3}x_{34}+x_{13}x_{27}-x_{12}x_{26}+x_{11}x_{25}-x_{10}x_{24}-x_{17}x_{23}+x_{16}x_{22}-x_{15}x_{21}-x_{14}x_{20},$

$-\chi_{3}x_{35}+x_{12}x_{27}+x_{13}x_{26}-x_{10}x_{25}-x_{11}x_{24}-x_{16}x_{23}-x_{17}x_{22}+x_{14}x_{21}-x_{15}x_{20},$

$-\chi_{3}x_{36}+x_{11}x_{27}-x_{10}x_{26}-x_{13}x_{25}+x_{12}x_{24}+x_{15}x_{23}-x_{14}x_{22}-x_{17}x_{21}-x_{16}x_{20},$

$-\chi_{3}x_{37}-x_{10}x_{27}-x_{11}x_{26}-x_{12}x_{25}-x_{13}x_{24}+x_{14}x_{23}+x_{15}x_{22}+x_{16}x_{21}-x_{17}x_{20}$))

In case of $D=Rd_{01}^{`}/d_{01}$:      

$DX=(0,0,0,$

\ \ \ \ \ \ \ \ \ $(x_{11},-x_{10},0,0,0,0,0,0),$

\ \ \ \ \ \ \ \ \ $(-\frac{x_{21}}{2},\frac{x_{20}}{2},-\frac{x_{23}}{2},\frac{x_{22}}{2},-\frac{x_{25}}{2},\frac{x_{24}}{2},-\frac{x_{27}}{2},\frac{x_{26}}{2}),$

\ \ \ \ \ \ \ \ \ $(-\frac{x_{31}}{2},\frac{x_{30}}{2},\frac{x_{33}}{2},-\frac{x_{32}}{2},\frac{x_{35}}{2},-\frac{x_{34}}{2},\frac{x_{37}}{2},-\frac{x_{36}}{2}))$.

\noindent
Therfore we have C$_{m}(DX$, \textrm{F}$_{m}(X,X))=0$.
Other cases can be confirmed in the same way.

For $R_{1},R_{2} \in  \{ Rd_{ij}^{`}/d_{ij}(0 \leq i < j \leq 7),Rm_{ij}^{`}/m_{ij}(1 \leq i \leq 3, 0 \leq j \leq 7) , i(R\tau_{1}^{`}/\tau_{1}-R\tau_{2}^{`}/\tau_{2}), i(R\tau_{1}^{`}/\tau_{1}+R\tau_{2}^{`}/\tau_{2}), iRt_{ij}^{`}/t_{ij}(1 \leq i \leq 3, 0 \leq j \leq 7)\}$,
let's  calculate $tr(R_{1}.^{t}R_{2})$.

\ \ \ \ \ \ \ \ \ \ \ \ \ \ \ \ $tr(R_{1}.^{t}R_{2})=0 \ \ \ (R_{1} \neq R_{2})$,

\ \ \ \ \ \ \ \ \ \ \ \ \ \ \ \ $tr(R_{1}.^{t}R_{1})=6 \ \ \ (R_{1} \in  \{ Rd_{ij}^{`}/d_{ij}(0 \leq i < j \leq 7)\}$,

\ \ \ \ \ \ \ \ \ \ \ \ \ \ \ \ $tr(R_{2}.^{t}R_{2})=\frac{13}{2} \ \ (R_{2} \in \{Rm_{ij}^{`}/m_{ij}(1 \leq i \leq 3, 0 \leq j \leq 7) \}$,

\ \ \ \ \ \ \ \ \ \ \ \ \ \ \ \ $tr(R_{1}.^{t}R_{1})=-6 \ \ \ (R_{1} =i(R\tau_{1}^{`}/\tau_{1}-R\tau_{2}^{`}/\tau_{2}) )$,

\ \ \ \ \ \ \ \ \ \ \ \ \ \ \ \ $tr(R_{1}.^{t}R_{1})=-18 \ \ (R_{1} =i(R\tau_{1}^{`}/\tau_{1}+R\tau_{2}^{`}/\tau_{2}) )$,

\ \ \ \ \ \ \ \ \ \ \ \ \ \ \ \ $tr(R_{2}.^{t}R_{2})=-\frac{13}{2} \ \ (R_{2} \in \{iRt_{ij}^{`}/t_{ij}(1 \leq i \leq 3, 0 \leq j \leq 7) \}$.

Since dim \ge$_{6}=78$, then $\{ Rd_{ij}^{`}(0 \leq i < j \leq 7),Rm_{ij}^{`}(1 \leq i \leq 3, 0 \leq j \leq 7), i(R\tau_{1}^{`}/\tau_{1}-R\tau_{2}^{`}/\tau_{2})\tau_{1}, i(R\tau_{1}^{`}/\tau_{1}+R\tau_{2}^{`}/\tau_{2})\tau_{2}, iRt_{ij}^{`}(1 \leq i \leq 3, 0 \leq j \leq 7) \}$ are 
orthogonal bases of \ge$_{6}$.
\ \ \ \ \ \emph{Q.E.D.}

\bigskip

\bigskip

\emph{Lemma 15.16.}  The simply connected compact Lie group $E_{6}$ is expressed by

$E_{6}=\{\alpha \in Iso_{\C}(\C^{27})%
 \mid \textrm{F}_{m}(\alpha X , \alpha Y)=\iota \alpha \iota \textrm{F}_{m}(X, Y),\langle \alpha X,\alpha Y \rangle= \langle X,Y \rangle \}$.

\bigskip

\emph{Proof.}  By \emph{Definition 3.3}, \emph{Proposition 3.4}, and \emph{Remark 6.5}, we have the above \emph{Lemma}.\ \ \ \ \ \emph{Q.E.D.}

\bigskip

\section{The exceptional simle Lie algebra \gr$_{7}^{\C}$ of type $E_{7}$}

\bigskip

\ \ \ \ \ \emph{Definition 16.1.} \ We define the followings for elements $%
Ra_{ij},R\alpha _{k}\in $\gR\ga$^{\C},Rb_{ij},$

\noindent
$R\beta _{k}\in $\gR\gb$^{\C},R\rho _{1}\in $\gR\gp$^{\C} .$

\ \ \ \ \ $\ \ \ \ \ \ \ \ \ \ \ \ \ Ua_{ij}=Ra_{ij}/a_{ij}\ \ (1\leq i\leq 3,0\leq
j\leq 7),$

\ \ \ \ \ \ \ \ \ \ \ \ \ \ \ \ $\ \ U\alpha _{k}=R\alpha _{k}/\alpha _{k}$ \ $(1\leq
k\leq 3),$

\ \ \ \ \ $\ \ \ \ \ \ \ \ \ \ \ \ \ Ub_{ij}=Rb_{ij}/b_{ij}\ \ (1\leq i\leq 3,0\leq
j\leq 7),$

\ \ \ \ \ \ \ \ \ \ \ \ \ \ \ \ $\ \ U\beta _{k}=R\beta _{k}/\beta _{k}$ \ $(1\leq
k\leq 3),$

\ \ \ \ \ \ \ \ \ \ \ \ \ \ \ \ $\ \ U\rho _{1}=R\rho _{1}/\rho _{1}$ $.$

\bigskip

\emph{Lemma 16.2. }\ We have the following Lie bracket operations.

\ \ \ \ $[Ua_{ij},Ua_{kl}]=0,(1\leq i,k\leq 3,0\leq j,l\leq 7),$

\ \ \ \ $[U\alpha _{i},U\alpha _{k}]=0,(1\leq i,k\leq 3),$

\ \ \ \ $[Ua_{ij},U\alpha _{k}]=0,(1\leq i,k\leq 3,0\leq j\leq 7),$

\ \ \ \ $[Ub_{ij},Ub_{kl}]=0,(1\leq i,k\leq 3,0\leq j,l\leq 7),$

\ \ \ \ $[U\beta _{i},U\beta _{k}]=0,(1\leq i,k\leq 3),$

\ \ \ \ $[Ub_{ij},U\beta _{k}]=0,(1\leq i,k\leq 3,0\leq j\leq 7),$

\ \ \ \ $[Ua_{1i},Ub_{1j}]=2Ud_{ij}\ ,(0\leq i\neq j\leq 7),$

\ \ \ \ $[Ua_{1i},Ub_{1i}]=-\frac{4}{3}U\tau _{1}+\frac{2}{3}u\tau2+2U\rho _{1},(0\leq i\leq 7),$

\ \ \ \ $[Ua_{1i},Ub_{2j}]=-Sn(i+1,j+1)Um_{3k}+Sn(i+1,j+1)Ut_{3k},$

\ \ \ \ \ \ \ \ \ \ \ \ \ \ \ \ \ \ \ \ \ $(k=Ca(i+1,j+1),\ 0\leq i,j\leq7),$

\ \ \ \ $[Ua_{1i},Ub_{3j}]=-Sn(j+1,i+1)Um_{2k}-Sn(j+1,i+1)Ut_{2k},$

\ \ \ \ \ \ \ \ \ \ \ \ \ \ \ \ \ \ \ \ \ $(k=Ca(j+1,i+1),\ 0\leq i,j\leq7),$

\ \ \ \ $[Ua_{1i},U\beta _{1}]=0,$

\ \ \ \ $[Ua_{1i},U\beta _{2}]=-Um_{1i}+Ut_{1i},(0\leq i\leq 7),$

\ \ \ \ $[Ua_{1i},U\beta _{3}]=Um_{1i}+Ut_{1i},(0\leq i\leq 7),$

\ \ \ \ $[Ua_{2i},Ub_{1j}]=-Sn(j+1,i+1)Um_{3k}-Sn(j+1,i+1)Ut_{3k},$

\ \ \ \ \ \ \ \ \ \ \ \ \ \ \ \ \ \ \ \ \ $(k=Ca(j+1,i+1),\ 0\leq i,j\leq7),$

\ \ \ \ $[Ua_{2i},Ub_{2j}]=2\sum\limits_{0\leq n<l\leq7}Mv(ki,kj)Ud_{nl},$

\ \ \ \ \ \ \ \ \ \ \ \ \ \ \ \ \ \ \ \ \ $(ki=Nu(i+1,j+1),kj=Nu(n+1,l+1),\ 0\leq i\neq j\leq 7),$

\ \ \ \ $[Ua_{2i},Ub_{2i}]=\frac{2}{3}U\tau _{1}-\frac{4}{3}u\tau2+2U\rho _{1},(0\leq i\leq 7),$

\ \ \ \ $[Ua_{2i},Ub_{3j}]=-Sn(i+1,j+1)Um_{1k}+Sn(i+1,j+1)Ut_{1k},$

\ \ \ \ \ \ \ \ \ \ \ \ \ \ \ \ \ \ \ \ \ $(k=Ca(i+1,j+1),\ 0\leq i,j\leq 7),$

\ \ \ \ $[Ua_{2i},U\beta _{1}]=Um_{2i}+Ut_{2i},(0\leq i\leq 7),$

\ \ \ \ $[Ua_{2i},U\beta _{2}]=0,$

\ \ \ \ $[Ua_{2i},U\beta _{3}]=-Um_{2i}+Ut_{2i},(0\leq i\leq 7),$

\ \ \ \ $[Ua_{3i},Ub_{1j}]=-Sn(i+1,j+1)Um_{2k}+Sn(i+1,j+1)Ut_{2k},$

\ \ \ \ \ \ \ \ \ \ \ \ \ \ \ \ \ \ \ \ \ $(k=Ca(i+1,j+1),\ 0\leq i,j\leq 7),$

\ \ \ \ $[Ua_{3i},Ub_{2j}]=-Sn(j+1,i+1)Um_{1k}-Sn(j+1,i+1)Ut_{1k},$

\ \ \ \ \ \ \ \ \ \ \ \ \ \ \ \ \ \ \ \ \ $(k=Ca(j+1,i+1),\ 0\leq i,j\leq 7),$

\ \ \ \ $[Ua_{3i},Ub_{3j}]=2\sum\limits_{0\leq n<l\leq 7}Mv^{2}(ki,kj)Ud_{nl},$

\ \ \ \ \ \ \ \ \ \ \ \ \ \ \ \ \ \ \ \ \ $(ki=Nu(i+1,j+1),kj=Nu(n+1,l+1),0\leq i\neq j\leq 7),$

\ \ \ \ $[Ua_{3i},Ub_{3i}]=\frac{2}{3}U\tau _{1}+\frac{2}{3}U\tau _{2}+2U\rho _{1},(0\leq i\leq 7),$

\ \ \ \ $[Ua_{3i},U\beta _{1}]=-Um_{3i}+Ut_{3i},(0\leq i\leq 7),$

\ \ \ \ $[Ua_{3i},U\beta _{2}]=Um_{3i}+Ut_{3i},(0\leq i\leq 7),$

\ \ \ \ $[Ua_{3i},U\beta _{3}]=0,$

\ \ \ \ $[U\alpha _{1},Ub_{1j}]=0,(0\leq j\leq 7),$

\ \ \ \ $[U\alpha _{1},Ub_{2j}]=-Um_{2j}+Ut_{2j},(0\leq j\leq 7),$

\ \ \ \ $[U\alpha _{1},Ub_{3j}]=Um_{3j}+Ut_{3j},(0\leq j\leq 7),$

\ \ \ \ $[U\alpha _{1},U\beta _{1}]=\frac{4}{3}U\tau _{1}-\frac{2}{3}U\tau _{2}+U\rho _{1},$

\ \ \ \ $[U\alpha _{1},U\beta _{2}]=0,$

\ \ \ \ $[U\alpha _{1},U\beta _{3}]=0,$

\ \ \ \ $[U\alpha _{2},Ub_{1j}]=Um_{1j}+Ut_{1j},(0\leq j\leq 7),$

\ \ \ \ $[U\alpha _{2},Ub_{2j}]=0,(0\leq j\leq 7),$

\ \ \ \ $[U\alpha _{2},Ub_{3j}]=-Um_{3j}+Ut_{3j},(0\leq j\leq 7),$

\ \ \ \ $[U\alpha _{2},U\beta _{1}]=0,$

\ \ \ \ $[U\alpha _{2},U\beta _{2}]=-\frac{2}{3}U\tau _{1}+\frac{4}{3}U\tau _{2}+U\rho _{1},$

\ \ \ \ $[U\alpha _{2},U\beta _{3}]=0,$

\ \ \ \ $[U\alpha _{3},Ub_{1j}]=-Um_{1j}+Ut_{1j},(0\leq j\leq 7),$

\ \ \ \ $[U\alpha _{3},Ub_{2j}]=Um_{2j}+Ut_{2j},(0\leq j\leq 7)$,

\ \ \ \ $[U\alpha _{3},Ub_{3j}]=0,(0\leq j\leq 7),$

\ \ \ \ $[U\alpha _{3},U\beta _{1}]=0,$

\ \ \ \ $[U\alpha _{3},U\beta _{2}]=0,$

\ \ \ \ $[U\alpha _{3},U\beta _{3}]=-\frac{2}{3}U\tau _{1}-\frac{2}{3}U\tau _{2}+U\rho _{1}.$

\bigskip

\emph{Proof. \ }We have the above Lie bracket operations with
calculations using Maxima.\ \ \ \ \emph{Q.E.D.}

\bigskip

\emph{Lemma 16.3.} \ We have the following Lie bracket operations.

\ \ \ \ $[Ud_{ij},Ua_{1k}]=-Ua_{1j}$ $($in case of $k=i),$

\ \ \ \  \ \ \ \ \ \ \ \ \ \ \ \ \ \ \ \ $=Ua_{1i}$ $($in case of $k=j),$

$\ \ \ \ \ \ \ \ \ \ \ \ \ \ \ \ \  \ \ \ =0$ $($in case of $%
k\neq i,j),(0\leq i<j\leq 7),$

\ \ \ \ $[Ud_{ij},Ua_{2k}]=-\sum\limits_{0\leq n<l\leq 7}Mv^{2}(ki,kj)Ua_{2l}$ $($where $k=n)$

\ \ \ \ \ \ \ \ \ \ \ \ \ \ \ \ \ \ \ \ \ \ \ $
+\sum\limits_{0\leq n<l\leq 7}Mv^{2}(ki,kj)Ua_{2n}$ $($where $k=l)$,

\ \ \ \ \ \ \ \ \ \ \ \ \ \ \ \ \ \ \ \ \ $(ki=Nu(i+1,j+1),kj=Nu(n+1,l+1),0\leq i<j\leq 7),$

\ \ \ \ $[Ud_{ij},Ua_{3k}]=-\sum\limits_{0\leq n<l\leq 7}Mv(ki,kj)Ua_{3l}$ $($where $n=k)$

\ \ \ \ \ \ \ \ \ \ \ \ \ \ \ \ \ \ \ \ \ \ \ $
+\sum\limits_{0\leq n<l\leq 7}Mv(ki,kj)Ua_{3n}$ $($where $l=k),$

\ \ \ \ \ \ \ \ \ \ \ \ \ \ \ \ \ \ \ \ \ $(ki=Nu(i+1,j+1),kj=Nu(n+1,l+1),0\leq i<j\leq 7),$

\ \ \ \ $[Ud_{ij},U\alpha _{k}]=0,(0\leq i<j\leq 7,k=1,2).$

\ \ \ \ $[Ud_{ij},Ub_{1k}]=-Ub_{1j}$ $($in case of $k=i),$

$\ \ \ \ \ \ \ \ \ \ \ \ \ \ \ \ \ \ \ \ =Ub_{1i}$ $($in case of $k=j),$

$\ \ \ \ \ \ \ \ \ \ \ \ \ \ \ \  \ \ \ \ =0$ $($in case of $k\neq i,j),(0\leq i<j\leq 7),$

\ \ \ \ $[Ud_{ij},Ub_{2k}]=-\sum\limits_{0\leq n<l\leq 7}Mv^{2}(ki,kj)Ub_{2l}$ $($where $k=n)$

\ \ \ \ \ \ \ \ \ \ \ \ \ \ \ \ \ \ \ \ \ \ \ $
+\sum\limits_{0\leq n<l\leq 7}Mv^{2}(ki,kj)Ub_{2n}$ $($where $k=l)$,

$\ \ \ \ \ \ \ \ \ \ \ \ \ \ \ \ \ \ \ \ \ \ \ \ \ \ \ \ \ \ \
(ki=Nu(i+1,j+1),kj=Nu(n+1,l+1),0\leq i<j\leq 7),$

\ \ \ \ $[Ud_{ij},Ub_{3k}]=-\sum\limits_{0\leq n<l\leq 7}Mv(ki,kj)Ub_{3l}$ $($where $n=k)$

\ \ \ \ \ \ \ \ \ \ \ \ \ \ \ \ \ \ \ \ \ \ \ $
+\sum\limits_{0\leq n<l\leq 7}Mv(ki,kj)Ub_{3n}$ $($where $l=k),$

$\ \ \ \ \ \ \ \ \ \ \ \ \ \ \ \ \ \ \ \ \ \ \ \ \ \ \ \ \ \ \
(ki=Nu(i+1,j+1),kj=Nu(n+1,l+1),0\leq i<j\leq 7),$

\ \ \ \ $[Ud_{ij},U\beta _{k}]=0,(0\leq i<j\leq 7,k=1,2).$

\bigskip

\emph{Proof. \ }We have the above Lie bracket operations with
calculations using Maxima.\ \ \ \ \emph{Q.E.D.}

\bigskip

\emph{Lemma 16.4.} \ We have the following Lie bracket operations.

\ \ \ \ $[Um_{1i},Ua_{1j}]=0\ ,(0\leq i\neq j\leq 7),$

\ \ \ \ $[Um_{1i},Ua_{1i}]=U\alpha _{2}-U\alpha _{3},(0\leq i\leq 7),$

\ \ \ \ $[Um_{1i},Ua_{2j}]=\frac{1}{2}Sn(i+1,j+1)Ua_{3k},(k=Ca(i+1,j+1),\ 0\leq i,j\leq 7),$

\ \ \ \ $[Um_{1i},Ua_{3j}]=-\frac{1}{2}Sn(j+1,i+1)Ua_{2k},(k=Ca(j+1,i+1),\ 0\leq i,j\leq 7),$

\ \ \ \ $[Um_{1i},U\alpha _{1}]=0\ ,(0\leq i\leq 7),$

\ \ \ \ $[Um_{1i},U\alpha _{2}]=-\frac{1}{2}Ua_{1i}\ ,(0\leq i\leq 7),$

\ \ \ \ $[Um_{1i},U\alpha _{3}]=\frac{1}{2}Ua_{1i}\ ,(0\leq i\leq 7),$

\ \ \ \ $[Um_{1i},Ub_{1j}]=0\ ,(0\leq i\neq j\leq 7),$

\ \ \ \ $[Um_{1i},Ub_{1i}]=U\beta _{2}-U\beta _{3},(0\leq i\leq 7),$

\ \ \ \ $[Um_{1i},Ub_{2j}]=\frac{1}{2}Sn(i+1,j+1)Ub_{3k},(k=Ca(i+1,j+1),\ 0\leq i,j\leq 7),$

\ \ \ \ $[Um_{1i},Ub_{3j}]=-\frac{1}{2}Sn(j+1,i+1)Ub_{2k},(k=Ca(j+1,i+1),\ 0\leq i,j\leq 7),$

\ \ \ \ $[Um_{1i},U\beta _{1}]=0\ ,(0\leq i\leq 7),$

\ \ \ \ $[Um_{1i},U\beta _{2}]=-\frac{1}{2}Ub_{1i}\ ,(0\leq i\leq 7),$

\ \ \ \ $[Um_{1i},U\beta _{3}]=\frac{1}{2}Ub_{1i}\ ,(0\leq i\leq 7),$

\ \ \ \ $[Um_{2i},Ua_{1j}]=-\frac{1}{2}Sn(j+1,i+1)Ua_{3k},(k=Ca(j+1,i+1),\ 0\leq i,j\leq 7),$

\ \ \ \ $[Um_{2i},Ua_{2i}]=-U\alpha _{1}+U\alpha _{3},(0\leq i\leq 7),$

\ \ \ \ $[Um_{2i},Ua_{2j}]=0\ ,(0\leq i\neq j\leq 7),$

\ \ \ \ $[Um_{2i},Ua_{3j}]=\frac{1}{2}Sn(i+1,j+1)Ua_{1k},(k=Ca(i+1,j+1),\ 0\leq i,j\leq 7),$

\ \ \ \ $[Um_{2i},U\alpha _{1}]=\frac{1}{2}Ua_{2i}\ ,(0\leq i\leq 7)$,

\ \ \ \ $[Um_{2i},U\alpha _{2}]=0\ ,(0\leq i\leq 7)$,

\ \ \ \ $[Um_{2i},U\alpha _{3}]=-\frac{1}{2}Ua_{2i}\ ,(0\leq i\leq 7),$

\ \ \ \ $[Um_{2i},Ub_{1j}]=-\frac{1}{2}Sn(j+1,i+1)Ub_{3k},(k=Ca(j+1,i+1),\ 0\leq i,j\leq 7),$

\ \ \ \ $[Um_{2i},Ub_{2i}]=-U\beta _{1}+U\beta _{3},(0\leq i\leq 7),$

\ \ \ \ $[Um_{2i},Ub_{2j}]=0\ ,(0\leq i\neq j\leq 7),$

\ \ \ \ $[Um_{2i},Ub_{3j}]=\frac{1}{2}Sn(i+1,j+1)Ub_{1k}, (k=Ca(i+1,j+1),\ 0\leq i,j\leq 7),$

\ \ \ \ $[Um_{2i},U\beta _{1}]=\frac{1}{2}Ub_{2i}\ ,(0\leq i\leq 7),$

\ \ \ \ $[Um_{2i},U\beta _{2}]=0\ ,(0\leq i\leq 7),$

\ \ \ \ $[Um_{2i},U\beta _{3}]=-\frac{1}{2}Ub_{2i}\ ,(0\leq i\leq 7),$

\ \ \ \ $[Um_{3i},Ua_{1j}]=\frac{1}{2}Sn(i+1,j+1)Ua_{2k},(k=Ca(i+1,j+1),\ 0\leq i,j\leq 7),$

\ \ \ \ $[Um_{3i},Ua_{2j}]=-\frac{1}{2}Sn(j+1,i+1)Ua_{1k},(k=Ca(j+1,i+1),\ 0\leq i,j\leq 7),$

\ \ \ \ $[Um_{3i},Ua_{3i}]=U\alpha _{1}-U\alpha _{2},(0\leq i\leq 7),$

\ \ \ \ $[Um_{3i},Ua_{3j}]=0\ ,(0\leq i\neq j\leq 7),$

\ \ \ \ $[Um_{3i},U\alpha _{1}]=-\frac{1}{2}Ua_{3i}\ ,(0\leq i\leq 7),$

\ \ \ \ $[Um_{3i},U\alpha _{2}]=\frac{1}{2}Ua_{3i}\ ,(0\leq i\leq 7),$

\ \ \ \ $[Um_{3i},U\alpha _{3}]=0\ ,(0\leq i\leq 7),$

\ \ \ \ $[Um_{3i},Ub_{1j}]=\frac{1}{2}Sn(i+1,j+1)Ub_{2k},(k=Ca(i+1,j+1),\ 0\leq i,j\leq 7),$

\ \ \ \ $[Um_{3i},Ub_{2j}]=-\frac{1}{2}Sn(j+1,i+1)Ub_{1k},(k=Ca(j+1,i+1),\ 0\leq i,j\leq 7),$

\ \ \ \ $[Um_{3i},Ub_{3i}]=U\beta _{1}-U\beta _{2},(0\leq i\leq 7),$

\ \ \ \ $[Um_{3i},Ub_{3j}]=0\ ,(0\leq i\neq j\leq 7),$

\ \ \ \ $[Um_{3i},U\beta _{1}]=-\frac{1}{2}Ub_{3i}\ ,(0\leq i\leq 7),$

\ \ \ \ $[Um_{3i},U\beta _{2}]=\frac{1}{2}Ub_{3i}\ ,(0\leq i\leq 7),$

\ \ \ \ $[Um_{3i},U\beta _{3}]=0\ ,(0\leq i\leq 7).$

\bigskip

\emph{Proof. \ }We have the above Lie bracket operations with
calculations using Maxima.\ \ \ \ \emph{Q.E.D.}

\bigskip

\emph{Lemma 16.5.} \ We have the following Lie bracket operations.

\ \ \ \ $[Ut_{1i},Ua_{1j}]=0\ ,(0\leq i\neq j\leq 7),$

\ \ \ \ $[Ut_{1i},Ua_{1i}]=U\alpha _{2}+U\alpha _{3},(0\leq i\leq 7),$

\ \ \ \ $[Ut_{1i},Ua_{2j}]=\frac{1}{2}Sn(i+1,j+1)Ua_{3k},(k=Ca(i+1,j+1),\ 0\leq i,j\leq 7),$

\ \ \ \ $[Ut_{1i},Ua_{3j}]=\frac{1}{2}Sn(j+1,i+1)Ua_{2k},(k=Ca(j+1,i+1),\ 0\leq i,j\leq 7),$

\ \ \ \ $[Ut_{1i},U\alpha _{1}]=0\ ,(0\leq i\leq 7),$

\ \ \ \ $[Ut_{1i},U\alpha _{2}]=\frac{1}{2}Ua_{1i}\ ,(0\leq i\leq 7),$

\ \ \ \ $[Ut_{1i},U\alpha _{3}]=\frac{1}{2}Ua_{1i}\ ,(0\leq i\leq 7),$

\ \ \ \ $[Ut_{1i},Ub_{1j}]=0\ ,(0\leq i\neq j\leq 7),$

\ \ \ \ $[Ut_{1i},Ub_{1i}]=-U\beta _{2}-U\beta _{3},(0\leq i\leq 7),$

\ \ \ \ $[Ut_{1i},Ub_{2j}]=-\frac{1}{2}Sn(i+1,j+1)Ub_{3k},(k=Ca(i+1,j+1),\ 0\leq i,j\leq 7),$

\ \ \ \ $[t_{1i},Ub_{3j}]=-\frac{1}{2}Sn(j+1,i+1)Ub_{2k},(k=Ca(j+1,i+1),\ 0\leq i,j\leq 7),$

\ \ \ \ $[Ut_{1i},U\beta _{1}]=0\ ,(0\leq i\leq 7),$

\ \ \ \ $[Ut_{1i},U\beta _{2}]=-\frac{1}{2}Ub_{1i}\ ,(0\leq i\leq 7),$

\ \ \ \ $[Ut_{1i},U\beta _{3}]=-\frac{1}{2}Ub_{1i}\ ,(0\leq i\leq 7),$

\ \ \ \ $Ut_{2i},Ua_{1j}]=\frac{1}{2}Sn(j+1,i+1)Ua_{3k},(k=Ca(j+1,i+1),\ 0\leq i,j\leq 7),$

\ \ \ \ $[Ut_{2i},Ua_{2i}]=U\alpha _{1}+U\alpha _{3},(0\leq i\leq 7),$

\ \ \ \ $[Ut_{2i},Ua_{2j}]=0\ ,(0\leq i\neq j\leq 7),$

\ \ \ \ $[Ut_{2i},Ua_{3j}]=\frac{1}{2}Sn(i+1,j+1)Ua_{1k},(k=Ca(i+1,j+1),\ 0\leq i,j\leq 7),$

\ \ \ \ $[Ut_{2i},U\alpha _{1}]=\frac{1}{2}Ua_{2i}\ ,(0\leq i\leq 7)$,

\ \ \ \ $[Ut_{2i},U\alpha _{2}]=0\ ,(0\leq i\leq 7)$,

\ \ \ \ $[Ut_{2i},U\alpha _{3}]=\frac{1}{2}Ua_{2i}\ ,(0\leq i\leq 7),$

\ \ \ \ $[Ut_{2i},Ub_{1j}]=-\frac{1}{2}Sn(j+1,i+1)Ub_{3k},(k=Ca(j+1,i+1),\ 0\leq i,j\leq 7),$

\ \ \ \ $[Ut_{2i},Ub_{2i}]=-U\beta _{1}-U\beta _{3},(0\leq i\leq 7),$

\ \ \ \ $[Ut_{2i},Ub_{2j}]=0\ ,(0\leq i\neq j\leq 7),$

\ \ \ \ $Ut_{2i},Ub_{3j}]=\frac{1}{2}Sn(i+1,j+1)Ub_{1k},(k=Ca(i+1,j+1),\ 0\leq i,j\leq 7),$

\ \ \ \ $[Ut_{2i},U\beta _{1}]=-\frac{1}{2}Ub_{2i}\ ,(0\leq i\leq 7),$

\ \ \ \ $[Ut_{2i},U\beta _{2}]=0\ ,(0\leq i\leq 7),$

\ \ \ \ $[Ut_{2i},U\beta _{3}]=-\frac{1}{2}Ub_{2i}\ ,(0\leq i\leq 7),$

\ \ \ \ $[Ut_{3i},Ua_{1j}]=\frac{1}{2}Sn(i+1,j+1)Ua_{2k},(k=Ca(i+1,j+1),\ 0\leq i,j\leq 7),$

\ \ \ \ $[Ut_{3i},Ua_{2j}]=\frac{1}{2}Sn(j+1,i+1)Ua_{1k},(k=Ca(j+1,i+1),\ 0\leq i,j\leq 7),$

\ \ \ \ $[Ut_{3i},Ua_{3i}]=U\alpha _{1}+U\alpha _{2},(0\leq i\leq 7),$

\ \ \ \ $[Ut_{3i},Ua_{3j}]=0\ ,(0\leq i\neq j\leq 7),$

\ \ \ \ $[Ut_{3i},U\alpha _{1}]=\frac{1}{2}Ua_{3i}\ ,(0\leq i\leq 7),$

\ \ \ \ $[Ut_{3i},U\alpha _{2}]=\frac{1}{2}Ua_{3i}\ ,(0\leq i\leq 7),$

\ \ \ \ $[Ut_{3i},U\alpha _{3}]=0\ ,(0\leq i\leq 7),$

\ \ \ \ $ [Ut_{3i},Ub_{1j}]=-\frac{1}{2}Sn(i+1,j+1)Ub_{2k},(k=Ca(i+1,j+1),\ 0\leq i,j\leq 7),$

\ \ \ \ $ [Ut_{3i},Ub_{2j}]=-\frac{1}{2}Sn(j+1,i+1)Ub_{1k},(k=Ca(j+1,i+1),\ 0\leq i,j\leq 7),$

\ \ \ \ $[Ut_{3i},Ub_{3i}]=-U\beta _{1}-U\beta _{2},(0\leq i\leq 7),$

\ \ \ \ $[Ut_{3i},Ub_{3j}]=0\ ,(0\leq i\neq j\leq 7),$

\ \ \ \ $[Ut_{3i},U\beta _{1}]=-\frac{1}{2}Ub_{3i}\ ,(0\leq i\leq 7),$

\ \ \ \ $[Ut_{3i},U\beta _{2}]=-\frac{1}{2}Ub_{3i}\ ,(0\leq i\leq 7),$

\ \ \ \ $ [Ut_{3i},U\beta _{3}]=0\ ,(0\leq i\leq 7),$

\ \ \ \ $[U\tau _{1},Ua_{1j}]=-\frac{1}{2}Ua_{1j}\ ,(0\leq j\leq 7),$

\ \ \ \ $ [U\tau _{1},Ua_{2j}]=0$ $,(0\leq j\leq 7),$

\ \ \ \ $[U\tau _{1},Ua_{3j}]=\frac{1}{2}Ua_{3j}\ ,(0\leq j\leq 7),$

\ \ \ \ $[U\tau _{1},U\alpha _{1}]=U\alpha _{1}\ ,$

\ \ \ \ $[U\tau _{1},U\alpha _{2}]=0\ ,$

\ \ \ \ $[U\tau _{1},U\alpha _{3}]=-U\alpha _{3}\ ,$

\ \ \ \ $[U\tau _{1},Ub_{1j}]=\frac{1}{2}Ub_{1j}\ ,(0\leq j\leq 7),$

\ \ \ \ $ [U\tau _{1},Ub_{2j}]=0$ $,(0\leq j\leq 7),$

\ \ \ \ $[U\tau _{1},Ub_{3j}]=-\frac{1}{2}Ub_{3j}\ ,(0\leq j\leq 7),$

\ \ \ \ $ [U\tau _{1},U\beta _{1}]=-U\beta _{1}$ $,$

\ \ \ \ $[U\tau _{1},U\beta _{2}]=0\ ,$

\ \ \ \ $[U\tau _{1},U\beta _{3}]=U\beta _{3}$ $,$

\ \ \ \ $[U\tau _{2},Ua_{1j}]=0,(0\leq j\leq 7),$

\ \ \ \ $[U\tau _{2},Ua_{2j}]=-\frac{1}{2}Ua_{2j}\ ,(0\leq j\leq 7),$

\ \ \ \ $[U\tau _{2},Ua_{3j}]=\frac{1}{2}Ua_{3j}\ ,(0\leq j\leq 7),$

\ \ \ \ $[U\tau _{2},U\alpha _{1}]=0\ ,$

\ \ \ \ $[U\tau _{2},U\alpha _{2}]=U\alpha _{2}\ ,$

\ \ \ \ $[U\tau _{2},U\alpha _{3}]=-U\alpha _{3}\ ,$

\ \ \ \ $[U\tau _{2},Ub_{1j}]=0,(0\leq j\leq 7),$

\ \ \ \ $[U\tau _{2},Ub_{2j}]=\frac{1}{2}Ub_{2j}\ ,(0\leq j\leq 7),$

\ \ \ \ $[U\tau _{2},Ub_{3j}]=-\frac{1}{2}Ub_{3j}\ ,(0\leq j\leq 7),$

\ \ \ \ $[U\tau _{2},U\beta _{1}]=0\ ,$

\ \ \ \ $[U\tau _{2},U\beta _{2}]=-U\beta _{2}$ $,$

\ \ \ \ $[U\tau _{2},U\beta _{3}]=U\beta _{3}$\ .

\bigskip

\emph{Proof. \ }We have the above Lie bracket operations with
calculations using Maxima.\ \ \ \ \emph{Q.E.D.}

\bigskip

\emph{Lemma 16.6.} \ We have the following Lie bracket operations.

\ \ \ \ $[Ud_{ij},U\rho _{1}]=0$ $,(0\leq i<j\leq 7),$

\ \ \ \ $[Um_{1j},U\rho _{1}]=0$ $,(0\leq j\leq 7),$

\ \ \ \ $[Um_{2j},U\rho _{1}]=0$ $,(0\leq j\leq 7),$

\ \ \ \ $[Um_{3j},U\rho _{1}]=0$ $,(0\leq j\leq 7),$

\ \ \ \ $[Ut_{1j},U\rho _{1}]=0$ $,(0\leq j\leq 7),$

\ \ \ \ $[Ut_{2j},U\rho _{1}]=0$ $,(0\leq j\leq 7),$

\ \ \ \ $[Ut_{3j},U\rho _{1}]=0$ $,(0\leq j\leq 7),$

\ \ \ \ $[U\tau _{1},U\rho _{1}]=0$ $,$

\ \ \ \ $[U\tau _{2},U\rho _{1}]=0$ $,$

\ \ \ \ $[Ua_{1j},U\rho _{1}]=-\frac{2}{3}Ua_{1j}$ $,(0\leq j\leq 7),$

\ \ \ \ $[Ua_{2j},U\rho _{1}]=-\frac{2}{3}Ua_{2j}$ $,(0\leq j\leq 7),$

\ \ \ \ $[Ua_{3j},U\rho _{1}]=-\frac{2}{3}Ua_{3j}$ $,(0\leq j\leq 7),$

\ \ \ \ $[U\alpha _{i},U\rho _{1}]=-\frac{2}{3}U\alpha _{i}$ $,(1\leq i\leq 3),$

\ \ \ \ $[Ub_{1j},U\rho _{1}]=\frac{2}{3}Ub_{1j}$ $,(0\leq j\leq 7),$

\ \ \ \ $[Ub_{2j},U\rho _{1}]=\frac{2}{3}Ub_{2j}$ $,(0\leq j\leq 7),$

\ \ \ \ $[Ub_{3j},U\rho _{1}]=\frac{2}{3}Ub_{3j}$ $,(0\leq j\leq 7),$

\ \ \ \ $[U\beta _{i},U\rho _{1}]=\frac{2}{3}U\beta _{i}$ $,(1\leq i\leq 3),$

\bigskip

\emph{Proof. \ }We have the above Lie bracket operations with
calculations using Maxima.\ \ \ \ \emph{Q.E.D.}

\bigskip

\emph{Lemma 16.7.} \ \gr$_{7}^{\C}$ is simple.

\bigskip

\emph{Proof. }Case 1\emph{\ }for $Ra_{10}\in $\gR\ga$^{\C}.$

\noindent
By \emph{Lemma 16.3}, we have

\ \ \ \ $\{[Rd_{0j},Ra_{10}] \mid 1\leq j\leq 7\}=\{Ra_{1i} \mid 1\leq i\leq 7\}.$

\noindent
By \emph{Lemma 16.4}, we have \ 

\ \ \ \ $\{[Rm_{3i},Ra_{1j}] \mid 0\leq i\neq j\leq 7\}=\{Ra_{2i} \mid 0\leq i\leq 7\},$

\ \ \ \ $\{[Rm_{2i},Ra_{1j}] \mid 0\leq i\neq j\leq 7\}=\{Ra_{3i} \mid 0\leq i\leq 7\}.$

\noindent
By \emph{Lemma 16.2}, we have \ 

\ \ \ \ $\{[Ra_{1i},Rb_{1j}] \mid 0\leq i<j\leq 7\}=\{Rd_{ij} \mid 0\leq i<j\leq 7\},$

\ \ \ \ $\{[Ra_{1i},Rb_{1i}],[Ra_{2i},Rb_{2i}],[Ra_{3i},Rb_{3i}] \mid 0\leq i\leq 7\}=\{R\tau _{1},R\tau _{2},R\rho _{1}\}.$

\noindent
By \emph{Lemma 16.5}, we have

\ \ \ \ $\{[R\tau _{1},R\alpha _{1}],[R\tau _{1},R\alpha _{3}],[R\tau _{2},R\alpha _{2}]\}=\{R\alpha _{1},R\alpha _{2},R\alpha _{3}\}.$

\noindent
By \emph{Lemma 16.2}, we have \ 

\ \ \ \ $\{[Ra_{1i},R\beta _{2}],[Ra_{1i},R\beta _{3}] 
\mid 0\leq i\leq 7\}=\{Rm_{1i},Rt_{1i} \mid 0\leq
i\leq 7\},$

\ \ \ \ $\{[Ra_{2i},R\beta _{1}],[Ra_{2i},R\beta _{3}] \mid 0\leq i\leq 7\}=\{Rm_{2i},Rt_{2i} \mid 0\leq i\leq 7\},$

\ \ \ \ $\{[Ra_{3i},R\beta _{1}],[Ra_{3i},R\beta _{2}] \mid 0\leq i\leq 7\}=\{Rm_{3i},Rt_{3i} \mid 0\leq i\leq 7\}$.

\noindent
By \emph{Lemma 16.4}, we have

\ \ \ \ $\{[Rm_{1i},R\beta _{2}],[Rm_{2i},R\beta _{3}],[Rm_{3i},R\beta _{1}] \mid 0\leq i\leq 7\}$

\ \ \ \ \ \ \ \ $=\{Rb_{1i},Rb_{2i},Rb_{3i} \mid 0\leq i\leq 7\}.$

\noindent
By \emph{Lemma 16.5}, we have

\ \ \ \ $\{[R\tau _{1},R\beta _{1}],[R\tau _{1},R\beta _{3}],[R\tau _{2},R\beta _{2}]\}=\{R\beta _{1},R\beta _{2},R\beta _{3}\}.$

\noindent
Hence we have

\ \ \ \ $\{[Ra10,x],[[Ra10,x],y] \mid x,y \in $\gr$_{7}^{\C} \}= \ $\gr$_{7}^{\C}.$

Case 2 for $Ra_{1i}\in $\gR\ga$^{\C},(1\leq i\leq 7).$

\noindent
By \emph{Lemma 16.3}, we have

\ \ \ \ $\{[Rd_{0j},Ra_{1j}]$ $|$ $1\leq j\leq 7\}=\{Ra10\}.$

\noindent
Therefore,this applies to Case 1.

Case 3 for $Ra_{2i}\in $\gR\ga$^{\C},(0\leq i\leq 7).$

\noindent
By \emph{Lemma 16.4}, we have \ 

\ \ \ \ $\{[Rm_{3i},Ra_{2j}] \mid 0\leq i\neq j\leq 7\}=\{Ra_{1i} \mid 0\leq i\leq 7\}.$

\noindent
Therefore,this applies to Case 2.

Case 4 for $Ra_{3i}\in $\gR\ga$^{\C},(1\leq i\leq 7).$

\noindent
By \emph{Lemma 16.4}, we have \ 

\ \ \ \ $\{[Rm_{2i},Ra_{3j}] \mid 0\leq i\neq j\leq 7\}=\{Ra_{1i} \mid 0\leq i\leq 7\}.$

\noindent
Therefore,this applies to Case 2.

Case 5 for $R\alpha _{1},R\alpha _{2},R\alpha _{3}\in $\gR\ga$^{\C}.$

\noindent
By \emph{Lemma 16.5}, we have

\ \ \ \ $\{[Rt_{2i},R\alpha _{1}] \mid 0\leq i\leq 7\}=\{Ra_{2i} \mid 0\leq i\leq 7\},$

\ \ \ \ $\{[Rt_{1i},R\alpha _{2}] \mid 0\leq i\leq 7\}=\{Ra_{1i} \mid 0\leq i\leq 7\},$

\ \ \ \ $\{[Rt_{1i},R\alpha _{3}] \mid 0\leq i\leq 7\}=\{Ra_{1i} \mid 0\leq i\leq 7\}.$

\noindent
Therefore,this applies to Case 2 and Case 1.

Case 6 for $Rb_{1i},Rb_{2i},Rb_{3i},R\beta _{1},R\beta _{2},R\beta _{3}\in $\gR\gb$^{\C}.$

\noindent
Same as Case 1,Case 2,Case 3, Case 4 and Case 5.

Case 7 for $R\rho _{1}\in $\gR\gp$^{\C}.$

\noindent
By \emph{Lemma 16.6}, we have

\ \ \ \ $\{[Ra_{1i},R\rho _{1}] \mid 0\leq i\leq 7\}=\{Ra_{1i} \mid 0\leq i\leq 7\}.$

\noindent
So, this applies to Case 1.

Since \gr$_{7}^{\C}=\ $\gr$%
_{6}^{\C} \bigoplus $\gR\ga$^{\C} \bigoplus $\gR\gb$^{\C} \bigoplus $\gR\gp$^{\C} \ $, therefore \gr$_{7}^{\C}$ is simple. \ \ \ \ 
\emph{Q.E.D.}

\bigskip

\emph{Lemma 16.8.} \ The Killing form $B_{7}$ of the Lie algebra \gr$_{7}^{\C}$ is given by

\ \ \ \ \ \ \ \ \ \ \ \ \ \ \ \ \ $B_{7}(R_{1},R_{2})=\frac{3}{5}tr(R_{1}R_{2}),R_{1},R_{2}\in $\gr$_{7}^{\C}.$ \ 

\bigskip

\emph{Proof}. \ Since \gr$_{7}^{\C}$ is simple ,there exist $\kappa \in \C$ such that

\ \ \ \ \ \ \ \ \ \ \ \ \ \ \ \ \ $B_{7}(R_{1},R_{2})=\kappa tr(R_{1}R_{2}).$

\noindent
To determine this $\kappa ,$let $R_{1}=U\rho _{1}.$ $(adR_{1})^{2}$ is calculated as
follows.

\noindent
For $R_{2}=(R_{62}+A_{2}+B_{2}+\rho _{2})\in $\gr$_{7}^{\C},R_{62}\in $\gr$%
_{6}^{\C},A_{2}\in $\gR\ga$,B_{2}\in $\gR\gb$,\rho _{2}\in $\gR\gp$ ,$

\noindent
we have

\ \ \ \ \ \ \ \ \ \ \ \ \ \ $[R_{1},[R_{1},R_{2}]]=[R_{1},\frac{2}{3}A_{2}-\frac{2}{3}B_{2}]=\frac{4}{9}A_{2}+\frac{4}{9}B_{2}.$

\noindent
Hence

\ \ \ \ \ \ \ \ \ \ \ \ \ \ \ \ $B_{7}(R_{1},R_{1})=tr((adR_{1})^{2})=\frac{4}{9}%
\times (3\times 8+3)+\frac{4}{9}\times (3\times 8+3)=24.$

\noindent
On the other hand

\ \ \ \ \ \ \ \ \ \ \ \ \ \ $\ tr(U\rho _{1}.U\rho _{1})=40.$

\noindent
Therefore $\kappa =\frac{24}{40}=\frac{3}{5}.$ \ \ \ \ Q.E.D.

\bigskip

\emph{Lemma 16.9.} The rank of the Lie algebra \gr$_{7}^{\C}$
is 7. The roots of \gr$_{7}^{\C}$ relative 

\noindent
to some Cartan subalgebra of \gr$_{7}^{\C}$ are given by

$\ \ \ \ \ \ \ \ \ \ \ \ \ \ \ \ \pm (\lambda _{k}-\lambda _{l}),\pm
(\lambda _{k}+\lambda _{l}),0\leq k<l\leq 3,$

\ \ \ \ \ \ \ \ \ \ $\ \ \ \ \ \ \pm \lambda _{k}\pm \frac{1}{2}(\mu
_{1}+2\mu _{2}),0\leq k\leq 3,$

$\ \ \ \ \ \ \ \ \ \ \ \ \ \ \ \ \pm \frac{1}{2}(-\lambda _{0}-\lambda
_{1}+\lambda _{2}-\lambda _{3})\pm \frac{1}{2}(-2\mu _{1}-\mu _{2}),$

$\ \ \ \ \ \ \ \ \ \ \ \ \ \ \ \ \pm \frac{1}{2}(\lambda _{0}+\lambda
_{1}+\lambda _{2}-\lambda _{3})\pm \frac{1}{2}(-2\mu _{1}-\mu _{2}),$

\ \ \ \ \ \ \ \ \ \ \ \ \ $\ \ \ \pm \frac{1}{2}(-\lambda _{0}+\lambda
_{1}+\lambda _{2}+\lambda _{3})\pm \frac{1}{2}(-2\mu _{1}-\mu _{2}),$

$\ \ \ \ \ \ \ \ \ \ \ \ \ \ \ \ \pm \frac{1}{2}(\lambda _{0}-\lambda
_{1}+\lambda _{2}+\lambda _{3})\pm \frac{1}{2}(-2\mu _{1}-\mu _{2}),$

\ \ \ \ \ \ \ \ \ \ \ \ \ \ $\ \ \pm \frac{1}{2}(\lambda _{0}-\lambda
_{1}+\lambda _{2}-\lambda _{3})\pm \frac{1}{2}(\mu _{1}-\mu _{2}),$

$\ \ \ \ \ \ \ \ \ \ \ \ \ \ \ \ \pm \frac{1}{2}(-\lambda _{0}+\lambda
_{1}+\lambda _{2}-\lambda _{3})\pm \frac{1}{2}(\mu _{1}-\mu _{2}),$

\ \ \ \ \ \ \ \ \ \ \ \ \ \ \ $\ \pm \frac{1}{2}(\lambda _{0}+\lambda
_{1}+\lambda _{2}+\lambda _{3})\pm \frac{1}{2}(\mu _{1}-\mu _{2}),$

$\ \ \ \ \ \ \ \ \ \ \ \ \ \ \ \ \pm \frac{1}{2}(-\lambda _{0}-\lambda
_{1}+\lambda _{2}+\lambda _{3})\pm \frac{1}{2}(\mu _{1}-\mu _{2}),$

\ \ \ \ \ \ \ \ \ \ \ \ \ $\ \ \ \pm (\mu _{j}+\frac{2}{3}\nu ),0\leq j\leq
2,\pm (-\mu _{1}-\mu _{2}+\frac{2}{3}\nu )$

\ \ \ \ \ \ \ \ \ \ \ \ \ \ \ $\ \pm \lambda _{k}\pm (\frac{1}{2}\mu _{1}-%
\frac{2}{3}\nu ),0\leq k\leq 3,$

\ \ \ \ \ \ \ \ \ \ \ \ \ \ \ $\ \pm \frac{1}{2}(-\lambda _{0}-\lambda
_{1}+\lambda _{2}-\lambda _{3})\pm (\frac{1}{2}\mu _{2}-\frac{2}{3}\nu ),$

$\ \ \ \ \ \ \ \ \ \ \ \ \ \ \ \ \pm \frac{1}{2}(\lambda _{0}+\lambda
_{1}+\lambda _{2}-\lambda _{3})\pm (\frac{1}{2}\mu _{2}-\frac{2}{3}\nu ),$

\ \ \ \ \ \ \ \ \ \ \ \ \ $\ \ \ \pm \frac{1}{2}(-\lambda _{0}+\lambda
_{1}+\lambda _{2}+\lambda _{3})\pm (\frac{1}{2}\mu _{2}-\frac{2}{3}\nu ),$

$\ \ \ \ \ \ \ \ \ \ \ \ \ \ \ \ \pm \frac{1}{2}(\lambda _{0}-\lambda
_{1}+\lambda _{2}+\lambda _{3})\pm (\frac{1}{2}\mu _{2}-\frac{2}{3}\nu ),$

\ \ \ \ \ \ \ \ \ \ \ \ \ \ $\ \ \pm \frac{1}{2}(-\lambda _{0}+\lambda
_{1}-\lambda _{2}+\lambda _{3})\pm (-\frac{1}{2}\mu _{1}-\frac{1}{2}\mu _{2}-%
\frac{2}{3}\nu ),$

$\ \ \ \ \ \ \ \ \ \ \ \ \ \ \ \ \pm \frac{1}{2}(-\lambda _{0}+\lambda
_{1}+\lambda _{2}-\lambda _{3})\pm (-\frac{1}{2}\mu _{1}-\frac{1}{2}\mu _{2}-%
\frac{2}{3}\nu ),$

\ \ \ \ \ \ \ \ \ \ \ \ \ \ \ $\ \pm \frac{1}{2}(\lambda _{0}+\lambda
_{1}+\lambda _{2}+\lambda _{3})\pm (-\frac{1}{2}\mu _{1}-\frac{1}{2}\mu _{2}-%
\frac{2}{3}\nu ),$

$\ \ \ \ \ \ \ \ \ \ \ \ \ \ \ \ \pm \frac{1}{2}(-\lambda _{0}-\lambda
_{1}+\lambda _{2}+\lambda _{3})\pm (-\frac{1}{2}\mu _{1}-\frac{1}{2}\mu _{2}-%
\frac{2}{3}\nu ),$

\bigskip

\emph{Proof. \ }Let

\ \ \ \gh$=\left\{ h=h_{\delta }+H+V\in \text{\gr}_{7}^{\C}%
\middle| %
\begin{array}{c}
h_{\delta }=\sum\limits_{k=0}^{3}\lambda_{k}H_{k}=\sum\limits_{k=0}^{3}-\lambda _{k}iUd_{k4+k},\\
H=\mu _{1}U\tau _{1}+\mu_{2}U\tau _{2}, \\ 
V=\nu U\rho _{1},\lambda _{k},\mu _{j},\nu \in \C%
\end{array}%
\right\} ,$

\noindent
then \gh \ is an abelian subalgebra of \gr$_{7}^{\C} ( $it will
be a Cartan subalgebra of \gr$_{7}^{\C} ).$
 
 \noindent
That \gh \ is abelian is clear from \emph{Lemma 16.6 .}

$(1)$ The roots of \gr$_{6}^{\C}$ are also roots of \gr$%
_{7}^{\C} $. Indeed,let $\alpha $ be a root of \gr$_{6}^{\C}$ and $S\in 
$\gr$_{6}^{\C}\subset $\gr$_{7}^{\C}$ be
a root vector belong to $\alpha .$ Then

$[h,S]=[h_{\delta }+H+V,S]=[h_{\delta }+H,S]+[V,S]=\alpha (h)S,($Since S$\in 
$\gr$_{6}^{\C},$so $[V,S]=0).$ Hence $\alpha $ is a root of \gr$_{7}^{\C}.$

$(2)$ Let we put $Sa_{0j}:$

\ \ \ \ \ \ \ \ \ \ \ \ \ \ \ \ \ \ \ \ \ \ $Sa_{0j}=U\alpha_{j},(1\leq j\leq 3).$

\noindent
Then we have by \emph{Lemma 16.6}

\ \ \ \ \ \ \ \ \ $[h_{\delta }+H+V,Sa_{0j}]=\mu_{j}U\alpha_{j}+\frac{2}{3}%
vU\alpha_{j}=(\mu_{j}+\frac{2}{3}v)Sa_{0j},$

\ \ \ \ \ \ \ \ \ \ \ \ \ \ \ \ \ \ 
\ \ \ \ \ \ \ \ \ \ \ \ \ \ \ \ \ \ \ $(1\leq j\leq 3,\mu _{3}=-\mu _{1}-\mu_{2}).$

\noindent
Let we put $Sb_{0j}:$

\ \ \ \ \ \ \ \ \ \ \ \ \ \ \ \ \ \ \ \ \ \ $Sb_{0j}=U\beta _{j},(1\leq j\leq 3).$

\noindent
Then we have by \emph{Lemma 16.6}

\ \ \ \ \ \ \ \ \ $[h_{\delta }+H+V,Sb_{0j}]=-\mu _{j}U\beta _{j}-\frac{2}{3}%
vU\beta _{j}=-(\mu _{j}+\frac{2}{3}v)Sb_{0j},$

\ \ \ \ \ \ \ \ \ \ \ \ \ \ \ \ \ \ 
\ \ \ \ \ \ \ \ \ \ \ \ \ \ \ \ \ \ \ $(1\leq j\leq 3,\mu _{3}=-\mu _{1}-\mu
_{2}).$

$(3)$ Let we put $Sa_{1k}:$

\ \ \ \ \ \ \ \ \ \ \ \ \ \ \ \ \ \ \ \ \ \ $Sa_{1k}=Ua_{1k}+iUa_{14+k},(0%
\leq k\leq 3).$

\noindent
Then we have by \emph{Lemma 16.3} and\emph{\ Lemma 16.5}

\ \ \ \ \ \ \ \ \ $[h_{\delta }+H+V,Sa_{1k}]=(\lambda _{k}-\frac{1}{2}\mu
_{1}+\frac{2}{3}v)Sa_{1k},(0\leq k\leq 3).$

\noindent
Let we put $Sa_{4k}:$

\ \ \ \ \ \ \ \ \ \ \ \ \ \ \ \ \ \ \ \ \ \ $Sa_{4k}=Ua_{1k}-iUa_{14},(0\leq
k\leq 3).$

\noindent
Then we have by \emph{Lemma 16.3} and\emph{\ Lemma 16.5}

\ \ \ \ \ \ \ \ \ $[h_{\delta }+H+V,Sa_{4k}]=(-\lambda _{k}-\frac{1}{2}\mu
_{1}+\frac{2}{3}v)Sa_{4k},(0\leq k\leq 3).$

$(4)$ Let we put $Sb_{1k}:$

\ \ \ \ \ \ \ \ \ \ \ \ \ \ \ \ \ \ \ \ \ \ $Sb_{1k}=Ub_{1k}+iUb_{14+k},(0%
\leq k\leq 3).$

\noindent
Then we have by \emph{Lemma 16.3} and\emph{\ Lemma 16.5}

\ \ \ \ \ \ \ \ \ $[h_{\delta }+H+V,Sb_{1k}]=(\lambda _{k}+\frac{1}{2}\mu
_{1}-\frac{2}{3}v)Sb_{1k},(0\leq k\leq 3).$

\noindent
Let we put $Sb_{1k+4}:$

\ \ \ \ \ \ \ \ \ \ \ \ \ \ \ \ \ \ \ \ \ \ $Sb_{4k}=Ub_{1k}-iUb_{14+k},(0%
\leq k\leq 3).$

\noindent
Then we have by \emph{Lemma 16.3} and\emph{\ Lemma 16.5}

\ \ \ \ \ \ \ \ \ $[h_{\delta }+H+V,Sb_{4k}]=(-\lambda _{k}+\frac{1}{2}\mu
_{1}-\frac{2}{3}v)Sb_{4k},(0\leq k\leq 3).$

$(5)$ Let we put $Sa_{2k},Sa_{5k}:$

\ \ \ \ \ \ \ \ \ \ \ \ \ \ \ \ \ \ \ \ \ \ $Sa_{2k}=Ua_{2k}+iUa_{24+k},(0%
\leq k\leq 3),$

\ \ \ \ \ \ \ \ \ \ \ \ \ \ \ \ \ \ \ \ \ \ $Sa_{5k}=Ua_{2k}-iUa_{24+k},(0%
\leq k\leq 3).$

\noindent
Then we have by \emph{Lemma 16.3} and\emph{\ Lemma 16.5}

\ \ \ \ \ \ \ \ \ $[h_{\delta }+H+V,Sa_{20}]=(\frac{1}{2}(-\lambda
_{0}-\lambda _{1}+\lambda _{2}-\lambda _{3})-\frac{1}{2}\mu _{2}+\frac{2}{3}%
v)Sa_{20},$

\ \ \ \ \ \ \ \ \ $[h_{\delta }+H+V,Sa_{21}]=(\frac{1}{2}(\lambda
_{0}+\lambda _{1}+\lambda _{2}-\lambda _{3})-\frac{1}{2}\mu _{2}+\frac{2}{3}%
v)Sa_{21},$

\ \ \ \ \ \ \ \ \ $[h_{\delta }+H+V,Sa_{22}]=(\frac{1}{2}(-\lambda
_{0}+\lambda _{1}+\lambda _{2}+\lambda _{3})-\frac{1}{2}\mu _{2}+\frac{2}{3}%
v)Sa_{22},$

\ \ \ \ \ \ \ \ \ $[h_{\delta }+H+V,Sa_{23}]=(\frac{1}{2}(\lambda
_{0}-\lambda _{1}+\lambda _{2}+\lambda _{3})-\frac{1}{2}\mu _{2}+\frac{2}{3}%
v)Sa_{23},$

\ \ \ \ \ \ \ \ \ $[h_{\delta }+H+V,Sa_{50}]=(-\frac{1}{2}(-\lambda
_{0}-\lambda _{1}+\lambda _{2}-\lambda _{3})-\frac{1}{2}\mu _{2}+\frac{2}{3}%
v)Sa_{50},$

\ \ \ \ \ \ \ \ \ $[h_{\delta }+H+V,Sa_{51}]=(-\frac{1}{2}(\lambda
_{0}+\lambda _{1}+\lambda _{2}-\lambda _{3})-\frac{1}{2}\mu _{2}+\frac{2}{3}%
v)Sa_{51},$

\ \ \ \ \ \ \ \ \ $[h_{\delta }+H+V,Sa_{52}]=(-\frac{1}{2}(-\lambda
_{0}+\lambda _{1}+\lambda _{2}+\lambda _{3})-\frac{1}{2}\mu _{2}+\frac{2}{3}%
v)Sa_{52},$

\ \ \ \ \ \ \ \ \ $[h_{\delta }+H+V,Sa_{53}]=(-\frac{1}{2}(\lambda
_{0}-\lambda _{1}+\lambda _{2}+\lambda _{3})-\frac{1}{2}\mu _{2}+\frac{2}{3}%
v)Sa_{53}.$

(6)Let we put $Sb_{2k},Sb_{5k}:$

\ \ \ \ \ \ \ \ \ \ \ \ \ \ \ \ \ \ \ \ \ \ $Sb_{2k}=Ub_{2k}+iUb_{24+k},(0%
\leq k\leq 3),$

\ \ \ \ \ \ \ \ \ \ \ \ \ \ \ \ \ \ \ \ \ \ $Sb_{5k}=Ub_{2k}-iUb_{24+k},(0%
\leq k\leq 3).$

\noindent
Then we have by \emph{Lemma 16.3} and\emph{\ Lemma 16.5}

\ \ \ \ \ \ \ \ \ $[h_{\delta }+H+V,Sb_{20}]=(\frac{1}{2}(-\lambda
_{0}-\lambda _{1}+\lambda _{2}-\lambda _{3})+\frac{1}{2}\mu _{2}-\frac{2}{3}%
v)Sb_{20},$

\ \ \ \ \ \ \ \ \ $[h_{\delta }+H+V,Sb_{21}]=(\frac{1}{2}(\lambda
_{0}+\lambda _{1}+\lambda _{2}-\lambda _{3})+\frac{1}{2}\mu _{2}-\frac{2}{3}%
v)Sb_{21},$

\ \ \ \ \ \ \ \ \ $[h_{\delta }+H+V,Sb_{22}]=(\frac{1}{2}(-\lambda
_{0}+\lambda _{1}+\lambda _{2}+\lambda _{3})+\frac{1}{2}\mu _{2}-\frac{2}{3}%
v)Sb_{22},$

\ \ \ \ \ \ \ \ \ $[h_{\delta }+H+V,Sb_{23}]=(\frac{1}{2}(\lambda
_{0}-\lambda _{1}+\lambda _{2}+\lambda _{3})+\frac{1}{2}\mu _{2}-\frac{2}{3}%
v)Sb_{23},$

\ \ \ \ \ \ \ \ \ $[h_{\delta }+H+V,Sb_{50}]=(-\frac{1}{2}(-\lambda
_{0}-\lambda _{1}+\lambda _{2}-\lambda _{3})+\frac{1}{2}\mu _{2}-\frac{2}{3}%
v)Sb_{50},$

\ \ \ \ \ \ \ \ \ $[h_{\delta }+H+V,Sb_{51}]=(-\frac{1}{2}(\lambda
_{0}+\lambda _{1}+\lambda _{2}-\lambda _{3})+\frac{1}{2}\mu _{2}-\frac{2}{3}%
v)Sb_{51},$

\ \ \ \ \ \ \ \ \ $[h_{\delta }+H+V,Sb_{52}]=(-\frac{1}{2}(-\lambda
_{0}+\lambda _{1}+\lambda _{2}+\lambda _{3})+\frac{1}{2}\mu _{2}-\frac{2}{3}%
v)Sb_{52},$

\ \ \ \ \ \ \ \ \ $[h_{\delta }+H+V,Sb_{53}]=(-\frac{1}{2}(\lambda
_{0}-\lambda _{1}+\lambda _{2}+\lambda _{3})+\frac{1}{2}\mu _{2}-\frac{2}{3}%
v)Sb_{53}.$

$(7)$ Let we put $Sa_{3k},Sa6_{k}:$

\ \ \ \ \ \ \ \ \ \ \ \ \ \ \ \ \ \ \ \ \ \ $Sa_{3k}=Ua_{3k}+iUa_{34+k},(0%
\leq k\leq 3),$

\ \ \ \ \ \ \ \ \ \ \ \ \ \ \ \ \ \ \ \ \ \ $Sa6_{k}=Ua_{3k}-iUa_{34+k},(0%
\leq k\leq 3).$

\noindent
Then we have by \emph{Lemma 16.3} and\emph{\ Lemma 16.5}

\ \ \ \ \ \ \ \ \ $[h_{\delta }+H+V,Sa_{30}]=(\frac{1}{2}(-\lambda
_{0}+\lambda _{1}-\lambda _{2}+\lambda _{3})+\frac{1}{2}\mu _{1}+\frac{1}{2}%
\mu _{2}+\frac{2}{3}\nu )Sa_{30},$

\ \ \ \ \ \ \ \ \ $[h_{\delta }+H+V,Sa_{31}]=(\frac{1}{2}(-\lambda
_{0}+\lambda _{1}+\lambda _{2}-\lambda _{3})+\frac{1}{2}\mu _{1}+\frac{1}{2}%
\mu _{2}+\frac{2}{3}\nu )Sa_{31},$

\ \ \ \ \ \ \ \ \ $[h_{\delta }+H+V,Sa_{32}]=(\frac{1}{2}(\lambda
_{0}+\lambda _{1}+\lambda _{2}+\lambda _{3})+\frac{1}{2}\mu _{1}+\frac{1}{2}%
\mu _{2}+\frac{2}{3}\nu )Sa_{32},$

\ \ \ \ \ \ \ \ \ $[h_{\delta }+H+V,Sa_{33}]=(\frac{1}{2}(-\lambda
_{0}-\lambda _{1}+\lambda _{2}+\lambda _{3})+\frac{1}{2}\mu _{1}+\frac{1}{2}%
\mu _{2}+\frac{2}{3}\nu )Sa_{33},$

\ \ \ \ \ \ \ \ \ $[h_{\delta }+H+V,Sa_{60}]=(-\frac{1}{2}(-\lambda
_{0}+\lambda _{1}-\lambda _{2}+\lambda _{3})+\frac{1}{2}\mu _{1}+\frac{1}{2}%
\mu _{2}+\frac{2}{3}\nu )Sa_{60},$

\ \ \ \ \ \ \ \ \ $[h_{\delta }+H+V,Sa_{61}]=(-\frac{1}{2}(-\lambda
_{0}+\lambda _{1}+\lambda _{2}-\lambda _{3})+\frac{1}{2}\mu _{1}+\frac{1}{2}%
\mu _{2}+\frac{2}{3}\nu )Sa_{61},$

\ \ \ \ \ \ \ \ \ $[h_{\delta }+H+V,Sa_{62}]=(-\frac{1}{2}(\lambda
_{0}+\lambda _{1}+\lambda _{2}+\lambda _{3})+\frac{1}{2}\mu _{1}+\frac{1}{2}%
\mu _{2}+\frac{2}{3}\nu )Sa_{62},$

\ \ \ \ \ \ \ \ \ $[h_{\delta }+H+V,Sa_{63}]=(-\frac{1}{2}(-\lambda
_{0}-\lambda _{1}+\lambda _{2}+\lambda _{3})+\frac{1}{2}\mu _{1}+\frac{1}{2}%
\mu _{2}+\frac{2}{3}\nu )Sa_{63},$

$(8)$ Let we put $Sb_{3k},Sb_{6k}:$

\ \ \ \ \ \ \ \ \ \ \ \ \ \ \ \ \ \ \ \ \ \ $Sb_{3k}=Ub_{3k}+iUb_{34+k},(0%
\leq k\leq 3),$

\ \ \ \ \ \ \ \ \ \ \ \ \ \ \ \ \ \ \ \ \ \ $Sb_{6k}=Ub_{3k}-iUb_{34+k},(0%
\leq k\leq 3).$

\noindent
Then we have by \emph{Lemma 16.3} and\emph{\ Lemma 16.5}

\ \ \ \ \ \ \ \ \ $[h_{\delta }+H+V,Sb_{30}]=(\frac{1}{2}(-\lambda
_{0}+\lambda _{1}-\lambda _{2}+\lambda _{3})-\frac{1}{2}\mu _{1}-\frac{1}{2}%
\mu _{2}-\frac{2}{3}\nu )Sb_{30},$

\ \ \ \ \ \ \ \ \ $[h_{\delta }+H+V,Sb_{31}]=(\frac{1}{2}(-\lambda
_{0}+\lambda _{1}+\lambda _{2}-\lambda _{3})-\frac{1}{2}\mu _{1}-\frac{1}{2}%
\mu _{2}-\frac{2}{3}\nu )Sb_{31},$

\ \ \ \ \ \ \ \ \ $[h_{\delta }+H+V,Sb_{32}]=(\frac{1}{2}(\lambda
_{0}+\lambda _{1}+\lambda _{2}+\lambda _{3})-\frac{1}{2}\mu _{1}-\frac{1}{2}%
\mu _{2}-\frac{2}{3}\nu )Sb_{32},$

\ \ \ \ \ \ \ \ \ $[h_{\delta }+H+V,Sb_{33}]=(\frac{1}{2}(-\lambda
_{0}-\lambda _{1}+\lambda _{2}+\lambda _{3})-\frac{1}{2}\mu _{1}-\frac{1}{2}%
\mu _{2}-\frac{2}{3}\nu )Sb_{33},$

\ \ \ \ \ \ \ \ \ $[h_{\delta }+H+V,Sb_{60}]=(-\frac{1}{2}(\lambda
_{0}-\lambda _{1}+\lambda _{2}-\lambda _{3})-\frac{1}{2}\mu _{1}-\frac{1}{2}%
\mu _{2}-\frac{2}{3}\nu )Sb_{60},$

\ \ \ \ \ \ \ \ \ $[h_{\delta }+H+V,Sb_{61}]=(-\frac{1}{2}(-\lambda
_{0}+\lambda _{1}+\lambda _{2}-\lambda _{3})-\frac{1}{2}\mu _{1}-\frac{1}{2}%
\mu _{2}-\frac{2}{3}\nu )Sb_{61},$

\ \ \ \ \ \ \ \ \ $[h_{\delta }+H+V,Sb_{62}]=(-\frac{1}{2}(\lambda
_{0}+\lambda _{1}+\lambda _{2}+\lambda _{3})-\frac{1}{2}\mu _{1}-\frac{1}{2}%
\mu _{2}-\frac{2}{3}\nu )Sb_{62},$

\ \ \ \ \ \ \ \ \ $[h_{\delta }+H+V,Sb_{63}]=(-\frac{1}{2}(-\lambda
_{0}-\lambda _{1}+\lambda _{2}+\lambda _{3})-\frac{1}{2}\mu _{1}-\frac{1}{2}%
\mu _{2}-\frac{2}{3}\nu )Sb_{63}.$

\noindent
Hence we have roots of \gr$_{7}^{\C}$ and these associated root
vectors.

By direct calculations using Maxima, we also have $(1),(2),(3),(4),(5),(6)$,$(7)$ and $(8)$. \ \ \ \ \emph{Q.E.D.}

\bigskip

\emph{Theorem 16.10. }\ In the root system of \emph{Lemma 16.9} ,

\ \ \ \ \ \ \ \ \ \ \ \ \ \ \ $\alpha _{1}=\lambda _{0}-\lambda _{1},\alpha
_{2}=\lambda _{1}-\lambda _{2},\alpha _{3}=\lambda _{2}-\lambda _{3},$

\ \ \ \ \ \ \ \ \ \ \ \ \ \ \ $\alpha _{4}=\frac{1}{2}(-\lambda _{0}-\lambda
_{1}-\lambda _{2}+\lambda _{3})-\frac{1}{2}(2\mu _{1}+\mu _{2}),$

\ \ \ \ \ \ \ \ \ \ \ \ \ \ \ $\alpha _{5}=\frac{1}{2}(\lambda _{0}+\lambda
_{1}+\lambda _{2}+\lambda _{3})+\frac{1}{2}(\mu _{1}-\mu _{2}),$

\ \ \ \ \ \ \ \ \ \ \ \ \ \ \ $\alpha _{6}=\mu _{2}+\frac{2}{3}v,\alpha
_{7}=\mu _{1}+\mu _{2}-\frac{2}{3}v$

\noindent
is a fundamental root system of the Lie algebra \gr$_{7}^{\C}$ and

$\ \ \ \ \ \ \ \ \ \ \ \ \ \ \ \ \ \mu =\alpha _{1}+2\alpha _{2}+3\alpha
_{3}+4\alpha _{4}+3\alpha _{5}+2\alpha _{6}+2\alpha _{7}$

\noindent
is the highest root. The Dynkin diagram and extended Dynkin diagram of 
\gr$_{7}^{\C}$ are respectively given by

\begin{figure}[H]
\centering
\includegraphics[width=12cm, height=2.0cm]{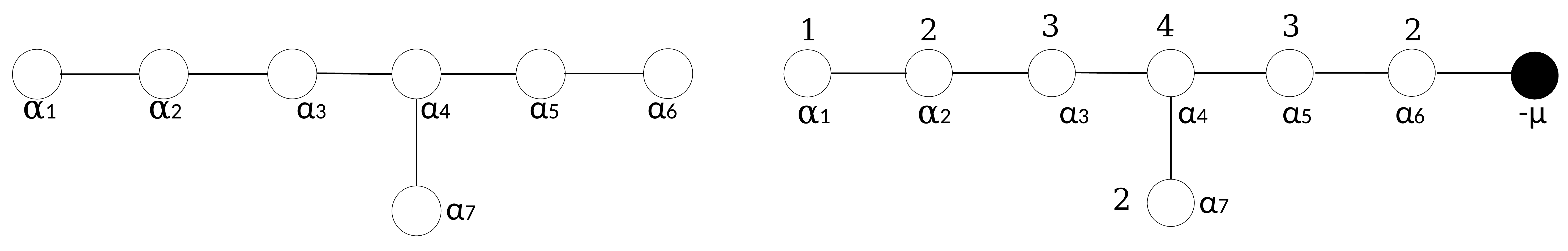}
\end{figure}

\ \ \ \ \emph{Proof}. \ In the following, the notation $n_{1}$ $n_{2}$ $%
n_{3} $ $n_{4}$ $n_{5}$ $n_{6}$ $n_{7}$ denotes the root
$n_{1}\alpha _{1}+n_{2}\alpha _{2}+n_{3}\alpha _{3}+n_{4}\alpha
_{4}+n_{5}\alpha _{5}+n_{6}\alpha _{6}+n_{7}\alpha _{7}.$ Now,all positive
roots of \gr$_{7}^{\C}$ are represented by

\begin{flushright}
$\ \ \ \ \ \ \ \ \ \lambda _{0}-\lambda _{1}=%
\begin{array}{lllllll}
1 & 0 & 0 & 0 & 0 & 0 & 0%
\end{array}%
,\lambda _{0}+\lambda _{1}=%
\begin{array}{lllllll}
1 & 2 & 2 & 2 & 2 & 1 & 1%
\end{array}%
,$

$\ \ \ \ \ \ \ \ \ \lambda _{0}-\lambda _{2}=%
\begin{array}{lllllll}
1 & 1 & 0 & 0 & 0 & 0 & 0%
\end{array}%
,\lambda _{0}+\lambda _{2}=%
\begin{array}{lllllll}
1 & 1 & 2 & 2 & 2 & 1 & 1%
\end{array}%
,$

$\ \ \ \ \ \ \ \ \ \lambda _{0}-\lambda _{3}=%
\begin{array}{lllllll}
1 & 1 & 1 & 0 & 0 & 0 & 0%
\end{array}%
,\lambda _{0}+\lambda _{3}=%
\begin{array}{lllllll}
1 & 1 & 1 & 2 & 2 & 1 & 1%
\end{array}%
,$

$\ \ \ \ \ \ \ \ \ \lambda _{1}-\lambda _{2}=%
\begin{array}{lllllll}
0 & 1 & 0 & 0 & 0 & 0 & 0%
\end{array}%
,\lambda _{1}+\lambda _{2}=%
\begin{array}{lllllll}
0 & 1 & 2 & 2 & 2 & 1 & 1%
\end{array}%
,$

$\ \ \ \ \ \ \ \ \ \lambda _{1}-\lambda _{3}=%
\begin{array}{lllllll}
0 & 1 & 1 & 0 & 0 & 0 & 0%
\end{array}%
,\lambda _{1}+\lambda _{3}=%
\begin{array}{lllllll}
0 & 1 & 1 & 2 & 2 & 1 & 1%
\end{array}%
,$

$\ \ \ \ \ \ \ \ \ \lambda _{2}-\lambda _{3}=%
\begin{array}{lllllll}
0 & 0 & 1 & 0 & 0 & 0 & 0%
\end{array}%
,\lambda _{2}+\lambda _{3}=%
\begin{array}{lllllll}
0 & 0 & 1 & 2 & 2 & 1 & 1%
\end{array}%
,$

$\ \ \ \ \ \ \ \ \ \lambda _{0}+\frac{1}{2}(\mu _{1}+2\mu _{2})=%
\begin{array}{lllllll}
1 & 1 & 1 & 1 & 1 & 1 & 1%
\end{array}%
,$

$\ \ \ \ \ \ \ \ \ \lambda _{1}+\frac{1}{2}(\mu _{1}+2\mu _{2})=%
\begin{array}{lllllll}
0 & 1 & 1 & 1 & 1 & 1 & 1%
\end{array}%
,$

$\ \ \ \ \ \ \ \ \ \lambda _{2}+\frac{1}{2}(\mu _{1}+2\mu _{2})=%
\begin{array}{lllllll}
0 & 0 & 1 & 1 & 1 & 1 & 1%
\end{array}%
,$

$\ \ \ \ \ \ \ \ \ \lambda _{3}+\frac{1}{2}(\mu _{1}+2\mu _{2})=%
\begin{array}{lllllll}
0 & 0 & 0 & 1 & 1 & 1 & 1%
\end{array}%
,$

$\ \ \ \ \ \ \ \ \ \lambda _{0}-\frac{1}{2}(\mu _{1}+2\mu _{2})=%
\begin{array}{lllllll}
1 & 1 & 1 & 1 & 1 & 0 & 0%
\end{array}%
,$

$\ \ \ \ \ \ \ \ \ \lambda _{1}-\frac{1}{2}(\mu _{1}+2\mu _{2})=%
\begin{array}{lllllll}
0 & 1 & 1 & 1 & 1 & 0 & 0%
\end{array}%
,$

$\ \ \ \ \ \ \ \ \ \lambda _{2}-\frac{1}{2}(\mu _{1}+2\mu _{2})=%
\begin{array}{lllllll}
0 & 0 & 1 & 1 & 1 & 0 & 0%
\end{array}%
,$

$\ \ \ \ \ \ \ \ \ \lambda _{3}-\frac{1}{2}(\mu _{1}+2\mu _{2})=%
\begin{array}{lllllll}
0 & 0 & 0 & 1 & 1 & 0 & 0%
\end{array}%
,$

$\ \ \ \ \ \ \ \ \ \ \frac{1}{2}(-\lambda _{0}-\lambda _{1}+\lambda
_{2}-\lambda _{3})+\frac{1}{2}(-2\mu _{1}-\mu _{2})=%
\begin{array}{lllllll}
0 & 0 & 1 & 1 & 0 & 0 & 0%
\end{array}%
,$

$\ \ \ \ \ \ \ \ \ \ \ \ \frac{1}{2}(\lambda _{0}+\lambda _{1}+\lambda
_{2}-\lambda _{3})+\frac{1}{2}(-2\mu _{1}-\mu _{2})=%
\begin{array}{lllllll}
1 & 2 & 3 & 3 & 2 & 1 & 1%
\end{array}%
,$

$\ \ \ \ \ \ \ \ \ \ \frac{1}{2}(-\lambda _{0}+\lambda _{1}+\lambda
_{2}+\lambda _{3})+\frac{1}{2}(-2\mu _{1}-\mu _{2})=%
\begin{array}{lllllll}
0 & 1 & 2 & 3 & 2 & 1 & 1%
\end{array}%
,$

$\ \ \ \ \ \ \ \ \ \ \ \frac{1}{2}(\lambda _{0}-\lambda _{1}+\lambda
_{2}+\lambda _{3})+\frac{1}{2}(-2\mu _{1}-\mu _{2})=%
\begin{array}{lllllll}
1 & 1 & 2 & 3 & 2 & 1 & 1%
\end{array}%
,$

\ \ $\ \ \ \ \ \ \ \ \frac{1}{2}(\lambda _{0}+\lambda _{1}-\lambda
_{2}+\lambda _{3})+\frac{1}{2}(-2\mu _{1}-\mu _{2})=%
\begin{array}{lllllll}
1 & 2 & 2 & 3 & 2 & 1 & 1%
\end{array}%
,$

$\ \ \ \ \ \ \ \ \ \ \frac{1}{2}(-\lambda _{0}-\lambda _{1}-\lambda
_{2}+\lambda _{3})+\frac{1}{2}(-2\mu _{1}-\mu _{2})=%
\begin{array}{lllllll}
0 & 0 & 0 & 1 & 0 & 0 & 0%
\end{array}%
,$

$\ \ \ \ \ \ \ \ \ \ \ \ \frac{1}{2}(\lambda _{0}-\lambda _{1}-\lambda
_{2}-\lambda _{3})+\frac{1}{2}(-2\mu _{1}-\mu _{2})=%
\begin{array}{lllllll}
1 & 1 & 1 & 1 & 0 & 0 & 0%
\end{array}%
,$

$\ \ \ \ \ \ \ \ \ \ \frac{1}{2}(-\lambda _{0}+\lambda _{1}-\lambda
_{2}-\lambda _{3})+\frac{1}{2}(-2\mu _{1}-\mu _{2})=%
\begin{array}{lllllll}
0 & 1 & 1 & 1 & 0 & 0 & 0%
\end{array}%
,$

$\ \ \ \ \ \ \ \ \ \ \ \ \ \ \ \frac{1}{2}(\lambda _{0}-\lambda _{1}+\lambda
_{2}-\lambda _{3})-\frac{1}{2}(\mu _{1}-\mu _{2})=%
\begin{array}{lllllll}
1 & 1 & 2 & 2 & 1 & 1 & 1%
\end{array}%
,$

$\ \ \ \ \ \ \ \ \ \ \ \ \ \ \ \frac{1}{2}(\lambda _{0}-\lambda _{1}-\lambda
_{2}+\lambda _{3})-\frac{1}{2}(\mu _{1}-\mu _{2})=%
\begin{array}{lllllll}
1 & 1 & 1 & 2 & 1 & 1 & 1%
\end{array}%
,$

$\ \ \ \ \ \ \ \ \ \ \ \ \ \ \ \frac{1}{2}(\lambda _{0}+\lambda _{1}+\lambda
_{2}+\lambda _{3})+\frac{1}{2}(\mu _{1}-\mu _{2})=%
\begin{array}{lllllll}
0 & 0 & 0 & 0 & 1 & 0 & 0%
\end{array}%
,$

$\ \ \ \ \ \ \ \ \ \ \ \ \ \ \ \frac{1}{2}(\lambda _{0}+\lambda _{1}-\lambda
_{2}-\lambda _{3})-\frac{1}{2}(\mu _{1}-\mu _{2})=%
\begin{array}{lllllll}
1 & 2 & 2 & 2 & 1 & 1 & 1%
\end{array}%
,$

$\ \ \ \ \ \ \ \ \ \ \ \ \ \frac{1}{2}(-\lambda _{0}+\lambda _{1}-\lambda
_{2}+\lambda _{3})-\frac{1}{2}(\mu _{1}-\mu _{2})=%
\begin{array}{lllllll}
0 & 1 & 1 & 2 & 1 & 1 & 1%
\end{array}%
,$

$\ \ \ \ \ \ \ \ \ \ \ \ \ \frac{1}{2}(-\lambda _{0}+\lambda _{1}+\lambda
_{2}-\lambda _{3})-\frac{1}{2}(\mu _{1}-\mu _{2})=%
\begin{array}{lllllll}
0 & 1 & 2 & 2 & 1 & 1 & 1%
\end{array}%
,$

$\ \ \ \ \ \ \ \ \ \ \ \ \ \ \ \frac{1}{2}(\lambda _{0}+\lambda _{1}+\lambda
_{2}+\lambda _{3})-\frac{1}{2}(\mu _{1}-\mu _{2})=%
\begin{array}{lllllll}
1 & 2 & 3 & 4 & 3 & 2 & 2%
\end{array}%
,$

$\ \ \ \ \ \ \ \ \ \ \ \ \ \frac{1}{2}(-\lambda _{0}-\lambda _{1}+\lambda
_{2}+\lambda _{3})-\frac{1}{2}(\mu _{1}-\mu _{2})=%
\begin{array}{lllllll}
0 & 0 & 1 & 2 & 1 & 1 & 1%
\end{array}%
,$

\ \ \ \ \ \ \ \ \ \ \ \ \ \ \ \ \ \ \ \ \ \ \ \ \ \ \ \ \ \ \ \ \ \ \ \ \ \
\ \ \ \ \ \ \ \ \ \ \ \ $-\mu _{1}-\frac{2}{3}v=%
\begin{array}{lllllll}
1 & 2 & 3 & 4 & 2 & 1 & 2%
\end{array}%
,$

\ \ \ \ \ \ \ \ \ \ \ \ \ \ \ \ \ \ \ \ \ \ \ \ \ \ \ \ \ \ \ \ \ \ \ \ \ \
\ \ \ \ \ \ \ \ \ \ \ \ \ \ $\mu _{2}+\frac{2}{3}v=%
\begin{array}{lllllll}
0 & 0 & 0 & 0 & 0 & 1 & 0%
\end{array}%
,$

\ \ \ \ \ \ \ \ \ \ \ \ \ \ \ \ \ \ \ \ \ \ \ \ \ \ \ \ \ \ \ \ \ \ \ \ \ \
\ \ \ \ \ \ \ $\mu _{1}+\mu _{2}-\frac{2}{3}v=%
\begin{array}{lllllll}
0 & 0 & 0 & 0 & 0 & 0 & 1%
\end{array}%
,$

\ \ \ \ \ \ \ \ \ \ \ \ \ \ \ \ \ \ \ \ \ \ \ \ \ \ \ \ \ \ \ \ \ \ \ \ \ \
\ \ \ $\lambda _{0}-\frac{1}{2}\mu _{1}+\frac{2}{3}v=%
\begin{array}{lllllll}
1 & 1 & 1 & 1 & 1 & 1 & 0%
\end{array}%
,$

\ \ \ \ \ \ \ \ \ \ \ \ \ \ \ \ \ \ \ \ \ \ \ \ \ \ \ \ \ \ \ \ \ \ \ \ \ \
\ \ \ $\lambda _{1}-\frac{1}{2}\mu _{1}+\frac{2}{3}v=%
\begin{array}{lllllll}
0 & 1 & 1 & 1 & 1 & 1 & 0%
\end{array}%
,$

\ \ \ \ \ \ \ \ \ \ \ \ \ \ \ \ \ \ \ \ \ \ \ \ \ \ \ \ \ \ \ \ \ \ \ \ \ \
\ \ \ $\lambda _{2}-\frac{1}{2}\mu _{1}+\frac{2}{3}v=%
\begin{array}{lllllll}
0 & 0 & 1 & 1 & 1 & 1 & 0%
\end{array}%
,$

\ \ \ \ \ \ \ \ \ \ \ \ \ \ \ \ \ \ \ \ \ \ \ \ \ \ \ \ \ \ \ \ \ \ \ \ \ \
\ \ \ $\lambda _{3}-\frac{1}{2}\mu _{1}+\frac{2}{3}v=%
\begin{array}{lllllll}
0 & 0 & 0 & 1 & 1 & 1 & 0%
\end{array}%
,$

\ \ \ \ \ \ \ \ \ \ \ \ \ \ \ \ \ \ \ \ \ \ \ \ \ \ \ \ \ \ \ \ \ \ \ \ \ \
\ \ \ $\lambda _{0}+\frac{1}{2}\mu _{1}-\frac{2}{3}v=%
\begin{array}{lllllll}
1 & 1 & 1 & 1 & 1 & 0 & 1%
\end{array}%
,$

\ \ \ \ \ \ \ \ \ \ \ \ \ \ \ \ \ \ \ \ \ \ \ \ \ \ \ \ \ \ \ \ \ \ \ \ \ \
\ \ \ $\lambda _{1}+\frac{1}{2}\mu _{1}-\frac{2}{3}v=%
\begin{array}{lllllll}
0 & 1 & 1 & 1 & 1 & 0 & 1%
\end{array}%
,$

\ \ \ \ \ \ \ \ \ \ \ \ \ \ \ \ \ \ \ \ \ \ \ \ \ \ \ \ \ \ \ \ \ \ \ \ \ \
\ \ \ $\lambda _{2}+\frac{1}{2}\mu _{1}-\frac{2}{3}v=%
\begin{array}{lllllll}
0 & 0 & 1 & 1 & 1 & 0 & 1%
\end{array}%
,$

\ \ \ \ \ \ \ \ \ \ \ \ \ \ \ \ \ \ \ \ \ \ \ \ \ \ \ \ \ \ \ \ \ \ \ \ \ \
\ \ \ $\lambda _{3}+\frac{1}{2}\mu _{1}-\frac{2}{3}v=%
\begin{array}{lllllll}
0 & 0 & 0 & 1 & 1 & 0 & 1%
\end{array}%
,$

$\ \ \ \ \ \ \ \ \ \ \ \ \ \ \ \frac{1}{2}(-\lambda _{0}-\lambda
_{1}+\lambda _{2}-\lambda _{3})+\frac{1}{2}\mu _{2}-\frac{2}{3}v=%
\begin{array}{lllllll}
0 & 0 & 1 & 1 & 0 & 0 & 1%
\end{array}%
,$

$\ \ \ \ \ \ \ \ \ \ \ \ \ \ \ \ \ \frac{1}{2}(\lambda _{0}+\lambda
_{1}+\lambda _{2}-\lambda _{3})+\frac{1}{2}\mu _{2}-\frac{2}{3}v=%
\begin{array}{lllllll}
1 & 2 & 3 & 3 & 2 & 1 & 2%
\end{array}%
,$

$\ \ \ \ \ \ \ \ \ \ \ \ \ \ \ \frac{1}{2}(-\lambda _{0}+\lambda
_{1}+\lambda _{2}+\lambda _{3})+\frac{1}{2}\mu _{2}-\frac{2}{3}v=%
\begin{array}{lllllll}
0 & 1 & 2 & 3 & 2 & 1 & 2%
\end{array}%
,$

$\ \ \ \ \ \ \ \ \ \ \ \ \ \ \ \ \ \frac{1}{2}(\lambda _{0}-\lambda
_{1}+\lambda _{2}+\lambda _{3})+\frac{1}{2}\mu _{2}-\frac{2}{3}v=%
\begin{array}{lllllll}
1 & 1 & 2 & 3 & 2 & 1 & 2%
\end{array}%
,$

$\ \ \ \ \ \ \ \ \ \ \ \ \ \ \ \ \ \frac{1}{2}(\lambda _{0}+\lambda
_{1}-\lambda _{2}+\lambda _{3})+\frac{1}{2}\mu _{2}-\frac{2}{3}v=%
\begin{array}{lllllll}
1 & 2 & 2 & 3 & 2 & 1 & 2%
\end{array}%
,$

$\ \ \ \ \ \ \ \ \ \ \ \ \ \ \ \frac{1}{2}(-\lambda _{0}-\lambda
_{1}-\lambda _{2}+\lambda _{3})+\frac{1}{2}\mu _{2}-\frac{2}{3}v=%
\begin{array}{lllllll}
0 & 0 & 0 & 1 & 0 & 0 & 1%
\end{array}%
,$

$\ \ \ \ \ \ \ \ \ \ \ \ \ \ \ \ \ \frac{1}{2}(\lambda _{0}-\lambda
_{1}-\lambda _{2}-\lambda _{3})+\frac{1}{2}\mu _{2}-\frac{2}{3}v=%
\begin{array}{lllllll}
1 & 1 & 1 & 1 & 0 & 0 & 1%
\end{array}%
,$

$\ \ \ \ \ \ \ \ \ \ \ \ \ \ \ \frac{1}{2}(-\lambda _{0}+\lambda
_{1}-\lambda _{2}-\lambda _{3})+\frac{1}{2}\mu _{2}-\frac{2}{3}v=%
\begin{array}{lllllll}
0 & 1 & 1 & 1 & 0 & 0 & 1%
\end{array}%
,$

$\ \ \ \ \ \ \frac{1}{2}(\lambda _{0}-\lambda _{1}+\lambda _{2}-\lambda
_{3})-\frac{1}{2}(\mu _{1}+\mu _{2})-\frac{2}{3}v=%
\begin{array}{lllllll}
1 & 1 & 2 & 2 & 1 & 0 & 1%
\end{array}%
,$

$\ \ \ \ \ \ \ \frac{1}{2}(\lambda _{0}-\lambda _{1}-\lambda _{2}+\lambda
_{3})-\frac{1}{2}(\mu _{1}+\mu _{2})-\frac{2}{3}v=%
\begin{array}{lllllll}
1 & 1 & 1 & 2 & 1 & 0 & 1%
\end{array}%
,$

$\ \ \ \ \ \ \ \frac{1}{2}(\lambda _{0}+\lambda _{1}+\lambda _{2}+\lambda
_{3})+\frac{1}{2}(\mu _{1}+\mu _{2})+\frac{2}{3}v=%
\begin{array}{lllllll}
0 & 0 & 0 & 0 & 1 & 1 & 0%
\end{array}%
,$

$\ \ \ \ \ \ \ \frac{1}{2}(\lambda _{0}+\lambda _{1}-\lambda _{2}-\lambda
_{3})-\frac{1}{2}(\mu _{1}+\mu _{2})-\frac{2}{3}v=%
\begin{array}{lllllll}
1 & 2 & 2 & 2 & 1 & 0 & 1%
\end{array}%
,$

$\ \ \ \ \frac{1}{2}(-\lambda _{0}+\lambda _{1}-\lambda _{2}+\lambda _{3})-%
\frac{1}{2}(\mu _{1}+\mu _{2})-\frac{2}{3}v=%
\begin{array}{lllllll}
0 & 1 & 1 & 2 & 1 & 0 & 1%
\end{array}%
,$

$\ \ \ \ \frac{1}{2}(-\lambda _{0}+\lambda _{1}+\lambda _{2}-\lambda _{3})-%
\frac{1}{2}(\mu _{1}+\mu _{2})-\frac{2}{3}v=%
\begin{array}{lllllll}
0 & 1 & 2 & 2 & 1 & 0 & 1%
\end{array}%
,$

$\ \ \ \ \ \ \ \frac{1}{2}(\lambda _{0}+\lambda _{1}+\lambda _{2}+\lambda
_{3})-\frac{1}{2}(\mu _{1}+\mu _{2})-\frac{2}{3}v=%
\begin{array}{lllllll}
1 & 2 & 3 & 4 & 3 & 1 & 2%
\end{array}%
,$

$\ \ \ \ \ \frac{1}{2}(-\lambda _{0}-\lambda _{1}+\lambda _{2}+\lambda _{3})-%
\frac{1}{2}(\mu _{1}+\mu _{2})-\frac{2}{3}v=%
\begin{array}{lllllll}
0 & 0 & 1 & 2 & 1 & 0 & 1%
\end{array}%
,$
\end{flushright}

\noindent
Hence $\Pi =\{\alpha _{1},\alpha _{2},\alpha _{3},\alpha _{4},\alpha
_{5},\alpha _{6},\alpha _{7}\}$ is a fundamental root system of \gr$%
_{7}^{\C}$ . The real part \gh$_{R}$ of \gh \ is

\ \ \ \gh$_{R}=\left\{ h=h_{\delta }+H+V\in \text{\gr}_{7}^{\C}%
\begin{tabular}{c|}
\\ 
\\ 
\end{tabular}%
\begin{array}{c}
h_{\delta }=\sum\limits_{k=0}^{3}\lambda
_{k}H_{k}=\sum\limits_{k=0}^{3}-\lambda _{k}iUd_{k4+k},\\
H=\mu _{1}U\tau _{1}+\mu_{2}U\tau _{2}, \\ 
V=\nu U\rho _{1},\lambda _{k},\mu _{j},\nu \in \R%
\end{array}%
\right\} .$

The Killing form $B_{7}$ of \gr$_{7}^{\C}$ is $B_{7}(R_{1},R_{2})=\frac{3}{5}%
tr(R_{1}R_{2})$ (\emph{Lemma 16.8} ), so that the Killing form $B_{7}$ of \gr$_{7}^{\C}$ on \gh$_{R}$ is given by

$\ \ \ \ \ \ \ \ B_{7}(h,h^{\prime
})=6(6\sum\limits_{k=0}^{3}\lambda _{k}\lambda _{k}^{\prime }+6\mu _{1}\mu
_{1}^{\prime }+3\mu _{1}\mu _{2}^{\prime }+3\mu _{2}\mu _{1}^{\prime }+6\mu
_{2}\mu _{2}^{\prime }+4vv^{\prime }),$

\noindent
for $h=\sum\limits_{k=0}^{3}\lambda _{k}H_{k}+\mu_{1}U\tau_{1}+\mu_{2}U\tau
2+vU\rho _{1},h^{\prime }=\sum\limits_{k=0}^{3}\lambda _{k}^{\prime }H_{k}+\mu
_{1}^{\prime }U\tau _{1}+\mu_{2}^{\prime }U\tau _{2}$

$+v^{\prime }U\rho _{1}\in $\gh$_{R}.$

\noindent
Indeed,by calculate with Maxima we have

\ \ \ \ $B_{7}(h,h^{\prime })=36\lambda _{0}\lambda _{0}^{\prime }+36\lambda
_{1}\lambda _{1}^{\prime }+36\lambda _{2}\lambda _{2}^{\prime }+36\lambda
_{3}\lambda _{3}^{\prime }$

\ \ \ \ \ \ \ \ \ \ \ \ \ \ \ \ \ $+36\mu _{1}\mu _{1}^{\prime }+18\mu _{1}\mu
_{2}^{\prime }+18\mu _{2}\mu _{1}^{\prime }+36\mu _{2}\mu _{2}^{\prime
}+24vv^{\prime }.$

Now,the canonical elements $H\alpha _{i}\in $\gh$_{R}$ corresponding
to $\alpha _{i}$ $(B_{7}(H\alpha ,H)=\alpha (H),H\in $\gh$_{R}$%
\textbf{) }are determined as follows.

\ \ \ \ $H_{\alpha _{1}}=\frac{1}{36}%
(H_{0}-H_{1}),H_{\alpha _{2}}=\frac{1}{36}(H_{1}-H_{2}),H_{\alpha _{3}}=\frac{1}{36}(H_{2}-H_{3}),$

\ \ \ \ $H_{\alpha 4}=\frac{1}{72}((-H_{0}-H_{1}-H_{2}+H_{3})-2U\tau _{1}),$

\ \ \ \ $H_{\alpha 5}=\frac{1}{72}((H_{0}+H_{1}+H_{2}+H_{3})+2U\tau _{1}-2U\tau _{2}),$

\ \ \ \ $H_{\alpha 6}=\frac{1}{54}(-U\tau_{1}+2U\tau _{2}+\frac{3}{2}vU\rho _{1}),$

\ \ \ \ $H_{\alpha 7}=\frac{1}{54}(U\tau_{1}+U\tau_{2}-\frac{3}{2}vU\rho _{1}).$

\noindent
Thus we have

\ \ \ \ $(\alpha _{1},\alpha _{1})=B_{7}(H_{\alpha1},H_{\alpha _{1}})=\frac{1}{18},$

\ \ \ \ $(\alpha _{i},\alpha _{i})=\frac{1}{18},(i=2,3,4,5,6,7),$

\ \ \ \ $(\alpha _{1},\alpha _{2})=(\alpha _{2},\alpha
_{3})=(\alpha _{3},\alpha _{4})=(\alpha _{4},\alpha _{5})=(\alpha
_{4},\alpha _{7})=(\alpha _{5},\alpha _{6})=-\frac{1}{36},$

\ \ \ \ $(\alpha _{i},\alpha _{j})=0$, otherwise,

\ \ \ \ $(-\mu ,-\mu )=\frac{1}{18},(-\mu ,\alpha
_{6})=-\frac{1}{36},(-\mu ,\alpha _{i})=0,(i=1,2,3,5,7),$

\noindent
using them,we can draw the Dynkin diagram and the extended Dynkin diagram of 
\gr$_{7}^{\C}.$\ \ \ \ \emph{Q.E.D.}

\bigskip

\emph{Corollary 16.11.} For a $248\times 248$ matrix $X$,let $X|_{7}$ be the matrix in which the $133\times 133$
elements in the upper left corner are clipped from  $X$. Furthermore, let \gr\gd$|_{7}=\{Rd|_{7} \mid Rd\in $\gr\gd$\}$ , 
\gR\gm$|_{7}=\{Rm|_{7} \mid Rm\in $\gR\gm$\}$,\gR\gt$|_{7}=\{Rt|_{7} \mid Rt\in $\gR\gt$\}$,\gR\ga$|_{7}=\{Ra|_{7} \mid Ra\in $\gR\ga$\}$,
\gR\gb$|_{7}=\{Rb|_{7} \mid Rb\in $\gR\gb$\}$ and \gR\gp$|_{7}=\{R\rho_{1}|_{7} \mid R\rho_{1}\in $\gR\gp$\}$.
\emph{Theorem 16.10} holds for \gr$_{7}^{\C}|_{7}=$\gr\gd$^{\C}|_{7}\oplus $\gR\gm$^{\C}|_{7}\oplus $\gR\gt$^{\C}|_{7}$

\noindent
$\oplus $\gR\ga$^{\C}|_{7}$
$\oplus $\gR\gb$^{\C}|_{7}\oplus $\gR\gp$^{\C}|_{7}$ 
as well. However, the Killing form $B_{7}(R_{1},R_{2})=tr(R_{1}R_{2}) \ \ (R_{1},R_{2} \in $\gr$_{7}^{\C}|_{7}).$

\bigskip

\emph{Proof} We have the above with calculations using Maxima.\ \ \ \ \emph{Q.E.D.}

\bigskip

\emph{Definition 16.12.} \ We define followings,

\ge\ga$^{0}=\{c_{k}^{0}(U\alpha _{k}-U\beta _{k})+ic_{k}^{1}(U\alpha _{k}+U\beta
_{k}) \mid c_{k}^{0},c_{k}^{1}\in \R,1\leq k\leq 3\},$

\ \ \ \ \ \ $=\{\alpha_{k}U\alpha_{k}-\iota\alpha_{k}U\beta_{k} \mid \alpha_{k}\in \C,1\leq k\leq 3\},$

\ \ \ \ \ \ \ \ where $\iota $ is the complex conjugation of $\C$,

\ge\ga$^{1}=%
\{c_{kj}^{0}(Ua_{kj}-Ub_{kj})+ic_{kj}^{1}(Ua_{kj}+Ub_{kj}) \mid c_{kj}^{0},c_{kj}^{1}\in \R,$

\ \ \ \ \ \ \ \ \ \ \ \ \ \ \ \ \ \ \ \ \ \ \ \ \ \ \ \ \ \ \ \ \ \ \ \ \ \ \ \ \ \ \ \ \ \ \ \ \ \ $1\leq k\leq 3,0\leq j\leq 7\}$,

\ \ \ \ \ \ $=\{a_{kj}Ua_{kj}-\iota a_{kj}Ub_{kj} \mid a_{kj}\in \C,1\leq k\leq 3,0\leq j\leq 7\}$,

\ge\ga$=$\ge\ga$^{0}\oplus $\ge\ga$^{1},$

\ge\gp$=\{\rho _{1}U\rho_{1} \mid \rho _{1}\in i\R\}.$

\noindent
Also we define as follows.

\ge\gd$_{7}=\{X|_{7} \mid X\in $\ge\gd$\},$

\ge\gm$_{7}=\{X|_{7} \mid X\in $\ge\gm$\},$

\ge\gt$^{0}_{7}=\{X|_{7} \mid X\in $\ge\gt$^{0}\}, $ \ge\gt$^{1}_{7}=\{X|_{7} \mid X\in $\ge\gt$^{1}\}, $ \ge\gt$_{7}=\{X|_{7} \mid X\in $\ge\gt$\},$

\ge\ga$^{0}_{7}=\{X|_{7} \mid X\in $\ge\ga$^{0}\}, $ \ge\ga$^{1}_{7}=\{X|_{7} \mid X\in $\ge\ga$^{1}\}, $ \ge\ga$_{7}=\{X|_{7} \mid X\in $\ge\ga$\},$

\ge\gp$_{7}=\{X|_{7} \mid X\in $\ge\gp$ \}.$

\bigskip

\emph{Theorem 16.13.}  Let we put

\ge$_{7}=$\ge\gd$_{7}\oplus $\ge\gm$_{7}\oplus $\ge\gt$_{7}\oplus $\ge\ga$_{7}\oplus $\ge\gp$ _{7}$ . \ 

\noindent
Then \ge$_{7}$ is a compact exceptional simple Lie algebra of type $E_{7}.$

\bigskip

\emph{Proof.} \ By \emph{Lemma 16.2,16.3, 16.4, 16.5}, and \emph{16.6}, \ge$_{7}$ is a Lie algebra.
And by calculations using Maxima,we have as follows.

\noindent
For $X_{1}\in $\ge\gd$_{7},X_{2}\in $\ge\gm$_{7},X_{3}\in $\ge\gt$^{0}_{7},X_{4}\in $\ge\gt$^{1}_{7},X_{5}\in
$\ge\ga$^{0}_{7},X_{6}\in $\ge\ga$^{1}_{7}$,
$X_{7}\in $\ge\gp$ _{7}$,

$B_{7}(X_{1},X_{1})$

$=-36(d_{67}^{2}+d_{57}^{2}+d_{56}^{2}+d_{47}^{2}+d_{46}%
^{2}+d_{45}^{2}+d_{37}^{2}$

\ \ \ \ \ \ $\ +d_{36}^{2}+d_{35}^{2}+d_{34}^{2}+d_{27}^{2}+d_{26}%
^{2}+d_{25}^{2}+d_{24}^{2}$

\ \ \ \ \ \ $\ +d_{23}^{2}+d_{17}^{2}+d_{16}^{2}+d_{15}^{2}+d_{14}%
^{2}+d_{13}^{2}+d_{12}^{2}$

\ \ \ \ \ \ $\ +d_{07}^{2}+d_{06}^{2}+d_{05}^{2}+d_{04}^{2}+d_{03}%
^{2}+d_{02}^{2}+d_{01}^{2}),$

$B_{7}(X_{2},X_{2})$

$=-36(m_{37}^{2}+m_{36}^{2}+m_{35}^{2}+m_{34}^{2}+m_{33}%
^{2}+m_{32}^{2}+m_{31}^{2}+m_{30}^{2}$

\ \ $\ +m_{27}^{2}+m_{26}^{2}+m_{25}^{2}+m_{24}^{2}+m_{23}%
^{2}+m_{22}^{2}+m_{21}^{2}+m_{20}^{2}$

\ \ $\ +m_{17}^{2}+m_{16}^{2}+m_{15}^{2}+m_{14}^{2}+m_{13}%
^{2}+m_{12}^{2}+m_{11}^{2}+m_{10}^{2}),$

$B_{7}(X_{3},X_{3})=-36(\tau _{1}^{2}+3\tau _{2}^{2})$

$B_{7}(X_{4},X_{4})$

$=-36(t_{37}^{2}+t_{36}^{2}+t_{35}^{2}+t_{34}^{2}+t_{33}%
^{2}+t_{32}^{2}+t_{31}^{2}+t_{30}^{2}$

\ \ \ \ \ \ $\ +t_{27}^{2}+t_{26}^{2}+t_{25}^{2}+t_{24}^{2}+t_{23}%
^{2}+t_{22}^{2}+t_{21}^{2}+t_{20}^{2}$

\ \ \ \ \ $\ \ $+$t_{17}^{2}+t_{16}^{2}+t_{15}^{2}+t_{14}^{2}+t_{13}%
^{2}+t_{12}^{2}+t_{11}^{2}+t_{10}^{2}),$

$B_{7}(X_{5},X_{5})=-72(\iota \alpha _{3}\alpha _{3}+\iota \alpha _{2}\alpha _{2}+\iota \alpha
_{1}\alpha _{1}),$

$B_{7}(X_{6},X_{6})$

=$-144(\iota a_{37}a_{37}+\iota a_{36}a_{36}+\iota a_{35}a_{35}+\iota a_{34}a_{34}+\iota
a_{33}a_{33}+\iota a_{32}a_{32}+\iota a_{31}a_{31}$

\ \ \ \ $+\iota a_{30}a_{30}+\iota a_{27}a_{27}+\iota a_{26}a_{26}+\iota a_{25}a_{25}+\iota a_{24}a_{24}+\iota a_{23}a_{23}+\iota a_{22}a_{22}$

\ \ \ \ $+\iota a_{21}a_{21}+\iota a_{20}a_{20}+\iota a_{17}a_{17}+\iota a_{16}a_{16}+\iota a_{15}a_{15}+\iota a_{14}a_{14}+\iota a_{13}a_{13}$

\ \ \ \ $+\iota a_{12}a_{12}+\iota a_{11}a_{11}+\iota a_{10}a_{10}),$

$B_{7}(X_{7},X_{7})=-24\rho _{1}^{2},$

$B_{7}(Xi,Xj)=0,(i\neq j).$

\noindent
Therfor we have $B_{7}(X,X)<0$,for $^{\forall }X\neq 0,X\in $\ge$_{7}.$ 
Then \ge$_{7}$ is compact.
\ \ \emph{Q.E.D.}

\bigskip

\section{The exceptional simle Lie algebra \gr$_{8}^{\C}$  of type $E_{8}$}

\bigskip

\ \ \ \ \emph{Definition 17.1.} \ We define the followings for elements $%
Rx_{ij},R\chi _{k}\in $\gR\gx$^{\C},$

\noindent
$Ry_{ij},R\gamma _{k}\in $\gR\gy$^{\C},R\xi _{1}\in $\gR\gk$^{\C} ,R\eta _{1}\in $\gR\gi$^{\C} ,Rz_{ij},R\mu  _{k}\in
$\gR\gz$^{\C},Rw_{ij},R\psi _{k}\in $\gR\gw$^{\C}$,

\noindent
$R\zeta _{1}\in $\gR\gl$^{\C} ,R\omega _{1}\in $\gR\go$^{\C} ,Rr_{1}\in $\gR\gr$^{\C},Rs_{1}\in $\gR\gs$^{\C},Ru_{1}\in $\gR\gu$^{\C}.$

\ \ \ \ \ $\ \ \ \ \ \ \ \ \ \ \ \ \ Ux_{ij}=Rx_{ij}/x_{ij}\ \ (1\leq i\leq 3,0\leq
j\leq 7),$

\ \ \ \ \ \ \ \ \ \ \ \ \ \ \ \ $\ \ U\chi _{k}=R\chi _{k}/\chi _{k}$ \ $(1\leq k\leq
3),$

\ \ \ \ \ $\ \ \ \ \ \ \ \ \ \ \ \ \ Uy_{ij}=Ry_{ij}/y_{ij}\ \ (1\leq i\leq 3,0\leq
j\leq 7),$

\ \ \ \ \ \ \ \ \ \ \ \ \ \ \ \ $\ \ U\gamma _{k}=R\gamma _{k}/\gamma _{k}$ \ $(1\leq
k\leq 3),$

\ \ \ \ \ \ \ \ \ \ \ \ \ \ \ \ $\ \ U\xi _{1}=R\xi _{1}/\xi _{1}$ $,$

\ \ \ \ \ \ \ \ \ \ \ \ \ \ \ \ $\ \ U\eta _{1}=R\eta _{1}/\eta _{1}$ $,$

\ \ \ \ \ $\ \ \ \ \ \ \ \ \ \ \ \ \ Uz_{ij}=Rz_{ij}/z_{ij}\ \ (1\leq i\leq 3,0\leq
j\leq 7),$

\ \ \ \ \ \ \ \ \ \ \ \ \ \ \ \ $\ \ U\mu_{k}=R\mu  _{k}/\mu  _{k}$ \ $(1\leq k\leq
3), $

\ \ \ \ \ $\ \ \ \ \ \ \ \ \ \ \ \ \ Uw_{ij}=Rw_{ij}/w_{ij}\ \ (1\leq i\leq 3,0\leq
j\leq 7),$

\ \ \ \ \ \ \ \ \ \ \ \ \ \ \ \ $\ \ U\psi_{k}=R\psi _{k}/\psi _{k}$ \ $(1\leq k\leq
3),$

\ \ \ \ \ \ \ \ \ \ \ \ \ \ \ \ $\ \ U\zeta_{1}=R\zeta _{1}/\zeta _{1}$ $,$

\ \ \ \ \ \ \ \ \ \ \ \ \ \ \ \ $\ \ U\omega_{1}=R\omega _{1}/\omega _{1}$ $,$

\ \ \ \ \ \ \ \ \ \ \ \ \ \ \ \ $\ \ Ur_{1}=Rr_{1}/r_{1}$ $,$

\ \ \ \ \ \ \ \ \ \ \ \ \ \ \ \ $\ \ Us_{1}=Rs_{1}/s_{1}$ $,$

\ \ \ \ \ \ \ \ \ \ \ \ \ \ \ \ $\ \ Uu_{1}=Ru_{1}/u_{1}$ $.$

\bigskip

\emph{Lemma 17.2. }\ We have the following Lie bracket operations.

\ \ \ \ $[Ud_{ij},Ux_{1k}]=-Ux_{1j}$ $($in case of $k=i),$

\ \ \ \ \ \ \ \ \ \ \ \ \ \ \ \ \ \ \ \ $=Ux_{1i}$ $($in case of $k=j),$

\ \ \ \ \ \ \ \ \ \ \ \ \ \ \ \ \ \ \ \ $=0$ $($in case of $k\neq i,j),(0\leq i<j\leq 7),$

\ \ \ \ $ [Ud_{ij},Ux_{2k}]=-\sum\limits_{0\leq n<l\leq 7}Mv^{2}(ki,kj)Ux_{2l}$ $($where $k=n)$

\ \ \ \ \ \ \ \ \ \ \ \ \ \ \ \ \ \ \ \ $\ \ \
+\sum\limits_{0\leq n<l\leq 7}Mv^{2}(ki,kj)Ux_{2n}$ $($where $k=l)$,

$\ \ \ \ \ \ \ \ \ \ \ \ \ \ \ \ \ \ \ \ \ \ \ \ \ \ \ \ \ \ \
(ki=Nu(i+1,j+1),kj=Nu(n+1,l+1),0\leq i<j\leq 7),$

\ \ \ \ $[Ud_{ij},Ux_{3k}]=-\sum\limits_{0\leq n<l\leq 7}Mv(ki,kj)Ux_{3l}$ $($where $n=k)$

\ \ \ \ \ \ \ \ \ \ \ \ \ \ \ \ \ \ \ \ \  $\ \ \
+\sum\limits_{0\leq n<l\leq 7}Mv(ki,kj)Ux_{3n}$ $($where $l=k),$

$\ \ \ \ \ \ \ \ \ \ \ \ \ \ \ \ \ \ \ \ \ \ \ \ \ \ \ \ \ \ \
(ki=Nu(i+1,j+1),kj=Nu(n+1,l+1),0\leq i<j\leq 7),$

\ \ \ \ $[Ud_{ij},U\chi _{k}]=0,(0\leq i<j\leq 7,k=1,2,3).$

\ \ \ \ $[Um_{1i},Ux_{1j}]=0,(0\leq i\neq j\leq 7),$

\ \ \ \ $[Um_{1i},Ux_{1i}]=U\chi _{2}-U\chi _{3},(0\leq i\leq 7),$

\ \ \ \ $[Um_{1i},Ux_{2j}]=\frac{1}{2}Sn(i+1,j+1)Ux_{3k},(k=Ca(i+1,j+1),\ 0\leq i,j\leq 7),$

\ \ \ \ $[Um_{1i},Ux_{3j}]=-\frac{1}{2}Sn(j+1,i+1)Ux_{2k},(k=Ca(j+1,i+1),\ 0\leq i,j\leq 7),$

\ \ \ \ $[Um_{1i},U\chi _{1}]=0\ ,(0\leq i\leq 7),$

\ \ \ \ $[Um_{1i},U\chi _{2}]=-\frac{1}{2}Ux_{1i}\ ,(0\leq i\leq 7),$

\ \ \ \ $[Um_{1i},U\chi _{3}]=\frac{1}{2}Ux_{1i}\ ,(0\leq i\leq 7),$

\ \ \ \ $[Um_{2i},Ux_{1j}]=-\frac{1}{2}Sn(j+1,i+1)Ux_{3k},(k=Ca(j+1,i+1),\ 0\leq i,j\leq 7),$

\ \ \ \ $[Um_{2i},Ux_{2i}]=-U\chi _{1}+U\chi _{3},(0\leq i\leq 7),$

\ \ \ \ $[Um_{2i},Ux_{2j}]=0\ ,(0\leq i\neq j\leq 7),$

\ \ \ \ $[Um_{2i},Ux_{3j}]=\frac{1}{2}Sn(i+1,j+1)Ux_{1k},(k=Ca(i+1,j+1),\ 0\leq i,j\leq 7),$

\ \ \ \ $[Um_{2i},U\chi _{1}]=\frac{1}{2}Ux_{2i}\ ,(0\leq i\leq 7)$,

\ \ \ \ $[Um_{2i},U\chi _{2}]=0\ ,(0\leq i\leq 7)$,

\ \ \ \ $[Um_{2i},U\chi _{3}]=-\frac{1}{2}Ux_{2i}\ ,(0\leq i\leq 7),$

\ \ \ \ $[Um_{3i},Ux_{1j}]=\frac{1}{2}Sn(i+1,j+1)Ux_{2k},(k=Ca(i+1,j+1),\ 0\leq i,j\leq 7),$

\ \ \ \ $[Um_{3i},Ux_{2j}]=-\frac{1}{2}Sn(j+1,i+1)Ux_{1k},(k=Ca(j+1,i+1),\ 0\leq i,j\leq 7),$

\ \ \ \ $[Um_{3i},Ux_{3i}]=U\chi _{1}-U\chi _{2},(0\leq i\leq 7),$

\ \ \ \ $[Um_{3i},Ux_{3j}]=0\ ,(0\leq i\neq j\leq 7),$

\ \ \ \ $[Um_{3i},U\chi _{1}]=-\frac{1}{2}Ux_{3i}\ ,(0\leq i\leq 7),$

\ \ \ \ $[Um_{3i},U\chi _{2}]=\frac{1}{2}Ux_{3i}\ ,(0\leq i\leq 7),$

\ \ \ \ $[Um_{3i},U\chi _{3}]=0\ ,(0\leq i\leq 7),$

\ \ \ \ $[Ut_{1i},Ux_{1j}]=0\ ,(0\leq i\neq j\leq 7),$

\ \ \ \ $[Ut_{1i},Ux_{1i}]=U\chi _{2}+U\chi _{3},(0\leq i\leq 7),$

\ \ \ \ $[Ut_{1i},Ux_{2j}]=\frac{1}{2}Sn(i+1,j+1)Ux_{3k},(k=Ca(i+1,j+1),\ 0\leq i,j\leq 7),$

\ \ \ \ $[Ut_{1i},Ux_{3j}]=\frac{1}{2}Sn(j+1,i+1)Ux_{2k},(k=Ca(j+1,i+1),\ 0\leq i,j\leq 7),$

\ \ \ \ $[Ut_{1i},U\chi _{1}]=0\ ,(0\leq i\leq 7),$

\ \ \ \ $[Ut_{1i},U\chi _{2}]=\frac{1}{2}Ux_{1i}\ ,(0\leq i\leq 7),$

\ \ \ \ $[Ut_{1i},U\chi _{3}]=\frac{1}{2}Ux_{1i}\ ,(0\leq i\leq 7),$

\ \ \ \ $[Ut_{2i},Ux_{1j}]=\frac{1}{2}Sn(i+1,j+1)Ux_{3k},(k=Ca(i+1,j+1),\ 0\leq i,j\leq 7),$

\ \ \ \ $[Ut_{2i},Ux_{2i}]=U\chi _{1}+U\chi _{3},(0\leq i\leq 7),$

\ \ \ \ $[Ut_{2i},Ux_{2j}]=0\ ,(0\leq i\neq j\leq 7),$

\ \ \ \ $[Ut_{2i},Ux_{3j}]=\frac{1}{2}Sn(i+1,j+1)Ux_{1k},(k=Ca(i+1,j+1),\ 0\leq i,j\leq 7),$

\ \ \ \ $[Ut_{2i},U\chi _{1}]=\frac{1}{2}Ux_{2i}\ ,(0\leq i\leq 7)$,

\ \ \ \ $[Ut_{2i},U\chi _{2}]=0\ ,(0\leq i\leq 7)$,

\ \ \ \ $[Ut_{2i},U\chi _{3}]=\frac{1}{2}Ux_{2i}\ ,(0\leq i\leq 7),$

\ \ \ \ $[Ut_{3i},Ux_{1j}]=\frac{1}{2}Sn(i+1,j+1)Ux_{2k},(k=Ca(i+1,j+1),\ 0\leq i,j\leq 7),$

\ \ \ \ $[Ut_{3i},Ux_{2j}]=\frac{1}{2}Sn(j+1,i+1)Ux_{1k},(k=Ca(j+1,i+1),\ 0\leq i,j\leq 7),$

\ \ \ \ $[Ut_{3i},Ux_{3i}]=U\chi _{1}+U\chi _{2},(0\leq i\leq 7),$

\ \ \ \ $[Ut_{3i},Ux_{3j}]=0\ ,(0\leq i\neq j\leq 7),$

\ \ \ \ $[Ut_{3i},U\chi _{1}]=\frac{1}{2}Ux_{3i}\ ,(0\leq i\leq 7),$

\ \ \ \ $[Ut_{3i},U\chi _{2}]=\frac{1}{2}Ux_{3i}\ ,(0\leq i\leq 7),$

\ \ \ \ $[Ut_{3i},U\chi _{3}]=0\ ,(0\leq i\leq 7),$

\ \ \ \ $[U\tau _{1},Ux_{1j}]=-\frac{1}{2}Ux_{1j}\ ,(0\leq j\leq 7),$

\ \ \ \ $[U\tau _{1},Ux_{2j}]=0$ $,(0\leq j\leq 7),$

\ \ \ \ $[U\tau _{1},Ux_{3j}]=\frac{1}{2}Ux_{3j}\ ,(0\leq j\leq 7),$

\ \ \ \ $[U\tau _{1},U\chi _{1}]=U\chi _{1}\ ,$

\ \ \ \ $[U\tau _{1},U\chi _{2}]=0\ ,$

\ \ \ \ $[U\tau _{1},U\chi _{3}]=-U\chi _{3}\ ,$

\ \ \ \ $[U\tau _{2},Ux_{1j}]=0,(0\leq j\leq 7),$

\ \ \ \ $[U\tau _{2},Ux_{2j}]=-\frac{1}{2}Ux_{2j}\ ,(0\leq j\leq 7),$

\ \ \ \ $[U\tau _{2},Ux_{3j}]=\frac{1}{2}Ux_{3j}\ ,(0\leq j\leq 7),$

\ \ \ \ $[U\tau _{2},U\chi _{1}]=0\ ,$

\ \ \ \ $[U\tau _{2},U\chi _{2}]=U\chi _{2}\ ,$

\ \ \ \ $[U\tau _{2},U\chi _{3}]=-U\chi _{3}\ .$

\bigskip

\emph{Proof. \ }We have the above Lie bracket operations with
calculations using Maxima.\ \ \ \ \emph{Q.E.D.}

\bigskip

\emph{Lemma 17.3. }\ We have the following Lie bracket operations.

\ \ \ \ $[Ua_{1i},Ux_{1i}]=-2U\gamma _{1},(0\leq i\leq 7),$

\ \ \ \ $[Ua_{1i},Ux_{1j}]=0,(\ 0\leq i\neq j\leq 7),$

\ \ \ \ $[Ua_{1i},Ux_{2j}]=Sn(i+1,j+1)Uy_{3k},(k=Ca(i+1,j+1),\ 0\leq i,j\leq 7),$

\ \ \ \ $[Ua_{1i},Ux_{3j}]=Sn(j+1,i+1)Uy_{2k},(k=Ca(j+1,i+1),\ 0\leq i,j\leq 7),$

\ \ \ \ $[Ua_{2i},Ux_{1j}]=Sn(j+1,i+1)Uy_{3k},(k=Ca(j+1,i+1),\ 0\leq i,j\leq 7),$

\ \ \ \ $[Ua_{2i},Ux_{2i}]=-2U\gamma _{2},(0\leq i\leq 7),$

\ \ \ \ $[Ua_{2i},Ux_{2j}]=0,(\ 0\leq i\neq j\leq 7),$

\ \ \ \ $[Ua_{2i},Ux_{3j}]=Sn(i+1,j+1)Uy_{1k},(k=Ca(i+1,j+1),\ 0\leq i,j\leq 7),$

\ \ \ \ $[Ua_{3i},Ux_{1j}]=Sn(i+1,j+1)Uy_{2k},(k=Ca(i+1,j+1),\ 0\leq i,j\leq 7),$

\ \ \ \ $[Ua_{3i},Ux_{2j}]=Sn(j+1,i+1)Uy_{1k},(k=Ca(j+1,i+1),\ 0\leq i,j\leq 7),$

\ \ \ \ $[Ua_{3i},Ux_{3i}]=-2U\gamma _{3},(0\leq i\leq 7),$

\ \ \ \ $[Ua_{3i},Ux_{3j}]=0,(\ 0\leq i\neq j\leq 7),$

\ \ \ \ $[Ua_{1i},U\chi _{1}]=-Uy_{1i},(\ 0\leq i\leq 7),$

\ \ \ \ $[Ua_{1i},U\chi _{2}]=0,$

\ \ \ \ $[Ua_{1i},U\chi _{3}]=0,$

\ \ \ \ $[Ua_{2i},U\chi _{1}]=0,$

\ \ \ \ $[Ua_{2i},U\chi _{2}]=-Uy_{2i},(\ 0\leq i\leq 7),$

\ \ \ \ $[Ua_{2i},U\chi _{3}]=0,$

\ \ \ \ $[Ua_{3i},U\chi _{1}]=0,$

\ \ \ \ $[Ua_{3i},U\chi _{2}]=0,$

\ \ \ \ $[Ua_{3i},U\chi _{3}]=-Uy_{3i},(\ 0\leq i\leq 7),$

\ \ \ \ $[U\alpha _{1},Ux_{1i}]=-Uy_{1i},(\ 0\leq i\leq 7),$

\ \ \ \ $[U\alpha _{1},Ux_{2i}]=0,(\ 0\leq i\leq 7),$

\ \ \ \ $[U\alpha _{1},Ux_{3i}]=0,(\ 0\leq i\leq 7),$

\ \ \ \ $[U\alpha _{2},Ux_{1i}]=0,(\ 0\leq i\leq 7),$

\ \ \ \ $[U\alpha _{2},Ux_{2i}]=-Uy_{2i},(\ 0\leq i\leq 7),$

\ \ \ \ $[U\alpha _{2},Ux_{3i}]=0,(\ 0\leq i\leq 7),$

\ \ \ \ $[U\alpha _{3},Ux_{1i}]=0,(\ 0\leq i\leq 7),$

\ \ \ \ $[U\alpha _{3},Ux_{2i}]=0,(\ 0\leq i\leq 7),$

\ \ \ \ $[U\alpha _{3},Ux_{3i}]=-Uy_{3i},(\ 0\leq i\leq 7),$

\ \ \ \ $[U\alpha _{1},U\chi _{1}]=0,$

\ \ \ \ $[U\alpha _{1},U\chi _{2}]=U\gamma _{3},$

\ \ \ \ $[U\alpha _{1},U\chi _{3}]=U\gamma _{2},$

\ \ \ \ $[U\alpha _{2},U\chi _{1}]=U\gamma _{3},$

\ \ \ \ $[U\alpha _{2},U\chi _{2}]=0,$

\ \ \ \ $[U\alpha _{2},U\chi _{3}]=U\gamma _{1},$

\ \ \ \ $[U\alpha _{3},U\chi _{1}]=U\gamma _{2},$

\ \ \ \ $[U\alpha _{3},U\chi _{2}]=U\gamma _{1},$

\ \ \ \ $[U\alpha _{3},U\chi _{3}]=0.$

\bigskip

\emph{Proof. \ }We have the above Lie bracket operations with
calculations using Maxima.\ \ \ \ \emph{Q.E.D.}

\bigskip

\emph{Lemma 17.4. }\ We have the following Lie bracket operations.

\ \ \ \ $[Ub_{1i},Ux_{1i}]=2U\eta _{1},(0\leq i\leq 7),$

\ \ \ \ $[Ub_{1i},Ux_{1j}]=0$,$(\ 0\leq i\neq j\leq 7),$

\ \ \ \ $[Ub_{1i},Ux_{2j}]=0$,$(\ 0\leq i,j\leq 7),$

\ \ \ \ $[Ub_{1i},Ux_{3j}]=0$,$(\ 0\leq i,j\leq 7),$

\ \ \ \ $[Ub_{2i},Ux_{1j}]=0$,$(\ 0\leq i,j\leq 7),$

\ \ \ \ $[Ub_{2i},Ux_{2i}]=2U\eta _{1},(0\leq i\leq 7),$

\ \ \ \ $[Ub_{2i},Ux_{2j}]=0$,$(\ 0\leq i\neq j\leq 7),$

\ \ \ \ $[Ub_{2i},Ux_{3j}]=0$,$(\ 0\leq i,j\leq 7),$

\ \ \ \ $[Ub_{3i},Ux_{1j}]=0$,$(\ 0\leq i,j\leq 7),$

\ \ \ \ $[Ub_{3i},Ux_{2j}]=0$,$(\ 0\leq i,j\leq 7),$

\ \ \ \ $[Ub_{3i},Ux_{3i}]=2U\eta _{1},(0\leq i\leq 7),$

\ \ \ \ $[Ub_{3i},Ux_{3j}]=0$,$(\ 0\leq i\neq j\leq 7),$

\ \ \ \ $[Ub_{ij},U\chi _{k}]=0$,$(\ 1\leq i,k\leq 3,0\leq j\leq 7),$

\ \ \ \ $[U\beta _{i},Ux_{kj}]=0$,$(\ 1\leq i,k\leq 3,0\leq j\leq 7), $

\ \ \ \ $[U\beta _{i},U\chi _{i}]=U\eta _{1},(1\leq i\leq 3),$

\ \ \ \ $[U\beta _{i},U\chi _{j}]=0,(1\leq i\neq j\leq 3).$

\bigskip

\emph{Proof. \ }We have the above Lie bracket operations with
calculations using Maxima.\ \ \ \ \emph{Q.E.D.}

\bigskip

\emph{Lemma 17.5. }\ We have the following Lie bracket operations.

\ \ \ \ $[U\rho _{1},Ux_{ij}]=-\frac{1}{3}Ux_{ij},(1\leq i\leq 3,0\leq j\leq 7),$

\ \ \ \ $[U\rho _{1},U\chi _{i}]=-\frac{1}{3}U\chi _{i},(1\leq i\leq 3). $

\bigskip

\emph{Proof. \ }We have the above Lie bracket operations with
calculations using Maxima.\ \ \ \ \emph{Q.E.D.}

\bigskip

\emph{Lemma 17.6. }\ We have the following Lie bracket operations.

\ \ \ \ $[Ud_{ij},Uy_{1k}]=-Uy_{1j}$ $($in case of $k=i),$

\ \ \ \ $\ \ \ \ \ \ \ \ \ \ \ \ \ \ \ \ =Uy_{1i}$ $($in case of $k=j),$

\ \ \ \ $\ \ \ \ \ \ \ \ \ \ \ \ \ \ \ \ =0$ $($in case of $k\neq i,j),(0\leq i<j\leq 7),$

\ \ \ \ $[Ud_{ij},Uy_{2k}]=-\sum\limits_{0\leq n<l\leq 7}Mv^{2}(ki,kj)Uy_{2l}$ $($where $k=n)$

\ \ \ \ \ \ \ \ \ \ \ \ \ \ \ \ \ \ \ \ $\ \ \
+\sum\limits_{0\leq n<l\leq 7}Mv^{2}(ki,kj)Uy_{2n}$ $($where $k=l)$,

\ \ \ \ \ \ \ \ $\ \ \ \ \ \ \ \ \
(ki=Nu(i+1,j+1),kj=Nu(n+1,l+1),0\leq i<j\leq 7),$

\ \ \ \ $[Ud_{ij},Uy_{3k}]=-\sum\limits_{0\leq n<l\leq 7}Mv(ki,kj)Uy_{3l}$ $($where $n=k)$

\ \ \ \ \ \ \ \ \ \ \ \ \ \ \ \ \ \ \ \ \ $\ \ \
+\sum\limits_{0\leq n<l\leq 7}Mv(ki,kj)Uy_{3n}$ $($where $l=k),$

\ \ \ \ $\ \ \ \ \ \ \ \ \ \ \ \ \ \ \ \ \ \ \ \
(ki=Nu(i+1,j+1),kj=Nu(n+1,l+1),0\leq i<j\leq 7),$

\ \ \ \ $[Ud_{ij},U\gamma _{k}]=0,(0\leq i<j\leq 7,k=1,2,3)$,

\ \ \ \ $[Um_{1i},Uy_{1i}]=U\gamma _{2}-U\gamma _{3},(0\leq i\leq 7),$

\ \ \ \ $[Um_{1i},Uy_{2j}]=\frac{1}{2}Sn(i+1,j+1)Uy_{3k},(k=Ca(i+1,j+1),\ 0\leq i,j\leq 7),$

\ \ \ \ $[Um_{1i},Uy_{3j}]=-\frac{1}{2}Sn(j+1,i+1)Uy_{2k},(k=Ca(j+1,i+1),\ 0\leq i,j\leq 7),$

\ \ \ \ $[Um_{1i},U\gamma _{1}]=0\ ,(0\leq i\leq 7),$

\ \ \ \ $[Um_{1i},U\gamma _{2}]=-\frac{1}{2}Uy_{1i}\ ,(0\leq i\leq 7),$

\ \ \ \ $[Um_{1i},U\gamma _{3}]=\frac{1}{2}Uy_{1i}\ ,(0\leq i\leq 7),$

\ \ \ \ $[Um_{2i},Uy_{1j}]=-\frac{1}{2}Sn(j+1,i+1)Uy_{3k},(k=Ca(j+1,i+1),\ 0\leq i,j\leq 7),$

\ \ \ \ $[Um_{2i},Uy_{2i}]=-U\gamma _{1}+U\gamma _{3},(0\leq i\leq 7),$

\ \ \ \ $[Um_{2i},Uy_{2j}]=0\ ,(0\leq i\neq j\leq 7),$

\ \ \ \ $[Um_{2i},Uy_{3j}]=\frac{1}{2}Sn(i+1,j+1)Uy_{1k},(k=Ca(i+1,j+1),\ 0\leq i,j\leq 7),$

\ \ \ \ $[Um_{2i},U\gamma _{1}]=\frac{1}{2}Uy_{2i}\ ,(0\leq i\leq 7)$,

\ \ \ \ $[Um_{2i},U\gamma _{2}]=0\ ,(0\leq i\leq 7)$,

\ \ \ \ $[Um_{2i},U\gamma _{3}]=-\frac{1}{2}Uy_{2i}\ ,(0\leq i\leq 7),$

\ \ \ \ $[Um_{3i},Uy_{1j}]=\frac{1}{2}Sn(i+1,j+1)Uy_{2k},(k=Ca(i+1,j+1),\ 0\leq i,j\leq 7),$

\ \ \ \ $[Um_{3i},Uy_{2j}]=-\frac{1}{2}Sn(j+1,i+1)Uy_{1k},(k=Ca(j+1,i+1),\ 0\leq i,j\leq 7),$

\ \ \ \ $[Um_{3i},Uy_{3i}]=U\gamma _{1}-U\gamma _{2},(0\leq i\leq 7),$

\ \ \ \ $[Um_{3i},Uy_{3j}]=0\ ,(0\leq i\neq j\leq 7),$

\ \ \ \ $[Um_{3i},U\gamma _{1}]=-\frac{1}{2}Uy_{3i}\ ,(0\leq i\leq 7),$

\ \ \ \ $[Um_{3i},U\gamma _{2}]=\frac{1}{2}Uy_{3i}\ ,(0\leq i\leq 7),$

\ \ \ \ $[Um_{3i},U\gamma _{3}]=0\ ,(0\leq i\leq 7),$

\ \ \ \ $[Ut_{1i},Uy_{1j}]=0\ ,(0\leq i\neq j\leq 7),$

\ \ \ \ $[Ut_{1i},Uy_{1i}]=-U\gamma _{2}-U\gamma _{3},(0\leq i\leq 7),$

\ \ \ \ $[Ut_{1i},Uy_{2j}]=-\frac{1}{2}Sn(i+1,j+1)Uy_{3k},(k=Ca(i+1,j+1),\ 0\leq i,j\leq 7),$

\ \ \ \ $[Ut_{1i},Uy_{3j}]=-\frac{1}{2}Sn(j+1,i+1)Uy_{2k},(k=Ca(j+1,i+1),\ 0\leq i,j\leq 7),$

\ \ \ \ $[Ut_{1i},U\gamma _{1}]=0\ ,(0\leq i\leq 7),$

\ \ \ \ $[Ut_{1i},U\gamma _{2}]=-\frac{1}{2}Uy_{1i}\ ,(0\leq i\leq 7),$

\ \ \ \ $[Ut_{1i},U\gamma _{3}]=-\frac{1}{2}Uy_{1i}\ ,(0\leq i\leq 7),$

\ \ \ \ $[Ut_{2i},Uy_{1j}]=-\frac{1}{2}Sn(i+1,j+1)Uy_{3k},(k=Ca(i+1,j+1),\ 0\leq i,j\leq 7),$

\ \ \ \ $[Ut_{2i},Uy_{2i}]=-U\gamma _{1}-U\gamma _{3},(0\leq i\leq 7),$

\ \ \ \ $[Ut_{2i},Uy_{2j}]=0\ ,(0\leq i\neq j\leq 7),$

\ \ \ \ $[Ut_{2i},Uy_{3j}]=-\frac{1}{2}Sn(i+1,j+1)Uy_{1k},(k=Ca(i+1,j+1),\ 0\leq i,j\leq 7),$

\ \ \ \ $[Ut_{2i},U\gamma _{1}]=-\frac{1}{2}Uy_{2i}\ ,(0\leq i\leq 7),$

\ \ \ \ $[Ut_{2i},U\gamma _{2}]=0\ ,(0\leq i\leq 7)$,

\ \ \ \ $[Ut_{2i},U\gamma _{3}]=-\frac{1}{2}Uy_{2i}\ ,(0\leq i\leq 7),$

\ \ \ \ $[Ut_{3i},Uy_{1j}]=-\frac{1}{2}Sn(i+1,j+1)Uy_{2k},(k=Ca(i+1,j+1),\ 0\leq i,j\leq 7),$

\ \ \ \ $[Ut_{3i},Uy_{2j}]=-\frac{1}{2}Sn(j+1,i+1)Uy_{1k},(k=Ca(j+1,i+1),\ 0\leq i,j\leq 7),$

\ \ \ \ $[Ut_{3i},Uy_{3i}]=-U\gamma _{1}-U\gamma _{2},(0\leq i\leq 7),$

\ \ \ \ $[Ut_{3i},Uy_{3j}]=0\ ,(0\leq i\neq j\leq 7),$

\ \ \ \ $[Ut_{3i},U\gamma _{1}]=-\frac{1}{2}Uy_{3i}\ ,(0\leq i\leq 7),$

\ \ \ \ $[Ut_{3i},U\gamma _{2}]=-\frac{1}{2}Uy_{3i}\ ,(0\leq i\leq 7),$

\ \ \ \ $[Ut_{3i},U\gamma _{3}]=0\ ,(0\leq i\leq 7),$

\ \ \ \ $[U\tau _{1},Uy_{1j}]=\frac{1}{2}Uy_{1j}\ ,(0\leq j\leq 7),$

\ \ \ \ $[U\tau _{1},Uy_{2j}]=0$ $,(0\leq j\leq 7),$

\ \ \ \ $[U\tau _{1},Uy_{3j}]=-\frac{1}{2}Uy_{3j}\ ,(0\leq j\leq 7),$

\ \ \ \ $[U\tau _{1},U\gamma _{1}]=-U\gamma _{1}\ ,$

\ \ \ \ $[U\tau _{1},U\gamma _{2}]=0\ ,$

\ \ \ \ $[U\tau _{1},U\gamma _{3}]=U\gamma _{3}\ ,$

\ \ \ \ $[U\tau _{2},Uy_{1j}]=0,(0\leq j\leq 7),$

\ \ \ \ $[U\tau _{2},Uy_{2j}]=\frac{1}{2}Uy_{2j}\ ,(0\leq j\leq 7),$

\ \ \ \ $[U\tau _{2},Uy_{3j}]=-\frac{1}{2}Uy_{3j}\ ,(0\leq j\leq 7),$

\ \ \ \ $[U\tau _{2},U\gamma _{1}]=0\ ,$

\ \ \ \ $[U\tau _{2},U\gamma _{2}]=-U\gamma _{2}\ ,$

\ \ \ \ $[U\tau _{2},U\gamma _{3}]=U\gamma _{3}\ .$

\bigskip

\emph{Proof. \ }We have the above Lie bracket operations with
calculations using Maxima.\ \ \ \ \emph{Q.E.D.}

\bigskip

\emph{Lemma 17.7. }\ We have the following Lie bracket operations.

\ \ \ \ $[Ub_{1i},Uy_{1i}]=-2U\chi _{1},(0\leq i\leq 7),$

\ \ \ \ $[Ub_{1i},Uy_{1j}]=0,(\ 0\leq i\neq j\leq 7),$

\ \ \ \ $[Ub_{1i},Uy_{2j}]=-Sn(i+1,j+1)Ux_{3k},(k=Ca(i+1,j+1),\ 0\leq i,j\leq 7),$

\ \ \ \ $[Ub_{1i},Uy_{3j}]=0,(\ 0\leq i,j\leq 7),$

\ \ \ \ $[Ub_{2i},Uy_{1j}]=0,(\ 0\leq i,j\leq 7),$

\ \ \ \ $[Ub_{2i},Uy_{2i}]=-2U\chi _{2},(0\leq i\leq 7),$

\ \ \ \ $[Ub_{2i},Uy_{2j}]=0,(\ 0\leq i\neq j\leq 7),$

\ \ \ \ $[Ub_{2i},Uy_{3j}]=-Sn(i+1,j+1)Ux_{1k},(k=Ca(i+1,j+1),\ 0\leq i,j\leq 7),$

\ \ \ \ $[Ub_{3i},Uy_{1j}]=-Sn(j+1,i+1)Ux_{2k},(k=Ca(j+1,i+1),\ 0\leq i,j\leq 7),$

\ \ \ \ $[Ub_{3i},Uy_{2j}]=0,(\ 0\leq i\neq j\leq 7),$

\ \ \ \ $[Ub_{3i},Uy_{3i}]=-2U\chi _{3},(0\leq i\leq 7),$

\ \ \ \ $[Ub_{3i},Uy_{3j}]=0,(\ 0\leq i\neq j\leq 7),$

\ \ \ \ $[Ub_{1i},U\gamma _{1}]=-Ux_{1i},(\ 0\leq i\leq 7),$

\ \ \ \ $[Ub_{1i},U\gamma _{2}]=0,$

\ \ \ \ $[Ub_{1i},U\gamma _{3}]=0,$

\ \ \ \ $[Ub_{2i},U\gamma _{1}]=0,$

\ \ \ \ $[Ub_{2i},U\gamma _{2}]=-Ux_{2i},(\ 0\leq i\leq 7),$

\ \ \ \ $[Ub_{2i},U\gamma _{3}]=0,$

\ \ \ \ $[Ub_{3i},U\gamma _{1}]=0,$

\ \ \ \ $[Ub_{3i},U\gamma _{2}]=0,$

\ \ \ \ $[Ub_{3i},U\gamma _{3}]=-Ux_{3i},(\ 0\leq i\leq 7),$

\ \ \ \ $[U\beta _{1},Uy_{1i}]=-Ux_{1i},(\ 0\leq i\leq 7),$

\ \ \ \ $[U\beta _{1},Uy_{2i}]=0,(\ 0\leq i\leq 7),$

\ \ \ \ $[U\beta _{1},Uy_{3i}]=0,(\ 0\leq i\leq 7),$

\ \ \ \ $[U\beta _{2},Uy_{1i}]=0,(\ 0\leq i\leq 7),$

\ \ \ \ $[U\beta _{2},Uy_{2i}]=-Ux_{2i},(\ 0\leq i\leq 7),$

\ \ \ \ $[U\beta _{2},Uy_{3i}]=0,(\ 0\leq i\leq 7),$

\ \ \ \ $[U\beta _{3},Uy_{1i}]=0,(\ 0\leq i\leq 7),$

\ \ \ \ $[U\beta _{3},Uy_{2i}]=0,(\ 0\leq i\leq 7),$

\ \ \ \ $[U\beta _{3},Uy_{3i}]=-Ux_{3i},(\ 0\leq i\leq 7),$

\ \ \ \ $[U\beta _{1},U\gamma _{1}]=0,$

\ \ \ \ $[U\beta _{1},U\gamma _{2}]=U\chi _{3},$

\ \ \ \ $[U\beta _{1},U\gamma _{3}]=U\chi _{2},$

\ \ \ \ $[U\beta _{2},U\gamma _{1}]=U\chi _{3},$

\ \ \ \ $[U\beta _{2},U\gamma _{2}]=0,$

\ \ \ \ $[U\beta _{2},U\gamma _{3}]=U\chi _{1},$

\ \ \ \ $[U\beta _{3},U\gamma _{1}]=U\chi _{2},$

\ \ \ \ $[U\beta _{3},U\gamma _{2}]=U\chi _{1},$

\ \ \ \ $[U\beta _{3},U\gamma _{3}]=0.$

\bigskip

\emph{Proof. \ }We have the above Lie bracket operations with
calculations using Maxima.\ \ \ \ \emph{Q.E.D.}

\bigskip

\emph{Lemma 17.8. }\ We have the following Lie bracket operations.

\ \ \ \ $[Ua_{1i},Uy_{1i}]=2U\xi _{1},(0\leq i\leq 7),$

\ \ \ \ $[Ua_{1i},Uy_{1j}]=0$,$(\ 0\leq i\neq j\leq 7),$

\ \ \ \ $[Ua_{1i},Uy_{2j}]=0$,$(\ 0\leq i,j\leq 7),$

\ \ \ \ $[Ua_{1i},Uy_{3j}]=0$,$(\ 0\leq i,j\leq 7),$

\ \ \ \ $[Ua_{2i},Uy_{1j}]=0$,$(\ 0\leq i,j\leq 7),$

\ \ \ \ $[Ua_{2i},Uy_{2i}]=2U\xi _{2},(0\leq i\leq 7),$

\ \ \ \ $[Ua_{2i},Uy_{2j}]=0$,$(\ 0\leq i\neq j\leq 7),$

\ \ \ \ $[Ua_{2i},Uy_{3j}]=0$,$(\ 0\leq i,j\leq 7),$

\ \ \ \ $[Ua_{3i},Uy_{1j}]=0$,$(\ 0\leq i,j\leq 7),$

\ \ \ \ $[Ua_{3i},Uy_{2j}]=0$,$(\ 0\leq i,j\leq 7),$

\ \ \ \ $[Ua_{3i},Uy_{3i}]=2U\xi _{3},(0\leq i\leq 7),$

\ \ \ \ $[Ua_{3i},Uy_{3j}]=0$,$(\ 0\leq i\neq j\leq 7),$

\ \ \ \ $[Ua_{ij},U\gamma _{k}]=0$,$(\ 1\leq i,k\leq 3,0\leq j\leq 7),$

\ \ \ \ $[U\alpha _{i},Uy_{kj}]=0$,$(\ 1\leq i,k\leq 3,0\leq j\leq 7),$

\ \ \ \ $[U\alpha _{i},U\gamma _{i}]=U\xi _{i},(1\leq i\leq 3),$

\ \ \ \ $[U\alpha _{i},u\gamma j]=0,(1\leq i\neq j\leq 3).$

\bigskip

\emph{Proof. \ }We have the above Lie bracket operations with
calculations using Maxima.\ \ \ \ \emph{Q.E.D.}

\bigskip

\emph{Lemma 17.9. }\ We have the following Lie bracket operations.

\ \ \ \ $[U\rho _{1},Uy_{ij}]=\frac{1}{3}Uy_{ij},(1\leq i\leq 3,0\leq j\leq 7),$

\ \ \ \ $[U\rho _{1},U\gamma _{i}]=\frac{1}{3}U\gamma _{i},(1\leq i\leq 3).$

\bigskip

\emph{Proof. \ }We have the above Lie bracket operations with
calculations using Maxima.\ \ \ \ \emph{Q.E.D.}

\bigskip

\emph{Lemma 17.10. }\ We have the following Lie bracket operations.

\ \ \ \ $[Ud_{ij},U\xi _{1}]=0,(0\leq i<j\leq 7),$

\ \ \ \ $[Um_{ij},U\xi _{1}]=0,(1\leq i\leq 3,0\leq j\leq 7),$

\ \ \ \ $[Ut_{ij},U\xi _{1}]=0,(1\leq i\leq 3,0\leq j\leq 7),$

\ \ \ \ $[U\tau _{i},U\xi _{1}]=0,(1\leq i\leq 2),$

\ \ \ \ $[Ua_{ij},U\xi _{1}]=0,(1\leq i\leq 3,0\leq j\leq 7),$

\ \ \ \ $[U\alpha _{i},U\xi _{1}]=0,(1\leq i\leq 3),$

\ \ \ \ $[Ub_{ij},U\xi _{1}]=Uy_{ij},(1\leq i\leq 3,0\leq j\leq 7),$

\ \ \ \ $[U\beta _{i},U\xi _{1}]=U\gamma _{i},(1\leq i\leq 3),$

\ \ \ \ $[U\rho _{1},U\xi _{1}]=U\xi _{1}$,

\ \ \ \ $[Ud_{ij},U\eta _{1}]=0,(0\leq i<j\leq 7),$

\ \ \ \ $[Um_{ij},U\eta _{1}]=0,(1\leq i\leq 3,0\leq j\leq 7),$

\ \ \ \ $[Ut_{ij},U\eta _{1}]=0,(1\leq i\leq 3,0\leq j\leq 7),$

\ \ \ \ $[U\tau _{i},U\eta _{1}]=0,(1\leq i\leq 2),$

\ \ \ \ $[Ub_{ij},U\eta _{1}]=0,(1\leq i\leq 3,0\leq j\leq 7),$

\ \ \ \ $[U\beta _{i},U\eta _{1}]=0,(1\leq i\leq 3),$

\ \ \ \ $[Ua_{ij},U\eta _{1}]=Ux_{ij},(1\leq i\leq 3,0\leq j\leq 7),$

\ \ \ \ $[U\alpha _{i},U\eta _{1}]=U\chi _{i},(1\leq i\leq 3),$

\ \ \ \ $[U\rho _{1},U\eta _{1}]=-U\eta _{1}$.

\bigskip

\emph{Proof. \ }We have the above Lie bracket operations with
calculations using Maxima.\ \ \ \ \emph{Q.E.D.}

\bigskip

\emph{Lemma 17.11. }\ We have the following Lie bracket operations.

\ \ \ \ $[Ud_{ij},Uz_{1k}]=-Uz_{1j}$ $($in case of $k=i),$

\ \ \ \ $\ \ \ \ \ \ \ \ \ \ \ \ \ \ \ \ =Uz_{1i}$ $($in case of $k=j),$

\ \ \ \ $\ \ \ \ \ \ \ \ \ \ \ \ \ \ \ \ =0$ $($in case of $k\neq i,j),(0\leq i<j\leq 7),$

\ \ \ \ $[Ud_{ij},Uz_{2k}]=-\sum\limits_{0\leq n<l\leq 7}Mv^{2}(ki,kj)Uz_{2l}$ $($where $k=n)$

\ \ \ \ \ \ \ \ \ \ \ \ \ \ \ \ \ \ \ \ $\ \ \
+\sum\limits_{0\leq n<l\leq 7}Mv^{2}(ki,kj)Uz_{2n}$ $($where $k=l)$,

\ \ \ \ $\ \ \ \ \ \ \ \ \ \ \ \ \ \ \ \ \ \ \ \
(ki=Nu(i+1,j+1),kj=Nu(n+1,l+1),0\leq i<j\leq 7),$

\ \ \ \ $[Ud_{ij},Uz_{3k}]=-\sum\limits_{0\leq n<l\leq 7}Mv(ki,kj)Uz_{3l}$ $($where $n=k)$

\ \ \ \ \ \ \ \ \ \ \ \ \ \ \ \ \ \ \ \ \ $\ \ \
+\sum\limits_{0\leq n<l\leq 7}Mv(ki,kj)Uz_{3n}$ $($where $l=k),$

\ \ \ \ $\ \ \ \ \ \ \ \ \ \ \ \ \ \ \ \ \ \ \ \
(ki=Nu(i+1,j+1),kj=Nu(n+1,l+1),0\leq i<j\leq 7),$

 \ \ \ \ $[Ud_{ij},U\mu_{k}]=0,(0\leq i<j\leq 7,k=1,2,3)$,

\ \ \ \ $[Um_{1i},Uz_{1j}]=0,(0\leq i\neq j\leq 7),$

\ \ \ \ $[Um_{1i},Uz_{1i}]=U\mu_{2}-U\mu_{3},(0\leq i\leq 7),$

\ \ \ \ $[Um_{1i},Uz_{2j}]=\frac{1}{2}Sn(i+1,j+1)Uz_{3k},(k=Ca(i+1,j+1),\ 0\leq i,j\leq 7),$

\ \ \ \ $[Um_{1i},Uz_{3j}]=-\frac{1}{2}Sn(j+1,i+1)Uz_{2k},(k=Ca(j+1,i+1),\ 0\leq i,j\leq 7),$

\ \ \ \ $[Um_{1i},U\mu_{1}]=0\ ,(0\leq i\leq 7),$

\ \ \ \ $[Um_{1i},U\mu_{2}]=-\frac{1}{2}Uz_{1i}\ ,(0\leq i\leq 7),$

\ \ \ \ $[Um_{1i},U\mu_{3}]=\frac{1}{2}Uz_{1i}\ ,(0\leq i\leq 7),$

\ \ \ \ $[Um_{2i},Uz_{1j}]=-\frac{1}{2}Sn(j+1,i+1)Uz_{3k},(k=Ca(j+1,i+1),\ 0\leq i,j\leq 7),$

\ \ \ \ $[Um_{2i},Uz_{2i}]=-U\mu_{1}+U\mu_{3},(0\leq i\leq 7),$

\ \ \ \ $[Um_{2i},Uz_{2j}]=0\ ,(0\leq i\neq j\leq 7),$

\ \ \ \ $[Um_{2i},Uz_{3j}]=\frac{1}{2}Sn(i+1,j+1)Uz_{1k},(k=Ca(i+1,j+1),\ 0\leq i,j\leq 7),$

\ \ \ \ $[Um_{2i},U\mu_{1}]=\frac{1}{2}Uz_{2i}\ ,(0\leq i\leq 7),$

\ \ \ \ $[Um_{2i},U\mu_{2}]=0\ ,(0\leq i\leq 7)$,

\ \ \ \ $[Um_{2i},U\mu_{3}]=-\frac{1}{2}Uz_{2i}\ ,(0\leq i\leq 7),$

\ \ \ \ $[Um_{3i},Uz_{1j}]=\frac{1}{2}Sn(i+1,j+1)Uz_{2k},(k=Ca(i+1,j+1),\ 0\leq i,j\leq 7),$

\ \ \ \ $[Um_{3i},Uz_{2j}]=-\frac{1}{2}Sn(j+1,i+1)Uz_{1k},(k=Ca(j+1,i+1),\ 0\leq i,j\leq 7),$

\ \ \ \ $[Um_{3i},Uz_{3i}]=U\mu_{1}-U\mu_{2},(0\leq i\leq 7),$

\ \ \ \ $[Um_{3i},Uz_{3j}]=0\ ,(0\leq i\neq j\leq 7),$

\ \ \ \ $[Um_{3i},U\mu_{1}]=-\frac{1}{2}Uz_{3i}\ ,(0\leq i\leq 7),$

\ \ \ \ $[Um_{3i},U\mu_{2}]=\frac{1}{2}Uz_{3i}\ ,(0\leq i\leq 7),$

\ \ \ \ $[Um_{3i},U\mu_{3}]=0\ ,(0\leq i\leq 7),$

\ \ \ \ $[Ut_{1i},Uz_{1j}]=0\ ,(0\leq i\neq j\leq 7),$

\ \ \ \ $[Ut_{1i},Uz_{1i}]=U\mu_{2}+U\mu_{3},(0\leq i\leq 7),$

\ \ \ \ $[Ut_{1i},Uz_{2j}]=\frac{1}{2}Sn(i+1,j+1)Uz_{3k},(k=Ca(i+1,j+1),\ 0\leq i,j\leq 7),$

\ \ \ \ $[Ut_{1i},Uz_{3j}]=\frac{1}{2}Sn(j+1,i+1)Uz_{2k},(k=Ca(j+1,i+1),\ 0\leq i,j\leq 7),$

\ \ \ \ $[Ut_{1i},U\mu_{1}]=0\ ,(0\leq i\leq 7),$

\ \ \ \ $[Ut_{1i},U\mu_{2}]=\frac{1}{2}Uz_{1i}\ ,(0\leq i\leq 7),$

\ \ \ \ $[Ut_{1i},U\mu_{3}]=\frac{1}{2}Uz_{1i}\ ,(0\leq i\leq 7),$

\ \ \ \ $[Ut_{2i},Uz_{1j}]=\frac{1}{2}Sn(i+1,j+1)Uz_{3k},(k=Ca(i+1,j+1),\ 0\leq i,j\leq 7),$

\ \ \ \ $[Ut_{2i},Uz_{2i}]=U\mu_{1}+U\mu_{3},(0\leq i\leq 7),$

\ \ \ \ $[Ut_{2i},Uz_{2j}]=0\ ,(0\leq i\neq j\leq 7),$

\ \ \ \ $[Ut_{2i},Uz_{3j}]=\frac{1}{2}Sn(i+1,j+1)Uz_{1k},(k=Ca(i+1,j+1),\ 0\leq i,j\leq 7),$

\ \ \ \ $[Ut_{2i},U\mu_{1}]=\frac{1}{2}Uz_{2i}\ ,(0\leq i\leq 7),$

\ \ \ \ $[Ut_{2i},U\mu_{2}]=0\ ,(0\leq i\leq 7)$,

\ \ \ \ $[Ut_{2i},U\mu_{3}]=\frac{1}{2}Uz_{2i}\ ,(0\leq i\leq 7),$

\ \ \ \ $[Ut_{3i},Uz_{1j}]=\frac{1}{2}Sn(i+1,j+1)Uz_{2k},(k=Ca(i+1,j+1),\ 0\leq i,j\leq 7),$

\ \ \ \ $[Ut_{3i},Uz_{2j}]=\frac{1}{2}Sn(j+1,i+1)Uz_{1k},(k=Ca(j+1,i+1),\ 0\leq i,j\leq 7),$

\ \ \ \ $[Ut_{3i},Uz_{3i}]=U\mu_{1}+U\mu_{2},(0\leq i\leq 7),$

\ \ \ \ $[Ut_{3i},Uz_{3j}]=0\ ,(0\leq i\neq j\leq 7),$

\ \ \ \ $[Ut_{3i},U\mu_{1}]=\frac{1}{2}Uz_{3i}\ ,(0\leq i\leq 7),$

\ \ \ \ $[Ut_{3i},U\mu_{2}]=\frac{1}{2}Uz_{3i}\ ,(0\leq i\leq 7),$

\ \ \ \ $[Ut_{3i},U\mu_{3}]=0\ ,(0\leq i\leq 7),$

\ \ \ \ $[U\tau _{1},Uz_{1j}]=-\frac{1}{2}Uz_{1j}\ ,(0\leq j\leq 7),$

\ \ \ \ $[U\tau _{1},Uz_{2j}]=0$ $,(0\leq j\leq 7),$

\ \ \ \ $[U\tau _{1},Uz_{3j}]=\frac{1}{2}Uz_{3j}\ ,(0\leq j\leq 7),$

\ \ \ \ $[U\tau _{1},U\mu_{1}]=U\mu_{1}\ ,$

\ \ \ \ $[U\tau _{1},U\mu_{2}]=0\ ,$

\ \ \ \ $[U\tau _{1},U\mu_{3}]=-U\mu_{3}\ ,$

\ \ \ \ $[U\tau _{2},Uz_{1j}]=0,(0\leq j\leq 7),$

\ \ \ \ $[U\tau _{2},Uz_{2j}]=-\frac{1}{2}Uz_{2j}\ ,(0\leq j\leq 7),$

\ \ \ \ $[U\tau _{2},Uz_{3j}]=\frac{1}{2}Uz_{3j}\ ,(0\leq j\leq 7),$

\ \ \ \ $[U\tau _{2},U\mu_{1}]=0\ ,$

\ \ \ \ $[U\tau _{2},U\mu_{2}]=U\mu_{2}\ ,$

\ \ \ \ $[U\tau _{2},U\mu_{3}]=-U\mu_{3}\ .$

\bigskip

\emph{Proof. \ }We have the above Lie bracket operations with
calculations using Maxima.\ \ \ \ \emph{Q.E.D.}

\bigskip

\emph{Lemma 17.12. }\ We have the following Lie bracket operations.

\ \ \ \ $[Ua_{1i},Uz_{1i}]=-2U\psi_{1},(0\leq i\leq 7),$

\ \ \ \ $[Ua_{1i},Uz_{1j}]=0,(\ 0\leq i\neq j\leq 7),$

\ \ \ \ $[Ua_{1i},Uz_{2j}]=Sn(i+1,j+1)Uw_{3k},(k=Ca(i+1,j+1),\ 0\leq i,j\leq 7),$

\ \ \ \ $[Ua_{1i},Uz_{3j}]=Sn(j+1,i+1)Uw_{2k},(k=Ca(j+1,i+1),\ 0\leq i,j\leq 7),$

\ \ \ \ $[Ua_{2i},Uz_{1j}]=Sn(j+1,i+1)Uw_{3k},(k=Ca(j+1,i+1),\ 0\leq i,j\leq 7),$

\ \ \ \ $[Ua_{2i},Uz_{2i}]=-2U\psi_{2},(0\leq i\leq 7),$

\ \ \ \ $[Ua_{2i},Uz_{2j}]=0,(\ 0\leq i\neq j\leq 7),$

\ \ \ \ $[Ua_{2i},Uz_{3j}]=Sn(i+1,j+1)Uw_{1k},(k=Ca(i+1,j+1),\ 0\leq i,j\leq 7),$

\ \ \ \ $[Ua_{3i},Uz_{1j}]=Sn(i+1,j+1)Uw_{2k},(k=Ca(i+1,j+1),\ 0\leq i,j\leq 7),$

\ \ \ \ $[Ua_{3i},Uz_{2j}]=Sn(j+1,i+1)Uw_{1k},(k=Ca(j+1,i+1),\ 0\leq i,j\leq 7),$

\ \ \ \ $[Ua_{3i},Uz_{3i}]=-2U\psi_{3},(0\leq i\leq 7),$

\ \ \ \ $[Ua_{3i},Uz_{3j}]=0,(\ 0\leq i\neq j\leq 7),$

\ \ \ \ $[Ua_{1i},U\mu_{1}]=-Uw_{1},(\ 0\leq i\leq 7),$

\ \ \ \ $[Ua_{1i},U\mu_{2}]=0,$

\ \ \ \ $[Ua_{1i},U\mu_{3}]=0,$

\ \ \ \ $[Ua_{2i},U\mu_{1}]=0,$

\ \ \ \ $[Ua_{2i},U\mu_{2}]=-Uw_{2},(\ 0\leq i\leq 7),$

\ \ \ \ $[Ua_{2i},U\mu_{3}]=0,$

\ \ \ \ $[Ua_{3i},U\mu_{1}]=0,$

\ \ \ \ $[Ua_{3i},U\mu_{2}]=0,$

\ \ \ \ $[Ua_{3i},U\mu_{3}]=-Uw_{3},(\ 0\leq i\leq 7),$

\ \ \ \ $[U\alpha _{1},Uz_{1i}]=-Uw_{1},(\ 0\leq i\leq 7),$

\ \ \ \ $[U\alpha _{1},Uz_{2i}]=0,(\ 0\leq i\leq 7),$

\ \ \ \ $[U\alpha _{1},Uz_{3i}]=0,(\ 0\leq i\leq 7),$

\ \ \ \ $[U\alpha _{2},Uz_{1i}]=0,(\ 0\leq i\leq 7),$

\ \ \ \ $[U\alpha _{2},Uz_{2i}]=-Uw_{2},(\ 0\leq i\leq 7),$

\ \ \ \ $[U\alpha _{2},Uz_{3i}]=0,(\ 0\leq i\leq 7),$

\ \ \ \ $[U\alpha _{3},Uz_{1i}]=0,(\ 0\leq i\leq 7),$

\ \ \ \ $[U\alpha _{3},Uz_{2i}]=0,(\ 0\leq i\leq 7),$

\ \ \ \ $[U\alpha _{3},Uz_{3i}]=-Uw_{3},(\ 0\leq i\leq 7),$

\ \ \ \ $[U\alpha _{1},U\mu_{1}]=0,$

\ \ \ \ $[U\alpha _{1},U\mu_{2}]=U\psi_{3},$

\ \ \ \ $[U\alpha _{1},U\mu_{3}]=U\psi_{2},$

\ \ \ \ $[U\alpha _{2},U\mu_{1}]=U\psi_{3},$

\ \ \ \ $[U\alpha _{2},U\mu_{2}]=0,$

\ \ \ \ $[U\alpha _{2},U\mu_{3}]=U\psi_{1},$

\ \ \ \ $[U\alpha _{3},U\mu_{1}]=U\psi_{2},$

\ \ \ \ $[U\alpha _{3},U\mu_{2}]=U\psi_{1},$

\ \ \ \ $[U\alpha _{3},U\mu_{3}]=0.$

\bigskip

\emph{Proof. \ }We have the above Lie bracket operations with
calculations using Maxima.\ \ \ \ \emph{Q.E.D.}

\bigskip

\emph{Lemma 17.13. }\ We have the following Lie bracket operations.

\ \ \ \ $[Ub_{1i},Uz_{1i}]=2U\omega_{1},(0\leq i\leq 7),$

\ \ \ \ $[Ub_{1i},Uz_{1j}]=0$,$(\ 0\leq i\neq j\leq 7),$

\ \ \ \ $[Ub_{1i},Uz_{2j}]=0$,$(\ 0\leq i,j\leq 7),$

\ \ \ \ $[Ub_{1i},Uz_{3j}]=0$,$(\ 0\leq i,j\leq 7),$

\ \ \ \ $[Ub_{2i},Uz_{1j}]=0$,$(\ 0\leq i,j\leq 7),$

\ \ \ \ $[Ub_{2i},Uz_{2i}]=2U\omega_{1},(0\leq i\leq 7),$

\ \ \ \ $[Ub_{2i},Uz_{2j}]=0$,$(\ 0\leq i\neq j\leq 7),$

\ \ \ \ $[Ub_{2i},Uz_{3j}]=0$,$(\ 0\leq i,j\leq 7),$

\ \ \ \ $[Ub_{3i},Uz_{1j}]=0$,$(\ 0\leq i,j\leq 7),$

\ \ \ \ $[Ub_{3i},Uz_{2j}]=0$,$(\ 0\leq i,j\leq 7),$

\ \ \ \ $[Ub_{3i},Uz_{3i}]=2U\omega_{1},(0\leq i\leq 7),$

\ \ \ \ $[Ub_{3i},Uz_{3j}]=0$,$(\ 0\leq i\neq j\leq 7),$

\ \ \ \ $[Ub_{ij},U\mu_{k}]=0$,$(\ 1\leq i,k\leq 3,0\leq j\leq 7),$

\ \ \ \ $[U\beta _{i},uzkj]=0$,$(\ 1\leq i,k\leq 3,0\leq j\leq 7), $

\ \ \ \ $[U\beta _{i},U\mu_{i}]=U\omega_{1},(1\leq i\leq 3),$

\ \ \ \ $[U\beta _{i},U\mu_{j}]=0,(1\leq i\neq j\leq 3).$

\bigskip

\emph{Proof. \ }We have the above Lie bracket operations with
calculations using Maxima.\ \ \ \ \emph{Q.E.D.}

\bigskip

\emph{Lemma 17.14. }\ We have the following Lie bracket operations.

\ \ \ \ $[U\rho _{1},Uz_{ij}]=-\frac{1}{3}Uz_{ij},(1\leq i\leq 3,0\leq j\leq 7),$

\ \ \ \ $[U\rho _{1},U\mu_{i}]=-\frac{1}{3}U\mu_{i},(1\leq i\leq 3).$

\bigskip

\emph{Proof. \ }We have the above Lie bracket operations with
calculations using Maxima.\ \ \ \ \emph{Q.E.D.}

\bigskip

\emph{Lemma 17.15. }\ We have the following Lie bracket operations.

\ \ \ \ $[Ud_{ij},Uw_{1k}]=-Uw_{1j}$ $($in case of $k=i),$

\ \ \ \ $\ \ \ \ \ \ \ \ \ \ \ \ \ \ \ \ =Uw_{1i}$ $($in case of $k=j),$

\ \ \ \ $\ \ \ \ \ \ \ \ \ \ \ \ \ \ \ \ =0$ $($in case of $k\neq i,j),(0\leq i<j\leq 7),$

\ \ \ \ $[Ud_{ij},Uw_{2k}]=-\sum\limits_{0\leq n<l\leq 7}Mv^{2}(ki,kj)Uw_{2l}$ $($where $k=n)$

\ \ \ \ \ \ \ \ \ \ \ \ \ \ \ \ \ \ \ \ \ $\ \ \
+\sum\limits_{0\leq n<l\leq 7}Mv^{2}(ki,kj)Uw_{2n}$ $($where $k=l)$

\ \ \ \ $\ \ \ \ \ \ \ \ \ \ \ \ \ \ \ \ \ \ \ \
(ki=Nu(i+1,j+1),kj=Nu(n+1,l+1),0\leq i<j\leq 7),$

\ \ \ \ $[Ud_{ij},Uw_{3k}]=-\sum\limits_{0\leq n<l\leq 7}Mv(ki,kj)Uw_{3l}$ $($where $n=k)$

\ \ \ \ \ \ \ \ \ \ \ \ \ \ \ \ \ \ \ \ \ $\ \ \
+\sum\limits_{0\leq n<l\leq 7}Mv(ki,kj)Uw_{3n}$ $($where $l=k),$

\ \ \ \ $\ \ \ \ \ \ \ \ \ \ \ \ \ \ \ \ \ \ \ \
(ki=Nu(i+1,j+1),kj=Nu(n+1,l+1),0\leq i<j\leq 7),$

\ \ \ \ $[Ud_{ij},U\psi_{k}]=0,(0\leq i<j\leq 7,k=1,2,3).$

\ \ \ \ $[Um_{1i},Uw_{1i}]=U\psi_{2}-U\psi_{3},(0\leq i\leq 7),$

\ \ \ \ $[Um_{1i},Uw_{2j}]=\frac{1}{2}Sn(i+1,j+1)Uw_{3k},(k=Ca(i+1,j+1),\ 0\leq i,j\leq 7),$

\ \ \ \ $[Um_{1i},Uw_{3j}]=-\frac{1}{2}Sn(j+1,i+1)Uw_{2k},(k=Ca(j+1,i+1),\ 0\leq i,j\leq 7),$

\ \ \ \ $[Um_{1i},U\psi_{1}]=0\ ,(0\leq i\leq 7),$

\ \ \ \ $[Um_{1i},U\psi_{2}]=-\frac{1}{2}Uw_{1i}\ ,(0\leq i\leq 7),$

\ \ \ \ $[Um_{1i},U\psi_{3}]=\frac{1}{2}Uw_{1i}\ ,(0\leq i\leq 7),$

\ \ \ \ $[Um_{2i},Uw_{1j}]=-\frac{1}{2}Sn(j+1,i+1)Uw_{3k},(k=Ca(j+1,i+1),\ 0\leq i,j\leq 7),$

\ \ \ \ $[Um_{2i},Uw_{2i}]=-U\psi_{1}+U\psi_{3},(0\leq i\leq 7),$

\ \ \ \ $[Um_{2i},Uw_{2j}]=0\ ,(0\leq i\neq j\leq 7),$

\ \ \ \ $[Um_{2i},Uw_{3j}]=\frac{1}{2}Sn(i+1,j+1)Uw_{1k},(k=Ca(i+1,j+1),\ 0\leq i,j\leq 7),$

\ \ \ \ $[Um_{2i},U\psi_{1}]=\frac{1}{2}Uw_{2i}\ ,(0\leq i\leq 7)$,

\ \ \ \ $[Um_{2i},U\psi_{2}]=0\ ,(0\leq i\leq 7)$,

\ \ \ \ $[Um_{2i},U\psi_{3}]=-\frac{1}{2}Uw_{2i}\ ,(0\leq i\leq 7),$

\ \ \ \ $[Um_{3i},Uw_{1j}]=\frac{1}{2}Sn(i+1,j+1)Uw_{2k},(k=Ca(i+1,j+1),\ 0\leq i,j\leq 7),$

\ \ \ \ $[Um_{3i},Uw_{2j}]=-\frac{1}{2}Sn(j+1,i+1)Uw_{1k},(k=Ca(j+1,i+1),\ 0\leq i,j\leq 7),$

\ \ \ \ $[Um_{3i},Uw_{3i}]=U\psi_{1}-U\psi_{2},(0\leq i\leq 7),$

\ \ \ \ $[Um_{3i},Uw_{3j}]=0\ ,(0\leq i\neq j\leq 7),$

\ \ \ \ $[Um_{3i},U\psi_{1}]=-\frac{1}{2}Uw_{3i}\ ,(0\leq i\leq 7),$

\ \ \ \ $[Um_{3i},U\psi_{2}]=\frac{1}{2}Uw_{3i}\ ,(0\leq i\leq 7),$

\ \ \ \ $[Um_{3i},U\psi_{3}]=0\ ,(0\leq i\leq 7),$

\ \ \ \ $[Ut_{1i},Uw_{1j}]=0\ ,(0\leq i\neq j\leq 7),$

\ \ \ \ $[Ut_{1i},Uw_{1i}]=-U\psi_{2}-U\psi_{3},(0\leq i\leq 7),$

\ \ \ \ $[Ut_{1i},Uw_{2j}]=-\frac{1}{2}Sn(i+1,j+1)Uw_{3k},(k=Ca(i+1,j+1),\ 0\leq i,j\leq 7),$

\ \ \ \ $[Ut_{1i},Uw_{3j}]=-\frac{1}{2}Sn(j+1,i+1)Uw_{2k},(k=Ca(j+1,i+1),\ 0\leq i,j\leq 7),$

\ \ \ \ $[Ut_{1i},U\psi_{1}]=0\ ,(0\leq i\leq 7),$

\ \ \ \ $[Ut_{1i},U\psi_{2}]=-\frac{1}{2}Uw_{1i}\ ,(0\leq i\leq 7),$

\ \ \ \ $[Ut_{1i},U\psi_{3}]=-\frac{1}{2}Uw_{1i}\ ,(0\leq i\leq 7),$

\ \ \ \ $[Ut_{2i},Uw_{1j}]=-\frac{1}{2}Sn(i+1,j+1)Uw_{3k},(k=Ca(i+1,j+1),\ 0\leq i,j\leq 7),$

\ \ \ \ $[Ut_{2i},Uw_{2i}]=-U\psi_{1}-U\psi_{3},(0\leq i\leq 7),$

\ \ \ \ $[Ut_{2i},Uw_{2j}]=0\ ,(0\leq i\neq j\leq 7),$

\ \ \ \ $[Ut_{2i},Uw_{3j}]=-\frac{1}{2}Sn(i+1,j+1)Uw_{1k},(k=Ca(i+1,j+1),\ 0\leq i,j\leq 7),$

\ \ \ \ $[Ut_{2i},U\psi_{1}]=-\frac{1}{2}Uw_{2i}\ ,(0\leq i\leq 7)$,

\ \ \ \ $[Ut_{2i},U\psi_{2}]=0\ ,(0\leq i\leq 7)$,

\ \ \ \ $[Ut_{2i},U\psi_{3}]=-\frac{1}{2}Uw_{2i}\ ,(0\leq i\leq 7),$

\ \ \ \ $[Ut_{3i},Uw_{1j}]=-\frac{1}{2}Sn(i+1,j+1)Uw_{2k},(k=Ca(i+1,j+1),\ 0\leq i,j\leq 7),$

\ \ \ \ $[Ut_{3i},Uw_{2j}]=-\frac{1}{2}Sn(j+1,i+1)Uw_{1k},(k=Ca(j+1,i+1),\ 0\leq i,j\leq 7),$

\ \ \ \ $[Ut_{3i},Uw_{3i}]=-U\psi_{1}-U\psi_{2},(0\leq i\leq 7),$

\ \ \ \ $[Ut_{3i},Uw_{3j}]=0\ ,(0\leq i\neq j\leq 7),$

\ \ \ \ $[Ut_{3i},U\psi_{1}]=-\frac{1}{2}Uw_{3i}\ ,(0\leq i\leq 7),$

\ \ \ \ $[Ut_{3i},U\psi_{2}]=-\frac{1}{2}Uw_{3i}\ ,(0\leq i\leq 7),$

\ \ \ \ $[Ut_{3i},U\psi_{3}]=0\ ,(0\leq i\leq 7),$

\ \ \ \ $[U\tau _{1},Uw_{1j}]=\frac{1}{2}Uw_{1j}\ ,(0\leq j\leq 7),$

\ \ \ \ $[U\tau _{1},Uw_{2j}]=0$ $,(0\leq j\leq 7),$

\ \ \ \ $[U\tau _{1},Uw_{3j}]=-\frac{1}{2}Uw_{3j}\ ,(0\leq j\leq 7),$

\ \ \ \ $[U\tau _{1},U\psi_{1}]=-U\psi_{1}\ ,$

\ \ \ \ $[U\tau _{1},U\psi_{2}]=0\ ,$

\ \ \ \ $[U\tau _{1},U\psi_{3}]=U\psi_{3}\ ,$

\ \ \ \ $[U\tau _{2},Uw_{1j}]=0,(0\leq j\leq 7),$

\ \ \ \ $[U\tau _{2},Uw_{2j}]=\frac{1}{2}Uw_{2j}\ ,(0\leq j\leq 7),$

\ \ \ \ $[U\tau _{2},Uw_{3j}]=-\frac{1}{2}Uw_{3j}\ ,(0\leq j\leq 7),$

\ \ \ \ $[U\tau _{2},U\psi_{1}]=0\ ,$

\ \ \ \ $[U\tau _{2},U\psi_{2}]=-U\psi_{2}\ ,$

\ \ \ \ $[U\tau _{2},U\psi_{3}]=U\psi_{3}\ .$

\bigskip

\emph{Proof. \ }We have the above Lie bracket operations with
calculations using Maxima.\ \ \ \ \emph{Q.E.D.}

\bigskip

\emph{Lemma 17.16. }\ We have the following Lie bracket operations.

\ \ \ \ $[Ub_{1i},Uw_{1i}]=-2U\mu_{1},(0\leq i\leq 7),$

\ \ \ \ $[Ub_{1i},Uw_{1j}]=0,(\ 0\leq i\neq j\leq 7),$

\ \ \ \ $[Ub_{1i},Uw_{2j}]=-Sn(i+1,j+1)Uz_{3k},(k=Ca(i+1,j+1),\ 0\leq i,j\leq 7),$

\ \ \ \ $[Ub_{1i},Uw_{3j}]=0,(\ 0\leq i,j\leq 7),$

\ \ \ \ $[Ub_{2i},Uw_{1j}]=0,(\ 0\leq i,j\leq 7),$

\ \ \ \ $[Ub_{2i},Uw_{2i}]=-2U\mu_{2},(0\leq i\leq 7),$

\ \ \ \ $[Ub_{2i},Uw_{2j}]=0,(\ 0\leq i\neq j\leq 7),$

\ \ \ \ $[Ub_{2i},Uw_{3j}]=-Sn(i+1,j+1)Uz_{1k},(k=Ca(i+1,j+1),\ 0\leq i,j\leq 7),$

\ \ \ \ $[Ub_{3i},Uw_{1j}]=-Sn(i+1,j+1)Uz_{2k},(k=Ca(i+1,j+1),\ 0\leq i,j\leq 7),$

\ \ \ \ $[Ub_{3i},Uw_{2j}]=0,(\ 0\leq i\neq j\leq 7),$

\ \ \ \ $[Ub_{3i},Uw_{3i}]=-2U\mu_{3},(0\leq i\leq 7),$

\ \ \ \ $[Ub_{3i},Uw_{3j}]=0,(\ 0\leq i\neq j\leq 7),$

\ \ \ \ $[Ub_{1i},U\psi_{1}]=-Uz_{1i},(\ 0\leq i\leq 7),$

\ \ \ \ $[Ub_{1i},U\psi_{2}]=0,$

\ \ \ \ $[Ub_{1i},U\psi_{3}]=0,$

\ \ \ \ $[Ub_{2i},U\psi_{1}]=0,$

\ \ \ \ $[Ub_{2i},U\psi_{2}]=-Uz_{2i},(\ 0\leq i\leq 7),$

\ \ \ \ $[Ub_{2i},U\psi_{3}]=0,$

\ \ \ \ $[Ub_{3i},U\psi_{1}]=0,$

\ \ \ \ $[Ub_{3i},U\psi_{2}]=0,$

\ \ \ \ $[Ub_{3i},U\psi_{3}]=-Uz_{3i},(\ 0\leq i\leq 7),$

\ \ \ \ $[U\beta _{1},Uw_{1i}]=-Uz_{1i},(\ 0\leq i\leq 7),$

\ \ \ \ $[U\beta _{1},Uw_{2i}]=0,(\ 0\leq i\leq 7),$

\ \ \ \ $[U\beta _{1},Uw_{3i}]=0,(\ 0\leq i\leq 7),$

\ \ \ \ $[U\beta _{2},Uw_{1i}]=0,(\ 0\leq i\leq 7),$

\ \ \ \ $[U\beta _{2},Uw_{2i}]=-Uz_{2i},(\ 0\leq i\leq 7),$

\ \ \ \ $[U\beta _{2},Uw_{3i}]=0,(\ 0\leq i\leq 7),$

\ \ \ \ $[U\beta _{3},Uw_{1i}]=0,(\ 0\leq i\leq 7),$

\ \ \ \ $[U\beta _{3},Uw_{2i}]=0,(\ 0\leq i\leq 7),$

\ \ \ \ $[U\beta _{3},Uw_{3i}]=-Uz_{3i},(\ 0\leq i\leq 7),$

\ \ \ \ $[U\beta _{1},U\psi_{1}]=0,$

\ \ \ \ $[U\beta _{1},U\psi_{2}]=U\mu_{3},$

\ \ \ \ $[U\beta _{1},U\psi_{3}]=U\mu_{2},$

\ \ \ \ $[U\beta _{2},U\psi_{1}]=U\mu_{3},$

\ \ \ \ $[U\beta _{2},U\psi_{2}]=0,$

\ \ \ \ $[U\beta _{2},U\psi_{3}]=U\mu_{1},$

\ \ \ \ $[U\beta _{3},U\psi_{1}]=U\mu_{2},$

\ \ \ \ $[U\beta _{3},U\psi_{2}]=U\mu_{1},$

\ \ \ \ $[U\beta _{3},U\psi_{3}]=0.$

\bigskip

\emph{Proof. \ }We have the above Lie bracket operations with
calculations using Maxima.\ \ \ \ \emph{Q.E.D.}

\bigskip

\emph{Lemma 17.17. }\ We have the following Lie bracket operations.

\ \ \ \ $[Ua_{1i},Uw_{1i}]=2U\zeta_{1},(0\leq i\leq 7),$

\ \ \ \ $[Ua_{1i},Uw_{1j}]=0$,$(\ 0\leq i\neq j\leq 7),$

\ \ \ \ $[Ua_{1i},Uw_{2j}]=0$,$(\ 0\leq i,j\leq 7),$

\ \ \ \ $[Ua_{1i},Uw_{3j}]=0$,$(\ 0\leq i,j\leq 7),$

\ \ \ \ $[Ua_{2i},Uw_{1j}]=0$,$(\ 0\leq i,j\leq 7),$

\ \ \ \ $[Ua_{2i},Uw_{2i}]=2U\zeta_{2},(0\leq i\leq 7),$

\ \ \ \ $[Ua_{2i},Uw_{2j}]=0$,$(\ 0\leq i\neq j\leq 7),$

\ \ \ \ $[Ua_{2i},Uw_{3j}]=0$,$(\ 0\leq i,j\leq 7),$

\ \ \ \ $[Ua_{3i},Uw_{1j}]=0$,$(\ 0\leq i,j\leq 7),$

\ \ \ \ $[Ua_{3i},Uw_{2j}]=0$,$(\ 0\leq i,j\leq 7),$

\ \ \ \ $[Ua_{3i},Uw_{3i}]=2U\zeta_{3},(0\leq i\leq 7),$

\ \ \ \ $[Ua_{3i},Uw_{3j}]=0$,$(\ 0\leq i\neq j\leq 7),$

\ \ \ \ $[Ua_{ij},U\psi_{k}]=0$,$(\ 1\leq i,k\leq 3,0\leq j\leq 7),$

\ \ \ \ $[U\alpha _{i},Uw_{kj}]=0$,$(\ 1\leq i,k\leq 3,0\leq j\leq 7),$

\ \ \ \ $[U\alpha _{i},U\psi_{i}]=U\zeta_{i},(1\leq i\leq 3),$

\ \ \ \ $[U\alpha _{i},U\psi_{j}]=0,(1\leq i\neq j\leq 3).$

\bigskip

\emph{Proof. \ }We have the above Lie bracket operations with
calculations using Maxima.\ \ \ \ \emph{Q.E.D.}

\bigskip

\emph{Lemma 17.18. }\ We have the following Lie bracket operations.

\ \ \ \ $[U\rho _{1},Uw_{ij}]=\frac{1}{3}Uw_{ij},(1\leq i\leq 3,0\leq
j\leq 7),$

\ \ \ \ $[U\rho _{1},U\psi_{i}]=\frac{1}{3}U\psi_{i},(1\leq i\leq 3).$

\bigskip

\emph{Proof. \ }We have the above Lie bracket operations with
calculations using Maxima.\ \ \ \ \emph{Q.E.D.}

\bigskip

\emph{Lemma 17.19. }\ We have the following Lie bracket operations.

\ \ \ \ $[Ud_{ij},U\zeta_{1}]=0,(0\leq i<j\leq 7),$

\ \ \ \ $[Um_{ij},U\zeta_{1}]=0,(1\leq i\leq 3,0\leq j\leq 7),$

\ \ \ \ $[Ut_{ij},U\zeta_{1}]=0,(1\leq i\leq 3,0\leq j\leq 7),$

\ \ \ \ $[U\tau _{i},U\zeta_{1}]=0,(1\leq i\leq 2),$

\ \ \ \ $[Ua_{ij},U\zeta_{1}]=0,(1\leq i\leq 3,0\leq j\leq 7),$

\ \ \ \ $[U\alpha _{i},U\zeta_{1}]=0,(1\leq i\leq 3),$

\ \ \ \ $[Ub_{ij},U\zeta_{1}]=Uw_{ij},(1\leq i\leq 3,0\leq j\leq 7),$

\ \ \ \ $[U\beta _{i},U\zeta_{1}]=U\psi_{i},(1\leq i\leq 3),$

\ \ \ \ $[U\rho _{1},U\zeta_{1}]=U\zeta_{1}$,

\ \ \ \ $[Ud_{ij},U\omega_{1}]=0,(0\leq i<j\leq 7),$

\ \ \ \ $[Um_{ij},U\omega_{1}]=0,(1\leq i\leq 3,0\leq j\leq 7),$

\ \ \ \ $[Ut_{ij},U\omega_{1}]=0,(1\leq i\leq 3,0\leq j\leq 7),$

\ \ \ \ $[U\tau _{i},U\omega_{1}]=0,(1\leq i\leq 2),$

\ \ \ \ $[Ub_{ij},U\omega_{1}]=0,(1\leq i\leq 3,0\leq j\leq 7),$

\ \ \ \ $[U\beta _{i},U\omega_{1}]=0,(1\leq i\leq 3),$

\ \ \ \ $[Ua_{ij},U\omega_{1}]=Uz_{ij},(1\leq i\leq 3,0\leq j\leq 7),$

\ \ \ \ $[U\alpha _{i},U\omega_{1}]=U\mu_{i},(1\leq i\leq 3),$

\ \ \ \ $[U\rho _{1},U\omega_{1}]=-U\omega_{1}$

\bigskip

\emph{Proof. \ }We have the above Lie bracket operations with
calculations using Maxima.\ \ \ \ \emph{Q.E.D.}

\bigskip

\emph{Lemma 17.20. }\ We have the following Lie bracket operations.

\ \ \ \ $[Ud_{ij},Ur_{1}]=0,(0\leq i<j\leq 7),$

\ \ \ \ $[Um_{ij},Ur_{1}]=0,(1\leq i\leq 3,0\leq j\leq 7),$

\ \ \ \ $[Ut_{ij},Ur_{1}]=0,(1\leq i\leq 3,0\leq j\leq 7),$

\ \ \ \ $[U\tau _{i},Ur_{1}]=0,(1\leq i\leq 2),$

\ \ \ \ $[Ua_{ij},Ur_{1}]=0,(1\leq i\leq 3,0\leq j\leq 7),$

\ \ \ \ $[U\alpha _{i},Ur_{1}]=0,(1\leq i\leq 3),$

\ \ \ \ $[Ub_{ij},Ur_{1}]=0,(1\leq i\leq 3,0\leq j\leq 7),$

\ \ \ \ $[U\beta _{i},Ur_{1}]=0,(1\leq i\leq 3),$

\ \ \ \ $[U\rho _{1},Ur_{1}]=0,$

\ \ \ \ $[Ud_{ij},Us_{1}]=0,(0\leq i<j\leq 7),$

\ \ \ \ $[Um_{ij},Us_{1}]=0,(1\leq i\leq 3,0\leq j\leq 7),$

\ \ \ \ $[Ut_{ij},Us_{1}]=0,(1\leq i\leq 3,0\leq j\leq 7),$

\ \ \ \ $[U\tau _{i},Us_{1}]=0,(1\leq i\leq 2),$

\ \ \ \ $[Ub_{ij},Us_{1}]=0,(1\leq i\leq 3,0\leq j\leq 7),$

\ \ \ \ $[U\beta _{i},Us_{1}]=0,(1\leq i\leq 3),$

\ \ \ \ $[Ua_{ij},Us_{1}]=0,(1\leq i\leq 3,0\leq j\leq 7),$

\ \ \ \ $[U\alpha _{i},Us_{1}]=0,(1\leq i\leq 3),$

\ \ \ \ $[U\rho _{1},Us_{1}]=0,$

\ \ \ \ $[Ud_{ij},Uu_{1}]=0,(0\leq i<j\leq 7),$

\ \ \ \ $[Um_{ij},Uu_{1}]=0,(1\leq i\leq 3,0\leq j\leq 7),$

\ \ \ \ $[Ut_{ij},Uu_{1}]=0,(1\leq i\leq 3,0\leq j\leq 7),$

\ \ \ \ $[U\tau _{i},Uu_{1}]=0,(1\leq i\leq 2),$

\ \ \ \ $[Ub_{ij},Uu_{1}]=0,(1\leq i\leq 3,0\leq j\leq 7),$

\ \ \ \ $[U\beta _{i},Uu_{1}]=0,(1\leq i\leq 3),$

\ \ \ \ $[Ua_{ij},Uu_{1}]=0,(1\leq i\leq 3,0\leq j\leq 7),$

\ \ \ \ $[U\alpha _{i},Uu_{1}]=0,(1\leq i\leq 3),$

\ \ \ \ $[U\rho _{1},Uu_{1}]=0,$

\bigskip

\emph{Proof. \ }We have the above Lie bracket operations with
calculations using Maxima.\ \ \ \ \emph{Q.E.D.}

\bigskip

\emph{Lemma 17.21.} \ We have the following Lie bracket operations.

\ \ \ \ $[Ux_{ij},Ux_{kl}]=0,(1\leq i,k\leq 3,0\leq j,l\leq 7),$

\ \ \ \ $[U\chi _{i},Ux_{kl}]=0,(1\leq i,k\leq 3,0\leq l\leq 7),$

\ \ \ \ $[U\chi _{i},U\chi _{k}]=0,(1\leq i,k\leq 3),$

\ \ \ \ $[Ux_{1i},Uy_{1i}]=\frac{1}{2}Us_{1},(0\leq i\leq 7),$

\ \ \ \ $[Ux_{1i},Uy_{1j}]=0,(0\leq i\neq j\leq 7),$

\ \ \ \ $[Ux_{1i},Uy_{2j}]=0,(0\leq i,j\leq 7),$

\ \ \ \ $[Ux_{1i},Uy_{3j}]=0,(0\leq i,j\leq 7),$

\ \ \ \ $[Ux_{2i},Uy_{1j}]=0,(0\leq i,j\leq 7),$

\ \ \ \ $[Ux_{2i},Uy_{2i}]=\frac{1}{2}Us_{1},(0\leq i\leq 7),$

\ \ \ \ $[Ux_{2i},Uy_{2j}]=0,(0\leq i\neq j\leq 7),$

\ \ \ \ $[Ux_{2i},Uy_{3j}]=0,(0\leq i,j\leq 7),$

\ \ \ \ $[Ux_{3i},Uy_{1j}]=0,(0\leq i,j\leq 7),$

\ \ \ \ $[Ux_{3i},Uy_{2j}]=0,(0\leq i,j\leq 7),$

\ \ \ \ $[Ux_{3i},Uy_{3i}]=\frac{1}{2}Us_{1},(0\leq i\leq 7),$

\ \ \ \ $[Ux_{3i},Uy_{3j}]=0,(0\leq i\neq j\leq 7),$

\ \ \ \ $[Ux_{ij},U\gamma _{k}]=0,(1\leq i,k\leq 3,0\leq j\leq 7),$

\ \ \ \ $[U\chi _{k},Uy_{ij}]=0,(1\leq i,k\leq 3,0\leq j\leq 7),$

\ \ \ \ $[U\chi _{i},U\gamma _{i}]=\frac{1}{4}Us_{1},(1\leq i\leq 3),$

\ \ \ \ $[U\chi _{i},U\gamma _{k}]=0,(1\leq i,k\leq 3),$

\ \ \ \ $[Ux_{ij},U\xi _{1}]=0,(1\leq i\leq 3,0\leq j\leq 7),$

\ \ \ \ $[U\chi _{i},U\xi _{1}]=0,(1\leq i\leq 3),$

\ \ \ \ $[Ux_{ij},U\eta _{1}]=0,(1\leq i\leq 3,0\leq j\leq 7),$

\ \ \ \ $[U\chi _{i},U\eta _{1}]=0,(1\leq i\leq 3),$

\ \ \ \ $[Uy_{ij},Uy_{kl}]=0,(1\leq i,k\leq 3,0\leq j,l\leq 7),$

\ \ \ \ $[U\gamma _{i},Uy_{kl}]=0,(1\leq i,k\leq 3,0\leq l\leq 7),$

\ \ \ \ $[U\gamma _{i},U\gamma _{k}]=0,(1\leq i,k\leq 3),$

\ \ \ \ $[Uy_{ij},U\xi _{1}]=0,(1\leq i\leq 3,0\leq j\leq 7),$

\ \ \ \ $[U\gamma _{i},U\xi _{1}]=0,(1\leq i\leq 3),$

\ \ \ \ $[Uy_{ij},U\eta _{1}]=0,(1\leq i\leq 3,0\leq j\leq 7),$

\ \ \ \ $[U\gamma _{i},U\eta _{1}]=0,(1\leq i\leq 3),$

\ \ \ \ $[U\xi _{1},U\eta _{1}]=\frac{1}{4}Us_{1},(1\leq i\leq 3),$

\bigskip

\emph{Proof. \ }We have the above Lie bracket operations with
calculations using Maxima.\ \ \ \ \emph{Q.E.D.}

\bigskip

\emph{Lemma 17.22.} \ We have the following Lie bracket operations.

\ \ \ \ $[Ux_{1i},Uz_{1j}]=0\ ,(0\leq i\neq j\leq 7),$

\ \ \ \ $[Ux_{1i},Uz_{1i}]=-\frac{1}{2}U\beta _{1},(0\leq i\leq 7),$

\ \ \ \ $[Ux_{1i},Uz_{2j}]=\frac{1}{4}Sn(i+1,j+1)Ub_{3k},(k=Ca(i+1,j+1),\ 0\leq i,j\leq 7),$

\ \ \ \ $[Ux_{1i},Uz_{3j}]=\frac{1}{4}Sn(j+1,i+1)Ub_{2k},(k=Ca(j+1,i+1),\ 0\leq i,j\leq 7),$

\ \ \ \ $[Ux_{1i},U\mu_{1}]=-\frac{1}{4}Ub_{1i}\ ,(0\leq i\leq 7),$

\ \ \ \ $[Ux_{1i},U\mu_{2}]=0\ ,(0\leq i\leq 7),$

\ \ \ \ $[Ux_{1i},U\mu_{3}]=0\ ,(0\leq i\leq 7),$

\ \ \ \ $[Ux_{2i},Uz_{1j}]=\frac{1}{4}Sn(j+1,i+1)Ub_{3k},(k=Ca(j+1,i+1),\ 0\leq i,j\leq 7),$

\ \ \ \ $[Ux_{2i},Uz_{2i}]=-\frac{1}{2}U\beta _{2},(0\leq i\leq 7),$

\ \ \ \ $[Ux_{2i},Uz_{2j}]=0\ ,(0\leq i\neq j\leq 7),$

\ \ \ \ $[Ux_{2i},Uz_{3j}]=\frac{1}{4}Sn(i+1,j+1)Ub_{1k},(k=Ca(i+1,j+1),\ 0\leq i,j\leq 7),$

\ \ \ \ $[Ux_{2i},U\mu_{1}]=0\ ,(0\leq i\leq 7),$

\ \ \ \ $[Ux_{2i},U\mu_{2}]=-\frac{1}{4}Ub_{2i}\ ,(0\leq i\leq 7),$

\ \ \ \ $[Ux_{2i},U\mu_{3}]=0\ ,(0\leq i\leq 7),$

\ \ \ \ $[Ux_{3i},Uz_{1j}]=\frac{1}{4}Sn(i+1,j+1)Ub_{2k},(k=Ca(i+1,j+1),\ 0\leq i,j\leq 7),$

\ \ \ \ $[Ux_{3i},Uz_{2j}]=\frac{1}{4}Sn(j+1,i+1)Ub_{1k},(k=Ca(j+1,i+1),\ 0\leq i,j\leq 7),$

\ \ \ \ $[Ux_{3i},Uz_{3i}]=-\frac{1}{2}U\beta _{3},(0\leq i\leq 7),$

\ \ \ \ $[Ux_{3i},Uz_{3j}]=0\ ,(0\leq i\neq j\leq 7),$

\ \ \ \ $[Ux_{3i},U\mu_{1}]=0\ ,(0\leq i\leq 7),$

\ \ \ \ $[Ux_{3i},U\mu_{2}]=0\ ,(0\leq i\leq 7),$

\ \ \ \ $[Ux_{3i},U\mu_{3}]=-\frac{1}{4}Ub_{3i}\ ,(0\leq i\leq 7),$

\ \ \ \ $[U\chi _{1},Uz_{1j}]=-\frac{1}{4}Ub_{1j}\ ,(0\leq j\leq 7),$

\ \ \ \ $[U\chi _{1},Uz_{2j}]=0$ $,(0\leq j\leq 7),$

\ \ \ \ $[U\chi _{1},Uz_{3j}]=0$ $,(0\leq j\leq 7),$

\ \ \ \ $[U\chi _{1},U\mu_{1}]=0$ $,$

\ \ \ \ $[U\chi _{1},U\mu_{2}]=\frac{1}{4}U\beta _{3}\ ,$

\ \ \ \ $[U\chi _{1},U\mu_{3}]=\frac{1}{4}U\beta _{2}\ ,$

\ \ \ \ $[U\chi _{2},Uz_{1j}]=0$ $,(0\leq j\leq 7),$

\ \ \ \ $[U\chi _{2},Uz_{2j}]=-\frac{1}{4}Ub_{2j},(0\leq j\leq 7),$

\ \ \ \ $[U\chi _{2},Uz_{3j}]=0$ $,(0\leq j\leq 7),$

\ \ \ \ $[U\chi _{2},U\mu_{1}]=\frac{1}{4}U\beta _{3}\ ,$

\ \ \ \ $[U\chi _{2},U\mu_{2}]=0,$

\ \ \ \ $[U\chi _{2},U\mu_{3}]=\frac{1}{4}U\beta _{1}\ ,$

\ \ \ \ $[U\chi _{3},Uz_{1j}]=0\ ,$

\ \ \ \ $[U\chi _{3},Uz_{2j}]=0\ ,$

\ \ \ \ $[U\chi _{3},Uz_{3j}]=-\frac{1}{4}Ub_{3j},(0\leq j\leq 7),$

\ \ \ \ $[U\chi _{3},U\mu_{1}]=\frac{1}{4}U\beta _{2}\ ,$

\ \ \ \ $[U\chi _{3},U\mu_{2}]=\frac{1}{4}U\beta _{1}\ ,$

\ \ \ \ $[U\chi _{3},U\mu_{3}]=0\ .$

\bigskip

\emph{Proof. \ }We have the above Lie bracket operations with
calculations using Maxima.\ \ \ \ \emph{Q.E.D.}

\bigskip

\emph{Lemma 17.23. }\ We have the following Lie bracket operations.

\ \ \ \ $[Ux_{1i},Uw_{1j}]=-\frac{1}{2}Ud_{ij}\ ,(0\leq i\neq j\leq 7),$

\ \ \ \ $[Ux_{1i},Uw_{1i}]=\frac{1}{3}U\tau _{1}-\frac{1}{6}U\tau _{2}+%
\frac{1}{4}U\rho _{1}-\frac{1}{4}Ur_{1},(0\leq i\leq 7),$

\ \ \ \ $[Ux_{1i},Uw_{2j}]=\frac{1}{4}Sn(i+1,j+1)Um_{3k}-\frac{1}{4}Sn(i+1,j+1)Ut_{3k},$

\ \ \ \ \ \ \ \ \ \ \ \ \ \ \ \ $(k=Ca(i+1,j+1),0\leq i,j\leq 7),$

\ \ \ \ $[Ux_{1i},Uw_{3j}]=-\frac{1}{4}Sn(j+1,i+1)Um_{2k}-\frac{1}{4}Sn(j+1,i+1)Ut_{2k},$

\ \ \ \ \ \ \ \ \ \ \ \ \ \ \ \ $(k=Ca(j+1,i+1),0\leq i,j\leq 7),$

\ \ \ \ $[Ux_{1i},U\psi_{1}]=0,$

\ \ \ \ $[Ux_{1i},U\psi_{2}]=\frac{1}{4}Um_{1i}-\frac{1}{4}Ut_{1i},(0\leq
i\leq 7),$

\ \ \ \ $[Ux_{1i},U\psi_{3}]=-\frac{1}{4}Um_{1i}-\frac{1}{4}Ut_{1i},(0\leq i\leq 7),$

\ \ \ \ $[Ux_{2i},Uw_{1j}]=-\frac{1}{4}Sn(j+1,i+1)Um_{3k}-\frac{1}{4}Sn(j+1,i+1)Ut_{3k},$

\ \ \ \ \ \ \ \ \ \ \ \ \ \ \ \ \ $(k=Ca(j+1,i+1),0\leq i,j\leq 7),$

\ \ \ \ $[Ux_{2i},Uw_{2j}]=-\frac{1}{2}\sum\limits_{0\leq n<l\leq 7}Mv(ki,kj)Ud_{nl},$

\ \ \ \ \ \ \ \ \ \ \ \ \ \ \ \ \ $(ki=Nu(i+1,j+1),kj=Nu(n+1,l+1),0\leq i\neq j\leq 7),$

\ \ \ \ $[Ux_{2i},Uw_{2i}]=-\frac{1}{6}U\tau _{1}+\frac{1}{3}U\tau _{2}+%
\frac{1}{4}U\rho _{1}-\frac{1}{4}Ur_{1},(0\leq i\leq 7),$

\ \ \ \ $[Ux_{2i},Uw_{3j}]=\frac{1}{4}Sn(i+1,j+1)Um_{1k}-\frac{1}{4}Sn(i+1,j+1)Ut_{1k},$

\ \ \ \ \ \ \ \ \ \ \ \ \ \ \ \ \ $(k=Ca(i+1,j+1),\ 0\leq i,j\leq 7),$

\ \ \ \ $[Ux_{2i},U\psi_{1}]=-\frac{1}{4}Um_{2i}-\frac{1}{4}Ut_{2i},(0\leq i\leq 7),$

\ \ \ \ $[Ux_{2i},U\psi_{2}]=0,$

\ \ \ \ $[Ux_{2i},U\psi_{3}]=\frac{1}{4}Um_{2i}-\frac{1}{4}Ut_{2i},(0\leq i\leq 7),$

\ \ \ \ $[Ux_{3i},Uw_{1j}]=\frac{1}{4}Sn(i+1,j+1)Um_{2k}-\frac{1}{4}Sn(i+1,j+1)Ut_{2k},$

\ \ \ \ \ \ \ \ \ \ \ \ \ \ \ \ \ $(k=Ca(i+1,j+1),\ 0\leq i,j\leq 7),$

\ \ \ \ $[Ux_{3i},Uw_{2j}]=-\frac{1}{4}Sn(j+1,i+1)Um_{1k}-\frac{1}{4}Sn(j+1,i+1)Ut_{1k},$

\ \ \ \ \ \ \ \ \ \ \ \ \ \ \ \ \ $(k=Ca(j+1,i+1),0\leq i,j\leq 7),$

\ \ \ \ $[Ux_{3i},Uw_{3j}]=-\frac{1}{2}\sum\limits_{0\leq n<l\leq 7}Mv^{2}(ki,kj)Ud_{nl},$

\ \ \ \ \ \ \ \ \ \ \ \ \ \ \ \ \ $(ki=Nu(i+1,j+1),kj=Nu(n+1,l+1),0\leq i,j\leq 7),$

\ \ \ \ $[Ux_{3i},Uw_{3i}]=-\frac{1}{6}U\tau _{1}-\frac{1}{6}U\tau _{2}+%
\frac{1}{4}U\rho _{1}-\frac{1}{4}Ur_{1},(0\leq i\leq 7),$

\ \ \ \ $[Ux_{3i},U\psi_{1}]=\frac{1}{4}Um_{3i}-\frac{1}{4}Ut_{3i},(0\leq i\leq 7),$

\ \ \ \ $[Ux_{3i},U\psi_{2}]=-\frac{1}{4}Um_{3i}-\frac{1}{4}Ut_{3i},(0\leq i\leq 7),$

\ \ \ \ $[Ux_{3i},U\psi_{3}]=0,$

\ \ \ \ $[U\chi _{1},Uw_{1j}]=0,(0\leq j\leq 7),$

\ \ \ \ $[U\chi _{1},Uw_{2j}]=\frac{1}{4}Um_{2j}-\frac{1}{4}Ut_{2j},(0\leq
j\leq 7),$

\ \ \ \ $[U\chi _{1},Uw_{3j}]=-\frac{1}{4}Um_{3j}-\frac{1}{4}Ut_{3j},(0\leq j\leq 7),$

\ \ \ \ $[U\chi _{1},U\psi_{1}]=-\frac{1}{3}U\tau _{1}+\frac{1}{6}%
U\tau _{2}+\frac{1}{8}U\rho _{1}-\frac{1}{8}Ur_{1},$

\ \ \ \ $[U\chi _{1},U\psi_{2}]=0,$

\ \ \ \ $[U\chi _{1},U\psi_{3}]=0,$

\ \ \ \ $[U\chi _{2},Uw_{1j}]=-\frac{1}{4}Um_{1j}-\frac{1}{4}Ut_{1j},(0\leq j\leq 7),$

\ \ \ \ $[U\chi _{2},Uw_{2j}]=0,(0\leq j\leq 7),$

\ \ \ \ $[U\chi _{2},Uw_{3j}]=\frac{1}{4}Um_{3j}-\frac{1}{4}Ut_{3j},(0\leq
j\leq 7),$

\ \ \ \ $[U\chi _{2},U\psi_{1}]=0,$

\ \ \ \ $[U\chi _{2},U\psi_{2}]=\frac{1}{6}U\tau _{1}-\frac{1}{3}U\tau
_{2}+\frac{1}{8}U\rho _{1}-\frac{1}{8}Ur_{1},$

\ \ \ \ $[U\chi _{2},U\psi_{3}]=0,$

\ \ \ \ $[U\chi _{3},Uw_{1j}]=\frac{1}{4}Um_{1j}-\frac{1}{4}Ut_{1j},(0\leq
j\leq 7),$

\ \ \ \ $[U\chi _{3},Uw_{2j}]=-\frac{1}{4}Um_{2j}-\frac{1}{4}Ut_{2j},(0\leq j\leq 7),,$

\ \ \ \ $[U\chi _{3},Uw_{3j}]=0,(0\leq j\leq 7),$

\ \ \ \ $[U\chi _{3},U\psi_{1}]=0,$

\ \ \ \ $[U\chi _{3},U\psi_{2}]=0,$

\ \ \ \ $[U\chi _{3},U\psi_{3}]=\frac{1}{6}U\tau _{1}+\frac{1}{6}U\tau
_{2}+\frac{1}{8}U\rho _{1}-\frac{1}{8}Ur_{1}.$

\bigskip

\emph{Proof. \ }We have the above Lie bracket operations with
calculations using Maxima.\ \ \ \ \emph{Q.E.D.}

\bigskip

\emph{Lemma 17.24. }\ We have the following Lie bracket operations.

\ \ \ \ $[Ux_{1i},U\zeta_{1}]=\frac{1}{4}Ua_{1i},(0\leq i\leq 7),$

\ \ \ \ $[Ux_{2i},U\zeta_{1}]=\frac{1}{4}Ua_{2i},(0\leq i\leq 7),$

\ \ \ \ $[Ux_{3i},U\zeta_{1}]=\frac{1}{4}Ua_{3i},(0\leq i\leq 7),$

\ \ \ \ $[U\chi _{1},U\zeta_{1}]=\frac{1}{4}U\alpha _{1},$

\ \ \ \ $[U\chi _{2},U\zeta_{1}]=\frac{1}{4}U\alpha _{2},$

\ \ \ \ $[U\chi _{3},U\zeta_{1}]=\frac{1}{4}U\alpha _{3},$

\ \ \ \ $[Ux_{1i},U\omega_{1}]=0,(0\leq i\leq 7),$

\ \ \ \ $[Ux_{2i},U\omega_{1}]=0,(0\leq i\leq 7),$

\ \ \ \ $[Ux_{3i},U\omega_{1}]=0,(0\leq i\leq 7),$

\ \ \ \ $[U\chi _{1},U\omega_{1}]=0,$

\ \ \ \ $[U\chi _{2},U\omega_{1}]=0,$

\ \ \ \ $[U\chi _{3},U\omega_{1}]=0.$

\bigskip

\emph{Proof. \ }We have the above Lie bracket operations with
calculations using Maxima.\ \ \ \ \emph{Q.E.D.}

\bigskip

\emph{Lemma 17.25. }\ We have the following Lie bracket operations.

\ \ \ \ $[Uz_{1i},Uy_{1j}]=\frac{1}{2}Ud_{ij}\ ,(0\leq i\neq j\leq 7),$

\ \ \ \ $[Uz_{1i},Uy_{1i}]=-\frac{1}{3}U\tau _{1}+\frac{1}{6}U\tau _{2}-%
\frac{1}{4}U\rho _{1}-\frac{1}{4}Ur_{1},(0\leq i\leq 7),$

\ \ \ \ $[Uz_{1i},Uy_{2j}]=-\frac{1}{4}Sn(i+1,j+1)Um_{3k}+\frac{1}{4}Sn(i+1,j+1)Ut_{3k},$

\ \ \ \ \ \ \ \ \ \ \ \ \ \ \ \ \ $(k=Ca(i+1,j+1),\ 0\leq i,j\leq 7),$

\ \ \ \ $[Uz_{1i},Uy_{3j}]=\frac{1}{4}Sn(j+1,i+1)Um_{2k}+\frac{1}{4}Sn(j+1,i+1)Ut_{2k},$

\ \ \ \ \ \ \ \ \ \ \ \ \ \ \ \ \ $(k=Ca(j+1,i+1),\ 0\leq i,j\leq 7),$

\ \ \ \ $[Uz_{1i},U\gamma _{1}]=0,$

\ \ \ \ $[Uz_{1i},U\gamma _{2}]=-\frac{1}{4}Um_{1i}+\frac{1}{4}Ut_{1i},(0\leq i\leq 7),$

\ \ \ \ $[Uz_{1i},U\gamma _{3}]=\frac{1}{4}Um_{1i}+\frac{1}{4}Ut_{1i},(0\leq i\leq 7),$

\ \ \ \ $[Uz_{2i},Uy_{1j}]=\frac{1}{4}Sn(j+1,i+1)Um_{3k}+\frac{1}{4}Sn(j+1,i+1)Ut_{3k},$

\ \ \ \ \ \ \ \ \ \ \ \ \ \ \ \ \ $(k=Ca(j+1,i+1),\ 0\leq i,j\leq 7),$

\ \ \ \ $[Uz_{2i},Uy_{2j}]=\frac{1}{2}\sum\limits_{0\leq n<l\leq 7}Mv(ki,kj)Ud_{nl},$

\ \ \ \ \ \ \ \ \ \ \ \ \ \ \ \ \ $(ki=Nu(i+1,j+1),kj=Nu(n+1,l+1),\ 0\leq i\neq j\leq 7),$

\ \ \ \ $[Uz_{2i},Uy_{2i}]=\frac{1}{6}U\tau _{1}-\frac{1}{3}U\tau _{2}-%
\frac{1}{4}U\rho _{1}-\frac{1}{4}Ur_{1},(0\leq i\leq 7),$

\ \ \ \ $[Uz_{2i},Uy_{3j}]=-\frac{1}{4}Sn(i+1,j+1)Um_{1k}+\frac{1}{4}Sn(i+1,j+1)Ut_{1k},$

\ \ \ \ \ \ \ \ \ \ \ \ \ \ \ \ \ $(k=Ca(i+1,j+1),\ 0\leq i,j\leq 7),$

\ \ \ \ $[Uz_{2i},U\gamma _{1}]=\frac{1}{4}Um_{2i}+\frac{1}{4}Ut_{2i},(0\leq i\leq 7),$

\ \ \ \ $[Uz_{2i},U\gamma _{2}]=0,$

\ \ \ \ $[Uz_{2i},U\gamma _{3}]=-\frac{1}{4}Um_{2i}+\frac{1}{4}Ut_{2i},(0\leq i\leq 7),$

\ \ \ \ $[Uz_{3i},Uy_{1j}]=-\frac{1}{4}Sn(i+1,j+1)Um_{2k}+\frac{1}{4}Sn(i+1,j+1)Ut_{2k},$

\ \ \ \ \ \ \ \ \ \ \ \ \ \ \ \ \ $(k=Ca(i+1,j+1),\ 0\leq i,j\leq 7),$

\ \ \ \ $[Uz_{3i},Uy_{2j}]=\frac{1}{4}Sn(j+1,i+1)Um_{1k}+\frac{1}{4}Sn(j+1,i+1)Ut_{1k},$

\ \ \ \ \ \ \ \ \ \ \ \ \ \ \ \ \ $(k=Ca(j+1,i+1),\ 0\leq i,j\leq 7),$

\ \ \ \ $[Uz_{3i},Uy_{3j}]=\frac{1}{2}\sum\limits_{0\leq n<l\leq 7}Mv^{2}(ki,kj)Ud_{nl},$

\ \ \ \ \ \ \ \ \ \ \ \ \ \ \ \ \ $(ki=Nu(i+1,j+1),kj=Nu(n+1,l+1),0\leq i\neq j\leq 7),$

\ \ \ \ $[Uz_{3i},Uy_{3i}]=\frac{1}{6}U\tau _{1}+\frac{1}{6}U\tau _{2}-%
\frac{1}{4}U\rho _{1}-\frac{1}{4}Ur_{1},(0\leq i\leq 7),$

\ \ \ \ $[Uz_{3i},U\gamma _{1}]=-\frac{1}{4}Um_{3i}+\frac{1}{4}Ut_{3i},(0\leq i\leq 7),$

\ \ \ \ $[Uz_{3i},U\gamma _{2}]=\frac{1}{4}Um_{3i}+\frac{1}{4}Ut_{3i},(0\leq i\leq 7),$

\ \ \ \ $[Uz_{3i},U\gamma _{3}]=0,$

\ \ \ \ $[U\mu_{1},Uy_{1j}]=0,(0\leq j\leq 7),$

\ \ \ \ $[U\mu_{1},Uy_{2j}]=-\frac{1}{4}Um_{2j}+\frac{1}{4}Ut_{2j},(0\leq
j\leq 7),$

\ \ \ \ $[U\mu_{1},Uy_{3j}]=\frac{1}{4}Um_{3j}+\frac{1}{4}Ut_{3j},(0\leq
j\leq 7),$

\ \ \ \ $[U\mu_{1},U\gamma _{1}]=\frac{1}{3}U\tau _{1}-\frac{1}{6}%
U\tau _{2}-\frac{1}{8}U\rho _{1}-\frac{1}{8}Ur_{1},$

\ \ \ \ $[U\mu_{1},U\gamma _{2}]=0,$

\ \ \ \ $[U\mu_{1},U\gamma _{3}]=0,$

\ \ \ \ $[U\mu_{2},Uy_{1j}]=\frac{1}{4}Um_{1j}+\frac{1}{4}Ut_{1j},(0\leq
j\leq 7),$

\ \ \ \ $[U\mu_{2},Uy_{2j}]=0,(0\leq j\leq 7),$

\ \ \ \ $[U\mu_{2},Uy_{3j}]=-\frac{1}{4}Um_{3j}+\frac{1}{4}Ut_{3j},(0\leq
j\leq 7),$

\ \ \ \ $[U\mu_{2},U\gamma _{1}]=0,$

\ \ \ \ $[U\mu_{2},U\gamma _{2}]=-\frac{1}{6}U\tau _{1}+\frac{1}{3}%
U\tau _{2}-\frac{1}{8}U\rho _{1}-\frac{1}{8}Ur_{1},$

\ \ \ \ $[U\mu_{2},U\gamma _{3}]=0,$

\ \ \ \ $[U\mu_{3},Uy_{1j}]=-\frac{1}{4}Um_{1j}+\frac{1}{4}Ut_{1j},(0\leq
j\leq 7),$

\ \ \ \ $[U\mu_{3},Uy_{2j}]=\frac{1}{4}Um_{2j}+\frac{1}{4}Ut_{2j},(0\leq
j\leq 7)$,

\ \ \ \ $[U\mu_{3},Uy_{3j}]=0,(0\leq j\leq 7),$

\ \ \ \ $[U\mu_{3},U\gamma _{1}]=0,$

\ \ \ \ $[U\mu_{3},U\gamma _{2}]=0,$

\ \ \ \ $[U\mu_{3},U\gamma _{3}]=-\frac{1}{6}U\tau _{1}-\frac{1}{6}%
U\tau _{2}-\frac{1}{8}U\rho _{1}-\frac{1}{8}Ur_{1}.$

\bigskip

\emph{Proof. \ }We have the above Lie bracket operations with
calculations using Maxima.\ \ \ \ \emph{Q.E.D.}

\bigskip

\emph{Lemma 17.26.} \ We have the following Lie bracket operations.

\ \ \ \ $[Uy_{1i},Uw_{1j}]=0\ ,(0\leq i\neq j\leq 7),$

\ \ \ \ $[Uy_{1i},Uw_{1i}]=\frac{1}{2}U\alpha _{1},(0\leq i\leq 7),$

\ \ \ \ $[Uy_{1i},Uw_{2j}]=-\frac{1}{4}Sn(i+1,j+1)Ua_{3k},(k=Ca(i+1,j+1),\ 0\leq i,j\leq 7),$

\ \ \ \ $[Uy_{1i},Uw_{3j}]=-\frac{1}{4}Sn(j+1,i+1)Ua_{2k},(k=Ca(j+1,i+1),\ 0\leq i,j\leq 7),$

\ \ \ \ $[Uy_{1i},U\psi_{1}]=\frac{1}{4}Ua_{1i}\ ,(0\leq i\leq 7),$

\ \ \ \ $[Uy_{1i},U\psi_{2}]=0\ ,(0\leq i\leq 7),$

\ \ \ \ $[Uy_{1i},U\psi_{3}]=0\ ,(0\leq i\leq 7),$

\ \ \ \ $[Uy_{2i},Uw_{1j}]=-\frac{1}{4}Sn(j+1,i+1)Ua_{3k},(k=Ca(j+1,i+1),\ 0\leq i,j\leq 7),$

\ \ \ \ $[Uy_{2i},Uw_{2i}]=\frac{1}{2}U\alpha _{2},(0\leq i\leq 7),$

\ \ \ \ $[Uy_{2i},Uw_{2j}]=0\ ,(0\leq i\neq j\leq 7),$

\ \ \ \ $[Uy_{2i},Uw_{3j}]=-\frac{1}{4}Sn(i+1,j+1)Ua_{1k},(k=Ca(i+1,j+1),\ 0\leq i,j\leq 7),$

\ \ \ \ $[Uy_{2i},U\psi_{1}]=0\ ,(0\leq i\leq 7),$

\ \ \ \ $[Uy_{2i},U\psi_{2}]=\frac{1}{4}Ua_{2i}\ ,(0\leq i\leq 7),$

\ \ \ \ $[Uy_{2i},U\psi_{3}]=0\ ,(0\leq i\leq 7),$

\ \ \ \ $[Uy_{3i},Uw_{1j}]=-\frac{1}{4}Sn(i+1,j+1)Ua_{2k},(k=Ca(i+1,j+1),\ 0\leq i,j\leq 7),$

\ \ \ \ $[Uy_{3i},Uw_{2j}]=-\frac{1}{4}Sn(j+1,i+1)Ua_{1k},(k=Ca(j+1,i+1),\ 0\leq i,j\leq 7),$

\ \ \ \ $[Uy_{3i},Uw_{3i}]=\frac{1}{2}U\alpha _{3},(0\leq i\leq 7),$

\ \ \ \ $[Uy_{3i},Uw_{3j}]=0\ ,(0\leq i\neq j\leq 7),$

\ \ \ \ $[Uy_{3i},U\psi_{1}]=0\ ,(0\leq i\leq 7),$

\ \ \ \ $[Uy_{3i},U\psi_{2}]=0\ ,(0\leq i\leq 7),$

\ \ \ \ $[Uy_{3i},U\psi_{3}]=\frac{1}{4}Ua_{3i}\ ,(0\leq i\leq 7),$

\ \ \ \ $[U\gamma _{1},Uw_{1j}]=\frac{1}{4}Ua_{1j}\ ,(0\leq j\leq 7),$

\ \ \ \ $[U\gamma _{1},Uw_{2j}]=0$ $,(0\leq j\leq 7),$

\ \ \ \ $[U\gamma _{1},Uw_{3j}]=0$ $,(0\leq j\leq 7),$

\ \ \ \ $[U\gamma _{1},U\psi_{1}]=0$ $,$

\ \ \ \ $[U\gamma _{1},U\psi_{2}]=-\frac{1}{4}U\alpha _{3}\ ,$

\ \ \ \ $[U\gamma _{1},U\psi_{3}]=-\frac{1}{4}U\alpha _{2}\ ,$

\ \ \ \ $[U\gamma _{2},Uw_{1j}]=0$ $,(0\leq j\leq 7),$

\ \ \ \ $[U\gamma _{2},Uw_{2i}]=\frac{1}{4}Ua_{2i},(0\leq i\leq 7),$

\ \ \ \ $[U\gamma _{2},Uw_{3j}]=0$ $,(0\leq j\leq 7),$

\ \ \ \ $[U\gamma _{2},U\psi_{1}]=-\frac{1}{4}U\alpha _{3}\ ,$

\ \ \ \ $[U\gamma _{2},U\psi_{2}]=0,$

\ \ \ \ $[U\gamma _{2},U\psi_{3}]=-\frac{1}{4}U\alpha _{1}\ ,$

\ \ \ \ $[U\gamma _{3},Uw_{1j}]=0\ ,$

\ \ \ \ $[U\gamma _{3},Uw_{2j}]=0\ ,$

\ \ \ \ $[U\gamma _{3},Uw_{3i}]=\frac{1}{4}Ua_{3i},(0\leq i\leq 7),$

\ \ \ \ $[U\gamma _{3},U\psi_{1}]=-\frac{1}{4}U\alpha _{2}\ ,$

\ \ \ \ $[U\gamma _{3},U\psi_{2}]=-\frac{1}{4}U\alpha _{1}\ ,$

\ \ \ \ $[U\gamma _{3},U\psi_{3}]=0\ .$

\bigskip

\emph{Proof. \ }We have the above Lie bracket operations with
calculations using Maxima.\ \ \ \ \emph{Q.E.D.}

\bigskip

\emph{Lemma 17.27. }\ We have the following Lie bracket operations.

\ \ \ \ $[Uy_{1i},U\zeta_{1}]=0,(0\leq i\leq 7),$

\ \ \ \ $[Uy_{2i},U\zeta_{1}]=0,(0\leq i\leq 7),$

\ \ \ \ $[Uy_{3i},U\zeta_{1}]=0,(0\leq i\leq 7),$

\ \ \ \ $[U\gamma _{1},U\zeta_{1}]=0,$

\ \ \ \ $[U\gamma _{2},U\zeta_{1}]=0,$

\ \ \ \ $[U\gamma _{3},U\zeta_{1}]=0,$

\ \ \ \ $[Uy_{1i},U\omega_{1}]=-\frac{1}{4}Ub_{1i},(0\leq i\leq 7),$

\ \ \ \ $[Uy_{2i},U\omega_{1}]=-\frac{1}{4}Ub_{2i}0,(0\leq i\leq 7),$

\ \ \ \ $[Uy_{3i},U\omega_{1}]=-\frac{1}{4}Ub_{3i},(0\leq i\leq 7),$

\ \ \ \ $[U\gamma _{1},U\omega_{1}]=-\frac{1}{4}U\beta _{1},$

\ \ \ \ $[U\gamma _{2},U\omega_{1}]=-\frac{1}{4}U\beta _{2},$

\ \ \ \ $[U\gamma _{3},U\omega_{1}]=-\frac{1}{4}U\beta _{3}.$

\bigskip

\emph{Proof. \ }We have the above Lie bracket operations with
calculations using Maxima.\ \ \ \ \emph{Q.E.D.}

\bigskip

\emph{Lemma 17.28. }\ We have the following Lie bracket operations.

\ \ \ \ $[U\xi _{1},Uz_{1i}]=\frac{1}{4}Ua_{1i},(0\leq i\leq 7),$

\ \ \ \ $[U\xi _{1},Uz_{2i}]=\frac{1}{4}Ua_{2i},(0\leq i\leq 7),$

\ \ \ \ $[U\xi _{1},Uz_{3i}]=\frac{1}{4}Ua_{3i},(0\leq i\leq 7),$

\ \ \ \ $[U\xi _{1},U\mu_{1}]=\frac{1}{4}U\alpha _{1},$

\ \ \ \ $[U\xi _{1},U\mu_{2}]=\frac{1}{4}U\alpha _{2},$

\ \ \ \ $[U\xi _{1},U\mu_{3}]=\frac{1}{4}U\alpha _{3},$

\ \ \ \ $[U\xi _{1},Uw_{1i}]=0,(0\leq i\leq 7),$

\ \ \ \ $[U\xi _{1},Uw_{2i}]=0,(0\leq i\leq 7),$

\ \ \ \ $[U\xi _{1},Uw_{3i}]=0,(0\leq i\leq 7),$

\ \ \ \ $[U\xi _{1},U\psi_{1}]=0,$

\ \ \ \ $[U\xi _{1},U\psi_{2}]=0,$

\ \ \ \ $[U\xi _{1},U\psi_{3}]=0,$

\ \ \ \ $[U\eta _{1},Uz_{1i}]=0,(0\leq i\leq 7),$

\ \ \ \ $[U\eta _{1},Uz_{2i}]=0,(0\leq i\leq 7),$

\ \ \ \ $[U\eta _{1},Uz_{3i}]=0,(0\leq i\leq 7),$

\ \ \ \ $[U\eta _{1},U\mu_{1}]=0,$

\ \ \ \ $[U\eta _{1},U\mu_{2}]=0,$

\ \ \ \ $[U\eta _{1},U\mu_{3}]=0,$

\ \ \ \ $[U\eta _{1},Uw_{1i}]=-\frac{1}{4}Ub_{1i},(0\leq i\leq 7),$

\ \ \ \ $[U\eta _{1},Uw_{2i}]=-\frac{1}{4}Ub_{2i},(0\leq i\leq 7),$

\ \ \ \ $[U\eta _{1},Uw_{3i}]=-\frac{1}{4}Ub_{3i},(0\leq i\leq 7),$

\ \ \ \ $[U\eta _{1},U\psi_{1}]=-\frac{1}{4}U\beta _{1},$

\ \ \ \ $[U\eta _{1},U\psi_{2}]=-\frac{1}{4}U\beta _{2},$

\ \ \ \ $[U\eta _{1},U\psi_{3}]=-\frac{1}{4}U\beta _{3}.$

\bigskip

\emph{Proof. \ }We have the above Lie bracket operations with
calculations using Maxima.\ \ \ \ \emph{Q.E.D.}

\bigskip

\emph{Lemma 17.29. }\ We have the following Lie bracket operations.

\ \ \ \ $[U\xi _{1},U\zeta_{1}]=0,$

\ \ \ \ $[U\xi _{1},U\omega_{1}]=-\frac{3}{8}U\rho _{1}-\frac{1}{8}Ur_{1}, $

\ \ \ \ $[U\eta _{1},U\zeta_{1}]=-\frac{3}{8}U\rho _{1}+\frac{1}{8}Ur_{1}, $

\ \ \ \ $[U\eta _{1},U\omega_{1}]=0.$

\bigskip

\emph{Proof. \ }We have the above Lie bracket operations with
calculations using Maxima.\ \ \ \ \emph{Q.E.D.}

\bigskip

\emph{Lemma 17.30.} \ We have the following Lie bracket operations.

\ \ \ \ $[Uz_{ij},Uz_{kl}]=0,(1\leq i,k\leq 3,0\leq j,l\leq 7),$

\ \ \ \ $[U\mu_{i},Uz_{kl}]=0,(1\leq i,k\leq 3,0\leq l\leq 7),$

\ \ \ \ $[U\mu_{i},U\mu_{k}]=0,(1\leq i,k\leq 3),$

\ \ \ \ $[Uz_{1i},Uw_{1i}]=-\frac{1}{2}Uu_{1},(0\leq i\leq 7),$

\ \ \ \ $[Uz_{1i},Uw_{1j}]=0,(0\leq i\neq j\leq 7),$

\ \ \ \ $[Uz_{1i},Uw_{2j}]=0,(0\leq i,j\leq 7),$

\ \ \ \ $[Uz_{1i},Uw_{3j}]=0,(0\leq i,j\leq 7),$

\ \ \ \ $[Uz_{2i},Uw_{1j}]=0,(0\leq i,j\leq 7),$

\ \ \ \ $[Uz_{2i},Uw_{2i}]=-\frac{1}{2}Uu_{1},(0\leq i\leq 7),$

\ \ \ \ $[Uz_{2i},Uw_{2j}]=0,(0\leq i\neq j\leq 7),$

\ \ \ \ $[Uz_{2i},Uw_{3j}]=0,(0\leq i,j\leq 7),$

\ \ \ \ $[Uz_{3i},Uw_{1j}]=0,(0\leq i,j\leq 7),$

\ \ \ \ $[Uz_{3i},Uw_{2j}]=0,(0\leq i,j\leq 7),$

\ \ \ \ $[Uz_{3i},Uw_{3i}]=-\frac{1}{2}Uu_{1},(0\leq i\leq 7),$

\ \ \ \ $[Uz_{3i},Uw_{3j}]=0,(0\leq i\neq j\leq 7),$

\ \ \ \ $[Uz_{ij},U\psi_{k}]=0,(1\leq i,k\leq 3,0\leq j\leq 7),$

\ \ \ \ $[U\mu_{k},Uw_{ij}]=0,(1\leq i,k\leq 3,0\leq j\leq 7),$

\ \ \ \ $[U\mu_{i},U\psi_{i}]=-\frac{1}{4}Uu_{1},(1\leq i\leq 3),$

\ \ \ \ $[U\mu_{i},U\psi_{k}]=0,(1\leq i,k\leq 3),$

\ \ \ \ $[Uz_{ij},U\zeta_{1}]=0,(1\leq i\leq 3,0\leq j\leq 7),$

\ \ \ \ $[U\mu_{i},U\zeta_{1}]=0,(1\leq i\leq 3),$

\ \ \ \ $[Uz_{ij},U\omega_{1}]=0,(1\leq i\leq 3,0\leq j\leq 7),$

\ \ \ \ $[U\mu_{i},U\omega_{1}]=0,(1\leq i\leq 3),$

\ \ \ \ $[Uw_{ij},Uw_{kl}]=0,(1\leq i,k\leq 3,0\leq j,l\leq 7),$

\ \ \ \ $[U\psi_{i},Uw_{kl}]=0,(1\leq i,k\leq 3,0\leq l\leq 7),$

\ \ \ \ $[U\psi_{i},U\psi_{k}]=0,(1\leq i,k\leq 3),$

\ \ \ \ $[Uw_{ij},U\zeta_{1}]=0,(1\leq i\leq 3,0\leq j\leq 7),$

\ \ \ \ $[U\psi_{i},U\zeta_{1}]=0,(1\leq i\leq 3),$

\ \ \ \ $[Uw_{ij},U\omega_{1}]=0,(1\leq i\leq 3,0\leq j\leq 7),$

\ \ \ \ $[U\psi_{i},U\omega_{1}]=0,(1\leq i\leq 3),$

\ \ \ \ $[U\zeta_{1},U\omega_{1}]=-\frac{1}{4}Uu_{1},(1\leq i\leq 3),$

\bigskip

\emph{Proof. \ }We have the above Lie bracket operations with
calculations using Maxima.\ \ \ \ \emph{Q.E.D.}

\bigskip

\emph{Lemma 17.31.} \ We have the following Lie bracket operations.

\ \ \ \ $[Ux_{ij},Ur_{1}]=-Ux_{ij},(1\leq i\leq 3,0\leq j\leq 7),$

\ \ \ \ $[U\chi _{i},Ur_{1}]=-U\chi _{i},(1\leq i\leq 3),$

\ \ \ \ $[Uy_{ij},Ur_{1}]=-Uy_{ij},(1\leq i\leq 3,0\leq j\leq 7),$

\ \ \ \ $[U\gamma _{i},Ur_{1}]=-U\gamma _{i},(1\leq i\leq 3),$

\ \ \ \ $[U\xi _{1},Ur_{1}]=-U\xi _{1},$

\ \ \ \ $[U\eta _{1},Ur_{1}]=-U\eta _{1},$

\ \ \ \ $[Ux_{ij},Us_{1}]=0,(1\leq i\leq 3,0\leq j\leq 7),$

\ \ \ \ $[U\chi _{i},Us_{1}]=0,(1\leq i\leq 3),$

\ \ \ \ $[Uy_{ij},Us_{1}]=0,(1\leq i\leq 3,0\leq j\leq 7),$

\ \ \ \ $[U\gamma _{i},Us_{1}]=0,(1\leq i\leq 3),$

\ \ \ \ $[U\xi _{1},Us_{1}]=0,$

\ \ \ \ $[U\eta _{1},Us_{1}]=0,$

\ \ \ \ $[Ux_{ij},Uu_{1}]=-Uz_{ij},(1\leq i\leq 3,0\leq j\leq 7),$

\ \ \ \ $[U\chi _{i},Uu_{1}]=-U\mu_{i},(1\leq i\leq 3),$

\ \ \ \ $[Uy_{ij},Uu_{1}]=-Uw_{ij},(1\leq i\leq 3,0\leq j\leq 7),$

\ \ \ \ $[U\gamma _{i},Uu_{1}]=-U\psi_{i},(1\leq i\leq 3),$

\ \ \ \ $[U\xi _{1},Uu_{1}]=-U\zeta_{1},$

\ \ \ \ $[U\eta _{1},Uu_{1}]=-U\omega_{1}.$

\bigskip

\emph{Proof. \ }We have the above Lie bracket operations with
calculations using Maxima.\ \ \ \ \emph{Q.E.D.}

\bigskip

\emph{Lemma 17.32.} \ We have the following Lie bracket operations.

\ \ \ \ $[Uz_{ij},Ur_{1}]=Uz_{ij},(1\leq i\leq 3,0\leq j\leq 7),$

\ \ \ \ $[U\mu_{i},Ur_{1}]=U\mu_{i},(1\leq i\leq 3),$

\ \ \ \ $[Uw_{ij},Ur_{1}]=Uw_{ij},(1\leq i\leq 3,0\leq j\leq 7),$

\ \ \ \ $[U\psi_{i},Ur_{1}]=U\psi_{i},(1\leq i\leq 3),$

\ \ \ \ $[U\zeta_{1},Ur_{1}]=U\zeta_{1},$

\ \ \ \ $[U\omega_{1},Ur_{1}]=U\omega_{1},$

\ \ \ \ $[Uz_{ij},Us_{1}]=-Ux_{ij},(1\leq i\leq 3,0\leq j\leq 7),$

\ \ \ \ $[U\mu_{i},Us_{1}]=-U\chi _{i},(1\leq i\leq 3),$

\ \ \ \ $[Uw_{ij},Us_{1}]=-Uy_{ij},(1\leq i\leq 3,0\leq j\leq 7),$

\ \ \ \ $[U\psi_{i},Us_{1}]=-U\gamma _{i},(1\leq i\leq 3),$

\ \ \ \ $[U\zeta_{1},Us_{1}]=-U\xi _{1},$

\ \ \ \ $[U\omega_{1},Us_{1}]=-U\eta _{1},$

\ \ \ \ $[Uz_{ij},Uu_{1}]=0,(1\leq i\leq 3,0\leq j\leq 7),$

\ \ \ \ $[U\mu_{i},Uu_{1}]=0,(1\leq i\leq 3),$

\ \ \ \ $[Uw_{ij},Uu_{1}]=0,(1\leq i\leq 3,0\leq j\leq 7),$

\ \ \ \ $[U\psi_{i},Uu_{1}]=0,(1\leq i\leq 3),$

\ \ \ \ $[U\zeta_{1},Uu_{1}]=0,$

\ \ \ \ $[U\omega_{1},Uu_{1}]=0.$

\bigskip

\emph{Proof. \ }We have the above Lie bracket operations with
calculations using Maxima.\ \ \ \ \emph{Q.E.D.}

\bigskip

\emph{Lemma 17.33.} \ We have the following Lie bracket operations.

\ \ \ \ $[Ur_{1},Us_{1}]=2Us_{1}.$

\ \ \ \ $[Ur_{1},Uu_{1}]=-2Uu_{1},$

\ \ \ \ $[Us_{1},Uu_{1}]=Ur_{1}.$

\bigskip

\emph{Proof. \ }We have the above Lie bracket operations with
calculations using Maxima.\ \ \ \ \emph{Q.E.D.}

\bigskip

\emph{Lemma 17.34.} \ \gr$_{8}^{\C}$ is simple.

\bigskip

\emph{Proof. \ }Case 1 for $Rx_{10}\in $\gR\gx$^{\C}.$

\noindent
By \emph{Lemma 17.2}, we have

\ \ \ \ $\{[Rd_{0j},Rx_{10}] \mid 1\leq j\leq 7\}=\{Rx_{1j}$ $|$ $1\leq j\leq 7\},$

\ \ \ \ $\{[Rm_{3j},Rx_{10}] \mid 0\leq j\leq 7\ =\{Rx_{2j}$ $|$ $0\leq j\leq 7\},$

\ \ \ \ $\{[Rm_{2j},Rx_{10}] \mid 0\leq j\leq 7\}=\{Rx_{1j}$ $|$ $0\leq j\leq 7\},$

\ \ \ \ $\{[Ut_{2j}-Um_{2j},Rx_{2j}] \mid 1\leq j\leq 7\}=\{R\chi _{1}\},$

\ \ \ \ $\{[Ut_{2j}+Um_{2j},Rx_{3j}] \mid 1\leq j\leq 7\}=\{R\chi _{3}\},$

\ \ \ \ $\{[Ut_{3j}-Um_{3j},Rx_{1j}] \mid 1\leq j\leq 7\}=\{R\chi _{2}\}.$

\noindent
Then we have

$(1)$\ \ \ \ $\{[x,Rx_{10}],[y,[x,Rx_{10}] ]\mid x,y\in
$\gr$_{4}^C\}\supset $\gR\gx$^{\C}.$

\noindent
By \emph{Lemma 17.3}, we have

\ \ \ \ \ $\{[Ra_{2i},Rx_{3j}] \mid 0\leq i,j\leq 7\}=\{Ry_{1k}$ $|$ $1\leq k\leq 7\},$

\ \ \ \ \ $\{[Ra_{3i},Rx_{1j}] \mid 0\leq i,j\leq 7\}=\{Ry_{2k}$ $|$ $1\leq k\leq 7\},$

\ \ \ \ \ $\{[Ra_{1i},Rx_{2j}] \mid 0\leq i,j\leq 7\}=\{Ry_{3k}$ $|$ $1\leq k\leq 7\},$

\ \ \ \ \ $\{[Ra_{1i},Rx_{1i}] \mid 0\leq i\leq 7\}=\{R\gamma _{1}\},$

\ \ \ \ \ $\{[Ra_{2i},Rx_{2i}] \mid 0\leq i\leq 7\}=\{R\gamma _{2}\},$

\ \ \ \ \ $\{[Ra_{3i},Rx_{3i}] \mid 0\leq i\leq 7\}=\{R\gamma _{3}\}.$

\noindent
Then we have

$(2)$\ \ \ \ \ $\{[a,x] \mid a\in $\gR\ga$^{\C},x\in $\gR\gx$^{\C}\}\supset $\gR\gy$^{\C}.$

\noindent
By \emph{Lemma 17.22}, we have

\ \ \ \ \ $\{[Rx_{2i},Rz_{3j}] \mid 0\leq i,j\leq 7\}=\{Rb_{1k}$ $|$ $%
1\leq k\leq 7\},$

\ \ \ \ \ $\{[Rx_{3i},Rz_{1j}] \mid 0\leq i,j\leq 7\}=\{Rb_{2k}$ $|$ $1\leq k\leq 7\},$

\ \ \ \ \ $\{[Rx_{1i},Rz_{2j}] \mid 0\leq i,j\leq 7\}=\{Rb_{3k}$ $|$ $1\leq k\leq 7\},$

\ \ \ \ \ $\{[Rx_{1i},Rz_{1i}] \mid 0\leq i\leq 7\}=\{R\beta _{1}\},$

\ \ \ \ \ $\{[Rx_{2i},Rz_{2i}] \mid 0\leq i\leq 7\}=\{R\beta _{2}\},$

\ \ \ \ \ $\{[Rx_{3i},Rz_{3i}] \mid 0\leq i\leq 7\}=\{R\beta _{3}\}.$

\noindent
Then we have

$(3)$\ \ \ \ \ $\{[x,z] \mid x\in $\gR\gx$^{\C},z\in $\gR\gz$^{\C}\}\supset $\gR\gb$^{\C}.$

\noindent
By \emph{Lemma 17.26}, we have

\ \ \ \ \ $\{[Ry_{2i},Rw_{3j}] \mid 0\leq i,j\leq 7\}=\{Ra_{1k}$ $|$ $1\leq k\leq 7\},$

\ \ \ \ \ $\{[Ry_{3i},Rw_{1j}] \mid 0\leq i,j\leq 7\}=\{Ra_{2k}$ $|$ $1\leq k\leq 7\},$

\ \ \ \ \ $\{[Ry_{1i},Rw_{2j}] \mid 0\leq i,j\leq 7\}=\{Ra_{3k}$ $|$ $%
1\leq k\leq 7\},$

\ \ \ \ \ $\{[Ry_{1i},Rw_{1i}] \mid 0\leq i\leq 7\}=\{R\alpha _{1}\},$

\ \ \ \ \ $\{[Ry_{2i},Rw_{2i}] \mid 0\leq i\leq 7\}=\{R\alpha _{2}\},$

\ \ \ \ \ $\{[Ry_{3i},Rw_{3i}] \mid 0\leq i\leq 7\}=\{R\alpha _{3}\}.$

\noindent
Then we have

$(4)$\ \ \ \ \ $\{[y,w] \mid y\in $\gR\gy$^{\C},w\in $\gR\gw$^{\C}\}\supset $\gR\ga$^{\C}.$

\noindent
By \emph{Lemma 17.23}, we have

\ \ \ \ $\{[Rx_{1i},Rw_{1j}] \mid 0\leq i\neq j\leq 7\}=\{Rd_{ij}$ $|$
$0\leq i\neq j\leq 7\},$

\ \ \ $\ \{[Rx_{1i},U\psi_{2}-U\psi_{3}] \mid 0\leq i\leq
7\}=\{Rm_{1i}  \mid 0\leq i\leq 7\},$

\ \ \ $\ \{[Rx_{2i},U\psi_{3}-U\psi_{1}] \mid 0\leq i\leq
7\}=\{Rm_{2i}  \mid 0\leq i\leq 7\},$

\ \ \ $\ \{[Rx_{3i},U\psi_{1}-U\psi_{2}] \mid 0\leq i\leq
7\}=\{Rm_{3i}  \mid 0\leq i\leq 7\},$

\ \ \ $\ \{[Rx_{1i},U\psi_{2}+U\psi_{3}] \mid 0\leq i\leq
7\}=\{Rt_{1i}  \mid 0\leq i\leq 7\},$

\ \ \ $\ \{[Rx_{2i},U\psi_{3}+U\psi_{1}] \mid 0\leq i\leq
7\}=\{Rt_{2i}  \mid 0\leq i\leq 7\},$

\ \ \ $\ \{[Rx_{3i},U\psi_{1}+U\psi_{2}] \mid 0\leq i\leq
7\}=\{Rt_{3i}  \mid 0\leq i\leq 7\},$

\ \ \ \ $[Ux_{1i},Uw_{1i}]-[Ux_{3i},Uw_{3i}]=\frac{1}{2}U\tau _{1},$

\ \ \ \ $[Ux_{2i},Uw_{2i}]-[Ux_{3i},Uw_{3i}]=\frac{1}{2}U\tau _{2}.$

\noindent
Then we have

$(5)$\ \ \ \ \ $\{[x,w] \mid y\in $\gR\gx$^{\C},w\in $\gR\gw$^{\C}\}\supset
\ $\gr\gd$^{\C} \bigoplus $\gR\gm$^{\C} \bigoplus $\gR\gt$^{\C}\ =$\gr$_{6}^{\C}.$

\noindent
By \emph{Lemma 17.12}, we have

\ \ \ \ \ $\{[Ra_{2i},Rz_{3j}] \mid 0\leq i,j\leq 7\}=\{Rw_{1k}$ $|$ $%
1\leq k\leq 7\},$

\ \ \ \ \ $\{[Ra_{3i},Rz_{1j}] \mid 0\leq i,j\leq 7\}=\{Rw_{2k}$ $|$ $%
1\leq k\leq 7\},$

\ \ \ \ \ $\{[Ra_{1i},Rz_{2j}] \mid 0\leq i,j\leq 7\}=\{Rw_{3k}$ $|$ $%
1\leq k\leq 7\},$

\ \ \ \ \ $\{[Ra_{1i},Rz_{1i}] \mid 0\leq i\leq 7\}=\{R\psi _{1}\},$

\ \ \ \ \ $\{[Ra_{2i},Rz_{2i}] \mid 0\leq i\leq 7\}=\{R\psi _{2}\},$

\ \ \ \ \ $\{[Ra_{3i},Rz_{3i}] \mid 0\leq i\leq 7\}=\{R\psi _{3}\}.$

\noindent
Then we have

$(6)$\ \ \ \ \ $\{[a,z] \mid a\in $\gR\ga$^{\C},z\in $\gR\gz$^{\C}\}\supset $\gR\gw$^{\C}.$

\noindent
By \emph{Lemma 17.16}, we have

\ \ \ \ \ $\{[Rb_{2i},Rw_{3j}] \mid 0\leq i,j\leq 7\}=\{Rz_{1k}$ $|$ $%
1\leq k\leq 7\},$

\ \ \ \ \ $\{[Rb_{3i},Rw_{1j}] \mid 0\leq i,j\leq 7\}=\{Rz_{2k}$ $|$ $%
1\leq k\leq 7\},$

\ \ \ \ \ $\{[Rb_{1i},Rw_{2j}] \mid 0\leq i,j\leq 7\}=\{Rz_{3k}$ $|$ $%
1\leq k\leq 7\},$

\ \ \ \ \ $\{[Rb_{1i},Rw_{1i}] \mid 0\leq i\leq 7\}=\{R\mu _{1}\},$

\ \ \ \ \ $\{[Rb_{2i},Rw_{2i}] \mid 0\leq i\leq 7\}=\{R\mu _{2}\},$

\ \ \ \ \ $\{[Rb_{3i},Rw_{3i}] \mid 0\leq i\leq 7\}=\{R\mu _{3}\}.$

\noindent
Then we have

$(7)$\ \ \ \ \ $\{[b,w] \mid b\in $\gR\gb$^{\C},w\in $\gR\gw$^{\C}\}\supset $\gR\gz$^{\C}.$

\noindent
By \emph{Lemma 17.4}, we have

\ \ \ \ \ $\{[Rb_{1i},Rx_{1i}] \mid 0\leq i\leq 7\}=\{R\eta _{1}\},$

\noindent
so we have

$(8)$ \ \ \ $\{[b,x] \mid b\in $\gR\gb$^{\C},x\in $\gR\gx$^{\C}\}\supset $\gR\gi$^{\C} .$

\noindent
By \emph{Lemma 17.8}, we have

\ \ \ \ \ $\{[Ra_{1i},Ry_{1i}] \mid 0\leq i\leq 7\}=\{R\xi _{1}\},$

\noindent
so we have

$(9)$ \ \ \ $\{[a,y] \mid a\in $\gR\ga$^{\C},y\in $\gR\gy$^{\C}\}\supset $\gR\gk$^{\C} .$

By \emph{Lemma 17.13}, we have

\ \ \ \ \ $\{[Rb_{1i},Rz_{1i}] \mid 0\leq i\leq 7\}=\{R\omega _{1}\},$

\noindent
so we have

$(10)$ \ \ \ $\{[b,z] \mid b\in $\gR\gb$^{\C},z\in $\gR\gz$^{\C}\}\supset $\gR\go$^{\C}. $

\noindent
By \emph{Lemma 17.17}, we have

\ \ \ \ \ $\{[Ra_{1i},Rw_{1i}] \mid 0\leq i\leq 7\}=\{R\zeta _{1}\},$

\noindent
so we have

$(11)$ \ \ \ $\{[a,w] \mid a\in $\gR\ga$^{\C},w\in $\gR\gw$^{\C}\}\supset $\gR\gl$^{\C} .$

\noindent
By \emph{Lemma 17.29}, we have

\ \ \ \ $[U\eta _{1},U\zeta_{1}]+[U\xi _{1},U\omega_{1}]=-\frac{6}{8}U\rho _{1},$

\ \ \ \ $[U\eta _{1},U\zeta_{1}]-[U\xi _{1},U\omega_{1}]=\frac{2}{8}Ur_{1}.$

\noindent
Then we have

$(12)$\ \ \ \ \ $\{[x,y] \mid x\in $\gR\gk$^{\C} \bigoplus $\gR\gi$^{\C}, y\in
$\gR\gl$^{\C}  \bigoplus $\gR\go$^{\C} \ $ $\}\supset \ $\gR\gp$^{\C}  \bigoplus $\gR\gr$^{\C}\ .$

\noindent
By \emph{Lemma 17.12}, we have

\ \ \ \ \ $\{[Rx_{1i},Ry_{1i}] \mid 0\leq i\leq 7\}=\{Rs_{1}\},$

\noindent
so we have

$(13)$ \ \ \ $\{[x,y] \mid x\in $\gR\gx$^{\C},y\in $\gR\gy$^{\C}\}\supset $\gR\gs$^{\C}.$

\noindent
By \emph{Lemma 17.30}, we have

\ \ \ \ \ $\{[Rz_{1i},Rw_{1i}] \mid 0\leq i\leq 7\}=\{Ru_{1}\},$

\noindent
so we have

$(14)$ \ \ \ $\{[z,w] \mid z\in $\gR\gz$^{\C},w\in $\gR\gw$^{\C}\}\supset $\gR\gu$^{\C}.$

\noindent
By $(1),(2),(3),(4),(5),(6),(7),(8),(9),(10),(11),(12),(13)$ and $(14)$,we
have

\ \ \ $\{[Rx_{10},x],[[Rx_{10},x],y] \mid x,y\in $\gr$%
_{8}^{\C}\}=$\gr$_{8}^{\C}.$

Case 2 for $Rx_{1j}\in $\gR\gx$^{\C}.$

\noindent
By \emph{Lemma 17.2}, we have

\ \ \ $\{[Rd_{0j},Rx_{1j}] \mid 0\leq j\leq 7\}=\{Rx_{10}\}.$

\noindent
Therefore,this applies to Case 1.

Case 3 for $Rx_{2j}\in $\gR\gx$^{\C},(0\leq i\leq 7).$

\noindent
By \emph{Lemma 17.2}, we have

\ \ \ $\{[Rm_{3i},Rx_{2j}] \mid 0\leq i\neq j\leq 7\}=\{Rx_{1k}$ $|$ $%
0\leq k\leq 7\}.$

\noindent
Therefore,this applies to Case 2.

Case 4 for $Rx_{3j}\in $\gR\gx$^{\C},(0\leq i\leq 7).$

\noindent
By \emph{Lemma 17.2}, we have

\ \ \ $\{[Rm_{2i},Rx_{3j}] \mid 0\leq i\neq j\leq 7\}=\{Rx_{1k}$ $|$ $%
0\leq k\leq 7\}.$

\noindent
Therefore,this applies to Case 2.

Case 5 for $R\chi _{1},R\chi _{2},R\chi _{3}\in $\gR\gx$^{\C}.$

\noindent
By \emph{Lemma 17.2}, we have

\ \ $\{[Rm_{2i},R\chi _{1}] \mid 0\leq i\leq 7\}=\{Rx_{2k}$ $|$ $%
0\leq k\leq 7\},$

\ \ $\{[Rm_{1i},R\chi _{2}] \mid 0\leq i\leq 7\}=\{Rx_{1k}$ $|$ $%
0\leq k\leq 7\},$

\ \ $\{[Rm_{1i},R\chi _{3}] \mid 0\leq i\leq 7\}=\{Rx_{1k}$ $|$ $%
0\leq k\leq 7\}.$

\noindent
Therefore,this applies to Case 1 and Case 2.

Case 6 for $y\in $\gR\gy$^{\C},z\in $\gR\gz$^{\C},w\in $\gR\gw$^{\C}.$

\noindent
Same as Case 1,Case 2,Case 3,Case 4 and Case 5.

Case 7 for $\xi \in $\gR\gk$^{\C} ,\eta \in $\gR\gi$^{\C} ,\zeta \in $\gR\gl$^{\C} ,\omega \in $\gR\go$^{\C} .$

\noindent
By \emph{Lemma 17.24}, we have

\ \ \ $\{[x,\zeta ] \mid x\in $\gR\gx$^{\C}\}=$\gR\ga$^{\C}.$

\noindent
From $(2),$ this applies to Case 6.

\noindent
By \emph{Lemma 17.27}, we have

\ \ \ $\{[y,\omega ] \mid y\in $\gR\gy$^{\C}\}=$\gR\gb$^{\C}.$

\noindent
From $(7),$ this applies to Case 6.

\noindent
By \emph{Lemma 17.28}, we have

\ \ \ $\{[z,\xi ] \mid z\in $\gR\gz$^{\C}\}=$\gR\ga$^{\C}.$

\ \ \ $\{[w,\eta ] \mid w\in $\gR\gw$^{\C}\}=$\gR\gb$^{\C}.$

\noindent
From $(2)$ and $(7),$ these applies to Case 6.

Case 8 for $r\in $\gR\gr$^{\C},s\in $\gR\gs$^{\C},u\in $\gR\gu$^{\C}.$

\noindent
By \emph{Lemma 17.31}, we have

\ \ \ $\{[Rx_{1i},Rr_{1}] \mid 0\leq i\leq 7\}=\{Rx_{1i}  \mid 0\leq i\leq 7\},$

\ \ \ $\{[Rx_{1i},Ru_{1}] \mid 0\leq i\leq 7\}=\{Rz_{1i}  \mid 0\leq i\leq 7\}.$

\noindent
By \emph{Lemma 17.32}, we have

\ \ \ $\{[Rz_{1i},Rs_{1}] \mid 0\leq i\leq 7\}=\{Rx_{1i} \mid 0\leq i\leq 7\}.$

\noindent
Therefore,these applies to Case 1 and Case 2.

By Case 1,Case 2,Case 3,Case 4,Case 5,Case 6, and Case 7, \gr$%
_{8}^{\C}$ is simple. \ \ \ \ \emph{Q.E.D.}

\bigskip

\emph{Lemma 17.35.} The Killing form $B_{8}$ of the Lie algebra %
\gr$_{8}^{\C}$ is given by

\ \ \ \ \ \ \ \ \ \ \ \ \ \ \ \ \ $B_{8}(R_{1},R_{2})=tr(R_{1}R_{2}),R_{1},R_{2}\in $\gr%
$_{8}^{\C}.$ \ 

\bigskip

\emph{Proof}. \ Since \gr$_{8}^{\C}$ is simple ,there exist $%
\kappa \in \C$ such that

\ \ \ \ \ \ \ \ \ \ \ \ \ \ \ \ \ $B_{8}(R_{1},R_{2})=\kappa tr(R_{1}R_{2}).$

\noindent
To determine this $\kappa ,$let $R_{1}=Ur_{1}.$ $(adR_{1})^{2}$ is calculated as
follows.

\noindent
For $R_{2}=(R_{72}+X_{2}+Y_{2}+\xi _{2}+\eta _{2}+Z_{2}+W_{2}+\zeta _{2}+\omega _{2}+r_{2}$
$+s_{2}+u_{2})\in $\gr$_{8}^{\C},R_{72}\in $\gr$_{7}^{\C},$

\noindent
$X_{2}\in $\gR\gx$^{\C},Y_{2}\in $\gR\gy$^{\C},\xi _{2}\in $\gR\gk$^{\C} ,\eta _{2}\in $\gR\gi$^{\C} ,$
$Z_{2}\in $\gR\gz$^{\C},W_{2}\in $\gR\gw$^{\C},\zeta _{2}\in $\gR\gl$^{\C} ,\omega _{2}\in $\gR\go$^{\C}$,

\noindent
$r_{2}\in $\gR\gr$^{\C},s_{2}\in $\gR\gs$^{\C},u_{2}\in $\gR\gu$^{\C},$ 
we have by \emph{Lemma 17.20}, \emph{Lemma 17.31}, \emph{Lemma 17.32} and 
\emph{Lemma 17.33},

$[R_{1},[R_{1},R_{2}]]=[R_{1},X_{2}+Y_{2}+\xi _{2}+\eta _{2}-Z_{2}-W_{2}-\zeta _{2}-\omega _{2}+s_{2}-u_{2}],$

\ \ \ \ \ \ \ \ \ \ \ \ \ \ \ \ \ \ \ \ $=X_{2}+Y_{2}+\xi _{2}+\eta _{2}+Z_{2}+W_{2}+\zeta
_{2}+\omega _{2}+4s_{2}+4u_{2}.$

\noindent
Hence

\ \ \ \ $B_{8}(R_{1},R_{1})=tr((adR_{1})^{2})=(27+27+1+1+27+27+1+1+4+4)=120.$

\noindent
On the other hand

\  \ \ \ $tr(Ur_{1}.Ur_{1})=120.$

\noindent
Therefore $\kappa =1.$ \ \ \ \ \emph{Q.E.D.}

\bigskip

\emph{Lemma 17.36.} The rank of the Lie algebra \gr$_{8}^{\C}$ is
8. The roots of \gr$_{8}^{\C}$ relative to some Cartan
subalgebra are given by

$\ \ \ \ \ \ \ \ \ \ \ \ \ \ \ \ \pm (\lambda _{k}-\lambda _{l}),\pm
(\lambda _{k}+\lambda _{l}),0\leq k<l\leq 3,$

\ \ \ \ \ \ \ \ \ \ $\ \ \ \ \ \ \pm \lambda _{k}\pm \frac{1}{2}(\mu
_{1}+2\mu _{2}),0\leq k\leq 3,$

$\ \ \ \ \ \ \ \ \ \ \ \ \ \ \ \ \pm \frac{1}{2}(-\lambda _{0}-\lambda
_{1}+\lambda _{2}-\lambda _{3})\pm \frac{1}{2}(-2\mu _{1}-\mu _{2}),$

$\ \ \ \ \ \ \ \ \ \ \ \ \ \ \ \ \pm \frac{1}{2}(\lambda _{0}+\lambda
_{1}+\lambda _{2}-\lambda _{3})\pm \frac{1}{2}(-2\mu _{1}-\mu _{2}),$

\ \ \ \ \ \ \ \ \ \ \ \ \ $\ \ \ \pm \frac{1}{2}(-\lambda _{0}+\lambda
_{1}+\lambda _{2}+\lambda _{3})\pm \frac{1}{2}(-2\mu _{1}-\mu _{2}),$

$\ \ \ \ \ \ \ \ \ \ \ \ \ \ \ \ \pm \frac{1}{2}(\lambda _{0}-\lambda
_{1}+\lambda _{2}+\lambda _{3})\pm \frac{1}{2}(-2\mu _{1}-\mu _{2}),$

\ \ \ \ \ \ \ \ \ \ \ \ \ \ $\ \ \pm \frac{1}{2}(\lambda _{0}-\lambda
_{1}+\lambda _{2}-\lambda _{3})\pm \frac{1}{2}(\mu _{1}-\mu _{2}),$

$\ \ \ \ \ \ \ \ \ \ \ \ \ \ \ \ \pm \frac{1}{2}(-\lambda _{0}+\lambda
_{1}+\lambda _{2}-\lambda _{3})\pm \frac{1}{2}(\mu _{1}-\mu _{2}),$

\ \ \ \ \ \ \ \ \ \ \ \ \ \ \ $\ \pm \frac{1}{2}(\lambda _{0}+\lambda
_{1}+\lambda _{2}+\lambda _{3})\pm \frac{1}{2}(\mu _{1}-\mu _{2}),$

$\ \ \ \ \ \ \ \ \ \ \ \ \ \ \ \ \pm \frac{1}{2}(-\lambda _{0}-\lambda
_{1}+\lambda _{2}+\lambda _{3})\pm \frac{1}{2}(\mu _{1}-\mu _{2}),$

\ \ \ \ \ \ \ \ \ \ \ \ \ $\ \ \ \pm (\mu _{k}+\frac{2}{3}\nu ),0\leq k\leq
2,\pm (-\mu _{1}-\mu _{2}+\frac{2}{3}\nu ),$

\ \ \ \ \ \ \ \ \ \ \ \ \ \ \ $\ \pm \lambda _{k}\pm (\frac{1}{2}\mu _{1}-%
\frac{2}{3}\nu ),0\leq k\leq 3,$

\ \ \ \ \ \ \ \ \ \ \ \ \ \ \ $\ \pm \frac{1}{2}(-\lambda _{0}-\lambda
_{1}+\lambda _{2}-\lambda _{3})\pm (\frac{1}{2}\mu _{2}-\frac{2}{3}\nu ),$

$\ \ \ \ \ \ \ \ \ \ \ \ \ \ \ \ \pm \frac{1}{2}(\lambda _{0}+\lambda
_{1}+\lambda _{2}-\lambda _{3})\pm (\frac{1}{2}\mu _{2}-\frac{2}{3}\nu ),$

\ \ \ \ \ \ \ \ \ \ \ \ \ $\ \ \ \pm \frac{1}{2}(-\lambda _{0}+\lambda
_{1}+\lambda _{2}+\lambda _{3})\pm (\frac{1}{2}\mu _{2}-\frac{2}{3}\nu ),$

$\ \ \ \ \ \ \ \ \ \ \ \ \ \ \ \ \pm \frac{1}{2}(\lambda _{0}-\lambda
_{1}+\lambda _{2}+\lambda _{3})\pm (\frac{1}{2}\mu _{2}-\frac{2}{3}\nu ),$

\ \ \ \ \ \ \ \ \ \ \ \ \ \ $\ \ \pm \frac{1}{2}(-\lambda _{0}+\lambda
_{1}-\lambda _{2}+\lambda _{3})\pm (-\frac{1}{2}\mu _{1}-\frac{1}{2}\mu _{2}-%
\frac{2}{3}\nu ),$

$\ \ \ \ \ \ \ \ \ \ \ \ \ \ \ \ \pm \frac{1}{2}(-\lambda _{0}+\lambda
_{1}+\lambda _{2}-\lambda _{3})\pm (-\frac{1}{2}\mu _{1}-\frac{1}{2}\mu _{2}-%
\frac{2}{3}\nu ),$

\ \ \ \ \ \ \ \ \ \ \ \ \ \ \ $\ \pm \frac{1}{2}(\lambda _{0}+\lambda
_{1}+\lambda _{2}+\lambda _{3})\pm (-\frac{1}{2}\mu _{1}-\frac{1}{2}\mu _{2}-%
\frac{2}{3}\nu ),$

$\ \ \ \ \ \ \ \ \ \ \ \ \ \ \ \ \pm \frac{1}{2}(-\lambda _{0}-\lambda
_{1}+\lambda _{2}+\lambda _{3})\pm (-\frac{1}{2}\mu _{1}-\frac{1}{2}\mu _{2}-%
\frac{2}{3}\nu ),$

\ \ \ \ \ \ \ \ \ \ \ \ \ \ \ \ $\pm (\mu _{j}-\frac{1}{3}v)\pm r,0\leq
j\leq 2,\pm (-\mu _{1}-\mu _{2}-\frac{1}{3}\nu )\pm r,$

\ \ \ \ \ \ \ \ \ \ \ \ \ \ \ $\ \pm \lambda _{k}\pm (\frac{1}{2}\mu _{1}+%
\frac{1}{3}\nu )\pm r,0\leq k\leq 3,$

\ \ \ \ \ \ \ \ \ \ \ \ \ \ \ $\ \pm \frac{1}{2}(-\lambda _{0}-\lambda
_{1}+\lambda _{2}-\lambda _{3})\pm (\frac{1}{2}\mu _{2}+\frac{1}{3}\nu )\pm
r,$

$\ \ \ \ \ \ \ \ \ \ \ \ \ \ \ \ \pm \frac{1}{2}(\lambda _{0}+\lambda
_{1}+\lambda _{2}-\lambda _{3})\pm (\frac{1}{2}\mu _{2}+\frac{1}{3}\nu )\pm
r,$

\ \ \ \ \ \ \ \ \ \ \ \ \ $\ \ \ \pm \frac{1}{2}(-\lambda _{0}+\lambda
_{1}+\lambda _{2}+\lambda _{3})\pm (\frac{1}{2}\mu _{2}+\frac{1}{3}\nu )\pm
r,$

$\ \ \ \ \ \ \ \ \ \ \ \ \ \ \ \ \pm \frac{1}{2}(\lambda _{0}-\lambda
_{1}+\lambda _{2}+\lambda _{3})\pm (\frac{1}{2}\mu _{2}+\frac{1}{3}\nu )\pm
r,$

\ \ \ \ \ \ \ \ \ \ \ \ \ \ $\ \ \pm \frac{1}{2}(\lambda _{0}-\lambda
_{1}+\lambda _{2}-\lambda _{3})\pm (-\frac{1}{2}\mu _{1}-\frac{1}{2}\mu _{2}+%
\frac{1}{3}\nu )\pm r,$

\ \ \ \ \ \ \ \ \ \ \ \ \ \ $\ \ \pm \frac{1}{2}(\lambda _{0}-\lambda
_{1}-\lambda _{2}+\lambda _{3})\pm (-\frac{1}{2}\mu _{1}-\frac{1}{2}\mu _{2}+%
\frac{1}{3}\nu )\pm r,$

\ \ \ \ \ \ \ \ \ \ \ \ \ \ \ $\ \pm \frac{1}{2}(\lambda _{0}+\lambda
_{1}+\lambda _{2}+\lambda _{3})\pm (-\frac{1}{2}\mu _{1}-\frac{1}{2}\mu _{2}+%
\frac{1}{3}\nu )\pm r,$

$\ \ \ \ \ \ \ \ \ \ \ \ \ \ \ \ \pm \frac{1}{2}(\lambda _{0}+\lambda
_{1}-\lambda _{2}-\lambda _{3})\pm (-\frac{1}{2}\mu _{1}-\frac{1}{2}\mu _{2}+%
\frac{1}{3}\nu )\pm r,$

\ \ \ \ \ \ \ \ \ \ \ \ \ \ $\ \ \pm 2r,$

\ \ \ \ \ \ \ \ \ \ \ \ \ $\ \ \ \pm v\pm r.$

\bigskip

\emph{Proof. \ }Let

\gh$=\left\{ h=h_{\delta }+H+V+R\in \text{\gr}_{8}^{\C}%
\middle| %
\begin{array}{c}
h_{\delta }=\sum\limits_{k=0}^{3}\lambda_{k}H_{k}=\sum\limits_{k=0}^{3}-\lambda _{k}iUd_{k4+k},\\
H=\mu _{1}U\tau _{1}+\mu_{2}U\tau _{2}, \\ 
V=\nu U\rho _{1},R=rUr_{1},\lambda _{k},\mu _{j},\nu ,r\in \C%
\end{array}%
\right\} ,$

\noindent
then \gh\  is an abelian subalgebra of \gr$_{8}^{\C}$ $($it will
be a Cartan subalgebra of \gr$_{8}^{\C}$ $)$. 

\noindent
That \gh\  is abelian is clear from \emph{Lemma 17.20 .}

$(1)$ The roots of \gr$_{7}^{\C}$ are also roots of \gr$%
_{8}^{\C} $. Indeed,let $\alpha $ be a root of \gr$_{7}^{\C}$ and $S\in 
$\gr$_{7}^{\C}\subset $\gr$_{8}^{\C}$
be a root vector belong to $\alpha .$ Then

$[h,S]=[h_{\delta }+H+V+R,S]=[h_{\delta }+H+V,S]+[R,S]=\alpha (h)S,$

$($Since S$\in $\gr$_{7}^{\C},$so $[R,S]=0).$

\noindent
Hence $\alpha $ is a root of \gr$_{8}^{\C}.$

$(2)$ Let we put $Sx_{0k}:$

\ \ \ \ \ \ \ \ \ \ \ \ \ \ \ \ \ \ \ \ \ \ $Sx_{0k}=U\chi _{k},(1\leq k\leq
3). $

\noindent
Then we have by \emph{Lemma 17.2 ,Lemma 17.5 }and\emph{\ Lemma 17.30 ,}

$[h_{\delta }+H+V+R,Sx_{0k}]=\mu _{k}U\chi _{k}-\frac{1}{3}vU\chi _{k}+rU\chi
k $

$\ \ \ \ \ \ \ \ \ \ \ \ \ \ \ \ \ \ \ \ \ \ \ \ \ \ \ \ \ \ \ \ =((\mu _{k}-%
\frac{1}{3}v)+r)Sx_{0k},(1\leq k\leq 3,\mu _{3}=-\mu _{1}-\mu _{2}).$

\noindent
Let we put $Sy_{0k}:$

\ \ \ \ \ \ \ \ \ \ \ \ \ \ \ \ \ \ \ \ \ \ $Sy_{0k}=U\gamma _{k},(1\leq k\leq
3).$

\noindent
Then we have by \emph{Lemma 17.6 ,Lemma 17.9 }and\emph{\ Lemma 17.30 ,}

$[h_{\delta }+H+V+R,Sy_{0k}]=-\mu _{k}U\gamma _{k}+\frac{1}{3}vU\gamma
k+rU\gamma _{k}$

$\ \ \ \ \ \ \ \ \ \ \ \ \ \ \ \ \ \ \ \ \ \ \ \ \ \ \ \ \ \ \ \ =(-(\mu
_{k}-\frac{1}{3}v)+r)Sy_{0k},(1\leq k\leq 3,\mu _{3}=-\mu _{1}-\mu _{2}).$

\noindent
Let we put $Sz_{0k}:$

\ \ \ \ \ \ \ \ \ \ \ \ \ \ \ \ \ \ \ \ \ \ $Sz_{0k}=U\mu_{k},(1\leq k\leq 3).$

\noindent
Then we have by \emph{Lemma 17.11 ,Lemma 17.14 }and\emph{\ Lemma 17.31 ,}

$[h_{\delta }+H+V+R,Sz_{0k}]=\mu _{k}U\mu_{k}-\frac{1}{3}vU\mu_{k}-rU\mu_{k}$

$\ \ \ \ \ \ \ \ \ \ \ \ \ \ \ \ \ \ \ \ \ \ \ \ \ \ \ \ \ \ \ \ =((\mu _{k}-%
\frac{1}{3}v)-r)Sz_{0k},(1\leq k\leq 3,\mu _{3}=-\mu _{1}-\mu _{2}).$

\noindent
Let we put $Sw_{0k}:$

\ \ \ \ \ \ \ \ \ \ \ \ \ \ \ \ \ \ \ \ \ \ $Sw_{0k}=U\psi_{k},(1\leq k\leq
3). $

\noindent
Then we have by \emph{Lemma 17.15 ,Lemma 17.18 }and\emph{\ Lemma 17.31 ,}

\ $[h_{\delta }+H+V+R,Sw_{0k}]=-\mu _{k}U\psi_{k}+\frac{1}{3}vU\psi_{k}-rU\psi_{k}$

$\ \ \ \ \ \ \ \ \ \ \ \ \ \ \ \ \ \ \ \ \ \ \ \ \ \ \ \ \ \ \ \ =(-(\mu
_{k}-\frac{1}{3}v)-r)Sw_{0k},(1\leq k\leq 3,\mu _{3}=-\mu _{1}-\mu _{2}).$

\noindent
Henc $(\mu _{k}-\frac{1}{3}v)+r,-(\mu _{k}-\frac{1}{3}v)+r,(1\leq k\leq
3),(\mu _{k}-\frac{1}{3}v)-r,-(\mu _{k}-\frac{1}{3}v)-r,(1\leq k\leq 3)$
are roots \gr$_{8}^{\C}$ and $U\chi _{j},u\gamma j,U\mu_{j},U\psi_{j},(1\leq
j\leq 3)$ are root vectors respectively.

$(3)$ Let we put $Sx_{1k}$ and $Sx_{4k}:$

\ \ \ \ \ \ \ \ \ \ \ \ \ \ \ \ \ \ \ \ \ \ $Sx_{1k}=Ux_{1k}+iUx_{14+k},(0%
\leq k\leq 3),$

\ \ \ \ \ \ \ \ \ \ \ \ \ \ \ \ \ \ \ \ \ \ $Sx_{4k}=Ux_{1k}-iUx_{14+k},(0%
\leq k\leq 3).$

\noindent
Then we have by \emph{Lemma 17.2 ,Lemma 17.5 }and\emph{\ Lemma 17.30 ,}

$[h_{\delta }+H+V+R,Sx_{1k}]=(\lambda _{k}-(\frac{1}{2}\mu _{1}+\frac{1}{%
3}v)+r)Sx_{1k},(0\leq k\leq 3),$

$[h_{\delta }+H+V+R,Sx_{4k}]=(-\lambda _{k}-(\frac{1}{2}\mu _{1}+\frac{1%
}{3}v)+r)Sx_{4k},(0\leq k\leq 3).$

\noindent
Hence $\pm \lambda _{k}-(\frac{1}{2}\mu _{1}+\frac{1}{3}v)+r$ are roots of 
\gr$_{8}^{\C}$ .

\noindent
Let we put $Sy_{1k}$ and $Sy_{4k}:$

\ \ \ \ \ \ \ \ \ \ \ \ \ \ \ \ \ \ \ \ \ \ $Sy_{1k}=Uy_{1k}+iUy_{14+k},(0%
\leq k\leq 3),$

\ \ \ \ \ \ \ \ \ \ \ \ \ \ \ \ \ \ \ \ \ \ $Sy_{4k}=Uy_{1k}-iUy_{14+k},(0%
\leq k\leq 3).$

\noindent
Then we have by \emph{Lemma 17.6 ,Lemma 17.9 }and\emph{\ Lemma 17.30 ,}

$[h_{\delta }+H+V+R,Sy_{1k}]=(\lambda _{k}+(\frac{1}{2}\mu _{1}+\frac{1}{%
3}v)+r)Sy_{1k},(0\leq k\leq 3),$

$[h_{\delta }+H+V+R,Sy_{4k}]=(-\lambda _{k}+(\frac{1}{2}\mu _{1}+\frac{1%
}{3}v)+r)Sy_{4k},(0\leq k\leq 3).$

\noindent
Hence $\pm \lambda _{k}+(\frac{1}{2}\mu _{1}+\frac{1}{3}v)+r$ are roots of 
\gr$_{8}^{\C}$ .

\noindent
Let we put $Sz_{1k}$ and $Sz_{4k}:$

\ \ \ \ \ \ \ \ \ \ \ \ \ \ \ \ \ \ \ \ \ \ $Sz_{1k}=Uz_{1k}+iUz_{14+k},(0%
\leq k\leq 3),$

\ \ \ \ \ \ \ \ \ \ \ \ \ \ \ \ \ \ \ \ \ \ $Sz_{4k}=Uz_{1k}-iUz_{14+k},(0%
\leq k\leq 3).$

\noindent
Then we have by \emph{Lemma 17.11 ,Lemma 17.14 }and\emph{\ Lemma 17.31 ,}

$[h_{\delta }+H+V+R,Sz_{1k}]=(\lambda _{k}-(\frac{1}{2}\mu _{1}+\frac{1}{%
3}v)-r)Sz_{1k},(0\leq k\leq 3),$

$[h_{\delta }+H+V+R,Sz_{4k}]=(-\lambda _{k}-(\frac{1}{2}\mu _{1}+\frac{1%
}{3}v)-r)Sz_{4k},(0\leq k\leq 3).$

\noindent
Hence $\pm \lambda _{k}-(\frac{1}{2}\mu _{1}+\frac{1}{3}v)-r$ are roots of 
\gr$_{8}^{\C}$ .

\noindent
Let we put $Sw_{1k}$ and $Sw_{4k}:$

\ \ \ \ \ \ \ \ \ \ \ \ \ \ \ \ \ \ \ \ \ \ $Sw_{1k}=Uw_{1k}+iUw_{14+k},(0%
\leq k\leq 3),$

\ \ \ \ \ \ \ \ \ \ \ \ \ \ \ \ \ \ \ \ \ \ $Sw_{4k}=Uw_{1k}-iUw_{14+k},(0%
\leq k\leq 3).$

\noindent
Then we have by \emph{Lemma 17.15 ,Lemma 17.18 }and\emph{\ Lemma 17.31 ,}

$[h_{\delta }+H+V+R,Sw_{1k}]=(\lambda _{k}+(\frac{1}{2}\mu _{1}+\frac{1}{%
3}v)-r)Sw_{1k},(0\leq k\leq 3),$

$[h_{\delta }+H+V+R,Sw_{4k}]=(-\lambda _{k}+(\frac{1}{2}\mu _{1}+\frac{1%
}{3}v)-r)Sw_{4k},(0\leq k\leq 3).$

\noindent
Hence $\pm \lambda _{k}+(\frac{1}{2}\mu _{1}+\frac{1}{3}v)-r$ are roots of 
\gr$_{8}^{\C}$ .

$(4)$ Let we put $Sx_{2k}$ and $Sx_{5k}:$

\ \ \ \ \ \ \ \ \ \ \ \ \ \ \ \ \ \ \ \ \ \ $Sx_{2k}=Ux_{2k}+iUx_{24+k},(0%
\leq k\leq 3),$

\ \ \ \ \ \ \ \ \ \ \ \ \ \ \ \ \ \ \ \ \ \ $Sx_{5k}=Ux_{2k}-iUx_{24+k},(0%
\leq k\leq 3).$

\noindent
Then we have by \emph{Lemma 17.2 ,Lemma 17.5 }and\emph{\ Lemma 17.30 },

$[h_{\delta }+H+V+R,Sx_{20}]=(\frac{1}{2}(-\lambda _{0}-\lambda
_{1}+\lambda _{2}-\lambda _{3})-(\frac{1}{2}\mu _{2}+\frac{1}{3}\nu
)+r)Sx_{20},$

\ \ \ \ \ \ \ \ \ \ \ \ \ \ \ \ \ \ \ \ \ \ \ \ \ \ \ \ \ \ \ \ \ \ $(0\leq k\leq 3),$

$[h_{\delta }+H+V+R,Sx_{21}]=(\frac{1}{2}(\lambda _{0}+\lambda
_{1}+\lambda _{2}-\lambda _{3})-(\frac{1}{2}\mu _{2}+\frac{1}{3}\nu
)+r)Sx_{21},$

\ \ \ \ \ \ \ \ \ \ \ \ \ \ \ \ \ \ \ \ \ \ \ \ \ \ \ \ \ \ \ \ \ \ $(0\leq k\leq 3),$

$[h_{\delta }+H+V+R,Sx_{22}]=(\frac{1}{2}(-\lambda _{0}+\lambda
_{1}+\lambda _{2}+\lambda _{3})-(\frac{1}{2}\mu _{2}+\frac{1}{3}\nu
)+r)Sx_{22},$

\ \ \ \ \ \ \ \ \ \ \ \ \ \ \ \ \ \ \ \ \ \ \ \ \ \ \ \ \ \ \ \ \ \ $(0\leq k\leq 3),$

$[h_{\delta }+H+V+R,Sx_{23}]=(\frac{1}{2}(\lambda _{0}-\lambda
_{1}+\lambda _{2}+\lambda _{3})-(\frac{1}{2}\mu _{2}+\frac{1}{3}\nu
)+r)Sx_{23},$

\ \ \ \ \ \ \ \ \ \ \ \ \ \ \ \ \ \ \ \ \ \ \ \ \ \ \ \ \ \ \ \ \ \ $(0\leq k\leq 3),$

$[h_{\delta }+H+V+R,Sx_{50}]=(-\frac{1}{2}(-\lambda _{0}-\lambda
_{1}+\lambda _{2}-\lambda _{3})-(\frac{1}{2}\mu _{2}+\frac{1}{3}\nu
)+r)Sx_{50},$

\ \ \ \ \ \ \ \ \ \ \ \ \ \ \ \ \ \ \ \ \ \ \ \ \ \ \ \ \ \ \ \ \ \ $(0\leq k\leq 3),$

$[h_{\delta }+H+V+R,Sx_{51}]=(-\frac{1}{2}(\lambda _{0}+\lambda
_{1}+\lambda _{2}-\lambda _{3})-(\frac{1}{2}\mu _{2}+\frac{1}{3}\nu
)+r)Sx_{51},$

\ \ \ \ \ \ \ \ \ \ \ \ \ \ \ \ \ \ \ \ \ \ \ \ \ \ \ \ \ \ \ \ \ \ $(0\leq k\leq 3),$

$[h_{\delta }+H+V+R,Sx_{52}]=(-\frac{1}{2}(-\lambda _{0}+\lambda
_{1}+\lambda _{2}+\lambda _{3})-(\frac{1}{2}\mu _{2}+\frac{1}{3}\nu
)+r)Sx_{52},$

\ \ \ \ \ \ \ \ \ \ \ \ \ \ \ \ \ \ \ \ \ \ \ \ \ \ \ \ \ \ \ \ \ \ $(0\leq k\leq 3),$

$[h_{\delta }+H+V+R,Sx_{53}]=(-\frac{1}{2}(\lambda _{0}-\lambda
_{1}+\lambda _{2}+\lambda _{3})-(\frac{1}{2}\mu _{2}+\frac{1}{3}\nu
)+r)Sx_{53},$

\ \ \ \ \ \ \ \ \ \ \ \ \ \ \ \ \ \ \ \ \ \ \ \ \ \ \ \ \ \ \ \ \ \ $(0\leq k\leq 3).$

\noindent
Let we put $Sy_{2k}$ and $Sy_{5k}:$

\ \ \ \ \ \ \ \ \ \ \ \ \ \ \ \ \ \ \ \ \ \ $Sy_{2k}=Uy_{2k}+iUy_{24+k},(0%
\leq k\leq 3),$

\ \ \ \ \ \ \ \ \ \ \ \ \ \ \ \ \ \ \ \ \ \ $Sy_{5k}=Uy_{2k}-iUy_{24+k},(0%
\leq k\leq 3).$

\noindent
Then we have by \emph{Lemma 17.6 ,Lemma 17.9 }and\emph{\ Lemma 17.30 ,}

$[h_{\delta }+H+V+R,Sy_{20}]=(\frac{1}{2}(-\lambda _{0}-\lambda
_{1}+\lambda _{2}-\lambda _{3})+(\frac{1}{2}\mu _{2}+\frac{1}{3}\nu
)+r)Sy_{20},$

\ \ \ \ \ \ \ \ \ \ \ \ \ \ \ \ \ \ \ \ \ \ \ \ \ \ \ \ \ \ \ \ \ \ $(0\leq k\leq 3),$

$[h_{\delta }+H+V+R,Sy_{21}]=(\frac{1}{2}(\lambda _{0}+\lambda
_{1}+\lambda _{2}-\lambda _{3})+(\frac{1}{2}\mu _{2}+\frac{1}{3}\nu
)+r)Sy_{21},$

\ \ \ \ \ \ \ \ \ \ \ \ \ \ \ \ \ \ \ \ \ \ \ \ \ \ \ \ \ \ \ \ \ \ $(0\leq k\leq 3),$

$[h_{\delta }+H+V+R,Sy_{22}]=(\frac{1}{2}(-\lambda _{0}+\lambda
_{1}+\lambda _{2}+\lambda _{3})+(\frac{1}{2}\mu _{2}+\frac{1}{3}\nu
)+r)Sy_{22},$

\ \ \ \ \ \ \ \ \ \ \ \ \ \ \ \ \ \ \ \ \ \ \ \ \ \ \ \ \ \ \ \ \ \ $(0\leq k\leq 3),$

$[h_{\delta }+H+V+R,Sy_{23}]=(\frac{1}{2}(\lambda _{0}-\lambda
_{1}+\lambda _{2}+\lambda _{3})+(\frac{1}{2}\mu _{2}+\frac{1}{3}\nu
)+r)Sy_{23},$

\ \ \ \ \ \ \ \ \ \ \ \ \ \ \ \ \ \ \ \ \ \ \ \ \ \ \ \ \ \ \ \ \ \ $(0\leq k\leq 3),$

$[h_{\delta }+H+V+R,Sy_{50}]=(-\frac{1}{2}(-\lambda _{0}-\lambda
_{1}+\lambda _{2}-\lambda _{3})+(\frac{1}{2}\mu _{2}+\frac{1}{3}\nu
)+r)Sy_{50},$

\ \ \ \ \ \ \ \ \ \ \ \ \ \ \ \ \ \ \ \ \ \ \ \ \ \ \ \ \ \ \ \ \ \ $(0\leq k\leq 3),$

$[h_{\delta }+H+V+R,Sy_{51}]=(-\frac{1}{2}(\lambda _{0}+\lambda
_{1}+\lambda _{2}-\lambda _{3})+(\frac{1}{2}\mu _{2}+\frac{1}{3}\nu
)+r)Sy_{51},$

\ \ \ \ \ \ \ \ \ \ \ \ \ \ \ \ \ \ \ \ \ \ \ \ \ \ \ \ \ \ \ \ \ \ $(0\leq k\leq 3),$

$[h_{\delta }+H+V+R,Sy_{52}]=(-\frac{1}{2}(-\lambda _{0}+\lambda
_{1}+\lambda _{2}+\lambda _{3})+(\frac{1}{2}\mu _{2}+\frac{1}{3}\nu
)+r)Sy_{52},$

\ \ \ \ \ \ \ \ \ \ \ \ \ \ \ \ \ \ \ \ \ \ \ \ \ \ \ \ \ \ \ \ \ \ $(0\leq k\leq 3),$

$[h_{\delta }+H+V+R,Sy_{53}]=(-\frac{1}{2}(\lambda _{0}-\lambda
_{1}+\lambda _{2}+\lambda _{3})+(\frac{1}{2}\mu _{2}+\frac{1}{3}\nu
)+r)Sy_{53},$

\ \ \ \ \ \ \ \ \ \ \ \ \ \ \ \ \ \ \ \ \ \ \ \ \ \ \ \ \ \ \ \ \ \ $(0\leq k\leq 3).$

\noindent
Let we put $Sz_{2k}$ and $Sz_{5k}:$

\ \ \ \ \ \ \ \ \ \ \ \ \ \ \ \ \ \ \ \ \ \ $Sz_{2k}=Uz_{2k}+iUz_{24+k},(0%
\leq k\leq 3),$

\ \ \ \ \ \ \ \ \ \ \ \ \ \ \ \ \ \ \ \ \ \ $Sz_{5k}=Uz_{2k}-iUz_{24+k},(0%
\leq k\leq 3).$

\noindent
Then we have by \emph{Lemma 17.11 ,Lemma 17.14 }and\emph{\ Lemma 17.31 ,}

$[h_{\delta }+H+V+R,Sx_{20}]=(\frac{1}{2}(-\lambda _{0}-\lambda
_{1}+\lambda _{2}-\lambda _{3})-(\frac{1}{2}\mu _{2}+\frac{1}{3}\nu
)-r)Sx_{20},$

\ \ \ \ \ \ \ \ \ \ \ \ \ \ \ \ \ \ \ \ \ \ \ \ \ \ \ \ \ \ \ \ \ \ $(0\leq k\leq 3),$

$[h_{\delta }+H+V+R,Sx_{21}]=(\frac{1}{2}(\lambda _{0}+\lambda
_{1}+\lambda _{2}-\lambda _{3})-(\frac{1}{2}\mu _{2}+\frac{1}{3}\nu
)-r)Sx_{21},$

\ \ \ \ \ \ \ \ \ \ \ \ \ \ \ \ \ \ \ \ \ \ \ \ \ \ \ \ \ \ \ \ \ \ $(0\leq k\leq 3),$

$[h_{\delta }+H+V+R,Sx_{22}]=(\frac{1}{2}(-\lambda _{0}+\lambda
_{1}+\lambda _{2}+\lambda _{3})-(\frac{1}{2}\mu _{2}+\frac{1}{3}\nu
)-r)Sx_{22},$

\ \ \ \ \ \ \ \ \ \ \ \ \ \ \ \ \ \ \ \ \ \ \ \ \ \ \ \ \ \ \ \ \ \ $(0\leq k\leq 3),$

$[h_{\delta }+H+V+R,Sx_{23}]=(\frac{1}{2}(\lambda _{0}-\lambda
_{1}+\lambda _{2}+\lambda _{3})-(\frac{1}{2}\mu _{2}+\frac{1}{3}\nu
)-r)Sx_{23},$

\ \ \ \ \ \ \ \ \ \ \ \ \ \ \ \ \ \ \ \ \ \ \ \ \ \ \ \ \ \ \ \ \ \ $(0\leq k\leq 3),$

$[h_{\delta }+H+V+R,Sx_{50}]=(-\frac{1}{2}(-\lambda _{0}-\lambda
_{1}+\lambda _{2}-\lambda _{3})-(\frac{1}{2}\mu _{2}+\frac{1}{3}\nu
)-r)Sx_{50},$

\ \ \ \ \ \ \ \ \ \ \ \ \ \ \ \ \ \ \ \ \ \ \ \ \ \ \ \ \ \ \ \ \ \ $(0\leq k\leq 3),$

$[h_{\delta }+H+V+R,Sx_{51}]=(-\frac{1}{2}(\lambda _{0}+\lambda
_{1}+\lambda _{2}-\lambda _{3})-(\frac{1}{2}\mu _{2}+\frac{1}{3}\nu
)-r)Sx_{51},$

\ \ \ \ \ \ \ \ \ \ \ \ \ \ \ \ \ \ \ \ \ \ \ \ \ \ \ \ \ \ \ \ \ \ $(0\leq k\leq 3),$

$[h_{\delta }+H+V+R,Sx_{52}]=(-\frac{1}{2}(-\lambda _{0}+\lambda
_{1}+\lambda _{2}+\lambda _{3})-(\frac{1}{2}\mu _{2}+\frac{1}{3}\nu
)-r)Sx_{52},$

\ \ \ \ \ \ \ \ \ \ \ \ \ \ \ \ \ \ \ \ \ \ \ \ \ \ \ \ \ \ \ \ \ \ $(0\leq k\leq 3),$

$[h_{\delta }+H+V+R,Sx_{53}]=(-\frac{1}{2}(\lambda _{0}-\lambda
_{1}+\lambda _{2}+\lambda _{3})-(\frac{1}{2}\mu _{2}+\frac{1}{3}\nu
)-r)Sx_{53},$

\ \ \ \ \ \ \ \ \ \ \ \ \ \ \ \ \ \ \ \ \ \ \ \ \ \ \ \ \ \ \ \ \ \ $(0\leq k\leq 3).$

\noindent
Let we put $Sw_{2k}$ and $Sw_{5k}:$

\ \ \ \ \ \ \ \ \ \ \ \ \ \ \ \ \ \ \ \ \ \ $Sw_{2k}=Uw_{2k}+iUw_{24+k},(0%
\leq k\leq 3),$

\ \ \ \ \ \ \ \ \ \ \ \ \ \ \ \ \ \ \ \ \ \ $Sw_{5k}=Uw_{2k}-iUw_{24+k},(0%
\leq k\leq 3).$

\noindent
Then we have by \emph{Lemma 17.15 ,Lemma 17.18 }and\emph{\ Lemma 17.31 ,}

$[h_{\delta }+H+V+R,Sw_{20}]=(\frac{1}{2}(-\lambda _{0}-\lambda
_{1}+\lambda _{2}-\lambda _{3})+(\frac{1}{2}\mu _{2}+\frac{1}{3}\nu
)-r)Sw_{20},$

\ \ \ \ \ \ \ \ \ \ \ \ \ \ \ \ \ \ \ \ \ \ \ \ \ \ \ \ \ \ \ \ \ \ $(0\leq k\leq 3),$

$[h_{\delta }+H+V+R,Sw_{21}]=(\frac{1}{2}(\lambda _{0}+\lambda
_{1}+\lambda _{2}-\lambda _{3})+(\frac{1}{2}\mu _{2}+\frac{1}{3}\nu
)-r)Sw_{21},$

\ \ \ \ \ \ \ \ \ \ \ \ \ \ \ \ \ \ \ \ \ \ \ \ \ \ \ \ \ \ \ \ \ \ $(0\leq k\leq 3),$

$[h_{\delta }+H+V+R,Sw_{22}]=(\frac{1}{2}(-\lambda _{0}+\lambda
_{1}+\lambda _{2}+\lambda _{3})+(\frac{1}{2}\mu _{2}+\frac{1}{3}\nu
)-r)Sw_{22},$

\ \ \ \ \ \ \ \ \ \ \ \ \ \ \ \ \ \ \ \ \ \ \ \ \ \ \ \ \ \ \ \ \ \ $(0\leq k\leq 3),$

$[h_{\delta }+H+V+R,Sw_{23}]=(\frac{1}{2}(\lambda _{0}-\lambda
_{1}+\lambda _{2}+\lambda _{3})+(\frac{1}{2}\mu _{2}+\frac{1}{3}\nu
)-r)Sw_{23},$

\ \ \ \ \ \ \ \ \ \ \ \ \ \ \ \ \ \ \ \ \ \ \ \ \ \ \ \ \ \ \ \ \ \ $(0\leq k\leq 3),$

$[h_{\delta }+H+V+R,Sw_{50}]=(-\frac{1}{2}(-\lambda _{0}-\lambda
_{1}+\lambda _{2}-\lambda _{3})+(\frac{1}{2}\mu _{2}+\frac{1}{3}\nu
)-r)Sw_{50},$

\ \ \ \ \ \ \ \ \ \ \ \ \ \ \ \ \ \ \ \ \ \ \ \ \ \ \ \ \ \ \ \ \ \ $(0\leq k\leq 3),$

$[h_{\delta }+H+V+R,Sw_{51}]=(-\frac{1}{2}(\lambda _{0}+\lambda
_{1}+\lambda _{2}-\lambda _{3})+(\frac{1}{2}\mu _{2}+\frac{1}{3}\nu
)-r)Sw_{51},$

\ \ \ \ \ \ \ \ \ \ \ \ \ \ \ \ \ \ \ \ \ \ \ \ \ \ \ \ \ \ \ \ \ \ $(0\leq k\leq 3),$

$[h_{\delta }+H+V+R,Sw_{52}]=(-\frac{1}{2}(-\lambda _{0}+\lambda
_{1}+\lambda _{2}+\lambda _{3})+(\frac{1}{2}\mu _{2}+\frac{1}{3}\nu
)-r)Sw_{52},$

\ \ \ \ \ \ \ \ \ \ \ \ \ \ \ \ \ \ \ \ \ \ \ \ \ \ \ \ \ \ \ \ \ \ $(0\leq k\leq 3),$

$[h_{\delta }+H+V+R,Sw_{53}]=(-\frac{1}{2}(\lambda _{0}-\lambda
_{1}+\lambda _{2}+\lambda _{3})+(\frac{1}{2}\mu _{2}+\frac{1}{3}\nu
)-r)Sw_{53},$

\ \ \ \ \ \ \ \ \ \ \ \ \ \ \ \ \ \ \ \ \ \ \ \ \ \ \ \ \ \ \ \ \ \ $(0\leq k\leq 3).$

\noindent
Hence $\pm \frac{1}{2}(-\lambda _{0}-\lambda _{1}+\lambda _{2}-\lambda
_{3})\pm (\frac{1}{2}\mu _{2}+\frac{1}{3}\nu )\pm r,$

$\pm \frac{1}{2}(\lambda _{0}+\lambda _{1}+\lambda _{2}-\lambda _{3})\pm (%
\frac{1}{2}\mu _{2}+\frac{1}{3}\nu )\pm r,$

$\pm \frac{1}{2}(-\lambda _{0}+\lambda _{1}+\lambda _{2}+\lambda _{3})\pm (%
\frac{1}{2}\mu _{2}+\frac{1}{3}\nu )\pm r,\pm \frac{1}{2}(\lambda
_{0}-\lambda _{1}+\lambda _{2}+\lambda _{3})\pm (\frac{1}{2}\mu _{2}+\frac{1%
}{3}\nu )\pm r$ are roots of \gr$_{8}^{\C}$ .

$(5)$ Let we put $Sx_{3k}$ and $Sx_{6k}:$

\ \ \ \ \ \ \ \ \ \ \ \ \ \ \ \ \ \ \ \ \ \ $Sx_{3k}=Ux_{3k}+iUx_{34+k},(0%
\leq k\leq 3),$

\ \ \ \ \ \ \ \ \ \ \ \ \ \ \ \ \ \ \ \ \ \ $Sx_{5k}=Ux_{3k}-iUx_{34+k},(0%
\leq k\leq 3).$

\noindent
Then we have by \emph{Lemma 17.2 ,Lemma 17.5 }and\emph{\ Lemma 17.30 }

$[h_{\delta }+H+V+R,Sx_{30}]=(-\frac{1}{2}(\lambda _{0}-\lambda _{1}+\lambda
_{2}-\lambda _{3})-(-\frac{1}{2}\mu _{1}-\frac{1}{2}\mu _{2}+\frac{1}{3}\nu
)+r)Sx_{30},$

\ \ \ \ \ \ \ \ \ \ \ \ \ \ \ \ \ \ \ \ \ \ \ \ \ \ \ \ \ \ \ \ \ \ $(0\leq k\leq 3),$

$[h_{\delta }+H+V+R,Sx_{31}]=(-\frac{1}{2}(\lambda _{0}-\lambda _{1}-\lambda
_{2}+\lambda _{3})-(-\frac{1}{2}\mu _{1}-\frac{1}{2}\mu _{2}+\frac{1}{3}\nu
)+r)Sx_{31},$

\ \ \ \ \ \ \ \ \ \ \ \ \ \ \ \ \ \ \ \ \ \ \ \ \ \ \ \ \ \ \ \ \ \ $(0\leq k\leq 3),$

$[h_{\delta }+H+V+R,Sx_{32}]=(\frac{1}{2}(\lambda _{0}+\lambda _{1}+\lambda
_{2}+\lambda _{3})-(-\frac{1}{2}\mu _{1}-\frac{1}{2}\mu _{2}+\frac{1}{3}\nu
)+r)Sx_{32},$

\ \ \ \ \ \ \ \ \ \ \ \ \ \ \ \ \ \ \ \ \ \ \ \ \ \ \ \ \ \ \ \ \ \ $(0\leq k\leq 3),$

$[h_{\delta }+H+V+R,Sx_{33}]=(-\frac{1}{2}(\lambda _{0}+\lambda _{1}-\lambda
_{2}-\lambda _{3})-(-\frac{1}{2}\mu _{1}-\frac{1}{2}\mu _{2}+\frac{1}{3}\nu
)+r)Sx_{33},$

\ \ \ \ \ \ \ \ \ \ \ \ \ \ \ \ \ \ \ \ \ \ \ \ \ \ \ \ \ \ \ \ \ \ $(0\leq k\leq 3),$

$[h_{\delta }+H+V+R,Sx_{60}]=(\frac{1}{2}(\lambda _{0}-\lambda _{1}+\lambda
_{2}-\lambda _{3})-(-\frac{1}{2}\mu _{1}-\frac{1}{2}\mu _{2}+\frac{1}{3}\nu
)+r)Sx_{60},$

\ \ \ \ \ \ \ \ \ \ \ \ \ \ \ \ \ \ \ \ \ \ \ \ \ \ \ \ \ \ \ \ \ \ $(0\leq k\leq 3),$

$[h_{\delta }+H+V+R,Sx_{61}]=(\frac{1}{2}(\lambda _{0}-\lambda _{1}-\lambda
_{2}+\lambda _{3})-(-\frac{1}{2}\mu _{1}-\frac{1}{2}\mu _{2}+\frac{1}{3}\nu
)+r)Sx_{61},$

\ \ \ \ \ \ \ \ \ \ \ \ \ \ \ \ \ \ \ \ \ \ \ \ \ \ \ \ \ \ \ \ \ \ $(0\leq k\leq 3),$

$[h_{\delta }+H+V+R,Sx_{62}]=(-\frac{1}{2}(\lambda _{0}+\lambda _{1}+\lambda
_{2}+\lambda _{3})-(-\frac{1}{2}\mu _{1}-\frac{1}{2}\mu _{2}+\frac{1}{3}\nu
)+r)Sx_{62},$

\ \ \ \ \ \ \ \ \ \ \ \ \ \ \ \ \ \ \ \ \ \ \ \ \ \ \ \ \ \ \ \ \ \ $(0\leq k\leq 3),$

$[h_{\delta }+H+V+R,Sx_{63}]=(\frac{1}{2}(\lambda _{0}+\lambda _{1}-\lambda
_{2}-\lambda _{3})-(-\frac{1}{2}\mu _{1}-\frac{1}{2}\mu _{2}+\frac{1}{3}\nu
)+r)Sx_{63},$

$\ \ \ \ \ \ \ \ \ \ \ \ \ \ \ \ \ \ \ \ \ \ \ \ \ \ \ \ \ \ \ \ \ \ (0\leq k\leq 3).$

\noindent
Let we put $Sy_{3k}$ and $Sy_{6k}:$

\ \ \ \ \ \ \ \ \ \ \ \ \ \ \ \ \ \ \ \ \ \ $Sy_{3k}=Uy_{3k}+iUy_{34+k},(0%
\leq k\leq 3),$

\ \ \ \ \ \ \ \ \ \ \ \ \ \ \ \ \ \ \ \ \ \ $Sy_{6k}=Uy_{3k}-iUy_{34+k},(0%
\leq k\leq 3).$

\noindent
Then we have by \emph{Lemma 17.6 ,Lemma 17.9 }and\emph{\ Lemma 17.30 ,}

$[h_{\delta }+H+V+R,Sy_{30}]=(-\frac{1}{2}(\lambda _{0}-\lambda _{1}+\lambda
_{2}-\lambda _{3})+(-\frac{1}{2}\mu _{1}-\frac{1}{2}\mu _{2}+\frac{1}{3}\nu
)+r)Sy_{30},$

\ \ \ \ \ \ \ \ \ \ \ \ \ \ \ \ \ \ \ \ \ \ \ \ \ \ \ \ \ \ \ \ \ \ $(0\leq k\leq 3),$

$[h_{\delta }+H+V+R,Sy_{31}]=(-\frac{1}{2}(\lambda _{0}-\lambda _{1}-\lambda
_{2}+\lambda _{3})+(-\frac{1}{2}\mu _{1}-\frac{1}{2}\mu _{2}+\frac{1}{3}\nu
)+r)Sy_{31},$

\ \ \ \ \ \ \ \ \ \ \ \ \ \ \ \ \ \ \ \ \ \ \ \ \ \ \ \ \ \ \ \ \ \ $(0\leq k\leq 3),$

$[h_{\delta }+H+V+R,Sy_{32}]=(\frac{1}{2}(\lambda _{0}+\lambda _{1}+\lambda
_{2}+\lambda _{3})+(-\frac{1}{2}\mu _{1}-\frac{1}{2}\mu _{2}+\frac{1}{3}\nu
)+r)Sy_{32},$

\ \ \ \ \ \ \ \ \ \ \ \ \ \ \ \ \ \ \ \ \ \ \ \ \ \ \ \ \ \ \ \ \ \ $(0\leq k\leq 3),$

$[h_{\delta }+H+V+R,Sy_{33}]=(-\frac{1}{2}(\lambda _{0}+\lambda _{1}-\lambda
_{2}-\lambda _{3})+(-\frac{1}{2}\mu _{1}-\frac{1}{2}\mu _{2}+\frac{1}{3}\nu
)+r)Sy_{33},$

\ \ \ \ \ \ \ \ \ \ \ \ \ \ \ \ \ \ \ \ \ \ \ \ \ \ \ \ \ \ \ \ \ \ $(0\leq k\leq 3),$

$[h_{\delta }+H+V+R,Sy_{60}]=(\frac{1}{2}(\lambda _{0}-\lambda _{1}+\lambda
_{2}-\lambda _{3})+(-\frac{1}{2}\mu _{1}-\frac{1}{2}\mu _{2}+\frac{1}{3}\nu
)+r)Sy_{60},$

\ \ \ \ \ \ \ \ \ \ \ \ \ \ \ \ \ \ \ \ \ \ \ \ \ \ \ \ \ \ \ \ \ \ $(0\leq k\leq 3),$

$[h_{\delta }+H+V+R,Sy_{61}]=(\frac{1}{2}(\lambda _{0}-\lambda _{1}-\lambda
_{2}+\lambda _{3})+(-\frac{1}{2}\mu _{1}-\frac{1}{2}\mu _{2}+\frac{1}{3}\nu
)+r)Sy_{61},$

\ \ \ \ \ \ \ \ \ \ \ \ \ \ \ \ \ \ \ \ \ \ \ \ \ \ \ \ \ \ \ \ \ \ $(0\leq k\leq 3),$

$[h_{\delta }+H+V+R,Sy_{62}]=(-\frac{1}{2}(\lambda _{0}+\lambda _{1}+\lambda
_{2}+\lambda _{3})+(-\frac{1}{2}\mu _{1}-\frac{1}{2}\mu _{2}+\frac{1}{3}\nu
)+r)Sy_{62},$

\ \ \ \ \ \ \ \ \ \ \ \ \ \ \ \ \ \ \ \ \ \ \ \ \ \ \ \ \ \ \ \ \ \ $(0\leq k\leq 3),$

$[h_{\delta }+H+V+R,Sy_{63}]=(\frac{1}{2}(\lambda _{0}+\lambda _{1}-\lambda
_{2}-\lambda _{3})+(-\frac{1}{2}\mu _{1}-\frac{1}{2}\mu _{2}+\frac{1}{3}\nu
)+r)Sy_{63},$

\ \ \ \ \ \ \ \ \ \ \ \ \ \ \ \ \ \ \ \ \ \ \ \ \ \ \ \ \ \ \ \ \ \ $(0\leq k\leq 3).$

\noindent
Let we put $Sz_{3k}$ and $Sz_{6k}:$

\ \ \ \ \ \ \ \ \ \ \ \ \ \ \ \ \ \ \ \ \ \ $Sz_{3k}=Uz_{3k}+iUz_{34+k},(0%
\leq k\leq 3),$

\ \ \ \ \ \ \ \ \ \ \ \ \ \ \ \ \ \ \ \ \ \ $Sz_{5k}=Uz_{3k}-iUz_{34+k},(0%
\leq k\leq 3).$

\noindent
Then we have by \emph{Lemma 17.11 ,Lemma 17.14 }and\emph{\ Lemma 17.31 ,}

$[h_{\delta }+H+V+R,Sz_{30}]=(-\frac{1}{2}(\lambda _{0}-\lambda _{1}+\lambda
_{2}-\lambda _{3})-(-\frac{1}{2}\mu _{1}-\frac{1}{2}\mu _{2}+\frac{1}{3}\nu
)-r)Sz_{30},$

\ \ \ \ \ \ \ \ \ \ \ \ \ \ \ \ \ \ \ \ \ \ \ \ \ \ \ \ \ \ \ \ \ \ $(0\leq k\leq 3),$

$[h_{\delta }+H+V+R,Sz_{31}]=(-\frac{1}{2}(\lambda _{0}-\lambda _{1}-\lambda
_{2}+\lambda _{3})-(-\frac{1}{2}\mu _{1}-\frac{1}{2}\mu _{2}+\frac{1}{3}\nu
)-r)Sz_{31},$

\ \ \ \ \ \ \ \ \ \ \ \ \ \ \ \ \ \ \ \ \ \ \ \ \ \ \ \ \ \ \ \ \ \ $(0\leq k\leq 3),$

$[h_{\delta }+H+V+R,Sz_{32}]=(\frac{1}{2}(\lambda _{0}+\lambda _{1}+\lambda
_{2}+\lambda _{3})-(-\frac{1}{2}\mu _{1}-\frac{1}{2}\mu _{2}+\frac{1}{3}\nu
)-r)Sz_{32},$

\ \ \ \ \ \ \ \ \ \ \ \ \ \ \ \ \ \ \ \ \ \ \ \ \ \ \ \ \ \ \ \ \ \ $(0\leq k\leq 3),$

$[h_{\delta }+H+V+R,Sz_{33}]=(-\frac{1}{2}(\lambda _{0}+\lambda _{1}-\lambda
_{2}-\lambda _{3})-(-\frac{1}{2}\mu _{1}-\frac{1}{2}\mu _{2}+\frac{1}{3}\nu
)-r)Sz_{33},$

\ \ \ \ \ \ \ \ \ \ \ \ \ \ \ \ \ \ \ \ \ \ \ \ \ \ \ \ \ \ \ \ \ \ $(0\leq k\leq 3),$

$[h_{\delta }+H+V+R,Sz_{60}]=(\frac{1}{2}(\lambda _{0}-\lambda _{1}+\lambda
_{2}-\lambda _{3})-(-\frac{1}{2}\mu _{1}-\frac{1}{2}\mu _{2}+\frac{1}{3}\nu
)-r)Sz_{60},$

\ \ \ \ \ \ \ \ \ \ \ \ \ \ \ \ \ \ \ \ \ \ \ \ \ \ \ \ \ \ \ \ \ \ $(0\leq k\leq 3),$

$[h_{\delta }+H+V+R,Sz_{61}]=(\frac{1}{2}(\lambda _{0}-\lambda _{1}-\lambda
_{2}+\lambda _{3})-(-\frac{1}{2}\mu _{1}-\frac{1}{2}\mu _{2}+\frac{1}{3}\nu
)-r)Sz_{61},$

\ \ \ \ \ \ \ \ \ \ \ \ \ \ \ \ \ \ \ \ \ \ \ \ \ \ \ \ \ \ \ \ \ \ $(0\leq k\leq 3),$

$[h_{\delta }+H+V+R,Sz_{62}]=(-\frac{1}{2}(\lambda _{0}+\lambda _{1}+\lambda
_{2}+\lambda _{3})-(-\frac{1}{2}\mu _{1}-\frac{1}{2}\mu _{2}+\frac{1}{3}\nu
)-r)Sz_{62},$

\ \ \ \ \ \ \ \ \ \ \ \ \ \ \ \ \ \ \ \ \ \ \ \ \ \ \ \ \ \ \ \ \ \ $(0\leq k\leq 3),$

$[h_{\delta }+H+V+R,Sz_{63}]=(\frac{1}{2}(\lambda _{0}+\lambda _{1}-\lambda
_{2}-\lambda _{3})-(-\frac{1}{2}\mu _{1}-\frac{1}{2}\mu _{2}+\frac{1}{3}\nu
)-r)Sz_{63},$

\ \ \ \ \ \ \ \ \ \ \ \ \ \ \ \ \ \ \ \ \ \ \ \ \ \ \ \ \ \ \ \ \ \ $(0\leq k\leq 3).$

\noindent
Let we put $Sw_{3k}$ and $Sw_{6k}:$

\ \ \ \ \ \ \ \ \ \ \ \ \ \ \ \ \ \ \ \ \ \ $Sw_{3k}=Uw_{3k}+iUw_{34+k},(0%
\leq k\leq 3),$

\ \ \ \ \ \ \ \ \ \ \ \ \ \ \ \ \ \ \ \ \ \ $Sw_{6k}=Uw_{3k}-iUw_{34+k},(0%
\leq k\leq 3).$

\noindent
Then we have by \emph{Lemma 17.15 ,Lemma 17.18 }and\emph{\ Lemma 17.31 ,}

$[h_{\delta }+H+V+R,Sw_{30}]=(-\frac{1}{2}(\lambda _{0}-\lambda _{1}+\lambda
_{2}-\lambda _{3})+(-\frac{1}{2}\mu _{1}-\frac{1}{2}\mu _{2}+\frac{1}{3}\nu
)-r)Sw_{30},$

\ \ \ \ \ \ \ \ \ \ \ \ \ \ \ \ \ \ \ \ \ \ \ \ \ \ \ \ \ \ \ \ \ \ $(0\leq k\leq 3),$

$[h_{\delta }+H+V+R,Sw_{31}]=(-\frac{1}{2}(\lambda _{0}-\lambda _{1}-\lambda
_{2}+\lambda _{3})+(-\frac{1}{2}\mu _{1}-\frac{1}{2}\mu _{2}+\frac{1}{3}\nu
)-r)Sw_{31},$

\ \ \ \ \ \ \ \ \ \ \ \ \ \ \ \ \ \ \ \ \ \ \ \ \ \ \ \ \ \ \ \ \ \ $(0\leq k\leq 3),$

$[h_{\delta }+H+V+R,Sw_{32}]=(\frac{1}{2}(\lambda _{0}+\lambda _{1}+\lambda
_{2}+\lambda _{3})+(-\frac{1}{2}\mu _{1}-\frac{1}{2}\mu _{2}+\frac{1}{3}\nu
)-r)Sw_{32},$

\ \ \ \ \ \ \ \ \ \ \ \ \ \ \ \ \ \ \ \ \ \ \ \ \ \ \ \ \ \ \ \ \ \ $(0\leq k\leq 3),$

$[h_{\delta }+H+V+R,Sw_{33}]=(-\frac{1}{2}(\lambda _{0}+\lambda _{1}-\lambda
_{2}-\lambda _{3})+(-\frac{1}{2}\mu _{1}-\frac{1}{2}\mu _{2}+\frac{1}{3}\nu
)-r)Sw_{33},$

\ \ \ \ \ \ \ \ \ \ \ \ \ \ \ \ \ \ \ \ \ \ \ \ \ \ \ \ \ \ \ \ \ \ $(0\leq k\leq 3),$

$[h_{\delta }+H+V+R,Sw_{60}]=(\frac{1}{2}(\lambda _{0}-\lambda _{1}+\lambda
_{2}-\lambda _{3})+(-\frac{1}{2}\mu _{1}-\frac{1}{2}\mu _{2}+\frac{1}{3}\nu
)-r)Sw_{60},$

\ \ \ \ \ \ \ \ \ \ \ \ \ \ \ \ \ \ \ \ \ \ \ \ \ \ \ \ \ \ \ \ \ \ $(0\leq k\leq 3),$

$[h_{\delta }+H+V+R,Sw_{61}]=(\frac{1}{2}(\lambda _{0}-\lambda _{1}-\lambda
_{2}+\lambda _{3})+(-\frac{1}{2}\mu _{1}-\frac{1}{2}\mu _{2}+\frac{1}{3}\nu
)-r)Sw_{61},$

\ \ \ \ \ \ \ \ \ \ \ \ \ \ \ \ \ \ \ \ \ \ \ \ \ \ \ \ \ \ \ \ \ \ $(0\leq k\leq 3),$

$[h_{\delta }+H+V+R,Sw_{62}]=(-\frac{1}{2}(\lambda _{0}+\lambda _{1}+\lambda
_{2}+\lambda _{3})+(-\frac{1}{2}\mu _{1}-\frac{1}{2}\mu _{2}+\frac{1}{3}\nu
)-r)Sw_{62},$

\ \ \ \ \ \ \ \ \ \ \ \ \ \ \ \ \ \ \ \ \ \ \ \ \ \ \ \ \ \ \ \ \ \ $(0\leq k\leq 3),$

$[h_{\delta }+H+V+R,Sw_{63}]=(\frac{1}{2}(\lambda _{0}+\lambda _{1}-\lambda
_{2}-\lambda _{3})+(-\frac{1}{2}\mu _{1}-\frac{1}{2}\mu _{2}+\frac{1}{3}\nu
)-r)Sw_{63},$

\ \ \ \ \ \ \ \ \ \ \ \ \ \ \ \ \ \ \ \ \ \ \ \ \ \ \ \ \ \ \ \ \ \ $(0\leq k\leq 3).$

\noindent
Hence $\pm \frac{1}{2}(\lambda _{0}-\lambda _{1}+\lambda _{2}-\lambda
_{3})\pm (-\frac{1}{2}\mu _{1}-\frac{1}{2}\mu _{2}+\frac{1}{3}\nu )\pm r,$

$\pm \frac{1}{2}(\lambda _{0}-\lambda _{1}-\lambda _{2}+\lambda _{3})\pm (-%
\frac{1}{2}\mu _{1}-\frac{1}{2}\mu _{2}+\frac{1}{3}\nu )\pm r,$

$\pm \frac{1}{2}(\lambda _{0}+\lambda _{1}+\lambda _{2}+\lambda _{3})\pm (-%
\frac{1}{2}\mu _{1}-\frac{1}{2}\mu _{2}+\frac{1}{3}\nu )\pm r,$

$\pm \frac{1}{2}(\lambda _{0}+\lambda _{1}-\lambda _{2}-\lambda _{3})\pm (-%
\frac{1}{2}\mu _{1}-\frac{1}{2}\mu _{2}+\frac{1}{3}\nu )\pm r$ are roots of 
\gr$_{8}^{\C}$ .

$(6)$ Let we put $Su=Uu_{1}$ and $Ss=Us_{1}$.

\noindent
Then we have by \emph{Lemma 17.10 }and\emph{\ Lemma 17.30 ,}

\ \ \ \ \ \ \ \ \ \ \ \ \ \ \ \ \ \ \ \ \ \ \ \ \ $[h_{\delta
}+H+V+R,Ss]=2rSs,$

\ \ \ \ \ \ \ \ \ \ \ \ \ \ \ \ \ \ \ \ \ \ \ \ \ $[h_{\delta
}+H+V+R,Su]=-2rSu.$

\noindent
Hence $\pm 2r$ are roots of \gr$_{8}^{\C}$ .

$(7)$ Let we put $S\xi =U\xi _{1},S\eta =U\eta _{1},S\zeta =U\zeta_{1}$ and $S\omega
=U\omega_{1}$.

\noindent
Then we have by \emph{Lemma 17.19 }and \emph{Lemma 17.31 ,}

\ \ \ \ \ \ \ \ \ \ \ \ \ \ \ \ \ \ \ \ \ \ \ \ \ $[h_{\delta }+H+V+R,S\xi
]=(v+r)S\xi ,$

\ \ \ \ \ \ \ \ \ \ \ \ \ \ \ \ \ \ \ \ \ \ \ \ \ $[h_{\delta }+H+V+R,S\eta
]=(-v+r)S\eta ,$

\ \ \ \ \ \ \ \ \ \ \ \ \ \ \ \ \ \ \ \ \ \ \ \ \ $[h_{\delta }+H+V+R,S\zeta
]=(v-r)S\zeta ,$

\ \ \ \ \ \ \ \ \ \ \ \ \ \ \ \ \ \ \ \ \ \ \ \ \ $[h_{\delta
}+H+V+R,S\omega ]=(-v-r)S\omega .$

\noindent
Hence $\pm v\pm r$ are roots of \gr$_{8}^{\C}$ .

By direct calculations using Maxima, we also have $(1),(2),(3),(4),(5),(6)$ and $(7)$.
\ \ \ \ \emph{Q.E.D.}

\bigskip

\emph{Theorem 17.37.} In the root system of Lemma 17.34 ,

\ \ \ \ \ \ \ \ \ \ \ \ \ \ \ \ \ \ \ \ \ \ $\alpha _{1}=\frac{1}{2}(\lambda
0-\lambda 1-\lambda 2-\lambda 3)+\frac{1}{2}(-2\mu _{1}-\mu _{2}),$

\ \ \ \ \ \ \ \ \ \ \ \ \ \ \ \ \ \ \ \ \ \ $\alpha _{2}=\mu _{1}-\frac{1}{3}v-r,$
\ \ $\alpha _{3}=2r,$

\ \ \ \ \ \ \ \ \ \ \ \ \ \ \ \ \ \ \ \ \ \ $\alpha _{4}=\mu _{2}-\frac{1}{3}v-r,$
\ \ $\alpha _{5}=\lambda _{3}-\frac{1}{2}(\mu _{1}+2\mu _{2}),$

\ \ \ \ \ \ \ \ \ \ \ \ \ \ \ \ \ \ \ \ \ \ $\alpha _{6}=\lambda _{2}-\lambda _{3},$
\ \ $\alpha _{7}=\lambda _{1}-\lambda _{2},$ \ \ \ $\alpha _{8}=v-r,$

\noindent
is a fundamental root of the Lie algebra \gr$_{8}^{\C}$ and

\ \ \ \ \ \ \ \ \ \ \ \ \ \ \ \ \ \ \ \ \ \ $\mu =2\alpha _{1}+4\alpha
_{2}+6\alpha _{3}+5\alpha _{4}+4\alpha _{5}+3\alpha _{6}+2\alpha _{7}+3\alpha _{8}$

\noindent
is the highest root. The Dynkin diagram and the extended Dynkin diagram of 
\gr$_{8}^{\C}$ are respectively given by

\begin{figure}[H]
\centering
\includegraphics[width=12cm, height=3.0cm]{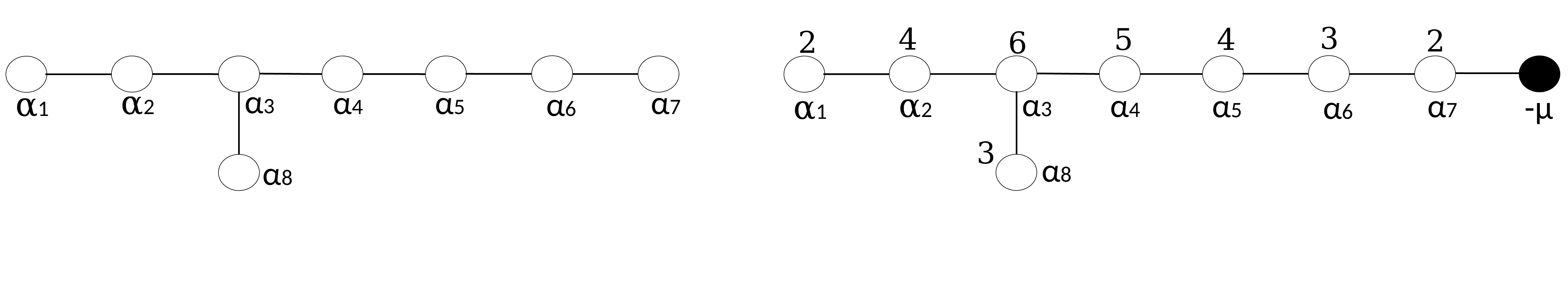}
\end{figure}

\emph{Proof.} In the following,the notation $n_{1}$ $n_{2}$ $n_{3}$ $%
n_{4}$ $n_{5}$ $n_{6}$ $n_{7}$ $n_{8}$ denotes the root
$n_{1}\alpha _{1}+n_{2}\alpha _{2}+n_{3}\alpha _{3}+n_{4}\alpha
_{4}+n_{5}\alpha _{5}+n_{6}\alpha _{6}+n_{7}\alpha _{7}+n_{8}\alpha _{8}.$
Now,all positive roots of \gr$_{8}^{\C}$ are represented by

\begin{flushright}
$\lambda _{0}-\lambda _{1}=%
\begin{array}{llllllll}
2 & 3 & 4 & 3 & 2 & 1 & 0 & 2%
\end{array},$

$\lambda _{0}+\lambda _{1}=%
\begin{array}{llllllll}
2 & 4 & 6 & 5 & 4 & 3 & 2 & 3%
\end{array},$

$\lambda _{0}-\lambda _{2}=%
\begin{array}{llllllll}
2 & 3 & 4 & 3 & 2 & 1 & 1 & 2%
\end{array},$

$\lambda _{0}+\lambda _{2}=%
\begin{array}{llllllll}
2 & 4 & 6 & 5 & 4 & 3 & 1 & 3%
\end{array},$

$\lambda _{0}-\lambda _{3}=%
\begin{array}{llllllll}
2 & 3 & 4 & 3 & 2 & 2 & 1 & 2%
\end{array},$

$\lambda _{0}+\lambda _{3}=%
\begin{array}{llllllll}
2 & 4 & 6 & 5 & 4 & 2 & 1 & 3%
\end{array},$

$\lambda _{1}-\lambda _{2}=%
\begin{array}{llllllll}
0 & 0 & 0 & 0 & 0 & 0 & 1 & 0%
\end{array},$

$\lambda _{1}+\lambda _{2}=%
\begin{array}{llllllll}
0 & 1 & 2 & 2 & 2 & 2 & 1 & 1%
\end{array},$

$\lambda _{1}-\lambda _{3}=%
\begin{array}{llllllll}
0 & 0 & 0 & 0 & 0 & 1 & 1 & 0%
\end{array},$

$\lambda _{1}+\lambda _{3}=%
\begin{array}{llllllll}
0 & 1 & 2 & 2 & 2 & 1 & 1 & 1%
\end{array},$

$\lambda _{2}-\lambda _{3}=%
\begin{array}{llllllll}
0 & 0 & 0 & 0 & 0 & 1 & 0 & 0%
\end{array},$

$\lambda _{2}+\lambda _{3}=%
\begin{array}{llllllll}
0 & 1 & 2 & 2 & 2 & 1 & 0 & 1%
\end{array},$

$\lambda _{0}+\frac{1}{2}(\mu _{1}+2\mu _{2})=%
\begin{array}{llllllll}
2 & 4 & 6 & 5 & 3 & 2 & 1 & 3%
\end{array},$

$\lambda _{1}+\frac{1}{2}(\mu _{1}+2\mu _{2})=%
\begin{array}{llllllll}
0 & 1 & 2 & 2 & 1 & 1 & 1 & 1%
\end{array},$

$\lambda _{2}+\frac{1}{2}(\mu _{1}+2\mu _{2})=%
\begin{array}{llllllll}
0 & 1 & 2 & 2 & 1 & 1 & 0 & 1%
\end{array},$

$\lambda _{3}+\frac{1}{2}(\mu _{1}+2\mu _{2})=%
\begin{array}{llllllll}
0 & 1 & 2 & 2 & 1 & 0 & 0 & 1%
\end{array},$

$\lambda _{0}-\frac{1}{2}(\mu _{1}+2\mu _{2})=%
\begin{array}{llllllll}
2 & 3 & 4 & 3 & 3 & 2 & 1 & 2%
\end{array},$

$\lambda _{1}-\frac{1}{2}(\mu _{1}+2\mu _{2})=%
\begin{array}{llllllll}
0 & 0 & 0 & 0 & 1 & 1 & 1 & 0%
\end{array},$

$\lambda _{2}-\frac{1}{2}(\mu _{1}+2\mu _{2})=%
\begin{array}{llllllll}
0 & 0 & 0 & 0 & 1 & 1 & 0 & 0%
\end{array},$

$\lambda _{3}-\frac{1}{2}(\mu _{1}+2\mu _{2})=%
\begin{array}{llllllll}
0 & 0 & 0 & 0 & 1 & 0 & 0 & 0%
\end{array},$

$\frac{1}{2}(\lambda _{0}+\lambda _{1}-\lambda
_{2}+\lambda _{3})+\frac{1}{2}(-2\mu _{1}-\mu _{2})=%
\begin{array}{llllllll}
1 & 1 & 2 & 2 & 2 & 1 & 1 & 1%
\end{array},$

$\frac{1}{2}(\lambda _{0}+\lambda _{1}+\lambda
_{2}-\lambda _{3})+\frac{1}{2}(-2\mu _{1}-\mu _{2})=%
\begin{array}{llllllll}
1 & 1 & 2 & 2 & 2 & 1 & 1 & 1%
\end{array},$

$\frac{1}{2}(\lambda _{0}-\lambda _{1}-\lambda
_{2}-\lambda _{3})+\frac{1}{2}(-2\mu _{1}-\mu _{2})=%
\begin{array}{llllllll}
1 & 0 & 0 & 0 & 0 & 0 & 0 & 0%
\end{array},$

$\frac{1}{2}(\lambda _{0}-\lambda _{1}+\lambda
_{2}+\lambda _{3})+\frac{1}{2}(-2\mu _{1}-\mu _{2})=%
\begin{array}{llllllll}
1 & 1 & 2 & 2 & 2 & 1 & 0 & 1%
\end{array},$

$\frac{1}{2}(\lambda _{0}+\lambda _{1}-\lambda
_{2}+\lambda _{3})-\frac{1}{2}(-2\mu _{1}-\mu _{2})=%
\begin{array}{llllllll}
1 & 3 & 4 & 3 & 2 & 1 & 1 & 2%
\end{array},$

$\frac{1}{2}(\lambda _{0}+\lambda _{1}+\lambda
_{2}-\lambda _{3})-\frac{1}{2}(-2\mu _{1}-\mu _{2})=%
\begin{array}{llllllll}
1 & 3 & 4 & 3 & 2 & 2 & 1 & 2%
\end{array},$

$\frac{1}{2}(\lambda _{0}-\lambda _{1}-\lambda
_{2}-\lambda _{3})-\frac{1}{2}(-2\mu _{1}-\mu _{2})=%
\begin{array}{llllllll}
1 & 2 & 2 & 1 & 0 & 0 & 0 & 1%
\end{array},$

$\frac{1}{2}(\lambda _{0}-\lambda _{1}+\lambda
_{2}+\lambda _{3})-\frac{1}{2}(-2\mu _{1}-\mu _{2})=%
\begin{array}{llllllll}
1 & 3 & 4 & 3 & 2 & 1 & 0 & 2%
\end{array},$

$\frac{1}{2}(\lambda _{0}-\lambda _{1}+\lambda
_{2}-\lambda _{3})+\frac{1}{2}(\mu _{1}-\mu _{2})=%
\begin{array}{llllllll}
1 & 2 & 2 & 1 & 1 & 1 & 0 & 1%
\end{array},$

$\frac{1}{2}(\lambda _{0}-\lambda _{1}-\lambda
_{2}+\lambda _{3})+\frac{1}{2}(\mu _{1}-\mu _{2})=%
\begin{array}{llllllll}
1 & 2 & 2 & 1 & 1 & 0 & 0 & 1%
\end{array},$

$\frac{1}{2}(\lambda _{0}+\lambda _{1}+\lambda
_{2}+\lambda _{3})+\frac{1}{2}(\mu _{1}-\mu _{2})=%
\begin{array}{llllllll}
1 & 3 & 4 & 3 & 3 & 2 & 1 & 2%
\end{array},$

$\frac{1}{2}(\lambda _{0}+\lambda _{1}-\lambda
_{2}-\lambda _{3})+\frac{1}{2}(\mu _{1}-\mu _{2})=%
\begin{array}{llllllll}
1 & 2 & 2 & 1 & 1 & 1 & 1 & 1%
\end{array},$

$\frac{1}{2}(\lambda _{0}-\lambda _{1}+\lambda
_{2}-\lambda _{3})-\frac{1}{2}(\mu _{1}-\mu _{2})=%
\begin{array}{llllllll}
1 & 1 & 2 & 2 & 1 & 1 & 0 & 1%
\end{array},$

$\frac{1}{2}(\lambda _{0}-\lambda _{1}-\lambda
_{2}+\lambda _{3})-\frac{1}{2}(\mu _{1}-\mu _{2})=%
\begin{array}{llllllll}
1 & 1 & 2 & 2 & 1 & 0 & 0 & 1%
\end{array},$

$\frac{1}{2}(\lambda _{0}+\lambda _{1}+\lambda
_{2}+\lambda _{3})-\frac{1}{2}(\mu _{1}-\mu _{2})=%
\begin{array}{llllllll}
1 & 2 & 4 & 4 & 3 & 2 & 1 & 2%
\end{array},$

$\frac{1}{2}(\lambda _{0}+\lambda _{1}-\lambda
_{2}-\lambda _{3})-\frac{1}{2}(\mu _{1}-\mu _{2})=%
\begin{array}{llllllll}
1 & 1 & 2 & 2 & 1 & 1 & 1 & 1%
\end{array},$

$\mu _{1}+\frac{2}{3}v=%
\begin{array}{llllllll}
0 & 1 & 1 & 0 & 0 & 0 & 0 & 1%
\end{array},$

$\mu _{2}+\frac{2}{3}v=%
\begin{array}{llllllll}
0 & 0 & 1 & 1 & 0 & 0 & 0 & 1%
\end{array},$

$\mu _{1}+\mu _{2}-\frac{2}{3}v=%
\begin{array}{llllllll}
0 & 1 & 1 & 1 & 0 & 0 & 0 & 0%
\end{array},$

$\lambda _{0}-\frac{1}{2}\mu _{1}+\frac{2}{3}v=%
\begin{array}{llllllll}
2 & 3 & 5 & 4 & 3 & 2 & 1 & 3%
\end{array},$

$\lambda _{1}-\frac{1}{2}\mu _{1}+\frac{2}{3}v=%
\begin{array}{llllllll}
0 & 0 & 1 & 1 & 1 & 1 & 1 & 1%
\end{array},$

$\lambda _{2}-\frac{1}{2}\mu _{1}+\frac{2}{3}v=%
\begin{array}{llllllll}
0 & 0 & 1 & 1 & 1 & 1 & 0 & 1%
\end{array},$

$\lambda _{3}-\frac{1}{2}\mu _{1}+\frac{2}{3}v=%
\begin{array}{llllllll}
0 & 0 & 1 & 1 & 1 & 0 & 0 & 1%
\end{array},$

$\lambda _{0}+\frac{1}{2}\mu _{1}-\frac{2}{3}v=%
\begin{array}{llllllll}
2 & 4 & 5 & 4 & 3 & 2 & 1 & 2%
\end{array},$

$\lambda _{1}+\frac{1}{2}\mu _{1}-\frac{2}{3}v=%
\begin{array}{llllllll}
0 & 1 & 1 & 1 & 1 & 1 & 1 & 0%
\end{array},$

$\lambda _{2}+\frac{1}{2}\mu _{1}-\frac{2}{3}v=%
\begin{array}{llllllll}
0 & 1 & 1 & 1 & 1 & 1 & 0 & 0%
\end{array},$

$\lambda _{3}+\frac{1}{2}\mu _{1}-\frac{2}{3}v=%
\begin{array}{llllllll}
0 & 1 & 1 & 1 & 1 & 0 & 0 & 0%
\end{array},$

$\frac{1}{2}(\lambda _{0}+\lambda
_{1}-\lambda _{2}+\lambda _{3})+\frac{1}{2}\mu _{2}-\frac{2}{3}v=%
\begin{array}{llllllll}
1 & 2 & 3 & 3 & 2 & 1 & 1 & 1%
\end{array},$

$\frac{1}{2}(\lambda _{0}+\lambda
_{1}+\lambda _{2}-\lambda _{3})+\frac{1}{2}\mu _{2}-\frac{2}{3}v=%
\begin{array}{llllllll}
1 & 2 & 3 & 3 & 2 & \text{2} & 1 & 1%
\end{array},$

$\frac{1}{2}(\lambda _{0}-\lambda
_{1}-\lambda _{2}-\lambda _{3})+\frac{1}{2}\mu _{2}-\frac{2}{3}v=%
\begin{array}{llllllll}
1 & 1 & 1 & 1 & 0 & \text{0} & 0 & 0%
\end{array}$

$\frac{1}{2}(\lambda _{0}-\lambda
_{1}+\lambda _{2}+\lambda _{3})+\frac{1}{2}\mu _{2}-\frac{2}{3}v=%
\begin{array}{llllllll}
1 & 2 & 3 & 3 & 2 & \text{1} & 0 & 1%
\end{array},$

$\frac{1}{2}(\lambda _{0}+\lambda
_{1}-\lambda _{2}+\lambda _{3})-\frac{1}{2}\mu _{2}+\frac{2}{3}v=%
\begin{array}{llllllll}
1 & 2 & 3 & 2 & 2 & \text{1} & 1 & 2%
\end{array},$

$\frac{1}{2}(\lambda _{0}+\lambda
_{1}+\lambda _{2}-\lambda _{3})-\frac{1}{2}\mu _{2}+\frac{2}{3}v=%
\begin{array}{llllllll}
1 & 2 & 3 & 2 & 2 & \text{2} & 1 & 2%
\end{array},$

$\frac{1}{2}(\lambda _{0}-\lambda
_{1}-\lambda _{2}-\lambda _{3})-\frac{1}{2}\mu _{2}+\frac{2}{3}v=%
\begin{array}{llllllll}
1 & 1 & 1 & 0 & 0 & \text{0} & 0 & 1%
\end{array},$

$\frac{1}{2}(\lambda _{0}-\lambda _{1}+\lambda
_{2}+\lambda _{3}))-\frac{1}{2}\mu _{2}+\frac{2}{3}v=%
\begin{array}{llllllll}
1 & 2 & 3 & 2 & 2 & \text{1} & 0 & 2%
\end{array},$

$\frac{1}{2}(\lambda _{0}-\lambda _{1}+\lambda _{2}-\lambda
_{3})-\frac{1}{2}(\mu _{1}+\mu _{2})-\frac{2}{3}v=%
\begin{array}{llllllll}
1 & 1 & 1 & 1 & 1 & \text{1} & 0 & 0%
\end{array},$

$\frac{1}{2}(\lambda _{0}-\lambda _{1}-\lambda _{2}+\lambda
_{3})-\frac{1}{2}(\mu _{1}+\mu _{2})-\frac{2}{3}v=%
\begin{array}{llllllll}
1 & 1 & 1 & 1 & 1 & \text{0} & 0 & 0%
\end{array},$

$\frac{1}{2}(\lambda _{0}+\lambda _{1}+\lambda _{2}+\lambda
_{3})-\frac{1}{2}(\mu _{1}+\mu _{2})-\frac{2}{3}v=%
\begin{array}{llllllll}
1 & 2 & 3 & 3 & 3 & 2 & 1 & 1%
\end{array},$

$\frac{1}{2}(\lambda _{0}+\lambda _{1}-\lambda _{2}-\lambda
_{3})-\frac{1}{2}(\mu _{1}+\mu _{2})-\frac{2}{3}v=%
\begin{array}{llllllll}
1 & 1 & 1 & 1 & 1 & \text{1} & 1 & 0%
\end{array},$

$\frac{1}{2}(\lambda _{0}-\lambda _{1}+\lambda _{2}-\lambda
_{3})+\frac{1}{2}(\mu _{1}+\mu _{2})+\frac{2}{3}v=%
\begin{array}{llllllll}
1 & 2 & 3 & 2 & 1 & 1 & 0 & 2%
\end{array},$

$\frac{1}{2}(\lambda _{0}-\lambda _{1}-\lambda _{2}+\lambda
_{3})+\frac{1}{2}(\mu _{1}+\mu _{2})+\frac{2}{3}v=%
\begin{array}{llllllll}
1 & 2 & 3 & 2 & 1 & 0 & 0 & 2%
\end{array},$

$\frac{1}{2}(\lambda _{0}+\lambda _{1}+\lambda _{2}+\lambda
_{3})+\frac{1}{2}(\mu _{1}+\mu _{2})+\frac{2}{3}v=%
\begin{array}{llllllll}
1 & 3 & 5 & 4 & 3 & 2 & 1 & 3%
\end{array},$

$\frac{1}{2}(\lambda _{0}+\lambda _{1}-\lambda _{2}-\lambda
_{3})+\frac{1}{2}(\mu _{1}+\mu _{2})+\frac{2}{3}v=%
\begin{array}{llllllll}
1 & 2 & 3 & 2 & 1 & 1 & 1 & 2%
\end{array},$

$\mu _{1}-\frac{1}{3}v+r=%
\begin{array}{llllllll}
0 & 1 & 1 & 0 & 0 & 0 & 0 & 0%
\end{array},$

$\mu _{2}-\frac{1}{3}v+r=%
\begin{array}{llllllll}
0 & 0 & 1 & 1 & 0 & 0 & 0 & 0%
\end{array},$

$\mu _{1}+\mu _{2}+\frac{1}{3}v+r=%
\begin{array}{llllllll}
0 & 1 & 2 & 1 & 0 & 0 & 0 & 0%
\end{array},$

$\mu _{1}-\frac{1}{3}v-r=%
\begin{array}{llllllll}
0 & 1 & 0 & 0 & 0 & 0 & 0 & 0%
\end{array},$

$\mu _{2}-\frac{1}{3}v-r=%
\begin{array}{llllllll}
0 & 0 & 0 & 1 & 0 & 0 & 0 & 0%
\end{array},$

$\mu _{1}+\mu _{2}+\frac{1}{3}v+r=%
\begin{array}{llllllll}
0 & 1 & 1 & 1 & 0 & 0 & 0 & 1%
\end{array},$

$\lambda _{0}+\frac{1}{2}\mu _{1}+\frac{1}{3}v+r=%
\begin{array}{llllllll}
2 & 4 & 6 & 4 & 3 & 2 & 1 & 3%
\end{array},$

$\lambda _{0}+\frac{1}{2}\mu _{1}+\frac{1}{3}v-r=%
\begin{array}{llllllll}
2 & 4 & 5 & 4 & 3 & 2 & 1 & 3%
\end{array},$

$\lambda _{0}-\frac{1}{2}\mu _{1}-\frac{1}{3}v+r=%
\begin{array}{llllllll}
2 & 3 & 5 & 4 & 3 & 2 & 1 & 2%
\end{array},$

$\lambda _{0}-\frac{1}{2}\mu _{1}-\frac{1}{3}v-r=%
\begin{array}{llllllll}
2 & 3 & 4 & 4 & 3 & 2 & 1 & 2%
\end{array},$

$\lambda _{1}+\frac{1}{2}\mu _{1}+\frac{1}{3}v+r=%
\begin{array}{llllllll}
0 & 1 & 2 & 1 & 1 & 1 & 1 & 1%
\end{array},$

$\lambda _{1}+\frac{1}{2}\mu _{1}+\frac{1}{3}v-r=%
\begin{array}{llllllll}
0 & 1 & 1 & 1 & 1 & 1 & 1 & 1%
\end{array},$

$\lambda _{1}-\frac{1}{2}\mu _{1}-\frac{1}{3}v+r=%
\begin{array}{llllllll}
0 & 0 & 1 & 1 & 1 & 1 & 1 & 0%
\end{array},$

$\lambda _{1}-\frac{1}{2}\mu _{1}-\frac{1}{3}v-r=%
\begin{array}{llllllll}
0 & 0 & 0 & 1 & 1 & 1 & 1 & 0%
\end{array},$

$\lambda _{2}+\frac{1}{2}\mu _{1}+\frac{1}{3}v+r=%
\begin{array}{llllllll}
0 & 1 & 2 & 1 & 1 & 1 & 0 & 1%
\end{array},$

$\lambda _{2}+\frac{1}{2}\mu _{1}+\frac{1}{3}v-r=%
\begin{array}{llllllll}
0 & 1 & 1 & 1 & 1 & 1 & 0 & 1%
\end{array},$

$\lambda _{2}-\frac{1}{2}\mu _{1}-\frac{1}{3}v+r=%
\begin{array}{llllllll}
0 & 0 & 1 & 1 & 1 & 1 & 0 & 0%
\end{array},$

$\lambda _{2}-\frac{1}{2}\mu _{1}-\frac{1}{3}v-r=%
\begin{array}{llllllll}
0 & 0 & 0 & 1 & 1 & 1 & 0 & 0%
\end{array},$

$\lambda _{3}+\frac{1}{2}\mu _{1}+\frac{1}{3}v+r=%
\begin{array}{llllllll}
0 & 1 & 2 & 1 & 1 & 0 & 0 & 1%
\end{array},$

$\lambda _{3}+\frac{1}{2}\mu _{1}+\frac{1}{3}v-r=%
\begin{array}{llllllll}
0 & 1 & 1 & 1 & 1 & 0 & 0 & 1%
\end{array},$

$\lambda _{3}-\frac{1}{2}\mu _{1}-\frac{1}{3}v+r=%
\begin{array}{llllllll}
0 & 0 & 1 & 1 & 1 & 0 & 0 & 0%
\end{array},$

$\lambda _{3}-\frac{1}{2}\mu _{1}-\frac{1}{3}v-r=%
\begin{array}{llllllll}
0 & 0 & 0 & 1 & 1 & 0 & 0 & 0%
\end{array},$

$\frac{1}{2}(\lambda _{0}+\lambda _{1}-\lambda _{2}+\lambda
_{3})+(\frac{1}{2}\mu _{2}+\frac{1}{3}\nu )+r=%
\begin{array}{llllllll}
1 & 2 & 4 & 3 & 2 & 1 & 1 & 2%
\end{array},$

$\frac{1}{2}(\lambda _{0}+\lambda _{1}-\lambda _{2}+\lambda
_{3})+(\frac{1}{2}\mu _{2}+\frac{1}{3}\nu )-r=%
\begin{array}{llllllll}
1 & 2 & 3 & 3 & 2 & 1 & 1 & 2%
\end{array},$

$\frac{1}{2}(\lambda _{0}+\lambda _{1}-\lambda _{2}+\lambda
_{3})-(\frac{1}{2}\mu _{2}+\frac{1}{3}\nu )+r=%
\begin{array}{llllllll}
1 & 2 & 3 & 2 & 2 & 1 & 1 & 1%
\end{array},$

$\frac{1}{2}(\lambda _{0}+\lambda _{1}-\lambda _{2}+\lambda
_{3})-(\frac{1}{2}\mu _{2}+\frac{1}{3}\nu )-r=%
\begin{array}{llllllll}
1 & 2 & 2 & 2 & 2 & 1 & 1 & 1%
\end{array},$

$\frac{1}{2}(\lambda _{0}+\lambda _{1}+\lambda _{2}-\lambda
_{3})+(\frac{1}{2}\mu _{2}+\frac{1}{3}\nu )+r=%
\begin{array}{llllllll}
1 & 2 & 4 & 3 & 2 & 2 & 1 & 2%
\end{array},$

$\frac{1}{2}(\lambda _{0}+\lambda _{1}+\lambda _{2}-\lambda
_{3})+(\frac{1}{2}\mu _{2}+\frac{1}{3}\nu )-r=%
\begin{array}{llllllll}
1 & 2 & 3 & 3 & 2 & 2 & 1 & 2%
\end{array},$

$\frac{1}{2}(\lambda _{0}+\lambda _{1}+\lambda _{2}-\lambda
_{3})-(\frac{1}{2}\mu _{2}+\frac{1}{3}\nu )+r=%
\begin{array}{llllllll}
1 & 2 & 3 & 2 & 2 & 2 & 1 & 1%
\end{array},$

$\ \ \ \ \ \ \ \frac{1}{2}(\lambda _{0}+\lambda _{1}+\lambda _{2}-\lambda
_{3})-(\frac{1}{2}\mu _{2}+\frac{1}{3}\nu )-r=%
\begin{array}{llllllll}
1 & 2 & 2 & 2 & 2 & 2 & 1 & 1%
\end{array}%
,$

$\frac{1}{2}(\lambda _{0}-\lambda _{1}-\lambda _{2}-\lambda
_{3})+(\frac{1}{2}\mu _{2}+\frac{1}{3}\nu )+r=%
\begin{array}{llllllll}
1 & 1 & 2 & 1 & 0 & 0 & 0 & 1%
\end{array},$

$\frac{1}{2}(\lambda _{0}-\lambda _{1}-\lambda _{2}-\lambda
_{3})+(\frac{1}{2}\mu _{2}+\frac{1}{3}\nu )-r=%
\begin{array}{llllllll}
1 & 1 & 1 & 1 & 0 & 0 & 0 & 1%
\end{array},$

$\frac{1}{2}(\lambda _{0}-\lambda _{1}-\lambda _{2}-\lambda
_{3})-(\frac{1}{2}\mu _{2}+\frac{1}{3}\nu )+r=%
\begin{array}{llllllll}
1 & 1 & 1 & 0 & 0 & 0 & 0 & 0%
\end{array},$

$\frac{1}{2}(\lambda _{0}-\lambda _{1}-\lambda _{2}-\lambda
_{3})-(\frac{1}{2}\mu _{2}+\frac{1}{3}\nu )-r=%
\begin{array}{llllllll}
1 & 1 & 0 & 0 & 0 & 0 & 0 & 0%
\end{array},$

$\frac{1}{2}(\lambda _{0}-\lambda _{1}+\lambda _{2}+\lambda
_{3})+(\frac{1}{2}\mu _{2}+\frac{1}{3}\nu )+r=%
\begin{array}{llllllll}
1 & 2 & 4 & 3 & 2 & 1 & 0 & 2%
\end{array},$

$\frac{1}{2}(\lambda _{0}-\lambda _{1}+\lambda _{2}+\lambda
_{3})+(\frac{1}{2}\mu _{2}+\frac{1}{3}\nu )-r=%
\begin{array}{llllllll}
1 & 2 & 3 & 3 & 2 & 1 & 0 & 2%
\end{array},$

$\frac{1}{2}(\lambda _{0}-\lambda _{1}+\lambda _{2}+\lambda
_{3})-(\frac{1}{2}\mu _{2}+\frac{1}{3}\nu )+r=%
\begin{array}{llllllll}
1 & 2 & 3 & 2 & 2 & 1 & 0 & 1%
\end{array},$

$\frac{1}{2}(\lambda _{0}-\lambda _{1}+\lambda _{2}+\lambda
_{3})-(\frac{1}{2}\mu _{2}+\frac{1}{3}\nu )-r=%
\begin{array}{llllllll}
1 & 2 & 2 & 2 & 2 & 1 & 0 & 1%
\end{array},$

$\frac{1}{2}(\lambda _{0}-\lambda _{1}+\lambda _{2}-\lambda _{3})+(-\frac{1%
}{2}\mu _{1}-\frac{1}{2}\mu _{2}+\frac{1}{3}\nu )+r=%
\begin{array}{llllllll}
1 & 1 & 2 & 1 & 1 & 1 & 0 & 1%
\end{array},$

$\frac{1}{2}(\lambda _{0}-\lambda _{1}+\lambda _{2}-\lambda _{3})+(-\frac{1%
}{2}\mu _{1}-\frac{1}{2}\mu _{2}+\frac{1}{3}\nu )-r=%
\begin{array}{llllllll}
1 & 1 & 1 & 1 & 1 & 1 & 0 & 1%
\end{array},$

$\frac{1}{2}(\lambda _{0}-\lambda _{1}+\lambda _{2}-\lambda _{3})-(-\frac{1%
}{2}\mu _{1}-\frac{1}{2}\mu _{2}+\frac{1}{3}\nu )+r=%
\begin{array}{llllllll}
1 & 2 & 3 & 2 & 1 & 1 & 0 & 1%
\end{array},$

$\frac{1}{2}(\lambda _{0}-\lambda _{1}+\lambda _{2}-\lambda _{3})-(-\frac{1%
}{2}\mu _{1}-\frac{1}{2}\mu _{2}+\frac{1}{3}\nu )-r=%
\begin{array}{llllllll}
1 & 2 & 2 & 2 & 1 & 1 & 0 & 1%
\end{array},$

$\frac{1}{2}(\lambda _{0}-\lambda _{1}-\lambda _{2}+\lambda _{3})+(-\frac{1%
}{2}\mu _{1}-\frac{1}{2}\mu _{2}+\frac{1}{3}\nu )+r=%
\begin{array}{llllllll}
1 & 1 & 2 & 1 & 1 & 0 & 0 & 1%
\end{array},$

$\frac{1}{2}(\lambda _{0}-\lambda _{1}-\lambda _{2}+\lambda _{3})+(-\frac{1%
}{2}\mu _{1}-\frac{1}{2}\mu _{2}+\frac{1}{3}\nu )-r=%
\begin{array}{llllllll}
1 & 1 & 1 & 1 & 1 & 0 & 0 & 1%
\end{array},$

$\frac{1}{2}(\lambda _{0}-\lambda _{1}-\lambda _{2}+\lambda _{3})-(-\frac{1%
}{2}\mu _{1}-\frac{1}{2}\mu _{2}+\frac{1}{3}\nu )+r=%
\begin{array}{llllllll}
1 & 2 & 3 & 2 & 1 & 0 & 0 & 1%
\end{array},$

$\frac{1}{2}(\lambda _{0}-\lambda _{1}-\lambda _{2}+\lambda _{3})-(-\frac{1%
}{2}\mu _{1}-\frac{1}{2}\mu _{2}+\frac{1}{3}\nu )-r=%
\begin{array}{llllllll}
1 & 2 & 2 & 2 & 1 & 0 & 0 & 1%
\end{array},$

$\frac{1}{2}(\lambda _{0}+\lambda _{1}+\lambda _{2}+\lambda _{3})+(-\frac{1%
}{2}\mu _{1}-\frac{1}{2}\mu _{2}+\frac{1}{3}\nu )+r=%
\begin{array}{llllllll}
1 & 2 & 4 & 3 & 3 & 2 & 1 & 2%
\end{array},$

$\frac{1}{2}(\lambda _{0}+\lambda _{1}+\lambda _{2}+\lambda _{3})+(-\frac{1%
}{2}\mu _{1}-\frac{1}{2}\mu _{2}+\frac{1}{3}\nu )-r=%
\begin{array}{llllllll}
1 & 2 & 3 & 3 & 3 & 2 & 1 & 2%
\end{array},$

$\frac{1}{2}(\lambda _{0}+\lambda _{1}+\lambda _{2}+\lambda _{3})-(-\frac{1%
}{2}\mu _{1}-\frac{1}{2}\mu _{2}+\frac{1}{3}\nu )+r=%
\begin{array}{llllllll}
1 & 3 & 5 & 4 & 3 & 2 & 1 & 2%
\end{array},$

$\frac{1}{2}(\lambda _{0}+\lambda _{1}+\lambda _{2}+\lambda _{3})-(-\frac{1%
}{2}\mu _{1}-\frac{1}{2}\mu _{2}+\frac{1}{3}\nu )-r=%
\begin{array}{llllllll}
1 & 3 & 4 & 4 & 3 & 2 & 1 & 2%
\end{array},$

$ \frac{1}{2}(\lambda _{0}+\lambda _{1}-\lambda _{2}-\lambda _{3})+(-\frac{1%
}{2}\mu _{1}-\frac{1}{2}\mu _{2}+\frac{1}{3}\nu )+r=%
\begin{array}{llllllll}
1 & 1 & 2 & 1 & 1 & 1 & 1 & 1%
\end{array},$

$\frac{1}{2}(\lambda _{0}+\lambda _{1}-\lambda _{2}-\lambda _{3})+(-\frac{1%
}{2}\mu _{1}-\frac{1}{2}\mu _{2}+\frac{1}{3}\nu )-r=%
\begin{array}{llllllll}
1 & 1 & 1 & 1 & 1 & 1 & 1 & 1%
\end{array},$

$ \frac{1}{2}(\lambda _{0}+\lambda _{1}-\lambda _{2}-\lambda _{3})-(-\frac{1%
}{2}\mu _{1}-\frac{1}{2}\mu _{2}+\frac{1}{3}\nu )+r=%
\begin{array}{llllllll}
1 & 2 & 3 & 2 & 1 & 1 & 1 & 1%
\end{array},$

$ \frac{1}{2}(\lambda _{0}+\lambda _{1}-\lambda _{2}-\lambda _{3})-(-\frac{1%
}{2}\mu _{1}-\frac{1}{2}\mu _{2}+\frac{1}{3}\nu )-r=%
\begin{array}{llllllll}
1 & 2 & 2 & 2 & 1 & 1 & 1 & 1%
\end{array},$

$2r$ \ \ \ \ $=%
\begin{array}{llllllll}
0 & 0 & 1 & 0 & 0 & 0 & 0 & 0%
\end{array},$

$v+r=%
\begin{array}{llllllll}
0 & 0 & 1 & 0 & 0 & 0 & 0 & 1%
\end{array},$

$\ v-r=%
\begin{array}{llllllll}
0 & 0 & 0 & 0 & 0 & 0 & 0 & 1%
\end{array}.$
\end{flushright}

\noindent
Hence $\Pi =\{\alpha _{1},\alpha _{2},\alpha _{3},\alpha _{4},\alpha
_{5},\alpha _{6},\alpha _{7},\alpha _{8}\}$ is a fundamental root system of 
\gr$_{8}^{\C}$ . The real part of \gh$_{R}$ of \gh\  is

\gh$_{R}=\left\{ h=h_{\delta }+H+V+R\in \text{\gr}_{8}^{\C}%
\middle| %
\begin{array}{c}
h_{\delta }=\sum\limits_{k=0}^{3}\lambda
_{k}H_{k}=\sum\limits_{k=0}^{3}-\lambda _{k}iUd_{k4+k},\\
H=\mu _{1}U\tau _{1}+\mu_{2}U\tau _{2}, \\ 
V=\nu U\rho _{1},R=rUr_{1},\\
\lambda _{k},\mu _{j},\nu ,r\in \R%
\end{array}%
\right\} .$

The Killing form $B_{8}$ of \gr$_{8}^{\C}$ is $B_{8}(R_{1},R_{2})=tr(R_{1}R_{2})$ (%
\emph{Lemma 17.33}), so that the Killing form $B_{8}$ of \gr$_{8}^{\C}$ on \gh$_{R}$ is given by

$B_{8}(h,h^{\prime })=60\sum\limits_{k=0}^{3}\lambda
_{k}\lambda _{k}^{\prime }+60\mu _{1}\mu _{1}^{\prime }+30\mu _{1}\mu
_{2}^{\prime }+30\mu _{2}\mu _{1}^{\prime }+60\mu _{2}\mu _{2}^{\prime
}+40vv^{\prime }+120rr^{\prime },$

\noindent
for $h=\sum\limits_{k=0}^{3}\lambda _{k}H_{k}+\mu _{1}U\tau _{1}+\mu _{2}U\tau
_{2}+vU\rho _{1}+rUr_{1},$

$h^{\prime }=\sum\limits_{k=0}^{3}\lambda _{k}^{\prime }H_{k}+\mu
_{1}^{\prime }U\tau _{1}+\mu _{2}^{\prime }U\tau _{2}+v^{\prime }U\rho _{1}+r^{\prime
}Ur_{1}\in $\gh$_{R}.$

\noindent
Indeed,by calculate with Maxima we have the above .

Now,the canonical elements $H\alpha _{i}\in $\gh$_{R}$ corresponding
to $\alpha _{i}$ $(B_{8}(H\alpha ,H)=\alpha (H),H\in $\gh$_{R}$%
\textbf{) }are determined as follows.

$\ \ \ \ \ \ \ \ \ \ \ \ \ \ \ \ \ \ H_{\alpha _{1}}=\frac{1}{120}%
(H_{0}-H_{1}-H_{2}-H_{3})-\frac{1}{60}U\tau _{1},$

$\ \ \ \ \ \ \ \ \ \ \ \ \ \ \ \ \ \ H_{\alpha _{2}}=\frac{1}{90}(2U\tau _{1}-U\tau
_{2})-\frac{1}{120}U\rho _{1}-\frac{1}{120}Ur_{1},$

$\ \ \ \ \ \ \ \ \ \ \ \ \ \ \ \ \ \ H_{\alpha _{3}}=\frac{1}{60}Ur_{1},$

$\ \ \ \ \ \ \ \ \ \ \ \ \ \ \ \ \ \ H_{\alpha 4}=\frac{1}{90}(-U\tau
_{1}+2U\tau _{2})-\frac{1}{120}U\rho _{1}-\frac{1}{120}Ur_{1},$

\ \ \ \ \ \ \ \ \ \ \ \ \ \ \ \ \ \ $H_{\alpha 5}=\frac{1}{60}H_{3}-\frac{1}{%
60}U\tau _{2},$

\ \ \ \ \ \ \ \ \ \ \ \ \ \ \ \ \ \ $H_{\alpha 6}=\frac{1}{60}(H_{2}-H_{3}),$

\ \ \ \ \ \ \ \ \ \ \ \ \ \ \ \ \ \ $H_{\alpha 7}=\frac{1}{60}(H_{1}-H_{2}),$

$\ \ \ \ \ \ \ \ \ \ \ \ \ \ \ \ \ \ H_{\alpha 8}=\frac{1}{40}U\rho _{1}-\frac{1%
}{120}Ur_{1},$

\noindent
Thus we have

\ \ \ \ \ \ \ \ \ \ \ \ \ $(\alpha _{1},\alpha _{1})=B_{8}(H_{\alpha
1},H_{\alpha _{1}})=\frac{1}{30},$

\ \ \ \ \ \ \ \ \ \ \ \ \ $\ (\alpha _{i},\alpha _{i})=\frac{1}{30}%
,(i=2,3,4,5,6,7,8),$

\ \ \ \ \ \ \ \ \ \ \ \ $\ (\alpha _{1},\alpha _{2})=(\alpha _{2},\alpha
_{3})=(\alpha _{3},\alpha _{4})=(\alpha _{4},\alpha _{5})=(\alpha
_{5},\alpha _{6})=(\alpha _{6},\alpha _{7})$

\ \ \ \ \ \ \ \ \ \ \ \ \ \ \ \ \ \ \ \ \ \ \ \ \ $=(\alpha _{3},\alpha _{8})=-\frac{1}{60},$

\ \ \ \ \ \ \ \ \ \ \ \ \ $\ (\alpha _{i},\alpha _{j})=0$, otherwise,

\ \ \ \ \ \ \ \ \ \ \ $\ (-\mu ,-\mu )=\frac{1}{30},(-\mu ,\alpha
_{7})=-\frac{1}{60},(-\mu ,\alpha _{i})=0,(i=1,2,3,4,5,6,8).$

\noindent
Using them,we can draw the Dynkin diagram and the extended Dynkin 
diagram of  \gr$_{8}^{\C}.$\ \ \ \ \emph{Q.E.D.}

\bigskip

\emph{Definition 17.38.} \ We define followings,

\ge\gx$^{0}=\{c_{k}^{0}(U\chi _{k}+U\psi_{k})+ic_{k}^{1}(U\chi _{k}-U\psi
_{k}) \mid c_{k}^{0},c_{k}^{1}\in \R,1\leq k\leq 3\},$

\ \ \ \ \ \ $=\{\chi_{k}U\chi_{k}+\iota\chi_{k}U\psi_{k} \mid \chi_{k}\in \C,1\leq k\leq 3\},$

\ \ \ \ \ \ \ \ where $\iota $ is the complex conjugation of $\C$,

\ge\gx$^{1}=%
\{c_{kj}^{0}(Ux_{kj}+Uw_{kj})+ic_{kj}^{1}(Ux_{kj}-Uw_{kj}) \mid c_{kj}^{0},c_{kj}^{1}\in \R,$

\ \ \ \ \ \ \ \ \ \ \ \ \ \ \ \ \ \ \ \ \ \ \ \ \ \ \ \ \ \ \ \ \ \ \ \ \ \ \ \ \ \ \ \ \ \ \ \ \ \ $1\leq k\leq 3,0\leq j\leq 7\}$,

\ \ \ \ \ \ $=\{x_{kj}Ux_{kj}+\iota x_{kj}Uw_{kj} \mid x_{kj}\in \C,1\leq k\leq 3,0\leq j\leq 7\}$,

\ge\gx$=$\ge\gx$^{0}\oplus $\ge\gx$^{1},$

\ge\gy$^{0}=\{c_{k}^{0}(U\gamma _{k}-U\mu_{k})+ic_{k}^{1}(U\gamma _{k}+U\mu
_{k}) \mid c_{k}^{0},c_{k}^{1}\in \R,1\leq k\leq 3\},$

\ \ \ \ \ \ $=\{\gamma_{k}U\gamma_{k}-\iota\gamma_{k}U\mu_{k} \mid \gamma_{k}\in \C,1\leq k\leq 3\},$

\ge\gy$^{1}=%
\{c_{kj}^{0}(Uy_{kj}-Uz_{kj})+ic_{kj}^{1}(Uy_{kj}+Uz_{kj}) \mid c_{kj}^{0},c_{kj}^{1}\in \R,$

\ \ \ \ \ \ \ \ \ \ \ \ \ \ \ \ \ \ \ \ \ \ \ \ \ \ \ \ \ \ \ \ \ \ \ \ \ \ \ \ \ \ \ \ \ \ \ \ \ \ $1\leq k\leq 3,0\leq j\leq 7\}$,

\ \ \ \ \ \ $=\{y_{kj}Uy_{kj}-\iota y_{kj}Uz_{kj} \mid y_{kj}\in \C,1\leq k\leq 3,0\leq j\leq 7\}$,

\ge\gy$=$\ge\gy$^{0}\oplus $\ge\gy$^{1},$

\ge\gk$=\{c^{0}(U\xi _{1}+U\omega_{1})+ic^{1}(U\xi _{1}-U\omega
_{1}) \mid c^{0},c^{1}\in \R\},$

 \ \ \ \ $=\{\xi_{1}U\xi_{1}+\iota\xi_{1}U\omega_{1} \mid \xi_{1}\in \C\},$
 
\ge\gi$=\{c^{0}(U\eta _{1}-U\zeta_{1})+ic^{1}(U\xi _{1}+U\zeta
_{1}) \mid c^{0},c^{1}\in \R\},$

 \ \ \ \ $=\{\eta_{1}U\eta_{1}-\iota\eta_{1}U\zeta_{1} \mid \eta_{1}\in \C\},$

\ge\gr$=\{r_{1}Ur_{1} \mid r_{1}\in i\R\},$

\ge\gs$=\{c^{0}(Us_{1}-Uu_{1})+ic^{1}(Us_{1}+Uu
_{1}) \mid c^{0},c^{1}\in \R\},$

 \ \ \ \ $=\{s_{1}Us_{1}-\iota s_{1}Uu_{1} \mid s_{1}\in \C\},$

\bigskip

\emph{Theorem 17.39.} Let we put

\ge$_{8}=$\ge\gd$\oplus $\ge\gm$\oplus $\ge\gt$\oplus $\ge\ga$\oplus $\ge\gp$ \oplus $\ge\gx$\oplus $\ge\gy$\oplus
$\ge\gk$ \oplus $\ge\gi$ \oplus $\ge\gr$\oplus $\ge\gs. \ 

\noindent
Then \ge$_{8}$ is a compact exceptional simple Lie algebra of type $E_{8}.$

\bigskip

\emph{Proof.}\ By from \emph{Lemma 17.2} to \emph{Lemma 17.33}, \ge$_{8}$ is a Lie algebra. 
And by calculations using Maxima,we have as follows.

\noindent
For $X_{1}\in $\ge\gd$,X_{2}\in $\ge\gm$,X_{3}\in $\ge\gt$^{0},X_{4}\in $\ge\gt$^{1},X_{5}\in $\ge\ga$^{0},X_{6}\in $\ge\ga$^{1},$

\noindent
$X_{7}\in $\ge\gp$ ,X_{8}\in $\ge\gx$^{0},X_{9}\in $\ge\gx$^{1},X_{10}\in $\ge\gy$^{0},X_{11}\in $\ge\gy$^{1},X_{12}\in $\ge\gk$ ,$

\noindent
$X_{13}\in $\ge\gi$ ,X_{14}\in $\ge\gr$,X_{15}\in $\ge\gs,

$B_{8}(X_{1},X_{1})$

$=-60(d_{67}^{2}+d_{57}^{2}+d_{56}^{2}+d_{47}^{2}+d_{46}%
^{2}+d_{45}^{2}+d_{37}^{2}$

\ \ \ \ \ \ $\ +d_{36}^{2}+d_{35}^{2}+d_{34}^{2}+d_{27}^{2}+d_{26}%
^{2}+d_{25}^{2}+d_{24}^{2}$

\ \ \ \ \ \ $\ +d_{23}^{2}+d_{17}^{2}+d_{16}^{2}+d_{15}^{2}+d_{14}%
^{2}+d_{13}^{2}+d_{12}^{2}$

\ \ \ \ \ \ $\ +d_{07}^{2}+d_{06}^{2}+d_{05}^{2}+d_{04}^{2}+d_{03}%
^{2}+d_{02}^{2}+d_{01}^{2}),$

$B_{8}(X_{2},X_{2})$

$=-60(m_{37}^{2}+m_{36}^{2}+m_{35}^{2}+m_{34}^{2}+m_{33}%
^{2}+m_{32}^{2}+m_{31}^{2}+m_{30}^{2}$

\ \ \ \ $\ +m_{27}^{2}+m_{26}^{2}+m_{25}^{2}+m_{24}^{2}+m_{23}%
^{2}+m_{22}^{2}+m_{21}^{2}+m_{20}^{2}$

\ \ \ \ $\ +m_{17}^{2}+m_{16}^{2}+m_{15}^{2}+m_{14}^{2}+m_{13}%
^{2}+m_{12}^{2}+m_{11}^{2}+m_{10}^{2}),$

$B_{8}(X_{3},X_{3})=-60(\tau _{1}^{2}+3\tau _{2}^{2})$

$B_{8}(X_{4},X_{4})$

$=-60(t_{37}^{2}+t_{36}^{2}+t_{35}^{2}+t_{34}^{2}+t_{33}%
^{2}+t_{32}^{2}+t_{31}^{2}+t_{30}^{2}$

\ \ \ \ \ \ $\ +t_{27}^{2}+t_{26}^{2}+t_{25}^{2}+t_{24}^{2}+t_{23}%
^{2}+t_{22}^{2}+t_{21}^{2}+t_{20}^{2}$

\ \ \ \ \ $\ \ $+$t_{17}^{2}+t_{16}^{2}+t_{15}^{2}+t_{14}^{2}+t_{13}%
^{2}+t_{12}^{2}+t_{11}^{2}+t_{10}^{2}),$

$B_{8}(X_{5},X_{5})=-120(\iota \alpha _{3}\alpha _{3}+\iota \alpha _{2}\alpha _{2}+\iota
\alpha _{1}\alpha _{1}),$

$B_{8}(X_{6},X_{6})$

=$-240(\iota a_{37}a_{37}+\iota a_{36}a_{36}+\iota a_{35}a_{35}+\iota a_{34}a_{34}+\iota
a_{33}a_{33}+\iota a_{32}a_{32}+\iota a_{31}a_{31}$

\ \ \ \ $+\iota a_{30}a_{30}$+$\iota a_{27}a_{27}+\iota a_{26}a_{26}+\iota a_{25}a_{25}+\iota a_{24}a_{24}+\iota
a_{23}a_{23}+\iota a_{22}a_{22}$

\ \ \ \ $+\iota a_{21}a_{21}+\iota a_{20}a_{20}$+$\iota a_{17}a_{17}+\iota a_{16}a_{16}+\iota a_{15}a_{15}+\iota
a_{14}a_{14}+\iota a_{13}a_{13}$

\ \ \ \ $+\iota a_{12}a_{12}+\iota a_{11}a_{11}+\iota a_{10}a_{10}),$

$B_{8}(X_{7},X_{7})=-40\rho _{1}^{2},$

$B_{8}(X_{8},X_{8})=-30(\iota \chi _{3}\chi _{3}+\iota \chi _{2}\chi _{2}+\iota \chi _{1}\chi _{1}),$

$B_{8}(X_{9},X_{9})$

$=-60(\iota x_{37}x_{37}+\iota x_{36}x_{36}+\iota x_{35}x_{35}+\iota x_{34}x_{34}+\iota x_{33}x_{33}+\iota
x_{32}x_{32}+\iota x_{31}x_{31}$

\ \ \ \ $+\iota x_{30}x_{30}+\iota x_{27}x_{27}+\iota x_{26}x_{26}+\iota x_{25}x_{25}+\iota x_{24}x_{24}+\iota
x_{23}x_{23}+\iota x_{22}x_{22}$

\ \ \ \ $+\iota x_{21}x_{21}+\iota x_{20}x_{20}+\iota x_{17}x_{17}+\iota x_{16}x_{16}+\iota x_{15}x_{15}+\iota
x_{14}x_{14}+\iota x_{13}x_{13}$

\ \ \ \ $+\iota x_{12}x_{12}+\iota x_{11}x_{11}+\iota x_{10}x_{10}),$

$B_{8}(X_{10},X_{10})=-30(\iota \gamma _{3}\gamma _{3}+\iota \gamma _{2}\gamma _{2}+\iota
\gamma _{1}\gamma _{1}),$

$B_{8}(X_{11},X_{11})$

$=-60(\iota y_{37}y_{37}+\iota y_{36}y_{36}+\iota y_{35}y_{35}+\iota y_{34}y_{34}+\iota y_{33}y_{33}+\iota
y_{32}y_{32}+\iota y_{31}y_{31}$

\ \ \ \ $+\iota y_{30}y_{30}+\iota y_{27}y_{27}+\iota y_{26}y_{26}+\iota y_{25}y_{25}+\iota y_{24}y_{24}+\iota
y_{23}y_{23}+\iota y_{22}y_{22}$

\ \ \ \ $+\iota y_{21}y_{21}+\iota y_{20}y_{20}+\iota y_{17}y_{17}+\iota y_{16}y_{16}+\iota y_{15}y_{15}+\iota
y_{14}y_{14}+\iota y_{13}y_{13}$

\ \ \ \ $+\iota y_{12}y_{12}+\iota y_{11}y_{11}+\iota y_{10}y_{10}),$

$B_{8}(X_{12},X_{12})=-30\iota \xi _{1} \xi _{1},$

$B_{8}(X_{13},X_{13})=-30\iota \eta _{1} \eta _{1},$

$B_{8}(X_{14},X_{14})=-120r_{1}^{2},$

$B_{8}(X_{15},X_{15})=-120\iota s_{1} s_{1},$

$B_{8}(X_{i},X_{j})=0,(i\neq j).$

\noindent
Therfor we have $B_{8}(X,X)<0$,for $^{\forall }X\neq 0,X\in $\ge$_{8}.$ Then \ge$_{8}$ is compact.\ \ \ $%
Q.E.D.$

\bigskip

\emph{Remark 17.40.}
For $R \in e_{8}$, $^{t}\overline R+R=0$ does not hold where $^{t}R$ means \emph{matrix} transpose and $\overline{R}$ means complex conjugate.

\bigskip 

\section{The exceptional simple Lie algebra \gr$_{2}^{\C}$ of type $G_{2}$}

\bigskip

\ \ \ \ \emph{Proposition 18.1. }(I.Yokota\cite[\emph{Theorem1.9.3.}]{Yokota1})

$G_{2} = \{ \alpha \in Iso_{\R}($\gC$) \mid \alpha (xy)=(\alpha x)(\alpha y) \}$ is a simply connected
compact Lie group.

\bigskip

\emph{Proposition 18.2. } (I.Yokota\cite[\emph{Theorem1.4.1.}]{Yokota1})

The Lie algebra \gg$_{2}$ of the Lie group $G_{2}$ is given by

\gg$_{2}=\{D \in Hom_{\R}$(\gC$) \mid D(xy)=(Dx)y+x(Dy)\}$.

\bigskip

\emph{Definition 18.3.} \ We define the following elements and a Lie subalgebra of \gr\gd$.$

\ \ \ \ \ $S_{1}=2Ud_{12}-Ud_{47}-Ud_{56}$,

\ \ \ \ \ $S_{2}=2Ud_{13}-Ud_{46}+Ud_{57}$,

\ \ \ \ \ $S_{3}=2Ud_{14}+Ud_{27}+Ud_{36}$,

\ \ \ \ \ $S_{4}=2Ud_{15}+Ud_{26}-Ud_{37}$,

\ \ \ \ \ $S_{5}=2Ud_{16}-Ud_{25}-Ud_{34}$,

\ \ \ \ \ $S_{6}=2Ud_{17}-Ud_{24}+Ud_{35}$,

\ \ \ \ \ $S_{7}=2Ud_{23}-Ud_{45}-Ud_{67}(=-2L_{1}-L_{6})$,

\ \ \ \ \ $S_{11}=Ud_{45}+Ud_{67}$,

\ \ \ \ \ $S_{12}=-Ud_{46}+Ud_{57}$,

\ \ \ \ \ $S_{13}=Ud_{47}+Ud_{56}$,

\ \ \ \ \ $S_{21}=Ud_{45}-Ud_{67}$

\ \ \ \ \ $S_{22}=-Ud_{46}-Ud_{57}$

\ \ \ \ \ $S_{23}=-Ud_{47}+Ud_{56}$

\ \ \ \ \ $L_{1}=-Ud_{23}+Ud_{45}$,

\ \ \ \ \ $L_{2}=Ud_{24}+Ud_{35}$,

\ \ \ \ \ $L_{3}=-Ud_{25}+Ud_{34}$,

\ \ \ \ \ $L_{4}=Ud_{26}+Ud_{37}$,

\ \ \ \ \ $L_{5}=-Ud_{27}+Ud_{36}$,

\ \ \ \ \ $L_{6}=-Ud_{45}+Ud_{67}$,

\ \ \ \ \ $L_{7}=Ud_{46}+Ud_{57}$,

\ \ \ \ \ $L_{8}=-Ud_{47}+Ud_{56}$,

\ \ \ \ \ $K_{1}=-Ud_{01}-Ud_{23}$,

\ \ \ \ \ $K_{2}=-Ud_{02}+Ud_{13}$,

\ \ \ \ \ $K_{3}=-Ud_{03}-Ud_{12}$,

\ \ \ \ \ $K_{4}=Ud_{01}-Ud_{23}$,

\ \ \ \ \ $K_{5}=Ud_{02}+Ud_{13}$,

\ \ \ \ \ $K_{6}=Ud_{03}-Ud_{12}$,

\ \ \ \ \ $K_{11}=Ud_{01}+Ud_{23}+Ud_{45}+Ud_{67}$,

\ \ \ \ \ $K_{12}=-Ud_{01}+Ud_{23}+Ud_{45}+Ud_{67}$,

\ \ \ \ \ \gr$_{2}=\{$Lie algebra over real numbers generated by

\ \ \ \ \ \ \ \ \ \ \ \ \ \ \ $S_{1},S_{2},S_{3},S_{4},S_{5},S_{6},S_{7},L_{2},L_{3},L_{4},L_{5},L_{6},L_{7},L_{8}\}$.

\bigskip

\emph{Lemma 18.4.} $\ \{S_{1},S_{2},S_{3},S_{4},S_{5},S_{6},S_{7},L_{2},L_{3},L_{4},L_{5},L_{6},L_{7},L_{8}\}$
are orthogonal bases of \gr$_{2}$.

\bigskip

\emph{Proof. \ }We have the following Lie bracket operations with
calculations using Maxima.
The following tables show the operation results of $[A,B]$.
\setlength{\tabcolsep}{1mm} 
\begin{flushleft}
{\fontsize{7pt}{10pt} \selectfont%
\begin{tabular}
[c]{@{}c@{}c@{}|@{}c@{}c@{}c@{}c@{}c@{}c@{}c|}\cline{3-9}
&  &  &  &  & B &  &  & \\\cline{3-9}
&  & $S_{1}$ & \multicolumn{1}{|c}{$S_{2}$} & \multicolumn{1}{|c}{$S_{3}$} &
\multicolumn{1}{|c}{$S_{4}$} & \multicolumn{1}{|c}{$S_{5}$} &
\multicolumn{1}{|c}{$S_{6}$} & \multicolumn{1}{|c|}{$S_{7}$}\\\hline
\multicolumn{1}{|c|}{} & \multicolumn{1}{|c|}{$S_{1}$} & $0$ &
\multicolumn{1}{|c}{$-2S_{7}$} & \multicolumn{1}{|c}{$2S_{6}-3L_{2}$} &
\multicolumn{1}{|c}{$2S_{5}+3L_{3}$} & \multicolumn{1}{|c}{$-2S_{4}-3L_{4}$} &
\multicolumn{1}{|c}{$-2S_{3}+3L_{5}$} & \multicolumn{1}{|c|}{$2S_{2}$%
}\\\cline{2-9}%
\multicolumn{1}{|c|}{} & \multicolumn{1}{|c|}{$S_{2}$} & $2S_{7}$ &
\multicolumn{1}{|c}{$0$} & \multicolumn{1}{|c}{$2S_{5}-3L_{3}$} &
\multicolumn{1}{|c}{$-2S_{6}-3L_{2}$} & \multicolumn{1}{|c}{$-2S_{3}-3L_{5}$}
& \multicolumn{1}{|c}{$2S_{4}-3L_{4}$} & \multicolumn{1}{|c|}{$-2S_{1}$%
}\\\cline{2-9}%
\multicolumn{1}{|c|}{A} & \multicolumn{1}{|c|}{$S_{3}$} & $-2S_{6}+3L_{2}$ &
\multicolumn{1}{|c}{$-2S_{5}+3L_{3}$} & \multicolumn{1}{|c}{$0$} &
\multicolumn{1}{|c}{$S_{7}+3L_{6}$} & \multicolumn{1}{|c}{$2S_{2}-3L_{7}$} &
\multicolumn{1}{|c}{$2S_{1}+3L_{8}$} & \multicolumn{1}{|c|}{$-S_{4}$%
}\\\cline{2-9}%
\multicolumn{1}{|c|}{} & \multicolumn{1}{|c|}{$S_{4}$} & $-2S_{5}-3L_{3}$ &
\multicolumn{1}{|c}{$2S_{6}+3L_{2}$} & \multicolumn{1}{|c}{$-S_{7}-3L_{6}$} &
\multicolumn{1}{|c}{$0$} & \multicolumn{1}{|c}{$2S_{1}-3L_{8}$} &
\multicolumn{1}{|c}{$-2S_{2}-3L_{7}$} & \multicolumn{1}{|c|}{$S_{3}$%
}\\\cline{2-9}%
\multicolumn{1}{|c|}{} & \multicolumn{1}{|c|}{$S_{5}$} & $2S_{4}+3L_{4}$ &
\multicolumn{1}{|c}{$2S_{3}+3L_{5}$} & \multicolumn{1}{|c}{$-2S_{2}+3L_{7}$} &
\multicolumn{1}{|c}{$-2S_{1}+3L_{8}$} & \multicolumn{1}{|c}{$0$} &
\multicolumn{1}{|c}{$S_{7}-3L_{6}$} & \multicolumn{1}{|c|}{$-S_{6}$%
}\\\cline{2-9}%
\multicolumn{1}{|c|}{} & \multicolumn{1}{|c|}{$S_{6}$} & $2S_{3}-3L_{5}$ &
\multicolumn{1}{|c}{$-2S_{4}+3L_{4}$} & \multicolumn{1}{|c}{$-2S_{1}-3L_{8}$}
& \multicolumn{1}{|c}{$2S_{2}+3L_{7}$} & \multicolumn{1}{|c}{$-S_{7}+3L_{6}$}
& \multicolumn{1}{|c}{$0$} & \multicolumn{1}{|c|}{$S_{5}$}\\\cline{2-9}%
\multicolumn{1}{|c|}{} & \multicolumn{1}{|c|}{$S_{7}$} & $-2S_{2}$ &
\multicolumn{1}{|c}{$2S_{1}$} & \multicolumn{1}{|c}{$S_{4}$} &
\multicolumn{1}{|c}{$-S_{3}$} & \multicolumn{1}{|c}{$S_{6}$} &
\multicolumn{1}{|c}{$-S_{5}$} & \multicolumn{1}{|c|}{$0$}\\\hline
\end{tabular}
}

{\fontsize{7pt}{10pt} \selectfont%
\begin{tabular}
[c]{cc|ccccccc|}\cline{3-9}
&  &  &  &  & B &  &  & \\\cline{3-9}
&  & $L_{2}$ & \multicolumn{1}{|c}{$L_{3}$} & \multicolumn{1}{|c}{$L_{4}$} &
\multicolumn{1}{|c}{$L_{5}$} & \multicolumn{1}{|c}{$L_{6}$} &
\multicolumn{1}{|c}{$L_{7}$} & \multicolumn{1}{|c|}{$L_{8}$}\\\hline
\multicolumn{1}{|c|}{} & \multicolumn{1}{|c|}{$L_{2}$} & $0$ &
\multicolumn{1}{|c}{$-S_{7}-L_{6}$} & \multicolumn{1}{|c}{$-L_{7}$} &
\multicolumn{1}{|c}{$-L_{8}$} & \multicolumn{1}{|c}{$L_{3}$} &
\multicolumn{1}{|c}{$L_{4}$} & \multicolumn{1}{|c|}{$L_{5}$}\\\cline{2-9}%
\multicolumn{1}{|c|}{} & \multicolumn{1}{|c|}{$L_{3}$} & $S_{7}+L_{6}$ &
\multicolumn{1}{|c}{$0$} & \multicolumn{1}{|c}{$L_{8}$} &
\multicolumn{1}{|c}{$-L_{7}$} & \multicolumn{1}{|c}{$-L_{2}$} &
\multicolumn{1}{|c}{$L_{5}$} & \multicolumn{1}{|c|}{$-L_{4}$}\\\cline{2-9}%
\multicolumn{1}{|c|}{} & \multicolumn{1}{|c|}{$L_{4}$} & $L_{7}$ &
\multicolumn{1}{|c}{$-L_{8}$} & \multicolumn{1}{|c}{$0$} &
\multicolumn{1}{|c}{$-S_{7}+L_{6}$} & \multicolumn{1}{|c}{$-L_{5}$} &
\multicolumn{1}{|c}{$-L_{2}$} & \multicolumn{1}{|c|}{$L_{3}$}\\\cline{2-9}%
\multicolumn{1}{|c|}{A} & \multicolumn{1}{|c|}{$L_{5}$} & $L_{8}$ &
\multicolumn{1}{|c}{$L_{7}$} & \multicolumn{1}{|c}{$S_{7}-L_{6}$} &
\multicolumn{1}{|c}{$0$} & \multicolumn{1}{|c}{$L_{4}$} &
\multicolumn{1}{|c}{$-L_{3}$} & \multicolumn{1}{|c|}{$-L_{2}$}\\\cline{2-9}%
\multicolumn{1}{|c|}{} & \multicolumn{1}{|c|}{$L_{6}$} & $-L_{3}$ &
\multicolumn{1}{|c}{$L_{2}$} & \multicolumn{1}{|c}{$L_{5}$} &
\multicolumn{1}{|c}{$-L_{4}$} & \multicolumn{1}{|c}{$0$} &
\multicolumn{1}{|c}{$2L_{8}$} & \multicolumn{1}{|c|}{$-2L_{7}$}\\\cline{2-9}%
\multicolumn{1}{|c|}{} & \multicolumn{1}{|c|}{$L_{7}$} & $-L_{4}$ &
\multicolumn{1}{|c}{$-L_{5}$} & \multicolumn{1}{|c}{$L_{2}$} &
\multicolumn{1}{|c}{$L_{3}$} & \multicolumn{1}{|c}{$-2L_{8}$} &
\multicolumn{1}{|c}{$0$} & \multicolumn{1}{|c|}{$2L_{6}$}\\\cline{2-9}%
\multicolumn{1}{|c|}{} & \multicolumn{1}{|c|}{$L_{8}$} & $-L_{5}$ &
\multicolumn{1}{|c}{$L_{4}$} & \multicolumn{1}{|c}{$-L_{3}$} &
\multicolumn{1}{|c}{$L_{2}$} & \multicolumn{1}{|c}{$2L_{7}$} &
\multicolumn{1}{|c}{$-2L_{6}$} & \multicolumn{1}{|c|}{$0$}\\\hline
\end{tabular}
}

{\fontsize{7pt}{10pt} \selectfont%
\begin{tabular}
[c]{cc|ccccccc|}\cline{3-9}
&  &  &  &  & B &  &  & \\\cline{3-9}
&  & $L_{2}$ & \multicolumn{1}{|c}{$L_{3}$} & \multicolumn{1}{|c}{$L_{4}$} &
\multicolumn{1}{|c}{$L_{5}$} & \multicolumn{1}{|c}{$L_{6}$} &
\multicolumn{1}{|c}{$L_{7}$} & \multicolumn{1}{|c|}{$L_{8}$}\\\hline
\multicolumn{1}{|c|}{} & \multicolumn{1}{|c|}{$S_{1}$} & $S_{3}$ &
\multicolumn{1}{|c}{$-S_{4}$} & \multicolumn{1}{|c}{$S_{5}$} &
\multicolumn{1}{|c}{$-S_{6}$} & \multicolumn{1}{|c}{$0$} &
\multicolumn{1}{|c}{$0$} & \multicolumn{1}{|c|}{$0$}\\\cline{2-9}%
\multicolumn{1}{|c|}{} & \multicolumn{1}{|c|}{$S_{2}$} & $S_{4}$ &
\multicolumn{1}{|c}{$S_{3}$} & \multicolumn{1}{|c}{$S_{6}$} &
\multicolumn{1}{|c}{$S_{5}$} & \multicolumn{1}{|c}{$0$} &
\multicolumn{1}{|c}{$0$} & \multicolumn{1}{|c|}{$0$}\\\cline{2-9}%
\multicolumn{1}{|c|}{A} & \multicolumn{1}{|c|}{$S_{3}$} & $-S_{1}$ &
\multicolumn{1}{|c}{$-S_{2}$} & \multicolumn{1}{|c}{$0$} &
\multicolumn{1}{|c}{$0$} & \multicolumn{1}{|c}{$-S_{4}$} &
\multicolumn{1}{|c}{$S_{5}$} & \multicolumn{1}{|c|}{$-S_{6}$}\\\cline{2-9}%
\multicolumn{1}{|c|}{} & \multicolumn{1}{|c|}{$S_{4}$} & $-S_{2}$ &
\multicolumn{1}{|c}{$S_{1}$} & \multicolumn{1}{|c}{$0$} &
\multicolumn{1}{|c}{$0$} & \multicolumn{1}{|c}{$S_{3}$} &
\multicolumn{1}{|c}{$S_{6}$} & \multicolumn{1}{|c|}{$S_{5}$}\\\cline{2-9}%
\multicolumn{1}{|c|}{} & \multicolumn{1}{|c|}{$S_{5}$} & $0$ &
\multicolumn{1}{|c}{$0$} & \multicolumn{1}{|c}{$-S_{1}$} &
\multicolumn{1}{|c}{$-S_{2}$} & \multicolumn{1}{|c}{$S_{6}$} &
\multicolumn{1}{|c}{$-S_{3}$} & \multicolumn{1}{|c|}{$-S_{4}$}\\\cline{2-9}%
\multicolumn{1}{|c|}{} & \multicolumn{1}{|c|}{$S_{6}$} & $0$ &
\multicolumn{1}{|c}{$0$} & \multicolumn{1}{|c}{$-S_{2}$} &
\multicolumn{1}{|c}{$S_{1}$} & \multicolumn{1}{|c}{$-S_{5}$} &
\multicolumn{1}{|c}{$-S_{4}$} & \multicolumn{1}{|c|}{$S_{3}$}\\\cline{2-9}%
\multicolumn{1}{|c|}{} & \multicolumn{1}{|c|}{$S_{7}$} & $-3L_{3}$ &
\multicolumn{1}{|c}{$3L_{2}$} & \multicolumn{1}{|c}{$-3L_{5}$} &
\multicolumn{1}{|c}{$3L_{4}$} & \multicolumn{1}{|c}{$0$} &
\multicolumn{1}{|c}{$0$} & \multicolumn{1}{|c|}{$0$}\\\hline
\end{tabular}
}
\end{flushleft}

\ \ \ \ \ 

Evidently, $\{S_{1},S_{2},S_{3},S_{4},S_{5},S_{6},S_{7},L_{2},L_{3},L_{4},L_{5},L_{6},L_{7},L_{8}\}$ are  independent.
So these are bases of \gr$_{2}$.

For $X,Y\in  \{ S_{1}, S_{2}, S_{3}, S_{4}, S_{5}, S_{6}, S_{7}, L_{2}, L_{3}, L_{4}, L_{5}, L_{6}, L_{7}, L_{8}\}$,
we calculate $tr(X^{t}Y)$:
 
$\ \ \ \ \ \ \ \ \ \ \ \ \ tr(X^{t}Y)=0,(X\neq Y),$

$\ \ \ \ \ \ \ \ \ \ \ \ \ tr(X^{t}X)=360,(X=S_{1},\cdot\cdot\cdot,S_{7}),$

$\ \ \ \ \ \ \ \ \ \ \ \ \ tr(X^{t}X)=120,(X=L_{2},\cdot\cdot\cdot,L_{8}).$
\ \ \ \ \ \emph{Q.E.D.}
 
 \bigskip

\emph{Lemma 18.5.} \ \gr$_{2}$ is simple.

\bigskip

\emph{Proof. }\ In case of $L_{2}\in$ \gr$_{2}$.

\noindent
\ By \emph{Lemma 18.4}, we have

\ \ \ \ \ $\{[[L_{2},x],y] \mid x,y=L_{2},\cdot \cdot \cdot ,L_{8}\}\ni L_{2},\cdot \cdot \cdot ,L_{8}$,

\ \ \ \ \ $\{[[L_{2},x],y] \mid x,y=S_{1},\cdot \cdot \cdot ,S_{7}\}\ni S_{1},\cdot \cdot \cdot ,S_{7}$.

\noindent
Then we have $\{[[L_{2},x],y] \mid x,y\in $\gr$_{2}\}=$\gr$_{2}$.

In cases of $L_{3}, \cdot \cdot \cdot ,L_{8}$ are the same as $L_{2}$.

In case of $S_{3}, \cdot \cdot \cdot ,S_{6}$.

\noindent
\ By \emph{Lemma 18.4}, we have

\ \ \ \ \ $[[S_{3},S_{4}],L_{8}]=-6L_{7},[[S_{5},S_{6}],L_{8}]=6L_{7}$.

\noindent
Then,these applies to case of $L_{7}$.

In case of $S_{1},S_{2},S_{7}$.

\ \ \ \ \ $[[S_{1},L_{2}],[S_{2},L_{2}]]=[S_{3},S_{4}],[S_{7},S_{3}]=S_{4}$.

\noindent
Then,these applies to case $S_{3},S_{4}$.

Therfore \gr$_{2}$ is simple. \ \ \emph{Q.E.D.}

\bigskip

\emph{Lemma 18.6.} \ The Killing form $B_{2}$ of the Lie algebra 
\gr$_{2}^{\C}$ is given by

\ \ \ \ \ \ \ \ \ \ \ \ \ \ \ \ \ $B_{2}(R_{1},R_{2})=\frac{2}{15}%
tr(R_{1}R_{2}),R_{1},R_{2}\in $\gr$_{2}^{\C}.$

\bigskip

\emph{\ Proof.} \ Since \gr$_{2}^{\C}$ is simple,there
exist $\kappa \in \C$ such that

\ \ \ \ \ \ \ \ \ \ \ \ \ \ \ \ $B_{2}(R_{1},R_{2})=\kappa tr(R_{1}R_{2}).$

\noindent
To determin this $\kappa $,we put $R=R_{1}=R_{2}=L_{1}\in $\gr$_{2}^{\C}.$
$(adR)^{2}$ is calculated as follows by \emph{Lemma 18.4}  .

\ \ \ \ \ $\ [L_{1},[L_{1},L_{2}]]=[L_{1},2L_{3}]=-4L_{2},$

\ \ \ \ \ $\ [L_{1},[L_{1},L_{3}]]=[L_{1},-2L_{2}]=-4L_{3},$

\ \ \ \ \ \ $[L_{1},[L_{1},L_{4}]]=[L_{1},L_{5}]=-L_{4},$

\ \ \ \ \ \ $[L_{1},[L_{1},L_{5}]]=[L_{1},-L_{4}]=-L_{5},$

\ \ \ \ \ \ $[L_{1},[L_{1},L_{6}]]=[L_{1},L_{7}]=-L_{6},$

\ \ \ \ \ \ $[L_{1},[L_{1},L_{7}]]=[L_{1},-L_{6}]=-L_{7},$

\ \ \ \ \ \ $[L_{1},[L_{1},S_{1}]]=[L_{1},S_{2}]=-S_{1},$

\ \ \ \ \ \ $[L_{1},[L_{1},S_{2}]]=[L_{1},-S_{1}]=-S_{2},$

\ \ \ \ \ \ $[L_{1},[L_{1},S_{3}]]=[L_{1},-S_{4}]=-S_{3},$

\ \ \ \ \ \ $[L_{1},[L_{1},S_{4}]]=[L_{1},S_{3}]=-S_{4},$

\ \ \ \ \ the others $=0$.

\noindent
Hence we have

\ \ \ \ $B_{2}(L_{1},L_{1})=tr((ad(L_{1}))^{2})=(-4)\times 2+(-1)\times 8=-16.$

\noindent
On the other hand ,by calculate with Maxima we have

\ \ \ \ $tr((L_{1})(L_{1}))=-120.$

\noindent
Hence we have $\kappa =\frac{-16}{-120}=\frac{2}{15}.$ \ \ \ \ \emph{Q.E.D.}

\bigskip

\emph{Lemma 18.7.} The rank of the Lie algebra \gr$_{2}^{\C}$
is 2. The roots of \gr$_{2}^{\C}$ relative to some Cartan
subalgebra of \gr$_{2}^{\C}$ are given by

$\ \ \ \ \ \ \ \ \ \ \ \ \ \ \ \ \ \ \pm (\lambda _{1}-\lambda _{2}),\pm
(\lambda _{1}-\lambda _{3}),\pm (\lambda _{2}-\lambda _{3}).$

\ \ \ \ \ \ \ \ \ \ \ \ \ \ \ \ \ \ $\pm \lambda _{1},\pm\lambda _{2},\pm \lambda _{3},$

\ \ \ \ \ \ \ \ \ \ \ \ \ \ \ \ \ \ with $\lambda _{1}+\lambda _{2}+\lambda _{3}=0$.

\bigskip

\emph{Proof.} Let \gh$=\{h=\sum\limits_{k=1}^{3}\lambda
_{k}H_{k}|H_{k}=-iUd_{2k2k+1},\lambda _{1}+\lambda _{2}+\lambda
_{3}=0,\lambda _{k}\in \C\}$,

\noindent
then \gh \ is an abelian subalgebra of \gr$_{2}^{\C}$(it will be
a Cartan subalgebra of \gr$_{2}^{\C}$).

We have followings with calculations using Maxima.

\ \ \ \ \ $[h,L_{2}+iL_{3}]=(\lambda _{1}-\lambda _{2})(L_{2}+iL_{3}),$

\ \ \ \ \ $[h,L_{4}+iL_{5}]=(\lambda _{1}-\lambda _{3})(L_{4}+iL_{5}),$

\ \ \ \ \ $[h,L_{7}+iL_{8}]=(\lambda _{2}-\lambda _{3})(L_{7}+iL_{8}),$

\ \ \ \ \ $[h,-L_{2}+iL_{3}]=-(\lambda _{1}-\lambda _{2})(-L_{2}+iL_{3}),$

\ \ \ \ \ $[h,-L_{4}+iL_{5}]=-(\lambda _{1}-\lambda _{3})(-L_{4}+iL_{5}),$

\ \ \ \ \ $[h,-L_{7}+iL_{8}]=-(\lambda _{2}-\lambda _{3})(-L_{7}+iL_{8}),$

\ \ \ \ \ $[h,S_{1}+iS_{2}]=\lambda _{1}(S_{1}+iS_{2}),$

\ \ \ \ \ $[h,S_{3}+iS_{4}]=\lambda _{2}(S_{3}+iS_{4}),$

\ \ \ \ \ $[h,S_{5}+iS_{6}]=\lambda _{3}(S_{5}+iS_{6}),$

\ \ \ \ \ $[h,S_{1}-iS_{2}]=-\lambda _{1}(S_{1}-iS_{2}),$

\ \ \ \ \ $[h,S_{3}-iS_{4}]=-\lambda _{2}(S_{3}-iS_{4}),$

\ \ \ \ \ $[h,S_{5}-iS_{6}]=-\lambda _{3}(S_{5}-iS_{6}).$
\ \ \ \ \ \ \ \ \ \ \ \ \ \emph{Q.E.D.}

\bigskip

\emph{Theorem 18.8.} \ In the root system of \emph{Lemma 18.7},

\ \ \ \ \ \ \ $\alpha _{1}=\lambda _{1}-\lambda _{2},\alpha _{2}=\lambda_{2}$

\noindent
is a fundamental root sysyten of the Lie algebra \gr$_{2}^{\C}$ and

$\ \ \ \ \ \ \ \mu =2\alpha _{1}+3\alpha _{2}$

\noindent
is the highest root. The Dynkin diagram and the extended Dynkin diagram 
of  \gr$_{2}^{\C}$ are respectively given by

\begin{figure}[H]
\centering
\includegraphics[width=10cm, height=2.0cm]{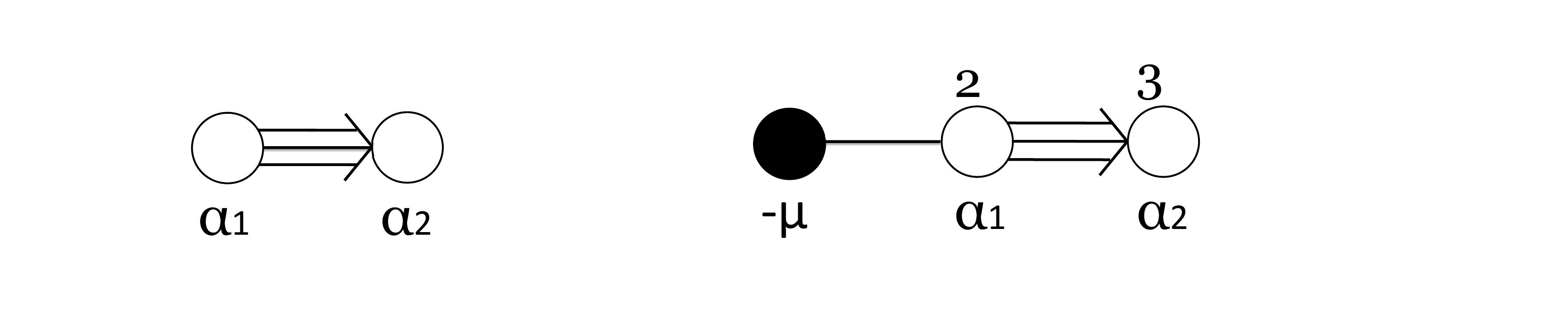}
\end{figure}

\emph{Proof.} In the following, the notation $n_{1}$ $n_{2}$ 
denotes the root $n_{1}\alpha _{1}+n_{2}\alpha _{2}$.

\noindent
For example, $\alpha =%
\begin{array}{@{}ll@{}}
n_{1} & n_{2} %
\end{array}%
$ means \ $\alpha =n_{1}\alpha _{1}+n_{2}\alpha _{2}$ .

Now,all positive roots of \gr$_{2}^{\C}$ are represented by
\begin{align*}
\lambda _{1} &= %
\begin{array}{ll}
1 & 1%
\end{array},\\
\lambda _{2} &= %
\begin{array}{ll}
0 & 1 %
\end{array},\\
-\lambda _{3}& = %
\begin{array}{ll}
1 & 2 %
\end{array},\\
\lambda _{1}-\lambda _{2} &= %
\begin{array}{ll}
1 & 0 %
\end{array},\\
 \lambda_{1}-\lambda _{3} &= %
\begin{array}{ll}
2 & 3%
\end{array},\\
\lambda _{2}-\lambda _{3} &= %
\begin{array}{ll}
1 & 3 %
\end{array},\\
\end{align*}

\noindent
Hence $\Pi =\{\alpha _{1},\alpha _{2} \}$ is a
fundamental root system of \gr$_{2}^{\C}$ and

\noindent
$\mu =2\alpha _{1}+3\alpha _{2}$ is the highest root. The real part of \gh$_{R}$
of \gh\  is

\ \ \ \ \ \ \ \ \ \ \ \ \ \ \ \ \ \ \ \ \gh$_{R}=\{H=\sum%
\limits_{k=0}^{3}\lambda _{k}H_{k} \mid \lambda _{k}\in \R\}.$

The Killing form $B_{2}$ of \gr$_{2}^{\C}$ is $B_{2}(R_{1},R_{2})=\frac{2}{15%
}tr(R_{1}R_{2})$ (\emph{Lemma 18.5.}),  so that

\ \ \ \ \ \ \ \ \ \ \ $B_{2}(H,H^{\prime} )=8\sum\limits_{k=1}^{3}\lambda
_{k}\lambda _{k}^{\prime },H=\sum\limits_{k=1}^{3}\lambda
_{k}H_{k},H^{\prime }=\sum\limits_{k=1}^{3}\lambda _{k}^{\prime
}H_{k}^{\prime }\in $\gh$_{R}.$

\noindent
Indeed,by calculate with Maxima we have

$\ \ \ \ \ \ \ \ \ \ \ B_{2}(H,H^{\prime} )=\frac{2}{15}tr(HH^{\prime} )=%
\frac{2}{15}(60\sum\limits_{k=1}^{3}\lambda _{k}\lambda _{k}^{\prime
})=8\sum\limits_{k=1}^{3}\lambda _{k}\lambda _{k}^{\prime }.$

Now,the canonical elements H$_{\alpha _{i}}\in $\gh$_{R}$ corresponding
to $\alpha _{i}(B_{2}(H_{\alpha _{i}},H)=\alpha _{i}(H),H\in $\gh$_{R})$
are determind by

$\ \ \ \ \ \ \ \ \ \ \ H_{\alpha _{1}}=\frac{1}{8}%
(H_{1}-H_{2}),H_{\alpha _{2}}=-\frac{1}{24}H_{1}+\frac{1}{12}H_{2}-\frac{1}{24}H_{\alpha _{3}},$

\noindent
Hence we have

$\ \ \ \ \ \ \ \ \ \ \ \ (\alpha _{1},\alpha _{1})=\frac{1}{4}$,
$\ \ \ \ \ \ \ \ \ \ (\alpha _{2},\alpha _{2})=\frac{1}{12}$,
$\ \ \ \ \ \ \ \ \ \ (\alpha _{1},\alpha _{2})=-\frac{1}{8}$,

\ \ \ \ \ \ \ \ \ \ \ \ $(-\mu ,-\mu )=\frac{1}{4}$,
\ \ \ \ \ \ \ \ $(-\mu ,\alpha_{1})=-\frac{1}{8}$,
\ \ \ \ \ \ \ \ $(-\mu ,\alpha _{2})=0$.

\noindent
using them,we can draw the Dynkin diagram and the extended Dynkin 
diagram of  \gr$_{2}^{\C}.$ \ \ \ \ \emph{Q.E.D.}

\bigskip

\emph{Corollary 18.9.} For a $248\times 248$ matrix $X$,let $X|_{2}$ be the matrix in which the $28\times 28$
elements in the upper left corner are clipped from  $X$. Furthermore, let \gr$_{2}^{\C}|_{2}=\{R|_{2} \mid R\in $\gr$_{2}^{\C} \}$ .
\emph{Theorem 18.7.} holds for \gr$_{2}^{\C}|_{2}$ as well. However, the
Killing form $B_{2}(R_{1},R_{2})=tr(R_{1}R_{2}) \ \ (R_{1},R_{2} \in $\gr$_{2}^{\C}|_{2}).$

\bigskip

Next, let digitalize the exceptional simple Lie algebra of type $G_{2}$ as
matrices of \gD$_{4} \subset M(8 \times 8,\R)$.

By \emph{Lemma 13.8}, \gr\gd \ is isomorphic to \gD$_{4}$ under the correspondence fd as follows:

\ \ \ \ \ \ \ \ \ \ \ \ \ \ fd\ :\ \gr\gd \ $\ni$ Ud$_{ij}\rightarrow D_{ij}\in $\ \gD$_{4}$.

By identifying $S_{i}$ with fd($S_{i}$) and $L_{i\text{ }}$with fd($L_{i}$)
,  we have followings:

\begin{flushleft}
{\fontsize{7pt}{8pt} \selectfont%
$S_{1}=\left(
\begin{array}
[c]{@{}cccccccc@{}}%
0 & 0 & 0 & 0 & 0 & 0 & 0 & 0\\
0 & 0 & 2 & 0 & 0 & 0 & 0 & 0\\
0 & -2 & 0 & 0 & 0 & 0 & 0 & 0\\
0 & 0 & 0 & 0 & 0 & 0 & 0 & 0\\
0 & 0 & 0 & 0 & 0 & 0 & 0 & -1\\
0 & 0 & 0 & 0 & 0 & 0 & -1 & 0\\
0 & 0 & 0 & 0 & 0 & 1 & 0 & 0\\
0 & 0 & 0 & 0 & 1 & 0 & 0 & 0
\end{array}
\right)  ,S_{2}=\left(
\begin{array}
[c]{@{}cccccccc@{}}%
0 & 0 & 0 & 0 & 0 & 0 & 0 & 0\\
0 & 0 & 0 & 2 & 0 & 0 & 0 & 0\\
0 & 0 & 0 & 0 & 0 & 0 & 0 & 0\\
0 & -2 & 0 & 0 & 0 & 0 & 0 & 0\\
0 & 0 & 0 & 0 & 0 & 0 & -1 & 0\\
0 & 0 & 0 & 0 & 0 & 0 & 0 & 1\\
0 & 0 & 0 & 0 & 1 & 0 & 0 & 0\\
0 & 0 & 0 & 0 & 0 & -1 & 0 & 0
\end{array}
\right)  $,

$S_{3}=\left(
\begin{array}
[c]{@{}cccccccc@{}}%
0 & 0 & 0 & 0 & 0 & 0 & 0 & 0\\
0 & 0 & 0 & 0 & 2 & 0 & 0 & 0\\
0 & 0 & 0 & 0 & 0 & 0 & 0 & 1\\
0 & 0 & 0 & 0 & 0 & 0 & 1 & 0\\
0 & -2 & 0 & 0 & 0 & 0 & 0 & 0\\
0 & 0 & 0 & 0 & 0 & 0 & 0 & 0\\
0 & 0 & 0 & -1 & 0 & 0 & 0 & 0\\
0 & 0 & -1 & 0 & 0 & 0 & 0 & 0
\end{array}
\right)  ,S_{4}=\left(
\begin{array}
[c]{@{}cccccccc@{}}%
0 & 0 & 0 & 0 & 0 & 0 & 0 & 0\\
0 & 0 & 0 & 0 & 0 & 2 & 0 & 0\\
0 & 0 & 0 & 0 & 0 & 0 & 1 & 0\\
0 & 0 & 0 & 0 & 0 & 0 & 0 & -1\\
0 & 0 & 0 & 0 & 0 & 0 & 0 & 0\\
0 & -2 & 0 & 0 & 0 & 0 & 0 & 0\\
0 & 0 & -1 & 0 & 0 & 0 & 0 & 0\\
0 & 0 & 0 & 1 & 0 & 0 & 0 & 0
\end{array}
\right)  $,

$S_{5}=\left(
\begin{array}
[c]{@{}cccccccc@{}}%
0 & 0 & 0 & 0 & 0 & 0 & 0 & 0\\
0 & 0 & 0 & 0 & 0 & 0 & 2 & 0\\
0 & 0 & 0 & 0 & 0 & -1 & 0 & 0\\
0 & 0 & 0 & 0 & -1 & 0 & 0 & 0\\
0 & 0 & 0 & 1 & 0 & 0 & 0 & 0\\
0 & 0 & 1 & 0 & 0 & 0 & 0 & 0\\
0 & -2 & 0 & 0 & 0 & 0 & 0 & 0\\
0 & 0 & 0 & 0 & 0 & 0 & 0 & 0
\end{array}
\right)  ,S_{6}=\left(
\begin{array}
[c]{@{}cccccccc@{}}%
0 & 0 & 0 & 0 & 0 & 0 & 0 & 0\\
0 & 0 & 0 & 0 & 0 & 0 & 0 & 2\\
0 & 0 & 0 & 0 & -1 & 0 & 0 & 0\\
0 & 0 & 0 & 0 & 0 & 1 & 0 & 0\\
0 & 0 & 1 & 0 & 0 & 0 & 0 & 0\\
0 & 0 & 0 & -1 & 0 & 0 & 0 & 0\\
0 & 0 & 0 & 0 & 0 & 0 & 0 & 0\\
0 & -2 & 0 & 0 & 0 & 0 & 0 & 0
\end{array}
\right)  $,

$S_{7}=\left(
\begin{array}
[c]{@{}cccccccc@{}}%
0 & 0 & 0 & 0 & 0 & 0 & 0 & 0\\
0 & 0 & 0 & 0 & 0 & 0 & 0 & 2\\
0 & 0 & 0 & 0 & -1 & 0 & 0 & 0\\
0 & 0 & 0 & 0 & 0 & 1 & 0 & 0\\
0 & 0 & 1 & 0 & 0 & 0 & 0 & 0\\
0 & 0 & 0 & -1 & 0 & 0 & 0 & 0\\
0 & 0 & 0 & 0 & 0 & 0 & 0 & 0\\
0 & -2 & 0 & 0 & 0 & 0 & 0 & 0
\end{array}
\right)  ,L_{2}=\left(
\begin{array}
[c]{@{}cccccccc@{}}%
0 & 0 & 0 & 0 & 0 & 0 & 0 & 0\\
0 & 0 & 0 & 0 & 0 & 0 & 0 & 0\\
0 & 0 & 0 & 0 & 1 & 0 & 0 & 0\\
0 & 0 & 0 & 0 & 0 & 1 & 0 & 0\\
0 & 0 & -1 & 0 & 0 & 0 & 0 & 0\\
0 & 0 & 0 & -1 & 0 & 0 & 0 & 0\\
0 & 0 & 0 & 0 & 0 & 0 & 0 & 0\\
0 & 0 & 0 & 0 & 0 & 0 & 0 & 0
\end{array}
\right)  ,$

$L_{3}=\left(
\begin{array}
[c]{@{}cccccccc@{}}%
0 & 0 & 0 & 0 & 0 & 0 & 0 & 0\\
0 & 0 & 0 & 0 & 0 & 0 & 0 & 0\\
0 & 0 & 0 & 0 & 0 & -1 & 0 & 0\\
0 & 0 & 0 & 0 & 1 & 0 & 0 & 0\\
0 & 0 & 0 & -1 & 0 & 0 & 0 & 0\\
0 & 0 & 1 & 0 & 0 & 0 & 0 & 0\\
0 & 0 & 0 & 0 & 0 & 0 & 0 & 0\\
0 & 0 & 0 & 0 & 0 & 0 & 0 & 0
\end{array}
\right)  ,L_{4}=\left(
\begin{array}
[c]{@{}cccccccc@{}}%
0 & 0 & 0 & 0 & 0 & 0 & 0 & 0\\
0 & 0 & 0 & 0 & 0 & 0 & 0 & 0\\
0 & 0 & 0 & 0 & 0 & 0 & 1 & 0\\
0 & 0 & 0 & 0 & 0 & 0 & 0 & 1\\
0 & 0 & 0 & 0 & 0 & 0 & 0 & 0\\
0 & 0 & 0 & 0 & 0 & 0 & 0 & 0\\
0 & 0 & -1 & 0 & 0 & 0 & 0 & 0\\
0 & 0 & 0 & -1 & 0 & 0 & 0 & 0
\end{array}
\right)  $,

$L_{5}=\left(
\begin{array}
[c]{@{}cccccccc@{}}%
0 & 0 & 0 & 0 & 0 & 0 & 0 & 0\\
0 & 0 & 0 & 0 & 0 & 0 & 0 & 0\\
0 & 0 & 0 & 0 & 0 & 0 & 0 & -1\\
0 & 0 & 0 & 0 & 0 & 0 & 1 & 0\\
0 & 0 & 0 & 0 & 0 & 0 & 0 & 0\\
0 & 0 & 0 & 0 & 0 & 0 & 0 & 0\\
0 & 0 & 0 & -1 & 0 & 0 & 0 & 0\\
0 & 0 & 1 & 0 & 0 & 0 & 0 & 0
\end{array}
\right)  ,L_{6}=\left(
\begin{array}
[c]{@{}cccccccc@{}}%
0 & 0 & 0 & 0 & 0 & 0 & 0 & 0\\
0 & 0 & 0 & 0 & 0 & 0 & 0 & 0\\
0 & 0 & 0 & 0 & 0 & 0 & 0 & 0\\
0 & 0 & 0 & 0 & 0 & 0 & 0 & 0\\
0 & 0 & 0 & 0 & 0 & -1 & 0 & 0\\
0 & 0 & 0 & 0 & 1 & 0 & 0 & 0\\
0 & 0 & 0 & 0 & 0 & 0 & 0 & 1\\
0 & 0 & 0 & 0 & 0 & 0 & -1 & 0
\end{array}
\right)  $,

$L_{7}=\left(
\begin{array}
[c]{@{}cccccccc@{}}%
0 & 0 & 0 & 0 & 0 & 0 & 0 & 0\\
0 & 0 & 0 & 0 & 0 & 0 & 0 & 0\\
0 & 0 & 0 & 0 & 0 & 0 & 0 & 0\\
0 & 0 & 0 & 0 & 0 & 0 & 0 & 0\\
0 & 0 & 0 & 0 & 0 & 0 & 1 & 0\\
0 & 0 & 0 & 0 & 0 & 0 & 0 & 1\\
0 & 0 & 0 & 0 & -1 & 0 & 0 & 0\\
0 & 0 & 0 & 0 & 0 & -1 & 0 & 0
\end{array}
\right)  ,L_{8}=\left(
\begin{array}
[c]{@{}cccccccc@{}}%
0 & 0 & 0 & 0 & 0 & 0 & 0 & 0\\
0 & 0 & 0 & 0 & 0 & 0 & 0 & 0\\
0 & 0 & 0 & 0 & 0 & 0 & 0 & 0\\
0 & 0 & 0 & 0 & 0 & 0 & 0 & 0\\
0 & 0 & 0 & 0 & 0 & 0 & 0 & -1\\
0 & 0 & 0 & 0 & 0 & 0 & 1 & 0\\
0 & 0 & 0 & 0 & 0 & -1 & 0 & 0\\
0 & 0 & 0 & 0 & 1 & 0 & 0 & 0
\end{array}
\right)  .$}
\end{flushleft}

Therfore we have \emph{Lemma 18.10}.

\bigskip

\emph{Lemma 18.10.} The above $\{S_{1},\cdot\cdot\cdot,S_{7},L_{2},\cdot\cdot
\cdot,L_{8}\}$ are bases of the exceptional simple Lie algebra of type $G_{2}$.

\bigskip

\emph{Lemma 18.11.}  \gg$_{2}$ is expressed  by 

$\{ D \in Hom_{\R}(\R^{8}) \mid D(\textrm{C}_{m}(x, y))=\textrm{C}_{m}(D(x), y)+\textrm{C}_{m}(x, D(y)), x,y\in \R^{8} \}.$

\noindent
And $\{S_{1},\cdot\cdot\cdot,S_{7},L_{2},\cdot\cdot
\cdot,L_{8}\}$ are orthogonal bases of \gg$_{2}$.

\bigskip

\emph{Proof} \ By \emph{Remark 6.3} and \emph{Proposition 18.2}, \gg$_{2}$ is expressed  by 

$\{ D \in Hom_{\R}(\R^{8}) \mid D(\textrm{C}_{m}(x, y))=\textrm{C}_{m}(D(x), y)+\textrm{C}_{m}(x, D(y)), x,y\in \R^{8} \}.$

\noindent
Let check each element of $\{S_{1},\cdot\cdot\cdot,S_{7},L_{2},\cdot\cdot\cdot,L_{8}\}$
satisfy the condition:

\ \ \ \ \ \ \ \ $D(\textrm{C}_{m}(x, y))=\textrm{C}_{m}(D(x), y)+\textrm{C}_{m}(x, D(y))$.

In case of $S_{1}:$

$S_{1}(\textrm{C}_{m}(x, y))=$

$(0,$

$2(x_{5}y_{7}-x_{4}y_{6}-x_{7}y_{5}+x_{6}y_{4}-x_{1}y_{3}+x_{0}y_{2}%
+x_{3}y_{1}+x_{2}y_{0}),$

$-2(x_{6}y_{7}-x_{7}y_{6}+x_{4}y_{5}-x_{5}y_{4}+x_{2}y_{3}-x_{3}y_{2}%
+x_{0}y_{1}+x_{1}y_{0}),$

$0,$

$-x_{0}y_{7}-x_{1}y_{6}-x_{2}y_{5}-x_{3}y_{4}+x_{4}y_{3}+x_{5}y_{2}+x_{6}%
y_{1}-x_{7}y_{0},$

$x_{1}y_{7}-x_{0}y_{6}-x_{3}y_{5}+x_{2}y_{4}+x_{5}y_{3}-x_{4}y_{2}-x_{7}%
y_{1}-x_{6}y_{0},$

$-x_{2}y_{7}-x_{3}y_{6}+x_{0}y_{5}+x_{1}y_{4}+x_{6}y_{3}+x_{7}y_{2}-x_{4}%
y_{1}+x_{5}y_{0},$

$-x_{3}y_{7}+x_{2}y_{6}-x_{1}y_{5}+x_{0}y_{4}+x_{7}y_{3}-x_{6}y_{2}+x_{5}%
y_{1}+x_{4}y_{0}).$

$\textrm{C}_{m}(S_{1}(x), y)=$

$(-x_{4}y_{7}-x_{5}y_{6}+x_{6}y_{5}+x_{7}y_{4}+2x_{1}y_{2}-2x_{2}y_{1},$

$x_{5}y_{7}-x_{4}y_{6}-x_{7}y_{5}+x_{6}y_{4}-2x_{1}y_{3}+2x_{2}y_{0},$

$-x_{6}y_{7}+x_{7}y_{6}-x_{4}y_{5}+x_{5}y_{4}-2x_{2}y_{3}-2x_{1}y_{0},$

$-x_{7}y_{7}-x_{6}y_{6}-x_{5}y_{5}-x_{4}y_{4}+2x_{2}y_{2}+2x_{1}y_{1},$

$-2x_{1}y_{6}-2x_{2}y_{5}+x_{4}y_{3}-x_{5}y_{2}-x_{6}y_{1}-x_{7}y_{0},$

$2x_{1}y_{7}+2x_{2}y_{4}+x_{5}y_{3}+x_{4}y_{2}+x_{7}y_{1}-x_{6}y_{0},$

$-2x_{2}y_{7}+_{2}x_{1}y_{4}+x_{6}y_{3}-x_{7}y_{2}+x_{4}y_{1}+x_{5}y_{0},$

$2x_{2}y_{6}-_{2}x_{1}y_{5}+x_{7}y_{3}+x_{6}y_{2}-x_{5}y_{1}+x_{4}y_{0}).$

$\textrm{C}_{m}(x, S_{1}(y))=$

$(x_{4}y_{7}+x_{5}y_{6}-x_{6}y_{5}-x_{7}y_{4}-2x_{1}y_{2}+2x_{2}y_{1},$

$x_{5}y_{7}-x_{4}y_{6}-x_{7}y_{5}+x_{6}y_{4}+2x_{0}y_{2}+2x_{3}y_{1},$

$-x_{6}y_{7}+x_{7}y_{6}-x_{4}y_{5}+x_{5}y_{4}+2x_{3}y_{2}-2x_{0}y_{1},$

$x_{7}y_{7}+x_{6}y_{6}+x_{5}y_{5}+x_{4}y_{4}-2x_{2}y_{2}-2x_{1}y_{1},$

$-x_{0}y_{7}+x_{1}y_{6}+x_{2}y_{5}-x_{3}y_{4}+2x_{5}y_{2}+2x_{6}y_{1},$

$-x_{1}y_{7}-x_{0}y_{6}-x_{3}y_{5}-x_{2}y_{4}-2x_{4}y_{2}-2x_{7}y_{1},$

$x_{2}y_{7}-x_{3}y_{6}+x_{0}y_{5}-x_{1}y_{4}+2x_{7}y_{2}-2x_{4}y_{1},$

$-x_{3}y_{7}-x_{2}y_{6}+x_{1}y_{5}+x_{0}y_{4}-2x_{6}y_{2}+2x_{5}y_{1}).$

\noindent
Therfore $S_{1}\textrm{C}_{m}(x, y)=\textrm{C}_{m}(S_{1}(x), y)+\textrm{C}_{m}(x, S_{1}(y)).$

Other cases can be confirmed in the same way.

For $R_{1},R_{2} \in \{S_{1},\cdot\cdot\cdot,S_{7},L_{2},\cdot\cdot\cdot,L_{8}\}$,
let's  calculate $tr(R_{1}.^{t}R_{2})$.

\ \ \ \ \ \ \ \ \ \ \ \ \ \ \ \ $tr(R_{1}.^{t}R_{2})=0 \ \ \ (R_{1} \neq R_{2})$,

\ \ \ \ \ \ \ \ \ \ \ \ \ \ \ \ $tr(R_{1}.^{t}R_{1})=12 \ \ (R_{1} \in  \{S_{1},\cdot\cdot\cdot,S_{7}\})$,

\ \ \ \ \ \ \ \ \ \ \ \ \ \ \ \ $tr(R_{2}.^{t}R_{2})=4 \  \ \ (R_{2} \in \{L_{2},\cdot\cdot\cdot,L_{8}\})$.

Since dim \gg$_{2}=14$, $\{S_{1},\cdot\cdot\cdot,S_{7},L_{2},\cdot\cdot\cdot,L_{8}\}$
are orthogonal bases of \gg$_{2}$.    
 \ \ \ \ \emph{Q.E.D.}

\bigskip

\emph{Lemma 18.12.}  The simply connected compact Lie group $G_{2}$ is expressed by

$G_{2}=\{\alpha \in Iso_{\R}(\R^{8}) \mid \alpha (\textrm{C}_{m}(X, Y))=\textrm{C}_{m}(\alpha X,\alpha Y) \}$ .

\bigskip

\emph{Proof.}  By \emph{Proposition 18.1} and \emph{Remark 6.3}, we have the above \emph{Lemma}.\ \ \ \ \ \emph{Q.E.D.}

\section{Some Lie subalgebras of \gr$_{8}^{\C}$}

\subsection{\gs\gp$(1)$ and \gs\gu$(3)$ of \gr$_{2}^{\C}$}
\ \ \ \emph{Definition 19.1.} \ We define the following matrices.

$sp_{1}=\left( 
\begin{array}{ll}
i & 0 \\ 
0 & -i%
\end{array}%
\right) ,sp_{2}=\left( 
\begin{array}{ll}
0 & i \\ 
i & 0%
\end{array}%
\right) ,sp_{3}=\left( 
\begin{array}{ll}
0 & -1 \\ 
1 & 0%
\end{array}%
\right) ,$

$su_{1}=\left( 
\begin{array}{lll}
i & 0 & 0 \\ 
0 & -i & 0 \\ 
0 & 0 & 0%
\end{array}%
\right) ,su_{2}=\left( 
\begin{array}{lll}
i & 0 & 0 \\ 
0 & 0 & 0 \\ 
0 & 0 & -i%
\end{array}%
\right) ,su_{3}=\left( 
\begin{array}{lll}
0 & 0 & 0 \\ 
0 & 0 & 1 \\ 
0 & -1 & 0%
\end{array}%
\right) ,$

$su_{4}=\left( 
\begin{array}{lll}
0 & 0 & 0 \\ 
0 & 0 & i \\ 
0 & i & 0%
\end{array}%
\right) ,$
$su_{5}=\left( 
\begin{array}{lll}
0 & 0 & -1 \\ 
0 & 0 & 0 \\ 
1 & 0 & 0%
\end{array}%
\right) ,su_{6}=\left( 
\begin{array}{lll}
0 & 0 & i \\ 
0 & 0 & 0 \\ 
i & 0 & 0%
\end{array}%
\right) ,$

$su_{7}=\left( 
\begin{array}{lll}
0 & 1 & 0 \\ 
-1 & 0 & 0 \\ 
0 & 0 & 0%
\end{array}%
\right) ,su_{8}=\left( 
\begin{array}{lll}
0 & i & 0 \\ 
i & 0 & 0 \\ 
0 & 0 & 0%
\end{array}%
\right) .$

$sp_{1},sp_{2}$ and $sp_{3}$ are bases of  \gs \gp $(1)=$%
$\{X\in M(2\times 2,\C) \mid \ ^{t}XJ+JX=0\}$ where $J=-sp_{3}.$

$su_{k}(k=1,2,3,4,5,6,7,8)$ are bases of \gs \gu $(3)=\{X\in M(3\times 3,\C) \mid \ X^{\ast}+X=0,tr(X)=0\}$ where $\ast$ means matrix transpose and complex conjugate.

\bigskip

\emph{Definition 19.2.}  We define the following Lie subalgebras of \gr$_{2}^{\C}$.

\gr$^{1}$\gs\gp$(1)=\{ c_{1}S_{1}+c_{2}S_{2}+c_{3}S_{7} \mid c_{1} ,c_{2},c_{3}\in \C \}$,

\gr$^{2}$\gs\gp$(1)=\{ c_{1}L_{6}+c_{2}L_{7}+c_{3}L_{8} \mid c_{1} ,c_{2},c_{3}\in \C \}$,

(cf. I.Yokota\cite[\emph{Theorem1.10.1. Remark 2.}]{Yokota1}), 

\gr$^{1}$\gs\gu$(3)=\{$Lie algebra generated by $L_{k}(1\leq k\leq 8)$ over complex numbers$\}$,

(cf. I.Yokota\cite[\emph{Theorem1.5.1. }]{Yokota1}), 

\bigskip

\emph{Lemma 19.3.} \ \gr$^{1}$\gs\gp$(1)$ and \gr$^{2}$\gs\gp$(1)$ are isomorphic to \gs \gp$(1)$. 

\gr$^{1}$\gs\gp$(1) \oplus $\gr$^{2}$\gs\gp$(1)$ is a direct sum Lie subalgebra of \gr$_{2}^{\C}$.
\bigskip

\emph{Proof. }\ We have following Lie bracket operations.

$[sp_{1},sp_{2}]=2sp_{3},[sp_{2},sp_{3}]=2sp_{1},[sp_{3},sp_{1}]=2sp_{2}$,

$[S_{2},S_{1}]=2S_{7},[S_{1},S_{7}]=2S_{2},[S_{7},S_{2}]=2S_{1}$,

$[L_{6},L_{7}]=2L_{8},[L_{7},L_{8}]=2L_{6},[L_{8},L_{6}]=2L_{7}$.

\noindent
By corresponding $S_{2}$ to $sp_{1},S_{1}$ to $sp_{2}$, and $S_{7}$ to $sp_{3}$, \gr$^{1}$\gs\gp$(1)$
is isomorphic to \gs \gp$(1)$.
Similarly by corresponding $L_{6}$ to $sp_{1},L_{7}$ to $sp_{2}$, and $L_{8}$ to $sp_{3}$, \gr$^{2}$\gs\gp$(1)$
is isomorphic to \gs \gp$(1)$. 

For $\forall x \in $\gr$^{1}$\gs\gp$(1)$,$\forall y \in $\gr$^{2}$\gs\gp$(1)$, we have $[x,y]=0$ by using Maxima.  
\ \ \ \ \emph{Q.E.D.}

\bigskip

\emph{Lemma 19.4.} \ \gr$^{1}$\gs\gu$(3)$ is isomorphic to \gs \gu$(3)$.

\bigskip

\emph{Proof. }\ We have following Lie bracket operations with direct calculation , by \emph{Lemma 18.2} and by using Maxima.

\noindent
The following tables show the operation results of $[A,B]$.

\setlength{\tabcolsep}{1.5mm} 
\begin{flushleft}
{\fontsize{8pt}{10pt} \selectfont%
\begin{tabular}{@{}c@{}c@{}|@{}c@{}c@{}c@{}c@{}c@{}c@{}c@{}c@{}|@{}}
\cline{3-10}
&  &  &  &  & B &  &  &  &  \\ \cline{3-10}
&  & $su_{1}$ & \multicolumn{1}{|c}{$su_{2}$} & \multicolumn{1}{|c}{$su_{3}$} & 
\multicolumn{1}{|c}{$su_{4}$} & \multicolumn{1}{|c}{$su_{5}$} & 
\multicolumn{1}{|c}{$su_{6}$} & \multicolumn{1}{|c}{$su_{7}$} & 
\multicolumn{1}{|c|}{$su_{8}$} \\ \hline
\multicolumn{1}{|c}{} & \multicolumn{1}{|c|}{$su_{1}$} & $0$ & 
\multicolumn{1}{|c}{$0$} & \multicolumn{1}{|c}{$-su_{4}$} & \multicolumn{1}{|c}{%
$su_{3}$} & \multicolumn{1}{|c}{$-su_{6}$} & \multicolumn{1}{|c}{$su_{5}$} & 
\multicolumn{1}{|c}{$2su_{8}$} & \multicolumn{1}{|c|}{$-2su_{7}$} \\ 
\cline{2-3}\cline{2-10}
\multicolumn{1}{|c}{} & \multicolumn{1}{|c|}{$su_{2}$} & $0$ & 
\multicolumn{1}{|c}{$0$} & \multicolumn{1}{|c}{$su_{4}$} & \multicolumn{1}{|c}{$%
-su_{3}$} & \multicolumn{1}{|c}{$-2su_{6}$} & \multicolumn{1}{|c}{$2su_{5}$} & 
\multicolumn{1}{|c}{$su_{8}$} & \multicolumn{1}{|c|}{$-su_{7}$} \\ \cline{2-10}
\multicolumn{1}{|c}{} & \multicolumn{1}{|c|}{$su_{3}$} & $su_{4}$ & 
\multicolumn{1}{|c}{$-su_{4}$} & \multicolumn{1}{|c}{$0$} & \multicolumn{1}{|c}{%
$2(su_{2}-su_{1})$} & \multicolumn{1}{|c}{$-su_{7}$} & \multicolumn{1}{|c}{$su_{8}$} & 
\multicolumn{1}{|c}{$su_{5}$} & \multicolumn{1}{|c|}{$-su_{6}$} \\ \cline{2-10}
\multicolumn{1}{|c}{A} & \multicolumn{1}{|c|}{$su_{4}$} & $-su_{3}$ & 
\multicolumn{1}{|c}{$su_{3}$} & \multicolumn{1}{|c}{-$2(su_{2}-su_{1})$} & 
\multicolumn{1}{|c}{$0$} & \multicolumn{1}{|c}{$su_{8}$} & \multicolumn{1}{|c}{$%
su_{7}$} & \multicolumn{1}{|c}{$-su_{6}$} & \multicolumn{1}{|c|}{$-su_{5}$} \\ 
\cline{2-10}
\multicolumn{1}{|c}{} & \multicolumn{1}{|c|}{$su_{5}$} & $su_{6}$ & 
\multicolumn{1}{|c}{$2su_{6}$} & \multicolumn{1}{|c}{$su_{7}$} & 
\multicolumn{1}{|c}{$-su_{8}$} & \multicolumn{1}{|c}{$0$} & \multicolumn{1}{|c}{
-$2su_{2}$} & \multicolumn{1}{|c}{$-su_{3}$} & \multicolumn{1}{|c|}{$su_{4}$} \\ 
\cline{2-10}
\multicolumn{1}{|c}{} & \multicolumn{1}{|c|}{$su_{6}$} & $-su_{5}$ & 
\multicolumn{1}{|c}{$-2su_{5}$} & \multicolumn{1}{|c}{$-su_{8}$} & 
\multicolumn{1}{|c}{$-su_{7}$} & \multicolumn{1}{|c}{$2su_{2}$} & 
\multicolumn{1}{|c}{$0$} & \multicolumn{1}{|c}{$su_{4}$} & \multicolumn{1}{|c|}{%
$su_{3}$} \\ \cline{2-10}
\multicolumn{1}{|c}{} & \multicolumn{1}{|c|}{$su_{7}$} & $-2su_{8}$ & 
\multicolumn{1}{|c}{$-su_{8}$} & \multicolumn{1}{|c}{$-su_{5}$} & 
\multicolumn{1}{|c}{$su_{6}$} & \multicolumn{1}{|c}{$su_{3}$} & 
\multicolumn{1}{|c}{$-su_{4}$} & \multicolumn{1}{|c}{$0$} & 
\multicolumn{1}{|c|}{$2su_{1}$} \\ \cline{2-10}
\multicolumn{1}{|c}{} & \multicolumn{1}{|c|}{$su_{8}$} & $2su_{7}$ & 
\multicolumn{1}{|c}{$su_{7}$} & \multicolumn{1}{|c}{$su_{6}$} & 
\multicolumn{1}{|c}{$su_{5}$} & \multicolumn{1}{|c}{$-su_{4}$} & 
\multicolumn{1}{|c}{$-su_{3}$} & \multicolumn{1}{|c}{$-2su_{1}$} & 
\multicolumn{1}{|c|}{$0$} \\ \hline
\end{tabular}
}
\end{flushleft}

\setlength{\tabcolsep}{0.6mm} 
\begin{flushleft}
{\fontsize{8pt}{10pt} \selectfont%
\begin{tabular}{@{}c@{}c@{}|@{}c@{}c@{}c@{}c@{}c@{}c@{}c@{}c@{}|@{}}
\cline{3-3}\cline{3-10}
&  &  &  &  & B &  &  &  &  \\ \cline{3-10}
&  & $L_{1}$ & \multicolumn{1}{|c}{$(-L_{6})$} & \multicolumn{1}{|c}{$L_{4}$} & 
\multicolumn{1}{|c}{$(-L_{5})$} & \multicolumn{1}{|c}{$L_{7}$} & 
\multicolumn{1}{|c}{$L_{8}$} & \multicolumn{1}{|c}{$L_{2}$} & \multicolumn{1}{|c|}{%
$L_{3}$} \\ \cline{1-3}\cline{3-10}
\multicolumn{1}{|c|}{} & \multicolumn{1}{|c|}{$L_{1}$} & $0$ & 
\multicolumn{1}{|c}{$0$} & \multicolumn{1}{|c}{$-(-L_{5})$} & 
\multicolumn{1}{|c}{$L_{4}$} & \multicolumn{1}{|c}{$-L_{8}$} & \multicolumn{1}{|c}{%
$L_{7}$} & \multicolumn{1}{|c}{$2L_{3}$} & \multicolumn{1}{|c|}{$-2L_{2}$} \\ 
\cline{2-10}
\multicolumn{1}{|c|}{} & \multicolumn{1}{|c|}{$(-L_{6})$} & $0$ & 
\multicolumn{1}{|c}{$0$} & \multicolumn{1}{|c}{$(-L_{5})$} & 
\multicolumn{1}{|c}{$-L_{4}$} & \multicolumn{1}{|c}{$-2L_{8}$} & 
\multicolumn{1}{|c}{$2L_{7}$} & \multicolumn{1}{|c}{$L_{3}$} & 
\multicolumn{1}{|c|}{$-L_{2}$} \\ \cline{2-10}
\multicolumn{1}{|c|}{} & \multicolumn{1}{|c|}{$L_{4}$} & $(-L_{5})$ & 
\multicolumn{1}{|c}{$-(-L_{5})$} & \multicolumn{1}{|c}{$0$} & 
\multicolumn{1}{|c}{$2((-L_{6})-L_{1})$} & \multicolumn{1}{|c}{$-L_{2}$} & 
\multicolumn{1}{|c}{$L_{3}$} & \multicolumn{1}{|c}{$L_{7}$} & \multicolumn{1}{|c|}{%
$-L_{8}$} \\ \cline{2-10}
\multicolumn{1}{|c|}{A} & \multicolumn{1}{|c|}{$(-L_{5})$} & $-L_{4}$ & 
\multicolumn{1}{|c}{$L_{4}$} & \multicolumn{1}{|c}{$-2((-L_{6})-L_{1})$} & 
\multicolumn{1}{|c}{$0$} & \multicolumn{1}{|c}{$L_{3}$} & \multicolumn{1}{|c}{$%
L_{2}$} & \multicolumn{1}{|c}{$-L_{8}$} & \multicolumn{1}{|c|}{$-L_{7}$} \\ 
\cline{2-10}
\multicolumn{1}{|c|}{} & \multicolumn{1}{|c|}{$L_{7}$} & $L_{8}$ & 
\multicolumn{1}{|c}{$2L_{8}$} & \multicolumn{1}{|c}{$L_{2}$} & \multicolumn{1}{|c}{%
$-L_{3}$} & \multicolumn{1}{|c}{$0$} & \multicolumn{1}{|c}{$-2(-L_{6})$} & 
\multicolumn{1}{|c}{$-L_{4}$} & \multicolumn{1}{|c|}{$(-L_{5})$} \\ \cline{2-10}
\multicolumn{1}{|c|}{} & \multicolumn{1}{|c|}{$L_{8}$} & $-L_{7}$ & 
\multicolumn{1}{|c}{$-2L_{7}$} & \multicolumn{1}{|c}{$-L_{3}$} & 
\multicolumn{1}{|c}{$-L_{2}$} & \multicolumn{1}{|c}{$2(-L_{6})$} & 
\multicolumn{1}{|c}{$0$} & \multicolumn{1}{|c}{$(-L_{5})$} & 
\multicolumn{1}{|c|}{$L_{4}$} \\ \cline{2-10}
\multicolumn{1}{|c|}{} & \multicolumn{1}{|c|}{$L_{2}$} & $-2L_{3}$ & 
\multicolumn{1}{|c}{$-L_{3}$} & \multicolumn{1}{|c}{$-L_{7}$} & 
\multicolumn{1}{|c}{$L_{8}$} & \multicolumn{1}{|c}{$L_{4}$} & \multicolumn{1}{|c}{$%
-(-L_{5})$} & \multicolumn{1}{|c}{$0$} & \multicolumn{1}{|c|}{$2L_{1}$} \\ 
\cline{2-10}
\multicolumn{1}{|c|}{} & \multicolumn{1}{|c|}{$L_{3}$} & $2L_{2}$ & 
\multicolumn{1}{|c}{$L_{2}$} & \multicolumn{1}{|c}{$L_{8}$} & \multicolumn{1}{|c}{$%
L_{7}$} & \multicolumn{1}{|c}{$-(-L_{5})$} & \multicolumn{1}{|c}{$-L_{4}$} & 
\multicolumn{1}{|c}{$-2L_{1}$} & \multicolumn{1}{|c|}{$0$} \\ \hline
\end{tabular}
}
\end{flushleft}

By corresponding $L_{1}$ to $su_{1},-L_{6}$ to $su_{2},L_{4}$ to $su_{3},-L_{5}$ to $su_{4},L_{7}$ to $%
su_{5},L_{8}$ to $su_{6}$, $L_{2}$ to $su_{7}$, and $L_{3}$ to $su_{8}$, \gr$^{1}$\gs\gu$(3)$ is isomorphic to \gs\gu$(3)$.
\ \ \ \ \emph{Q.E.D.}

\bigskip

\subsection{\gs\gp$(3)$, \gs\gu$(3)$, and \gs\go$(9)$ of \gr$_{4}^{\C}$}

\ \ \ \ \emph{Definition 19.5.}  We denote the quaternion by $H$, and let we put 

\ \ \ \ \ \ \gs\gp$(3)=\left\{ X\in M(3\times 3,H)|X^{\ast }+X=0\right\} $, 

\noindent
where * is a matrix transpose and quaternion conjugate.

\noindent
Then $X\in$ \gs\gp$(3)$ is expressed by

$\left( 
\begin{array}{ccc}
\chi _{1} & x_{3} & -\overline{x_{2}} \\ 
-\overline{x_{3}} & \chi _{2} & x_{1} \\ 
x_{2} & -\overline{x_{1}} & \chi _{3}%
\end{array}%
\right) $ ($\overline{x}$ means quaternion conjugate of $x$).

\noindent
Let we express X and the elements of X as follows:

$X=\chi _{1}E_{1}+\chi _{2}E_{2}+\chi _{3}E_{3}+A_{1}(x_{1})+A_{2}(x_{2})+A_{3}(x_{3}),$

$\chi _{1}=a_{1}e_{1}+b_{1}e_{2}+c_{1}e_{3},$

$\chi _{2}=a_{2}e_{1}+b_{2}e_{2}+c_{2}e_{3},$

$\chi _{3}=a_{3}e_{1}+b_{3}e_{2}+c_{3}e_{3},$

$x_{1}=x_{10}e_{0}+x_{11}e_{1}+x_{12}e_{2}+x_{13}e_{3},$

$x_{2}=x_{20}e_{0}+x_{21}e_{1}+x_{22}e_{2}+x_{23}e_{3},$

$x_{3}=x_{30}e_{0}+x_{31}e_{1}+x_{32}e_{2}+x_{33}e_{3}$, ($\{e_{0},e_{1},e_{2},e_{3}\} $ are bases of $H$).

\ \ \ \ \ $E_{1}=\left( 
\begin{array}{ccc}
1 & 0 & 0 \\ 
0 & 0 & 0 \\ 
0 & 0 & 0%
\end{array}%
\right) ,$ \ \ \ $E_{2}=\left( 
\begin{array}{ccc}
0 & 0 & 0 \\ 
0 & 1 & 0 \\ 
0 & 0 & 0%
\end{array}%
\right) ,$ \ \ \ \ $E_{3}=\left( 
\begin{array}{ccc}
0 & 0 & 0 \\ 
0 & 0 & 0 \\ 
0 & 0 & 1%
\end{array}%
\right) ,$

{\fontsize{9pt}{10pt} \selectfont%
$A_{1}(a)=\left( 
\begin{array}{ccc}
0 & 0 & 0 \\ 
0 & 0 & a \\ 
0 & -\overline{a} & 0%
\end{array}%
\right) ,$ $A_{2}(a)=\left( 
\begin{array}{ccc}
0 & 0 & -\overline{a} \\ 
0 & 0 & 0 \\ 
a & 0 & 0%
\end{array}%
\right) ,$ $A_{3}(a)=\left( 
\begin{array}{ccc}
0 & a & 0 \\ 
-\overline{a} & 0 & 0 \\ 
0 & 0 & 0%
\end{array}%
\right) $,
}

$a_{k},b_{k},c_{k},x_{kl}\in \R,(k=1,2,3,l=0,1,2,3), a \in H.$
\bigskip

\emph{Definition 19.6.}  We define the following Lie subalgebras of \gr$_{4}^{\C}$.

\gr\gs\gp$(3)=\{ $Lie algebra generated by $K_{n}(1\leq n\leq
6),S_{11},S_{12},S_{13},$

\ \ \ \ \ \ \ \ \ \ \ \ \ \ $Um_{kl}(1\leq k\leq 3,0\leq l \leq 3)$ over real numbers $\} $.

(cf. I.Yokota\cite[\emph{Theorem2.11.2. Remark.}]{Yokota1}), 

\gr$^{2}$\gs\gu$(3)=\{$Lie algebra generated by  $K_{11},K_{12},Um10,Um11,Um20,Um21$,

\ \ \ \ \ \ \ \ \ \ \ \ \ \ \ \ $Um30,Um31$ over complex numbers$\}$.

(cf. I.Yokota\cite[\emph{Theorem2.12.2. Remark 1.}]{Yokota1}).

\bigskip

\emph{Lemma 19.7.} \ We have the following Lie bracket operations.

$[\chi _{1}E_{1},A_{1}(x_{1})]=0,$

$[\chi _{1}E_{1},A_{2}(x_{2})]=-A_{2}(x_{2}\chi _{1}),$

\ \ $x_{2}\chi _{1}=-(-a_{1}x_{21}-b_{1}x_{22}-c_{1}x_{23})e_{0}-(a_{1}x_{20}-b_{1}x_{23}+c_{1}x_{22})e_{1}$

\ \ \ \ \ \ \ \ \ \ \ \ \ $-(a_{1}x_{23}+b_{1}x_{20}-c_{1}x_{21})e_{2}-(-a_{1}x_{22}+b_{1}x_{21}+c_{1}x_{20})e_{3},$

$[\chi _{1}E_{1},A_{3}(x_{3})]=A_{3}(\chi _{1}x_{3}),$

$\ \ \ \chi _{1}x_{3}=(-a_{1}x_{31}-b_{1}x_{32}-c_{1}x_{33})e_{0}+(a_{1}x_{30}+b_{1}x_{33}-c_{1}x_{32})e_{1}$

\ \ \ \ \ \ \ \ \ \ \ \ $+(-a_{1}x_{33}+b_{1}x_{30}+c_{1}x_{31})e_{2}+(a_{1}x_{32}-b_{1}x_{31}+c_{1}x_{30})e_{3},$

$[\chi _{2}E_{2},A_{1}(x_{1})]=A_{1}(\chi _{2}x_{1}),$

\ \ \ \ $\chi _{2}x_{1}=(-a_{2}x_{11}-b_{2}x_{12}-c_{2}x_{13})e_{0}+(a_{2}x_{10}+b_{2}x_{13}-c_{2}x_{12})e_{1}$

\ \ \ \ \ \ \ \ \ \ \ \ \ \ $+(-a_{2}x_{13}+b_{2}x_{10}+c_{2}x_{11})e_{2}+(a_{2}x_{12}-b_{2}x_{11}+c_{2}x_{10})e_{3},$

$[\chi _{2}E_{2},A_{2}(x_{2})]=0,$

$[\chi _{2}E_{2},A_{3}(x_{3})]=-A_{3}(x_{3}\chi _{2}),$

$\ \ \ x_{3}\chi _{2}=-(-a_{2}x_{31}-b_{2}x_{32}-c_{2}x_{33})e_{0}-(a_{2}x_{30}-b_{2}x_{33}+c_{2}x_{32})e_{1}$

\ \ \ \ \ \ \ \ \ \ \ \ \ \ $-(a_{2}x_{33}+b_{2}x_{30}-c_{2}x_{31})e_{2}-(-a_{2}x_{32}+b_{2}x_{31}+c_{2}x_{30})e_{3},$

$[\chi _{3}E_{3},A_{1}(x_{1})]=-A_{1}(x_{1}\chi _{3}),$

$\ \ x_{1}\chi _{3}=-(-a_{3}x_{11}-b_{3}x_{12}-c_{3}x_{13})e_{0}-(a_{3}x_{10}-b_{3}x_{13}+c_{3}x_{12})e_{1}$

\ \ \ \ \ \ \ \ \ \ \ \ \ $-(a_{3}x_{13}+b_{3}x_{10}-c_{3}x_{11})e_{2}-(-a_{3}x_{12}+b_{3}x_{11}+c_{3}x_{10})e_{3},$

$[\chi _{3}E_{3},A_{2}(x_{2})]=A_{2}(\chi _{3}x_{2}),$

$\ \ \chi _{3}x_{2}=(-a_{3}x_{21}-b_{3}x_{22}-c_{3}x_{23})e_{0}+(a_{3}x_{20}+b_{3}x_{23}-c_{3}x_{22})e_{1}$

\ \ \ \ \ \ \ \ \ \ \ \ $+(-a_{3}x_{23}+b_{3}x_{20}+c_{3}x_{21})e_{2}+(a_{3}x_{22}-b_{3}x_{21}+c_{3}x_{20})e_{3},$

$[\chi _{3}E_{3},A_{3}(x_{3})]=0,$

$[A_{1}(x_{1}),A_{2}(x_{2})]=A_{3}(-\overline{x_{1}x_{2}})$,

$[A_{2}(x_{2}),A_{3}(x_{3})]=A_{1}(-\overline{x_{2}x_{3}})$,

$[A_{3}(x_{3}),A_{1}(x_{1})]=A_{2}(-\overline{x_{3}x_{1}})$,

$\ \ -\overline{x_{i}x_{j}}=(-x_{i0}x_{j0}+x_{i1}x_{j1}+x_{i2}x_{j2}+x_{i3}x_{j3})e_{0}$

\ \ \ \ \ \ \ \ \ \ \ \ \ \ $+(x_{i0}x_{j1}+x_{i1}x_{j0}+x_{i2}x_{j3}-x_{i3}x_{j2})e_{1}$

\ \ \ \ \ \ \ \ \ \ \ \ \ \ $+(x_{i0}x_{j2}-x_{i1}x_{j3}+x_{i2}x_{j0}+x_{i3}x_{j1})e_{2}$

\ \ \ \ \ \ \ \ \ \ \ \ \ \ $+(x_{i0}x_{j3}+x_{i1}x_{j2}-x_{i2}x_{j1}+x_{i3}x_{j0})e_{3}$.

\bigskip

\emph{Proof. \ }We have the above Lie bracket operations by
direct calculation.

\emph{Q.E.D.}

\bigskip

\emph{Lemma 19.8.} \ We have the following Lie bracket operations.

(1) Let we put $N_{1},N_{2}$, and $N_{3}$:

\ \ \ \ $N_{1}=a_{1}S_{11}+b_{1}S_{12}+c_{1}S_{13},$
$N_{2}=a_{2}K_{1}+b_{2}K_{2}+c_{2}K_{3},$

\ \ \ \ $N_{3}=a_{3}K_{4}+b_{3}K_{5}+c_{3}K_{6},$

$[N_{i},N_{j}]=0,\ (i,j=1,2,3)$.

(2)The following tables show the operation results of $[A,B]$.

(2.1) tables of bracket operations $[K_{i},Um_{kl}]$

\begin{tabular}{|cc|cccc|}
\hline
& \multicolumn{1}{c|}{} &  & $B$ &  &  \\ \cline{3-3}\cline{3-6}
& \multicolumn{1}{c|}{} & $Um_{10}$ & \multicolumn{1}{|c}{$Um_{11}$} & 
\multicolumn{1}{|c}{$Um_{12}$} & \multicolumn{1}{|c|}{$Um_{13}$} \\ \hline
& \multicolumn{1}{|c|}{$K_{1}$} & $Um_{11}$ & \multicolumn{1}{|c}{$-Um_{10}$} & 
\multicolumn{1}{|c}{$Um_{13}$} & \multicolumn{1}{|c|}{$-Um_{12}$} \\ \cline{2-6}
& \multicolumn{1}{|c|}{$K_{2}$} & $Um_{12}$ & \multicolumn{1}{|c}{$-Um_{13}$} & 
\multicolumn{1}{|c}{$-Um_{10}$} & \multicolumn{1}{|c|}{$Um_{11}$} \\ \cline{2-6}
$A$ & \multicolumn{1}{|c|}{$K_{3}$} & $Um_{13}$ & \multicolumn{1}{|c}{$Um_{12}$} & 
\multicolumn{1}{|c}{$-Um_{11}$} & \multicolumn{1}{|c|}{$-Um_{10}$} \\ \cline{2-6}
& \multicolumn{1}{|c|}{$K_{4}$} & $-Um_{11}$ & \multicolumn{1}{|c}{$Um_{10}$} & 
\multicolumn{1}{|c}{$Um_{13}$} & \multicolumn{1}{|c|}{$-Um_{12}$} \\ \cline{2-6}
& \multicolumn{1}{|c|}{$K_{5}$} & $-Um_{12}$ & \multicolumn{1}{|c}{$-Um_{13}$} & 
\multicolumn{1}{|c}{$Um_{10}$} & \multicolumn{1}{|c|}{$Um_{11}$} \\ \cline{2-6}
& \multicolumn{1}{|c|}{$K_{6}$} & $-Um_{13}$ & \multicolumn{1}{|c}{$Um_{12}$} & 
\multicolumn{1}{|c}{$-Um_{11}$} & \multicolumn{1}{|c|}{$Um_{10}$} \\ \hline
\end{tabular}

\begin{tabular}{|cc|cccc|}
\hline
& \multicolumn{1}{c|}{} &  & $B$ &  &  \\ \cline{3-3}\cline{3-6}
& \multicolumn{1}{c|}{} & $Um_{20}$ & \multicolumn{1}{|c}{$Um_{21}$} & 
\multicolumn{1}{|c}{$Um_{22}$} & \multicolumn{1}{|c|}{$Um_{23}$} \\ \hline
& \multicolumn{1}{|c|}{$K_{1}$} & $0$ & \multicolumn{1}{|c}{$0$} & 
\multicolumn{1}{|c}{0} & \multicolumn{1}{|c|}{0} \\ \cline{2-6}
& \multicolumn{1}{|c|}{$K_{2}$} & 0 & \multicolumn{1}{|c}{0} & 
\multicolumn{1}{|c}{0} & \multicolumn{1}{|c|}{0} \\ \cline{2-6}
$A$ & \multicolumn{1}{|c|}{$K_{3}$} & 0 & \multicolumn{1}{|c}{0} & 
\multicolumn{1}{|c}{0} & \multicolumn{1}{|c|}{0} \\ \cline{2-6}
& \multicolumn{1}{|c|}{$K_{4}$} & $Um_{21}$ & \multicolumn{1}{|c}{$-Um_{20}$} & 
\multicolumn{1}{|c}{$Um_{23}$} & \multicolumn{1}{|c|}{$-Um_{22}$} \\ \cline{2-6}
& \multicolumn{1}{|c|}{$K_{5}$} & $Um_{22}$ & \multicolumn{1}{|c}{$-Um_{23}$} & 
\multicolumn{1}{|c}{$-Um_{20}$} & \multicolumn{1}{|c|}{$Um_{21}$} \\ \cline{2-6}
& \multicolumn{1}{|c|}{$K_{6}$} & $Um_{23}$ & \multicolumn{1}{|c}{$Um_{22}$} & 
\multicolumn{1}{|c}{$-Um_{21}$} & \multicolumn{1}{|c|}{$-Um_{20}$} \\ \hline
\end{tabular}

\begin{tabular}{|cc|cccc|}
\hline
& \multicolumn{1}{c|}{} &  & $B$ &  &  \\ \cline{3-6}
& \multicolumn{1}{c|}{} & $Um_{30}$ & \multicolumn{1}{|c}{$Um_{31}$} & 
\multicolumn{1}{|c}{$Um_{32}$} & \multicolumn{1}{|c|}{$Um_{33}$} \\ \hline
& \multicolumn{1}{|c|}{$K_{1}$} & $-Um_{31}$ & \multicolumn{1}{|c}{$Um_{30}$} & 
\multicolumn{1}{|c}{$Um_{33}$} & \multicolumn{1}{|c|}{$-Um_{32}$} \\ \cline{2-6}
& \multicolumn{1}{|c|}{$K_{2}$} & $-Um_{32}$ & \multicolumn{1}{|c}{$-Um_{33}$} & 
\multicolumn{1}{|c}{$Um_{30}$} & \multicolumn{1}{|c|}{$Um_{31}$} \\ \cline{2-6}
$A$ & \multicolumn{1}{|c|}{$K_{3}$} & $-Um_{33}$ & \multicolumn{1}{|c}{$Um_{32}$} & 
\multicolumn{1}{|c}{$-Um_{31}$} & \multicolumn{1}{|c|}{$Um_{30}$} \\ \cline{2-6}
& \multicolumn{1}{|c|}{$K_{4}$} & $0$ & \multicolumn{1}{|c}{$0$} & 
\multicolumn{1}{|c}{$0$} & \multicolumn{1}{|c|}{$0$} \\ \cline{2-6}
& \multicolumn{1}{|c|}{$K_{5}$} & $0$ & \multicolumn{1}{|c}{$0$} & 
\multicolumn{1}{|c}{$0$} & \multicolumn{1}{|c|}{$0$} \\ \cline{2-6}
& \multicolumn{1}{|c|}{$K_{6}$} & $0$ & \multicolumn{1}{|c}{$0$} & 
\multicolumn{1}{|c}{$0$} & \multicolumn{1}{|c|}{$0$} \\ \hline
\end{tabular}

\bigskip

(2.2)tables of blacket operations $[S_{i},Um_{kl}]$

\begin{tabular}{|cc|cccc|}
\hline
& \multicolumn{1}{c|}{} &  & $B$ &  &  \\ \cline{3-6}
& \multicolumn{1}{c|}{} & $Um_{10}$ & \multicolumn{1}{|c}{$Um_{11}$} & 
\multicolumn{1}{|c}{$Um_{12}$} & \multicolumn{1}{|c|}{$Um_{13}$} \\ \hline
& \multicolumn{1}{|c|}{$S_{11}$} & $0$ & \multicolumn{1}{|c}{$0$} & 
\multicolumn{1}{|c}{$0$} & \multicolumn{1}{|c|}{$0$} \\ \cline{2-6}
$A$ & \multicolumn{1}{|c|}{$S_{12}$} & $0$ & \multicolumn{1}{|c}{$0$} & 
\multicolumn{1}{|c}{$0$} & \multicolumn{1}{|c|}{$0$} \\ \cline{2-6}
& \multicolumn{1}{|c|}{$S_{13}$} & $0$ & \multicolumn{1}{|c}{$0$} & 
\multicolumn{1}{|c}{$0$} & \multicolumn{1}{|c|}{$0$} \\ \hline
\end{tabular}

\begin{tabular}{|cc|cccc|}
\hline
& \multicolumn{1}{c|}{} &  & $B$ &  &  \\ \cline{3-6}
& \multicolumn{1}{c|}{} & $Um_{20}$ & \multicolumn{1}{|c}{$Um_{21}$} & 
\multicolumn{1}{|c}{$Um_{22}$} & \multicolumn{1}{|c|}{$Um_{23}$} \\ \hline
& \multicolumn{1}{|c|}{$S_{11}$} & $-Um_{23}$ & \multicolumn{1}{|c}{$Um_{22}$} & 
\multicolumn{1}{|c}{$-Um_{21}$} & \multicolumn{1}{|c|}{$Um_{20}$} \\ \cline{2-6}
$A$ & \multicolumn{1}{|c|}{$S_{12}$} & $-Um_{22}$ & \multicolumn{1}{|c}{$-Um_{23}$} & 
\multicolumn{1}{|c}{$Um_{20}$} & \multicolumn{1}{|c|}{$Um_{21}$} \\ \cline{2-6}
& \multicolumn{1}{|c|}{$S_{13}$} & $-Um_{21}$ & \multicolumn{1}{|c}{$Um_{20}$} & 
\multicolumn{1}{|c}{$Um_{23}$} & \multicolumn{1}{|c|}{$-Um_{22}$} \\ \hline
\end{tabular}

\begin{tabular}{|cc|cccc|}
\hline
& \multicolumn{1}{c|}{} &  & $B$ &  &  \\ \cline{3-6}
& \multicolumn{1}{c|}{} & $Um_{30}$ & \multicolumn{1}{|c}{$Um_{31}$} & 
\multicolumn{1}{|c}{$Um_{32}$} & \multicolumn{1}{|c|}{$Um_{33}$} \\ \hline
& \multicolumn{1}{|c|}{$S_{11}$} & $Um_{33}$ & \multicolumn{1}{|c}{$Um_{32}$} & 
\multicolumn{1}{|c}{$-Um_{31}$} & \multicolumn{1}{|c|}{$-Um_{30}$} \\ \cline{2-6}
$A$ & \multicolumn{1}{|c|}{$S_{12}$} & $Um_{32}$ & \multicolumn{1}{|c}{$-Um_{33}$} & 
\multicolumn{1}{|c}{$-Um_{30}$} & \multicolumn{1}{|c|}{$Um_{31}$} \\ \cline{2-6}
& \multicolumn{1}{|c|}{$S_{13}$} & $Um_{31}$ & \multicolumn{1}{|c}{$-Um_{30}$} & 
\multicolumn{1}{|c}{$Um_{33}$} & \multicolumn{1}{|c|}{$-Um_{32}$} \\ \hline
\end{tabular}

\bigskip

\emph{Proof. \ }We have the above Lie bracket operations by \emph{Lemma 13.5}, \emph{Lemma 14.3} and \emph{Lemma 14.4} . 
\ \ \ \emph{Q.E.D.}

\bigskip

\emph{Lemma 19.9.} \ \gr\gs\gp$(3)$ is isomorphic to \gs \gp$(3)$. 

\noindent
Moreover \gr\gs\gp$(3) \oplus $\gr$^{2}$\gs\gp$(1)$ is a direct sum Lie subalgebra of \gr$_{4}^{\C}$.
\bigskip

\emph{Proof. }\ We define a $\R$-linear mapping  $f:$\gs\gp$(3)\rightarrow $\gr\gs\gp$(3)$ by

$f(\chi _{1}E_{1})=a_{1}S_{11}+b_{1}S_{12}+c_{1}S_{13},$

$f(\chi _{2}E_{2})=a_{2}K_{1}+b_{2}K_{2}+c_{2}K_{3},$

$f(\chi _{3}E_{3})=a_{3}K_{4}+b_{3}K_{5}+c_{3}K_{6},$

$f(A_{1}(x_{1}))=\sum\limits_{0\leq l\leq 3}(x_{1l})Um_{1l}$,

$f(A_{2}(x_{2}))=\sum\limits_{0\leq l\leq 3}(x_{2l})Um_{2l}$,

$f(A_{3}(x_{3}))=\sum\limits_{0\leq l\leq 3}(x_{3l})Um_{3l}$.

\noindent
By \emph{Lemma 19.7}  and \emph{Lemma 19.8}, $f$ gives an isomorphism.

For $\forall x \in $\gr\gs\gp$(3)$,$\forall y \in $\gr$^{2}$\gs\gp$(1)$, we have $[x,y]=0$ by using Maxima.  
\ \ \ \ \ \ \emph{Q.E.D.}

\bigskip

\emph{Lemma 19.10.} \ \gr$^{2}$\gs\gu$(3)$ is isomorphic to \gs \gu$(3)$.

\noindent
Moreover \gr$^{1}$\gs\gu$(3) \oplus $\gr$^{2}$\gs\gu$(3)$ is a direct sum Lie subalgebra of \gr$_{4}^{\C}$.

\bigskip

\emph{Proof. }\ We have following Lie bracket operations by using Maxima.

\noindent
The following tables show the operation results of $[A,B]$.

\setlength{\tabcolsep}{0.6mm} 
\begin{flushleft}
{\fontsize{6pt}{10pt} \selectfont%
\begin{tabular}{@{}c@{}c@{}|@{}c@{}c@{}c@{}c@{}c@{}c@{}c@{}c@{}|@{}}
\cline{3-10}
&  &  &  &  & B &  &  &  &  \\ \cline{3-10}
&  & $K_{11}$ & \multicolumn{1}{|c}{$K_{12}$} & \multicolumn{1}{|c}{$2Um_{10}$} & 
\multicolumn{1}{|c}{$2Um_{11}$} & \multicolumn{1}{|c}{$2Um_{20}$} & 
\multicolumn{1}{|c}{$2Um_{21}$} & \multicolumn{1}{|c}{$2Um_{30}$} & 
\multicolumn{1}{|c|}{$2Um_{31}$} \\ \hline
\multicolumn{1}{|c}{} & \multicolumn{1}{|c|}{$K_{11}$} & $0$ & 
\multicolumn{1}{|c}{$0$} & \multicolumn{1}{|c}{$-2Um_{11}$} & \multicolumn{1}{|c}{%
$2Um_{10}$} & \multicolumn{1}{|c}{$-2Um_{21}$} & \multicolumn{1}{|c}{$2Um_{20}$} & 
\multicolumn{1}{|c}{$4Um_{31}$} & \multicolumn{1}{|c|}{$-4Um_{30}$} \\ 
\cline{2-3}\cline{2-10}
\multicolumn{1}{|c}{} & \multicolumn{1}{|c|}{$K_{12}$} & $0$ & 
\multicolumn{1}{|c}{$0$} & \multicolumn{1}{|c}{$2Um_{11}$} & \multicolumn{1}{|c}{$%
-2Um_{10}$} & \multicolumn{1}{|c}{$-4Um_{21}$} & \multicolumn{1}{|c}{$4Um_{20}$} & 
\multicolumn{1}{|c}{$2Um_{31}$} & \multicolumn{1}{|c|}{$-2Um_{30}$} \\ \cline{2-10}
\multicolumn{1}{|c}{} & \multicolumn{1}{|c|}{$2Um_{10}$} & $2Um_{11}$ & 
\multicolumn{1}{|c}{$-2Um_{11}$} & \multicolumn{1}{|c}{$0$} & \multicolumn{1}{|c}{%
$2(K_{12}-K_{11})$} & \multicolumn{1}{|c}{$-2Um_{30}$} & \multicolumn{1}{|c}{$2Um_{31}$} & 
\multicolumn{1}{|c}{$2Um_{20}$} & \multicolumn{1}{|c|}{$-2Um_{21}$} \\ \cline{2-10}
\multicolumn{1}{|c}{A} & \multicolumn{1}{|c|}{$2Um_{11}$} & $-2Um_{10}$ & 
\multicolumn{1}{|c}{$2Um_{10}$} & \multicolumn{1}{|c}{-$2(K_{12}-K_{11})$} & 
\multicolumn{1}{|c}{$0$} & \multicolumn{1}{|c}{$2Um_{31}$} & \multicolumn{1}{|c}{$%
2Um_{30}$} & \multicolumn{1}{|c}{$-2Um_{21}$} & \multicolumn{1}{|c|}{$-2Um_{20}$} \\ 
\cline{2-10}
\multicolumn{1}{|c}{} & \multicolumn{1}{|c|}{$2Um_{20}$} & $2Um_{21}$ & 
\multicolumn{1}{|c}{$4Um_{21}$} & \multicolumn{1}{|c}{$2Um_{30}$} & 
\multicolumn{1}{|c}{$-2Um_{31}$} & \multicolumn{1}{|c}{$0$} & \multicolumn{1}{|c}{
-$2K_{12}$} & \multicolumn{1}{|c}{$-2Um_{10}$} & \multicolumn{1}{|c|}{$2Um_{11}$} \\ 
\cline{2-10}
\multicolumn{1}{|c}{} & \multicolumn{1}{|c|}{$2Um_{21}$} & $-2Um_{20}$ & 
\multicolumn{1}{|c}{$-4Um_{20}$} & \multicolumn{1}{|c}{$-2Um_{31}$} & 
\multicolumn{1}{|c}{$-2Um_{30}$} & \multicolumn{1}{|c}{$2K_{12}$} & 
\multicolumn{1}{|c}{$0$} & \multicolumn{1}{|c}{$2Um_{11}$} & \multicolumn{1}{|c|}{%
$2Um_{10}$} \\ \cline{2-10}
\multicolumn{1}{|c}{} & \multicolumn{1}{|c|}{$2Um_{30}$} & $-4Um_{31}$ & 
\multicolumn{1}{|c}{$-2Um_{31}$} & \multicolumn{1}{|c}{$-2Um_{20}$} & 
\multicolumn{1}{|c}{$2Um_{21}$} & \multicolumn{1}{|c}{$2Um_{10}$} & 
\multicolumn{1}{|c}{$-2Um_{11}$} & \multicolumn{1}{|c}{$0$} & 
\multicolumn{1}{|c|}{$2K_{11}$} \\ \cline{2-10}
\multicolumn{1}{|c}{} & \multicolumn{1}{|c|}{$2Um_{31}$} & $4Um_{30}$ & 
\multicolumn{1}{|c}{$2Um_{30}$} & \multicolumn{1}{|c}{$2Um_{21}$} & 
\multicolumn{1}{|c}{$2Um_{20}$} & \multicolumn{1}{|c}{$-2Um_{11}$} & 
\multicolumn{1}{|c}{$-2Um_{10}$} & \multicolumn{1}{|c}{$-2K_{11}$} & 
\multicolumn{1}{|c|}{$0$} \\ \hline
\end{tabular}
}
\end{flushleft}

By  corresponding $K_{11}$ to $su_{1},K_{12}$ to $su_{2},2Um_{10}$ to $su_{3},2Um_{11}$ to $su_{4}$,
$2Um_{20}$ to $su_{5},2Um_{21}$ to $su_{6}$, $2Um_{30}$ to $su_{7}$, and $2Um_{31}$ to $su_{8}$, \gr$^{2}$\gs\gu$(3)$ is isomorphic to \gs\gu$(3)$.

For $\forall x \in $\gr$^{1}$\gs\gu$(3)$,$\forall y \in $\gr$^{2}$\gs\gu$(3)$, we have $[x,y]=0$ by using Maxima.  
\ \ \ \ \ \ \emph{Q.E.D.}

\bigskip

\emph{Lemma 19.11.}  Let's define the following vector space of \gr$_{4}^{\C}$.

\gr\gs\go$(9)= $\gr\gd$\bigoplus $\gR\gm$_{1}$,

\noindent
where \gR\gm$_{1}=\{Rm_{1j} \in $\gR\gm$ \mid  m_{1j} \in \R,(0\leq j\leq 7)\}$.
Then \gr\gs\go$(9)$ is a  Lie algebra.

\bigskip

\emph{Proof}  By the use of Maxima for Lie bracket operations, \gr\gs\go$(9)$ is a Lie algebra.
\ \ \ \ \ \ \emph{Q.E.D.}

\bigskip

\emph{Lemma 19.12.}  \gr\gs\go$(9)$ is isomorphic to \gs\go$(9)=\{X \in M(9\times9,\R) \mid$ $^{t}X + X = 0 \}$.

\bigskip

\emph{Proof.}  We define a $\R$-linear mapping $f:$\gr\gs\go$(9)\rightarrow $\gs\go$(9)$ by

$f(d_{ij}Ud_{ij})=d_{ij}D_{ij}$, $(0\leq i <j \leq 7)$,

$f(m_{1j}Um_{1j})=m_{1j}D_{j8}$, $(0 \leq j \leq 7)$,

\noindent
where $D_{ij}\in M(9\times9,\R)$, (see \emph{Definition 1.12}).

We have following Lie bracket opearatins by using Maxima.

$f([\sum\limits_{0\leq i<j \leq 7}d_{ij}Ud_{ij}+\sum\limits_{0\leq j\leq 7}m_{1j}Um_{1j},%
\sum\limits_{0\leq i<j \leq 7}h_{ij}Ud_{ij}+\sum\limits_{0\leq j\leq 7}n_{1j}Um_{1j}])$

$=[\sum\limits_{0\leq i<j \leq 7}d_{ij}D_{ij}+\sum\limits_{0\leq j\leq 7}m_{1j}D_{j8},%
\sum\limits_{0\leq i<j \leq 7}h_{ij}D_{ij}+\sum\limits_{0\leq j\leq 7}n_{1j}D_{j8}])$,

\noindent
$(d_{ij},m_{1j},h_{ij},n_{1j} \in \R)$.
So $f$ gives an isomorphism.\ \ \ \ \ \ \emph{Q.E.D.}

\bigskip

\subsection{\gs\go$(10)$ , \gu$(1)$ , \gs\gu$(6)$,  \gs\gu$(3)$, and \gs\gp$(4)$ of \gr$_{6}^{\C}$ }

\ \ \ \ \emph{Lemma 19.13.}  Let's define the following vector spaces of \gr$_{6}^{\C}$.

\gr\gs\go$(10)= $\gr\gd$\bigoplus $\gR\gm$_{1} \bigoplus (i*$\gR\gt$_{1}) \bigoplus (-i*$\gR\gt$_{4})$,

\gr\gu$(1)=\{ i\theta(2U\tau_{1}-U\tau_{2}) \mid \theta \in \R \}$,

\noindent
where \gR\gt$_{1}=\{Rt_{1j} \in $\gR\gt$ \mid  t_{1j} \in \R,(0\leq j\leq 7)\}$,
\gR\gt$_{4}=\{R\tau_{2} \in $\gR\gt$ \mid  \tau_{2} \in \R\}$.
Then \gr\gs\go$(10)$ and \gr\gu$(1)$ are Lie algebras.
\bigskip

\emph{Proof}  By the use of Maxima for Lie bracket operations, \gr\gs\go$(10)$ and \gr\gu$(1)$ are Lie algebras.
\ \ \ \ \ \ \emph{Q.E.D.}

\bigskip

\emph{Lemma 19.14.}  \gr\gs\go$(10)$ is isomorphic to \gs\go$(10)=\{X \in M(10\times10,\R) \mid$ $^{t}X + X = 0 \}$,
 and \gr\gu$(1)$ is isomorphic to \gu$(1)=\{i\theta \mid \theta \in \R \}$.
 
 \noindent
 And \gr\gs\go$(10) \oplus $\gr\gu$(1)$ is a direct sum Lie subalgebra of \gr$_{6}^{\C}$.

\bigskip

\emph{Proof.}  We define a $\R$-linear mapping $f:$\gr\gs\go$(10)\rightarrow $\gs\go$(10)$ by

$f(d_{ij}Ud_{ij})=d_{ij}D_{ij}$, $(0\leq i <j \leq 7)$,

$f(m_{1j}Um_{1j})=m_{1j}D_{j8}$, $(0 \leq j \leq 7)$,

$f(i*t_{1j}Ut_{1j})=t_{1j}D_{j9}$, $(0 \leq j \leq 7)$, $f(-i*\tau_{2}U\tau_{2})=\tau_{2}D_{89}$,

\noindent
where $D_{ij}\in M(10\times10,\R)$, (see \emph{Definition 1.12}).

We have following Lie bracket operations by using Maxima.

$f([\sum\limits_{0\leq i<j \leq 7}d_{ij}Ud_{ij}+\sum\limits_{0\leq j\leq 7}m_{1j}Um_{1j}%
+\ i*\sum\limits_{0\leq j\leq 7}t_{1j}Ut_{1j}-i*\tau_{2}U\tau_{2}$ ,

\ \ \ \ $\sum\limits_{0\leq i<j \leq 7}h_{ij}Ud_{ij}+\sum\limits_{0\leq j\leq 7}n_{1j}Um_{1j}%
+\ i*\sum\limits_{0\leq j\leq 7}s_{1j}Ut_{1j}-i*\tau_{3}U\tau_{2}])$

$=[\sum\limits_{0\leq i<j \leq 7}d_{ij}D_{ij}+\sum\limits_{0\leq j\leq 7}m_{1j}D_{j8}%
+\sum\limits_{0\leq j\leq 7}t_{1j}D_{j9}+\tau_{2}D_{89}$ ,

\ \ \ \ $\sum\limits_{0\leq i<j \leq 7}h_{ij}D_{ij}+\sum\limits_{0\leq j\leq 7}n_{1j}D_{j8}%
+\sum\limits_{0\leq j\leq 7}s_{1j}D_{j9}+\tau_{3}D_{89}])$,

\noindent
$(d_{ij},m_{1j},h_{ij},n_{1j},t_{1j},s_{1j},\tau_{2},\tau_{3} \in \R)$.
So $f$ gives an isomorphism.

Evidently \gr\gu$(1)$ is isomorphic to \gu$(1)$ .

For $\forall x \in $\gr\gs\go$(10)$,$\forall y \in $\gr\gu$(1)$, we have $[x,y]=0$ by using Maxima.  

\emph{Q.E.D.}

\bigskip

\emph{Lemma 19.15.}  Let's define the following vector space of \gr$_{6}^{\C}$.

\gr\gs\gu$(6)=$\gr\gn$  \bigoplus $\gR\gm$_{f} \bigoplus (i*$\gR\gt$_{\tau})\bigoplus (i*$\gR\gt$_{f})$,

\noindent
where 
\gr\gn$=\{ a_{1}S_{11}+b_{1}S_{12}+c_{1}S_{13}+a_{2}K_{1}+b_{2}K_{2}+c_{2}K_{3}+a_{3}K_{4}+b_{3}K_{5}+c_{3}K_{6} \mid a_{k},b_{k},c_{k} \in \R,%
(1\leq k\leq 3) \}$,
\gR\gm$_{f}=\{Rm_{kj} \in $\gR\gm$\ \mid  m_{kj} \in \R,(1\leq k\leq 3, 0\leq j\leq 3) \}$,
\gR\gt$_{\tau}=\{R\tau_{1}, R\tau_{2}\in $\gR\gt$\ \mid  \tau_{1}, \tau_{2}\in \R\}$,
\gR\gt$_{f}=\{Rt_{kj} \in $\gR\gt$\ \mid t_{kj} \in \R, (1\leq k\leq 3, 0\leq j\leq 3 \}$.
Then \gr\gs\gu$(6)$ is a 35-dimensional Lie algebra.
\bigskip

\emph{Proof}  By the use of Maxima for Lie bracket operations, \gr\gs\gu$(6)$ is a 35-dimensional Lie algebra.
\ \ \ \ \ \ \emph{Q.E.D.}

\bigskip

\emph{Definition 19.16.}  We denote $X \in $ \gs\gu$(6)=\{ X \in M(6 \times 6,\C) \mid \ X^{*}+X=0,tr(X)=0 \}$
($\ast$ means matrix transpose and complex conjugate) by following.

\noindent
$X=$
{\fontsize{8pt}{10pt} \selectfont%
$\left(%
\begin{array}{@{}llllll@{}}
ia_{1}$+$i\tau _{1} & $-$b_{1}$-$ic_{1} & m_{m5}$+$it_{t5} & $-$m_{m6}$-$it_{t6} & $-$%
\overline{m_{m3}}$+$i\overline{t_{t3}} & $-$m_{m4}$-$it_{t4} \\ 
b_{1}$-$ic_{1} & $-$ia_{1}$+$i\tau _{1} & \overline{m_{m6}}$+$i\overline{t_{t6}} & 
\overline{m_{m5}}$+$i\overline{t_{t5}} & \overline{m_{m4}}$-$i\overline{t_{t4}}
& m_{m3}$-$it_{t3} \\ 
$-$\overline{m_{m5}}$+$i\overline{t_{t5}} & $-$m_{m6}$-$it_{t6} & ia_{2}$+$i\tau _{2}
& $-$b_{2}$-$ic_{2} & m_{m1}$+$it_{t1} & $-$m_{m2}$-$it_{t2} \\ 
\overline{m_{m6}}$-$i\overline{t_{t6}} & m_{m5}$-$it_{t5} & b_{2}$-$ic_{2} & 
$-$ia_{2}$+$i\tau _{2} & \overline{m_{m2}}$+$i\overline{t_{t2}} & \overline{m_{m1}}%
$+$i\overline{t_{t1}} \\ 
m_{m3}$+$it_{t3} & $-$m_{m4}$-$it_{t4} & $-$\overline{m_{m1}}$+$i\overline{t_{t1}} & 
$-$m_{m2}$-$it_{t2} & ia_{3}$-$i\tau _{1}$-$i\tau _{2} & $-$b_{3}$-$ic_{3} \\ 
\overline{m_{m4}}$+$i\overline{t_{t4}} & \overline{m_{m3}}$+$i\overline{t_{t3}}
& \overline{m_{m2}}$-$i\overline{t_{t2}} & m_{m1}$-$it_{t1} & b_{3}$-$ic_{3} & 
$-$ia_{3}$-$i\tau _{1}$-$i\tau _{2}%
\end{array}%
\right)$
}

\noindent
($\overline{m}$ means complex conjugate of $m$) 

\noindent
where $a_{k},b_{k},c_{k},\tau_{1},\tau_{2}\in R,(1\leq k\leq 3),m_{mk},t_{tk}\in C,(1\leq k\leq 6)$.

(cf. I.Yokota\cite[\emph{Section 3.11.}]{Yokota1}).

\bigskip

\emph{Lemma 19.17.}  \gr\gs\gu$(6)$ is isomorphic to \gs\gu$(6)$. 

\noindent
Moreover \gr\gs\gu$(6) \oplus $\gr$^{2}$\gs\gp$(1)$ is a direct sum Lie subalgebra of \gr$_{6}^{\C}$. 
\bigskip

\emph{Proof}  We correspond the elements of \gr\gs\gu$(6)$ with the elements of \gs\gu$(6)$ as follows: 

$\left( a_{1}S_{11}+b_{1}S_{12}+c_{1}S_{13}\right) /\sqrt{2}\rightarrow 
\begin{tabular}{|c|c|}
\hline
$ia_{1}$ & $-b_{1}-ic_{1}$ \\ \hline
$b_{1}-ic_{1}$ & $-ia_{1}$ \\ \hline
\end{tabular}%
,$

$\left( a_{2}K_{1}+b_{2}K_{2}+c_{2}K_{3}\right) /\sqrt{2}\rightarrow 
\begin{tabular}{|c|c|}
\hline
$ia_{2}$ & $-b_{2}-ic_{2}$ \\ \hline
$b_{2}-ic_{2}$ & $-ia_{2}$ \\ \hline
\end{tabular}%
,$

$\left( a_{3}K_{4}+b_{3}K_{5}+c_{3}K_{6}\right) /\sqrt{2}\rightarrow 
\begin{tabular}{|c|c|}
\hline
$ia_{3}$ & $-b_{3}-ic_{3}$ \\ \hline
$b_{3}-ic_{3}$ & $-ia_{3}$ \\ \hline
\end{tabular}%
,$

$\left( Rm_{10}+Rm_{11}+Rm_{12}+Rm_{13}\right) \sqrt{2}\rightarrow 
m_{m1}=m_{10}+im_{11}$,

\ \ \ \ \ \ \ \ \ \ \ \ \ \ \ \ \ \ \ \ \ \ \ \ \ \ \ \ \ \ \ \ \ \ \ \ \ \ \ \ \ \ \ \ \ \ \ \ \ \ \ \ $m_{m2}=m_{12}+im_{13}$,

$\left( Rm_{20}+Rm_{21}+Rm_{22}+Rm_{23}\right) \sqrt{2}\rightarrow 
m_{m3}=m_{20}+im_{21}$,

\ \ \ \ \ \ \ \ \ \ \ \ \ \ \ \ \ \ \ \ \ \ \ \ \ \ \ \ \ \ \ \ \ \ \ \ \ \ \ \ \ \ \ \ \ \ \ \ \ \ \ \ $m_{m4}=m_{22}+im_{23}$,

$\left( Rm_{30}+Rm_{31}+Rm_{32}+Rm_{33}\right) \sqrt{2}\rightarrow 
m_{m5}=m_{30}+im_{31}$,

\ \ \ \ \ \ \ \ \ \ \ \ \ \ \ \ \ \ \ \ \ \ \ \ \ \ \ \ \ \ \ \ \ \ \ \ \ \ \ \ \ \ \ \ \ \ \ \ \ \ \ \ $m_{m6}=m_{32}+im_{33}$,

$i\left( Rt_{10}+Rt_{11}+Rt_{12}+Rt_{13}\right) \sqrt{2}\rightarrow 
\ \ \ \ it_{t1}=it_{10}-t_{11}$,

\ \ \ \ \ \ \ \ \ \ \ \ \ \ \ \ \ \ \ \ \ \ \ \ \ \ \ \ \ \ \ \ \ \ \ \ \ \ \ \ \ \ \ \ \ \ \ \ \ \ \ \ $it_{t2}=it_{12}-t_{13}$,

$i\left( Rt_{20}+Rt_{21}+Rt_{22}+Rt_{23}\right) \sqrt{2}\rightarrow 
\ \ \ \ it_{t3}=it_{20}-t_{21}$,

\ \ \ \ \ \ \ \ \ \ \ \ \ \ \ \ \ \ \ \ \ \ \ \ \ \ \ \ \ \ \ \ \ \ \ \ \ \ \ \ \ \ \ \ \ \ \ \ \ \ \ \ $it_{t4}=it_{22}-t_{23}$,

$i\left( Rt_{30}+Rt_{31}+Rt_{32}+Rt_{33}\right) \sqrt{2}\rightarrow 
\ \ \ \ it_{t5}=it_{30}-t_{31}$,

\ \ \ \ \ \ \ \ \ \ \ \ \ \ \ \ \ \ \ \ \ \ \ \ \ \ \ \ \ \ \ \ \ \ \ \ \ \ \ \ \ \ \ \ \ \ \ \ \ \ \ \ $it_{t6}=it_{32}-t_{33}$,

$i(R\tau _{1}+R\tau _{2})\sqrt{2}\rightarrow 
\begin{tabular}{|c|c|c|c|c|c|}
\hline
$i\tau _{1}$ &  &  &  &  &  \\ \hline
& $i\tau _{1}$ &  &  &  &  \\ \hline
&  & $i\tau _{2}$ &  &  &  \\ \hline
&  &  & $i\tau _{2}$ &  &  \\ \hline
&  &  &  & $-i\tau _{1}-i\tau _{2}$ &  \\ \hline
&  &  &  &  & $-i\tau _{1}-i\tau _{2}$ \\ \hline
\end{tabular}
$.

By the use of Maxima for Lie bracket operations, these  correspondences give an isomorphism of \gr\gs\gu$(6)$ to \gs\gu$(6)$.

For $\forall x \in $\gr\gs\gu$(6)$,$\forall y \in $\gr$^{2}$\gs\gp$(1)$, we have $[x,y]=0$ by using Maxima.  
\ \ \ \ \ \ \emph{Q.E.D.}

\bigskip

\emph{Definition 19.18.} \ We define the following elements and a Lie subalgebra of \gr$_{6}^{\C}.$

$T_{1}=iU\tau _{1}+(Ud_{01}-Ud_{23}-Ud_{45}-Ud_{67})/2$

$T_{2}=iU\tau _{2}+Ud_{01}$

$T_{3}=-iUt_{31}-Um_{30}$

$T_{4}=iUt_{30}-Um_{31}$

$T_{5}=iUt_{10}-Um_{11}$

$T_{6}=iUt_{11}+Um_{10}$

$T_{7}=iUt_{20}-Um_{21}$

$T_{8}=iUt_{21}+Um_{20}$

\gr$^{3}$\gs\gu$(3)=\{$Lie algebra generated by $T_{k}(1\leq k\leq 8)$ over complex numbers$\}$,

(cf. I.Yokota\cite[\emph{Theorem3.13.5. Remark 1. }]{Yokota1}).

\bigskip

\emph{Lemma 19.19.} \ \gr$^{3}$\gs\gu$(3)$ is isomorphic to \gs \gu$(3)$. Moreover \gr$^{1}$\gs\gu$(3)$ and \gr$^{3}$\gs\gu$(3)$ are commutative.  
\bigskip

\emph{Proof. }\ We have following Lie bracket operations by using Maxima.

\noindent
The following tables show the operation results of $[A,B]$.

\setlength{\tabcolsep}{1.5mm} 
\begin{flushleft}
{\fontsize{8pt}{10pt} \selectfont%
\begin{tabular}{@{}c@{}c@{}|@{}c@{}c@{}c@{}c@{}c@{}c@{}c@{}c@{}|@{}}
\cline{3-10}
&  &  &  &  & B &  &  &  &  \\ \cline{3-10}
&  & $T_{1}$ & \multicolumn{1}{|c}{$T_{2}$} & \multicolumn{1}{|c}{$T_{3}$} & 
\multicolumn{1}{|c}{$T_{4}$} & \multicolumn{1}{|c}{$T_{5}$} & 
\multicolumn{1}{|c}{$T_{6}$} & \multicolumn{1}{|c}{$T_{7}$} & 
\multicolumn{1}{|c|}{$T_{8}$} \\ \hline
\multicolumn{1}{|c}{} & \multicolumn{1}{|c|}{$T_{1}$} & $0$ & 
\multicolumn{1}{|c}{$0$} & \multicolumn{1}{|c}{$-T_{4}$} & \multicolumn{1}{|c}{%
$T_{3}$} & \multicolumn{1}{|c}{$-T_{6}$} & \multicolumn{1}{|c}{$T_{5}$} & 
\multicolumn{1}{|c}{$2T_{8}$} & \multicolumn{1}{|c|}{$-2T_{7}$} \\ 
\cline{2-3}\cline{2-10}
\multicolumn{1}{|c}{} & \multicolumn{1}{|c|}{$T_{2}$} & $0$ & 
\multicolumn{1}{|c}{$0$} & \multicolumn{1}{|c}{$T_{4}$} & \multicolumn{1}{|c}{$%
-T_{3}$} & \multicolumn{1}{|c}{$-2T_{6}$} & \multicolumn{1}{|c}{$2T_{5}$} & 
\multicolumn{1}{|c}{$T_{8}$} & \multicolumn{1}{|c|}{$-T_{7}$} \\ \cline{2-10}
\multicolumn{1}{|c}{} & \multicolumn{1}{|c|}{$T_{3}$} & $T_{4}$ & 
\multicolumn{1}{|c}{$-T_{4}$} & \multicolumn{1}{|c}{$0$} & \multicolumn{1}{|c}{%
$2(T_{2}-T_{1})$} & \multicolumn{1}{|c}{$-T_{7}$} & \multicolumn{1}{|c}{$T_{8}$} & 
\multicolumn{1}{|c}{$T_{5}$} & \multicolumn{1}{|c|}{$-T_{6}$} \\ \cline{2-10}
\multicolumn{1}{|c}{A} & \multicolumn{1}{|c|}{$T_{4}$} & $-T_{3}$ & 
\multicolumn{1}{|c}{$T_{3}$} & \multicolumn{1}{|c}{-$2(T_{2}-T_{1})$} & 
\multicolumn{1}{|c}{$0$} & \multicolumn{1}{|c}{$T_{8}$} & \multicolumn{1}{|c}{$%
T_{7}$} & \multicolumn{1}{|c}{$-T_{6}$} & \multicolumn{1}{|c|}{$-T_{5}$} \\ 
\cline{2-10}
\multicolumn{1}{|c}{} & \multicolumn{1}{|c|}{$T_{5}$} & $T_{6}$ & 
\multicolumn{1}{|c}{$2T_{6}$} & \multicolumn{1}{|c}{$T_{7}$} & 
\multicolumn{1}{|c}{$-T_{8}$} & \multicolumn{1}{|c}{$0$} & \multicolumn{1}{|c}{
-$2T_{2}$} & \multicolumn{1}{|c}{$-T_{3}$} & \multicolumn{1}{|c|}{$T_{4}$} \\ 
\cline{2-10}
\multicolumn{1}{|c}{} & \multicolumn{1}{|c|}{$T_{6}$} & $-T_{5}$ & 
\multicolumn{1}{|c}{$-2T_{5}$} & \multicolumn{1}{|c}{$-T_{8}$} & 
\multicolumn{1}{|c}{$-T_{7}$} & \multicolumn{1}{|c}{$2T_{2}$} & 
\multicolumn{1}{|c}{$0$} & \multicolumn{1}{|c}{$T_{4}$} & \multicolumn{1}{|c|}{%
$T_{3}$} \\ \cline{2-10}
\multicolumn{1}{|c}{} & \multicolumn{1}{|c|}{$T_{7}$} & $-2T_{8}$ & 
\multicolumn{1}{|c}{$-T_{8}$} & \multicolumn{1}{|c}{$-T_{5}$} & 
\multicolumn{1}{|c}{$T_{6}$} & \multicolumn{1}{|c}{$T_{3}$} & 
\multicolumn{1}{|c}{$-T_{4}$} & \multicolumn{1}{|c}{$0$} & 
\multicolumn{1}{|c|}{$2T_{1}$} \\ \cline{2-10}
\multicolumn{1}{|c}{} & \multicolumn{1}{|c|}{$T_{8}$} & $2T_{7}$ & 
\multicolumn{1}{|c}{$T_{7}$} & \multicolumn{1}{|c}{$T_{6}$} & 
\multicolumn{1}{|c}{$T_{5}$} & \multicolumn{1}{|c}{$-T_{4}$} & 
\multicolumn{1}{|c}{$-T_{3}$} & \multicolumn{1}{|c}{$-2T_{1}$} & 
\multicolumn{1}{|c|}{$0$} \\ \hline
\end{tabular}
}
\end{flushleft}

Similar to \emph{Lemma 19.4}, by corresponding $T_{1}$ to $su_{1},T_{2}$ to $su_{2},T_{3}$ to $su_{3},T_{4}$ to $su_{4},T_{5}$ to $%
su_{5},T_{6}$ to $su_{6}$, $T_{7}$ to $su_{7}$, and $T_{8}$ to $su_{8}$, \gr$^{3}$\gs\gu$(3)$ is isomorphic to \gs\gu$(3)$.

For $\forall x \in $\gr$^{1}$\gs\gu$(3)$,$\forall y \in $\gr$^{3}$\gs\gu$(3)$, we have $[x,y]=0$ by using Maxima.  
\ \ \ \ \ \ \emph{Q.E.D.}

\bigskip

\emph{Definition 19.20.}  We denote the quaternion by $H$, and let we put 

\ \ \ \ \ \ \gs\gp$(4)=\left\{ X\in M(4\times 4,H)|X^{\ast }+X=0\right\} $,

\noindent
where * is a matrix transpose and quaternion conjugate.

\noindent
Then $X\in$ \gs\gp$(4)$ is expressed by

\ \ \ \ \ $X=\left(
\begin{array}
[c]{cccc}%
\chi_{1} & x_{7} & -\overline{x_{6}} & x_{8}\\
-\overline{x_{7}} & \chi_{2} & x_{5} & x_{9}\\
x_{6} & -\overline{x_{5}} & \chi_{3} & x_{10}\\
-\overline{x_{8}} & -\overline{x_{9}} & -\overline{x_{10}} & \chi_{4}%
\end{array}
\right)  $

\noindent
$X=\chi _{1}E_{1}+\chi _{2}E_{2}+\chi _{3}E_{3}+\chi_{4}E_{4}+A_{5}(x_{5})+A_{6}(x_{6})+A_{7}(x_{7})+A_{8}(x_{8})+A_{9}(x_{9})+A_{10}(x_{10}),$
where
$\chi_{1}=a_{1}e_{1}+b_{1}e_{2}+c_{1}e_{3}$,
$\chi_{2}=a_{2}e_{1}+b_{2}e_{2}+c_{2}e_{3}$,

\ \ \ \ $\chi_{3}=a_{3}e_{1}+b_{3}e_{2}+c_{3}e_{3}$,
$\chi_{4}=a_{4}e_{1}+b_{4}e_{2}+c_{4}e_{3}$,

\ \ \ \ $x_{5}=x_{10}e_{0}+x_{11}e_{1}+x_{12}e_{2}+x_{13}e_{3}$,
$x_{6}=x_{20}e_{0}+x_{21}e_{1}+x_{22}e_{2}+x_{23}e_{3}$,

\ \ \ \ $x_{7}=x_{30}e_{0}+x_{31}e_{1}+x_{32}e_{2}+x_{33}e_{3}$,
$x_{8}=z_{10}e_{0}+z_{11}e_{1}+z_{12}e_{2}+z_{13}e_{3}$,

\ \ \ \ $x_{9}=z_{20}e_{0}+z_{21}e_{1}+z_{22}e_{2}+z_{23}e_{3}$,
$x_{10}=z_{30}e_{0}+z_{31}e_{1}+z_{32}e_{2}+z_{33}e_{3}$,

\begin{flushleft}
{\fontsize{7pt}{8pt} \selectfont%
$E_{1}=\left(
\begin{array}
[c]{cccc}%
1 & 0 & 0 & 0\\
0 & 0 & 0 & 0\\
0 & 0 & 0 & 0\\
0 & 0 & 0 & 0
\end{array}
\right)  ,E_{2}=\left(
\begin{array}
[c]{cccc}%
0 & 0 & 0 & 0\\
0 & 1 & 0 & 0\\
0 & 0 & 0 & 0\\
0 & 0 & 0 & 0
\end{array}
\right)  ,$

$E_{3}=\left(
\begin{array}
[c]{cccc}%
0 & 0 & 0 & 0\\
0 & 0 & 0 & 0\\
0 & 0 & 1 & 0\\
0 & 0 & 0 & 0
\end{array}
\right)  ,E_{4}=\left(
\begin{array}
[c]{cccc}%
0 & 0 & 0 & 0\\
0 & 0 & 0 & 0\\
0 & 0 & 0 & 0\\
0 & 0 & 0 & 1
\end{array}
\right)  ,$

$A_{5}(x)=\left(
\begin{array}
[c]{cccc}%
0 & 0 & 0 & 0\\
0 & 0 & x & 0\\
0 & -\overline{x} & 0 & 0\\
0 & 0 & 0 & 0
\end{array}
\right)  ,A_{6}(x)=\left(
\begin{array}
[c]{cccc}%
0 & 0 & -\overline{x} & 0\\
0 & 0 & 0 & 0\\
x & 0 & 0 & 0\\
0 & 0 & 0 & 0
\end{array}
\right)  ,A_{7}(x)=\left(
\begin{array}
[c]{cccc}%
0 & x & 0 & 0\\
-\overline{x} & 0 & 0 & 0\\
0 & 0 & 0 & 0\\
0 & 0 & 0 & 0
\end{array}
\right)  ,$

$A_{8}(x)=\left(
\begin{array}
[c]{cccc}%
0 & 0 & 0 & x\\
0 & 0 & 0 & 0\\
0 & 0 & 0 & 0\\
-\overline{x} & 0 & 0 & 0
\end{array}
\right)  ,A_{9}(x)=\left(
\begin{array}
[c]{cccc}%
0 & 0 & 0 & 0\\
0 & 0 & 0 & x\\
0 & 0 & 0 & 0\\
0 & -\overline{x} & 0 & 0
\end{array}
\right)  ,A_{10}(x)=\left(
\begin{array}
[c]{cccc}%
0 & 0 & 0 & 0\\
0 & 0 & 0 & 0\\
0 & 0 & 0 & x\\
0 & 0 & -\overline{x} & 0
\end{array}
\right)$  ,}
\end{flushleft}

\noindent
$\{e_{0},e_{1},e_{2},e_{3}\} $ are bases of $H$, $a_{k},b_{k},c_{k},x_{ij},z_{ij} \in \R, (1 \le k \le 4, 1 \le i \le 3, 0 \le j \le 3), x \in H$.

\bigskip

\emph{Definition 19.21.}  We define the following elements and a Lie subalgebra of \gr$_{6}^{\C}$.

$M_{1}=x_{10}Um_{10}+x_{11}Um_{11}+x_{12}Um_{12}+x_{13}Um_{13}$,

$M_{2}=x_{20}Um_{20}+x_{21}Um_{21}+x_{22}Um_{22}+x_{23}Um_{23}$,

$M_{3}=x_{30}Um_{30}+x_{31}Um_{31}+x_{32}Um_{32}+x_{33}Um_{33}$,

$T_{11}=z_{10}iUt_{17}+z_{11}iUt_{16}+z_{12}iUt_{15}+z_{13}iUt_{14}$,

$T_{12}=z_{20}iUt_{27}+z_{21}iUt_{26}+z_{22}iUt_{25}+z_{23}iUt_{24}$,

$T_{13}=z_{30}iUt_{37}+z_{31}iUt_{36}+z_{32}iUt_{35}+z_{33}iUt_{34}$,

$rsp_{1}=a_{1}S_{11}+b_{1}S_{12}+c_{1}S_{13}$,
$rsp_{2}=a_{2}K_{1}+b_{2}K_{2}+c_{2}K_{3}$,

$rsp_{3}=a_{3}K_{4}+b_{3}K_{5}+c_{3}K_{6}$,
$rsp_{4}=a_{4}S_{21}+b_{4}S_{22}+c_{4}S_{23}$,

\noindent
$rsp_{5}=2M_{1}$,
$rsp_{6}=2M_{2}$,
$rsp_{7}=2M_{3}$,
$rsp_{8}=2T_{11}$,
$rsp_{9}=2T_{12}$,
$rsp_{10}=2T_{13}$,

\noindent
where $a_{k},b_{k},c_{k},x_{ij},z_{ij} \in \R, (1 \le k \le 4, 1 \le i \le 3, 0 \le j \le 3)$.

\gr\gs\gp$(4)=\{ $Lie algebra generated by $K_{n}(1\leq n\leq 6),S_{11},S_{12},S_{13},S_{21},S_{22}$,

\ \ \ \ \ \ \ \ \ \ \ \ $S_{23},Um_{kj}(1\leq k\leq 3,0\leq j \leq 3)$,$iUt_{lj}(4 \le l \le 7,0\le j \le 3)$ 

\ \ \ \ \ \ \ \ \ \ \ \ over real numbers $\} $.

\bigskip

\emph{Lemma 19.22.} \ \gr\gs\gp$(4)$ is isomorphic to \gs \gp$(4)$. 

\bigskip

\emph{Proof.}  For $X, Y \in \mathfrak{sp}(4)$, we denote $X, Y$ as follows:

 $X=\chi _{1}E_{1}+\chi _{2}E_{2}+\chi _{3}E_{3}+\chi_{4}E_{4}+A_{5}(x_{5})+A_{6}(x_{6})+A_{7}(x_{7})+A_{8}(x_{8})+A_{9}(x_{9})+A_{10}(x_{10}),$
where
$\chi_{1}=a_{1}e_{1}+b_{1}e_{2}+c_{1}e_{3}$,
$\chi_{2}=a_{2}e_{1}+b_{2}e_{2}+c_{2}e_{3}$,

\ \ \ \ $\chi_{3}=a_{3}e_{1}+b_{3}e_{2}+c_{3}e_{3}$,
$\chi_{4}=a_{4}e_{1}+b_{4}e_{2}+c_{4}e_{3}$,

\ \ \ \ $x_{5}=x_{10}e_{0}+x_{11}e_{1}+x_{12}e_{2}+x_{13}e_{3}$,
$x_{6}=x_{20}e_{0}+x_{21}e_{1}+x_{22}e_{2}+x_{23}e_{3}$,

\ \ \ \ $x_{7}=x_{30}e_{0}+x_{31}e_{1}+x_{32}e_{2}+x_{33}e_{3}$,
$x_{8}=z_{10}e_{0}+z_{11}e_{1}+z_{12}e_{2}+z_{13}e_{3}$,

\ \ \ \ $x_{9}=z_{20}e_{0}+z_{21}e_{1}+z_{22}e_{2}+z_{23}e_{3}$,
$x_{10}=z_{30}e_{0}+z_{31}e_{1}+z_{32}e_{2}+z_{33}e_{3}$,

 $Y=\gamma _{1}E_{1}+\gamma _{2}E_{2}+\gamma _{3}E_{3}+\gamma_{4}E_{4}+A_{5}(y_{5})+A_{6}(y_{6})+A_{7}(y_{7})+A_{8}(y_{8})+A_{9}(y_{9})+A_{10}(y_{10}),$
where
$\gamma_{1}=a_{5}e_{1}+b_{5}e_{2}+c_{5}e_{3}$,
$\gamma_{2}=a_{6}e_{1}+b_{6}e_{2}+c_{6}e_{3}$,

\ \ \ \ $\gamma_{3}=a_{7}e_{1}+b_{7}e_{2}+c_{7}e_{3}$,
$\gamma_{4}=a_{8}e_{1}+b_{8}e_{2}+c_{8}e_{3}$,

\ \ \ \ $y_{5}=y_{10}e_{0}+y_{11}e_{1}+y_{12}e_{2}+y_{13}e_{3}$,
$y_{6}=y_{20}e_{0}+y_{21}e_{1}+y_{22}e_{2}+y_{23}e_{3}$,

\ \ \ \ $y_{7}=y_{30}e_{0}+y_{31}e_{1}+y_{32}e_{2}+y_{33}e_{3}$,
$y_{8}=w_{10}e_{0}+w_{11}e_{1}+w_{12}e_{2}+w_{13}e_{3}$,

\ \ \ \ $y_{9}=w_{20}e_{0}+w_{21}e_{1}+w_{22}e_{2}+w_{23}e_{3}$,
$y_{10}=w_{30}e_{0}+w_{31}e_{1}+w_{32}e_{2}+w_{33}e_{3}$.

\noindent
And for $rsp_{1}^{`},rsp_{2}^{`},rsp_{3}^{`},rsp_{4}^{`},rsp_{5}^{`},rsp_{6}^{`},rsp_{7}^{`},rsp_{8}^{`},rsp_{9}^{`},rsp_{10}^{`} \in \mathfrak{rsp}(4)$,
we denote them as follows:

$rsp_{1}^{`}=a_{5}S_{11}+b_{5}S_{12}+c_{5}S_{13}$,
$rsp_{2}^{`}=a_{6}K_{1}+b_{6}K_{2}+c_{6}K_{3}$,

$rsp_{3}^{`}=a_{7}K_{4}+b_{7}K_{5}+c_{7}K_{6}$,
$rsp_{4}^{`}=a_{8}S_{21}+b_{8}S_{22}+c_{8}S_{23}$,

\noindent
$rsp_{5}^{`}=2M_{1}^{`}$,
$rsp_{6}^{`}=2M_{2}^{`}$,
$rsp_{7}^{`}=2M_{3}^{`}$,
$rsp_{8}^{`}=2T_{11}^{`}$,
$rsp_{9}^{`}=2T_{12}^{`}$,
$rsp_{10}^{`}=2T_{13}^{`}$,

$M_{1}^{`}=y_{10}Um_{10}+y_{11}Um_{11}+y_{12}Um_{12}+y_{13}Um_{13}$,

$M_{2}^{`}=y_{20}Um_{20}+y_{21}Um_{21}+y_{22}Um_{22}+y_{23}Um_{23}$,

$M_{3}^{`}=y_{30}Um_{30}+y_{31}Um_{31}+y_{32}Um_{32}+y_{33}Um_{33}$,

$T_{11}^{`}=w_{10}iUt_{17}+w_{11}iUt_{16}+w_{12}iUt_{15}+w_{13}iUt_{14}$,

$T_{12}^{`}=w_{20}iUt_{27}+w_{21}iUt_{26}+w_{22}iUt_{25}+w_{23}iUt_{24}$,

$T_{13}^{`}=w_{30}iUt_{37}+w_{31}iUt_{36}+w_{32}iUt_{35}+w_{33}iUt_{34}$.

We consider the following mapping f: $\mathfrak{sp}(4) \rightarrow \mathfrak{rsp}(4)$,

f($\chi_{i}E_{i})=rsp_{i}, (i=1,2,3,4)$, and f($A_{k}(x_{k}))=rsp_{k}, (k=5,6,7,8,9,10)$.

Let's calculate bracket operations in $\mathfrak{sp}(4)$ and $\mathfrak{rsp}(4)$.

\noindent
Case1:

$[\chi_{1}E_{1}, \gamma_{1}E_{1}]=(2(b_{1}c_{5}-b_{5}c_{1})e_{1}+2(a_{5}c_{1}-a_{1}c_{5})e_{2}+2(a_{1}b_{5}-a_{5}b_{1})e_{3})E_{1}$,

$[rsp_{1},rsp_{1}^{`}]=Ud_{67}(2b_{1}c_{5}-2b_{5}c_{1})+Ud_{45}(2b_{1}c_{5}-2b_{5}c_{1})+Ud_{46}(2a_{1}c_{5}-2a_{5}c_{1})+Ud_{57}(2a_{5}c_{1}-2a_{1}c_{5})+Ud_{56}(2a_{1}b_{5}-2a_{5}b_{1})+Ud_{47}(2a_{1}b_{5}-2a_{5}b_{1})$,

\noindent
therefore f($[\chi_{1}E_{1}, \gamma_{1}E_{1}])=[rsp_{1},rsp_{1}^{`}]$.

$[\chi_{1}E_{1}, \gamma_{2}E_{2}]=0$, $[rsp_{1},rsp_{2}^{`}]=0$.

$[\chi_{1}E_{1}, \gamma_{3}E_{3}]=0$, $[rsp_{1},rsp_{3}^{`}]=0$.

$[\chi_{1}E_{1}, \gamma_{4}E_{4}]=0$, $[rsp_{1},rsp_{4}^{`}]=0$.

\noindent
Case2:

$[\chi_{1}E_{1}, A_{5}(y_{5})]=0$, $[rsp_{1},rsp_{5}^{`}]=0$.

$[\chi_{1}E_{1}, A_{6}(y_{6})]=A_{6}((c_{1}y_{23}+b_{1}y_{22}+a_{1}y_{21})e_{0}+(b_{1}y_{23}-c_{1}y_{22}-a_{1}y_{20})e_{1}+(-a_{1}y_{23}+c_{1}y_{21}-b_{1}y_{20})e_{2}+(a_{1}y_{22}-b_{1}y_{21}-c_{1}y_{20})e_{3})$,

$[rsp_{1},rsp_{6}^{`}]=2Um_{20}(c_{1}y_{23}+b_{1}y_{22}+a_{1}y_{21})-2Um_{21}(-b_{1}y_{23}+c_{1}y_{22}+a_{1}y_{20})-2Um_{22}(a_{1}y_{23}-c_{1}y_{21}+b_{1}y_{20})-2Um_{23}(-a_{1}y_{22}+b_{1}y_{21}+c_{1}y_{20})$,

\noindent
therefore f($[\chi_{1}E_{1}, A_{6}(y_{6})])=[rsp_{1},rsp_{6}^{`}]$.

$[\chi_{1}E_{1}, A_{7}(y_{7})]=A_{7}((-c_{1}y_{33}-b_{1}y_{32}-a_{1}y_{31})e_{0}+(b_{1}y_{33}-c_{1}y_{32}+a_{1}y_{30})e_{1}+(-a_{1}y_{33}+c_{1}y_{31}+b_{1}y_{30})e_{2}+(a_{1}y_{32}-b_{1}y_{31}+c_{1}y_{30})e_{3})$,

$[rsp_{1},rsp_{7}^{`}]=2Um_{30}(-c_{1}y_{33}-b_{1}y_{32}-a_{1}y_{31})-2Um_{31}(-b_{1}y_{33}+c_{1}y_{32}-a_{1}y_{30})-2Um_{32}(a_{1}y_{33}-c_{1}y_{31}-b_{1}y_{30})-2Um_{33}(-a_{1}y_{32}+b_{1}y_{31}-c_{1}y_{30})$,

\noindent
therefore f($[\chi_{1}E_{1}, A_{7}(y_{7})])=[rsp_{1},rsp_{7}^{`}]$.

\noindent
Case3:

$[\chi_{1}E_{1}, A_{8}(y_{8})]=A_{8}((-c_{1}w_{13}-b_{1}w_{12}-a_{1}w_{11})e_{0}+(b_{1}w_{13}-c_{1}w_{12}+a_{1}w_{10})e_{1}+(-a_{1}w_{13}+c_{1}w_{11}+b_{1}w_{10})e_{2}+(a_{1}w_{12}-b_{1}w_{11}+c_{1}w_{10})e_{3})$,

$[rsp_{1},rsp_{8}^{`}]=-2iUt_{17}(c_{1}w_{13}+b_{1}w_{12}+a_{1}w_{11})-2iUt_{16}(-b_{1}w_{13}+c_{1}w_{12}-a_{1}w_{10})-2iUt_{15}(a_{1}w_{13}-c_{1}w_{11}-b_{1}w_{10})-2iUt_{14}(-a_{1}w_{12}+b_{1}w_{11}-c_{1}w_{10})$,

\noindent
therefore f($[\chi_{1}E_{1}, A_{8}(y_{8}])=[rsp_{1},rsp_{8}^{`}]$.

$[\chi_{1}E_{1}, A_{9}(y_{9})]=0$, $[rsp_{1},rsp_{9}^{`}]=0$.

$[\chi_{1}E_{1}, A_{10}(y_{10})]=0$, $[rsp_{1},rsp_{10}^{`}]=0$.

\noindent
Case4:

$[\chi_{2}E_{2}, \gamma_{2}E_{2}]=(2(b_{2}c_{6}-b_{6}c_{2})e_{1}+2(a_{6}c_{2}-a_{2}c_{6})e_{2}+2(a_{2}b_{6}-2a_{6}b_{2})e_{3})E_{2}$,

$[rsp_{2},rsp_{2}^{`}]=Ud_{23}(2b_{6}c_{2}-2b_{2}c_{6})+Ud_{01}(2b_{6}c_{2}-2b_{2}c_{6})+Ud_{02}(2a_{2}c_{6}-2a_{6}c_{2})+Ud_{13}(2a_{6}c_{2}-2a_{2}c_{6})+Ud_{12}(2a_{6}b_{2}-2a_{2}b_{6})+Ud_{03}(2a_{6}b_{2}-2a_{2}b_{6})$,

\noindent
therefore f($[\chi_{2}E_{2}, \gamma_{2}E_{2}])=[rsp_{2},rsp_{2}^{`}]$.

$[\chi_{2}E_{2}, \gamma_{3}E_{3}]=0$, $[rsp_{2},rsp_{3}^{`}]=0$.

$[\chi_{2}E_{2}, \gamma_{4}E_{4}]=0$, $[rsp_{2},rsp_{4}^{`}]=0$.

\noindent
Case5:

$[\chi_{2}E_{2}, A_{5}(y_{5})]=A_{5}((-c_{2}y_{13}-b_{2}y_{12}-a_{2}y_{11})e_{0}+(b_{2}y_{13}-c_{2}y_{12}+a_{2}y_{10})e_{1}+(-a_{2}y_{13}+c_{2}y_{11}+b_{2}y_{10})e_{2}+(a_{2}y_{12}-b_{2}y_{11}+c_{2}y_{10})e_{3})$,

$[rsp_{2},rsp_{5}^{`}]=2Um_{10}(-c_{2}y_{13}-b_{2}y_{12}-a_{2}y_{11})-2Um_{11}(-b_{2}y_{13}+c_{2}y_{12}-a_{2}y_{10})-2Um_{12}(a_{2}y_{13}-c_{2}y_{11}-b_{2}y_{10})-2Um_{13}(-a_{2}y_{12}+b_{2}y_{11}-c_{2}y_{10})$,

\noindent
therefore f($[\chi_{2}E_{2}, A_{5}(y_{5})])=[rsp_{2},rsp_{5}^{`}]$.

$[\chi_{2}E_{2}, A_{6}(y_{6})]=0$,$[rsp_{2},rsp_{6}^{`}]=0$.

$[\chi_{2}E_{2}, A_{7}(y_{7})]=A_{7}((c_{2}y_{33}+b_{2}y_{32}+a_{2}y_{31})e_{0}+(b_{2}y_{33}-c_{2}y_{32}-a_{2}y_{30})e_{1}+(-a_{2}y_{33}+c_{2}y_{31}-b_{2}y_{30})e_{2}+(a_{2}y_{32}-b_{2}y_{31}-c_{2}y_{30})e_{3})$,

$[rsp_{2},rsp_{7}^{`}]=2Um_{30}(c_{2}y_{33}+b_{2}y_{32}+a_{2}y_{31})-2Um_{31}(-b_{2}y_{33}+c_{2}y_{32}+a_{2}y_{30})-2Um_{32}(a_{2}y_{33}-c_{2}y_{31}+b_{2}y_{30})-2Um_{33}(-a_{2}y_{32}+b_{2}y_{31}+c_{2}y_{30})$,

\noindent
therefore f($[\chi_{2}E_{2}, A_{7}(y_{7})])=[rsp_{2},rsp_{7}^{`}]$.

\noindent
Case6:

$[\chi_{2}E_{2}, A_{8}(y_{8})]=0$,$[rsp_{2},rsp_{8}^{`}]=0$.

$[\chi_{2}E_{2}, A_{9}(y_{9})]=A_{9}((-c_{2}w_{23}-b_{2}w_{22}-a_{2}w_{21})e_{0}+(b_{2}w_{23}-c_{2}w_{22}+a_{2}w_{20})e_{1}+(-a_{2}w_{23}+c_{2}w_{21}+b_{2}w_{20})e_{2}+(a_{2}w_{22}-b_{2}w_{21}+c_{2}w_{20})e_{3})$

$[rsp_{2},rsp_{9}^{`}]=-2iUt_{27}(c_{2}w_{23}+b_{2}w_{22}+a_{2}w_{21})-2iUt_{26}(-b_{2}w_{23}+c_{2}w_{22}-a_{2}w_{20})-2iUt_{25}(a_{2}w_{23}-c_{2}w_{21}-b_{2}w_{20})-2iUt_{24}(-a_{2}w_{22}+b_{2}w_{21}-c_{2}w_{20})$,

\noindent
therefore f($[\chi_{2}E_{2}, A_{9}(y_{9})])=[rsp_{2},rsp_{9}^{`}]$.

$[\chi_{2}E_{2}, A_{10}(y_{10})]=0$, $[rsp_{2},rsp_{10}^{`}]=0$.

\noindent
Case7:

$[\chi_{3}E_{3}, \gamma_{3}E_{3}]=((2b_{3}c_{7}-2b_{7}c_{3})e_{1}+(2a_{7}c_{3}-2a_{3}c_{7})e_{2}+(2a_{3}b_{7}-2a_{7}b_{3})e_{3})E_{3}$,

$[rsp_{3},rsp_{3}^{`}]=Ud_{01}(2b_{3}c_{7}-2b_{7}c_{3})+Ud_{23}(2b_{7}c_{3}-2b_{3}c_{7})+Ud_{13}(2a_{7}c_{3}-2a_{3}c_{7})+Ud_{02}(2a_{7}c_{3}-2a_{3}c_{7})+Ud_{03}(2a_{3}b_{7}-2a_{7}b_{3})+Ud_{12}(2a_{7}b_{3}-2a_{3}b_{7})$,

\noindent
therefore f($[\chi_{3}E_{3}, \gamma_{3}E_{3}])=[rsp_{3},rsp_{3}^{`}]$.

$[\chi_{3}E_{3}, \gamma_{4}E_{4}]=0$,$[rsp_{3},rsp_{4}^{`}]=0$.

\noindent
Case8:

$[\chi_{3}E_{3}, A_{5}(y_{5})]=A_{5}((c_{3}y_{13}+b_{3}y_{12}+a_{3}y_{11})e_{0}+(b_{3}y_{13}-c_{3}y_{12}-a_{3}y_{10})e_{1}+(-a_{3}y_{13}+c_{3}y_{11}-b_{3}y_{10})e_{2}+(a_{3}y_{12}-b_{3}y_{11}-c_{3}y_{10})e_{3})$,

$[rsp_{3},rsp_{5}^{`}]=2Um_{10}(c_{3}y_{13}+b_{3}y_{12}+a_{3}y_{11})-2Um_{11}(-b_{3}y_{13}+c_{3}y_{12}+a_{3}y_{10})-2Um_{12}(a_{3}y_{13}-c_{3}y_{11}+b_{3}y_{10})-2Um_{13}(-a_{3}y_{12}+b_{3}y_{11}+c_{3}y_{10})$,

\noindent
therefore f($[\chi_{3}E_{3}, A_{5}(y_{5})])=[rsp_{3},rsp_{5}^{`}]$.

$[\chi_{3}E_{3}, A_{6}(y_{6})]=A_{6}((-c_{3}y_{23}-b_{3}y_{22}-a_{3}y_{21})e_{0}+(b_{3}y_{23}-c_{3}y_{22}+a_{3}y_{20})e_{1}+(-a_{3}y_{23}+c_{3}y_{21}+b_{3}y_{20})e_{2}+(a_{3}y_{22}-b_{3}y_{21}+c_{3}y_{20})e_{3})$,

$[rsp_{3},rsp_{6}^{`}]=2Um_{20}(-c_{3}y_{23}-b_{3}y_{22}-a_{3}y_{21})-2Um_{21}(-b_{3}y_{23}+c_{3}y_{22}-a_{3}y_{20})-2Um_{22}(a_{3}y_{23}-c_{3}y_{21}-b_{3}y_{20})-2Um_{23}(-a_{3}y_{22}+b_{3}y_{21}-c_{3}y_{20})$,

\noindent
therefore f($[\chi_{3}E_{3}, A_{6}(y_{6})])=[rsp_{3},rsp_{6}^{`}]$.

$[\chi_{3}E_{3}, A_{7}(y_{7})]=0$, $[rsp_{3},rsp_{7}^{`}]=0$.

\noindent
Case9:

$[\chi_{3}E_{3}, A_{8}(y_{8})]=0$, $[rsp_{3},rsp_{8}^{`}]=0$.

$[\chi_{3}E_{3}, A_{9}(y_{9})]=0$, $[rsp_{3},rsp_{9}^{`}]=0$.

$[\chi_{3}E_{3}, A_{10}(y_{10})]=A_{10}((-c_{3}w_{33}-b_{3}w_{32}-a_{3}w_{31})e_{0}+(b_{3}w_{33}-c_{3}w_{32}+a_{3}w_{30})e_{1}+(-a_{3}w_{33}+c_{3}w_{31}+b_{3}w_{30})e_{2}+(a_{3}w_{32}-b_{3}w_{31}+c_{3}w_{30})e_{3})$,

$[rsp_{3},rsp_{10}^{`}]=-2iUt_{37}(c_{3}w_{33}+b_{3}w_{32}+a_{3}w_{31})-2iUt_{36}(-b_{3}w_{33}+c_{3}w_{32}-a_{3}w_{30})-2iUt_{35}(a_{3}w_{33}-c_{3}w_{31}-b_{3}w_{30})-2iUt_{34}(-a_{3}w_{32}+b_{3}w_{31}-c_{3}w_{30})$,

\noindent
therefore f($[\chi_{3}E_{3}, A_{10}(y_{10})])=[rsp_{3},rsp_{10}^{`}]$.

\noindent
Case10:

$[\chi_{4}E_{4}, \gamma_{4}E_{4}]=((2b_{4}c_{8}-2b_{8}c_{4})e_{1}+(2a_{8}c_{4}-2a_{4}c_{8})e_{2}+(2a_{4}b_{8}-2a_{8}b_{4})e_{3})E_{4}$,

$[rsp_{4},rsp_{4}^{`}]=Ud_{45}(2b_{4}c_{8}-2b_{8}c_{4})+Ud_{67}(2b_{8}c_{4}-2b_{4}c_{8})+Ud_{57}(2a_{4}c_{8}-2a_{8}c_{4})+Ud_{46}(2a_{4}c_{8}-2a_{8}c_{4})+Ud_{56}(2a_{4}b_{8}-2a_{8}b_{4})+Ud_{47}(2a_{8}b_{4}-2a_{4}b_{8})$,

\noindent
therefore f($[\chi_{4}E_{4}, \gamma_{4}E_{4}])=[rsp_{4},rsp_{4}^{`}]$.

\noindent
Case11:

$[\chi_{4}E_{4}, A_{5}(y_{5})]=0$, $[rsp_{4},rsp_{5}^{`}]=0$.

$[\chi_{4}E_{4}, A_{6}(y_{6})]=0$, $[rsp_{4},rsp_{6}^{`}]=0$.

$[\chi_{4}E_{4}, A_{7}(y_{7})]=0$, $[rsp_{4},rsp_{7}^{`}]=0$.

\noindent
Case12:

$[\chi_{4}E_{4}, A_{8}(y_{8})]=A_{8}((c_{4}w_{13}+b_{4}w_{12}+a_{4}w_{11})e_{0}+(b_{4}w_{13}-c_{4}w_{12}-a_{4}w_{10})e_{1}+(-a_{4}w_{13}+c_{4}w_{11}-b_{4}w_{10})e_{2}+(a_{4}w_{12}-b_{4}w_{11}-c_{4}w_{10})e_{3})$,

$[rsp_{4},rsp_{8}^{`}]=-2iUt_{17}(-c_{4}w_{13}-b_{4}w_{12}-a_{4}w_{11})-2iUt_{16}(-b_{4}w_{13}+c_{4}w_{12}+a_{4}w_{10})-2iUt_{15}(a_{4}w_{13}-c_{4}w_{11}+b_{4}w_{10})-2iUt_{14}(-a_{4}w_{12}+b_{4}w_{11}+c_{4}w_{10})$,

\noindent
therefore f($[\chi_{4}E_{4}, A_{8}(y_{8})])=[rsp_{4},rsp_{8}^{`}]$.

$[\chi_{4}E_{4}, A_{9}(y_{9})]=A_{9}((c_{4}w_{23}+b_{4}w_{22}+a_{4}w_{21})e_{0}+(b_{4}w_{23}-c_{4}w_{22}-a_{4}w_{20})e_{1}+(-a_{4}w_{23}+c_{4}w_{21}-b_{4}w_{20})e_{2}+(a_{4}w_{22}-b_{4}w_{21}-c_{4}w_{20})e_{3})$,

$[rsp_{4},rsp_{9}^{`}]=-2iUt_{27}(-c_{4}w_{23}-b_{4}w_{22}-a_{4}w_{21})-2iUt_{26}(-b_{4}w_{23}+c_{4}w_{22}+a_{4}w_{20})-2iUt_{25}(a_{4}w_{23}-c_{4}w_{21}+b_{4}w_{20})-2iUt_{24}(-a_{4}w_{22}+b_{4}w_{21}+c_{4}w_{20})$,

\noindent
therefore f($[\chi_{4}E_{4}, A_{9}(y_{9})])=[rsp_{4},rsp_{9}^{`}]$.

$[\chi_{4}E_{4}, A_{10}(y_{10})]=A_{10}((c_{4}w_{33}+b_{4}w_{32}+a_{4}w_{31})e_{0}+(b_{4}w_{33}-c_{4}w_{32}-a_{4}w_{30})e_{1}$

\noindent
$+(-a_{4}w_{33}+c_{4}w_{31}-b_{4}w_{30})e_{2}+(a_{4}w_{32}-b_{4}w_{31}-c_{4}w_{30})e_{3})$,

$[rsp_{4},rsp_{10}^{`}]=-2iUt_{37}(-c_{4}w_{33}-b_{4}w_{32}-a_{4}w_{31})-2iUt_{36}(-b_{4}w_{33}+c_{4}w_{32}+a_{4}w_{30})-2iUt_{35}(a_{4}w_{33}-c_{4}w_{31}+b_{4}w_{30})-2iUt_{34}(-a_{4}w_{32}+b_{4}w_{31}+c_{4}w_{30})$,

\noindent
therefore f($[\chi_{4}E_{4}, A_{10}(y_{10})])=[rsp_{4},rsp_{10}^{`}]$.

\noindent
Case13:

$[A_{5}(x_{5}),A_{5}(y_{5})]=((2x_{12}y_{13}-2x_{13}y_{12}+2x_{10}y_{11}-2x_{11}y_{10})e_{1}+(-2x_{11}y_{13}+2x_{10}y_{12}+2x_{13}y_{11}-2x_{12}y_{10})e_{2}+(2x_{10}y_{13}+2x_{11}y_{12}-2x_{12}y_{11}-2x_{13}y_{10})e_{3})E_{2}+((2x_{12}y_{13}-2x_{13}y_{12}-2x_{10}y_{11}+2x_{11}y_{10})e_{1}+(-2x_{11}y_{13}-2x_{10}y_{12}+2x_{13}y_{11}+2x_{12}y_{10})e_{2}+(-2x_{10}y_{13}+2x_{11}y_{12}-2x_{12}y_{11}+2x_{13}y_{10})e_{3})E_{3}$,

\noindent
then  f($[A_{5}(x_{5}),A_{5}(y_{5})])=$
$(2x_{12}y_{13}-2x_{13}y_{12}+2x_{10}y_{11}-2x_{11}y_{10})(-Ud_{01}-Ud_{23})
+(-2x_{11}y_{13}+2x_{10}y_{12}+2x_{13}y_{11}-2x_{12}y_{10})(-Ud_{02}+Ud_{13})
+(2x_{10}y_{13}+2x_{11}y_{12}-2x_{12}y_{11}-2x_{13}y_{10})(-Ud_{03}-Ud_{12})
+((2x_{12}y_{13}-2x_{13}y_{12}-2x_{10}y_{11}+2x_{11}y_{10})(Ud_{01}-Ud_{23})
+(-2x_{11}y_{13}-2x_{10}y_{12}+2x_{13}y_{11}+2x_{12}y_{10})(Ud_{02}+Ud_{13})
+(-2x_{10}y_{13}+2x_{11}y_{12}-2x_{12}y_{11}+2x_{13}y_{10})(Ud_{03}-Ud_{12})$,

$[rsp_{5},rsp_{5}^{`}]=Ud_{23}(4x_{13}y_{12}-4x_{12}y_{13})+Ud_{13}(4x_{13}y_{11}-4x_{11}y_{13})$

\noindent
$+Ud_{03}(4x_{13}y_{10}-4x_{10}y_{13})+Ud_{12}(4x_{12}y_{11}-4x_{11}y_{12})+Ud_{02}(4x_{12}y_{10}-4x_{10}y_{12})+Ud_{01}(4x_{11}y_{10}-4x_{10}y_{11})$,

\noindent
therefore f($[A_{5}(x_{5}),A_{5}(y_{5})])=[rsp_{5},rsp_{5}^{`}]$.

$[A_{5}(x_{5}),A_{6}(y_{6})]=A_{7}((x_{13}y_{23}+x_{12}y_{22}+x_{11}y_{21}-x_{10}y_{20})e_{0}+(x_{12}y_{23}-x_{13}y_{22}+x_{10}y_{21}+x_{11}y_{20})e_{1}+(-x_{11}y_{23}+x_{10}y_{22}+x_{13}y_{21}+x_{12}y_{20})e_{2}+(x_{10}y_{23}+x_{11}y_{22}-x_{12}y_{21}+x_{13}y_{20})e_{3}))$,

$[rsp_{5},rsp_{6}^{`}]=2Um_{30}(x_{13}y_{23}+x_{12}y_{22}+x_{11}y_{21}-x_{10}y_{20})-2Um_{31}(-x_{12}y_{23}+x_{13}y_{22}-x_{10}y_{21}-x_{11}y_{20})-2Um_{32}(x_{11}y_{23}-x_{10}y_{22}-x_{13}y_{21}-x_{12}y_{20})$

\noindent
$-2Um_{33}(-x_{10}y_{23}-x_{11}y_{22}+x_{12}y_{21}-x_{13}y_{20})$,

\noindent
therefore f($[A_{5}(x_{5}),A_{6}(y_{6})])=[rsp_{5},rsp_{6}^{`}]$.

$[A_{5}(x_{5}),A_{7}(y_{7})]=A_{6}((-x_{13}y_{33}-x_{12}y_{32}-x_{11}y_{31}+x_{10}y_{30})e_{0}+(x_{12}y_{33}-x_{13}y_{32}-x_{10}y_{31}-x_{11}y_{30})e_{1}+(-x_{11}y_{33}-x_{10}y_{32}+x_{13}y_{31}-x_{12}y_{30})e_{2}+(-x_{10}y_{33}+x_{11}y_{32}-x_{12}y_{31}-x_{13}y_{30})e_{3})$,

$[rsp_{5},rsp_{7}^{`}]=2Um_{20}(-x_{13}y_{33}-x_{12}y_{32}-x_{11}y_{31}+x_{10}y_{30})-2Um_{21}(-x_{12}y_{33}+x_{13}y_{32}+x_{10}y_{31}+x_{11}y_{30})-2Um_{22}(x_{11}y_{33}+x_{10}y_{32}-x_{13}y_{31}+x_{12}y_{30})$

\noindent
$-2Um_{23}(x_{10}y_{33}-x_{11}y_{32}+x_{12}y_{31}+x_{13}y_{30})$,

\noindent
therefore f($[A_{5}(x_{5}),A_{7}(y_{7})])=[rsp_{5},rsp_{7}^{`}]$.

\noindent
Case14:

$[A_{5}(x_{5}),A_{8}(y_{8})]=0$, $[rsp_{5},rsp_{8}^{`}]=0$.

$[A_{5}(x_{5}),A_{9}(y_{9})]=A_{10}((-w_{23}x_{13}-w_{22}x_{12}-w_{21}x_{11}-w_{20}x_{10})e_{0}+(-w_{22}x_{13}+w_{23}x_{12}+w_{20}x_{11}-w_{21}x_{10})e_{1}+(w_{21}x_{13}+w_{20}x_{12}-w_{23}x_{11}-w_{22}x_{10})e_{2}+(w_{20}x_{13}-w_{21}x_{12}+w_{22}x_{11}-w_{23}x_{10})e_{3}))$,

$[rsp_{5},rsp_{9}^{`}]=-2iUt_{37}(w_{23}x_{13}+w_{22}x_{12}+w_{21}x_{11}+w_{20}x_{10})-2iUt_{36}(w_{22}x_{13}-w_{23}x_{12}-w_{20}x_{11}+w_{21}x_{10})-2iUt_{35}(-w_{21}x_{13}-w_{20}x_{12}+w_{23}x_{11}+w_{22}x_{10})-2iUt_{34}(-w_{20}x_{13}+w_{21}x_{12}-w_{22}x_{11}+w_{23}x_{10})$,

\noindent
therefore f($[A_{5}(x_{5}),A_{9}(y_{9})])=[rsp_{5},rsp_{9}^{`}]$.

$[A_{5}(x_{5}),A_{10}(y_{10})]=A_{9}((-w_{33}x_{13}-w_{32}x_{12}-w_{31}x_{11}+w_{30}x_{10})e_{0}$

\noindent
$+(-w_{32}x_{13}+w_{33}x_{12}+w_{30}x_{11}+w_{31}x_{10})e_{1}+(w_{31}x_{13}+w_{30}x_{12}-w_{33}x_{11}+w_{32}x_{10})e_{2}+(w_{30}x_{13}-w_{31}x_{12}+w_{32}x_{11}+w_{33}x_{10})e_{3})$,

$[rsp_{5},rsp_{10}^{`}]=-2iUt_{27}(w_{33}x_{13}+w_{32}x_{12}+w_{31}x_{11}-w_{30}x_{10})-2iUt_{26}(w_{32}x_{13}-w_{33}x_{12}-w_{30}x_{11}-w_{31}x_{10})-2iUt_{25}(-w_{31}x_{13}-w_{30}x_{12}+w_{33}x_{11}-w_{32}x_{10})-2iUt_{24}(-w_{30}x_{13}+w_{31}x_{12}-w_{32}x_{11}-w_{33}x_{10})$,

\noindent
therefore f($[A_{5}(x_{5}),A_{10}(y_{10})])=[rsp_{5},rsp_{10}^{`}]$.

\noindent
Case15:

$[A_{6}(x_{6}),A_{6}(y_{6})]=((2x_{22}y_{23}-2x_{23}y_{22}-2x_{20}y_{21}+2x_{21}y_{20})e_{1}+(-2x_{21}y_{23}-2x_{20}y_{22}+2x_{23}y_{21}+2x_{22}y_{20})e_{2}+(-2x_{20}y_{23}+2x_{21}y_{22}-2x_{22}y_{21}+2x_{23}y_{20})e_{3})E_{1}+((2x_{22}y_{23}-2x_{23}y_{22}+2x_{20}y_{21}-2x_{21}y_{20})e_{1}+(-2x_{21}y_{23}+2x_{20}y_{22}+2x_{23}y_{21}-2x_{22}y_{20})e_{2}+(2x_{20}y_{23}+2x_{21}y_{22}-2x_{22}y_{21}-2x_{23}y_{20})e_{3})E_{3}$,

$[rsp_{6},rsp_{6}^{`}]=Ud_{01}(2x_{22}y_{23}-2x_{23}y_{22}+2x_{20}y_{21}-2x_{21}y_{20})+Ud_{67}(2x_{22}y_{23}-2x_{23}y_{22}-2x_{20}y_{21}+2x_{21}y_{20})+Ud_{45}(2x_{22}y_{23}-2x_{23}y_{22}-2x_{20}y_{21}+2x_{21}y_{20})+Ud_{23}(-2x_{22}y_{23}+2x_{23}y_{22}-2x_{20}y_{21}+2x_{21}y_{20})+Ud_{46}(2x_{21}y_{23}+2x_{20}y_{22}-2x_{23}y_{21}-2x_{22}y_{20})+Ud_{13}(-2x_{21}y_{23}+2x_{20}y_{22}+2x_{23}y_{21}-2x_{22}y_{20})+Ud_{02}(-2x_{21}y_{23}+2x_{20}y_{22}+2x_{23}y_{21}-2x_{22}y_{20})+Ud_{57}(-2x_{21}y_{23}-2x_{20}y_{22}+2x_{23}y_{21}+2x_{22}y_{20})+Ud_{03}(2x_{20}y_{23}+2x_{21}y_{22}-2x_{22}y_{21}-2x_{23}y_{20})+Ud_{56}(-2x_{20}y_{23}+2x_{21}y_{22}-2x_{22}y_{21}+2x_{23}y_{20})+Ud_{47}(-2x_{20}y_{23}+2x_{21}y_{22}-2x_{22}y_{21}+2x_{23}y_{20})+Ud_{12}(-2x_{20}y_{23}-2x_{21}y_{22}+2x_{22}y_{21}+2x_{23}y_{20})$,

\noindent
therefore f($[A_{6}(x_{6}),A_{6}(y_{6})])=[rsp_{6},rsp_{6}^{`}]$.

$[A_{6}(x_{6}),A_{7}(y_{7})]=A_{5}((x_{23}y_{33}+x_{22}y_{32}+x_{21}y_{31}-x_{20}y_{30})e_{0}+(x_{22}y_{33}-x_{23}y_{32}+x_{20}y_{31}+x_{21}y_{30})e_{1}+(-x_{21}y_{33}+x_{20}y_{32}+x_{23}y_{31}+x_{22}y_{30})e_{2}+(x_{20}y_{33}+x_{21}y_{32}-x_{22}y_{31}+x_{23}y_{30})e_{3})$,

$[rsp_{6},rsp_{7}^{`}]=2Um_{10}(x_{23}y_{33}+x_{22}y_{32}+x_{21}y_{31}-x_{20}y_{30})-2Um_{11}(-x_{22}y_{33}+x_{23}y_{32}-x_{20}y_{31}-x_{21}y_{30})-2Um_{12}(x_{21}y_{33}-x_{20}y_{32}-x_{23}y_{31}-x_{22}y_{30})$

\noindent
$-2Um_{13}(-x_{20}y_{33}-x_{21}y_{32}+x_{22}y_{31}-x_{23}y_{30})$,

\noindent
therefore f($[A_{6}(x_{6}),A_{7}(y_{7})])=[rsp_{6},rsp_{7}^{`}]$.

\noindent
Case16:

$[A_{6}(x_{6}),A_{8}(y_{8})]=A_{10}((-w_{13}x_{23}-w_{12}x_{22}-w_{11}x_{21}+w_{10}x_{20})e_{0}+(-w_{12}x_{23}+w_{13}x_{22}+w_{10}x_{21}+w_{11}x_{20})e_{1}+(w_{11}x_{23}+w_{10}x_{22}-w_{13}x_{21}+w_{12}x_{20})e_{2}+(w_{10}x_{23}-w_{11}x_{22}+w_{12}x_{21}+w_{13}x_{20})e_{3})$,

$[rsp_{6},rsp_{8}^{`}]=-2iUt_{37}(w_{13}x_{23}+w_{12}x_{22}+w_{11}x_{21}-w_{10}x_{20})-2iUt_{36}(w_{12}x_{23}-w_{13}x_{22}-w_{10}x_{21}-w_{11}x_{20})-2iUt_{35}(-w_{11}x_{23}-w_{10}x_{22}+w_{13}x_{21}-w_{12}x_{20})-2iUt_{34}(-w_{10}x_{23}+w_{11}x_{22}-w_{12}x_{21}-w_{13}x_{20})$,

\noindent
therefore f($[A_{6}(x_{6}),A_{8}(y_{8})])=[rsp_{6},rsp_{8}^{`}]$.

$[A_{6}(x_{6}),A_{9}(y_{9})]=0$, $[rsp_{6},rsp_{9}^{`}]=0$.

$[A_{6}(x_{6}),A_{10}(y_{10})]=A_{8}((-w_{33}x_{23}-w_{32}x_{22}-w_{31}x_{21}-w_{30}x_{20})e_{0}$

\noindent
$+(-w_{32}x_{23}+w_{33}x_{22}+w_{30}x_{21}-w_{31}x_{20})e_{1}+(w_{31}x_{23}+w_{30}x_{22}-w_{33}x_{21}-w_{32}x_{20})e_{2}+(w_{30}x_{23}-w_{31}x_{22}+w_{32}x_{21}-w_{33}x_{20})e_{3})$,

$[rsp_{6},rsp_{10}^{`}]=-2iUt_{17}(w_{33}x_{23}+w_{32}x_{22}+w_{31}x_{21}+w_{30}x_{20})-2iUt_{16}(w_{32}x_{23}-w_{33}x_{22}-w_{30}x_{21}+w_{31}x_{20})-2iUt_{15}(-w_{31}x_{23}-w_{30}x_{22}+w_{33}x_{21}+w_{32}x_{20})-2iUt_{14}(-w_{30}x_{23}+w_{31}x_{22}-w_{32}x_{21}+w_{33}x_{20})$,

\noindent
therefore f($[A_{6}(x_{6}),A_{10}(y_{10})])=[rsp_{6},rsp_{10}^{`}]$.

\noindent
Case17:

$[A_{7}(x_{7}),A_{7}(y_{7})]=((2x_{32}y_{33}-2x_{33}y_{32}+2x_{30}y_{31}-2x_{31}y_{30})e_{1}+(-2x_{31}y_{33}+2x_{30}y_{32}+2x_{33}y_{31}-2x_{32}y_{30})e_{2}+(2x_{30}y_{33}+2x_{31}y_{32}-2x_{32}y_{31}-2x_{33}y_{30})e_{3})E_{1}+((2x_{32}y_{33}-2x_{33}y_{32}-2x_{30}y_{31}+2x_{31}y_{30})e_{1}+(-2x_{31}y_{33}-2x_{30}y_{32}+2x_{33}y_{31}+2x_{32}y_{30})e_{2}+(-2x_{30}y_{33}+2x_{31}y_{32}-2x_{32}y_{31}+2x_{33}y_{30})e_{3})E_{2}$,

$[rsp_{7},rsp_{7}^{`}]=Ud_{67}(2x_{32}y_{33}-2x_{33}y_{32}+2x_{30}y_{31}-2x_{31}y_{30})+Ud_{45}(2x_{32}y_{33}-2x_{33}y_{32}+2x_{30}y_{31}-2x_{31}y_{30})+Ud_{23}(-2x_{32}y_{33}+2x_{33}y_{32}+2x_{30}y_{31}-2x_{31}y_{30})+Ud_{01}(-2x_{32}y_{33}+2x_{33}y_{32}+2x_{30}y_{31}-2x_{31}y_{30})+Ud_{02}(2x_{31}y_{33}+2x_{30}y_{32}-2x_{33}y_{31}-2x_{32}y_{30})+Ud_{46}(2x_{31}y_{33}-2x_{30}y_{32}-2x_{33}y_{31}+2x_{32}y_{30})+Ud_{57}(-2x_{31}y_{33}+2x_{30}y_{32}+2x_{33}y_{31}-2x_{32}y_{30})+Ud_{13}(-2x_{31}y_{33}-2x_{30}y_{32}+2x_{33}y_{31}+2x_{32}y_{30})+Ud_{56}(2x_{30}y_{33}+2x_{31}y_{32}-2x_{32}y_{31}-2x_{33}y_{30})+Ud_{47}(2x_{30}y_{33}+2x_{31}y_{32}-2x_{32}y_{31}-2x_{33}y_{30})+Ud_{12}(2x_{30}y_{33}-2x_{31}y_{32}+2x_{32}y_{31}-2x_{33}y_{30})+Ud_{03}(2x_{30}y_{33}-2x_{31}y_{32}+2x_{32}y_{31}-2x_{33}y_{30})$,

\noindent
therefore f($[A_{7}(x_{7}),A_{7}(y_{7})])=[rsp_{7},rsp_{7}^{`}]$.

\noindent
Case18:

$[A_{7}(x_{7}),A_{8}(y_{8})]=A_{9}((-w_{13}x_{33}-w_{12}x_{32}-w_{11}x_{31}-w_{10}x_{30})e_{0}+(-w_{12}x_{33}+w_{13}x_{32}+w_{10}x_{31}-w_{11}x_{30})e_{1}+(w_{11}x_{33}+w_{10}x_{32}-w_{13}x_{31}-w_{12}x_{30})e_{2}+(w_{10}x_{33}-w_{11}x_{32}+w_{12}x_{31}-w_{13}x_{30})e_{3})$,

$[rsp_{7},rsp_{8}^{`}]=-2iUt_{27}(w_{13}x_{33}+w_{12}x_{32}+w_{11}x_{31}+w_{10}x_{30})-2iUt_{26}(w_{12}x_{33}-w_{13}x_{32}-w_{10}x_{31}+w_{11}x_{30})-2iUt_{25}(-w_{11}x_{33}-w_{10}x_{32}+w_{13}x_{31}+w_{12}x_{30})-2iUt_{24}(-w_{10}x_{33}+w_{11}x_{32}-w_{12}x_{31}+w_{13}x_{30})$,

\noindent
therefore f($[A_{7}(x_{7}),A_{8}(y_{8})])=[rsp_{7},rsp_{8}^{`}]$.

$[A_{7}(x_{7}),A_{9}(y_{9})]=A_{8}((-w_{23}x_{33}-w_{22}x_{32}-w_{21}x_{31}+w_{20}x_{30})e_{0}+(-w_{22}x_{33}+w_{23}x_{32}+w_{20}x_{31}+w_{21}x_{30})e_{1}+(w_{21}x_{33}+w_{20}x_{32}-w_{23}x_{31}+w_{22}x_{30})e_{2}+(w_{20}x_{33}-w_{21}x_{32}+w_{22}x_{31}+w_{23}x_{30})e_{3})$,

$[rsp_{7},rsp_{9}^{`}]=-2iUt_{17}(w_{23}x_{33}+w_{22}x_{32}+w_{21}x_{31}-w_{20}x_{30})-2iUt_{16}(w_{22}x_{33}-w_{23}x_{32}-w_{20}x_{31}-w_{21}x_{30})-2iUt_{15}(-w_{21}x_{33}-w_{20}x_{32}+w_{23}x_{31}-w_{22}x_{30})-2iUt_{14}(-w_{20}x_{33}+w_{21}x_{32}-w_{22}x_{31}-w_{23}x_{30})$,

\noindent
therefore f($[A_{7}(x_{7}),A_{9}(y_{9})])=[rsp_{7},rsp_{9}^{`}]$.

$[A_{7}(x_{7}),A_{10}(y_{10})]=0$, $[rsp_{7},rsp_{10}^{`}]=0$.

\noindent
Case19:

$[A_{8}(x_{8}),A_{8}(y_{8})]=((-2w_{12}z_{13}+2w_{13}z_{12}-2w_{10}z_{11}+2w_{11}z_{10})e_{1}+(2w_{11}z_{13}-2w_{10}z_{12}-2w_{13}z_{11}+2w_{12}z_{10})e_{2}+(-2w_{10}z_{13}-2w_{11}z_{12}+2w_{12}z_{11}+2w_{13}z_{10})e_{3})E_{1}+((-2w_{12}z_{13}+2w_{13}z_{12}+2w_{10}z_{11}-2w_{11}z_{10})e_{1}+(2w_{11}z_{13}+2w_{10}z_{12}-2w_{13}z_{11}-2w_{12}z_{10})e_{2}+(2w_{10}z_{13}-2w_{11}z_{12}+2w_{12}z_{11}-2w_{13}z_{10})e_{3})E_{4}$

$[rsp_{8},rsp_{8}^{`}]=Ud_{45}(4w_{13}z_{12}-4w_{12}z_{13})+Ud_{46}(4w_{13}z_{11}-4w_{11}z_{13})$

\noindent
$+Ud_{47}(4w_{13}z_{10}-4w_{10}z_{13})+Ud_{56}(4w_{12}z_{11}-4w_{11}z_{12})+Ud_{57}(4w_{12}z_{10}-4w_{10}z_{12})+Ud_{67}(4w_{11}z_{10}-4w_{10}z_{11})$,

\noindent
therefore f($[A_{8}(x_{8}),A_{8}(y_{8})])=[rsp_{8},rsp_{8}^{`}]$.

$[A_{8}(x_{8}),A_{9}(y_{9})]=A_{7}((-w_{23}z_{13}-w_{22}z_{12}-w_{21}z_{11}-w_{20}z_{10})e_{0}+(-w_{22}z_{13}+w_{23}z_{12}-w_{20}z_{11}+w_{21}z_{10})e_{1}+(w_{21}z_{13}-w_{20}z_{12}-w_{23}z_{11}+w_{22}z_{10})e_{2}+(-w_{20}z_{13}-w_{21}z_{12}+w_{22}z_{11}+w_{23}z_{10})e_{3})$,

$[rsp_{8},rsp_{9}^{`}]=2Um_{30}(-w_{23}z_{13}-w_{22}z_{12}-w_{21}z_{11}-w_{20}z_{10})-2Um_{31}(w_{22}z_{13}-w_{23}z_{12}+w_{20}z_{11}-w_{21}z_{10})-2Um_{32}(-w_{21}z_{13}+w_{20}z_{12}+w_{23}z_{11}-w_{22}z_{10})-2Um_{33}(w_{20}z_{13}+w_{21}z_{12}-w_{22}z_{11}-w_{23}z_{10})$,

\noindent
therefore f($[A_{8}(x_{8}),A_{9}(y_{9})])=[rsp_{8},rsp_{9}^{`}]$.

$[A_{8}(x_{8}),A_{10}(y_{10})]=A_{6}((w_{33}z_{13}+w_{32}z_{12}+w_{31}z_{11}+w_{30}z_{10})e_{0}+(-w_{32}z_{13}+w_{33}z_{12}-w_{30}z_{11}+w_{31}z_{10})e_{1}+(w_{31}z_{13}-w_{30}z_{12}-w_{33}z_{11}+w_{32}z_{10})e_{2}+(-w_{30}z_{13}-w_{31}z_{12}+w_{32}z_{11}+w_{33}z_{10})e_{3})$,

$[rsp_{8},rsp_{10}^{`}]=2Um_{20}(w_{33}z_{13}+w_{32}z_{12}+w_{31}z_{11}+w_{30}z_{10})-2Um_{21}(w_{32}z_{13}-w_{33}z_{12}+w_{30}z_{11}-w_{31}z_{10})-2Um_{22}(-w_{31}z_{13}+w_{30}z_{12}+w_{33}z_{11}-w_{32}z_{10})-2Um_{23}(w_{30}z_{13}+w_{31}z_{12}-w_{32}z_{11}-w_{33}z_{10})$,

\noindent
therefore f($[A_{8}(x_{8}),A_{10}(y_{10})])=[rsp_{8},rsp_{10}^{`}]$.

\noindent
Case20:

$[A_{9}(x_{9}),A_{9}(y_{9})]=((-2w_{22}z_{23}+2w_{23}z_{22}-2w_{20}z_{21}+2w_{21}z_{20})e_{1}+(2w_{21}z_{23}-2w_{20}z_{22}-2w_{23}z_{21}+2w_{22}z_{20})e_{2}+(-2w_{20}z_{23}-2w_{21}z_{22}+2w_{22}z_{21}+2w_{23}z_{20})e_{3})E_{2}+((-2w_{22}z_{23}+2w_{23}z_{22}+2w_{20}z_{21}-2w_{21}z_{20})e_{1}+(2w_{21}z_{23}+2w_{20}z_{22}-2w_{23}z_{21}-2w_{22}z_{20})e_{2}+(2w_{20}z_{23}-2w_{21}z_{22}+2w_{22}z_{21}-2w_{23}z_{20})e_{3})E_{4}$,

$[rsp_{9},rsp_{9}^{`}]=Ud_{23}(2w_{22}z_{23}-2w_{23}z_{22}+2w_{20}z_{21}-2w_{21}z_{20})+Ud_{01}(2w_{22}z_{23}-2w_{23}z_{22}+2w_{20}z_{21}-2w_{21}z_{20})+Ud_{67}(2w_{22}z_{23}-2w_{23}z_{22}-2w_{20}z_{21}+2w_{21}z_{20})+Ud_{45}(-2w_{22}z_{23}+2w_{23}z_{22}+2w_{20}z_{21}-2w_{21}z_{20})+Ud_{13}(2w_{21}z_{23}-2w_{20}z_{22}-2w_{23}z_{21}+2w_{22}z_{20})+Ud_{02}(-2w_{21}z_{23}+2w_{20}z_{22}+2w_{23}z_{21}-2w_{22}z_{20})$

\noindent
$+Ud_{57}(-2w_{21}z_{23}-2w_{20}z_{22}+2w_{23}z_{21}+2w_{22}z_{20})+Ud_{46}(-2w_{21}z_{23}-2w_{20}z_{22}+2w_{23}z_{21}+2w_{22}z_{20})+Ud_{12}(2w_{20}z_{23}+2w_{21}z_{22}-2w_{22}z_{21}-2w_{23}z_{20})$

\noindent
$+Ud_{03}(2w_{20}z_{23}+2w_{21}z_{22}-2w_{22}z_{21}-2w_{23}z_{20})+Ud_{56}(2w_{20}z_{23}-2w_{21}z_{22}+2w_{22}z_{21}-2w_{23}z_{20})+Ud_{47}(-2w_{20}z_{23}+2w_{21}z_{22}-2w_{22}z_{21}+2w_{23}z_{20})$,

\noindent
therefore f($[A_{9}(x_{9}),A_{9}(y_{9})])=[rsp_{9},rsp_{9}^{`}]$.

$[A_{9}(x_{9}),A_{10}(y_{10})]=A_{5}((-w_{33}z_{23}-w_{32}z_{22}-w_{31}z_{21}-w_{30}z_{20})e_{0}+(-w_{32}z_{23}+w_{33}z_{22}-w_{30}z_{21}+w_{31}z_{20})e_{1}+(w_{31}z_{23}-w_{30}z_{22}-w_{33}z_{21}+w_{32}z_{20})e_{2}+(-w_{30}z_{23}-w_{31}z_{22}+w_{32}z_{21}+w_{33}z_{20})e_{3})$,

$[rsp_{9},rsp_{10}^{`}]=2Um_{10}(-w_{33}z_{23}-w_{32}z_{22}-w_{31}z_{21}-w_{30}z_{20})-2Um_{11}(w_{32}z_{23}-w_{33}z_{22}+w_{30}z_{21}-w_{31}z_{20})-2Um_{12}(-w_{31}z_{23}+w_{30}z_{22}+w_{33}z_{21}-w_{32}z_{20})-2Um_{13}(w_{30}z_{23}+w_{31}z_{22}-w_{32}z_{21}-w_{33}z_{20})$,

\noindent
therefore f($[A_{9}(x_{9}),A_{10}(y_{10})])=[rsp_{9},rsp_{10}^{`}]$.

\noindent
Case21:

$[A_{10}(x_{10}),A_{10}(y_{10})]=((-2w_{32}z_{33}+2w_{33}z_{32}-2w_{30}z_{31}+2w_{31}z_{30})e_{1}$

\noindent
$+(2w_{31}z_{33}-2w_{30}z_{32}-2w_{33}z_{31}+2w_{32}z_{30})e_{2}+(-2w_{30}z_{33}-2w_{31}z_{32}+2w_{32}z_{31}+2w_{33}z_{30})e_{3})E_{3}+((-2w_{32}z_{33}+2w_{33}z_{32}+2w_{30}z_{31}-2w_{31}z_{30})e_{1}+(2w_{31}z_{33}+2w_{30}z_{32}-2w_{33}z_{31}-2w_{32}z_{30})e_{2}+(2w_{30}z_{33}-2w_{31}z_{32}+2w_{32}z_{31}-2w_{33}z_{30})e_{3})E_{4}$,

$[rsp_{10},rsp_{10}^{`}]=Ud_{23}(2w_{32}z_{33}-2w_{33}z_{32}+2w_{30}z_{31}-2w_{31}z_{30})$

\noindent
$+Ud_{67}(2w_{32}z_{33}-2w_{33}z_{32}-2w_{30}z_{31}+2w_{31}z_{30})+Ud_{45}(-2w_{32}z_{33}+2w_{33}z_{32}+2w_{30}z_{31}-2w_{31}z_{30})+Ud_{01}(-2w_{32}z_{33}+2w_{33}z_{32}-2w_{30}z_{31}+2w_{31}z_{30})$

\noindent
$+Ud_{13}(2w_{31}z_{33}-2w_{30}z_{32}-2w_{33}z_{31}+2w_{32}z_{30})+Ud_{02}(2w_{31}z_{33}-2w_{30}z_{32}-2w_{33}z_{31}+2w_{32}z_{30})+Ud_{57}(-2w_{31}z_{33}-2w_{30}z_{32}+2w_{33}z_{31}+2w_{32}z_{30})$

\noindent
$+Ud_{46}(-2w_{31}z_{33}-2w_{30}z_{32}+2w_{33}z_{31}+2w_{32}z_{30})+Ud_{12}(2w_{30}z_{33}+2w_{31}z_{32}-2w_{32}z_{31}-2w_{33}z_{30})+Ud_{56}(2w_{30}z_{33}-2w_{31}z_{32}+2w_{32}z_{31}-2w_{33}z_{30})$

\noindent
$+Ud_{47}(-2w_{30}z_{33}+2w_{31}z_{32}-2w_{32}z_{31}+2w_{33}z_{30})+Ud_{03}(-2w_{30}z_{33}-2w_{31}z_{32}+2w_{32}z_{31}+2w_{33}z_{30})$,

\noindent
therefore f($[A_{10}(x_{10}),A_{10}(y_{10})])=[rsp_{10},rsp_{10}^{`}]$.

By the above bracket calculations, $\mathfrak{rsp}(4)$ is isomorphic  to $\mathfrak{sp}(4)$.
\ \ \ \ \ \emph{Q.E.D.}

%
%
%
%
%
%
%
%
%
%
\bigskip

\newpage
\section*{Appendix A }
\addcontentsline{toc}{section}{A  Corrections of I.Yokota\cite{Yokota1} and M.Sato and T.Kimura\cite{SatoKimura1}}

\

We fix some printing mistakes in I.Yokota\cite{Yokota1} and M.Sato and T.Kimura\cite{SatoKimura1}
of the root system.

1.  In I.Yokota\cite[\emph{Theorem 2.6.1}]{Yokota1}, ( mistake $\rightarrow$ correction ) 

$(1)\lambda _{3}=\alpha _{3}\rightarrow \lambda _{2}=\alpha _{3},$

$(2)\lambda _{2}+\lambda _{3}=\alpha _{1}+2\alpha _{2}+4\alpha _{3}+4\alpha
_{4}\rightarrow \lambda _{2}+\lambda _{3}=\alpha _{1}+2\alpha _{2}+4\alpha
_{3}+2\alpha _{4}.$

2.  In I.Yokota\cite[\emph{Theorem 4.6.2}]{Yokota1},  ( mistake $\rightarrow$ correction ) 

$(1)\alpha _{6}=\mu _{2}+\frac{3}{2}\nu \rightarrow \alpha _{6}=\mu _{2}+%
\frac{2}{3}\nu ,$

$(2)\alpha _{7}=-\mu _{3}-\frac{3}{2}\nu \rightarrow \alpha _{7}=-\mu _{3}-%
\frac{2}{3}\nu ,$

$(3)\frac{1}{2}(\lambda _{0}-\lambda _{1}-\lambda _{2}+\lambda _{3})+\frac{1%
}{2}\mu _{3}+\frac{2}{3}\nu =1\ \ 1\ \ 1\ \ 2\ \ 1\ \ 0\ \ 1$

$\rightarrow \frac{1}{2}(\lambda
_{0}-\lambda _{1}-\lambda _{2}+\lambda _{3})+\frac{1}{2}\mu _{3}-\frac{2}{3}%
\nu =1\ \ 1\ \ 1\ \ 2\ \ 1\ \ 0\ \ 1,$

$(4)\frac{1}{2}(\lambda _{0}+\lambda _{1}+\lambda _{2}+\lambda _{3})-\frac{1%
}{2}\mu _{3}-\frac{2}{3}\nu =0\ \ 0\ \ 0\ \ 0\ \ 1\ \ 1\ \ 0$

$\rightarrow \frac{1}{2}(\lambda
_{0}+\lambda _{1}+\lambda _{2}+\lambda _{3})-\frac{1}{2}\mu _{3}+\frac{2}{3}%
\nu =0\ \ 0\ \ 0\ \ 0\ \ 1\ \ 1\ \ 0.$

3. In I.Yokota\cite[\emph{Theorem 5.6.1}]{Yokota1},  ( mistake $\rightarrow$ correction ) 

$(1)\pm (\mu _{j}-\frac{1}{3}\nu +r)\rightarrow \pm (\mu _{j}-\frac{1}{3}\nu
)\pm r,$

$(2)\pm \frac{1}{2}(\lambda _{0}+\lambda _{1}-\lambda _{2}-\lambda _{3})\pm (%
\frac{1}{2}\mu _{2}+\frac{1}{3}\nu )\pm r$

\ $\rightarrow \pm \frac{1}{2}(\lambda
_{0}-\lambda _{1}-\lambda _{2}-\lambda _{3})\pm (\frac{1}{2}\mu _{2}+\frac{1%
}{3}\nu )\pm r,$

$(3)\pm \frac{1}{2}(\lambda _{0}+\lambda _{1}-\lambda _{2}-\lambda _{3})\pm (%
\frac{1}{2}\mu _{2}+\frac{1}{3}\nu )\pm r$

\ $\rightarrow \pm \frac{1}{2}%
(-\lambda _{0}+\lambda _{1}-\lambda _{2}-\lambda _{3})\pm (\frac{1}{2}\mu
_{2}+\frac{1}{3}\nu )\pm r,$

$(4)\mu _{2}-\frac{1}{3}\nu +r=0\ \ 0\ \ 0\ \ 1\ \ 0\ \ 0\ \ 0\ \ 0$

$\rightarrow \mu _{2}-\frac{1}{3}\nu
-r=0\ \ 0\ \ 0\ \ 1\ \ 0\ \ 0\ \ 0\ \ 0,$

$(5)-\mu _{3}-\frac{1}{3}\nu -r=0\ \ 1\ \ 1\ \ 1\ \ 0\ \ 0\ \ 0\ \ 1$

$\rightarrow \ -\mu _{3}+\frac{1}{3}\nu
-r=0\ \ 1\ \ 1\ \ 1\ \ 0\ \ 0\ \ 0\ \ 1.$

$(6)H_{\alpha _{1}}=(\Phi (\frac{1}{120}%
(H_{0}-H_{1}-H_{2}-H_{3})+2(E_{3}-E_{1})^{{}},0,0,0),0,0,0,0,0),$

$\rightarrow H_{\alpha _{1}}=(\Phi (\frac{1}{120}(H_{0}-H_{1}-H_{2}-H_{3})+%
\frac{1}{60}(E_{3}-E_{1})^{{}},0,0,0),0,0,0,0,0)$.

4.  In M.Sato and T.Kimura\cite[\emph{Example 41}]{SatoKimura1}, there is a mistake
as follows:

In the root system of $\Delta$, $\Lambda^{'}$ and $\Lambda^{*}$ are opposite to each other.

\bigskip

\section*{Appendix B}
\addcontentsline{toc}{section}{B  Matrix expression of elements of Lie groups $F_{4}$, $E_{6}$, and $G_{2}$}

\

\emph{Remark 20.1.} \ \ We denote elements of the Lie group corresponding to the bases of the Lie algebras \gf$_{4}$ and \ge$_{6}$.
Let's put  followings:

$Ud_{ij}^{`}=Rd_{ij}^{`}/d_{ij}$, $U m_{ij}^{`}=Rm_{ij}^{`}/m_{ij}$, $Ut_{ij}^{`}=Rt_{ij}^{`}/t_{ij}$, $U\tau_{i}^{`}=R\tau_{i}^{`}/\tau_{i}$.

\noindent
And let's put cos$(x)$  as $c_{2}$, cos$(\frac{x}{2})$ as $c$, sin$(x)$ as $s_{2}$, sin$(\frac{x}{2})$ as $s$, exp$(xi)$ as $r_{2}$, and exp$(\frac{x}{2}i)$ as $r$.
Then we have the following expressions:

\noindent
$exp(xUd_{01}^{`})=$
{\fontsize{6pt}{8pt} \selectfont%
$\left( 
%
\right)$ }

\bigskip

\emph{Remark 20.2.} \ \ $\{$exp$(x_{ij}Ud_{ij}^{`}) \ (0 \le i<j \le 7)$,\ exp$(x_{ij}Um_{ij}^{`}) \ (1 \le i \le 3, 0 \le j \le 7) \}$
are make up bases of a simply connected compact Lie group of type $F_{4}$.

\noindent
And $\{$exp$(x_{ij}Ud_{ij}^{`}) \ (0 \le i<j \le 7)$,\ exp$(x_{ij}Um_{ij}^{`}) \ (1 \le i \le 3, 0 \le j \le 7)$,\ exp$(x_{ij}Ut_{ij}^{`}) \ (1 \le i \le 3, 0 \le j \le 7)$,
\ exp$(x_{1}(U\tau_{1}^{`}-U\tau_{2}^{`}))$,\ exp$(x_{2}(U\tau_{1}^{`}+U\tau_{2}^{`}))$ $ \}$
are make up bases of a simply connected compact Lie group of type $E_{6}$.

\bigskip

\emph{Remark 20.3.} \ \ We denote elements of the Lie group corresponding to the bases of the Lie algebras \gg$_{2}$.

\begin{flushleft}
{\fontsize{7pt}{8pt} \selectfont%
$s_{1}(x):=\exp(xS_{1})=\left(
\begin{array}
[c]{cccccccc}%
1 & 0 & 0 & 0 & 0 & 0 & 0 & 0\\
0 & \cos2x & \sin2x & 0 & 0 & 0 & 0 & 0\\
0 & -\sin2x & \cos2x & 0 & 0 & 0 & 0 & 0\\
0 & 0 & 0 & 1 & 0 & 0 & 0 & 0\\
0 & 0 & 0 & 0 & \cos x & 0 & 0 & -\sin x\\
0 & 0 & 0 & 0 & 0 & \cos x & -\sin x & 0\\
0 & 0 & 0 & 0 & 0 & \sin x & \cos x & 0\\
0 & 0 & 0 & 0 & \sin x & 0 & 0 & \cos x
\end{array}
\right)  ,$

$s_{2}(x):=\exp(xS_{2})=\left(
\begin{array}
[c]{cccccccc}%
1 & 0 & 0 & 0 & 0 & 0 & 0 & 0\\
0 & \cos2x & 0 & \sin2x & 0 & 0 & 0 & 0\\
0 & 0 & 1 & 0 & 0 & 0 & 0 & 0\\
0 & -\sin2x & 0 & \cos2x & 0 & 0 & 0 & 0\\
0 & 0 & 0 & 0 & \cos x & 0 & -\sin x & 0\\
0 & 0 & 0 & 0 & 0 & \cos x & 0 & \sin x\\
0 & 0 & 0 & 0 & \sin x & 0 & \cos x & 0\\
0 & 0 & 0 & 0 & 0 & -\sin x & 0 & \cos x
\end{array}
\right)  ,$

$s_{3}(x):=\exp(xS_{3})=\left(
\begin{array}
[c]{cccccccc}%
1 & 0 & 0 & 0 & 0 & 0 & 0 & 0\\
0 & \cos2x & 0 & 0 & \sin2x & 0 & 0 & 0\\
0 & 0 & \cos x & 0 & 0 & 0 & 0 & \sin x\\
0 & 0 & 0 & \cos x & 0 & 0 & \sin x & 0\\
0 & -\sin2x & 0 & 0 & \cos2x & 0 & 0 & 0\\
0 & 0 & 0 & 0 & 0 & 1 & 0 & 0\\
0 & 0 & 0 & -\sin x & 0 & 0 & \cos x & 0\\
0 & 0 & -\sin x & 0 & 0 & 0 & 0 & \cos x
\end{array}
\right)  ,$

$s_{4}(x):=\exp(xS_{4})=\left(
\begin{array}
[c]{cccccccc}%
1 & 0 & 0 & 0 & 0 & 0 & 0 & 0\\
0 & \cos2x & 0 & 0 & 0 & \sin2x & 0 & 0\\
0 & 0 & \cos x & 0 & 0 & 0 & \sin x & 0\\
0 & 0 & 0 & \cos x & 0 & 0 & 0 & -\sin x\\
0 & 0 & 0 & 0 & 1 & 0 & 0 & 0\\
0 & -\sin2x & 0 & 0 & 0 & \cos2x & 0 & 0\\
0 & 0 & -\sin x & 0 & 0 & 0 & \cos x & 0\\
0 & 0 & 0 & \sin x & 0 & 0 & 0 & \cos x
\end{array}
\right)  ,$

$s_{5}(x):=\exp(xS_{5})=\left(
\begin{array}
[c]{cccccccc}%
1 & 0 & 0 & 0 & 0 & 0 & 0 & 0\\
0 & \cos2x & 0 & 0 & 0 & 0 & \sin2x & 0\\
0 & 0 & \cos x & 0 & 0 & -\sin x & 0 & 0\\
0 & 0 & 0 & \cos x & -\sin x & 0 & 0 & 0\\
0 & 0 & 0 & \sin x & \cos x & 0 & 0 & 0\\
0 & 0 & \sin x & 0 & 0 & \cos x & 0 & 0\\
0 & -\sin2x & 0 & 0 & 0 & 0 & \cos2x & 0\\
0 & 0 & 0 & 0 & 0 & 0 & 0 & 1
\end{array}
\right)  ,$

$s_{6}(x):=\exp(xS_{6})=\left(
\begin{array}
[c]{cccccccc}%
1 & 0 & 0 & 0 & 0 & 0 & 0 & 0\\
0 & \cos2x & 0 & 0 & 0 & 0 & 0 & \sin2x\\
0 & 0 & \cos x & 0 & -\sin x & 0 & 0 & 0\\
0 & 0 & 0 & \cos x & 0 & \sin x & 0 & 0\\
0 & 0 & \sin x & 0 & \cos x & 0 & 0 & 0\\
0 & 0 & 0 & -\sin x & 0 & \cos x & 0 & 0\\
0 & 0 & 0 & 0 & 0 & 0 & 1 & 0\\
0 & -\sin2x & 0 & 0 & 0 & 0 & 0 & \cos2x
\end{array}
\right)  ,$

$s_{7}(x):=\exp(xS_{7})=\left(
\begin{array}
[c]{cccccccc}%
1 & 0 & 0 & 0 & 0 & 0 & 0 & 0\\
0 & 1 & 0 & 0 & 0 & 0 & 0 & 0\\
0 & 0 & \cos2x & \sin2x & 0 & 0 & 0 & 0\\
0 & 0 & -\sin2x & \cos2x & 0 & 0 & 0 & 0\\
0 & 0 & 0 & 0 & \cos x & -\sin x & 0 & 0\\
0 & 0 & 0 & 0 & \sin x & \cos x & 0 & 0\\
0 & 0 & 0 & 0 & 0 & 0 & \cos x & -\sin x\\
0 & 0 & 0 & 0 & 0 & 0 & \sin x & \cos x
\end{array}
\right)  ,$

$l_{2}(x):=\exp(xL_{2})=\left(
\begin{array}
[c]{cccccccc}%
1 & 0 & 0 & 0 & 0 & 0 & 0 & 0\\
0 & 1 & 0 & 0 & 0 & 0 & 0 & 0\\
0 & 0 & \cos x & 0 & \sin x & 0 & 0 & 0\\
0 & 0 & 0 & \cos x & 0 & \sin x & 0 & 0\\
0 & 0 & -\sin x & 0 & \cos x & 0 & 0 & 0\\
0 & 0 & 0 & -\sin x & 0 & \cos x & 0 & 0\\
0 & 0 & 0 & 0 & 0 & 0 & 1 & 0\\
0 & 0 & 0 & 0 & 0 & 0 & 0 & 1
\end{array}
\right)  ,$

$l_{3}(x):=\exp(xL_{3})=\left(
\begin{array}
[c]{cccccccc}%
1 & 0 & 0 & 0 & 0 & 0 & 0 & 0\\
0 & 1 & 0 & 0 & 0 & 0 & 0 & 0\\
0 & 0 & \cos x & 0 & 0 & -\sin x & 0 & 0\\
0 & 0 & 0 & \cos x & \sin x & 0 & 0 & 0\\
0 & 0 & 0 & -\sin x & \cos x & 0 & 0 & 0\\
0 & 0 & \sin x & 0 & 0 & \cos x & 0 & 0\\
0 & 0 & 0 & 0 & 0 & 0 & 1 & 0\\
0 & 0 & 0 & 0 & 0 & 0 & 0 & 1
\end{array}
\right)  ,$

$l_{4}(x):=\exp(xL_{4})=\left(
\begin{array}
[c]{cccccccc}%
1 & 0 & 0 & 0 & 0 & 0 & 0 & 0\\
0 & 1 & 0 & 0 & 0 & 0 & 0 & 0\\
0 & 0 & \cos x & 0 & 0 & 0 & \sin x & 0\\
0 & 0 & 0 & \cos x & 0 & 0 & 0 & \sin x\\
0 & 0 & 0 & 0 & 1 & 0 & 0 & 0\\
0 & 0 & 0 & 0 & 0 & 1 & 0 & 0\\
0 & 0 & -\sin x & 0 & 0 & 0 & \cos x & 0\\
0 & 0 & 0 & -\sin x & 0 & 0 & 0 & \cos x
\end{array}
\right)  ,$

$l_{5}(x):=\exp(xL_{5})=\left(
\begin{array}
[c]{cccccccc}%
1 & 0 & 0 & 0 & 0 & 0 & 0 & 0\\
0 & 1 & 0 & 0 & 0 & 0 & 0 & 0\\
0 & 0 & \cos x & 0 & 0 & 0 & 0 & -\sin x\\
0 & 0 & 0 & \cos x & 0 & 0 & \sin x & 0\\
0 & 0 & 0 & 0 & 1 & 0 & 0 & 0\\
0 & 0 & 0 & 0 & 0 & 1 & 0 & 0\\
0 & 0 & 0 & -\sin x & 0 & 0 & \cos x & 0\\
0 & 0 & \sin x & 0 & 0 & 0 & 0 & \cos x
\end{array}
\right)  ,$

$l_{6}(x):=\exp(xL_{6})=\left(
\begin{array}
[c]{cccccccc}%
1 & 0 & 0 & 0 & 0 & 0 & 0 & 0\\
0 & 1 & 0 & 0 & 0 & 0 & 0 & 0\\
0 & 0 & 1 & 0 & 0 & 0 & 0 & 0\\
0 & 0 & 0 & 1 & 0 & 0 & 0 & 0\\
0 & 0 & 0 & 0 & \cos x & -\sin x & 0 & 0\\
0 & 0 & 0 & 0 & \sin x & \cos x & 0 & 0\\
0 & 0 & 0 & 0 & 0 & 0 & \cos x & \sin x\\
0 & 0 & 0 & 0 & 0 & 0 & -\sin x & \cos x
\end{array}
\right)  ,$

$l_{7}(x):=\exp(xL_{7})=\left(
\begin{array}
[c]{cccccccc}%
1 & 0 & 0 & 0 & 0 & 0 & 0 & 0\\
0 & 1 & 0 & 0 & 0 & 0 & 0 & 0\\
0 & 0 & 1 & 0 & 0 & 0 & 0 & 0\\
0 & 0 & 0 & 1 & 0 & 0 & 0 & 0\\
0 & 0 & 0 & 0 & \cos x & 0 & \sin x & 0\\
0 & 0 & 0 & 0 & 0 & \cos x & 0 & \sin x\\
0 & 0 & 0 & 0 & -\sin x & 0 & \cos x & 0\\
0 & 0 & 0 & 0 & 0 & -\sin x & 0 & \cos x
\end{array}
\right)  ,$

$l_{8}(x):=\exp(xL_{8})=\left(
\begin{array}
[c]{cccccccc}%
1 & 0 & 0 & 0 & 0 & 0 & 0 & 0\\
0 & 1 & 0 & 0 & 0 & 0 & 0 & 0\\
0 & 0 & 1 & 0 & 0 & 0 & 0 & 0\\
0 & 0 & 0 & 1 & 0 & 0 & 0 & 0\\
0 & 0 & 0 & 0 & \cos x & 0 & 0 & -\sin x\\
0 & 0 & 0 & 0 & 0 & \cos x & \sin x & 0\\
0 & 0 & 0 & 0 & 0 & -\sin x & \cos x & 0\\
0 & 0 & 0 & 0 & \sin x & 0 & 0 & \cos x
\end{array}
\right)  $}.
\end{flushleft}

\emph{Remark 20.4.} \ \ \ $\{s_{1}(x_{1}),s_{2}(x_{2}),s_{3}(x_{3}),s_{4}(x_{4}%
),s_{5}(x_{5}),s_{6}(x_{6})$,

\noindent
$s_{7}(x_{7}),l_{2}(x_{8}),l_{3}(x_{9}%
),l_{4}(x_{10}),l_{5}(x_{11}),l_{6}(x_{12}),l_{7}(x_{13}),l_{8}(x_{14}), x_{i} \in \R \}$

\noindent
are make up bases of a simply connected compact Lie group of type $G_{2}$.

\end{document}